\documentclass[10pt]{amsart}
\usepackage{latexsym,amssymb,verbatim}
\usepackage[all]{xy}
\usepackage{shuffle} % Used for the shuffle product symbol.
\usepackage{comment} % Used to comment out and toggle the visibility of whole swaths of text.
\usepackage{tabu} % Used for tables with advanced features.
\usepackage{mathdots}
\usepackage{color}
\definecolor{mygreen}{rgb}{0,.4,0}
\definecolor{myblue}{rgb}{0,0,.5}
\usepackage[bookmarks=true,colorlinks,linkcolor=myblue,citecolor=mygreen]{hyperref}
\usepackage[type={CC}, modifier={by}, version={4.0},]{doclicense}

\usepackage[margin=1in]{geometry}

\newtheoremstyle{plainsl}% <name>
  {8pt plus 2pt minus 4pt}% <Space above>
  {8pt plus 2pt minus 4pt}% <Space below>
  {\slshape}% <Body font>
  {0pt}% <Indent amount>
  {\bfseries}% <Theorem head font>
  {.}% <Punctuation after theorem head>
  {5pt plus 1pt minus 1pt}% <Space after theorem headi>
  {}% <Theorem head spec (can be left empty, meaning `normal')>

\theoremstyle{plainsl}
  \newtheorem{theorem}{Theorem}[subsection]
  \newtheorem{proposition}[theorem]{Proposition}
  \newtheorem{lemma}[theorem]{Lemma}
  \newtheorem{corollary}[theorem]{Corollary}
  \newtheorem{conjecture}[theorem]{Conjecture}
\theoremstyle{definition}
  \newtheorem{definition}[theorem]{Definition}
  \newtheorem{example}[theorem]{Example}
  \newtheorem{exercise}[theorem]{Exercise}
  
  \newtheorem{algorithm}[theorem]{Algorithm}
  \newtheorem{question}[theorem]{Question}
 \theoremstyle{remark}
  \newtheorem{remark}[theorem]{Remark}

\newenvironment{statement}{\begin{quote}}{\end{quote}}

\newenvironment{verlong}{}{}
\newenvironment{vershort}{}{}
\newenvironment{nosolutions}{}{}
\excludecomment{verlong}
\includecomment{vershort}
\includecomment{nosolutions}

\newenvironment{commentedout}{}{}
\excludecomment{commentedout}

\newcommand\horiline{\noindent\makebox[\linewidth]{\rule{\textwidth}{1pt}} \newline} % Draw a horizontal line from left to right margin.

\newcommand\arxiv[1]{\href{http://www.arxiv.org/abs/#1}{\texttt{arXiv:#1}}}

\newcommand\qbin[3]{\left[\begin{matrix} #1 \\ #2 \end{matrix} \right]_{#3}}
\newcommand\shuf[2]{{#1}\, \shuffle \,{#2}}

\newcommand\ncps[2]{{#1}\left<\left< {#2} \right>\right>}
\newcommand\powser[2]{{#1}\left[\left[ {#2} \right]\right]}

\newcommand\id{\operatorname{id}} % Identity map/morphism.
\newcommand\St{\operatorname{St}}

\newcommand\Hom{\operatorname{Hom}}
\newcommand\Ind{\operatorname{Ind}}
\newcommand\Inf{\operatorname{Infl}}
\newcommand\Res{\operatorname{Res}}
\newcommand\End{\operatorname{End}}
\newcommand\Sym{\operatorname{Sym}}
\newcommand\Qsym{\operatorname{QSym}}
\newcommand\Nsym{\operatorname{NSym}}
\newcommand\FQsym{\operatorname{FQSym}}
\newcommand\Irr{\operatorname{Irr}}
\newcommand\ch{\operatorname{ch}}
\newcommand\wt{\operatorname{wt}}
\newcommand\mods{\operatorname{mods}}
\newcommand\trace{\operatorname{trace}}
\newcommand\assoc{\operatorname{assoc}}
\newcommand\loops{\operatorname{loops}}
\newcommand\coloops{\operatorname{coloops}}

\newcommand\fin{\operatorname{fin}}
\newcommand\opp{\operatorname{opp}}
\newcommand\res{\operatorname{res}}
\newcommand\ind{\operatorname{ind}}

\newcommand\lex{\operatorname{lex}}
\newcommand\rev{\operatorname{rev}}
\newcommand\inc{\operatorname{inc}}
\newcommand\std{\operatorname{std}}
\newcommand\cont{\operatorname{cont}}
\newcommand\rank{\operatorname{rank}}
\newcommand\acyc{\operatorname{acyc}}
\newcommand\strict{\operatorname{strict}}
\newcommand\weak{\operatorname{natural}}
\newcommand\cols{\operatorname{cols}}
\newcommand\sgn{\operatorname{sgn}}
\newcommand\Des{\operatorname{Des}}
\newcommand\Par{\operatorname{Par}}
\newcommand\Comp{\operatorname{Comp}}
\newcommand\Rib{\operatorname{Rib}}
\newcommand\ps{{\operatorname*{ps}}}
\newcommand\shufmult{\mathbin{\underline{\shuffle}}} % Shuffle product in an algebra
\newcommand\colstrut{\vphantom{\dfrac{a}{a}}} % Vertical buffer space.
\newcommand{\arinj}{\ar@{_{(}->}} % Monomorphism.
\newcommand{\arinjrev}{\ar@{^{(}->}} % Monomorphism.
\newcommand{\arsurj}{\ar@{->>}} % Epimorphism.
\newcommand{\arelem}{\ar@{|->}} % Arrow to be used in describing a map on elements.
\newcommand{\arback}{\ar@{<-}} % Backward arrow.

\newcommand\sq{\square}
\newcommand\bsq{\blacksquare}
\newcommand{\squarebin}{\mathbin{\square}}
\newcommand\one{1} % Multiplicative identity of group or ring.
\newcommand\triv{\underline{1}} % Trivial character of a group or indicator function of a subset.

\newcommand\Symm{\mathfrak{S}} % Symmetric group.
\newcommand\Lieg{\mathfrak{g}}
\newcommand\Liep{\mathfrak{p}}
\newcommand\ooo{\mathfrak{o}}
\newcommand\mmm{\mathfrak{m}}

\newcommand\kk{\mathbf{k}}

\newcommand\starstar{{\circledast}}

\newcommand\ii{{\mathbf{i}}}
\newcommand\jj{{\mathbf{j}}}
\newcommand\xx{{\mathbf{x}}}
\newcommand\XX{{\mathbf{X}}}
\newcommand\yy{{\mathbf{y}}}
\newcommand\YY{{\mathbf{Y}}}
\newcommand\zz{{\mathbf{z}}}
\newcommand\RRR{{\mathbf{R}}}
\newcommand\pt{\mathbf{pt}}

\newcommand\AAA{{\mathcal{A}}}
\newcommand\BBB{{\mathcal{B}}}
\newcommand\CCC{{\mathcal{C}}}
\newcommand\FFF{{\mathcal{F}}}
\newcommand\GGG{{\mathcal{G}}}
\newcommand\HHH{{\mathcal{H}}}
\newcommand\JJJ{{\mathcal{J}}}
\newcommand\LLL{\mathcal{L}}
\newcommand\MMM{{\mathcal{M}}}
\newcommand\PPP{{\mathcal{P}}}
\newcommand\UUU{{\mathcal{U}}}

\newcommand\ZZ{{\mathbb Z}}
\newcommand\FF{{\mathbb F}}
\newcommand\QQ{{\mathbb Q}}
\newcommand\NN{{\mathbb N}}
\newcommand\CC{{\mathbb C}}
\newcommand\RR{{\mathbb R}}
\newcommand\Proj{{\mathbb P}}

\newcommand{\dfn}[1]{\emph{{#1}}\index{#1}} % Emphasizes a word and adds it to the index.
\newcommand{\idx}[1]{{#1}\index{#1}} % Prints a word and adds it to the index.

\usepackage{ifthen}
\numberwithin{equation}{subsection}

\begin{verlong}
\title[Hopf algebras in combinatorics (detailed version containing solutions)]{Hopf algebras in combinatorics (detailed version containing solutions)}
\end{verlong}
\begin{vershort}
\title[Hopf algebras in combinatorics (version containing solutions)]{Hopf algebras in combinatorics (version containing solutions)}
\end{vershort}
\begin{nosolutions}
\title[Hopf algebras in combinatorics]{Hopf algebras in combinatorics}
\end{nosolutions}

\author{Darij Grinberg}
\email{darijgrinberg@gmail.com}
\address{Drexel University,
Korman Center, Room 263,
15 S 33rd Street,
Philadelphia PA, 19104,
USA
// (temporary)
Mathematisches Forschungsinstitut Oberwolfach,
Schwarzwaldstrasse 9--11, 77709 Oberwolfach,
Germany
}
\author{Victor Reiner}
\email{reiner@math.umn.edu}
\address{School of Mathematics\\
University of  Minnesota\\
Minneapolis, MN 55455\\
USA}

\date{July 27, 2020 (with minor corrections July 28, 2020)}

\makeindex

\begin{document}

\maketitle

\tableofcontents

\doclicenseThis

\horiline

\newpage

%%%%%%%%%%%%%%%%%%%%%%%%%%%%%%%%%%%%%
\section*{Introduction}
%%%%%%%%%%%%%%%%%%%%%%%%%%%%%%%%%%%%%
%Certain Hopf algebras arise in combinatorics because they have bases
%naturally parametrized by combinatorial objects (partitions, compositions,
%permutations, tableaux, graphs, trees, posets, polytopes, etc).
%The rigidity in the structure of
%a Hopf algebra can lead to enlightening proofs, and many
%interesting invariants of combinatorial objects turn out to be evaluations
%of Hopf morphisms.

The concept of a Hopf algebra % is an algebraic notion which
crystallized out of algebraic topology and the study of
algebraic groups in the 1940s and 1950s (see
\cite{AndruskiewitschSantos-history} and \cite{Cartier} for
its history). Being a fairly elementary algebraic notion
itself, it subsequently found applications in other
mathematical disciplines, and is now particularly commonplace in
representation theory\footnote{where it provides explanations
for similarities between group representations and Lie
algebra representations}.

These notes concern themselves
(after a brief introduction into the algebraic foundations of
Hopf algebra theory in Chapter~\ref{Hopf-intro-section})
with the Hopf algebras that appear in combinatorics.
These Hopf algebras tend to have bases naturally parametrized
by combinatorial objects (partitions, compositions,
permutations, tableaux, graphs, trees, posets, polytopes, etc.),
and their Hopf-algebraic operations often encode basic
operations on these objects\footnote{such as concatenating two
compositions, or taking the disjoint union of two graphs
-- but, more often, operations which return a multiset
of results, such as cutting a composition into two pieces
at all possible places, or partitioning a poset into two
subposets in every way that satisfies a certain axiom}.
Combinatorial results can then be seen as particular cases of
general algebraic properties of Hopf algebras
(e.g., the multiplicativity
of the M\"obius function can be recovered from the fact that
the antipode of a Hopf algebra is an algebra
anti-endomorphism), and many interesting invariants of
combinatorial objects turn out to be evaluations
of Hopf morphisms. In some cases (particularly that of
symmetric functions), the rigidity in the
structure of a Hopf algebra can lead to enlightening proofs.

One of the most elementary interesting examples of a
combinatorial Hopf algebra is that of the symmetric functions.
We will devote all of Chapter~\ref{Sym-section} to studying
it, deviating from the usual treatments (such as in
Stanley \cite[Ch. 7]{Stanley}, Sagan \cite{Sagan} and
Macdonald \cite{Macdonald}) by introducing the Hopf-algebraic
structure early on and using it to obtain combinatorial
results. Chapter~\ref{PSH-section} will underpin the
importance of this algebra by proving Zelevinsky's
main theorem of PSH theory, which (roughly) claims that a
Hopf algebra over $\ZZ$ satisfying a certain set of axioms
must be a tensor product of copies of the Hopf algebra
of symmetric functions. These axioms are fairly restrictive,
so this result is far from curtailing the diversity of
combinatorial Hopf algebras; but they are natural enough
that, as we will see in
Chapter~\ref{representations-section}, they are satisfied
for a Hopf algebra of representations of symmetric groups.
As a consequence, this Hopf algebra will be revealed
isomorphic to the symmetric functions -- this is the
famous Frobenius correspondence between symmetric
functions and characters of symmetric groups, usually
obtained through other ways (\cite[\S 7.3]{Fulton},
\cite[\S 4.7]{Sagan}). We will further elaborate on the
representation theories of wreath products and general
linear groups over finite fields; while Zelevinsky's PSH
theory does not fully explain the latter, it illuminates
it significantly.

In the next chapters, we will study further examples of
combinatorial Hopf algebras: the quasisymmetric functions
and the noncommutative symmetric functions in
Chapter~\ref{Qsym-section}, various other algebras (of
graphs, posets, matroids, etc.) in
Chapter~\ref{sect.ABS}, and the Malvenuto-Reutenauer
Hopf algebra of permutations in Chapter~\ref{sect.MR}.

The main prerequisite for reading these notes is a good
understanding of graduate algebra\footnote{William
Schmitt's expositions \cite{Schmitt-exp} are
tailored to a reader interested in combinatorial Hopf
algebras; his notes on modules and algebras cover a
significant part of what we need from abstract algebra,
whereas those on categories cover all category theory
we will use and much more.}, in particular multilinear
algebra (tensor products, symmetric powers and exterior
powers)\footnote{Keith Conrad's expository notes
\cite{Conrad-blurbs} are useful, even if not
comprehensive, sources for the latter.} and basic
categorical language\footnote{We also will use a few
nonstandard notions from linear algebra that are explained
in the Appendix (Chapter~\ref{chp.appendix}).}. In
Chapter~\ref{representations-section}, familiarity with
representation theory of finite groups (over $\CC$)
is assumed, along
with the theory of finite fields and (at some places) the
rational canonical form of a matrix. Only basic knowledge
of combinatorics is required (except for a few spots in
Chapter~\ref{sect.ABS}), and familiarity with
geometry and topology is needed only to understand some
tangential remarks. The concepts of Hopf algebras and
coalgebras and the basics of symmetric function theory will
be introduced as needed. We will work over a commutative
base ring most of the time, but no commutative algebra
(besides, occasionally, properties of modules over a PID)
will be used.

% [DG][v41] Added footnote with Schmitt reference.

These notes began as an accompanying text for
Fall 2012 Math 8680 Topics in Combinatorics,
a graduate class taught by the second author at
the University of Minnesota.
The first author has since added many exercises
(and solutions%
\begin{nosolutions}
\footnote{The version of the notes you are reading does not contain said solutions. The version that does can be downloaded from \url{http://www.cip.ifi.lmu.de/~grinberg/algebra/HopfComb-sols.pdf} or compiled from the sourcecode.}
\end{nosolutions}
), as well as Chapter~\ref{sect.QSym.lyndon}
on Lyndon words and the polynomiality of $\Qsym$.
The notes might still grow, and any
comments, corrections and complaints are welcome!

The course was an attempt to focus on examples that
we find interesting, but which are hard to find
fully explained currently in books or in one paper.
Much of the subject of combinatorial Hopf algebras is
fairly recent (1990s onwards) and still spread over
research papers, although sets of lecture notes do exist,
such as Foissy's \cite{Foissy}.
A reference which we discovered late, having a great deal of overlap with these notes
is Hazewinkel, Gubareni, and Kirichenko \cite{HazewinkelGubareniKirichenko}.
References for the purely algebraic theory of Hopf
algebras are much more frequent (see the beginning of
Chapter~\ref{Hopf-intro-section} for a list).
Another recent text that has a significant amount of
material in common with ours (but focuses on
representation theory and probability applications)
is M\'eliot's \cite{Meliot}.

% [DG][v71] Added \cite{Meliot} reference.

Be warned that our notes are highly idiosyncratic in choice of topics,
and they steal heavily from the sources in the bibliography.

% [DG][v26] I wrote the above introduction (partly using
% your syllabus page http://www-users.math.umn.edu/~reiner/Classes/8680_Fall2012_syllabus.html ).
% Do you know a better source for tensor algebra than
% Keith Conrad's blurbs?

\vskip.2in
\noindent
\textbf{Warnings:} Unless otherwise specified ...
\begin{enumerate}
\item[$\bullet$]
\dfn{$\kk$} here usually denotes a commutative ring\footnote{As explained below,
``ring'' means ``associative ring with $1$''.
The most important cases are when $\kk$ is a field or when $\kk = \ZZ$.}.
\item[$\bullet$]
all maps between $\kk$-modules are $\kk$-linear.
\item[$\bullet$]
every ring or $\kk$-algebra is associative and has a $1$, and every ring morphism or $\kk$-algebra morphism preserves the $1$'s.
\item[$\bullet$]
all $\kk$-algebras $A$ have the property that $\left(\lambda 1_A\right) a = a \left(\lambda 1_A\right) = \lambda a$ for all $\lambda \in \kk$ and $a \in A$.
\item[$\bullet$]
all tensor products are over $\kk$ (unless a subscript specifies a different base ring).\index{$\otimes$}
\item[$\bullet$]
\dfn{$\one$} will denote the multiplicative identity in some ring like $\kk$ or in some
$\kk$-algebra (sometimes also the identity of a group written multiplicatively).
\item[$\bullet$]
for any set $S$, we denote by \dfn{$\id_S$} (or by $\id$) the identity map on $S$.
\item[$\bullet$]
The symbols \dfn{$\subset$} (for ``subset'') and \dfn{$<$}
(for ``subgroup'') don't imply properness
(so $\ZZ \subset \ZZ$ and $\ZZ < \ZZ$).
\item[$\bullet$]
the $n$-th symmetric group (i.e., the group of all permutations of
$\left\{1,2,\ldots,n\right\}$) is denoted \dfn{$\Symm_n$}.
\item[$\bullet$]
A permutation $\sigma \in \Symm_n$ will often be identified with
the $n$-tuple $\left(\sigma\left(1\right), \sigma\left(2\right),
\ldots, \sigma\left(n\right)\right)$, which will occasionally be
written without commas and parentheses (i.e., as follows:
$\sigma\left(1\right)\sigma\left(2\right)\cdots\sigma\left(n\right)$).
This is called the \dfn{one-line notation} for permutations.
\item[$\bullet$]
The product of permutations $a \in \Symm_n$ and $b \in \Symm_n$ is defined by
$(ab)(i) = a(b(i))$ for all $i$.
\item[$\bullet$]
\emph{Words}\index{word} over (or in) an \dfn{alphabet} $I$ simply mean finite tuples of elements
of a set $I$. It is customary to write such a word $\left(a_1,a_2,\ldots,a_k\right)$ as
$a_1a_2\ldots a_k$ when this is not likely to be confused for multiplication.
\item[$\bullet$]
$\NN := \left\{0,1,2,\ldots\right\}$.\index{$\NN$}
\item[$\bullet$]
if $i$ and $j$ are any two objects, then \dfn{$\delta_{i,j}$} denotes the \dfn{Kronecker
delta} of $i$ and $j$; this is the integer $1$ if $i = j$ and $0$ otherwise.
\item[$\bullet$]
a \dfn{family} of objects indexed by a set $I$ means a choice of an object $f_i$ for each element $i \in I$; this family will be denoted either by \dfn{$\left(f_i\right)_{i \in I}$} or by \dfn{$\left\{f_i\right\}_{i \in I}$} (and sometimes the ``$i \in I$'' will be omitted when the context makes it obvious -- so we just write \dfn{$\left\{f_i\right\}$}).
% [DG][v72] Added this. Note that the $\left(...\right)$ notation is mine while the $\left\{...\right\}$ notation is yours; I'd make them uniform but I'm afraid there are just too many places to edit.
\item[$\bullet$]
several objects $s_1, s_2, \ldots, s_k$ are said to be
\dfn{distinct} if every $i \neq j$ satisfy $s_i \neq s_j$.
\item[$\bullet$]
similarly,
several sets $S_1, S_2, \ldots, S_k$ are said to be \dfn{disjoint}
if every $i \neq j$ satisfy $S_i \cap S_j = \varnothing$.
\item[$\bullet$]
the symbol \dfn{$\sqcup$} (and the corresponding quantifier $\bigsqcup$)
denotes a disjoint union of sets or posets.
For example, if $S_1, S_2, \ldots, S_k$ are $k$ sets, then
$\bigsqcup_{i=1}^k S_i$ is their disjoint union.
This \idx{disjoint union} can mean either of the following two things:
\begin{itemize}
\item It can mean the union $\bigcup_{i=1}^k S_i$ in the case when
      the sets $S_1, S_2, \ldots, S_k$ are disjoint.
      This is called an ``\idx{internal disjoint union}'', and is simply
      a way to refer to the union of sets while simultaneously
      claiming that these sets are disjoint.
      Thus, of course, it is only well-defined if the sets are
      disjoint.
\item It can also mean the union
      $\bigcup_{i=1}^k \left\{i\right\} \times S_i$.
      This is called an ``\idx{external disjoint union}'', and is
      well-defined whether or not the sets
      $S_1, S_2, \ldots, S_k$ are disjoint; it is a way to
      assemble the sets $S_1, S_2, \ldots, S_k$ into a larger
      set which contains a copy of each of their elements
      that ``remembers'' which set this element comes from.
\end{itemize}
The two meanings are different, but in the case when
$S_1, S_2, \ldots, S_k$ are disjoint, they are isomorphic.
We hope the reader will not have a hard time telling which
of them we are trying to evoke.
\par
Similarly, the notion of a direct sum of $\kk$-modules
has two meanings (``\idx{internal direct sum}'' and
``\idx{external direct sum}'').
\item[$\bullet$] A sequence $\left(w_1, w_2, \ldots, w_k\right)$
of numbers (or, more generally, of elements of a poset) is said
to be \dfn{strictly increasing} (or, for short, \dfn{increasing})
if it satisfies $w_1 < w_2 < \cdots < w_k$.
A sequence $\left(w_1, w_2, \ldots, w_k\right)$
of numbers (or, more generally, of elements of a poset) is said
to be \dfn{weakly increasing} (or \dfn{nondecreasing})
if it satisfies $w_1 \leq w_2 \leq \cdots \leq w_k$.
Reversing the inequalities, we obtain the definitions of a
\dfn{strictly decreasing} (a.k.a. \dfn{decreasing}) and of a
\dfn{weakly decreasing} (a.k.a. \dfn{nonincreasing})
sequence.
All these definitions extend in an obvious way to infinite sequences.
Note that ``nondecreasing'' is not the same as ``not decreasing'';
for example, any sequence having at most one entry is both
decreasing and nondecreasing, whereas the sequence $\left(1, 3, 1\right)$
is neither.
\end{enumerate}

% [DG] I added the last three bullets. I couldn't believe my eyes when I first
% saw that Sage and GAP (and a fair share of British literature) disagree with
% me on the product of permutations.

% [DG][v17] Changed default type of $\kk$ from field to commutative ring.

% [DG][v19] The three different notions we used to denote by $1$ now are
% called differently:
% -- $1$ or $\one$ is used (interchangeably) for multiplicative identities of
%    rings and monoids.
% -- $\id$ is used for identity maps.
% -- $\triv$ is used for trivial characters of groups and for indicator
%    functions of subsets. This currently resolves to $\underline{1}$,
%    but can be changed easily.

% [DG][v64] Added explanation for $\sqcup$ and $\bigsqcup$.

% [DG][v71] Added bullet point about ``nondecreasing'' vs ``not decreasing''.

Hopefully context will resolve some of the ambiguities.

% [DG][v26] Moved the above conventions and warnings to the introduction,
% as they are global to the text.

% [DG][v37] The outermost layer of sectioning is now called Chapters,
% the next-outermost Sections, etc. (This used to be the case in some
% places but not others; now it should be homogeneous.)

\newpage

%%%%%%%%%%%%%%%%%%%%%%%%%%%%%%%%%%%%%
\section{What is a Hopf algebra?}
\label{Hopf-intro-section}
%%%%%%%%%%%%%%%%%%%%%%%%%%%%%%%%%%%%%

The standard references for Hopf algebras are Abe \cite{Abe} and Sweedler \cite{Sweedler},
and some other good ones are
\cite{CaenepeelVercruysse, ChariPressley, DascalescuNastasescuRaianu,
HazewinkelGubareniKirichenko, Kassel, Kytola-hopf, Montgomery,
Radford-Hopf, Schweigert-Hopf, Wisbauer}.
See also Foissy \cite{Foissy} and Manchon \cite{Manchon} for introductions to Hopf algebras
tailored to combinatorial applications.
Most texts only study Hopf algebras over fields (with exceptions such as
\cite{ChariPressley, CaenepeelVercruysse, Wisbauer}). We will work over arbitrary commutative
rings\footnote{and we will profit from this generality in Chapters
\ref{PSH-section} and \ref{representations-section}, where we will
be applying the theory of Hopf algebras to $\kk = \ZZ$ in a way that would
not be possible over $\kk = \QQ$}, which requires some more care at certain
points (but we will not go deep enough into the algebraic theory to witness
the situation over commutative rings diverge seriously from that over
fields).

% [DG][v17] Added reference to Caenepeel-Vercruysse (just discovered today)
% and a couple lines on fields vs. rings.

% [DG][v14] Added Foissy and Manchon references. The combinatorics they do
% is mostly different from ours, but the Hopf algebra theory intersects a
% lot (unavoidably).
% Foissy's work is among the most readable I have seen on this subject;
% we should reference it more...

Let's build up the definition of Hopf algebra
structure bit-by-bit, starting with the more familiar definition of algebras.

\subsection{Algebras}

Recall that an \emph{associative $\kk$-algebra}\index{algebra}\index{associative}
is defined to be
a $\kk$-module $A$ equipped with an associative $\kk$-bilinear map
$\operatorname{mult} : A \times A \to A$ (the
\dfn{multiplication map} of $A$) and an element $\one \in A$
(the \emph{(multiplicative) unity}\index{unity of an algebra}
or \emph{identity}\index{identity of an algebra} of $A$)
that is neutral for this map $\operatorname{mult}$ (that is, it
satisfies $\operatorname{mult} \left(a, \one\right) =
\operatorname{mult} \left(\one, a\right) = a$ for all $a \in A$).
If we recall that
\begin{itemize}
\item $\kk$-bilinear maps $A \times A \to A$
      are in 1-to-1 correspondence with $\kk$-linear maps
      $A \otimes A \to A$ (by the universal property of the tensor
      product), and
\item elements of $A$ are in 1-to-1 correspondence with $\kk$-linear
      maps $\kk \to A$,
\end{itemize}
then we can restate this classical definition of associative
$\kk$-algebras
as follows in terms of $\kk$-linear maps\footnote{Explicitly speaking,
we are replacing the $\kk$-bilinear multiplication map
$\operatorname{mult} : A \times A \to A$ by the $\kk$-linear map
$m : A \otimes A \to A, \  a \otimes b \mapsto \operatorname{mult}\left(a, b\right)$,
and we are replacing the element $\one \in A$ by the $\kk$-linear map
$u : \kk \to A, \  1_\kk \mapsto \one$.}:

% [DG][v14] Added preceding sentence as well as the convention that
% all $\kk$-algebras are central (Damien Calaque, who has some papers on
% Hopf algebras, is one of the minority who disagrees with this). Also
% added $\kk$-linearity assumption on $m$ and $u$ in definition below.

% [DG][v77] Added preceding paragraph and slightly clarified the
% definition below (as well as that of a coalgebra).

\begin{definition}
\label{def.algebra}
An \emph{associative $\kk$-algebra}\index{algebra}\index{associative}
is a $\kk$-module $A$
equipped with a $\kk$-linear \dfn{associative operation}\index{$m$}
$A \otimes A \overset{m}{\rightarrow} A$,
and a $\kk$-linear \dfn{unit}\index{$u$}
$\kk \overset{u}{\rightarrow} A$,
for which the following two diagrams are commutative:

\begin{equation}
\label{associativity-diagram}
\xymatrix{
& A  \otimes A \otimes A \ar[dl]_{m \otimes \id} \ar[dr]^{\id \otimes m} &  \\
A \otimes A \ar[dr]^{m} & &  A \otimes A \ar[dl]_{m} \\
& A &
}
\end{equation}

\begin{equation}
\label{unit-diagram}
\xymatrix{
A \otimes \kk \ar[d]_{\id \otimes u} & A \ar[l] \ar[d]_{\id} \ar[r] & \kk \otimes A \ar[d]^{u \otimes \id} \\
A \otimes A \ar[r]^-{m} & A & A \otimes A \ar[l]_-{m}
}
\end{equation}
where the maps $A \rightarrow A \otimes \kk$ and
$A \rightarrow \kk \otimes A$ are the isomorphisms
sending $a \mapsto  a \otimes \one$ and
$a \mapsto \one \otimes a$.

We abbreviate ``associative $\kk$-algebra'' as ``$\kk$-algebra''
(associativity is assumed unless otherwise specified) or as
``algebra'' (when $\kk$ is clear from the context).
We sometimes refer to $m$ as the ``multiplication map'' of $A$
as well.
\end{definition}

As we said, the multiplication map $m : A \otimes A \to A$ sends
each $a \otimes b$ to the product $ab$, and the unit map $u : \kk \to A$
sends the identity $1_\kk$ of $\kk$ to the identity $1_A$ of $A$.

% [DG][v77] Added the preceding paragraph and the last sentence
% of the definition above.

Well-known examples of $\kk$-algebras are
\emph{tensor} and \emph{symmetric algebras},
which we can think of as algebras of \emph{words} and 
\emph{multisets}, respectively.

\begin{example}
\label{exa.tensor-alg}
If $V$ is a $\kk$-module and $n \in \NN$, then the \emph{$n$-fold
tensor power $V^{\otimes n}$}\index{tensor power}
of $V$ is the $\kk$-module
$\underbrace{V \otimes V \otimes \cdots \otimes V}_{n\text{ times}}$.
(For $n = 0$, this is the $\kk$-module $\kk$, spanned by the
``empty tensor'' $1_\kk$.)

% [DG][v19] Introduced tensor powers here, since you are using them
% here.

The \dfn{tensor algebra} $T(V)=\bigoplus_{n \geq 0} V^{\otimes n}$
on a $\kk$-module $V$ is an associative $\kk$-algebra
spanned (as $\kk$-module) by decomposable tensors
$v_1 v_2 \cdots v_k := v_1 \otimes v_2 \otimes \cdots \otimes v_k$ with
$k\in\NN$ and $v_1, v_2, \ldots, v_k \in V$. Its multiplication is
defined $\kk$-linearly by
\[
m \left( v_1 v_2 \cdots v_k \otimes w_1 w_2 \cdots w_\ell \right)
:= v_1 v_2 \cdots v_k w_1 w_2 \cdots w_\ell
\]
\footnote{Some remarks about our notation (which we are using here
and throughout these notes) are in order.

Since we are working with tensor products of $\kk$-modules
like $T\left(V\right)$ -- which themselves are made of tensors -- here,
we must specify what the $\otimes$ sign means in expressions like
$a \otimes b$ where $a$ and $b$ are elements of $T\left(V\right)$.
Our convention is the following:
When $a$ and $b$ are elements of a tensor algebra $T\left(V\right)$,
we always understand $a \otimes b$ to mean the pure tensor
$a \otimes b \in T\left(V\right) \otimes T\left(V\right)$
rather than the product of $a$ and $b$ inside the tensor algebra
$T\left(V\right)$. The latter product will plainly be written $ab$.

The operator
precedence between $\otimes$ and multiplication in $T\left(V\right)$
is such that multiplication in $T\left(V\right)$ binds more
tightly than the $\otimes$ sign; e.g., the term $ab \otimes cd$ means
$\left(ab\right) \otimes \left(cd\right)$.
The same convention applies to any algebra instead of
$T\left(V\right)$.
}
for all $k,\ell\in\NN$ and
$v_1, v_2, \ldots, v_k, w_1, w_2, \ldots, w_\ell$ in $V$.
The unit map $u : \kk \to T\left(V\right)$ sends $1_\kk$ to the
empty tensor $1_{T\left(V\right)} = 1_\kk \in \kk = V^{\otimes 0}$.

If $V$ is a free $\kk$-module, say with $\kk$-basis $\{x_i\}_{i \in I}$,
then $T(V)$ has a $\kk$-basis of decomposable tensors
$x_{i_1} \cdots x_{i_k} := x_{i_1} \otimes \cdots \otimes x_{i_k}$
indexed by \emph{words} $(i_1,\ldots,i_k)$ in the alphabet $I$,
and the multiplication on this basis is given by concatenation of words:
\[
m( x_{i_1} \cdots x_{i_k} \otimes  x_{j_1} \cdots x_{j_\ell}) =
x_{i_1} \cdots x_{i_k} x_{j_1} \cdots x_{j_\ell}.
\]
\end{example}

% [DG][v17] Separated general and free cases in the above definition.

% [DG][v27] I merged three footnotes about notations into a single one,
% which defines how notation like $v_1 v_2 \cdots v_k \otimes
% w_1 w_2 \cdots w_\ell$ is to be understood. (It was a mistake of mine
% to have three footnotes for this to begin with!)

Recall that a \dfn{two-sided ideal} of a $\kk$-algebra $A$
is defined to be a $\kk$-submodule $J$ of $A$ such that
all $j \in J$ and $a \in A$ satisfy $ja \in J$ and $aj \in J$.
Using tensors, we can restate this as follows:
A \dfn{two-sided ideal} of a $\kk$-algebra $A$
means a $\kk$-submodule $J$ of $A$
satisfying $m(J \otimes A) \subset J$
and $m(A \otimes J) \subset J$.
Often, the word ``two-sided'' is omitted and one just speaks of
an ideal.

It is well-known that if $J$ is a two-sided ideal of a
$\kk$-algebra $A$,
then one can form a \emph{quotient algebra} $A/J$.

% [DG][v77] Expanded the above two paragraphs (previously a
% single paragraph).

\begin{example}
\label{exa.SymV}
Let $V$ be a $\kk$-module.
The \dfn{symmetric algebra} $\Sym(V)=\bigoplus_{n \geq 0} \Sym^n(V)$ is the
quotient of $T(V)$ by the two-sided ideal generated by
all elements $x y - y x$ with $x,y$ in $V$.
When $V$ is a free $\kk$-module with basis $\left\{x_i\right\}_{i \in I}$,
this symmetric algebra $S\left(V\right)$ can be identified
with a (commutative) polynomial algebra $\kk[x_i]_{i \in I}$,
having a $\kk$-basis of (commutative) monomials
$x_{i_1} \cdots x_{i_k}$ as
$\left\{ i_1,\ldots,i_k \right\}_{\operatorname{multiset}}$
runs through all finite multisubsets\footnote{By
a \dfn{multisubset} of a set $S$, we mean a multiset each of whose
elements belongs to $S$ (but can appear arbitrarily often).}
of $I$, and with multiplication defined $\kk$-linearly via
multiset union\footnote{The \dfn{multiset union} of two finite
multisets $A$ and $B$ is defined to be the multiset $C$ with the
property that every $x$ satisfies
\[
\left(\text{multiplicity of $x$ in $C$}\right)
= \left(\text{multiplicity of $x$ in $A$}\right)
+ \left(\text{multiplicity of $x$ in $B$}\right) .
\]
Equivalently, the multiset union of
$\left\{a_1, a_2, \ldots, a_k\right\}_{\operatorname{multiset}}$
and
$\left\{b_1, b_2, \ldots, b_\ell\right\}_{\operatorname{multiset}}$
is
$\left\{a_1, a_2, \ldots, a_k, b_1, b_2, \ldots, b_\ell\right\}_{\operatorname{multiset}}$.
The multiset union is also known as the
\emph{disjoint union}\index{disjoint union of multisets}
of multisets.
}.
\end{example}

% [DG][v25] Added clarifications about basis and footnote about
% the meaning of "multisubsets".

% [DG][v77] Added footnote about multiset union. Unfortunately,
% the definition is NOT standard in the literature; to some
% people, the multiset union takes maxima instead of sums of
% multiplicities.

Note that the $\kk$-module $\kk$ itself canonically becomes a
$\kk$-algebra.
Its associative operation $m : \kk \otimes \kk \to \kk$ is the canonical
isomorphism $\kk \otimes \kk \to \kk$, and its unit $u : \kk \to \kk$
is the identity map.

% [DG][v13] I added the above nugget of pedantry, and a similar one
% about the $\kk$-coalgebra $\kk$, to explain your use of labels like
% $\Delta$ on a map $\kk \otimes \kk \to \kk$ in commutative diagrams.

Topology and group theory give more examples.

\begin{example}
The \emph{cohomology algebra} $H^*(X;\kk)=\bigoplus_{i \geq 0}H^i(X;\kk)$
with coefficients in $\kk$ for
a topological space $X$ has an associative \emph{cup product}.  Its unit
$\kk =H^*(\pt;\kk) \overset{u}{\rightarrow} H^*(X;\kk)$ is induced
from the unique (continuous) map $X \rightarrow \pt$, where $\pt$ is a one-point space.
\end{example}

\begin{example}
\label{group-algebra-example}
For a group $G$, the \dfn{group algebra} \dfn{$\kk G$} has $\kk$-basis
$\{t_g\}_{g \in G}$ and multiplication defined $\kk$-linearly by
$t_g t_h = t_{gh}$, and
unit defined by $u(1) = t_e$, where $e$ is the identity element of $G$.
\end{example}

\subsection{Coalgebras}

In Definition~\ref{def.algebra}, we have defined the notion of an
algebra entirely in terms of linear maps; thus, by reversing
all arrows, we can define a dual notion, which is called a
\emph{coalgebra}.
If we are to think of the multiplication $A \otimes A  \rightarrow A$
in an algebra as \emph{putting together} two basis elements of $A$ to get a
sum of basis elements of $A$, then coalgebra structure should be thought of
as \emph{taking basis elements apart}.

% [DG][v77] Added the first sentence in the above paragraph.

\begin{definition}
\label{def.coalgebra}
A \emph{co-associative $\kk$-coalgebra}\index{coalgebra}\index{coassociative}
is a $\kk$-module $C$ equipped with
a \dfn{comultiplication}\index{$\Delta$}, that is, a $\kk$-linear map
$C \overset{\Delta}{\rightarrow} C \otimes C$,
and a $\kk$-linear \dfn{counit}\index{$\epsilon$}
$C \overset{\epsilon}{\rightarrow} \kk$
for which the following diagrams (which are exactly the diagrams
in \eqref{associativity-diagram} and \eqref{unit-diagram} but with
\emph{all arrows reversed}) are commutative:

% [DG][v77] Reworded the preceding paragraph somewhat.

\begin{equation}
\label{coassociativity-diagram}
\xymatrix{
& C  \otimes C \otimes C &  \\
C \otimes C \ar[ur]^{\Delta \otimes \id} & &  C \otimes C \ar[ul]_{\id \otimes \Delta} \\
& C \ar[ur]^\Delta \ar[ul]_\Delta&
}
\end{equation}

\begin{equation}
\label{counit-diagram}
\xymatrix{
C \otimes \kk \ar[r] & C & \kk \otimes C \ar[l] \\
C \otimes C \ar[u]^{\id \otimes \epsilon} & C \ar[l]^-{\Delta} \ar[u]^\id \ar[r]_-{\Delta} & C \otimes C \ar[u]^{\epsilon \otimes \id}
}
\end{equation}
Here the maps $C \otimes \kk \rightarrow C$ and
$\kk \otimes C \rightarrow C$ are the isomorphisms
sending $c \otimes \one \mapsto c$ and
$\one \otimes c \mapsto c$.

We abbreviate ``co-associative $\kk$-coalgebra'' as ``$\kk$-coalgebra''
(co-associativity, i.e., the commutativity of the diagram
\eqref{coassociativity-diagram}, is assumed unless otherwise specified) or as
``coalgebra'' (when $\kk$ is clear from the context).

Sometimes, the word ``coproduct'' is used as a synonym for
``comultiplication''\footnote{although the word ``coproduct''
already has a different meaning in algebra}.
\end{definition}

% [DG][v25] Added the last sentence about "coproduct". (You do
% use it in these notes.)

One often uses the \dfn{Sweedler notation}
\begin{equation}
\Delta(c) = \sum_{(c)} c_1 \otimes c_2 = \sum c_1 \otimes c_2
\label{eq.sweedler-not.def1}
\end{equation}
to abbreviate formulas involving $\Delta$.
This means that an expression of the form
$\sum_{(c)} f\left(c_1, c_2\right)$ (where $f : C \times C \to M$
is some $\kk$-bilinear map from $C \times C$ to some
$\kk$-module $M$)
has to be understood to mean $\sum_{k=1}^m f\left(d_k, e_k\right)$,
where $k \in \NN$ and $d_1, d_2, \ldots, d_k \in C$ and
$e_1, e_2, \ldots, e_k \in C$ are chosen such that
$\Delta\left(c\right) = \sum_{k=1}^m d_k \otimes e_k$.
(There are many ways to choose such $k$, $d_i$ and $e_i$, but
they all produce the same result
$\sum_{k=1}^m f\left(d_k, e_k\right)$.
Indeed, the result they produce is
$F\left(\Delta\left(c\right)\right)$, where
$F : C \otimes C \to M$ is the $\kk$-linear map induced by
the bilinear map $f$.)
For example, commutativity of the
left square in \eqref{counit-diagram} asserts that
$\sum_{(c)} c_1 \epsilon(c_2) = c$ for each $c \in C$.
Likewise, commutativity of the right square in
\eqref{counit-diagram} asserts that
$\sum_{(c)} \epsilon(c_1) c_2 = c$ for each $c \in C$.
The commutativity of \eqref{coassociativity-diagram} can
be written as
$\sum_{(c)} \Delta(c_1) \otimes c_2 = \sum_{(c)} c_1 \otimes \Delta(c_2)$,
or (using nested Sweeedler notation to unravel the two remaining $\Delta$'s)
as
\[
\sum_{(c)} \sum_{(c_1)} (c_1)_1 \otimes (c_1)_2 \otimes c_2
=
\sum_{(c)} \sum_{(c_2)} c_1 \otimes (c_2)_1 \otimes (c_2)_2 .
\]

% [DG][v79] Added a bit of explanation for the Sweedler notation.

The $\kk$-module $\kk$ itself canonically becomes a
$\kk$-coalgebra, with its comultiplication $\Delta : \kk \to \kk \otimes \kk$
being the canonical isomorphism $\kk \to \kk \otimes \kk$, and its counit
$\epsilon : \kk \to \kk$ being the identity map.

\begin{example}
Let $\kk$ be a field.
The \emph{homology} $H_*(X;\kk)=\bigoplus_{i \geq 0} H_i(X;\kk)$
for a topological space $X$ is naturally a coalgebra: the
(continuous) \dfn{diagonal embedding}
$X \rightarrow X \times X$ sending $x \mapsto (x,x)$ induces
a coassociative map
\[
H_*(X; \kk) \rightarrow H_*(X \times X;\kk) \cong H_*(X;\kk) \otimes H_*(X;\kk)
\]
in which the last isomorphism comes from the \emph{K\"unneth theorem}
with field coefficients $\kk$.
As before, the unique (continuous) map $X \rightarrow \pt$ induces the
counit $H_*(X;\kk) \overset{\epsilon}{\rightarrow} H_*(\pt;\kk) \cong \kk$.
\end{example}

\begin{exercise}
\label{exe.counit.unique}
Let $C$ be a $\kk$-module, and let
$\Delta : C \rightarrow C \otimes C$ be a $\kk$-linear map.
Prove that there exists
\emph{at most one} $\kk$-linear map
$\epsilon : C \rightarrow \kk$
such that the diagram \eqref{counit-diagram} commutes.
\end{exercise}

% [DG][v3] Added exercise, 19 Apr 2014.

For us, the notion of a coalgebra serves mostly as a stepping
stone towards that of a Hopf algebra, which will be the focus
of these notes.
However, coalgebras have interesting properties of their own
(see, e.g., \cite{Manetti-voyage}).

\subsection{Morphisms, tensor products, and bialgebras}

Just as we rewrote the definition of an algebra in terms of linear
maps (in Definition~\ref{def.algebra}), we can likewise rephrase
the standard definition of a morphism of algebras:

% [DG][v77] Added the preceding paragraph.
% Also slightly rewrote the first sentence of the following
% definition:

\begin{definition} \label{def.morphism-alg-coalg}
A \dfn{morphism of algebras}
is a $\kk$-linear map $A \overset{\varphi}{\rightarrow} B$
between two $\kk$-algebras $A$ and $B$
that makes the following two diagrams commute:

\begin{equation}
\qquad
\xymatrix{
A  \ar[r]^{\varphi} & B  \\
A \otimes A \ar[u]^{m_A} \ar[r]^{\varphi \otimes \varphi} & B \otimes B \ar[u]_{m_B}
}
\qquad
\xymatrix{
A \ar[rr]^{\varphi}& & B  \\
&\kk \ar[ul]^{u_A} \ar[ur]_{u_B}&
}
\end{equation}
Here the subscripts on $m_A,m_B,u_A,u_B$ indicate for which algebra they are
part of the structure (e.g., the map $u_A$ is the map $u$ of the algebra $A$);
we will occasionally use such conventions from now on.

Similarly, a \dfn{morphism of coalgebras}
is a $\kk$-linear map $C \overset{\varphi}{\rightarrow} D$
between two $\kk$-coalgebras $C$ and $D$
that makes the reverse diagrams commute:

\begin{equation}
\label{coalgebra-morphism-diagrams}
\xymatrix{
C \ar[d]_{\Delta_C} \ar[r]^{\varphi} & D \ar[d]^{\Delta_D} \\
C \otimes C \ar[r]^{\varphi \otimes \varphi} & D \otimes D
}
\qquad
\xymatrix{
C \ar[dr]_{\epsilon_C} \ar[rr]^{\varphi}& & D \ar[dl]^{\epsilon_D}  \\
&\kk&
}
\end{equation}

As usual, we shall use the word ``\emph{homomorphism}''
as a synonym for ``morphism'', and we will say
``$\kk$-coalgebra homomorphism'' for ``homomorphism
of coalgebras'' (and similarly for algebras and other
structures).

As usual, the word ``\dfn{isomorphism}'' (of algebras,
of coalgebras, or of other structures that we will define
further below) means ``invertible morphism whose inverse is
a morphism as well''.
Two algebras (or coalgebras, or other structures) are said
to be \dfn{isomorphic} if there exists an isomorphism
between them.
\end{definition}

% [DG][v80] Added the above paragraph about isomorphisms.

\begin{example}
Let $\kk$ be a field.
Continuous maps $X \overset{f}{\rightarrow} Y$ of topological spaces
induce algebra morphisms $H^*(Y;\kk) \rightarrow H^*(X;\kk)$,
and coalgebra morphisms $H_*(X;\kk) \rightarrow H_*(Y;\kk)$.
\end{example}

Coalgebra morphisms behave similarly to algebra morphisms in
many regards:
For example, the inverse of an invertible coalgebra morphism
is again a coalgebra morphism\footnote{The easy proof
of this fact is left to the reader.}. Thus, the invertible
coalgebra morphisms are precisely the coalgebra isomorphisms.

% [DG][v70] Added the above paragraph.

\begin{definition}
\label{def.tensor-prod-algs}
Given two $\kk$-algebras $A, B$, their tensor product $A \otimes B$ also
becomes a $\kk$-algebra defining the multiplication bilinearly via\index{tensor product of algebras}
\[
m((a \otimes b) \otimes (a' \otimes b')) := aa' \otimes bb' ,
\]
or, in other words, $m_{A \otimes B}$ is the composite map
\[
\xymatrixcolsep{4pc}
\xymatrix{
A \otimes B \otimes A \otimes B
\ar[r]^-{\id \otimes T \otimes \id}
&
A \otimes A \otimes B \otimes B \ar[r]^-{m_A \otimes m_B}
&
A \otimes B
}
\]
where $T$ is the \dfn{twist map} $B \otimes A \rightarrow A \otimes B$
that sends $b \otimes a \mapsto a \otimes b$.
(See Exercise~\ref{exe.tensor-prod-algs}(a) below for a proof that
this $\kk$-algebra $A \otimes B$ is well-defined.)

Here we are omitting the topologist's sign\index{topologist's sign convention}
in the twist map which should be present for graded algebras
and coalgebras that come from cohomology and homology:
For homogeneous elements $a$ and $b$ of two graded modules
$A$ and $B$, the topologist's twist map
$T : B \otimes A \to A \otimes B$ sends
\begin{equation}
\label{topologist-twist}
b \otimes a \longmapsto  (-1)^{\deg(b)\deg(a)} a \otimes b
\end{equation}
instead of $b \otimes a \mapsto a \otimes b$.
This means that, if one is using the topologists' convention,
most of our examples which we later call \emph{graded}
should actually be considered to live in only \emph{even} degrees
(which can be achieved, e.g., by artificially doubling all
degrees in their grading).  We will, however, keep to our own
definitions (so that our twist map $T$ will always send
$b \otimes a \mapsto a \otimes b$) unless otherwise
noted.  Only in parts of Exercise~\ref{exe.Sym.wedge} will we use
the topologist's sign.  Readers interested in the wide world of
algebras defined using the topologist's sign convention
(which is also known as the \dfn{Koszul sign rule})
can consult \cite[Appendix A2]{Eisenbud}; see also
\cite{GrosshansRotaStein} for applications to algebraic combinatorics%
\footnote{To be precise, \cite{GrosshansRotaStein} works with the
related concept of \dfn{superalgebras}, which are graded by elements
of $\ZZ / 2 \ZZ$ rather than $\NN$ but use the same sign convention
as the topologists have for algebras.}.

% old wording:
% We will ignore this issue, and hope that it causes no confusion later!

% [DG][v25] I changed the above paragraph to be less scary to
% the reader and less apologetic towards topologists
% -- it used to sound as though we were cheating,
% while all it is supposed to say (I hope) is that we are missing
% out on some of the variety of Hopf algebras. To compensate, I
% pointed out that the bialgebra is super in the homological
% example below and in the exercise on the exterior algebra.
% I hope I caught all places where this matters.
% 
% Note that if we wanted to deal with the super setting for real,
% I think we would also have to change the definitions of
% (co)commutativity as well, unless we require $2$ invertible.
% (And the Jacobi identity for Lie algebras becomes even worse.)

% [DG][v48] Replaced rightarrows/longrightarrows with oversets by
% xymatrix diagrams. This should prevent the labels of the arrows
% from protruding beyond the arrows.

% [DG][v80] Added more sentences (including references) to the
% above semi-apologetic paragraph about the topologist's sign.

The unit element of $A \otimes B$ is $\one_A \otimes \one_B$,
meaning that the
unit map $\kk \overset{u_{A \otimes B}}{\rightarrow} A \otimes B$ is the composite
\[
\xymatrixcolsep{4pc}
\xymatrix{
\kk \ar[r] &
\kk \otimes \kk \ar[r]^-{u_A \otimes u_B} &
A \otimes B
} .
\]

Similarly, given two coalgebras $C, D$, one can make $C \otimes D$
a coalgebra in which the comultiplication and counit maps are the composites of\index{tensor product of coalgebras}
\[
\xymatrixcolsep{4pc}
\xymatrix{
C \otimes D \ar[r]^-{\Delta_C \otimes \Delta_D} &
C \otimes C \otimes D \otimes D
\ar[r]^-{\id \otimes T \otimes \id} &
C \otimes D \otimes C \otimes D
}
\]
and
\[
\xymatrixcolsep{4pc}
\xymatrix{
C \otimes D \ar[r]^-{\epsilon_C \otimes \epsilon_D} &
\kk \otimes \kk \ar[r] & \kk
} .
\]
(See Exercise~\ref{exe.tensor-prod-algs}(b) below for a proof that
this $\kk$-coalgebra $C \otimes D$ is well-defined.)
\end{definition}

\begin{exercise}
\phantomsection\label{exe.tensor-prod-algs}

\begin{enumerate}
\item[(a)] Let $A$ and $B$ be two $\kk$-algebras. Show that the
$\kk$-algebra $A\otimes B$ introduced in
Definition~\ref{def.tensor-prod-algs} is actually
well-defined (i.e., its multiplication
and unit satisfy the axioms of a $\kk$-algebra).

\item[(b)] Let $C$ and $D$ be two $\kk$-coalgebras. Show that the
$\kk$-coalgebra $C\otimes D$ introduced in
Definition~\ref{def.tensor-prod-algs} is actually
well-defined (i.e., its comultiplication and counit satisfy
the axioms of a $\kk$-coalgebra).
\end{enumerate}
\end{exercise}

% [DG][v68] Added the above exercise for the sake of
% completeness.

It is straightforward to show that the concept of tensor
products of algebras and of coalgebras satisfy the properties
one would expect:

\begin{itemize}
\item For any three $\kk$-coalgebras $C$, $D$ and $E$,
the $\kk$-linear map
\[
\left(C \otimes D\right) \otimes E
\to C \otimes \left(D \otimes E\right),
\qquad \left( c \otimes d \right) \otimes e
\mapsto c \otimes \left( d \otimes e \right)
\]
is a coalgebra isomorphism.
This allows us to speak of the $\kk$-coalgebra
$C \otimes D \otimes E$ without worrying about the
parenthesization.
\item For any two $\kk$-coalgebras $C$ and $D$, the
$\kk$-linear map
\[
T : C \otimes D \to D \otimes C,
\qquad c \otimes d \mapsto d \otimes c
\]
is a coalgebra isomorphism.
\item For any $\kk$-coalgebra $C$, the $\kk$-linear maps
\begin{align*}
C &\to \kk \otimes C, \qquad c \mapsto 1 \otimes c
\qquad \text{and} \\
C &\to C \otimes \kk, \qquad c \mapsto c \otimes 1
\end{align*}
are coalgebra isomorphisms.
\item Similar properties hold for algebras instead of
coalgebras.
\end{itemize}

% [DG][v70] Added the above list.

\noindent
One of the first signs that these definitions interact nicely is the following straightforward proposition.

\begin{proposition}
\label{bialgebra-prop}
When $A$ is both a $\kk$-algebra and a
$\kk$-coalgebra, the following are equivalent:
\begin{itemize}
\item The maps $\Delta$ and $\epsilon$ are morphisms for the algebra structure $(A,m,u)$.
\item The maps $m$ and $u$ are morphisms for the coalgebra structure $(A,\Delta,\epsilon)$.
\item These four diagrams commute:
\begin{equation}
\label{bialgebra-diagrams}
\begin{array}{ccc}
& \xymatrix{
 & A \otimes A \ar[dl]_{\Delta \otimes \Delta} \ar[ddr]_m & \\
A \otimes A \otimes A \otimes A \ar[dd]_{\id \otimes T \otimes \id} &  & \\
 & & A \ar[ddl]_\Delta\\
A \otimes A \otimes A \otimes A \ar[dr]_{m \otimes m} &  & \\
 & A \otimes A &
} & \\
\xymatrix{
A \otimes A \ar[r]^{\epsilon \otimes \epsilon} \ar[d]_{m} &\kk \otimes \kk \ar[d]^m \\
A \ar[r]_{\epsilon}  &\kk
} & &
\xymatrix{
\kk \ar[r]^{u} \ar[d]_{\Delta} &A \ar[d]^\Delta \\
\kk \otimes \kk \ar[r]_{u \otimes u}  &A \otimes A
}\\
& \xymatrix{
\kk \ar[dr]_{u} \ar[rr]^{\id} & & \kk \\
 &A \ar[ur]_{\epsilon} &
} &
\end{array}
\end{equation}
\end{itemize}
\end{proposition}

% [DG][v78] Slight rewording of the first two equivalent
% conditions above.

\begin{exercise} \phantomsection
\label{exe.tensorprod.morphism}
\begin{itemize}
\item[(a)] If $A$, $A^\prime$, $B$ and $B^\prime$ are four
$\kk$-algebras, and $f : A \to A^\prime$ and $g : B \to B^\prime$
are two $\kk$-algebra homomorphisms, then show that $f \otimes g :
A \otimes B \to A^\prime \otimes B^\prime$ is a $\kk$-algebra
homomorphism.
\item[(b)] If $C$, $C^\prime$, $D$ and $D^\prime$ are four
$\kk$-coalgebras, and $f : C \to C^\prime$ and $g : D \to D^\prime$
are two $\kk$-coalgebra homomorphisms, then show that
$f \otimes g : C \otimes D \to C^\prime \otimes D^\prime$ is a
$\kk$-coalgebra homomorphism.
\end{itemize}
\end{exercise}

% [DG][v25] Added the above exercise (as a lemma for another
% exercise further down).

\begin{definition}
Call the $\kk$-module $A$ a \emph{$\kk$-bialgebra}\index{bialgebra} if it is
a $\kk$-algebra and $\kk$-coalgebra satisfying the three
equivalent conditions in Proposition~\ref{bialgebra-prop}.
\end{definition}

\begin{example}
For a group $G$, one can make the \idx{group algebra} \idx{$\kk G$} a coalgebra with
counit $\kk G \overset{\epsilon}{\rightarrow} \kk$ mapping
$t_g \mapsto 1$ for all $g$ in $G$, and with comultiplication
$\kk G \overset{\Delta}{\rightarrow} \kk G \otimes \kk G$ given by
$\Delta(t_g) := t_g \otimes t_g$.  Checking the various diagrams in
\eqref{bialgebra-diagrams} commute is easy.  For example, one
can check the pentagonal diagram on each basis element $t_g \otimes t_h$:
\[
\xymatrix{
 & t_g \otimes t_h \ar[dl]_{\Delta \otimes \Delta} \ar[ddr]_m & \\
t_g \otimes t_g \otimes t_h \otimes t_h \ar[dd]_{\id \otimes T \otimes \id} &  & \\
 & & t_{gh} \ar[ddl]_\Delta\\
t_g \otimes t_h \otimes t_g \otimes t_h \ar[dr]_{m \otimes m} &  & \\
 & t_{gh} \otimes t_{gh} &
}
\]
\end{example}

\begin{remark}
\label{comultiplication-as-tensor-action}
In fact, one can think of adding a bialgebra structure to a $\kk$-algebra $A$
as a way of making $A$-modules $M,N$ have an $A$-module structure on their
tensor product $M \otimes N$:  the algebra $A \otimes A$ already acts naturally
on $M \otimes N$, so one can let $a$ in $A$ act via $\Delta(a)$ in $A \otimes A$. 
In the theory of group representations over $\kk$, that is, $\kk G$-modules $M$,
this is how one defines the \dfn{diagonal action} of $G$ on $M \otimes N$,
namely $t_g$ acts as $t_g \otimes t_g$.
\end{remark}

\begin{definition}
An element $x$ in a coalgebra for which $\Delta(x) = x \otimes x$ and
$\epsilon(x)=1$ is called \dfn{group-like}.

% [DG] The "coalgebra" above used to say "bialgebra", but I've made it a bit
% less wasteful. (You don't ever use it in either generality, though.)

An element $x$ in a bialgebra for which $\Delta(x) = 1 \otimes x+x\otimes 1$
is called \dfn{primitive}.
We shall also sometimes abbreviate
``primitive element'' as ``primitive''.
\end{definition}

% [DG][v64] Added the preceding sentence.

\begin{example}
Let $V$ be a $\kk$-module.
The \dfn{tensor algebra} $T(V)=\bigoplus_{n \geq 0} V^{\otimes n}$ is
a coalgebra, with counit $\epsilon$ equal to
the identity on $V^{\otimes 0} = \kk$ and the zero map on
$V^{\otimes n}$ for $n > 0$, and with comultiplication defined
to make the elements $x$ in $V^{\otimes 1}=V$ all primitive:
\[
\Delta(x) := \one \otimes x + x \otimes \one \text{ for }x \in V^{\otimes 1}.
\]
Since the elements of $V$
generate $T(V)$ as a $\kk$-algebra, and since $T(V) \otimes T(V)$
is also an associative $\kk$-algebra, the universal property of $T(V)$ as
the free associative $\kk$-algebra on the generators $V$
allows one to define $T(V) \overset{\Delta}{\rightarrow} T(V) \otimes T(V)$ arbitrarily on $V$, and extend it as
an algebra morphism.

It may not be obvious that this $\Delta$ is coassociative,
but one can prove this as follows. Note that
\[
\left(
(\id \otimes \Delta) \circ \Delta
\right)(x) =
x \otimes \one \otimes \one +
\one \otimes x \otimes \one +
\one \otimes \one \otimes x
=\left( (\Delta \otimes \id) \circ \Delta \right)(x)
\]
for every $x$ in $V$.  Hence the two
maps $(\id \otimes \Delta) \circ \Delta$ and
$(\Delta \otimes \id) \circ \Delta$, considered
as algebra morphisms
$T(V) \rightarrow T(V) \otimes T(V) \otimes T(V)$,
must coincide on every element of $T(V)$ since they coincide on $V$.
We leave it as an exercise to check the map $\epsilon$ defined as
above satisfies the counit axioms \eqref{counit-diagram}.

Here is a sample calculation in $T(V)$ when $x,y,z$
are three elements of $V$:
\begin{align*}
\Delta(xyz) &= \Delta(x) \Delta(y) \Delta(z) \\
&=(1 \otimes x+x\otimes 1)(1 \otimes y+y\otimes 1)(1 \otimes z+z\otimes 1) \\
&=(1 \otimes xy + x \otimes y+ y \otimes x + xy \otimes 1)(1 \otimes z+z\otimes 1) \\
&=1 \otimes xyz
+ x \otimes yz + y \otimes xz + z \otimes xy\\
&\qquad + xy \otimes z + xz \otimes y + yz \otimes x
+xyz \otimes 1 .
\end{align*}
This illustrates the idea that comultiplication ``takes basis elements apart''
(and, in the case of $T(V)$, not just basis elements, but any decomposable tensors).
Here for any $v_1, v_2, \ldots, v_n$ in $V$ one has
\[
\Delta\left( v_1 v_2 \cdots v_n \right) =
\sum
v_{j_1} \cdots v_{j_r}
\otimes v_{k_1} \cdots v_{k_{n-r}}
\]
where the sum is over ordered pairs
$\left( j_1, j_2, \ldots, j_r \right) , \left( k_1, k_2, \ldots, k_{n-r} \right)$
of complementary subwords of the word $\left( 1, 2, \ldots, n \right)$.
\footnote{More formally speaking, the sum is over all permutations
$\left(j_1, j_2, \ldots, j_r, k_1, k_2, \ldots, k_{n-r}\right)$ of
$\left(1,2,\ldots,n\right)$ satisfying $j_1 < j_2 < \cdots < j_r$
and $k_1 < k_2 < \cdots < k_{n-r}$.} Equivalently (and in a more
familiar language),
\begin{equation}
\label{eq.tensor-alg.m}
\Delta\left( v_1 v_2 \cdots v_n \right) =
\sum_{I \subset \left\{1,2,\ldots,n\right\}}
v_I \otimes v_{\left\{1,2,\ldots,n\right\} \setminus I},
\end{equation}
where $v_J$ (for $J$ a subset of $\left\{1,2,\ldots,n\right\}$)
denotes the product of all $v_j$ with $j \in J$ in the order of
increasing $j$.
\end{example}

% [DG][v25] Added footnote clarifying the use of $\otimes$, and
% added \eqref{eq.tensor-alg.m}.

% [DG][v27] Moved footnote into an earlier example.

We can rewrite the axioms of a $\kk$-bialgebra $A$ using Sweedler
notation.
Indeed, asking for $\Delta : A \to A \otimes A$ to be a
$\kk$-algebra morphism is equivalent to requiring that
\begin{equation}
\sum_{(ab)} (ab)_1 \otimes (ab)_2
= \sum_{(a)} \sum_{(b)} a_1 b_1 \otimes a_2 b_2
\qquad \text{ for all $a, b \in A$}
\label{eq.bialgebra.sweedler-Deltaab}
\end{equation}
and $\sum_{(1)} 1_1 \otimes 1_2 = 1_A \otimes 1_A$.
(The other axioms have already been rewritten or don't
need Sweedler notation.)

% [DG][v79] Added the paragraph above.

Recall one can quotient a $\kk$-algebra $A$ by a two-sided ideal $J$ to
obtain a quotient algebra $A/J$.
An analogous construction can be done for coalgebras using the
following concept, which is dual to that of a two-sided ideal:

\begin{definition}
In a coalgebra $C$, a \dfn{two-sided coideal}\index{coideal} is a
$\kk$-submodule $J \subset C$ for which
\begin{align*}
\Delta(J) &\subset J \otimes C + C \otimes J , \\
\epsilon(J) &= 0 .
\end{align*}
The quotient $\kk$-module $C/J$ then inherits a coalgebra
structure\footnote{Indeed, $J \otimes C + C \otimes J$ is contained in
the kernel of the canonical map $C \otimes C \to
\left(C/J\right) \otimes \left(C/J\right)$;
therefore, the condition $\Delta(J) \subset J \otimes C + C \otimes J$
shows that the map $C \overset{\Delta}{\rightarrow} C \otimes C
\twoheadrightarrow \left(C/J\right) \otimes \left(C/J\right)$ factors through
a map $\overline{\Delta} : C/J \to \left(C/J\right) \otimes \left(C/J\right)$.
Likewise, $\epsilon(J) = 0$ shows that the map $\epsilon : C \to \kk$
factors through a map $\overline{\epsilon} : C/J \to \kk$.
Equipping $C/J$ with these maps $\overline{\Delta}$ and $\overline{\epsilon}$,
we obtain a coalgebra (as the commutativity of the required diagrams
follows from the corresponding property of $C$).}.
Similarly, in a bialgebra $A$, a subset $J \subset A$ which is both
a two-sided ideal and two-sided coideal gives rise to a quotient
bialgebra $A/J$.
\end{definition}

% [DG][v77] More details added to the above footnote.

\begin{exercise}
\label{exe.kernel-is-coideal}
Let $A$ and $C$ be two $\kk$-coalgebras, and $f : A \to C$ a
surjective coalgebra homomorphism.
\begin{itemize}
\item[(a)] If $f$ is surjective, then show that $\ker f$
is a two-sided coideal of $A$.
\item[(b)] If $\kk$ is a field, then show that $\ker f$
is a two-sided coideal of $A$.
\end{itemize}
\end{exercise}

% [DG][v13] Added the above exercise.

% [DG][v23] Extended it by part (b).

\begin{example}
\label{symmetric-algebra-as-bialgebra-example}
Let $V$ be a $\kk$-module.
The \emph{symmetric algebra} $\Sym(V)$ was defined as the
quotient of the tensor algebra $T(V)$
by the two-sided ideal $J$ generated by
all \dfn{commutators} $[x,y]=x y - y x$ for $x,y$ in $V$
(see Example~\ref{exa.SymV}).
Note that $x,y$ are primitive elements in $T(V)$,
and the following very reusable calculation
shows that \emph{the commutator of two primitives is primitive}:
\begin{align}
\Delta [x,y] &= \Delta( xy-yx )
= \Delta\left(x\right) \Delta\left(y\right) - \Delta\left(y\right) \Delta\left(x\right) \nonumber\\
& \qquad \left(\text{since $\Delta$ is an algebra homomorphism}\right) \nonumber\\
&=(\one \otimes x+x\otimes \one)(\one \otimes y+y\otimes \one)
-(\one \otimes y+y\otimes \one)(\one \otimes x+x\otimes \one)\nonumber\\
&=\one \otimes xy - \one \otimes yx +
  xy \otimes \one - yx \otimes \one \nonumber\\
&\qquad +
  x \otimes y + y \otimes x
 - x \otimes y - y \otimes x \nonumber\\
&=\one \otimes (xy-yx) + (xy-yx) \otimes \one \nonumber\\
&=\one \otimes [x,y] + [x,y] \otimes \one .
\label{commutator-of-primitives}
\end{align}
In particular, the commutators $[x,y]$ have $\Delta[x,y]$ in
$J \otimes T(V) + T(V) \otimes J$.  They also satisfy
$\epsilon([x,y])=0$.  Since they are generators for
$J$ as a two-sided ideal, it is not hard to see this implies
$\Delta (J) \subset J \otimes T(V) + T(V) \otimes J$,
and $\epsilon(J) = 0$.  Thus $J$ is
also a two-sided coideal, and
$\Sym(V)=T(V)/J$ inherits a bialgebra structure.
\end{example}

In fact we will see in Section~\ref{PSH-implies-polynomial-section}
that symmetric algebras are the universal example of bialgebras which are
\emph{graded, connected, commutative, cocommutative}.  But first we should
define some of these concepts.

\begin{definition} \phantomsection
\begin{enumerate}

\item[(a)]
A \emph{graded $\kk$-module}\index{graded module}\footnote{%
also known as an ``\emph{$\NN$-graded $\kk$-module}''}
is a $\kk$-module $V$ equipped with a $\kk$-module
direct sum decomposition $V=\bigoplus_{n \geq 0} V_n$.
In this case, the addend $V_n$ (for any given $n \in \NN$)
is called the
\emph{$n$-th homogeneous component}\index{homogeneous component}
(or the \emph{$n$-th graded component}\index{graded component})
of the graded $\kk$-module $V$.
Furthermore, elements $x$ in $V_n$
are said to be \dfn{homogeneous} of degree $n$; occasionally, the
notation $\deg(x)=n$ is used to signify this\footnote{This
notation should not be taken too literally, as it would absurdly
imply that $\deg(0)$ ``equals'' every $n \in \NN$ at the same time,
since $0 \in V_n$ for all $n$.}.
The decomposition $\bigoplus_{n \geq 0} V_n$ of $V$ (that is,
the family of submodules $\left( V_n \right)_{n \in \NN}$)
is called the \dfn{grading} of $V$.

% [DG][v14] Where do you use the freeness assumption? I'd much like to
% get rid of it. (In PSH-theory, freeness follows from the axioms
% anyway. In the section about dual spaces, I think we should be
% maximally explicit about its use.)

% [DG][v17] Moved the freeness assumption to the definition of
% "finite type".

\item[(b)]
The tensor product $V \otimes W$ of two graded $\kk$-modules $V$
and $W$ is, by default, endowed
with the graded module structure in which%
\index{tensor product of graded modules}
\[
(V \otimes W)_n:=\bigoplus_{i+j=n} V_i \otimes W_j .
\]

\item[(c)]
A $\kk$-linear map $V \overset{\varphi}{\rightarrow} W$ between
two graded $\kk$-modules is called
\emph{graded}\index{graded linear map} if $\varphi(V_n) \subset W_n$ for all $n$.
Graded $\kk$-linear maps are also called
\emph{homomorphisms of graded $\kk$-modules}\index{homomorphism of graded $\kk$-modules}.
An \dfn{isomorphism of graded $\kk$-modules} means an
invertible graded $\kk$-linear map whose inverse is also graded.%
\footnote{We shall see in Exercise~\ref{exe.graded.iso} that
the ``whose inverse is also graded'' requirement is actually
superfluous (i.e., it is automatically satisfied for an
invertible graded $\kk$-linear map); we are imposing it only in
order to stick to our tradition of defining ``isomorphisms''
as invertible morphisms whose inverses are morphisms as well.}

\item[(d)]
Say that a $\kk$-algebra (or coalgebra, or bialgebra) is
\emph{graded}\index{graded algebra}\index{graded coalgebra}\index{graded bialgebra}
if it is a graded $\kk$-module and
all of the relevant structure maps ($u,\epsilon,m,\Delta$) are graded.

\item[(e)]
Say that a graded $\kk$-module $V$ is
\emph{connected}\index{connected graded module}
if $V_0 \cong \kk$.

\item[(f)]
Let $V$ be a graded $\kk$-module. Then, a
\emph{graded $\kk$-submodule of $V$}\index{graded submodule}
(sometimes also called a
\dfn{homogeneous $\kk$-submodule of $V$}\index{homogeneous submodule})
means a graded $\kk$-module $W$ such that $W \subset V$
as sets, and such that the inclusion map
$W \hookrightarrow V$ is a graded $\kk$-linear map. \\
Note that if $W$ is a graded $\kk$-submodule of $V$,
then the grading of $W$ is uniquely determined by
the underlying set of $W$ and the grading of $V$ --
namely, the $n$-th graded component $W_n$ of $W$
is $W_n = W \cap V_n$ for each $n \in \NN$.
Thus, we can specify a graded $\kk$-submodule of $V$
without explicitly specifying its grading.
From this point of view, a graded $\kk$-submodule of
$V$ can also be defined as a $\kk$-submodule $W$ of $V$
satisfying
$W = \sum_{n \in \NN} \left(W \cap V_n\right)$.
(This sum is automatically a direct sum, and thus
defines a grading on $W$.)

\end{enumerate}
\end{definition}

\begin{example}
Let $\kk$ be a field.
A path-connected space $X$ has its homology and cohomology
\begin{align*}
H_*(X;\kk)&=\bigoplus_{i \geq 0} H_i(X;\kk) ,\\
H^*(X;\kk)&=\bigoplus_{i \geq 0} H^i(X;\kk)
\end{align*}
carrying the structure of connected graded coalgebras and algebras, respectively.  If in addition, $X$ is a topological group,
or even less strongly, a \emph{homotopy-associative $H$-space} (e.g. the \emph{loop space} $\Omega Y$ on some other space $Y$),
the continuous multiplication map
$
X \times X \rightarrow X
$
induces an algebra structure on $H_*(X;\kk)$ and a coalgebra structure on
$H^*(X;\kk)$, so that each become bialgebras in the topologist's sense
(i.e., with the twist as in \eqref{topologist-twist}),
and these bialgebras are dual to each other in a sense soon to be discussed.
This was Hopf's motivation:  the (co-)homology of a compact Lie group carries bialgebra structure that explains
why it takes a certain form;  see Cartier \cite[\S 2]{Cartier}.
\end{example}

% [DG][v25] Isn't the Pontryagin product non-associative in general
% for $H$-spaces? Or do field coefficients somehow spare us this?

\begin{example}
Let $V$ be a graded $\kk$-module.
Then, its tensor algebra $T(V)$ and its
symmetric algebra $\Sym(V)$ are graded Hopf algebras.
The grading is given as follows:
If $v_1, v_2, \ldots, v_k$ are homogeneous elements of $V$
having degrees $i_1, i_2, \ldots, i_k$, respectively,
then the elements $v_1 v_2 \cdots v_k$ of $T(V)$ and $\Sym(V)$
are homogeneous of degree $i_1 + i_2 + \cdots + i_k$.
That is, we have
\[
\deg\left(v_1 v_2 \cdots v_k\right)
= \deg\left(v_1\right) + \deg\left(v_2\right) + \cdots + \deg\left(v_k\right)
\]
for any homogeneous elements $v_1, v_2, \ldots, v_k$ of $V$.

Assuming that $V_0=0$, the graded algebras $T(V)$ and $\Sym(V)$ are connected.
This is a fairly common situation in combinatorics.
For example, we will often turn a (non-graded) $\kk$-module $V$
into a graded $\kk$-module by declaring that all elements of
$V$ are homogeneous of degree $1$, but at other times, it will
make sense to have $V$ live in different (positive) degrees.
\end{example}

% [DG][v25] Inserted parenthesized "more general" statement.

\begin{exercise}
\label{exe.graded.iso}
Let $V$ and $W$ be two graded $\kk$-modules. Prove that
if $f : V \to W$ is an invertible graded $\kk$-linear map,
then its inverse $f^{-1} : W \to V$ is also graded.
\end{exercise}

% [DG][v80] Added the above exercise (I think it is tacitly
% used in a few places).

\begin{exercise}
\label{exe.primitives.graded-coideal}
Let $A = \bigoplus_{n \geq 0} A_n$ be a graded $\kk$-bialgebra. We denote
by $\Liep$ the set of all primitive elements of $A$.

\begin{itemize}

\item[(a)]
Show that
$\Liep$ is a graded $\kk$-submodule of $A$
(that is, we have
$\Liep = \bigoplus_{n \geq 0} \left(\Liep \cap A_n\right)$).

\item[(b)]
Show that $\Liep$ is a two-sided coideal of $A$.

\end{itemize}
\end{exercise}

% [DG][v14] Added the above nearly trivial exercise.

\begin{exercise}
\label{graded-connected-exercise}
Let $A$ be a connected graded $\kk$-bialgebra. Show the following:
\begin{itemize}
\item[(a)] The $\kk$-submodule $\kk = \kk\cdot 1_A$ of $A$ lies in $A_0$.
\item[(b)] The map $u$ is an isomorphism $\kk \overset{u}{\rightarrow} A_0$.
\item[(c)] We have $A_0 = \kk \cdot 1_A$.
\item[(d)] The two-sided ideal
$\ker \epsilon$
is the $\kk$-module of positive degree elements
$I = \bigoplus_{n > 0} A_n$.
\item[(e)] The map $\epsilon$ restricted to $A_0$ is the inverse isomorphism
$A_0 \overset{\epsilon}{\rightarrow} \kk$ to $u$.
\item[(f)] For every $x \in A$, we have
\[
\Delta(x) \in x \otimes 1 + A \otimes I.
\]
\item[(g)] Every $x$ in $I$ satisfies % has comultiplication of the form
\[
\Delta(x) = 1 \otimes x + x \otimes 1 + \Delta_+(x) ,
\qquad \text{where $\Delta_+(x)$ lies in $I \otimes I$.}
\]
\item[(h)] Every $n > 0$ and every $x \in A_n$ satisfy
\[
\Delta \left( x \right)
= 1 \otimes x + x \otimes 1 + \Delta_+ \left( x \right) ,
\qquad
\text{where $\Delta_+ \left( x \right)$ lies in
$\sum_{k=1}^{n-1} A_k \otimes A_{n-k}$.}
\]
\end{itemize}

% [DG][v77] I'm not fond of talking about the "comultiplication of $x$",
% so I've edited part (g) accordingly.

(Use only the gradedness of the unit $u$ and
counit $\epsilon$ maps, along with commutativity of diagrams
\eqref{counit-diagram}, and
\eqref{bialgebra-diagrams} and the connectedness of $A$.)
\end{exercise}

% [DG][v13] I have added part (e) (a weaker version of (f) which
% holds for all $x$, not just for $x \in I$ -- comes in handy in
% proofs --, and also serves as a stepping stone for (f)). I also
% tried hard to formalize the proof of this exercise down to every
% detail, and never found myself using \eqref{unit-diagram}, so
% I have removed the reference to this diagram (mentioning the
% connectedness instead).
% I also clarified what $\kk$ means in part (a) and rewrote some
% of the exercise.

% [DG][v66] Inserted parts (c) and (h), which follow quickly from
% (b) and (c) but are nevertheless worth spelling out (since they
% are used in this form).

Having discussed graded $\kk$-modules, let us also define
the concept of a \emph{graded basis}, which is the analogue of the
notion of a basis in the graded context.
Roughly speaking, a graded basis of a graded $\kk$-module is
a basis that comprises bases of all its homogeneous components.
More formally:

\begin{definition}
\label{def.graded-basis}
Let $V=\bigoplus_{n \geq 0} V_n$ be a graded $\kk$-module.
A \dfn{graded basis} of the graded $\kk$-module $V$
means a basis
$\left\{v_i\right\}_{i \in I}$ of the $\kk$-module $V$ whose
indexing set $I$ is partitioned into subsets
$I_0, I_1, I_2, \ldots$ (which are allowed to be empty) with the
property that, for every $n \in \NN$, the subfamily
$\left\{v_i\right\}_{i \in I_n}$ is a basis of the $\kk$-module
$V_n$.
\end{definition}

\begin{example}
Consider the polynomial ring $\kk\left[x\right]$ in one
variable $x$ over $\kk$.
This is a graded $\kk$-module (graded by the degree of a
polynomial; thus, each $x^n$ is homogeneous of degree $n$).
Then, the family
$\left(x^n\right)_{n \in \NN} = \left(x^0, x^1, x^2, \ldots\right)$
is a graded basis of $\kk\left[x\right]$
(presuming that its indexing set $\NN$ is partitioned into
the one-element subsets
$\left\{0\right\}, \left\{1\right\}, \left\{2\right\}, \ldots$).
The family
$\left(\left(-x\right)^n\right)_{n \in \NN}
= \left(x^0, -x^1, x^2, -x^3, \ldots\right)$
is a graded basis of $\kk\left[x\right]$ as well.
But the family
$\left(\left(1+x\right)^n\right)_{n \in \NN}$
is not, since it contains non-homogeneous elements.
\end{example}

% [DG][v80] I have extracted the definition of a graded
% basis from a footnote (it's used a few times) and added
% the example following it.

We end this section by discussing morphisms between
bialgebras. They are defined as one would expect:

\begin{definition}
A \dfn{morphism of bialgebras} (also known as a
\dfn{$\kk$-bialgebra homomorphism})
is a $\kk$-linear map $A \overset{\varphi}{\rightarrow} B$
between two $\kk$-bialgebras $A$ and $B$
that is simultaneously a $\kk$-algebra homomorphism
and a $\kk$-coalgebra homomorphism.
\end{definition}

For example, any $\kk$-linear map $f : V \to W$
between two $\kk$-modules $V$ and $W$ induces a $\kk$-linear map
$T\left(f\right) : T\left(V\right) \to T\left(W\right)$
between their tensor algebras
(which sends each $v_1 v_2 \cdots v_k \in T\left(V\right)$ to
$f\left(v_1\right) f\left(v_2\right) \cdots f\left(v_k\right) \in T\left(W\right)$)
as well as a $\kk$-linear map
$\Sym\left(f\right) : \Sym\left(V\right) \to \Sym\left(W\right)$
between their symmetric algebras;
both of these maps $T\left(f\right)$ and $\Sym\left(f\right)$
are morphisms of bialgebras.

% [DG][v80] Added above definition and example.

Graded bialgebras come with a special family of endomorphisms,
as the following exercise shows:

\begin{exercise}
\label{exe.graded.Dq}
Fix $q \in \kk$. Let $A = \bigoplus_{n\in \NN} A_n$ be a graded
$\kk$-bialgebra (where the $A_n$ are the homogeneous components of $A$).
Let $D_q : A \rightarrow A$ be the $\kk$-module endomorphism of $A$
defined by setting
\[
D_q \left(  a\right) = q^n a
\qquad \text{for each } n \in \NN \text{ and each } a \in A_n.
\]
(It is easy to see that this is well-defined; equivalently,
$D_q$ could be defined as the direct sum
$\bigoplus_{n\in \NN } \left(  q^n \cdot \id_{A_n}\right)
: \bigoplus_{n\in \NN } A_n \rightarrow \bigoplus_{n\in \NN } A_n$
of the maps $q^n \cdot \id_{A_n} : A_n \rightarrow A_n$.)

Prove that $D_q$ is a $\kk$-bialgebra homomorphism.
\end{exercise}

% [DG][v80] Added the above exercise, as it is tacitly used when
% transferring properties from the antipode of $\Lambda$ to the
% fundamental involution $\omega$ and vice versa.

The tensor product of two bialgebras is canonically a
bialgebra, as the following proposition shows:

\begin{proposition}
\label{prop.bialg.tensor}
Let $A$ and $B$ be two $\kk$-bialgebras.
Then, $A \otimes B$ is both a $\kk$-algebra and a
$\kk$-coalgebra (by Definition~\ref{def.tensor-prod-algs}).
These two structures, combined, turn $A \otimes B$ into
a $\kk$-bialgebra.
\end{proposition}

\begin{exercise} \phantomsection
\label{exe.bialg.tensor}
\begin{enumerate}
\item[(a)] Prove Proposition~\ref{prop.bialg.tensor}.

\item[(b)] Let $G$ and $H$ be two groups.
Show that the $\kk$-bialgebra $\kk G \otimes \kk H$
(defined as in Proposition~\ref{prop.bialg.tensor})
is isomorphic to the $\kk$-bialgebra
$\kk \left[G \times H\right]$.
(The notation $\kk \left[S\right]$ is a synonym for
$\kk S$.)
\end{enumerate}
\end{exercise}

% [DG][v70] Added the above proposition and exercise.

\subsection{Antipodes and Hopf algebras}

There is one more piece of structure needed to make a bialgebra
a Hopf algebra, although it will come for free in the connected graded case.

\begin{definition}
\label{convolution-algebra}
For any coalgebra $C$ and algebra $A$, one can endow the $\kk$-module
$\Hom(C,A)$ (which consists of all $\kk$-linear maps from $C$ to $A$)
% [DG][v77] Minor rewording here.
with an associative algebra structure called
the \dfn{convolution algebra}:  Define the product \dfn{$f \star g$}\index{$\star$}
of two maps $f,g$ in $\Hom(C,A)$ by
% [DG][v62] Replaced "send $f,g$ in $\Hom(C,A)$ to $f \star g$ defined by"
% by "Define the product $f \star g$ of two maps $f,g$ in $\Hom(C,A)$ by".
$(f \star g)(c) = \sum f(c_1) g(c_2)$, using
the Sweedler notation\footnote{See the paragraph around
\eqref{eq.sweedler-not.def1} for the meaning of this
notation.}
$\Delta(c) = \sum c_1 \otimes c_2$.
Equivalently, $f \star g$ is the composite
\[
\xymatrixcolsep{4pc}
\xymatrix{
C \ar[r]^-{\Delta}
& C \otimes C \ar[r]^-{f \otimes g}
& A \otimes A \ar[r]^-{m}
& A
} .
\]
The associativity of this multiplication $\star$ is easy to
check (see Exercise~\ref{exe.convolution.assoc} below).

The map $u \circ \epsilon$ is a two-sided
identity element for $\star$, meaning that every $f \in \Hom(C,A)$ satisfies
\[
\sum f(c_1) \epsilon(c_2) = f(c) = \sum \epsilon(c_1) f(c_2)
\]
for all $c \in C$.
One sees this by adding a top row to \eqref{counit-diagram}:

\begin{equation}
\label{added-to-counit-diagram}
\xymatrix{
A \otimes \kk \ar[r] & A & \kk \otimes A \ar[l] \\
C \otimes \kk \ar[u]^{f \otimes \id} \ar[r] & C \ar[u]^{f} & \kk \otimes C \ar[u]^{\id \otimes f} \ar[l] \\
C \otimes C \ar[u]^{\id \otimes \epsilon} & C \ar[l]^-{\Delta} \ar[u]^\id \ar[r]_-{\Delta} & C \otimes C \ar[u]^{\epsilon \otimes \id}
}
\end{equation}
\end{definition}

In particular, when one has a bialgebra $A$, the convolution product
$\star$ gives an associative algebra structure on $\End(A):=\Hom(A,A)$.

\begin{exercise}
\label{exe.convolution.assoc}
Let $C$ be a $\kk$-coalgebra and $A$ be a $\kk$-algebra. Show that the
binary operation $\star$ on $\Hom\left(C, A\right)$ is associative.
\end{exercise}

% [DG][v32] Added above exercise.

The product $f \star g$ of two elements $f$ and $g$ in a
convolution algebra $\Hom(C, A)$ is often called their
\dfn{convolution}.

The following simple (but useful) property of convolution algebras
says essentially that the $\kk$-algebra
$\left( \Hom\left(C,A\right), \star \right)$ is a covariant
functor in $A$ and a contravariant functor in $C$,
acting on morphisms by pre- and post-composition:

\begin{proposition}
\label{prop.convolution.functor}
Let $C$ and $C'$ be two $\kk$-coalgebras, and
let $A$ and $A'$ be two $\kk$-algebras.
Let $\gamma : C \to C'$
be a $\kk$-coalgebra morphism.
Let $\alpha : A \to A'$
be a $\kk$-algebra morphism.

The map
\[
\Hom\left(C',A\right) \to \Hom\left(C,A'\right),
\qquad
f \mapsto \alpha \circ f \circ \gamma
\]
is a $\kk$-algebra homomorphism from the convolution
algebra $\left(\Hom\left(C', A\right), \star\right)$ to the
convolution algebra
$\left(\Hom\left(C, A'\right), \star\right)$.
\end{proposition}

\begin{proof}[Proof of Proposition~\ref{prop.convolution.functor}.]
Denote this map by $\varphi$.
We must show that $\varphi$ is a $\kk$-algebra homomorphism.

Recall that $\alpha$ is an algebra morphism; thus,
$\alpha \circ m_A = m_{A'} \circ \left(\alpha \otimes \alpha\right)$
and
$\alpha \circ u_A = u_{A'}$.
Also, $\gamma$ is a coalgebra morphism; thus,
$\Delta_{C'} \circ \gamma = \left(\gamma \otimes \gamma\right) \circ \Delta_C$
and
$\epsilon_{C'} \circ \gamma = \epsilon_C$.

Now, the definition of $\varphi$ yields
$
\varphi(u_A \circ \epsilon_{C'})
= \underbrace{\alpha \circ u_A}_{= u_{A'}}
    \circ \underbrace{\epsilon_{C'} \circ \gamma}_{= \epsilon_C}
= u_{A'} \circ \epsilon_{C}
$;
in other words, $\varphi$ sends the unity of the algebra
$\left(\Hom\left(C', A\right), \star\right)$ to the
unity of the algebra
$\left(\Hom\left(C, A'\right), \star\right)$.

Furthermore, every $f \in \Hom \left(C',A\right)$ and
$g \in \Hom\left(C',A\right)$
satisfy
\begin{align}
\varphi(f \star g)
&=\alpha \circ \underbrace{(f \star g)}_{= m_A \circ \left( f \otimes g \right) \circ \Delta_{C'}} \circ \gamma \nonumber\\
&= \underbrace{\alpha \circ m_A}_{= m_{A'} \circ \left(\alpha \otimes \alpha\right)}
   \circ (f \otimes g)
   \circ \underbrace{\Delta_{C'} \circ \gamma}_{= \left(\gamma \otimes \gamma\right) \circ \Delta_C}
   \nonumber\\
&=m_{A'}
   \circ \underbrace{(\alpha \otimes \alpha) \circ (f \otimes g) \circ (\gamma \otimes \gamma)}_{= (\alpha \circ f \circ \gamma) \otimes
          (\alpha \circ g \circ \gamma)} \circ \Delta_C \nonumber\\
&=m_{A'} \circ \left( (\alpha \circ f \circ \gamma) \otimes
          (\alpha \circ g \circ \gamma) \right) \circ \Delta_C \nonumber\\
&= \underbrace{(\alpha \circ f \circ \gamma)}_{= \varphi(f)} \star \underbrace{(\alpha \circ g \circ \gamma)}_{= \varphi(g)}
= \varphi(f) \star \varphi(g) .
\label{pre-and-post-composition-in-convolution}
\end{align}
Thus, $\varphi$ is a $\kk$-algebra homomorphism
(since $\varphi$ is a $\kk$-linear map and
sends the unity of the algebra
$\left(\Hom\left(C', A\right), \star\right)$ to the
unity of the algebra
$\left(\Hom\left(C, A'\right), \star\right)$).
\end{proof}

% [DG][v65] Created the above proposition, which used to be a step
% in the proof of Proposition \ref{antipodes-give-convolution-inverses}.
% Its uses range far beyond that proof, so it makes much sense to me
% to have it separate.
% Extended the proof with more details, too.

\begin{exercise}
\label{exe.convolution.tensor}
Let $C$ and $D$ be two $\kk$-coalgebras, and let $A$ and $B$ be two
$\kk$-algebras. Prove that:
\begin{itemize}
\item[(a)] If $f:C\rightarrow A$,
$f^{\prime}:C\rightarrow A$, $g:D\rightarrow B$ and
$g^{\prime}:D\rightarrow B$ are four $\kk$-linear maps, then
\[
\left(  f\otimes g\right)  \star\left(  f^{\prime}\otimes g^{\prime}\right)
=\left(  f\star f^{\prime}\right)  \otimes\left(  g\star g^{\prime}\right)
\]
in the convolution algebra $\Hom\left(  C\otimes D,A\otimes
B\right)  $.
\item[(b)] Let $R$ be the $\kk$-linear map
$\left(\Hom\left(C,A\right),\star\right) \otimes
 \left(\Hom\left(D,B\right),\star\right) \to
 \left(\Hom\left(C \otimes D, A \otimes B\right),\star\right)$
which sends every tensor
$f \otimes g \in 
 \left(\Hom\left(C,A\right),\star\right) \otimes
 \left(\Hom\left(D,B\right),\star\right)$
to the map $f \otimes g : C \otimes D \to A \otimes B$. (Notice
that the tensor $f \otimes g$ and the map $f \otimes g$ are different
things which happen to be written in the same way.) Then, $R$ is a
$\kk$-algebra homomorphism.
\end{itemize}
\end{exercise}

% [DG][v14] Added the above exercise (simple computation and useful
% lemma for Dynkin idempotent exercise).

\begin{exercise}
\label{exe.convolution.tensor-curry}
Let $C$ and $D$ be two $\kk$-coalgebras.
Let $A$ be a $\kk$-algebra.
Let $\Phi$ be the canonical $\kk$-module isomorphism
$\Hom \left( C \otimes D, A \right)
\to \Hom \left( C, \Hom \left( D, A \right) \right)$
(defined by
$\left( \left( \Phi \left( f \right) \right) \left( c \right)
\right) \left( d \right)
= f \left( c \otimes d \right)$
for all $f \in \Hom \left( C \otimes D, A \right)$,
$c \in C$ and $d \in D$).
Prove that $\Phi$ is a $\kk$-algebra isomorphism
\[
\left( \Hom \left( C \otimes D, A \right) , \star \right)
\to \left( \Hom \left( C, \left( \Hom \left( D, A \right)
, \star\right) \right) , \star \right)  .
\]
\end{exercise}

% [DG][v68] Added the above simple but (to me) surprising
% exercise.

\begin{definition}
A bialgebra $A$ is called a \dfn{Hopf algebra} if there is an
element \dfn{$S$} (called an \dfn{antipode} for $A$) in $\End(A)$
which is a $2$-sided inverse under $\star$ for the identity map $\id_A$.
In other words, this diagram commutes:
\begin{equation}
\label{antipode-diagram}
\xymatrix{
&A \otimes A \ar[rr]^{S \otimes \id_A}& &A \otimes A \ar[dr]^m& \\
A \ar[ur]^\Delta \ar[rr]^{\epsilon} \ar[dr]_\Delta& &\kk \ar[rr]^{u} & & A\\
&A \otimes A \ar[rr]_{\id_A \otimes S}& &A \otimes A \ar[ur]_m& \\
}
\end{equation}
Or equivalently, if we follow the Sweedler notation
in writing $\Delta(a) = \sum a_1 \otimes a_2$,
then
\begin{equation}
\label{antipode-defining-relation}
\sum_{(a)} S(a_1) a_2 = u(\epsilon(a)) = \sum_{(a)} a_1 S(a_2).
\end{equation}
\end{definition}

\begin{example}
For a group algebra $\kk G$, one can define an antipode $\kk$-linearly
via $S(t_g)=t_{g^{-1}}$.  The top pentagon in the above diagram commutes because
\[
(S \star \id)(t_g) = m ((S \otimes \id)(t_g \otimes t_g)) =
S(t_g) t_g = t_{g^{-1}} t_g = t_e = (u \circ \epsilon)(t_g).
\]
\end{example}

Note that when it exists, the antipode $S$ is unique, as with all 2-sided inverses in associative algebras: if $S, S'$ are both 2-sided $\star$-inverses to $\id_A$ then
\[
S'=(u\circ \epsilon) \star S'
=(S \star \id_A) \star S'
=S \star (\id_A \star S')
=S \star (u \circ \epsilon)= S.
\]
Thus, we can speak of ``\emph{the antipode}'' of a Hopf algebra.

Unlike the comultiplication $\Delta$, the antipode $S$ of a Hopf algebra is not always an algebra homomorphism. It is instead an algebra \emph{anti-homomorphism}, a notion we shall now introduce:

\begin{definition} \phantomsection
\label{def.anti-hom}

\begin{itemize}

\item[(a)] For any two $\kk$-modules $U$ and $V$, we let
$T_{U,V} : U \otimes V \rightarrow V \otimes U$\index{$T_{U,V}$}
be the $\kk$-linear map
$U \otimes V \rightarrow V \otimes U$ sending every
$u \otimes v$ to $v \otimes u$.
This map $T_{U,V}$ is called the \dfn{twist map} for $U$ and $V$.

\item[(b)] A
\emph{$\kk$-algebra anti-homomorphism}\index{anti-homomorphism of algebras}
means a $\kk$-linear
map $f : A \rightarrow B$ between two $\kk$-algebras $A$ and $B$
which satisfies
$f \circ m_A
= m_B \circ \left( f \otimes f \right) \circ T_{A, A}$
and $f \circ u_A = u_B$.

\item[(c)] A
\emph{$\kk$-coalgebra anti-homomorphism}\index{anti-homomorphism of coalgebras}
means a $\kk$-linear
map $f : C \rightarrow D$ between two $\kk$-coalgebras $C$ and $D$
which satisfies
$\Delta_D \circ f
= T_{D,D} \circ \left( f \otimes f \right) \circ \Delta_C$
and $\epsilon_D \circ f = \epsilon_C$.

\item[(d)] A
\emph{$\kk$-algebra anti-endomorphism}\index{anti-endomorphism of an algebra}
of a $\kk$-algebra $A$ means a $\kk$-algebra anti-homomorphism from $A$
to $A$.

\item[(e)] A
\emph{$\kk$-coalgebra anti-endomorphism}\index{anti-endomorphism of a coalgebra}
of a $\kk$-coalgebra $C$ means a $\kk$-coalgebra anti-homomorphism from $C$
to $C$.

\end{itemize}
\end{definition}

Parts (b) and (c) of Definition~\ref{def.anti-hom} can be restated in terms of elements:
\begin{itemize}
\item
A $\kk$-linear map $f : A \rightarrow B$ between two $\kk$-algebras $A$ and $B$
is a $\kk$-algebra anti-homomorphism if and only if it satisfies
$f \left(ab\right) = f \left(b\right) f \left(a\right)$ for all $a, b \in A$
as well as $f \left(1\right) = 1$.
\item
A $\kk$-linear map $f : C \to D$ between two $\kk$-coalgebras $C$  and $D$
is a $\kk$-coalgebra anti-homomorphism if and only if it satisfies
$\sum_{\left(f\left(c\right)\right)} \left(f\left(c\right)\right)_1 \otimes \left(f\left(c\right)\right)_2
= \sum_{\left(c\right)} f\left(c_2\right) \otimes f\left(c_1\right)$
and $\epsilon\left(f\left(c\right)\right) = \epsilon\left(c\right)$
for all $c \in C$.
\end{itemize}

\begin{example}
Let $n \in \NN$, and consider the $\kk$-algebra $\kk^{n \times n}$
of $n \times n$-matrices over $\kk$.
The map $\kk^{n \times n} \to \kk^{n \times n}$ that sends
each matrix $A$ to its transpose $A^T$ is a $\kk$-algebra
anti-endomorphism of $\kk^{n \times n}$.
\end{example}

We warn the reader that the composition of two $\kk$-algebra anti-homomorphisms is not generally a $\kk$-algebra anti-homomorphism again, but rather a $\kk$-algebra homomorphism. The same applies to coalgebra anti-homomorphisms.
Other than that, however, anti-homomorphisms share many of the helpful properties of homomorphisms.
In particular, two $\kk$-algebra anti-homomorphisms are identical if they agree on a generating set of their domain.
Thus, the next proposition is useful when one wants to check that a certain map \emph{is} the antipode in a particular Hopf algebra, by checking it on an algebra generating set.

% [DG][v80] Added the above definition (previously partly included
% in footnotes and exercises) and the paragraph after it.

\begin{proposition}
\label{antipodes-are-antiendomorphisms}
The antipode $S$ in a Hopf algebra $A$ is an algebra anti-endomorphism:
$S(\one)=\one$, and $S(ab)=S(b)S(a)$ for all $a,b$ in $A$.
\end{proposition}
\begin{proof}
This is surprisingly nontrivial; the following argument comes from
\cite[proof of Proposition 4.0.1]{Sweedler}.

Since $\Delta$ is an algebra morphism, one has
$\Delta(1)=1\otimes 1$, and therefore
$
\one= u \epsilon (\one) = S(\one) \cdot \one = S(\one).
$

To show $S(ab)=S(b)S(a)$,
consider $A \otimes A$ as a coalgebra and $A$ as an algebra.
Then $\Hom(A \otimes A,A)$ is an associative algebra with a convolution
product $\starstar$ (to be distinguished from the convolution $\star$
on $\End(A)$), having two-sided identity element $u_A \epsilon_{A \otimes A}$.
We define three elements $f$, $g$, $h$
of $\Hom(A \otimes A,A)$ by
\begin{align*}
f(a \otimes b) &= ab, \\
g(a \otimes b) &=S(b) S(a), \\
h(a \otimes b) &=S(ab).
\end{align*}
We will show that these three elements have the property that
\begin{equation}
\label{three-antipode-calculations}
h \starstar f
= u_A \epsilon_{A \otimes A}
= f \starstar g ,
\end{equation}
which would then show the desired equality $h=g$ via associativity:
\[
h=h \starstar (u_A \epsilon_{A \otimes A})
=h \starstar (f \starstar g)
=(h \starstar f) \starstar g
=(u_A \epsilon_{A \otimes A}) \starstar g=g.
\]
So we evaluate the three elements in \eqref{three-antipode-calculations}
on $a \otimes b$. To do so, we use Sweedler notation -- i.e.,
we assume $\Delta(a)=\sum_{(a)} a_1 \otimes a_2$ and
$\Delta(b) = \sum_{(b)} b_1 \otimes b_2$, and hence
$\Delta(ab)= \sum_{(a),(b)} a_1 b_1 \otimes a_2 b_2$
(by \eqref{eq.bialgebra.sweedler-Deltaab});
then,
\begin{align*}
(u_A \epsilon_{A \otimes A})(a \otimes b)
&=u_A( \epsilon_A(a) \epsilon_A(b))= u_A(\epsilon_A (ab)).\\
& \\
(h \starstar f)(a \otimes b)
&=\sum_{(a),(b)} h(a_1 \otimes b_1) f(a_2 \otimes b_2) \\
&=\sum_{(a),(b)} S(a_1 b_1) a_2 b_2 \\
&= (S \star \id_A) (ab)
= u_A (\epsilon_A (ab)).\\
& \\
(f \starstar g)(a \otimes b)
&=\sum_{(a),(b)} f(a_1 \otimes b_1) g(a_2 \otimes b_2) \\
&=\sum_{(a),(b)} a_1 b_1 S(b_2) S(a_2) \\
&=\sum_{(a)} a_1 \cdot (\id_A \star S)(b) \cdot S(a_2) \\
&= u_A(\epsilon_A(b)) \sum_{(a)} a_1 S(a_2)
= u_A (\epsilon_A (b)) u_A (\epsilon_A (a))
= u_A (\epsilon_A (ab)) .
\end{align*}
These results are equal, so that \eqref{three-antipode-calculations}
holds, and we conclude that $h = g$ as explained above.
\end{proof}

% [DG][v77] Slightly uncrowded a couple sentences and added
% one at the end. The proof didn't get longer because
% the blacksquare at the end of the proof already took a
% whole line for itself :)

\begin{remark}
\label{antipode-as-dual-action}
Recall from Remark~\ref{comultiplication-as-tensor-action}
that the comultiplication on a bialgebra $A$ allows one
to define an $A$-module structure on the tensor product $M \otimes N$
of two $A$-modules $M,N$.  Similarly, the anti-endomorphism $S$
in a Hopf algebra allows one to turn \emph{left} $A$-modules into
\emph{right} $A$-modules, or vice-versa.\footnote{Be warned that
these two transformations are not mutually inverse! Turning a left
$A$-module into a right one and then again into a left one using the
antipode might lead to a non-isomorphic $A$-module, unless the
antipode $S$ satisfies $S^2=\id$.}
 E.g., left $A$-modules
$M$ naturally have a right $A$-module structure on the
dual $\kk$-module $M^*:=\Hom(M,\kk)$, defined via $(fa)(m):=f(am)$
for $f$ in $M^*$ and $a$ in $A$.  The antipode $S$
can be used to turn this back
into a left $A$-module $M^*$, via $(af)(m) = f(S(a)m)$.

% [DG][v17] Added footnote. Is the exclamation mark too alarmist?

For groups $G$ and left $\kk G$-modules
(group representations) $M$, this is
how one defines the \dfn{contragredient action} of $G$ on $M^*$,
namely $t_g$ acts as $(t_g f)(m) = f( t_{g^{-1}} m)$.

More generally, if $A$ is a Hopf algebra and $M$ and $N$ are
two left $A$-modules, then $\Hom\left(M, N\right)$ (the $\Hom$
here means $\Hom_\kk$, not $\Hom_A$) canonically
becomes a left $A$-module by setting
\[
\left(af\right)\left(m\right) = \sum_{\left(a\right)} a_1 f\left(S\left(a_2\right) m\right)
\qquad \qquad \text{for all } a \in A, \ f \in \Hom\left(M, N\right) \text{ and } m \in M.
\]
\footnote{In more abstract terms, this $A$-module structure is given
by the composition
\[
\xymatrixcolsep{3pc}
\xymatrix{
A \ar[r]^-{\Delta} & A \otimes A \ar[r]^{\id_A \otimes S} &
A \otimes A^{\operatorname{op}} \ar[r] & \End\left(\Hom\left(M, N\right)\right)
},
\]
where the last arrow is the morphism
\begin{align*}
\nonumber
A \otimes A^{\operatorname{op}} &\longrightarrow \End\left(\Hom\left(M, N\right)\right), \\
a \otimes b &\longmapsto \left(f \mapsto \left(M \to N, \ m \mapsto a f\left(bm\right)\right)\right).
\end{align*}
Here, $A^{\operatorname{op}}$ denotes the \dfn{opposite algebra}
of $A$, which is the $\kk$-algebra differing from $A$ only in the
multiplication being twisted (the product of $a$ and $b$ in
$A^{\operatorname{op}}$ is defined to be the product of $b$ and $a$
in $A$). As $\kk$-modules, $A^{\operatorname{op}} = A$, but we prefer
to use $A^{\operatorname{op}}$ instead of $A$ here to ensure that all
morphisms in the above composition are algebra morphisms.}
When $A$ is the group algebra $\kk G$ of a group $G$, this leads to
\[
\left(t_g f\right)\left(m\right) = t_g f\left(t_{g^{-1}} m\right)
\qquad \qquad \text{for all } g \in G, \ f \in \Hom\left(M, N\right) \text{ and } m \in M.
\]
This is
precisely how one commonly makes $\Hom\left(M, N\right)$ a
representation of $G$ for two representations $M$ and $N$.

Along the same lines, whenever $A$ is a $\kk$-bialgebra,
we are supposed to think of the counit
$A \overset{\epsilon}{\rightarrow} \kk$ as giving a way to
make $\kk$ into a \emph{trivial}\index{trivial module}
$A$-module. This $A$-module $\kk$
behaves as one would expect: the canonical isomorphisms
$\kk \otimes M \to M$, $M \otimes \kk \to M$ and
(if $A$ is a Hopf algebra)
$\Hom\left(M, \kk\right) \to M^*$ are $A$-module isomorphisms
for any $A$-module $M$.
\end{remark}

% [DG][v23] Added definition of $A$-module $\Hom\left(M, N\right)$.

\begin{corollary}
\label{commutative-implies-involutive-antipode-cor}
Let $A$ be a commutative Hopf algebra.
Then, its antipode is an involution: $S^2=\id_A$.
\end{corollary}
\begin{proof}
One checks that $S^2=S \circ S$ is a right $\star$-inverse to $S$, as follows:
\begin{align*}
(S \star S^2)(a)
  &= \sum_{(a)} S(a_1) S^2(a_2) \\
  &= S \left( \sum_{(a)} S(a_2) a_1 \right)
    \qquad \left( \text{by Proposition~\ref{antipodes-are-antiendomorphisms}} \right) \\
  &= S \left( \sum_{(a)} a_1 S(a_2) \right)
    \qquad \left( \text{by commutativity of } A \right) \\
  &= S \left( u (\epsilon(a)) \right)  \\
  & = u (\epsilon(a))
    \qquad \left( \text{since }S(\one)=\one \text{ by Proposition~\ref{antipodes-are-antiendomorphisms}} \right) .
\end{align*}
Since $S$ itself is the $\star$-inverse to $\id_A$, this shows that
$S^2 = \id_A$.
\end{proof}

% [DG][v63] Added a couple more details to the above proof and
% made the annotations follow my usual style
% (\qquad \left( \text{ ... } \right)).

\begin{remark}
\label{rmk.commutative-implies-involutive-antipode-cor.cocomm}
We won't need it, but it is easy to adapt the above proof to show that
$S^2=\id_A$ also holds for \emph{cocommutative} Hopf algebras
(the dual notion to commutativity; see Definition~\ref{def.cocomm} below for the precise definition);
see \cite[Corollary 1.5.12]{Montgomery} or \cite[Proposition 4.0.1 6)]{Sweedler}
or Exercise~\ref{exe.cocomm-antipode-S2} below.
For a general Hopf algebra which is not finite-dimensional over a field
$\kk$, the antipode
$S$ may not even have finite order, even in the connected graded setting.
E.g., Aguiar and Sottile \cite{AguiarSottile}
show that the Malvenuto-Reutenauer Hopf algebra of permutations has antipode of infinite
order.  In general, antipodes need not even be invertible \cite{Takeuchi}.
\end{remark}

\begin{proposition}
\label{prop.hopf.tensor}
Let $A$ and $B$ be two Hopf algebras.
Then, the $\kk$-bialgebra $A \otimes B$ (defined as in
Proposition~\ref{prop.bialg.tensor}) is a Hopf algebra.
The antipode of this Hopf algebra $A \otimes B$ is the
map $S_A \otimes S_B : A \otimes B \to A \otimes B$,
where $S_A$ and $S_B$ are the antipodes of the Hopf
algebras $A$ and $B$.
\end{proposition}

\begin{exercise}
\label{exe.hopf.tensor}
Prove Proposition~\ref{prop.hopf.tensor}.
\end{exercise}

% [DG][v70] Added the above proposition & exercise.

In our frequent setting of connected graded bialgebras,
antipodes come for free.

\begin{proposition}
\label{graded-connected-bialgebras-have-antipodes}
A connected graded bialgebra $A$ has a unique antipode $S$,
which is a graded map $A \overset{S}{\longrightarrow} A$, endowing it with a Hopf structure.
\end{proposition}
\begin{proof}
Let us try to define a ($\kk$-linear) left $\star$-inverse $S$ to $\id_A$
on each homogeneous component $A_n$, via induction on $n$.

In the base case $n=0$, Proposition~\ref{antipodes-are-antiendomorphisms} and its proof show that one
must define $S(\one)=\one$ so $S$ is the identity on $A_0=\kk$.

In the inductive step, recall from
Exercise~\ref{graded-connected-exercise}(h) that
a homogeneous element $a$ of degree $n > 0$
has $\Delta(a) = a \otimes \one  + \sum a'_1 \otimes a'_2$,
with each $\deg(a'_1) < n$.
(Here $\sum a'_1 \otimes a'_2$ stands for a sum of
tensors $a'_{1,k} \otimes a'_{2,k}$, with each $a'_{1,k}$
being homogeneous of degree $\deg(a'_{1,k}) < n$.
This is a slight variation on Sweedler notation.)
Hence in order to have $S \star \id_A = u \epsilon$,
one must define $S(a)$ in such a way that
$
S(a) \cdot \one + \sum S(a'_1) a'_2 = u \epsilon (a) = 0
$
and hence
$
S(a) := - \sum S(a'_1) a'_2,
$
where $S(a'_1)$ have already been uniquely defined
by induction (since $\deg(a'_{1,k}) < n$).
This does indeed define such a left $\star$-inverse $S$ to
$\id_A$, by induction.  It is also a graded map by induction.

The same argument shows how to define a right $\star$-inverse
$S'$ to $\id_A$.  Then $S=S'$ is a
two-sided $\star$-inverse to $\id_A$ by the associativity of $\star$.
\end{proof}

Here is another consequence of the fact that $S(\one)=\one$.

\begin{proposition}
\label{antipode-on-primitives}
In bialgebras, primitive elements $x$ have $\epsilon(x)=0$,
and in Hopf algebras, they have $S(x)=-x$.
\end{proposition}
\begin{proof}
In a bialgebra, $\epsilon(1)=1$.  Hence
$\Delta(x) = \one \otimes x + x \otimes \one$ implies
via \eqref{counit-diagram} that
$1 \cdot \epsilon(x) + \epsilon(1) x = x$, so $\epsilon(x) = 0$.
It also implies via \eqref{antipode-diagram} that
$S(x) \one + S(\one) x=u \epsilon(x)=u(0)=0$, so $S(x)=-x$.
\end{proof}

\noindent
Thus, whenever $A$ is a Hopf algebra generated as an algebra by its
primitive elements, $S$ is its unique $\kk$-algebra
anti-endomorphism that negates all primitive elements.

\begin{example}
\label{exa.tensor-alg.antipode}
The tensor and symmetric algebras
$T(V)$ and $\Sym(V)$ are each generated by $V$,
and each element of $V$ is primitive when regarded as an element
of either of them.  Hence one
has in $T(V)$ that
\begin{equation}
\label{antipode-in-tensor-algebra}
S(x_{i_1} x_{i_2} \cdots x_{i_k})=(-x_{i_k}) \cdots (-x_{i_2})(-x_{i_1}) =
(-1)^k x_{i_k} \cdots x_{i_2}x_{i_1}
\end{equation}
for each word $(i_1,\ldots,i_k)$ in the alphabet $I$
if $V$ is a free $\kk$-module with basis $\left\{x_i\right\}_{i\in I}$.
The same holds in $\Sym(V)$ for each multiset
$\left\{ i_1,\ldots,i_k \right\}_{\operatorname{multiset}}$,
recalling that the monomials are now commutative.  In other words,
for a commutative polynomial $f(x_1, x_2, \ldots, x_n)$ in
$\Sym(V)$, the antipode $S$ sends $f(x_1, x_2, \ldots, x_n)$
to $f(-x_1, -x_2, \ldots, -x_n)$,
negating all the variables.
\end{example}

The antipode for a connected graded Hopf algebra
has an interesting formula due to
Takeuchi \cite{Takeuchi}, reminiscent of P. Hall's formula for the M\"obius function of a
poset\footnote{In fact, for incidence Hopf algebras, Takeuchi's formula generalizes
Hall's formula-- see Corollary~\ref{P.Hall-formula}.}.  For the sake of stating
this, consider (for every $k \in \NN$) the $k$-fold \emph{tensor power}
$
A^{\otimes k}=A \otimes \cdots \otimes A
$
(defined in Example~\ref{exa.tensor-alg})
and define
\dfn{iterated multiplication and comultiplication} maps
\[
\xymatrixcolsep{3pc}
\xymatrix{
A^{\otimes k} \ar[r]^{m^{(k-1)}} & A
}
\qquad \text{and} \qquad
\xymatrix{
A \ar[r]^{\Delta^{(k-1)}} & A^{\otimes k}
}
\]
by induction over $k$, setting
$m^{(-1)}=u$, $\Delta^{(-1)}=\epsilon$, $m^{(0)}=\Delta^{(0)}=\id_A$,
and
\[
\begin{array}{rclr}
m^{(k)}&=m \circ (\id_A \otimes m^{(k-1)})
%    &=m \circ (m^{(k-1)} \otimes \id_A)
    & \qquad \text{ for every } k \geq 1; \\
\Delta^{(k)}&=(\id_A \otimes \Delta^{(k-1)}) \circ \Delta
%    &=(\Delta^{(k-1)} \otimes \id_A) \circ \Delta
    & \qquad \text{ for every } k \geq 1.
\end{array}
\]
Using associativity and coassociativity, one can see that for $k \geq 1$
these maps also satisfy
\[
\begin{array}{rclr}
m^{(k)} % &=m \circ (\id_A \otimes m^{(k-1)})
    &=m \circ (m^{(k-1)} \otimes \id_A)
    & \qquad \text{ for every } k \geq 1; \\
\Delta^{(k)} % &=(\id_A \otimes \Delta^{(k-1)}) \circ \Delta
    &=(\Delta^{(k-1)} \otimes \id_A) \circ \Delta
    & \qquad \text{ for every } k \geq 1
\end{array}
\]
(so we could just as well have used $\id_A \otimes m^{(k-1)}$ instead of
$m^{(k-1)} \otimes \id_A$ in defining them) and further symmetry properties
(see Exercise~\ref{exe.iterated.m} and Exercise~\ref{exe.iterated.Delta}).
They are how
one gives meaning to the right sides of these equations:
\begin{align*}
m^{(k)}(a^{(1)} \otimes \cdots \otimes a^{(k+1)})
&= a^{(1)} \cdots a^{(k+1)} ;\\
\Delta^{(k)}(b)
&= \sum b_1 \otimes \cdots \otimes b_{k+1} \text{ in Sweedler notation.}
\end{align*}

\begin{exercise}
\label{exe.iterated.m}
Let $A$ be a $\kk$-algebra. Let us define, for every $k \in \NN$, a
$\kk$-linear map
$m^{\left(k\right)} : A^{\otimes \left(k+1\right)} \to A$. Namely, we
define these maps by induction over $k$, with the induction base
$m^{\left(0\right)} = \id_A$, and with the induction step
$m^{\left(k\right)} = m \circ \left(\id_A \otimes m^{\left(k-1\right)}\right)$
for every $k \geq 1$. (This generalizes our definition of
$m^{\left(k\right)}$ for Hopf algebras $A$ given above, except for
$m^{\left(-1\right)}$ which we have omitted.)
\begin{itemize}
\item[(a)] Show that
$m^{\left(k\right)} = m \circ \left(m^{\left(i\right)} \otimes m^{\left(k-1-i\right)}\right)$
for every $k \geq 0$ and $0 \leq i \leq k-1$.
\item[(b)] Show that
$m^{\left(k\right)} = m \circ \left(m^{\left(k-1\right)} \otimes \id_A\right)$
for every $k \geq 1$.
\item[(c)] Show that
$m^{\left(k\right)} = m^{\left(k-1\right)} \circ
 \left(\id_{A^{\otimes i}} \otimes m \otimes \id_{A^{\otimes \left(k-1-i\right)}}\right)$
for every $k \geq 0$ and $0 \leq i \leq k-1$.
\item[(d)] Show that
$m^{\left(k\right)}
= m^{\left(k-1\right)} \circ \left(\id_{A^{\otimes \left(k-1\right)}} \otimes m\right)
= m^{\left(k-1\right)} \circ \left(m \otimes \id_{A^{\otimes \left(k-1\right)}}\right)$
for every $k \geq 1$.
\end{itemize}
\end{exercise}

\begin{exercise}
\label{exe.iterated.Delta}
Let $C$ be a $\kk$-coalgebra. Let us define, for every $k \in \NN$, a
$\kk$-linear map
$\Delta^{\left(k\right)} : C \to C^{\otimes \left(k+1\right)}$. Namely, we
define these maps by induction over $k$, with the induction base
$\Delta^{\left(0\right)} = \id_C$, and with the induction step
$\Delta^{\left(k\right)} = \left(\id_C \otimes \Delta^{\left(k-1\right)}\right) \circ \Delta$
for every $k \geq 1$. (This generalizes our definition of
$\Delta^{\left(k\right)}$ for Hopf algebras $A$ given above, except for
$\Delta^{\left(-1\right)}$ which we have omitted.)
\begin{itemize}
\item[(a)] Show that
$\Delta^{\left(k\right)} = \left(\Delta^{\left(i\right)} \otimes \Delta^{\left(k-1-i\right)}\right) \circ \Delta$
for every $k \geq 0$ and $0 \leq i \leq k-1$.
\item[(b)] Show that
$\Delta^{\left(k\right)} = \left(\Delta^{\left(k-1\right)} \otimes \id_C\right) \circ \Delta$
for every $k \geq 1$.
\item[(c)] Show that
$\Delta^{\left(k\right)} =
 \left(\id_{C^{\otimes i}} \otimes \Delta \otimes \id_{C^{\otimes \left(k-1-i\right)}}\right)
 \circ \Delta^{\left(k-1\right)}$
for every $k \geq 0$ and $0 \leq i \leq k-1$.
\item[(d)] Show that
$\Delta^{\left(k\right)}
= \left(\id_{C^{\otimes \left(k-1\right)}} \otimes \Delta\right) \circ \Delta^{\left(k-1\right)}
= \left(\Delta \otimes \id_{C^{\otimes \left(k-1\right)}}\right) \circ \Delta^{\left(k-1\right)}$
for every $k \geq 1$.
\end{itemize}
\end{exercise}

\begin{remark}
Exercise~\ref{exe.iterated.m} holds more generally for nonunital
associative algebras $A$ (that is, $\kk$-modules $A$ equipped with a
$\kk$-linear map $m : A \otimes A \to A$ such that the diagram
\eqref{associativity-diagram} is commutative, but not necessarily admitting
a unit map $u$). Similarly, Exercise~\ref{exe.iterated.Delta} holds for
non-counital coassociative coalgebras $C$. The existence of a unit in $A$,
respectively a counit in $C$, allows slightly extending these two exercises
by additionally introducing maps $m^{\left(-1\right)} = u : \kk \to A$ and
$\Delta^{\left(-1\right)} = \epsilon : C \to \kk$; however, not much is
gained from this extension.\footnote{The identity
$m^{\left(k\right)} = m \circ \left(\id_A \otimes m^{\left(k-1\right)}\right)$
for a $\kk$-algebra $A$ still holds when $k=0$ if it is interpreted in the
right way (viz., if $A$ is identified with $A \otimes \kk$ using the
canonical homomorphism).}
\end{remark}

% [DG][v19] Detailed the introduction of iterated multiplication and
% comultiplication maps and added the two exercises above to prove their
% properties.

\begin{exercise}
\label{exe.iterated.patras}
For every $k \in \NN$ and every $\kk$-bialgebra $H$, consider the map
$\Delta_H^{\left(k\right)} : H \to H^{\otimes \left(k+1\right)}$ (this
is the map $\Delta^{\left(k\right)}$ defined
as in Exercise~\ref{exe.iterated.Delta} for $C = H$), and the map
$m_H^{\left(k\right)} : H^{\otimes \left(k+1\right)} \to H$ (this is the
map $m^{\left(k\right)}$ defined as in Exercise~\ref{exe.iterated.m}
for $A = H$).

Let $H$ be a $\kk$-bialgebra. Let $k \in \NN$. Show that:\footnote{The
following statements are taken from \cite{Patras-descent}; specifically,
part (c) is \cite[Lem. 1.8]{Patras-descent}.}
\begin{itemize}
\item[(a)] The map
$m_H^{\left(k\right)} : H^{\otimes \left(k+1\right)} \to H$ is a
$\kk$-coalgebra homomorphism.
\item[(b)] The map
$\Delta_H^{\left(k\right)} : H \to H^{\otimes \left(k+1\right)}$ is a
$\kk$-algebra homomorphism.
\item[(c)] We have $m_{H^{\otimes \left(k+1\right)}}^{\left(\ell\right)}
\circ \left(\Delta^{\left(k\right)}_H\right)^{\otimes \left(\ell + 1\right)}
= \Delta^{\left(k\right)}_H \circ m^{\left(\ell\right)}_H$ for
every $\ell \in \NN$.
\item[(d)] We have
$\left(m_H^{\left(\ell\right)}\right)^{\otimes \left(k+1\right)}
\circ \Delta_{H^{\otimes \left(\ell + 1\right)}}^{\left(k\right)}
= \Delta^{\left(k\right)}_H \circ m^{\left(\ell\right)}_H$ for
every $\ell \in \NN$.
\end{itemize}
\end{exercise}

% [DG][v25] Added above exercise (which should chip away a little bit
% at the internal product theory, cf. its use in \cite{Patras-descent},
% although it seems to me that it is the easiest part of the proof).

The iterated multiplication and comultiplication maps allow explicitly
computing the convolution of multiple maps; the following formula will
often be used without explicit mention:

\begin{exercise}
\label{exe.convolution.k}
Let $C$ be a $\kk$-coalgebra, and $A$ be a $\kk$-algebra. Let $k\in\NN$.
Let $f_1, f_2, \ldots, f_k$ be $k$ elements of
$\Hom\left(C,A\right)$. Show that
\[
f_1 \star f_2 \star \cdots \star f_k
= m^{\left(k-1\right)}_A
  \circ \left(f_1 \otimes f_2 \otimes \cdots \otimes f_k\right)
  \circ \Delta^{\left(k-1\right)}_C.
\]
\end{exercise}

% [DG][v19] Added the above exercise and referenced it from the proof of
% Proposition \ref{Takeuchi-prop}.

We are now ready to state Takeuchi's formula for the antipode:

\begin{proposition}
\label{Takeuchi-prop}
In a connected graded Hopf algebra $A$, the
antipode has formula
\begin{equation}
\label{Takeuchi-formula}
\begin{aligned}
S&=\sum_{k \geq 0} (-1)^k m^{(k-1)} f^{\otimes k} \Delta^{(k-1)} \\
 &=u\epsilon - f + m \circ f^{\otimes 2} \circ \Delta
    - m^{(2)} \circ f^{\otimes 3} \circ \Delta^{(2)} + \cdots
\end{aligned}
\end{equation}
where $f:=\id_A - u\epsilon$ in $\End(A)$.
\end{proposition}
\begin{proof}
We argue as in \cite[proof of Lemma 14]{Takeuchi} or \cite[\S 5]{AguiarSottile}.
For any $f$ in $\End(A)$, the following explicit formula expresses
its $k$-fold convolution power $f^{\star k}:=f \star \cdots \star f$
in terms of its tensor powers $f^{\otimes k}:=f \otimes \cdots \otimes f$
(according to Exercise~\ref{exe.convolution.k}):
\[
f^{\star k} = m^{(k-1)} \circ f^{\otimes k} \circ \Delta^{(k-1)}.
\]
Therefore any $f$ annihilating $A_0$ will be
\dfn{locally $\star$-nilpotent} on $A$, meaning that for each $n$ one
has that $A_n$ is annihilated by $f^{\star m}$ for every $m > n$:  homogeneity forces that
for $a$ in $A_n$, every summand of $\Delta^{(m-1)}(a)$ must contain among its
$m$ tensor factors at least one factor lying in $A_0$,
so each summand is annihilated by $f^{\otimes m}$, and
$f^{\star m}(a)=0$.

In particular such $f$ have the property that
$u \epsilon + f$ has as two-sided $\star$-inverse
\begin{align*}
(u\epsilon + f)^{\star (-1)}
&= u\epsilon - f + f \star f - f \star f \star f + \cdots \\
& = \sum_{k \geq 0} (-1)^k f^{\star k}
= \sum_{k \geq 0} (-1)^k m^{(k-1)} \circ f^{\otimes k} \circ \Delta^{(k-1)}.
\end{align*}
The proposition follows upon taking
$f:=\id_A - u\epsilon$, which annihilates $A_0$.
\end{proof}

\begin{remark}
\label{non-graded-Takeuchi-remark}
In fact, one can see that Takeuchi's formula applies more generally to
define an antipode $A \overset{S}{\longrightarrow} A$
in any (not necessarily graded)
bialgebra $A$ where the map $\id_A - u\epsilon$ is locally
$\star$-nilpotent.

It is also worth noting that the proof of Proposition~\ref{Takeuchi-prop}
gives an alternate proof of
Proposition~\ref{graded-connected-bialgebras-have-antipodes}.
\end{remark}

To finish our discussion of antipodes, we mention some properties
(taken from \cite[Lemma 4.0.3]{Sweedler}) relating antipodes
to convolutional inverses.

\begin{proposition}
\label{antipodes-give-convolution-inverses}
Let $H$ be a Hopf algebra with antipode $S$.
\begin{enumerate}
\item[(a)] For any algebra $A$ and algebra morphism $H \overset{\alpha}{\rightarrow} A$,
one has $\alpha \circ S = \alpha^{\star -1}$, the convolutional inverse to $\alpha$ in $\Hom(H,A)$.
\item[(b)] For any coalgebra $C$ and coalgebra morphism $C \overset{\gamma}{\rightarrow} H$,
one has $S \circ \gamma = \gamma^{\star -1}$, the convolutional inverse to $\gamma$ in $\Hom(C,H)$.
\end{enumerate}
\end{proposition}
\begin{proof}
We prove (a); the proof of (b) is similar.

For assertion (a), note that 
Proposition~\ref{prop.convolution.functor}
(applied to $H$, $H$, $H$, $A$, $\id_H$ and
$\alpha$ instead of $C$, $C'$, $A$, $A'$, $\gamma$
and $\alpha$) shows that the map
\[
\Hom \left(H, H\right) \to \Hom \left(H, A\right),
\qquad
f \mapsto \alpha \circ f
\]
is a $\kk$-algebra homomorphism from the convolution
algebra $\left( \Hom \left(H, H\right), \star \right)$
to the convolution algebra
$\left( \Hom \left(H, A\right), \star \right)$.
Denoting this homomorphism by $\varphi$, we thus have
$\varphi\left((\id_H)^{\star -1}\right)
= \left(\varphi(\id_H)\right)^{\star -1}$
(since $\kk$-algebra homomorphisms preserve inverses).
Now,
\[
\alpha \circ S = \varphi(S)
               = \varphi\left((\id_H)^{\star -1}\right)
               = \left(\varphi(\id_H)\right)^{\star -1}
               = \left(\alpha \circ \id_H\right)^{\star -1}
               = \alpha^{\star -1}.
\]
\end{proof}

A rather useful consequence of Proposition~\ref{antipodes-give-convolution-inverses}
is the fact (\cite[Lemma 4.0.4]{Sweedler})
that a bialgebra morphism between Hopf algebras automatically respects the
antipodes:

\begin{corollary} \label{cor.bialg-mor-is-Hopf}
Let $H_1$ and $H_2$ be Hopf algebras with antipodes $S_1$ and $S_2$, respectively.
Then, any bialgebra morphism $H_1 \overset{\beta}{\rightarrow} H_2$ is a Hopf morphism\footnote{A \dfn{Hopf morphism} (or, more officially, a \dfn{Hopf algebra morphism}, or \dfn{homomorphism of Hopf algebras}) between two Hopf algebras $A$ and $B$ is defined to be a bialgebra morphism $f : A \to B$ that satisfies $f \circ S_A = S_B \circ f$.}, that is, it commutes with the antipodes (i.e., we have $\beta \circ S_1 = S_2 \circ \beta$).
\end{corollary}

\begin{proof}
Proposition~\ref{antipodes-give-convolution-inverses}(a) (applied to $H = H_1$, $S = S_1$, $A = H_2$ and $\alpha = \beta$)
yields $\beta \circ S_1 = \beta^{\star -1}$.
Proposition~\ref{antipodes-give-convolution-inverses}(b) (applied to $H = H_2$, $S = S_2$, $C = H_1$ and $\gamma = \beta$)
yields $S_2 \circ \beta = \beta^{\star -1}$.
Comparing these equalities shows that $\beta \circ S_1 = S_2 \circ \beta$, qed.
\end{proof}

% [DG][v80] Moved the above corollary out of the proposition.

\begin{exercise}
\label{exe.coalg-anti}
Prove that the antipode $S$ of a Hopf algebra $A$ is a
coalgebra anti-endomorphism, i.e., that it satisfies
$\epsilon \circ S = \epsilon$ and
$\Delta \circ S = T \circ \left( S \otimes S \right) \circ \Delta$,
where $T : A\otimes A \rightarrow A \otimes A$ is the
twist map sending every $a\otimes b$ to $b\otimes a$.
\end{exercise}

\begin{exercise}
\label{exe.coalg-anti.gen}
If $C$ is a $\kk$-coalgebra and if $A$ is a $\kk$-algebra, then a
$\kk$-linear map $f : C \rightarrow A$ is
said to be \dfn{$\star$-invertible} if it is invertible as an
element of the $\kk$-algebra
$\left( \Hom \left( C, A \right) , \star \right)$.
In this case, the multiplicative inverse $f^{\star\left( -1\right) }$
of $f$ in $\left( \Hom \left( C, A \right) , \star \right)$
is called the \dfn{$\star$-inverse} of $f$.

Recall the concepts introduced in Definition~\ref{def.anti-hom}.

\begin{itemize}
\item[(a)] If $C$ is a $\kk$-bialgebra, if $A$ is a $\kk$-algebra,
and if $r : C \rightarrow A$ is a $\star$-invertible $\kk$-algebra
homomorphism, then prove that the $\star$-inverse $r^{\star\left(
-1\right)  }$ of $r$ is a $\kk$-algebra anti-homomorphism.

\item[(b)] If $C$ is a $\kk$-bialgebra, if $A$ is a $\kk$-coalgebra,
and if $r : A\rightarrow C$ is a $\star$-invertible $\kk$-coalgebra
homomorphism, then prove that the $\star$-inverse $r^{\star\left(
-1\right)  }$ of $r$ is a $\kk$-coalgebra anti-homomorphism.

\item[(c)] Derive Proposition~\ref{antipodes-are-antiendomorphisms}
from Exercise~\ref{exe.coalg-anti.gen}(a), and derive
Exercise~\ref{exe.coalg-anti} from
Exercise~\ref{exe.coalg-anti.gen}(b).

\item[(d)] Prove
Corollary~\ref{commutative-implies-involutive-antipode-cor} again
using Proposition~\ref{antipodes-give-convolution-inverses}.

\item[(e)] If $C$ is a graded $\kk$-coalgebra, if $A$ is a graded
$\kk$-algebra, and if $r : C \to A$ is a $\star$-invertible
$\kk$-linear map that is graded, then prove that the
$\star$-inverse $r^{\star\left(  -1\right)  }$ of $r$ is also
graded.

\end{itemize}
\end{exercise}

% [DG][v48] Added preceding exercise (nothing particularly interesting,
% but can help in the rare cases when an antipode does not exist).

% [DG][v49] Added part (d).

% [DG][v80] Added part (e), since I realized I do use it tacitly
% in another solution.

\begin{exercise} \phantomsection
\label{exe.antipode.inverse}
\begin{itemize}

\item[(a)]
Let $A$ be a Hopf algebra. If
$P:A\rightarrow A$ is a $\kk$-linear map such that every $a\in A$
satisfies
\[
\sum_{\left( a\right)} P\left( a_2 \right) \cdot a_1
= u\left( \epsilon\left( a\right) \right) ,
\]
then prove that the antipode $S$ of $A$ is invertible and its inverse is $P$.

\item[(b)]
Let $A$ be a Hopf algebra. If
$P:A\rightarrow A$ is a $\kk$-linear map such that every $a\in A$
satisfies
\[
\sum_{\left( a\right)} a_2 \cdot P\left( a_1 \right)
= u\left( \epsilon\left( a\right) \right) ,
\]
then prove that the antipode $S$ of $A$ is invertible and its inverse is $P$.

\item[(c)]
Show that the antipode of a connected graded Hopf algebra is invertible.

\end{itemize}
\end{exercise}

(Compare this exercise to \cite[Lemma 1.5.11]{Montgomery}.)

% [DG][v3] Added two exercises above.

\begin{definition}
Let $C$ be a $\kk$-coalgebra. A \dfn{subcoalgebra} of $C$ means a
$\kk$-coalgebra $D$ such that $D \subset C$ and such that the canonical
inclusion map $D \to C$ is a $\kk$-coalgebra homomorphism\footnote{In
this definition, we follow \cite[p. 55]{Morris} and
\cite[\S 6.7]{Wisbauer}; other authors may use
other definitions.}. When $\kk$
is a field, we can equivalently define a subcoalgebra of $C$ as a
$\kk$-submodule $D$ of $C$ such that $\Delta_C\left(D\right)$ is
a subset of the $\kk$-submodule $D \otimes D$ of $C \otimes C$;
however, this might no longer be equivalent when $\kk$ is not a
field\footnote{This is because the $\kk$-submodule $D \otimes D$ of
$C \otimes C$ is generally not isomorphic to the $\kk$-module
$D \otimes D$. See \cite[p. 56]{Morris} for specific counterexamples
for the non-equivalence of the two notions of a subcoalgebra. Notice
that the equivalence is salvaged if
$D$ is a direct summand of $C$ as a $\kk$-module (see
Exercise~\ref{exe.subcoalgebra.summand} for this).}.

Similarly, a \dfn{subbialgebra} of a bialgebra $C$ is a
$\kk$-bialgebra $D$ such that $D \subset C$ and such that the canonical
inclusion map $D \to C$ is a $\kk$-bialgebra homomorphism.
Also, a \dfn{Hopf subalgebra} of a Hopf algebra $C$ is a
$\kk$-Hopf algebra $D$ such that $D \subset C$ and such that the canonical
inclusion map $D \to C$ is a $\kk$-Hopf algebra homomorphism.\footnote{By
Corollary~\ref{cor.bialg-mor-is-Hopf}, we can also
define it as a subbialgebra of $C$ that happens to be a Hopf algebra.}
\end{definition}

% [DG][v17] Added definition above.
% I am referencing the proceedings volume edited by Morris. Should
% I reference the Nichols-Sweedler paper instead?

\begin{exercise}
\label{exe.subcoalgebra.summand}
Let $C$ be a $\kk$-coalgebra. Let $D$ be a $\kk$-submodule of $C$
such that $D$ is a direct summand of $C$ as a $\kk$-module (i.e.,
there exists a $\kk$-submodule $E$ of $C$ such that $C = D \oplus E$).
(This is automatically satisfied if $\kk$ is a field.)
Assume that $\Delta\left(D\right) \subset C \otimes D$ and
$\Delta\left(D\right) \subset D \otimes C$. (Here, we are abusing
the notation $C \otimes D$ to denote the $\kk$-submodule of $C \otimes C$
spanned by tensors of the form $c \otimes d$ with $c \in C$ and
$d \in D$; similarly, $D \otimes C$ should be understood.)
Show that there is a
canonically defined $\kk$-coalgebra structure on $D$ which makes $D$
a subcoalgebra of $C$.
\end{exercise}

The next exercise is implicit in \cite[\S 5]{AguiarBergeronSottile}:

\begin{exercise}
\label{exe.subcoalgebra.ABS}
Let $\kk$ be a field. Let $C$ be a $\kk$-coalgebra, and let $U$ be any
$\kk$-module. Let $f : C \to U$ be a $\kk$-linear map. Recall the map
$\Delta^{\left(2\right)} : C \to C^{\otimes 3}$ from
Exercise~\ref{exe.iterated.Delta}. Let
$K = \ker\left(\left(\id_C \otimes f \otimes \id_C\right) \circ
\Delta^{\left(2\right)}\right)$.
\begin{itemize}
\item[(a)] Show that $K$ is a $\kk$-subcoalgebra of $C$.
\item[(b)] Show that every $\kk$-subcoalgebra of $C$ which is a subset
of $\ker f$ must be a subset of $K$.
\end{itemize}
\end{exercise}

% [DG][v19] Added the above two exercises. Feel free to make the second
% one into a theorem/proof if you write the Aguiar-Bergeron-Sottile II
% section.

\begin{exercise} \phantomsection
\label{exe.antipode.connected.general}
\begin{itemize}
\item[(a)]
Let $C = \bigoplus_{n\geq 0} C_n$ be a graded $\kk$-coalgebra, and $A$ be
any $\kk$-algebra. Notice that $C_0$ itself is a $\kk$-subcoalgebra of
$C$. Let $h : C \to A$ be a $\kk$-linear map such that the restriction
$h\mid_{C_0}$ is a $\star$-invertible map in $\Hom\left(C_0,A\right)$.
Prove that $h$ is a $\star$-invertible map in $\Hom\left(C,A\right)$. (This is
a weaker version of Takeuchi's \cite[Lemma 14]{Takeuchi}.)

\item[(b)]
Let $A = \bigoplus_{n\geq 0} A_n$ be a graded $\kk$-bialgebra. Notice
that $A_0$ is a subbialgebra of $A$. Assume that $A_0$ is a Hopf algebra.
Show that $A$ is a Hopf algebra.

\item[(c)]
Obtain yet another proof of
Proposition~\ref{graded-connected-bialgebras-have-antipodes}.
\end{itemize}
\end{exercise}

\begin{exercise}
\label{exe.I=0}
Let $A = \bigoplus_{n\geq 0} A_n$ be a connected graded
$\kk$-bialgebra. Let $\Liep$ be the $\kk$-submodule of $A$
consisting of the primitive elements of $A$.
\begin{itemize}
\item[(a)] If $I$ is a two-sided coideal of $A$ such that $I \cap \Liep = 0$ and
such that $I = \bigoplus_{n \geq 0}\left( I\cap A_n \right)$, then prove that
$I=0$.
\item[(b)] Let $f:A \to C$ be a graded surjective coalgebra homomorphism from
$A$ to a graded $\kk$-coalgebra $C$. If $f\mid_{\Liep}$ is
injective, then prove that $f$ is injective.
\item[(c)] Assume that $\kk$ is a field. Show that the claim of
Exercise~\ref{exe.I=0}(b) is valid even without requiring $f$ to be surjective.
\end{itemize}
\end{exercise}

\begin{remark}
Exercise~\ref{exe.I=0} (b) and (c) are often used in order to prove that certain
coalgebra homomorphisms are injective.

The word ``bialgebra'' can be replaced by ``coalgebra'' in Exercise~\ref{exe.I=0},
provided that the
notion of a connected graded coalgebra is defined correctly (namely, as a
graded coalgebra such that the restriction of $\epsilon$ to the $0$-th graded
component is an isomorphism), and the notion of the element $1$ of a
connected graded coalgebra is defined accordingly (namely, as the preimage of
$1\in\kk$ under the restriction of $\epsilon$ to the $0$-th graded component).
\end{remark}

% [DG][v13] Added this exercise. The solution comes from an old mail
% of mine and is really a mess; sorry.

% [DG][v23] Added part (c) to this exercise.

\subsection{Commutativity, cocommutativity}

Recall that a $\kk$-algebra $A$ is \dfn{commutative}
if and only if all $a, b \in A$ satisfy $ab=ba$.
Here is a way to restate this classical definition using
tensors instead of pairs of elements:

\begin{definition} \label{def.comm}
A $\kk$-algebra $A$ is said to be \dfn{commutative} if
% $ab=ba$ for all $a, b \in A$, that is,
the following diagram commutes:
\begin{equation}
\label{commutativity-diagram}
\xymatrix{
A \otimes A \ar[rr]^T \ar[dr]_m& &A \otimes A \ar[dl]^m\\
 &A&
}
\end{equation}
where $T$ is the twist map $T_{A,A}$ (see
Definition~\ref{def.anti-hom}(a) for its definition).
% This, of course, is a mere restatement of the classical definition
% of a commutative $\kk$-algebra using tensors instead of pairs of
% elements.
% % [DG][v62] Added the preceding sentence.
\end{definition}

Having thus redefined commutative algebras in terms of
tensors and linear maps, we can dualize this definition
(reversing all arrows)
and obtain the notion of \emph{cocommutative coalgebras}:

\begin{definition} \label{def.cocomm}
A $\kk$-coalgebra $C$ is said to be \dfn{cocommutative} if
the following diagram commutes:
\begin{equation}
\label{cocommutativity-diagram}
\xymatrix{
C \otimes C \ar[rr]^T & &C \otimes C \\
 &C \ar[ul]^\Delta \ar[ur]_\Delta&
}
\end{equation}
where $T$ is the twist map $T_{C,C}$ (see
Definition~\ref{def.anti-hom}(a) for its definition).
\end{definition}

% [DG][v80] Split the above definition into two.

\begin{example}
Group algebras $\kk G$ are always cocommutative.
They are commutative if and only if $G$ is abelian or $\kk = 0$.

Tensor algebras $T(V)$ are always cocommutative,
but not generally commutative\footnote{If $\kk$ is a field,
then $T(V)$ is commutative if and only if $\dim_\kk V \leq 1$.}.

% [DG][v17] Weakened the "if and only if" since $\kk$ is an arbitrary
% commutative ring now.

Symmetric algebras $\Sym(V)$ are always cocommutative and commutative.

Homology and cohomology of $H$-spaces are always cocommutative and commutative
\emph{in the topologist's sense} where one reinterprets that
twist map $A \otimes A \overset{T}{\rightarrow} A \otimes A$ to have
the extra sign as in \eqref{topologist-twist}.
\end{example}

% [DG][v25] Is the kind of (co)commutativity given by
% your definition (reinterpreted using the topologist's twist)
% really what passes for (co)commutativity in the homological
% setting? I remember people asking for more (for commutativity,
% there is an axiom $xx = 0$ for $x$ odd; for cocommutativity,
% don't know what). Is this still true in the homological
% setting?

Note how the cocommutative Hopf algebras $T(V), \Sym(V)$
have much of their structure controlled by their $\kk$-submodules $V$,
which consist of primitive elements only (although, in general, not of
all their primitive elements).  This is not far from the truth in general,
and closely related to Lie algebras.

% [DG][v17] Replaced "their subspace $V$ of primitive elements" by
% "their $\kk$-submodules $V$, which consist of primitive elements only
% (although, in general, not of all their primitive elements)" to get rid
% of ambiguity.

\begin{exercise}
\label{exe.liealg}
Recall that a \dfn{Lie algebra} over $\kk$ is a $\kk$-module
$\Lieg$ with a $\kk$-bilinear map $[\cdot,\cdot]:\Lieg\times\Lieg\to\Lieg$
that satisfies $[x,x]=0$ for $x$ in $\Lieg$, and the \dfn{Jacobi identity}
\begin{align*}
[x,[y,z]]&=[[x,y],z] + [y,[x,z]], \text{ or equivalently }\\
[x,[y,z]]&+[z,[x,y]]+[y,[z,x]]=0
\end{align*}
for all $x, y, z \in \Lieg$.
This $\kk$-bilinear map $[\cdot,\cdot]$ is called the \dfn{Lie
bracket} of $\Lieg$.
\begin{itemize}
\item[(a)]
Check that any associative algebra $A$ gives rise to a
Lie algebra by means of the commutator operation $[a,b]:=ab-ba$.
\item[(b)]
If $A$ is also a bialgebra, show that the $\kk$-submodule
of primitive elements $\Liep \subset A$ is closed under the Lie
bracket, that is, $[\Liep, \Liep] \subset \Liep$, and hence forms
a Lie subalgebra.
\end{itemize}

% [DG][v9] Replaced "$[\cdot,\cdot]$" by
% "$[\cdot,\cdot]:\Lieg\times\Lieg\to\Lieg$".

Conversely, given a Lie algebra $\Liep$,
one constructs the \dfn{universal enveloping algebra}
$\UUU(\Liep):=T(\Liep)/J$\index{$\UUU(\Liep)$}
as the quotient of the tensor algebra
$T(\Liep)$ by the two-sided ideal $J$ generated by all elements
$xy-yx-[x,y]$ for $x,y$ in $\Liep$.
\begin{itemize}
\item[(c)]
Show that $J$ is also a two-sided coideal in $T(\Liep)$ for
its usual coalgebra structure, and hence the quotient
$\UUU(\Liep)$ inherits the structure of a cocommutative bialgebra.
\item[(d)]
Show that the antipode $S$ on $T(\Liep)$ preserves $J$, meaning
that $S(J) \subset J$, and hence $\UUU(\Liep)$ inherits the
structure of a (cocommutative) Hopf algebra.
\end{itemize}
\end{exercise}

There are theorems,
discussed in \cite[\S3.8]{Cartier}, \cite[Chap. 5]{Montgomery},
\cite[\S 3.2]{DuchampMinhTolluChienNghia}
giving various mild hypotheses in addition to  cocommutativity
which imply that the inclusion of the $\kk$-module $\Liep$ of primitives in
a Hopf algebra $A$ extends to a Hopf isomorphism
$\UUU(\Liep) \cong A$.

% [DG][v13] Added {DuchampMinhTolluChienNghia} reference since it is
% one of the few to cover base rings which are not fields. (Another
% is Fresse's operads book, but it is still in the process of being
% written and the proof of Cartier-Milnor-Moore will be rewritten.)

\begin{exercise}
\label{exe.commutative.convolution}
Let $C$ be a cocommutative $\kk$-coalgebra.
Let $A$ be a commutative $\kk$-algebra.
Show that the convolution algebra
$\left( \Hom \left( C, A \right) , \star \right)$ is commutative
(i.e., every $f, g \in \Hom \left( C, A \right)$
satisfy $f \star g = g \star f$).
\end{exercise}

% [DG][v62] Added the simple exercise above.

\begin{exercise} \phantomsection
\label{exe.cocommutative.coalghom}
\begin{itemize}
\item[(a)] Let $C$ be a $\kk$-coalgebra. Show that $C$ is cocommutative if and
only if its comultiplication $\Delta_C : C \to C \otimes C$ is
a $\kk$-coalgebra homomorphism.
\item[(b)] Let $A$ be a $\kk$-algebra. Show that $A$ is commutative if and
only if its multiplication $m_A : A \otimes A \to A$ is
a $\kk$-algebra homomorphism.
\end{itemize}
\end{exercise}

% [DG][v14] Added this exercise.

\begin{remark}
If $C$ is a $\kk$-coalgebra, then $\epsilon_C : C \to \kk$ is always a
$\kk$-coalgebra homomorphism. Similarly, $u_A : \kk \to A$ is a
$\kk$-algebra homomorphism whenever $A$ is a $\kk$-algebra.
\end{remark}

% [DG][v25] Added above remark to accompany the exercise just above it.

\begin{exercise} \phantomsection
\label{exe.comm-cocomm.anti}
\begin{itemize}

\item[(a)] Let $A$ and $B$ be two $\kk$-algebras, at least
one of which is commutative.
Prove that the $\kk$-algebra anti-homomorphisms from $A$ to $B$
are the same as the $\kk$-algebra homomorphisms from $A$ to $B$.

\item[(b)] State and prove the dual of this result.

\end{itemize}
\end{exercise}

% [DG][v80] Added the above near-trivial exercise as its result has
% been tacitly used a few times.

\begin{exercise}
\label{exe.iterated.comm}
Let $A$ be a commutative $\kk$-algebra, and let $k \in \NN$. The
symmetric group $\Symm_k$ acts on the $k$-fold tensor power $A^{\otimes k}$
by permuting the tensor factors:
$\sigma\left(v_1 \otimes v_2 \otimes \cdots \otimes v_k\right)
= v_{\sigma^{-1}\left(1\right)} \otimes v_{\sigma^{-1}\left(2\right)}
  \otimes \cdots \otimes v_{\sigma^{-1}\left(k\right)}$
for all $v_1, v_2, \ldots, v_k \in A$ and $\sigma \in \Symm_k$.
For every $\pi \in \Symm_k$, denote
by $\rho\left(\pi\right)$ the action of $\pi$ on $A^{\otimes k}$ (this is
an endomorphism of $A^{\otimes k}$).
Show that every $\pi \in \Symm_k$ satisfies
$m^{\left(k-1\right)} \circ \left(\rho\left(\pi\right)\right)
= m^{\left(k-1\right)}$. (Recall that
$m^{\left(k-1\right)} : A^{\otimes k} \to A$
is defined as in Exercise~\ref{exe.iterated.m} for $k \geq 1$, and by
$m^{\left(-1\right)} = u : \kk \to A$ for $k = 0$.)
\end{exercise}

\begin{exercise}
\label{exe.iterated.cocomm}
State and solve the analogue of Exercise~\ref{exe.iterated.comm} for
cocommutative $\kk$-coalgebras.
\end{exercise}

% [DG][v19] Added above two exercises.

\begin{exercise} \phantomsection
\label{exe.alghoms.comm}
\begin{itemize}
\item[(a)] If $H$ is a $\kk$-bialgebra and $A$ is a commutative
$\kk$-algebra, and if $f$ and $g$ are two $\kk$-algebra homomorphisms
$H \to A$, then prove that $f \star g$ also is a $\kk$-algebra
homomorphism $H \to A$.
\item[(b)] If $H$ is a $\kk$-bialgebra and $A$ is a commutative
$\kk$-algebra, and if $f_1, f_2, \ldots, f_k$ are several
$\kk$-algebra homomorphisms $H \to A$, then prove that
$f_1 \star f_2 \star \cdots \star f_k$ also is a $\kk$-algebra
homomorphism $H \to A$.
\item[(c)] If $H$ is a Hopf algebra and $A$ is a commutative
$\kk$-algebra, and if $f : H \to A$ is a $\kk$-algebra
homomorphism, then prove that $f \circ S : H \to A$ (where $S$
is the antipode of $H$) is again a $\kk$-algebra homomorphism,
and is a $\star$-inverse to $f$.
\item[(d)] If $A$ is a commutative $\kk$-algebra, then show that
$m^{\left(k\right)}$ is a $\kk$-algebra homomorphism for every
$k \in \NN$. (The map
$m^{\left(k\right)} : A^{\otimes \left(k+1\right)} \to A$
is defined as in Exercise~\ref{exe.iterated.m}.)
\item[(e)] If $C^\prime$ and $C$ are two $\kk$-coalgebras, if
$\gamma : C \to C^\prime$ is a $\kk$-coalgebra homomorphism,
if $A$ and $A^\prime$ are two $\kk$-algebras, if
$\alpha : A \to A^\prime$ is a $\kk$-algebra homomorphism, and
if $f_1, f_2, \ldots, f_k$ are several $\kk$-linear maps
$C^\prime \to A$, then prove that
\[
\alpha \circ \left(f_1 \star f_2 \star \cdots \star f_k\right)
\circ \gamma
= \left(\alpha \circ f_1 \circ \gamma\right)
\star \left(\alpha \circ f_2 \circ \gamma\right)
\star \cdots
\star \left(\alpha \circ f_k \circ \gamma\right) .
\]
\item[(f)] If $H$ is a commutative $\kk$-bialgebra, and $k$
and $\ell$ are two nonnegative integers, then prove that
$\id_H^{\star k} \circ \id_H^{\star \ell}
= \id_H^{\star \left(k\ell\right)}$.
\item[(g)] If $H$ is a commutative $\kk$-Hopf algebra, and $k$
and $\ell$ are two integers, then prove that
$\id_H^{\star k} \circ \id_H^{\star \ell}
= \id_H^{\star \left(k\ell\right)}$. (These powers
$\id_H^{\star k}$, $\id_H^{\star \ell}$ and 
$\id_H^{\star \left(k\ell\right)}$ are well-defined
since $\id_H$ is $\star$-invertible.)
\item[(h)] State and prove the duals of parts (a)--(g) of this
exercise.
\end{itemize}
\end{exercise}

\begin{remark}
The maps $\id_H^{\star k}$ for $k \in \NN$ are known as the
\dfn{Adams operators} of the bialgebra $H$; they are studied,
inter alia, in \cite{AguiarLauve-2014}. Particular cases
(and variants) of
Exercise~\ref{exe.alghoms.comm}(f) appear in
\cite[Corollaire II.9]{Patras-descent} and
\cite[Theorem 1]{GerstenhaberSchack-shuffle}.
Exercise~\ref{exe.alghoms.comm}(f) and its dual are
\cite[Prop. 1.6]{Loday-serie-Hausdorff}.
\end{remark}

% [DG][v25] Added this exercise (part (a) of which I use later).
% 
% The convolution powers of the identity in part (f) are called the
% Adams operators/-ions; you might have some words to say about why
% (I don't know enough geometry).

\begin{exercise}
\label{exe.cocomm-antipode-S2}
Prove that the antipode $S$ of a cocommutative Hopf algebra
$A$ satisfies $S^2 = \id_A$. (This was a statement made in
Remark~\ref{rmk.commutative-implies-involutive-antipode-cor.cocomm}.)
\end{exercise}

% [DG][v80] Added this exercise. (It used to be part (e) in
% Exercise~\ref{exe.coalg-anti.gen} before, but then I realized
% that "cocommutative" was not defined at the point where
% that exercise was posed, so this part was premature.)

\begin{exercise}
\label{exe.dynkin}
Let $A$ be a cocommutative graded Hopf algebra with antipode
$S$. Define a $\kk$-linear map $E : A\rightarrow A$ by having
$E\left( a \right) = \left( \deg a \right) \cdot a$
for every homogeneous element $a$ of $A$.

\begin{itemize}
\item[(a)] Prove that for every $a\in A$, the elements $\left(  S\star
E\right)  \left(  a\right)  $ and $\left(  E\star S\right)  \left(  a\right)
$ (where $\star$ denotes convolution in $\Hom\left(A,A\right)  $) are
primitive.

\item[(b)] Prove that for every primitive $p\in A$, we have
$ \left(  S\star E\right)  \left(  p\right) 
= \left(  E\star S\right)  \left(  p\right)
= E\left(  p\right)  $.

\item[(c)] Prove that for every $a\in A$ and every primitive $p\in A$, we
have $\left(  S\star E\right)  \left(  ap\right)  =\left[  \left(  S\star
E\right)  \left(  a\right)  ,p\right]  +\epsilon\left(  a\right)  E\left(
p\right)  $, where $\left[  u,v\right] $ denotes the commutator $uv-vu$ of
$u$ and $v$.

\item[(d)] If $A$ is connected and $\QQ$ is a subring of $\kk$,
prove that the $\kk$-algebra $A$
is generated by the $\kk$-submodule $\Liep$ consisting of the
primitive elements of $A$.

\item[(e)] Assume that $A$ is the tensor algebra $T\left( V \right)$
of a $\kk$-module $V$, and that the $\kk$-submodule
$V = V^{\otimes 1}$ of $T\left( V \right)$ is the degree-$1$
homogeneous component of $A$.
Show that
$\left( S\star E\right)  \left( x_1 x_2 \ldots x_n \right)
= \left[ \ldots\left[ \left[ x_1, x_2 \right] , x_3 \right] ,
\ldots, x_n \right]  $ for any $n \geq 1$ and any
$x_1, x_2, \ldots , x_n \in V$.
\end{itemize}
\end{exercise}

% [DG][v13] I added part (d) to the above exercise.

% [DG][v37] Added grading condition in (e).

\begin{remark}
\label{rmk.dynkin}
Exercise~\ref{exe.dynkin} gives rise to a certain idempotent
map $A\rightarrow A$ when $\kk$ is a commutative $\QQ$-algebra
and $A$ is a cocommutative connected graded $\kk$-Hopf algebra.
Namely, the $\kk$-linear map
$A \rightarrow A$ sending every homogeneous $a \in A$ to
$\frac{1}{\deg a} \left( S\star E \right) \left( a \right)$
(or $0$ if $\deg a = 0$) is
idempotent and is a projection on the $\kk$-module of primitive
elements of $A$. It is
called the \dfn{Dynkin idempotent}; see
\cite{PatrasReutenauer} for more of its properties.\footnote{We
will see another such idempotent in
Exercise~\ref{exe.eulerian-idp}.}
Part (c) of the exercise is more or less Baker's identity.
\end{remark}

% [DG][v3] Added exercise above.

\subsection{\label{subsect.duals}Duals}

Recall that for $\kk$-modules $V$, taking the dual $\kk$-module
$V^* :=\Hom(V,\kk)$ reverses $\kk$-linear maps.  That is, every
$\kk$-linear map $V \overset{\varphi}{\rightarrow} W$ induces
an \dfn{adjoint map}
$W^* \overset{\varphi^*}{\rightarrow} V^*$ defined uniquely by
\[
( f, \varphi(v) ) = ( \varphi^*(f), v )
\]
in which $(f,v)$ is the bilinear pairing $V^* \times V \rightarrow \kk$
sending $(f,v) \mapsto f(v)$.  If $V$ and $W$ are finite free
$\kk$-modules\footnote{A $\kk$-module is said to be \dfn{finite
free} if it has a finite basis. If $\kk$ is a field, then a
finite free $\kk$-module is the same as a finite-dimensional
$\kk$-vector space.}, more can be said: When $\varphi$ is
expressed in terms of a basis $\{v_i\}_{i \in I}$ for $V$
and a basis $\{w_j\}_{j \in J}$ for $W$ by some matrix,
the map $\varphi^*$ is expressed by the transpose matrix in terms
of the dual bases of these two bases\footnote{If
$\{v_i\}_{i \in I}$ is a basis of a finite free $\kk$-module $V$,
then the \dfn{dual basis} of this basis is defined as the basis
$\{f_i\}_{i \in I}$ of $V^*$ that satisfies
$(f_i,v_j)=\delta_{i,j}$ for all $i$ and $j$. (Recall that
$\delta_{i,j}$ is the Kronecker delta: $\delta_{i,j}= 1$ if $i=j$
and $0$ else.)}.

% [DG][v19] Revamped the last sentence.

The correspondence $\varphi \mapsto \varphi^*$ between
$\kk$-linear maps $V \overset{\varphi}{\rightarrow} W$ and
$\kk$-linear maps $W^* \overset{\varphi^*}{\rightarrow} V^*$ is
one-to-one when $W$ is finite free. However, this is not the case
in many combinatorial situations (in which $W$ is usually free but
not finite free). Fortunately, many of the good properties of
finite free
modules carry over to a certain class of graded modules as long
as the dual $V^*$ is replaced by a smaller module $V^o$ called the
graded dual. Let us first introduce the latter:

When $V=\bigoplus_{n \geq 0} V_n$ is a graded $\kk$-module,
note that the dual $V^* = \prod_{n\geq 0} (V_n)^*$ can contain
functionals $f$ supported on infinitely many $V_n$. However, we
can consider the
$\kk$-submodule $V^o:=\bigoplus_{n \geq 0} (V_n)^* \subset
\prod_{n\geq 0} (V_n)^* = V^*$\index{$V^o$},
sometimes called the \dfn{graded dual}\footnote{Do not mistake
this for the coalgebraic restricted dual $A^\circ$ of
\cite[\S 6.0]{Sweedler}.},
consisting of the functions $f$ that vanish on all but
finitely many $V_n$.
Notice that $V^o$ is graded, whereas $V^*$ (in general) is not.
If $V \overset{\varphi}{\rightarrow} W$ is a graded $\kk$-linear
map, then the adjoint map
$W^* \overset{\varphi^*}{\rightarrow} V^*$
restricts to a graded $\kk$-linear map
$W^o \rightarrow V^o$, which we (abusively) still denote by
$\varphi^*$.

A graded $\kk$-module $V=\bigoplus_{n \geq 0} V_n$ is said to be
\dfn{of finite type} if each $V_n$ is a finite free
$\kk$-module\footnote{This meaning of ``finite type'' can differ
from the standard one.}.
When the graded $\kk$-module $V$ is of finite type,
the graded $\kk$-module $V^o$ is
again of finite type\footnote{More precisely: Let
$V=\bigoplus_{n \geq 0} V_n$ be of finite type, and let
$\left\{v_i\right\}_{i \in I}$ be a \dfn{graded basis} of $V$,
that is, a basis of the $\kk$-module $V$ such that
the indexing set $I$ is partitioned into subsets
$I_0, I_1, I_2, \ldots$ (which are allowed to be empty) with the
property that, for every $n \in \NN$, the subfamily
$\left\{v_i\right\}_{i \in I_n}$ is a basis of the $\kk$-module
$V_n$. Then, we can define a family $\left\{f_i\right\}_{i \in I}$
of elements of $V^o$ by setting
$\left(f_i, v_j\right) = \delta_{i,j}$ for all $i, j \in I$.
This family $\left\{f_i\right\}_{i \in I}$ is a graded basis of
the graded $\kk$-module $V^o$. (Actually, for every $n \in \NN$,
the subfamily $\left\{f_i\right\}_{i \in I_n}$ is a basis of
the $\kk$-submodule $\left(V_n\right)^\ast$ of $V^o$ -- indeed
the dual basis to the basis $\left\{v_i\right\}_{i \in I_n}$ of
$V_n$.) This basis $\left\{f_i\right\}_{i \in I}$ is said to be
the \dfn{dual basis} to the basis $\left\{v_i\right\}_{i \in I}$
of $V$.}
and satisfies $\left(V^o\right)^o \cong V$.
Many other properties of finite free modules are salvaged in this
situation; most importantly:
The correspondence $\varphi \mapsto \varphi^*$ between graded
$\kk$-linear maps $V \rightarrow W$ and graded
$\kk$-linear maps $W^o \rightarrow V^o$ is
one-to-one when $W$ is of finite type\footnote{Only $W$ has
to be of finite type here; $V$ can be any graded $\kk$-module.}.

% [DG][v7] Replaced "which is again a graded vector space" by
% "This $V^o$ is again a graded vector space" to disambiguate grammar.

% [DG][v17] The above part of this section rewritten in a
% ring-friendly and more explicit (and leisurely?) way.
% 
% Also, renamed "restricted dual" by "graded dual" because I have
% seen the Sweedler "zero-dual" being referred to as restricted
% dual, and this is a different beast.
% 
% Shouldn't the statement about matrices require bases for both
% $V$ and $W$?

% [DG][v48] Introduced dual basis of $V^o$ in the finite-type
% setting in a footnote. (The concept is used a lot.)

Reversing the diagrams should then make it clear that, in the
finite free or finite-type situation, duals of algebras are
coalgebras, and vice-versa,
and duals of bialgebras or Hopf algebras are bialgebras or Hopf algebras.
For example, the product in a Hopf algebra $A$ of finite type uniquely defines
the coproduct of $A^o$ via adjointness:
\[
\left(\Delta_{A^o}(f) , a \otimes b\right)_{A \otimes A}
 = (f , a b)_{A}.
\]
Thus if $A$ has a basis $\{a_i\}_{i \in I}$ with \dfn{product structure
constants} $\{c^i_{j,k}\}$, meaning
\[
a_j a_k = \sum_{i \in I} c^i_{j,k} a_i,
\]
then the dual basis $\{f_i\}_{i \in I}$ has the same  $\{c^i_{j,k}\}$ as its
\dfn{coproduct structure constants}:
\[
\Delta_{A^o}(f_i) = \sum_{(j,k) \in I \times I} c^i_{j,k} f_j \otimes f_k.
\]
The assumption that $A$ be of finite type was indispensable here; in general,
the dual of a $\kk$-algebra does not become a $\kk$-coalgebra. However,
the dual of a $\kk$-coalgebra still becomes a $\kk$-algebra, as shown in the
following exercise:

\begin{exercise}
\label{exe.dual-algebra}
For any two $\kk$-modules $U$ and $V$, let $\rho_{U,V} :
U^{\ast} \otimes V^{\ast} \to \left(U \otimes V\right)^{\ast}$ be the
$\kk$-linear map which sends every tensor
$f \otimes g \in U^{\ast} \otimes V^{\ast}$ to the composition
$U \otimes V \overset{f \otimes g} {\longrightarrow}
\kk \otimes \kk \overset{m_\kk} {\longrightarrow} \kk$ of the
map\footnote{Keep in mind that the \emph{tensor}
$f \otimes g \in U^{\ast} \otimes V^{\ast}$ is not the same as the
\emph{map} $U \otimes V \overset{f \otimes g} {\longrightarrow}
\kk \otimes \kk$.} $f \otimes g$ with the canonical isomorphism
$\kk \otimes \kk \overset{m_\kk} {\longrightarrow} \kk$. When $\kk$ is a
field and $U$ is finite-dimensional, this map
$\rho_{U,V}$ is a $\kk$-vector space isomorphism (and usually regarded
as the identity); more generally, it is injective whenever $\kk$ is a
field\footnote{Over arbitrary rings it does not have to be even that!}.
% http://mathoverflow.net/questions/56255/duals-and-tensor-products has
% an example for the footnote.
Also, let $s : \kk \to \kk^*$ be the canonical isomorphism.
Prove that:
\begin{itemize}
\item[(a)] If $C$ is a $\kk$-coalgebra, then $C^*$ becomes a
$\kk$-algebra if we define its associative operation by
$m_{C^*} = \Delta_C^* \circ \rho_{C,C} : C^* \otimes C^* \to C^*$ and
its unit map to be $\epsilon_C^* \circ s : \kk \to C^*$.
\ \ \ \ \footnote{If $C$ is a finite free $\kk$-module, then this
$\kk$-algebra structure is the same as the one defined above by
adjointness. But the advantage of the new definition is that it
works even if $C$ is not a finite free $\kk$-module.}
\item[(b)] The $\kk$-algebra structure defined on $C^*$ in
part (a) is precisely the one defined on $\Hom\left(C,\kk\right) = C^*$
in Definition~\ref{convolution-algebra} applied to $A=\kk$.
\item[(c)] If $C$ is a graded $\kk$-coalgebra, then $C^o$ is
a $\kk$-subalgebra of the $\kk$-algebra $C^*$ defined in part (a).
\item[(d)] If $f : C \to D$ is a homomorphism of $\kk$-coalgebras, then
$f^* : D^* \to C^*$ is a homomorphism of $\kk$-algebras.
\item[(e)] Let $U$ be a graded $\kk$-module (not necessarily of
finite type), and let $V$ be a graded $\kk$-module of finite type.
Then, there is a 1-to-1 correspondence between graded $\kk$-linear maps
$U \to V$ and graded $\kk$-linear maps $V^o \to U^o$ given by
$f \mapsto f^*$.
\item[(f)] Let $C$ be a graded $\kk$-coalgebra (not necessarily of
finite type), and let $D$ be a graded $\kk$-coalgebra of finite type.
Part (e) of this exercise shows that
there is a 1-to-1 correspondence between graded $\kk$-linear maps
$C \to D$ and graded $\kk$-linear maps $D^o \to C^o$ given by
$f \mapsto f^*$. This correspondence has the property that a given
graded $\kk$-linear map $f : C \to D$ is a $\kk$-coalgebra morphism
if and only if $f^* : D^o \to C^o$ is a $\kk$-algebra morphism.
\end{itemize}
\end{exercise}

% [DG][v14] Added the above exercise, and referenced it from Section 6
% in which it *now* is used.

Another example of a Hopf algebra is provided by the so-called shuffle
algebra. Before we introduce it, let us define the \emph{shuffles} of
two words:

\begin{definition}
\label{shuffles}
Given two words $a=\left(  a_{1},a_{2}, \ldots ,a_{n}\right)  $ and
$b=\left(  b_{1},b_{2}, \ldots ,b_{m}\right)  $, the \emph{multiset of shuffles
of $a$ and $b$}\index{shuffle of words} is defined as the multiset
\[
\left\{  \left(  c_{w\left(  1\right)  },c_{w\left(  2\right)
}, \ldots ,c_{w\left(  n+m\right)  }\right)  \ : 
\ w\in\operatorname{Sh}_{n,m}\right\}_{\text{multiset}}  ,
\]
where $\left(  c_{1},c_{2}, \ldots ,c_{n+m}\right)  $ is the concatenation $a\cdot
b=\left(  a_{1},a_{2}, \ldots ,a_{n},b_{1},b_{2}, \ldots ,b_{m}\right)  $, and where
\dfn{$\operatorname{Sh}_{n,m}$} is the subset\footnote{\textbf{Warning:} This
definition of $\operatorname{Sh}_{n,m}$ is highly nonstandard, and many
authors define $\operatorname{Sh}_{n,m}$ to be the set of the inverses
of the permutations belonging to what we call $\operatorname{Sh}_{n,m}$.}
\[
\left\{  w\in\Symm_{n+m}\ : \ w^{-1}\left(  1\right)  <w^{-1}\left(  2\right)
<\cdots<w^{-1}\left(  n\right)  ;\ w^{-1}\left(  n+1\right)  <w^{-1}\left(  n+2\right)
<\cdots<w^{-1}\left(  n+m\right)  \right\}
\]
of the symmetric group $\Symm_{n+m}$. Informally speaking, the shuffles of
the two words $a$ and $b$ are the words obtained by overlaying the words $a$
and $b$, after first moving their letters apart so that no letters get
superimposed when the words are overlayed\footnote{For instance, if
$a = \left(1,3,2,1\right)$ and $b = \left(2,4\right)$, then
the shuffle $\left(1,2,3,2,4,1\right)$ of $a$ and $b$ can be obtained by
moving the letters of $a$ and $b$ apart as follows:
\[
\begin{array}{ccccccc}
a = & 1 &   & 3 & 2 &   & 1 \\
b = &   & 2 &   &   & 4 &   
\end{array}
\]
and then overlaying them to obtain
$\begin{array}{cccccc} 1 & 2 & 3 & 2 & 4 & 1 \end{array}$. Other ways of
moving letters apart lead to further shuffles (not always distinct).}.
In particular, any shuffle of
$a$ and $b$ contains $a$ and $b$ as subsequences. The multiset of shuffles of $a$ and
$b$ has $\binom{m+n}{n}$ elements
(counted with multiplicity) and is denoted by $\shuf{a}{b}$.
For instance, the shuffles of $\left(1,2,1\right)$
and $\left(3,2\right)$ are
\begin{align*}
\left({\underline 1},{\underline 2},{\underline 1},3,2\right),
\left({\underline 1},{\underline 2},3,{\underline 1},2\right),
\left({\underline 1},{\underline 2},3,2,{\underline 1}\right),
\left({\underline 1},3,{\underline 2},{\underline 1},2\right),
\left({\underline 1},3,{\underline 2},2,{\underline 1}\right), \\
\left({\underline 1},3,2,{\underline 2},{\underline 1}\right),
\left(3,{\underline 1},{\underline 2},{\underline 1},2\right),
\left(3,{\underline 1},{\underline 2},2,{\underline 1}\right),
\left(3,{\underline 1},2,{\underline 2},{\underline 1}\right),
\left(3,2,{\underline 1},{\underline 2},{\underline 1}\right),
\end{align*}
listed here as often as they appear in the multiset
$\shuf{ \left(1,2,1\right) }{ \left(3,2\right) }$. Here we have
underlined the letters taken from $a$ -- that is,
the letters at positions $w^{-1}\left(  1\right)$, $w^{-1}\left(  2\right)$,
$\ldots$, $w^{-1}\left(  n\right)$.
\end{definition}

% [DG][v3] I added the definition above.

% [DG][v25] Added footnote with an example of the overlaying procedure.
% And another footnote warning readers against my nonstandard
% notation (my shuffles seem to be most people's unshuffles).
% Hmm... maybe I should change it?

\begin{example}
\label{exa.shuffle-alg}
When $A=T(V)$ is the tensor algebra for a finite free $\kk$-module
$V$, having $\kk$-basis
$\{x_i\}_{i \in I}$, its graded dual $A^o$ is another Hopf
algebra whose basis $\left\{ y_{(i_1,\ldots,i_\ell)} \right\}$ (the
dual basis of the basis $\left\{ x_{i_1} \cdots x_{i_\ell} \right\}$
of $A = T\left(V\right)$) is indexed by
words in the alphabet $I$. This Hopf algebra $A^o$ could be called the
\emph{shuffle algebra} of $V^*$.
(To be more precise, it is isomorphic to the shuffle algebra of $V^*$
introduced in Proposition~\ref{prop.shuffle-alg} further below; we
prefer not to call $A^o$ itself the shuffle algebra of $V^*$, since
$A^o$ has several disadvantages\footnote{Specifically, $A^o$ has
the disadvantages of being defined only when $V^*$ is the
dual of a finite free $\kk$-module $V$, and depending on a choice of
basis,
whereas Proposition~\ref{prop.shuffle-alg} will define shuffle algebras
in full generality and canonically.}.)
Duality shows that the \emph{cut} coproduct in $A^o$ is defined by
\begin{equation}
\label{eq.exa.shuffle-alg.Delta}
\Delta y_{(i_1,\ldots,i_\ell)} =
\sum_{j=0}^{\ell} y_{(i_1,\ldots,i_j)} \otimes y_{(i_{j+1},i_{j+2},\ldots,i_\ell)}.
\end{equation}
For example,
\[
\Delta y_{abcb}  = y_{\varnothing} \otimes y_{abcb}
+y_{a} \otimes y_{bcb}
+y_{ab} \otimes y_{cb}
+y_{abc} \otimes y_{b}
+y_{abcb} \otimes y_{\varnothing} .
\]
Duality also shows that the
\emph{shuffle} product in $A^o$ will be given by
\begin{equation}
\label{eq.exa.shuffle-alg.m}
y_{(i_1,\ldots,i_\ell)} y_{(j_1,\ldots,j_m)} =
\sum_{\kk=(k_1,\ldots,k_{\ell+m})
        \in \shuf{\ii}{\jj} } y_{(k_1,\ldots,k_{\ell+m})}
\end{equation}
where $\shuf{\ii}{\jj}$ (as in Definition~\ref{shuffles})
denotes the multiset of the $\binom{\ell+m}{\ell}$
words obtained as \emph{shuffles} of the two
words $\ii=(i_1,\ldots,i_\ell)$ and $\jj=(j_1,\ldots,j_m)$.
For example,
\begin{align*}
y_{ab} y_{cb} &= y_{abcb} + y_{acbb} + y_{cabb} + y_{cabb}+ y_{acbb} + y_{cbab}\\
&= y_{abcb} + 2 y_{acbb} + 2 y_{cabb} +  y_{cbab} .
\end{align*}
Equivalently, one has
\begin{align}
\label{eq.exa.shuffle-alg.m2}
y_{(i_1,i_2,\ldots,i_\ell)} y_{(i_{\ell+1},i_{\ell+2},\ldots,i_{\ell+m})}
&=
\sum\limits_{\substack{w \in \Symm_{\ell+m}:\\w(1)<\cdots<w(\ell),\\
       w(\ell+1) < \cdots < w(\ell+m)}}
       y_{\left(i_{w^{-1}(1)},i_{w^{-1}(2)},\ldots,i_{w^{-1}(\ell+m)}\right)} \\
&=
\label{eq.exa.shuffle-alg.m3}
\sum\limits_{\substack{\sigma \in \operatorname{Sh}_{\ell,m}}}
       y_{\left(i_{\sigma\left(1\right)},i_{\sigma\left(2\right)},\ldots,i_{\sigma\left(\ell+m\right)}\right)}
\end{align}
(using the notations of Definition~\ref{shuffles} again).
Lastly, the antipode $S$ of $A^o$ is the adjoint of the antipode of $A=T(V)$ described in \eqref{antipode-in-tensor-algebra}:
\[
S y_{(i_1,i_2,\ldots,i_\ell)}
=(-1)^{\ell} y_{(i_\ell,\ldots,i_2,i_1)}.
\]
Since the coalgebra $T\left(V\right)$ is cocommutative, its graded dual
$T\left(V\right)^o$ is commutative.

\end{example}

% [DG][v25] Added parenthetical remark on more general shuffle algebras,
% and yet one more way to write down the multiplication formula.

% [DG][v31] Fixed the use of two different notations
% $\operatorname{Sh}_{\ell,m}$ and
% $\operatorname{Sh}(\ell,m)$ for one and the same thing.
% This is still a heavily nonstandard notation :/

\begin{exercise}
\label{exe.Sym.1dim}
Let $V$ be a $1$-dimensional free $\kk$-module with basis element $x$,
so $\Sym(V) \cong \kk[x]$, with
$\kk$-basis $\{1=x^0,x^1,x^2,\ldots\}$.

\begin{itemize}

\item[(a)]
Check that the powers $x^i$ satisfy
\begin{align*}
x^i \cdot x^j &= x^{i+j} , \\
\Delta(x^n) &= \sum_{i+j=n} \binom{n}{i} x^i \otimes x^j , \\
S(x^n) &= (-1)^n x^n .
\end{align*}

\item[(b)]
Check that the dual basis elements
$\{f^{(0)},f^{(1)},f^{(2)},\ldots\}$ for $\Sym(V)^o$,
defined by $f^{(i)}(x^j)=\delta_{i,j}$, satisfy
\begin{align*}
f^{(i)} f^{(j)} &= \binom{i+j}{i} f^{(i+j)} , \\
\Delta(f^{(n)}) &= \sum_{i+j=n} f^{(i)} \otimes f^{(j)} , \\
S(f^{(n)}) &= (-1)^n f^{(n)} .
\end{align*}

\item[(c)]
Show that if $\QQ$ is a subring of $\kk$, then the $\kk$-linear
map $\Sym(V)^{o} \rightarrow \Sym(V)$ sending
$f^{(n)} \mapsto \frac{x^n}{n!}$ is a graded Hopf isomorphism.

For this reason, the Hopf structure on $\Sym(V)^o$ is
called a \dfn{divided power algebra}.

\item[(d)]
Show that when $\kk$ is a field of characteristic $p >0$,
one has $(f^{(1)})^p = 0$, and hence why there can
be no Hopf isomorphism $\Sym(V)^{o} \rightarrow \Sym(V)$.
\end{itemize}
\end{exercise}

\begin{exercise}
\label{exe.Sym.wedge}
Let $V$ have $\kk$-basis $\{x_1,\ldots,x_n\}$,  and let $V \oplus V$ have
$\kk$-basis $\{x_1,\ldots,x_n, y_1,\ldots,y_n\}$, so that one has
isomorphisms
\[
\Sym(V \oplus V) \cong \kk[\xx,\yy] \cong \kk[\xx] \otimes \kk[\yy] \cong
\Sym(V) \otimes \Sym(V).
\]
Here we are using the abbreviations
$\xx = \left(x_1, x_2, \ldots, x_n\right)$
and $\yy = \left(y_1, y_2, \ldots, y_n\right)$.

\begin{itemize}

\item[(a)]
Show that our usual coproduct on $\Sym(V)$ can be re-expressed
as follows:
\[
\begin{array}{ccc}
\Sym(V) & &  \Sym(V) \otimes \Sym(V) \\
\Vert & & \Vert \\
\kk[\xx] & \overset{\Delta}{\longrightarrow} & \kk[\xx,\yy] , \\
 f(x_1,\ldots,x_n) & \longmapsto & f(x_1+y_1,\ldots,x_n+y_n) .\\
\end{array}
\]
In other words, it is induced from the diagonal map
\begin{equation}
\label{vector-space-diagonal-map}
\begin{array}{rcl}
V & \longrightarrow & V \oplus V ,\\
x_i & \longmapsto & x_i+y_i .\\
\end{array}
\end{equation}

\item[(b)]
One can similarly define a coproduct
on the \dfn{exterior algebra} $\wedge V$, which is the quotient $T(V)/J$
where $J$ is the two-sided ideal generated by the elements
$\{x^2(=x \otimes x)\}_{x \in V}$ in $T^2(V)$.
The ideal $J$ is a graded $\kk$-submodule of $T(V)$
(this is not obvious!), and
the quotient $T(V) / J$ becomes a graded commutative algebra
\[
\wedge V
 = \bigoplus_{d=0}^n \wedge^d V \left( = \bigoplus_{d=0}^\infty \wedge^d V \right),
\]
if one views the elements of $V=\wedge^1 V$ as having \emph{odd} degree, and
uses the topologist's sign convention (as in
\eqref{topologist-twist}).  One again has
$
\wedge(V \oplus V) = \wedge V \otimes \wedge V
$
as graded algebras.  Show that one can again let the diagonal map
\eqref{vector-space-diagonal-map} induce a map
\begin{equation}
\label{eq.exe.Sym.wedge.b.map}
\begin{array}{ccc}
\wedge(V) & \overset{\Delta}{\longrightarrow}
 & \wedge V \otimes \wedge V , \\
  f(x_1,\ldots,x_n) & \longmapsto & f(x_1+y_1,\ldots,x_n+y_n) \\
 \Vert & & \Vert \\
 \sum c_{i_1,\ldots,i_d} \cdot x_{i_1} \wedge \cdots \wedge x_{i_d} &
  & \sum c_{i_1,\ldots,i_d} \cdot (x_{i_1}+y_{i_1} )\wedge \cdots \wedge (x_{i_d}+y_{i_d}) ,\\
\end{array}
\end{equation}
which makes $\wedge V$ into a connected graded Hopf algebra.

\item[(c)]
Show that in the tensor algebra $T(V)$, if
one views the elements of $V=V^{\otimes 1}$ as having odd degree,
and uses the topologist's sign convention
\eqref{topologist-twist} in
the twist map when defining $T(V)$,
then for any $x$ in $V$ one has
$
\Delta(x^2) = 1 \otimes x^2 + x^2 \otimes 1.
$

\item[(d)]
Let us use the convention \eqref{topologist-twist} as in
part (c).
Show that the two-sided ideal $J \subset T(V)$
generated by $\{x^2\}_{x \in V}$ is also a two-sided coideal
and a graded $\kk$-submodule of $T(V)$,
and hence the quotient $\wedge V = T(V)/J$ inherits the structure
of a graded bialgebra.
Check that the coproduct on $\wedge V$
inherited from $T(V)$
is the same as the one defined in part (b).

% [DG][v80] Added the "graded" part to parts (a) and (d),
% because we cannot avoid keeping track of the grading
% when we are using the topologist's sign convention.
% (I mean, we could, if we defined superalgebras, but we
% haven't.)
% Anyway, this makes the exercise more interesting.

\end{itemize}

[\textbf{Hint:} The ideal $J$ in part (b) is a graded
$\kk$-submodule of $T(V)$, but this is not completely
obvious (not all elements of $V$ have to be
homogeneous!).]

% [DG][v25] Added this hint, instead of the old hint after (c)
% that I've commented out.
\end{exercise}

\begin{exercise}
\label{exe.convolution.dual}
Let $C$ be a $\kk$-coalgebra. As we know
from Exercise~\ref{exe.dual-algebra}(a), this makes $C^\ast$
into a $\kk$-algebra.

Let $A$ be a $\kk$-algebra which is finite free as
$\kk$-module. This makes $A^\ast$ into
a $\kk$-coalgebra.

Let $f : C \to A$ and $g : C \to A$ be two $\kk$-linear maps.
Show that $f^\ast \star g^\ast = \left(f \star g\right)^\ast$.
\end{exercise}

% [DG][v19] Added exercise above.

The above arguments might have created the impression that duals
of bialgebras have good properties only under certain restrictive
conditions (e.g., the dual of a bialgebra $H$ does not generally
become a bialgebra unless $H$ is of finite type), and so they
cannot be used in proofs and constructions unless one is willing
to sacrifice some generality (e.g., we had to require
$V$ to be finite free in Example~\ref{exa.shuffle-alg}). While
the first part of this impression is true, the second is not
always; often there is a way to gain back the generality lost
from using duals. As an example of this, let us define the
shuffle algebra of an arbitrary $\kk$-module (not just of a
dual of a finite free $\kk$-module as in
Example~\ref{exa.shuffle-alg}):

\begin{proposition}
\label{prop.shuffle-alg}
Let $V$ be a $\kk$-module. Define a $\kk$-linear map
$\Delta_\shuffle : T\left(V\right) \to T\left(V\right) \otimes T\left(V\right)$
by setting
\begin{align*}
\Delta_\shuffle \left( v_1 v_2 \cdots v_n \right)
= \sum_{k=0}^n \left( v_1 v_2 \cdots v_k \right) \otimes
\left( v_{k+1} v_{k+2} \cdots v_n \right)
\qquad \qquad
\text{ for all } n \in \NN \text{ and }
v_1, v_2, \ldots, v_n \in V .
\end{align*}
\footnote{This is well-defined, because the right hand side is
$n$-multilinear in $v_1, v_2, \ldots, v_n$, and because
any $n$-multilinear map
$V^{\times n} \to M$ into a $\kk$-module $M$ gives rise to a
unique $\kk$-linear map $V^{\otimes n} \to M$.} Define a
$\kk$-bilinear map $\shufmult : T\left(V\right)
\times T\left(V\right) \to T\left(V\right)$, which will be written
in infix notation (that is, we will write $a \shufmult b$
instead of $\shufmult \left(a, b\right)$), by
setting\footnote{Many authors use the symbol
$\shuffle$ instead of $\shufmult$ here, but we prefer to
reserve the former notation for the shuffle product of words.}
\begin{align*}
\left( v_1 v_2 \cdots v_\ell \right) \shufmult
\left( v_{\ell + 1} v_{\ell + 2} \cdots v_{\ell + m} \right)
&= \sum_{\sigma \in \operatorname{Sh}_{\ell, m}}
v_{\sigma\left(1\right)} v_{\sigma\left(2\right)} \cdots
v_{\sigma\left(\ell + m\right)} \\
& \qquad \qquad
\text{ for all } \ell, m \in \NN \text{ and }
v_1, v_2, \ldots, v_{\ell + m} \in V .
\end{align*}
\footnote{Again, this is well-defined by the $\ell + m$-multilinearity of
the right hand side.}
Consider also the comultiplication
$\epsilon$ of the Hopf algebra $T\left(V\right)$.

Then, the $\kk$-module $T\left(V\right)$, endowed with the
multiplication $\shufmult$, the unit
$1_{T\left(V\right)} \in V^{\otimes 0} \subset T\left(V\right)$,
the comultiplication
$\Delta_\shuffle$ and the counit $\epsilon$, becomes a commutative
Hopf algebra.
This Hopf algebra is called the \dfn{shuffle algebra} of $V$, and
denoted by $\operatorname{Sh}\left(V\right)$.
The antipode of the Hopf algebra $\operatorname{Sh}\left(V\right)$
is precisely the antipode $S$ of $T\left(V\right)$.
\end{proposition}

% [DG][v27] Changed the formatting of this definition.

\begin{exercise}
\label{exe.shuffle-alg}
Prove Proposition~\ref{prop.shuffle-alg}.

[\textbf{Hint:} When $V$ is a finite free $\kk$-module,
Proposition~\ref{prop.shuffle-alg} follows from
Example~\ref{exa.shuffle-alg}. The trick is to derive the general
case from this specific one. Every $\kk$-linear map $f : W \to V$
between two $\kk$-modules $W$ and $V$ induces a map
$T\left(f\right) : T\left(W\right) \to T\left(V\right)$ which
preserves $\Delta_\shuffle$, $\shufmult$,
$1_{T\left(W\right)}$, $\epsilon$ and $S$ (in the appropriate
meanings -- e.g., preserving $\Delta_\shuffle$ means
$\Delta_\shuffle \circ T\left(f\right)
= \left(T\left(f\right) \otimes T\left(f\right)\right)
\circ \Delta_\shuffle$). Show that each of the equalities that need to
be proven in order to verify Proposition~\ref{prop.shuffle-alg} can be
``transported'' along such a map $T\left(f\right)$ from a
$T\left(W\right)$ for a suitably chosen finite free $\kk$-module
$W$.]
\end{exercise}

% [DG][v25] Added above proposition & proof. (I might eventually
% elaborate on shuffle algebras as an aside in the Hazewinkel
% basis section; in fact, one of the simplest consequences of
% the Lyndon shuffle lemma is Radford's theorem on shuffle
% algebras, which proves that $\Qsym$ is a polynomial algebra
% over a base ring of characteristic $0$. This is a warmup to the
% harder result by Hazewinkel that it is always a polynomial
% algebra.)
% 
% I have no idea how this became 5 pages :(

% [DG][v26] This is now an exercise, and the proof is in the
% solutions. Does the hint make sense?

It is also possible to prove Proposition~\ref{prop.shuffle-alg}
``by foot'', as long as one is ready to make combinatorial
arguments about cutting shuffles.

\begin{remark} \phantomsection
\label{rmk.shuffle-alg.dual}

\begin{itemize}

\item[(a)] Let $V$ be a finite free $\kk$-module. The
Hopf algebra $T\left(V\right)^o$ (studied in
Example~\ref{exa.shuffle-alg}) is naturally isomorphic
to the shuffle algebra $\operatorname{Sh}\left(V^\ast\right)$
(defined as in Proposition~\ref{prop.shuffle-alg} but for
$V^\ast$ instead of $V$) as Hopf algebras, by the
obvious isomorphism (namely, the direct sum of the
isomorphisms $\left(V^{\otimes n}\right)^{\ast} \to
\left(V^\ast\right)^{\otimes n}$ over all $n \in \NN$).
\ \ \ \ \footnote{This can be verified
by comparing \eqref{eq.exa.shuffle-alg.Delta} with
the definition of $\Delta_\shuffle$, and comparing
\eqref{eq.exa.shuffle-alg.m3} with the definition of
$\shufmult$.}

\item[(b)] The same statement applies to the case when
$V$ is a graded $\kk$-module of finite type
% (i.e., a graded $\kk$-module such that each
% $V_n$ is finite free)
satisfying $V_0 = 0$
rather than a finite free
$\kk$-module, provided that $V^\ast$
and $\left(V^{\otimes n}\right)^{\ast}$ are replaced by
$V^o$ and $\left(V^{\otimes n}\right)^o$.

\end{itemize}

\end{remark}

% [DG][v31] Added remark above.

We shall return to shuffle algebras in
Section~\ref{subsect.shuffle.radford}, where we will show that
under certain conditions ($\QQ$ being a subring of $\kk$, and
$V$ being a free $\kk$-module) the algebra structure on a
shuffle algebra $\operatorname{Sh}(V)$
is a polynomial algebra in an appropriately chosen set
of generators\footnote{This says nothing about the coalgebra
structure on $\operatorname{Sh}(V)$ -- which is much more
complicated in these generators.}.

% [DG][v36] Added above sentence.

\subsection{\label{subsect.leray}Infinite sums and Leray's theorem}

In this section (which can be skipped, as it will not be used except in a few
exercises), we will see how a Hopf algebra structure on a $\kk$-algebra
reveals knowledge about the $\kk$-algebra itself. Specifically, we will
show that if $\kk$ is a commutative $\QQ$-algebra, and if $A$ is
any commutative connected graded $\kk$-Hopf algebra, then $A$ as a
$\kk$-algebra must be (isomorphic to) a symmetric algebra of a
$\kk$-module\footnote{If $\kk$ is a field, then this simply
means that $A$ as a $\kk$-algebra must be a polynomial ring over
$\kk$.}. This is a specimen of a class of facts which are commonly
called \emph{Leray theorems}; for different specimens, see \cite[Theorem
7.5]{MilnorMoore} or \cite[p. 17, ``Hopf's theorem'']{Cartier} or
\cite[\S 2.5, A, B, C]{Cartier} or
\cite[Theorem 3.8.3]{Cartier}.\footnote{Notice that many of these sources
assume $\kk$ to be a field; some of their proofs rely on this
assumption.} In a sense, these facts foreshadow Zelevinsky's theory of
positive self-dual Hopf algebras, which we shall encounter in
Chapter~\ref{PSH-section}; however, the latter theory works in a much less
general setting (and makes much stronger claims).

We shall first explore the possibilities of applying a formal power series $v$
to a linear map $f:C\rightarrow A$ from a coalgebra $C$ to an algebra $A$.
We have already seen an example of this in the proof of
Proposition~\ref{Takeuchi-formula} above (where the power series
$\sum_{k\geq0}\left(  -1\right)  ^{k}T^{k}
\in\kk \left[  \left[  T\right] \right]  $
was applied to the locally $\star$-nilpotent map
$\id_{A}-u_{A}\epsilon_{A}:A\rightarrow A$); we shall now take a more systematic
approach and establish general criteria for when such applications are
possible. First, we will have to make sense of infinite sums of maps from a
coalgebra to an algebra. This is somewhat technical, but the effort will pay off.

\begin{definition}
\label{def.fin-supp}
Let $A$ be an abelian group (written additively).

We say that a family $\left(  a_{q}\right)  _{q\in Q}\in A^{Q}$ of elements of
$A$ is \dfn{finitely supported} if all but finitely many $q\in Q$ satisfy
$a_{q}=0$. Clearly, if $\left(  a_{q}\right)  _{q\in Q}\in A^{Q}$ is a
finitely supported family, then the sum
$\sum_{q\in Q}a_{q}$\index{infinite sums}
is well-defined
(since all but finitely many of its addends are $0$). Sums like this satisfy
the usual rules for sums, even though their indexing set $Q$ may be infinite.
(For example, if $\left(  a_{q}\right)  _{q\in Q}$ and $\left(  b_{q}\right)
_{q\in Q}$ are two finitely supported families in $A^{Q}$, then the family
$\left(  a_{q}+b_{q}\right)  _{q\in Q}$ is also finitely supported, and we
have $\sum_{q\in Q}a_{q}+\sum_{q\in Q}b_{q}=\sum_{q\in Q}\left(  a_{q}%
+b_{q}\right)  $.)
\end{definition}

\begin{definition}
\label{def.pw-fin-supp}
Let $C$ and $A$ be two $\kk$-modules.

We say that a family $\left(  f_{q}\right)  _{q\in Q} \in \left(
\Hom \left(  C,A\right)  \right)  ^{Q}$ of maps
$f_{q} \in \Hom \left(  C,A\right)  $ is
\dfn{pointwise finitely supported} if for each $x\in C$, the family
$\left(  f_{q}\left(  x\right) \right)  _{q\in Q}\in A^{Q}$%
\index{infinite sums}
of elements of $A$ is finitely supported.\footnote{Here are some
examples of pointwise finitely supported families:
\par
\begin{itemize}
\item If $Q$ is a finite set, then any family $\left(  f_{q}\right)  _{q\in
Q}\in\left(  \Hom \left(  C,A\right)  \right)  ^{Q}$ is
pointwise finitely supported.
\par
\item More generally, any finitely supported family $\left(  f_{q}\right)
_{q\in Q}\in\left(  \Hom \left(  C,A\right)  \right)  ^{Q}$ is
pointwise finitely supported.
\par
\item If $C$ is a graded $\kk$-module, and if $\left(  f_{n}\right)
_{n\in\NN}\in\left(  \Hom \left(  C,A\right)  \right)
^{\NN}$ is a family of maps such that $f_{n}\left(  C_{m}\right)  =0$
whenever $n\neq m$, then the family $\left(  f_{n}\right)  _{n\in\NN}$
is pointwise finitely supported.
\par
\item If $C$ is a graded $\kk$-coalgebra and $A$ is any
$\kk$-algebra, and if $f\in\Hom \left(  C,A\right)  $ satisfies
$f\left(  C_{0}\right)  =0$, then the family $\left(  f^{\star n}\right)
_{n\in\NN}\in\left(  \Hom \left(  C,A\right)  \right)
^{\NN}$ is pointwise finitely supported. (This will be proven in
Proposition~\ref{prop.convolution-series}(h).)
\end{itemize}
} If $\left(  f_{q}\right)  _{q\in Q}\in\left(  \Hom \left(
C,A\right)  \right)  ^{Q}$ is a pointwise finitely supported family, then the
sum $\sum_{q\in Q}f_{q}$ is defined to be the map $C\rightarrow A$ sending
each $x\in C$ to $\sum_{q\in Q}f_{q}\left(  x\right)  $.\ \ \ \ \footnote{This
definition of $\sum_{q\in Q}f_{q}$ generalizes the usual definition of
$\sum_{q\in Q}f_{q}$ when $Q$ is a finite set (because if $Q$ is a finite set,
then any family $\left(  f_{q}\right)  _{q\in Q}\in\left(  \Hom %
\left(  C,A\right)  \right)  ^{Q}$ is pointwise finitely supported).}
\end{definition}

Note that the concept of a ``pointwise finitely supported'' family
$\left( f_q \right) _{q\in Q}
\in \left( \Hom \left( C, A \right) \right)^Q$
is precisely the concept of a ``summable'' family in
\cite[Definition 1]{DuchampMinhTolluChienNghia}.

\begin{definition}
\label{def.subsect.leray.standing1}
For the rest of Section~\ref{subsect.leray},
we shall use the following conventions:

\begin{itemize}
\item Let $C$ be a $\kk$-coalgebra. Let $A$ be a $\kk$-algebra.

\item We shall avoid our standard practice of denoting the unit map
$u_{A}:\kk\rightarrow A$ of a $\kk$-algebra $A$ by $u$; instead,
we will use the letter $u$ (without the subscript $A$) for other purposes.
\end{itemize}
\end{definition}

Definition~\ref{def.pw-fin-supp} allows us to work with infinite sums in
$\Hom \left(  C,A\right)  $, provided that we are summing a
pointwise finitely supported family. We shall next state some properties of
such sums:\footnote{See Exercise~\ref{exe.pw-fin-supp.P1-5} below for the
proofs of these properties.}

\begin{proposition}
\label{prop.pw-fin-supp.P1}
Let $\left(  f_{q}\right)  _{q\in Q}\in\left(
\Hom \left(  C,A\right)  \right)  ^{Q}$ be a pointwise finitely
supported family. Then, the map $\sum_{q\in Q}f_{q}$ belongs to
$\Hom \left(  C,A\right)  $.
\end{proposition}

\begin{proposition}
\label{prop.pw-fin-supp.P2}
Let $\left(  f_{q}\right)  _{q\in Q}$ and $\left(
g_{q}\right)  _{q\in Q}$ be two pointwise finitely supported families in
$\left(  \Hom \left(  C,A\right)  \right)  ^{Q}$. Then, the
family $\left(  f_{q}+g_{q}\right)  _{q\in Q}\in\left(  \Hom %
\left(  C,A\right)  \right)  ^{Q}$ is also pointwise finitely supported, and
satisfies
\[
\sum_{q\in Q}f_{q}+\sum_{q\in Q}g_{q}=\sum_{q\in Q}\left(  f_{q}+g_{q}\right)
.
\]

\end{proposition}

\begin{proposition}
\label{prop.pw-fin-supp.P3}
Let $\left(  f_{q}\right)  _{q\in Q}\in\left(
\Hom \left(  C,A\right)  \right)  ^{Q}$ and $\left(
g_{r}\right)  _{r\in R}\in\left(  \Hom \left(  C,A\right)
\right)  ^{R}$ be two pointwise finitely supported families. Then, the family
$\left(  f_{q}\star g_{r}\right)  _{\left(  q,r\right)  \in Q\times R}%
\in\left(  \Hom \left(  C,A\right)  \right)  ^{Q\times R}$ is
pointwise finitely supported, and satisfies%
\[
\sum_{\left(  q,r\right)  \in Q\times R}\left(  f_{q}\star g_{r}\right)
=\left(  \sum_{q\in Q}f_{q}\right)  \star\left(  \sum_{r\in R}g_{r}\right)  .
\]

\end{proposition}

Roughly speaking, the above three propositions say that sums of the form
$\sum_{q\in Q}f_{q}$ (where $\left(  f_{q}\right)  _{q\in Q}$ is a pointwise
finitely supported family) satisfy the usual rules for finite sums.
Furthermore, the following properties of pointwise finitely supported families hold:

\begin{proposition}
\label{prop.pw-fin-supp.P4}
Let $\left(  f_{q}\right)  _{q\in Q}\in\left(
\Hom \left(  C,A\right)  \right)  ^{Q}$ be a pointwise finitely
supported family. Let $\left(  \lambda_{q}\right)  _{q\in Q}\in\kk^{Q}$
be any family of elements of $\kk$. Then, the family $\left(
\lambda_{q}f_{q}\right)  _{q\in Q}\in\left(  \Hom \left(
C,A\right)  \right)  ^{Q}$ is pointwise finitely supported.
\end{proposition}

\begin{proposition}
\label{prop.pw-fin-supp.P5}
Let $\left(  f_{q}\right)  _{q\in Q}\in\left(
\Hom \left(  C,A\right)  \right)  ^{Q}$ and $\left(
g_{q}\right)  _{q\in Q}\in\left(  \Hom \left(  C,A\right)
\right)  ^{Q}$ be two families such that $\left(  f_{q}\right)  _{q\in Q}$ is
pointwise finitely supported. Then, the family $\left(  f_{q}\star
g_{q}\right)  _{q\in Q}\in\left(  \Hom \left(  C,A\right)
\right)  ^{Q}$ is also pointwise finitely supported.
\end{proposition}

\begin{exercise}
\label{exe.pw-fin-supp.P1-5}
Prove Propositions~\ref{prop.pw-fin-supp.P1}%
,~\ref{prop.pw-fin-supp.P2},~\ref{prop.pw-fin-supp.P3}%
,~\ref{prop.pw-fin-supp.P4}~and~\ref{prop.pw-fin-supp.P5}.
\end{exercise}

We can now define the notion of a ``pointwise $\star$-nilpotent'' map.
Roughly speaking, this will mean an element of
$\left(  \Hom \left(  C,A\right)  ,\star\right)  $
that can be substituted into any power series because its powers (with
respect to the convolution $\star$) form a pointwise finitely supported
family. Here is the definition:

\begin{definition}
\phantomsection\label{def.leray.pwnilp}

\begin{enumerate}
\item[(a)] A map $f\in\Hom \left(  C,A\right)  $ is said to be
\dfn{pointwise $\star$-nilpotent} if and only if the family
$\left(  f^{\star n}\right)  _{n\in\NN}\in\left(  \Hom
\left(  C,A\right)  \right)  ^{\NN}$ is pointwise finitely supported.
Equivalently, a map $f\in\Hom \left(  C,A\right)  $ is pointwise
$\star$-nilpotent if and only if for each $x\in C$, the family $\left(
f^{\star n}\left(  x\right)  \right)  _{n\in\NN}$ of elements of $A$ is
finitely supported.

\item[(b)] If $f\in\Hom \left(  C,A\right)  $ is a pointwise
$\star$-nilpotent map, and if $\left(  \lambda_{n}\right)  _{n\in\NN}
\in\kk^{\NN}$ is any family of scalars, then the family
$\left(  \lambda_{n}f^{\star n}\right)  _{n\in\NN}\in\left(
\Hom \left(  C,A\right)  \right)  ^{\NN}$ is pointwise
finitely supported\footnote{This follows easily from
Proposition~\ref{prop.pw-fin-supp.P4} above. (In fact, the map $f$ is
pointwise $\star$-nilpotent, and thus the family $\left(  f^{\star n}\right)
_{n\in\NN}\in\left(  \Hom \left(  C,A\right)  \right)
^{\NN}$ is pointwise finitely supported (by the definition of
``pointwise $\star$-nilpotent''). Hence,
Proposition~\ref{prop.pw-fin-supp.P4} (applied to $Q = \NN$ and $\left(
f_{q}\right)  _{q\in Q}=\left(  f^{\star n}\right)  _{n\in\NN}$ and
$\left(  \lambda_{q}\right)  _{q\in Q}=\left(  \lambda_{n}\right)
_{n\in\NN}$) shows that the family $\left(  \lambda_{n}f^{\star
n}\right)  _{n\in\NN}\in\left(  \Hom \left(  C,A\right)
\right)  ^{\NN}$ is pointwise finitely supported.)}, and thus the
infinite sum $\sum_{n\geq0}\lambda_{n}f^{\star n}=\sum_{n\in\NN}%
\lambda_{n}f^{\star n}$ is well-defined and belongs to $\Hom
\left(  C,A\right)  $ (by Proposition~\ref{prop.pw-fin-supp.P1}%
).\footnote{Notice that the concept of ``local $\star$-nilpotence''
we used in the proof of
Proposition~\ref{Takeuchi-prop} serves the same function (viz., ensuring that
the sum $\sum_{n\in\NN}\lambda_{n}f^{\star n}$ is well-defined). But
local $\star$-nilpotence is only defined when a grading is present, whereas
pointwise $\star$-nilpotence is defined in the general case. Also, local
$\star$-nilpotence is more restrictive (i.e., a locally $\star$-nilpotent map
is always pointwise $\star$-nilpotent, but the converse does not always
hold).}

\item[(c)] We let \dfn{$\mathfrak{n}\left(  C,A\right)  $}
be the set of all pointwise $\star$-nilpotent maps
$f\in\Hom \left(  C,A\right) $.
Note that this is not necessarily a $\kk$-submodule of
$\Hom \left(  C,A\right)  $.

\item[(d)] Consider the ring $\kk\left[  \left[  T\right]  \right]  $
of formal power series in an indeterminate $T$ over $\kk$. For any
power series $u\in\kk\left[  \left[  T\right]  \right]  $ and any
$f\in\mathfrak{n}\left(  C,A\right)  $, we define a map
$u^{\star}\left(  f\right)  \in\Hom \left(  C,A\right)  $%
\index{$u^{\star}\left(  f\right)$}
by $u^{\star}\left(  f\right)
=\sum_{n\geq0}u_{n}f^{\star n}$, where $u$ is written in the form
$u=\sum_{n\geq0}u_{n}T^{n}$ with $\left(  u_{n}\right)  _{n\geq0}\in
\kk^{\NN}$. (This sum $\sum_{n\geq0}u_{n}f^{\star n}$ is
well-defined in $\Hom \left(  C,A\right)  $, since $f$ is
pointwise $\star$-nilpotent.)
\end{enumerate}
\end{definition}

The following proposition gathers some properties of pointwise $\star
$-nilpotent maps\footnote{See Exercise~\ref{exe.convolution-series} below for
the proofs of these properties.}:

\begin{proposition}
\phantomsection\label{prop.convolution-series}

\begin{enumerate}
\item[(a)] For any $f\in\mathfrak{n}\left(  C,A\right)  $ and $k\in\NN$,
we have%
\begin{equation}
\left(  T^{k}\right)  ^{\star}\left(  f\right)  =f^{\star k}.
\label{eq.exe.convolution-series.a.monom}%
\end{equation}

\item[(b)] For any $f\in\mathfrak{n}\left(  C,A\right)  $ and $u,v\in
\kk\left[  \left[  T\right]  \right]  $, we have%
\begin{align}
\left(  u+v\right)  ^{\star}\left(  f\right)   &  =u^{\star}\left(  f\right)
+v^{\star}\left(  f\right)  \ \ \ \ \ \ \ \ \ \ \text{and}%
\label{eq.exe.convolution-series.b.1}\\
\left(  uv\right)  ^{\star}\left(  f\right)   &  =u^{\star}\left(  f\right)
\star v^{\star}\left(  f\right)  . \label{eq.exe.convolution-series.b.2}%
\end{align}
Also, for any $f\in\mathfrak{n}\left(  C,A\right)  $ and $u\in\kk%
\left[  \left[  T\right]  \right]  $ and $\lambda\in\kk$, we have%
\begin{equation}
\left(  \lambda u\right)  ^{\star}\left(  f\right)  =\lambda u^{\star}\left(
f\right)  . \label{eq.exe.convolution-series.b.3}%
\end{equation}
Also, for any $f\in\mathfrak{n}\left(  C,A\right)  $, we have%
\begin{align}
0^{\star}\left(  f\right)   &  =0\ \ \ \ \ \ \ \ \ \ \text{and}%
\label{eq.exe.convolution-series.b.4}\\
1^{\star}\left(  f\right)   &  =u_{A}\epsilon_{C}.
\label{eq.exe.convolution-series.b.5}%
\end{align}

\item[(c)] If $f,g\in\mathfrak{n}\left(  C,A\right)  $ satisfy $f\star
g=g\star f$, then $f+g\in\mathfrak{n}\left(  C,A\right)  $.

\item[(d)] For any $\lambda\in\kk$ and $f\in\mathfrak{n}\left(
C,A\right)  $, we have $\lambda f\in\mathfrak{n}\left(  C,A\right)  $.

\item[(e)] If $f\in\mathfrak{n}\left(  C,A\right)  $ and $g\in
\Hom \left(  C,A\right)  $ satisfy $f\star g=g\star f$, then
$f\star g\in\mathfrak{n}\left(  C,A\right)  $.

\item[(f)] If $v\in\kk\left[  \left[  T\right]  \right]  $ is a power
series whose constant term is $0$, then $v^{\star}\left(  f\right)
\in\mathfrak{n}\left(  C,A\right)  $ for each $f\in\mathfrak{n}\left(
C,A\right)  $.

\item[(g)] If $u,v\in\kk\left[  \left[  T\right]  \right]  $ are two
power series such that the constant term of $v$ is $0$, and if $f\in
\mathfrak{n}\left(  C,A\right)  $ is arbitrary, then%
\begin{equation}
\left(  u\left[  v\right]  \right)  ^{\star}\left(  f\right)  =u^{\star
}\left(  v^{\star}\left(  f\right)  \right)  .
\label{eq.exe.convolution-series.g.1}
\end{equation}
Here, $u\left[  v\right]  $ denotes the
\emph{composition}\index{composition of power series}
of $u$ with $v$;
this is the power series obtained by substituting $v$ for $T$ in $u$. (This
power series is well-defined, since $v$ has constant term $0$.) Furthermore,
notice that the right hand side of (\ref{eq.exe.convolution-series.g.1}) is
well-defined, since Proposition~\ref{prop.convolution-series}(f) shows that
$v^{\star}\left(  f\right)  \in\mathfrak{n}\left(  C,A\right)  $.

\item[(h)] If $C$ is a graded $\kk$-coalgebra, and if $f\in
\Hom \left(  C,A\right)  $ satisfies $f\left(  C_{0}\right)
=0$, then $f\in\mathfrak{n}\left(  C,A\right)  $.

\item[(i)] If $B$ is any $\kk$-algebra, and if $s:A\rightarrow B$ is
any $\kk$-algebra homomorphism, then every $u\in\kk\left[
\left[  T\right]  \right]  $ and $f\in\mathfrak{n}\left(  C,A\right)  $
satisfy%
\[
s\circ f\in\mathfrak{n}\left(  C,B\right)  \ \ \ \ \ \ \ \ \ \ \text{and}%
\ \ \ \ \ \ \ \ \ \ u^{\star}\left(  s\circ f\right)  =s\circ\left(  u^{\star
}\left(  f\right)  \right)  .
\]

\item[(j)] If $C$ is a connected graded $\kk$-bialgebra, and if
$F:C\rightarrow A$ is a $\kk$-algebra homomorphism, then $F-u_{A}%
\epsilon_{C}\in\mathfrak{n}\left(  C,A\right)  $.
\end{enumerate}
\end{proposition}

\begin{example}
Let $C$ be a graded $\kk$-coalgebra. Let $f\in\Hom \left(  C,A\right)  $
be such that $f\left(  C_{0}\right)  =0$. Then, we claim
that the map $u_{A}\epsilon_{C}+f:C\rightarrow A$ is $\star$-invertible. (This
observation has already been made in the proof of Proposition
\ref{Takeuchi-prop}, at least in the particular case when $C=A$.)

Let us see how this claim follows from
Proposition~\ref{prop.convolution-series}. First,
Proposition~\ref{prop.convolution-series}(h) shows that $f\in\mathfrak{n}%
\left(  C,A\right)  $. Now, define a power series $u\in\kk\left[
\left[  T\right]  \right]  $ by $u=1+T$. Then, the power series $u$ has
constant term $1$, and thus has a multiplicative inverse $v=u^{-1}%
\in\kk\left[  \left[  T\right]  \right]  $. Consider this $v$.
(Explicitly, $v=\sum_{n\geq0}\left(  -1\right)  ^{n}T^{n}$, but this does not
matter for us.) Now, (\ref{eq.exe.convolution-series.b.2}) yields $\left(
uv\right)  ^{\star}\left(  f\right)  =u^{\star}\left(  f\right)  \star
v^{\star}\left(  f\right)  $. Since $uv=1$ (because $v=u^{-1}$), we have
$\left(  uv\right)  ^{\star}\left(  f\right)  =1^{\star}\left(  f\right)
=u_{A}\epsilon_{C}$ (by (\ref{eq.exe.convolution-series.b.5})). Thus,
$u^{\star}\left(  f\right)  \star v^{\star}\left(  f\right)  =\left(
uv\right)  ^{\star}\left(  f\right)  =u_{A}\epsilon_{C}$. Hence, the map
$u^{\star}\left(  f\right)  $ has a right $\star$-inverse.

Also, from $u=1+T$, we obtain%
\begin{align*}
u^{\star}\left(  f\right)   &  =\left(  1+T\right)  ^{\star}\left(  f\right)
=\underbrace{1^{\star}\left(  f\right)  }_{=u_{A}\epsilon_{C}}%
+\underbrace{T^{\star}\left(  f\right)  }_{\substack{=f^{\star1}\\\text{(by
(\ref{eq.exe.convolution-series.a.monom}), applied to }k=1\text{)}%
}}\ \ \ \ \ \ \ \ \ \ \left(  \text{by (\ref{eq.exe.convolution-series.b.1}%
)}\right) \\
&  =u_{A}\epsilon_{C}+\underbrace{f^{\star1}}_{=f}=u_{A}\epsilon_{C}+f.
\end{align*}
Thus, the map $u_{A}\epsilon_{C}+f$ has a right $\star$-inverse (since the map
$u^{\star}\left(  f\right)  $ has a right $\star$-inverse). A similar argument
shows that this map $u_{A}\epsilon_{C}+f$ has a left $\star$-inverse.
Consequently, the map $u_{A}\epsilon_{C}+f$ is $\star$-invertible.
\end{example}

\begin{exercise}
\label{exe.convolution-series}
Prove Proposition~\ref{prop.convolution-series}.
\end{exercise}

\begin{definition}
\phantomsection\label{def.subsect.leray.standing2}

\begin{enumerate}
\item[(a)] For the rest of Section~\ref{subsect.leray}, we assume that
$\kk$ is a commutative $\QQ$-algebra. Thus, the two formal power
series $\exp=\sum_{n\geq0}\dfrac{1}{n!}T^{n}\in\kk\left[  \left[
T\right]  \right]  $\index{$\exp$} and
$\log\left(  1+T\right)  =\sum_{n\geq1}\dfrac{\left(
-1\right)  ^{n-1}}{n}T^{n}\in\kk\left[  \left[  T\right]  \right]  $%
\index{$\log$}
are well-defined.

\item[(b)] Define two power series $\overline{\exp}\in\kk\left[
\left[  T\right]  \right]  $\index{$\overline{\exp}$}
and $\overline{\log}\in\kk\left[  \left[
T\right]  \right]  $\index{$\overline{\log}$}
by $\overline{\exp} = \exp-1$ and
$\overline{\log} = \log\left(  1+T\right)  $.

\item[(c)] If $u$ and $v$ are two power series in $\kk\left[  \left[
T\right]  \right]  $ such that $v$ has constant term $0$, then
\dfn{$u\left[ v \right]$} denotes the
\emph{composition}\index{composition of power series}
of $u$ with $v$; this is the
power series obtained by substituting $v$ for $T$ in $u$.
\end{enumerate}
\end{definition}

The following proposition is just a formal analogue of the well-known fact
that the exponential function and the logarithm are mutually inverse (on their
domains of definition):\footnote{See Exercise~\ref{exe.exp-log} below for the
proof of this proposition, as well as of the lemma and proposition that follow
afterwards.}

\begin{proposition}
\label{prop.exp-log-inv}
Both power series $\overline{\exp}$ and $\overline{\log}$
have constant term $0$ and satisfy $\overline{\exp}\left[
\overline{\log}\right]  =T$ and $\overline{\log}\left[  \overline{\exp
}\right]  =T$.
\end{proposition}

For any map $f\in\mathfrak{n}\left(  C,A\right)  $, the power series $\exp$,
$\overline{\exp}$ and $\overline{\log}$ give rise to three further maps
$\exp^{\star}f$, $\overline{\exp}^{\star}f$ and $\overline{\log}^{\star}f$. We
can also define a map $\log^{\star}g$ whenever $g$ is a map in
$\Hom \left(  C,A\right)  $ satisfying $g-u_{A}\epsilon_{C}%
\in\mathfrak{n}\left(  C,A\right)  $ (but we cannot define $\log^{\star}f$ for
$f\in\mathfrak{n}\left(  C,A\right)  $, since $\log$ is not per se a power
series); in order to do this, we need a simple lemma:

\begin{lemma}
\label{lem.exp-log-log}Let $g\in\Hom \left(  C,A\right)  $ be
such that $g-u_{A}\epsilon_{C}\in\mathfrak{n}\left(  C,A\right)  $. Then,
$\overline{\log}^{\star}\left(  g-u_{A}\epsilon_{C}\right)  $ is a
well-defined element of $\mathfrak{n}\left(  C,A\right)  $.
\end{lemma}

\begin{definition}
\label{def.exp-log.log}
If $g\in\Hom \left(  C,A\right)  $ is a map satisfying
$g-u_{A}\epsilon_{C}\in\mathfrak{n}\left(  C,A\right)  $, then we define a map
$\log^{\star}g\in\mathfrak{n}\left(  C,A\right)  $\index{$\log^{\star}g$}
by $\log^{\star} g
=\overline{\log}^{\star}\left(  g-u_{A}\epsilon_{C}\right)  $. (This is
well-defined, according to Lemma~\ref{lem.exp-log-log}.)
\end{definition}

\begin{proposition}
\phantomsection\label{prop.exp-log}

\begin{enumerate}
\item[(a)] Each $f\in\mathfrak{n}\left(  C,A\right)  $ satisfies $\exp^{\star
}f-u_{A}\epsilon_{C}\in\mathfrak{n}\left(  C,A\right)  $ and%
\[
\log^{\star}\left(  \exp^{\star}f\right)  =f.
\]

\item[(b)] Each $g\in\Hom \left(  C,A\right)  $ satisfying
$g-u_{A}\epsilon_{C}\in\mathfrak{n}\left(  C,A\right)  $ satisfies%
\[
\exp^{\star}\left(  \log^{\star}g\right)  =g.
\]

\item[(c)] If $f,g\in\mathfrak{n}\left(  C,A\right)  $ satisfy $f\star
g=g\star f$, then $f+g\in\mathfrak{n}\left(  C,A\right)  $ and $\exp^{\star
}\left(  f+g\right)  =\left(  \exp^{\star}f\right)  \star\left(  \exp^{\star
}g\right)  $.

\item[(d)] The $\kk$-linear map $0:C\rightarrow A$ satisfies
$0\in\mathfrak{n}\left(  C,A\right)  $ and $\exp^{\star}0=u_{A}\epsilon_{C}$.

\item[(e)] If $f\in\mathfrak{n}\left(  C,A\right)  $ and $n\in\NN$,
then $nf\in\mathfrak{n}\left(  C,A\right)  $ and $\exp^{\star}\left(
nf\right)  =\left(  \exp^{\star}f\right)  ^{\star n}$.

\item[(f)] If $f\in\mathfrak{n}\left(  C,A\right)  $, then%
\begin{equation}
\log^{\star}\left(  f+u_{A}\epsilon_{C}\right)  =\sum_{n\geq1}\dfrac{\left(
-1\right)  ^{n-1}}{n}f^{\star n}. \label{eq.exe.convolution-series.logstar}%
\end{equation}

\end{enumerate}
\end{proposition}

\begin{example}
\label{exa.convolution-series.polynom-d0}Consider again the Hopf algebra
$\kk\left[  x\right]  $ from Exercise \ref{exe.Sym.1dim}. Let
$c_{1}:\kk\left[  x\right]  \rightarrow\kk$ be the $\kk%
$-linear map sending each polynomial $p\in\kk\left[  x\right]  $ to the
coefficient of $x^{1}$ in $p$. (In other words, $c_{1}$ sends each polynomial
$p\in\kk\left[  x\right]  $ to its derivative at $0$.)

Then, $c_{1}\left(  \left(  \kk\left[  x\right]  \right)  _{0}\right)
=0$ (as can easily be seen). Hence, Proposition~\ref{prop.convolution-series}%
(h) shows that $c_{1}\in\mathfrak{n}\left(  \kk\left[  x\right]
,\kk\right)  $. Thus, a map $\exp^{\star}\left(  c_{1}\right)
:\kk\left[  x\right]  \rightarrow\kk$ is well-defined. It is not
hard to see that this map is explicitly given by%
\[
\left(  \exp^{\star}\left(  c_{1}\right)  \right)  \left(  p\right)  =p\left(
1\right)  \ \ \ \ \ \ \ \ \ \ \text{for every }p\in\kk\left[  x\right]
.
\]
(In fact, this follows easily after showing that each $n\in\NN$
satisfies
\[
\left(  c_{1}\right)  ^{\star n}\left(  p\right)  =n!\cdot\left(  \text{the
coefficient of }x^{n}\text{ in }p\right)  \ \ \ \ \ \ \ \ \ \ \text{for every
}p\in\kk\left[  x\right]  ,
\]
which in turn is easily seen by induction.)

Note that the equality $\left(  \exp^{\star}\left(  c_{1}\right)  \right)
\left(  p\right)  =p\left(  1\right)  $ shows that the map $\exp^{\star
}\left(  c_{1}\right)  $ is a $\kk$-algebra homomorphism. This is a
particular case of a fact that we will soon see (Proposition
\ref{prop.leray.exp-algh1}).
\end{example}

\begin{exercise}
\label{exe.exp-log}Prove Proposition~\ref{prop.exp-log-inv},
Lemma~\ref{lem.exp-log-log} and Proposition~\ref{prop.exp-log}.
\end{exercise}

Next, we state another sequence of facts (some of which have nothing to do
with Hopf algebras), beginning with a fact about convolutions which is similar
to Proposition~\ref{prop.convolution.functor}:\footnote{See
Exercise~\ref{exe.leray.exp-alg} below for their proofs.}

\begin{proposition}
\label{prop.convolution.nat} Let $C$ and $C^{\prime}$ be two $\kk%
$-coalgebras, and let $A$ and $A^{\prime}$ be two $\kk$-algebras. Let
$\gamma:C\rightarrow C^{\prime}$ be a $\kk$-coalgebra morphism. Let
$\alpha:A\rightarrow A^{\prime}$ be a $\kk$-algebra morphism.

\begin{enumerate}
\item[(a)] If $f\in\Hom \left(  C,A\right)  $, $g\in
\Hom \left(  C,A\right)  $, $f^{\prime}\in\Hom
\left(  C^{\prime},A^{\prime}\right)  $ and $g^{\prime}\in\Hom
\left(  C^{\prime},A^{\prime}\right)  $ satisfy $f^{\prime}\circ\gamma
=\alpha\circ f$ and $g^{\prime}\circ\gamma=\alpha\circ g$, then $\left(
f^{\prime}\star g^{\prime}\right)  \circ\gamma=\alpha\circ\left(  f\star
g\right)  $.

\item[(b)] If $f\in\Hom \left(  C,A\right)  $ and $f^{\prime}%
\in\Hom \left(  C^{\prime},A^{\prime}\right)  $ satisfy
$f^{\prime}\circ\gamma=\alpha\circ f$, then each $n\in \NN$ satisfies
$\left(  f^{\prime}\right)  ^{\star n}\circ\gamma=\alpha\circ f^{\star n}$.
\end{enumerate}
\end{proposition}

\begin{proposition}
\label{prop.leray.f*nxy}Let $C$ be a $\kk$-bialgebra. Let $A$ be a
commutative $\kk$-algebra. Let $f\in\Hom \left(
C,A\right)  $ be such that $f\left(  \left(  \ker\epsilon\right)  ^{2}\right)
=0$ and $f\left(  1\right)  =0$. Then, any $x,y\in C$ and $n\in \NN$
satisfy%
\[
f^{\star n}\left(  xy\right)  =\sum_{i=0}^{n}\dbinom{n}{i}f^{\star i}\left(
x\right)  f^{\star\left(  n-i\right)  }\left(  y\right)  .
\]

\end{proposition}

\begin{proposition}
\label{prop.leray.exp-algh1}Let $C$ be a $\kk$-bialgebra. Let $A$ be a
commutative $\kk$-algebra. Let $f\in\mathfrak{n}\left(  C,A\right)  $
be such that $f\left(  \left(  \ker\epsilon\right)  ^{2}\right)  =0$ and
$f\left(  1\right)  =0$. Then, $\exp^{\star}f:C\rightarrow A$ is a
$\kk$-algebra homomorphism.
\end{proposition}

\begin{lemma}
\label{lem.sol.eulerian-idp.2.poly1} Let $V$ be any torsionfree abelian group
(written additively). Let $N\in{\NN}$. For every $k\in\left\{
0,1,\ldots,N\right\}  $, let $w_{k}$ be an element of $V$. Assume that
\begin{equation}
\sum_{k=0}^{N}w_{k}n^{k}=0\ \ \ \ \ \ \ \ \ \ \text{for all }n\in{\NN}.
\label{eq.lem.sol.eulerian-idp.2.poly1.ass}%
\end{equation}
Then, $w_{k}=0$ for every $k\in\left\{  0,1,\ldots,N\right\}  $.
\end{lemma}

\begin{lemma}
\label{lem.leray.poly1-inf}Let $V$ be a torsionfree abelian group (written
additively). Let $\left(  w_{k}\right)  _{k\in\NN}\in V^{\NN}$
be a finitely supported family of elements of $V$. Assume that%
\[
\sum_{k\in\NN}w_{k}n^{k}=0\ \ \ \ \ \ \ \ \ \ \text{for all }%
n\in \NN.
\]
Then, $w_{k}=0$ for every $k\in \NN$.
\end{lemma}

\begin{proposition}
\label{prop.leray.exp-algh2}Let $C$ be a graded $\kk$-bialgebra. Let
$A$ be a commutative $\kk$-algebra. Let $f\in\Hom \left(
C,A\right)  $ be such that $f\left(  C_{0}\right)  =0$. Assume
that\footnote{Notice that $\exp^{\star}f$ is well-defined, since
Proposition~\ref{prop.convolution-series}(h) yields $f\in\mathfrak{n}\left(
C,A\right)  $.} $\exp^{\star}f:C\rightarrow A$ is a $\kk$-algebra
homomorphism. Then, $f\left(  \left(  \ker\epsilon\right)  ^{2}\right)  =0$.
\end{proposition}

\begin{proposition}
\label{prop.leray.exp-algh3}Let $C$ be a connected graded $\kk%
$-bialgebra. Let $A$ be a commutative $\kk$-algebra. Let $f\in
\mathfrak{n}\left(  C,A\right)  $ be such that $f\left(  \left(  \ker
\epsilon\right)  ^{2}\right)  =0$ and $f\left(  1\right)  =0$. Assume further
that $f\left(  C\right)  $ generates the $\kk$-algebra $A$. Then,
$\exp^{\star}f:C\rightarrow A$ is a surjective $\kk$-algebra homomorphism.
\end{proposition}

\begin{exercise}
\label{exe.leray.exp-alg}Prove Lemmas~\ref{lem.sol.eulerian-idp.2.poly1}%
~and~\ref{lem.leray.poly1-inf} and Propositions~\ref{prop.convolution.nat}%
,~\ref{prop.leray.f*nxy},~\ref{prop.leray.exp-algh1}%
,~\ref{prop.leray.exp-algh2}~and~\ref{prop.leray.exp-algh3}.

[\textbf{Hint:} For Proposition~\ref{prop.leray.exp-algh2}, show first that
$\exp^{\star}\left(  nf\right)  =\left(  \exp^{\star}f\right)  ^{\star n}$ is
a $\kk$-algebra homomorphism for each $n\in \NN$. Turn this into
an equality between polynomials in $n$, and use
Lemma~\ref{lem.leray.poly1-inf}.]
\end{exercise}

With these preparations, we can state our version of \idx{Leray's theorem}:

\begin{theorem}
\label{thm.leray.leray-e}Let $A$ be a commutative connected graded
$\kk$-bialgebra.\footnote{Keep in mind that $\kk$ is assumed to
be a commutative $\QQ$-algebra.}

\begin{enumerate}
\item[(a)] We have $\id_{A}-u_{A}\epsilon_{A}%
\in\mathfrak{n}\left(  A,A\right)  $; thus, the map $\log^{\star}\left(
\id_{A}\right)  \in\mathfrak{n}\left(  A,A\right)  $
is well-defined. We denote this map
$\log^{\star}\left(  \id_{A}\right)  $ by $\mathfrak{e}$.

\item[(b)] We have $\ker\mathfrak{e}=\kk\cdot1_{A}+\left(  \ker
\epsilon\right)  ^{2}$ and $\mathfrak{e}\left(  A\right)  \cong\left(
\ker\epsilon\right)  /\left(  \ker\epsilon\right)  ^{2}$ (as $\kk$-modules).

\item[(c)] For each $\kk$-module $V$, let $\iota_{V}$ be the canonical
inclusion $V\rightarrow\Sym V$. Let $\mathfrak{q}$ be the map%
\[
A\overset{\mathfrak{e}}{\longrightarrow}\mathfrak{e}\left(  A\right)
\overset{\iota_{\mathfrak{e}\left(  A\right)  }}{\longrightarrow
}\Sym \left(  \mathfrak{e}\left(  A\right)  \right)  .
\]
Then, $\mathfrak{q}\in\mathfrak{n}\left(  A,\Sym \left(
\mathfrak{e}\left(  A\right)  \right)  \right)  $\ \ \ \ \footnote{Do not
mistake the map $\mathfrak{q}$ for $\mathfrak{e}$. While every $a\in A$
satisfies $\mathfrak{q}\left(  a\right)  =\mathfrak{e}\left(  a\right)  $, the
two maps $\mathfrak{q}$ and $\mathfrak{e}$ have different target sets, and
thus we do \textbf{not} have $\left(  \exp^{\star}\mathfrak{q}\right)  \left(
a\right) = \left(  \exp^{\star}\mathfrak{e}\right)  \left(  a\right)  $ for
every $a\in A$.}.

\item[(d)] Let $\mathbf{i}$ be the canonical inclusion $\mathfrak{e}\left(
A\right)  \rightarrow A$. Recall the universal property of the symmetric
algebra: If $V$ is a $\kk$-module, if $W$ is a commutative $\kk%
$-algebra, and if $\varphi:V\rightarrow W$ is any $\kk$-linear map,
then there exists a unique $\kk$-algebra homomorphism $\Phi
:\Sym V\rightarrow W$ satisfying $\varphi=\Phi\circ\iota_{V}$.
Applying this to $V=\mathfrak{e}\left(  A\right)  $, $W=A$ and $\varphi
=\mathbf{i}$, we conclude that there exists a unique $\kk$-algebra
homomorphism $\Phi:\Sym \left(  \mathfrak{e}\left(  A\right)
\right)  \rightarrow A$ satisfying $\mathbf{i}=\Phi\circ\iota_{\mathfrak{e}%
\left(  A\right)  }$. Denote this $\Phi$ by $\mathfrak{s}$. Then, the maps
$\exp^{\star}\mathfrak{q}:A\rightarrow\Sym \left(
\mathfrak{e}\left(  A\right)  \right)  $ and $\mathfrak{s}:\Sym %
\left(  \mathfrak{e}\left(  A\right)  \right)  \rightarrow A$ are mutually
inverse $\kk$-algebra isomorphisms.

\item[(e)] We have $A\cong\Sym \left(  \left(  \ker
\epsilon\right)  /\left(  \ker\epsilon\right)  ^{2}\right)  $ as $\kk$-algebras.

\item[(f)] The map $\mathfrak{e}:A\rightarrow A$ is a projection (i.e., it
satisfies $\mathfrak{e}\circ\mathfrak{e}=\mathfrak{e}$).
\end{enumerate}
\end{theorem}

\begin{remark}
\phantomsection\label{rmk.leray.leray-e}

\begin{enumerate}
\item[(a)] The main upshot of Theorem~\ref{thm.leray.leray-e} is that any
commutative connected graded $\kk$-bialgebra $A$ (where $\kk$ is
a commutative $\QQ$-algebra) is isomorphic \textbf{as a $\kk$-algebra}
to the symmetric algebra $\Sym W$ of some
$\kk$-module $W$. (Specifically, Theorem~\ref{thm.leray.leray-e}(e)
claims this for $W=\left(  \ker\epsilon\right)  /\left(  \ker\epsilon\right)
^{2}$, whereas Theorem~\ref{thm.leray.leray-e}(d) claims this for
$W=\mathfrak{e}\left(  A\right)  $; these two modules $W$ are isomorphic by
Theorem~\ref{thm.leray.leray-e}(b).) This is a useful statement even without
any specific knowledge about $W$, since symmetric algebras are a far tamer
class of algebras than arbitrary commutative algebras. For example, if
$\kk$ is a field, then symmetric algebras are just polynomial algebras
(up to isomorphism). This can be applied, for example, to the case of the
shuffle algebra $\operatorname*{Sh}\left(  V\right)  $ of a $\kk%
$-module $V$. The consequence is that the shuffle algebra $\operatorname*{Sh}%
\left(  V\right)  $ of any $\kk$-module $V$ (where $\kk$ is a
commutative $\QQ$-algebra) is isomorphic \textbf{as a $\kk$-algebra}
to a symmetric algebra $\Sym W$. When $V$ is
a free $\kk$-module, one can actually show that $\operatorname*{Sh}%
\left(  V\right)  $ is isomorphic \textbf{as a $\kk$-algebra}
to the symmetric algebra of a \textbf{free} $\kk$-module $W$ (that is,
to a polynomial ring over $\kk$); however, this $W$ is not easy to
characterize. Such a characterization is given by \emph{Radford's theorem}
(Theorem~\ref{thm.shuffle.radford} below) using the concept of \emph{Lyndon
words}. Notice that if $V$ has rank $\geq2$, then $W$ is not finitely generated.

\item[(b)] The isomorphism in Theorem~\ref{thm.leray.leray-e}(e) is generally
not an isomorphism of Hopf algebras. However, with a little (rather
straightforward) work, it reveals to be an isomorphism of \textbf{graded}
$\kk$-algebras. Actually, all maps mentioned in
Theorem~\ref{thm.leray.leray-e} are graded, provided that we use the
appropriate gradings for $\mathfrak{e}\left(  A\right)  $ and
$\Sym \left(  \mathfrak{e}\left(  A\right)  \right)  $. (To
define the appropriate grading for $\mathfrak{e}\left(  A\right)  $, we must
show that $\mathfrak{e}$ is a graded map, whence $\mathfrak{e}\left(
A\right)  $ is a homogeneous submodule of $A$; this provides $\mathfrak{e}%
\left(  A\right)  $ with the grading we seek. The grading on
$\Sym \left(  \mathfrak{e}\left(  A\right)  \right)  $ then
follows from the usual definition of the grading on the symmetric algebra
$\Sym V$ of a graded $\kk$-module $V$: Namely, if $V$ is
a graded $\kk$-module, then the $n$-th graded component of
$\Sym V$ is defined to be the span of all products of the form
$v_{1}v_{2}\cdots v_{k}\in\Sym V$, where $v_{1},v_{2}%
,\ldots,v_{k}\in V$ are homogeneous elements satisfying $\deg\left(
v_{1}\right)  +\deg\left(  v_{2}\right)  +\cdots+\deg\left(  v_{k}\right)  =n$.)

\item[(c)] The map $\mathfrak{e}:A\rightarrow A$ from
Theorem~\ref{thm.leray.leray-e} is called the \dfn{Eulerian idempotent} of
$A$.

\item[(d)] Theorem~\ref{thm.leray.leray-e} is concerned with commutative
bialgebras. Most of its claims have a ``dual version'',
concerning cocommutative bialgebras. Again, the
Eulerian idempotent plays a crucial role; but the result characterizes not the
$\kk$-algebra structure on $A$, but the $\kk$-coalgebra
structure on $A$. This leads to the Cartier-Milnor-Moore theorem; see
\cite[\S 3.8]{Cartier} and \cite[\S 3.2]{DuchampMinhTolluChienNghia}. We shall
say a bit about the Eulerian idempotent for a cocommutative bialgebra in
Exercises~\ref{exe.eulerian-idp}~and~\ref{exe.eulerian-idp.2}.
\end{enumerate}
\end{remark}

\begin{example}
\label{exa.leray.leray-e.Sym}Consider the symmetric algebra
$\Sym V$ of a $\kk$-module $V$. Then,
$\Sym V$ is a commutative connected graded $\kk%
$-bialgebra, and thus Theorem~\ref{thm.leray.leray-e} can be applied to
$A=\Sym V$. What is the projection $\mathfrak{e}:A\rightarrow A$
obtained in this case?

Theorem~\ref{thm.leray.leray-e}(b) shows that its kernel is%
\begin{equation}
\operatorname*{Ker}\mathfrak{e}=\underbrace{\kk\cdot1_{A}%
}_{=\Sym^{0}V}+\underbrace{\left(  \ker\epsilon
\right)  ^{2}}_{=\sum_{n\geq2}\Sym^{n}V}%
=\Sym^{0}V+\sum_{n\geq2}\Sym^{n}V
=\sum_{n\neq1}\Sym^{n}V.
\label{eq.exa.leray.leray-e.Sym.Ker}
\end{equation}
This does not yet characterize $\mathfrak{e}$ completely, because we have yet
to determine the action of $\mathfrak{e}$ on $\Sym^{1}V$.
Fortunately, the elements of $\Sym^{1}V$ are all primitive
(recall that
$\Delta_{\Sym V}\left(  v\right)  =1\otimes v+v\otimes1$
for each $v\in V$), and it can
easily be shown that the map $\mathfrak{e}$ fixes any primitive element of
$A$\ \ \ \ \footnote{See Exercise~\ref{exe.eulerian-idp}(f) further below for
this proof. (While Exercise~\ref{exe.eulerian-idp} requires $A$ to be
cocommutative, this requirement is not used in the solution to
Exercise~\ref{exe.eulerian-idp}(f). That said, this requirement is actually
satisfied for $A=\Sym V$, so we do not even need to avoid it
here.)}. Therefore, the map $\mathfrak{e}$ fixes all elements of
$\Sym^{1}V$. Since we also know that $\mathfrak{e}$
annihilates all elements of $\sum_{n\neq1}\Sym^{n}V$
(by (\ref{eq.exa.leray.leray-e.Sym.Ker})), we thus conclude that
$\mathfrak{e}$ is the canonical projection from the direct sum
$\Sym V=\bigoplus_{n\in\NN}\Sym^{n}V$ onto its addend
$\Sym^{1}V$.
\end{example}

\begin{example}
\label{exa.leray.leray-e.Sh}For this example, let $A$ be the shuffle algebra
$\operatorname*{Sh}\left(  V\right)  $ of a $\kk$-module $V$. (See
Proposition \ref{prop.shuffle-alg} for its definition, and keep in mind that
its product is being denoted by $\shufmult$, whereas the notation
$uv$ is still being used for the product of two elements $u$ and $v$ in the
\textbf{tensor} algebra $T\left(  V\right)  $.)

Theorem~\ref{thm.leray.leray-e} can be applied to $A=\operatorname*{Sh}\left(
V\right)  $. What is the projection $\mathfrak{e}:A\rightarrow A$ obtained in
this case?

Let us compute $\mathfrak{e}\left(  v_{1}v_{2}\right)  $ for two elements
$v_{1},v_{2}\in V$. Indeed, define a map $\widetilde{\id}:A\rightarrow A$
by $\widetilde{\id}=\id_{A}-u_{A}\epsilon_{A}$. Then, $\widetilde{\id}%
\in\mathfrak{n}\left(  A,A\right)  $ and $\log^{\star}\left(
\underbrace{\widetilde{\id}+ u_A \epsilon_A}%
_{=\id_{A}}\right)  =\log^{\star}\left(
\id_{A}\right)  =\mathfrak{e}$. Hence,
(\ref{eq.exe.convolution-series.logstar}) (applied to $C=A$ and
$f=\widetilde{\id}$) shows that
\begin{equation}
\mathfrak{e}=\sum_{n\geq1}\dfrac{\left(  -1\right)  ^{n-1}}{n}%
\widetilde{\id}^{\star n}. \label{eq.exa.leray.leray-e.Sh.e=}%
\end{equation}
Thus, we need to compute $\widetilde{\id}^{\star n}\left(
v_{1}v_{2}\right)  $ for each $n\geq1$.

Notice that the map $\widetilde{\id}$ annihilates $A_{0}$, but
fixes any element of $A_{k}$ for $k>0$. Thus,
\[
\widetilde{\id}\left(  w_{1}w_{2}\cdots w_{k}\right)  =%
\begin{cases}
w_{1}w_{2}\cdots w_{k}, & \text{if }k>0;\\
0, & \text{if }k=0
\end{cases}
\ \ \ \ \ \ \ \ \ \ \text{for any }w_{1},w_{2},\ldots,w_{k}\in V.
\]

But it is easy to see that the map $\widetilde{\id}^{\star
n}:A\rightarrow A$ annihilates $A_{k}$ whenever $n>k$. In particular, for
every $n>2$, the map $\widetilde{\id}^{\star n}:A\rightarrow A$
annihilates $A_{2}$, and therefore satisfies%
\begin{equation}
\widetilde{\id}^{\star n}\left(  v_{1}v_{2}\right)
=0\ \ \ \ \ \ \ \ \ \ \left(  \text{since }v_{1}v_{2}\in A_{2}\right)  .
\label{eq.exa.leray.leray-e.Sh.2}%
\end{equation}
It remains to find $\widetilde{\id}^{\star n}\left(  v_{1}%
v_{2}\right)  $ for $n\in\left\{  1,2\right\}  $.

We have $\widetilde{\id}^{\star1}=\widetilde{\id}$ and thus
\[
\widetilde{\id}^{\star1}\left(  v_{1}v_{2}\right)
=\widetilde{\id}\left(  v_{1}v_{2}\right)  =v_{1}v_{2}%
\]
and%
\begin{align*}
\widetilde{\id}^{\star2}\left(  v_{1}v_{2}\right)   &
=\underbrace{\widetilde{\id}\left(  1\right)  }_{=0}
\shufmult \underbrace{\widetilde{\id}\left(
v_{1}v_{2}\right)  }_{=v_{1}v_{2}}+\underbrace{\widetilde{\id}
\left(  v_{1}\right)  }_{=v_{1}} \shufmult
\underbrace{\widetilde{\id}\left(  v_{2}\right)  }_{=v_{2}%
}+\underbrace{\widetilde{\id}\left(  v_{1}v_{2}\right)
}_{=v_{1}v_{2}} \shufmult \underbrace{\widetilde{\id}
\left(  1\right)  }_{=0}\\
&  \ \ \ \ \ \ \ \ \ \ \left(  \text{since }\Delta_{\operatorname*{Sh}%
V}\left(  v_{1}v_{2}\right)  =1\otimes v_{1}v_{2}+v_{1}\otimes v_{2}%
+v_{1}v_{2}\otimes1\right) \\
&  =\underbrace{0 \shufmult \left(  v_{1}v_{2}\right)  }%
_{=0}+\underbrace{v_{1} \shufmult v_{2}}_{=v_{1}v_{2}+v_{2}v_{1}%
}+\underbrace{\left(  v_{1}v_{2}\right)  \shufmult 0}_{=0}\\
&  =v_{1}v_{2}+v_{2}v_{1}.
\end{align*}

Now, applying both sides of (\ref{eq.exa.leray.leray-e.Sh.e=}) to $v_{1}v_{2}%
$, we find
\begin{align*}
&  \mathfrak{e}\left(  v_{1}v_{2}\right) \\
&  =\sum_{n\geq1}\dfrac{\left(  -1\right)  ^{n-1}}{n}%
\widetilde{\id}^{\star n}\left(  v_{1}v_{2}\right)
=\underbrace{\dfrac{\left(  -1\right)  ^{1-1}}{1}}_{=1}%
\underbrace{\widetilde{\id}^{\star1}\left(  v_{1}v_{2}\right)
}_{=v_{1}v_{2}}+\underbrace{\dfrac{\left(  -1\right)  ^{2-1}}{2}}_{=\dfrac
{-1}{2}}\underbrace{\widetilde{\id}^{\star2}\left(  v_{1}%
v_{2}\right)  }_{=v_{1}v_{2}+v_{2}v_{1}}+\sum_{n\geq3}\dfrac{\left(
-1\right)  ^{n-1}}{n}\underbrace{\widetilde{\id}^{\star
n}\left(  v_{1}v_{2}\right)  }_{\substack{=0\\\text{(by
(\ref{eq.exa.leray.leray-e.Sh.2}))}}}\\
&  =v_{1}v_{2}+\dfrac{-1}{2}\left(  v_{1}v_{2}+v_{2}v_{1}\right)
+\underbrace{\sum_{n\geq3}\dfrac{\left(  -1\right)  ^{n-1}}{n}0}_{=0}%
=\dfrac{1}{2}\left(  v_{1}v_{2}-v_{2}v_{1}\right)  .
\end{align*}
This describes the action of $\mathfrak{e}$ on the graded component $A_{2}$ of
$A=\operatorname*{Sh}\left(  V\right)  $.

Similarly, we can describe $\mathfrak{e}$ acting on any other graded
component:%
\begin{align*}
\mathfrak{e}\left(  1\right)   &  =0;\\
\mathfrak{e}\left(  v_{1}\right)   &  =v_{1}\ \ \ \ \ \ \ \ \ \ \text{for each
}v_{1}\in V;\\
\mathfrak{e}\left(  v_{1}v_{2}\right)   &  =\dfrac{1}{2}\left(  v_{1}%
v_{2}-v_{2}v_{1}\right)  \ \ \ \ \ \ \ \ \ \ \text{for any }v_{1},v_{2}\in
V;\\
\mathfrak{e}\left(  v_{1}v_{2}v_{3}\right)   &  =\dfrac{1}{6}\left(
2v_{1}v_{2}v_{3}-v_{1}v_{3}v_{2}-v_{2}v_{1}v_{3}-v_{2}v_{3}v_{1}-v_{3}%
v_{1}v_{2}+2v_{3}v_{2}v_{1}\right)  \ \ \ \ \ \ \ \ \ \ \text{for any }%
v_{1},v_{2},v_{3}\in V,\\
&  \ldots.
\end{align*}
With some more work, one can show the following formula for the action of
$\mathfrak{e}$ on any nontrivial pure tensor:%
\begin{align*}
\mathfrak{e}\left(  v_{1}v_{2}\cdots v_{n}\right)   &  =\sum_{\sigma
\in\Symm_{n}}\left(  \sum_{k=1+\operatorname*{des}\left(  \sigma
^{-1}\right)  }^{n}\dfrac{\left(  -1\right)  ^{k-1}}{k}\dbinom
{n-1-\operatorname*{des}\left(  \sigma^{-1}\right)  }{k-1-\operatorname*{des}%
\left(  \sigma^{-1}\right)  }\right)  v_{\sigma\left(  1\right)  }%
v_{\sigma\left(  2\right)  }\cdots v_{\sigma\left(  n\right)  }\\
&= \sum_{\sigma \in \Symm_n}
   \dfrac{\left(-1\right)^{\operatorname{des}\left(\sigma^{-1}\right)}}{\operatorname{des}\left(\sigma^{-1}\right)+1}
   \dbinom{n}{\operatorname{des}\left(\sigma^{-1}\right)+1}^{-1}
v_{\sigma\left(  1\right)  }%
v_{\sigma\left(  2\right)  }\cdots v_{\sigma\left(  n\right)  }\\
&  \ \ \ \ \ \ \ \ \ \ \text{for any }n\geq1\text{ and }v_{1},v_{2}%
,\ldots,v_{n}\in V,
\end{align*}
where we use the notation $\operatorname*{des}\pi$ for the number of
descents\footnote{A \emph{descent}\index{descent of a permutation}
of a permutation $\pi\in\Symm_n$
means an $i\in\left\{  1,2,\ldots,n-1\right\}  $ satisfying $\pi\left(
i\right)  >\pi\left(  i+1\right)  $.} of any permutation $\pi\in
\Symm_{n}$.
(A statement essentially dual to this appears in \cite[Theorem 9.5]{Schmitt}.)

Theorem~\ref{thm.leray.leray-e}(b) yields $\ker\mathfrak{e}=\kk
\cdot1_{A}+\left(  \ker\epsilon\right)  ^{2}$. Notice, however, that $\left(
\ker\epsilon\right)  ^{2}$ means the square of the ideal $\ker\epsilon$ with
respect to the shuffle multiplication $\shufmult$; thus, $\left(
\ker\epsilon\right)  ^{2}$ is the $\kk$-linear span of all shuffle
products of the form $a \shufmult b$ with $a\in\ker\epsilon$ and
$b\in\ker\epsilon$.
\end{example}

\begin{exercise}
\label{exe.leray.leray-e}Prove Theorem~\ref{thm.leray.leray-e}.

[\textbf{Hint:} (a) is easy. For (b), define an element
$\widetilde{\id}$ of $\mathfrak{n}\left(  A,A\right)  $ by
$\widetilde{\id}=\id_{A}-u_{A}%
\epsilon_{A}$. Observe that $\mathfrak{e}=\sum_{n\geq1}\dfrac{\left(
-1\right)  ^{n-1}}{n}\widetilde{\id}^{\star n}$, and draw the
conclusions that $\mathfrak{e}\left(  1_{A}\right)  =0$ and that each $x\in A$
satisfies $\widetilde{\id}\left(  x\right)  -\mathfrak{e}%
\left(  x\right)  \in\left(  \ker\epsilon\right)  ^{2}$ (because
$\widetilde{\id}^{\star n}\left(  x\right)  \in\left(
\ker\epsilon\right)  ^{2}$ for every $n\geq2$). Use this to prove
$\ker\mathfrak{e}\subset\kk\cdot1_{A}+\left(  \ker\epsilon\right)
^{2}$. On the other hand, prove $\mathfrak{e}\left(  \left(  \ker
\epsilon\right)  ^{2}\right)  =0$ by applying
Proposition~\ref{prop.leray.exp-algh2}. Combine to obtain $\ker\mathfrak{e}%
=\kk\cdot1_{A}+\left(  \ker\epsilon\right)  ^{2}$. Finish (b) by
showing that $A/\left(  \kk\cdot1_{A}+\left(  \ker\epsilon\right)
^{2}\right)  \cong\left(  \ker\epsilon\right)  /\left(  \ker\epsilon\right)
^{2}$ as $\kk$-modules. Part (c) is easy again. For (d), first apply
Proposition~\ref{prop.convolution-series}(i) to show that $\exp^{\star}\left(
\mathfrak{s}\circ\mathfrak{q}\right)  =\mathfrak{s}\circ\left(  \exp^{\star
}\mathfrak{q}\right)  $. In light of $\mathfrak{s}\circ\mathfrak{q}%
=\mathfrak{e}$ and $\exp^{\star}\mathfrak{e}=\id_{A}$,
this becomes $\id_{A}=\mathfrak{s}\circ\left(
\exp^{\star}\mathfrak{q}\right)  $. To obtain part (d), it remains to show
that $\exp^{\star}\mathfrak{q}$ is a surjective $\kk$-algebra
homomorphism; but this follows from Proposition~\ref{prop.leray.exp-algh3}.
For (e), combine (d) and (b). For (f), use once again the observation that
each $x\in A$ satisfies $\widetilde{\id}\left(  x\right)
-\mathfrak{e}\left(  x\right)  \in\left(  \ker\epsilon\right)  ^{2}$.]
\end{exercise}

% [DG][v67] Added the above section.

\newpage

%%%%%%%%%%%%%%%%%%%%%%%%%%%%%%%%%%%%%
\section{Review of symmetric functions $\Lambda$ as Hopf algebra}
\label{Sym-section}
%%%%%%%%%%%%%%%%%%%%%%%%%%%%%%%%%%%%%

Here we review the ring of symmetric functions, borrowing heavily from
standard treatments, such as Macdonald \cite[Chap. I]{Macdonald},
Sagan \cite[Chap. 4]{Sagan}, Stanley \cite[Chap. 7]{Stanley},
and Mendes and Remmel \cite{MendesRemmel},
but emphasizing the Hopf structure early on.
Other recent references for this subject are \cite{Wildon2016},
\cite{Sam-symf}, \cite{Egge}, \cite[Chapters 2--3]{Meliot}
and \cite[Chapter 7]{Sagan2020}.

\subsection{\label{subsect.Sym-section.def}Definition of $\Lambda$}

As before, $\kk$ here is a commutative ring (hence could be a field or
the integers $\ZZ$; these are the usual choices).

Given an infinite variable set $\xx=(x_1,x_2,\ldots)$\index{$\xx$},
a monomial
$\xx^\alpha:=x_1^{\alpha_1} x_2^{\alpha_2} \cdots$\index{$\xx^\alpha$}
is indexed by a sequence
$\alpha=(\alpha_1,\alpha_2,\ldots)$ in $\NN^\infty$ having finite
support\footnote{The \dfn{support} of a sequence
$\alpha = \left(\alpha_1, \alpha_2, \alpha_3, \ldots\right)
\in \NN^\infty$ is defined to be the set of all positive integers
$i$ for which $\alpha_i \neq 0$.};
such sequences $\alpha$ are called \dfn{weak compositions}.
The nonzero entries of the sequence
$\alpha = \left(\alpha_1,\alpha_2,\ldots\right)$ are called the
\emph{parts}\index{part of a weak composition}
of the weak composition $\alpha$.

The sum $\alpha_1+\alpha_2+\alpha_3+\cdots$ of all entries of
a weak composition
$\alpha = \left(\alpha_1, \alpha_2, \alpha_3, \ldots\right)$
(or, equivalently, the sum of all parts of $\alpha$)
is called the \emph{size}\index{size of a weak composition}
of $\alpha$ and denoted by $\left|\alpha\right|$.\index{$\mid\alpha\mid$}

Consider the $\kk$-algebra $\kk\left[\left[\xx\right]\right]
:= \kk\left[\left[x_1, x_2, x_3, \ldots\right]\right]$\index{$\kk\left[\left[\xx\right]\right]$}
of all formal power series in the indeterminates $x_1, x_2, x_3, \ldots$
over $\kk$; these series are infinite $\kk$-linear combinations
$\sum_{\alpha} c_\alpha \xx^\alpha$ (with $c_\alpha$ in $\kk$)
of the monomials $\xx^\alpha$ where $\alpha$ ranges over all
weak compositions.
The product of two such formal power series is well-defined
by the usual multiplication rule.

The \emph{degree}\index{degree of a monomial}
of a monomial $\xx^{\alpha}$ is defined to be the number
$\deg(\xx^{\alpha}) := \sum_i \alpha_i \in \NN$\index{$\deg(\xx^{\alpha})$}.
Given a number $d \in \NN$,
we say that a formal power series
$f(\xx) = \sum_{\alpha} c_\alpha \xx^\alpha \in \kk\left[\left[\xx\right]\right]$
(with $c_\alpha$ in $\kk$)
is \dfn{homogeneous of degree $d$}\index{homogeneous power series}
if every weak composition
$\alpha$ satisfying $\deg(\xx^{\alpha}) \neq d$
must satisfy $c_\alpha = 0$.
In other words, a formal power series is homogeneous of degree $d$
if it is an infinite $\kk$-linear combination of monomials
of degree $d$.
Every formal power series $f \in \kk\left[\left[\xx\right]\right]$
can be uniquely represented as an infinite sum
$f_0 + f_1 + f_2 + \cdots$, where each $f_d$ is homogeneous
of degree $d$; in this case, we refer to each $f_d$ as the
\emph{$d$-th homogeneous component of $f$}\index{homogeneous component}.
Note that this does not make $\kk\left[\left[\xx\right]\right]$ into
a graded $\kk$-module, since these sums $f_0 + f_1 + f_2 + \cdots$
can have infinitely many nonzero addends.
Nevertheless, if $f$ and $g$ are homogeneous power series of
degrees $d$ and $e$, then $fg$ is homogeneous of degree $d+e$.

% [DG][v79] Added the above paragraph.

A formal power series
$f(\xx) = \sum_{\alpha} c_\alpha \xx^\alpha \in \kk\left[\left[\xx\right]\right]$
(with $c_\alpha$ in $\kk$)
is said to be \dfn{of bounded degree}\index{bounded degree}
if there exists some
bound $d=d(f) \in \NN$ such that every weak composition
$\alpha = \left(\alpha_1, \alpha_2, \alpha_3, \ldots\right)$
satisfying $\deg(\xx^{\alpha}) > d$
must satisfy $c_\alpha = 0$.
Equivalently, a formal power series $f \in \kk\left[\left[\xx\right]\right]$
is of bounded degree if all but finitely many of its
homogeneous components are zero.
(For example, $x_1^2 + x_2^2 + x_3^2 + \cdots$ and
$1 + x_1 + x_2 + x_3 + \cdots$ are of bounded
degree, while $x_1 + x_1 x_2 + x_1 x_2 x_3 + \cdots$ and
$1 + x_1 + x_1^2 + x_1^3 + \cdots$ are not.)
It is easy to see that the sum and the product of two power
series of bounded degree also have bounded degree.
Thus, the formal power series of bounded degree form
a $\kk$-subalgebra of $\kk\left[\left[\xx\right]\right]$,
which we call \dfn{$R(\xx)$}.
This subalgebra $R(\xx)$ is graded (by degree).

% [DG] I added the definition of "parts", requiring them to be nonzero
% since otherwise "partition with at most $n$ parts" doesn't make much
% sense. On the other hand, you speak of "nonzero parts" sometimes...

% [DG][v79] Parts are now nonzero everywhere. If zeroes are desired,
% I speak of entries instead.

The symmetric group $\Symm_n$ permuting the first $n$ variables
$x_1,\ldots,x_n$ acts as a group of automorphisms on $R(\xx)$, as does
the union $\Symm_{(\infty)}=\bigcup_{n \geq 0} \Symm_n$%
\index{$\Symm_{(\infty)}$}
of the infinite ascending chain
$
\Symm_0 \subset \Symm_1 \subset \Symm_2 \subset \cdots
$
of symmetric groups\footnote{This ascending chain is constructed as
follows: For every $n \in \NN$, there is an injective group
homomorphism $\iota_n : \Symm_n \to \Symm_{n+1}$ which sends every
permutation $\sigma \in \Symm_n$ to the permutation
$\iota_n \left(\sigma\right) = \tau \in \Symm_{n+1}$
defined by
\[
\tau\left(i\right) =
\begin{cases}
\sigma\left(i\right), &\text{if } i \leq n ; \\
i, &\text{if } i = n+1
\end{cases}
\qquad \text{ for all } i \in \left\{1,2,\ldots,n+1\right\} .
\]
These homomorphisms
$\iota_n$ for all $n$ form a chain
$\Symm_0 \overset{\iota_0}{\longrightarrow}
\Symm_1 \overset{\iota_1}{\longrightarrow}
\Symm_2 \overset{\iota_2}{\longrightarrow} \cdots$,
which is often regarded as a chain of inclusions.}.
This group $\Symm_{(\infty)}$ can also be
described as the group of all permutations of the set
$\left\{1,2,3,\ldots\right\}$ which leave all but finitely many
elements invariant.
It is known as the \dfn{finitary symmetric group}
on $\left\{1,2,3,\ldots\right\}$.

The group $\Symm_{(\infty)}$ also acts on the set of
all weak compositions by permuting their entries:
\begin{align*}
\sigma \left(\alpha_1, \alpha_2, \alpha_3, \ldots\right)
&= \left(\alpha_{\sigma^{-1}\left(1\right)},
\alpha_{\sigma^{-1}\left(2\right)},
\alpha_{\sigma^{-1}\left(3\right)}, \ldots\right)
\\
&\qquad
\text{for any weak composition $\left(\alpha_1, \alpha_2, \alpha_3, \ldots\right)$
and any $\sigma \in \Symm_{(\infty)}$}.
\end{align*}

These two actions are connected by the equality
$\sigma\left(\xx^\alpha\right) = \xx^{\sigma\alpha}$
for any
weak composition $\alpha$ and any $\sigma \in \Symm_{(\infty)}$.

% [DG][v14] Again, added $\Symm_0$. (You already used
% $\bigcup_{n \geq 0} \Symm_n$ in the Malvenuto-Reutenauer chapter,
% so this makes things more consistent.)

% [DG][v25] Added last sentence. (I define $\Symm_{(\infty,\infty)}$
% later, in the same way; the "union of a chain of finite symmetric
% groups" definition would be more complicated in this case.)

% [DG][v54] Added footnotes defining "support" and the chain of
% inclusions of the symmetric group.

% [DG][v75] Added lots of details to the above (e.g., examples
% of bounded-degree and not-bounded-degree power series); split
% some sentences.

\begin{definition}
The \emph{ring of symmetric functions in $\xx$ with coefficients in $\kk$}%
\index{symmetric function}\index{ring of symmetric functions}\index{$\Lambda$},
denoted
$
\Lambda = \Lambda_\kk = \Lambda(\xx) = \Lambda_\kk(\xx),
$
is the $\Symm_{(\infty)}$-invariant subalgebra $R(\xx)^{\Symm_{(\infty)}}$
of $R(\xx)$:
\begin{align*}
\Lambda & :=
\left\{ f \in R(\xx) : \sigma \left(f\right) = f
\text{ for all } \sigma \in \Symm_{(\infty)} \right\}
\\
&=
\left\{ f=\sum_\alpha c_\alpha \xx^\alpha \in R(\xx):
c_\alpha=c_\beta\text{ if }\alpha, \beta\text{ lie in the same }
\Symm_{(\infty)}\text{-orbit} \right\}.
\end{align*}
\end{definition}

We refer to the elements of $\Lambda$ as
\emph{symmetric functions}\index{symmetric function}
(over $\kk$); however, despite this terminology, they are not
functions in the usual sense.\footnote{Being power
series, they can be evaluated at appropriate families of
variables. But this does not make them functions (no more than
polynomials are functions).
The terminology ``symmetric function'' is thus not
well-chosen; but it is standard.}

% [DG][v49] Added above paragraph.

Note that $\Lambda$ is a graded $\kk$-algebra, since
$\Lambda = \bigoplus_{n \geq 0} \Lambda_n$ where $\Lambda_n$ are
the symmetric functions $f=\sum_\alpha c_\alpha \xx^\alpha$
which are \emph{homogeneous of degree $n$},
meaning $\deg(\xx^\alpha)=n$ for all $c_\alpha \neq 0$.

\begin{exercise}
\label{exe.Lambda.subs.fin}
Let $f \in R\left(\xx\right)$.
Let $A$ be a commutative $\kk$-algebra, and
$a_1, a_2, \ldots, a_k$ be finitely many elements of $A$. Show
that substituting $a_1, a_2, \ldots, a_k, 0, 0, 0, \ldots$ for
$x_1, x_2, x_3, \ldots$ in $f$ yields an infinite sum
in which all but finitely many addends are zero. Hence, this
sum has a value in $A$, which is commonly denoted by
$f\left(a_1, a_2, \ldots, a_k\right)$.
\end{exercise}

% [DG][v25] Added above exercise.

\begin{definition}
A \dfn{partition} $\lambda=(\lambda_1,\lambda_2,\ldots,\lambda_\ell,0,0,\ldots)$
is a weak composition whose entries weakly decrease:
$\lambda_1 \geq \cdots \geq \lambda_\ell > 0$.
The (uniquely defined) $\ell$ is said to be the
\emph{length}\index{length of a partition} of the partition $\lambda$
and denoted by \dfn{$\ell\left(\lambda\right)$}.
Thus, $\ell\left(\lambda\right)$ is the number of
parts\footnote{Recall that a \emph{part} of a partition means a
nonzero entry of the partition.} of $\lambda$.
One sometimes omits trailing zeroes from a partition:
e.g., one can write the partition
$\left(3,1,0,0,0,\ldots\right)$ as $\left(3,1\right)$.
We will often (but not always) write \dfn{$\lambda_i$}
for the $i$-th entry of the partition $\lambda$ (for instance,
if $\lambda = \left(5,3,1,1\right)$, then
$\lambda_2 = 3$ and $\lambda_5 = 0$).
If $\lambda_i$ is nonzero, we will also call it the
\dfn{$i$-th part} of $\lambda$.
The sum
$\lambda_1+\lambda_2+\cdots+\lambda_\ell = \lambda_1+\lambda_2+\cdots$ (where
$\ell = \ell\left(\lambda\right)$) of
all entries of $\lambda$ (or, equivalently, of all
parts of $\lambda$)
is the size $\left|\lambda\right|$ of $\lambda$\index{size of a partition}%
\index{$\mid\lambda\mid$}. % because | signs in index entries seem to be broken.
For a given integer $n$, the partitions of size $n$
are referred to as the \emph{partitions of $n$}\index{partition of $n$}.
The \idx{empty partition} $() = (0,0,0,\ldots)$
is denoted by \dfn{$\varnothing$}.

% [DG] I added the definitions of "length" and "partition of $n$" since you
% use them later.

% [DG][v25] Added more pedantry to the above definition :P

% [DG][v54] Replaced "whose parts weakly
% decrease" by "whose entries weakly decrease" in the definition of
% a partition, since $(3,0,2,0,0,\ldots)$ is not a partition.

Partitions (as defined above) are sometimes called
\emph{integer partitions}\index{integer partition}
in order to distinguish them from set partitions.

Every weak composition $\alpha$ lies in the  $\Symm_{(\infty)}$-orbit
of a unique partition
$\lambda=(\lambda_1,\lambda_2,\ldots,\lambda_\ell,0,0,\ldots)$
with $\lambda_1 \geq \cdots \geq \lambda_\ell > 0$.
For any partition $\lambda$, define the
\dfn{monomial symmetric function}\index{$m_\lambda$}
\begin{equation}
\label{monomial-symmetric-function-definition}
m_\lambda:=\sum_{\alpha \in \Symm_{(\infty)}\lambda}  \xx^\alpha.
\end{equation}
Letting $\lambda$ run through the set \dfn{$\Par$} of all partitions,
this gives the \emph{monomial $\kk$-basis}\index{monomial basis of $\Lambda$}
$\{m_{\lambda}\}$ of $\Lambda$.
Letting $\lambda$ run only through the set \dfn{$\Par_n$} of partitions of $n$
gives the monomial $\kk$-basis for $\Lambda_n$.
\end{definition}

\begin{example}
For $n=3$, one has
\begin{align*}
m_{(3)} &= x_1^3+x_2^3+x_3^3+ \cdots ,\\
m_{(2,1)} &= x_1^2 x_2 + x_1 x_2^2 + x_1^2 x_3 +x_1 x_3^2 + \cdots ,\\
m_{(1,1,1)} &= x_1 x_2 x_3 +x_1 x_2 x_4 + x_1 x_3 x_4 + x_2 x_3 x_4 + x_1 x_2 x_5 + \cdots .
\end{align*}
\end{example}

The monomial basis $\{m_{\lambda}\}_{\lambda \in \Par}$ of $\Lambda$
is thus a graded basis\footnote{See Definition~\ref{def.graded-basis}
for the meaning of ``graded basis''.}
of the graded $\kk$-module $\Lambda$.
(Here and in the following, when we say that a basis
$\left\{u_\lambda\right\}_{\lambda \in \Par}$ indexed by $\Par$
is a \idx{graded basis of $\Lambda$},
we tacitly understand that $\Par$ is partitioned into
$\Par_0, \Par_1, \Par_2, \ldots$, so that for each $n \in \NN$,
the subfamily
$\left\{u_\lambda\right\}_{\lambda \in \Par_n}$ should be a
basis for $\Lambda_n$.)

% [DG][v75] Added the above (stupid but necessary) paragraph.

\begin{remark}
\label{rmk.Lambda.symmetric}
We have defined the symmetric functions as the elements of
$R\left(\xx\right)$ invariant under the group $\Symm_{(\infty)}$.
However, they also are the elements of $R\left(\xx\right)$
invariant under the group $\Symm_\infty$ of \emph{all}
permutations of the set $\left\{1, 2, 3, \ldots\right\}$
(which acts on $R\left(\xx\right)$ in the same way as its
subgroup $\Symm_{(\infty)}$ does).%
\footnote{\textit{Proof.} We need to show that
$\Lambda = R\left(\xx\right)^{\Symm_\infty}$. Since
\[
\Lambda
= \left\{ f=\sum_\alpha c_\alpha \xx^\alpha \in R\left(\xx\right):
c_\alpha=c_\beta\text{ if }\alpha, \beta\text{ lie in the same }
\Symm_{(\infty)}\text{-orbit} \right\}
\]
and
\[
R\left(\xx\right)^{\Symm_\infty}
= \left\{ f=\sum_\alpha c_\alpha \xx^\alpha \in R\left(\xx\right):
c_\alpha=c_\beta\text{ if }\alpha, \beta\text{ lie in the same }
\Symm_{\infty}\text{-orbit} \right\} ,
\]
this will follow immediately if we can show that two weak
compositions $\alpha$ and $\beta$ lie in the same
$\Symm_{(\infty)}$-orbit if and only if they lie in the same
$\Symm_{\infty}$-orbit. But this is straightforward to check
(in fact, two weak compositions $\alpha$ and $\beta$ lie in the
same orbit under either group if and only if they have the same
multiset of nonzero entries).}
\end{remark}

% [DG][v39] Added above remark.

% [DG][v40] Replaced proof, since the previous argument wasn't very
% easy to follow.

\begin{remark}
\label{finite-variable-set-remark}
It is sometimes convenient to work with finite variable sets $x_1,\ldots,x_n$,
which one justifies as follows.  Note that the algebra homomorphism
\[
R(\xx) \rightarrow R(x_1,\ldots,x_n)=\kk[x_1,\ldots,x_n]
\]
which sends $x_{n+1},x_{n+2},\ldots$ to $0$ restricts to
an algebra homomorphism
\[
\Lambda_\kk(\xx) \rightarrow \Lambda_\kk(x_1,\ldots,x_n)=\kk[x_1,\ldots,x_n]^{\Symm_n}.
\]
Furthermore, this last homomorphism is a $\kk$-module isomorphism when restricted
to $\Lambda_i$ for $0 \leq i \leq n$, since it sends the monomial basis
elements $m_\lambda(\xx)$ to the monomial basis elements
$m_\lambda(x_1,\ldots,x_n)$.
Thus, when one proves identities in $\Lambda_\kk(x_1,\ldots,x_n)$ for all $n$,
they are valid in $\Lambda$, that is, $\Lambda$ is the inverse limit of
the $\Lambda(x_1,\ldots,x_n)$ in the category of graded
$\kk$-algebras.\footnote{\textit{Warning:} The word ``graded'' here
is crucial. Indeed, $\Lambda$ is \textbf{not}
the inverse limit of the $\Lambda(x_1,\ldots,x_n)$ in the category of
$\kk$-algebras. In fact, the latter limit is the $\kk$-algebra of
all symmetric power series $f$ in $\kk\left[\xx\right]$ with the
following property: For each $g \in \NN$, there exists a
$d \in \NN$ such that every monomial in $f$ that involves exactly $g$
distinct indeterminates has degree at most $d$. For example, the power
series
$\left(1+x_1\right)\left(1+x_2\right)\left(1+x_3\right)\cdots$
and
$m_{\left(1\right)} + m_{\left(2,2\right)} + m_{\left(3,3,3\right)}
+ \cdots$ satisfy this property, although they do not lie in $\Lambda$
(unless $\kk$ is a trivial ring).}

This characterization of $\Lambda$ as an inverse limit of
the graded $\kk$-algebras $\Lambda(x_1,\ldots,x_n)$ can be used as
an alternative definition of $\Lambda$. The definitions used by
Macdonald \cite{Macdonald} and Wildon \cite{Wildon2016} are closely
related (see \cite[\S 1.2, p. 19, Remark 1]{Macdonald},
\cite[\S A.11]{Hazewinkel2} and \cite[\S 1.7]{Wildon2016} for
discussions of this definition). It also suggests that much of the
theory of symmetric functions can be rewritten in terms of the
$\Lambda(x_1,\ldots,x_n)$ (at the cost of extra complexity); and this
indeed is possible\footnote{See, for example,
\cite[Chapter SYM]{LaksovLascouxPragaczThorup},
\cite{Prasad-schur} and
\cite[Chapters 10--11]{Loehr-bij} for various results
of this present chapter rewritten in terms of symmetric polynomials
in finitely many variables.}.
\end{remark}

% [DG][v58] Added stuff to the above remark. Also, added the Wildon
% reference.

One can also define a comultiplication on $\Lambda$ as follows.

Consider the countably infinite set of variables
$(\xx,\yy)=(x_1,x_2,\ldots,y_1,y_2,\ldots)$.\index{$(\xx,\yy)$}
Although it properly contains $\xx$, there are nevertheless
bijections between $\xx$ and $(\xx, \yy)$, since these two
variable sets have the same cardinality.

Let $R(\xx, \yy)$ denote the $\kk$-algebra of formal power series
in $(\xx,\yy)$ of bounded degree.
Let $\Symm_{(\infty,\infty)}$
be the group of all permutations of $\left\{ x_1, x_2, \ldots,
y_1, y_2, \ldots\right\}$ leaving
all but finitely many variables invariant.
Then, $\Symm_{(\infty,\infty)}$ acts on $R(\xx, \yy)$ by
permuting variables, in the same way as $\Symm_{(\infty)}$ acts on
$R(\xx)$.
The fixed space $R(\xx, \yy)^{\Symm_{(\infty,\infty)}}$
is a $\kk$-algebra, which we denote by $\Lambda(\xx, \yy)$.
This $\kk$-algebra $\Lambda(\xx, \yy)$
is isomorphic to $\Lambda = \Lambda(\xx)$, since there is
a bijection between the two sets of variables $(\xx, \yy)$ and
$\xx$.
More explicitly: The map
\begin{equation}
\begin{array}{rcl}
\Lambda=\Lambda(\xx) & \overset{\Delta}{\longrightarrow}
   & \Lambda(\xx,\yy) , \\
f(\xx)=f(x_1,x_2,\ldots) & \longmapsto & f(\xx,\yy)=f(x_1,x_2,\ldots,y_1,y_2,\ldots)
\end{array}
\label{symmetric-functions-in-two-sets-fxxyy}
\end{equation}
is a graded $\kk$-algebra isomorphism.
Here, $f \left( x_1, x_2, \ldots, y_1, y_2, \ldots \right)$ means
the result of choosing some bijection
\newline
$\phi : \left\{ x_1, x_2, x_3, \ldots \right\}
\to \left\{ x_1, x_2, \ldots, y_1, y_2, \ldots \right\}$ and
substituting $\phi\left(x_i\right)$ for every $x_i$ in $f$.
(The choice of $\phi$ is irrelevant since $f$ is
symmetric.\footnote{To be more precise, the choice of $\phi$ is
irrelevant because $f$ is $\Symm_\infty$-invariant, with the
notations of Remark~\ref{rmk.Lambda.symmetric}.})

The group $\Symm_{(\infty)} \times \Symm_{(\infty)}$ is a subgroup of
the group $\Symm_{(\infty,\infty)}$ (via the obvious injection,
which lets each
$\left(\sigma, \tau\right) \in \Symm_{(\infty)} \times \Symm_{(\infty)}$
act by separately
permuting the $x_1, x_2, x_3, \ldots$ using $\sigma$
and permuting the $y_1, y_2, y_3, \ldots$ using $\tau$),
% where $\Symm_{(\infty)} \times \Symm_{(\infty)}$ denotes
% permutations of (finite subsets of) the
% $\xx$  and separate permutations of (finite subsets of)
% the $\yy$, 
and thus also acts on $R(\xx, \yy)$.
Hence, we have an inclusion of $\kk$-algebras
$\Lambda(\xx,\yy)
= R(\xx, \yy)^{\Symm_{(\infty,\infty)}}
\subset R(\xx, \yy)^{\Symm_{(\infty)} \times \Symm_{(\infty)}}
\subset R(\xx, \yy)$.
The $\kk$-module $R(\xx,\yy)^{\Symm_{(\infty)} \times \Symm_{(\infty)}}$
has $\kk$-basis $\{ m_\lambda(\xx) m_\mu(\yy) \}_{\lambda, \mu \in \Par}$,
since $m_\lambda(\xx) m_\mu(\yy)$ is just the sum of
all monomials in the
$\Symm_{(\infty)} \times \Symm_{(\infty)}$-orbit of
$\xx^\lambda \yy^\mu$ (and since any
$\Symm_{(\infty)} \times \Symm_{(\infty)}$-orbit of
monomials has exactly one representative of the form
$\xx^\lambda \yy^\mu$ with $\lambda, \mu \in \Par$).
Here, of course, $\yy$ stands for the set of variables
$\left(y_1, y_2, y_3, \ldots\right)$, and we define
$\yy^\mu$ to be $y_1^{\mu_1} y_2^{\mu_2} \cdots$.

On the other hand, the map
% Note that when one decomposes the variables into two sets
% $(\xx,\yy)=(x_1,x_2,\ldots,y_1,y_2,\ldots)$, one has a ring homomorphism
\[
\begin{array}{rcl}
R(\xx) \otimes R(\xx) & \longrightarrow& R(\xx,\yy) , \\
f(\xx) \otimes g(\xx) & \longmapsto& f(\xx)g(\yy)
\end{array}
\]
is a $\kk$-algebra homomorphism.
Restricting it to $R(\xx)^{\Symm_{(\infty)}} \otimes R(\xx)^{\Symm_{(\infty)}}$,
we obtain a $\kk$-algebra homomorphism
\begin{equation}
\label{symmetric-functions-in-two-sets-iso}
\Lambda \otimes \Lambda
=R(\xx)^{\Symm_{(\infty)}} \otimes R(\xx)^{\Symm_{(\infty)}}
\longrightarrow R(\xx,\yy)^{\Symm_{(\infty)} \times \Symm_{(\infty)} } ,
\end{equation}
which is an isomorphism because it sends the basis
$\{ m_\lambda \otimes m_\mu \}_{\lambda, \mu \in \Par}$
of the $\kk$-module $\Lambda \otimes \Lambda$
to the basis $\{ m_\lambda(\xx) m_\mu(\yy) \}_{\lambda, \mu \in \Par}$
of the $\kk$-module $R(\xx,\yy)^{\Symm_{(\infty)} \times \Symm_{(\infty)}}$.
% As $\Symm_{(\infty)} \times \Symm_{(\infty)}$ is a subgroup of the
% group $\Symm_{(\infty,\infty)}$
% (the group of all permutations of $\left\{ x_1, x_2, \ldots,
% y_1, y_2, \ldots\right\}$ leaving
% all but finitely many variables invariant)
% acting on all of $(\xx,\yy)$,
Thus, we get an inclusion of graded $\kk$-algebras
\[
\Lambda(\xx,\yy) =  R(\xx,\yy)^{\Symm_{(\infty,\infty)}} \hookrightarrow
     R(\xx,\yy)^{\Symm_{(\infty)} \times \Symm_{(\infty)} } \cong \Lambda \otimes \Lambda
\]
where the last isomorphism is the inverse of the one in \eqref{symmetric-functions-in-two-sets-iso}.
This gives a comultiplication
\[
\begin{array}{rcl}
\Lambda=\Lambda(\xx) & \overset{\Delta}{\longrightarrow}
   & \Lambda(\xx,\yy) \hookrightarrow \Lambda \otimes \Lambda ,\\
f(\xx)=f(x_1,x_2,\ldots) & \longmapsto & f(\xx,\yy)=f(x_1,x_2,\ldots,y_1,y_2,\ldots).
\end{array}
\]
Here, $f \left( x_1, x_2, \ldots, y_1, y_2, \ldots \right)$
is understood as in \eqref{symmetric-functions-in-two-sets-fxxyy}.

% [DG][v25] Added the last two sentences about what
% $f \left( x_1, x_2, \ldots, y_1, y_2, \ldots \right)$ means.
% Or am I blind and this has already been dealt with?
% Also, replaced the "$\Symm_{(\infty)}$ acting on
% all of $(\xx,\yy)$" by $\Symm_{(\infty,\infty)}$ and
% explained what it means (as there is not really an obviously
% intended way to identify these two groups).

% [DG][v40] Added the footnote, as the permutation that takes
% one bijection to another might not have finite support.

% [DG][v75] Expanded the above explanation to give more details
% and to do fewer things in a single sentence.

\begin{example}
One has
\begin{align*}
\Delta m_{(2,1)} &= m_{(2,1)}(x_1,x_2,\ldots,y_1,y_2,\ldots) \\
&= x_1^2 x_2 + x_1 x_2^2 + \cdots \\
 &\quad + x_1^2 y_1 + x_1^2 y_2 + \cdots \\
 &\quad + x_1 y_1^2 + x_1 y_2^2 + \cdots \\
 &\quad + y_1^2 y_2 + y_1 y_2^2 + \cdots\\
&= m_{(2,1)}(\xx) + m_{(2)}(\xx) m_{(1)}(\yy) +
 m_{(1)}(\xx) m_{(2)}(\yy) + m_{(2,1)}(\yy) \\
&= m_{(2,1)} \otimes \one + m_{(2)} \otimes m_{(1)} +
m_{(1)}\otimes m_{(2)} + \one \otimes m_{(2,1)}.
\end{align*}
This example generalizes easily to the following formula:
\begin{equation}
\label{comultiplication-of-monomial-symmetric-function}
\Delta m_\lambda =
  \sum\limits_{\substack{ (\mu, \nu):\\ \mu \sqcup \nu = \lambda} } m_\mu \otimes m_\nu,
\end{equation}
in which \dfn{$\mu \sqcup \nu$}
is the partition obtained by taking the
multiset union  of the parts of $\mu$ and $\nu$, and then
reordering them to make them weakly decreasing.

\end{example}

Checking that $\Delta$ is coassociative amounts to checking that
\[
(\Delta \otimes \id) \circ \Delta f =
f(\xx,\yy,\zz) =
(\id \otimes \Delta) \circ \Delta f
\]
inside $\Lambda(\xx,\yy,\zz)$ as a subring of
$\Lambda \otimes \Lambda \otimes \Lambda$.

\begin{verlong}
Here is this argument in more detail:
\begin{proof}[First proof of the coassociativity of $\Delta$ on $\Lambda$.] 
Just as we identified the ring
$\Lambda\left(  \xx ,\yy \right)  $ with a subring of $\Lambda\otimes\Lambda$, we can
identify the ring $\Lambda\left(  \xx ,\yy ,\zz \right)  $
with a subring of $\Lambda\otimes\Lambda\otimes\Lambda$ (where $\left(
\xx ,\yy ,\zz \right)  $ denotes the family $\left(
x_{1},x_{2},...,y_{1},y_{2},...,z_{1},z_{2},...\right)  $ of indeterminates).
Let $f\in\Lambda$. We write $\Delta f\in\Lambda\otimes\Lambda$ in the form
$\Delta f=\sum\limits_{i\in I}g_{i}\otimes h_{i}$. Then, due to our
identification of the ring $\Lambda\left(  \xx ,\yy \right)  $
with a subring of $\Lambda\otimes\Lambda$, we have $\sum\limits_{i\in I}%
g_{i}\left(  \xx \right)  h_{i}\left(  \yy \right)
=\sum\limits_{i\in I}g_{i}\otimes h_{i}
=\Delta f=f\left(  \xx ,\yy \right)  $.
If we substitute $\left(  \xx ,\yy \right)  $ and
$\zz $ for $\xx $ and $\yy $ in this identity of formal
power series, then we
obtain $\sum\limits_{i\in I}g_{i}\left(  \xx ,\yy \right)
h_{i}\left(  \zz \right)  =f\left(  \left(  \xx ,\yy %
\right) , \zz\right)  =f\left(  \xx ,\yy ,\zz \right)  $. Now,
\begin{align*}
\left(  \Delta\otimes\id \right)  \underbrace{\left(  \Delta
f\right)  }_{=\sum\limits_{i\in I}g_{i}\otimes h_{i}}  &  =\left(
\Delta\otimes\id \right)  \left(  \sum\limits_{i\in I}%
g_{i}\otimes h_{i}\right)  =\sum\limits_{i\in I}\underbrace{\Delta\left(
g_{i}\right)  }_{=g_{i}\left(  \xx ,\yy \right)  }\otimes
h_{i}=\sum\limits_{i\in I}g_{i}\left(  \xx ,\yy \right)  \otimes
h_{i}\\
&  =\sum\limits_{i\in I}g_{i}\left(  \xx ,\yy \right)
h_{i}\left(  \zz \right) \\
&  \ \ \ \ \ \ \ \ \ \ \left(  \text{due to our identification of }%
\Lambda\left(  \xx ,\yy ,\zz \right)  \text{ with a subring
of }\Lambda\otimes\Lambda\otimes\Lambda\right) \\
&  =f\left(  \xx ,\yy ,\zz \right)
\end{align*}
and similarly $\left(  \id \otimes\Delta\right)  \left(  \Delta
f\right)  =f\left(  \xx ,\yy ,\zz \right)  $. Thus, $\left(
\Delta\otimes\id \right)  \left(  \Delta f\right)  =f\left(
\xx ,\yy ,\zz \right)  =\left(  \id \otimes\Delta\right)
\left(  \Delta f\right)  $. Since we have proven this
for every $f\in\Lambda$, we thus obtain $\left(  \Delta\otimes
\id \right)  \circ\Delta=\left(  \id \otimes
\Delta\right)  \circ\Delta$. In other words, the coassociativity of $\Delta$
on $\Lambda$ is proved.
\end{proof}

Alternatively, we can also prove the coassociativity of $\Delta$
using \eqref{comultiplication-of-monomial-symmetric-function}:

\begin{proof}[Second proof of the coassociativity of $\Delta$ on $\Lambda$.]
Every $\lambda\in\Par$ satisfies
\begin{align*}
&  \left(  \Delta\otimes\id \right)  \underbrace{\left(  \Delta
m_{\lambda}\right)  }_{\substack{=\sum\limits_{\substack{\left(  \mu,\gamma\right)
\in\Par \times\Par ;\\\mu\sqcup\gamma=\lambda}}
m_{\mu}\otimes m_{\gamma}
\\ \text{(by \eqref{comultiplication-of-monomial-symmetric-function})}
}}\\
&  =\left(  \Delta\otimes\id \right)  \left(
\sum\limits_{\substack{\left(  \mu,\gamma\right)  \in\Par %
\times\Par ;\\\mu\sqcup\gamma=\lambda}}m_{\mu}\otimes m_{\gamma
}\right)  =\underbrace{\sum\limits_{\substack{\left(  \mu,\gamma\right)
\in\Par \times\Par ;\\\mu\sqcup\gamma=\lambda}%
}}_{=\sum\limits_{\gamma\in\Par }\sum\limits_{\substack{\mu
\in\Par ;\\\mu\sqcup\gamma=\lambda}}}\underbrace{\Delta\left(
m_{\mu}\right)  }_{\substack{
=\sum\limits_{\substack{\left(  \alpha,\beta\right)
\in\Par \times\Par ;\\\alpha\sqcup\beta=\mu
}}m_{\alpha}\otimes m_{\beta}
\\ \text{(by \eqref{comultiplication-of-monomial-symmetric-function})}}}
\otimes m_{\gamma}\\
&  =\sum\limits_{\gamma\in\Par }\underbrace{
\sum\limits_{\substack{\mu\in\Par ;\\\mu\sqcup\gamma=\lambda}%
}\sum\limits_{\substack{\left(  \alpha,\beta\right)  \in\Par %
\times\Par ;\\\alpha\sqcup\beta=\mu}}}_{=\sum\limits_{\mu
\in\Par }\sum\limits_{\substack{\left(  \alpha,\beta\right)
\in\Par \times\Par ;\\\alpha\sqcup\beta=\mu
;\ \mu\sqcup\gamma=\lambda}}}m_{\alpha}\otimes m_{\beta}\otimes m_{\gamma}\\
&  =\sum\limits_{\gamma\in\Par }\sum\limits_{\mu\in
\Par }\underbrace{\sum\limits_{\substack{\left(  \alpha
,\beta\right)  \in\Par \times\Par ;\\\alpha
\sqcup\beta=\mu;\ \mu\sqcup\gamma=\lambda}}}_{\substack{=
\sum\limits_{\substack{\left(  \alpha,\beta\right)  \in\Par %
\times\Par ;\\\alpha\sqcup\beta=\mu;\ \alpha\sqcup\beta
\sqcup\gamma=\lambda}}\\\text{(because the assertion }\left(  \alpha
\sqcup\beta=\mu\text{ and }\mu\sqcup\gamma=\lambda\right)  \\\text{is
equivalent to }\left(  \alpha\sqcup\beta=\mu\text{ and}\ \alpha\sqcup
\beta\sqcup\gamma=\lambda\right)  \text{)}}}
m_{\alpha}\otimes m_{\beta}\otimes m_{\gamma}\\
&  =\sum\limits_{\gamma\in\Par }\underbrace{\sum\limits_{\mu
\in\Par }\sum\limits_{\substack{\left(  \alpha,\beta\right)
\in\Par \times\Par ;\\\alpha\sqcup\beta
=\mu;\ \alpha\sqcup\beta\sqcup\gamma=\lambda}}}_{=
\sum\limits_{\substack{\left(  \alpha,\beta\right)  \in\Par %
\times\Par ;\\\alpha\sqcup\beta\sqcup\gamma=\lambda}%
}\sum\limits_{\substack{\mu\in\Par ;\\\alpha\sqcup\beta=\mu}%
}}m_{\alpha}\otimes m_{\beta}\otimes m_{\gamma}=\sum\limits_{\gamma
\in\Par }\sum\limits_{\substack{\left(  \alpha,\beta\right)
\in\Par \times\Par ;\\\alpha\sqcup\beta
\sqcup\gamma=\lambda}}\underbrace{\sum\limits_{\substack{\mu\in
\Par ;\\\alpha\sqcup\beta=\mu}}m_{\alpha}\otimes m_{\beta
}\otimes m_{\gamma}}_{=m_{\alpha}\otimes m_{\beta}\otimes m_{\gamma}}\\
&  =\underbrace{\sum\limits_{\gamma\in\Par }
\sum\limits_{\substack{\left(  \alpha,\beta\right)  \in\Par %
\times\Par ;\\\alpha\sqcup\beta\sqcup\gamma=\lambda}}}%
_{=\sum\limits_{\substack{\left(  \alpha,\beta,\gamma\right)  \in
\Par \times\Par \times\Par %
;\\\alpha\sqcup\beta\sqcup\gamma=\lambda}}}m_{\alpha}\otimes m_{\beta}\otimes
m_{\gamma}=\sum\limits_{\substack{\left(  \alpha,\beta,\gamma\right)
\in\Par \times\Par \times\Par %
;\\\alpha\sqcup\beta\sqcup\gamma=\lambda}}m_{\alpha}\otimes m_{\beta}\otimes
m_{\gamma}
\end{align*}
and (similarly) $\left(  \id \otimes\Delta\right)  \left(
\Delta m_{\lambda}\right)  =\sum\limits_{\substack{\left(  \alpha,\beta
,\gamma\right)  \in\Par \times\Par \times
\Par ;\\\alpha\sqcup\beta\sqcup\gamma=\lambda}}m_{\alpha}\otimes
m_{\beta}\otimes m_{\gamma}$. Thus, every $\lambda\in\Par$
satisfies
\[
\left(  \Delta\otimes\id \right)  \left(  \Delta m_{\lambda
}\right)  =\sum\limits_{\substack{\left(  \alpha,\beta,\gamma\right)
\in\Par \times\Par \times\Par ;\\\alpha\sqcup\beta\sqcup\gamma=\lambda}}
m_{\alpha}\otimes m_{\beta}\otimes m_{\gamma}
=\left(  \id \otimes\Delta\right)  \left(  \Delta m_{\lambda}\right)  .
\]
Since $\left(  m_{\lambda}\right)  _{\lambda\in\Par }$ is a
basis of $\Lambda$, this yields that the maps $\left(  \Delta\otimes
\id \right)  \circ\Delta$ and
$\left(  \id \otimes\Delta\right)  \circ\Delta$ are equal to each
other on a basis of $\Lambda$. Thus,
$\left(  \Delta\otimes\id \right)  \circ
\Delta=\left(  \id \otimes\Delta\right)  \circ\Delta$. In other
words, the coassociativity of $\Delta$ on $\Lambda$ is proven again.
\end{proof}
\end{verlong}

The counit $\Lambda \overset{\epsilon}{\rightarrow} \kk$
is defined in the usual fashion for connected graded coalgebras,
namely $\epsilon$ annihilates $I=\bigoplus_{n > 0} \Lambda_n$, and
$\epsilon$ is the identity on $\Lambda_0=\kk$;  alternatively $\epsilon$
sends a symmetric function $f(\xx)$ to its constant term $f(0,0,\ldots)$.

Note that $\Delta$ is an algebra morphism
$\Lambda \rightarrow \Lambda \otimes \Lambda$
because it is a composition of maps which are all algebra morphisms.
As the unit and counit axioms are easily checked,
$\Lambda$ becomes a connected graded $\kk$-bialgebra of finite type,
and hence also a Hopf algebra by
Proposition~\ref{graded-connected-bialgebras-have-antipodes}.
We will identify its antipode more explicitly in Section~\ref{symm-antipode-section} below.

\subsection{Other Bases}

We introduce the usual other bases of $\Lambda$,
and explain their significance later.

\begin{definition}
\label{def.Lambda.peh}
Define the families of \dfn{power sum symmetric functions} $p_n$,
\dfn{elementary symmetric functions} $e_n$,
and \dfn{complete homogeneous symmetric functions} $h_n$,
for $n=1,2,3,\ldots$ by\index{$p_n$}\index{$e_n$}\index{$h_n$}
\begin{align}
p_n &:= x_1^n +x_2^n + \cdots = m_{(n)} ,
\label{def.pn} \\
e_n &:= \sum_{i_1 < \cdots < i_n}x_{i_1} \cdots x_{i_n} = m_{(1^n)} ,
\label{def.en} \\
h_n &:= \sum_{i_1 \leq \cdots \leq i_n}x_{i_1} \cdots x_{i_n} = \sum_{\lambda \in \Par_n} m_\lambda .
\label{def.hn}
\end{align}
Here, we are using the \dfn{multiplicative notation} for
partitions: whenever $\left(m_1, m_2, m_3, \ldots\right)$ is a
weak composition, $\left(1^{m_1}2^{m_2}3^{m_3}\cdots\right)$
denotes the partition $\lambda$ such that for every $i$, the
multiplicity of the part $i$ in $\lambda$ is $m_i$. The
$i^{m_i}$ satisfying $m_i = 0$ are often omitted from this notation,
and so the $\left(1^n\right)$ in \eqref{def.en} means
$\left(\underbrace{1, 1, \ldots, 1}_{n \text{ ones}}\right)$.
(For another example, $\left(1^2 3^1 4^3\right) =
\left(1^2 2^0 3^1 4^3 5^0 6^0 7^0 \cdots \right)$
means the partition $\left(4,4,4,3,1,1\right)$.)
% where $(1^n)=(1,1,\ldots,1)$, using a multiplicative notation
% $\lambda=(1^{m_1}2^{m_2}\cdots)$ if the multiplicity of the part $i$ in $\lambda$ is $m_i$.
By convention, also define $h_0=e_0=1$, and $h_n=e_n=0$ if $n < 0$.
Extend these multiplicatively to
partitions $\lambda=(\lambda_1,\lambda_2,\ldots,\lambda_\ell)$
with $\lambda_1 \geq \cdots \geq \lambda_\ell > 0$
by setting\index{$p_\lambda$}\index{$e_\lambda$}\index{$h_\lambda$}
\begin{align*}
p_\lambda&:=p_{\lambda_1} p_{\lambda_2} \cdots p_{\lambda_\ell} ,\\
e_\lambda&:=e_{\lambda_1} e_{\lambda_2} \cdots e_{\lambda_\ell} ,\\
h_\lambda&:=h_{\lambda_1} h_{\lambda_2} \cdots h_{\lambda_\ell} .\\
\end{align*}
Also define the \dfn{Schur function}\index{$s_\lambda$}
\begin{equation}
\label{Schur-function-definition}
s_\lambda:=\sum_T \xx^{\cont(T)}
\end{equation}
where $T$ runs through all \emph{column-strict tableaux}\index{column-strict tableau}
of shape $\lambda$, that is,
$T$ is an assignment of entries in $\{1,2,3,\ldots\}$
to the cells of the \emph{Ferrers diagram}\footnote{The
\dfn{Ferrers diagram} of a partition $\lambda$ is defined as
the set of all pairs
$\left(i, j\right) \in \left\{ 1, 2, 3, \ldots \right\}^2$
satisfying $j \leq \lambda_i$.
This is a set of cardinality $\left|\lambda\right|$.
Usually, one visually represents a Ferrers diagram by drawing its
elements $\left(i, j\right)$ as points on the plane, although
(unlike the standard convention for drawing points on the plane)
one lets the x-axis go top-to-bottom (i.e., the point
$\left(i+1, j\right)$ is one step below the point
$\left(i, j\right)$), and the y-axis go left-to-right (i.e.,
the point $\left(i, j+1\right)$ is one step to the right of the
point $\left(i, j\right)$). (This is the so-called
\dfn{English notation}, also known as the
\dfn{matrix notation} because it is precisely the way one
labels the entries of a matrix. Other notations appear in
literature, such as the French notation used, e.g., in
Malvenuto's \cite{Malvenuto}, and the Russian notation used,
e.g., in parts of Kerov's \cite{Kerov}.) These points are
drawn either as dots or as square boxes; in the latter case, the
boxes are centered at the points they represent, and they have
sidelength $1$ so that the boxes centered around
$\left(i, j\right)$ and $\left(i, j+1\right)$ touch each other
along a sideline. For example, the Ferrers diagram of the
partition $\left(3, 2, 2\right)$ is represented as
\[
\begin{matrix}
\bullet & \bullet & \bullet \\
\bullet & \bullet & \\
\bullet & \bullet & \\
\end{matrix}
\text{ (using dots) }
\qquad \qquad \text{ or as } \qquad \qquad
\begin{tabular}{ccc}
\hline
\multicolumn{1}{|c|}{} & \multicolumn{1}{c|}{} & \multicolumn{1}{c|}{}  \\
\hline
\multicolumn{1}{|c|}{} & \multicolumn{1}{c|}{} &  \\
\cline{1-2}
\multicolumn{1}{|c|}{} & \multicolumn{1}{c|}{} &  \\
\cline{1-2}
\end{tabular}
\text{ (using boxes).}
\]
The Ferrers diagram of a partition $\lambda$ uniquely determines
$\lambda$. One refers to the elements of the Ferrers diagram of
$\lambda$ as the \emph{cells}\index{cell of a Ferrers diagram}
(or \emph{boxes}\index{box of a Ferrers diagram}) of this
diagram (which is particularly natural when one represents them
by boxes) or, briefly, as the cells of $\lambda$. Notation like
``west'', ``north'', ``left'', ``right'', ``row'' and
``column'' concerning cells of Ferrers diagrams normally refers
to their visual representation.

Ferrers diagrams are also known as \emph{Young diagrams}\index{Young diagram}.

One can characterize the Ferrers diagrams of partitions as
follows: A finite subset $S$ of $\left\{ 1, 2, 3, \ldots \right\}^2$
is the Ferrers diagram of some partition if and only if for
every $\left(i, j\right) \in S$ and every $\left(i', j'\right)
\in \left\{ 1, 2, 3, \ldots \right\}^2$ satisfying $i' \leq i$
and $j' \leq j$, we have $\left(i', j'\right) \in S$. In other
words, a finite subset $S$ of $\left\{ 1, 2, 3, \ldots \right\}^2$
is the Ferrers diagram of some partition if and only if it is
a lower set of the poset $\left\{ 1, 2, 3, \ldots \right\}^2$
with respect to the componentwise order.}
 for $\lambda$, weakly
increasing left-to-right in rows, and strictly increasing top-to-bottom
in columns.  Here \dfn{$\cont(T)$} denotes the weak composition
$\left(|T^{-1}(1)|, |T^{-1}(2)|, |T^{-1}(3)|, \ldots\right)$, so that
$\xx^{\cont(T)}=\prod_i x_i^{|T^{-1}(i)|}$.  For
example,\footnote{To visually represent a column-strict tableau
$T$ of shape $\lambda$, we draw the same picture as when
representing the Ferrers diagram of $\lambda$, but with a little
difference: a cell $\left(i, j\right)$ is no longer represented
by a dot or box, but instead is represented by the entry
of $T$ assigned to this cell. Accordingly, the entry of $T$
assigned to a given cell $c$ is often referred to as
\emph{the entry of $T$ in $c$}\index{entry of a tableau}.}
\[
T=\begin{matrix}
1 & 1 & 1 & 4 & 7\\
2 & 3 & 3 &   & \\
4 & 4 & 6 &   & \\
6 & 7
\end{matrix}
\]
is a column-strict tableau of shape $\lambda=(5,3,3,2)$
with $\xx^{\cont(T)}=x_1^3 x_2^1 x_3^2 x_4^3 x_5^0 x_6^2 x_7^2$.
If $T$ is a column-strict tableau, then the weak composition $\cont(T)$
is called the \emph{content}\index{content of a tableau} of $T$.
\end{definition}

% [DG][v19] Added explicit definition of $\cont(T)$ (not just through its
% monomial).

% [DG][v25] Made the definition of multiplicative notation clearer.
% Unfortunately, switching from {array} to {align} in the definitions
% of p_n, e_n, h_n took the alignment of the second equality signs with
% it; do you have an idea how to get it back? (I want {align} due to
% the possibility of labelling each equation.)

% [DG][v32] Added two footnotes defining Ferrers diagrams and
% notations around tableaux.

Column-strict tableaux are also known as
\emph{semistandard tableaux}\index{semistandard tableau},
and some authors even
omit the adjective and just call them \emph{tableaux}\index{tableau}
(e.g., Fulton in \cite{Fulton}, a book entirely devoted
to them).

\begin{example}
One has
\begin{align*}
m_{(1)}&=p_{(1)}=e_{(1)}=h_{(1)}=s_{(1)}=x_1+x_2+x_3+\cdots , \\
s_{(n)}&=h_n ,\\
s_{(1^n)}&=e_n .
\end{align*}
\end{example}

\begin{example}
\label{transition-matrices-example}
One has for $\lambda=(2,1)$ that
\[
\begin{array}{rl}
p_{(2,1)}&=p_2 p_1 = (x_1^2+x_2^2+\cdots)(x_1+x_2+\cdots) \\
         &          = m_{(2,1)}+m_{(3)} ,\\
& \\
e_{(2,1)}&=e_2 e_1 = (x_1 x_2+x_1 x_3 + \cdots)(x_1+x_2+\cdots)\\
         &          = m_{(2,1)}+3m_{(1,1,1)} ,\\
& \\
h_{(2,1)}&=h_2 h_1 = (x_1^2 +x_2^2+\cdots + x_1x_2+x_1 x_3 + \cdots)(x_1+x_2+\cdots)\\
         &          = m_{(3)}+2m_{(2,1)}+3m_{(1,1,1)} ,\\

\end{array}
\]
and
\begin{align*}
&\begin{matrix}
s_{(2,1)}&=x_1^2x_2 &+x_1^2x_3&+x_1x_2^2&+x_1x_3^2&+x_1x_2x_3&+x_1x_2x_3&+x_1x_2x_4&+\cdots\\
&11            &11       &12            &13       &12        &13        &12 & \\
&2\phantom{1}        &3\phantom{1}   &2\phantom{1}        &3\phantom{1}   &3\phantom{1}    &2\phantom{1}    &4\phantom{1}  & \\
\end{matrix} \\
& \qquad \quad = m_{(2,1)}+2m_{(1,1,1)} .
\end{align*}
In fact, one has these transition matrices for $n=3$ expressing elements in terms of the monomial basis $m_\lambda$:
\[
\bordermatrix{
~          & p_{(3)} & p_{(2,1)} & p_{(1,1,1)} \cr
m_{(3)}    & 1       & 1         & 1  \cr
m_{(2,1)}  & 0       & 1         & 3\cr
m_{(1,1,1)}& 0       & 0         & 6
}
\ ,
\qquad
\bordermatrix{
~          & e_{(3)} & e_{(2,1)} & e_{(1,1,1)} \cr
m_{(3)}    & 0       & 0         & 1  \cr
m_{(2,1)}  & 0       & 1         & 3\cr
m_{(1,1,1)}& 1       & 3         & 6
}
\ ,
\]
\vskip.1in
\[
\bordermatrix{
~          & h_{(3)} & h_{(2,1)} & h_{(1,1,1)} \cr
m_{(3)}    & 1       & 1         & 1  \cr
m_{(2,1)}  & 1       & 2         & 3\cr
m_{(1,1,1)}& 1       & 3         & 6
}
\ ,
\qquad
\bordermatrix{
~          & s_{(3)} & s_{(2,1)} & s_{(1,1,1)} \cr
m_{(3)}    & 1       & 0         & 0  \cr
m_{(2,1)}  & 1       & 1         & 0\cr
m_{(1,1,1)}& 1       & 2         & 1
}
\ .
\]
\end{example}

Our next goal is to show that $e_\lambda, s_\lambda, h_\lambda$
(and, under some conditions, the $p_\lambda$ as well)
all give bases for $\Lambda$.  However at the moment it is not yet
even clear that $s_\lambda$ are symmetric!

% [DG][v48] Removed the $p_\lambda$ from the
% $p_\lambda, e_\lambda, s_\lambda, h_\lambda$ list and put them
% into parentheses.

\begin{proposition}
\label{Schur-functions-are-symmetric-prop}
Schur functions $s_\lambda$ are symmetric, that is, they lie
in $\Lambda$.
\end{proposition}
\begin{proof}
It suffices to show $s_\lambda$ is symmetric under swapping the variables
$x_i, x_{i+1}$, by providing an involution $\iota$ on the set
of all column-strict tableaux $T$ of shape $\lambda$ which switches
the $\cont(T)$ for $(i,i+1) \cont(T)$.
Restrict attention to the entries $i,i+1$
in $T$, which must look something like this:
\[
\begin{matrix}
   &   & & & & &   &   &i  &i  &i  &i&i+1&i+1 \\
   &i  &i&i&i&i&i+1&i+1&i+1&i+1&i+1& &   & \\
i+1&i+1&i+1& & &   &   &   &   &   & &   &
\end{matrix}
\]
One finds several vertically aligned pairs
$\begin{matrix} i\\ i+1 \end{matrix}$.  If one were to remove all such pairs, the remaining entries would be a sequence of rows, each looking like this:
\begin{equation}
\label{i-i+1-row-to-be-swapped}
\underbrace{i,i,\ldots,i}_{r\text{ occurrences}},
\underbrace{i+1,i+1,\ldots,i+1}_{s\text{ occurrences}} \, .
\end{equation}
An involution due to Bender and Knuth\index{Bender-Knuth involution}
tells us to leave fixed all the vertically aligned pairs
$\begin{matrix} i\\ i+1 \end{matrix}$,
but change each sequence of remaining entries
as in \eqref{i-i+1-row-to-be-swapped} to this:
\[
\underbrace{i,i,\ldots,i}_{s\text{ occurrences}},
\underbrace{i+1,i+1,\ldots,i+1}_{r\text{ occurrences}} \, .
\]
For example, the above configuration in $T$ would change
to
\[
\begin{matrix}
   &   & & & &   &   &   &i  &i  &i  &i  &i  &i+1 \\
   &i  &i&i&i&i+1&i+1&i+1&i+1&i+1&i+1&   &   & \\
i  &i+1&i+1& &   &   &   &   &   &   &   &   &
\end{matrix}
\]
It is easily checked that this map is an involution,
and that it has the effect of swapping $(i,i+1)$ in $\cont(T)$.
\end{proof}

\begin{remark}
\label{other-linear-orders-remark}
The symmetry of Schur functions allows one to
reformulate them via column-strict tableaux defined with respect to
\emph{any} total ordering $\LLL$ on the positive integers, rather than the usual
$1< 2 < 3<\cdots$.  For example, one can use
the \emph{reverse order}\footnote{This reverse order is what one uses when one defines a Schur function as a generating function for \emph{reverse semistandard tableaux} or \emph{column-strict plane partitions}; see Stanley \cite[Proposition 7.10.4]{Stanley}.}
$
\cdots <3<2<1,
$
or even more exotic orders, such as
\[
1<3<5<7<\cdots<2<4<6<8<\cdots.
\]
Say that an assignment $T$ of entries in $\{1,2,3,\ldots\}$ to the cells of
the Ferrers diagram of $\lambda$ is an \dfn{$\LLL$-column-strict tableau} if
it is weakly $\LLL$-increasing left-to-right in rows, and strictly
$\LLL$-increasing top-to-bottom in columns.

\begin{proposition}
\label{prop.other-linear-orders.slambda}
For any total order  $\LLL$ on the positive integers,
\begin{equation}
\label{different-total-order-Schur-function}
s_\lambda = \sum_T \xx^{\cont(T)}
\end{equation}
as $T$ runs through all $\LLL$-column-strict tableaux of shape $\lambda$.
\end{proposition}
\begin{proof}
Given a weak composition $\alpha=(\alpha_1,\alpha_2,\ldots)$  with
$\alpha_{n+1}=\alpha_{n+2}=\cdots=0$,
assume that the integers $1,2,\ldots,n$ are totally ordered by $\LLL$ as
$w(1) <_{\LLL} \cdots <_{\LLL} w(n)$ for some $w$ in $\Symm_n$.
Then the coefficient of
$\xx^{\alpha}=x_1^{\alpha_1} \cdots x_{n}^{\alpha_n}$
on the right side of \eqref{different-total-order-Schur-function}
is the same as the coefficient of $\xx^{w^{-1}(\alpha)}$
on the right side of \eqref{Schur-function-definition}
defining $s_\lambda$,
which by symmetry of $s_\lambda$ is the same
as the coefficient of $\xx^{\alpha}$ on the right side of \eqref{Schur-function-definition}.
\end{proof}
\end{remark}

It is now not hard to show that
$p_\lambda, e_\lambda, s_\lambda$ give bases by a triangularity
argument\footnote{See Section~\ref{sect.STmat} for some notions and
notations that will be used in this argument.}.  For this purpose,
let us introduce a useful partial order on partitions.

\begin{definition}
\label{def.dominance}
The \dfn{dominance} or \dfn{majorization} order on $\Par_n$
is the partial order on the set $\Par_n$ whose
greater-or-equal relation \dfn{$\triangleright$} is defined as follows:
For two partitions $\lambda$ and $\mu$ of $n$, we set
$\lambda \triangleright \mu$ (and say that $\lambda$
\emph{dominates}\index{dominate}, or \emph{majorizes}\index{majorize}, $\mu$)
if and only if
\[
\lambda_1+\lambda_2+ \cdots+\lambda_k \geq \mu_1 + \mu_2 + \cdots+\mu_k
\quad \text{ for } k=1,2,\ldots,n.
\]
\end{definition}

(The definition of dominance would not change if we would
replace ``for $k=1,2,\ldots,n$'' by ``for every positive
integer $k$'' or by ``for every $k \in \NN$''.)

% [DG][v29] Added the previous sentence.

% [DG][v56] Restated definition slightly.

\begin{definition}
\label{def.partition.conjugate}
For a partition $\lambda$, its
\emph{conjugate}\index{conjugate of a partition}
or \emph{transpose}\index{transpose of a partition}
partition \dfn{$\lambda^t$} is the one whose
Ferrers diagram is obtained from that of $\lambda$
by exchanging rows for columns
(i.e., by flipping the diagram across the ``main'', i.e.,
top-right-to-bottom-left, diagonal)\footnote{In more
rigorous terms:
The cells of the Ferrers diagram of $\lambda^t$ are the
pairs $\left(j, i\right)$, where $\left(i, j\right)$ ranges
over all cells of $\lambda$.
It is easy to see that this indeed uniquely determines
a partition $\lambda^t$.}.
% [DG][v62] Added the preceding footnote.
Alternatively, one has this formula
for its $i$-th entry:
\begin{equation}
\label{eq.partition.transpose}
(\lambda^t)_i := |\{ j: \lambda_j \geq i \}|.
\end{equation}
\end{definition}

% [DG][v23] Removed the comma before "is the one"; please let me know if
% it was intentional.

For example, $\left( 4,3,1 \right)^t = \left( 3,2,2,1 \right)$,
which can be easily verified by flipping the Ferrers diagram
of $\left( 4,3,1 \right)$ across the ``main diagonal'':
\[
\underbrace{
\begin{matrix}
\bullet & \bullet & \bullet & \bullet \\
\bullet & \bullet & \bullet \\
\bullet & \\
\end{matrix}
}_{\text{Ferrers diagram of $\left( 4,3,1 \right)$}}
\qquad \longmapsto \qquad
\underbrace{
\begin{matrix}
\bullet & \bullet & \bullet \\
\bullet & \bullet & \\
\bullet & \bullet & \\
\bullet
\end{matrix}
}_{\text{Ferrers diagram of $\left( 4,2,2,1 \right)$}}
\]
(or simply counting the boxes in each column of this diagram).

% [DG][v78] Added above example. Also added the "flipping"
% description of conjugation to the definition.

\begin{exercise}
\label{exe.partition.dominance}
Let $\lambda , \mu \in \Par_n$. Show that
$\lambda \triangleright \mu$ if and only if
$\mu^t \triangleright \lambda^t$.
\end{exercise}

% [DG][v23] Made the above paragraph into an actual exercise.

\begin{proposition}
\label{many-symmetric-function-bases-prop}
The families $\{e_\lambda\}$ and $\{s_\lambda\}$, as $\lambda$ runs through
all partitions, are graded bases for the
graded $\kk$-module $\Lambda_\kk$
whenever $\kk$ is a commutative ring.
The same holds for the family $\{p_\lambda\}$ when $\QQ$ is a
subring of $\kk$.
\end{proposition}

% [DG][v80] Replaced "bases" by "graded bases", since we
% prove this stronger fact in fact (and it is more useful).

Our proof of this proposition will involve three separate
arguments, one for each of the three alleged bases
$\{s_\lambda\}$, $\{e_\lambda\}$ and $\{p_\lambda\}$; however,
all these three arguments fit the same mold: Each one shows
that the alleged basis expands invertibly
triangularly\footnote{i.e., triangularly, with all diagonal
coefficients being invertible} in the
basis $\{m_\lambda\}$ (possibly after reindexing), with an
appropriately chosen partial order on the indexing set.
We will simplify our life by restricting ourselves to $\Par_n$
for a given $n \in \NN$, and by stating the common part of the
three arguments in a greater generality (so that we won't have
to repeat it thrice):

\begin{lemma}
\label{lem.many-symmetric-function-bases-prop.last-step}
Let $S$ be a finite poset. We write $\leq$ for the
smaller-or-equal relation of $S$.

Let $M$ be a free $\kk$-module with
a basis $\left( b_\lambda \right)_{\lambda \in S}$.
Let $\left( a_\lambda \right)_{\lambda \in S}$ be a further
family of elements of $M$.

For each $\lambda \in S$, let
$\left( g_{\lambda, \mu} \right)_{\mu \in S}$
be the family of the coefficients in the
expansion of $a_\lambda \in M$ in the basis
$\left( b_\mu \right)_{\mu \in S}$; in other words,
let $\left( g_{\lambda, \mu} \right)_{\mu \in S} \in \kk^S$
be such that
$a_\lambda = \sum\limits_{\mu \in S} g_{\lambda, \mu} b_\mu$.
Assume that:
\begin{itemize}
\item \textit{Assumption A1:} Any $\lambda \in S$ and $\mu \in S$ satisfy $g_{\lambda, \mu} = 0$ unless $\mu \leq \lambda$.
\item \textit{Assumption A2:} For any $\lambda \in S$, the element $g_{\lambda, \lambda}$ of $\kk$ is invertible.
\end{itemize}

Then, the family $\left( a_\lambda \right)_{\lambda \in S}$
is a basis of the $\kk$-module $M$.
\end{lemma}

\begin{proof}[Proof of Lemma~\ref{lem.many-symmetric-function-bases-prop.last-step}.]
Use the notations of Section~\ref{sect.STmat}.
Assumptions A1 and A2 yield that the
$S \times S$-matrix
$\left(g_{\lambda, \mu}\right)_{\left(\lambda, \mu\right) \in S \times S} \in \kk^{S \times S}$
is invertibly triangular. But the definition of the
$g_{\lambda, \mu}$ yields that the
family $\left( a_\lambda \right)_{\lambda \in S}$ expands in the family
$\left( b_\lambda \right)_{\lambda \in S}$ through this matrix
$\left(g_{\lambda, \mu}\right)_{\left(\lambda, \mu\right) \in S \times S}$.
Since the latter matrix is invertibly triangular, this shows
that the family $\left( a_\lambda \right)_{\lambda \in S}$ expands
invertibly triangularly in the family $\left( b_\lambda \right)_{\lambda \in S}$.
Therefore, Corollary~\ref{cor.STmat.expansion-tria-inv}(e) (applied
to
$\left( e_s \right)_{s \in S} = \left( a_\lambda \right)_{\lambda \in S}$
and
$\left( f_s \right)_{s \in S} = \left( b_\lambda \right)_{\lambda \in S}$)
shows that
$\left( a_\lambda \right)_{\lambda \in S}$ is a basis of the $\kk$-module $M$
(since $\left( b_\lambda \right)_{\lambda \in S}$ is a basis of
the $\kk$-module $M$).
\end{proof}

\begin{proof}[Proof of Proposition~\ref{many-symmetric-function-bases-prop}.]
We can restrict our attention to each homogeneous
component $\Lambda_n$
and partitions $\lambda$ of $n$.  Thus, we have to prove that,
for each $n \in \NN$, the families
$\left( e_\lambda \right)_{\lambda \in \Par_n}$ and
$\left( s_\lambda \right)_{\lambda \in \Par_n}$
are bases of the $\kk$-module $\Lambda_n$, and that the same
holds for
$\left( p_\lambda \right)_{\lambda \in \Par_n}$
if $\QQ$ is a subring of $\kk$.

% We check that in each case,
% Lemma~\ref{lem.many-symmetric-function-bases-prop.last-step}
% applies (with $S = \Par_n$ for some partial order on $\Par_n$,
% with $M = \Lambda_n$ and with $\left( b_s \right)_{s \in S}
% = \left( m_\lambda \right)_{\lambda \in \Par_n}$). This will
% then yield that the proposed basis is indeed a basis.
% the proposed basis expands invertibly triangularly\footnote{See
% Definition~\ref{def.STmat.expansion-tria}(b) for what this means.}
% in the basis $\left\{ m_\lambda \right\}_{\lambda \in \Par_n}$
% with some choice of orderings on $\Par_n$ indexing the rows and columns,
% as illustrated in Example~\ref{transition-matrices-example}.
% This will then yield (by Corollary~\ref{cor.STmat.expansion-tria-inv}(e))
% that the proposed basis is indeed a basis.

Fix $n \in \NN$. We already know that
$\left( m_\lambda \right)_{\lambda \in \Par_n}$ is a
basis of the $\kk$-module $\Lambda_n$.

\begin{enumerate}

\item[1.] We shall first show that
the family $\left( s_\lambda \right)_{\lambda \in \Par_n}$ is a
basis of the $\kk$-module $\Lambda_n$.

For every partition $\lambda$, we have
$s_\lambda = \sum_{\mu \in \Par} K_{\lambda,\mu} m_\mu$,
where the coefficient $K_{\lambda,\mu}$ is the \dfn{Kostka number} counting
the column-strict tableaux $T$ of shape $\lambda$ having $\cont(T)=\mu$;
this follows because both sides are symmetric functions,
and $K_{\lambda,\mu}$ is the coefficient of $\xx^\mu$ on both
sides\footnote{In general, in order to prove that two
symmetric functions $f$ and $g$ are equal, it suffices
to show that, for every $\mu \in \Par$, the coefficients
of $\xx^\mu$ in $f$ and in $g$ are equal. (Indeed, all
other coefficients are determined by these coefficients
because of the symmetry.)}.
Thus, for every $\lambda \in \Par_n$, one has
\begin{equation}
s_\lambda = \sum_{\mu \in \Par_n} K_{\lambda,\mu} m_\mu
\label{pf.many-symmetric-function-bases-prop.s-through-m}
\end{equation}
(since $s_\lambda$ is homogeneous of degree $n$).
\ \ \ \ \footnote{See Exercise~\ref{exe.Lambda.triangular}(c)
below for a detailed proof of
\eqref{pf.many-symmetric-function-bases-prop.s-through-m}.}
But if $\lambda$ and $\mu$ are partitions satisfying
$K_{\lambda, \mu} \neq 0$, then there exists a column-strict
tableau $T$ of shape $\lambda$ having $\cont(T) = \mu$
(since $K_{\lambda, \mu}$ counts such tableaux), and
therefore we must have
$\lambda_1+\lambda_2+ \cdots+\lambda_k \geq \mu_1 + \mu_2 + \cdots+\mu_k$
for each positive integer $k$
(since the entries $1,2,\ldots,k$ in $T$ must all lie within the
first $k$ rows of $\lambda$); in other words,
$\lambda \triangleright \mu$ (if $K_{\lambda,\mu} \neq 0$)
\ \ \ \ \footnote{See Exercise~\ref{exe.Lambda.triangular}(d)
below for a detailed proof of this fact.}.
In other words,
\begin{equation}
\text{any } \lambda \in \Par_n \text{ and }\mu \in \Par_n
\text{ satisfy }K_{\lambda, \mu} = 0
\text{ unless }\lambda \triangleright \mu .
\label{pf.many-symmetric-function-bases-prop.s-through-m.0}
\end{equation}
One can also check that $K_{\lambda, \lambda}=1$ for any
$\lambda \in \Par_n$
\ \ \ \ \footnote{See Exercise~\ref{exe.Lambda.triangular}(e)
below for a proof of this.}. Hence,
\begin{equation}
\text{for any } \lambda \in \Par_n \text{, the element }
K_{\lambda, \lambda} \text{ of } \kk \text{ is invertible.}
\label{pf.many-symmetric-function-bases-prop.s-through-m.1}
\end{equation}

Now, let us regard
the set $\Par_n$ as a poset, whose greater-or-equal relation
is $\triangleright$.
Lemma~\ref{lem.many-symmetric-function-bases-prop.last-step}
(applied to $S = \Par_n$, $M = \Lambda_n$,
$a_\lambda = s_\lambda$, $b_\lambda = m_\lambda$ and
$g_{\lambda, \mu} = K_{\lambda, \mu}$) shows that
the family $\left( s_\lambda \right)_{\lambda \in \Par_n}$ is
a basis of the $\kk$-module $\Lambda_n$
(because the Assumptions A1 and A2 of
Lemma~\ref{lem.many-symmetric-function-bases-prop.last-step}
are satisfied\footnote{Indeed, they follow from
\eqref{pf.many-symmetric-function-bases-prop.s-through-m.0}
and
\eqref{pf.many-symmetric-function-bases-prop.s-through-m.1},
respectively.}).

% Hence, the expansion \eqref{pf.many-symmetric-function-bases-prop.s-through-m}
% becomes
% $s_\lambda=\sum_{\mu: \lambda \triangleright \mu} K_{\lambda,\mu} m_\mu$.

\item[2.]
Before we show that $\left(e_\lambda\right)_{\lambda \in \Par_n}$
is a basis, we define a few
notations regarding integer matrices. A \dfn{$\{0,1\}$-matrix}
means a matrix whose entries belong to the set $\{0,1\}$. If
$A \in \NN^{\ell \times m}$ is a matrix, then the \dfn{row sums}
of $A$ means the $\ell$-tuple $\left(r_1, r_2, \ldots, r_\ell\right)$,
where each $r_i$ is the sum of all entries in the $i$-th row of $A$;
similarly, the \dfn{column sums} of $A$ means the $m$-tuple
$\left(c_1, c_2, \ldots, c_m\right)$, where each $c_j$ is the sum
of all entries in the $j$-th column of $A$.
(For instance, the row sums of the $\left\{0,1\right\}$-matrix
$\left( \begin{matrix} 0 & 1 & 1 & 0 & 0 \\
					   1 & 1 & 0 & 1 & 0 \end{matrix} \right)$
is $\left(2, 3\right)$, whereas its column sums is
$\left(1, 2, 1, 1, 0\right)$.)
We identify any $k$-tuple of nonnegative integers
$\left(a_1, a_2, \ldots, a_k\right)$ with the weak composition
$\left(a_1, a_2, \ldots, a_k, 0, 0, 0, \ldots\right)$;
thus, the row sums and the column sums of a matrix in
$\NN^{\ell \times m}$ can be viewed as weak compositions.
(For example, the column sums of the matrix
$\left( \begin{matrix} 0 & 1 & 1 & 0 & 0 \\
					   1 & 1 & 0 & 1 & 0 \end{matrix} \right)$
is the $5$-tuple $\left(1, 2, 1, 1, 0\right)$, and can be
viewed as the weak composition
$\left(1, 2, 1, 1, 0, 0, 0, \ldots \right)$.)

% [DG][v23] Added preceding paragraph, also replacing "row sum" by
% "row sums" below. I feel that this pedantry is helpful because
% "row sum" could also be understood as "entrywise sum of all rows",
% which is precisely what you call "column sums". (But I also feel
% it could have been done better.)

% [DG][v70] Added example and clarifications
% about how a partition can be the row/column sums of a
% finite matrix.

For every $\lambda \in \Par_n$, one has
\begin{equation}
e_\lambda = \sum_{\mu \in \Par_n} a_{\lambda,\mu} m_\mu ,
\label{pf.many-symmetric-function-bases-prop.e-through-m}
\end{equation}
where $a_{\lambda,\mu}$ counts $\{0,1\}$-matrices (of size
$\ell\left(\lambda\right) \times \ell\left(\mu\right)$) having
row sums $\lambda$ and column sums $\mu$: indeed, when one
expands $e_{\lambda_1} e_{\lambda_2} \cdots $,
choosing the monomial $x_{j_1} \ldots x_{j_{\lambda_i}}$ in the
$e_{\lambda_i}$ factor corresponds to putting
$1$'s in the $i$-th row and columns $j_1,\ldots,j_{\lambda_i}$ of the
$\{0,1\}$-matrix
\footnote{See Exercise~\ref{exe.Lambda.triangular}(g)
below for a detailed proof of
\eqref{pf.many-symmetric-function-bases-prop.e-through-m}.}.
Applying \eqref{pf.many-symmetric-function-bases-prop.e-through-m}
to $\lambda^t$ instead of $\lambda$, we see that
\begin{equation}
e_{\lambda^t} = \sum_{\mu \in \Par_n} a_{\lambda^t,\mu} m_\mu
\label{pf.many-symmetric-function-bases-prop.e-through-m.!}
\end{equation}
for every $\lambda \in \Par_n$. 

It is not hard to check\footnote{See
Exercise~\ref{exe.Lambda.triangular}(h)
below for a proof of this.
This is the easy implication in the \dfn{Gale-Ryser Theorem}.
(The hard implication is the converse: It
says that if $\lambda, \mu \in \Par_n$ satisfy
$\lambda^t \triangleright \mu$, then
there exists a $\left\{0,1\right\}$-matrix having
row sums $\lambda$ and column sums $\mu$, so that
$a_{\lambda,\mu}$ is a positive integer.
This is proven, e.g., in \cite{Krause-Gale}, in
\cite[Theorem 2.4]{Dahl-ComMat}
and in \cite[Section 5.2]{Wildon2016}.)}
that $a_{\lambda,\mu}$ vanishes unless
$\lambda^t \triangleright \mu$. Applying this to $\lambda^t$
instead of $\lambda$, we conclude that
\begin{equation}
\text{any } \lambda \in \Par_n \text{ and }\mu \in \Par_n
\text{ satisfy }a_{\lambda^t, \mu} = 0
\text{ unless }\lambda \triangleright \mu .
\label{pf.many-symmetric-function-bases-prop.e-through-m.0}
\end{equation}
Moreover, one can show that
$a_{\lambda^t, \lambda} = 1$ for each $\lambda \in \Par_n$
\ \ \ \ \footnote{See Exercise~\ref{exe.Lambda.triangular}(i)
below for a proof of this.}.
Hence,
\begin{equation}
\text{for any } \lambda \in \Par_n \text{, the element }
a_{\lambda^t, \lambda} \text{ of } \kk \text{ is invertible.}
\label{pf.many-symmetric-function-bases-prop.e-through-m.1}
\end{equation}

Now, let us regard
the set $\Par_n$ as a poset, whose greater-or-equal relation
is $\triangleright$.
Lemma~\ref{lem.many-symmetric-function-bases-prop.last-step}
(applied to $S = \Par_n$, $M = \Lambda_n$,
$a_\lambda = e_{\lambda^t}$, $b_\lambda = m_\lambda$ and
$g_{\lambda, \mu} = a_{\lambda^t, \mu}$) shows that
the family $\left( e_{\lambda^t} \right)_{\lambda \in \Par_n}$ is
a basis of the $\kk$-module $\Lambda_n$
(because the Assumptions A1 and A2 of
Lemma~\ref{lem.many-symmetric-function-bases-prop.last-step}
are satisfied\footnote{Indeed, they follow from
\eqref{pf.many-symmetric-function-bases-prop.e-through-m.0}
and
\eqref{pf.many-symmetric-function-bases-prop.e-through-m.1},
respectively.}).
Hence,
$\left( e_{\lambda} \right)_{\lambda \in \Par_n}$ is a
basis of $\Lambda_n$.

\item[3.]
Assume now that $\QQ$ is a subring of $\kk$.
For every $\lambda \in \Par_n$, one has
\begin{equation}
p_\lambda = \sum_{\mu \in \Par_n} b_{\lambda,\mu} m_\mu ,
\label{pf.many-symmetric-function-bases-prop.p-through-m}
\end{equation}
where $b_{\lambda,\mu}$ counts the ways to partition the nonzero parts
$\lambda_1,\ldots,\lambda_\ell$ (where $\ell = \ell\left(\lambda\right)$)
into blocks such that the sums of the blocks
give $\mu$; more formally, $b_{\lambda,\mu}$ is the number of maps
$\varphi: \{1,2,\ldots,\ell\} \rightarrow \{1,2,3,\ldots\}$ having
\[
\mu_j = \sum_{i: \varphi(i)=j} \lambda_i
\qquad \text{ for all }j=1,2,\ldots
\]
\footnote{See Exercise~\ref{exe.Lambda.triangular}(k)
below for a detailed proof of
\eqref{pf.many-symmetric-function-bases-prop.p-through-m}
(and see Exercise~\ref{exe.Lambda.triangular}(j) for a proof
that the numbers $b_{\lambda,\mu}$ are well-defined).}.
Again it is not hard to check that
\begin{equation}
\text{any } \lambda \in \Par_n \text{ and }\mu \in \Par_n
\text{ satisfy }b_{\lambda, \mu} = 0
\text{ unless }\mu \triangleright \lambda .
\label{pf.many-symmetric-function-bases-prop.p-through-m.0}
\end{equation}
\footnote{See Exercise~\ref{exe.Lambda.triangular}(l)
below for a proof of this.}
Furthermore, for any $\lambda \in \Par_n$, the element
$b_{\lambda, \lambda}$ is a positive integer\footnote{This
is proven in Exercise~\ref{exe.Lambda.triangular}(m)
below.}, and thus
invertible in $\kk$ (since $\QQ$ is a subring of $\kk$).
Thus,
\begin{equation}
\text{for any } \lambda \in \Par_n \text{, the element }
b_{\lambda, \lambda} \text{ of } \kk \text{ is invertible}
\label{pf.many-symmetric-function-bases-prop.p-through-m.1}
\end{equation}
(although we don't always have
$b_{\lambda, \lambda} = 1$ this time).

Now, let us regard
the set $\Par_n$ as a poset, whose smaller-or-equal relation
is $\triangleright$.
Lemma~\ref{lem.many-symmetric-function-bases-prop.last-step}
(applied to $S = \Par_n$, $M = \Lambda_n$,
$a_\lambda = p_\lambda$, $b_\lambda = m_\lambda$ and
$g_{\lambda, \mu} = b_{\lambda, \mu}$) shows that
the family $\left( p_\lambda \right)_{\lambda \in \Par_n}$ is
a basis of the $\kk$-module $\Lambda_n$
(because the Assumptions A1 and A2 of
Lemma~\ref{lem.many-symmetric-function-bases-prop.last-step}
are satisfied\footnote{Indeed, they follow from
\eqref{pf.many-symmetric-function-bases-prop.p-through-m.0}
and
\eqref{pf.many-symmetric-function-bases-prop.p-through-m.1},
respectively.}).

\end{enumerate}
\end{proof}

\begin{remark}
When $\QQ$ is not a subring of $\kk$, the family
$\left\{ p_\lambda \right\}$ is not (in general) a basis of $\Lambda_{\kk}$;
for instance, $e_2 = \frac12 \left(p_{(1,1)}-p_2\right) \in \Lambda_{\QQ}$
is not in the $\ZZ$-span of this family. However, if we define
$b_{\lambda,\mu}$ as in the above proof, then the $\ZZ$-linear span
of all $p_{\lambda}$ equals the $\ZZ$-linear span of all $b_{\lambda,
\lambda} m_{\lambda}$. Indeed, if
$\mu = \left(\mu_1, \mu_2, \ldots, \mu_k\right)$ with $k = \ell(\mu)$,
then $b_{\mu, \mu}$ is the size of the subgroup of $\Symm_k$
consisting of all permutations $\sigma \in \Symm_k$ having each $i$
satisfy $\mu_{\sigma(i)} = \mu_i$\ \ \ \ \footnote{See
Exercise~\ref{exe.Lambda.triangular}(n)
below for a proof of this.}. As a consequence, $b_{\mu, \mu}$
divides $b_{\lambda, \mu}$ for every partition $\mu$ of the same size as
$\lambda$ (because this group acts\footnote{Specifically, an element
$\sigma$ of the group takes
$\varphi: \{1,2,\ldots,\ell\} \rightarrow \{1,2,3,\ldots\}$ to
$\sigma \circ \varphi$.} freely on the set which is
enumerated by $b_{\lambda, \mu}$)\ \ \ \ \footnote{See
Exercise~\ref{exe.Lambda.triangular}(o)
below for a detailed proof of this.}. Hence, the
$\Par_n \times \Par_n$-matrix
$\left( \dfrac{b_{\lambda, \mu}}{b_{\mu, \mu}} \right)_{\left(\lambda, \mu\right) \in \Par_n \times \Par_n}$
has integer entries. Furthermore, this matrix is
unitriangular\footnote{Here, we are using the terminology defined
in Section~\ref{sect.STmat}, and we are regarding $\Par_n$ as a poset
whose smaller-or-equal relation is $\triangleright$.}
(indeed, \eqref{pf.many-symmetric-function-bases-prop.p-through-m.0}
shows that it is triangular, but its diagonal entries are clearly
$1$) and thus invertibly triangular. But
\eqref{pf.many-symmetric-function-bases-prop.p-through-m} shows
that the family $\left( p_\lambda \right)_{\lambda \in \Par_n}$
expands in the family
$\left( b_{\lambda, \lambda} m_\lambda \right)_{\lambda\in\Par_n}$
through this matrix. Hence, the family
$\left( p_\lambda \right)_{\lambda \in \Par_n}$
expands invertibly triangularly in the family
$\left( b_{\lambda, \lambda} m_\lambda \right)_{\lambda\in\Par_n}$.
Thus, Corollary~\ref{cor.STmat.expansion-tria-inv}(b) (applied
to $\ZZ$, $\Lambda_n$, $\Par_n$, $\left( p_\lambda \right)_{\lambda \in \Par_n}$
and $\left( b_{\lambda, \lambda} m_\lambda \right)_{\lambda \in \Par_n}$
instead of $\kk$,
$M$, $S$, $\left(e_s\right)_{s \in S}$ and $\left(f_s\right)_{s \in S}$)
shows that the $\ZZ$-submodule of $\Lambda_n$ spanned by
$\left( p_\lambda \right)_{\lambda \in \Par_n}$ is the
$\ZZ$-submodule of $\Lambda_n$ spanned by
$\left( b_{\lambda, \lambda} m_\lambda \right)_{\lambda \in \Par_n}$.
\end{remark}

% [DG] I added the remark above for no particular reason -- it just looks
% in place here. If you disagree, just remove it.

% [DG][v55] Rewrote the above proof and remark with less informal
% jargon and more structure and precision. (And more length,
% unfortunately.)

% [DG][v56] Rewrote the above proof according to our discussion.
% It has grown longer again, but hopefully it also reads better
% due to less cramping of different things in the same
% sentence.

The purpose of the following exercise is to fill in some details omitted from
the proof of Proposition~\ref{many-symmetric-function-bases-prop}.

\begin{exercise}
\phantomsection\label{exe.Lambda.triangular}
Let $n\in \NN$.

\begin{enumerate}
\item[(a)] Show that every $f \in \Lambda_n$ satisfies
\[
f = \sum_{\mu \in \Par_n}
\left(  \left[  \xx^{\mu}\right]  f\right)  m_{\mu}.
\]
Here, $\left[  \xx^{\mu}\right]  f$ denotes the coefficient of the
monomial $\xx^{\mu}$ in the power series $f$.
\end{enumerate}

Now, we introduce a notation (which generalizes the notation $K_{\lambda,\mu}$
from the proof of Proposition~\ref{many-symmetric-function-bases-prop}): For
any partition $\lambda$ and any weak composition $\mu$, we let $K_{\lambda
,\mu}$ denote the number of all column-strict tableaux $T$ of shape $\lambda$
having $\cont \left(  T\right)  =\mu$.

\begin{enumerate}
\item[(b)] Prove that this number $K_{\lambda,\mu}$ is well-defined (i.e.,
there are only finitely many column-strict tableaux $T$ of shape $\lambda$
having $\cont \left(  T\right)  =\mu$).

\item[(c)] Show that $s_{\lambda}=\sum_{\mu\in\Par_n}
K_{\lambda,\mu}m_{\mu}$ for every $\lambda\in\Par_n$.

\item[(d)] Show that $K_{\lambda,\mu}=0$ for any partitions $\lambda
\in\Par_{n}$ and $\mu\in\Par_{n}$
that don't satisfy $\lambda\triangleright\mu$.

\item[(e)] Show that $K_{\lambda,\lambda}=1$ for any $\lambda\in
\Par_{n}$.
\end{enumerate}

Next, we recall a further notation: For any two partitions $\lambda$ and $\mu
$, we let $a_{\lambda,\mu}$ denote the number of all $\left\{  0,1\right\}
$-matrices of size $\ell\left(  \lambda\right)  \times\ell\left(  \mu\right)
$ having row sums $\lambda$ and column sums $\mu$. (See the proof of
Proposition~\ref{many-symmetric-function-bases-prop} for the concepts of
$\left\{  0,1\right\}  $-matrices and of row sums and column sums.)

\begin{enumerate}
\item[(f)] Prove that this number $a_{\lambda,\mu}$ is well-defined (i.e.,
there are only finitely many $\left\{  0,1\right\}  $-matrices of size
$\ell\left(  \lambda\right)  \times\ell\left(  \mu\right)  $ having row sums
$\lambda$ and column sums $\mu$).

\item[(g)] Show that $e_{\lambda}=\sum_{\mu\in\Par_n}
a_{\lambda,\mu}m_{\mu}$ for every $\lambda\in\Par_{n}$.

\item[(h)] Show that $a_{\lambda,\mu}=0$ for any partitions $\lambda
\in\Par_{n}$ and $\mu\in\Par_{n}$
that don't satisfy $\lambda^{t}\triangleright\mu$.

\item[(i)] Show that $a_{\lambda^{t},\lambda}=1$ for any $\lambda
\in\Par_{n}$.

\end{enumerate}

Next, we introduce a further notation (which generalizes the notation
$b_{\lambda,\mu}$ from the proof of
Proposition~\ref{many-symmetric-function-bases-prop}): For any partition
$\lambda$ and any weak composition $\mu$, we let $b_{\lambda,\mu}$ be the
number of all maps $\varphi:\left\{  1,2,\ldots,\ell\right\}  \rightarrow
\left\{  1,2,3,\ldots\right\}  $ satisfying $\left(  \mu_{j}=
\sum\limits_{\substack{i\in\left\{  1,2,\ldots,\ell\right\}  ;\\\varphi
(i)=j}}\lambda_{i}\text{ for all }j\geq1\right)  $, where $\ell=\ell\left(
\lambda\right)  $.

\begin{enumerate}
\item[(j)] Prove that this number $b_{\lambda,\mu}$ is well-defined (i.e.,
there are only finitely many maps $\varphi:\left\{  1,2,\ldots,\ell\right\}
\rightarrow\left\{  1,2,3,\ldots\right\}  $ satisfying $\left(  \mu_{j}%
=\sum\limits_{\substack{i\in\left\{  1,2,\ldots,\ell\right\}  ;\\\varphi
(i)=j}}\lambda_{i}\text{ for all }j\geq1\right)  $).

\item[(k)] Show that $p_{\lambda}=\sum_{\mu\in\Par_n}
b_{\lambda,\mu}m_{\mu}$ for every $\lambda\in\Par_{n}$.

\item[(l)] Show that $b_{\lambda,\mu}=0$ for any partitions $\lambda
\in\Par_{n}$ and $\mu\in\Par_{n}$
that don't satisfy $\mu\triangleright\lambda$.

\item[(m)] Show that $b_{\lambda,\lambda}$ is a positive integer for any
$\lambda\in\Par_{n}$.

\item[(n)] Show that for any partition $\mu=\left(  \mu_{1},\mu_{2},\ldots
,\mu_{k}\right)  \in\Par_{n}$ with $k=\ell\left(
\mu\right)  $, the integer $b_{\mu,\mu}$ is the size of the subgroup of
$\Symm_k$ consisting of all permutations $\sigma\in\Symm_k$
having each $i$ satisfy $\mu_{\sigma(i)}=\mu_{i}$. (In particular, show that
this subgroup is indeed a subgroup.)

\item[(o)] Show that $b_{\mu,\mu}\mid b_{\lambda,\mu}$ for every $\lambda
\in\Par_{n}$ and $\mu\in\Par_n$.
\end{enumerate}
\end{exercise}

% [DG][v71] Added the above exercise, which fills some DYI
% parts of the proof of Proposition~\ref{many-symmetric-function-bases-prop}.

The bases $\left\{p_\lambda\right\}$ and $\left\{e_\lambda\right\}$
of $\Lambda$ are two examples of \emph{multiplicative bases}\index{multiplicative basis}:
these are bases constructed from a sequence $v_1, v_2, v_3, \ldots$
of symmetric functions by taking all possible finite products.
We will soon encounter another example.
First, let us observe that the finite products of a sequence
$v_1, v_2, v_3, \ldots$ of symmetric functions form a basis of
$\Lambda$ if and only if the sequence is an algebraically independent
generating set of $\Lambda$.
This holds more generally for any commutative algebra, as the following
simple exercise shows:

\begin{exercise}
\label{exe.alggen-via-span}
Let $A$ be a commutative $\kk$-algebra.
Let $v_1, v_2, v_3, \ldots$ be some elements of $A$.

For every partition $\lambda$, define an element $v_{\lambda} \in A$ by
$v_{\lambda} = v_{\lambda_1} v_{\lambda_2} \cdots v_{\lambda_{\ell\left( \lambda\right)  }}$.
Prove the following:

\begin{itemize}
\item[(a)] The $\kk$-subalgebra of $A$ generated by
$v_1, v_2, v_3, \ldots$ is the $\kk$-submodule of $A$ spanned by
the family $\left(  v_{\lambda}\right)  _{\lambda\in\Par}$.

\item[(b)] The elements $v_1, v_2, v_3, \ldots$ generate the
$\kk$-algebra $A$ if and only if the family
$\left(  v_{\lambda}\right) _{\lambda\in\Par}$
spans the $\kk$-module $A$.

\item[(c)] The elements $v_1, v_2, v_3, \ldots$ are algebraically
independent over $\kk$ if and only if the family
$\left(  v_{\lambda}\right)  _{\lambda\in\Par}$
is $\kk$-linearly independent.
\end{itemize}
\end{exercise}

% [DG][v75] Added the above exercise (as it is used tacitly
% several times below) and the paragraph above it.

The next exercise states two well-known identities
for the \emph{generating functions} of the sequences
$\left(e_0, e_1, e_2, \ldots\right)$ and
$\left(h_0, h_1, h_2, \ldots\right)$, which will be used
several times further below:

\begin{exercise}
\label{exe.H(t)-and-E(t)}
In the ring of formal power series
$\powser{\left(\powser{\kk}{\xx}\right)}{t}$, prove the two
identities
\begin{equation}
\prod_{i=1}^{\infty} \left(1 - x_i t\right)^{-1}
= 1 + h_1\left(\xx\right) t + h_2\left(\xx\right) t^2 + \cdots
= \sum_{n \geq 0} h_n\left(\xx\right) t^n
\label{h-generating-function-no-HE}
\end{equation}
and
\begin{equation}
\prod_{i=1}^{\infty} \left(1 + x_i t\right)
= 1 + e_1\left(\xx\right) t + e_2\left(\xx\right) t^2 + \cdots
= \sum_{n \geq 0} e_n\left(\xx\right) t^n .
\label{e-generating-function-no-HE}
\end{equation}
\end{exercise}

% [DG][v75] Added the above exercise.

\subsection{Comultiplications}

Thinking about comultiplication $\Lambda \overset{\Delta}{\rightarrow} \Lambda \otimes \Lambda$
on Schur functions forces us to immediately confront the following.

\begin{definition}
\label{skew-Schur-function-definition}
For partitions $\mu$ and $\lambda$ say that $\mu \subseteq \lambda$ if
$\mu_i \leq \lambda_i$ for $i=1,2,\ldots$.
In other words, two partitions $\mu$ and $\lambda$ satisfy
$\mu \subseteq \lambda$ if and only if the Ferrers diagram for
$\mu$ is a subset of the Ferrers diagram of $\lambda$.
In this case, define the \dfn{skew (Ferrers) diagram} $\lambda/\mu$ to
be their set difference.\footnote{In other words, the skew
Ferrers diagram $\lambda / \mu$ is the set of all
$\left(i, j\right) \in \left\{ 1, 2, 3, \ldots \right\}^2$
satisfying $\mu_i < j \leq \lambda_i$.

While the Ferrers diagram for a
single partition $\lambda$ uniquely determines $\lambda$, the
skew Ferrers diagram $\lambda/\mu$ does not uniquely determine
$\mu$ and $\lambda$. (For instance, it is empty whenever $\lambda
= \mu$.) When one wants to keep $\mu$ and $\lambda$ in memory, one
speaks of the \dfn{skew shape} $\lambda / \mu$; this simply means
the pair $\left(\mu, \lambda\right)$. Every notion defined for
skew Ferrers diagrams also makes sense for skew shapes, because to
any skew shape $\lambda / \mu$ we can assign the skew Ferrers diagram
$\lambda / \mu$ (even if not injectively). For instance, the
\emph{cells}\index{cell of a skew shape}
of the skew shape $\lambda / \mu$ are the cells of the
skew Ferrers diagram $\lambda / \mu$.

One can characterize the skew Ferrers diagrams as
follows: A finite subset $S$ of $\left\{ 1, 2, 3, \ldots \right\}^2$
is a skew Ferrers diagram (i.e., there exist two partitions
$\lambda$ and $\mu$ such that $\mu \subseteq \lambda$ and such
that $S$ is the skew Ferrers diagram $\lambda / \mu$)
if and only if for every $\left(i, j\right) \in S$, every
$\left(i', j'\right) \in \left\{ 1, 2, 3, \ldots \right\}^2$ and
every $\left(i'', j''\right) \in S$ satisfying $i'' \leq i' \leq i$
and $j'' \leq j' \leq j$, we have $\left(i', j'\right) \in S$.}

% [DG][v32] Added the previous footnote because you occasionally speak
% of skew shapes (not just in the context of "tableau of skew shape
% \lambda/\mu").

Then define the \dfn{skew Schur function} $s_{\lambda/\mu}(\xx)$\index{$s_{\lambda/\mu}$}
to be the sum $s_{\lambda/\mu}:=\sum_{T} \xx^{\cont(T)}$, where the sum ranges
over all \emph{column-strict tableaux}\index{column-strict tableau}
$T$ of shape $\lambda/\mu$, that is,
assignments of a value in $\{1,2,3,\ldots\}$ to each cell of $\lambda/\mu$,
weakly increasing left-to-right in rows, and strictly increasing top-to-bottom
in columns.
\end{definition}

\begin{example}
\label{skew-tableau-example}
Let $\lambda = \left(5,3,3,2\right)$ and
$\mu = \left(3,1,1,0\right)$. Then, $\mu \subseteq \lambda$.
The Ferrers diagrams for $\lambda$ and
$\mu$ and the skew Ferrers diagram for $\lambda / \mu$ look as follows:
\[
\underbrace{
\begin{matrix}
\bullet & \bullet & \bullet & \bullet & \bullet \\
\bullet & \bullet & \bullet \\
\bullet & \bullet & \bullet \\
\bullet & \bullet
\end{matrix}
}_{\text{Ferrers diagram of $\lambda$}}
\qquad
\underbrace{
\begin{matrix}
\bullet & \bullet & \bullet \\
\bullet \\
\bullet \\
\vphantom{\bullet}
\end{matrix}
}_{\text{Ferrers diagram of $\mu$}}
\qquad
\underbrace{
\begin{matrix}
\cdot   & \cdot   &   \cdot & \bullet & \bullet \\
\cdot   & \bullet & \bullet \\
\cdot   & \bullet & \bullet \\
\bullet & \bullet
\end{matrix}
}_{\text{skew Ferrers diagram of $\lambda / \mu$}}
\qquad
\]
(where the small dots represent boxes removed from the
diagram).
The filling
\[
T=\begin{matrix}
\cdot & \cdot & \cdot & 2 & 5\\
\cdot & 1     & 1     &   & \\
2     & 2     & 4     &   & \\
4     & 5
\end{matrix}
\]
is a column-strict tableau of shape $\lambda/\mu=(5,3,3,2)/(3,1,0,0)$
and it has $\xx^{\cont(T)}=x_1^2 x_2^3 x_3^0 x_4^2 x_5^2$.

On the other hand, if we took $\lambda = \left(5, 3, 1\right)$
and $\mu = \left(1, 1, 1, 1\right)$, then we wouldn't
have $\mu \subseteq \lambda$, since
$\mu_4 = 1 > 0 = \lambda_4$.
\end{example}

% [DG][v78] Expanded the above example to show not just the
% tableau, but also the diagram.

\begin{remark}
If $\mu$ and $\lambda$ are partitions such that $\mu \subseteq \lambda$,
then $s_{\lambda / \mu} \in \Lambda$. (This is proven similarly as
Proposition~\ref{Schur-functions-are-symmetric-prop}.) Actually,
if $\mu \subseteq \lambda$, then
$s_{\lambda / \mu} \in \Lambda_{\left|\lambda / \mu\right|}$,
where $\left|\lambda / \mu\right|$ denotes the number of cells of
the skew shape $\lambda / \mu$ (so
$\left|\lambda / \mu\right| = \left|\lambda\right| - \left|\mu\right|$).

It is customary to define $s_{\lambda / \mu}$ to be $0$ if we don't
have $\mu \subseteq \lambda$. This can also be seen by a literal
reading of the definition $s_{\lambda/\mu}:=\sum_{T} \xx^{\cont(T)}$,
as long as we understand that there are no column-strict tableaux of
shape $\lambda / \mu$ when $\lambda / \mu$ is not defined.

Clearly, every partition $\lambda$ satisfies $s_\lambda =
s_{\lambda / \varnothing}$.

It is easy to see that two partitions $\lambda$ and $\mu$ satisfy
$\mu \subseteq \lambda$ if and only if they satisfy
$\mu^t \subseteq \lambda^t$.
\end{remark}

% [DG][v25] Added preceding two paragraphs to this remark.

% [DG][v78] Added last sentence to the above remark.

\begin{exercise} \phantomsection
\label{exe.schur.rotate180}

\begin{itemize}

\item[(a)]
State and prove an analogue of
Proposition~\ref{prop.other-linear-orders.slambda} for skew Schur
functions.

\item[(b)]
Let $\lambda$, $\mu$, $\lambda'$ and $\mu'$ be partitions such that
$\mu \subseteq \lambda$ and $\mu' \subseteq \lambda'$. Assume that the
skew Ferrers diagram $\lambda' / \mu'$ can be obtained from the skew
Ferrers diagram $\lambda / \mu$ by a $180^{\circ}$ rotation.\footnote{For
example, this happens when $\lambda = \left(3,2\right)$,
$\mu = \left(1\right)$, $\lambda' = \left(5,4\right)$ and
$\mu' = \left(3,1\right)$.} Prove that
$s_{\lambda / \mu} = s_{\lambda' / \mu'}$.

\end{itemize}
\end{exercise}

% [DG][v14] Inserted the exercise above.

\begin{exercise}
\label{exe.schur.disconnected}
Let $\lambda$ and $\mu$ be two partitions, and let $k \in \NN$
be such that\footnote{As usual, we write $\nu_k$ for the
$k$-th entry of a partition $\nu$.}
$\mu_k \geq \lambda_{k+1}$. Let $F$ be the skew Ferrers
diagram $\lambda / \mu$. Let
$F_{\operatorname{rows}\leq k}$ denote the subset of $F$
consisting of all $\left(i, j\right) \in F$ satisfying
$i \leq k$. Let
$F_{\operatorname{rows} > k}$ denote the subset of $F$
consisting of all $\left(i, j\right) \in F$ satisfying
$i > k$.
Let $\alpha$ and $\beta$ be
two partitions such that $\beta \subseteq \alpha$ and such
that the skew Ferrers diagram
$\alpha / \beta$ can be obtained from
$F_{\operatorname{rows}\leq k}$
by parallel translation.
Let $\gamma$ and $\delta$ be two partitions such that
$\delta \subseteq \gamma$ and such that the skew
Ferrers diagram $\gamma / \delta$ can be obtained from
$F_{\operatorname{rows} > k}$ by parallel
translation.\footnote{Here is an example of the situation:
$\lambda = \left(6, 5, 5, 2, 2\right)$,
$\mu = \left(4, 4, 3, 1\right)$, $k = 3$ (satisfying
$\mu_k = \mu_3 = 3 \geq 2 = \lambda_4 = \lambda_{k+1}$),
$\alpha = \left(3, 2, 2\right)$,
$\beta = \left(1, 1\right)$,
$\gamma = \left(2, 2\right)$, and
$\delta = \left(1\right)$.} Prove that
$s_{\lambda / \mu}
= s_{\alpha / \beta} s_{\gamma / \delta}$.
\end{exercise}

% [DG][v32] Added the exercise above.

\begin{proposition}
\label{symm-comultiplication-formulas-prop}
The comultiplication $\Lambda \overset{\Delta}{\rightarrow} \Lambda \otimes \Lambda$
has the following effect on the symmetric functions discussed so far\footnote{
The abbreviated summation indexing $\sum_{i+j=n} t_{i,j}$ used here is intended to mean
\[
\sum_{ \substack{ (i,j) \in \NN^2; \\ i+j=n} } t_{i,j}.
\]}:
\begin{enumerate}
\item[ (i) ] $\Delta p_n = \one \otimes p_n + p_n \otimes \one$ for every $n \geq 1$,
that is, the power sums $p_n$ are primitive.
\item[ (ii) ] $\Delta e_n = \sum_{i+j=n} e_i \otimes e_j$ for every $n \in \NN$.
\item[ (iii) ] $\Delta h_n = \sum_{i+j=n} h_i \otimes h_j$ for every $n \in \NN$.
\item[ (iv) ] $\Delta s_{\lambda} = \sum_{\mu \subseteq \lambda} s_{\mu} \otimes s_{\lambda/\mu}$
for any partition $\lambda$.
\vspace{0.2pc}
\item[ (v) ] $\Delta s_{\lambda / \nu} = \sum\limits_{\substack{\mu \in \Par: \\ \nu \subseteq \mu \subseteq \lambda}} s_{\mu/\nu} \otimes s_{\lambda/\mu}$ for any partitions $\lambda$ and $\nu$.
\end{enumerate}

% [DG][v25] Added part (v).

\end{proposition}
\begin{proof}
Recall that $\Delta$ sends $f(\xx) \mapsto f(\xx,\yy)$, and one can easily check that
\begin{enumerate}
\item[ (i) ] $p_n(\xx,\yy) = \sum_i x_i^n + \sum_i y_i^n
= p_n(\xx) \cdot 1 + 1 \cdot p_n(\yy)$ for every $n \geq 1$;
\item[ (ii) ] $e_n(\xx,\yy) = \sum_{i+j=n} e_i(\xx)e_j(\yy)$ for every $n \in \NN$;
\item[ (iii) ] $h_n(\xx,\yy) = \sum_{i+j=n} h_i(\xx)h_j(\yy)$ for every $n \in \NN$.
\end{enumerate}
For assertion (iv), note that by \eqref{different-total-order-Schur-function}, one has
\begin{equation}
s_\lambda(\xx,\yy)=\sum_T (\xx,\yy)^{\cont(T)},
\label{pf.symm-comultiplication-formulas-prop.1}
\end{equation}
where the sum is over column-strict tableaux $T$ of shape $\lambda$ having entries
in the linearly ordered alphabet
%\footnote{
%There is a subtle point here:
%the linear ordering \eqref{two-variable-set-ordering} is not equivalent under
%permutations in $\Symm_{(infty)}$ to some other orderings, such as,
%$x_1 < y_1 < x_2 < y_2 < \cdots$.   Given $\lambda$ in $\Par_n$,
%it would be better to first assume one has
%\emph{finite} variable sets $\xx=(x_1,\ldots,x_n),
%\yy=(y_1,\ldots,y_n)$ in proving the identity
%$s_\lambda(\xx,\yy)=\sum_{\mu \subseteq \lambda} s_\mu(\xx) s_{\lambda/\mu}(\yy%)$
%as above;  the result then follows for infinite variable sets by
%the usual considerations.}
\begin{equation}
\label{two-variable-set-ordering}
x_1 < x_2 < \cdots < y_1 < y_2 < \cdots .
\end{equation}
\footnote{Here, $(\xx,\yy)^{\cont(T)}$ means the monomial
$\prod_{a \in \mathfrak{A}} a^{\left|T^{-1}\left(a\right)\right|}$,
where $\mathfrak{A}$ denotes the totally ordered alphabet
$x_1 < x_2 < \cdots < y_1 < y_2 < \cdots $. In other words,
$(\xx,\yy)^{\cont(T)}$ is the product of all entries of the tableau
$T$ (which is a monomial, since the entries of $T$ are not numbers
but variables).

The following rather formal argument should allay any doubts as to
why \eqref{pf.symm-comultiplication-formulas-prop.1} holds:
Let $\LLL$ denote the totally ordered set which is given by the
set $\left\{1,2,3,\ldots\right\}$ of positive integers, equipped
with the total order
$1 <_\LLL 3 <_\LLL 5 <_\LLL 7 <_\LLL \cdots <_\LLL 2 <_\LLL 4 <_\LLL 6 <_\LLL 8 <_\LLL \cdots$.
Then, \eqref{different-total-order-Schur-function} yields
$s_\lambda = \sum_T \xx^{\cont(T)}$
as $T$ runs through all $\LLL$-column-strict tableaux of shape $\lambda$.
Substituting the variables $x_1, y_1, x_2, y_2, x_3, y_3, \ldots$
for $x_1, x_2, x_3, x_4, x_5, x_6, \ldots$ (that is, substituting
$x_i$ for $x_{2i-1}$ and $y_i$ for $x_{2i}$) in this equality, we
obtain \eqref{pf.symm-comultiplication-formulas-prop.1}.}
For example,
\[
T=\begin{matrix}
x_1   & x_1 & x_1 & y_2 & y_5\\
x_2   & y_1 & y_1 &   & \\
y_2   & y_2 & y_4 &   & \\
y_4   & y_5
\end{matrix}
\]
is such a tableau of shape $\lambda=(5,3,3,2)$.  Note that the restriction of
$T$ to the alphabet $\xx$ gives a column-strict tableau $T_\xx$ of some shape
$\mu \subseteq \lambda$, and the restriction of $T$ to the alphabet $\yy$ gives
a column-strict tableau $T_\yy$ of shape $\lambda/\mu$
(e.g. for $T$ in the example above,
the tableau $T_\yy$ appeared in Example~\ref{skew-tableau-example}).
Consequently, one has
\begin{align}
s_{\lambda}(\xx,\yy) =& \sum_T \xx^{\cont(T_\xx)} \cdot \yy^{\cont(T_\yy)}
\nonumber \\
       =& \sum_{\mu \subseteq \lambda} \left( \sum_{T_\xx} \xx^{\cont(T_\xx)} \right)
                                      \left( \sum_{T_\yy} \yy^{\cont(T_\yy)} \right)
       = \sum_{\mu \subseteq \lambda} s_{\mu}(\xx)  s_{\lambda/\mu}(\yy).
\label{eq.Lambda.s.xy}
\end{align}

Assertion (v) is obvious in the case when we don't have
$\nu \subseteq \lambda$ (in fact, in this case, both
$s_{\lambda / \nu}$ and
$\sum\limits_{\substack{\mu \in \Par: \\ \nu \subseteq \mu \subseteq \lambda}}
s_{\mu/\nu} \otimes s_{\lambda/\mu}$
are clearly zero). In the remaining case, the proof of assertion (v)
is similar to that of (iv). (Of course,
the tableaux $T$ and $T_{\xx}$ now have skew shapes $\lambda/\nu$
and $\mu/\nu$, and instead of
\eqref{different-total-order-Schur-function},
we need to use the answer to Exercise~\ref{exe.schur.rotate180}(a).)
\end{proof}

% [DG][v25] Added footnote and proof of (v).

Notice that parts (ii) and (iii) of
Proposition~\ref{symm-comultiplication-formulas-prop} are particular
cases of part (iv), since $h_n = s_{(n)}$ and $e_n = s_{(1^n)}$.

% [DG][v46] Added previous sentence.

\begin{exercise} \phantomsection
\label{exe.Lambda.cocommutative}
\begin{itemize}
\item[(a)] Show that the Hopf algebra $\Lambda$ is cocommutative.
\item[(b)] Show that $\Delta s_{\lambda / \nu} = \sum\limits_{\substack{\mu \in \Par: \\ \nu \subseteq \mu \subseteq \lambda}} s_{\lambda/\mu} \otimes s_{\mu/\nu}$ for any partitions $\lambda$ and $\nu$.
\end{itemize}
\end{exercise}

% [DG][v25] Added the preceding exercise.

\begin{exercise}
\label{exe.schur.finite}
Let $n \in \NN$. Consider the finite variable set $\left(x_1, x_2, \ldots,
x_n\right)$ as a subset of $\xx = \left(x_1, x_2, x_3, \ldots\right)$.
Recall that $f \left( x_1, x_2, \ldots , x_n \right)$ is a
well-defined element of $\kk\left[x_1, x_2, \ldots, x_n\right]$
for every $f \in R\left(\xx\right)$ (and therefore also for
every $f \in \Lambda$, since $\Lambda \subset R\left(\xx\right)$),
according to Exercise~\ref{exe.Lambda.subs.fin}.
\begin{itemize}
\item[(a)] Show that any two partitions $\lambda$ and $\mu$ satisfy
\begin{vershort}
\[
s_{\lambda / \mu} \left(x_1, x_2, \ldots, x_n\right)
= \sum\limits_{\substack{T\text{ is a column-strict}\\ \text{tableau of shape } \lambda / \mu \text{;} \\ \text{all entries of } T \text{ belong} \\ \text{to } \left\{1,2,\ldots,n\right\}}}
\xx^{\cont\left(T\right)}.
\]
\end{vershort}
\begin{verlong}
\[
s_{\lambda / \mu} \left(x_1, x_2, \ldots, x_n\right)
= \sum\limits_{\substack{T\text{ is a column-strict}\\ \text{tableau of shape } \lambda / \mu \text{;} \\ \text{all entries of } T \text{ belong} \\ \text{to } \left\{1,2,\ldots,n\right\}}}
\left(x_1,x_2,\ldots,x_n\right)^{\cont\left(T\right)}.
\]
Here, $\left(x_1,x_2,\ldots,x_n\right)^{\cont\left(T\right)}$
is defined as
$\prod_{i=1}^{n} x_i ^{\left\vert T^{-1}\left( i\right) \right\vert }$
for any column-strict
tableau $T$ having the property that all entries of $T$ belong to $\left\{
1,2,\ldots,n\right\}  $.
\end{verlong}
\item[(b)] If $\lambda$ is a partition having more than $n$ parts%
\footnote{Recall that the \emph{parts} of a partition
are its nonzero entries.}, then
show that $s_\lambda \left(x_1, x_2, \ldots, x_n\right) = 0$.
\end{itemize}
\end{exercise}

% [DG][v25] Added the exercise above, as preparation for the bialternant
% section. More interesting parts could be added (characterizing the
% ideal of symmetric functions that vanish upon specializing to $n$
% variables in several way); do you think that would be good?

% [DG][v26] Replaced $\left(x_1,x_2,\ldots,x_n\right)^{\cont\left(T\right)}$
% by $\xx^{\cont\left(T\right)}$.

\begin{remark}
\label{rmk.many-symmetric-function-bases-prop.finite}
An analogue of Proposition~\ref{many-symmetric-function-bases-prop}
holds for symmetric polynomials in finitely many variables: Let
$N \in \NN$. Then, we have
\begin{itemize}
\item[(a)] The family
$\left\{ m_\lambda \left(x_1, x_2, \ldots, x_N\right) \right\}$,
as $\lambda$ runs through all partitions having length $\leq N$,
is a graded basis of the graded $\kk$-module
$\Lambda \left(x_1, x_2, \ldots, x_N\right)
= \kk \left[x_1, x_2, \ldots, x_N\right]^{\Symm_N}$.
\item[(b)] For any
partition $\lambda$ having length $> N$, we have
$m_\lambda \left(x_1, x_2, \ldots, x_N\right) = 0$.
\item[(c)] The family
$\left\{ e_\lambda \left(x_1, x_2, \ldots, x_N\right) \right\}$,
as $\lambda$ runs through all partitions whose parts are all
$\leq N$, is a graded basis of the graded $\kk$-module
$\Lambda \left(x_1, x_2, \ldots, x_N\right)$.
\item[(d)] The family
$\left\{ s_\lambda \left(x_1, x_2, \ldots, x_N\right) \right\}$,
as $\lambda$ runs through all partitions having length
$\leq N$, is a graded basis of the graded $\kk$-module
$\Lambda \left(x_1, x_2, \ldots, x_N\right)$.
\item[(e)] If $\QQ$ is a subring of $\kk$, then the family
$\left\{ p_\lambda \left(x_1, x_2, \ldots, x_N\right) \right\}$,
as $\lambda$ runs through all partitions having length
$\leq N$, is a graded basis of the graded $\kk$-module
$\Lambda \left(x_1, x_2, \ldots, x_N\right)$.
\item[(f)] If $\QQ$ is a subring of $\kk$, then the family
$\left\{ p_\lambda \left(x_1, x_2, \ldots, x_N\right) \right\}$,
as $\lambda$ runs through all partitions whose parts are all
$\leq N$, is a graded basis of the graded $\kk$-module
$\Lambda \left(x_1, x_2, \ldots, x_N\right)$.
\end{itemize}
Indeed, the claims (a) and (b) are obvious, while the
claims (c), (d) and (e) are proven similarly to our proof of
Proposition~\ref{many-symmetric-function-bases-prop}. We leave
the proof of (f) to the reader; this proof can also be found in
\cite[Theorem 10.86]{Loehr-bij}\footnote{See
\cite[Remark 10.76]{Loehr-bij} for why \cite[Theorem 10.86]{Loehr-bij}
is equivalent to our claim (f).}.

Claim (c) can be rewritten as follows: The elementary symmetric
polynomials $e_i \left(x_1, x_2, \ldots, x_N\right)$, for
$i \in \left\{1, 2, \ldots, N\right\}$, form an algebraically
independent generating set of
$\Lambda \left(x_1, x_2, \ldots, x_N\right)$. This is precisely
the well-known theorem (due to Gauss)\footnote{See, e.g.,
\cite[\textit{Symmetric Polynomials}, Theorem 5 and Remark 17]{Conrad-blurbs}
or \cite[\S 5.3]{vanderWaerden-1}
or \cite[Theorem 1]{BlumSmithCoskey}. In a slightly different
form, it also appears in \cite[Theorem (5.10)]{LaksovLascouxPragaczThorup}.}
that every symmetric polynomial in $N$
variables $x_1, x_2, \ldots, x_N$ can be written uniquely as a
polynomial in the $N$ elementary symmetric polynomials.
\end{remark}

% [DG][v39] Added remark above (it is strange to have a text on
% symmetric functions not mentioning the place every undergrad
% encounters them in). I have put it here since the meaning of
% evaluating at finitely many indeterminates is explained in the
% preceding exercise.

% [DG][v55] Added two more references for Gauss's theorem.

\subsection{The antipode, the involution $\omega$, and algebra generators}
\label{symm-antipode-section}

Since $\Lambda$ is a connected graded $\kk$-bialgebra, it will have an
antipode $\Lambda \overset{S}{\rightarrow} \Lambda$ making it a Hopf algebra
by Proposition~\ref{graded-connected-bialgebras-have-antipodes}.
However, we can identify $S$ more explicitly now.

% [DG][v62] "several issues will be resolved by identifying $S$ more explicitly now"
% --> "we can identify $S$ more explicitly now".

\begin{proposition}
\label{symm-antipode-and-h-basis}
Each of the families $\{e_n\}_{n=1,2,\ldots}$ and $\{h_n\}_{n=1,2,\ldots}$
are algebraically independent, and generate $\Lambda_\kk$ as a polynomial
algebra for any commutative ring $\kk$.  The same holds
for $\{p_n\}_{n=1,2,\ldots}$ when $\QQ$ is a subring of $\kk$.

Furthermore, the antipode $S$ acts as follows:
\begin{enumerate}
\item[ (i) ] $S(p_n) = -p_n$ for every positive integer $n$.
\item[ (ii) ] $S(e_n) = (-1)^n h_n$ for every $n \in \NN$.
\item[ (iii) ] $S(h_n) = (-1)^n e_n$ for every $n \in \NN$.
\end{enumerate}
\end{proposition}

\begin{proof}
The assertion that $\{e_n\}_{n \geq 1}$ are algebraically independent and generate
$\Lambda$ is equivalent to Proposition~\ref{many-symmetric-function-bases-prop}
asserting that $\{e_\lambda\}_{\lambda \in \Par}$ is a basis for $\Lambda$.
(Indeed, this equivalence follows from parts (b) and (c) of
Exercise~\ref{exe.alggen-via-span}, applied to $v_n = e_n$
and $v_\lambda = e_\lambda$.)
Thus, the former assertion is true.
If $\QQ$ is a subring of $\kk$, then a similar argument (using $p_n$ and
$p_\lambda$ instead of $e_n$ and $e_\lambda$)
shows that $\{p_n\}_{n \geq 1}$ are algebraically independent and generate
$\Lambda$.

The assertion $S(p_n) = -p_n$ follows from Proposition~\ref{antipode-on-primitives}
since $p_n$ is primitive by Proposition~\ref{symm-comultiplication-formulas-prop}(i).

For the remaining assertions, start with the easy generating function
identities\footnote{See the solution to Exercise~\ref{exe.H(t)-and-E(t)}
for the proofs of the identities.}
\begin{align}
\label{h-generating-function}
H(t)&:=\prod_{i=1}^\infty (1-x_it)^{-1}
     = 1+h_1(\xx)t+h_2(\xx)t^2+\cdots = \sum_{n \geq 0} h_n(\xx) t^n ; \\
\label{e-generating-function}
E(t)&:=\prod_{i=1}^\infty (1+x_it)
     = 1+e_1(\xx)t+e_2(\xx)t^2+\cdots = \sum_{n \geq 0} e_n(\xx) t^n .
\end{align}
These show that
\begin{equation}
\label{E(t)-H(t)-relation}
1=E(-t)H(t)=\left( \sum_{n \geq 0} e_n(\xx) (-t)^n \right)
            \left( \sum_{n \geq 0} h_n(\xx) t^n \right) .
\end{equation}
Hence, equating coefficients of powers of $t$, we see that for $n=0,1,2,\ldots$
we have
\begin{equation}
\label{e-h-relation}
\sum_{i+j=n} (-1)^i e_i h_j = \delta_{0,n}.
\end{equation}
This lets us recursively express the $e_n$ in terms of $h_n$ and vice-versa:
\begin{align}
e_0&=1=h_0 ; \label{e-h-recursions.1}\\
e_n&= e_{n-1} h_1 - e_{n-2} h_2 + e_{n-3} h_3 - \cdots ; \label{e-h-recursions.2}\\
h_n&= h_{n-1} e_1 - h_{n-2} e_2 + h_{n-3} e_3 - \cdots \label{e-h-recursions.3}
\end{align}
for $n =1,2,3,\ldots$
Now, let us use the algebraic independence of the
generators $\{e_n\}$ for $\Lambda$
to define a $\kk$-algebra endomorphism
\[
\begin{array}{rcl}
\Lambda &\overset{\omega}{\rightarrow} &\Lambda ,\\
e_n&\longmapsto&h_n \qquad \text{ (for positive integers $n$)}.
\end{array}
\]
Then,
\begin{align}
\omega\left(e_n\right) = h_n \qquad \text{ for each }
n \geq 0
\label{the-fundamental-involution-eh}
\end{align}
(indeed, this holds for $n > 0$ by definition, and for $n = 0$
because $\omega\left(e_0\right) = \omega\left(1\right) = 1 = h_0$).
Hence, the identical form of the two recursions \eqref{e-h-recursions.2}
and \eqref{e-h-recursions.3}
shows that
\begin{align}
\omega\left(h_n\right) = e_n \qquad \text{ for each }
n \geq 0
\label{the-fundamental-involution-he}
\end{align}
\footnote{%
Here is this argument in more detail:
We must show that $\omega\left(h_n\right) = e_n$ for each $n \geq 0$.
We shall prove this by strong induction on $n$.
Thus, we fix an $n \geq 0$, and assume as induction hypothesis that
$\omega\left(h_m\right) = e_m$ for each $m < n$.
We must then prove that $\omega\left(h_n\right) = e_n$.
If $n = 0$, then this is obvious; thus, assume WLOG that $n > 0$.
Hence,
\begin{align*}
\omega\left(h_n\right)
&= \omega\left( h_{n-1} e_1 - h_{n-2} e_2 + h_{n-3} e_3 - \cdots \right)
\qquad \left(\text{by \eqref{e-h-recursions.3}}\right) \\
&= \omega\left(h_{n-1}\right) \omega\left(e_1\right) - \omega\left(h_{n-2}\right) \omega\left(e_2\right) + \omega\left(h_{n-3}\right) \omega\left(e_3\right) - \cdots
\qquad \left(\text{since $\omega$ is a $\kk$-algebra homomorphism}\right) \\
&= e_{n-1} \omega\left(e_1\right) - e_{n-2} \omega\left(e_2\right) + e_{n-3} \omega\left(e_3\right) - \cdots
\qquad
\left(\text{since $\omega\left(h_m\right) = e_m$ for each $m < n$}\right) \\
&= e_{n-1} h_1 - e_{n-2} h_2 + e_{n-3} h_3 - \cdots
\qquad
\left(\text{since \eqref{the-fundamental-involution-eh} shows that
$\omega\left(e_m\right) = h_m$ for each $m \geq 0$}\right) \\
&= e_n \qquad \left(\text{by \eqref{e-h-recursions.2}}\right) ,
\end{align*}
as desired. This completes the induction step.}.
Combining this with \eqref{the-fundamental-involution-eh}, we
conclude that $\left(\omega\circ\omega\right) \left(e_n\right)
= e_n$ for each $n \geq 0$.
Therefore, the two $\kk$-algebra homomorphisms
$\omega \circ \omega : \Lambda \to \Lambda$ and $\id : \Lambda \to \Lambda$
agree on each element of the generating set $\{e_n\}$ of $\Lambda$.
Hence, they are equal, i.e., we have $\omega \circ \omega = \id$.
Therefore $\omega$ is an involution and therefore a $\kk$-algebra
automorphism of $\Lambda$.
This, in turn, yields that
the $\{h_n\}$ (being the images of the $\{e_n\}$ under this
automorphism) are another algebraically independent generating set
for $\Lambda$.

% [DG][v80] Added details to the above paragraph.

For the assertion about the antipode $S$ applied to $e_n$ or $h_n$,
note that the coproduct formulas for $e_n, h_n$ in
Proposition~\ref{symm-comultiplication-formulas-prop}(ii),(iii)
show that the defining relations for their
antipodes \eqref{antipode-defining-relation}
will in this case be
\begin{align*}
\sum_{i+j=n} S(e_i) e_j &= \delta_{0,n} = \sum_{i+j=n} e_i S(e_j) , \\
\sum_{i+j=n} S(h_i) h_j &= \delta_{0,n} = \sum_{i+j=n} h_i S(h_j)
\end{align*}
because $u \epsilon(e_n)=u \epsilon(h_n)=\delta_{0,n}$.
Comparing these to \eqref{e-h-relation}, one concludes via induction on $n$
that $S(e_n)=(-1)^n h_n$ and $S(h_n)=(-1)^n e_n$.
\end{proof}

The $\kk$-algebra endomorphism $\omega$ of $\Lambda$ defined in the
proof of Proposition~\ref{symm-antipode-and-h-basis} is
sufficiently important that we record its definition and a selection
of fundamental properties:

\begin{definition} \label{def.Lambda.omega}
Let \dfn{$\omega$} be the $\kk$-algebra homomorphism
\begin{equation}
\label{the-fundamental-involution-definition}
\begin{array}{rcl}
\Lambda &\to &\Lambda ,\\
e_n&\longmapsto&h_n \qquad \text{ (for positive integers $n$)}.
\end{array}
\end{equation}
This homomorphism $\omega$ is known as the
\dfn{fundamental involution} of $\Lambda$.
\end{definition}

\begin{proposition} \label{prop.Lambda.omega-and-S}
Consider the fundamental involution $\omega$ and the
antipode $S$ of the Hopf algebra $\Lambda$.

\begin{itemize}
\item[(a)] We have
\[
\omega \left( e_n \right) = h_n
\qquad \text{ for each $n \in \ZZ$.}
\]

\item[(b)] We have
\[
\omega \left( h_n \right) = e_n
\qquad \text{ for each $n \in \ZZ$.}
\]

\item[(c)] We have
\[
\omega \left( p_n \right) = \left(-1\right)^{n-1} p_n
\qquad \text{ for each positive integer $n$.}
\]

\item[(d)] The map $\omega$ is a $\kk$-algebra automorphism
of $\Lambda$ and an involution.

\item[(e)] If $n \in \NN$, then
\begin{equation}
\label{antipode-versus-omega-in-Sym}
S \left(f\right) = \left(-1\right)^n \omega \left(f\right)
\qquad \text{ for all } f \in \Lambda_n .
\end{equation}

\item[(f)] The map $\omega$ is a Hopf algebra automorphism
of $\Lambda$.

\item[(g)] The map $S$ is a Hopf algebra automorphism
of $\Lambda$.

\item[(h)] Every partition $\lambda$ satisfies the three
equalities
\begin{align}
\omega\left(  h_{\lambda}\right)   &  =e_{\lambda};
\label{eq.exe.Lambda.omegahlam.h}
\\
\omega\left(  e_{\lambda}\right)   &  =h_{\lambda};
\label{eq.exe.Lambda.omegahlam.e}
\\
\omega\left(  p_{\lambda}\right)   &  =\left(  -1\right)  ^{\left\vert
\lambda\right\vert -\ell\left(  \lambda\right)  }p_{\lambda}.
\label{eq.exe.Lambda.omegahlam.p}
\end{align}

\item[(i)] The map $\omega$ is an isomorphism of graded
$\kk$-modules.

\item[(j)] The family $\left(h_\lambda\right)_{\lambda \in \Par}$
is a graded basis of the graded $\kk$-module $\Lambda$.

\end{itemize}
\end{proposition}

\begin{exercise}
\label{exe.Lambda.omega-and-S}
Prove Proposition~\ref{prop.Lambda.omega-and-S}.

[\textbf{Hint:} Parts (a), (b) and (d) have been shown in
the proof of Proposition~\ref{symm-antipode-and-h-basis} above.
For part (e), let $D_{-1} : \Lambda \to \Lambda$ be the
$\kk$-algebra morphism sending each homogeneous $f \in \Lambda_n$
to $\left(-1\right)^n f$; then
argue that $\omega \circ D_{-1}$ and $S$ are two $\kk$-algebra
morphisms that agree on all elements
of the generating set $\{e_n\}$.
Derive part (c) from (d) and
Proposition~\ref{symm-antipode-and-h-basis}.
Part (h) then follows by multiplicativity.
For parts (f) and (g), check the coalgebra homomorphism
axioms on the $e_n$.
Parts (i) and (j) are easy consequences.]
\end{exercise}

% [DG][v80] Added the above proposition and exercise. All
% of it was already been said in some way, but now it is
% explicit and referenceable.

Proposition~\ref{prop.Lambda.omega-and-S}(e) shows that the
antipode $S$ on $\Lambda$ is, up to sign,
the same as the fundamental involution $\omega$.
Thus, studying $\omega$ is essentially equivalent to
studying $S$.

\begin{remark}
Up to now we have not yet derived how the involution $\omega$ and the
antipode $S$ act on (skew) Schur functions, which is quite beautiful:
If $\lambda$ and $\mu$ are partitions satisfying $\mu \subseteq \lambda$,
then
\begin{equation}
\label{skew-Schur-antipodes}
\begin{aligned}
\omega(s_{\lambda/\mu}) &= s_{\lambda^t/\mu^t} , \\
S(s_{\lambda/\mu})      &= (-1)^{|\lambda/\mu|} s_{\lambda^t/\mu^t} \\
\end{aligned}
\end{equation}
where recall that $\lambda^t$ is the transpose or conjugate partition
to $\lambda$, and $|\lambda/\mu|$ is the number of squares in
the skew diagram $\lambda/\mu$, that is, $|\lambda/\mu|=n-k$
if $\lambda, \mu$ lie in $\Par_n, \Par_k$ respectively.

We will deduce this later in three ways (once as an exercise using
the Pieri rules in Exercise~\ref{exe.pieri.omega},
once again using
skewing operators in Exercise~\ref{Darij's-skewing-exercise},
and for the third time from the action of
the antipode in $\Qsym$ on $P$-partition enumerators
in Corollary~\ref{antipode-on-skew-Schur-corollary}).
However, one could also deduce it immediately from our knowledge of the action of
$\omega$ and $S$ on $e_n, h_n$, if we were to prove the following
famous \emph{Jacobi-Trudi}\index{Jacobi-Trudi formula}
and \emph{dual Jacobi-Trudi}\index{dual Jacobi-Trudi formula}
formulas\footnote{The second of the following identities is also
known as the \dfn{von N\"agelsbach-Kostka identity}.}:

\begin{theorem}
\label{Jacobi-Trudi-formulae}
Skew Schur functions are the following polynomials in $\{h_n\}, \{e_n\}$:
\begin{align}
s_{\lambda/\mu} &=
\det( h_{\lambda_i - \mu_j - i + j} )_{i,j=1,2,\ldots,\ell} ,
\label{eq.jacobi-trudi.h} \\
s_{\lambda^t/\mu^t} &=
\det( e_{\lambda_i - \mu_j - i + j} )_{i,j=1,2,\ldots,\ell}
\label{eq.jacobi-trudi.e}
\end{align}
for any two partitions $\lambda$ and $\mu$ and any
$\ell \in \NN$ satisfying $\ell\left(\lambda\right) \leq \ell$ and
$\ell\left(\mu\right) \leq \ell$.
\end{theorem}

% [DG][v25] Replaced "$\lambda$ has at most $\ell$ nonzero parts"
% by "$\ell\left(\lambda\right) \leq \ell$ and
% $\ell\left(\mu\right) \leq \ell$" to avoid making a subtly false
% statement if $\mu$ is not contained in $\lambda$.

\noindent
Since we appear not to need these formulas in the sequel, we will
not prove them right away.
\begin{vershort}
However, a proof is sketched in the
solution to Exercise~\ref{exe.jacobi-trudi.macdonald},
\end{vershort}
\begin{verlong}
However, two proofs are sketched in the
solutions to Exercise~\ref{exe.jacobi-trudi} and to
Exercise~\ref{exe.jacobi-trudi.macdonald},
\end{verlong}
and various proofs are well-explained in
\cite[(39) and (41)]{Leeuwen-altSchur}, \cite[\S I.5]{Macdonald},
\cite[Thm. 7.1]{RosasSagan},
\cite[\S 4.5]{Sagan}, \cite[\S 7.16]{Stanley},
\cite[Thms. 3.5 and 3.5$^\ast$]{Wachs}; also, a
simultaneous generalization of both formulas is shown in
\cite[Theorem 11]{GesselViennot}, and three others in
\cite[1.9]{Remmel}, \cite[Thm. 3.1]{HamelGoulden}
and \cite{Jin-outsidenested}.
An elegant treatment of Schur polynomials taking the
Jacobi-Trudi formula \eqref{eq.jacobi-trudi.h}
as the \emph{definition} of $s_\lambda$
is given by Tamvakis \cite{Tamvakis}.
\end{remark}

% [DG][v14] I added the Gessel-Viennot reference -- while I didn't read it,
% it looks very interesting. The generalization there replaces the
% "weakly increasing" and "strictly increasing" conditions of semistandard
% Young tableaux by two relations satisfying fairly liberal constraints.
% Also, I added the van Leeuwen and Rosas-Sagan references; each of them
% proves more than just Jacobi-Trudi.

% [DG][v19] Added Wachs reference. She uses divided-difference operators
% to prove it. This might be a good thing to look into when writing the
% chapter about 0-Hecke algebras!

\subsection{Cauchy product, Hall inner product, self-duality}

The Schur functions, although a bit unmotivated right now, have
special properties with regard to the Hopf structure.  One property is
intimately connected with the following \dfn{Cauchy identity}.

\begin{theorem}
\label{thm.cauchy-identity}
In the power series ring $\kk\left[\left[\xx, \yy\right]\right] :=
\kk \left[\left[x_1, x_2, \ldots, y_1, y_2, \ldots\right]\right]$, one has
the following expansion:
\begin{equation}
\label{Cauchy-identity}
\prod_{i,j=1}^\infty (1-x_i y_j)^{-1}
=\sum_{\lambda \in \Par} s_{\lambda}(\xx) s_\lambda(\yy).
\end{equation}
\end{theorem}

\begin{remark}
The left hand side of \eqref{Cauchy-identity} is known as the
\dfn{Cauchy product}, or \dfn{Cauchy kernel}.

An equivalent version of the equality \eqref{Cauchy-identity}
is obtained by replacing each $x_i$ by $x_i t$, and writing the
resulting identity in the power series ring $R(\xx, \yy)[[t]]$:
\begin{equation}
\label{Cauchy-identity-with-t}
\prod_{i,j=1}^\infty (1-tx_i y_j)^{-1}
=\sum_{\lambda \in \Par} t^{|\lambda|} s_{\lambda}(\xx) s_\lambda(\yy).
\end{equation}
(Recall that $|\lambda|=\lambda_1 + \lambda_2 + \cdots +\lambda_\ell$
for any partition $\lambda = \left(\lambda_1,\lambda_2,\ldots,\lambda_\ell\right)$.)
\end{remark}

% [DG][v17] Replaced "a suitable completion of $R(\xx) \otimes R(\yy)$"
% by "$\kk[[\xx, \yy]]$" due to the elusiveness of $R(\xx) \otimes R(\yy)$.
% (I don't have a feeling for this tensor product, and I am not sure
% if anyone has; when $\kk$ is not a field, it may not inject into
% $R(\xx, \yy)$ by the canonical map, as shown in
% http://projecteuclid.org/DPubS?verb=Display&version=1.0&service=UI&handle=euclid.pjm/1102959646 .
% A completion *might* deal with the issues, but there are too many
% possible candidates.)

\begin{proof}[Proof of Theorem~\ref{thm.cauchy-identity}]
We follow the standard combinatorial proof
(see \cite[\S 4.8]{Sagan},\cite[\S 7.11,7.12]{Stanley}), which rewrites the
left and right sides of \eqref{Cauchy-identity-with-t}, and then compares them with the
\emph{Robinson-Schensted-Knuth} (RSK) bijection.\footnote{The
RSK bijection has been introduced by Knuth \cite{Knuth-RSK},
where what we call ``biletters'' is referred to as
``two-line arrays''. The most important ingredient of this
algorithm -- the RS-insertion operation -- however goes back to
Schensted. The special case of the RSK algorithm where the
biword has to be a permutation (written in two-line notation)
and the two tableaux have to be \emph{standard} (i.e., each
of them has content $\left(1^n\right)$, where $n$ is the size
of their shape)
is the famous \dfn{Robinson-Schensted correspondence}
\cite{Leeuwen-rsk}. More about these algorithms can be found
in \cite[Chapter 3]{Sagan}, \cite[Chapter 5]{MendesRemmel},
\cite[\S 7.11-7.12]{Stanley}, \cite[Sections 10.9--10.22]{Loehr-bij},
\cite[Chapters 1 and A]{Fulton}, \cite[\S 3, \S 6]{BritzFomin}
and various other places.}
On the left side, expanding out each geometric series
\[
(1-tx_i y_j)^{-1}=1+tx_iy_j+(tx_iy_j)^2 +(tx_iy_j)^3+\cdots
\]
and thinking of $(x_i y_j)^m$ as $m$ occurrences of a
\emph{biletter}\footnote{A \dfn{biletter} here simply means a pair
of letters, written as a column vector. A \dfn{letter} means a
positive integer.} $\binom{i}{j}$,
we see that the left hand side can be rewritten as the sum of
$t^\ell \left(x_{i_1} y_{j_1}\right) \left(x_{i_2} y_{j_2}\right)
\cdots \left(x_{i_\ell} y_{j_\ell}\right)$
over all multisets
$\left\{\binom{i_1}{j_1}, \ldots, \binom{i_\ell}{j_\ell}\right\}_{\operatorname{multiset}}$
of biletters.
Order the biletters in such a multiset in the
lexicographic order $\leq_{lex}$, which is the total order on the
set of all biletters defined by
\[
\binom{i_1}{j_1} \leq_{lex} \binom{i_2}{j_2}
\quad \Longleftrightarrow \quad
\left( \text{we have } i_1 \leq i_2 \text{, and if }
i_1 = i_2 \text{, then } j_1 \leq j_2 \right).
\]
% which first checks if $i_1 \leq i_2$ and
% then if $i_1 = i_2$ checks if $j_1 \leq j_2$.
Defining a \dfn{biword} to be an array
$\binom{\ii}{\jj} = \binom{i_1 \cdots i_\ell}{j_1 \cdots j_\ell}$
in which the biletters are ordered
$\binom{i_1}{j_1} \leq_{lex} \cdots \leq_{lex} \binom{i_\ell}{j_\ell}$,
then the left side of \eqref{Cauchy-identity-with-t}
is the sum $\sum t^{\ell} \xx^{\cont(\ii)} \yy^{\cont(\jj)}$ over
all biwords $\binom{\ii}{\jj}$, where $\ell$ stands for the number of biletters in the biword.
On the right side, expanding out the Schur functions as sums of tableaux
gives $\sum_{(P,Q)} t^{\ell} \xx^{\cont(Q)} \yy^{\cont(P)}$ in which the sum is over all ordered
pairs $(P,Q)$ of column-strict tableaux \emph{having the same shape}\footnote{And
this shape should be the Ferrers diagram of a partition (not just a skew
diagram).}, with $\ell$ cells. (We shall refer to such pairs
as \emph{tableau pairs} from now on.)

The \dfn{Robinson-Schensted-Knuth algorithm} gives us a bijection between
the biwords $\binom{\ii}{\jj}$ and the tableau pairs $(P,Q)$, which has
the property that
\begin{align*}
\cont(\ii) &= \cont(Q),\\
\cont(\jj) &= \cont(P)
\end{align*}
(and that the length $\ell$ of the biword $\binom{\ii}{\jj}$
equals the size $\left|\lambda\right|$ of the common shape of
$P$ and $Q$; but this follows automatically from
$\cont(\ii) = \cont(Q)$). Clearly, once such a bijection is
constructed, the equality \eqref{Cauchy-identity-with-t} will
follow.

Before we define this algorithm, we introduce a simpler operation
known as \dfn{RS-insertion} (short for Robinson-Schensted insertion).
RS-insertion takes as input a column-strict tableau $P$ and a letter
$j$, and returns a new column-strict tableau $P'$ along with a corner
cell\footnote{A \dfn{corner cell} of a tableau or of a Ferrers
diagram is defined to be a cell $c$ which belongs to the tableau
(resp. diagram) but whose immediate neighbors to the east and to the
south don't. For example, the cell $\left(3,2\right)$ is a corner
cell of the Ferrers diagram of the partition $\left(3,2,2,1\right)$,
and thus also of any tableau whose shape is this partition.
But the cell $\left(2,2\right)$ is not a corner cell of this
Ferrers diagram, since its immediate neighbor to the south is
still in the diagram.}
$c$ of $P'$, which is constructed as follows:
Start out by setting $P' = P$. The letter $j$ tries to insert itself
into the first row of $P'$ by either bumping out the leftmost letter in the first row
strictly larger than $j$, or else placing itself at the right end of the row if
no such larger letter exists.  If a letter was bumped from the first row, this letter follows
the same rules to insert itself into the second row, and so on\footnote{Here,
rows are allowed to be empty -- so it is possible that a letter is bumped
from the last nonempty row of $P'$ and settles in the next, initially
empty, row.}.  This series of bumps must eventually come to an
end\footnote{since we can
only bump out entries from nonempty rows}.  At the end of the bumping,
the tableau $P'$ created has an extra corner cell not present in $P$.
If we call this corner cell $c$, then $P'$ (in its final form) and $c$
are what the RS-insertion operation returns. One says that $P'$ is the
result of \emph{inserting}\footnote{This terminology is
reminiscent of insertion into binary search trees, a basic operation
in theoretical computer science. This is more than superficial similarity;
there are, in fact, various analogies between Ferrers diagrams (and
their fillings) and unlabelled plane binary trees (resp. their labellings),
and one of them is the analogy between RS-insertion and binary search
tree insertion. See
\cite[\S 4.1]{Hivert}.} $j$ into the tableau $P$. It is straightforward
to see that this resulting filling $P'$ is a column-strict
tableau\footnote{Indeed, the reader can check that $P'$ remains a
column-strict tableau throughout the algorithm that defines
RS-insertion. (The only part of this that isn't obvious is showing that
when a letter $t$ bumped out of some row $k$ is inserted into row $k+1$,
the property that the letters increase strictly down columns is
preserved. Argue that the bumping-out of $t$ from row
$k$ was caused by the insertion of another letter $u < t$, and that the
cell of row $k+1$ into which $t$ is then being inserted is in the same
column as this $u$, or in a column further left than it.)}.

\begin{example}
To give an example
of this operation, let us insert the letter $j = 3$ into the
column-strict tableau
$\begin{matrix} 1 & 1 & 3 & 3 & 4 \\ 2 & 2 & 4 & 6 & \\
3 & 4 & 7 & & \\ 5 & & & & \end{matrix}$ (we are showing all
intermediate states of $P'$; the underlined letter is always
the one that is going to be bumped out at the next step):
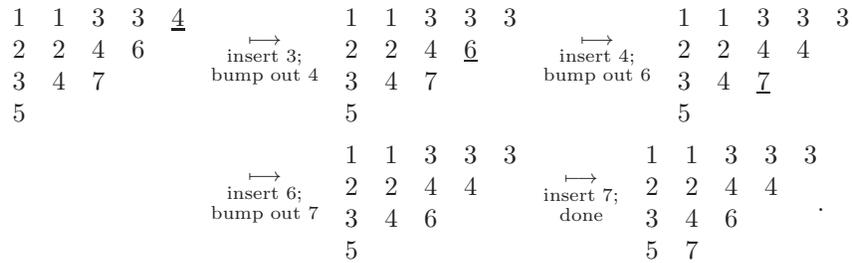
\begin{align*}
\begin{matrix} 1 & 1 & 3 & 3 & \underline{4} \\ 2 & 2 & 4 & 6 & \\
3 & 4 & 7 & & \\ 5 & & & & \end{matrix}
\quad & \overset{\longmapsto}{\substack{\text{insert 3;} \\ \text{bump out 4}}} \quad
\begin{matrix} 1 & 1 & 3 & 3 & 3 \\ 2 & 2 & 4 & \underline{6} & \\
3 & 4 & 7 & & \\ 5 & & & & \end{matrix}
\quad \overset{\longmapsto}{\substack{\text{insert 4;} \\ \text{bump out 6}}} \quad
\begin{matrix} 1 & 1 & 3 & 3 & 3 \\ 2 & 2 & 4 & 4 & \\
3 & 4 & \underline{7} & & \\ 5 & & & & \end{matrix} \\
\quad & \overset{\longmapsto}{\substack{\text{insert 6;} \\ \text{bump out 7}}} \quad
\begin{matrix} 1 & 1 & 3 & 3 & 3 \\ 2 & 2 & 4 & 4 & \\
3 & 4 & 6 & & \\ 5 & & & & \end{matrix}
\quad \overset{\longmapsto}{\substack{\text{insert 7;} \\ \text{done}}} \quad
\begin{matrix} 1 & 1 & 3 & 3 & 3 \\ 2 & 2 & 4 & 4 & \\
3 & 4 & 6 & & \\ 5 & 7 & & & \end{matrix} .
\end{align*}
The last tableau in this sequence is the column-strict
tableau that is returned. The corner cell that is returned
is the second cell of the fourth row (the one containing $7$).
\end{example}

RS-insertion will be used as a step in the RSK algorithm; the
construction will rely on a simple fact known as the
\emph{row bumping lemma}. Let us first define the notion of
a \dfn{bumping path} (or \dfn{bumping route}):
If $P$ is a column-strict tableau, and $j$ is a letter, then
some letters are inserted into some cells when RS-insertion
is applied to $P$ and $j$. The sequence of these cells (in
the order in which they see letters inserted into them) is called
the \emph{bumping path} for $P$ and $j$. This bumping path
always ends with the corner cell $c$ which is returned by
RS-insertion. As an example, when
$j=1$ is inserted into the tableau $P$
shown below, the result $P'$ is shown with all entries on
the bumping path underlined:
\[
P=\begin{matrix}
1&1&2&2&3\\
2&2&4&4&\\
3&4&5& &\\
4&6&6& &
\end{matrix}
\qquad
\underset{j=1}{\overset{\text{insert}}{\longmapsto}}
\qquad
P'=\begin{matrix}
1&1&\underline{1}&2&3\\
2&2&\underline{2}&4&\\
3&4&\underline{4}& &\\
4&\underline{5}&6& & \\
\underline{6}& & & &
\end{matrix}
\]
A first simple observation about bumping paths is that
bumping paths \emph{trend weakly left} -- that is, if the
bumping path of $P$ and $j$ is $\left(c_1,c_2,\ldots,c_k\right)$,
then, for each $1\leq i < k$, the cell $c_{i+1}$ lies in the
same column as $c_i$ or in a column further left.\footnote{This
follows easily from the preservation of column-strictness during
RS-insertion.} A subtler property of bumping paths is the following
\dfn{row bumping lemma} (\cite[p. 9]{Fulton}):

\begin{statement}
\textbf{Row bumping lemma:} Let $P$ be a column-strict tableau,
and let $j$ and $j'$ be two letters. Applying RS-insertion to the
tableau $P$ and the letter $j$ yields a new column-strict tableau
$P'$ and a corner cell $c$. Applying RS-insertion to the tableau
$P'$ and the letter $j'$ yields a new column-strict tableau $P''$
and a corner cell $c'$.
\begin{enumerate}
\item[(a)] Assume that $j \leq j'$. Then, the bumping path
for $P'$ and $j'$
stays strictly to the right, within each row, of the bumping path
for $P$ and $j$.
The cell $c'$ (in which the bumping path for $P'$ and $j'$ ends) is
in the same row as the cell $c$ (in which the bumping path for $P$
and $j$ ends) or in a row further up; it is also in a column further
right than $c$.
\item[(b)] Assume instead that $j > j'$. Then, the bumping path
for $P'$ and $j'$
stays weakly to the left, within each row, of the bumping path
for $P$ and $j$.
The cell $c'$ (in which the bumping path for $P'$ and $j'$ ends) is
in a row further down than the cell $c$ (in which the bumping
path for $P$ and $j$ ends); it is also in the same column as $c$ or
in a column further left.
\end{enumerate}
\end{statement}
This lemma can be easily proven by induction over the
row.\footnote{We leave the details to the reader, only giving the
main idea for (a) (the proof of (b) is similar). To prove the
first claim of (a), it is enough to show that for every $i$,
if any letter is inserted into row $i$ during RS-insertion for
$P'$ and $j'$, then some letter is also inserted into row $i$
during RS-insertion for $P$ and $j$, and the former insertion
happens in a cell strictly to the right of the cell where the
latter insertion happens. This follows by induction over $i$.
In the induction step, we need to show that if, for a positive
integer $i$, we try to consecutively insert two letters $k$ and $k'$, in this order,
into the $i$-th row of a column-strict tableau, possibly bumping out
existing letters in the process, and if we have $k \leq k'$, then the
cell into which $k$ is inserted is strictly to the left of the cell
into which $k'$ is inserted, and the letter bumped out by the insertion
of $k$ is $\leq$ to the letter bumped out by the insertion of $k'$ (or
else the insertion of $k'$ bumps out no letter at all -- but it cannot
happen that $k'$ bumps out a letter but $k$ does not). This statement is
completely straightforward to check (by only studying the $i$-th
row). This way, the first claim of (a) is proven, and this entails
that the cell $c'$ (being the last cell of the bumping path for $P'$
and $j'$) is in the same row as the cell $c$ or in a row further up.
It only remains to show that $c'$ is in a column further right than $c$.
This follows by noticing that, if $k$ is the row in which the cell
$c'$ lies, then $c'$ is in a column further right than
the entry of the bumping path for $P$ and $j$ in row $k$ (by the
first claim of (a)), and this latter entry is further right than
or in the same column as the ultimate entry $c$ of this bumping path
(since bumping paths trend weakly left).}

We can now define the actual RSK algorithm.
Let $\binom{\ii}{\jj}$ be a biword.
Starting with the pair $(P_0,Q_0)=(\varnothing,\varnothing)$ and
$m=0$, the algorithm applies the following steps
(see Example~\ref{RSK-example} below):
\begin{itemize}
\item If $i_{m+1}$ does not exist (that is, $m$ is the length of
$\ii$), stop.
\item Apply RS-insertion to the column-strict tableau $P_m$
and the letter $j_{m+1}$ (the bottom letter of
$\binom{i_{m+1}}{j_{m+1}}$). Let $P_{m+1}$ be the resulting
column-strict tableau, and let $c_{m+1}$ be the resulting
corner cell.
\item Create $Q_{m+1}$ from $Q_m$ by adding the top letter
$i_{m+1}$ of $\binom{i_{m+1}}{j_{m+1}}$ to $Q_m$ in the cell
$c_{m+1}$ (which, as we recall, is the extra corner cell of
$P_{m+1}$ not present in $P_m$).
\item Set $m$ to $m+1$.
\end{itemize}
After all of the biletters have been thus processed,
the result of the RSK algorithm
is $(P_\ell,Q_\ell)=:(P,Q)$.

\begin{example}
\label{RSK-example}
The term in the expansion of the left side of \eqref{Cauchy-identity} corresponding to
\[
(x_1 y_2)^1(x_1 y_4)^1(x_2 y_1)^1(x_4 y_1)^1(x_4 y_3)^2 (x_5 y_2)^1
\]
is the biword
$\binom{\ii}{\jj}=\binom{1124445}{2411332}$, whose RSK algorithm goes
as follows:
\[
\begin{array}{rcl|rcl}
P_0&= &\varnothing & Q_0 &=& \varnothing \\
 & & & & &\\
P_1&= &\begin{matrix} 2 \end{matrix}
 & Q_1 &= &\begin{matrix} 1 \end{matrix} \\
 & & & & &\\
P_2&= &\begin{matrix} 2&4 \end{matrix}
 & Q_2 &= &\begin{matrix} 1&1  \end{matrix} \\
 & & & & &\\
P_3&= &\begin{matrix} 1&4\\ 2 \end{matrix}
 & Q_3 &= &\begin{matrix} 1&1\\ 2 \end{matrix} \\
 & & & & &\\
P_4&= &\begin{matrix} 1&1\\ 2&4 \end{matrix}
 & Q_4 &= &\begin{matrix} 1&1\\ 2&4 \end{matrix} \\
 & & & & &\\
P_5&= &\begin{matrix} 1&1&3\\ 2&4& \end{matrix}
 & Q_5 &= &\begin{matrix} 1&1&4\\ 2&4& \end{matrix} \\
 & & & & &\\
P_6&= &\begin{matrix} 1&1&3&3\\ 2&4& & \end{matrix}
 & Q_6 &= &\begin{matrix} 1&1&4&4\\ 2&4& & \end{matrix} \\
 & & & & &\\
P:=P_7&= &\begin{matrix} 1&1&2&3\\ 2&3& & \\ 4& & & \end{matrix}
 & Q:=Q_7 &= &\begin{matrix} 1&1&4&4\\ 2&4& & \\5& & & \end{matrix} \\
\end{array}
\]
\end{example}

The bumping rule obviously maintains the property
that $P_m$ is a column-strict tableau of some Ferrers shape throughout.
It should be clear that $(P_m,Q_m)$ have the same shape at each stage.
Also, the construction of $Q_m$ shows that it is at least weakly
increasing in rows and weakly increasing in columns throughout.
What is perhaps least clear is
that $Q_m$ remains strictly increasing down columns.
That is, when one has a string of equal
letters on top $i_m=i_{m+1}=\cdots=i_{m+r}$,
so that on bottom one bumps in $j_m \leq j_{m+1} \leq \cdots \leq j_{m+r}$,
one needs to know that the new cells form a \dfn{horizontal strip}, that is,
no two of them lie in the same column\footnote{Actually, each of these
new cells (except for the first one) is in a column further right
than the previous one. We will use this stronger fact further below.}.
This follows from (the last claim of) part (a) of the row bumping lemma.
Hence, the result $\left(P, Q\right)$ of the RSK algorithm
is a tableau pair.

To see that the RSK map is a bijection, we show how to
recover $\binom{\ii}{\jj}$ from $(P,Q)$.  This is done by
\dfn{reverse bumping} from $(P_{m+1},Q_{m+1})$ to recover
both the biletter $\binom{i_{m+1}}{j_{m+1}}$ and the tableaux $(P_m,Q_m)$, as follows.
Firstly, $i_{m+1}$ is the maximum entry of $Q_{m+1}$, and $Q_m$ is obtained by removing
the rightmost occurrence of this letter $i_{m+1}$ from $Q_{m+1}$.
\footnote{It necessarily has to be the rightmost occurrence, since
(according to the previous footnote) the cell into which
$i_{m+1}$ was filled at the step from $Q_m$ to $Q_{m+1}$ lies further
right than any existing cell of $Q_m$ containing the letter $i_{m+1}$.}
To produce $P_m$ and $j_{m+1}$, find the position of the
rightmost occurrence of $i_{m+1}$ in $Q_{m+1}$,
and start \emph{reverse bumping} in $P_{m+1}$ from the entry in
this same position, where reverse bumping an
entry means inserting it into one row higher by having it bump
out the rightmost entry which is strictly smaller.\footnote{Let
us give a few more details on this ``reverse bumping'' procedure.
Reverse bumping (also known as \dfn{RS-deletion} or
\dfn{reverse RS-insertion}) is an operation which takes a
column-strict tableau $P'$ and a corner cell $c$ of $P'$,
and constructs a column-strict tableau $P$ and a letter $j$
such that RS-insertion for $P$ and $j$ yields $P'$ and $c$.
It starts by setting $P = P'$, and removing the entry in the
cell $c$ from $P$. This removed entry is then denoted by $k$,
and is inserted into the row of $P$ above $c$, bumping out the
rightmost entry which is smaller than $k$. The letter which is
bumped out -- say, $\ell$ --, in turn, is inserted into the row
above it, bumping out the rightmost entry which is smaller than
$\ell$. This procedure continues in the same way until an entry
is bumped out of the first row (which will eventually happen).
The reverse bumping operation returns the
resulting tableau $P$ and the entry which is bumped out of the
first row.
\par
It is straightforward to check that the reverse bumping
operation is well-defined (i.e., $P$ does stay a column-strict
tableau throughout the procedure) and is the inverse of the
RS-insertion operation. (In fact, these two operations undo each
other step by step.)}
The entry bumped out of the
first row is $j_{m+1}$, and the resulting tableau is $P_{m}$.

% [DG][v46] More details on reverse bumping.

Finally, to see that the RSK map is surjective, one needs to show
that the reverse bumping procedure can be applied to any pair
$(P,Q)$ of column-strict tableaux of the same shape, and will
result in a (lexicographically ordered) biword $\binom{\ii}{\jj}$.
We leave this verification to the reader.\footnote{It is easy to
see that repeatedly applying reverse bumping to $(P,Q)$ will
result in a sequence
$\binom{i_\ell}{j_\ell}, \binom{i_{\ell-1}}{j_{\ell-1}},
\ldots, \binom{i_1}{j_1}$ of biletters such that applying the
RSK algorithm to
$\binom{i_1 \cdots i_\ell}{j_1 \cdots j_\ell}$ gives back
$(P,Q)$. The question is why we have
$\binom{i_1}{j_1} \leq_{lex} \cdots \leq_{lex} \binom{i_\ell}{j_\ell}$.
Since the chain of inequalities $i_1 \leq i_2 \leq \cdots \leq i_\ell$
is clear from the choice of entry to reverse-bump, it only remains
to show that for every string $i_m=i_{m+1}=\cdots=i_{m+r}$
of equal top letters, the corresponding bottom letters
weakly increase (that is,
$j_m \leq j_{m+1} \leq \cdots \leq j_{m+r}$). One way to see
this is the following:
\par
Assume the contrary; i.e., assume that
the bottom letters corresponding
to some string $i_m=i_{m+1}=\cdots=i_{m+r}$ of equal top
letters do not weakly increase.
Thus, $j_{m+p} > j_{m+p+1}$ for some
$p \in \left\{0, 1, \ldots, r-1\right\}$.
Consider this $p$.
\par
Let us consider the cells containing the equal letters
$i_m=i_{m+1}=\cdots=i_{m+r}$ in the tableau $Q_{m+r}$.
Label these cells as $c_m, c_{m+1}, \ldots, c_{m+r}$
from left to right (noticing that no two of them lie
in the same column, since $Q_{m+r}$ is column-strict).
By the definition of reverse bumping, the first entry
to be reverse bumped from $P_{m+r}$ is the entry in
position $c_{m+r}$ (since this is the rightmost
occurrence of the letter $i_{m+r}$ in $Q_{m+r}$); then, the next entry
to be reverse bumped is the one in position $c_{m+r-1}$,
etc., moving further and further left.
Thus, for each $q \in \left\{0, 1, \ldots, r\right\}$,
the tableau $P_{m+q-1}$ is obtained from $P_{m+q}$ by
reverse bumping the entry in position $c_{m+q}$.
Hence, conversely, the tableau $P_{m+q}$ is obtained
from $P_{m+q-1}$ by RS-inserting the entry $j_{m+q}$,
which creates the corner cell $c_{m+q}$.
\par
But recall that $j_{m+p} > j_{m+p+1}$.
Hence, part (b) of the row bumping lemma (applied to
$P_{m+p-1}$, $j_{m+p}$, $j_{m+p+1}$, $P_{m+p}$,
$c_{m+p}$, $P_{m+p+1}$ and $c_{m+p+1}$
instead of $P$, $j$, $j'$, $P'$, $c$, $P''$ and $c'$)
shows that the cell $c_{m+p+1}$ is in the same
column as the cell $c_{m+p}$ or in a column further
left.
But this contradicts the fact that the cell $c_{m+p+1}$ is
in a column further right than the cell $c_{m+p}$ (since
we have labeled our cells as
$c_m, c_{m+1}, \ldots, c_{m+r}$ from left to right,
and no two of them lied in the same column).
This contradiction completes our proof.}
\end{proof}

% [DG][v45] Revamped the above proof, adding significant amounts of
% detail. What triggered the rewriting was my realization that both
% parts of the row bumping lemma are needed (as opposed to only the
% first one); but now that I'm doing it I've also added an example
% for RS-insertion, separated RS-insertion from the whole RSK
% algorithm (which is no longer called "RSK insertion" to remove
% confusion), inserted some references, and moved out the row
% bumping lemma from the niche it occupied in the proof of RSK.
% I hope this is clearer now...

This is by far not the only known proof of Theorem~\ref{thm.cauchy-identity}.
Two further proofs will be sketched in Exercise~\ref{exe.pieri.cauchy}
and Exercise~\ref{exe.cauchy.cauchy}.

Before we move on to extracting identities in $\Lambda$ from Theorem
\ref{thm.cauchy-identity}, let us state (as an exercise) a simple technical
fact that will be useful:

\begin{exercise}
\label{exe.Lambda.infinite-indep}
Let $\left( q_\lambda \right) _{\lambda \in \Par}$
be a basis of the $\kk$-module $\Lambda$.
Assume that for each partition $\lambda$, the element
$q_\lambda \in \Lambda$
is homogeneous of degree $\left\vert \lambda\right\vert $.

\begin{enumerate}
\item[(a)] If two families $\left(  a_{\lambda}\right)  _{\lambda
\in \Par }\in\kk^{ \Par }$ and $\left(
b_{\lambda}\right)  _{\lambda\in \Par }\in\kk%
^{ \Par }$ satisfy%
\begin{equation}
\sum_{\lambda\in \Par }a_{\lambda}q_{\lambda}\left(
\xx\right)  =\sum_{\lambda\in \Par }b_{\lambda}q_{\lambda
}\left(  \xx\right)  \label{eq.exe.Lambda.infinite-indep.a.ass}%
\end{equation}
in $\kk\left[  \left[  \xx\right]  \right]  $, then $\left(
a_{\lambda}\right)  _{\lambda\in \Par }=\left(  b_{\lambda
}\right)  _{\lambda\in \Par }$.\ \ \ \ \footnote{Note that this
does not immediately follow from the linear independence of the basis $\left(
q_{\lambda}\right)  _{\lambda\in \Par }$. Indeed, linear
independence would help if the sums in
(\ref{eq.exe.Lambda.infinite-indep.a.ass}) were finite, but they are not. A
subtler argument (involving the homogeneity of the $q_{\lambda}$) thus has to
be used.}

\item[(b)] Consider a further infinite family $\mathbf{y}=\left(  y_{1}%
,y_{2},y_{3},\ldots\right)  $ of indeterminates (disjoint from $\xx$).
If two families $\left(  a_{\mu,\nu}\right)  _{\left(  \mu,\nu\right)
\in\operatorname{Par}^{2}}\in\kk^{\operatorname{Par}^{2}}$ and $\left(
b_{\mu,\nu}\right)  _{\left(  \mu,\nu\right)  \in\operatorname{Par}^{2}}%
\in\kk^{\operatorname{Par}^{2}}$ satisfy
\begin{equation}
\sum_{\left(  \mu,\nu\right)  \in\operatorname{Par}^{2}}a_{\mu,\nu}q_{\mu
}\left( \xx \right)  q_{\nu}\left( \yy \right)
=\sum_{\left(  \mu,\nu\right)  \in\operatorname{Par}^{2}}b_{\mu,\nu}q_{\mu
}\left( \xx \right)  q_{\nu}\left( \yy \right)
\label{eq.exe.Lambda.infinite-indep.b.ass}
\end{equation}
in $\kk\left[  \left[  \xx,\mathbf{y}\right]  \right]  $, then
$\left(  a_{\mu,\nu}\right)  _{\left(  \mu,\nu\right)  \in\operatorname{Par}%
^{2}}=\left(  b_{\mu,\nu}\right)  _{\left(  \mu,\nu\right)  \in
\operatorname{Par}^{2}}$.

\item[(c)] Consider a further infinite family $\mathbf{z}=\left(  z_{1}%
,z_{2},z_{3},\ldots\right)  $ of indeterminates (disjoint from $\xx$
and $\mathbf{y}$). If two families $\left(  a_{\lambda,\mu,\nu}\right)
_{\left(  \mu,\nu,\lambda\right)  \in\operatorname{Par}^{3}}\in\kk%
^{\operatorname{Par}^{3}}$ and $\left(  b_{\lambda,\mu,\nu}\right)  _{\left(
\mu,\nu,\lambda\right)  \in\operatorname{Par}^{3}}\in\kk%
^{\operatorname{Par}^{3}}$ satisfy
\begin{equation}
\sum_{\left(  \mu,\nu,\lambda\right)  \in\operatorname{Par}^{3}}a_{\lambda
,\mu,\nu}q_{\mu}\left(  \xx\right)  q_{\nu}\left(  {\mathbf{y}%
}\right)  q_{\lambda}\left(  {\mathbf{z}}\right)  =\sum_{\left(  \mu
,\nu,\lambda\right)  \in\operatorname{Par}^{3}}b_{\lambda,\mu,\nu}q_{\mu
}\left(  \xx\right)  q_{\nu}\left(  {\mathbf{y}}\right)  q_{\lambda
}\left(  {\mathbf{z}}\right)  \label{eq.exe.Lambda.infinite-indep.c.ass}%
\end{equation}
in $\kk\left[  \left[  \xx,\mathbf{y},\mathbf{z}\right]
\right]  $, then $\left(  a_{\lambda,\mu,\nu}\right)  _{\left(  \mu
,\nu,\lambda\right)  \in\operatorname{Par}^{3}}=\left(  b_{\lambda,\mu,\nu
}\right)  _{\left(  \mu,\nu,\lambda\right)  \in\operatorname{Par}^{3}}$.
\end{enumerate}
\end{exercise}

\begin{remark}
Clearly, for any $n\in \NN$, we can state an analogue of Exercise
\ref{exe.Lambda.infinite-indep} for $n$ infinite families $\xx%
_{i}=\left(  x_{i,1},x_{i,2},x_{i,3},\ldots\right)  $ of indeterminates (with
$i\in\left\{  1,2,\ldots,n\right\}  $). The three parts of Exercise
\ref{exe.Lambda.infinite-indep} are the particular cases of this analogue for
$n=1$, for $n=2$ and for $n=3$. We have shied away from stating this analogue
in full generality
because these particular cases are the only ones we will need.
\end{remark}

\begin{corollary}
\label{self-duality-corollary}
In the Schur function basis $\{s_\lambda\}$ for $\Lambda$, the
structure constants for multiplication and comultiplication are the same, that is,
if one defines scalars
$c^\lambda_{\mu,\nu}, \hat{c}^\lambda_{\mu,\nu}$ via the unique expansions
\begin{align}
s_\mu s_\nu &= \sum_{\lambda} c^\lambda_{\mu,\nu} s_\lambda ,
\label{symm-structure-constants.1} \\
\Delta(s_\lambda) &= \sum_{\mu, \nu} \hat{c}^\lambda_{\mu,\nu} s_\mu \otimes s_\nu ,
\label{symm-structure-constants.2}
\end{align}
then $c^\lambda_{\mu,\nu}=\hat{c}^\lambda_{\mu,\nu}$.
\end{corollary}
\begin{proof}
Work in the ring $\kk\left[\left[\xx , \yy , \zz\right]\right]$, where
$\yy = \left(y_1, y_2, y_3, \ldots\right)$ and $\zz = \left(z_1, z_2,
z_3, \ldots\right)$ are two new sets of variables.
The identity \eqref{Cauchy-identity}
lets one interpret both $c^\lambda_{\mu,\nu}, \hat{c}^\lambda_{\mu,\nu}$
as the coefficient\footnote{Let us explain why speaking of
coefficients makes sense here:

We want to use the fact that if a power series
$f \in \kk\left[\left[\xx, \yy, \zz\right]\right]$ is written
in the form $f = \sum_{\left(\mu, \nu, \lambda\right) \in \Par^3}
a_{\lambda, \mu, \nu} s_\mu\left(\xx\right) s_\nu\left(\yy\right)
s_\lambda\left(\zz\right)$ for some coefficients
$a_{\lambda, \mu, \nu} \in \kk$, then these coefficients
$a_{\lambda, \mu, \nu}$ are uniquely determined by $f$.
But this fact is precisely the claim of
Exercise~\ref{exe.Lambda.infinite-indep}(c) above
(applied to $q_\lambda = s_\lambda$).}
 of $s_\mu(\xx) s_\nu(\yy) s_\lambda(\zz)$ in
the product
\begin{align*}
\prod_{i,j=1}^\infty (1-x_i z_j)^{-1}
\prod_{i,j=1}^\infty (1-y_i z_j)^{-1}
&\overset{\eqref{Cauchy-identity}}{=} \left( \sum_\mu s_\mu(\xx) s_\mu(\zz) \right)
\left( \sum_\nu s_\nu(\yy) s_\nu(\zz) \right) \\
&=\sum_{\mu,\nu} s_\mu(\xx) s_\nu(\yy) \cdot s_\mu(\zz) s_\nu(\zz) \\
&=
\sum_{\mu,\nu} s_\mu(\xx) s_\nu(\yy) \left( \sum_\lambda c^\lambda_{\mu,\nu} s_\lambda(\zz) \right)
\end{align*}
since, regarding $x_1,x_2,\ldots,y_1,y_2,\ldots$ as lying in a single variable set
$(\xx,\yy)$, separate from the variables $\zz$, the Cauchy identity \eqref{Cauchy-identity}
expands the same product as
\begin{align*}
\prod_{i,j=1}^\infty (1-x_i z_j)^{-1}
\prod_{i,j=1}^\infty (1-y_i z_j)^{-1}
&=\sum_{\lambda} s_\lambda(\xx,\yy) s_\lambda(\zz)\\
&= \sum_\lambda \left( \sum_{\mu,\nu} \hat{c}^\lambda_{\mu,\nu}
     s_\mu(\xx) s_\nu(\yy) \right) s_\lambda(\zz).
\end{align*}
\end{proof}

% [DG][v17] Added footnote.

% [DG][v25] Added first sentence to the proof.

\begin{definition}
The coefficients $c^\lambda_{\mu,\nu}=\hat{c}^\lambda_{\mu,\nu}$ appearing in
the expansions \eqref{symm-structure-constants.1}
and \eqref{symm-structure-constants.2}
are called \dfn{Littlewood-Richardson coefficients}.
\end{definition}

\begin{remark}
\label{third-L-R-interpretation-remark}
We will interpret $c^\lambda_{\mu,\nu}$
combinatorially in Section~\ref{bialternant-section}. By now,
however, we can already prove some properties of these
coefficients:

We have
\begin{equation}
\label{eq.third-L-R-interpretation-remark.commutativity}
c^\lambda_{\mu,\nu} = c^\lambda_{\nu,\mu}
\qquad \qquad \text{ for all } \lambda, \mu, \nu \in \Par
\end{equation}
(by comparing coefficients in
$\sum_{\lambda} c^\lambda_{\mu,\nu} s_\lambda = s_\mu s_\nu
= s_\nu s_\mu = \sum_{\lambda} c^\lambda_{\nu,\mu} s_\lambda$).
Furthermore, let
$\lambda$ and $\mu$ be two partitions (not necessarily
satisfying $\mu \subseteq \lambda$).
Comparing the expansion
\[
s_\lambda(\xx,\yy) = \Delta(s_\lambda)
   = \sum_{\mu,\nu} c^\lambda_{\mu,\nu} s_\mu(\xx) s_\nu(\yy)
   = \sum_{\mu \in \Par} \left( \sum_{\nu \in \Par} c^\lambda_{\mu,\nu} s_\nu(\yy) \right) s_\mu(\xx)
\]
with
\[
s_\lambda(\xx,\yy) = \sum_{\mu \subseteq \lambda} s_\mu(\xx) s_{\lambda/\mu}(\yy)
= \sum_{\mu \in \Par} s_\mu(\xx) s_{\lambda/\mu}(\yy)
\]
\footnote{In the last equality, we removed the condition
$\mu \subseteq \lambda$ on the addends of the sum; this does not
change the value of the sum (because we have $s_{\lambda/\mu} = 0$
whenever we don't have $\mu \subseteq \lambda$).},
one concludes that
\[
\sum_{\mu \in \Par} \left( \sum_{\nu \in \Par} c^\lambda_{\mu,\nu} s_\nu(\yy) \right) s_\mu(\xx)
= \sum_{\mu \in \Par} s_\mu(\xx) s_{\lambda/\mu}(\yy)
= \sum_{\mu \in \Par} s_{\lambda/\mu}(\yy) s_\mu(\xx) .
\]
Treating the indeterminates $\yy$ as constants, and
comparing coefficients before $s_\mu(\xx)$ on both sides
of this equality\footnote{``Comparing coefficients''
means applying Exercise~\ref{exe.Lambda.infinite-indep}(a)
to $q_\lambda = s_\lambda$
in this case (although the base ring $\kk$ is now replaced
by $\kk\left[\left[\yy\right]\right]$, and the index
$\mu$ is used instead of $\lambda$, since $\lambda$ is
already taken).}, we arrive at
another standard interpretation for $c^\lambda_{\mu,\nu}$:
\[
s_{\lambda/\mu} = \sum_{\nu} c^\lambda_{\mu,\nu} s_\nu .
\]
% whenever $\lambda$ and $\mu$ are two partitions (not necessarily
% satisfying $\mu \subseteq \lambda$).
In particular, $c^\lambda_{\mu,\nu}$ vanishes unless
$\mu \subseteq \lambda$. Consequently,
$c^\lambda_{\mu,\nu}$ vanishes unless $\nu \subseteq \lambda$
as well (since $c^\lambda_{\mu,\nu} = c^\lambda_{\nu,\mu}$)
and furthermore vanishes unless the equality
$\left|\mu\right|+\left|\nu\right| = \left|\lambda\right|$
holds\footnote{In fact, this is clear when we don't have
$\mu \subseteq \lambda$. When we do have $\mu \subseteq \lambda$,
this follows from observing that
$s_{\lambda/\mu} \in \Lambda_{\left|\lambda/\mu\right|}$
has zero coefficient before $s_\nu$ whenever
$\left|\mu\right|+\left|\nu\right| \neq \left|\lambda\right|$.}.
Altogether, we conclude that $c^\lambda_{\mu,\nu}$ vanishes unless
$\mu, \nu \subseteq \lambda$ and $|\mu|+|\nu|=|\lambda|$.
\end{remark}

% [DG][v25] Added many details to the above remark.

\begin{exercise}
\label{exe.Lambda.bialg-c}
Show that any four partitions $\kappa$, $\lambda$, $\varphi$ and
$\psi$ satisfy
\begin{vershort}
\[
\sum_{\rho \in \Par} c^\rho_{\kappa, \lambda} c^\rho_{\varphi, \psi}
= \sum_{\left(\alpha , \beta , \gamma , \delta\right) \in \Par^4}
c^\lambda_{\beta, \delta} c^\varphi_{\alpha, \beta}
c^\psi_{\gamma, \delta}.
\]
\end{vershort}
\begin{verlong}
\[
\sum_{\rho \in \Par} c^\rho_{\kappa, \lambda} c^\rho_{\varphi, \psi}
= \sum_{\alpha \in \Par} \sum_{\beta \in \Par} \sum_{\gamma \in \Par}
\sum_{\delta \in \Par} c^\kappa_{\alpha, \gamma}
c^\lambda_{\beta, \delta} c^\varphi_{\alpha, \beta}
c^\psi_{\gamma, \delta}.
\]
\end{verlong}
\end{exercise}

% [DG][v25] Added the above exercise. Is it too scary? (It is the
% bialgebra axiom, written in Schur coordinates, of course.)

\begin{exercise} \phantomsection
\label{exe.cauchy.skew}
\begin{itemize}
\item[(a)] For any partition $\mu$, prove that
\[
\sum_{\lambda \in \Par} s_\lambda\left(\xx\right) s_{\lambda/\mu}\left(\yy\right)
= s_\mu\left(\xx\right) \cdot \prod_{i,j=1}^\infty \left(1-x_iy_j\right)^{-1}
\]
in the power series ring $\kk\left[\left[\xx, \yy\right]\right]
= \kk\left[\left[x_1, x_2, x_3, \ldots, y_1, y_2, y_3, \ldots\right]\right]$.
\item[(b)] Let $\alpha$ and $\beta$ be two partitions. Show that
\[
\sum_{\lambda \in \Par} s_{\lambda/\alpha}\left(\xx\right) s_{\lambda/\beta}\left(\yy\right)
= \left( \sum_{\rho \in \Par} s_{\beta/\rho}\left(\xx\right) s_{\alpha/\rho}\left(\yy\right) \right)
 \cdot \prod_{i,j=1}^\infty \left(1-x_iy_j\right)^{-1}
\]
in the power series ring $\kk\left[\left[\xx, \yy\right]\right]
= \kk\left[\left[x_1, x_2, x_3, \ldots, y_1, y_2, y_3, \ldots\right]\right]$.
\end{itemize}

[\textbf{Hint:} For (b), expand the product
\[
\prod_{i,j=1}^\infty \left(1-x_i y_j\right)^{-1}
\prod_{i,j=1}^\infty \left(1-x_i w_j\right)^{-1}
\prod_{i,j=1}^\infty \left(1-z_i y_j\right)^{-1}
\prod_{i,j=1}^\infty \left(1-z_i w_j\right)^{-1}
\]
in the power series ring
$\kk\left[\left[x_1,x_2,x_3, \ldots ,y_1,y_2,y_3, \ldots ,z_1,z_2,z_3, \ldots ,w_1,w_2,w_3, \ldots \right]\right]$
in two ways: once by applying Theorem~\ref{thm.cauchy-identity} to the
two variable sets $\left(\zz, \xx\right)$ and $\left(\mathbf{w}, \yy\right)$
and then using \eqref{eq.Lambda.s.xy}; once again by applying
\eqref{Cauchy-identity} to the two variable sets $\zz$ and $\mathbf{w}$ and
then applying Exercise~\ref{exe.cauchy.skew}(a) twice.]
\end{exercise}

The statement of Exercise~\ref{exe.cauchy.skew}(b) is known as the
\dfn{skew Cauchy identity}, and appears in
Sagan-Stanley \cite[Cor. 6.12]{SaganStanley},
Stanley \cite[exercise 7.27(c)]{Stanley}
and Macdonald \cite[\S I.5, example 26]{Macdonald}; it seems to be due to
Zelevinsky. It generalizes the statement of Exercise~\ref{exe.cauchy.skew}(a),
which in turn is a generalization of Theorem~\ref{thm.cauchy-identity}.

% [DG][v25] Added the above exercise (following Macdonald mostly).

\begin{definition}
\label{def.hall-inner-product}
Define the \dfn{Hall inner product} on $\Lambda$ to be the $\kk$-bilinear
form $(\cdot,\cdot)$ which makes $\{s_\lambda\}$ an orthonormal basis,
that is, $(s_\lambda, s_\nu)=\delta_{\lambda,\nu}$.
\end{definition}

\begin{exercise} \phantomsection
\label{exe.Lambda.hall.gr}
\begin{itemize}
\item[(a)] If $n$ and $m$ are two distinct nonnegative integers, and if
$f \in \Lambda_n$ and $g \in \Lambda_m$, then show that
$\left(f, g\right) = 0$.
\item[(b)] If $n \in \NN$ and $f \in \Lambda_n$, then prove that
$\left(h_n, f\right) = f\left(1\right)$ (where $f\left(1\right)$ is
defined as in Exercise~\ref{exe.Lambda.subs.fin}).
\item[(c)] Show that $\left(f, g\right) = \left(g, f\right)$ for all
$f \in \Lambda$ and $g \in \Lambda$.
(In other words, the Hall inner product is symmetric.)
\end{itemize}
\end{exercise}

% [DG][v25] Added the preceding (trivial) exercise.

% [DG][v78] Added the trivial part (c) (for ease of referencing).

The Hall inner product induces a $\kk$-module homomorphism
$\Lambda \to \Lambda^o$ (sending every $f \in \Lambda$ to
the $\kk$-linear map $\Lambda \to \kk,\ g \mapsto \left(f, g\right)$).
This homomorphism is invertible (since the Hall inner product
has an orthonormal basis), so that $\Lambda^o \cong \Lambda$
as $\kk$-modules. But in fact, more can be said:

% [DG][v46] Added the above paragraph to clarify how the isomorphism
% below is obtained.

\begin{corollary}
\label{cor.Lambda.selfdual}
The isomorphism $\Lambda^o \cong \Lambda$ induced by the
Hall inner product is an isomorphism of Hopf algebras.
\end{corollary}

\begin{proof}
We have seen that the orthonormal basis  $\{s_\lambda\}$ of Schur functions is
\emph{self-dual}, in the sense that its
multiplication and comultiplication structure constants
are the same.  Thus the isomorphism $\Lambda^o \cong \Lambda$ induced by the
Hall inner product is an isomorphism of bialgebras\footnote{Here are
some details on the proof:
\par
Let $\gamma : \Lambda \to \Lambda^o$ be the $\kk$-module
isomorphism $\Lambda \to \Lambda^o$ induced by the Hall
inner product.
We want to show that $\gamma$ is an isomorphism of bialgebras.
\par
Let $\left\{s_\lambda^*\right\}$ be the basis of $\Lambda^o$
dual to the basis
$\left\{s_\lambda\right\}$ of $\Lambda$.
Thus, for any partition $\lambda$, we have
\begin{equation}
\gamma\left(s_\lambda\right) = s_\lambda^*
\label{pf.cor.Lambda.selfdual.fn1.0}
\end{equation}
(since any partition $\mu$ satisfies
$\left(\gamma\left(s_\lambda\right)\right)\left(s_\mu\right)
= \left(s_\lambda, s_\mu\right)
= \delta_{\lambda, \mu}
= s_\lambda^* \left(s_\mu\right)$,
and thus the two $\kk$-linear maps
$\gamma\left(s_\lambda\right) : \Lambda \to \kk$
and $s_\lambda^* : \Lambda \to \kk$
are equal to each other on the basis
$\left\{s_\mu\right\}$ of $\Lambda$,
which forces them to be identical).
\par
The coproduct structure constants of the basis $\left\{s_\lambda^*\right\}$
of $\Lambda^o$ equal the product structure constants of the basis
$\left\{s_\lambda\right\}$ of $\Lambda$
(according to our discussion of duals in Section~\ref{subsect.duals}).
Since the latter are the Littlewood-Richardson numbers
$c^\lambda_{\mu,\nu}$ (because of
\eqref{symm-structure-constants.1}), we thus conclude that the former
are $c^\lambda_{\mu,\nu}$ as well.
In other words, every $\lambda \in \Par$ satisfies
\begin{equation}
\Delta_{\Lambda^o} s_\lambda^*
= \sum_{\mu, \nu} c^\lambda_{\mu, \nu} s_\mu^* \otimes s_\nu^*
\label{pf.cor.Lambda.selfdual.fn1.1}
\end{equation}
(where the sum is over all pairs $\left(\mu, \nu\right)$ of partitions).
On the other hand, applying the map
$\gamma \otimes \gamma : \Lambda \otimes \Lambda \to \Lambda^o \otimes \Lambda^o$
to the equality
\eqref{symm-structure-constants.2} yields
\begin{align*}
\left(\gamma \otimes \gamma\right) \left( \Delta(s_\lambda) \right)
&= \left(\gamma \otimes \gamma\right) \left(\sum_{\mu, \nu} \hat{c}^\lambda_{\mu,\nu} s_\mu \otimes s_\nu \right)
= \sum_{\mu, \nu} \underbrace{\hat{c}^\lambda_{\mu,\nu}}_{=c^\lambda_{\mu,\nu}} \underbrace{\gamma\left(s_\mu\right)}_{\substack{= s_\mu^* \\ \text{(by \eqref{pf.cor.Lambda.selfdual.fn1.0})}}} \otimes \underbrace{\gamma\left(s_\nu \right)}_{\substack{= s_\nu^* \\ \text{(by \eqref{pf.cor.Lambda.selfdual.fn1.0})}}}
= \sum_{\mu, \nu} c^\lambda_{\mu, \nu} s_\mu^* \otimes s_\nu^*
\\
&= \Delta_{\Lambda^o} \underbrace{s_\lambda^*}_{\substack{= \gamma\left(s_\lambda\right) \\ \text{(by \eqref{pf.cor.Lambda.selfdual.fn1.0})}}}
\qquad \left(\text{by \eqref{pf.cor.Lambda.selfdual.fn1.1}}\right) \\
&= \Delta_{\Lambda^o} \left( \gamma\left(s_\lambda\right) \right)
\end{align*}
for each $\lambda \in \Par$.
In other words, the two $\kk$-linear maps
$\left(\gamma \otimes \gamma\right) \circ \Delta$ and
$\Delta_{\Lambda^o} \circ \gamma$ are equal to each other on
each $s_\lambda$ with $\lambda \in \Par$.
Hence, these two maps must be identical (since the $s_\lambda$ form a
basis of $\Lambda$). Hence,
$\Delta_{\Lambda^o} \circ \gamma = \left(\gamma \otimes \gamma\right) \circ \Delta$.
\par
Our next goal is to show that $\epsilon_{\Lambda^o} \circ \gamma = \epsilon$.
% Indeed, every partition $\lambda \in \Par$ satisfies
% $s_\varnothing^*\left(s_\lambda\right)
% = \left(s_\varnothing, s_\lambda\right)
% = \delta_{\varnothing, \lambda}
% = \epsilon\left(s_\lambda\right)$.
% Thus, $s_\varnothing^* = \epsilon$ as $\kk$-linear maps
% $\Lambda \to \kk$ (since the $s_\lambda$ form a basis of $\Lambda$).
% \par
Indeed, each $\lambda \in \Par$ satisfies
\begin{align*}
\left( \epsilon_{\Lambda^o} \circ \gamma \right) \left(s_\lambda\right)
&= \epsilon_{\Lambda^o} \left(\gamma\left(s_\lambda\right)\right)
= \left(\gamma\left(s_\lambda\right)\right) \left(1\right)
\qquad \left(\text{by the definition of } \epsilon_{\Lambda^o}\right) \\
&= \left(s_\lambda, \underbrace{1}_{= s_\varnothing}\right)
= \left(s_\lambda, s_\varnothing\right)
= \delta_{\lambda, \varnothing}
= \epsilon \left(s_\lambda\right) .
\end{align*}
Hence, $\epsilon_{\Lambda^o} \circ \gamma = \epsilon$.
Combined with
$\Delta_{\Lambda^o} \circ \gamma = \left(\gamma \otimes \gamma\right) \circ \Delta$,
this shows that $\gamma$ is a $\kk$-coalgebra homomorphism.
Similar reasoning can be used to prove that $\gamma$ is a $\kk$-algebra homomorphism.
Altogether, we thus conclude that $\gamma$ is a bialgebra homomorphism.
Since $\gamma$ is a $\kk$-module isomorphism,
this yields that $\gamma$ is an isomorphism of bialgebras.
Qed.
}, and hence also a Hopf
algebra isomorphism by Corollary~\ref{cor.bialg-mor-is-Hopf}.
\end{proof}

% [DG][v64] Added the long footnote above, since I remember tripping up
% there when I first read the notes.

We next identify two other dual pairs of bases, by expanding
the Cauchy product in two other ways.

\begin{proposition}
\label{prop-Two-other-Cauchy-expansions}
One can also expand
\begin{equation}
\label{Two-other-Cauchy-expansions}
\prod_{i,j=1}^\infty (1-x_i y_j)^{-1}
=\sum_{\lambda \in \Par} h_{\lambda}(\xx) m_\lambda(\yy)
=\sum_{\lambda \in \Par} z_\lambda^{-1} p_{\lambda}(\xx) p_\lambda(\yy)
\end{equation}
where $z_\lambda:=m_1! \cdot 1^{m_1} \cdot m_2! \cdot 2^{m_2} \cdots$%
\index{$z_\lambda$}
if $\lambda$ is written in multiplicative notation as
$\lambda=(1^{m_1},2^{m_2},\ldots)$ with multiplicity $m_i$ for the part $i$.
(Here, we assume that $\QQ$ is a subring of $\kk$ for the last equality.)
\end{proposition}

% [DG][v17] Added last sentence.

\begin{remark}
\label{rmk.zlambda.centralizer}
It is relevant later (and explains the notation) that $z_\lambda$ is
the size of the $\Symm_n$-centralizer subgroup for a permutation having
cycle type\footnote{If $\sigma$ is a permutation of a finite set $X$,
then the \dfn{cycle type} of $\sigma$ is defined
as the list of the lengths of all cycles of $\sigma$ (that is,
of all orbits of $\sigma$ acting on $X$) written in decreasing order.
This is clearly a partition of $\left|X\right|$. (Some other authors
write it in increasing order instead, or treat it as a multiset.)

For instance, the permutation of the set $\left\{0, 3, 6, 9, 12\right\}$
that sends $0$ to $3$, $3$ to $9$, $6$ to $6$, $9$ to $0$, and
$12$ to $12$ has cycle type $\left(3, 1, 1\right)$, since
the cycles of this permutation have lengths $3$, $1$ and $1$.

It is known that two permutations in $\Symm_n$ have the same cycle
type if and only if they are conjugate. Thus, for a given partition
$\lambda$ with $\left|\lambda\right| = n$, any two permutations in
$\Symm_n$ having cycle type $\lambda$ are conjugate and therefore
their $\Symm_n$-centralizer subgroups have the same size.}
$\lambda$ with $|\lambda|=n$. This is a classical (and
fairly easy) result (see, e.g., \cite[Prop. 1.1.1]{Sagan} or
\cite[Prop. 7.7.3]{Stanley} for a proof).
\end{remark}

% [DG][v45] Added references to Sagan & Stanley.

% [DG][v46] Added footnote defining cycle type, as most authors do not
% encode it as a partition. (Annoyingly, Keith Conrad of all people
% writes it increasingly; but worse are the sources that leave out
% 1-cycles.)

\begin{proof}[Proof of Proposition~\ref{prop-Two-other-Cauchy-expansions}.]
For the first expansion, note that \eqref{h-generating-function-no-HE} shows
\begin{align*}
\prod_{i,j=1}^\infty (1-x_i y_j)^{-1}
&=\prod_{j=1}^\infty \sum_{n \geq 0} h_n(\xx) y_j^n\\
&=\sum\limits_{\substack{ \text{weak }\\\text{ compositions } \\ (n_1,n_2,\ldots)}}
     ( h_{n_1}(\xx)h_{n_2}(\xx)\cdots  ) (y_1^{n_1} y_2^{n_2}\cdots )\\
&=\sum_{\lambda \in \Par}
     \sum\limits_{\substack{ \text{weak }\\\text{ compositions } \\ (n_1,n_2,\ldots)
              \\ \text{ satisfying } \\ (n_1,n_2,\ldots) \in \Symm_{(\infty)}\lambda }}
     \underbrace{( h_{n_1}(\xx)h_{n_2}(\xx)\cdots  )}_{\substack{=h_\lambda(\xx) \\
                            \text{ (since } (n_1,n_2,\ldots) \in \Symm_{(\infty)}\lambda \text{)}}}
     \underbrace{(y_1^{n_1} y_2^{n_2}\cdots )}_{=\yy^{(n_1, n_2, \ldots)}} \\
&=\sum_{\lambda \in \Par}
     h_\lambda(\xx)
     \underbrace{\sum\limits_{\substack{ \text{weak }\\\text{ compositions } \\ (n_1,n_2,\ldots)
              \\ \text{ satisfying } \\ (n_1,n_2,\ldots) \in \Symm_{(\infty)}\lambda }}
     \yy^{(n_1, n_2, \ldots)}}_{= m_\lambda(\yy)} \\
&=\sum_{\lambda \in \Par} h_{\lambda}(\xx) m_\lambda(\yy).
\end{align*}
For the second expansion (and for later use in
the proof of Theorem~\ref{Sym-to-Hall-isomorphism-thm})
note that
\begin{equation}
\label{H(t)-logarithm}
\log H(t)
= \log \prod_{i=1}^\infty (1-x_i t)^{-1}
=\sum_{i=1}^\infty -\log(1-x_i t)
=\sum_{i=1}^\infty \sum_{m=1}^\infty \frac{(x_it)^m}{m}
=\sum_{m=1}^\infty \frac{1}{m} p_m(\xx) t^m ,
\end{equation}
so that taking $\frac{d}{dt}$ then shows that
\begin{equation}
\label{E(t)-H(t)-P(t)-relation}
P(t):=\sum_{m \geq 0} p_{m+1}t^m = \frac{H'(t)}{H(t)}= H'(t) E(-t).
\end{equation}
A similar calculation shows that
\begin{equation}
\label{eq.cauchy-via-simple-power-sums}
%\begin{aligned}
\log \prod_{i,j=1}^\infty (1-x_i y_j)^{-1}
%&=\sum_{i,j=1}^\infty -\log(1-x_i y_j) \\
%&=\sum_{i,j=1}^\infty \left( \frac{(x_iy_j)^1}{1}+\frac{(x_iy_j)^2}{2}+\frac{(x_iy_j)^3}{3}+\cdots \right) \\
%&=\sum_{m=1}^\infty \frac{1}{m} \sum_{i,j=1}^\infty (x_iy_j)^m
=\sum_{m=1}^\infty \frac{1}{m} p_m(\xx) p_m(\yy)
%\end{aligned}
\end{equation}
and hence
\begin{align*}
&\prod_{i,j=1}^\infty (1-x_i y_j)^{-1}
=\exp\left( \sum_{m=1}^\infty \frac{1}{m} p_m(\xx) p_m(\yy) \right)
= \prod_{m=1}^\infty \exp \left( \frac{1}{m} p_m(\xx) p_m(\yy) \right) \\
&= \prod_{m=1}^\infty \sum_{k=0}^\infty \frac{1}{k!} \left( \frac{1}{m} p_m(\xx) p_m(\yy) \right)^k
= \sum\limits_{\substack{ \text{weak compositions } \\ \left(k_1,k_2,k_3,\ldots\right)}}
  \prod_{m=1}^\infty \left(\frac{1}{k_m!} \left( \frac{1}{m} p_m(\xx) p_m(\yy) \right)^{k_m}\right) \\
&  \qquad \qquad \left(\text{by the product rule}\right) \\
&= \sum\limits_{\substack{ \text{weak compositions } \\ \left(k_1,k_2,k_3,\ldots\right)}}
  \prod_{m=1}^\infty \frac{\left(p_m(\xx) p_m(\yy)\right)^{k_m}}{k_m! m^{k_m}}
= \sum\limits_{\substack{ \text{weak compositions } \\ \left(k_1,k_2,k_3,\ldots\right)}}
  \frac{\prod_{m=1}^\infty \left(p_m(\xx)\right)^{k_m} \prod_{m=1}^\infty \left(p_m(\yy)\right)^{k_m}}{\prod_{m=1}^\infty \left(k_m! m^{k_m}\right)} \\
&= \sum\limits_{\substack{ \text{weak compositions } \\ \left(k_1,k_2,k_3,\ldots\right)}}
  \frac{p_{\left(1^{k_1}2^{k_2}3^{k_3}\cdots\right)}\left(\xx\right) p_{\left(1^{k_1}2^{k_2}3^{k_3}\cdots\right)}\left(\yy\right)}{z_{\left(1^{k_1}2^{k_2}3^{k_3}\cdots\right)}}
=\sum_{\lambda \in \Par} \frac{ p_\lambda(\xx) p_\lambda(\yy) }{z_\lambda}
\end{align*}
due to the fact that every partition can be uniquely written in the
form $\left(1^{k_1}2^{k_2}3^{k_3}\cdots\right)$ with
$\left(k_1, k_2, k_3, \ldots\right)$ a weak composition.
\begin{commentedout}
[\textbf{Note:} The following alternative proof has been commented
out since it uses the multinomial formula for infinite sums,
which might be a bit less elementary than what the above proof
does.]
\begin{align*}
&\prod_{i,j=1}^\infty (1-x_i y_j)^{-1}
=\exp\left( \sum_{m=1}^\infty \frac{1}{m} p_m(\xx) p_m(\yy) \right)
=\sum_{k=0}^\infty \frac{1}{k!}
      \left( \sum_{m=1}^\infty \frac{1}{m} p_m(\xx) p_m(\yy) \right)^k\\
&=\sum_{k=0}^\infty \frac{1}{k!}
           \sum\limits_{\substack{(m_1,m_2,\ldots) \in \NN^\infty:\\ m_1+m_2+\cdots=k}}
             \binom{k}{m_1,m_2,\ldots}
               \left( \frac{p_1(\xx) p_1(\yy)}{1} \right)^{m_1}
               \left( \frac{p_2(\xx) p_2(\yy)}{2} \right)^{m_2} \cdots \\
& \qquad \qquad \left(\text{by the multinomial formula}\right) \\
&=\sum\limits_{\substack{ \text{weak }\\\text{ compositions } \\ (m_1,m_2,\ldots)}}
             \frac{ \left( p_1(\xx) p_1(\yy) \right)^{m_1} }{ 1^{m_1} m_1! }
             \cdot \frac{ \left( p_2(\xx) p_2(\yy)\right)^{m_2} }{ 2^{m_2} m_2! }
             \cdots
=\sum_{\lambda \in \Par} \frac{ p_\lambda(\xx) p_\lambda(\yy) }{z_\lambda}.
\end{align*}
\end{commentedout}
\end{proof}

% [DG][v25] Added some more detail to the proof of the first expansion.

% [DG][v37] Added "by the multinomial formula".
% Actually, there is an alternative way to do the last computation,
% which proceeds by simplifying
% $\exp\left( \sum_{m=1}^\infty \frac{1}{m} p_m(\xx) p_m(\yy) \right)$
% to
% $\prod_{m=1}^\infty \exp\left( \frac{1}{m} p_m(\xx) p_m(\yy) \right)$
% and then expanding the $\exp$ into an infinite sum and expanding the
% product by the product rule. This avoids the multinomial formula.
% Do you also consider it simpler?

% [DG][v38] Proof replaced by the one avoiding the multinomial
% formula (and in much more detail).

\begin{corollary} \phantomsection
\label{dual-bases-h-m-p}
\begin{itemize}
\item[(a)] With respect to the Hall inner product on $\Lambda$, one also has
dual bases $\{h_\lambda\}$ and $\{m_\lambda\}$.
\item[(b)] If $\QQ$ is a subring of $\kk$, then $\{p_\lambda\}$ and
$\left\{z_\lambda^{-1} p_\lambda\right\}$ are also dual bases with
respect to the Hall inner product on $\Lambda$.
\item[(c)] If $\RR$ is a subring of $\kk$, then
$\left\{\dfrac{p_\lambda}{\sqrt{z_\lambda}} \right\}$ is an orthonormal
basis of $\Lambda$ with respect to the Hall inner product.
\end{itemize}
\end{corollary}

% [DG][v7] I added the "in $\Lambda_\RR$" part, since it is not defined
% over $\ZZ$. IMHO, $\RR$ and the square roots are a red herring here,
% but I trust the reader to translate the claim into a rational form
% if needed...
% 
% I also made the fraction a \dfrac, because an inline-sized
% $\frac{p_\lambda}{\sqrt{z_\lambda}}$ didn't look nice. If you disagree,
% just change it back.

% [DG][v25] Split the corollary into (a), (b) and (c), adding part (b) to
% avoid unnecessary use of reals.

\begin{proof}
Since \eqref{Cauchy-identity} and \eqref{Two-other-Cauchy-expansions}
showed
\[
\prod_{i,j=1}^\infty (1-x_i y_j)^{-1}
=\sum_{\lambda \in \Par} s_{\lambda}(\xx) s_\lambda(\yy)\\
=\sum_{\lambda \in \Par} h_{\lambda}(\xx) m_\lambda(\yy)\\
=\sum_{\lambda \in \Par} p_{\lambda}(\xx)
                       z_\lambda^{-1} p_{\lambda}(\yy)
=\sum_{\lambda \in \Par} \frac{p_{\lambda}(\xx)}{\sqrt{z_\lambda}}
                       \frac{p_{\lambda}(\yy)}{\sqrt{z_\lambda}} ,
\]
it suffices to show that any pair of graded
bases\footnote{See Definition~\ref{def.graded-basis} for the
concept of a ``graded basis'', and recall our convention that
a graded basis of $\Lambda$
is tacitly assumed to have its indexing set $\Par$ partitioned
into $\Par_0, \Par_1, \Par_2, \ldots$.
Thus, a graded basis of $\Lambda$ means a basis
$\left\{w_\lambda\right\}_{\lambda \in \Par}$ of the $\kk$-module
$\Lambda$ (indexed by the partitions $\lambda \in \Par$)
with the property that, for every $n \in \NN$, the subfamily
$\left\{w_\lambda\right\}_{\lambda \in \Par_n}$ is a basis
of the $\kk$-module $\Lambda_n$.}
 $\{u_\lambda\},\{v_\lambda\}$ of $\Lambda$ having
\[
\sum_{\lambda \in \Par} s_{\lambda}(\xx) s_\lambda(\yy)=
\sum_{\lambda \in \Par} u_{\lambda}(\xx) v_\lambda(\yy)
\]
will be dual with respect to $(\cdot,\cdot)$.
To show this, consider such a pair of graded bases.
Write transition matrices
$A = \left(a_{\nu,\lambda}\right)_{\left(\nu, \lambda\right) \in \Par \times \Par}$
and
$B = \left(b_{\nu,\lambda}\right)_{\left(\nu, \lambda\right) \in \Par \times \Par}$
uniquely expressing
\begin{align}
u_\lambda&=\sum_\nu a_{\nu,\lambda} s_\nu, \label{eq.pf.dual-bases-h-m-p.a} \\
v_\lambda&=\sum_\nu b_{\nu,\lambda} s_\nu. \label{eq.pf.dual-bases-h-m-p.b}
\end{align}
Recall that $\Par = \bigsqcup_{r \in \NN} \Par_r$. Hence, we
can view $A$ as a block matrix, where the blocks are indexed by
pairs of nonnegative integers, and the $\left(r, s\right)$-th block
is $\left(a_{\nu, \lambda}\right)_{\left(\nu, \lambda\right) \in
\Par_r \times \Par_s}$.
For reasons of homogeneity\footnote{More precisely: The power series
$u_\lambda$ is homogeneous of degree $\left|\lambda\right|$, and the
power series $s_\nu$ is homogeneous of degree $\left|\nu\right|$.},
we have $a_{\nu, \lambda} = 0$ for any
$\left(\nu,\lambda\right) \in \Par^2$ satisfying
$\left|\nu\right| \neq \left|\lambda\right|$. Therefore, the
$\left(r, s\right)$-th block of $A$ is zero whenever $r \neq s$.
In other words, the block
matrix $A$ is block-diagonal. Similarly, $B$ can be viewed as a
block-diagonal matrix. The diagonal blocks of $A$ and $B$ are finite
square matrices (since $\Par_r$ is a finite set for each $r \in \NN$);
therefore, products such as $A^t B$, $B^t A$ and $AB^t$ are
well-defined (since all sums involved in their definition have only
finitely many nonzero addends)
and subject to the law of associativity. Moreover, the
matrix $A$ is invertible (being a transition matrix between two bases),
and its inverse is again block-diagonal (because $A$ is
block-diagonal).

The equalities \eqref{eq.pf.dual-bases-h-m-p.a} and
\eqref{eq.pf.dual-bases-h-m-p.b} show that
$(u_\alpha,v_\beta) = \sum_{\nu} a_{\nu,\alpha} b_{\nu,\beta}$
(by the orthonormality of the $s_\lambda$).
Hence, we want to prove that
$\sum_{\nu} a_{\nu,\alpha} b_{\nu,\beta} = \delta_{\alpha, \beta}$.
In other words, we want to prove that
$A^t B = I$, that is, $B^{-1}=A^t$.
  On the other hand, one has
\[
\sum_{\lambda } s_{\lambda}(\xx) s_\lambda(\yy)
=\sum_{\lambda } u_{\lambda}(\xx) v_\lambda(\yy)
=\sum_\lambda \sum_\nu a_{\nu,\lambda} s_\nu(\xx)
  \sum_\rho b_{\rho,\lambda} s_\rho(\yy).
\]
Comparing coefficients\footnote{Comparing coefficients is legitimate
because if a power series $f \in \kk\left[\left[\xx, \yy\right]\right]$
is written
in the form $f = \sum_{\left(\nu, \rho\right) \in \Par^2}
a_{\rho, \nu} s_\nu\left(\xx\right)
s_\rho\left(\yy\right)$ for some coefficients
$a_{\rho, \nu} \in \kk$, then these coefficients
$a_{\rho, \nu}$ are uniquely determined by $f$. This is
just a restatement of Exercise \ref{exe.Lambda.infinite-indep}(b).}
of $s_\nu(\xx) s_\rho(\yy)$ forces
$\sum_\lambda a_{\nu,\lambda} b_{\rho,\lambda} = \delta_{\nu,\rho}$,
or in other words, $AB^t=I$. Since $A$ is invertible, this
yields $B^tA = I$, and hence $A^tB=I$, as desired.\footnote{In our
argument above, we have obtained the invertibility of $A$ from the
fact that $A$ is a transition matrix between two bases. Here is an
alternative way to prove that $A$ is invertible: \par
Recall that $A$ and $B^t$ are block-diagonal matrices. Hence,
the equality $AB^t = I$ rewrites as
$A_{r,r} \left(B^t\right)_{r,r} = I$ for all $r \in \NN$, where we are
using the notation $C_{r,s}$ for the $\left(r, s\right)$-th block of a
block matrix $C$.
But this shows that each diagonal block $A_{r, r}$ of $A$ is
right-invertible. Therefore, each diagonal block $A_{r, r}$ of $A$ is
invertible (because $A_{r, r}$ is a square matrix of finite size,
and such matrices are always invertible when they are
right-invertible). Consequently, the block-diagonal matrix $A$ is
invertible, and its inverse is again a block-diagonal matrix (whose
diagonal blocks are the inverses of the $A_{r,r}$).}
\end{proof}

% [DG][v17] Added footnote.

% [DG][v23] Added clarification that the two bases should be graded;
% also, changed the use of finiteness to a use of the invertibility of
% $B$ (which is presumed known since $\{u_\lambda\}$ and $\{v_\lambda\}$
% are assumed to be bases). The old argument is now in the footnote.

% [DG][v58] Expanded footnote, and added more details. The concept of
% inverses of infinite matrices is somewhat slippery: Infinite
% matrices are not in 1-to-1 correspondence with homomorphisms between
% modules (because a matrix can only correspond to a homomorphism if
% each of its columns has only finitely many nonzero entries); thus,
% we would have to distinguish between invertibility-as-morphism and
% invertibility-as-matrix, and barring the former to make sure that
% associativity holds. I have introduced the block structure early on
% to get rid of these problems.

Corollary~\ref{dual-bases-h-m-p} is a known and fundamental
fact\footnote{For example, Corollary~\ref{dual-bases-h-m-p}(a)
appears in \cite[Corollary 3.3]{Leeuwen-altSchur} (though the
definition of Schur functions in \cite{Leeuwen-altSchur} is different
from ours; we will meet this alternative definition later on),
and parts (b) and (c) of Corollary~\ref{dual-bases-h-m-p} are
equivalent to \cite[\S I.4, (4.7)]{Macdonald} (though Macdonald
defines the Hall inner product using
Corollary~\ref{dual-bases-h-m-p}(a)).}.
However, our definition of the Hall inner product is unusual;
most authors (e.g., Macdonald in \cite[\S I.4, (4.5)]{Macdonald},
Hazewinkel/Gubareni/Kirichenko in
\cite[Def. 4.1.21]{HazewinkelGubareniKirichenko}, and Stanley in
\cite[(7.30)]{Stanley}) \emph{define} the Hall inner product
as the bilinear form satisfying
$\left(h_\lambda, m_\mu\right) = \delta_{\lambda, \mu}$ (or,
alternatively,
$\left(m_\lambda, h_\mu\right) = \delta_{\lambda, \mu}$), and
only later prove that the basis $\left\{s_\lambda\right\}$ is
orthonormal with respect to this scalar product. (Of course,
the fact that this definition is equivalent to our
Definition~\ref{def.hall-inner-product} follows either from
this orthonormality, or from our
Corollary~\ref{dual-bases-h-m-p}(a).)

% [DG][v45] Added preceding paragraph (closing a gap in the van Leeuwen
% paper and connecting the notes with some older literature).

The tactic applied in the proof of Corollary~\ref{dual-bases-h-m-p} can
not only be used to show that certain bases of $\Lambda$ are dual, but
also, with a little help from linear algebra over rings
(Exercise~\ref{exe.linalg.sur-bij}), it can be strengthened to show that
certain families of symmetric functions are bases to begin with, as we
will see in Exercise~\ref{exe.Lambda.basis-vlam} and
Exercise~\ref{exe.Lambda.bases}.

\begin{exercise} \phantomsection
\label{exe.linalg.sur-bij}
\begin{itemize}
\item[(a)] Prove that if an endomorphism of a finitely generated $\kk$-module
is surjective, then this endomorphism is a $\kk$-module isomorphism.
\item[(b)] Let $A$ be a finite free $\kk$-module with finite
basis $\left(  \gamma_{i}\right)  _{i\in I}$. Let $\left(  \beta_{i}\right)
_{i\in I}$ be a family of elements of $A$ which spans the $\kk$-module
$A$. Prove that $\left(  \beta_{i}\right)  _{i\in I}$ is a $\kk$-basis of
$A$.
\end{itemize}
\end{exercise}

\begin{exercise}
\label{exe.Lambda.basis-vlam}
For each partition $\lambda$, let $v_{\lambda}$
be an element of $\Lambda_{\left\vert \lambda\right\vert }$. Assume that the
family $\left(  v_{\lambda}\right)  _{\lambda \in \Par}$ spans
the $\kk$-module $\Lambda$.
Prove that the family $\left(  v_{\lambda
}\right)  _{\lambda \in \Par}$ is a graded basis of the
graded $\kk$-module $\Lambda$.
\end{exercise}

\begin{exercise} \phantomsection
\label{exe.Lambda.bases}
\begin{itemize}
\item[(a)] Assume that for every partition $\lambda$, two homogeneous elements
$u_\lambda$ and $v_\lambda$ of $\Lambda$, both having degree
$\left|\lambda\right|$, are given. Assume further that
\[
\sum_{\lambda \in \Par} s_\lambda\left(\xx\right) s_\lambda\left(\yy\right)
= \sum_{\lambda \in \Par} u_\lambda\left(\xx\right) v_\lambda\left(\yy\right)
\]
in $\kk\left[\left[\xx, \yy\right]\right]
= \kk\left[\left[x_1, x_2, x_3, \ldots, y_1, y_2, y_3, \ldots\right]\right]$.
Show that $\left(u_\lambda\right)_{\lambda \in \Par}$ and
$\left(v_\lambda\right)_{\lambda \in \Par}$ are $\kk$-bases of
$\Lambda$, and actually are dual bases with respect to the
Hall inner product on $\Lambda$.
\item[(b)] Use this to give a new proof of the fact that
$\left(h_\lambda\right)_{\lambda \in \Par}$ is a $\kk$-basis of
$\Lambda$.
\end{itemize}
\end{exercise}

% [DG][v23] Added the above two exercises, to factor out some steps in
% other solutions (e.g., of Exercise \ref{exe.witt}(h)).

\begin{exercise}
\label{exe.E(t)-H(t)-P(t)-relation}
Prove that
$\sum_{m \geq 0} p_{m+1} t^m = \dfrac{H'\left(t\right)} {H\left(t\right)}$.
(This was proven in \eqref{E(t)-H(t)-P(t)-relation} in the case when
$\QQ$ is a subring of $\kk$, but here we make no requirements on $\kk$.)
\end{exercise}

% [DG][v25] Added above exercise (as it was subtly -- sufficiently so
% that I missed it myself -- used in another's solution).

The following exercises give some useful criteria for algebraic independence
of families of symmetric functions:

\begin{exercise}
\label{exe.Lambda.alggen-ind}
Let $v_1, v_2, v_3, \ldots$ be elements of $\Lambda$.
Assume that $v_n \in \Lambda_n$ for each positive integer $n$.
Assume further that $v_1, v_2, v_3, \ldots$ generate the
$\kk$-algebra $\Lambda$.
Then:

\begin{enumerate}
\item[(a)] Prove that $v_1, v_2, v_3, \ldots$ are algebraically independent
over $\kk$.

\item[(b)] For every partition $\lambda$, define an element $v_{\lambda}
\in\Lambda$ by $v_{\lambda}=v_{\lambda_{1}}v_{\lambda_{2}}\cdots
v_{\lambda_{\ell\left(  \lambda\right)  }}$. Prove that the family $\left(
v_{\lambda}\right)  _{\lambda\in\Par}$ is a graded basis of the
graded $\kk$-module $\Lambda$.
\end{enumerate}
\end{exercise}

\begin{exercise}
\label{exe.Lambda.alggen-crit1}
For each partition $\lambda$, let $a_{\lambda} \in \kk$.
Assume that the element $a_{\left(  n\right)  } \in \kk$
is invertible for each positive integer $n$.
Let $v_1, v_2, v_3, \ldots$ be elements of $\Lambda$ such that
each positive integer $n$ satisfies
$v_n = \sum_{\lambda \in \Par_n} a_{\lambda} h_{\lambda}$.
Prove that the elements $v_1, v_2, v_3, \ldots$ generate the
$\kk$-algebra $\Lambda$ and are algebraically independent over
$\kk$.
\end{exercise}

\begin{exercise}
\label{exe.Lambda.alggen-crit1hall}
Let $v_1, v_2, v_3, \ldots$ be elements of $\Lambda$.
Assume that $v_n \in \Lambda_n$ for each positive integer $n$.
Assume further that $\left( p_n, v_n \right) \in \kk$ is invertible
for each positive integer $n$.
Prove that the elements $v_1, v_2, v_3, \ldots$ generate the
$\kk$-algebra $\Lambda$ and are
algebraically independent over $\kk$.
\end{exercise}

% [DG][v75] Added the above three exercises, as well as
% Exercise \ref{exe.Lambda.basis-vlam} further above.

\begin{exercise}
\label{exe.Lambda.hall-coefficient}
Let $f \in \Lambda$, and let $\beta$ be a weak composition.
Let $\mu \in \Par$ be the partition consisting of the nonzero
entries of $\beta$ (sorted in decreasing order).%
\footnote{For example, if $\beta=\left(  1,0,3,1,2,3,0,0,0,\ldots
\right)  $, then $\mu=\left(  3,3,2,1,1\right)  $.}
Prove that
\[
\left( f, h_\mu \right)
= \left( h_\mu, f \right)
= \left( \text{the coefficient of $\xx^{\beta}$ in $f$} \right)  .
\]

\end{exercise}

% [DG][v78] Added the above simple exercise, in order to cite it
% in a later proof.

\begin{exercise}
\label{exe.Lambda.h-e-through-p}
Assume that $\QQ$ is a subring of $\kk$.
Define a positive integer $z_{\lambda}$ for each $\lambda \in \Par$
as in Proposition~\ref{prop-Two-other-Cauchy-expansions}.
Prove that every $n \in \NN$ satisfies the two equalities
\begin{equation}
h_n = \sum_{\lambda \in \Par_n} z_{\lambda}^{-1} p_{\lambda}
\label{eq.exe.Lambda.h-e-through-p.h}
\end{equation}
and
\begin{equation}
e_n = \sum_{\lambda \in \Par_n}
\left(  -1\right) ^{\left\vert \lambda\right\vert -\ell\left(  \lambda\right)  }
z_{\lambda}^{-1} p_{\lambda}.
\label{eq.exe.Lambda.h-e-through-p.e}
\end{equation}
\end{exercise}

% [DG][v79] Added the above exercise out of completionism.

\subsection{Bialternants, Littlewood-Richardson:  Stembridge's concise proof}
\label{bialternant-section}

There is a more natural way in which Schur functions arise as a $\kk$-basis
for $\Lambda$, coming from consideration of polynomials in a finite variable set,
and the relation between those which are symmetric and those which are \emph{alternating}.

For the remainder of this section,
fix a nonnegative integer $n$, and let $\xx=(x_1,\ldots,x_n)$ be a finite variable set.
This means that $s_{\lambda/\mu}=s_{\lambda/\mu}(\xx)=\sum_T \xx^{\cont(T)}$ is a
generating function for column-strict tableaux $T$ as in
Definition~\ref{skew-Schur-function-definition}, but with the extra condition
that $T$ have entries in $\{1,2,\ldots,n\}$. \ \ \ \ 
\footnote{See Exercise~\ref{exe.schur.finite}(a) for this.}
As a consequence, $s_{\lambda / \mu}$ is a polynomial in
$\kk \left[x_1, x_2, \ldots, x_n\right]$ (not just a
power series), since there are only finitely many column-strict
tableaux $T$ of shape $\lambda / \mu$ having all their entries in
$\left\{ 1, 2, \ldots , n \right\}$.
We will assume without
further mention that all partitions appearing in the section
have at most $n$ parts.

% [DG][v44] Added footnote and sentence about polynomiality.

\begin{definition}
Let $\kk$ be the ring $\ZZ$ or a field of characteristic not equal to $2$.
(We require this to avoid certain annoyances in the discussion of
alternating polynomials in characteristic $2$.)

Say that a polynomial $f(\xx)=f(x_1,\ldots,x_n)$ is
\emph{alternating}\index{alternating polynomial} if
for every permutation $w$ in $\Symm_n$ one has that
\[
(wf)(\xx) = f(x_{w(1)},\ldots,x_{w(n)}) = \sgn(w) f(\xx).
\]
Let $\Lambda^{\sgn} \subset \kk[x_1,\ldots,x_n]$
denote the subset of alternating polynomials\footnote{When $\kk$ has characteristic $2$
(or, more generally, is an arbitrary commutative ring),
it is probably best to define the alternating polynomials
$\Lambda^{\sgn}_\kk$ as the $\kk$-submodule
$\Lambda^{\sgn} \otimes_{\ZZ} \kk$
of $\ZZ[x_1,\ldots,x_n] \otimes_{\ZZ} \kk  \cong \kk[x_1,\ldots,x_n]$.}.
\end{definition}

As with $\Lambda$ and its monomial basis $\{m_\lambda\}$, there is
an obvious $\kk$-basis for $\Lambda^{\sgn}$, coming from the fact that
a polynomial $f=\sum_\alpha c_\alpha \xx^\alpha$
is alternating if and only if $c_{w(\alpha)} = \sgn(w) c_\alpha$ for every $w$ in $\Symm_n$
and every $\alpha \in \NN^n$.  This means that every alternating $f$
is a $\kk$-linear combination of the following elements.

% [DG][v25] Replaced "weak composition $\alpha$" by "$\alpha \in \NN^n$"
% to avoid having to explain that weak compositions are assumed to
% have length $\leq n$ here.

\begin{definition}
\label{def.alternant}
For $\alpha=(\alpha_1,\ldots,\alpha_n)$ in $\NN^n$, define
the \dfn{alternant}
\[
a_{\alpha}:= \sum_{ w \in \Symm_n} \sgn(w) w(\xx^{\alpha})
= \det \left[ \begin{matrix}
x_1^{\alpha_1} & \cdots & x_1^{\alpha_n} \\
x_2^{\alpha_1} & \cdots & x_2^{\alpha_n} \\
\vdots & \ddots & \vdots \\
x_n^{\alpha_1} & \cdots & x_n^{\alpha_n}
\end{matrix}
\right].
\]
\end{definition}

\begin{example}
One has
\[
a_{(1,5,0)} = x_1^1 x_2^5 x_3^0
               -  x_1^5 x_2^1 x_3^0 -  x_1^1 x_2^0 x_3^5 -  x_1^0 x_2^5 x_3^1
               + x_1^0 x_2^1 x_3^5 + x_1^5 x_2^0 x_3^1
            = - a_{(5,1,0)}.
\]
Similarly, $a_{w\left(\alpha\right)} = \sgn(w) a_\alpha$ for every
$w \in \Symm_n$ and every $\alpha \in \NN^n$.

Meanwhile, $a_{(5,2,2)} = 0$ since the transposition $t=\binom{123}{132}$ fixes
$(5,2,2)$ and hence
\[
a_{(5,2,2)} = t(a_{(5,2,2)}) = \sgn(t) a_{(5,2,2)} = -a_{(5,2,2)}.
\]
\footnote{One subtlety should be addressed:
We want to prove that $a_{(5,2,2)} = 0$ in $\kk\left[x_1,\ldots,x_n\right]$
for every commutative ring $\kk$. It is clearly enough to prove that
$a_{(5,2,2)} = 0$ in $\ZZ\left[x_1,\ldots,x_n\right]$. Since $2$ is not a
zero-divisor in $\ZZ\left[x_1,\ldots,x_n\right]$, we can achieve this by showing that
$a_{(5,2,2)} = -a_{(5,2,2)}$. We would not be able to make this argument
directly over an arbitrary commutative ring $\kk$.}
Alternatively, $a_{(5,2,2)}=0$ as it is a determinant of a matrix with two equal columns.
Similarly, $a_\alpha = 0$ for every $n$-tuple $\alpha \in \NN^n$ having two
equal entries.
\end{example}

\noindent
This example illustrates that, for a $\kk$-basis for  $\Lambda^{\sgn}$,
one can restrict attention to alternants $a_\alpha$ in which $\alpha$ is a
\dfn{strict partition}, i.e., in which $\alpha$ satisfies $\alpha_1 > \alpha_2 > \cdots > \alpha_n$.
One can therefore uniquely express $\alpha = \lambda + \rho$,
where $\lambda$ is a (weak) partition
$\lambda_1 \geq \cdots \geq \lambda_n \geq 0$ and where
$
\rho:=(n-1,n-2,\ldots,2,1,0)
$
is sometimes called the \dfn{staircase partition}\footnote{The
name is owed to its Ferrers shape. For instance, if $n = 5$, then
the Ferrers diagram of $\rho$ (represented using dots) has
the form
\[
\begin{matrix}
\bullet & \bullet & \bullet & \bullet \\
\bullet & \bullet & \bullet & \\
\bullet & \bullet & & \\
\bullet & & & 
\end{matrix}
.
\]
}.
For example $\alpha=(5,1,0)=(3,0,0) + (2,1,0)=\lambda+\rho$.

% [DG][v32] Made a footnote out of your brief remark about staircase
% shapes.

\begin{proposition}
\label{prop.alternants.space}
Let $\kk$ be the ring $\ZZ$ or a field of characteristic not equal to $2$.

The alternants $\{ a_{\lambda +\rho} \}$ as $\lambda$ runs through the partitions
with at most $n$ parts form a $\kk$-basis for $\Lambda^{\sgn}$.
In addition, the \emph{bialternants} $\{ \frac{ a_{\lambda +\rho} }{ a_\rho } \}$
as $\lambda$ runs through the same set form a $\kk$-basis for
$\Lambda(x_1,\ldots,x_n)=\kk[x_1,\ldots,x_n]^{\Symm_n}$.
\end{proposition}
\begin{proof}
The first assertion should be clear from our previous discussion:  the
alternants $\{ a_{\lambda +\rho} \}$ span $\Lambda^{\sgn}$ by definition, and
they are $\kk$-linearly independent because they are supported on disjoint
sets of monomials $\xx^\alpha$.

The second assertion follows from the first,
after proving the following \textbf{Claim}: $f(\xx)$ lies in
$\Lambda^{\sgn}$ if and only if $f(\xx)=a_\rho \cdot g(\xx)$ where
$g(\xx)$ lies in $\kk[\xx]^{\Symm_n}$ and where
\[
a_{\rho} = \det( x_i^{n-j} )_{i,j=1,2,\ldots,n} =\prod_{1 \leq i < j \leq n} (x_i - x_j)
\]
is the \dfn{Vandermonde determinant/product}.
In other words,
\[
\Lambda^{\sgn}= a_{\rho}  \cdot \kk[\xx]^{\Symm_n}
\]
is a free $\kk[\xx]^{\Symm_n}$-module of rank $1$, with $a_\rho$ as its
$\kk[\xx]^{\Symm_n}$-basis element.

To see the Claim, first note the inclusion
\[
\Lambda^{\sgn} \supset a_{\rho} \cdot \kk[\xx]^{\Symm_n}
\]
since the product of a symmetric polynomial and an alternating polynomial is
an alternating polynomial.  For the reverse inclusion, note that since an alternating polynomial $f(\xx)$ changes
sign whenever one exchanges two distinct variables $x_i,x_j$, it must vanish upon setting
$x_i = x_j$, and therefore be divisible by $x_i - x_j$, so divisible by
the entire product $\prod_{1 \leq i < j \leq n} (x_i - x_j)=a_\rho$.  But then
the quotient $g(\xx) = \frac{f(\xx)}{a_\rho}$ is symmetric, as
it is a quotient of two alternating polynomials.
\end{proof}

% [DG][v3] Replaced $\supseteq$ by $\supset$ to match notation in the rest
% of the notes. (Also one more time further below.)

Let us now return to the general setting, where $\kk$ is an arbitrary
commutative ring. We are not requiring that the assumptions
of Proposition~\ref{prop.alternants.space} be valid; we can still study
the $a_\alpha$ of Definition~\ref{def.alternant}, but we cannot use
Proposition~\ref{prop.alternants.space} anymore.
We will show that the fraction
$\frac{ a_{\lambda +\rho} }{ a_\rho }$ is nevertheless a well-defined
polynomial in $\Lambda\left(x_1, \ldots, x_n\right)$ whenever $\lambda$
is a partition\footnote{This
can also be deduced by base change from the $\kk = \ZZ$ case of
Proposition~\ref{prop.alternants.space}.}, and in fact equals the
Schur function $s_\lambda (\xx)$. As a consequence, the mysterious
bialternant basis $\{ \frac{ a_{\lambda +\rho} }{ a_\rho } \}$
of $\Lambda\left(x_1, \ldots, x_n\right)$ defined in
Proposition~\ref{prop.alternants.space} still exists in the general
setting, and is plainly the Schur functions $\{ s_\lambda(\xx) \}$.
Stembridge \cite{Stembridge} noted
that one could give a remarkably concise proof of an even stronger
assertion, which simultaneously gives one of the standard combinatorial
interpretations for the Littlewood-Richardson coefficients $c^{\lambda}_{\mu,\nu}$.
For the purposes of stating it,
we introduce for a tableau $T$ the notation $T|_{\cols \geq j}$ (resp. $T|_{\cols \leq j}$)
to indicate the subtableau which is the restriction
of $T$ to the union of its columns $j,j+1,j+2,\ldots$ (resp. columns $1,2,\ldots,j$).

\begin{example}
If $T = \begin{matrix} & & 1 & 2 \\ & 2 & 2 & 3 \\ 3 & 5 \end{matrix}$,
then
\[
T|_{\cols \geq 3} = \begin{matrix} 1 & 2 \\ 2 & 3 \end{matrix}
\qquad \text{and} \qquad
T|_{\cols \leq 2} = \begin{matrix} & \\ & 2 \\ 3 & 5 \end{matrix}
\]
(note that $T|_{\cols \leq 2}$ has an empty first row).
\end{example}

% [DG][v77] Added the above example.

\begin{theorem}
\label{Stembridge-theorem}
For partitions $\lambda,\mu,\nu$ with $\mu \subseteq \lambda$, one has\footnote{Again,
we can drop the requirement that $\mu \subseteq \lambda$, provided that we understand
that there are no column-strict tableaux of shape $\lambda / \mu$ unless $\mu
\subseteq \lambda$.}
\[
a_{\nu+\rho} s_{\lambda/\mu} = \sum_T a_{\nu+\cont(T)+\rho}
\]
where $T$ runs through all column-strict tableaux with entries in $\{1,2,\ldots,n\}$
of shape $\lambda/\mu$ with
the property that for each $j=1,2,3,\ldots$, the weak composition
$\nu + \cont(T|_{\cols\geq j})$ is a partition.
\end{theorem}

% [DG][v77] Replaced "one has $\nu + \cont(T|_{\cols\geq j})$ a partition"
% by "the weak composition $\nu + \cont(T|_{\cols\geq j})$ is a partition".

\noindent
Before proving Theorem~\ref{Stembridge-theorem}, let us see some of its consequences.

\begin{corollary}
\label{cor.stembridge.slambda}
For any partition $\lambda$, we have\footnote{Notice that division
by $a_\rho$ is unambiguous in the ring $\kk\left[x_1,\ldots,x_n\right]$,
since $a_\rho$ is not a zero-divisor (in fact,
$a_\rho = \prod_{1 \leq i < j \leq n} (x_i - x_j)$ is the product of the
binomials $x_i - x_j$, none of which is a zero-divisor). \vspace{0.2pc}}
\[
s_\lambda(\xx) = \frac{ a_{\lambda +\rho} }{ a_\rho }.
\]
\end{corollary}
\begin{proof}
Fix a partition $\lambda$.
Take $\nu=\mu=\varnothing$ in Theorem~\ref{Stembridge-theorem}.  Note
that there is only one column-strict tableau $T$
of shape $\lambda$ such that each $\cont(T|_{\cols\geq j})$ is a partition,
namely the tableau having every entry in row $i$ equal to $i$:
\[
\begin{matrix}
1 & 1 & 1 & 1 & 1\\
2 & 2 & 2 &   & \\
3 & 3 & 3 &   & \\
4 & 4
\end{matrix}
\]
\footnote{\textit{Proof.} It is clear that the tableau
having every entry in row $i$ equal to $i$ indeed satisfies the condition
that each $\cont(T|_{\cols\geq j})$ is a partition.
It remains to show that it is the only column-strict tableau
(of shape $\lambda$) satisfying this condition.

Let $T$ be a column-strict tableau of shape $\lambda$ satisfying
the condition that each $\cont(T|_{\cols\geq j})$ is a partition.
We must show that for each $i$, every entry in row $i$ of $T$ is equal
to $i$.
Assume the contrary.
Thus, there exists some $i$ such that row $i$ of $T$ contains an entry
distinct from $i$.
Consider the smallest such $i$.
Hence, rows $1, 2, \ldots, i-1$ of $T$ are filled with entries
$1, 2, \ldots, i-1$, whereas row $i$ has some entry distinct from $i$.
Choose some $j$ such that the $j$-th entry of row $i$ of $T$ is
distinct from $i$.
This entry cannot be smaller than $i$ (since it has $i-1$ entries
above it in its column, and the entries of $T$ increase strictly
down columns); thus, it has to be larger than $i$.
Therefore, all entries in rows $i, i+1, i+2, \ldots$ of $T|_{\cols\geq j}$
are larger than $i$ as well (since they lie southeast of this entry).
Hence, each entry of $T|_{\cols\geq j}$ is either smaller than $i$
(if it is in one of rows $1, 2, \ldots, i-1$)
or larger than $i$ (if it is in row $i$ or further down).
Thus, $i$ is not an entry of $T|_{\cols\geq j}$.
In other words, $\cont_i (T|_{\cols\geq j}) = 0$.
Since $\cont (T|_{\cols\geq j})$ is a partition, we thus conclude that
$\cont_k (T|_{\cols\geq j}) = 0$ for all $k > i$.
In other words, $T|_{\cols\geq j}$ has no entries larger than $i$.
But this contradicts the fact that the $j$-th entry of row $i$ of $T$
is larger than $i$.
This contradiction completes our proof.}.
Furthermore, this $T$ has $\cont(T)=\lambda$, so the theorem says $a_{\rho} s_\lambda= a_{\lambda+\rho}$.
\end{proof}

% [DG][v77] Added footnote to the above proof.

\begin{example}
For $n=2$, so that $\rho=(1,0)$, if we take $\lambda=(4,2)$, then one
has
\begin{align*}
\frac{a_{\lambda+\rho}}{a_\rho}
&=\frac{a_{(4,2)+(1,0)}}{a_{(1,0)}}
=\frac{a_{(5,2)}}{a_{(1,0)}} \\
&=\frac{x_1^5 x_2^2 - x_1^2 x_2^5}{x_1 - x_2}\\
&= x_1^4 x_2^2 + x_1^3 x_2^3 + x_1^2 x_2^4\\
&= \xx^{\cont{
\left(\begin{matrix}
1111\\
22\phantom{1}\phantom{1}
\end{matrix}\right)
}}
+
\xx^{\cont{
\left(\begin{matrix}
1112\\
22\phantom{1}\phantom{1}
\end{matrix}\right)
}}
+
\xx^{\cont{
\left(\begin{matrix}
1122\\
22\phantom{1}\phantom{1}
\end{matrix}\right)
}}\\
&=s_{(4,2)}
=s_\lambda.
\end{align*}
\end{example}

Some authors use the equality in Corollary~\ref{cor.stembridge.slambda}
to \emph{define} the Schur polynomial
$s_\lambda \left(x_1, x_2, \ldots, x_n\right)$ in $n$ variables;
this definition, however, has the drawback of not generalizing
easily to infinitely many variables or to skew Schur
functions\footnote{With some effort, it is possible to use
Corollary~\ref{cor.stembridge.slambda} in order to define
the Schur function $s_\lambda$ in infinitely many variables.
Indeed, one can define this Schur function as the unique
element of $\Lambda$ whose evaluation at $\left(x_1, x_2,
\ldots, x_n\right)$ equals $\frac{ a_{\lambda +\rho} }{ a_\rho }$
for every $n \in \NN$. If one wants to use such a definition,
however, one needs to check that such an element exists. This is
the approach to defining $s_\lambda$ taken in
\cite[Definition 1.4.2]{Leeuwen-altSchur} and in
\cite[\S I.3]{Macdonald}.}.

% [DG][v45] Added preceding paragraph. (Its purpose is to
% advertise the van Leeuwen paper and explain how it ties in
% with these notes.)

Next divide through by $a_\rho$ on both sides of Theorem~\ref{Stembridge-theorem}
(and use Corollary~\ref{cor.stembridge.slambda}) to give the following.
\begin{corollary}
\label{finite-L-R-rule}
For partitions $\lambda, \mu, \nu$ having at most $n$ parts, one has
\begin{equation}
\label{eq.finite-L-R-rule.1}
s_{\nu} s_{\lambda/\mu} = \sum_T s_{\nu+\cont(T)}
\end{equation}
where $T$ runs through the same set as in Theorem~\ref{Stembridge-theorem}.
In particular, taking $\nu=\varnothing$, we obtain
\begin{equation}
\label{eq.finite-L-R-rule.2}
s_{\lambda/\mu} = \sum_T s_{\cont(T)}
\end{equation}
where in the sum $T$ runs through all column-strict tableaux of shape $\lambda/\mu$
for which each $\cont(T|_{\cols\geq j})$ is a partition.

%Rephrasing, this,  whenever, $\lambda, \mu, \nu$ all have at most $n$ parts
%the Littlewood-Richardson coefficient $c^\lambda_{\mu,\nu}$ is the
%number of such tableaux $T$ having $\cont(T)=\nu$.
\end{corollary}

% [DG][v49] Moved the actual LR rule (Corollary~\ref{L-R-rule})
% further down (meaning after the proof of
% Theorem~\ref{Stembridge-theorem}) and split it into two corollaries,
% the first of which sets the number of variables to infinity and the
% second of which passes from equalities between functions to a formula
% for the coefficients. This adds detail to the proof and also creates
% a referenceable skew-LR rule (Corollary~\ref{infinite-L-R-rule}),
% which I will eventually reference when I do Zelevinsky's LR rule.

\begin{proof}[Proof of Theorem~\ref{Stembridge-theorem}]
Start by rewriting the left side of the theorem:
% , and using the fact that
% $w(s_{\lambda/\mu})=s_{\lambda/\mu}$ for any $w$ in $\Symm_n$:
\begin{align*}
a_{\nu+\rho} s_{\lambda/\mu}
&=\sum_{w \in \Symm_n} \sgn(w) \xx^{w(\nu+\rho)}
     s_{\lambda/\mu}
 =\sum_{w \in \Symm_n} \sgn(w) \xx^{w(\nu+\rho)}
     w(s_{\lambda/\mu}) \\
& \qquad \qquad
    \left(\text{since $w(s_{\lambda/\mu})=s_{\lambda/\mu}$ for any $w \in \Symm_n$}\right) \\
&=\sum_{w \in \Symm_n} \sgn(w) \xx^{w(\nu+\rho)}
     \sum\limits_{\substack{\text{column-strict }T\\\text{ of shape }\lambda/\mu}} \xx^{w(\cont(T))}  \\
&=\sum\limits_{\substack{\text{column-strict }T\\\text{ of shape }\lambda/\mu}}
         \qquad \sum_{w \in \Symm_n} \sgn(w) \xx^{w(\nu+\cont(T)+\rho)}  \\
&= \sum\limits_{\substack{\text{column-strict }T\\\text{ of shape }\lambda/\mu}}
      a_{\nu+\cont(T)+\rho} .
\end{align*}
We wish to cancel out all the summands indexed by column-strict tableaux $T$ which
fail any of the conditions that $\nu + \cont(T|_{\cols\geq j})$ be a partition.
Given such a $T$, find the maximal $j$ for which it fails this
condition\footnote{Such a $j$ exists because $\nu + \cont(T|_{\cols\geq j})$ is a
partition for all sufficiently high $j$ (in fact, $\nu$ itself is a partition).},
and
then find the minimal $k$ for which
\[
\nu_k+\cont_k(T|_{\cols\geq j}) <
\nu_{k+1}+\cont_{k+1}(T|_{\cols\geq j}).
\]
Maximality of $j$ forces
\[
\nu_k+\cont_k(T|_{\cols\geq j+1}) \geq
\nu_{k+1}+\cont_{k+1}(T|_{\cols\geq j+1}).
\]
Since column-strictness implies that column $j$ of $T$ can contain at most one
occurrence of $k$ or of $k+1$ (or neither or both),
the previous two inequalities imply that
column $j$ must contain an occurrence of $k+1$ and no occurrence of $k$, so
that
\[
\nu_k+\cont_k(T|_{\cols\geq j}) + 1
= \nu_{k+1}+\cont_{k+1}(T|_{\cols\geq j}).
\]
This implies that the adjacent transposition $t_{k,k+1}$ swapping $k$ and $k+1$
fixes the vector $\nu+\cont(T|_{\cols\geq j})+\rho$.

Now create a new tableau $T^*$ from $T$ by
applying the Bender-Knuth involution (from the proof of
Proposition~\ref{Schur-functions-are-symmetric-prop})
on letters $k,k+1$,
but \emph{only to columns} $1,2,\ldots,j-1$ of $T$,
leaving columns $j,j+1,j+2,\ldots$ unchanged.%
\footnote{See Example~\ref{exa.pf.Stembridge-theorem.invol}
below for an example of this construction.}
One should check that $T^*$ is still column-strict, but this holds
because column $j$ of $T$ has no occurrences of letter $k$.
Note that
\[
t_{k,k+1} \cont(T|_{\cols\leq j-1}) =
              \cont(T^*|_{\cols\leq j-1})
\]
and hence
\[
t_{k,k+1} ( \nu+\cont(T)+\rho ) = \nu+\cont(T^*)+\rho ,
\]
so that $a_{\nu+\cont(T)+\rho} = - a_{\nu+\cont(T^*)+\rho}$.

Because $T$ and $T^*$ have exactly the same columns  $j,j+1,j+2,\ldots$,
the tableau $T^*$ is also a violator of at least one of the conditions that
$\nu + \cont(T^*|_{\cols\geq j})$ be a partition, and
has the same choice of maximal $j$ and minimal $k$ as did $T$.
Hence the map $T \mapsto T^*$ is an involution on the violators
that lets one cancel their summands
$a_{\nu+\cont(T)+\rho}$ and $a_{\nu+\cont(T^*)+\rho}$ in pairs.%
\footnote{One remark is in order: The tableaux $T$ and $T^*$ may
be equal. In this case, the summands
$a_{\nu+\cont(T)+\rho}$ and $a_{\nu+\cont(T^*)+\rho}$ do not cancel,
as they are the same summand.
However, this summand is zero (because
$t_{k,k+1} ( \nu+\cont(T)+\rho )
= \nu+\cont\left(\underbrace{T^*}_{=T}\right)+\rho
= \nu+\cont(T)+\rho$ shows that the $n$-tuple $\nu+\cont(T)+\rho$
has two equal entries, and thus $a_{\nu+\cont(T)+\rho} = 0$),
and thus does not affect the sum.}
\end{proof}

% [DG][v17] I've done various changes to the preceding (bialternant)
% section. The case of $\kk$ being a commutative ring is now the
% new normal, whereas all requirements on $\kk$ to be $\ZZ$ or not
% to have characteristic $2$ are now explicit. The corollary which
% gives the Schur function as a quotient of two alternants is stated
% in the general case, and a footnote now serves to explain why
% the quotient is well-defined. In contrast, the proposition which
% states that the alternants are a basis for $\Lambda^\sgn$ is only
% stated in the restricted setting (as it was before, but now
% explicitly). The argument that $a_{(5,2,2)} = 0$ in the Example
% has gotten a footnote to explain why it works despite $2$ not
% being invertible.

% [DG][v19] Added footnote in the above proof about why a maximal $j$
% exists. I stumbled upon this subtlety when I tried to generalize
% the proof.

% [DG][v77] Added the footnote at the end of the proof, about the
% case $T^* = T$ (thanks Jang Soo Kim for the pointer).

\begin{example}
\label{exa.pf.Stembridge-theorem.invol}
Here is an example of the construction of $T^*$ in the above
proof.
Let $n = 6$ and
$\lambda = \left(5,4,4\right)$ and $\mu = \left(2,2\right)$
and $\nu = \left(1\right)$.
Let $T$ be the column-strict tableau
\[
\begin{matrix}
  &   & 1 & 2 & 2 \\
  &   & 2 & 3 \\
2 & 2 & 3 & 4
\end{matrix}
\quad \text{ of shape } \lambda / \mu .
\]
Then, % $\cont(T) = \left(1,5,2,1,0,0,0,\ldots\right)$ and
\begin{align*}
\cont(T|_{\cols \geq 5}) &= \left(0,1,0,0,0,\ldots\right)
\quad \text{(since $T|_{\cols \geq 5}$ has a single entry, which is $2$)}, \\
\text{so that} \quad
\nu + \cont(T|_{\cols \geq 5}) &= \left(1,1,0,0,0,\ldots\right)
\text{ is a partition.}
\end{align*}
But
\begin{align*}
\cont(T|_{\cols \geq 4}) &= \left(0,2,1,1,0,0,0,\ldots\right) , \\
\text{and thus} \quad
\nu + \cont(T|_{\cols \geq 4}) &= \left(1,2,1,1,0,0,0,\ldots\right)
\text{ is not a partition.}
\end{align*}
Thus, the $j$ in the above proof of Theorem~\ref{Stembridge-theorem}
is $4$.
Furthermore, the $k$ in the proof is $1$, since
$\nu_1 + \cont_1 (T|_{\cols \geq 4}) = 1 + 0
= 1 < 2 = 0 + 2 = \nu_2 + \cont_2 (T|_{\cols \geq 4})$.
Thus, $T^*$ is obtained from $T$ by applying the Bender-Knuth
involution on letters $1,2$ to columns $1,2,3$ only, leaving
columns $4,5$ unchanged. The result is
\[
T^* =
\begin{matrix}
  &   & 1 & 2 & 2 \\
  &   & 2 & 3 \\
1 & 1 & 3 & 4
\end{matrix}
\quad .
\]
\end{example}

% [DG][v77] Added the above example.

So far (in this section) we have worked with a finite set of
variables $x_1,x_2,\ldots,x_n$ (where $n$ is a fixed nonnegative
integer) and with partitions having at most $n$ parts. We now drop
these conventions and restrictions; thus, partitions again mean
arbitrary partitions, and $\xx$ again means the infinite family
$\left(x_1,x_2,x_3,\ldots\right)$ of variables. In this setting,
we have the following analogue of Corollary~\ref{finite-L-R-rule}:

\begin{corollary}
\label{infinite-L-R-rule}
For partitions $\lambda, \mu, \nu$ (of any lengths), one has
\begin{equation}
\label{eq.infinite-L-R-rule.1}
s_{\nu} s_{\lambda/\mu} = \sum_T s_{\nu+\cont(T)}
\end{equation}
where $T$ runs through all column-strict tableaux
of shape $\lambda/\mu$ with
the property that for each $j=1,2,3,\ldots$, the weak
composition $\nu + \cont(T|_{\cols\geq j})$ is a partition.
In particular, taking $\nu=\varnothing$, we obtain
\begin{equation}
\label{eq.infinite-L-R-rule.2}
s_{\lambda/\mu} = \sum_T s_{\cont(T)}
\end{equation}
where in the sum $T$ runs through all column-strict
tableaux of shape $\lambda/\mu$
for which each $\cont(T|_{\cols\geq j})$ is a partition.
\end{corollary}

\begin{proof}[Proof of Corollary~\ref{infinite-L-R-rule}.]
Essentially, Corollary~\ref{infinite-L-R-rule} is obtained from
Corollary~\ref{finite-L-R-rule} by ``letting $n$ (that is,
the number of variables) tend to $\infty$''. This can be formalized
in different ways: One way is to endow the ring
of power series $\kk\left[\left[\xx\right]\right]
= \kk\left[\left[x_1,x_2,x_3,\ldots\right]\right]$ with the
coefficientwise topology\footnote{This topology is defined as
follows:
\par
We endow the ring $\kk$ with the discrete topology. Then, we can regard
the $\kk$-module $\kk \left[ \left[ \xx \right] \right]$ as a direct
product of infinitely many copies of $\kk$ (by identifying every power
series in $\kk \left[ \left[ \xx \right] \right]$ with the
family of its coefficients). Hence, the product topology is a well-defined
topology on $\kk \left[ \left[ \xx \right] \right]$; this
topology is denoted as the \dfn{coefficientwise topology}.
Its name is due to the fact that a sequence
$\left( a_n \right)_{n\in \NN}$ of power series converges
to a power series $a$ with respect to this topology if and
only if for every monomial $\mathfrak{m}$, all sufficiently
high $n\in \NN$ satisfy
\[
\left(  \text{the coefficient of }\mathfrak{m}\text{ in }a_{n}\right)
=\left(  \text{the coefficient of }\mathfrak{m}\text{ in }a\right)  .
\]
}, and to show that the left hand side of
\eqref{eq.finite-L-R-rule.1} tends to the left hand side of
\eqref{eq.infinite-L-R-rule.1} when $n \to \infty$, and the same
holds for the right hand sides. A different approach proceeds by
regarding $\Lambda$ as the inverse limit of the
$\Lambda\left(x_1,x_2,\ldots,x_n\right)$.
\end{proof}

Comparing coefficients of a given Schur function $s_\nu$
in \eqref{eq.infinite-L-R-rule.2}, we obtain
the following version of the Littlewood-Richardson rule.
\begin{corollary}
\label{L-R-rule}
For partitions $\lambda, \mu, \nu$ (of any lengths),
the Littlewood-Richardson coefficient $c^\lambda_{\mu,\nu}$ counts
column-strict tableaux $T$ of shape $\lambda/\mu$ with $\cont(T)=\nu$ having
the property that each $\cont(T|_{\cols\geq j})$ is a partition.
\end{corollary}

% [DG][v19] Last time, I suggested to add
% http://mathoverflow.net/questions/118225/how-to-show-a-certain-determinant-is-non-zero
% as an exercise. Turns out that the intended solution was false...

\subsection{The Pieri and Assaf-McNamara skew Pieri rule}

The classical \dfn{Pieri rule} refers to two special cases of the Littlewood-Richardson rule.
To state them, recall that a skew shape
is called a \emph{horizontal (resp. vertical) strip}\index{horizontal strip}\index{vertical strip}
if no two of its cells lie in the same column (resp. row).
A \emph{horizontal (resp. vertical) $n$-strip}\index{horizontal $n$-strip}\index{vertical $n$-strip}
(for $n \in \NN$)
shall mean a horizontal (resp. vertical) strip of size $n$
(that is, having exactly $n$ cells).

% [DG][v25] Added definition of "strip"s (before, only "$n$-strips"
% were defined, but "strips" were referred to -- but more importantly
% I want to use "strip"s in exercises).

\begin{theorem}
\label{Pieri-rules}
For every partition $\lambda$ and any $n \in \NN$, we have
\begin{align}
s_\lambda h_n =&
\sum\limits_{\substack{\lambda^+: \lambda^+/\lambda \text{ is a }\\\text{horizontal }n\text{-strip}}}  s_{\lambda^+} ;
\label{eq.pieri.h} \\
s_\lambda e_n =&
\sum\limits_{\substack{\lambda^+: \lambda^+/\lambda \text{ is a }\\\text{vertical }n\text{-strip}}}  s_{\lambda^+} .
\label{eq.pieri.e}
\end{align}
\end{theorem}

\begin{example}
\label{exam.Pieri-rules.1}
In the following equality, we are representing each partition
by its Ferrers diagram\footnote{And we are drawing each Ferrers
diagram with its boxes spaced out, in order to facilitate
counting the boxes.}.
\[
\begingroup % keep the change local
\setlength\arraycolsep{0.2pc} % Make the spacing equal in rows and in columns.
\begin{aligned}
&\begin{matrix}
s &   &   & \\
  &\sq&\sq&\sq\\
  &\sq&\sq&   \\
  &\sq&\sq&
\end{matrix}
\qquad \bullet \qquad
\begin{matrix}
h_2 &   &   \\
%  &\sq&\sq \\
  &   & \\
  &   & \\
  &   & \\
\end{matrix}
\\
 & \\
&= \ \ \ 
%\qquad
\begin{matrix}
s &   &   & \\
  &\sq&\sq&\sq\\
  &\sq&\sq&   \\
  &\sq&\sq& \\
  &\bsq&\bsq&
\end{matrix}
\quad
+
\quad
\begin{matrix}
s &   &   & \\
  &\sq&\sq&\sq\\
  &\sq&\sq&\bsq\\
  &\sq&\sq& \\
  &\bsq&  &
\end{matrix}
\quad 
+
\quad
\begin{matrix}
s &   &   &   &\\
  &\sq&\sq&\sq&\bsq\\
  &\sq&\sq&   &\\
  &\sq&\sq&   &\\
  &\bsq&  &   &
\end{matrix}
\\
& \\
&
\qquad \qquad \qquad
\ 
+
\quad 
\begin{matrix}
s &   &   &   &\\
  &\sq&\sq&\sq&\bsq\\
  &\sq&\sq&\bsq &\\
  &\sq&\sq&   &\\
\end{matrix}
\quad
+
\quad
\begin{matrix}
s &   &   &   &    &\\
  &\sq&\sq&\sq&\bsq&\bsq\\
  &\sq&\sq&   &    &\\
  &\sq&\sq&   &    &\\
\end{matrix}
\end{aligned}
\endgroup
\]
If $\lambda$ is the partition $\left(3,2,2\right)$ on
the left hand side, then all partitions $\lambda^+$ on the
right hand side visibly have the property that
$\lambda^+ / \lambda$ is a horizontal $2$-strip\footnote{We
have colored the boxes of $\lambda^+ / \lambda$ black.},
as \eqref{eq.pieri.h} predicts.
\end{example}

% [DG][v77] Added explanations to the above example.
% Made the space between columns (\arraycolsep) in the
% Ferrers diagrams smaller so it approximately matches
% the space between rows (is there an exact way to do
% this? LaTeX doesn't know \arrayrowsep).
% Also, replaced h_{[][]} by h_2, since I don't see
% much of a point in drawing the diagram under an h.

\begin{vershort}
\begin{proof}[Proof of Theorem~\ref{Pieri-rules}]
For the first Pieri formula involving $h_n$, as $h_n=s_{(n)}$ one has
\[
s_\lambda h_n = \sum_{\lambda^+} c^{\lambda^+}_{\lambda,(n)} s_{\lambda^+}.
\]
Corollary~\ref{L-R-rule} says $c^{\lambda^+}_{\lambda,(n)}$ counts
column-strict tableaux $T$ of shape $\lambda^+/\lambda$
having $\cont(T)=(n)$ (i.e. all entries of $T$ are $1$'s),
with an extra condition.  Since its entries are all equal,
such a $T$ must certainly have shape being a horizontal strip, and
more precisely a horizontal $n$-strip (since it has $n$ cells).
Conversely,
for any horizontal $n$-strip, there is a unique such filling,
and it will trivially satisfy the extra condition
that $\cont(T|_{\cols\geq j})$ is a partition for each $j$.
Hence $c^{\lambda^+}_{\lambda,(n)}$ is $1$ if $\lambda^+/\lambda$ is
a horizontal $n$-strip, and $0$ else.

For the second Pieri formula involving $e_n$, using $e_n=s_{(1^n)}$ one has
\[
s_\lambda e_n = \sum_{\lambda^+} c^{\lambda^+}_{\lambda,(1^n)} s_{\lambda^+}.
\]
Corollary~\ref{L-R-rule} says $c^{\lambda^+}_{\lambda,(1^n)}$
counts column-strict tableaux $T$ of shape $\lambda^+/\lambda$
having $\cont(T)=(1^n)$, so its entries are $1,2,\ldots,n$ each
occurring once, with the extra condition that $1,2,\ldots,n$
appear from right to left.  Together with
the tableau condition, this forces at most one entry in each row, that is
$\lambda^+/\lambda$ is a vertical strip, and then there is a unique
way to fill it  (maintaining column-strictness and the extra condition
that $1,2,\ldots,n$ appear from right to left).  Thus $c^{\lambda^+}_{\lambda,(1^n)}$
is $1$ if $\lambda^+/\lambda$ is a vertical $n$-strip, and $0$ else.
\end{proof}
\end{vershort}

\begin{verlong}
\begin{proof}[Proof of Theorem~\ref{Pieri-rules}]
For the first Pieri formula involving $h_n$, as $h_n=s_{(n)}$ one has
\begin{equation}
\label{pf.Pieri-rules.h}
s_\lambda h_n = \sum_{\lambda^+} c^{\lambda^+}_{\lambda,(n)} s_{\lambda^+}.
\end{equation}
Corollary~\ref{L-R-rule} says $c^{\lambda^+}_{\lambda,(n)}$ counts
column-strict tableaux $T$ of shape $\lambda^+/\lambda$
having $\cont(T)=(n)$ (i.e. all entries of $T$ are $1$'s),
with an extra condition (stating that
$\cont(T|_{\cols\geq j})$ is a partition for each $j$).  Since
its entries are all equal,
such a $T$ must certainly have shape being a horizontal strip, and
more precisely a horizontal $n$-strip (since it has $n$ cells).
Conversely,
for any horizontal $n$-strip, there is a unique such filling,
and it will trivially satisfy the extra condition
that $\cont(T|_{\cols\geq j})$ is a partition for each $j$.
Hence $c^{\lambda^+}_{\lambda,(n)}$ is $1$ if $\lambda^+/\lambda$ is
a horizontal $n$-strip, and $0$ else. Thus,
\eqref{pf.Pieri-rules.h} becomes
\[
s_\lambda h_n =
\sum\limits_{\substack{\lambda^+: \lambda^+/\lambda \text{ is a }\\\text{horizontal }n\text{-strip}}}  s_{\lambda^+} ,
\]
and \eqref{eq.pieri.h} is proven.

For the second Pieri formula involving $e_n$, using $e_n=s_{(1^n)}$ one has
\begin{equation}
\label{pf.Pieri-rules.e}
s_\lambda e_n = \sum_{\lambda^+} c^{\lambda^+}_{\lambda,(1^n)} s_{\lambda^+}.
\end{equation}
Corollary~\ref{L-R-rule} says $c^{\lambda^+}_{\lambda,(1^n)}$
counts column-strict tableaux $T$ of shape $\lambda^+/\lambda$
having $\cont(T)=(1^n)$, so its entries are $1,2,\ldots,n$ each
occurring once, with an extra condition (stating that
$\cont(T|_{\cols\geq j})$ is a partition for each $j$).
This extra condition forces the entries $1,2,\ldots,n$ to
appear in the tableau $T$ from right to left (more
precisely, if $i$ and $i'$ are two entries of $T$
satisfying $i' \leq i$, then the column in which $i'$ appears
is weakly to the right of the column in which $i$
appears\footnote{\textit{Proof:} Assume the
contrary. Then, there are two entries $i$ and $i'$ of $T$
satisfying $i' \leq i$, such that $i'$ appears in a column
strictly to the left of the column in which $i$ appears.
Consider these $i$ and $i'$. Let $j$ be the column in which
$i$ appears, and $j'$ be the column in which $i'$ appears;
thus, $j' < j$. We have $i' \leq i$, so that
$\left(\cont(T|_{\cols\geq j})\right)_{i'} \geq
\left(\cont(T|_{\cols\geq j})\right)_i$ (due to the extra
condition saying that $\cont(T|_{\cols\geq j})$ is a partition).
Hence,
$\left(\cont(T|_{\cols\geq j})\right)_{i'} \geq
\left(\cont(T|_{\cols\geq j})\right)_i \geq 1$ (since $i$ does
appear in $T|_{\cols\geq j}$, being an entry of column $j$ of
$T$), which means that at least one entry of
$T|_{\cols\geq j}$ must be equal to $i'$. In other words, $i'$
appears in one of the columns $j, j+1, j+2, \ldots$ of the
tableau $T$. Since $i'$ also appears in column $j'$ (by the
definition of $j'$), which is \textbf{not} one of the
columns $j, j+1, j+2, \ldots$ (since $j' < j$), this shows
that $i'$ appears
(at least) twice in the tableau $T$. But this is absurd, since
the entries of $T$ are $1,2,\ldots,n$ each occurring once. This
contradiction proves that our assumption was wrong, qed.});
it also forces the tableau $T$ to have at
most one entry in each row\footnote{\textit{Proof:} Assume the
contrary. Thus, some row of $T$ has (at least) two entries
-- say, an entry $i$ in a column $j$, and an entry $i'$ in a
different column $j'$. Assume WLOG that $j' < j$. Thus,
$i' \leq i$ (since the entries of $T$ increase weakly
left-to-right along rows). Now, a contradiction can be obtained
as in the previous footnote, and thus our assumption was wrong,
qed.}. Thus, for such a tableau $T$ to exist,
$\lambda^+/\lambda$ has to be a vertical strip (because $T$
must have at most one entry in each row), more precisely
a vertical $n$-strip. Conversely, given a partition $\lambda^+$
for which $\lambda^+/\lambda$ is a vertical
$n$-strip, there is a unique way to fill it to obtain a
column-strict tableau $T$ satisfying our extra
condition\footnote{Indeed, one must fill $\lambda^+/\lambda$
with the entries $1,2,\ldots,n$ from top to bottom; this is
the only way to maintain column-strictness and the extra condition
that $1,2,\ldots,n$ appear from right to left.}.
Thus $c^{\lambda^+}_{\lambda,(1^n)}$
is $1$ if $\lambda^+/\lambda$ is a vertical $n$-strip, and $0$ else.
Hence, \eqref{pf.Pieri-rules.e} becomes
\[
s_\lambda e_n =
\sum\limits_{\substack{\lambda^+: \lambda^+/\lambda \text{ is a }\\\text{vertical }n\text{-strip}}}  s_{\lambda^+} ,
\]
and \eqref{eq.pieri.e} is proven.
\end{proof}
\end{verlong}

% [DG][v14] Replaced "maintaining column-strictness" by
% "(maintaining column-strictness and the extra condition
% that $1,2,\ldots,n$ appear from right to left)". Or am I
% wrong here?

In 2009, Assaf and McNamara \cite{AssafMcNamara} proved an elegant
generalization.\index{Assaf-McNamara skew Pieri rule}

\begin{theorem}
\label{Assaf-McNamara-skew-Pieri}
For any partitions $\lambda$ and $\mu$ and any $n \in \NN$,
we have\footnote{Note that $\mu \subseteq \lambda$ is not required.
(The left hand sides are $0$ otherwise, but this does not
trivialize the equalities.)}
\begin{align}
s_{\lambda/\mu} h_n &=
\sum\limits_{\substack{\lambda^+,\mu^-:\\
                  \lambda^+/\lambda \text{ a horizontal strip;}\\
                   \mu/\mu^{-}\text{ a vertical strip;}\\
                   |\lambda^+/\lambda|+|\mu/\mu^-|=n}}
           (-1)^{|\mu/\mu^{-}|} s_{\lambda^+/\mu^-} ;
\label{eq.Assaf-McNamara-skew-Pieri.1}
\\
s_{\lambda/\mu} e_n &=
\sum\limits_{\substack{\lambda^+,\mu^-: \\
                 \lambda^+/\lambda \text{ a vertical strip;}\\
                    \mu/\mu^{-}\text{ a horizontal strip;}\\
                    |\lambda^+/\lambda|+|\mu/\mu^-|=n}}
           (-1)^{|\mu/\mu^{-}|} s_{\lambda^+/\mu^-} .
\label{eq.Assaf-McNamara-skew-Pieri.2}
\end{align}
\end{theorem}

\begin{example}
With the same conventions as in Example~\ref{exam.Pieri-rules.1}%
\footnote{but this time coloring both the boxes in $\lambda^+ / \lambda$
and the boxes in $\mu / \mu^-$ black},
we have
\[
\begingroup % keep the change local
\setlength\arraycolsep{0.2pc} % Make the spacing equal in rows and in columns.
\begin{aligned}
&\begin{matrix}
s &   &   & \\
  &   &\sq&\sq\\
  &   &\sq&   \\
  &\sq&\sq&
\end{matrix}
\qquad \bullet \qquad
\begin{matrix}
h_2 &   &   \\
%  &\sq&\sq \\
  &   & \\
  &   & \\
  &   & \\
\end{matrix}
\\
 & \\
&=\qquad
\begin{matrix}
s &   &   & \\
  &   &\sq&\sq\\
  &   &\sq&   \\
  &\sq&\sq& \\
  &\bsq&\bsq&
\end{matrix}
\qquad+\qquad
\begin{matrix}
s &   &   & \\
  &   &\sq&\sq\\
  &   &\sq&\bsq\\
  &\sq&\sq& \\
  &\bsq&  &
\end{matrix}
\qquad+\qquad
\begin{matrix}
s &   &   &   &\\
  &   &\sq&\sq&\bsq\\
  &   &\sq&   &\\
  &\sq&\sq&   &\\
  &\bsq&  &   &
\end{matrix}
\\
& \\
&\qquad \qquad \qquad
+\qquad
\begin{matrix}
s &   &   &   &\\
  &   &\sq&\sq&\bsq\\
  &   &\sq&\bsq &\\
  &\sq&\sq&   &\\
\end{matrix}
\qquad+\qquad\begin{matrix}
s &   &   &   &    &\\
  &   &\sq&\sq&\bsq&\bsq\\
  &   &\sq&   &    &\\
  &\sq&\sq&   &    &\\
\end{matrix}
\\
 & \\
&\qquad \qquad
-\quad
\begin{matrix}
s &    &   &   \\
  &    &\sq&\sq\\
  &\bsq&\sq&   \\
  &\sq &\sq&   \\
  &\bsq&   &
\end{matrix}
\qquad-\qquad
\begin{matrix}
s &    &   &   \\
  &    &\sq&\sq\\
  &\bsq&\sq&\bsq\\
  &\sq &\sq&   \\
\end{matrix}
\qquad-\qquad
\begin{matrix}
s &    &   &   &\\
  &    &\sq&\sq&\bsq\\
  &\bsq&\sq&   &\\
  &\sq &\sq&   &
\end{matrix}
\\
 & \\
&\qquad\qquad\qquad\qquad
+\qquad
\begin{matrix}
s &    &   &   \\
  &\bsq&\sq&\sq\\
  &\bsq&\sq&   \\
  &\sq &\sq&
\end{matrix}
\end{aligned}
\endgroup
\]
which illustrates the first equality of
Theorem~\ref{Assaf-McNamara-skew-Pieri}.
\end{example}

% [DG][v77] Same changes to this example as above.

\noindent
Theorem~\ref{Assaf-McNamara-skew-Pieri} is proven
in the next section, using an important Hopf algebra tool.

\begin{exercise}
\label{exe.horistrip.postnikov}
Let $\lambda = \left(\lambda_1, \lambda_2, \lambda_3, \ldots\right)$ and
$\mu = \left(\mu_1, \mu_2, \mu_3, \ldots\right)$
be two partitions such that $\mu \subseteq \lambda$.

\begin{enumerate}

\item[(a)] Show
that $\lambda / \mu$ is a horizontal strip if and only if every
$i \in \left\{ 1,2,3,\ldots \right\}$ satisfies
$\mu_i \geq \lambda_{i+1}$.\ \ \ \ \footnote{In
other words, $\lambda / \mu$ is a horizontal strip if and only if
$\left( \lambda_2, \lambda_3, \lambda_4, \ldots \right)
\subseteq \mu$.
This simple observation has been used by Pak and Postnikov
\cite[\S 10]{PakPostnikov-oscil} for a new approach to RSK-type
algorithms.}

\item[(b)] Show
that $\lambda / \mu$ is a vertical strip if and only if every
$i \in \left\{ 1,2,3,\ldots \right\}$ satisfies
$\lambda_i \leq \mu_i + 1$.

\end{enumerate}
\end{exercise}

% [DG][v25] Added this (near-trivial) exercise.

% [DG][v48] Added part (b) (for the sake of "symmetry" and possibly
% of later reference; the proof is a one-liner).

\begin{exercise} \phantomsection
\label{exe.pieri.h2}
\begin{itemize}
\item[(a)] Let $\lambda$ and $\mu$ be two partitions such that $\mu
\subseteq \lambda$. Let $n \in \NN$.
Show that $\left(h_n, s_{\lambda / \mu}\right)$ equals $1$ if $\lambda / \mu$
is a horizontal $n$-strip, and equals $0$ otherwise.
\item[(b)] Use part (a) to give a new proof of
\eqref{eq.pieri.h}.
\end{itemize}
\end{exercise}

% [DG][v25] Added above exercise showing that the first Pieri rule is
% simpler than Littlewood-Richardson (provided one knows Cauchy).
% (I don't see a similar simple argument for the second Pieri rule,
% though.)

\begin{exercise}
\label{exe.pieri.RSK}
Prove Theorem~\ref{Pieri-rules} again using the ideas of
the proof of Theorem~\ref{thm.cauchy-identity}.
\end{exercise}

% [DG][v47] Added previous exercise (it comes almost for free once RSK is
% introduced).

\begin{verlong}
\begin{exercise}
\label{exe.pieri.alt-e}
Obtain a new proof for \eqref{eq.pieri.e} using
Corollary~\ref{cor.stembridge.slambda} and the formula
\eqref{e-generating-function}.\footnote{Such a proof is
given in \cite[\S I.5, Rmk. 1, pp. 73--74]{Macdonald} and in
\cite[Proposition 2.3]{Leeuwen-altSchur}. A similar, but more
complicated, proof can be given for \eqref{eq.pieri.h}; this is
done, for example, in \cite[Proposition 2.5]{Leeuwen-altSchur}.}
\end{exercise}
\end{verlong}

\begin{exercise}
\label{exe.cauchy.cauchy}
Let $A$ be a commutative ring, and $n \in \NN$.
\begin{itemize}
\item[(a)] Let
$a_1, a_2, \ldots, a_n$ be $n$ elements of $A$.
Let $b_1, b_2, \ldots, b_n$ be $n$ further elements
of $A$.
If $a_i - b_j$ is an invertible element of
$A$ for every $i \in \left\{1,2,\ldots,n\right\}$ and
$j \in \left\{1,2,\ldots,n\right\}$, then prove that
\[
\det\left(\left(\frac{1}{a_i - b_j}\right)_{i,j = 1,2,\ldots,n}\right)
= \frac{\prod_{1 \leq j < i \leq n} \left(\left(a_i - a_j\right)
\left(b_j - b_i\right) \right)}{\prod_{\left(i, j\right)
\in \left\{1,2,\ldots,n\right\}^2}
\left(a_i - b_j\right)}.
\]
\item[(b)] Let
$a_1, a_2, \ldots, a_n$ be $n$ elements of $A$.
Let $b_1, b_2, \ldots, b_n$ be $n$ further elements
of $A$.
If $1 - a_i b_j$ is an invertible element
of $A$ for every $i \in \left\{1,2,\ldots,n\right\}$ and
$j \in \left\{1,2,\ldots,n\right\}$, then prove that
\[
\det\left(\left(\frac{1}{1 - a_i b_j}\right)_{i,j = 1,2,\ldots,n}\right)
= \frac{\prod_{1 \leq j < i \leq n} \left(\left(a_i - a_j\right)
\left(b_i - b_j\right) \right)}{\prod_{\left(i, j\right)
\in \left\{1,2,\ldots,n\right\}^2}
\left(1 - a_i b_j\right)}.
\]
\item[(c)] Use the result of part (b) to give a new
proof for Theorem~\ref{thm.cauchy-identity}.\footnote{This
approach to Theorem~\ref{thm.cauchy-identity} is taken in
\cite[\S 4]{CrawleyBoevey} (except that
\cite{CrawleyBoevey} only works with finitely many
variables).}
\end{itemize}
\end{exercise}

The determinant on the left hand side of
Exercise~\ref{exe.cauchy.cauchy}(a) is known as the \dfn{Cauchy
determinant}.

% [DG][v25] Added above exercise.

\begin{exercise}
\label{exe.jacobi-trudi.len2.0}
Prove that $s_{\left(a,b\right)} = h_a h_b - h_{a+1} h_{b-1}$ for any
two integers $a \geq b \geq 0$ (where we set $h_{-1} = 0$ as usual).

(Note that this is precisely the Jacobi-Trudi formula
\eqref{eq.jacobi-trudi.h} in the case when
$\lambda = \left(a,b\right)$ is a partition with at most two
entries and $\mu = \varnothing$.)
\end{exercise}

\begin{exercise}
\label{exe.pieri.cauchy}
If $\lambda$ is a partition and $\mu$ is a weak composition, let
$K_{\lambda, \mu}$ denote the number of column-strict tableaux $T$ of
shape $\lambda$ having $\cont\left(T\right) = \mu$. (This
$K_{\lambda, \mu}$ is called the
\dfn{$\left(\lambda, \mu\right)$-Kostka number}.)
\begin{itemize}
\item[(a)] Use Theorem~\ref{Pieri-rules} to show that every
partition $\mu$ satisfies
$h_\mu = \sum_\lambda K_{\lambda, \mu} s_\lambda$, where the sum
ranges over all partitions $\lambda$.
\item[(b)] Use this to give a new proof for
Theorem~\ref{thm.cauchy-identity}.\footnote{Of course, this gives
a new proof of Theorem~\ref{thm.cauchy-identity} only when coupled
with a proof of Theorem~\ref{Pieri-rules} which does not rely on
Theorem~\ref{thm.cauchy-identity}. The proof of
Theorem~\ref{Pieri-rules} we gave in the text above did not rely
on Theorem~\ref{thm.cauchy-identity}, whereas the proof of
\eqref{eq.pieri.h} given in Exercise~\ref{exe.pieri.h2}(b) did.}
\item[(c)] Give a new proof of the fact (previously shown as
Proposition~\ref{prop.Lambda.omega-and-S}(j)) that
$\left(h_\lambda\right)_{\lambda\in\Par}$ is a graded basis of
the graded $\kk$-module $\Lambda$.
\end{itemize}
\end{exercise}

% [DG][v19] Added above exercise.

% [DG][v23] Added part (c).

\begin{exercise} \phantomsection
\label{exe.pieri.omega}
\begin{itemize}
\item[(a)] Define a $\kk$-linear map $\mathfrak{Z} : \Lambda \to \Lambda$
by having it send $s_\lambda$ to $s_{\lambda^t}$ for every partition
$\lambda$. (This is clearly well-defined, since $\left(s_\lambda\right)_{\lambda \in \Par}$
is a $\kk$-basis of $\Lambda$.) Show that
\[
\mathfrak{Z}\left(f h_n\right) = \mathfrak{Z}\left(f\right) \cdot \mathfrak{Z}\left(h_n\right)
\qquad \qquad \text{ for every } f \in \Lambda \text{ and every } n \in \NN .
\]
\item[(b)] Show that $\mathfrak{Z} = \omega$.
\item[(c)] Show that $c_{\mu, \nu}^{\lambda} = c_{\mu^t, \nu^t}^{\lambda^t}$ for
any three partitions $\lambda$, $\mu$ and $\nu$.
\item[(d)] Use this to prove \eqref{skew-Schur-antipodes}.\footnote{The
first author learned this approach to \eqref{skew-Schur-antipodes} from
Alexander Postnikov.}
\end{itemize}
\end{exercise}

% [DG][v23] Added the preceding exercise.

\begin{exercise} \phantomsection
\label{exe.cauchy.dual}
\begin{itemize}
\item[(a)] Show that
\[
\prod_{i,j=1}^\infty \left(1+x_iy_j\right)
= \sum_{\lambda \in \Par} s_\lambda\left(\xx\right) s_{\lambda^t}\left(\yy\right)
= \sum_{\lambda \in \Par} e_\lambda\left(\xx\right) m_\lambda\left(\yy\right)
\]
in the power series ring $\kk\left[\left[\xx, \yy\right]\right]
= \kk\left[\left[x_1, x_2, x_3, \ldots, y_1, y_2, y_3, \ldots\right]\right]$.
\item[(b)] Assume that $\QQ$ is a subring of $\kk$. Show that
\[
\prod_{i,j=1}^\infty \left(1+x_iy_j\right)
= \sum_{\lambda \in \Par}
\left(-1\right)^{\left|\lambda\right| - \ell\left(\lambda\right)}
z_\lambda^{-1} p_\lambda\left(\xx\right) p_\lambda\left(\yy\right)
\]
in the power series ring $\kk\left[\left[\xx, \yy\right]\right]
= \kk\left[\left[x_1, x_2, x_3, \ldots, y_1, y_2, y_3, \ldots\right]\right]$,
where $z_\lambda$ is defined as in
Proposition~\ref{prop-Two-other-Cauchy-expansions}.
\end{itemize}
\end{exercise}

The first equality of Exercise~\ref{exe.cauchy.dual}(a) appears in
\cite[Thm. 7.14.3]{Stanley}, \cite[Thm. 4.8.6]{Sagan} and several other
references under the name of the \dfn{dual Cauchy identity}, and is commonly
proven using a ``dual'' analogue of the Robinson-Schensted-Knuth algorithm.

% [DG][v25] Added the above exercise (on the dual Cauchy identity).

% [DG][v45] Added a solution to its part (a), using a variation on RSK
% which uses the same RS-insertion step as usual RSK but processes the
% biletters in a different order.
% Curiously it seems to be unpopular in literature (instead, everyone
% uses dual RSK).

% [DG][v46] Added the dual RSK solution too (mainly to nail down the
% properties of dual RSK, which I need to know for a different reason).

\begin{verlong}
The following two exercises (Exercise~\ref{exe.jacobi-trudi} and
Exercise~\ref{exe.jacobi-trudi.macdonald}) give two different proofs of
the Jacobi-Trudi identities.

\begin{exercise}
\label{exe.jacobi-trudi}
For any two partitions $\lambda$ and $\mu$ and any $k \in \NN$, define a
symmetric function $\widetilde{s}_{\lambda / \mu}^{[k]} \in \Lambda$
by
\[
\widetilde{s}_{\lambda / \mu}^{[k]}
= \det \left(\left( h_{\lambda_i - \mu_j - i + j} \right)_{i,j=1,2,\ldots,k}  \right).
\]
(As usual, we are denoting by $\nu_i$ the $i$-th entry of a partition $\nu$
here.)
\begin{itemize}
\item[(a)] For any partitions $\lambda$ and $\mu$ and any $k \in \NN$
satisfying $k \geq \ell\left(\lambda\right)$ and
$k \geq \ell\left(\mu\right)$, show that the symmetric function
$\widetilde{s}_{\lambda / \mu}^{[k]}$ does not depend on $k$. We thus
can define a symmetric function $\widetilde{s}_{\lambda / \mu}$ for
any two partitions $\lambda$ and $\mu$ by setting
$\widetilde{s}_{\lambda / \mu} = \widetilde{s}_{\lambda / \mu}^{[k]}$,
where $k$ is any element of $\NN$ satisfying
$k \geq \ell\left(\lambda\right)$ and $k \geq \ell\left(\mu\right)$.
In the rest of this exercise, we will study this
$\widetilde{s}_{\lambda / \mu}$ (and show that it equals
$s_{\lambda / \mu}$).
\item[(b)] If two partitions $\lambda$ and $\mu$ do not satisfy
$\mu \subseteq \lambda$, then show that
$\widetilde{s}_{\lambda / \mu} = 0$.
\item[(c)] For any two partitions $\lambda$ and $\nu$, show that
\[
\Delta \widetilde{s}_{\lambda / \nu}
= \sum\limits_{\substack{\mu \in \Par: \\ \nu \subseteq \mu \subseteq \lambda}}
  \widetilde{s}_{\lambda / \mu} \otimes \widetilde{s}_{\mu / \nu}.
\]
\item[(d)] Show that
$\widetilde{s}_{\lambda / \mu} \in \Lambda_{\left|\lambda/\mu\right|}$
for any partitions $\lambda$ and $\mu$ satisfying
$\mu \subseteq \lambda$.
\item[(e)] For any two partitions $\lambda$ and $\mu$ and any
$n \in \NN$, show that
\[
\left( h_n, \widetilde{s}_{\lambda / \mu} \right)
= \left( h_n, s_{\lambda / \mu} \right).
\]
\item[(f)] If $x$ is a primitive element of $\Lambda$ such
that every $n \in \NN$ satisfies $\left( h_n, x \right) = 0$,
then show that $x = 0$.
\item[(g)] Prove that
$\widetilde{s}_{\lambda / \mu} = s_{\lambda / \mu}$ for any
two partitions $\lambda$ and $\mu$ satisfying
$\mu \subseteq \lambda$.
\item[(h)] Use this to obtain a proof of
Theorem~\ref{Jacobi-Trudi-formulae}.
\end{itemize}
\end{exercise}

% [DG][v75] Added the preceding exercise, but only to the verlong
% version. (Not the best proof of Jacobi-Trudi, but Hopfish and
% somewhat new. But it's become a mess when I wrote it up a few
% years ago, so I never felt like including it in the notes.
% But this is the verlong version, so no one will find it here
% except for those who really care.
% Macdonald uses similar arguments for a different result; I
% wasn't aware of % that when I started writing this up.)
\end{verlong}

\begin{exercise}
\label{exe.jacobi-trudi.macdonald}
Prove Theorem~\ref{Jacobi-Trudi-formulae}.

[\textbf{Hint:}\footnote{This is the proof given in
Stanley \cite[\S 7.16, Second Proof of Thm. 7.16.1]{Stanley} and
Macdonald \cite[proof of (5.4)]{Macdonald}.}
Switch $\xx$ and $\yy$ in the formula of
Exercise~\ref{exe.cauchy.skew}(a), and specialize the resulting
equality by replacing $\yy$ by a finite set of variables
$\left(y_1, y_2, \ldots, y_\ell\right)$;
then, set $n = \ell$ and $\rho = \left(n-1,n-2,\ldots,0\right)$,
and multiply with the alternant
$a_\rho\left(y_1, y_2, \ldots, y_\ell\right)$, using
Corollary~\ref{cor.stembridge.slambda} to simplify the result;
finally, extract the coefficient of $\yy^{\lambda + \rho}$.]
\end{exercise}

% [DG][v25] Added the preceding exercise (proving Jacobi-Trudi).

\begin{exercise}
\label{exe.hall-e-f}Prove the following:

\begin{enumerate}
\item[(a)] We have $\left(  S\left(  f\right)  ,S\left(  g\right)  \right)
=\left(  f,g\right)  $ for all $f\in\Lambda$ and $g\in\Lambda$.

\item[(b)] We have $\left(  e_{n},f\right)  =\left(  -1\right)  ^{n}%
\cdot\left(  S\left(  f\right)  \right)  \left(  1\right)  $ for any
$n \in \NN$ and $f \in \Lambda_n$. (See Exercise~\ref{exe.Lambda.subs.fin}
for the meaning of $\left(  S\left(  f\right)  \right)  \left(  1\right)  $.)
\end{enumerate}
\end{exercise}

% [DG][v75] Added the above exercise.

\subsection{Skewing and Lam's proof of the skew Pieri rule}

We codify here the operation $s_\mu^\perp$ of \emph{skewing by $s_\mu$},
acting on Schur functions via
\[
s_\mu^\perp(s_\lambda) = s_{\lambda/\mu}
\]
(where, as before, one defines $s_{\lambda/\mu}=0$ if $\mu \not\subseteq \lambda$).
These operations play a crucial role
\begin{itemize}
\item
in Lam's proof of the
skew Pieri rule,
\item in Lam, Lauve, and Sottile's proof \cite{LamLauveSottile}
of a more general  skew Littlewood-Richardson rule
that had been conjectured by Assaf and McNamara, and
\item in Zelevinsky's structure theory of PSH's to be developed in the next chapter.
\end{itemize}
We are going to define them in the general setting of any graded Hopf
algebra.

% [DG][v9] Here and in the following, "PSH-algebra" replaced by
% "PSH".

\begin{definition}
\label{def.skewing}
Given a graded Hopf algebra $A$, and its (graded) dual $A^o$, let
$(\cdot,\cdot)=(\cdot,\cdot)_A : A^o \times A \to \kk$ be
the pairing defined by
$(f,a):=f(a)$ for $f$ in $A^o$ and $a$ in $A$.
Then define for each
$f$ in $A^o$ an operator
$A \overset{f^\perp}{\rightarrow} A$ as follows\footnote{This $f^\perp(a)$ is
called $a \leftharpoonup f$ in Montgomery \cite[Example 1.6.5]{Montgomery}.}:
for $a$ in $A$
with $\Delta(a) = \sum a_1 \otimes a_2$, let
\[
f^\perp(a) = \sum (f,a_1) a_2.
\]
In other words, $f^\perp$ is the composition
\[
\xymatrix{
A \ar[r]^-{\Delta} & A \otimes A \ar[r]^-{f \otimes \id} & \kk \otimes A \ar[r]^-{\cong} & A ,
}
\]
where the rightmost arrow is the canonical isomorphism $\kk \otimes A \to A$.
This operator $f^\perp$ is called \dfn{skewing by $f$}.
\end{definition}

% [DG][v64] Removed the "finite type" condition from the above
% definition (I shall instead require it whenever it is necessary).
% Also, added the alternative definition through the composition
% of arrows.

\noindent
Now, recall that the Hall inner product induces an isomorphism
$\Lambda^o \cong \Lambda$ (by Corollary~\ref{cor.Lambda.selfdual}).
Hence, we can regard any element $f \in \Lambda$ as an element
of $\Lambda^o$; this allows us to define an operator
$f^\perp : \Lambda \to \Lambda$ for each $f \in \Lambda$
(by regarding $f$ as an element of $\Lambda^o$, and applying
Definition~\ref{def.skewing} to $A = \Lambda$).
Explicitly, this operator is given by
\begin{equation}
f^\perp(a) = \sum (f, a_1) a_2 \qquad \text{ whenever } \qquad
\Delta(a) = \sum a_1 \otimes a_2 ,
\label{eq.skewing.Lambda.general}
\end{equation}
where the inner product
$(f, a_1)$ is now understood as a Hall inner product.

Recall that each partition $\lambda$ satisfies
\[
\Delta s_\lambda
= \sum_{\mu \subseteq \lambda} s_\mu \otimes s_{\lambda/\mu}
= \sum_{\nu \subseteq \lambda} s_\nu \otimes s_{\lambda/\nu}
= \sum_{\nu} s_\nu \otimes s_{\lambda/\nu}
\]
(since $s_{\lambda/\nu} = 0$ unless $\nu \subseteq \lambda$).
Hence, for any two partitions $\lambda$ and $\mu$, we have
\begin{align}
s_\mu^\perp \left(s_\lambda\right)
&= \sum_\nu \underbrace{\left(s_\mu, s_\nu\right)}_{= \delta_{\mu, \nu}} s_{\lambda/\nu}
\qquad
\left(\text{by \eqref{eq.skewing.Lambda.general}, applied to }
      f = s_\mu \text{ and } a = s_\lambda \right)
\nonumber\\
&= \sum_\nu \delta_{\mu, \nu} s_{\lambda/\nu}
= s_{\lambda/\mu} .
\label{eq.skewing.Lambda.s}
\end{align}
Thus, skewing acts on the Schur functions exactly as desired.

% [DG][v46] Added footnote on how $\Lambda$ becomes selfdual.

% [DG][v64] Uncrowded the above paragraph, hopefully making it
% a lot clearer.

\begin{proposition}
\label{skewing-properties-prop}
Let $A$ be a graded Hopf algebra.
The $f^\perp$ operators $A \rightarrow A$ have the following properties.
\begin{enumerate}
\item[(i)] For every $f \in A^o$, the map $f^\perp$ is adjoint to left multiplication $A^o \overset{f \cdot}{\rightarrow} A^o$ in the sense that
\[
(g, f^\perp(a) ) = (fg, a).
\]
\item[(ii)] For every $f, g \in A^o$, we have $(fg)^\perp(a) = g^\perp(f^\perp(a))$, that is,
$A$ becomes a right $A^o$-module via the $f^\perp$ action.\footnote{This
makes sense, since $A^o$ is a $\kk$-algebra (by
Exercise~\ref{exe.dual-algebra}(c), applied to $C=A$).}
\item[(iii)] The unity $1_{A^o}$ of the $\kk$-algebra $A^o$
satisfies $\left(1_{A^o}\right)^\perp = \id_A$.
\item[(iv)] Assume that $A$ is of finite type (so $A^o$ becomes a
Hopf algebra, not just an algebra).
If an $f \in A^o$ satisfies $\Delta(f) = \sum f_1 \otimes f_2$, then
\[
f^\perp (a b) = \sum f_1^\perp(a) f_2^\perp(b).
\]
In particular, if $f$ is primitive in $A^o$,
so that $\Delta(f) = f \otimes \one + \one \otimes f$,
then $f^\perp$ is a \emph{derivation}:
\[
f^\perp (a b) = f^\perp(a) \cdot b + a \cdot f^\perp(b).
\]
\end{enumerate}
\end{proposition}

% [DG][v75] Added part (iii) to this proposition (we
% are already using it tacitly).

\begin{proof}
For (i), note that
\[
(g, f^\perp(a) )
=\sum (f,a_1) (g,a_2)
=( f \otimes g , \Delta_A (a) )
=( m_{A^o}(f \otimes g) , a )
= (fg, a).
\]

For (ii), using (i) and considering any $h$ in $A^o$, one has that
\[
(h, (fg)^\perp(a))
= (fgh, a)
= (gh, f^\perp(a))
= (h, g^\perp(f^\perp(a))).
\]

For (iii), we recall that the unity $1_{A^o}$ of $A^o$ is the counit $\epsilon$ of $A$,
and thus every $a \in A$ satisfies
\begin{align*}
\left(1_{A^o}\right)^\perp \left(a\right)
&= \epsilon^\perp \left(a\right)
= \sum_{\left(a\right)} \underbrace{\left(\epsilon, a_1\right)}_{=\epsilon\left(a_1\right)} a_2
\qquad \left(\text{by the definition of $\epsilon^\perp$}\right) \\
&= \sum_{\left(a\right)} \epsilon\left(a_1\right) a_2
= a \qquad \left(\text{by the axioms of a coalgebra}\right) ,
\end{align*}
so that $\left(1_{A^o}\right)^\perp = \id_A$.

For (iv), noting that
\[
\Delta(ab)= \Delta(a)\Delta(b)
= \left( \sum_{(a)} a_1 \otimes a_2 \right)
\left( \sum_{(b)} b_1 \otimes b_2 \right)
= \sum_{(a),(b)} a_1 b_1 \otimes a_2 b_2,
\]
one has that
\begin{align*}
f^\perp(ab)
&= \sum_{(a),(b)} (f,a_1 b_1)_A \,\, a_2 b_2
= \sum_{(a),(b)} (\Delta(f),a_1 \otimes b_1)_{A \otimes A} \,\, a_2 b_2 \\
&= \sum_{(f),(a),(b)} (f_1,a_1)_A (f_2,b_1)_A \,\, a_2 b_2 \\
&= \sum_{(f)} \left( \sum_{(a)} (f_1,a_1)_A a_2 \right)  \left( \sum_{(b)} (f_2,b_1)_A  b_2 \right)
=\sum_{(f)} f_1^\perp(a) f_2^\perp(b).
\end{align*}

%considering any $g$ in $A^o$, one has that
%\begin{align*}
%(g,f^\perp(ab))
%&= (fg, ab)_A \\
%&= (fg, m_A(a \otimes b))_A \\
%&= ( \Delta_{A^o}(fg), a \otimes b )_{A \otimes A}\\
%&= ( \Delta_{A^o}(f) \Delta_{A^o}(g), a \otimes b )_{A \otimes A}\\
%&= \left( \left( \sum f_1 \otimes f_2 \right) \Delta_{A^o}(g), a \otimes b \right)_{A \otimes A}\\
%&= \left(  \Delta_{A^o}(g), \sum f_1^\perp(a) \otimes f_2^\perp(b) \right)_{A \otimes A}\\
%&= \left(  g, \sum f_1^\perp(a) f_2^\perp(b) \right)_A.\\
%\end{align*}
\end{proof}

The Pieri rules (Theorem~\ref{Pieri-rules}) expressed
multiplication by $h_n$ or by $e_n$ in the basis
$\left(s_\lambda\right)_{\lambda \in \Par}$ of $\Lambda$.
We can similarly express skewing by $h_n$ or by $e_n$:

\begin{proposition}
\label{prop.coPieri-rules}
For every partition $\lambda$ and any $n \in \NN$, we have
\begin{align}
h_n^\perp s_\lambda =&
\sum\limits_{\substack{\lambda^-: \lambda/\lambda^- \text{ is a }\\\text{horizontal }n\text{-strip}}}  s_{\lambda^-} ;
\label{eq.copieri.h} \\
e_n^\perp s_\lambda =&
\sum\limits_{\substack{\lambda^-: \lambda/\lambda^- \text{ is a }\\\text{vertical }n\text{-strip}}}  s_{\lambda^-} .
\label{eq.copieri.e}
\end{align}
\end{proposition}

\begin{exercise}
\label{exe.Lambda.coPieri-rules}
Prove Proposition~\ref{prop.coPieri-rules}.

[\textbf{Hint:} Use Theorem~\ref{Pieri-rules}
and $(s_{\mu^-},e_n^\perp s_\mu) = (e_n s_{\mu^-}, s_\mu)$.]
\end{exercise}

% [DG][v79] Added the above proposition and exercise
% (or, rather, moved them hither from the solutions)
% to better bring the point across that skewing is
% like "multiplication backwards", and to make the
% proof of the skew Pieri rules below a little bit
% clearer.

\noindent
The following interaction between multiplication and $h^\perp$ is the key to
deducing the skew Pieri formula from the usual Pieri formulas.
\begin{lemma}
\label{Lam's-lemma}
For any $f,g$ in $\Lambda$ and any $n \in \NN$, one has
\[
f \cdot h_n^\perp(g)
= \sum_{k=0}^n (-1)^k h_{n-k}^\perp(e_k^\perp(f) \cdot g).
\]
\end{lemma}
\begin{proof}
Starting with the right side, first apply
Proposition~\ref{skewing-properties-prop}(iv):
\begin{align*}
&\sum_{k=0}^n (-1)^k \underbrace{h_{n-k}^\perp(e_k^\perp(f) \cdot g)}_{\substack{= \sum_{j=0}^{n-k} h_j^\perp(e_k^\perp(f)) \cdot h_{n-k-j}^\perp(g) \\ \text{(by Proposition~\ref{skewing-properties-prop}(iv), applied }\\ \text{to } h_{n-k} \text{, } e_k^\perp(f) \text{ and } g \text{ instead of } f \text{, } a \text{ and } b \text{)}}}
\\
&=\sum_{k=0}^n (-1)^{k} \sum_{j=0}^{n-k}
                       h_j^\perp(e_k^\perp(f))
                      \cdot h_{n-k-j}^\perp(g)\\
&=\sum_{j=0}^n \sum_{k=0}^{n-j} (-1)^{k}
                       h_j^\perp(e_k^\perp(f))
                      \cdot h_{n-k-j}^\perp(g)\\
&=\sum_{j=0}^n
       \sum_{i=0}^{n-j} (-1)^{n-i-j}
                      h_j^\perp(e_{n-i-j}^\perp(f))
                      \cdot h_i^\perp(g) \qquad (\text{reindexing }i:=n-k-j \text{ in the inner sum}\,\,) \\
&=\sum_{i=0}^n (-1)^{n-i}
       \left( \sum_{j=0}^{n-i}
                      (-1)^j h_j^\perp(e_{n-i-j}^\perp(f)) \right)
                      \cdot h_i^\perp(g) \\
&=\sum_{i=0}^n (-1)^{n-i}
       \left( \sum_{j=0}^{n-i} (-1)^j e_{n-i-j} h_j\right)^\perp \left( f \right)
                      \cdot h_i^\perp(g) \qquad
                       (\text{by Proposition~\ref{skewing-properties-prop}(ii)}\,\,) \\
&=1^\perp(f) \cdot h_n^\perp(g) = f \cdot h_n^\perp(g)
\end{align*}
where the second-to-last equality used \eqref{e-h-relation}.
\end{proof}

\begin{proof}[Proof of Theorem~\ref{Assaf-McNamara-skew-Pieri}]
We prove \eqref{eq.Assaf-McNamara-skew-Pieri.1}; the
equality \eqref{eq.Assaf-McNamara-skew-Pieri.2} is analogous, swapping
$h_i \leftrightarrow e_i$ and swapping the words ``vertical'' $\leftrightarrow$ ``horizontal''.
%Symmetry of $(\cdot,\cdot)_\Lambda$, commutativity of $\Lambda$, and
%Proposition~\ref{skewing-properties-prop}(i) imply for $f$ in $\Lambda$,
For any $f \in \Lambda$, we have
\begin{align}
\left( s_{\lambda/\mu}, f \right)
&= \left( s_\mu^\perp(s_\lambda),f \right)
\qquad \left(\text{by \eqref{eq.skewing.Lambda.s}}\right) \nonumber\\
&= \left( f, s_\mu^\perp(s_\lambda) \right)
\qquad \left(\text{by symmetry of } (\cdot,\cdot)_\Lambda \right) \nonumber\\
&= \left( s_\mu f , s_\lambda \right)
\qquad \left(\text{by Proposition~\ref{skewing-properties-prop}(i)} \right) \nonumber\\
&= \left( s_\lambda ,s_\mu f \right)
\qquad \left(\text{by symmetry of } (\cdot,\cdot)_\Lambda \right) .
\label{skewing-versus-multiplication-in-Sym}
\end{align}
Hence for any $g$ in $\Lambda$, one can compute that
\begin{align}
( h_n s_{\lambda/\mu} \,\, ,\,\, g )
&\underset{~\ref{skewing-properties-prop}(i)}{\overset{Prop.}{=}}
( s_{\lambda/\mu} \,\, ,\,\, h_n^\perp g )
\overset{\eqref{skewing-versus-multiplication-in-Sym}}{=}
 ( s_{\lambda}\,\, , \,\, s_\mu \cdot h_n^\perp g ) \nonumber \\
&\underset{\ref{Lam's-lemma}}{\overset{Lemma}{=}}
 \sum_{k=0}^n (-1)^k ( s_\lambda \,\, , \,\,
                 h_{n-k}^\perp (e_k^\perp(s_\mu) \cdot g ) ) \nonumber \\
&\underset{~\ref{skewing-properties-prop}(i)}{\overset{Prop.}{=}}
\sum_{k=0}^n (-1)^k ( h_{n-k} s_\lambda \,\, , e_k^\perp(s_\mu) \cdot g ) .
\label{Lam-proof-calculation}
\end{align}
The first Pieri rule in Theorem~\ref{Pieri-rules} lets one rewrite
$h_{n-k} s_\lambda = \sum_{\lambda^+} s_{\lambda^+}$, with the sum running through
$\lambda^+$ for which $\lambda^+/\lambda$ is a horizontal $(n-k)$-strip.
Meanwhile, \eqref{eq.copieri.e} lets one rewrite
$e_k^\perp s_\mu = \sum_{\mu^-} s_{\mu^-}$, with the sum running through
$\mu^-$ for which $\mu/\mu^-$ is a vertical $k$-strip.
% since $(s_{\mu^-},e_k^\perp s_\mu) = (e_k s_{\mu^-}, s_\mu)$.
Thus the right hand side of \eqref{Lam-proof-calculation} becomes
\[
\sum_{k=0}^n (-1)^k \left( \sum_{\lambda^+}  s_{\lambda^+} \,\, , \,\,
                         \sum_{\mu^-} s_{\mu^-} \cdot g \right) \\
\overset{\eqref{skewing-versus-multiplication-in-Sym}}{=}
\left( \sum_{k=0}^n  (-1)^k \sum_{(\lambda^+,\mu^-)}  s_{\lambda^+/\mu^-} , g \right)
\]
where the sum is over the pairs $(\lambda^+, \mu^-)$
for which $\lambda^+/\lambda$ is a horizontal $(n-k)$-strip and
$\mu/\mu^-$ is a vertical $k$-strip.
This proves \eqref{eq.Assaf-McNamara-skew-Pieri.1}.
\end{proof}

\begin{exercise}
\label{exe.aguiar.appli-omega}
Let $n \in \NN$.

\begin{itemize}

\item[(a)]
For every $k \in \NN$,
let $p\left(n,k\right)$ denote the number of partitions of $n$ of
length $k$. Let $c\left(n\right)$ denote the number of
\emph{self-conjugate}\index{self-conjugate partition}
partitions of $n$ (that is, partitions $\lambda$ of $n$ satisfying
$\lambda^t = \lambda$). Show that
\[
\left(-1\right)^n c\left(n\right) = \sum_{k=0}^n \left(-1\right)^k p\left(n,k\right).
\]
(This application of Hopf algebras was found by Aguiar and Lauve,
\cite[\S 5.1]{AguiarLauve-2014}. See also
\cite[Chapter 1, Exercise 22(b)]{Stanley} for an elementary proof.)

\item[(b)]
For every partition $\lambda$, let
$C\left(\lambda\right)$ denote the number of corner cells of the
Ferrers diagram of $\lambda$ (these are the cells of the Ferrers
diagram whose neighbors to the east and to the south both lie
outside of the Ferrers diagram). For every partition $\lambda$, let
$\mu_1\left(\lambda\right)$ denote the number of parts of $\lambda$
equal to $1$. Show that
\[
\sum_{\lambda \in \Par_n} C\left(\lambda\right)
= \sum_{\lambda \in \Par_n} \mu_1\left(\lambda\right).
\]
(This is also due to Stanley.)

\end{itemize}
\end{exercise}

% [DG][v9] Aguiar/Lauve give some references for part (a) which don't use
% Hopf algebras. Part (b) works in greater generality: For any given
% integer $k$, the sum over all $\lambda \in \Par_n$ of the number
% of length-$k$ hooks in $\lambda$ equals the sum over all
% $\lambda \in \Par_n$ of the number of parts of $\lambda$ equal to
% $k$. This was an exercise Stanley gave in his class last term which
% I couldn't solve; yesterday I rediscovered it (with a 5-line
% solution) in looking for applications of Hopf algebras...
% (How) Should we reference these sources? (The more general version
% of (b) could be introduced if we ever make a section on the
% Murnaghan-Nakayama rule.)

% [VR][v10] I think a good reference for part (a) of the exercise is 
% Stanley, Enum. Comb. Vol 1 Chapter 1 Exercise 9(b), 
% which he solves with some easy generating functions.
%
% For part (b) of the exercise, what about Proposition 4.7 on page 943 of  
% Stanley's "Differential posets", J. Amer. Math. Soc. 1 (1988) ?

\begin{exercise}
\label{Darij's-skewing-exercise}
The goal of this exercise is to prove
\eqref{skew-Schur-antipodes} using the skewing operators that
we have developed.\footnote{Make sure not to use the results
of Exercise~\ref{exe.pieri.omega} or Exercise~\ref{exe.cauchy.dual}
or Exercise~\ref{exe.hall-e-f} here, or anything else that relied
on \eqref{skew-Schur-antipodes}, in order to avoid circular
reasoning.}
Recall the involution $\omega:\Lambda\rightarrow\Lambda$ defined in
\eqref{the-fundamental-involution-definition}.
\begin{itemize}
\item[(a)]
Show that $\omega\left( p_{\lambda} \right)
 = \left( -1\right)^{\left\vert \lambda\right\vert
       -\ell\left( \lambda\right) } p_{\lambda}$
for any $\lambda\in\Par$, where
$\ell\left( \lambda\right)$ denotes the length of the partition $\lambda$.
\item[(b)]
Show that $\omega$ is an isometry.
\item[(c)]
Show that this same map $\omega:\Lambda\rightarrow\Lambda$ is a
Hopf automorphism.
\item[(d)]
Prove that $\omega \left( a^{\perp} b \right) = \left(
\omega \left( a\right) \right) ^{\perp} \left( \omega \left( b \right)
\right)$ for every $a\in\Lambda$ and $b\in\Lambda$.
\item[(e)] For any partition
$\lambda = \left( \lambda_1, \ldots, \lambda_\ell \right)$
with length $\ell \left( \lambda \right) = \ell$, prove that
\[
e_{\ell}^{\perp} s_{\lambda}
= s_{\left( \lambda_1 - 1, \lambda_2 - 1, \ldots, \lambda_{\ell} - 1 \right)} .
\]
\item[(f)] For any partition
$\lambda = \left( \lambda_1, \lambda_2, \ldots \right)$, prove that
\[
h_{\lambda_1}^{\perp} s_{\lambda}
= s_{\left( \lambda_2 , \lambda_3 , \lambda_4 , \ldots \right)}.
\]
\item[(g)]
Prove \eqref{skew-Schur-antipodes}.
\end{itemize}
\end{exercise}

\begin{exercise}
\label{exe.skew-p-by-e}Let $n$ be a positive integer. Prove the following:

\begin{enumerate}
\item[(a)] We have $\left(  e_{n},p_{n}\right)  =\left(  -1\right)  ^{n-1}$.

\item[(b)] We have $\left(  e_{m},p_{n}\right)  =0$ for each $m\in \NN$
satisfying $m\neq n$.

\item[(c)] We have $e_{n}^{\perp}p_{n}=\left(  -1\right)  ^{n-1}$.

\item[(d)] We have $e_{m}^{\perp}p_{n}=0$ for each positive integer $m$
satisfying $m\neq n$.
\end{enumerate}
\end{exercise}

% [DG][v75] Added the above exercise.

\subsection{Assorted exercises on symmetric functions}

% [DG][v25] Made this a separate subsection. It doesn't really live
% up to its name so far, but it might grow.

Over a hundred exercises on symmetric functions are collected in
Stanley's \cite[chapter 7]{Stanley}, and even more (but without any
hints or references) on his website\footnote{\ 
\url{http://math.mit.edu/~rstan/ec/ch7supp.pdf}}.
Further sources for results related to symmetric functions are
Macdonald's work, including his monograph \cite{Macdonald}
and his expository \cite{Macdonald-variations}. In this section,
we gather a few exercises that are not too difficult to handle
with the material given above.

\begin{exercise} \phantomsection
\label{exe.bernstein}
\begin{itemize}

\item[(a)] Let $m \in \ZZ$. Prove that, for every $f \in \Lambda$,
the infinite sum
$\sum_{i \in \NN} \left(-1\right)^i h_{m+i} e_i^\perp f$ is
convergent in the discrete topology (i.e., all but finitely many
addends of this sum are zero). Hence, we can define a map
$\mathbf{B}_m : \Lambda \to \Lambda$ by setting
\[
\mathbf{B}_m\left(f\right)
= \sum_{i \in \NN} \left(-1\right)^i h_{m+i} e_i^\perp f
\qquad \qquad \text{ for all } f \in \Lambda .
\]
Show that this map $\mathbf{B}_m$ is $\kk$-linear.

\item[(b)] Let $\lambda = \left(\lambda_1, \lambda_2, \lambda_3,
\ldots\right)$ be a partition, and let $m \in \ZZ$ be such
that $m \geq \lambda_1$. Show that
\[
\sum_{i \in \NN} \left(-1\right)^i h_{m+i} e_i^\perp s_\lambda
= s_{\left(m, \lambda_1, \lambda_2, \lambda_3, \ldots\right)} .
\]

\item[(c)] Let $n \in \NN$. For every $n$-tuple
$\left(\alpha_1, \alpha_2, \ldots, \alpha_n\right) \in \ZZ^n$,
we define an element
$\overline{s}_{\left(\alpha_1, \alpha_2, \ldots, \alpha_n\right)}
\in \Lambda$ by
\[
\overline{s}_{\left(\alpha_1, \alpha_2, \ldots, \alpha_n\right)}
= \det\left( \left( h_{\alpha_i - i + j} \right)_{i, j = 1, 2, \ldots, n} \right) .
\]
Show that
\begin{equation}
\label{eq.exe.bernstein.c.1}
s_\lambda
= \overline{s}_{\left(\lambda_1, \lambda_2, \ldots, \lambda_n\right)}
\end{equation}
for every partition
$\lambda = \left(\lambda_1, \lambda_2, \lambda_3, \ldots\right)$
having at most $n$ parts\footnote{Recall that a \emph{part}
of a partition means a nonzero entry of the partition.}.

Furthermore, show that for
every $n$-tuple
$\left(\alpha_1, \alpha_2, \ldots, \alpha_n\right) \in \ZZ^n$,
the symmetric function
$\overline{s}_{\left(\alpha_1, \alpha_2, \ldots, \alpha_n\right)}$
either is $0$ or equals $\pm s_\nu$ for some partition $\nu$ having
at most $n$ parts.

Finally, show that for any $n$-tuples
$\left(\alpha_1, \alpha_2, \ldots, \alpha_n\right) \in \ZZ^n$
and
$\left(\beta_1, \beta_2, \ldots, \beta_n\right) \in \NN^n$,
we have
\begin{equation}
\label{eq.exe.bernstein.c.3}
\overline{s}_{\left(\beta_1, \beta_2, \ldots, \beta_n\right)}^\perp
\overline{s}_{\left(\alpha_1, \alpha_2, \ldots, \alpha_n\right)}
= \det\left( \left( h_{\alpha_i - \beta_j - i + j} \right)_{i, j = 1, 2, \ldots, n} \right) .
\end{equation}

\item[(d)] For every $n \in \NN$, every $m \in \ZZ$ and every
$n$-tuple
$\left(\alpha_1, \alpha_2, \ldots, \alpha_n\right) \in \ZZ^n$,
prove that
\begin{equation}
\label{eq.exe.bernstein.d}
\sum_{i \in \NN} \left(-1\right)^i h_{m+i} e_i^\perp
\overline{s}_{\left(\alpha_1, \alpha_2, \ldots, \alpha_n\right)}
= \overline{s}_{\left(m, \alpha_1, \alpha_2, \ldots,
\alpha_n\right)} ,
\end{equation}
where we are using the notations of
Exercise~\ref{exe.bernstein}(c).

\item[(e)] For every $n \in \NN$ and every $n$-tuple
$\left(\alpha_1, \alpha_2, \ldots, \alpha_n\right) \in \ZZ^n$,
prove that
\[
\overline{s}_{\left(\alpha_1, \alpha_2, \ldots, \alpha_n\right)}
= \left(\mathbf{B}_{\alpha_1} \circ \mathbf{B}_{\alpha_2}
\circ \cdots \circ \mathbf{B}_{\alpha_n}\right) \left(1\right) ,
\]
where we are using the notations of
Exercise~\ref{exe.bernstein}(c) and Exercise~\ref{exe.bernstein}(a).

\item[(f)] For every $m \in \ZZ$ and every positive integer $n$,
prove that
$\mathbf{B}_{m} \left( p_n \right) = h_m p_n - h_{m+n}$.
Here, we are using the notations of Exercise~\ref{exe.bernstein}(a).

\end{itemize}
\end{exercise}

\begin{remark}
The map $\mathbf{B}_m$ defined in Exercise~\ref{exe.bernstein}(a) is
the so-called
\emph{$m$-th Bernstein creation operator}\index{Bernstein creation operator};
it appears in
Zelevinsky \cite[\S 4.20(a)]{Zelevinsky} and has been introduced by
J.N. Bernstein, who found the result of Exercise~\ref{exe.bernstein}(b).
It is called a ``Schur row adder'' in \cite{Garsia-adder}.
Exercise~\ref{exe.bernstein}(e) appears in
Berg/Bergeron/Saliola/Serrano/Zabrocki
\cite[Theorem 2.3]{BergBergeronSaliolaSerranoZabrocki-immac1}, where it
is used as a prototype for defining \emph{noncommutative} analogues
of Schur functions, the so-called \emph{immaculate functions}.
The particular case of Exercise~\ref{exe.bernstein}(e) for
$\left(\alpha_1, \alpha_2, \ldots, \alpha_n\right)$ a partition
of length $n$ (a restatement of Exercise~\ref{exe.bernstein}(b))
is proven in \cite[\S I.5, example 29]{Macdonald}.
\end{remark}

% [DG][v25] Added the above exercise, an occasion to cite the
% immaculate-functions paper (and plug a minor expository gap
% created by its authors referring to Zelevinsky for a theorem
% which is more general than what Zelevinsky states).
% 
% When I started writing the solution to this problem, I was
% expecting that the Jacobi-Trudi solution to part (b) would
% immediately generalize to the less restrictive setting of
% (d). It does, but this proved messier than I expected it to be;
% cases need to be distinguished and the claims of (c) need to
% be used.

% [DG][v75] Added part (f) to the above exercise.

\begin{exercise} \phantomsection
\label{exe.witt}
\begin{itemize}

\item[(a)]
Prove that there exists a unique family
$\left( x_n \right) _{n \geq 1}$ of elements of $\Lambda$
such that
\[
H \left( t \right) = \prod_{n=1}^{\infty} \left( 1 - x_n t^n \right) ^{-1}.
\]
Denote this family $\left( x_n \right) _{n \geq 1}$ by
$\left( w_n \right) _{n \geq 1}$.\index{$w_n$}
For instance,
\begin{align*}
w_{1} &=    s_{\left( 1 \right)},\ \ \ \ \ \ \ \ \ \ 
w_{2}  =  - s_{\left( 1,1 \right)},\ \ \ \ \ \ \ \ \ \ 
w_{3}  =  - s_{\left( 2,1 \right)},\\
w_{4} &=  - s_{\left( 1,1,1,1 \right)} - s_{\left( 2,1,1 \right)}
          - s_{\left( 2,2 \right)} - s_{\left( 3,1 \right)},\ \ \ \ \ \ \ \ \ \ 
w_{5}  =  - s_{\left( 2,1,1,1 \right)} - s_{\left( 2,2,1 \right)}
          - s_{\left( 3,1,1 \right)} - s_{\left( 3,2 \right)}
          - s_{\left( 4,1 \right)}.
\end{align*}

\item[(b)]
Show that $w_n$ is homogeneous of degree $n$ for every positive
integer $n$.

\item[(c)]
For every partition $\lambda$, define $w_{\lambda} \in \Lambda$ by
$w_{\lambda} = w_{\lambda_1} w_{\lambda_2} \cdots w_{\lambda_{\ell}}$
(where $\lambda = \left( \lambda_1, \lambda_2, \ldots,
\lambda_{\ell} \right)$ with $\ell = \ell\left( \lambda \right)$).
Notice that $w_{\lambda}$ is homogeneous of degree
$\left\vert \lambda\right\vert$. Prove that $\sum_{\lambda \in \Par_n}
w_{\lambda} = h_n$ for every $n\in\NN$.

\item[(d)]
Show that $\left\{ w_{\lambda} \right\} _{\lambda\in \Par}$
is a $\kk$-basis of $\Lambda$. (This basis is
called the \dfn{Witt basis}\footnote{This is due to its relation
with Witt vectors in the appropriate sense. Most
of the work on this basis has been done by Reutenauer and
Hazewinkel.}; it is studied in
\cite[\S 9-\S 10]{Hazewinkel2}.\footnote{It also implicitly
appears in \cite[\S 5]{BakerRichter}. Indeed, the $q_n$ of
\cite{BakerRichter} are our $w_n$ (for $\kk = R$).})

\item[(e)]
Prove that $p_n = \sum_{d \mid n} d w_d^{n/d}$ for every positive
integer $n$. (Here, the summation sign $\sum_{d \mid n}$ means a sum
over all positive divisors $d$ of $n$.)

\item[(f)]
We are going to show that $-w_n$ is a sum of Schur functions
(possibly with repetitions, but without signs!) for every $n \geq 2$. (For
$n=1$, the opposite is true: $w_{1}$ is a single Schur function.) This proof
goes back to Doran \cite{Doran}\footnote{See also Stanley
\cite[Exercise 7.46]{Stanley}.}.

For any positive integers $n$ and $k$, define $f_{n,k} \in \Lambda$ by
$f_{n,k} = \sum\limits_{\substack{\lambda\in \Par_n,\\\min\lambda\geq k}}w_{\lambda}$,
where $\min\lambda$ denotes the smallest
part\footnote{Recall that a \emph{part} of a partition means
a nonzero entry of the partition.} of $\lambda$. Show that
\[
-f_{n,k} = s_{\left( n-1,1 \right)} + \sum_{i=2}^{k-1} f_{i,i} f_{n-i,i}
\ \ \ \ \ \ \ \ \ \ \text{for every } n\geq k\geq 2.
\]
Conclude that $-f_{n,k}$ is a sum of Schur functions for every $n\in
\NN$ and $k\geq 2$. Conclude that $-w_n$ is a sum of Schur functions
for every $n\geq 2$.

\item[(g)]
For every partition $\lambda$, define $r_{\lambda}\in\Lambda$ by
$r_{\lambda}=\prod_{i\geq 1}h_{v_{i}}\left(  x_{1}^{i},x_{2}^{i},x_{3}
^{i},\ldots\right)  $, where $v_{i}$ is the number of occurrences of $i$ in
$\lambda$. Show that $\sum_{\lambda\in\Par}w_{\lambda}\left(
\xx\right)  r_{\lambda}\left(  \yy\right)  =\prod_{i,j=1}
^{\infty}\left(  1-x_{i}y_{j}\right)  ^{-1}$.

\item[(h)]
Show that $\left\{ r_\lambda \right\} _{\lambda \in \Par}$ and
$\left\{ w_\lambda \right\} _{\lambda \in \Par}$ are dual bases of
$\Lambda$.

\end{itemize}
\end{exercise}

% [DG][v7] Added exercise above.

\begin{exercise}
\label{exe.Lambda.maps-on-pn}
For this exercise, set $\kk = \ZZ$, and consider $\Lambda =
\Lambda_\ZZ$ as a subring of $\Lambda_\QQ$. Also,
consider $\Lambda \otimes_\ZZ \Lambda$ as a subring of
$\Lambda_\QQ \otimes_\QQ \Lambda_\QQ$. \ \ \ \ \footnote{Here
is how this works: We have $\Lambda_\QQ \cong \QQ \otimes_\ZZ
\Lambda$. But fundamental properties of tensor products yield
\begin{equation}
\label{eq.exe.Lambda.maps-on-pn.inj}
\QQ \otimes_\ZZ \left(\Lambda \otimes_\ZZ \Lambda\right)
\cong
\underbrace{\left(\QQ \otimes_\ZZ \Lambda\right)}_{\cong \Lambda_\QQ}
\otimes_\QQ
\underbrace{\left(\QQ \otimes_\ZZ \Lambda\right)}_{\cong \Lambda_\QQ}
\cong \Lambda_\QQ \otimes_\QQ \Lambda_\QQ
\end{equation}
as $\QQ$-algebras.
But $\Lambda \otimes_\ZZ \Lambda$ is a free $\ZZ$-module
(since $\Lambda$ is a free $\ZZ$-module), and so the canonical
ring homomorphism $\Lambda \otimes_\ZZ \Lambda \to
\QQ \otimes_\ZZ \left(\Lambda \otimes_\ZZ \Lambda\right)$
sending every $u$ to $1_\QQ \otimes_\ZZ u$ is injective.
Composing this ring homomorphism with the $\QQ$-algebra
isomorphism of \eqref{eq.exe.Lambda.maps-on-pn.inj} gives
an injective ring homomorphism
$\Lambda \otimes_\ZZ \Lambda \to
\Lambda_\QQ \otimes_\QQ \Lambda_\QQ$. We use this latter
homomorphism
to identify $\Lambda \otimes_\ZZ \Lambda$ with a subring of
$\Lambda_\QQ \otimes_\QQ \Lambda_\QQ$.}
Recall that the family
$\left(p_n\right)_{n \geq 1}$ generates the $\QQ$-algebra
$\Lambda_\QQ$, but does not generate the $\ZZ$-algebra
$\Lambda$.
\begin{itemize}

\item[(a)] Define a $\QQ$-linear map $Z : \Lambda_\QQ \to
\Lambda_\QQ$ by setting
\[
Z\left(p_\lambda\right) = z_\lambda p_\lambda
\qquad \qquad \text{for every partition } \lambda ,
\]
where $z_\lambda$ is defined as in
Proposition~\ref{prop-Two-other-Cauchy-expansions}.\footnote{This is
well-defined, since $\left(p_\lambda\right)_{\lambda \in \Par}$ is
a $\QQ$-module basis of $\Lambda_\QQ$.} Show that
$Z \left( \Lambda \right) \subset \Lambda$.

\item[(b)] Define a $\QQ$-algebra homomorphism
$\Delta_\times : \Lambda_\QQ \to \Lambda_\QQ \otimes_\QQ \Lambda_\QQ$
by setting
\[
\Delta_\times \left(p_n\right) = p_n \otimes p_n
\qquad \qquad \text{for every positive integer } n.
\]
\footnote{This is well-defined, since the family
$\left(p_n\right)_{n \geq 1}$ generates the $\QQ$-algebra
$\Lambda_\QQ$ and is algebraically independent.} Show that
$\Delta_\times \left( \Lambda \right) \subset
\Lambda \otimes_\ZZ \Lambda$.

\item[(c)] Let $r \in \ZZ$. Define a $\QQ$-algebra homomorphism
$\epsilon_r : \Lambda_\QQ \to \QQ$
by setting
\[
\epsilon_r \left(p_n\right) = r
\qquad \qquad \text{for every positive integer } n.
\]
\footnote{This is well-defined, since the family
$\left(p_n\right)_{n \geq 1}$ generates the $\QQ$-algebra
$\Lambda_\QQ$ and is algebraically independent.} Show that
$\epsilon_r \left( \Lambda \right) \subset \ZZ$.

\item[(d)] Let $r \in \ZZ$. Define a $\QQ$-algebra homomorphism
$\mathbf{i}_r : \Lambda_\QQ \to \Lambda_\QQ$
by setting
\[
\mathbf{i}_r \left(p_n\right) = r p_n
\qquad \qquad \text{for every positive integer } n.
\]
\footnote{This is well-defined, since the family
$\left(p_n\right)_{n \geq 1}$ generates the $\QQ$-algebra
$\Lambda_\QQ$ and is algebraically independent.} Show that
$\mathbf{i}_r \left( \Lambda \right) \subset \Lambda$.

\item[(e)] Define a $\QQ$-linear map $\operatorname{Sq} : \Lambda_\QQ \to
\Lambda_\QQ$ by setting
\[
\operatorname{Sq}\left(p_\lambda\right) = p_\lambda^2
\qquad \qquad \text{for every partition } \lambda .
\]
\footnote{This is
well-defined, since $\left(p_\lambda\right)_{\lambda \in \Par}$ is
a $\QQ$-module basis of $\Lambda_\QQ$.} Show that
$\operatorname{Sq} \left( \Lambda \right) \subset \Lambda$.

\item[(f)] Let $r \in \ZZ$. Define a $\QQ$-algebra homomorphism
$\Delta_r : \Lambda_\QQ \to \Lambda_\QQ \otimes_\QQ \Lambda_\QQ$
by setting
\[
\Delta_r \left(p_n\right) = \sum_{i=1}^{n-1} \binom{n}{i}
p_i \otimes p_{n-i} + r \otimes p_n + p_n \otimes r
\qquad \qquad \text{for every positive integer } n.
\]
\footnote{This is well-defined, since the family
$\left(p_n\right)_{n \geq 1}$ generates the $\QQ$-algebra
$\Lambda_\QQ$ and is algebraically independent.} Show that
$\Delta_r \left( \Lambda \right) \subset
\Lambda \otimes_\ZZ \Lambda$.

\item[(g)] Consider the map $\Delta_\times$ introduced in
Exercise~\ref{exe.Lambda.maps-on-pn}(b) and the map
$\epsilon_1$ introduced in
Exercise~\ref{exe.Lambda.maps-on-pn}(c). Show that the
$\QQ$-algebra $\Lambda_\QQ$, endowed with the
comultiplication $\Delta_\times$ and the counit
$\epsilon_1$, becomes a cocommutative
$\QQ$-bialgebra.\footnote{But unlike $\Lambda_\QQ$ with the
usual coalgebra structure, it is neither graded nor a Hopf
algebra.}

\item[(h)] Define a $\QQ$-bilinear map $* : \Lambda_\QQ
\times \Lambda_\QQ \to \Lambda_\QQ$, which will be written
in infix notation (that is, we will write $a * b$
instead of $*\left(a, b\right)$), by setting
\[
p_\lambda * p_\mu = \delta_{\lambda, \mu} z_\lambda p_\lambda
\qquad \qquad \text{for any partitions } \lambda \text{ and }
\mu
\]
(where $z_\lambda$ is defined as in
Proposition~\ref{prop-Two-other-Cauchy-expansions}).
\footnote{This is
well-defined, since $\left(p_\lambda\right)_{\lambda \in \Par}$ is
a $\QQ$-module basis of $\Lambda_\QQ$.}
Show that $f * g \in \Lambda$ for any $f \in \Lambda$
and $g \in \Lambda$.

\item[(i)] Show that $\epsilon_1\left(f\right) =
f\left(1\right)$ for every $f \in \Lambda_\QQ$ (where
we are using the notation $\epsilon_r$ defined in
Exercise~\ref{exe.Lambda.maps-on-pn}(c)).
\end{itemize}

[\textbf{Hint:}
\begin{itemize}
\item For (b), show that, for every $f \in \Lambda_\QQ$,
the tensor $\Delta_\times \left(f\right)$ is the
preimage of
$f\left( \left(x_i y_j\right)_{\left(i,j\right) \in
\left\{1,2,3,\ldots\right\}^2} \right)
= f\left( x_1 y_1, x_1 y_2, x_1 y_3, \ldots,
x_2 y_1, x_2 y_2, x_2 y_3, \ldots,
\ldots \right) \in \QQ\left[\left[\xx,\yy\right]\right]$
under the canonical injection
$\Lambda_\QQ \otimes_\QQ \Lambda_\QQ \to
\QQ\left[\left[\xx,\yy\right]\right]$ which maps every
$f \otimes g$ to $f\left(\xx\right) g\left(\yy\right)$.
(This requires making sure that the evaluation
$f\left( \left(x_i y_j\right)_{\left(i,j\right) \in
\left\{1,2,3,\ldots\right\}^2} \right)$ is well-defined
to begin with, i.e., converges as a formal power series.)
\par For an alternative solution to (b), compute
$\Delta_\times \left(h_n\right)$ or
$\Delta_\times \left(e_n\right)$.
\item For (c), compute $\epsilon_r \left(e_n\right)$
or $\epsilon_r \left(h_n\right)$.
\item Reduce (d) to (b) and (c) using
Exercise~\ref{exe.tensorprod.morphism}.
\item Reduce (e) to (b).
\item (f) is the hardest part. It is tempting to try and
interpret the definition of $\Delta_r$ as a convoluted way
of saying that $\Delta_r \left(f\right)$ is the preimage
of
$f\left( \left(x_i + y_j\right)_{\left(i,j\right) \in
\left\{1,2,3,\ldots\right\}^2} \right)$ under the
canonical injection
$\Lambda_\QQ \otimes_\QQ \Lambda_\QQ \to
\QQ\left[\left[\xx,\yy\right]\right]$ which maps every
$f \otimes g$ to $f\left(\xx\right) g\left(\yy\right)$.
However, this does not make sense since the evaluation
$f\left( \left(x_i + y_j\right)_{\left(i,j\right) \in
\left\{1,2,3,\ldots\right\}^2} \right)$ is (in general)
not well-defined\footnote{e.g., it involves summing
infinitely many $x_1$'s if $f = e_1$} (and even if it
was, it would fail to explain the $r$). So we need to get down
to finitely many variables.
For every $N \in \NN$, define a
$\QQ$-algebra homomorphism $\mathcal{E}_N :
\Lambda_\QQ \otimes_\QQ \Lambda_\QQ \to
\QQ\left[x_1, x_2, \ldots, x_N, y_1, y_2, \ldots, y_N\right]$
by sending each $f \otimes g$ to
$f\left(x_1, x_2, \ldots, x_N\right)
g\left(y_1, y_2, \ldots, y_N\right)$. Show that
$\Delta_N \left(\Lambda\right) \subset
\mathcal{E}_N^{-1}\left(
\ZZ\left[x_1, x_2, \ldots, x_N, y_1, y_2, \ldots, y_N\right]
\right)$. This shows that, at least, the coefficients
of $\Delta_r \left(f\right)$ in front of the
$m_\lambda \otimes m_\mu$ with $\ell\left(\lambda\right)
\leq r$ and $\ell\left(\mu\right) \leq r$ (in the
$\QQ$-basis
$\left(m_\lambda \otimes m_\mu\right)_{\lambda, \mu \in \Par}$
of $\Lambda_\QQ \otimes_\QQ \Lambda_\QQ$)
are integral for $f \in \Lambda$. Of course, we want
all coefficients. Show that
$\Delta_a = \Delta_b \star \left(\Delta_{\Lambda_\QQ}
\circ \mathbf{i}_{a-b}\right)$ in
$\Hom\left(\Lambda_\QQ, \Lambda_\QQ \otimes_\QQ
\Lambda_\QQ\right)$ for any integers $a$ and $b$. This
allows ``moving'' the $r$. This approach to (f) was partly
suggested to the first author by Richard Stanley.
\item For (h), notice that Definition~\ref{def.self-duality}(b)
(below) allows us to construct a bilinear form
$\left(\cdot, \cdot\right)_{\Lambda_\QQ \otimes_\QQ \Lambda_\QQ}
: \left(\Lambda_\QQ \otimes_\QQ \Lambda_\QQ\right)
\times \left(\Lambda_\QQ \otimes_\QQ \Lambda_\QQ\right)
\to \QQ$ from the Hall inner product
$\left(\cdot, \cdot\right) : \Lambda_\QQ \times \Lambda_\QQ \to
\QQ$. Show that
\begin{equation}
\label{eq.exe.Lambda.maps-on-pn.adjointness}
\left(a * b, c\right)
= \left(a \otimes b, \Delta_\times \left(c\right) \right)_{\Lambda_\QQ \otimes_\QQ \Lambda_\QQ}
\qquad \qquad \text{ for all } a, b, c \in \Lambda_\QQ ,
\end{equation}
and then use (b).
\end{itemize}
]
\end{exercise}

\begin{remark}
The map $\Delta_\times$ defined in
Exercise~\ref{exe.Lambda.maps-on-pn}(b) is known as the
\dfn{internal comultiplication} (or \dfn{Kronecker
comultiplication}) on $\Lambda_\QQ$. Unlike the standard
comultiplication $\Delta_{\Lambda_\QQ}$, it is not a graded
map, but rather sends every homogeneous component
$\left(\Lambda_\QQ\right)_n$ into
$\left(\Lambda_\QQ\right)_n \otimes \left(\Lambda_\QQ\right)_n$.
The bilinear map $*$ from
Exercise~\ref{exe.Lambda.maps-on-pn}(h) is the so-called
\dfn{internal multiplication} (or \dfn{Kronecker
multiplication}), and is similarly not graded but rather
takes
$\left(\Lambda_\QQ\right)_n \times \left(\Lambda_\QQ\right)_m$
to $\left(\Lambda_\QQ\right)_n$ if $n = m$ and to $0$ otherwise.

The analogy between the two internal structures is not perfect:
While we saw in Exercise~\ref{exe.Lambda.maps-on-pn}(g) how
the internal comultiplication yields another bialgebra structure
on $\Lambda_\QQ$, it is not true that the internal
multiplication (combined with the usual coalgebra structure of
$\Lambda_\QQ$) forms a bialgebra structure as well. What is
missing is a multiplicative unity; if we would take the closure
of $\Lambda_\QQ$ with respect to the grading, then
$1 + h_1 + h_2 + h_3 + \cdots$ would be such a unity.

The structure constants of the internal comultiplication
on the Schur basis $\left(s_\lambda\right)_{\lambda \in \Par}$
are equal to the structure constants
of the internal multiplication on the Schur basis\footnote{This
can be obtained, e.g., from
\eqref{eq.exe.Lambda.maps-on-pn.adjointness}.}, and are
commonly referred to as the
\dfn{Kronecker coefficients}. They are known to be
nonnegative integers (this follows from
Exercise~\ref{exe.repSn.inner-tensor.intprod}(c)\footnote{Their
integrality can also be easily deduced from
Exercise~\ref{exe.Lambda.maps-on-pn}(b).}), but no
combinatorial proof is known for their
nonnegativity. Combinatorial
interpretations for these coefficients akin to the
Littlewood-Richardson rule have been found only in special
cases (cf., e.g., \cite{Rosas-kronprod} and
\cite{Blasiak-kronhook} and \cite{Liu-simplkron}).

The map $\Delta_r$ of
Exercise~\ref{exe.Lambda.maps-on-pn}(f) also has some
classical theory behind it, relating to Chern classes of
tensor products (\cite{Manivel-chern},
\cite[\S I.4, example 5]{Macdonald}).
\end{remark}

% [DG][v25] Added exercise above. Part (f) is the statement of
% http://mathoverflow.net/questions/120924/is-the-renormalized-third-comultiplication-on-mathbfsymm-integral?rq=1
% (but the proof is not the one I mentioned in my answer).
% It took me a while to find this solution and I am not sure
% my hint is anywhere near good...

% [DG][v46] Added part (i) (more or less a triviality).

Parts (b), (c), (d), (e) and (f) of Exercise~\ref{exe.Lambda.maps-on-pn} are
instances of a general phenomenon: Many $\ZZ$-algebra homomorphisms
$\Lambda\rightarrow A$ (with $A$ a commutative ring, usually torsionfree) are
easiest to define by first defining a $\QQ$-algebra homomorphism
$\Lambda_\QQ \rightarrow A \otimes \QQ$ and then showing that
this homomorphism restricts to a $\ZZ$-algebra homomorphism
$\Lambda\rightarrow A$. One might ask for general criteria when this is
possible; specifically, for what choices of $\left( b_n \right)  _{n\geq 1}
\in A^{\left\{  1,2,3,\ldots\right\}  }$ does there exist a
$\ZZ$-algebra homomorphism $\Lambda\rightarrow A$ sending the $p_n$ to
$b_n$ ? Such choices are called \dfn{ghost-Witt vectors} in Hazewinkel
\cite{Hazewinkel2}, and we can give various equivalent conditions for a family
$\left( b_n \right) _{n \geq 1}$ to be a ghost-Witt vector:

\begin{exercise}
\label{exe.witt.ghost-equiv}
Let $A$ be a commutative ring.

For every $n\in\left\{  1,2,3,\ldots\right\} $, let $\varphi_{n}:A\rightarrow
A$ be a ring endomorphism of $A$. Assume that the following properties hold:

\begin{itemize}
\item We have $\varphi_{n}\circ\varphi_{m}=\varphi_{nm}$ for any two positive
integers $n$ and $m$.

\item We have $\varphi_{1}= \id$.

\item We have $\varphi_{p}\left(  a\right)  \equiv a^p \operatorname{mod}pA$
for every $a\in A$ and every prime number $p$.
\end{itemize}

\noindent (For example, when $A = \ZZ$, one can set $\varphi_n = \id$
for all $n$; this simplifies the exercise somewhat. More generally, setting
$\varphi_{n} = \id$ works whenever $A$ is a binomial
ring\footnote{A \dfn{binomial ring} is defined to be a torsionfree (as an
additive group) commutative ring $A$ which has one of the following equivalent
properties:
\par
\begin{itemize}
\item For every $n\in \NN$ and $a\in A$, we have $a\left(  a-1\right)
\cdots\left(  a-n+1\right)  \in n!\cdot A$. (That is, binomial coefficients
$\dbinom{a}{n}$ with $a\in A$ and $n\in \NN$ are defined in $A$.)
\par
\item We have $a^{p}\equiv a\operatorname{mod}pA$ for every $a\in A$ and every
prime number $p$.
\end{itemize}
\par
See \cite{Xantcha} and the references therein for studies of these rings. It
is not hard to check that $\ZZ$ and every localization of $\ZZ$
are binomial rings, and so is any commutative $\QQ$-algebra
as well as the ring
\[
\left\{  P\in \QQ \left[  X\right]  \ \mid\ P\left(  n\right)
\in \ZZ \text{ for every }n\in \ZZ \right\}
\]
(but not the ring $ \ZZ \left[  X\right] $ itself).}. However, the results of
this exercise are at their most useful when $A$ is a multivariate polynomial
ring $ \ZZ \left[  x_{1},x_{2},x_{3},\ldots\right]  $ over $\ZZ$
and the homomorphism $\varphi_{n}$ sends every $P\in A$ to
$P\left(  x_{1}^{n},x_{2}^{n},x_{3}^{n},\ldots\right)  $.)

Let $\mu$ denote the
\dfn{number-theoretic M\"obius function}\index{M\"obius function}\index{$\mu$};
this is
the function $\left\{  1,2,3,\ldots\right\}  \rightarrow \ZZ$ defined by
\[
\mu\left(  m\right)  = \begin{cases}
0,& \text{if }m\text{ is not squarefree;}\\
\left(  -1\right)  ^{\left(  \text{number of prime factors of }m\right)
},& \text{if }m\text{ is squarefree}
\end{cases}
\qquad \qquad \text{for every positive integer } m.
\]

Let \dfn{$\phi$} denote the \dfn{Euler totient function};
this is the function
$\left\{  1,2,3,\ldots\right\}  \rightarrow \NN$ which sends every
positive integer $m$ to the number of elements of $\left\{  1,2,\ldots
,m\right\}  $ coprime to $m$.

Let $\left(  b_{n}\right)  _{n\geq1}\in A^{\left\{  1,2,3,\ldots\right\}  }$ be a
family of elements of $A$. Prove that the following seven assertions are equivalent:

\begin{itemize}
\item \textit{Assertion $\mathcal{C}$:} For every positive integer
$n$ and every prime factor $p$ of $n$, we have
\[
\varphi_{p}\left(  b_{n/p}\right)  \equiv b_{n}
\operatorname{mod} p^{v_{p}\left(  n\right)  }A.
\]
Here, $v_{p}\left(  n\right)  $ denotes the exponent of $p$ in the prime
factorization of $n$.

\item \textit{Assertion $\mathcal{D}$:} There exists a family
$\left(  \alpha_{n}\right)  _{n\geq1}\in A^{\left\{  1,2,3,\ldots\right\}  }$
of elements of $A$ such that every positive integer $n$ satisfies
\[
b_{n}=\sum_{d\mid n}d\alpha_{d}^{n/d}.
\]
\footnote{Here and in the following, summations of the form $\sum_{d\mid n}$
range over all \textbf{positive} divisors of $n$.}

\item \textit{Assertion $\mathcal{E}$:} There exists a family
$\left(  \beta_{n}\right)  _{n\geq1}\in A^{\left\{  1,2,3,\ldots\right\}  }$
of elements of $A$ such that every positive integer $n$ satisfies
\[
b_{n}=\sum_{d\mid n}d\varphi_{n/d}\left(  \beta_{d}\right)  .
\]

\item \textit{Assertion $\mathcal{F}$:} Every positive integer $n$
satisfies
\[
\sum_{d\mid n}\mu\left(  d\right)  \varphi_{d}\left(  b_{n/d}\right)  \in nA.
\]

\item \textit{Assertion $\mathcal{G}$:} Every positive integer $n$
satisfies
\[
\sum_{d\mid n}\phi\left(  d\right)  \varphi_{d}\left(  b_{n/d}\right)  \in
nA.
\]

\item \textit{Assertion $\mathcal{H}$:} Every positive integer $n$
satisfies
\[
\sum_{i=1}^{n}\varphi_{n/\gcd\left(  i,n\right)  }\left(  b_{\gcd\left(
i,n\right)  }\right)  \in nA.
\]

\item \textit{Assertion $\mathcal{J}$:} There exists a ring
homomorphism $\Lambda_{ \ZZ }\rightarrow A$ which, for every positive
integer $n$, sends $p_{n}$ to $b_{n}$.
\end{itemize}

[\textbf{Hint:} The following identities hold for every positive integer $n$:
\begin{align}
\sum_{d\mid n}\phi\left(  d\right)   &
=n;\label{eq.exe.witt.ghost-equiv.phi*1}\\
\sum_{d\mid n}\mu\left(  d\right)   &  =\delta_{n,1}
;\label{eq.exe.witt.ghost-equiv.mu*1}\\
\sum_{d\mid n}\mu\left(  d\right)  \dfrac{n}{d}  &  =\phi\left(  n\right)
;\label{eq.exe.witt.ghost-equiv.mu*n}\\
\sum_{d\mid n}d\mu\left(  d\right)  \phi\left(  \dfrac{n}{d}\right)   &
=\mu\left(  n\right)  . \label{eq.exe.witt.ghost-equiv.nmu*phi}
\end{align}
Furthermore, the following simple lemma is useful: If $k$ is a positive
integer, and if $p\in \NN$, $a\in A$ and $b\in A$ are such that $a\equiv
b\operatorname{mod}p^{k}A$, then $a^{p^{\ell}}\equiv b^{p^{\ell}}
\operatorname{mod}p^{k+\ell}A$ for every $\ell\in \NN$.]
\end{exercise}

\begin{remark}
Much of Exercise \ref{exe.witt.ghost-equiv} is folklore, but it is hard to
pinpoint concrete appearances in literature. The equivalence $\mathcal{C}
\Longleftrightarrow\mathcal{D}$ appears in Hesselholt \cite[Lemma
1]{Hesselholt-witt} and \cite[Lemma 1.1]{Hesselholt-deRhamWitt}
(in slightly greater generality), where it is referred to
as Dwork's lemma and used in the construction of the Witt vector functor. This
equivalence is also \cite[Lemma 9.93]{Hazewinkel2}. The equivalence
$\mathcal{D}\Longleftrightarrow\mathcal{F}\Longleftrightarrow\mathcal{G}
\Longleftrightarrow\mathcal{H}$ in the case $A = \ZZ$ is \cite[Corollary
on p. 10]{DressSiebeneicher-burnside}, where it is put into the context of
Burnside rings and necklace counting. The equivalence $\mathcal{C}
\Longleftrightarrow\mathcal{F}$ for finite families $\left(  b_{n}\right)
_{n\in\left\{  1,2,\ldots,m\right\}  }$ in lieu of $\left(  b_{n}\right)
_{n\geq1}$ is \cite[Exercise 5.2 \textbf{a}]{Stanley}. One of the likely
oldest relevant sources is Schur's \cite{Schur-1937}, which proves the
equivalence $\mathcal{C}\Longleftrightarrow\mathcal{D}\Longleftrightarrow
\mathcal{F}$ for finite families $\left(  b_{n}\right)  _{n\in\left\{
1,2,\ldots,m\right\}  }$, as well as a ``finite
version'' of $\mathcal{C}\Longleftrightarrow\mathcal{J}$
(Schur did not have $\Lambda$, but was working with actual power sums of roots
of polynomials).
\end{remark}

\begin{exercise}
\label{exe.witt.ghost-exa}
Let $A$ denote the ring $\ZZ$. For every
$n\in\left\{  1,2,3,\ldots\right\}  $, let $\varphi_{n}$ denote the identity
endomorphism $\id$ of $A$. Prove that the seven equivalent
assertions $\mathcal{C}$, $\mathcal{D}$, $\mathcal{E}$, $\mathcal{F}$,
$\mathcal{G}$, $\mathcal{H}$ and $\mathcal{J}$ of Exercise
\ref{exe.witt.ghost-equiv} are satisfied for each of the following families
$\left(  b_{n}\right)  _{n\geq1}\in \ZZ ^{\left\{  1,2,3,\ldots\right\}  }$:

\begin{itemize}
\item the family $\left(  b_{n}\right)  _{n\geq1}=\left(  q^{n}\right)
_{n\geq1}$, where $q$ is a given integer.

\item the family
$\left(  b_{n}\right)  _{n\geq1}=\left(  q\right)  _{n\geq1}$, where
$q$ is a given integer.

\item the family $\left(  b_{n}\right)  _{n\geq1}
= \left(  \dbinom{qn}{rn}\right)  _{n\geq1}$, where $r\in \QQ $ and
$q\in \ZZ$ are
given. (Here, a binomial coefficient $\dbinom{a}{b}$ has to be interpreted as
$0$ when $b\notin \NN$.)

\item the family $\left(  b_{n}\right)  _{n\geq1}
= \left(  \dbinom{qn-1}{rn-1}\right)  _{n\geq1}$, where $r\in \ZZ$
and $q\in \ZZ$ are given.
\end{itemize}
\end{exercise}

\begin{exercise}
\label{exe.witt.ghost-app}
For every $n\in\left\{  1,2,3,\ldots\right\}  $,
define a map $\mathbf{f}_{n}:\Lambda\rightarrow\Lambda$ by setting
\[
\mathbf{f}_{n}\left(  a\right) =
a\left(  x_{1}^{n},x_{2}^{n},x_{3}^{n},\ldots\right)
\qquad \qquad \text{for every }a\in\Lambda.
\]
(So what $\mathbf{f}_{n}$ does to a symmetric function is replacing all
variables $x_{1},x_{2},x_{3},\ldots$ by their $n$-th powers.)

\begin{enumerate}
\item[(a)] Show that $\mathbf{f}_{n}:\Lambda\rightarrow\Lambda$ is a
$\kk$-algebra homomorphism for every $n\in\left\{  1,2,3,\ldots
\right\}  $.

\item[(b)] Show that $\mathbf{f}_{n}\circ\mathbf{f}_{m}=\mathbf{f}_{nm}$ for
any two positive integers $n$ and $m$.

\item[(c)] Show that $\mathbf{f}_{1}= \id$.

\item[(d)] Prove that $\mathbf{f}_n : \Lambda \to \Lambda$ is a Hopf
algebra homomorphism for every $n \in \left\{ 1, 2, 3, \ldots \right\}$.

\item[(e)] Prove that $\mathbf{f}_2 \left( h_m \right)
= \sum_{i=0}^{2m} \left( -1\right) ^i h_i h_{2m-i}$ for every $m \in \NN$.

\item[(f)] Assume that $\kk = \ZZ$. Prove that $\mathbf{f}_{p}\left(
a\right)  \equiv a^{p}\operatorname{mod}p\Lambda$ for every $a\in\Lambda$ and
every prime number $p$.

\item[(g)] Use Exercise \ref{exe.witt.ghost-equiv} to obtain new solutions to
parts (b), (c), (d), (e) and (f) of Exercise~\ref{exe.Lambda.maps-on-pn}.
\end{enumerate}
\end{exercise}

The maps $\mathbf{f}_n$ constructed in Exercise~\ref{exe.witt.ghost-app} are
known as the \dfn{Frobenius endomorphisms} of $\Lambda$. They are a
(deceptively) simple particular case of the notion of \emph{plethysm}
(\cite[Chapter 7, Appendix 2]{Stanley} and \cite[Section I.8]{Macdonald}),
and are often used as intermediate
steps in computing more complicated plethysms\footnote{In the notations of
\cite[(A2.160)]{Stanley}, the value $\mathbf{f}_n\left(a\right)$ for an
$a \in \Lambda$ can be written as $a\left[p_n\right]$ or (when $\kk = \ZZ$)
as $p_n\left[a\right]$.}.

% [DG][v34] Added the above three exercises. They are somewhat
% tangentially related to symmetric functions, but I wanted a
% place I could reference, as the result is useful and too hard
% to find explicitly in any single point in the literature.

% [DG][v38] Added part (d) to the exercise, and comments about the
% relation to plethysm.

% [DG][v40] Added part (e) and solution to part (d).

\begin{exercise}
\label{exe.Lambda.verschiebung}
For every $n \in \left\{ 1, 2, 3, \ldots \right\}$, define a
$\kk$-algebra homomorphism
$\mathbf{v}_n : \Lambda \rightarrow \Lambda$ by
\[
\mathbf{v}_n \left( h_m \right)
=  \begin{cases}
h_{m/n},& \text{if }n\mid m;\\
0,& \text{if }n\nmid m
\end{cases}
\qquad\qquad\text{for every positive integer }m
\]
\footnote{This is well-defined, since the family
$\left( h_m \right) _{m \geq 1}$ generates the $\kk$-algebra
$\Lambda$ and is algebraically independent.}.

\begin{enumerate}
\item[(a)] Show that any positive integers $n$ and $m$ satisfy
\[
\mathbf{v}_n \left( p_m \right)
= \begin{cases}
n p_{m/n},& \text{if }n\mid m;\\
0,& \text{if }n\nmid m
\end{cases} \quad .
\]

\item[(b)] Show that any positive integers $n$ and $m$ satisfy
\[
\mathbf{v}_n \left( e_m \right)
= \begin{cases}
\left( -1 \right) ^{m - m/n} e_{m/n},& \text{if }n\mid m;\\
0,& \text{if }n\nmid m
\end{cases} \quad .
\]

\item[(c)] Prove that $\mathbf{v}_n \circ \mathbf{v}_m = \mathbf{v}_{nm}$ for
any two positive integers $n$ and $m$.

\item[(d)] Prove that $\mathbf{v}_1 = \id$.

\item[(e)] Prove that $\mathbf{v}_n : \Lambda \rightarrow \Lambda$ is a Hopf
algebra homomorphism for every $n \in \left\{ 1, 2, 3, \ldots \right\}$.
\end{enumerate}

Now, consider also the maps $\mathbf{f}_n : \Lambda \rightarrow \Lambda$ defined
in Exercise~\ref{exe.witt.ghost-app}. Fix a positive integer $n$.

\begin{enumerate}
\item[(f)] Prove that the maps $\mathbf{f}_n : \Lambda \rightarrow \Lambda$ and
$\mathbf{v}_n : \Lambda \rightarrow \Lambda$ are adjoint with respect to the
Hall inner product on $\Lambda$.

\item[(g)] Show that $\mathbf{v}_n \circ \mathbf{f}_n
= \id_{\Lambda}^{\star n}$.

\item[(h)] Prove that
$\mathbf{f}_n \circ \mathbf{v}_m = \mathbf{v}_m \circ \mathbf{f}_n$
whenever $m$ is a positive integer coprime to $n$.
\end{enumerate}

Finally, recall the $w_m \in \Lambda$ defined in Exercise~\ref{exe.witt}.

\begin{enumerate}
\item[(i)] Show that any positive integer $m$ satisfies
\[
\mathbf{v}_n \left( w_m \right)
= \begin{cases}
w_{m/n},& \text{if }n\mid m;\\
0,& \text{if }n\nmid m
\end{cases} \quad .
\]

\end{enumerate}
\end{exercise}

The homomorphisms $\mathbf{v}_n : \Lambda \rightarrow \Lambda$ defined in
Exercise~\ref{exe.Lambda.verschiebung} are called the \dfn{Verschiebung
endomorphisms} of $\Lambda$; this name comes from German, where
``Verschiebung'' means ``shift''.
This terminology, as well as that of Frobenius
endomorphisms, originates in the theory of Witt vectors, and the connection
between the Frobenius and Verschiebung endomorphisms of $\Lambda$ and the
identically named operators on Witt vectors is elucidated in \cite[Chapter
13]{Hazewinkel2}\footnote{which is also where most of the statements of
Exercises \ref{exe.witt.ghost-app} and \ref{exe.Lambda.verschiebung} come
from}.

% [DG][v39] Added above exercise.

% [DG][v75] Corrected egregious error in part (b).

\begin{exercise}
\label{exe.shareshian-wachs.path}
Fix $n \in \NN$. For any $n$-tuple $w = \left(w_1, w_2, \ldots, w_n\right)$ of
integers, define the \emph{descent set}\index{descent set of a permutation}
$\Des\left(w\right)$ of $w$ to be
the set $\left\{i \in \left\{1, 2, \ldots, n-1\right\} : w_i > w_{i+1}
\right\}$.

\begin{itemize}
\item[(a)] We say that an $n$-tuple $\left(w_1, w_2, \ldots, w_n\right)$ is
\emph{Smirnov}\index{Smirnov tuple}
if every $i \in \left\{1, 2, \ldots, n-1\right\}$
satisfies $w_i \neq w_{i+1}$.

Fix $k \in \NN$, and let $X_{n, k} \in \kk\left[\left[\xx\right]\right]$
denote the sum of the monomials $x_{w_1} x_{w_2} \cdots x_{w_n}$ over
all Smirnov $n$-tuples $w = \left(w_1, w_2, \ldots, w_n\right) \in
\left\{1, 2, 3, \ldots\right\}^n$ satisfying
$\left|\Des\left(w\right)\right| = k$. Prove that
$X_{n, k} \in \Lambda$.

\item[(b)] For any $n$-tuple $w = \left(w_1, w_2, \ldots, w_n\right)$,
define the \dfn{stagnation set}
$\operatorname{Stag}\left(w\right)$ of $w$ to be the set
\newline % I hate having to do that :/
$\left\{i \in \left\{1, 2, \ldots, n-1\right\} : w_i = w_{i+1}
\right\}$. (Thus, an $n$-tuple is Smirnov if and only if its stagnation
set is empty.)

For any $d \in \NN$ and $s \in \NN$, define a power
series $X_{n, d, s} \in \kk\left[\left[\xx\right]\right]$
as the sum of the monomials $x_{w_1} x_{w_2} \cdots x_{w_n}$ over
all $n$-tuples $w = \left(w_1, w_2, \ldots, w_n\right)
\in \left\{1, 2, 3, \ldots\right\}^n$ satisfying
$\left|\Des\left(w\right)\right| = d$ and
$\left|\operatorname{Stag}\left(w\right)\right| = s$. Prove that
$X_{n, d, s} \in \Lambda$ for any nonnegative integers $d$ and $s$.

\item[(c)] Assume that $n$ is positive.
For any $d \in \NN$ and $s \in \NN$, define three further
power series $U_{n, d, s}$, $V_{n, d, s}$ and $W_{n, d, s}$ in
$\kk\left[\left[\xx\right]\right]$ by the following formulas:
\begin{align}
\label{eq.exe.shareshian-wachs.path.c.U}
U_{n, d, s}
&= \sum\limits_{\substack{w = \left(w_1, w_2, \ldots, w_n\right)
    \in \left\{1, 2, 3, \ldots\right\}^n ; \\
  \left|\Des\left(w\right)\right| = d ; \ 
  \left|\operatorname{Stag}\left(w\right)\right| = s ; \\
  w_1 < w_n}} x_{w_1} x_{w_2} \cdots x_{w_n} ; \\
\label{eq.exe.shareshian-wachs.path.c.V}
V_{n, d, s}
&= \sum\limits_{\substack{w = \left(w_1, w_2, \ldots, w_n\right)
    \in \left\{1, 2, 3, \ldots\right\}^n ; \\
  \left|\Des\left(w\right)\right| = d ; \ 
  \left|\operatorname{Stag}\left(w\right)\right| = s ; \\
  w_1 = w_n}} x_{w_1} x_{w_2} \cdots x_{w_n} ; \\
\label{eq.exe.shareshian-wachs.path.c.W}
W_{n, d, s}
&= \sum\limits_{\substack{w = \left(w_1, w_2, \ldots, w_n\right)
    \in \left\{1, 2, 3, \ldots\right\}^n ; \\
  \left|\Des\left(w\right)\right| = d ; \ 
  \left|\operatorname{Stag}\left(w\right)\right| = s ; \\
  w_1 > w_n}} x_{w_1} x_{w_2} \cdots x_{w_n} .
\end{align}
Prove that these three power series
$U_{n, d, s}$, $V_{n, d, s}$ and $W_{n, d, s}$ belong to
$\Lambda$.
\end{itemize}
\end{exercise}

\begin{remark}
The function $X_{n, k}$ in Exercise~\ref{exe.shareshian-wachs.path}(a)
is a simple example
(\cite[Example 2.5, Theorem C.3]{ShareshianWachs-2014})
of a chromatic quasisymmetric function that happens to be
symmetric. See Shareshian/Wachs \cite{ShareshianWachs-2014} for more
general criteria for such functions to be symmetric, as well as deeper
results. For example, \cite[Theorem 6.3]{ShareshianWachs-2014} gives
an expansion for a wide class of chromatic quasisymmetric functions in
the Schur basis of $\Lambda$, which, in particular, shows that our
$X_{n, k}$ satisfies
\[
X_{n, k} = \sum_{\lambda \in \Par_n} a_{\lambda, k} s_\lambda ,
\]
where $a_{\lambda, k}$ is the number of all assignments $T$ of entries
in $\left\{1, 2, \ldots, n\right\}$ to the cells of the Ferrers
diagram of $\lambda$ such that the following four conditions are
satisfied:
\begin{itemize}
\item Every element of $\left\{1, 2, \ldots, n\right\}$ is used
precisely once in the assignment (i.e., we have $\cont\left(T\right)
= \left(1^n\right)$).
\item Whenever a cell $y$ of the Ferrers diagram lies immediately to the right of
a cell $x$, we have $T\left(y\right) - T\left(x\right) \geq 2$.
\item Whenever a cell $y$ of the Ferrers diagram lies immediately below
a cell $x$, we have $T\left(y\right) - T\left(x\right) \geq -1$.
\item There exist precisely $k$ elements
$i \in \left\{1, 2, \ldots, n-1\right\}$ such that the cell
$T^{-1}\left(i\right)$ lies in a row below $T^{-1}\left(i+1\right)$.
\end{itemize}
Are there any such rules for the $X_{n, d, s}$ of part (b)?

Smirnov $n$-tuples are more usually called Smirnov words, or
(occasionally) \idx{Carlitz words}.

See \cite[Chapter 6]{Ellzey-thesis} for further properties of
the symmetric functions $U_{n, d, 0}$, $V_{n, d, 0}$ and
$W_{n, d, 0}$ from Exercise~\ref{exe.shareshian-wachs.path}(c)
(or, more precisely, of their generating functions
$\sum_d U_{n, d, 0} t^d$ etc.).
\end{remark}

% [DG][v25] Added above exercise. The second solution to part
% (a) (which I would not have written up if I had known of the other
% solution by then) is a fine example of my inability to write
% combinatorial proofs.
% 
% Part (b) as well as the name "stagnation set" is mine (or I have just
% not skimmed enough of the Shareshian-Wachs paper); do you know
% the "right" notion?

% [DG][v77] Added the last paragraph of the above remark.

\begin{exercise} \phantomsection
\label{exe.pe-ep-determinants}

\begin{itemize}

\item[(a)] Let $n \in \NN$. Define a matrix
$A_n = \left(a_{i,j}\right)_{i,j = 1, 2, \ldots, n} \in \Lambda^{n \times n}$
by
\[
a_{i,j} =
  \begin{cases}
   p_{i-j+1}, & \text{if } i \geq j ; \\
   i,         & \text{if } i = j - 1 ; \\
   0,         & \text{if } i < j - 1
  \end{cases}
\qquad \qquad \text{ for all } \left( i, j \right) \in \left\{ 1, 2, \ldots , n \right\}^2 .
\]
This matrix $A_n$ looks as follows:
\[
A_n
=
\left(
\begin{array}[c]{cccccc}
p_1   & 1 & 0 & \cdots & 0 & 0\\
p_2   & p_1   & 2 & \cdots & 0 & 0\\
p_3   & p_2   & p_1 & \cdots & 0 & 0\\
\vdots & \vdots & \vdots & \ddots & \vdots & \vdots\\
p_{n-1} & p_{n-2} & p_{n-3} & \cdots & p_1   & n-1\\
p_n & p_{n-1} & p_{n-2} & \cdots & p_2   & p_1 
\end{array}
\right) .
\]
Show that $\det\left(A_n\right) = n! e_n$.

\item[(b)] Let $n$ be a positive integer. Define a matrix
$B_n = \left(b_{i,j}\right)_{i,j = 1, 2, \ldots, n} \in \Lambda^{n \times n}$
by
\[
b_{i,j} =
  \begin{cases}
   i e_i,     & \text{if } j = 1 ; \\
   e_{i-j+1}, & \text{if } j > 1
% The following is redundant because we have $e_k = 0$ for negative $k$, and $e_0 = 1$:
%   e_{i-j+1}, & \text{if } i \geq j > 1 ; \\
%   1,         & \text{if } i = j - 1 ; \\
%   0,         & \text{if } i < j - 1
  \end{cases}
\qquad \qquad \text{ for all } \left( i, j \right) \in \left\{ 1, 2, \ldots , n \right\}^2 .
\]
The matrix $B_n$ looks as follows:
\begin{align*}
B_n
&=
\left(
\begin{array}[c]{cccccc}
e_1  & e_0 & e_{-1} & \cdots & e_{-n+3} & e_{-n+2} \\
2e_2  & e_1  & e_0 & \cdots & e_{-n+4} & e_{-n+3}\\
3e_3  & e_2  & e_1 & \cdots & e_{-n+5} & e_{-n+4}\\
\vdots & \vdots & \vdots & \ddots & \vdots & \vdots\\
\left(  n-1\right)  e_{n-1}  & e_{n-2} & e_{n-3}  & \cdots & e_1 & e_0\\
n e_n & e_{n-1}  & e_{n-2} & \cdots & e_2  & e_1
\end{array}
\right) \\
&=
\left(
\begin{array}[c]{cccccc}
e_1  & 1 & 0 & \cdots & 0 & 0\\
2e_2  & e_1  & 1 & \cdots & 0 & 0\\
3e_3  & e_2  & e_1 & \cdots & 0 & 0\\
\vdots & \vdots & \vdots & \ddots & \vdots & \vdots\\
\left(  n-1\right)  e_{n-1}  & e_{n-2} & e_{n-3}  & \cdots & e_1 & 1\\
n e_n & e_{n-1}  & e_{n-2} & \cdots & e_2  & e_1
\end{array}
\right) .
\end{align*}
Show that $\det\left(B_n\right) = p_n$.

\end{itemize}

\end{exercise}

The formulas of Exercise~\ref{exe.pe-ep-determinants}, for
finitely many variables, appear in
Prasolov's \cite[\S 4.1]{Prasolov-linalg}\footnote{where our
symmetric functions $e_k, h_k, p_k$, evaluated in finitely many
indeterminates, are denoted $\sigma_k, p_k, s_k$, respectively}.
In \cite[\S 4.2]{Prasolov-linalg}, Prasolov gives four more
formulas, which express $e_n$ as a polynomial in the
$h_1,h_2,h_3,\ldots$, or $h_n$ as a polynomial in the
$e_1,e_2,e_3,\ldots$, or $p_n$ as a polynomial in the
$h_1,h_2,h_3,\ldots$, or $n!h_n$ as a polynomial in the
$p_1,p_2,p_3,\ldots$. These are not novel for us, since the first
two of them are particular cases of
Theorem~\ref{Jacobi-Trudi-formulae}, whereas the latter two
can be derived from Exercise~\ref{exe.pe-ep-determinants} by
applying $\omega$. (Note that $\omega$ is only well-defined
on symmetric functions in infinitely many indeterminates, so
we need to apply $\omega$ \textbf{before} evaluating at
finitely many indeterminates; this explains why Prasolov has
to prove the latter two identities separately.)

% [DG][v30] Added this exercise. I tried not to copy the solution of
% Exercise 9.3 in https://github.com/darijgr/lambda too
% much, but there is only so much I can do when the solution is
% essentially the same.

% [DG][v32] Added paragraph with Prasolov reference.

\begin{exercise}
\label{exe.hookschur}
In the following, if $k \in \NN$, we shall use the notation
$1^k$ for $\underbrace{1,1,\ldots,1}_{k\text{ times}}$ (in
contexts such as $\left( n, 1^m \right)$). So, for example,
$\left( 3, 1^4 \right)$ is the partition
$\left( 3, 1, 1, 1, 1 \right)$.

\begin{itemize}

\item[(a)] Show that
$e_n h_m = s_{\left( m+1, 1^{n-1} \right)} + s_{\left( m, 1^n \right)}$
for any two positive integers $n$ and $m$.

\item[(b)] Show that
\[
\sum_{i=0}^{b} \left(-1\right)^i h_{a+i+1} e_{b-i}
= s_{\left( a+1, 1^b \right)}
\]
for any $a \in \NN$ and $b \in \NN$.

\item[(c)] Show that
\[
\sum_{i=0}^{b} \left(-1\right)^i h_{a+i+1} e_{b-i}
= \left(-1\right)^b \delta_{a+b,-1}
\]
for any negative integer $a$ and every $b \in \NN$.
(As usual, we set $h_j = 0$ for $j < 0$ here.)

\item[(d)] Show that
\begin{align*}
\Delta s_{\left( a+1, 1^b \right)}
&= 1 \otimes s_{\left( a+1, 1^b \right)}
 + s_{\left( a+1, 1^b \right)} \otimes 1 \\
& \ \ \ \ \ \ \ \ \ \
 + \sum\limits_{\substack{\left( c, d, e, f \right) \in \NN^4; \\
         c+e=a-1; \\ d+f=b}}
   s_{\left( c+1, 1^d \right)} \otimes s_{\left( e+1, 1^f \right)}
 + \sum\limits_{\substack{\left( c, d, e, f \right) \in \NN^4; \\
         c+e=a; \\ d+f=b-1}}
   s_{\left( c+1, 1^d \right)} \otimes s_{\left( e+1, 1^f \right)}
\end{align*}
for any $a \in \NN$ and $b \in \NN$.

\end{itemize}
\end{exercise}

% [DG][v30] Added above exercise.

Our next few exercises survey some results on
Littlewood-Richardson coefficients.

\begin{exercise}
\label{exe.liri.180}
Let $m \in \NN$ and $k \in \NN$. Let
$\lambda$ and $\mu$ be two partitions such that $\ell\left(  \lambda\right)
\leq k$ and $\ell\left(  \mu\right)  \leq k$. Assume that all parts of
$\lambda$ and all parts of $\mu$ are $\leq m$. (It is easy to see
that this assumption is equivalent to requiring $\lambda_i \leq m$
and $\mu_i \leq m$ for every positive integer $i$.\ \ \ \ \footnote{As usual,
we are denoting by $\nu_i$ the $i$-th entry of a partition $\nu$ here.}).
Let
$\lambda^{\vee}$ and $\mu^{\vee}$ denote the $k$-tuples $\left(  m-\lambda
_{k},m-\lambda_{k-1},\ldots,m-\lambda_{1}\right)  $ and $\left(  m-\mu
_{k},m-\mu_{k-1},\ldots,m-\mu_{1}\right)  $, respectively.

\begin{enumerate}
\item[(a)] Show that $\lambda^{\vee}$ and $\mu^{\vee}$ are partitions, and
that $s_{\lambda/\mu} = s_{\mu^{\vee}/\lambda^{\vee}}$.

\item[(b)] Show that
$c_{\mu,\nu}^{\lambda} = c_{\lambda^{\vee},\nu}^{\mu^{\vee}}$
for any partition $\nu$.

\item[(c)] Let $\nu$ be a partition such that $\ell\left(  \nu\right)  \leq
k$, and such that all parts of $\nu$ are $\leq m$. Let $\nu^{\vee}$ denote the
$k$-tuple $\left(  m-\nu_{k},m-\nu_{k-1},\ldots,m-\nu_{1}\right)  $. Show that
$\nu^{\vee}$ is a partition, and satisfies
\[
c_{\mu,\nu}^{\lambda} = c_{\nu,\mu}^{\lambda}
= c_{\lambda^{\vee},\nu}^{\mu^{\vee}} = c_{\nu,\lambda^{\vee}}^{\mu^{\vee}}
= c_{\mu,\lambda^{\vee}}^{\nu^{\vee}} = c_{\lambda^{\vee},\mu}^{\nu^{\vee}} .
\]

\item[(d)] Show that
\[
s_{\lambda^{\vee}}\left(  x_{1},x_{2},\ldots,x_{k}\right)
=\left(  x_{1} x_{2}\cdots x_{k}\right)  ^{m}
 \cdot s_{\lambda}\left(  x_{1}^{-1},x_{2}^{-1},\ldots,x_{k}^{-1}\right)
\]
in the Laurent polynomial ring $\kk \left[  x_{1},x_{2},\ldots
,x_{k},x_{1}^{-1},x_{2}^{-1},\ldots,x_{k}^{-1}\right]  $.

\item[(e)] Let $r$ be a nonnegative integer. Show that
$\left( r+\lambda_{1}, r+\lambda_{2}, \ldots, r+\lambda_{k} \right)$ is
a partition and satisfies
\[
s_{\left( r+\lambda_{1}, r+\lambda_{2}, \ldots, r+\lambda_{k} \right)}
\left( x_{1},x_{2},\ldots,x_{k} \right)
= \left(  x_{1}x_{2}\cdots x_{k}\right)^{r}
\cdot s_{\lambda}\left(  x_{1},x_{2},\ldots,x_{k}\right)
\]
in the polynomial ring $\kk \left[  x_{1},x_{2},\ldots,x_{k}\right]  $.
\end{enumerate}
\end{exercise}

\begin{exercise}
\label{exe.liri.180ab}
Let $m \in \NN$, $n \in \NN$ and $k \in \NN$.
Let $\mu$ and $\nu$ be two partitions such that $\ell\left(
\mu\right)  \leq k$ and $\ell\left(  \nu\right)  \leq k$. Assume that all
parts of $\mu$ are $\leq m$ (that is, $\mu_i \leq m$ for every positive
integer $i$)\ \ \ \ \footnote{As usual, we are denoting by $\nu_i$ the
$i$-th entry of a partition $\nu$ here.},
and that all parts of $\nu$ are $\leq n$ (that is, $\nu_i \leq n$
for every positive integer $i$). Let $\mu
^{\vee\left\{  m\right\}  }$ denote the $k$-tuple $\left(  m-\mu_{k}
,m-\mu_{k-1},\ldots,m-\mu_{1}\right)  $, and let $\nu^{\vee\left\{  n\right\}
}$ denote the $k$-tuple $\left(  n-\nu_{k},n-\nu_{k-1},\ldots,n-\nu
_{1}\right)  $.

\begin{enumerate}
\item[(a)] Show that $\mu^{\vee\left\{  m\right\}  }$ and
$\nu^{\vee\left\{ n\right\}  }$ are partitions.
\end{enumerate}

Now, let $\lambda$ be a further partition such that $\ell\left(
\lambda\right)  \leq k$.

\begin{enumerate}
\item[(b)] If not all parts of $\lambda$ are $\leq m+n$, then show that
$c_{\mu,\nu}^{\lambda} = 0$.

\item[(c)] If all parts of $\lambda$ are $\leq m+n$, then show that
$c_{\mu,\nu}^{\lambda} = c_{\mu^{\vee\left\{  m\right\}  },\nu^{\vee\left\{
n\right\}  }}^{\lambda^{\vee\left\{  m+n\right\}  }}$,
where $\lambda ^{\vee\left\{  m+n\right\}  }$ denotes the $k$-tuple
$\left(  m+n-\lambda_{k},m+n-\lambda_{k-1},\ldots,m+n-\lambda_{1}\right)  $.
\end{enumerate}
\end{exercise}

The results of Exercise~\ref{exe.pieri.omega}(c) and
Exercise~\ref{exe.liri.180}(c) are two
\dfn{symmetries of Littlewood-Richardson
coefficients}\footnote{The result of Exercise~\ref{exe.liri.180ab}(c) can also
be regarded as a symmetry of Littlewood-Richardson coefficients; see
\cite[\S 3.3]{Azenhas-LRsymm24}.}; combining them yields further such
symmetries. While these symmetries were relatively easy consequences of our
algebraic definition of the Littlewood-Richardson coefficients, it is a much
more challenging task to derive them bijectively from a combinatorial
definition of these coefficients (such as the one given in
Corollary~\ref{L-R-rule}). Some such derivations appear in \cite{ThomasYong},
in \cite{AzenhasKingTerada}, in
\cite[Example 3.6, Proposition 5.11 and references therein]{BenkartSottileStroomer},
\cite[\S 5.1, \S A.1, \S A.4]{Fulton} and
\cite[(2.12)]{KirillovBerenstein} (though a different combinatorial
interpretation of $c_{\mu,\nu}^{\lambda}$ is used in the latter three).

% [DG][v44] Added the preceding two exercises.

\begin{exercise}
\label{exe.liri.simple}
Recall our usual notations: For every partition $\lambda$
and every positive integer $i$, the $i$-th entry of
$\lambda$ is denoted by $\lambda_i$. The sign $\triangleright$
stands for dominance order. We let $\lambda^t$ denote
the conjugate partition of a partition $\lambda$.

For any two partitions $\mu$ and $\nu$, we define two new
partitions $\mu + \nu$ and $\mu \sqcup \nu$ of
$\left|\mu\right| + \left|\nu\right|$ as follows:
\begin{itemize}
\item The partition $\mu + \nu$ is defined as
$\left(\mu_1 + \nu_1, \mu_2 + \nu_2, \mu_3 + \nu_3, \ldots\right)$.
\item The partition $\mu \sqcup \nu$ is defined as the
result of sorting the list
$\left(\mu_1, \mu_2, \ldots, \mu_{\ell\left(\mu\right)},
\nu_1, \nu_2, \ldots, \nu_{\ell\left(\nu\right)}\right)$
in decreasing order.
\end{itemize}

\begin{itemize}

\item[(a)] Show that any two partitions $\mu$ and $\nu$
satisfy
$\left(\mu + \nu\right)^t = \mu^t \sqcup \nu^t$
and
$\left(\mu \sqcup \nu\right)^t = \mu^t + \nu^t$.

\item[(b)] Show that any two partitions $\mu$ and $\nu$
satisfy $c^{\mu + \nu}_{\mu, \nu} = 1$ and
$c^{\mu \sqcup \nu}_{\mu, \nu} = 1$.

\item[(c)] If $k \in \NN$ and $n \in \NN$ satisfy
$k \leq n$, and if $\mu \in \Par_k$,
$\nu \in \Par_{n-k}$ and $\lambda \in \Par_n$ are such that
$c^\lambda_{\mu, \nu} \neq 0$, then prove that
$\mu + \nu \triangleright \lambda \triangleright \mu \sqcup \nu$.

\item[(d)] If $n \in \NN$ and $m \in \NN$ and $\alpha, \beta
\in \Par_n$ and $\gamma, \delta \in \Par_m$ are such that
$\alpha \triangleright \beta$ and
$\gamma \triangleright \delta$, then show that
$\alpha + \gamma \triangleright \beta + \delta$ and
$\alpha \sqcup \gamma \triangleright \beta \sqcup \delta$.

\item[(e)] Let $m \in \NN$ and $k \in \NN$, and let $\lambda$
be the partition
$\left(m^k\right) = \left(\underbrace{m, m, \ldots, m}_{k
\text{ times}}\right)$.
Show that any two partitions $\mu$ and $\nu$ satisfy
$c^\lambda_{\mu, \nu} \in \left\{0, 1\right\}$.

\item[(f)] Let $a \in \NN$ and $b \in \NN$, and let $\lambda$
be the partition $\left(a+1, 1^b\right)$ (using the notation
of Exercise~\ref{exe.hookschur}).
Show that any two partitions $\mu$ and $\nu$ satisfy
$c^\lambda_{\mu, \nu} \in \left\{0, 1\right\}$.

\item[(g)] If $\lambda$ is any partition, and if $\mu$ and
$\nu$ are two rectangular partitions\footnote{A partition
is called \emph{rectangular}\index{rectangular partition}
if it has the form
$\left(m^k\right) = \left(\underbrace{m, m, \ldots, m}_{k
\text{ times}}\right)$ for some $m \in \NN$ and $k \in \NN$.},
then show that
$c^\lambda_{\mu, \nu} \in \left\{0, 1\right\}$.

\end{itemize}
\end{exercise}

Exercise~\ref{exe.liri.simple}(g) is part of Stembridge's
\cite[Thm. 2.1]{Stembridge-multfree}; we refer to that
article for further results of its kind.

% [DG][v29] Added this exercise.

% [DG][v30] Added part (f).

% [DG][v32] Added part (g).

The Littlewood-Richardson rule comes in many different
forms, whose equivalence is not always immediate. Our version
(Corollary~\ref{L-R-rule}) has the advantage of being the
simplest to prove and one of the simplest to state. Other
versions can be found in \cite[appendix 1 to Ch. 7]{Stanley},
Fulton's \cite[Ch. 5]{Fulton}
and van Leeuwen's \cite{Leeuwen-lrr}. We restrict ourselves
to proving some very basic equivalences that allow us to
restate parts of Corollary~\ref{L-R-rule}:

\begin{exercise}
\label{exe.liri.yamanouchi}
We shall use the following notations:

\begin{itemize}

\item If $T$ is a column-strict tableau and $j$ is a positive
integer, then we use the notation $T|_{\cols \geq j}$ for
the restriction of $T$ to the union of its columns
$j,j+1,j+2,\ldots$.
%, and we use the notation
%$T|_{\operatorname{rows} \geq j}$
%for the restriction of $T$ to the union of its rows
%$j,j+1,j+2,\ldots$.
(This notation has
already been used in Section~\ref{bialternant-section}.)

\item If $T$ is a column-strict tableau and $S$ is a set of
cells of $T$, then we write $T|_S$ for the
restriction of $T$ to the set $S$ of cells.\footnote{This
restriction $T|_S$
is not necessarily a tableau of skew shape; it is just a map
from $S$ to $\left\{1, 2, 3, \ldots\right\}$. The content
$\cont\left(T|_S\right)$ is nevertheless well-defined (in the
usual way: $\left(\cont\left(T|_S\right)\right)_i =
\left| \left(T|_S\right)^{-1} \left(i\right) \right|$).}

\item If $T$ is a column-strict tableau, then an
\emph{NE-set}\index{NE-set of a tableau}
of $T$ means a set $S$ of cells of $T$ such
that whenever $s \in S$, every cell of $T$ which lies
northeast\footnote{A cell $\left(r,c\right)$ is said to lie
\dfn{northeast} of a cell $\left(r',c'\right)$ if and
only if we have $r \leq r'$ and $c \geq c'$.} of $s$ must
also belong to $S$.

\item The \dfn{Semitic reading word}\footnote{The
notation comes from \cite{Leeuwen-lrr} and is a reference
to the Arabic and Hebrew way of writing.} of a
column-strict tableau $T$ is the concatenation\footnote{If
$s_1, s_2, s_3, \ldots$ are several words (finitely or
infinitely many), then the \dfn{concatenation}
$s_1 s_2 s_3 \cdots$ is defined as the word which is
obtained by starting with the empty word, then appending
$s_1$ to its end, then appending $s_2$ to the end of the
result, then appending $s_3$ to the end of the result, etc.}
$r_1 r_2 r_3 \cdots$, where $r_i$ is the word obtained
by reading the $i$-th row of $T$ from right to
left.\footnote{For example, the Semitic reading word of
the tableau
\[
\begin{matrix}
&   & 3 & 4 & 4 & 5 \\
& 1 & 4 & 6 &   &   \\
& 3 & 5 &   &   &
\end{matrix}
\]
is $5443 641 53$.
\par
The Semitic reading word of a tableau $T$ is what
is called the \dfn{reverse reading word} of $T$ in
\cite[\S A.1.3]{Stanley}.}

\item If $w = \left(w_1, w_2, \ldots, w_n\right)$ is a
word, then a \dfn{prefix} of $w$ means a word of the
form $\left(w_1, w_2, \ldots, w_i\right)$ for some
$i \in \left\{0, 1, \ldots, n\right\}$. (In particular,
both $w$ and the empty word are prefixes of $w$.)

A word $w$ over the set of positive integers is said to
be \emph{Yamanouchi}\index{Yamanouchi word}
if for any prefix $v$ of $w$ and
any positive integer $i$, there are at least as many
$i$'s among the letters of $v$ as there are
$\left(i+1\right)$'s among them.\footnote{For instance,
the words $11213223132$ and $1213$ are Yamanouchi, while
the words $132$, $21$ and $1121322332111$ are not. The
Dyck words (defined as in \cite[Example 6.6.6]{Stanley},
and written using $1$'s and $2$'s instead of $x$'s and
$y$'s) are precisely
the Yamanouchi words whose letters are $1$'s and $2$'s
and in which the letter $1$ appears as often as the
letter $2$.

Yamanouchi words are often called lattice permutations.}

\end{itemize}

Prove the following two statements:

\begin{itemize}

\item[(a)] Let $\mu$ be a partition.
Let $b_{i, j}$ be a nonnegative integer for every two
positive integers $i$ and $j$. Assume that $b_{i, j} = 0$
for all but finitely many pairs $\left(i, j\right)$.

The following two assertions are equivalent:

\begin{itemize}

\item
\textit{Assertion $\mathcal{A}$:}
There exist a partition $\lambda$
and a column-strict tableau $T$ of shape $\lambda / \mu$
such that all
$\left(i, j\right) \in
\left\{1, 2, 3, \ldots \right\}^2$ satisfy
\begin{equation}
\label{eq.exe.liri.yamanouchi.A.bij}
b_{i, j} = \left(\text{the number of all entries } i
\text{ in the } j \text{-th row of } T\right) .
\end{equation}

\item
\textit{Assertion $\mathcal{B}$:}
The inequality
\begin{equation}
\label{eq.exe.liri.yamanouchi.B.bij}
\mu_{j+1} + \left(b_{1, j+1} + b_{2, j+1} + \cdots + b_{i+1, j+1}\right)
\leq
\mu_j + \left(b_{1, j} + b_{2, j} + \cdots + b_{i, j}\right)
\end{equation}
holds for all
$\left(i, j\right) \in \NN \times
\left\{1, 2, 3, \ldots \right\}$.

\end{itemize}

\item[(b)] Let $\lambda$ and $\mu$ be two partitions,
and let $T$ be a column-strict tableau of shape
$\lambda / \mu$. Then, the following five assertions are
equivalent:

\begin{itemize}

\item
\textit{Assertion $\mathcal{C}$:} For every positive integer
$j$, the weak composition
$\cont\left(T|_{\cols \geq j}\right)$ is a partition.

\item
\textit{Assertion $\mathcal{D}$:} For every positive integers
$j$ and $i$, the number of entries $i+1$ in the first $j$
rows\footnote{The ``first $j$ rows'' mean the $1$-st row,
the $2$-nd row, etc., the $j$-th row (even if some of these
rows are empty).} of $T$ is $\leq$ to the number of
entries $i$ in the first $j-1$ rows of $T$.

\item
\textit{Assertion $\mathcal{E}$:} For every NE-set $S$
of $T$, the weak composition $\cont\left(T|_S\right)$ is a
partition.

\item
\textit{Assertion $\mathcal{F}$:} The Semitic reading word
of $T$ is Yamanouchi.

\item
\textit{Assertion $\mathcal{G}$:} There exists a column-strict
tableau $S$ whose shape is a partition and which satisfies
the following property: For any positive integers $i$ and $j$,
the number of entries $i$ in the $j$-th row of $T$ equals the
number of entries $j$ in the $i$-th row of $S$.

\end{itemize}
\end{itemize}

\end{exercise}

\begin{remark}
The equivalence of Assertions $\mathcal{C}$ and $\mathcal{F}$
in Exercise~\ref{exe.liri.yamanouchi}(b) is the
``not-too-difficult exercise'' mentioned in
\cite{Stembridge}. It yields the equivalence between our version
of the Littlewood-Richardson rule (Corollary~\ref{L-R-rule}) and
that in \cite[A1.3.3]{Stanley}.
\end{remark}

% [DG][v29] Added this exercise. It (or rather, the equivalences
% between Assertions A and B and between Assertions C and D) is
% needed for another exercise I plan to add; but it also provides
% bridges to some literature. It is still far from proving any
% nontrivial equivalence results between different LR rules
% (partly because we don't define jeu de taquin, hives or
% Gelfand-Tsetlin triangles in these notes and I don't feel like
% opening these cans of worms right now).
% 
% The solution bears silent witness to my inability to write
% combinatorial proofs; it gained in pedantry and length every
% time I discovered a mistake, and that was quite often...

In the next exercises, we shall restate Corollary~\ref{infinite-L-R-rule} in a
different form. While Corollary~\ref{infinite-L-R-rule} provided a
decomposition of the product of a skew Schur function with a Schur function
into a sum of Schur functions, the different form that we will encounter in
Exercise~\ref{exe.liri.skew.1}(b) will give a combinatorial interpretation for
the Hall inner product between two skew Schur functions. Let us first
generalize Exercise~\ref{exe.liri.yamanouchi}(b):

\begin{exercise}
\label{exe.liri.yamanouchi.kappa}
Let us use the notations of Exercise~\ref{exe.liri.yamanouchi}. Let
$\kappa$, $\lambda$ and $\mu$ be three partitions, and let $T$ be a
column-strict tableau of shape $\lambda/\mu$.

\begin{itemize}
\item[(a)] Prove that the following five assertions are equivalent:

\begin{itemize}
\item \textit{Assertion $\mathcal{C}^{\left(  \kappa\right)  }$:}
For every positive integer $j$, the weak composition
$\kappa + \cont \left( T|_{\operatorname{cols}\geq j}\right)$
is a partition.

\item \textit{Assertion $\mathcal{D}^{\left(  \kappa\right)  }$:}
For every positive integers $j$ and $i$, we have
\begin{align*}
&  \kappa_{i+1}+\left(  \text{the number of entries }i+1\text{ in the first
}j\text{ rows of }T\right) \\
&  \leq\kappa_{i}+\left(  \text{the number of entries }i\text{ in the first
}j-1\text{ rows of }T\right)  .
\end{align*}

\item \textit{Assertion $\mathcal{E}^{\left(  \kappa\right)  }$:}
For every NE-set $S$ of $T$, the weak composition $\kappa+ \cont
\left(  T|_{S}\right)  $ is a partition.

\item \textit{Assertion $\mathcal{F}^{\left(  \kappa\right)  }$:}
For every prefix $v$ of the Semitic reading word of $T$, and for every
positive integer $i$, we have
\begin{align*}
&  \kappa_{i}+\left(  \text{the number of }i\text{'s among the letters of
}v\right) \\
&  \geq\kappa_{i+1}+\left(  \text{the number of }\left(  i+1\right)  \text{'s
among the letters of }v\right)  .
\end{align*}

\item \textit{Assertion $\mathcal{G}^{\left(  \kappa\right)  }$:}
There exist a partition $\zeta$ and a column-strict tableau $S$ of shape
$\zeta/\kappa$ which satisfies the following property: For any positive
integers $i$ and $j$, the number of entries $i$ in the $j$-th row of $T$
equals the number of entries $j$ in the $i$-th row of $S$.
\end{itemize}

\item[(b)] Let $\tau$ be a partition such that $\tau=\kappa
+ \cont T$. Consider the five assertions $\mathcal{C}^{\left(
\kappa\right)  }$, $\mathcal{D}^{\left(  \kappa\right)  }$,
$\mathcal{E}^{\left(  \kappa\right)  }$,
$\mathcal{F}^{\left(  \kappa\right)  }$ and
$\mathcal{G}^{\left(  \kappa\right)  }$ introduced in
Exercise~\ref{exe.liri.yamanouchi.kappa}(a). Let us also consider the
following assertion:

\begin{itemize}
\item \textit{Assertion $\mathcal{H}^{\left(  \kappa\right)  }$:}
There exists a column-strict tableau $S$ of shape $\tau/\kappa$ which
satisfies the following property: For any positive integers $i$ and $j$, the
number of entries $i$ in the $j$-th row of $T$ equals the number of entries
$j$ in the $i$-th row of $S$.
\end{itemize}

Prove that the six assertions $\mathcal{C}^{\left(  \kappa\right)  }$,
$\mathcal{D}^{\left(  \kappa\right)  }$,
$\mathcal{E}^{\left(  \kappa\right)  }$,
$\mathcal{F}^{\left(  \kappa\right)  }$, $\mathcal{G}^{\left(
\kappa\right)  }$ and $\mathcal{H}^{\left(  \kappa\right)  }$ are equivalent.
\end{itemize}
\end{exercise}

Clearly, Exercise~\ref{exe.liri.yamanouchi}(b) is the particular case of
Exercise~\ref{exe.liri.yamanouchi.kappa} when $\kappa=\varnothing$.

Using Exercise~\ref{exe.liri.yamanouchi.kappa}, we can restate
Corollary~\ref{infinite-L-R-rule} in several ways:

\begin{exercise}
\label{exe.liri.skew.1}
Let $\lambda$, $\mu$ and $\kappa$ be three partitions.

\begin{enumerate}
\item[(a)] Show that
\[
s_{\kappa}s_{\lambda/\mu}=\sum_{T}s_{\kappa+ \cont T},
\]
where the sum ranges over all column-strict tableaux $T$ of shape $\lambda
/\mu$ satisfying the five equivalent assertions $\mathcal{C}^{\left(
\kappa\right)  }$, $\mathcal{D}^{\left(  \kappa\right)  }$,
$\mathcal{E}^{\left(  \kappa\right)  }$,
$\mathcal{F}^{\left(  \kappa\right)  }$ and
$\mathcal{G}^{\left(  \kappa\right)  }$ introduced in
Exercise~\ref{exe.liri.yamanouchi.kappa}(a).

\item[(b)] Let $\tau$ be a partition. Show that $\left(  s_{\lambda/\mu
},s_{\tau/\kappa}\right)  _{\Lambda}$ is the number of all column-strict
tableaux $T$ of shape $\lambda/\mu$ satisfying $\tau=\kappa
+ \cont T$ and also satisfying the six equivalent assertions
$\mathcal{C}^{\left(  \kappa\right)  }$, $\mathcal{D}^{\left(  \kappa\right)
}$, $\mathcal{E}^{\left(  \kappa\right)  }$, $\mathcal{F}^{\left(
\kappa\right)  }$, $\mathcal{G}^{\left(  \kappa\right)  }$ and
$\mathcal{H}^{\left(  \kappa\right)  }$ introduced in
Exercise~\ref{exe.liri.yamanouchi.kappa}.
\end{enumerate}
\end{exercise}

Exercise~\ref{exe.liri.skew.1}(a) is merely Corollary~\ref{infinite-L-R-rule},
rewritten in light of Exercise~\ref{exe.liri.yamanouchi.kappa}. Various parts
of it appear in the literature. For instance, \cite[(53)]{Leeuwen-altSchur}
easily reveals to be a restatement of the fact that $s_{\kappa}s_{\lambda/\mu
}=\sum_{T}s_{\nu+ \cont T}$, where the sum ranges over all
column-strict tableaux $T$ of shape $\lambda/\mu$ satisfying Assertion
$\mathcal{D}^{\left(  \kappa\right)  }$.

Exercise~\ref{exe.liri.skew.1}(b) is one version of a ``skew
Littlewood-Richardson rule'' that goes back to Zelevinsky
\cite{Zelevinsky-skewLR} (although Zelevinsky's version uses both a different
language and a combinatorial interpretation which is not obviously equivalent
to ours). It appears in various sources; for instance, \cite[Theorem 5.2,
second formula]{Leeuwen-altSchur} says that $\left(  s_{\lambda/\mu}%
,s_{\tau/\kappa}\right)  _{\Lambda}$ is the number of all column-strict
tableaux $T$ of shape $\lambda/\mu$ satisfying $\tau=\kappa
+ \cont T$ and the assertion $\mathcal{H}^{\left(  \kappa\right)
}$, whereas \cite[Theorem 1.2]{Gasharov} says that $\left(  s_{\lambda/\mu
},s_{\tau/\kappa}\right)  _{\Lambda}$ is the number of all all column-strict
tableaux $T$ of shape $\lambda/\mu$ satisfying $\tau=\kappa
+ \cont T$ and the assertion $\mathcal{F}^{\left(  \kappa\right)
}$. (Notice that Gasharov's proof of \cite[Theorem 1.2]{Gasharov} uses the
same involutions as Stembridge's proof of Theorem~\ref{Stembridge-theorem}; it
can thus be regarded as a close precursor to Stembridge's proof. However, it
uses the Jacobi-Trudi identities, while Stembridge's does not.)

% [DG][v54] Added the above two exercises. As exercises they are
% worthless (straightforward reformulations and generalizations of
% things already done), but I think their statements are worth
% making. If I ever come around to writing up the Zelevinsky
% version of the skew LR rule in terms of pictures, these
% exercises will be the springboard.

\begin{exercise}
\label{exe.liri.jordan}

Let $\mathbb{K}$ be a field.\footnote{This field has no relation
to the ring $\kk$, over which our symmetric functions are defined.}
If $N \in \mathbb{K}^{n \times n}$ is
a nilpotent matrix, then the \dfn{Jordan type} of $N$ is
defined to be the list of the sizes of the Jordan blocks in the
Jordan normal form of $N$, sorted in decreasing order\footnote{The
Jordan normal form of $N$ is well-defined even if $\mathbb{K}$ is
not algebraically closed, because $N$ is nilpotent (so the
characteristic polynomial of $N$ is $X^n$).}. This Jordan type
is a partition of $n$, and uniquely determines $N$ up to
similarity (i.e., two nilpotent $n \times n$-matrices
$N$ and $N'$ are similar if and only if the Jordan types of $N$
and $N'$ are equal). If $f$ is a nilpotent endomorphism of a
finite-dimensional $\mathbb{K}$-vector space $V$, then we
define the \dfn{Jordan type} of $f$ as the Jordan type of any
matrix representing $f$ (the choice of the matrix does not matter,
since the Jordan type of a matrix remains unchanged under
conjugation).

\begin{itemize}

\item[(a)] Let $n \in \NN$. Let $N \in \mathbb{K}^{n \times n}$
be a nilpotent matrix. Let $\lambda \in \Par_n$.
Show that the matrix $N$ has Jordan type $\lambda$ if and only
if every $k \in \NN$ satisfies
\[
\dim \left(\ker \left(N^k\right) \right)
= \left(\lambda^t\right)_1 + \left(\lambda^t\right)_2 +
  \ldots + \left(\lambda^t\right)_k .
\]
(Here, we are using the notation $\lambda^t$ for the transpose
of a partition $\lambda$, and the notation $\nu_i$ for the
$i$-th entry of a partition $\nu$.)

\item[(b)] Let $f$ be a nilpotent endomorphism of a
finite-dimensional $\mathbb{K}$-vector space $V$. Let $U$
be an $f$-stable $\mathbb{K}$-vector subspace of $V$ (that is, a
$\mathbb{K}$-vector subspace of $V$ satisfying
$f\left(U\right) \subset U$). Then, restricting $f$ to $U$ gives
a nilpotent endomorphism $f \mid U$ of $U$, and the endomorphism
$f$ also induces a nilpotent endomorphism $\overline f$ of the
quotient space $V/U$. Let $\lambda$, $\mu$ and $\nu$ be the
Jordan types of $f$, $f \mid U$ and $\overline f$, respectively.
Show that $c^\lambda_{\mu, \nu} \neq 0$ (if $\mathbb Z$ is a
subring of $\kk$).

\end{itemize}

[\textbf{Hint:} For (b), Exercise~\ref{exe.pieri.omega}(c)
shows that it is enough to prove that
$c^{\lambda^t}_{\mu^t, \nu^t} \neq 0$. Due to
Corollary~\ref{L-R-rule}, this only requires constructing a
column-strict tableau $T$ of shape $\lambda^t / \mu^t$ with
$\cont T = \nu^t$ which has the property that each
$\cont\left(T | _{\cols \geq j}\right)$ is a partition.
Construct this tableau by defining
$a_{i, j} = \dim \left( \left(f^i\right)^{-1}\left(U\right)
   \cap \ker\left(f^j\right) \right)$
for all $\left( i, j \right) \in \NN^2$,
and requiring that the number of entries $i$ in the $j$-th
row of $T$ be
$a_{i, j} - a_{i, j-1} - a_{i-1, j} + a_{i-1, j-1}$
for all $\left(i, j\right) \in \left\{1,2,3,\ldots\right\}^2$.
Use Exercise~\ref{exe.liri.yamanouchi}(a) to prove that this
indeed defines a column-strict tableau, and
Exercise~\ref{exe.liri.yamanouchi}(b) to verify that it
satisfies the condition on
$\cont\left(T | _{\cols \geq j}\right)$.]

\end{exercise}

\begin{remark}
Exercise~\ref{exe.liri.jordan} is a taste of the connections
between the combinatorics of partitions and the Jordan normal
form. Much more can, and has, been said. Marc van Leeuwen's
\cite{Leeuwen-geometry} is dedicated to some of these
connections; in particular, our
Exercise~\ref{exe.liri.jordan}(a) is
\cite[Proposition 1.1]{Leeuwen-geometry}, and a far stronger
version of Exercise~\ref{exe.liri.jordan}(b) appears in
\cite[Theorem 4.3 (2)]{Leeuwen-geometry}, albeit only for the
case of an infinite $\mathbb{K}$. One can prove a converse to
Exercise~\ref{exe.liri.jordan}(b) as well: If
$c^\lambda_{\mu, \nu} \neq 0$, then there exist $V$, $f$ and
$U$ satisfying the premises of
Exercise~\ref{exe.liri.jordan}(b). When $\mathbb{K}$ is a
finite field, we can ask enumerative questions, such as how
many $U$'s are there for given $V$, $f$, $\lambda$, $\mu$
and $\nu$; we will see a few answers in
Section~\ref{Hall-algebra-section} (specifically,
Proposition~\ref{prop.hall-alg.product}), and a more
detailed treatment is given in \cite[Ch. 2]{Macdonald}.

The relationship between partitions and Jordan normal forms
can be exploited to provide linear-algebraic proofs of purely
combinatorial facts. See \cite[Sections 6 and 9]{BritzFomin}
for some examples. Note that \cite[Lemma 9.10]{BritzFomin} is
the statement that, under the conditions of
Exercise~\ref{exe.liri.jordan}(b), we have
$\nu \subseteq \lambda$. This is a direct consequence of
Exercise~\ref{exe.liri.jordan}(b) (since
$c^\lambda_{\mu, \nu} \neq 0$ can happen only if
$\nu \subseteq \lambda$).
\end{remark}

% [DG][v29] My proof of the converse of
% Exercise~\ref{exe.liri.jordan}(b) is really ugly: For
% finite $\mathbb{K}$, I derive it from the positivity of
% the Hall polynomial's value at the size of $\mathbb{K}$
% (which follows from Theorem 4.2 in Klein's
% "The multiplication of Schur-functions and extensions
% of $p$-modules", J. London Math. Soc., 43 (1968),
% 280--284), whereas for infinite $\mathbb{K}$ it follows
% from Leeuwen's theorem I have cited (I have not read
% either of these papers, sorry). I still hope there
% is a combinatorial solution, some kind of construction
% of the nilpotent matrix from the Littlewood-Richardson
% tableau, but I wasn't able to make it work.
% 
% The remark likely shows off more ignorance than it
% does knowledge. If you can improve it, please do.

\begin{exercise}
\label{exe.Lambda.staircase-skew-1}
Let $a \in \Lambda$. Prove the following:

\begin{enumerate}
\item[(a)] The set $\left\{ g \in \Lambda \ \mid\ g^{\perp} a = \left(
\omega \left( g \right) \right) ^{\perp} a \right\}  $ is a
$\kk$-subalgebra of $\Lambda$.

\item[(b)] Assume that $e_k^{\perp} a = h_k^{\perp} a$ for each positive
integer $k$. Then, $g^{\perp} a = \left( \omega \left( g \right) \right)
^{\perp} a$ for each $g \in \Lambda$.
\end{enumerate}
\end{exercise}

\begin{exercise}
\label{exe.Lambda.staircase-skew-2}
Let $n \in \NN$. Let $\rho$ be the partition
$\left( n-1, n-2, \ldots, 1 \right)$. Prove that
$s_{\rho / \mu} = s_{\rho / \mu^{t}}$ for every $\mu \in \Par$.
\end{exercise}

\begin{remark}
Exercise~\ref{exe.Lambda.staircase-skew-2} appears in \cite[Corollary
7.32]{ReinerShawVanWilligenburg2007}, and is due to John Stembridge. Using
Remark~\ref{third-L-R-interpretation-remark}, we can rewrite it as yet another
equality between Littlewood-Richardson coefficients: Namely, $c_{\mu,\nu
}^{\rho} = c_{\mu^{t},\nu}^{\rho}$ for any $\mu \in \Par$ and
$\nu \in \Par$.
\end{remark}

% [DG][v59] Added the above two exercises and remark.

\newpage

%%%%%%%%%%%%%%%%%%%%%%%%%%%%%%%%%%%%%
\section{Zelevinsky's structure theory of positive self-dual Hopf algebras}
\label{PSH-section}
%%%%%%%%%%%%%%%%%%%%%%%%%%%%%%%%%%%%%

Chapter~\ref{Sym-section} showed that,
as a $\ZZ$-basis for the Hopf algebra $\Lambda=\Lambda_\ZZ$,
the Schur functions $\{s_\lambda\}$ have two special properties:
they have the \emph{same} structure constants $ c^{\lambda}_{\mu,\nu}$
for their multiplication as for their comultiplication
(Corollary~\ref{self-duality-corollary}), and these structure constants
are all \emph{nonnegative} integers (Corollary~\ref{L-R-rule}).
Zelevinsky \cite[\S 2,3]{Zelevinsky} isolated these two properties as crucial.

\begin{definition}
Say that a connected graded Hopf algebra $A$ over $\kk=\ZZ$ with a
distinguished $\ZZ$-basis  $\{\sigma_\lambda\}$ consisting of homogeneous
elements\footnote{not necessarily indexed by partitions}
is a \dfn{positive self-dual Hopf algebra} (or \dfn{PSH})
if it satisfies the two further axioms
\begin{itemize}
\item{\textbf{(self-duality)}} \index{self-duality}
The same structure constants
$a^{\lambda}_{\mu,\nu}$ appear for the product
$\sigma_\mu \sigma_\nu = \sum_{\lambda} a^{\lambda}_{\mu,\nu} \sigma_\lambda$ and
the coproduct
$\Delta \sigma_\lambda = \sum_{\mu,\nu} a^{\lambda}_{\mu,\nu} \sigma_\mu \otimes \sigma_\nu$.
\item{\textbf{(positivity)}} \index{positivity}
The $a^\lambda_{\mu,\nu}$ are all nonnegative (integers).
\end{itemize}
Call $\{\sigma_\lambda\}$ the \dfn{PSH-basis} of $A$.
\end{definition}
He then developed a beautiful structure theory for PSH's, explaining
how they can be uniquely expressed as tensor products of copies of PSH's each
isomorphic to $\Lambda$ after rescaling their grading.
The next few sections explain this, following his exposition closely.

\subsection{Self-duality implies polynomiality}
\label{PSH-implies-polynomial-section}

We begin with a property that forces a Hopf algebra to have
algebra structure which is a \emph{polynomial} algebra, specifically the
symmetric algebra $\Sym(\Liep)$, where $\Liep$ is the $\kk$-submodule
of primitive elements.

Recall from Exercise~\ref{graded-connected-exercise}(g)
that for a connected graded Hopf algebra $A = \bigoplus_{n=0}^\infty A_n$,
every $x$ in the two-sided ideal
$
I:=\ker \epsilon = \bigoplus_{n > 0} A_n
$
has the property that its comultiplication takes the form
\[
\Delta(x) = 1 \otimes x + x \otimes 1 + \Delta_+(x)
\]
where $\Delta_+(x)$ lies in $I \otimes I$.
Recall also that the elements $x$ for which $\Delta_+(x)=0$ are
called the \dfn{primitives}.  Denote by $\Liep$ the $\kk$-submodule
of primitive elements inside $A$.

Given a PSH $A$ (over $\kk = \ZZ$) with a PSH-basis $\{\sigma_\lambda\}$,
we consider the bilinear form $(\cdot,\cdot)_A : A \times A \to \ZZ$ on $A$
that makes this basis orthonormal.
Similarly, the elements $\{\sigma_\lambda \otimes\sigma_\mu\}$ give
an orthonormal basis for a form $(\cdot,\cdot)_{A \otimes A}$ on $A \otimes A$.
The bilinear form $(\cdot,\cdot)_A$ on the PSH $A$ gives rise
to a $\ZZ$-linear map $A \to A^o$, which is easily seen to be
injective and a $\ZZ$-algebra homomorphism. We thus identify
$A$ with a subalgebra of $A^o$. When $A$ is of finite type,
this map is a Hopf algebra isomorphism, thus allowing us to
identify $A$ with $A^o$.
This is an instance of the following notion of self-duality.

% [DG][v17] Rewrote the above paragraph. It used to speak of
% "a graded Hopf algebra $A$ of finite type over $\kk = \ZZ$ with a
% PSH-basis", which made it nebulous whether the basis was supposed
% to be graded. More importantly, you were regarding the injection
% of $A$ into $A^o$ as an inclusion even if $A$ is not finite-type.

\begin{definition} \phantomsection
\label{def.self-duality}
\begin{itemize}
\item[(a)] If $\left(\cdot, \cdot\right) : V \times W \to \kk$
is a bilinear form on the
product $V \times W$ of two graded $\kk$-modules
$V = \bigoplus_{n \geq 0} V_n$ and $W = \bigoplus_{n \geq 0} W_n$,
then we say that this form $\left(\cdot, \cdot\right)$ is
\emph{graded}\index{graded bilinear form}
if every two distinct nonnegative integers $n$ and $m$ satisfy
$\left(V_n, W_m\right) = 0$ (that is,
if every two homogeneous elements $v \in V$ and $w \in W$ having
distinct degrees satisfy $\left(v, w\right) = 0$).
\item[(b)] If $\left(\cdot, \cdot\right)_V : V \times V \to \kk$ and
$\left(\cdot, \cdot\right)_W : W \times W \to \kk$ are two
symmetric bilinear forms on
some $\kk$-modules $V$ and $W$, then we can canonically define a
symmetric bilinear form $\left(\cdot, \cdot\right)_{V\otimes W}$
on the $\kk$-module $V \otimes W$ by letting
\[
\left(v \otimes w, v' \otimes w'\right)_{V\otimes W}
= \left(v, v'\right)_V \left(w, w'\right)_W
\qquad \qquad \text{for all } v, v' \in V \text{ and }
w, w' \in W.
\]
This new bilinear form is graded if the original two forms
$\left(\cdot, \cdot\right)_V$ and $\left(\cdot, \cdot\right)_W$
were graded (presuming that $V$ and $W$ are graded).
\item[(c)] Say that a bialgebra $A$ is
\emph{self-dual}\index{self-dual bialgebra}
with respect to a
given symmetric bilinear form $(\cdot, \cdot) : A \times A \to \kk$ if
one has $(a, m(b \otimes c))_A = (\Delta(a) , b \otimes c)_{A \otimes A}$
and $(\one_A,a)=\epsilon(a)$ for $a,b,c$ in $A$.
If $A$ is a graded Hopf algebra of finite type, and this
form $\left(\cdot, \cdot\right)$ is graded,
then this is equivalent to the $\kk$-module map $A \rightarrow A^o$
induced by $(\cdot, \cdot)_A$ giving a Hopf algebra homomorphism.
\end{itemize}
\end{definition}

% [DG][v17] Revamped this definition, which was subtly wrong (gradedness
% condition on the form was missing, and it was claimed to yield a Hopf
% algebra isomorphism rather than homomorphism). Parts (a) and (b) are
% not very exciting, but probably necessary (linear algebra texts are
% too lazy to deal with these things).

Thus, any PSH $A$ is self-dual with respect to the bilinear form
$\left(\cdot, \cdot\right)_A$ that makes its PSH-basis orthonormal.

Notice also that the injective $\ZZ$-algebra homomorphism $A \to A^o$
obtained from the bilinear form $\left(\cdot, \cdot\right)_A$
on a PSH $A$
allows us to regard each $f \in A$ as an element of $A^o$.
Thus, for any PSH $A$ and any
$f \in A$, an operator $f^\perp : A \to A$
is well-defined (indeed, regard $f$ as an element of $A^o$, and
apply Definition~\ref{def.skewing}).

% [DG][v64] Inserted the above paragraph.

\begin{proposition}
\label{primitives-split-off-prop}
Let $A$ be a Hopf algebra over $\kk=\ZZ$ or $\kk=\QQ$ which is graded, connected, and self-dual
with respect to a positive definite graded\footnote{That is,
$(A_i,A_j)=0$ for $i \neq j$.}  bilinear form. Then:
\begin{itemize}
\item[(a)] Within the ideal $I$,
the $\kk$-submodule of primitives $\Liep$ is the orthogonal complement to the $\kk$-submodule $I^2$.
\item[(b)] In particular, $\Liep \cap I^2=0$.
\item[(c)] When $\kk=\QQ$, one has $I = \Liep \oplus I^2$.
\end{itemize}
\end{proposition}
\begin{proof}
(a) Note that $I^2=m(I \otimes I)$.  Hence an element $x$ in $I$ lies in the perpendicular
space to $I^2$ if and only if one has for all $y$ in $I \otimes I$ that
\[
0= ( x , m(y) )_A = (\Delta(x), y)_{A \otimes A} = (\Delta_+(x), y)_{A \otimes A}
\]
where the second equality uses self-duality, while the third equality uses
the fact that $y$ lies in $I \otimes I$ and the form $(\cdot,\cdot)_{A \otimes A}$
makes distinct homogeneous components orthogonal.
Since $y$ was arbitrary, this means $x$ is perpendicular
to $I^2$ if and only if $\Delta_+(x)=0$, that is, $x$ lies in $\Liep$.

(b) This follows from (a), since the form $\left(\cdot, \cdot\right)_A$ is positive definite.

(c) This follows from (a) using some basic linear algebra\footnote{Specifically,
either the
existence of an orthogonal projection on a subspace of a finite-dimensional
inner-product space over $\QQ$, or the fact that $\dim\left(W^\perp\right) = \dim V
- \dim W$ for a subspace $W$ of a finite-dimensional inner-product space $V$
over $\QQ$ can be used.}
when $A$ is of finite type
(which is the only case we will ever encounter in practice). See
Exercise~\ref{exe.primitive-split-off-prop.infinite} for the general proof.
\end{proof}

% [DG][v13] I added the three sentences under (b) and (c) to this proof
% because I think that any use of positivity is a nontrivial step in
% algebra and that (c) is not trivial at all in the infinite-type case.
% Am I overthinking things? (Here is why I believe (c) is not trivial:
% In a vector space with infinite basis $(a_1, a_2, ...)$, take the
% positive-definite bilinear form given by $(a_i, a_j) = \delta_{i,j}$,
% and consider the subspace spanned by $a_i - a_{i+1}$ for all $i$.
% The orthogonal complement of this subspace is $0$, and there is *no*
% direct sum decomposition anymore. Is there a good reason why it is
% clear that such situations don't happen in our case?)

% [DG][v17] Added footnote in proof of (c) to stress that everything must
% be done over $\QQ$.

\begin{remark}
One might wonder why we didn't just say $I = \Liep \oplus I^2$ even when $\kk=\ZZ$
in Proposition~\ref{primitives-split-off-prop}(c).
However, this is false
%\footnote{Thanks to Darij Grinberg for pointing this out to me.},
even for $A=\Lambda_\ZZ$:
the second homogeneous component $(\Liep \oplus I^2)_2$ is the
index $2$ sublattice of $\Lambda_2$ which is $\ZZ$-spanned by
$\{p_2,e_1^2\}$, containing $2e_2$, but not containing $e_2$ itself.
\end{remark}

Already the fact that $\Liep \cap I^2=0$ has a strong implication.

\begin{lemma}
\label{general-commutativity-lemma}
A connected graded Hopf algebra $A$ over any ring $\kk$ having $\Liep \cap I^2=0$
must necessarily be commutative (as an algebra).
\end{lemma}
% [DG][v62] "has algebra structure which is commutative"
% --> "must necessarily be commutative (as an algebra)".
\begin{proof}
The component $A_0 = \kk$ commutes with all of $A$.  This forms the base case for
an induction on $i+j$ in which one shows that any elements $x$ in $A_i$ and $y$ in $A_j$
with $i,j>0$ will have $[x,y]:=xy-yx = 0$.
Since $[x,y]$ lies in $I^2$, it suffices to show that
$[x,y]$ also lies in $\Liep$:
\begin{align*}
\Delta[x,y]
&= [\Delta(x),\Delta(y)]\\
&= [1 \otimes x + x \otimes 1+\Delta_+(x),
    1 \otimes y + y \otimes 1+\Delta_+(y)]\\
&= [1 \otimes x + x \otimes 1,1 \otimes y + y \otimes 1] \\
&\qquad   +[1 \otimes x + x \otimes 1,\Delta_+(y)]
   +[\Delta_+(x), 1 \otimes y + y \otimes 1]
   +[\Delta_+(x),\Delta_+(y)]\\
&=[1 \otimes x + x \otimes 1,1 \otimes y + y \otimes 1] \\
&=1 \otimes [x,y] + [x,y] \otimes 1
\end{align*}
showing that $[x,y]$ lies in $\Liep$.
Here the second-to-last equality used the inductive hypotheses:
 homogeneity implies that $\Delta_+(x)$ is a sum of
homogeneous tensors of the form $z_1 \otimes z_2$
satisfying $\deg (z_1 ),\deg (z_2 ) < i$, so that by
induction they will commute with $\one \otimes y, y \otimes \one$,
thus proving that $[\Delta_+ (x) , \one \otimes y + y \otimes \one] = 0$;
a symmetric argument shows $[\one \otimes x+x \otimes \one,\Delta_+ (y)] = 0$,
and a similar argument shows $[\Delta_+ (x) , \Delta_+ (y)] = 0$.
The last equality is an easy calculation, and was done already in
the process of proving
\eqref{commutator-of-primitives}.
\end{proof}

\begin{remark}
Zelevinsky actually shows \cite[Proof of A.1.3, p. 150]{Zelevinsky} that the
assumption of $\Liep \cap I^2 = 0$
(along with hypotheses of unit,
counit, graded, connected, and $\Delta$ being a morphism for multiplication)
already implies the \emph{associativity} of the multiplication in $A$ !
One shows by induction on $i+j+k$ that any $x,y,z$ in $A_i,A_j,A_k$
with $i,j,k>0$ have vanishing \dfn{associator}
$\assoc(x,y,z) := x(yz) - (xy)z$.  In the inductive step, one first
notes that $\assoc(x,y,z)$ lies in $I^2$, and then checks
that $\assoc(x,y,z)$ also lies in $\Liep$, by a calculation very similar to the one above,
repeatedly using the fact that $\assoc(x,y,z)$ is multilinear in its three arguments.
\end{remark}

% [DG][v17] Replaced "the assumption of self-duality" by
% "the assumption of $\Liep \cap I^2 = 0$", since the former
% would only be true if the form is positive definite.

\begin{exercise}
\label{exe.primitive-split-off-prop.infinite}
Prove Proposition~\ref{primitives-split-off-prop}(c) in the general case.
\end{exercise}

% [DG][v13] Added this exercise.

This leads to a general structure theorem.
\begin{theorem}
\label{general-structure-theorem}
If a connected graded Hopf algebra $A$ over
a field $\kk$ of characteristic zero has
$I=\Liep \oplus I^2$,
then the inclusion $\Liep \hookrightarrow A$ extends to a
Hopf algebra isomorphism from the symmetric algebra
$\Sym_\kk(\Liep) \rightarrow A$.  In particular, $A$ is
both commutative and cocommutative.
\end{theorem}

\noindent
Note that the hypotheses of Theorem~\ref{general-structure-theorem}
% [DG][v62] "these hypotheses" -> "the hypotheses of Theorem~\ref{general-structure-theorem}".
are valid, using Proposition~\ref{primitives-split-off-prop}(c),
whenever $A$ is obtained from a PSH (over $\ZZ$) by tensoring with $\QQ$.

\begin{proof}[Proof of Theorem~\ref{general-structure-theorem}]
Since Lemma~\ref{general-commutativity-lemma} implies that $A$ is commutative,
the universal property of $\Sym_\kk(\Liep)$ as a free commutative algebra
on generators $\Liep$ shows that the inclusion
$\Liep \hookrightarrow A$ at least extends to an algebra morphism
$\Sym_\kk(\Liep) \overset{\varphi}{\rightarrow} A$.
Since the Hopf structure on $\Sym_\kk(\Liep)$
makes the elements of $\Liep$ primitive (see
Example~\ref{symmetric-algebra-as-bialgebra-example}),
this $\varphi$ is actually a coalgebra morphism
(since $\Delta \circ \varphi = (\varphi \otimes \varphi) \circ \Delta$
and $\epsilon \circ \varphi =\epsilon$
need only to be checked on algebra generators),
hence a bialgebra morphism,
hence a Hopf algebra morphism
(by Corollary~\ref{cor.bialg-mor-is-Hopf}).
It remains to show that $\varphi$ is surjective, and injective.

For the surjectivity of $\varphi$, note that the hypothesis
$I=\Liep \oplus I^2$ implies that the composite $\Liep \hookrightarrow I \rightarrow I/I^2$
gives a $\kk$-vector space isomorphism.
What follows is a standard argument
to deduce that $\Liep$ generates $A$ as a commutative graded $\kk$-algebra.
One shows by induction
on $n$ that any homogeneous element $a$ in $A_n$ lies in the
$\kk$-subalgebra generated
by $\Liep$.  The base case $n=0$ is trivial as $a$ lies in $A_0 = \kk \cdot 1_A$.
% [DG][v62] Replaced "$a$ lies in $\kk$" by "$a$ lies in $A_0 = \kk \cdot 1_A$".
In the inductive step where $a$ lies in $I$, write $a \equiv p \bmod{I^2}$
for some $p$ in $\Liep$.  Thus $a=p + \sum_i b_i c_i$,
where $b_i,c_i$ lie in $I$ but have strictly smaller degree,
so that by induction they lie in the
subalgebra generated by $\Liep$, and hence so does $a$.

Note that the surjectivity argument did not use the assumption that
$\kk$ has characteristic zero, but we will now use it in the injectivity
argument for $\varphi$, to establish the following
\begin{equation}
\label{primitives-in-symmetric-algebra}
\text{\textbf{Claim:} Every primitive element of $\Sym(\Liep)$ lies in $\Liep=\Sym^1(\Liep)$.}
\end{equation}
Note that this claim fails in positive characteristic, e.g. if $\kk$ has characteristic $2$
then $x^2$ lies in $\Sym^2(\Liep)$, however
\[
\Delta (x^2) = 1 \otimes x^2 + 2 x \otimes x + x^2 \otimes 1
= 1 \otimes x^2  + x^2 \otimes 1.
\]
To prove the claim \eqref{primitives-in-symmetric-algebra},
assume not, so that by gradedness, there must exist some
primitive element $y \neq 0$ lying in some $\Sym^n(\Liep)$ with $n \geq 2$.
This would mean that $f(y) = 0$, where the map $f$ is defined
as the composition
\[
\Sym^n(\Liep) \overset{\Delta}{\longrightarrow}
 \bigoplus_{i+j=n} \Sym^i(\Liep) \otimes \Sym^j(\Liep)
   \overset{\text{projection}}{\longrightarrow}
   \Sym^1(\Liep) \otimes \Sym^{n-1}(\Liep)
\]
of the coproduct $\Delta$ with the component projection
of $\bigoplus_{i+j=n} \Sym^i(\Liep) \otimes \Sym^j(\Liep)$
onto $\Sym^1(\Liep) \otimes \Sym^{n-1}(\Liep)$.
However, one can check on a basis
that the multiplication backward
$\Sym^1(\Liep) \otimes \Sym^{n-1}(\Liep) \overset{m}{\rightarrow} \Sym^n(\Liep)$
has the property that $m \circ f = n \cdot \id_{\Sym^n(\Liep)}$:
Indeed,
\[
(m \circ f)(x_1 \cdots x_n)
 = m \left(\sum_{j=1}^n x_j \otimes x_1 \cdots \widehat{x_{j}} \cdots x_n \right)
 = n \cdot x_1 \cdots x_n
\]
for $x_1,\ldots,x_n$ in $\Liep$.
Then $n \cdot y = m(f(y)) = m(0) = 0$ leads to the contradiction that $y=0$,
since $\kk$ has characteristic zero.
Thus, \eqref{primitives-in-symmetric-algebra} is proven.

Now one can argue the injectivity of the (graded) map\footnote{The grading on $\Sym(\Liep)$ is induced from the grading on $\Liep$, a homogeneous subspace of $I \subset A$ as it is the kernel of the graded map $I \overset{\Delta_+}{\longrightarrow} A \otimes A$.} $\varphi$ by assuming that one has a nonzero homogeneous
element $u$ in $\ker(\varphi)$ of minimum degree.
In particular, $\deg(u) \geq 1$.  Also since $\Liep \hookrightarrow A$, one has that $u$ is not in $\Sym^1(\Liep) = \Liep$, and hence $u$ is not primitive
by \eqref{primitives-in-symmetric-algebra}.  Consequently $\Delta_+(u) \neq 0$, and one can find
a nonzero component $u^{(i,j)}$ of $\Delta_+(u)$ lying in
$\Sym(\Liep)_i \otimes \Sym(\Liep)_j$ for some $i,j > 0$.
Since this forces $i,j < \deg(u)$,
one has that $\varphi$ maps both $\Sym(\Liep)_i,\Sym(\Liep)_j$ injectively into $A_i,A_j$.
Hence the tensor product map
\[
\Sym(\Liep)_i \otimes \Sym(\Liep)_j
\overset{\varphi \otimes \varphi}{\longrightarrow}
A_i \otimes A_j
\]
is also injective\footnote{One needs to know that for two
injective maps $V_i \overset{\varphi_i}{\rightarrow} W_i$
of $\kk$-vector spaces $V_i,W_i$ with $i=1,2$, the
tensor product $\varphi_1 \otimes \varphi_2$ is also injective.
Factoring it as
$
\varphi_1 \otimes \varphi_2 =
\left( \id \otimes \varphi_2 \right) \circ
\left( \varphi_1 \otimes \id \right),
$
one sees that it suffices to show that for an injective
map $V \overset{\varphi}{\hookrightarrow} W$ of free $\kk$-modules,
and any free $\kk$-module $U$, the map
$V \otimes U \overset{\varphi \otimes \id}{\longrightarrow} W\otimes U$ is
also injective.  Since tensor products commute with direct sums,
and $U$ is (isomorphic to) a direct sum of copies of $\kk$, this
reduces to the easy-to-check case where $U=\kk$.

Note that some kind of freeness or flatness hypothesis on $U$ is needed here since, e.g. the injective $\ZZ$-module maps
$\ZZ \overset{\varphi_1 = (\cdot \times 2)}{\longrightarrow} \ZZ$ and
$\ZZ/2\ZZ \overset{\varphi_2= \id}{\longrightarrow} \ZZ/2\ZZ$
have $\varphi_1 \otimes \varphi_2=0$ on
$\ZZ \otimes_\ZZ \ZZ/2\ZZ \cong \ZZ/2\ZZ \neq 0$.}.
This implies $(\varphi \otimes \varphi)(u^{(i,j)}) \neq 0$, giving the contradiction that
\[
0=\Delta^A_+(0) = \Delta^A_+(\varphi(u)) = (\varphi \otimes \varphi)(\Delta^{\Sym(\Liep)}_+(u))
\]
contains the nonzero $A_i \otimes A_j$-component
$(\varphi\otimes \varphi)(u^{(i,j)})$.

(An alternative proof of the injectivity of $\varphi$ proceeds as
follows: By \eqref{primitives-in-symmetric-algebra}, the subspace
of primitive elements of $\Sym(\Liep)$ is $\Liep$, and clearly
$\varphi\mid_{\Liep}$ is injective. Hence,
Exercise~\ref{exe.I=0}(b) (applied to the homomorphism $\varphi$)
shows that $\varphi$ is injective.)
\end{proof}

% [DG][v13] I added the paragraph above; feel free to remove it if it
% doesn't fit.

Before closing this section, we mention one nonobvious corollary of
the Claim \eqref{primitives-in-symmetric-algebra}, when applied to
the ring of symmetric functions $\Lambda_\QQ$ with $\QQ$-coefficients,
since Proposition~\ref{symm-antipode-and-h-basis} says that
$\Lambda_\QQ=\QQ[p_1,p_2,\ldots]=\Sym(V)$ where $V=\QQ\{p_1,p_2,\ldots\}$.

\begin{corollary}
\label{primitives-in-Lambda}
The subspace $\Liep$ of primitives in $\Lambda_\QQ$ is
one-dimensional in each degree $n=1,2,\ldots$, and
spanned by $\{p_1,p_2,\ldots\}$.
\end{corollary}

We note in passing that this corollary can also be obtained in a
simpler fashion and a greater generality:

\begin{exercise}
\label{exe.Lambda.primitives}
Let $\kk$ be any commutative ring. Show that the primitive
elements of $\Lambda$ are precisely the elements of the
$\kk$-linear span of $p_1, p_2, p_3, \ldots$.
\end{exercise}

% [DG][v28] Added this exercise (one of the few exercises I can
% find for this section).

\subsection{The decomposition theorem}

Our goal here is Zelevinsky's theorem \cite[Theorem 2.2]{Zelevinsky} giving a canonical
decomposition of any PSH as a tensor product into PSH's that each
have only one primitive element in their PSH-basis.  For the sake
of stating it, we introduce some notation.

\begin{definition}
\label{indecomposable-PSH-factors-defn}
Given a PSH $A$ with
% [DG][v62] "Given $A$ a PSH" -> "Given a PSH $A$".
PSH-basis $\Sigma$,
let $\CCC:=\Sigma \cap \Liep$\index{$\CCC$} be the
primitive elements in  $\Sigma$.
For each $\rho$ in $\CCC$, let
$A(\rho) \subset A$ be the $\ZZ$-span of
\[
\Sigma(\rho):=\{\sigma \in \Sigma:
              \text{ there exists }n \geq 0\text{ with }
                (\sigma, \rho^n) \neq 0\}.
\]
\end{definition}

\begin{definition}
The
tensor product of two PSHs $A_1$ and $A_2$ with PSH-bases
$\Sigma_1$ and $\Sigma_2$ is defined as the graded Hopf algebra
$A_1 \otimes A_2$ with PSH-basis
$\left\{\sigma_1 \otimes \sigma_2\right\}
_{\left(\sigma_1, \sigma_2\right) \in \Sigma_1 \times \Sigma_2}$.
It is easy to see that this is again a PSH.
The tensor product of any finite family of PSHs is defined
similarly\footnote{For the empty family, it is the connected graded
Hopf algebra $\ZZ$ with PSH-basis $\left\{1\right\}$.}.
\end{definition}

% [DG][v62] Moved the above definition into an actual
% definition environment (it used to be in a footnote).

\begin{theorem}
\label{Zelevinsky-decomposition-theorem}
Any PSH $A$ has a canonical tensor product decomposition
\[
A = \bigotimes_{ \rho \in \CCC } A(\rho)
\]
with $A(\rho)$ a PSH, and $\rho$ the only
primitive element in its PSH-basis $\Sigma(\rho)$.
\end{theorem}
\noindent
Although in all the applications, $\CCC$ will be finite,
when $\CCC$ is infinite one should interpret the tensor product
in the theorem as the inductive limit of tensor products over finite subsets of $\CCC$,
that is, linear combinations of basic tensors $\bigotimes_\rho {a_\rho}$
in which there are only finitely many factors $a_\rho \neq 1$.

% [DG][v14] Added footnote about tensor products of PSHs.

The first step toward the theorem uses a certain unique
factorization property.

\begin{lemma}
\label{unique-factorization-lemma}
Let $\PPP$ be a set of pairwise orthogonal primitives in a PSH $A$.
Then,
\[
( \rho_1 \cdots \rho_r , \pi_1 \cdots \pi_s ) = 0
\]
for $\rho_i, \pi_j$ in $\PPP$ unless $r=s$ and
one can reindex so that $\rho_i=\pi_i$.
\end{lemma}
\begin{proof}
%Induct on $\min(r,s)$. 
Induct on $r$.
For $r > 0$, one has
\begin{align*}
(\rho_1 \cdots \rho_r , \pi_1 \cdots \pi_s )
&= (\rho_2 \cdots \rho_r , \rho_1^{\perp}(\pi_1 \cdots \pi_s) ) \\
&= (\rho_2 \cdots \rho_r ,
  \sum_{j=1}^s (\pi_1 \cdots \pi_{j-1} \cdot \rho_1^{\perp}(\pi_j) \cdot \pi_{j+1} \cdots \pi_s) )
\end{align*}
from Proposition~\ref{skewing-properties-prop}(iv)
because $\rho_1$ is primitive\footnote{Strictly speaking, this argument needs
further justification since $A$ might not be of finite type (and if it is
not, Proposition~\ref{skewing-properties-prop}(iv) cannot be applied).
It is more adequate to refer to the proof of
Proposition~\ref{skewing-properties-prop}(iv), which indeed goes through
with $\rho_1$ taking the role of $f$.}. 
 On the other hand, since each $\pi_j$ is
primitive, one has
$
\rho_1^\perp(\pi_j) = (\rho_1,\one) \cdot \pi_j + (\rho_1,\pi_j) \cdot \one
=(\rho_1,\pi_j)
$
which vanishes unless $\rho_1 = \pi_j$.
Hence $(\rho_1 \cdots \rho_r , \pi_1 \cdots \pi_s )=0$ unless $\rho_1 \in
\{\pi_1,\ldots,\pi_s\}$, in which case after reindexing so that
$\pi_1=\rho_1$, it equals
\[
n \cdot (\rho_1,\rho_1) \cdot
(\rho_2 \cdots \rho_r , \pi_2 \cdots \pi_s )
\]
if there are exactly $n$ occurrences of $\rho_1$ among
$\pi_1,\ldots,\pi_s$.  Now apply induction.
\end{proof}

% [DG][v14] Replaced "Induct on $\min(r,s)$" by "Induct on $r$" for minimal
% simplification.

% [DG][v17] Added footnote above.

So far the positivity hypothesis for a PSH has played little role.
Now we use it to introduce a certain partial order on the PSH $A$,
and then a semigroup grading.

\begin{definition}
For a subset $S$ of an abelian group,
let $\ZZ S$ (resp. $\NN S$) denote the subgroup of $\ZZ$-linear combinations
(resp. submonoid of $\NN$-linear combinations\footnote{Recall that
$\NN:=\{0,1,2,\ldots\}$.}) of the elements of $S$.

In a PSH $A$ with PSH-basis $\Sigma$, the subset $\NN \Sigma$ forms
a submonoid, and lets one define a partial order on $A$ via $a \leq b$
if $b-a$ lies in $\NN\Sigma$.
\end{definition}

We note a few trivial properties of this partial order:
\begin{itemize}
\item
The positivity hypothesis implies that
$
\NN \Sigma \cdot \NN \Sigma \subset \NN \Sigma.
$
\item
Hence multiplication by an element $c \geq 0$
(meaning $c$ lies in $\NN \Sigma$) preserves the order:
$a \leq b$ implies $ac \leq bc$ since
$(b-a)c$ lies in $\NN \Sigma$.
\item
Thus $0 \leq c \leq d$ and $0 \leq a \leq b$ together imply
$
ac \leq bc \leq bd
$.
\end{itemize}

\noindent
This allows one to introduce a semigroup grading on $A$.

\begin{definition}
Let $\NN^{\CCC}_{\fin}$ denote the additive submonoid
of $\NN^{\CCC}$ consisting of those $\alpha=(\alpha_\rho)_{\rho \in \CCC}$
with finite support.

Note that for any $\alpha$ in $\NN^{\CCC}_{\fin}$, one has that the product
$\prod_{ \rho \in \CCC } \rho^{\alpha_\rho} \geq 0$.  Define
\[
\Sigma(\alpha):=\{ \sigma \in \Sigma: \sigma \leq \prod_{\rho \in \CCC} \rho^{\alpha_\rho} \},
\]
that is, the subset of $\Sigma$ on which $\prod_{ \rho \in \CCC } \rho^{\alpha_\rho}$
has support.  Also define
\[
A_{(\alpha)} := \ZZ \Sigma(\alpha) \subset A.
\]
\end{definition}

\begin{proposition}
\label{semigroup-grading-prop}
The PSH $A$ has an $\NN^{\CCC}_{\fin}$-semigroup-grading:
one has an orthogonal direct sum decomposition
\[
A = \bigoplus_{ \alpha \in \NN^{\CCC}_{\fin}} A_{(\alpha)}
\]
for which
\begin{align}
\label{product-semigroup-graded}
A_{(\alpha)} A_{(\beta)}  &\subset A_{(\alpha+\beta)} , \\
\label{comultiplication-semigroup-graded}
\Delta A_{(\alpha)} &\subset
  \bigoplus_{\alpha=\beta+\gamma}A_{(\beta)} \otimes A_{(\gamma)}.
\end{align}
\end{proposition}
\begin{proof}
We will make free use of the fact that a PSH $A$ is commutative,
since it embeds in $A \otimes_\ZZ \QQ$, which is commutative by
Theorem~\ref{general-structure-theorem}.

Note that the orthogonality $\left( A_{(\alpha)}, A_{(\beta)} \right) = 0$
for $\alpha \neq \beta$ is equivalent to the assertion that
\[
\left(\prod_{ \rho \in \CCC } \rho^{\alpha_\rho},
\prod_{ \rho \in \CCC } \rho^{\beta_\rho} \right)=0 ,
\]
which follows from Lemma~\ref{unique-factorization-lemma}.

Next let us deal with the assertion \eqref{product-semigroup-graded}.
It suffices to check that when $\tau,\omega$ in $\Sigma$
lie in $A_{(\alpha)}, A_{(\beta)}$, respectively, then $\tau\omega$
lies in $A_{(\alpha+\beta)}$.  But note that any $\sigma$ in $\Sigma$
having $\sigma \leq \tau \omega$ will then have
\[
\sigma \leq \tau \omega
\leq
\prod_{ \rho \in \CCC } \rho^{\alpha_\rho} \cdot
\prod_{ \rho \in \CCC } \rho^{\beta_\rho}
=\prod_{ \rho \in \CCC } \rho^{\alpha_\rho+\beta_\rho}
\]
so that $\sigma$ lies in $A_{(\alpha+\beta)}$.  This means
that $\tau\omega$ lies in $A_{(\alpha+\beta)}$.

This lets us check that $\bigoplus_{ \alpha \in \NN^{\CCC}_{\fin}} A_{(\alpha)}$ exhaust $A$.  It suffices
to check that any $\sigma$ in $\Sigma$ lies in some $A_{(\alpha)}$.  Proceed
by induction on $\deg(\sigma)$, with the case $\sigma=\one$ being
trivial;  the element $\one$ always lies in $\Sigma$, and hence lies
in $A_{(\alpha)}$ for $\alpha=0$.  For $\sigma$ lying in $I$,
one either has
$(\sigma,a) \neq 0$ for some $a$ in $I^2$, or
else $\sigma$ lies in $(I^2)^\perp=\Liep$
(by Proposition~\ref{primitives-split-off-prop}(a)), so that
$\sigma$ is in $\CCC$
and we are done.
If $(\sigma,a) \neq 0$ with $a$ in $I^2$, then
$\sigma$ appears in the support of some
$\ZZ$-linear combination of elements $\tau\omega$ where $\tau,\omega$
lie in $\Sigma$ and have strictly smaller degree than $\sigma$ has.  
There exists at least
one such pair $\tau, \omega$ for which $(\sigma, \tau\omega) \neq 0$,
and therefore $\sigma \leq \tau\omega$.  Then by induction $\tau, \omega$
lie in some $A_{(\alpha)}, A_{(\beta)}$, respectively, so $\tau\omega$ lies in $A_{(\alpha+\beta)}$, and hence $\sigma$ lies
in $A_{(\alpha+\beta)}$ also.

% [DG][v14] Added "than $\sigma$ has".

Self-duality shows that
\eqref{product-semigroup-graded}
implies \eqref{comultiplication-semigroup-graded}:
if $a,b,c$ lie in $A_{(\alpha)}, A_{(\beta)}, A_{(\gamma)}$,
respectively, then
$
(\Delta a, b \otimes c)_{A \otimes A} = (a, bc)_A = 0
$
unless $\alpha=\beta+\gamma$.
\end{proof}

\begin{proposition}
\label{disjoint-support-prop}
For $\alpha, \beta$ in $\NN^\CCC_{\fin}$ with disjoint support, one has a bijection
\[
\begin{array}{rcl}
\Sigma(\alpha) \times \Sigma(\beta) & \longrightarrow & \Sigma(\alpha+\beta) ,\\
(\sigma, \tau) & \longmapsto & \sigma \tau.
\end{array}
\]
Thus, the multiplication map
$
A_{(\alpha)} \otimes A_{(\beta)} \rightarrow A_{(\alpha+\beta)}$
is an isomorphism.
\end{proposition}
\begin{proof}
We first check that for $\sigma_1,\sigma_2$ in  $\Sigma(\alpha)$
and $\tau_1,\tau_2$ in $\Sigma(\beta)$, one has
\begin{equation}
\label{disjoint-support-orthonormality}
(\sigma_1 \tau_1, \sigma_2 \tau_2) = \delta_{(\sigma_1,\tau_1),(\sigma_2,\tau_2)}.
\end{equation}
Note that this is equivalent to showing both
\begin{enumerate}
\item[$\bullet$]
that $\sigma \tau$ lie in $\Sigma(\alpha+\beta)$
so that the map is well-defined, since it shows
$(\sigma \tau, \sigma \tau)=1$, and
\item[$\bullet$]
that the map is injective.
\end{enumerate}
One calculates
\begin{align*}
(\sigma_1 \tau_1, \sigma_2 \tau_2)_A
&=(\sigma_1 \tau_1, m(\sigma_2 \otimes \tau_2))_A \\
&=(\Delta(\sigma_1 \tau_1), \sigma_2 \otimes \tau_2)_{A \otimes A}\\
&=(\Delta(\sigma_1) \Delta(\tau_1), \sigma_2 \otimes \tau_2)_{A \otimes A} .
\end{align*}
Note that due to \eqref{comultiplication-semigroup-graded}, $\Delta(\sigma_1) \Delta(\tau_1)$ lies in $\sum A_{(\alpha'+\beta')} \otimes A_{(\alpha'' + \beta'')}$ where
\begin{align*}
\alpha'+\alpha''&=\alpha ,\\
\beta'+\beta''&=\beta.
\end{align*}
Since $\sigma_2 \otimes \tau_2$
lies in $A_{(\alpha)} \otimes A_{(\beta)}$, the only nonvanishing terms in the inner product come
from those with
\begin{align*}
\alpha'+\beta'&=\alpha ,\\
\alpha''+\beta''&=\beta.
\end{align*}
As $\alpha,\beta$ have disjoint support, this can only happen if
\[
\alpha'=\alpha, \,\, \alpha''=0, \,\, \beta'=0,\,\, \beta''=\beta ;
\]
that is, the only nonvanishing term comes from
$(\sigma_1 \otimes \one)(\one \otimes \tau_1)=\sigma_1 \otimes \tau_1$.
Hence
\[
(\sigma_1 \tau_1, \sigma_2 \tau_2)_A
=(\sigma_1 \otimes \tau_1,\sigma_2 \otimes \tau_2)_{A \otimes A}\\
= \delta_{(\sigma_1,\tau_1),(\sigma_2,\tau_2)}.
\]

To see that the map is surjective, express
\begin{align*}
\prod_{\rho \in \CCC} \rho^{\alpha_\rho}&=\sum_i \sigma_i ,\\
\prod_{\rho \in \CCC} \rho^{\beta_\rho}&=\sum_j \tau_j
\end{align*}
with $\sigma_i \in \Sigma(\alpha)$ and $\tau_j \in \Sigma(\beta)$.
Then each product $\sigma_i \tau_j$ is in $\Sigma(\alpha +\beta)$ by
\eqref{disjoint-support-orthonormality}, and
\[
\prod_{\rho \in \CCC} \rho^{\alpha_\rho+\beta_\rho}=\sum_{i,j} \sigma_i \tau_j
\]
shows that $\{\sigma_i \tau_j\}$ exhausts $\Sigma(\alpha+\beta)$.  This
gives surjectivity.
\end{proof}

\begin{proof}[Proof of Theorem~\ref{Zelevinsky-decomposition-theorem}]
Recall from Definition~\ref{indecomposable-PSH-factors-defn} that
for each $\rho$ in $\CCC$, one defines $A(\rho) \subset A$ to
be the $\ZZ$-span of
\[
\Sigma(\rho):=\{\sigma \in \Sigma:
              \text{ there exists }n \geq 0\text{ with }
                (\sigma, \rho^n) \neq 0\}.
\]
In other words, $A(\rho):=\bigoplus_{n \geq 0} A_{(n \cdot e_\rho)}$
where $e_\rho$ in $\NN^{\CCC}_{\fin}$ is the standard basis element indexed by $\rho$.
Proposition~\ref{semigroup-grading-prop} then shows that $A(\rho)$ is a Hopf subalgebra of $A$.
Since every $\alpha$ in $\NN^{\CCC}_{\fin}$ can be expressed as the
(finite) sum $\sum_{\rho} \alpha_{\rho} e_\rho$, and the $e_\rho$ have disjoint support,
iterating Proposition~\ref{disjoint-support-prop} shows that
$A = \bigotimes_{ \rho \in \CCC } A(\rho)$.  Lastly, $\Sigma(\rho)$ is clearly a PSH-basis
for $A(\rho)$, and if $\sigma$ is any primitive
element in $\Sigma(\rho)$ then $(\sigma,\rho^n) \neq 0$ lets one conclude
via Lemma~\ref{unique-factorization-lemma} that $\sigma=\rho$ (and $n=1$).
\end{proof}

\subsection{$\Lambda$ is the unique indecomposable PSH}

The goal here is to prove the rest of Zelevinsky's structure theory for PSH's.
Namely, if $A$ has only one primitive element $\rho$
in its PSH-basis $\Sigma$, then $A$ must be isomorphic as a PSH to
the ring of symmetric functions $\Lambda$, after one rescales the grading of $A$.
Note that every $\sigma$ in $\Sigma$ has
$\sigma \leq \rho^n$ for some $n$, and hence has degree divisible by
the degree of $\rho$.  Thus one can divide all degrees by that of $\rho$
and assume $\rho$ has degree $1$.

The idea is to find within $A$ and $\Sigma$ a set of elements that
play the role of
\[
\{ h_n=s_{(n)} \}_{n=0,1,2,\ldots}, \qquad
\{ e_n=s_{(1^n)} \}_{n=0,1,2,\ldots}
\]
within $A=\Lambda$ and its PSH-basis of Schur
functions $\Sigma=\{s_\lambda\}$.  Zelevinsky's argument does this by
isolating some properties that turn out to characterize these elements:

\begin{enumerate}
\item[(a)] $h_0=e_0=1$, and $h_1=e_1=:\rho$ has $\rho^2$ a sum of two elements of $\Sigma$, namely
\[
\rho^2=h_2+e_2.
\]
\item[(b)] For all $n=0,1,2,\ldots$, there exist unique elements
$h_n, e_n$ in $A_n \cap \Sigma$ that satisfy
\begin{align*}
h_2^\perp e_n &= 0,\\
e_2^\perp h_n &= 0
\end{align*}
with $h_2,e_2$ being the two elements of $\Sigma$ introduced in (a).
\item[(c)] For $k=0,1,2,\ldots,n$ one has
\begin{align*}
h_k^\perp h_n & = h_{n-k} \text{ and }
  \sigma^\perp h_n = 0 \text{ for } \sigma \in \Sigma \setminus \{h_0,h_1,\ldots,h_n\}, \\
e_k^\perp e_n & = e_{n-k} \text{ and }
  \sigma^\perp e_n = 0 \text{ for } \sigma \in \Sigma \setminus \{e_0,e_1,\ldots,e_n\}. \\
\end{align*}
In particular, $e_k^\perp h_n = 0=h_k^\perp e_n$ for $k \geq 2$.
\item[(d)] Their coproducts are
\begin{align*}
\Delta(h_n) &= \sum_{i+j=n} h_i \otimes h_j,\\
\Delta(e_n) &= \sum_{i+j=n} e_i \otimes e_j.\\
\end{align*}
\end{enumerate}

We will prove Zelevinsky's result \cite[Theorem 3.1]{Zelevinsky}
as a combination of the following two theorems.

\begin{theorem}
\label{uniqueness-of-sym-parts-a-b-c-d}
Let $A$ be a PSH with PSH-basis $\Sigma$ containing only one primitive
$\rho$, and assume that the grading has been rescaled so that
$\rho$ has degree $1$.
Then, after renaming $\rho=e_1=h_1$, one can find unique sequences
$
\{ h_n\}_{n=0,1,2,\ldots}, \{ e_n \}_{n=0,1,2,\ldots}
$
of elements of $\Sigma$ having properties (a),(b),(c),(d) listed above.
\end{theorem}

The second theorem uses the following notion.

\begin{definition}
A \dfn{PSH-morphism} $A \overset{\varphi}{\rightarrow} A'$
between two PSH's $A,A'$ having PSH-bases $\Sigma, \Sigma'$
is a graded Hopf algebra morphism for which
$\varphi(\NN \Sigma) \subset \NN \Sigma'$.  If $A=A'$ and
$\Sigma=\Sigma'$ it will be called a \dfn{PSH-endomorphism}.
If $\varphi$ is an isomorphism and restricts to a bijection
$\Sigma \rightarrow \Sigma'$, it will be called a \dfn{PSH-isomorphism}%
\footnote{This definition is easily seen to be equivalent to
saying that a PSH-isomorphism is an invertible PSH-morphism
whose inverse is again a PSH-morphism.};
if it is both a PSH-isomorphism and an endomorphism, it is
a \dfn{PSH-automorphism}.\footnote{The reader should be warned
that not every invertible PSH-endomorphism is necessarily a
PSH-automorphism. For instance, it is an easy exercise to
check that
$\Lambda \otimes \Lambda \to \Lambda \otimes \Lambda,
\ f \otimes g \mapsto \sum_{(f)} f_1 \otimes f_2 g$
is a well-defined invertible PSH-endomorphism of the PSH
$\Lambda \otimes \Lambda$ with PSH-basis
$\left(s_\lambda \otimes s_\mu\right)_{\left(\lambda,
  \mu\right) \in \Par \times \Par}$,
but not a PSH-automorphism.}
\end{definition}

% [DG][v25] Added above footnote (a counterexample I was looking
% for long ago).

\begin{theorem}
\label{uniqueness-of-sym-parts-e-f-g}
The elements
$
\{ h_n\}_{n=0,1,2,\ldots}, \{ e_n \}_{n=0,1,2,\ldots}
$
in Theorem~\ref{uniqueness-of-sym-parts-a-b-c-d} also satisfy the following.

\begin{enumerate}
\item[(e)] The elements $h_n,e_n$ in $A$ satisfy the same relation \eqref{e-h-relation}
\[
\sum_{i+j=n} (-1)^i e_i h_j = \delta_{0,n}
\]
as their counterparts in $\Lambda$, along with the
property that
\[
A=\ZZ[h_1,h_2,\ldots]=\ZZ[e_1,e_2,\ldots].
\]
\item[(f)]  There is exactly one nontrivial automorphism $A \overset{\omega}{\rightarrow} A$ as
a PSH, swapping $h_n \leftrightarrow e_n$.

\item[(g)]  There are exactly two PSH-isomorphisms $A \rightarrow \Lambda$:
\begin{itemize}
\item one sending $h_n$ to the complete homogeneous symmetric functions $h_n(\xx)$, while sending $e_n$ to the elementary symmetric functions $e_n(\xx)$,
\item the second one (obtained by composing the first with $\omega$) sending
$h_n \mapsto e_n(\xx)$ and $e_n \mapsto h_n(\xx)$.
\end{itemize}
\end{enumerate}
\end{theorem}

Before embarking on the proof, we mention one more bit of convenient
terminology:  say that an element $\sigma$ in $\Sigma$ is a \dfn{constituent}
of $a$ in $\NN \Sigma$ when $\sigma \leq a$, that is, $\sigma$ appears with nonzero coefficient $c_\sigma$
in the  unique expansion $a=\sum_{\tau \in \Sigma} c_\tau \tau$.

\begin{proof}[Proof of Theorem~\ref{uniqueness-of-sym-parts-a-b-c-d}]
One fact that occurs frequently is this:
\begin{equation}
\label{every-sigma-is-a-constituent-of-powers-of-rho}
\text{Every }\sigma\text{ in }\Sigma \cap A_n\text{ is a constituent of }\rho^n.
\end{equation}
This follows from Theorem~\ref{Zelevinsky-decomposition-theorem},
since $\rho$ is the only primitive element of $\Sigma$:  one has
$A=A(\rho)$ and $\Sigma=\Sigma(\rho)$, so that $\sigma$ is a constituent
of some $\rho^m$, and homogeneity considerations force $m=n$.

Notice that $A$ is of finite type (due to
\eqref{every-sigma-is-a-constituent-of-powers-of-rho}). Thus, $A^o$ is
a graded Hopf algebra isomorphic to $A$.

% [DG][v17] Added previous paragraph. I ended up quoting this in three
% further proofs in Chapter 3 as it's subtly used in several places.
% Not the best solution, I guess.

\vskip.1in
\noindent
{\sf Assertion (a).} Note that
\[
(\rho^2,\rho^2)=(\rho^\perp(\rho^2),\rho)=(2\rho,\rho)=2
\]
using the fact that $\rho^\perp$ is a derivation since $\rho$ is primitive
(Proposition~\ref{skewing-properties-prop}(iv)).  On the other hand, expressing
$\rho^2 = \sum_{ \sigma \in \Sigma } c_\sigma \sigma$ with $c_\sigma$ in $\NN$,
one has
$
(\rho^2,\rho^2)=\sum_\sigma c_\sigma^2.
$
Hence exactly two of the $c_\sigma=1$, so $\rho^2$ has exactly two distinct
constituents.   Denote them by $h_2$ and $e_2$.
One concludes that $\Sigma \cap A_2=\{h_2,e_2\}$ from
\eqref{every-sigma-is-a-constituent-of-powers-of-rho}.

Note also that the same argument shows $\Sigma \cap A_1 = \{\rho\}$,
so that $A_1=\ZZ \rho$.  Since $\rho^\perp h_2$ lies in $A_1=\ZZ\rho$
and $(\rho^\perp h_2,\rho)=(h_2,\rho^2)=1$, we have $\rho^\perp h_2=\rho$.
Similarly $\rho^\perp e_2=\rho$.
\vskip.1in
\noindent
{\sf Assertion (b).}
We will show via induction on $n$ the following three assertions for $n \geq 1$:
\begin{equation}
\label{h_n-inductive-properties}
\begin{aligned}
& \bullet \ \text{There exists an
element }h_n\text{ in }\Sigma \cap A_n \text{ with }e_2^\perp h_n=0.\\
& \bullet \ \text{This element }h_n \text{ is unique.} \\
& \bullet \ \text{Furthermore }\rho^\perp h_n = h_{n-1}.\\
\end{aligned}
\end{equation}
In the base cases $n=1,2$, it is not hard to check that
our previously labelled elements, $h_1, h_2$ (namely $h_1:=\rho$,
and $h_2$ as named in part (a)) really \emph{are}
the unique elements satisfying these hypotheses.

In the inductive step, it turns out that
we will find $h_n$ as a constituent of $\rho h_{n-1}$.
Thus we again use the derivation property of $\rho^\perp$ to compute
that $\rho h_{n-1}$ has exactly two constituents:
\begin{align*}
( \rho h_{n-1}, \rho h_{n-1} )
&=( \rho^\perp (\rho h_{n-1}), h_{n-1} ) \\
&=( h_{n-1} + \rho \cdot \rho^\perp h_{n-1}, h_{n-1} )\\
&=( h_{n-1} + \rho h_{n-2}, h_{n-1} )\\
&=1+ ( h_{n-2}, \rho^\perp h_{n-1} )\\
&=1+ ( h_{n-2}, h_{n-2} )=1+1=2
\end{align*}
where the inductive hypothesis $\rho^\perp h_{n-1}=h_{n-2}$ was used twice.
We next show that exactly one of the two constituents of $\rho h_{n-1}$
is annihilated by $e_2^\perp$.  Note that since $e_2$ lies in $A_2$,
and $A_1$ has $\ZZ$-basis element $\rho$, there is a constant $c$ in $\ZZ$
such that
\begin{equation}
\label{coproduct-of-e2}
\Delta(e_2) = e_2 \otimes \one + c \rho \otimes \rho + \one \otimes e_2.
\end{equation}
On the other hand, (a) showed
\[
1=( e_2, \rho^2 )_A = (\Delta (e_2), \rho \otimes \rho )_{A \otimes A}
\]
so one must have $c=1$.
Therefore by Proposition~\ref{skewing-properties-prop}{(iv)} again,
\begin{equation}
\label{e_2-skewing-calculation}
\begin{array}{rcccccl}
e_2^\perp ( \rho h_{n-1} )
&=& e_2^\perp(\rho) h_{n-1} &+& \rho^\perp(\rho) \rho^\perp(h_{n-1}) &+& \rho e_2^{\perp}(h_{n-1}) \\
&=& 0&+&h_{n-2}&+&0\\
&=&h_{n-2} ,
\end{array}
\end{equation}
where the first term vanished due to degree considerations and the last term
vanished by the inductive hypothesis.
Bearing in mind that $\rho h_{n-1}$ lies in $\NN \Sigma$, and
in a PSH with PSH-basis $\Sigma$, any skewing operator
$\sigma^\perp$ for $\sigma$ in $\Sigma$ will
preserve $\NN \Sigma$, one concludes from
\eqref{e_2-skewing-calculation} that
\begin{enumerate}
\item[$\bullet$]
one of the two distinct constituents of the element $\rho h_{n-1}$ must be sent by $e_2^\perp$ to $h_{n-2}$, and
\item[$\bullet$]
the other constituent of $\rho h_{n-1}$ must be annihilated by
$e_2^\perp$; call this second constituent $h_n$.
\end{enumerate}

Lastly, to see that this $h_n$ is unique, it suffices to show that
any element $\sigma$ of $\Sigma \cap A_n$ which is killed by $e_2^\perp$ must be a constituent of
$\rho h_{n-1}$.  This holds for the following reason.
We know $\sigma \leq \rho^n$ by
\eqref{every-sigma-is-a-constituent-of-powers-of-rho},
and hence $0 \neq (\rho^n,\sigma) = (\rho^{n-1},\rho^\perp \sigma)$,
implying that $\rho^\perp \sigma \neq 0$.
On the other hand, since
$0= \rho^\perp e_2^\perp \sigma = e_2^\perp \rho^\perp \sigma$,
one has that $\rho^\perp \sigma$ is annihilated by $e_2^\perp$, and
hence $\rho^\perp \sigma$ must be a (positive) multiple of $h_{n-1}$ by part of
our inductive hypothesis.  Therefore $(\sigma, \rho h_{n-1}) = (\rho^\perp \sigma, h_{n-1})$
is positive, that is, $\sigma$ is a constituent of $\rho h_{n-1}$.

The preceding argument, applied to $\sigma=h_n$, shows that
$\rho^\perp h_n=ch_{n-1}$ for some $c$ in $\{1,2,\ldots\}$.  Since
$(\rho^\perp h_n,h_{n-1})=(h_n,\rho h_{n-1})=1$, this $c$ must be $1$, so
that $\rho^\perp h_n=h_{n-1}$.  This completes the induction step in
the proof of \eqref{h_n-inductive-properties}.

One can then argue, swapping the roles of $e_n, h_n$ in the
above argument, the existence and uniqueness of
a sequence $\{e_n\}_{n=0}^\infty$ in $\Sigma$ satisfying
the properties analogous to
\eqref{h_n-inductive-properties}, with $e_0:=1, e_1:=\rho$.

\vskip.1in
\noindent
{\sf Assertion (c).}
Iterating the property from (b) that $\rho^\perp h_n = h_{n-1}$ shows
that $(\rho^k)^\perp h_n = h_{n-k}$ for $0 \leq k \leq n$.  However one
also has an expansion
\[
\rho^k = c h_k + \sum\limits_{\substack{\sigma \in \Sigma \cap A_k:\\\sigma \neq h_k}} c_ \sigma \sigma
\]
for some integers $c, c_\sigma > 0$, since every $\sigma$ in $\Sigma \cap A_k$
is a constituent of $\rho^k$.  Hence
\[
1 = (h_{n-k}, h_{n-k}) =
((\rho^k)^\perp h_n,(\rho^k)^\perp h_n)
\geq c^2 (h_k^\perp h_n, h_k^\perp h_n)
\]
using Proposition~\ref{skewing-properties-prop}(ii).
Hence if we knew that $h_k^\perp h_n \neq 0$ this would force
\[
h_k^\perp h_n = (\rho^k)^\perp h_n = h_{n-k}
\]
as well as $\sigma^\perp h_n=0$ for all $\sigma \not\in \{h_0,h_1,\ldots,h_n\}$.  But
\[
(\rho^{n-k})^\perp h_k^\perp h_n =
h_k^\perp (\rho^{n-k})^\perp h_n =
h_k^\perp h_k = 1 \neq 0
\]
so $h_k^\perp h_n \neq 0$, as desired.
The argument for $e_k^\perp e_n = e_{n-k}$ is symmetric.

The last assertion in (c) follows if one checks that
$e_n \neq h_n$ for each $n \geq 2$, but this holds since
$e_2^\perp(h_n)=0$ but $e_2^\perp(e_n) = e_{n-2}$.

\vskip.1in
\noindent
{\sf Assertion (d).}
Part (c) implies that
\[
(\Delta h_n, \sigma \otimes \tau)_{A \otimes A}
=(h_n , \sigma \tau)_A
=(\sigma^\perp h_n, \tau)_A = 0
\]
unless $\sigma = h_k$ for some $k=0,1,2,\ldots,n$
and $\tau = h_{n-k}$.  Also one can compute
\[
(\Delta h_n, h_k \otimes h_{n-k})
= (h_n, h_k h_{n-k})
= (h_k^{\perp} h_n, h_{n-k})
\overset{(c)}{=}(h_{n-k}, h_{n-k}) =1.
\]
This is equivalent to the assertion
for $\Delta h_n$ in (d).  The argument for $\Delta e_n$ is symmetric.
\end{proof}

Before proving Theorem~\ref{uniqueness-of-sym-parts-e-f-g}, we note
some consequences of Theorem~\ref{uniqueness-of-sym-parts-a-b-c-d}.
Define for each partition
$\lambda=(\lambda_1 \geq \lambda_2 \geq \cdots \geq \lambda_\ell)$
the following two elements of $A$:
\begin{align*}
h_\lambda & = h_{\lambda_1} h_{\lambda_2} \cdots h_{\lambda_\ell} = h_{\lambda_1} h_{\lambda_2} \cdots, \\
e_\lambda & = e_{\lambda_1} e_{\lambda_2} \cdots e_{\lambda_\ell} = e_{\lambda_1} e_{\lambda_2} \cdots.
\end{align*}
Also, define the
\emph{lexicographic order}\index{lexicographic order on partitions}
on $\Par_n$
by saying $\lambda <_{\lex} \mu$ if $\lambda \neq \mu$ and
the smallest index $i$ for which $\lambda_i \neq \mu_i$ has
$\lambda_i < \mu_i$.  Recall also that $\lambda^t$ denotes the
\emph{conjugate} or \emph{transpose} partition to $\lambda$, obtained by
swapping rows and columns in the Ferrers diagram.

The following unitriangularity lemma will play a role in the proof
of Theorem~\ref{uniqueness-of-sym-parts-e-f-g}(e).
\begin{lemma}
\label{triangularity-lemma}
Under the hypotheses of Theorem~\ref{uniqueness-of-sym-parts-a-b-c-d},
for $\lambda, \mu$ in $\Par_n$, one has
\begin{equation}
e_\mu^\perp h_\lambda =
\begin{cases}
1, & \text{ if }\mu=\lambda^t ; \\
0, & \text{ if }\mu >_{\lex} \lambda^t.
\end{cases}
\label{eq.triangularity-lemma.perp}
\end{equation}
Consequently
\begin{equation}
\det\left[ (e_{\mu^t}, h_\lambda) \right]_{\lambda,\mu \in \Par_n}= 1.
\label{eq.triangularity-lemma.det=1}
\end{equation}
\end{lemma}

% [DG][v64] Replaced "=\pm 1" by "= 1" in the last equation, because
% this is a honest square matrix (with rows and columns both indexed
% by the same set) and is unitriangular, so how could its det be -1?

\begin{proof}
Notice that $A$ is of finite type (as shown in the
proof of Theorem~\ref{uniqueness-of-sym-parts-a-b-c-d}). Thus, $A^o$ is
a graded Hopf algebra isomorphic to $A$.

Also, notice that any $m \in \NN$ and any
$a_1, a_2, \ldots, a_\ell \in A$ satisfy
\begin{equation}
e_m^\perp \left(a_1 a_2 \cdots a_\ell\right)
= \sum_{i_1+\cdots+i_\ell = m} e_{i_1}^\perp\left(a_1\right) \cdots e_{i_\ell}^\perp\left(a_\ell\right) .
\label{pf.triangularity-lemma.1}
\end{equation}
Indeed, this follows by induction over $\ell$ using
Proposition~\ref{skewing-properties-prop}(iv)
(and the coproduct formula for $\Delta(e_n)$ in
Theorem~\ref{uniqueness-of-sym-parts-a-b-c-d}(d)).

In order to prove \eqref{eq.triangularity-lemma.perp},
induct on the length of $\mu$.  If $\lambda$ has length $\ell$, so that $\lambda_1^t=\ell$, then
\begin{align*}
e_\mu^\perp h_\lambda
&= e_{(\mu_2,\mu_3,\ldots)}^\perp
\left(
e_{\mu_1}^\perp (h_{\lambda_1} \cdots h_{\lambda_\ell})
\right)
\qquad \left(\text{since $e_\mu = e_{\mu_1} e_{(\mu_2,\mu_3,\ldots)}$
and thus $e_\mu^\perp = e_{(\mu_2,\mu_3,\ldots)}^\perp \circ e_{\mu_1}^\perp$}\right)
\\
&=e_{(\mu_2,\mu_3,\ldots)}^\perp
\sum_{i_1+\cdots+i_\ell=\mu_1}
e_{i_1}^\perp(h_{\lambda_1}) \cdots
e_{i_\ell}^\perp(h_{\lambda_\ell})
\qquad \left(\text{by \eqref{pf.triangularity-lemma.1}}\right) \\
&=e_{(\mu_2,\mu_3,\ldots)}^\perp
\sum\limits_{\substack{i_1+\cdots+i_\ell=\mu_1; \\ \text{each of } i_1, \ldots, i_\ell \text{ is } \leq 1}}
e_{i_1}^\perp(h_{\lambda_1}) \cdots
e_{i_\ell}^\perp(h_{\lambda_\ell})
\qquad \left(\text{since } e_k^\perp h_n = 0 \text{ for } k \geq 2 \right) \\
&=\begin{cases}
0, & \text{ if }\mu_1 > \ell=\lambda_1^t ;\\
e_{(\mu_2,\mu_3,\ldots)}^\perp  h_{(\lambda_1-1, \ldots,\lambda_\ell-1)}, & \text{ if }\mu_1 = \ell=\lambda_1^t
\end{cases}
\end{align*}
where the last equality used
\[
e_k^\perp(h_n)
=
\begin{cases}
h_{n-1}, & \text{ if }k=1;\\
0, &\text{ if }k \geq 2 .
\end{cases}
\]
Now apply the induction hypothesis, since
$(\lambda_1-1, \ldots,\lambda_\ell-1)^t = (\lambda_2^t,\lambda_3^t,\ldots)$.

To prove \eqref{eq.triangularity-lemma.det=1}, note that
any $\lambda,\mu$ in $\Par_n$ satisfy
$(e_{\mu^t}, h_\lambda)=(e_{\mu^t}^\perp(h_\lambda), 1)=e_{\mu^t}^\perp(h_\lambda)$
(since degree considerations enforce
$e_{\mu^t}^\perp(h_\lambda) \in A_0 = \kk \cdot 1$),
and thus
\[
(e_{\mu^t}, h_\lambda)
= e_{\mu^t}^\perp(h_\lambda)
= \begin{cases}
1, & \text{ if }\mu^t = \lambda^t ; \\
0, & \text{ if }\mu^t >_{\lex} \lambda^t
\end{cases}
\]
(by \eqref{eq.triangularity-lemma.perp}).
This means that the matrix
$\left[ (e_{\mu^t}, h_\lambda) \right]_{\lambda,\mu \in \Par_n}$
is unitriangular with respect to some total order on $\Par_n$
(namely, the lexicographic order on the conjugate partitions),
and hence has determinant $1$.
\end{proof}

% [DG][v64] Made the above proof more detailed.

The following proposition will be the crux of the
proof of Theorem~\ref{uniqueness-of-sym-parts-e-f-g}(f) and (g),
and turns out to be closely related to Kerov's \emph{asymptotic theory of
characters of the symmetric groups} \cite{Kerov}.

\begin{proposition}
\label{normalized-character-prop}
Given a PSH $A$ with PSH-basis $\Sigma$ containing only one primitive $\rho$,
the two maps $A \rightarrow \ZZ$
defined on $A=\bigoplus_{n \geq 0} A_n$ via
\begin{align*}
\delta_h&=\bigoplus_n h_n^\perp,\\
\delta_e&=\bigoplus_n e_n^\perp
\end{align*}
are characterized as the only two
$\ZZ$-linear maps $A \overset{\delta}{\rightarrow} \ZZ$ with the three
properties of being
\begin{itemize}
\item \textbf{positive}: $\delta(\NN \Sigma) \subset \NN$,
\item \textbf{multiplicative}: $\delta(a_1 a_2) = \delta(a_1) \delta(a_2)$ for all $a_1, a_2 \in A$, and
\item \textbf{normalized}: $\delta(\rho)=1$.
\end{itemize}
\end{proposition}
\begin{proof}
Notice that $A$ is of finite type (as shown in the
proof of Theorem~\ref{uniqueness-of-sym-parts-a-b-c-d}). Thus, $A^o$ is
a graded Hopf algebra isomorphic to $A$.

It should be clear from their definitions that $\delta_h, \delta_e$ are $\ZZ$-linear,
positive and normalized.  To see that $\delta_h$ is multiplicative, by $\ZZ$-linearity,
it suffices to check that for $a_1,a_2$ in $A_{n_1}, A_{n_2}$ with $n_1+n_2=n$, one has
\[
\delta_h(a_1 a_2) = h_n^\perp(a_1 a_2)
= \sum_{i_1 + i_2=n} h_{i_1}^\perp(a_1) h_{i_2}^\perp(a_2)
= h_{n_1}^\perp(a_1) h_{n_2}^\perp(a_2) = \delta_h(a_1) \delta_h(a_2)
\]
in which the second equality used Proposition~\ref{skewing-properties-prop}(iv)
and Theorem~\ref{uniqueness-of-sym-parts-a-b-c-d}(d).  The argument for $\delta_e$
is symmetric.

Conversely, given $A \overset{\delta}{\rightarrow} \ZZ$ which is $\ZZ$-linear,
positive, multiplicative, and normalized, note that
\[
\delta(h_2)+\delta(e_2) = \delta(h_2+e_2)
= \delta(\rho^2)
= \delta(\rho)^2
= 1^2 = 1
\]
and hence positivity implies that either $\delta(h_2)=0$ or $\delta(e_2)=0$.
Assume the latter holds, and we will show that $\delta=\delta_h$.

Given any $\sigma$ in $\Sigma \cap A_n \setminus \{h_n\}$, note that
$e_2^\perp \sigma \neq 0$ by Theorem~\ref{uniqueness-of-sym-parts-a-b-c-d}(b), and hence
$
0 \neq (e_2^\perp \sigma, \rho^{n-2}) =( \sigma, e_2 \rho^{n-2}).
$
Thus $\sigma$ is a constituent of $e_2 \rho^{n-2}$, so positivity implies
\[
0 \leq \delta(\sigma) \leq \delta(e_2 \rho^{n-2}) = \delta(e_2) \delta(\rho^{n-2}) = 0.
\]
Thus $\delta(\sigma)=0$ for $\sigma$ in $\Sigma \cap A_n \setminus \{h_n\}$.  Since
$\delta(\rho^n)=\delta(\rho)^n =1^n=1$, this forces $\delta(h_n)=1$, for each $n \geq 0$
(including $n=0$, as
$1=\delta(\rho)=\delta(\rho \cdot 1)=\delta(\rho)\delta(1)=1 \cdot \delta(1)=\delta(1)$).
Thus $\delta=\delta_h$.  The argument when $\delta(h_2)=0$ showing $\delta=\delta_e$ is
symmetric.
\end{proof}

\begin{proof}[Proof of Theorem~\ref{uniqueness-of-sym-parts-e-f-g}]
Many of the assertions of parts (e) and (f) will come from constructing
the unique nontrivial PSH-automorphism $\omega$ of $A$ from the antipode $S$:
for homogeneous $a$ in $A_n$, define $\omega(a):=(-1)^n S(a)$.  We now study some
of the properties of $S$ and $\omega$.

Notice that $A$ is of finite type (as shown in the
proof of Theorem~\ref{uniqueness-of-sym-parts-a-b-c-d}). Thus, $A^o$ is
a graded Hopf algebra isomorphic to $A$.

Since $A$ is a PSH, it is commutative by Theorem~\ref{general-structure-theorem}
(applied to $A \otimes_\ZZ \QQ$).  This implies both that
$S$ is an algebra endomorphism by
Proposition~\ref{antipodes-are-antiendomorphisms} (since
Exercise~\ref{exe.comm-cocomm.anti}(a) shows that the
algebra anti-endomorphisms of a commutative algebra are
the same as its algebra endomorphisms), and that
$S^2=\id_A$ by Corollary~\ref{commutative-implies-involutive-antipode-cor}.
Thus, $\omega$ is an algebra endomorphism and satisfies
$\omega^2 = \id_A$.

Since $A$ is self-dual and the defining diagram \eqref{antipode-diagram}
satisfied by the antipode $S$ is sent to itself when one replaces $A$ by $A^o$
and all maps by their adjoints, one concludes that $S=S^*$ (where $S^*$
means the restricted adjoint $S^* : A^o \to A^o$), i.e., $S$ is self-adjoint.
Since $S$ is an algebra endomorphism, and $S=S^*$, in fact $S$ is also a
coalgebra endomorphism, a bialgebra endomorphism, and a Hopf endomorphism
(by Corollary~\ref{cor.bialg-mor-is-Hopf}).  The same properties
are shared by $\omega$.

Since $\id_A = S^2 = S S^*$, one concludes that $S$ is an isometry, and hence so
is $\omega$.

Since $\rho$ is primitive, one has $S(\rho)=-\rho$ and $\omega(\rho)=\rho$.
Therefore $\omega(\rho^n)=\rho^n$ for $n=1,2,\ldots$.  Use this as follows
to check that $\omega$ is a PSH-automorphism, which amounts to checking that
every $\sigma$ in $\Sigma$ has $\omega(\sigma)$ in $\Sigma$:
\[
( \omega(\sigma),\omega(\sigma) )= (\sigma,\sigma)=1
\]
so that $\pm \omega(\sigma)$ lies in $\Sigma$, but also if $\sigma$ lies in $A_n$, then
\[
( \omega(\sigma), \rho^n )= (\sigma , \omega(\rho^n) )= (\sigma , \rho^n )>0.
\]
In summary, $\omega$ is a PSH-automorphism of $A$, an isometry, and an involution.

Let us try to determine the action of $\omega$ on the $\{h_n\}$.
By similar reasoning as in \eqref{coproduct-of-e2}, one has
\[
\Delta(h_2) = h_2 \otimes \one + \rho \otimes \rho + \one \otimes h_2.
\]
Thus $0 = S(h_2) + S(\rho)\rho + h_2$, and combining this with $S(\rho)=-\rho$,
one has $S(h_2) = e_2$.  Thus also $\omega(h_2)=(-1)^2 S(h_2)=e_2$.

We claim that this forces $\omega(h_n)=e_n$, because $h_2^\perp \omega(h_n)= 0$ via
the following calculation:  for any $a$ in $A$ one has
\begin{align*}
(h_2^\perp \omega(h_n), a)
&=  (\omega(h_n), h_2 a) \\
&=  ( h_n, \omega(h_2 a) ) \\
&=  ( h_n, e_2 \omega(a) ) \\
&=  ( e_2^\perp h_n, \omega(a) ) =  ( 0 , \omega(a) )  =0.
\end{align*}

Consequently the involution $\omega$ swaps $h_n$ and $e_n$, while the
antipode $S$ has $S(h_n)=(-1)^ne_n$ and $S(e_n)=(-1)^n h_n$.  Thus
the coproduct formulas in (d) and definition of the  antipode $S$ imply
the relation \eqref{e-h-relation} between $\{h_n\}$ and $\{e_n\}$.

This relation  \eqref{e-h-relation} also lets one recursively express the
$h_n$ as polynomials with integer coefficients in the $\{e_n\}$,
and vice-versa, so that $\{h_n\}$ and $\{e_n\}$  each generate the
same $\ZZ$-subalgebra $A'$ of $A$.  We wish to show that $A'$ exhausts $A$.

We argue that Lemma~\ref{triangularity-lemma} implies that the \emph{Gram matrix}
$\left[ (h_{\mu},h_{\lambda}) \right]_{\mu,\lambda \in \Par_n}$ has determinant $\pm 1$ as
follows.  Since $\{h_n\}$ and $\{e_n\}$  both generate $A'$,
there exists a $\ZZ$-matrix $(a_{\mu,\lambda})$ expressing
$e_{\mu^t} = \sum_{\lambda} a_{\mu,\lambda} h_\lambda$, and one has
\[
 \left[ (e_{\mu^t}, h_{\lambda}) \right]
= \left[ a_{\mu,\lambda} \right] \cdot
\left[ (h_{\mu},h_{\lambda}) \right].
\]
Taking determinants of these three $\ZZ$-matrices, and
using the fact that the determinant on the left is $1$
(by \eqref{eq.triangularity-lemma.det=1}),
both determinants on the right must also be $\pm 1$.

Now we will show that every $\sigma \in \Sigma \cap A_n$ lies in $A_n'$.
Uniquely express $\sigma=\sigma'+\sigma''$ in which $\sigma'$ lies
in the $\RR$-span $\RR A_n'$ and  $\sigma''$ lies in the real perpendicular
space  $(\RR A_n')^\perp$ inside $\RR \otimes_\ZZ A_n$.
One can compute $\RR$-coefficients $(c_\mu)_{\mu \in \Par_n}$
that express $\sigma' = \sum_\mu c_\mu h_\mu$ by solving the system
\[
\left( \sum_\mu c_\mu h_\mu, h_\lambda \right)
 = (\sigma, h_\lambda) \text{ for }\lambda \in \Par_n.
\]
This linear system is governed by the Gram matrix
$\left[ (h_{\mu},h_{\lambda}) \right]_{\mu,\lambda \in \Par_n}$
with determinant $\pm 1$, and its right side has $\ZZ$-entries since $\sigma, h_\lambda$
lie in $A$.  Hence the solution $(c_\mu)_{\mu \in \Par_n}$ will have $\ZZ$-entries,
so $\sigma'$ lies in $A'$.  Furthermore, $\sigma''= \sigma - \sigma'$ will lie in $A$,
and hence by the orthogonality of $\sigma',\sigma''$,
\[
1=(\sigma,\sigma) = (\sigma',\sigma')+(\sigma'',\sigma'').
\]
One concludes that either $\sigma''=0$, or $\sigma'=0$.  The latter cannot
occur since it would mean that $\sigma=\sigma''$ is perpendicular to all of $A'$.
But $\rho^n=h_1^n$ lies in $A'$, and $(\sigma,\rho^n) \neq 0$.  Thus $\sigma''=0$,
meaning $\sigma=\sigma'$ lies in $A'$.  This completes the proof of assertion (e).
Note that in the process, having shown
$\det(h_{\mu},h_{\lambda})_{\lambda,\mu \in \Par_n}=\pm 1$, one also knows
that $\{h_\lambda\}_{\lambda \in \Par_n}$ are $\ZZ$-linearly
independent, so that $\{h_1,h_2,\ldots\}$ are algebraically independent%
\footnote{by Exercise~\ref{exe.alggen-via-span}(c)},
and $A=\ZZ[h_1,h_2,\ldots]$ is the polynomial algebra generated by
$\{h_1,h_2,\ldots\}$.

For assertion (f), we have seen that $\omega$ gives such a
PSH-automorphism $A \rightarrow A$, swapping
$h_n \leftrightarrow e_n$.  Conversely, given a
PSH-automorphism $A \overset{\varphi}{\rightarrow} A$, consider the
positive, multiplicative, normalized $\ZZ$-linear map
$\delta:=\delta_h \circ \varphi: A \rightarrow \ZZ$.
Proposition~\ref{normalized-character-prop} shows that either
\begin{itemize}
\item $\delta=\delta_h$, which then forces $\varphi(h_n)=h_n$ for all $n$,
so $\varphi=\id_A$, or
\item
$\delta=\delta_e$,
which then forces $\varphi(e_n)=h_n$ for all $n$,
so $\varphi=\omega$.
\end{itemize}

For assertion (g), given a PSH $A$ with PSH-basis $\Sigma$ having
exactly one primitive $\rho$, since we have seen $A=\ZZ[h_1,h_2,\ldots]$,
where $h_n$ in $A$ is as defined in
Theorem~\ref{uniqueness-of-sym-parts-a-b-c-d},
one can uniquely define an algebra morphism
$A \overset{\varphi}{\rightarrow} \Lambda$ that sends the element
$h_n$ to the
complete homogeneous symmetric function $h_n(\xx)$.  Assertions (d) and (e)
show that $\varphi$ is a bialgebra isomorphism, and hence it is a Hopf
isomorphism.  To show that it is a PSH-isomorphism, we first note that
it is an isometry because one can iterate
Proposition~\ref{skewing-properties-prop}(iv) together with assertions (c) and
(d) to compute all inner products
\[
( h_\mu, h_\lambda )_A
=( \one , h_\mu^\perp h_\lambda )_A
=(1, h_{\mu_1}^\perp h_{\mu_2}^\perp\cdots \left( h_{\lambda_1} h_{\lambda_2} \cdots \right) )_A
\]
for $\mu, \lambda$ in $\Par_n$.  Hence
\[
( h_\mu, h_\lambda )_A
= ( h_\mu(\xx), h_\lambda(\xx) )_\Lambda
= ( \varphi(h_\mu), \varphi(h_\lambda) )_\Lambda .
\]
Once one knows $\varphi$ is an isometry, then elements
$\omega$ in $\Sigma \cap A_n$
are characterized in terms of the form $(\cdot,\cdot)$
by $(\omega,\omega)=1$ and $(\omega, \rho^n) > 0$.
Hence $\varphi$ sends each $\sigma$ in $\Sigma$ to a
Schur function $s_\lambda$, and is a PSH-isomorphism.
\end{proof}

\newpage

%%%%%%%%%%%%%%%%%%%%%%%%%%%%%%%%%%%%%
\section{Complex representations
for $\Symm_n$, wreath products, $GL_n(\FF_q)$}
\label{representations-section}
%%%%%%%%%%%%%%%%%%%%%%%%%%%%%%%%%%%%%

After reviewing the basics that we will need from representation and character theory of finite groups, we give Zelevinsky's three main examples of PSH's arising as spaces of virtual characters for three towers of finite groups:
\begin{enumerate}
\item[$\bullet$]
\emph{symmetric} groups,
\item[$\bullet$]
their \emph{wreath products} with any finite group, and
\item[$\bullet$]
the finite \emph{general linear} groups.
\end{enumerate}

Much in this chapter traces its roots to Zelevinsky's book \cite{Zelevinsky}. The results concerning the symmetric groups, however, are significantly older and spread across the literature: see, e.g., \cite[\S 7.18]{Stanley}, \cite[\S 7.3]{Fulton}, \cite[\S I.7]{Macdonald}, \cite[\S 4.7]{Sagan}, \cite{Knutson}, for proofs using different tools.

% [DG][v14] Added previous two sentences.

\subsection{Review of complex character theory}
\label{character-theory-review-section}

We shall now briefly discuss some basics of representation (and character) theory that will be used below. A good source for this material, including the crucial Mackey formula, is Serre \cite[Chaps. 1-7]{Serre}.\footnote{More advanced treatments of representation theory can be found in \cite{Webb} and \cite{Etingof-et-al}.}

% [DG][v58] Added footnote.

\subsubsection{Basic definitions, Maschke, Schur}

For a group $G$, a \emph{representation of $G$}\index{representation of a group} is a homomorphism $G \overset{\varphi}{\rightarrow} GL(V)$ for some vector space $V$ over a field.  We will take the field to be $\CC$ from now on, and we will also assume that $V$ is finite-dimensional over $\CC$.  Thus a representation of $G$ is the same as a finite-dimensional (left) $\CC G$-module $V$. (We use the notations $\CC G$ and $\CC \left[G\right]$ synonymously for the group algebra of $G$ over $\CC$. More generally, if $S$ is a set, then $\CC S = \CC \left[S\right]$ denotes the free $\CC$-module with basis $S$.)

% [DG][v34] Added preceding two sentences.

We also assume that $G$ is finite, so that
Maschke's Theorem\footnote{... which has a beautiful generalization to finite-dimensional Hopf algebras due to Larson and Sweedler; see Montgomery \cite[\S 2.2]{Montgomery}.} says that $\CC G$ is semisimple, meaning that every $\CC G$-module $U \subset V$ has a $\CC G$-module complement $U'$ with $V = U \oplus U'$.  Equivalently, \emph{indecomposable}\index{indecomposable module} $\CC G$-modules are the same thing as \emph{simple}\index{simple module} (=\emph{irreducible}\index{irreducible module}) $\CC G$-modules.

Schur's Lemma implies that for two simple $\CC G$-modules $V_1,V_2$, one has
\[
\Hom_{\CC G}(V_1,V_2) \cong
\begin{cases}
\CC , & \text{ if }V_1 \cong V_2;\\
0, & \text{ if } V_1 \not\cong V_2.
\end{cases}
\]

\subsubsection{Characters and Hom spaces}

A $\CC G$-module $V$ is completely determined up to isomorphism by its
\emph{character}\index{character of a module}
\[
\begin{array}{rcl}
G & \overset{\chi_V}{\longrightarrow} &\CC , \\
g & \longmapsto &\chi_V(g):=\trace(g: V \rightarrow V).
\end{array}
\]
This character $\chi_V$ is a \dfn{class function}, meaning it is
constant on $G$-conjugacy classes.  The space $R_\CC(G)$ of
class functions $G \rightarrow \CC$ has a Hermitian, positive definite
form
\[
(f_1,f_2)_G:=\frac{1}{|G|} \sum_{g \in G} f_1(g) \overline{f_2(g)}.
\]
For any two $\CC G$-modules $V_1, V_2$,
\begin{equation}
\label{character-scalprod}
( \chi_{V_1}, \chi_{V_2} )_G=
\dim_\CC \Hom_{\CC G}(V_1,V_2).
\end{equation}
The set of all \emph{irreducible characters}\index{$\Irr(G)$}\index{irreducible character}
\[
\Irr(G)=\{ \chi_V : V \text{ is a simple }\CC G\text{-module}\}
\]
forms an orthonormal basis of $R_\CC(G)$
with respect to this form, and spans a $\ZZ$-sublattice
\[
R(G):=\ZZ \Irr(G) \subset R_\CC(G)
\]
sometimes called the \emph{virtual characters}\index{$R(G)$}\index{virtual character}
of $G$.
For every $\CC G$-module $V$, the character $\chi_V$
belongs to $R(G)$.

% [DG][v25] Added the previous sentence.

Instead of working with the Hermitian form $\left(\cdot, \cdot\right)_G$
on $G$, we could also (and some authors do) define a $\CC$-bilinear
form $\left<\cdot,\cdot\right>_G$ on $R_\CC (G)$ by
\[
\left<f_1,f_2\right>_G := \frac{1}{|G|} \sum_{g \in G} f_1(g) f_2(g^{-1}).
\]
This form is not identical with $\left(\cdot, \cdot\right)_G$ (indeed,
$\left<\cdot,\cdot\right>_G$ is bilinear while $\left(\cdot,\cdot\right)_G$
is Hermitian), but it still satisfies \eqref{character-scalprod}, and
thus is identical with $\left(\cdot, \cdot\right)_G$ on $R(G) \times R(G)$.
Hence, for all we are going to do until
Section~\ref{Hall-algebra-section}, we could just as well use the form
$\left<\cdot,\cdot\right>_G$ instead of $\left(\cdot,\cdot\right)_G$.

% [DG] I added the previous paragraph for ideological reasons; feel
% free to take it out.

% [DG][v14] Added clarification that the bilinear-vs.-Hermitian
% question becomes important in the Hall algebra subsection. I am not
% fully sure what exactly is a better choice for that section.

\subsubsection{Tensor products}

Given two groups $G_1,G_2$ and $\CC G_i$-modules $V_i$ for $i=1,2$,
their tensor product $V_1 \otimes_\CC V_2$ becomes a $\CC[G_1 \times G_2]$-module
via $(g_1,g_2)(v_1 \otimes v_2)=g_1(v_1) \otimes g_2(v_2)$.  This module
is called the
\emph{(outer) tensor product}\index{outer tensor product}\index{tensor product of representations}
of $V_1$ and $V_2$.
When $V_1, V_2$ are
both simple, then so is $V_1 \otimes V_2$, and every
simple $\CC[G_1 \times G_2]$-module arises this way
(with $V_1$ and $V_2$ determined uniquely up to
isomorphism).\footnote{This is proven in
\cite[\S 3.2, Thm. 10]{Serre}. The fact that $\CC$ is
algebraically closed is essential for this!}
Thus one has identifications and isomorphisms
\begin{align*}
\Irr(G_1 \times G_2) &= \Irr(G_1) \times \Irr(G_2) ,\\
R(G_1 \times G_2) &\cong R(G_1) \otimes_{\ZZ} R(G_2) ;
\end{align*}
here, $\chi_{V_1}\otimes\chi_{V_2}\in R(G_1)\otimes_{\ZZ} R(G_2)$ is being
identified with $\chi_{V_1\otimes V_2}\in R(G_1\times G_2)$ for all
$\CC G_1$-modules $V_1$ and all $\CC G_2$-modules $V_2$.
The latter isomorphism is actually a restriction of the
isomorphism
$R_\CC(G_1 \times G_2) \cong R_\CC(G_1) \otimes_\CC R_\CC(G_2)$
under which every pure tensor
$\phi_1 \otimes \phi_2 \in R_\CC(G_1) \otimes_\CC R_\CC(G_2)$
corresponds to the class function
$G_1 \times G_2 \to \CC, \ \left(g_1, g_2\right) \mapsto
\phi_1\left(g_1\right) \otimes \phi_2\left(g_2\right)$.

% [DG][v45] Added the preceding sentence.

Given two $\CC G_1$-modules $V_1$ and $W_1$ and two $\CC G_2$-modules
$V_2$ and $W_2$, we have
\begin{equation}
\label{character-tensorprod}
\left(\chi_{V_1 \otimes V_2}, \chi_{W_1 \otimes W_2}\right)_{G_1 \times G_2}
= \left(\chi_{V_1}, \chi_{W_1}\right)_{G_1} \left(\chi_{V_2}, \chi_{W_2}\right)_{G_2}.
\end{equation}

% [DG] I added the previous statement for later use. (It is a corollary of
% Proposition 4.0a (b) in reiner-errata12.pdf, but of course it is much
% easier. Maybe Proposition 4.0a would make a nice exercise?)

\subsubsection{Induction and restriction}

Given a subgroup $H < G$ and $\CC H$-module $U$, one can use the fact that $\CC G$ is a $(\CC G, \CC H)$-bimodule to form the \dfn{induced $\CC G$-module}\index{$\Ind_H^G U$}
\[
\Ind_H^G U : = \CC G \otimes_{\CC H} U.
\]
The fact that $\CC G$ is free as a (right-)$\CC H$-module\footnote{... which also has a beautiful generalization to finite-dimensional Hopf algebras due to Nichols and Zoeller; see \cite[\S 3.1]{Montgomery}.}
on basis elements $\{t_{g}\}_{gH \in G/H}$
makes this tensor product easy to analyze.  For example one can
compute its character
\begin{equation}
\label{induced-character-formula}
\chi_{\Ind_H^G U}(g) %= \Ind_H^G(\chi_U)(g)
 =\frac{1}{|H|} \sum\limits_{\substack{k \in G:\\ k g k^{-1} \in H}} \chi_U(k g k^{-1}).
\end{equation}
\footnote{See \cite[\S 7.2, Prop. 20(ii)]{Serre} for the
proof of this equality. (Another proof is given in
\cite[Remark 5.9.2 (the Remark after Theorem 4.32 in the arXiv version)]{Etingof-et-al},
but \cite{Etingof-et-al} uses a different definition of
$\Ind_H^G U$; see Remark~\ref{rmk.induction.ind-as-hom}
for why it is equivalent to ours. Yet another proof of
\eqref{induced-character-formula} is given in
Exercise~\ref{exe.Indrho}(k).)}
One can also recognize when a $\CC G$-module
$V$ is isomorphic to $\Ind_H^G U$ for some $\CC H$-module $U$:  this happens
if and only if there is an $H$-stable subspace $U \subset V$ having the
property that $V = \bigoplus_{gH \in G/H} gU$.

The above construction of a $\CC G$-module $\Ind_H^G U$
corresponding to any $\CC H$-module $U$ is part of a
functor $\Ind_H^G$ from the category of $\CC H$-modules to the
category of $\CC G$-modules\footnote{On morphisms, it sends
any $f : U \to U'$ to $\id_{\CC G} \otimes_{\CC H} f :
\CC G \otimes_{\CC H} U \to \CC G \otimes_{\CC H} U'$.};
this functor is called \emph{induction}\index{induction of a $\CC H$-module}.

Besides induction on $\CC H$-modules, one can define induction
on class functions of $H$:

\begin{exercise}
\label{exe.induction.characters}
Let $G$ be a finite group, and $H$ a subgroup of $G$. Let $f
\in R_\CC \left(H\right)$ be a class function. We define the
\emph{induction}\index{induction of a class function}
$\Ind_H^G f$ of $f$ to be the
function $G \to \CC$ given by
\begin{equation}
\left(\Ind_H^G f\right) \left(g\right)
= \frac{1}{\left|H\right|} \sum\limits_{\substack{k \in G:\\ k g k^{-1} \in H}}
  f \left(k g k^{-1}\right)
\qquad \text{for all } g \in G.
\label{eq.exe.induction.characters.eq}
\end{equation}
%(This definition extends the second equality sign in
%\eqref{induced-character-formula} to arbitrary class functions.)

\begin{itemize}
\item[(a)]
Prove that this induction $\Ind_H^G f$ is a class
function on $G$, hence belongs to $R_\CC \left(G\right)$.

\item[(b)]
Let $J$ be a system of right coset\footnote{A \dfn{right coset}
of a subgroup $H$ in a group $G$ is defined to be a subset of $G$
having the form $Hj$ for some $j \in G$. Similarly, a
\dfn{left coset} has the form $jH$ for some $j \in G$.}
representatives for
$H \backslash G$, so that $G = \bigsqcup_{j \in J} H j$. Prove that
\[
\left(\Ind_H^G f\right)\left(g\right) =
\sum\limits_{\substack{j\in J:\\jgj^{-1}\in H}} f\left(jgj^{-1}\right)
\qquad \text{for all } g \in G.
\]

\end{itemize}
\end{exercise}

% [DG][v14] Added the exercise above.

The induction $\Ind_H^G$ defined in
Exercise~\ref{exe.induction.characters} is a $\CC$-linear map
$R_\CC \left(H\right) \to R_\CC \left(G\right)$. Since every
$\CC H$-module $U$ satisfies
\begin{equation}
\chi_{\Ind_H^G U} = \Ind_H^G(\chi_U)
\label{eq.induction.mods-and-class-functs}
\end{equation}
\footnote{This follows by comparing the value of
$\chi_{\Ind_H^G U}(g)$ obtained from
\eqref{induced-character-formula} with the value of
$\left(\Ind_H^G(\chi_U)\right)(g)$ found using
\eqref{eq.exe.induction.characters.eq}.}, this $\CC$-linear map
$\Ind_H^G$ restricts to a $\ZZ$-linear map $R\left(H\right) \to
R\left(G\right)$ (also denoted $\Ind_H^G$) which sends the
character $\chi_U$ of any
$\CC H$-module $U$ to the character $\chi_{\Ind_H^G U}$ of the
induced $\CC G$-module $\Ind_H^G U$.

% [DG][v46] Separated the definitions of induced modules and
% induced characters, and explicitly stated the relation between
% these two notions. I think identifying characters with
% modules can be rather confusing.
% Similar changes have been done below for restriction,
% inflation and $K$-fixed space construction. (Sometimes this
% boiled down to defining them for characters.)

\begin{exercise}
\label{exe.induction.transitivity}
Let $G$, $H$ and $I$ be three finite groups such that
$I < H < G$. Let $U$ be a $\CC I$-module. Prove that
$\Ind^G_H \Ind^H_I U \cong \Ind^G_I U$. (This fact is
often referred to as the \dfn{transitivity of induction}.)
\end{exercise}

% [DG][v27] Added above exercise.

\begin{exercise}
\label{exe.ind.tensor}
Let $G_1$ and $G_2$ be two groups. Let $H_1 < G_1$ and $H_2 < G_2$ be
two subgroups. Let $U_1$ be a $\CC H_1$-module, and $U_2$ be a
$\CC H_2$-module. Show that
\begin{equation}
\label{induction-vs-tensor-product}
\Ind_{H_1\times H_2}^{G_1\times G_2} \left( U_1 \otimes U_2 \right)
\cong \left(\Ind_{H_1}^{G_1} U_1\right) \otimes \left(\Ind_{H_2}^{G_2} U_2\right)
\end{equation}
as $\CC \left[ G_1 \times G_2 \right]$-modules.
\end{exercise}

% [DG][v3] This is used in the proof of the bialgebra identity.

% [DG][v28] Made this into an exercise.

The \emph{restriction}\index{restriction of a $\CC G$-module}\index{$\Res^G_H V$}
operation $V \mapsto \Res^G_H V$ restricts a
$\CC G$-module $V$ to a $\CC H$-module.
\dfn{Frobenius reciprocity} asserts the adjointness between
$\Ind_H^G$ and $\Res^G_H$
\begin{equation}
\label{Frobenius-reciprocity}
\Hom_{\CC G}( \Ind_H^G U, V ) \cong
\Hom_{\CC H}( U, \Res^G_H V ),
\end{equation}
as a special case ($S=A=\CC G, R=\CC H, B=U,C=V$)
of the general \dfn{adjoint associativity}
\begin{equation}
\label{Hom-tensor-adj}
\Hom_S(A \otimes_R B,C) \cong \Hom_R (B,\Hom_S(A,C))
\end{equation}
for $S,R$ two rings,
$A$ an $(S,R)$-bimodule,
$B$ a left $R$-module,
$C$ a left $S$-module.

We can define not just the restriction of a $\CC G$-module,
but also the \dfn{restriction of a class function} $f \in
R_\CC (G)$. When $H$ is a subgroup of $G$, the
restriction \dfn{$\Res^G_H f$} of an $f \in R_{\CC} (G)$ is defined
as the result of restricting the map $f : G \to \CC$ to $H$.
This $\Res^G_H f$ is easily seen to belong to $R_\CC (H)$, and
so $\Res^G_H$ is a $\CC$-linear map $R_\CC (G) \to R_\CC (H)$.
This map restricts to a $\ZZ$-linear map $R(G) \to R(H)$,
since we have $\Res^G_H \chi_V = \chi_{\Res^G_H V}$
for any $\CC G$-module $V$. Taking characters in
\eqref{Frobenius-reciprocity} (and recalling
$\Res^G_H \chi_V = \chi_{\Res^G_H V}$ and
\eqref{eq.induction.mods-and-class-functs}), we obtain
\begin{equation}
\label{Frobenius-reciprocity-chars}
( \Ind_H^G \chi_U, \chi_V )_G
= ( \chi_U, \Res^G_H \chi_V )_H .
\end{equation}
By bilinearity, this yields the equality
\[
\left( \Ind_H^G \alpha , \beta \right)_G
= \left( \alpha , \Res_H^G \beta \right)_H
\]
for any class functions $\alpha \in R_\CC (H)$ and
$\beta \in R_\CC (G)$ (since $R(G)$ spans $R_\CC(G)$
as a $\CC$-vector space).

\begin{exercise}
\label{exe.induction.ind-as-hom}
Let $G$ be a finite group, and let $H < G$. Let $U$ be a
$\CC H$-module. If $A$ and $B$ are two algebras, $P$ is a
$\left(B,A\right)$-bimodule and $Q$ is a left $B$-module,
then $\Hom_B\left(P, Q\right)$ is a left $A$-module
(since $\CC G$ is a $\left(\CC H, \CC G\right)$-bimodule).
As a
consequence, $\Hom_{\CC H}\left(\CC G, U\right)$ is a
$\CC G$-module. Prove that this $\CC G$-module is
isomorphic to $\Ind^G_H U$.
\end{exercise}

\begin{remark}
\label{rmk.induction.ind-as-hom}
Some texts \emph{define} the induction $\Ind^G_H U$ of a
$\CC H$-module $U$
to be $\Hom_{\CC H}\left(\CC G, U\right)$ (rather than to be
$\CC G \otimes_{\CC H} U$, as we did).\footnote{Or they define
it as a set of morphisms of $H$-sets from $G$ to $U$ (this is
how \cite[Def. 5.8.1 (Def. 4.28 in the arXiv version)]{Etingof-et-al}
defines it); this is easily
seen to be equivalent to $\Hom_{\CC H}\left(\CC G, U\right)$.}
As Exercise~\ref{exe.induction.ind-as-hom} shows, this
definition is equivalent to ours as long as $G$ is finite
(but not otherwise).
\end{remark}

% [DG][v14] Added preceding exercise and remark. (This exercise also
% appears in the Frobenius reciprocity section of
% \cite{Etingof-et-al}.)

Exercise~\ref{exe.induction.ind-as-hom} yields the following
``wrong-way'' version of Frobenius reciprocity:

\begin{exercise}
\label{exe.induction.wrong-way-frob}
Let $G$ be a finite group; let $H < G$. Let $U$ be a $\CC G$-module,
and let $V$ be a $\CC H$-module. Prove that
$\Hom_{\CC G}\left(U, \Ind^G_H V\right)
\cong \Hom_{\CC H}\left(\Res^G_H U, V\right)$.
\end{exercise}

% [DG][v14] Added the above exercise.

\subsubsection{Mackey's formula}

Mackey gave an alternate description of a
module which has been induced and then restricted.  To state it, for a subgroup $H < G$ and $g$ in $G$,
let $H^g:=g^{-1}Hg$ and ${}^gH:=gHg^{-1}$.
Given a $\CC H$-module $U$, say defined by a homomorphism
$H \overset{\varphi}{\rightarrow} GL(U)$,
let $U^{g}$ denote the $\CC[gHg^{-1}]$-module
on the same $\CC$-vector space $U$ defined by
the composite homomorphism
\[
\begin{array}{rcccl}
{}^gH &\longrightarrow&H &\overset{\varphi}{\longrightarrow}&GL(U),\\
h &\longmapsto &g^{-1}hg .& &
\end{array}
\]

\begin{theorem}(Mackey's formula)
\label{Mackey-formula}
Consider subgroups $H, K < G$,  and any $\CC H$-module $U$.
If $\{g_1,\ldots,g_t\}$ are double coset representatives for
$K\backslash G /H$, then
\[
\Res^G_K \Ind_H^G U
\cong \bigoplus_{i=1}^t
\Ind^K_{{}^{g_i}H \cap K} \left( \left( \Res^H_{H \cap K^{g_i}} U \right)^{g_i} \right) .
\]
\end{theorem}
\begin{proof}
In this proof, all tensor product symbols $\otimes$
should be interpreted as $\otimes_{\CC H}$.
Recall $\CC G$ has $\CC$-basis $\{t_g\}_{g \in G}$.
For subsets $S \subset G$,
let $\CC[S]$ denote the $\CC$-span of $\{t_g\}_{g \in S}$ in $\CC G$.

Note that each double coset $K g H$ gives rise to a sub-$(K,H)$-bimodule
$\CC[K g H]$ within $\CC G$, and one has
a $\CC K$-module direct sum decomposition
\[
\Ind_H^G U = \CC G \otimes U
= \bigoplus_{i=1}^t \CC[K g_i H] \otimes U.
\]
Hence it suffices to check for any element $g$ in $G$ that
\[
\CC[K g H] \otimes U
\cong
\Ind^K_{{}^{g}H \cap K} \left( \left( \Res^H_{H \cap K^{g}} U \right)^{g} \right).
\]
Note that ${}^{g}H \cap K$ is the subgroup of $K$ consisting of
the elements $k$ in $K$ for which $kg H =gH$.  Hence by picking $\{k_1,\ldots,k_s\}$ to be coset representatives for
$K/({}^gH \cap K)$,
one disjointly decomposes the double coset
\[
KgH = \bigsqcup_{j=1}^s k_j({}^gH \cap K)gH,
\]
giving a $\CC$-vector space direct sum decomposition
\begin{align*}
\CC[K g H] \otimes U
&=\bigoplus_{j=1}^s \CC\left[k_j \left( {}^gH \cap K\right) g H\right]  \otimes U \\
&\cong
\Ind_{{}^gH \cap K}^K
\left(\CC [\left( {}^gH \cap K \right) g H]  \otimes U\right).
\end{align*}
So it remains to check that one has a $\CC [{}^gH \cap K]$-module isomorphism
\[
\CC [\left( {}^gH \cap K \right) g H]  \otimes U
\cong \left( \Res^H_{H \cap K^g} U \right)^g.
\]
Bearing in mind that, for each $k$ in ${}^g H \cap K$ and $h$ in $H$, one has
$g^{-1}kg$ in $H$ and hence
\[
t_{kgh} \otimes u
= t_g \cdot t_{g^{-1}kg \cdot h}  \otimes u
=t_g \otimes g^{-1}kgh \cdot u ,
\]
one sees that this isomorphism can be defined by mapping
\[
t_{kgh} \otimes u \longmapsto  g^{-1}kgh \cdot u.
\]
\end{proof}

\subsubsection{Inflation and fixed points}

There are two (adjoint) constructions on representations
that apply when one has a normal subgroup $K \triangleleft G$.
Given a $\CC[G/K]$-module $U$, say defined by the
homomorphism $G/K \overset{\varphi}{\rightarrow} GL(U)$,
the \emph{inflation}\index{inflation of a $\CC[G/K]$-module}
of $U$ to a
$\CC G$-module \dfn{$\Inf_{G/K}^G U$} has the same underlying space $U$,
and is defined by the composite
homomorphism $G \rightarrow G/K \overset{\varphi}{\rightarrow} GL(U)$.
We will later use the easily-checked
fact that when $H < G$ is any other subgroup, one has
\begin{equation}
\label{restriction-inflation-commutation}
\Res^G_H \Inf^G_{G/K} U =
\Inf^H_{H/H \cap K} \Res^{G/K}_{H/H \cap K} U.
\end{equation}
(We regard $H / H \cap K$ as a subgroup of $G / K$, since the
canonical homomorphism $H / H \cap K \to G / K$ is injective.)

% [DG][v14] Added previous sentence.

Inflation turns out to be adjoint to the
\dfn{$K$-fixed space construction}\index{fixed space of a $\CC G$-module}
sending a $\CC G$-module $V$ to the $\CC[G/K]$-module\index{$V^K$}
\[
V^K:=\{v \in V: k(v)=v \text{ for }k \in K\} .
\]
Note that $V^K$ is indeed a $G$-stable subspace:
for any $v$ in $V^K$ and $g$ in $G$, one
has that $g(v)$ lies in $V^K$ since an element $k$ in $K$ satisfies
$kg(v) = g \cdot g^{-1}kg(v) = g(v)$ as $g^{-1}kg$ lies in $K$.
One has this adjointness
\begin{equation}
\label{inflation-reciprocity}
\Hom_{\CC G}(\Inf_{G/K}^G U, V) = \Hom_{\CC[G/K]}(U, V^K) ,
\end{equation}
because any $\CC G$-module homomorphism $\varphi$ on the left must
have the property that $k \varphi(u) = \varphi(k(u)) = \varphi(u)$
for all $k$ in $K$, so that $\varphi$ actually lies on the right.

We will also need the following formula for the
character $\chi_{V^K}$ in terms of the character $\chi_V$:
\begin{equation}
\label{deflation-character-formula}
\chi_{V^K}(gK) = \frac{1}{|K|} \sum_{k \in K} \chi_V(gk).
\end{equation}
To see this, note that when one has a $\CC$-linear
endomorphism $\varphi$ on a space $V$ that preserves some
$\CC$-subspace $W \subset V$, if $V \overset{\pi}{\rightarrow} W$
is any idempotent projection onto $W$, then the trace of the restriction
$\varphi|_W$ equals the trace of $\varphi \circ \pi$ on $V$.
Applying this to $W=V^K$ and $\varphi = g$, with
$\pi=\frac{1}{|K|} \sum_{k \in K} k$, gives
\eqref{deflation-character-formula}.\footnote{For
another proof of \eqref{deflation-character-formula}, see
Exercise~\ref{exe.Indrho}(l).}

Another way to restate \eqref{deflation-character-formula} is:
\begin{equation}
\label{deflation-character-formula-v2}
\chi_{V^K}(gK) = \frac{1}{|K|} \sum_{h \in gK} \chi_V(h).
\end{equation}

% [DG][v4] Trivial rewriting, only stated to factor out a very
% simple step from the proof of the bialgebra property below.

Inflation and $K$-fixed space construction can also be defined on
class functions. For inflation, this is particularly easy:
Inflation $\Inf_{G/K}^G f$ of an $f \in R_{\CC} (G/K)$ is defined as
the composition $\xymatrix{G \arsurj[r] & G/K \ar[r]^{f} & \CC}$.
This is a class function of $G$ and thus lies in $R_\CC (G)$.
Thus, inflation $\Inf_{G/K}^G$ is a $\CC$-linear map
$R_{\CC} (G/K) \to R_\CC (G)$. It restricts to a $\ZZ$-linear map
$R(G/K) \to R(G)$, since it is clear that every $\CC (G/K)$-module
$U$ satisfies $\Inf_{G/K}^G \chi_U = \chi_{\Inf_{G/K}^G U}$.

We can also use \eqref{deflation-character-formula}
(or \eqref{deflation-character-formula-v2}) as inspiration for
defining a ``$K$-fixed space construction'' on class
functions. Explicitly, for every
class function $f \in R_\CC (G)$, we define a class function
$f^K \in R_\CC (G/K)$ by
\[
f^K (gK) = \frac{1}{|K|} \sum_{k \in K} f(gk) = \frac{1}{|K|} \sum_{h \in gK} f(h).
\]
The map $\left(\cdot\right)^K : R_\CC (G) \to R_\CC (G/K),
\ f \mapsto f^K$ is $\CC$-linear, and restricts to a $\ZZ$-linear
map $R(G) \to R(G/K)$. Again, we have a compatibility with the
$K$-fixed point construction on modules: We have
$\chi_{V^K} = \left(\chi_V\right)^K$ for every $\CC G$-module
$V$.

Taking characters in \eqref{inflation-reciprocity}, we obtain
\begin{equation}
\label{inflation-reciprocity-chars}
( \Inf_{G/K}^G \chi_U, \chi_V )_G
= ( \chi_U, \chi_V^K )_{G/K}
\end{equation}
for any $\CC \left[G/K\right]$-module $U$ and any $\CC G$-module
$V$ (since $\chi_{\Inf_{G/K}^G U} = \Inf_{G/K}^G \chi_U$
and $\chi_{V^K} = \left(\chi_V\right)^K$). By $\ZZ$-linearity,
this implies that
\[
\left( \Inf_{G/K}^G \alpha , \beta \right)_G
= \left( \alpha, \beta^K \right)_{G/K}
\]
for any class functions $\alpha \in R_\CC \left(G/K\right)$
and $\beta \in R_\CC \left(G\right)$.

There is also an analogue of~\eqref{induction-vs-tensor-product}:

\begin{lemma}
Let $G_1$ and $G_2$ be two groups, and $K_1 < G_1$ and $K_2 < G_2$
be two respective subgroups. Let $U_i$ be a $\CC G_i$-module
for each $i \in \left\{1,2\right\}$. Then,
\begin{equation}
\label{invariants-vs-tensor-product}
\left(U_1 \otimes U_2\right)^{K_1 \times K_2} = U_1^{K_1} \otimes U_2^{K_2}
\end{equation}
(as subspaces of $U_1 \otimes U_2$).
\end{lemma}

\begin{proof}
The subgroup $K_1 = K_1 \times 1$ of $G_1 \times G_2$ acts on $U_1 \otimes U_2$,
and its fixed points are
$\left(U_1 \otimes U_2\right)^{K_1} = U_1^{K_1} \otimes U_2$
(because for a $\CC K_1$-module, tensoring with $U_2$ is the same as
taking a direct power, which clearly commutes with taking fixed points).
Similarly, $\left(U_1 \otimes U_2\right)^{K_2} = U_1 \otimes U_2^{K_2}$.
Now,
\[
\left(U_1 \otimes U_2\right)^{K_1 \times K_2}
= \left(U_1 \otimes U_2\right)^{K_1} \cap \left(U_1 \otimes U_2\right)^{K_2}
= \left(U_1^{K_1} \otimes U_2\right) \cap \left(U_1 \otimes U_2^{K_2}\right)
= U_1^{K_1} \otimes U_2^{K_2}
\]
according to the known linear-algebraic fact stating that if $P$ and $Q$ are
subspaces of two vector spaces $U$ and $V$, respectively, then
$\left(P \otimes V\right) \cap \left(U \otimes Q\right) = P \otimes Q$.
\end{proof}

% [DG][v3] This lemma is used (no longer) tacitly in the proof of
% the bialgebra property of A(GL).

\begin{exercise} \phantomsection
\label{exe.HomxHom}
\begin{itemize}
\item[(a)]
Let $G_1$ and $G_2$ be two groups. Let $V_i$ and $W_i$ be finite-dimensional
$\CC G_i$-modules for every $i \in \left\{1,2\right\}$. Prove that
the $\CC$-linear map
\[
\Hom_{\CC G_1} \left(V_1,W_1\right) \otimes \Hom_{\CC G_2} \left(V_2,W_2\right)
\to \Hom_{\CC \left[G_1\times G_2\right]} \left(V_1 \otimes V_2, W_1 \otimes W_2\right)
\]
sending each tensor $f\otimes g$ to the tensor product $f\otimes g$ of
homomorphisms is a vector space isomorphism.
\item[(b)] Use part (a) to give a new proof
of \eqref{character-tensorprod}.
\end{itemize}
\end{exercise}

% [DG][v23] I have moved this exercise here because its solution uses
% \eqref{invariants-vs-tensor-product}. It used to stand right after
% \eqref{character-tensorprod}, but there is a much simpler proof of
% \eqref{character-tensorprod} by direct computation.

As an aside, \eqref{restriction-inflation-commutation} has a
``dual'' analogue:

\begin{exercise}
\label{exe.ind-fixed-points}
Let $G$ be a finite group, and let $K \triangleleft G$ and $H < G$.
Let $U$ be a $\CC H$-module. As usual, regard $H/\left(H\cap K\right)$
as a subgroup of $G/K$. Show that
$\left(\Ind^G_H U\right)^K
\cong \Ind^{G/K}_{H/\left(H\cap K\right)} \left(U^{H\cap K}\right)$
as $\CC \left[G/K\right]$-modules.
\end{exercise}

% [DG][v14] Added above exercise. (I haven't seen this one anywhere.
% The particular case when $K = G$ is the Proposition on p. 121 of
% \cite{Knutson}, but proven there in a rather ugly way that
% depends on the ground field being a field of char $0$. I assume
% the general case can be found somewhere in group cohomology
% literature...)

Inflation also ``commutes'' with induction:

\begin{exercise}
\label{exe.ind-infl}
Let $G$ be a finite group, and let $K < H < G$ be such that
$K \triangleleft G$. Thus, automatically, $K \triangleleft H$,
and we regard the quotient $H/K$ as a subgroup of $G/K$.
Let $V$ be a $\CC \left[H/K\right]$-module.
Show that
$\Inf^G_{G/K} \Ind^{G/K}_{H/K} V
\cong \Ind^G_H \Inf^H_{H/K} V$
as $\CC G$-modules.
\end{exercise}

% [DG][v27] Added above exercise (to use it later in solving another).
% Note: unlike in Exercise~\ref{exe.ind-fixed-points}, the
% $K < H$ condition is important.

\begin{exercise}
\label{exe.fixed-vs-coinvariants}
Let $G$ be a finite group, and let $K \triangleleft G$.
Let $V$ be a $\CC G$-module. Let $I_{V,K}$ denote the $\CC$-vector
subspace of $V$ spanned by all elements of the form $v-kv$ for
$k \in K$ and $v \in V$.
\begin{itemize}
\item[(a)] Show that $I_{V,K}$ is a $\CC G$-submodule of $V$.
\item[(b)] Let $V_K$ denote the quotient $\CC G$-module
$V / I_{V,K}$. (This module is occasionally called the
\dfn{$K$-coinvariant module of $V$}\index{coinvariant module},
a name it sadly shares with
at least two other non-equivalent constructions in algebra.)
Show that $V_K \cong \Inf^G_{G/K} \left(V^K\right)$ as
$\CC G$-modules. (Use $\operatorname{char} \CC = 0$.)
\end{itemize}
\end{exercise}

% [DG][v14] Added the above exercise, not out of need but because it
% is a piece of culture people should see. The reminder to use
% $\operatorname{char} \CC = 0$ is more of an allusion to the fact
% that many other exercises here do not care about the characteristic
% of the ground field.

In the remainder of this subsection, we shall briefly survey
generalized notions of induction and restriction, defined in terms of a
group homomorphism $\rho$ rather than in terms of a group $G$ and a
subgroup $H$. These generalized notions (defined by van Leeuwen in
\cite[\S 2.2]{Leeuwen-hopf}) will not be used in the rest of
these notes, but they shed some new light on the facts about induction,
restriction, inflation and fixed point construction discussed above. (In
particular, they reveal that some of said facts have common generalizations.)

The reader might have noticed that the definitions of inflation and of
restriction (both for characters and for modules) are similar. In fact, they
both are particular cases of the following construction:

\begin{remark}
\label{rmk.Resrho}
Let $G$ and $H$ be two finite groups, and let $\rho : H \rightarrow G$ be
a group homomorphism.

\begin{itemize}
\item If $f \in R_\CC \left( G \right)$, then the \emph{$\rho$-restriction}\index{$\rho$-restriction of a class function}
$\Res_{\rho} f$ of $f$ is defined as the map
$f \circ \rho : H \rightarrow \CC$. This map is easily
seen to belong to $R_\CC \left( H \right)$.

\item If $V$ is a $\CC G$-module, then the \emph{$\rho$-restriction}\index{$\rho$-restriction of a $\CC G$-module}
$\Res_{\rho} V$ of $V$ is the $\CC H$-module with ground
space $V$ and action given by
\[
h \cdot v = \rho \left( h \right) \cdot v
\ \ \ \ \ \ \ \ \ \ \text{for every } h \in H \text{ and } v\in V .
\]

\end{itemize}

This construction generalizes both inflation and restriction:
If $H$ is a subgroup of $G$, and if $\rho : H \rightarrow G$ is the
inclusion map, then $\Res_{\rho} f = \Res_{H}^{G} f$
(for any $f \in R_\CC \left( G \right)$)
and $\Res_{\rho} V = \Res_{H}^{G} V$
(for any $\CC G$-module $V$). If, instead, we have
$G = H/K$ for a normal subgroup $K$ of $H$, and if
$\rho : H \rightarrow G$ is the projection map, then
$\Res_{\rho} f = \Inf_{H/K}^H f$
(for any $f \in R_\CC \left( H/K \right)$) and
$\Res_{\rho} V = \Inf_{H/K}^H V$
(for any $\CC \left[H/K\right]$-module $V$).
\end{remark}

A subtler observation is that induction and fixed point construction
can be generalized by a common notion. This is the subject of
Exercise~\ref{exe.Indrho} below.

\begin{exercise}
\label{exe.Indrho}
Let $G$ and $H$ be two finite groups, and let $\rho : H \rightarrow G$
be a group homomorphism. We introduce the following notations:

\begin{itemize}
\item If $f \in R_\CC \left( H \right)$, then the
\emph{$\rho$-induction}\index{$\rho$-induction of a class function}
$\Ind_{\rho} f$ of $f$ is a map
$G \rightarrow \CC$ which is defined as follows:
\[
\left( \Ind_{\rho} f \right) \left( g \right)
= \dfrac{1}{\left\vert H \right\vert}
  \sum\limits_{\substack{\left( h, k \right) \in H \times G; \\
                  k \rho\left( h \right) k^{-1} = g}}
  f \left( h \right)
\ \ \ \ \ \ \ \ \ \ \text{for every } g \in G.
\]

\item If $U$ is a $\CC H$-module, then the
\emph{$\rho$-induction}\index{$\rho$-induction of a $\CC H$-module}
$\Ind_{\rho} U$ of $U$ is defined as the $\CC G$-module
$\CC G \otimes_{\CC H} U$, where $\CC G$ is regarded as a
$\left( \CC G , \CC H \right)$-bimodule according to the following
rule: The left $\CC G$-module structure on $\CC G$ is plain
multiplication inside $\CC G$; the right $\CC H$-module structure
on $\CC G$ is induced by the $\CC$-algebra homomorphism
$\CC \left[ \rho \right] : \CC H \rightarrow \CC G$ (thus, it
is explicitly given by
$\gamma \eta
= \gamma \cdot \left( \CC \left[ \rho \right] \right) \eta$
for all $\gamma \in \CC G$ and $\eta \in \CC H$).
\end{itemize}

Prove the following properties of this construction:

\begin{enumerate}
\item[(a)] For every $f \in R_\CC \left( H \right)$, we have
$\Ind_{\rho} f \in R_\CC \left( G \right)$.

\item[(b)] For any finite-dimensional $\CC H$-module $U$, we have
$\chi_{\Ind_{\rho} U} = \Ind_{\rho} \chi_{U}$.

\item[(c)] If $H$ is a subgroup of $G$, and if
$\rho : H \rightarrow G$ is the inclusion map, then
$\Ind_{\rho} f = \Ind_H^G f$ for every
$f\in R_\CC \left( H \right)$.

\item[(d)] If $H$ is a subgroup of $G$, and if
$\rho : H \rightarrow G$ is the inclusion map, then
$\Ind_{\rho} U = \Ind_H^G U$ for every $\CC H$-module $U$.

\item[(e)] If $G = H/K$ for some normal subgroup $K$ of $H$, and if
$\rho : H \rightarrow G$ is the projection map, then
$\Ind_{\rho} f = f^K$ for every $f \in R_\CC \left( H \right)$.

\item[(f)] If $G = H/K$ for some normal subgroup $K$ of $H$, and if
$\rho : H \rightarrow G$ is the projection map, then
$\Ind_{\rho} U \cong U^K$ for every $\CC H$-module $U$.

\item[(g)] Any class functions $\alpha \in R_\CC \left( H \right)$
and $\beta \in R_\CC \left( G \right)$ satisfy
\begin{equation}
\left( \Ind_{\rho} \alpha, \beta \right) _G
= \left( \alpha, \Res_{\rho} \beta \right) _H
\label{eq.exe.Indrho.g.herm}
\end{equation}
and
\begin{equation}
\left\langle \Ind_{\rho} \alpha, \beta \right\rangle _G
= \left\langle \alpha, \Res_{\rho} \beta \right\rangle _H .
\label{eq.exe.Indrho.g.bil}
\end{equation}
(See Remark~\ref{rmk.Resrho} for the definition of
$\Res_{\rho} \beta$.)

\item[(h)] We have
$\Hom_{\CC G} \left( \Ind_{\rho} U, V \right)
\cong \Hom_{\CC H} \left( U, \Res_{\rho} V \right)$ for every
$\CC H$-module $U$ and every $\CC G$-module $V$. (See
Remark~\ref{rmk.Resrho} for the definition of $\Res_{\rho} V$.)

\item[(i)] Similarly to how we made $\CC G$ into a
$\left( \CC G, \CC H \right)$-bimodule, let us make $\CC G$ into a
$\left( \CC H, \CC G \right)$-bimodule (so the right
$\CC G$-module structure is plain multiplication inside $\CC G$,
whereas the left $\CC H$-module structure is induced by the
$\CC$-algebra homomorphism
$\CC \left[ \rho \right] : \CC H \rightarrow \CC G$). If $U$ is any
$\CC H$-module, then the $\CC G$-module
$\Hom_{\CC H} \left( \CC G, U \right)$ (defined as in
Exercise~\ref{exe.induction.ind-as-hom} using the
$\left( \CC H, \CC G \right)$-bimodule structure on
$\CC G$) is isomorphic to $\Ind_{\rho} U$.

\item[(j)] We have
$\Hom_{\CC G} \left( U, \Ind_{\rho} V \right)
\cong \Hom_{\CC H} \left( \Res_{\rho} U, V \right)$
for every $\CC G$-module $U$ and every $\CC H$-module $V$.
(See Remark~\ref{rmk.Resrho} for the definition of $\Res_{\rho} V$.)
\end{enumerate}

Furthermore:

\begin{enumerate}
\item[(k)] Use the above to prove the formula
\eqref{induced-character-formula}.

\item[(l)] Use the above to prove the formula
\eqref{deflation-character-formula}.
\end{enumerate}

[\textbf{Hint:} Part (b) of this exercise is hard. To solve it, it is useful
to have a way of computing the trace of a linear operator without knowing a
basis of the vector space it is acting on. There is a way to do this using a
``finite dual generating system'', which is a
somewhat less restricted notion than that of a basis\footnote{More precisely:
Let $\mathbb{K}$ be a field, and $V$ be a $\mathbb{K}$-vector space. A
\dfn{finite dual generating system} for $V$ means a triple
$\left( I, \left( a_i \right)_{i \in I},
        \left( f_i \right)_{i \in I} \right)$,
where
\par
\begin{itemize}
\item $I$ is a finite set;
\par
\item $\left( a_i \right)_{i \in I}$ is a family of elements of $V$;
\par
\item $\left( f_i \right)_{i \in I}$ is a family of elements of $V^{\ast}$
(where $V^{\ast}$ means
$\Hom_{\mathbb{K}} \left( V, \mathbb{K} \right)$)
\end{itemize}
\par
such that every $v\in V$ satisfies
$v = \sum_{i \in I} f_i \left( v \right) a_i$.
For example, if $\left( e_j \right) _{j \in J}$ is a
finite basis of the vector space $V$, and if
$\left( e_j^{\ast} \right) _{j \in J}$ is the basis of $V^{\ast}$ dual
to this basis $\left( e_j \right) _{j \in J}$, then
$\left( J, \left( e_j \right) _{j \in J},
\left( e_j^{\ast} \right) _{j \in J} \right)$ is a finite dual
generating system for $V$; however, most finite dual generating
systems are not obtained this way.
\par
The crucial observation is now that if
$\left( I, \left( a_i \right)_{i \in I},
        \left( f_i \right)_{i \in I} \right)$
is a finite dual generating system for a vector space $V$, and if $T$
is an endomorphism of $V$, then
\[
\trace T = \sum_{i \in I} f_i \left( T a_i \right) .
\]
Prove this!}. Try to create a finite dual generating system for
$\Ind_{\rho} U$ from one for $U$ (and from the group
$G$), and then use it to compute $\chi_{\Ind_{\rho} U}$.

The solution of part (i) is a modification of the solution of
Exercise~\ref{exe.induction.ind-as-hom}, but complicated by the fact
that $H$ is no longer (necessarily) a subgroup of $G$. Part (f) can be
solved by similar arguments, or using part (i), or using
Exercise~\ref{exe.fixed-vs-coinvariants}(b).]
\end{exercise}

The result of Exercise~\ref{exe.Indrho}(h) generalizes
\eqref{Frobenius-reciprocity} (because of Exercise~\ref{exe.Indrho}(d)), but
also generalizes \eqref{inflation-reciprocity} (due to
Exercise~\ref{exe.Indrho}(f)).
Similarly, Exercise~\ref{exe.Indrho}(g) generalizes both
\eqref{Frobenius-reciprocity-chars} and \eqref{inflation-reciprocity-chars}.
Similarly, Exercise~\ref{exe.Indrho}(i) generalizes
Exercise~\ref{exe.induction.ind-as-hom}, and
Exercise~\ref{exe.Indrho}(j) generalizes
Exercise~\ref{exe.induction.wrong-way-frob}.

Similarly, Exercise~\ref{exe.ind.tensor} is generalized by the
following exercise:

\begin{exercise}
\label{exe.Indrho.tensor}
Let $G_1$, $G_2$, $H_1$ and $H_2$ be four finite groups.
Let $\rho_1 : H_1 \rightarrow G_1$ and $\rho_2 : H_2 \rightarrow G_2$
be two group homomorphisms. These two homomorphisms clearly induce a
group homomorphism
$\rho_1 \times \rho_2 : H_1 \times H_2 \rightarrow G_1 \times G_2$.
Let $U_1$ be a $\CC H_1$-module, and $U_2$ be a $\CC H_2$-module.
Show that
\[
\Ind_{\rho_1 \times \rho_2} \left( U_1 \otimes U_2 \right)
\cong \left( \Ind_{\rho_1} U_1 \right)
      \otimes \left( \Ind_{\rho_2} U_2 \right)
\]
as $\CC \left[ G_1 \times G_2 \right]$-modules.
\end{exercise}

The $\Ind_{\rho}$ and $\Res_{\rho}$ operators behave ``functorially''
with respect to composition. Here is what this means:

\begin{exercise}
\label{exe.Indrho.transitivity}
Let $G$, $H$ and $I$ be three finite groups.
Let $\rho : H \rightarrow G$ and $\tau : I \rightarrow H$ be two
group homomorphisms.

\begin{enumerate}
\item[(a)] We have
$\Ind_{\rho} \Ind_{\tau} U \cong \Ind_{\rho \circ \tau} U$ for every
$\CC I$-module $U$.

\item[(b)] We have
$\Ind_{\rho} \Ind_{\tau} f = \Ind_{\rho \circ \tau} f$ for every
$f \in R_\CC \left( I \right)$.

\item[(c)] We have
$\Res_{\tau} \Res_{\rho} V = \Res_{\rho \circ \tau} V$ for every
$\CC G$-module $V$.

\item[(d)] We have
$\Res_{\tau} \Res_{\rho} f = \Res_{\rho \circ \tau} f$ for every
$f \in R_\CC \left( G \right)$.
\end{enumerate}
\end{exercise}

Exercise~\ref{exe.Indrho.transitivity}(a), of course, generalizes
Exercise~\ref{exe.induction.transitivity}.

% [DG][v52] Added the above discussion (beginning with the two
% paragraphs above Remark~\ref{rmk.Resrho} and ending before
% this very comment) of $\rho$-restriction and $\rho$-induction.
% We don't really need these operations (though one solution
% would have been simpler if I had had them at my disposal back
% when I was writing it), but they just feel morally relevant
% (particularly because they link induction with fixed point
% construction, and restriction with inflation, thus reducing
% the number of fundamental operations to two). Also I haven't
% seen my solution to Exercise~\ref{exe.Indrho}(b) anywhere
% and I thought it deserved writing up.

\subsubsection{Semidirect products}
Recall that a
\emph{semidirect product}\index{semidirect product of groups}
is a group $G \ltimes K$
having two subgroups $G,K$ with
\begin{enumerate}
\item[$\bullet$]
$K \triangleleft (G \ltimes K)$ is a normal subgroup,
\item[$\bullet$]
$G \ltimes K=GK=KG$, and
\item[$\bullet$]
$G \cap K=\{e\}$.
\end{enumerate}
In this setting one has two interesting adjoint constructions,
applied in Section~\ref{wreath-product-section}.

\begin{proposition}
\label{semidirect-product-adjoint-functors}
Fix a $\CC[G \ltimes K]$-module $V$.
\begin{enumerate}
\item[(i)]
For any $\CC G$-module $U$, one has $\CC[G \ltimes K]$-module
structure
\[
\Phi(U):=U \otimes V,
\]
determined via
\begin{align*}
k(u \otimes v)&= u \otimes k(v),\\
g(u \otimes v)&= g(u) \otimes g(v).\\
\end{align*}
\item[(ii)]
For any $\CC[G \ltimes K]$-module $W$, one has $\CC G$-module
structure
\[
\Psi(W):=
\Hom_{\CC K}(\Res^{G \ltimes K}_K V, \Res^{G \ltimes K}_K W),
\]
determined via
$
g(\varphi)=g \circ \varphi \circ g^{-1}.
$
\item[(iii)]
The maps
\[
\CC G-\mods
\underset{\Psi}{\overset{\Phi}{\rightleftharpoons}}
\CC[G \ltimes K]-\mods
\]
are adjoint in the sense that one has an isomorphism
\[
\begin{array}{ccc}
\Hom_{\CC G}(U,\Psi(W))
  &\longrightarrow &
      \Hom_{\CC[G\ltimes K]}(\Phi(U),W) \\
\Vert &  & \Vert\\
\Hom_{\CC G}(U,\Hom_{\CC K}(\Res^{G \ltimes K}_K V, \Res^{G \ltimes K}_K W)) & &
\Hom_{\CC[G\ltimes K]}(U \otimes V,W) , \\
& & \\
\varphi & \longmapsto & \overline{\varphi}(u \otimes v):=\varphi(u)(v) .
\end{array}
\]
\item[(iv)]
One has a $\CC G$-module isomorphism
\[
(\Psi \circ \Phi)(U) \cong U \otimes \End_{\CC K}(\Res^{G \ltimes K}_K  V).
\]
In particular, if $\Res^{G \ltimes K}_K  V$
is a simple $\CC K$-module, then $(\Psi \circ \Phi)(U) \cong U$.

\end{enumerate}
\end{proposition}
\begin{proof}
These are mostly straightforward exercises in the definitions.
To check assertion (iv), for example, note that $K$ acts only
in the right tensor factor in
$\Res^{G \ltimes K}_K(U \otimes V)$,
and hence as $\CC G$-modules one has
\begin{align*}
(\Psi \circ \Phi)(U)
        &= \Hom_{\CC K}(\Res^{G \ltimes K}_K V, \,\,
             \Res^{G \ltimes K}_K(U \otimes V)) \\
        &= \Hom_{\CC K}(\Res^{G \ltimes K}_K V, \,\,
             U \otimes \Res^{G \ltimes K}_K V) \\
        &= U \otimes \Hom_{\CC K}(\Res^{G \ltimes K}_K V, \,\,
                                  \Res^{G \ltimes K}_K V) \\
        &= U \otimes \End_{\CC K}(\Res^{G \ltimes K}_K V) .
\end{align*}
\end{proof}

\subsection{Three towers of groups}
\label{three-towers-section}

Here we consider three \idx{towers of groups}\index{$G_*$}\index{three towers}
\[
G_*= (G_0 < G_1 < G_2 < G_3 < \cdots)
\]
where either
\begin{enumerate}
\item[$\bullet$] $G_n=\Symm_n$, the \emph{symmetric group}\footnote{The symmetric group $\Symm_0$ is the group of all permutations of the empty set $\left\{1,2,\ldots,0\right\}=\varnothing$. It is a trivial group. (Note that $\Symm_1$ is also a trivial group.)}, or
\item[$\bullet$] $G_n=\Symm_n[\Gamma]$, the \emph{wreath product} of the
symmetric group with some arbitrary finite group $\Gamma$, or
\item[$\bullet$] $G_n=GL_n(\FF_q)$, the finite general linear group\footnote{The group $GL_0(\FF_q)$ is a trivial group, consisting of the empty $0\times 0$ matrix.}.
\end{enumerate}
Here the wreath product $\Symm_n[\Gamma]$ can be thought of informally as the group of \emph{monomial}\index{monomial matrix} $n \times n$ matrices whose nonzero entries lie in $\Gamma$, that is,
$n \times n$ matrices having exactly one nonzero entry in each row and column, and that entry is an element of $\Gamma$.  E.g.
\[
\left[
\begin{matrix}
0 & g_2 & 0 \\
g_1 & 0 & 0 \\
0   & 0 & g_3
\end{matrix}
\right]
\left[
\begin{matrix}
0   & 0   & g_6 \\
0   & g_5 & 0 \\
g_4   & 0 & 0
\end{matrix}
\right]
=
\left[
\begin{matrix}
0        & g_2g_5   & 0 \\
0        & 0        & g_1 g_6 \\
g_3g_4   & 0 & 0
\end{matrix}
\right].
\]
More formally, $\Symm_n[\Gamma]$ is
the semidirect product $\Symm_n \ltimes \Gamma^n$
in which $\Symm_n$ acts on $\Gamma^n$ via
$\sigma(\gamma_1,\ldots,\gamma_n)=
(\gamma_{\sigma^{-1}(1)},\ldots,\gamma_{\sigma^{-1}(n)}).
$

For each of the three towers $G_*$,
there are embeddings $G_i \times G_j \hookrightarrow G_{i+j}$
and we introduce maps \dfn{$\ind^{i+j}_{i,j}$}
taking $\CC[G_i \times G_j]$-modules to
$\CC G_{i+j}$-modules, as well as maps
\dfn{$\res^{i+j}_{i,j}$} carrying modules
in the reverse direction which are adjoint:
\begin{equation}
\label{ind-res-adjointness}
\Hom_{\CC G_{i+j}}( \ind^{i+j}_{i,j} U , V ) = \Hom_{\CC[G_i \times G_j]}( U, \res^{i+j}_{i,j} V ).
\end{equation}

\begin{definition}
\label{ind-res-defn}
For $G_n=\Symm_n$, one embeds $\Symm_i \times \Symm_j$ into $\Symm_{i+j}$
as the permutations that permute $\{1,2,\ldots,i\}$ and $\{i+1,i+2,\ldots,i+j\}$
separately.  Here one defines
\begin{align*}
\ind_{i,j}^{i+j}&:=\Ind_{\Symm_i \times \Symm_j}^{\Symm_{i+j}},\\
\res_{i,j}^{i+j}&:=\Res_{\Symm_i \times \Symm_j}^{\Symm_{i+j}}.
\end{align*}
For $G_n=\Symm_n[\Gamma]$, similarly embed
$\Symm_i[\Gamma] \times \Symm_j[\Gamma]$ into
$\Symm_{i+j}[\Gamma]$ as block monomial matrices whose two diagonal blocks
have sizes $i,j$ respectively, and define
\begin{align*}
\ind_{i,j}^{i+j}&:=\Ind_{\Symm_i[\Gamma] \times \Symm_j[\Gamma]}^{\Symm_{i+j}[\Gamma]},\\
\res_{i,j}^{i+j}&:=\Res_{\Symm_i[\Gamma] \times \Symm_j[\Gamma]}^{\Symm_{i+j}[\Gamma]}.
\end{align*}
For $G_n=GL_n(\FF_q)$, which we will denote just $GL_n$,
similarly embed $GL_i \times GL_j$ into
$GL_{i+j}$ as block diagonal matrices whose two diagonal blocks
have sizes $i,j$ respectively.  However, one also introduces
as an intermediate the
\emph{parabolic subgroup}\index{parabolic subgroup of $GL_n$}
$P_{i,j}$ consisting
of the block upper-triangular matrices of the form
\[
\left[
\begin{matrix}
g_i & \ell\\
0   & g_j
\end{matrix}
\right]
\]
where $g_i, g_j$ lie in $GL_i,GL_j$, respectively, and
$\ell$ in $\FF_q^{i \times j}$ is arbitrary.  One has
a quotient map $P_{i,j} \rightarrow GL_i \times GL_j$ whose kernel $K_{i,j}$
is the set of matrices of the form
\[
\left[
\begin{matrix}
I_i & \ell\\
0   & I_j
\end{matrix}
\right]
\]
with $\ell$ again arbitrary.
Here one defines
\begin{align*}
\ind_{i,j}^{i+j}&:=\Ind_{P_{i,j}}^{GL_{i+j}} \Inf^{P_{i,j}}_{GL_i \times GL_j},\\
\res_{i,j}^{i+j}&:=\left(\Res^{GL_{i+j}}_{P_{i,j}}(-)\right)^{K_{i,j}}.
\end{align*}
\end{definition}

In the case $G_n = GL_n$, the operation $\ind_{i,j}^{i+j}$ is sometimes
called \dfn{parabolic induction} or \dfn{Harish-Chandra induction}.
The operation $\res_{i,j}^{i+j}$ is essentially just the $K_{i,j}$-fixed point construction $V \mapsto V^{K_{i,j}}$.
However writing it as the above two-step composite makes it more obvious,
(via \eqref{Frobenius-reciprocity} and \eqref{inflation-reciprocity})
that $\res_{i,j}^{i+j}$ is again adjoint to $\ind_{i,j}^{i+j}$.

\begin{definition}
For each of the three towers $G_*$, define a graded $\ZZ$-module\index{$A(G_*)$}
\[
A:=A(G_*)=\bigoplus_{n \geq 0} R(G_n)
\]
with a bilinear form $(\cdot,\cdot)_A$ whose restriction to $A_n:=R(G_n)$
is the usual form $(\cdot,\cdot)_{G_n}$, and such that
$\Sigma := \bigsqcup_{n \geq 0} \Irr(G_n)$ gives an orthonormal
$\ZZ$-basis.  Notice that $A_0=\ZZ$ has
its basis element $\one$ equal to the unique irreducible character
of the trivial group $G_0$.
\end{definition}

% [DG] I have removed the $\{\chi\}$ from this definition (it used to
% say $\Sigma=\{\chi\}=\bigsqcup_{n \geq 0} \Irr(G_n)$) as I didn't find
% it very helpful (it would only make the impression that $\Sigma$ is
% a one-element set). Sorry if I missed its purpose!
%
% Also, this is one of the places I've replaced "so that" by "such
% that", as I do believe there is a difference between the two idioms.
% (I left wordings like "reindex so that $i$ becomes $1$" intact;
% I only changed cases like "there exists an [object] so that
% [condition]".)

Bearing in mind that $A_n=R(G_n)$ and
\[
A_i \otimes A_j =
R(G_i) \otimes R(G_j) \cong
R(G_i \times G_j) ,
\]
one then has candidates for product and coproduct
defined by
\[
\begin{array}{rrcl}
m:=\ind_{i,j}^{i+j}:&A_i \otimes A_j & \longrightarrow & A_{i+j} , \\
\Delta:=\bigoplus_{i+j=n} \res_{i,j}^{i+j}: &A_n & \longrightarrow & \bigoplus_{i+j=n} A_i \otimes A_j. \\
\end{array}
\]
The coassociativity of $\Delta$ is an easy consequence of transitivity of the
constructions of restriction and fixed points\footnote{More precisely,
using this transitivity, it is easily reduced to proving that
$K_{i+j, k} \cdot \left(K_{i,j} \times \left\{I_k\right\}\right)
= K_{i, j+k} \cdot \left(\left\{I_i\right\} \times K_{j,k}\right)$
(an equality between subgroups of $GL_{i+j+k}$) for any three
nonnegative integers $i, j, k$. But this equality
can be proven by realizing that both of its sides equal the
set of all block matrices of the form $\left(
\begin{array}{ccc}
I_i & \ell & \ell' \\
0 & I_j & \ell'' \\
0 & 0 & I_k
\end{array} \right)$
with $\ell$, $\ell'$ and $\ell''$ being matrices of sizes $i \times j$,
$i \times k$ and $j \times k$, respectively.}.
We could derive the associativity of $m$ from the transitivity of
induction and inflation, but this would be more
complicated\footnote{See Exercise~\ref{exe.towers.indres}(c) for
such a derivation.}; we will instead prove it differently.

We first show that the maps $m$ and $\Delta$ are adjoint with respect
to the forms $\left(\cdot,\cdot\right)_A$ and $\left(\cdot,\cdot\right)_{A\otimes A}$.
In fact, if $U$, $V$, $W$
are modules over $\CC G_i$, $\CC G_j$, $\CC G_{i+j}$, respectively, then
we can write the $\CC [G_i\times G_j]$-module $\res_{i,j}^{i+j} W$ as a
direct sum $\bigoplus_k X_k \otimes Y_k$ with $X_k$ being $\CC G_i$-modules
and $Y_k$ being $\CC G_j$-modules; we then have
\begin{equation}
\label{ind-res-adjointness-pf1}
\res_{i,j}^{i+j} \chi_W = \sum_k \chi_{X_k} \otimes \chi_{Y_k}
\end{equation}
and
\begin{align*}
\left(m\left(\chi_U \otimes \chi_V\right), \chi_W\right)_A
&= \left(\ind_{i,j}^{i+j} \left(\chi_{U\otimes V}\right), \chi_W\right)_A
= \left(\ind_{i,j}^{i+j} \left(\chi_{U\otimes V}\right), \chi_W\right)_{G_{i+j}} \\
&= \left(\chi_{U\otimes V}, \res_{i,j}^{i+j} \chi_W\right)_{G_i\times G_j}
= \left(\chi_{U\otimes V}, \sum_k \chi_{X_k} \otimes \chi_{Y_k}\right)_{G_i\times G_j} \\
&= \sum_k \left(\chi_{U\otimes V}, \chi_{X_k\otimes Y_k}\right)_{G_i\times G_j}
= \sum_k \left(\chi_U, \chi_{X_k}\right)_{G_i} \left(\chi_V, \chi_{Y_k}\right)_{G_j}
\end{align*}
(the third equality sign follows by taking dimensions in
\eqref{ind-res-adjointness} and recalling \eqref{character-scalprod};
the fourth equality sign follows from \eqref{ind-res-adjointness-pf1}; the
sixth one follows from \eqref{character-tensorprod}) and
\begin{align*}
\left(\chi_U \otimes \chi_V, \Delta\left(\chi_W\right)\right)_{A\otimes A}
&= \left(\chi_U \otimes \chi_V, \res_{i,j}^{i+j} \chi_W\right)_{A\otimes A}
= \left(\chi_U \otimes \chi_V, \sum_k \chi_{X_k} \otimes \chi_{Y_k}\right)_{A\otimes A} \\
&= \sum_k \left(\chi_U, \chi_{X_k}\right)_A \left(\chi_V, \chi_{Y_k}\right)_A
= \sum_k \left(\chi_U, \chi_{X_k}\right)_{G_i} \left(\chi_V, \chi_{Y_k}\right)_{G_j}
\end{align*}
(the first equality sign follows by removing all terms in $\Delta\left(\chi_W\right)$
whose scalar product with $\chi_U \otimes \chi_V$ vanishes for reasons of gradedness;
the second equality sign follows from
\eqref{ind-res-adjointness-pf1}), which in comparison yield
$\left(m\left(\chi_U \otimes \chi_V\right), \chi_W\right)_A
= \left(\chi_U \otimes \chi_V, \Delta\left(\chi_W\right)\right)_{A\otimes A}$,
thus showing that $m$ and $\Delta$ are adjoint maps. Therefore, $m$ is associative
(since $\Delta$ is coassociative).

% [DG] I didn't expect the details of the proof of adjointness to take this
% long. Maybe I am confusing myself here.
%
% The argument for associativity of $m$ using adjointness and coassociativity
% is an alternative to the 3-pages proof of the associativity of $m$ that I
% did in reiner-errata12.pdf. Of course, it doesn't show any functorial
% isomorphisms.

Endowing $A=\bigoplus_{n \geq 0} R(G_n)$ with the
obvious unit and counit maps, it thus becomes a graded, finite-type
$\ZZ$-algebra and $\ZZ$-coalgebra.

The next section addresses the issue of why
they form a bialgebra.  However, assuming this for the moment,
it should be clear that each of these algebras $A$ is a PSH having
$\Sigma=\bigsqcup_{n \geq 0} \Irr(G_n)$ as its PSH-basis.
$\Sigma$ is self-dual because $m, \Delta$ are defined by adjoint maps,
and it is positive because $m, \Delta$ take irreducible representations
to genuine representations not just virtual ones, and
hence have characters which are nonnegative sums of irreducible characters.

\begin{exercise}
\label{exe.Ind-ass}
Let $i$, $j$ and $k$ be three nonnegative integers. Let $U$ be a
$\CC \Symm_i$-module, let $V$ be a $\CC \Symm_j$-module, and let $W$ be a
$\CC \Symm_k$-module. Show that there are canonical
$\CC \left[\Symm_i \times \Symm_j \times \Symm_k\right]$-module
isomorphisms
\begin{align*}
\Ind^{\Symm_{i+j+k}}_{\Symm_{i+j} \times \Symm_k}\left(\Ind^{\Symm_{i+j}}_{\Symm_i \times \Symm_j}\left(U\otimes V\right)\otimes W\right)
&\cong \Ind^{\Symm_{i+j+k}}_{\Symm_i \times \Symm_j \times \Symm_k}\left(U\otimes V\otimes W\right) \\
&\cong \Ind^{\Symm_{i+j+k}}_{\Symm_i \times \Symm_{j+k}}\left(U \otimes \Ind^{\Symm_{j+k}}_{\Symm_j \times \Symm_k}\left(V\otimes W\right)\right).
\end{align*}

(Similar statements hold for the other two towers of groups and their
respective $\ind$ functors, although the one for the $GL_*$ tower is
harder to prove. See Exercise~\ref{exe.towers.indres}(a) for a
more general result.)
\end{exercise}

% [DG][v14] I recycled part of reiner-errata12.pdf in the
% above solution.

\subsection{Bialgebra and double cosets}
\label{double-coset-section}

To show that the algebra and coalgebras $A=A(G_*)$ are bialgebras,
the central issue is checking the pentagonal diagram in \eqref{bialgebra-diagrams}, that is, as maps $A \otimes A \rightarrow A \otimes A$, one has
\begin{equation}
\label{pentagonal-diagram-written-out}
\Delta \circ m = (m \otimes m) \circ (\id \otimes T \otimes \id) \circ (\Delta \otimes \Delta).
\end{equation}

In checking this, it is convenient to have a lighter notation for
various subgroups of the groups $G_n$ corresponding to compositions $\alpha$.

\begin{definition} \phantomsection
\label{def.composition-almost-composition}
\begin{itemize}

\item[(a)]
An \dfn{almost-composition} is a (finite) tuple
$\alpha = (\alpha_1,\alpha_2,\ldots,\alpha_\ell)$
of nonnegative integers. Its
\emph{length}\index{length of an almost-composition}
is defined to be $\ell$ and denoted by \dfn{$\ell(\alpha)$};
its \emph{size}\index{size of an almost-composition} is defined to be
$\alpha_1 + \alpha_2 + \cdots + \alpha_\ell$ and denoted by
$\left|\alpha\right|$;
its \emph{parts}\index{parts of an almost-composition}
are its entries $\alpha_1,\alpha_2,\ldots,\alpha_\ell$.
The almost-compositions of size $n$ are called the
\emph{almost-compositions of $n$}\index{almost-composition of $n$}.

\item[(b)] A \dfn{composition} is a finite tuple of positive
integers. Of course, any composition is an almost-composition,
and so all notions defined for almost-compositions (like
size and length) make sense for compositions.%
\index{length of a composition}\index{size of a composition}\index{parts of a composition}

Note that any partition of $n$ (written without trailing zeroes) is a
composition of $n$. We write $\varnothing$ (and sometimes, sloppily,
$(0)$, when there is no danger of mistaking it for the
almost-composition $(0)$) for the empty composition $()$.
\end{itemize}
\end{definition}

% [DG][v28] Introduced almost-compositions, and replaced the use of
% compositions by almost-compositions below. Indeed, this is needed in
% the proof of the bialgebra axiom (though one could avoid dealing with
% that, but that would be less organic).
% 
% I don't call them weak compositions, because the name is already taken
% for infinite sequences.

% [DG] This doesn't look like a great place to introduce compositions,
% but it's the best one I found, as they're used right below for the
% first time. In Chapter 2 it would probably be distracting.

\begin{definition}
Given an almost-composition $\alpha=(\alpha_1,\ldots,\alpha_\ell)$ of $n$,
define a subgroup
\[
G_\alpha \cong G_{\alpha_1} \times \cdots \times G_{\alpha_\ell}
< G_n
\]
via the block-diagonal embedding with diagonal blocks of
sizes $(\alpha_1,\ldots,\alpha_\ell)$.
This $G_\alpha$ is called a \dfn{Young subgroup} $\Symm_\alpha$ when
$G_n =\Symm_n$,
and a \dfn{Levi subgroup} when $G_n=GL_n$.
In the case when $G_n = \Symm_n[\Gamma]$, we also denote $G_\alpha$ by
$\Symm_\alpha[\Gamma]$.
In the case where $G_n=GL_n$, also define the
\emph{parabolic subgroup}\index{parabolic subgroup of $GL_n$}
$
P_\alpha
$
to be the subgroup of $G_n$ consisting of block-upper triangular
matrices whose diagonal blocks have sizes $(\alpha_1,\ldots,\alpha_\ell)$,
and let
$
K_\alpha
$
be the kernel of the obvious surjection $P_\alpha \rightarrow G_\alpha$ which sends
a block upper-triangular matrix to the tuple of its diagonal blocks
whose sizes are $\alpha_1,\alpha_2,\ldots,\alpha_\ell$.  Notice that
$P_{(i,j)} = P_{i,j}$ for any $i$ and $j$ with $i+j = n$; similarly,
$K_{(i,j)} = K_{i,j}$ for any $i$ and $j$ with $i+j = n$. We will also
abbreviate $G_{(i,j)} = G_i \times G_j$ by $G_{i,j}$.
%Also define operators $\ind^n_\alpha,
%\res^n_\alpha$ for the three towers analogous to those defined in
%Definition~\ref{ind-res-defn} when $\alpha=(i,j)$ has only two parts.

When $\left(\alpha_1, \alpha_2, \ldots, \alpha_\ell\right)$ is an
almost-composition, we abbreviate
$G_{\left(\alpha_1, \alpha_2, \ldots, \alpha_\ell\right)}$ by
$G_{\alpha_1, \alpha_2, \ldots, \alpha_\ell}$ (and similarly for the
$P$'s).
\end{definition}

% [DG][v28] Commented out the three lines above because these notations
% are not yet needed (and when they are needed, I'll define them in
% more detail).

\begin{definition}
\label{general-restriction}
Let $K$ and $H$ be two groups, $\tau : K \to H$ a group homomorphism,
and $U$ a $\CC H$-module. Then, $U^\tau$ is defined as the
$\CC K$-module with ground space $U$ and action given by
$ k \cdot u = \tau(k) \cdot u $ for all $k \in K$ and $u \in U$.
\ \ \ \ \footnote{We have already met this $\CC K$-module $U^\tau$
in Remark~\ref{rmk.Resrho}, where it was called $\Res_\tau U$.}
This very simple construction generalizes the definition of $U^g$
for an element $g \in G$, where $G$ is a group containing $H$ as
a subgroup; in fact, in this situation we have $U^g = U^\tau$,
where $K = {}^g H$ and $\tau : K \to H$ is the map $k \mapsto g^{-1}kg$.
\end{definition}

Using homogeneity, checking the bialgebra condition
\eqref{pentagonal-diagram-written-out}
in the homogeneous component $(A \otimes A)_n$
amounts to the following:
for each pair of representations $U_1,U_2$ of
$G_{r_1}, G_{r_2}$ with $r_1+r_2=n$, and
for each $(c_1,c_2)$ with $c_1+c_2=n$, one
must verify that
\begin{equation}
\label{bialgebra-as-Mackey}
\begin{aligned}
&\res_{c_1,c_2}^n\left(  \ind_{r_1,r_2}^n \left( U_1 \otimes U_2 \right) \right) \\
&\quad \cong \bigoplus_A
  \left( \ind_{a_{11},a_{21}}^{c_1} \otimes \ind^{c_2}_{a_{12},a_{22}}\right)
   \left( \left( \res^{r_1}_{a_{11},a_{12}} U_1 \otimes \res^{r_2}_{a_{21},a_{22}} U_2 \right)^{\tau_A^{-1}} \right)
\end{aligned}
\end{equation}
where the direct sum is over all matrices
$
A=\left[
\begin{matrix}
a_{11} & a_{12} \\
a_{21} & a_{22}
\end{matrix}
\right]
$
in $\NN^{2 \times 2}$ with row sums $(r_1,r_2)$ and column sums $(c_1,c_2)$,
and where $\tau_A$ is the obvious isomorphism
between the subgroups
\begin{equation}
\label{two-by-two-matrix-Mackey-subgroups}
\begin{aligned}
 G_{a_{11},a_{12},a_{21},a_{22}}
&\quad \left( < G_{r_1, r_2} \right) \qquad \text{and}\\
 G_{a_{11},a_{21},a_{12},a_{22}}
&\quad \left( < G_{c_1,c_2} \right)\\
\end{aligned}
\end{equation}
(we are using the inverse $\tau_A^{-1}$ of this isomorphism $\tau_A$
to identify modules for the first subgroup with modules
for the second subgroup, according to Definition~\ref{general-restriction}).

% [DG][v3] I replaced $\tau_A$ by $\tau_A^{-1}$ in the formula so as to match
% the definition \eqref{general-restriction} I introduced. (The definition
% defines $U^\tau$ as a pullback rather than a pushforward because the
% pushforward is defined only for isomorphisms $\tau$, while the pullback
% makes sense for any morphism.)

As one might guess, \eqref{bialgebra-as-Mackey} comes from the Mackey formula
(Theorem~\ref{Mackey-formula}),
once one identifies the appropriate double coset representatives.
This is just as easy to do in a slightly more general setting.

\begin{definition}
\label{def.compositions.intervaldec}
Given almost-compositions $\alpha,\beta$ of $n$ having
lengths $\ell,m$
and a matrix $A$ in $\NN^{\ell \times m}$ with
row sums $\alpha$ and column sums $\beta$, define a permutation $w_A$
in $\Symm_n$ as follows.  Disjointly decompose
$[n]=\{1,2,\ldots,n\}$ into consecutive intervals of numbers
\[
[n]=I_1 \sqcup \cdots \sqcup I_\ell \qquad \text{such that } |I_i|=\alpha_i
\]
(so the smallest $\alpha_1$ elements of $[n]$ go into $I_1$,
the next-smallest $\alpha_2$ elements of $[n]$ go into $I_2$,
and so on).
Likewise, disjointly decompose $[n]$ into consecutive intervals of numbers
\[
[n]=J_1 \sqcup \cdots \sqcup J_m \qquad \text{such that } |J_j|=\beta_j.
\]
For every $j\in [m]$, disjointly decompose $J_{j}$
into consecutive intervals of numbers $J_{j}=J_{j,1}\sqcup J_{j,2}\sqcup
\cdots\sqcup J_{j,\ell}$ such that every $i\in [\ell]$
satisfies $\left\vert J_{j,i}\right\vert =a_{ij}$.
For every $i\in [\ell]$, disjointly decompose $I_{i}$
into consecutive intervals of numbers $I_{i}=I_{i,1}\sqcup I_{i,2}\sqcup
\cdots\sqcup I_{i,m}$ such that every $j\in [m]$
satisfies $\left\vert I_{i,j}\right\vert =a_{ij}$.
Now, for every $i\in [\ell]$ and $j\in [m]$, let
$\pi_{i,j}$ be the increasing bijection
from $J_{j,i}$ to $I_{i,j}$ (this is well-defined since these two sets both
have cardinality $a_{ij}$). The disjoint union of these bijections $\pi_{i,j}$
over all $i$ and $j$ is a bijection $[n]\to [n]$ (since the disjoint union
of the sets $J_{j,i}$ over all $i$ and $j$ is $[n]$, and so is the disjoint
union of the sets $I_{i,j}$), that is, a permutation of $[n]$; this
permutation is what we call $w_A$.
%[Old version:]
%Then let $w_A$ map a consecutive string of $a_{ij}$ of the letters of
%$J_j$ into a consecutive string of $a_{ij}$ of the letters of $I_i$,
%in an increasing fashion.
\end{definition}

% [DG] Shortened my suggestion and still clumsy :(

% [DG][v80] Detailed the "disjointly decompose" a bit.

\begin{example}
Taking
$n=9$ and $\alpha=(4,5), \beta=(3,4,2)$, one has
\[
\begin{array}{ccc}
I_1=\{1,2,3,4\},&I_2=\{5,6,7,8,9\},& \\
J_1=\{1,2,3\}, &J_2=\{4,5,6,7\},&J_3=\{8,9\}.
\end{array}
\]
Then one possible matrix $A$ having row and column sums $\alpha, \beta$ is
$
A=\left[
\begin{matrix}
2 & 2 & 0\\
1 & 2 & 2
\end{matrix}
\right],
$
and its associated permutation $w_A$ written in two-line notation is
\[
\left(
\begin{matrix}
1 & 2 & 3 &|& 4 & 5 & 6 & 7&|& 8 & 9 \\
\underline{1} & \underline{2} & \underline{\underline{5}}&|
&\underline{3} & \underline{4} & \underline{\underline{6}}
&\underline{\underline{7}} &|
&\underline{\underline{8}}
&\underline{\underline{9}}
\end{matrix}
\right)
\]
with vertical lines dividing the sets $J_j$ on top, and
with elements of $I_i$ underlined $i$ times on the bottom.
\end{example}

\begin{remark}
Given almost-compositions $\alpha$ and $\beta$ of $n$
having lengths $\ell$ and $m$,
and a permutation $w \in \Symm_n$. It is easy to see that there exists a
matrix $A \in \NN^{\ell\times m}$ satisfying $w_A = w$ if and only if
the restriction of $w$ to each $J_j$ and the restriction of $w^{-1}$ to
each $I_i$ are increasing. In this case, the matrix $A$ is determined by
$a_{ij} = \left| w(J_j) \cap I_i \right|$.
\end{remark}

% [DG] I added this remark since it is (tacitly) used in the proof of the
% next proposition.

Among our three towers $G_*$ of groups, the symmetric group tower
($G_n = \Symm_n$) is the simplest one. We will now see that it also embeds
into the two others, in the sense that $\Symm_n$ embeds into $\Symm_n[\Gamma]$
for every $\Gamma$ and into $GL_n(\FF_q)$ for every $q$.

First, for every $n\in\NN$ and any group $\Gamma$, we embed
the group $\Symm_n$ into $\Symm_n[\Gamma]$ by means of the canonical embedding
$\Symm_n\rightarrow\Symm_n\ltimes\Gamma^{n}=\Symm_n[\Gamma] $.
If we regard elements of $\Symm_n[\Gamma] $ as $n\times n$ monomial matrices
with nonzero entries in $\Gamma$, then this boils down to identifying every
$\pi\in\Symm_n$ with the permutation matrix of $\pi$ (in which the $1$'s are
read as the neutral element of $\Gamma$). If $\alpha$ is an
almost-composition of $n$,
then this embedding $\Symm_n \rightarrow \Symm_n\left[  \Gamma\right]$
makes the subgroup $\Symm_{\alpha}$ of $\Symm_n$ become a subgroup of
$\Symm_n\left[ \Gamma\right]$, more precisely a subgroup of
$\Symm_{\alpha} \left[ \Gamma\right] < \Symm_n \left[ \Gamma\right]$.

For every $n\in\NN$ and every $q$, we embed the group $\Symm_n$ into
$GL_n \left(  \FF_{q}\right)  $ by identifying every permutation
$\pi \in \Symm_n$ with its permutation matrix in $GL_n \left(
\FF_{q}\right)  $. If $\alpha$ is an
almost-composition of $n$, then this embedding
makes the subgroup $\Symm_{\alpha}$ of $\Symm_n$ become a subgroup of
$GL_n \left(  \FF_{q}\right) $. If we let $G_n = GL_n \left(
\FF_{q}\right)  $, then $\Symm_{\alpha} < G_{\alpha} < P_{\alpha}$.

The embeddings we have just defined commute with the group embeddings
$G_n < G_{n+1}$ on both sides.

% [DG] Shortened my definitions; still too tedious...

\begin{proposition}
\label{double-coset-prop}
The permutations $\{w_A\}$,
as $A$ runs over all matrices in $\NN^{\ell\times m}$
having row sums $\alpha$ and column sums $\beta$,
give
\begin{itemize}
\item[(a)] a system of double coset representatives
           for $\Symm_\alpha \backslash \Symm_n / \Symm_\beta$;
\item[(b)] a system of double coset representatives
           for $\Symm_\alpha[\Gamma] \backslash \Symm_n[\Gamma] / \Symm_\beta[\Gamma]$;
\item[(b)] a system of double coset representatives
           for $P_\alpha \backslash GL_n / P_\beta$.
\end{itemize}
\end{proposition}

% [DG][v75] Broke this proposition up into (a), (b), (c).

\begin{proof}
(a) We give an algorithm to show
that every double coset $\Symm_\alpha w \Symm_\beta$
contains some $w_A$.  Start by altering
$w$ within its coset $w\Symm_\beta$, that is, by permuting the
\emph{positions} within each set $J_j$,
to obtain a representative $w'$ for $w\Symm_\beta$ in
which each set $w'(J_j)$ appears in increasing order in the second line of the
two-line notation for $w'$.  Then alter $w'$ within its coset
$\Symm_\alpha w'$, that is, by permuting the \emph{values} within each set
$I_i$, to obtain a representative
$w_A$ having the elements of each set $I_i$ appearing
in increasing order in the second line;  because the values within
each set $I_i$ are consecutive, this alteration will not ruin the property
that one had each set $w'(J_j)$ appearing in increasing order.
For example, one might have
\begin{align*}
w&=\left(
\begin{matrix}
1 & 2 & 3 &|& 4 & 5 & 6 & 7&|& 8 & 9 \\
\underline{4} & \underline{\underline{8}} & \underline{2}&|
&\underline{\underline{5}} & \underline{3} & \underline{\underline{9}}
&\underline{1} &|
&\underline{\underline{7}}
&\underline{\underline{6}}
\end{matrix}
\right) ,\\
w'&=\left(
\begin{matrix}
1 & 2 & 3 &|& 4 & 5 & 6 & 7&|& 8 & 9 \\
\underline{2} & \underline{4} & \underline{\underline{8}}&|
&\underline{1} & \underline{3} & \underline{\underline{5}}
&\underline{\underline{9}} &|
&\underline{\underline{6}}
&\underline{\underline{7}}
\end{matrix}
\right) \in w \Symm_\beta ,\\
w_A&=\left(
\begin{matrix}
1 & 2 & 3 &|& 4 & 5 & 6 & 7&|& 8 & 9 \\
\underline{1} & \underline{2} & \underline{\underline{5}}&|
&\underline{3} & \underline{4} & \underline{\underline{6}}
&\underline{\underline{7}} &|
&\underline{\underline{8}}
&\underline{\underline{9}}
\end{matrix}
\right) \in \Symm_\alpha w' \subset \Symm_\alpha w' \Symm_\beta = \Symm_\alpha w \Symm_\beta .\\
\end{align*}

Next note that
$\Symm_\alpha w_A \Symm_\beta=\Symm_\alpha w_B \Symm_\beta$ implies
$A=B$, since the quantities
\[
a_{i,j}(w):=|w(J_j) \cap I_i|
\]
are easily seen to be constant on double cosets
$\Symm_\alpha w \Symm_\beta$.

(b) Double coset representatives for
$\Symm_\alpha \backslash \Symm_n / \Symm_\beta$ should also provide
double coset representatives for
$\Symm_\alpha[\Gamma] \backslash \Symm_n[\Gamma] / \Symm_\beta[\Gamma]$,
since
\[
\Symm_\alpha[\Gamma]=\Symm_\alpha  \Gamma^n=\Gamma^n \Symm_\alpha.
\]
Thus, part (b) follows from part (a).

(c) In our proof of part (a) above, we showed that
$\Symm_\alpha w_A \Symm_\beta=\Symm_\alpha w_B \Symm_\beta$ implies
$A=B$.
A similar argument shows that
$P_\alpha w_A P_\beta=P_\alpha w_B P_\beta$ implies
$A=B$:  for $g$ in $GL_n$,
the rank $r_{ij}(g)$ of the matrix obtained by restricting $g$ to
rows $I_i \sqcup I_{i+1} \sqcup \cdots \sqcup I_\ell$
and columns $J_1 \sqcup J_2 \sqcup \cdots \sqcup J_j$ is
constant on double cosets $P_\alpha g P_\beta$, and
for a permutation matrix $w$ one can recover $a_{i,j}(w)$ from the formula
\[
a_{i,j}(w) =
r_{i,j}(w)-r_{i,j-1}(w)-r_{i+1,j}(w) + r_{i+1,j-1}(w).
\]

Thus it only remains to show that every double coset $P_\alpha g P_\beta$ contains
some $w_A$.  Since $\Symm_\alpha < P_\alpha$, and we have seen already that
every double coset $\Symm_\alpha w \Symm_\beta$ contains some $w_A$,
it suffices to show that every double coset
$P_\alpha g P_\beta$ contains some permutation $w$.  However, we claim that
this is  already true for the smaller double cosets $BgB$ where $B=P_{1^n}$ is
the \emph{Borel subgroup}\index{Borel subgroup of $GL_n$}
of upper triangular invertible matrices, that is,
one has the usual \emph{Bruhat decomposition}\index{Bruhat decomposition of $GL_n$}
\[
GL_n = \bigsqcup_{w \in \Symm_n} BwB.
\]
To prove this decomposition, we show how to
find a permutation $w$ in each double coset $BgB$.
The freedom to alter $g$ within its coset $gB$ allows one to scale columns and
add scalar multiples of earlier columns to later columns.  We claim that
using such column operations, one can
always find a representative $g'$ for coset $gB$ in which
\begin{enumerate}
\item[$\bullet$]
the bottommost nonzero entry of each column is $1$ (call this entry a \emph{pivot}),
\item[$\bullet$]
the entries to right of each pivot within its row are all $0$, and
\item[$\bullet$]
there is one pivot in each row and each column, so that
their positions are the positions of the $1$'s in some
permutation matrix $w$.
% [DG][v49] Replaced the line below by the two lines above.
% they lie in the positions of some permutation matrix $w$.
\end{enumerate}
In fact, we will see below that $BgB=BwB$ in this case.
The algorithm which produces $g'$ from $g$ is simple:  starting with the
leftmost column, find its bottommost nonzero entry, and scale the column
to make this entry a $1$, creating the pivot in this column.  Now use this
pivot to clear out all entries in its row to its right, using
column operations that subtract multiples of this column from later columns.
Having done this, move on to the next column to the right, and repeat, scaling to create a pivot, and using it to eliminate entries to its right.\footnote{To see that this works, we need to check three facts:
\begin{itemize}
\item[\textbf{(a)}] We will find a nonzero entry in every column during our
algorithm.
\item[\textbf{(b)}] Our column operations preserve
the zeroes lying to the right of already existing pivots.
\item[\textbf{(c)}] Every row contains exactly one pivot at the end
of the algorithm.
\end{itemize}
But fact \textbf{(a)} simply says that our matrix can never have an
all-zero column during the algorithm; this is clear (since the rank of
the matrix remains constant during the algorithm and was $n$ at its
beginning). Fact \textbf{(b)} holds because all our operations either
scale columns (which clearly preserves zero entries) or subtract a
multiple of the column $c$ containing the current pivot from a later
column $d$ (which will preserve every zero lying to the right of an
already existing pivot, because any already existing pivot must lie
in a column $b < c$ and therefore both
columns $c$ and $d$ have zeroes in its row). Fact \textbf{(c)} follows
from noticing that there are $n$ pivots altogether at the end of the
algorithm, but no row can contain two of them (since the entries to
the right of a pivot in its row are $0$).}

% [DG][v34] Added footnote. Replaced "has been scaled to $1$" by
% "is $1$" in the properties of $g'$ because we are describing $g'$
% there, not the algorithm that gives us $g'$.

For example, the typical matrix $g$ lying in the double coset $BwB$ where
\[
w=\left(
\begin{matrix}
1 & 2 & 3 &|& 4 & 5 & 6 & 7&|& 8 & 9 \\
\underline{4} & \underline{\underline{8}} & \underline{2}&|
&\underline{\underline{5}} & \underline{3} & \underline{\underline{9}}
&\underline{1} &|
&\underline{\underline{7}}
&\underline{\underline{6}}
\end{matrix}
\right)\\
\]
from before is one that can be altered within its coset $gB$ to look like this:
\[
g'=\left[
\begin{matrix}
* & * & * & * & * & * & 1 & 0 & 0 \\
* & * & 1 & 0 & 0 & 0 & 0 & 0 & 0 \\
* & * & 0 & * & 1 & 0 & 0 & 0 & 0 \\
1 & 0 & 0 & 0 & 0 & 0 & 0 & 0 & 0 \\
0 & * & 0 & 1 & 0 & 0 & 0 & 0 & 0 \\
0 & * & 0 & 0 & 0 & * & 0 & * & 1 \\
0 & * & 0 & 0 & 0 & * & 0 & 1 & 0 \\
0 & 1 & 0 & 0 & 0 & 0 & 0 & 0 & 0 \\
0 & 0 & 0 & 0 & 0 & 1 & 0 & 0 & 0 \\
\end{matrix}
\right] \in gB.
\]
Having found this $g'$ in $gB$, a similar algorithm using left multiplication
by $B$ shows that $w$ lies in
$Bg' \subset Bg'B=BgB$.  This time no scalings are required
to create the pivot entries:  starting with the bottom row, one
uses its pivot to eliminate all the entries above it in the same column
(shown by stars $*$ above)
by adding multiples of the bottom row to higher rows. Then do the same using the pivot in the next-to-bottom row, etc.  The result is the permutation matrix for $w$.
\end{proof}

\begin{remark}
The Bruhat decomposition $GL_n = \bigsqcup_{w \in \Symm_n} BwB$ is
related to the so-called \emph{LPU factorization} -- one of a
myriad of matrix factorizations appearing in linear
algebra.\footnote{Specifically, an \dfn{LPU factorization} of a
matrix $A \in GL_n(\FF)$ (for an arbitrary field $\FF$)
means a way to write $A$ as a product
$A = LPU$ with $L \in GL_n(\FF)$ being lower-triangular,
$U \in GL_n(\FF)$ being upper-triangular, and
$P \in \Symm_n \subset GL_n(\FF)$ being a permutation
matrix. Such a factorization always exists (although it is generally
not unique).
This can be derived from the Bruhat decomposition (see
Exercise~\ref{exe.bruhat.variations}(a) for a proof).
See also \cite{Strang-elim} for related discussion.}
It is actually a fairly general phenomenon, and requires neither the
finiteness of $\FF$, nor the invertibility, nor even the squareness
of the matrices (see Exercise~\ref{exe.bruhat.variations}(b) for an
analogue holding in a more general setup).
\end{remark}

\begin{exercise}
\label{exe.bruhat.variations}
Let $\FF$ be any field.

\begin{enumerate}
\item[(a)]
For any $n \in \NN$ and any $A \in GL_n(\FF)$, prove that there exist
a lower-triangular matrix $L \in GL_n(\FF)$, an upper-triangular
matrix $U \in GL_n(\FF)$ and a permutation matrix
$P \in \Symm_n \subset GL_n(\FF)$ (here, we identify permutations with
the corresponding permutation matrices) such that $A = LPU$.

\item[(b)]
Let $n \in \NN$ and $m \in \NN$. Let $F_{n, m}$ denote the set of all
$n \times m$-matrices $B \in \left\{0, 1\right\}^{n \times m}$ such that
each row of $B$ contains at most one $1$ and each column of $B$
contains at most one $1$. We regard $F_{n, m}$ as a subset of
$\FF^{n \times m}$ by means of regarding $\left\{0, 1\right\}$ as a
subset of $\FF$.

For every $k \in \NN$, we let $B_k$ denote the subgroup of
$GL_k(\FF)$ consisting of all upper-triangular matrices.

Prove that
\[
\FF^{n \times m}
= \bigsqcup_{f \in F_{n, m}} B_n f B_m .
\]
\end{enumerate}
\end{exercise}

% [DG][v49] Added the above remark and exercise (a side effect of TAing
% a Linear Algebra class). I have left out any comments on Bruhat
% decompositions for other algebraic groups, as that is beyond my
% paycheck (actually, probably not -- but beyond my competence).

\begin{corollary}
\label{three-towers-give-bialgebras}
For each of the three towers of groups $G_*$, the product and coproduct
structures on $A=A(G_*)$ endow it with a bialgebra structure, and hence
they form PSH's.
\end{corollary}
\begin{proof}
The first two towers $G_n = \Symm_n$ and $G_n=\Symm_n[\Gamma]$ have
product, coproduct defined by induction, restriction along
embeddings $G_i \times G_j < G_{i+j}$.  Hence the desired bialgebra equality
\eqref{bialgebra-as-Mackey} follows from Mackey's Theorem~\ref{Mackey-formula},
taking
$
G=G_n,
H=G_{(r_1,r_2)},
K=G_{(c_1,c_2)},
U=U_1\otimes U_2
$
with double coset representatives\footnote{Proposition~\ref{double-coset-prop} gives
as a system of double coset representatives for $G_{(c_1,c_2)} \backslash G_n / G_{(r_1,r_2)}$
the elements
\begin{align*}
&\left\{w_A \ : \  A\in\NN^{2\times 2}, \ A\text{ has row sums }
(c_1,c_2)\text{ and column sums }(r_1,r_2)\right\} \\
&= \left\{w_{A^t} \ : \  A\in\NN^{2\times 2}, \ A\text{ has row sums }
(r_1,r_2)\text{ and column sums }(c_1,c_2)\right\}
\end{align*}
where $A^t$ denotes the transpose matrix of $A$.}
\[
\{g_1,\ldots,g_t\} = \left\{w_{A^t} \ : \  A\in\NN^{2\times 2}, \ A\text{ has row sums }
(r_1,r_2)\text{ and column sums }(c_1,c_2)\right\}
\]
and checking for a given double coset
\[
KgH=(G_{c_1,c_2})w_{A^t}(G_{r_1,r_2})
\]
indexed by a matrix $A$ in $\NN^{2 \times 2}$ with row sums $(r_1,r_2)$ and
column sums $(c_1,c_2)$, that the two subgroups appearing on the left in
\eqref{two-by-two-matrix-Mackey-subgroups} are exactly
\begin{align*}
H \cap K^{w_{A^t}} &= G_{r_1,r_2} \cap (G_{c_1,c_2})^{w_{A^t}}, \\
{}^{w_{A^t}} H \cap K &= {}^{w_{A^t}}(G_{r_1,r_2}) \cap G_{c_1,c_2},
\end{align*}
respectively.
One should also apply \eqref{induction-vs-tensor-product} and check that
the isomorphism $\tau_A$ between the two subgroups
in \eqref{two-by-two-matrix-Mackey-subgroups} is the conjugation isomorphism
by $w_{A^t}$ (that is, $\tau_A(g) = w_{A^t} g w_{A^t}^{-1}$ for every
$g \in H \cap K^{w_{A^t}}$).  We leave all of these bookkeeping details to the reader to check.
\footnote{It helps to recognize $w_{A^t}$ as the permutation written in
two-line notation as
\[
\left(
\begin{matrix}
1 & 2 & \ldots & a_{11} &|& a_{11}+1 & a_{11}+2 & \ldots & r_1 &|& r_1+1 & r_1+2 & \ldots & a'_{22} &|& a'_{22}+1 & a'_{22}+2 & \ldots & n \\
1 & 2 & \ldots & a_{11} &|& c_1+1 & c_1+2 & \ldots & a'_{22} &|& a_{11}+1 & a_{11}+2 & \ldots & c_1 &|& a'_{22}+1 & a'_{22}+2 & \ldots & n
\end{matrix}
\right),
\]
where $a'_{22} = r_1 + a_{21} = c_1 + a_{12} = n - a_{22}$. In matrix form,
$w_{A^t}$ is the block matrix
$\left[
\begin{array}{cccc}
I_{a_{11}} & 0 & 0 & 0 \\
0 & 0 & I_{a_{21}} & 0 \\
0 & I_{a_{12}} & 0 & 0 \\
0 & 0 & 0 & I_{a_{22}}
\end{array}
\right]$.
}

% [DG] The two footnotes are probably way on the pedantic side. I wrote them in
% order to keep track of my computations, but I'm not sure if the reader will
% have any need for them.

For the tower with $G_n=GL_n$, there is slightly more work to be done
to check the equality \eqref{bialgebra-as-Mackey}.
Via Mackey's Theorem~\ref{Mackey-formula} and
Proposition~\ref{double-coset-prop}(c), the left side is
\begin{align}
\nonumber
&\res_{c_1,c_2}^n\left(  \ind_{r_1,r_2}^n \left( U_1 \otimes U_2 \right) \right) \\
\nonumber
&= \left(
 \Res^{G_n}_{P_{c_1,c_2}} \Ind^{G_n}_{P_{r_1,r_2}} \Inf^{P_{r_1,r_2}}_{G_{r_1,r_2}} \left( U_1 \otimes U_2 \right)
\right)^{K_{c_1,c_2}} \\
&= \bigoplus_A
\left(
 \Ind^{P_{c_1,c_2}}_{{}^{w_{A^t}} P_{r_1,r_2} \cap P_{c_1,c_2}}
   \left(\left( \Res^{P_{r_1,r_2}}_{P_{r_1,r_2} \cap P_{c_1,c_2}^{w_{A^t}}}
      \Inf^{P_{r_1,r_2}}_{G_{r_1,r_2}} \left( U_1 \otimes U_2 \right)
   \right)^{\tau_A^{-1}}\right)
\right)^{K_{c_1,c_2}}
\label{left-side-GL-bialgebra}
\end{align}
where $A$ runs over the usual $2 \times 2$ matrices.
The right side is a direct sum over this same set of matrices $A$:
\begin{align}
\nonumber
&\bigoplus_A
  \left( \ind_{a_{11},a_{21}}^{c_1} \otimes \ind^{c_2}_{a_{12},a_{22}}\right)
   \left( \left( \res^{r_1}_{a_{11},a_{12}} U_1 \otimes \res^{r_2}_{a_{21},a_{22}} U_2 \right)^{\tau_A^{-1}} \right) \\
\nonumber
&= \bigoplus_A
 \left(\Ind^{G_{c_1}}_{P_{a_{11},a_{21}}} \otimes \Ind^{G_{c_2}}_{P_{a_{12},a_{22}}}\right) \circ
  \left(\Inf^{P_{a_{11},a_{21}}}_{G_{a_{11},a_{21}}} \otimes \Inf^{P_{a_{12},a_{22}}}_{G_{a_{12},a_{22}}}\right) \\
\nonumber
&\qquad \qquad
    \left(
     \left(
      \left(\Res^{G_{r_1}}_{P_{a_{11},a_{12}}} U_1\right)^{K_{a_{11},a_{12}}}
      \otimes \left(\Res^{G_{r_2}}_{P_{a_{21},a_{22}}} U_2\right)^{K_{a_{21},a_{22}}}
     \right)^{\tau_A^{-1}} \right)\\
\nonumber
&= \bigoplus_A
 \Ind^{G_{c_1,c_2}}_{P_{a_{11},a_{21}} \times P_{a_{12},a_{22}}}  \\
&\qquad \qquad
  \Inf^{P_{a_{11},a_{21}} \times P_{a_{12},a_{22}}}_{G_{a_{11},a_{21},a_{12},a_{22}}}
    \left(
     \left(\left(
      \Res^{G_{r_1,r_2}}_{P_{a_{11},a_{12}} \times P_{a_{21},a_{22}}}  \left( U_1 \otimes U_2 \right)
     \right)^{K_{a_{11},a_{12}} \times K_{a_{21},a_{22}}}
    \right)^{\tau_A^{-1}}\right)
\label{right-side-GL-bialgebra}
\end{align}
(by \eqref{induction-vs-tensor-product}, \eqref{invariants-vs-tensor-product} and their
obvious analogues for restriction and inflation).
Thus it suffices to check for each $2 \times 2$ matrix $A$
that any $\CC G_{c_1,c_2}$-module of the
form $V_1 \otimes V_2$ has the same inner product with the $A$-summands of
\eqref{left-side-GL-bialgebra} and \eqref{right-side-GL-bialgebra}.
Abbreviate $w:=w_{A^t}$ and $\tau:=\tau_A^{-1}$.  

Notice that ${}^w P_{r_1,r_2}$ is the
group of all matrices having the block form
\begin{equation}
\left[
\begin{matrix}
g_{11} & h    & i & j \\
0     &g_{21} & 0 & k \\
d     & e    & g_{12} & \ell \\
0     & f    & 0 & g_{22}
\end{matrix}
\right]
\end{equation}
in which the diagonal blocks $g_{ij}$ for $i,j=1,2$ are invertible of
size $a_{ij} \times a_{ij}$,
while the blocks $h,i,j,k,\ell,d,e,f$ are all arbitrary matrices\footnote{The
blocks $i$ and $j$ have nothing to do with the indices $i,j$ in $g_{ij}$.}
of the appropriate (rectangular) block sizes.
Hence, ${}^{w}P_{r_1,r_2} \cap P_{c_1,c_2}$ is
the group of all matrices having the block form
\begin{equation}
\label{GL-bialgebra-matrixform-1}
\left[
\begin{matrix}
g_{11} & h    & i & j \\
0     &g_{21} & 0 & k \\
0     & 0    & g_{12} & \ell \\
0     & 0    & 0 & g_{22}
\end{matrix}
\right]
\end{equation}
in which the diagonal blocks $g_{ij}$ for $i,j=1,2$ are invertible of
size $a_{ij} \times a_{ij}$,
while the blocks $h,i,j,k,\ell$ are all arbitrary matrices
of the appropriate (rectangular) block sizes;  then
${}^{w}P_{r_1,r_2} \cap G_{c_1,c_2}$ is the subgroup where the blocks $i,j,k$
all vanish. The canonical projection
${}^{w}P_{r_1,r_2} \cap P_{c_1,c_2} \to {}^{w}P_{r_1,r_2} \cap G_{c_1,c_2}$
(obtained by restricting the projection $P_{c_1,c_2} \to G_{c_1,c_2}$)
has kernel ${}^{w}P_{r_1,r_2} \cap P_{c_1,c_2} \cap K_{c_1,c_2}$.
Consequently,
\begin{equation}
\label{GL-bialgebra-quotientgroup-1}
\left({}^wP_{r_1,r_2} \cap P_{c_1,c_2}\right)
/ \left({}^wP_{r_1,r_2} \cap P_{c_1,c_2} \cap K_{c_1,c_2}\right)
= {}^wP_{r_1,r_2} \cap G_{c_1,c_2}.
\end{equation}

Similarly,
\begin{equation}
\label{GL-bialgebra-quotientgroup-2}
\left(P_{r_1,r_2} \cap P_{c_1,c_2}^{w}\right)
/ \left(P_{r_1,r_2} \cap P_{c_1,c_2}^{w} \cap K_{r_1,r_2}\right)
= G_{r_1,r_2} \cap P_{c_1,c_2}^{w}.
\end{equation}

Computing first the inner product
of $V_1 \otimes V_2$ with the $A$-summand of
\eqref{left-side-GL-bialgebra}, and using adjointness properties, one gets
\begin{align*}
&
\left(
\left( \Res^{P_{r_1,r_2}}_{P_{r_1,r_2} \cap P_{c_1,c_2}^{w}}
      \Inf^{P_{r_1,r_2}}_{G_{r_1,r_2}} \left( U_1 \otimes U_2 \right)
   \right)^{\tau} , \right. \\
& \qquad \qquad \qquad \left.
\Res^{P_{c_1,c_2}}_{{}^{w}P_{r_1,r_2} \cap P_{c_1,c_2}}
      \Inf^{P_{c_1,c_2}}_{G_{c_1,c_2}} \left( V_1 \otimes V_2 \right)
\right)_{{}^{w}P_{r_1,r_2} \cap P_{c_1,c_2}} \\
&\overset{\eqref{restriction-inflation-commutation}}{=}
\left(
\left(
  \Inf^{P_{r_1,r_2} \cap P_{c_1,c_2}^{w}}_{G_{r_1,r_2} \cap P_{c_1,c_2}^{w}}
   \Res^{G_{r_1,r_2}}_{G_{r_1,r_2} \cap P_{c_1,c_2}^{w}} \left( U_1 \otimes U_2 \right)
   \right)^{\tau} , \right. \\
& \qquad \qquad \qquad \left.
  \Inf^{{}^wP_{r_1,r_2} \cap P_{c_1,c_2}}_{{}^wP_{r_1,r_2} \cap G_{c_1,c_2}}
   \Res^{G_{c_1,c_2}}_{{}^wP_{r_1,r_2} \cap G_{c_1,c_2}} \left( V_1 \otimes V_2 \right)
\right)_{{}^{w}P_{r_1,r_2} \cap P_{c_1,c_2}}
\end{align*}
(by \eqref{GL-bialgebra-quotientgroup-2} and \eqref{GL-bialgebra-quotientgroup-1}).
One can compute this inner product
by first recalling that ${}^{w}P_{r_1,r_2} \cap P_{c_1,c_2}$ is
the group of matrices having the block form \eqref{GL-bialgebra-matrixform-1}
%\[
%\left[
%\begin{matrix}
%g_{11} & h    & i & j \\
%0     &g_{21} & 0 & k \\
%0     & 0    & g_{12} & \ell \\
%0     & 0    & 0 & g_{22}
%\end{matrix}
%\right]
%\]
in which the diagonal blocks $g_{ij}$ for $i,j=1,2$ are invertible of
size $a_{ij} \times a_{ij}$,
while the blocks $h,i,j,k,\ell$ are all arbitrary matrices
of the appropriate (rectangular) block sizes;  then ${}^{w}P_{r_1,r_2} \cap G_{c_1,c_2}$ is the subgroup where the blocks $i,j,k$ all vanish.  The inner product above then becomes
\begin{equation}
\label{first-big-inner-product-sum}
\begin{aligned}
&\frac{1}{|{}^{w}P_{r_1,r_2} \cap P_{c_1,c_2}|}
\sum\limits_{\substack{(g_{ij})\\(h,i,j,k,\ell)}}
  \chi_{U_1}\left(\begin{matrix} g_{11}& i \\ 0 & g_{12} \end{matrix}\right)
  \chi_{U_2}\left(\begin{matrix} g_{21}& k \\ 0 & g_{22} \end{matrix}\right) \\
&\qquad \qquad \qquad \qquad \qquad \qquad \qquad
  \overline{\chi}_{V_1}\left(\begin{matrix} g_{11}& h \\ 0 & g_{21} \end{matrix}\right)
  \overline{\chi}_{V_2}\left(\begin{matrix} g_{12}& \ell \\ 0 & g_{22} \end{matrix}\right).
\end{aligned}
\end{equation}
If one instead computes the inner product of $V_1 \otimes V_2$ with the $A$-summand of
\eqref{right-side-GL-bialgebra}, using adjointness properties and \eqref{deflation-character-formula-v2} one gets
\begin{align*}
&\left(
\left(
 \left(
  \Res^{G_{r_1,r_2}}_{P_{a_{11},a_{12}} \times P_{a_{21},a_{22}}} \left( U_1 \otimes U_2 \right)
 \right)^{K_{a_{11},a_{12}} \times K_{a_{21},a_{22}}}
\right)^{\tau} , \right.\\
&\qquad\qquad
 \left.
 \left(
  \Res^{G_{c_1,c_2}}_{P_{a_{11},a_{21}} \times P_{a_{12},a_{22}}} \left( V_1 \otimes V_2 \right)
 \right)^{K_{a_{11},a_{21}} \times K_{a_{12},a_{22}}}
\right)_{G_{a_{11},a_{21},a_{12},a_{22}}} \\
&=\frac{1}{|G_{a_{11},a_{21},a_{12},a_{22}}|}
   \sum_{(g_{ij})}
     \frac{1}{|K_{a_{11},a_{12}} \times K_{a_{21},a_{22}}|}
       \sum_{(i,k)}
          \chi_{U_1}\left(\begin{matrix} g_{11}& i \\ 0 & g_{12} \end{matrix}\right)
          \chi_{U_2}\left(\begin{matrix} g_{21}& k \\ 0 & g_{22} \end{matrix}\right)\\
&\qquad \qquad \qquad \qquad \qquad
     \frac{1}{|K_{a_{11},a_{21}} \times K_{a_{12},a_{22}}|} \sum_{(h,\ell)}
  \overline{\chi}_{V_1}\left(\begin{matrix} g_{11}& h \\ 0 & g_{21} \end{matrix}\right)
  \overline{\chi}_{V_2}\left(\begin{matrix} g_{12}& \ell \\ 0 & g_{22} \end{matrix}\right).
\end{align*}
But this right hand side can be seen to equal
\eqref{first-big-inner-product-sum}, after one notes that
\[
|{}^{w}P_{r_1,r_2} \cap P_{c_1,c_2}|
=
|G_{a_{11},a_{21},a_{12},a_{22}}|
\cdot
|K_{a_{11},a_{12}} \times K_{a_{21},a_{22}}|
\cdot
|K_{a_{11},a_{21}} \times K_{a_{12},a_{22}}|
\cdot \#\{j \in \FF_q^{a_{11} \times a_{22}} \}
\]
and that the summands in  \eqref{first-big-inner-product-sum} are independent of the
matrix $j$ in the summation.
\end{proof}

We can also define a $\CC$-vector space $A_\CC$ as the direct sum
$\bigoplus_{n \geq 0} R_\CC(G_n)$.
In the same way as we have made $A = \bigoplus_{n \geq 0} R(G_n)$ into
a $\ZZ$-bialgebra, we can turn
$A_\CC = \bigoplus_{n \geq 0} R_\CC(G_n)$ into
a $\CC$-bialgebra\footnote{The definitions of $m$ and $\Delta$ for
this $\CC$-bialgebra look the same as for $A$: For instance, $m$ is
still defined to be $\ind^{i+j}_{i,j}$ on $\left(A_\CC\right)_i
\otimes \left(A_\CC\right)_j$, where
$\ind^{i+j}_{i,j}$ is defined by the same formulas as in
Definition~\ref{ind-res-defn}. However, the operators of induction,
restriction, inflation and $K$-fixed space construction appearing in
these formulas now act on class functions as opposed to modules.
\par
The fact that these maps $m$ and $\Delta$ satisfy the axioms of a
$\CC$-bialgebra is easy to check: they are merely the $\CC$-linear
extensions of the maps $m$ and $\Delta$ of the $\ZZ$-bialgebra
$A$ (this is because, for instance, induction of class functions
and induction of modules are related by the identity
\eqref{eq.induction.mods-and-class-functs}), and thus satisfy the
same axioms as the latter.}.
There is a $\CC$-bilinear form $\left(\cdot,\cdot\right)_{A_\CC}$
on $A_\CC$ which can be defined either as the $\CC$-bilinear
extension of the $\ZZ$-bilinear form
$\left(\cdot,\cdot\right)_A : A \times A \to \ZZ$ to $A_\CC$, or
(equivalently) as the $\CC$-bilinear form on $A_\CC$ which
restricts to
$\left<\cdot,\cdot\right>_{\Symm_n}$ on every homogeneous component
$R_\CC(G_n)$ and makes different homogeneous components mutually
orthogonal.
The obvious embedding of $A$ into the $\CC$-bialgebra $A_\CC$ (obtained
from the embeddings $R(G_n) \to R_\CC(G_n)$ for all $n$) respects the
bialgebra operations\footnote{This is because, for example, induction
of class functions harmonizes with induction of modules (i.e.,
the equality \eqref{eq.induction.mods-and-class-functs} holds).}, and
the $\CC$-bialgebra $A_\CC$ can be identified
with $A \otimes_\ZZ \CC$ (the result of extending scalars to $\CC$ in
$A$), because every finite group $G$ satisfies
$R_\CC (G) \cong R(G) \otimes_\ZZ \CC$.
The embedding of $A$ into $A_\CC$ also respects the bilinear
forms.

% [DG][v46] Added above short paragraph. Not sure if it makes anything
% clearer, but I felt like it should be here to define some things
% used later (hopefully not earlier).

\begin{exercise}
\label{exe.towers.indres}
Let $G_*$ be one of the three towers.

For every almost-composition
$\alpha = \left(\alpha_1,\alpha_2,\ldots,\alpha_\ell\right)$ of
$n \in \NN$, let us define a map \dfn{$\ind^n_\alpha$} which takes
$\CC G_\alpha$-modules to $\CC G_n$-modules as follows:
If $G_* = \Symm_*$ or $G_* = \Symm_*\left[\Gamma\right]$, we
set
\[
\ind^n_\alpha := \Ind^{G_n}_{G_\alpha}.
\]
If $G_* = GL_*$, then we set
\[
\ind^n_\alpha := \Ind^{G_n}_{P_\alpha} \Inf^{P_\alpha}_{G_\alpha}.
\]
(Note that $\ind^n_\alpha = \ind^n_{i,j}$ if $\alpha$ has the
form $\left(i, j\right)$.)

Similarly, for every almost-composition
$\alpha = \left(\alpha_1,\alpha_2,\ldots,\alpha_\ell\right)$ of
$n \in \NN$, let us define a map \dfn{$\res^n_\alpha$} which takes
$\CC G_n$-modules to $\CC G_\alpha$-modules as follows:
If $G_* = \Symm_*$ or $G_* = \Symm_*\left[\Gamma\right]$, we
set
\[
\res^n_\alpha := \Res^{G_n}_{G_\alpha}.
\]
If $G_* = GL_*$, then we set
\[
\res^n_\alpha := \left(\Res^{G_n}_{P_\alpha} \left( - \right)\right)^{K_\alpha}.
\]
(Note that $\res^n_\alpha = \res^n_{i,j}$ if $\alpha$ has the
form $\left(i, j\right)$.)

\begin{itemize}

\item[(a)] If
$\alpha = \left(\alpha_1,\alpha_2,\ldots,\alpha_\ell\right)$
is an almost-composition of an integer $n \in \NN$ satisfying
$\ell \geq 1$, and if $V_i$ is a
$\CC G_{\alpha_i}$-module for every $i \in \left\{1,2,\ldots,\ell\right\}$,
then show that
\begin{align*}
& \ind^n_{\alpha_1+\alpha_2+\cdots+\alpha_{\ell-1}, \alpha_\ell}
\left(
 \ind^{\alpha_1+\alpha_2+\cdots+\alpha_{\ell-1}}_{\left(\alpha_1,\alpha_2,\ldots,\alpha_{\ell-1}\right)}
  \left(V_1 \otimes V_2 \otimes \cdots \otimes V_{\ell-1}\right)
 \otimes V_\ell
\right) \\
&\cong \ind^n_{\alpha} \left(V_1 \otimes V_2 \otimes \cdots \otimes V_\ell\right) \\
&\cong \ind^n_{\alpha_1, \alpha_2+\alpha_3+\cdots+\alpha_\ell}
\left(
 V_1 \otimes
 \ind^{\alpha_2+\alpha_3+\cdots+\alpha_\ell}_{\left(\alpha_2,\alpha_3,\ldots,\alpha_\ell\right)}
  \left(V_2 \otimes V_3 \otimes \cdots \otimes V_\ell\right)
\right) .
\end{align*}

\item[(b)] Solve Exercise~\ref{exe.Ind-ass} again
using Exercise~\ref{exe.towers.indres}(a).

\item[(c)] We proved above that the map $m : A \otimes A \to A$
(where $A = A\left(G_*\right)$) is associative, by using the
adjointness of $m$ and $\Delta$. Give a new proof of this fact,
which makes no use of $\Delta$.

\item[(d)] If
$\alpha = \left(\alpha_1,\alpha_2,\ldots,\alpha_\ell\right)$
is an almost-composition of an $n \in \NN$, and if
$\chi_i \in R\left(G_{\alpha_i}\right)$ for every
$i \in \left\{1,2,\ldots,\ell\right\}$, then show that
\[
\chi_1 \chi_2 \cdots \chi_\ell
= \ind^n_\alpha \left(\chi_1 \otimes \chi_2 \otimes \cdots \otimes \chi_\ell\right)
\]
in $A = A\left(G_*\right)$.

\item[(e)] If $n \in \NN$, $\ell \in \NN$ and
$\chi \in R\left(G_n\right)$, then show that
\[
\Delta^{\left(\ell-1\right)} \chi
= \sum \res^n_\alpha \chi
\]
in $A^{\otimes \ell}$, where $A = A\left(G_*\right)$. Here,
the sum on the right hand side runs over all
almost-compositions $\alpha$ of $n$ having length $\ell$.
\end{itemize}
\end{exercise}

% [DG][v28] Added this exercise.

\subsection{Symmetric groups}
\label{symmeric-group-section}

Finally, some payoff.
Consider the tower of symmetric groups
$G_n=\Symm_n$, and $A=A(G_*)=:A(\Symm)$\index{$A(\Symm)$}.  Denote by
$\triv_{\Symm_n}, \sgn_{\Symm_n}$ the \emph{trivial} and \emph{sign} characters on
$\Symm_n$.  For a partition $\lambda$ of $n$, denote by
$\triv_{\Symm_\lambda}, \sgn_{\Symm_\lambda}$ the trivial and sign characters
restricted to the Young subgroup
$
\Symm_\lambda = \Symm_{\lambda_1} \times \Symm_{\lambda_2 } \times \cdots,
$
and denote by
$\triv_\lambda$ the class function
which is the characteristic function
for the $\Symm_n$-conjugacy class of permutations of cycle type $\lambda$.

\begin{theorem} \phantomsection
\label{symmetric-group-Frobenius-map-theorem}
\begin{itemize}

\item[(a)]
Irreducible complex characters $\{ \chi^\lambda \}$ of
$\Symm_n$ are indexed by partitions $\lambda$ in $\Par_n$,
and one has a PSH-isomorphism, the
\dfn{Frobenius characteristic map}\index{$\ch$}\footnote{It is unrelated
to the Frobenius endomorphisms from
Exercise~\ref{exe.witt.ghost-app}.},
\[
A=A(\Symm) \overset{\ch}{\longrightarrow} \Lambda
\]
that for $n \geq 0$ and $\lambda \in \Par_n$ sends
\[
\begin{array}{rcl}
\colstrut \triv_{\Symm_n}& \longmapsto & h_n , \\
%\qquad \text{(for }n \geq 0\text{)}\\
\colstrut \sgn_{\Symm_n}& \longmapsto & e_n , \\
%\qquad \text{(for }n \geq 0\text{)}\\
\colstrut \chi^\lambda & \longmapsto & s_\lambda , \\
%\qquad \text{(for }n \geq 0\text{ and }\lambda \in \Par_n\text{)}\\
\colstrut \Ind_{\Symm_\lambda}^{\Symm_n} \triv_{\Symm_\lambda} & \longmapsto &h_\lambda , \\
%\qquad \text{(for }n \geq 0\text{ and }\lambda \in \Par_n\text{)}\\
\colstrut \Ind_{\Symm_\lambda}^{\Symm_n} \sgn_{\Symm_\lambda} & \longmapsto &e_\lambda , \\
%\qquad \text{(for }n \geq 0\text{ and }\lambda \in \Par_n\text{)}\\
\colstrut \triv_\lambda & \longmapsto &\frac{p_\lambda}{z_\lambda} 
\end{array}
\]
(where $\ch$ is extended to a $\CC$-linear map
$A_\CC \to \Lambda_\CC$),
and for $n \geq 1$ sends
\[
\begin{array}{rcl}
\triv_{(n)}& \longmapsto &  \frac{p_n}{n} .
\end{array}
\]
Here, $z_\lambda$ is defined as in
Proposition~\ref{prop-Two-other-Cauchy-expansions}.
% [DG][v70] Inserted \colstruts into the above array to
% prevent the subscripts of one line from overlapping with
% the superscripts of the next.

% [DG][v13] I added the bit about $\ch$ being extended; otherwise, the
% image of $\triv_\lambda$ is not defined (neither $\triv_\lambda$ nor
% $\triv_{(n)}$ generally belong to $A$). Similar change in the proof
% below.

\item[(b)]
For each $n \geq 0$,
the involution on class functions $f: \Symm_n \rightarrow \CC$ 
sending $f \longmapsto \sgn_{\Symm_n} * f$ 
%For every $n \geq 0$, let $\omega_n : R(\Symm_n) \to R(\Symm_n)$ be the involution of the virtual character space $R(\Symm_n)$ given by
where 
\[
(\sgn_{\Symm_n} * f)(g):=\sgn(g) f(g)
\]
%\[
%\left(\omega_n f\right)\left(g\right) = \sgn\left(g\right) f\left(g\right) \qquad \text{for }g\in \Symm_n\text{ and }f\in R(\Symm_n).
%\]
%The direct sum $\bigoplus_{n\geq 0} \omega_n : \bigoplus_{n\geq 0} R(\Symm_n) \to \bigoplus_{n\geq 0} R(\Symm_n) $ of these maps $\omega_n$ is an involution of $\bigoplus_{n\geq 0} R(\Symm_n) = A(\Symm)$. This involution $\bigoplus_{n\geq 0} \omega_n$
preserves the $\ZZ$-sublattice $R(\Symm_n)$ of genuine characters.
The direct sum of these involutions induces an involution 
on $A=A(\Symm) = \bigoplus_{n\geq 0} R(\Symm_n)$
that corresponds under $\ch$ to the involution $\omega$ on $\Lambda$.

\end{itemize}
\end{theorem}

% [DG][v9] I replaced $\sgn_{\Symm_n} \cdot f$ by $\sgn_{\Symm_n} * f$ to
% avoid confusion with the usual product.

\begin{proof}
(a) Corollary~\ref{three-towers-give-bialgebras}
implies that the set $\Sigma=\bigsqcup_{n \geq 0} \Irr(\Symm_n)$
gives a PSH-basis for $A$.
Since a character $\chi$ of $\Symm_n$ has
\begin{equation}
\label{character-coproduct-in-symm}
\Delta(\chi) = \bigoplus_{i+j=n} \Res^{\Symm_n}_{\Symm_i \times \Symm_j} \chi,
\end{equation}
such an element $\chi \in \Sigma \cap A_n$ is never primitive for $n \geq 2$.
Hence the unique irreducible character $\rho=\triv_{\Symm_1}$ of $\Symm_1$
is the only element of $\CCC = \Sigma \cap \Liep$.

% [DG] I again removed $\{\chi\}$.

Thus Theorem~\ref{uniqueness-of-sym-parts-e-f-g}(g)
tells us that there are two PSH-isomorphisms
$A \rightarrow \Lambda$, each of which sends $\Sigma$ to the
PSH-basis of Schur functions $\{s_\lambda\}$ for $\Lambda$.
It also tells us that we can pin down one of the two isomorphisms
to call $\ch$, by insisting that it map the two characters
$\triv_{\Symm_2},\sgn_{\Symm_2}$ in $\Irr(\Symm_2)$ to $h_2,e_2$ (and not
$e_2, h_2$).

Bearing in mind the coproduct formula \eqref{character-coproduct-in-symm},
and the fact that $\triv_{\Symm_n}, \sgn_{\Symm_n}$ restrict, respectively,
to trivial and sign characters of $\Symm_i \times \Symm_j$ for
$i+j=n$, one finds that for $n \geq 2$ one has
$\sgn_{\Symm_2}^\perp$ annihilating $\triv_{\Symm_n}$,
and $\triv_{\Symm_2}^\perp$ annihilating $\sgn_{\Symm_n}$.
Therefore Theorem~\ref{uniqueness-of-sym-parts-a-b-c-d}(b)
(applied to $\Lambda$) implies
$
\triv_{\Symm_n}, \sgn_{\Symm_n}
$
are sent under $\ch$ to $h_n, e_n$.
Then the fact that
$\Ind_{\Symm_\lambda}^{\Symm_n} \triv_{\Symm_\lambda}, \Ind_{\Symm_\lambda}^{\Symm_n} \sgn_{\Symm_\lambda}$
are sent to $h_{\lambda}, e_\lambda$ follows via induction products.

% [DG][v9] Added "(applied to $\Lambda$)".

Recall that the $\CC$-vector space
$A_\CC = \bigoplus_{n\geq 0} R_\CC (\Symm_n)$ is a
$\CC$-bialgebra, and can be identified with $A \otimes_\ZZ \CC$.
The multiplication and the comultiplication of $A_\CC$ are
$\CC$-linear extensions of those of $A$,
and are still given by the same formulas
$m = \ind^{i+j}_{i,j}$ and $\Delta = \bigoplus_{i+j=n} \res^{i+j}_{i,j}$
as those of $A$ (but now, induction and restriction are defined for
class functions, not just for representations).
The $\CC$-bilinear form $\left(\cdot,\cdot\right)_{A_\CC}$ on $A_\CC$
extends both the $\ZZ$-bilinear form
$\left(\cdot,\cdot\right)_A$ on $A$ and the $\CC$-bilinear forms
$\left<\cdot,\cdot\right>_{\Symm_n}$ on all $R_\CC(\Symm_n)$.

% [DG][v45] Added the preceding pedantic paragraph since we talk
% of $\triv_{(n)}$ as being primitive in $A_\CC$ in the next
% paragraph, thus requiring a Hopf structure on $A_\CC$.

% [DG][v46] Rewrote it, as part of now comes earlier.

For the assertion about $\triv_{(n)}$, note that it is primitive in
$A_\CC$ for $n \geq 1$, because as a class function, the indicator function
of $n$-cycles vanishes upon restriction to $\Symm_i \times \Symm_j$
for $i+j=n$ if both $i,j \geq 1$;  these subgroups contain no $n$-cycles.
Hence Corollary~\ref{primitives-in-Lambda} implies that
$\ch(\triv_{(n)})$ is a scalar multiple of $p_n$.  To pin down the scalar,
note $p_n=m_{(n)}$ so
$
(h_n,p_n)_{\Lambda}=(h_n,m_n)_{\Lambda}=1,
$
while $\ch^{-1}(h_n)=\triv_{\Symm_n}$ has
\[
(\triv_{\Symm_n},\triv_{(n)})=\frac{1}{n!} \cdot (n-1)!=\frac{1}{n}.
\]
\footnote{The first equality sign in this computation uses
the fact that the number of all $n$-cycles in $\Symm_n$ is
$(n-1)!$. This is because any $n$-cycle in $\Symm_n$ can
be uniquely written in the form
$\left(i_1,i_2,\ldots,i_{n-1},n\right)$ (in cycle notation)
with $\left(i_1,i_2,\ldots,i_{n-1}\right)$ being a
permutation in $\Symm_{n-1}$ (written in one-line notation).}
Thus $\ch(\triv_{(n)})=\frac{p_n}{n}$.  The fact that
$\ch(\triv_{\lambda})=\frac{p_\lambda}{z_\lambda}$ then follows via induction
product calculations\footnote{For instance, one can use
\eqref{induced-character-formula} to show that
$z_\lambda \triv_\lambda
= \lambda_1 \lambda_2 \cdots \lambda_\ell
\cdot \triv_{\left(\lambda_1\right)} \triv_{\left(\lambda_2\right)} \cdots \triv_{\left(\lambda_\ell\right)}$
if $\lambda = \left(\lambda_1, \lambda_2, \ldots, \lambda_\ell\right)$
with $\ell = \ell\left(\lambda\right)$.
See Exercise~\ref{exe.frob-char.plambda}(d) for the details.}.
Part (b) follows from 
Exercise~\ref{internal-tensor-product-with-sgn-is-omega-exercise} below.
\end{proof}

% [DG][v14] Added the preceding footnote because IMHO that's a lot of
% proof that is being skipped here, and at least it should be
% explicitly left to the reader.

% [DG][v45] Added yet another footnote.

\begin{remark}
The paper of Liulevicius \cite{Liulevicius} gives
a very elegant alternate approach to the Frobenius map
as a Hopf isomorphism $A(\Symm) \overset{\ch}{\longrightarrow} \Lambda$,
inspired by equivariant $K$-theory and vector bundles over spaces which
are finite sets of points!
\end{remark}

\begin{exercise}
\label{exe.frob-char.plambda}

If $P$ is a subset of a group $G$, we denote by $\triv_{P}$ the map
$G \to \CC$ which sends every element of $P$ to $1$ and all
remaining elements of $G$ to $0$. \ \ \ \ \footnote{This is not in conflict
with the notation $\triv_G$ for the trivial character of $G$, since
$\triv_P = \triv_G$ for $P=G$. Note that $\triv_P$
is a class function when $P$ is a union of conjugacy classes of $G$.} For any
finite group $G$ and any $h\in G$, we introduce the following notations:

\begin{itemize}
\item Let $Z_{G}\left(  h\right)$ denote the centralizer of $h$ in $G$.

\item Let $\operatorname{Conj}_{G}\left(  h\right)  $ denote the conjugacy
class of $h$ in $G$.

\item Define a map $\alpha_{G,h} : G \to \CC$ by $\alpha
_{G,h}=\left\vert Z_{G}\left(  h\right)  \right\vert
\triv_{\operatorname{Conj}_{G}\left(  h\right)  }$. This map
$\alpha_{G,h}$ is a
class function\footnote{In fact, $\triv_{\operatorname{Conj}%
_{G}\left(  h\right)  }$ is a class function (since $\operatorname{Conj}%
_{G}\left(  h\right)  $ is a conjugacy class), and so $\alpha_{G,h}$ (being
the scalar multiple $\left\vert Z_{G}\left(  h\right)  \right\vert
\triv_{\operatorname{Conj}_{G}\left(  h\right)  }$ of
$\triv_{\operatorname{Conj}_{G}\left(  h\right)  }$) must also be a
class function.}.
\end{itemize}

\begin{enumerate}
\item[(a)] Prove that $\alpha_{G,h}\left(  g\right)  =\sum_{k\in G}\left[
khk^{-1}=g\right]  $ for every finite group $G$ and any $h\in G$ and $g\in G$.
Here, we are using the Iverson bracket notation (that is, for any statement
$\mathcal{A}$, we define $\left[  \mathcal{A}\right]  $ to be the integer $1$
if $\mathcal{A}$ is true, and $0$ otherwise).

\item[(b)] Prove that if $H$ is a subgroup of a finite group $G$, and if $h\in
H$, then $\Ind_{H}^{G}\alpha_{H,h}=\alpha_{G,h}$.

\item[(c)] Prove that if $G_{1}$ and $G_{2}$ are finite groups, and if
$h_{1}\in G_{1}$ and $h_{2}\in G_{2}$, then the canonical isomorphism
$R_\CC \left(  G_{1}\right)  \otimes R_\CC \left(
G_{2}\right)  \rightarrow R_\CC \left(  G_{1}\times G_{2}\right)  $
sends $\alpha_{G_{1},h_{1}}\otimes\alpha_{G_{2},h_{2}}$ to $\alpha
_{G_{1}\times G_{2},\left(  h_{1},h_{2}\right)  }$.

\item[(d)] Fill in the details of the proof of $\ch
(\underline{1}_{\lambda})=\frac{p_{\lambda}}{z_{\lambda}}$ in the proof of
Theorem~\ref{symmetric-group-Frobenius-map-theorem}.

\item[(e)] Obtain an alternative proof of Remark~\ref{rmk.zlambda.centralizer}.

\item[(f)] If $G$ and $H$ are two finite groups, and if
$\rho : H \rightarrow G$ is a group homomorphism, then prove that
$\Ind_{\rho} \alpha_{H, h} = \alpha_{G, \rho \left( h \right)}$ for
every $h \in H$, where $\Ind_{\rho} \alpha_{H, h}$ is defined as in
Exercise~\ref{exe.Indrho}.

\end{enumerate}
\end{exercise}

% [DG][v47] Added the preceding exercise. (I use parts (b) and (c) in
% the solution of another.)

\begin{exercise}
\label{internal-tensor-product-with-sgn-is-omega-exercise}
If $G$ is a group and $U_1$ and $U_2$ are two $\CC G$-modules, then the
tensor product $U_1 \otimes U_2$ is a $\CC\left[G\times G\right]$-module,
which can be made into a $\CC G$-module by letting $g \in G$ act as
$\left(g,g\right) \in G\times G$. This $\CC G$-module $U_1 \otimes U_2$
is called the \dfn{inner tensor product}\index{tensor product of representations}%
\footnote{Do not confuse this
with the inner product of characters.} of $U_1$ and $U_2$, and is a
restriction of the outer tensor product $U_1 \otimes U_2$ using the
inclusion map $G \to G\times G,\ g\mapsto\left(g,g\right)$.

Let $n\geq 0$, and let $\sgn_{\Symm_n}$ be the $1$-dimensional
$\CC \Symm_n$-module $\CC$ on which every $g \in \Symm_n$ acts as
multiplication by $\sgn(g)$. If $V$ is a $\CC \Symm_n$-module, show
that the involution on $A(\Symm)=\bigoplus_{n \geq 0} R(\Symm_n)$
defined in Theorem~\ref{symmetric-group-Frobenius-map-theorem}(b)
sends $\chi_V \mapsto \chi_{\sgn_{\Symm_n} \otimes V}$ 
%$\omega_n \left(\chi_V\right) = \chi_{\sgn_{\Symm_n} \otimes V}$,
where $\sgn_{\Symm_n} \otimes V$ is the inner tensor product of
$\sgn_{\Symm_n}$ and $V$.
Use this to show that this involution 
% $\bigoplus_{n\geq 0} \omega_n$ 
is a nontrivial PSH-automorphism of $A(\Symm)$, and deduce
Theorem~\ref{symmetric-group-Frobenius-map-theorem}(b).
%The fact that it corresponds to $\omega$ on $\Lambda$ now follows from
%Theorem~\ref{uniqueness-of-sym-parts-e-f-g}(f), once one notes that the latter
%operation induces a (nontrivial) PSH-automorphism of $A$.
\end{exercise}

\begin{exercise}
\label{exe.frob-char.ptype}
Let $n \in \NN$. For every permutation
$\sigma \in \Symm_n$, we let $\operatorname*{type}\sigma$ denote the
cycle type of $\sigma$. Extend $\operatorname*{ch}:A=A\left(  \Symm
\right)  \rightarrow\Lambda$ to a $\CC$-linear map $A_{\CC}
\rightarrow\Lambda_{\CC}$. We shall call the latter map
$\operatorname*{ch}$, too.

\begin{enumerate}
\item[(a)] Prove that every class function $f\in R_{\CC}\left(
 \Symm_n\right)  $ satisfies
\[
\operatorname*{ch}\left(  f\right)  =\dfrac{1}{n!}\sum_{\sigma\in
 \Symm_n} f\left(  \sigma\right)  p_{\operatorname*{type}\sigma}.
\]

\item[(b)] Let $H$ be a subgroup of $\Symm_n$. Prove that every class
function $f\in R_{\CC}\left(  H\right)  $ satisfies%
\[
\operatorname*{ch}\left(  \operatorname*{Ind}\nolimits_{H}^{ \Symm_{n}%
}f\right)  =\dfrac{1}{\left\vert H\right\vert }\sum_{h\in H}f\left(  h\right)
p_{\operatorname*{type}h}.
\]

\end{enumerate}
\end{exercise}

\begin{exercise} \phantomsection
\label{exe.symmetric-group-reps-rational}
\begin{itemize}
\item[(a)] Show that for every $n \geq 0$, every $g \in \Symm_n$ and every
finite-dimensional $\CC \Symm_n$-module $V$, we have
$\chi_V\left(g\right) \in \ZZ$.
\item[(b)] Show that for every $n \geq 0$ and every
finite-dimensional $\CC \Symm_n$-module $V$, there exists a
$\QQ \Symm_n$-module $W$ such that $V \cong \CC \otimes_{\QQ} W$.
(In the representation theorists' parlance, this says that all
representations of $\Symm_n$ are \emph{defined over $\QQ$}. This
part of the exercise requires some familiarity with representation
theory.)
\end{itemize}
\end{exercise}

\begin{remark}
Parts (a) and (b) of Exercise~\ref{exe.symmetric-group-reps-rational} both
follow from an even stronger result: For every $n\geq 0$ and every
finite-dimensional $\CC \Symm_{n}$-module $V$, there exists a
$\ZZ \Symm_n$-module $W$ which is finitely generated and free
as a $\ZZ$-module and satisfies $V\cong \CC\otimes_\ZZ W$
as $\CC \Symm_n$-modules. This follows from the combinatorial approach to the
representation theory of $\Symm_n$, in which the irreducible
representations of $\CC \Symm_n$ (the \emph{Specht modules}) are constructed
using Young tableaux and tabloids. See the literature on the symmetric group,
e.g., \cite{Sagan}, \cite[\S 7]{Fulton}, \cite{Wildon}
or \cite[Section 2.2]{Krob} for this approach.
\end{remark}

% [DG][v13] Added the above exercise and remark. Haven't done the
% wreath-product analogue so far.

The connection between $\Lambda$ and $A\left( \Symm \right)$ as
established in Theorem~\ref{symmetric-group-Frobenius-map-theorem} benefits
both the study of $\Lambda$ and that of $A\left( \Symm \right)  $. The
following two exercises show some applications to $\Lambda$:

\begin{exercise}
\label{exe.repSn.inner-tensor.intprod}
If $G$ is a group and $U_1$ and $U_2$ are two $\CC G$-modules,
then let $U_1 \boxtimes U_2$ denote the inner tensor product of
$U_1$ and $U_2$ (as defined in
Exercise~\ref{internal-tensor-product-with-sgn-is-omega-exercise}).
Consider also the binary operation $*$ on $\Lambda_\QQ$ defined in
Exercise~\ref{exe.Lambda.maps-on-pn}(h).

\begin{itemize}

\item[(a)] Show that
$\ch \left(\chi_{U_1 \boxtimes U_2}\right)
= \ch\left(\chi_{U_1}\right) * \ch\left(\chi_{U_2}\right)$
for any $n \in \NN$ and any two $\CC \Symm_n$-modules $U_1$ and
$U_2$.

\item[(b)] Use this to obtain a new solution for
Exercise~\ref{exe.Lambda.maps-on-pn}(h).

\item[(c)] Show that $s_\mu * s_\nu \in \sum_{\lambda \in \Par}
\NN s_\lambda$ for any two partitions $\mu$ and $\nu$.

\end{itemize}

[\textbf{Hint:} For any group $G$,
introduce a binary operation $*$ on $R_\CC \left(G\right)$
which satisfies
$\chi_{U_1 \boxtimes U_2} = \chi_{U_1} * \chi_{U_2}$ for any
two $\CC G$-modules $U_1$ and $U_2$.]
\end{exercise}

% [DG][v25] Added the above exercise.

\begin{exercise}
\label{exe.repSn.arithmetic-prod}
Define a $\QQ$-bilinear map
$\boxdot:\Lambda_{\QQ}\times\Lambda_{\QQ}\rightarrow
\Lambda_{\QQ}$, which will be written in infix notation (that is, we
will write $a\boxdot b$ instead of $\boxdot\left(  a,b\right)  $), by setting
\[
p_{\lambda}\boxdot p_{\mu}
= \prod_{i=1}^{\ell\left(  \lambda\right)  }
\prod_{j=1}^{\ell\left(  \mu\right)  }
p_{\operatorname{lcm}\left(  \lambda_{i},\mu_{j}\right)  }^{
\gcd\left(  \lambda_{i},\mu_{j}\right)  }
\qquad\qquad \text{for any partitions }\lambda\text{ and }\mu.
\]
\footnote{This is well-defined, since $\left(  p_{\lambda}\right)
_{\lambda\in\Par}$ is a $\QQ$-module basis of $\Lambda_{\QQ}$.}

\begin{enumerate}
\item[(a)] Show that $\Lambda_\QQ$, equipped with the binary
operation $\boxdot$, becomes a commutative $\QQ$-algebra with unity
$p_1$.

\item[(b)] For every $r\in\ZZ$, define the $\QQ$-algebra
homomorphism $\epsilon_{r}:\Lambda_\QQ\rightarrow\QQ$ as in
Exercise~\ref{exe.Lambda.maps-on-pn}(c). Show that $1\boxdot f=\epsilon
_{1}\left(  f\right)  1$ for every $f\in\Lambda_\QQ$ (where $1$
denotes the unity of $\Lambda$).

\item[(c)] Show that $s_{\mu}\boxdot s_{\nu}\in\sum_{\lambda\in
\Par} \NN s_{\lambda}$ for any two partitions $\mu$ and
$\nu$.

\item[(d)] Show that $f\boxdot g\in\Lambda$ for any $f\in\Lambda$ and
$g\in\Lambda$.
\end{enumerate}

[\textbf{Hint:} For every set $X$, let $\Symm_{X}$ denote the group of
all permutations of $X$. For two sets $X$ and $Y$, there is a canonical group
homomorphism $\Symm_{X}\times\Symm_{Y}\rightarrow
\Symm_{X\times Y}$, which is injective if $X$ and $Y$ are
nonempty. For positive
integers $n$ and $m$, this yields an embedding $\Symm_{n}
\times\Symm_{m}\rightarrow\Symm_{\left\{  1,2,\ldots,n\right\}
\times\left\{  1,2,\ldots,m\right\}  }$, which, once $\Symm_{\left\{
1,2,\ldots,n\right\}  \times\left\{  1,2,\ldots,m\right\}  }$ is identified
with $\Symm_{nm}$ (using an arbitrary but fixed bijection $\left\{
1,2,\ldots,n\right\}  \times\left\{  1,2,\ldots,m\right\}  \rightarrow\left\{
1,2,\ldots,nm\right\}  $), can be regarded as an embedding
$\Symm_{n}\times\Symm_{m}\rightarrow\Symm_{nm}$ and thus allows
defining a $\CC \Symm_{nm}$-module
$\Ind_{\Symm_{n}\times\Symm_{m}}^{\Symm_{nm}}\left(
U\otimes V\right)  $ for any $\CC \Symm_{n}$-module $U$ and any
$\CC \Symm_{m}$-module $V$. This gives a binary operation on
$A\left( \Symm \right)$. Show that this operation corresponds to
$\boxdot$ under the PSH-isomorphism $\ch:A\left(
\Symm\right)  \rightarrow\Lambda$.]
\end{exercise}

\begin{remark}
The statements (and the idea of the solution) of Exercise
\ref{exe.repSn.arithmetic-prod} are due to Manuel Maia and Miguel M\'{e}ndez
(see \cite{MaiaMendez} and, more explicitly, \cite{Mendez-MO}), who call the
operation $\boxdot$ the \dfn{arithmetic product}. Li \cite[Thm.
3.5]{Li-primegraphs} denotes it by $\boxtimes$ and relates it to the
enumeration of unlabelled graphs.
\end{remark}

% [DG][v47] Added the above exercise and remark.

\subsection{Wreath products}
\label{wreath-product-section}

Next consider the tower of groups
$G_n=\Symm_n[\Gamma]$ for a finite group $\Gamma$,
and the Hopf algebra $A=A(G_*)=:A(\Symm[\Gamma])$.  Recall
%from the previous section
(from Theorem~\ref{symmetric-group-Frobenius-map-theorem})
that irreducible complex representations $\chi^\lambda$ of $\Symm_n$ are
indexed by partitions $\lambda$ in $\Par_n$.  Index the
irreducible complex representations of $\Gamma$ as
$
\Irr(\Gamma) = \{ \rho_1,\ldots,\rho_d\}.
$

% [DG][v37] Replaced section reference by theorem
% reference above.

\begin{definition}
Define for a partition $\lambda$ in $\Par_n$
and $\rho$ in $\Irr(\Gamma)$
a representation $\chi^{\lambda,\rho}$
of $\Symm_n[\Gamma]$ in which $\sigma$ in $\Symm_n$ and
$\gamma=(\gamma_1,\ldots,\gamma_n)$ in $\Gamma^n$ act
on the space
$\chi^\lambda \otimes \left( \rho^{\otimes n} \right)$
as follows:
\begin{equation}
\label{wreath-irreducible-construction}
\begin{aligned}
\sigma( u \otimes (v_1 \otimes \cdots \otimes v_n) )
&= \sigma(u) \otimes (v_{\sigma^{-1}(1)} \otimes \cdots \otimes v_{\sigma^{-1}(n)}) ; \\
\gamma( u \otimes (v_1 \otimes \cdots \otimes v_n) )
&= u \otimes (\gamma_1 v_1 \otimes \cdots \otimes \gamma_n v_n) . \\
\end{aligned}
\end{equation}
\end{definition}

\begin{theorem}
\label{wreath-product-irreducibles}
The irreducible $\CC \Symm_n[\Gamma]$-modules
are the induced characters
\[
 \chi^{\underline{\lambda}}: =
\Ind_{\Symm_{\operatorname{degs}(\underline{\lambda})}[\Gamma]}^{\Symm_n[\Gamma]}
  \left(
   \chi^{\lambda^{(1)},\rho_1} \otimes \cdots
   \otimes \chi^{\lambda^{(d)},\rho_d}
  \right)
\]
as $\underline{\lambda}$ runs through all functions
\[
\begin{array}{rcl}
\Irr(\Gamma) &\overset{\underline{\lambda}}{\longrightarrow}&   \Par , \\
\rho_i &\longmapsto&  \lambda^{(i)}
\end{array}
\]
with the property that $\sum_{i=1}^d |\lambda^{(i)}|=n$.
Here, $\operatorname{degs}(\underline{\lambda})$
denotes the $d$-tuple
$\left(\left|\lambda^{(1)}\right|,\left|\lambda^{(2)}\right|,
\ldots,\left|\lambda^{(d)}\right|\right) \in \NN^d$, and
$\Symm_{\operatorname{degs}(\underline{\lambda})}$ is
defined as the subgroup
$\Symm_{\left|\lambda^{(1)}\right|} \times
\Symm_{\left|\lambda^{(2)}\right|} \times \cdots \times
\Symm_{\left|\lambda^{(d)}\right|}$ of $\Symm_n$.

% [DG][v14] Introduced $\operatorname{degs}(\underline{\lambda})$ to
% correct the definition of $\chi^{\underline{\lambda}}$ (it used to
% have a $\lambda$ instead, which made no sense).
% Also fixed a typo below ($\Lambda^{\otimes n}$ replaced by
% $\Lambda^{\otimes d}$).

Furthermore, one has a PSH-isomorphism
\[
\begin{array}{rcl}
A(\Symm[\Gamma]) &\longrightarrow &\Lambda^{\otimes d} ,\\
\chi^{\underline{\lambda}} & \longmapsto &
 s_{\lambda^{(1)}} \otimes \cdots \otimes s_{\lambda^{(d)}}.
\end{array}
\]
\end{theorem}
\begin{proof}
We know from Corollary~\ref{three-towers-give-bialgebras} that
$A(\Symm[\Gamma])$ is a PSH, with PSH-basis $\Sigma$
given by the union of all irreducible characters of all
groups $\Symm_n[\Gamma]$.  Therefore
Theorem~\ref{Zelevinsky-decomposition-theorem}
tells us that
$A(\Symm[\Gamma]) \cong \bigotimes_{\rho \in \CCC} A(\Symm[\Gamma])(\rho)$
where $\CCC$ is the set of irreducible characters which are
also primitive.  Just as in the case of $\Symm_n$,
it is clear from the definition of the coproduct that
an irreducible character $\rho$ of $\Symm_n[\Gamma]$ is
primitive if and only if $n=1$, that in this case
$\Symm_n[\Gamma]=\Gamma$,
and $\rho$ lies in $\Irr(\Gamma)=\{\rho_1,\ldots,\rho_d\}$.

The remaining assertions of the theorem will then follow from the definition
of the induction product algebra structure on $A(\Symm[\Gamma])$, once
we have shown that, for every $\rho \in \Irr(\Gamma)$,
 there is a PSH-isomorphism sending
\begin{equation}
\label{desired-wreath-PSH-iso}
\begin{array}{rcl}
A(\Symm) &\longrightarrow &A(\Symm[\Gamma])(\rho) ,\\
\chi^{\lambda} &\longmapsto &\chi^{\lambda,\rho}.
\end{array}
\end{equation}
Such an isomorphism comes from applying
Proposition~\ref{semidirect-product-adjoint-functors}
to the semidirect product $\Symm_n[\Gamma] = \Symm_n \ltimes \Gamma^n$,
so that $K=\Gamma^n, G=\Symm_n$, and fixing $V=\rho^{\otimes n}$ as
$\CC \Symm_n[\Gamma]$-module with structure as defined
in \eqref{wreath-irreducible-construction} (but with $\lambda$ set to
$\left(n\right)$, so that $\chi^\lambda$ is the trivial $1$-dimensional
$\CC \Symm_n$-module).
One obtains for each $n$, maps
\[
R(\Symm_n) \underset{\Psi}{\overset{\Phi}{\rightleftharpoons}}
R(\Symm_n[\Gamma])
\]
where
\[
\begin{array}{rcl}
\chi &\overset{\Phi}{\longmapsto} &\chi \otimes (\rho^{\otimes n}) ,\\
\alpha&\overset{\Psi}{\longmapsto}&
  \Hom_{\CC \Gamma^n}(\rho^{\otimes n},\alpha).
\end{array}
\]
Taking the direct sum of these maps for all $n$ gives maps
$
A(\Symm) \underset{\Psi}{\overset{\Phi}{\rightleftharpoons}}
A(\Symm[\Gamma]).
$

% [DG][v14] Changed $A(\Symm_n)$ to $R(\Symm_n)$ since you seem to
% reserve the notation $A$ for the direct sum of the representation
% rings rather than its specific addends (and, alas, also for a matrix).
% 
% Added "(but with $\lambda$ set to
% $\left(n\right)$, so that $\chi^\lambda$ is the trivial $1$-dimensional
% $\CC \Symm_n$-module)".

These maps are coalgebra morphisms because of
their interaction with restriction to $\Symm_i \times \Symm_j$.
Since Proposition~\ref{semidirect-product-adjoint-functors}(iii)
gives the adjointness property that
\[
(\chi,\Psi(\alpha))_{A(\Symm)} = (\Phi(\chi),\alpha)_{A(\Symm[\Gamma])},
\]
one concludes from the self-duality of $A(\Symm), A(\Symm[\Gamma])$
that $\Phi,\Psi$ are also algebra morphisms.  Since they take genuine
characters to genuine characters, they are PSH-morphisms.
Since $\rho$ being a simple $\CC \Gamma$-module implies that
$V=\rho^{\otimes n}$ is a simple $\CC \Gamma^n$-module,
Proposition~\ref{semidirect-product-adjoint-functors}(iv) shows that
\begin{equation}
\label{Psi-Phi-composite}
(\Psi \circ \Phi)(\chi) = \chi
\end{equation}
for all $\Symm_n$-characters $\chi$.   Hence $\Phi$ is an injective
PSH-morphism.   Using adjointness, \eqref{Psi-Phi-composite} also shows that
$\Phi$ sends $\CC \Symm_n$-simples $\chi$
to $\CC[\Symm_n[\Gamma]]$-simples $\Phi(\chi)$:
\[
(\Phi(\chi),\Phi(\chi))_{A(\Symm[\Gamma])}
 = ((\Psi \circ \Phi)(\chi),\chi)_{A(\Symm)}
 = (\chi,\chi)_{A(\Symm)} = 1.
\]
Since $\Phi(\chi)= \chi \otimes (\rho^{\otimes n})$
has $V=\rho^{\otimes n}$ as a constituent upon restriction to $\Gamma^n$,
Frobenius Reciprocity shows that the irreducible character
$\Phi(\chi)$ is a constituent of
$\Ind_{\Gamma^n}^{\Symm_n[\Gamma]} \rho^{\otimes n}= \rho^n$.  Hence the
entire image of $\Phi$ lies in $A(\Symm[\Gamma])(\rho)$ (due to how
we defined $A(\rho)$ in the proof of
Theorem~\ref{Zelevinsky-decomposition-theorem}), and so $\Phi$
must restrict to an isomorphism as desired in \eqref{desired-wreath-PSH-iso}.
\end{proof}

One of Zelevinsky's sample applications of the theorem is this branching rule.

\begin{corollary}
\label{wreath-product-branching-rule}
Given $\underline{\lambda}=(\lambda^{(1)},\ldots,\lambda^{(d)})$ with $\sum_{i=1}^d |\lambda^{(i)}|=n$,
one has
\[
\Res^{\Symm_n[\Gamma]}_{\Symm_{n-1}[\Gamma] \times \Gamma}  \left( \chi^{\underline{\lambda}} \right)
= \sum_{i=1}^d \sum\limits_{\substack{ \lambda^{(i)}_- \subseteq \lambda^{(i)}: \\ |\lambda^{(i)}/\lambda^{(i)}_-|=1}}
\chi^{(\lambda^{(1)}, \ldots,\lambda^{(i)}_-, \ldots,\lambda^{(d)})} \otimes \rho_i.
\]
(We are identifying functions $\underline{\lambda} : \Irr(\Gamma)
\to \Par$ with the corresponding $d$-tuples
$\left(\lambda^{(1)}, \lambda^{(2)}, \ldots, \lambda^{(d)}\right)$
here.)
\end{corollary}

\begin{example}
For $\Gamma$ a two-element group, so $\Irr(\Gamma)=\{\rho_1,\rho_2\}$ and
$d=2$, then
\[
\Res^{\Symm_6[\Gamma]}_{\Symm_{5}[\Gamma] \times \Gamma}
\left( \chi^{((3,1),(1,1))} \right)
=
\chi^{((3),(1,1))} \otimes \rho_1
+\chi^{((2,1),(1,1))} \otimes \rho_1
+\chi^{((3,1),(1))} \otimes \rho_2.
\]
\end{example}

\begin{proof}[Proof of Corollary~\ref{wreath-product-branching-rule}]
By Theorem~\ref{wreath-product-irreducibles}, this is equivalent to
computing in the Hopf algebra $A:=\Lambda^{\otimes d}$ the component of the
coproduct of $s_{\lambda^{(1)}} \otimes \cdots \otimes s_{\lambda^{(d)}}$
that lies in $A_{n-1} \otimes A_1$.  Working within each tensor factor $\Lambda$,
we conclude from Proposition~\ref{symm-comultiplication-formulas-prop}(iv)
% the Pieri formula implies
that the
$\Lambda_{\left|\lambda\right|-1} \otimes \Lambda_1$-component
of $\Delta(s_\lambda)$ is
\[
\sum\limits_{\substack{\lambda_- \subseteq \lambda: \\|\lambda/\lambda_-|=1}} s_{\lambda_-} \otimes \rho.
\]
One must apply this in each of the $d$ tensor factors of $A=\Lambda^{\otimes d}$, then sum on $i$.
\end{proof}

\subsection{General linear groups}
\label{general-linear-section}

We now consider the tower of finite general linear groups
$G_n=GL_n=GL_n(\FF_q)$ and $A=A(G_*)=:A(GL)$.
Corollary~\ref{three-towers-give-bialgebras} tells us that
$A(GL)$ is a PSH, with PSH-basis $\Sigma$
given by the union of all irreducible characters of all
groups $GL_n$.  Therefore Theorem~\ref{Zelevinsky-decomposition-theorem}
tells us that
\begin{equation}
\label{GL-Zelevinsky-decomposition}
A(GL) \cong \bigotimes_{\rho \in \CCC} A(GL)(\rho)
\end{equation}
where $\CCC=\Sigma \cap \Liep$ is the set of primitive irreducible characters.

\begin{definition}
Call an irreducible representation $\rho$ of $GL_n$
\emph{cuspidal}\index{cuspidal representation of $GL_n$}
for $n \geq 1$ if it lies in $\CCC$, that is, its restriction to proper parabolic subgroups $P_{i,j}$ with $i+j=n$
and $i,j>0$ contain no nonzero vectors which are $K_{i,j}$-invariant.
Given an irreducible character $\sigma$ of $GL_n$, say that $d(\sigma)=n$,
and let $\CCC_n:=\{ \rho \in \CCC: d(\rho)=n\}$ for $n \geq 1$
denote the subset of cuspidal characters of $GL_n$.
\end{definition}

Just as was the case for $\Symm_1$ and $\Symm_1[\Gamma]=\Gamma$, \emph{every}
irreducible character $\rho$ of $GL_1(\FF_q)=\FF_q^\times$ is cuspidal.
However, this does not exhaust the cuspidal characters.  In fact, one
can predict the number of cuspidal characters in $\CCC_n$, using
knowledge of the number of conjugacy classes in $GL_n$.
Let $\FFF$ denote the set of all nonconstant monic irreducible polynomials $f(x) \neq x$ in $\FF_q[x]$.
Let $\FFF_n:=\{ f \in \FFF: \deg(f)=n\}$ for $n \geq 1$.

\begin{proposition}
The number $|\CCC_n|$ of cuspidal characters of $GL_n(\FF_q)$
is the number of $|\FFF_n|$ of irreducible monic degree $n$ polynomials $f(x) \neq x$
in $\FF_q[x]$ with nonzero constant term.
\end{proposition}
\begin{proof}
We show  $|\CCC_n|=|\FFF_n|$ for $n \geq 1$ by strong induction on $n$.
For the base case\footnote{Actually, we don't need any base case for our
strong induction. We nevertheless handle the case $n=1$ as a warmup.}
$n=1$, just as with the families $G_n=\Symm_n$ and
$G_n=\Symm_n[\Gamma]$, when $n=1$ any irreducible character
$\chi$ of $G_1=GL_1(\FF_q)$
gives a primitive element of $A=A(GL)$, and hence is cuspidal.
Since $GL_1(\FF_q)=\FF_q^\times$ is abelian, there are $|\FF_q^\times|=q-1$
such cuspidal characters in $\CCC_1$, which agrees with the fact that
there are $q-1$ monic (irreducible)
linear polynomials $f(x) \neq x$ in $\FF_q[x]$, namely
$
\FFF_1:=\{ f(x)=x-c: c \in \FF_q^\times\}.
$

In the inductive step, use the fact that
the number $|\Sigma_n|$ of irreducible complex characters $\chi$
of $GL_n(\FF_q)$ equals its number of conjugacy classes.
These conjugacy classes are
uniquely represented by \emph{rational canonical forms},
which are parametrized
by functions $\underline{\lambda}: \FFF \rightarrow \Par$ with
the property that $\sum_{f \in \FFF} \deg(f) |\underline{\lambda}(f)|=n$.
On the other hand, \eqref{GL-Zelevinsky-decomposition} tells us that
$|\Sigma_n|$ is similarly parametrized by the functions
$\underline{\lambda}: \CCC \rightarrow \Par$ having
the property that $\sum_{\rho \in \CCC} d(\rho) |\underline{\lambda}(\rho)|=n$.
Thus we have parallel disjoint decompositions
\[
\begin{array}{rll}
\CCC &= \bigsqcup_{ n \geq 1} \CCC_n
         & \text{ where }\CCC_n=\left\{\rho \in \CCC:d(\rho)=n\right\},\\
\FFF &= \bigsqcup_{ n \geq 1} \FFF_n
         & \text{ where }\FFF_n=\left\{f \in \FFF:\deg(f)=n\right\},
\end{array}
\]
and hence an equality for all $n \geq 1$
\[
\left|\left\{ \CCC\overset{\underline{\lambda}}{\longrightarrow} \Par:
\quad
\sum_{\rho \in \CCC} d(\rho) |\underline{\lambda}(\rho)| =n\right\}\right| \\
=|\Sigma_n|
=\left|\left\{ \FFF\overset{\underline{\lambda}}{\longrightarrow} \Par:
\quad
\sum_{f \in \FFF} \deg(f) |\underline{\lambda}(f)| =n\right\}\right|.
\]
Since there is only one partition $\lambda$ having $|\lambda|=1$
(namely, $\lambda=(1)$),
this leads to parallel recursions
\begin{align*}
|\CCC_n|
&=|\Sigma_n|-
\left|\left\{ \bigsqcup_{i=1}^{n-1} \CCC_i \overset{\underline{\lambda}}{\longrightarrow} \Par: \quad
\sum_{\rho \in \CCC} d(\rho) |\underline{\lambda}(\rho)| =n\right\}\right|, \\
|\FFF_n|
&=|\Sigma_n|-
\left|\left\{ \bigsqcup_{i=1}^{n-1} \FFF_i \overset{\underline{\lambda}}{\longrightarrow} \Par: \quad
\sum_{f \in \FFF} \deg(f) |\underline{\lambda}(f)| =n\right\}\right|,
\end{align*}
and induction implies that $|\CCC_n|=|\FFF_n|$.
\end{proof}

We shall use the notation $\triv_H$ for the trivial character of a group
$H$ whenever $H$ is a finite group. This generalizes the notations
$\triv_{\Symm_n}$ and $\triv_{\Symm_\lambda}$ introduced above.

% [DG][v19] Added the above paragraph.

\begin{example}
\label{low-degree-cuspidal-labelling-example}
Taking $q=2$, let us list the sets $\FFF_n$ of monic
irreducible polynomials $f(x) \neq x$ in $\FF_2[x]$ of degree $n$
for $n \leq 3$, so that we know how many cuspidal characters
of $GL_n(\FF_q)$ in $\CCC_n$ to expect:
\begin{align*}
\FFF_1&=\{ x+1\};\\
\FFF_2&=\{ x^2+x+1 \};\\
\FFF_3&=\{ x^3+x+1, x^3+x^2+1\}.
\end{align*}
Thus we expect
\begin{enumerate}
\item[$\bullet$]
one cuspidal character of $GL_1(\FF_2)$, namely
$\rho_1(=\triv_{GL_1(\FF_2)})$,
\item[$\bullet$]
one cuspidal character $\rho_2$ of $GL_2(\FF_2)$, and
\item[$\bullet$]
two cuspidal characters $\rho_3,\rho_3'$ of $GL_3(\FF_2)$.
\end{enumerate}
We will say more about $\rho_2, \rho_3,\rho_3'$ in the next section.
\end{example}

\begin{exercise}
\label{exe.polynomials-and-necklaces}
Let $\mu : \left\{1,2,3,\ldots\right\} \to \ZZ$ denote the
\dfn{number-theoretic M\"obius function}\index{M\"obius function}\index{$\mu$},
defined by setting
$\mu(m)=(-1)^d$
if $m=p_1 \cdots p_d$ for $d$ distinct primes $p_1, p_2, \ldots, p_d$,
and $\mu(m)=0$ if $m$ is not squarefree.

\begin{itemize}

\item[(a)]
Show that for $n \geq 2$, we have
\begin{equation}
\label{primitive-necklace-number}
|\CCC_n| (=|\FFF_n|)=
\frac{1}{n} \sum_{d \mid n}
\mu\left(\frac{n}{d} \right) q^{d} .
\end{equation}
(Here, the summation
sign $\sum_{d \mid n}$ means a sum over all positive divisors
$d$ of $n$.)

% [DG][v34] Replaced the "d \text{ dividing } n" by a
% "d \mid n" in the summation subscript.

\item[(b)]
Show that \eqref{primitive-necklace-number} also counts
the \emph{necklaces}\index{necklace} with $n$ beads of $q$ colors
(= the equivalence classes under the $\ZZ/n\ZZ$-action of
cyclic rotation on sequences $(a_1,\ldots,a_n)$ in $\FF_q^n$)
which are \emph{primitive}\index{primitive necklace}
in the sense that no nontrivial rotation
fixes any of the sequences within the equivalence class.
For example, when $q=2$, here are systems of distinct
representatives of these primitive necklaces for $n=2,3,4$:
\begin{align*}
n=2:& \quad \{(0,1)\};\\
n=3:& \quad \{(0,0,1),(0,1,1)\};\\
n=4:& \quad \{(0,0,0,1),(0,0,1,1),(0,1,1,1)\}.
\end{align*}

\end{itemize}
\end{exercise}

The result of Exercise~\ref{exe.polynomials-and-necklaces}(a)
was stated by Gauss for prime $q$, and by Witt
for general $q$; it is discussed in \cite{CheboluMinac},
\cite[Section 7.6.2]{Reutenauer} and (for prime $q$)
\cite[(4.12.3)]{Granville-NTmasterclass}.
Exercise~\ref{exe.polynomials-and-necklaces}(b) is also
well-known.
See \cite[Section 7.6.2]{Reutenauer}
for a bijection explaining why the answers to both parts of
Exercise~\ref{exe.polynomials-and-necklaces} are the same.

\subsection{Steinberg's unipotent characters}

Not surprisingly, the (cuspidal) character $\iota:=\triv_{GL_1}$ of $GL_1(\FF_q)$
plays a distinguished role.  The parabolic subgroup $P_{(1^n)}$ of
$GL_n(\FF_q)$ is the Borel subgroup
$B$ of upper triangular matrices,
and we have
$
\iota^n = \Ind_B^{GL_n} \triv_B = \CC[GL_n/B]
$
(identifying representations with their characters as
usual)\footnote{\textit{Proof.} Exercise~\ref{exe.towers.indres}(d)
(applied to $G_* = GL_*$, $\ell = n$, $\alpha = \left(1^n\right)
= \left(\underbrace{1,1,\ldots,1}_{n \text{ times}}\right)$ and
$\chi_i = \iota$) gives
\[
\iota^n = \ind^n_{\left(1^n\right)} \iota^{\otimes n}
= \underbrace{\Ind^{G_n}_{P_{\left(1^n\right)}}}_{= \Ind^{GL_n}_B}
\underbrace{\Inf^{P_{\left(1^n\right)}}_{G_{\left(1^n\right)}}
   \iota^{\otimes n}}_{= \triv_{P_{\left(1^n\right)}} = \triv_B}
= \Ind^{GL_n}_B \triv_B = \CC \left[ GL_n / B \right] ,
\]
where the last equality follows from the general fact that if
$G$ is a finite group and $H$ is a subgroup of $G$, then
$\Ind^G_H \triv_H \cong \CC \left[G/H\right]$ as
$\CC G$-modules.}.
The subalgebra $A(GL)(\iota)$ of $A(GL)$ is the $\ZZ$-span of
the irreducible characters $\sigma$ that appear as constituents of
$
\iota^n = \Ind_B^{GL_n} \triv_B = \CC[GL_n/B]
$
for some $n$.

% [DG][v25] This $\iota^n = \Ind_B^{GL_n} \triv_B$ equality looks
% anything but trivial to me. I guess I see how to prove it with
% a lot of work, but shouldn't there be more details or an exercise?

% [DG][v28] Done.

\begin{definition}
An irreducible character $\sigma$ of $GL_n$ appearing as a constituent
of $\Ind_B^{GL_n} \triv_B = \CC[GL_n/B]$ is called a
\emph{unipotent character}\index{unipotent character of $GL_n$}.
Equivalently, by Frobenius reciprocity, $\sigma$ is unipotent if it contains a nonzero
$B$-invariant vector.
\end{definition}

In particular, $\triv_{GL_n}$ is a unipotent character of $GL_n$ for each $n$.

\begin{proposition}
\label{Sym-to-Steinberg-isomorphism}
One can choose $\Lambda \cong A(GL)(\iota)$ in
Theorem~\ref{uniqueness-of-sym-parts-e-f-g}(g) so that
$h_n \longmapsto \triv_{GL_n}$.
\end{proposition}
\begin{proof}
Theorem~\ref{uniqueness-of-sym-parts-a-b-c-d}(a)
tells us $\iota^2 = \Ind_B^{GL_2} \triv_B$ must have exactly two irreducible
constituents, one of which is $\triv_{GL_2}$;  call the other one $\St_2$.
Choose the isomorphism so as to send $h_2 \longmapsto \triv_{GL_2}$.
Then $h_n \mapsto \triv_{GL_n}$ follows from the
claim that $\St_2^\perp(\triv_{GL_n})=0$ for $n \geq 2$:
one has
\[
\Delta(\triv_{GL_n})
= \sum_{i+j=n} \left(\Res^{G_n}_{P_{i,j}} \triv_{GL_n}\right)^{K_{i,j}}
=\sum_{i+j=n}\triv_{GL_i} \otimes \triv_{GL_j}
\]
so that $\St_2 ^\perp(\triv_{GL_n})=(\St_2,\triv_{GL_2}) \triv_{GL_{n-2}} = 0$ since
$\St_2 \neq \triv_{GL_2}$.

\end{proof}

This subalgebra $A(GL)(\iota)$, and the unipotent characters $\chi^\lambda_q$ corresponding
under this isomorphism to the Schur functions $s_\lambda$, were introduced
by Steinberg \cite{Steinberg}.  He wrote down $\chi^\lambda_q$ as a virtual sum of induced
characters $\Ind_{P_\alpha}^{GL_n} \triv_{P_\alpha} (=\triv_{G_{\alpha_1}} \cdots \triv_{G_{\alpha_\ell}})$,
modelled on the
Jacobi-Trudi determinantal expression for $s_\lambda=\det(h_{\lambda_i-i+j})$.
Note that  $\Ind_{P_\alpha}^{GL_n} \triv_{P_\alpha}$ is the transitive
permutation representation $\CC[G/P_\alpha]$ for $GL_n$
permuting the \dfn{finite partial flag variety} $G/P_\alpha$, that is, the set
of \emph{$\alpha$-flags of subspaces}\index{flag of subspaces}
\[
\{0\} \subset V_{\alpha_1} \subset V_{\alpha_1+\alpha_2}
 \subset \cdots \subset  V_{\alpha_1+\alpha_2+\cdots+\alpha_{\ell-1}}
\subset \FF_q^n
\]
where $\dim_{\FF_q} V_d = d$ in each case.  This character has dimension equal
to $|G/P_\alpha|$, with formula given by the \dfn{$q$-multinomial coefficient}
(see e.g. Stanley \cite[\S 1.7]{Stanley}):
\[
\qbin{n}{\alpha}{q}=\frac{[n]!_q}{[\alpha_1]!_q \cdots [\alpha_\ell]!_q}
\]
where $[n]!_q:=[n]_q [n-1]_q \cdots [2]_q [1]_q$ and
$[n]_q:=1+q+\cdots+q^{n-1}=\frac{q^n-1}{q-1}$.

Our terminology $\St_2$ is motivated by the $n=2$ special case of
the \dfn{Steinberg character} $\St_n$, which is the
unipotent character corresponding under the isomorphism in
Proposition~\ref{Sym-to-Steinberg-isomorphism} to $e_n = s_{(1^n)}$.
It can be defined by the virtual sum
\[
\St_n:=
\chi^{(1^n)}_q = \sum_{\alpha} (-1)^{n-\ell(\alpha)} \Ind_{P_\alpha}^{GL_n} \triv_{P_\alpha}
\]
in which the sum runs through all compositions $\alpha$ of $n$.
This turns out to be the genuine character for
$GL_n(\FF_q)$ acting on the top homology group of its \dfn{Tits building}:
the simplicial complex whose vertices are
nonzero proper subspaces $V$ of $\FF_q^n$, and whose simplices correspond
to flags of nested subspaces.  One needs to know that this Tits building has
only top homology, so that one can deduce the above character formula from
the Hopf trace formula; see Bj\"orner \cite{Bjorner}.

\subsection{Examples: $GL_2(\FF_2)$ and $GL_3(\FF_2)$}

Let's get our hands dirty.

\begin{example}
For $n=2$, there are two unipotent characters, $\chi^{(2)}_q = \triv_{GL_2}$
and
\begin{equation}
\label{Steinberg-on-lines-description}
\St_2:=\chi^{(1,1)}_q
= \triv_{GL_1}^2 - \triv_{GL_2}
= \Ind_B^{GL_2} \triv_B -   \triv_{GL_2}
\end{equation}
since the Jacobi-Trudi formula \eqref{eq.jacobi-trudi.h} gives
$
s_{(1,1)}= \det \left[
\begin{matrix}
h_1 & h_2 \\
1  & h_1
\end{matrix} \right]
=h_1^2-h_2.
$
The description \eqref{Steinberg-on-lines-description} for
this Steinberg character $\St_2$ shows that it
has dimension
\[
|GL_2/B|-1 = (q+1)-1 = q
\]
and that one can think of it as follows:
consider the permutation action of $GL_2$ on the $q+1$ lines
$\{\ell_0,\ell_1,\ldots,\ell_q\}$
in the projective space $\Proj^1_{\FF_q}=GL_2(\FF_q)/B$, and
take the invariant subspace perpendicular to the sum of basis
elements $e_{\ell_0} + \cdots +e_{\ell_q}$.
\end{example}

\begin{example}
Continuing the previous example, but taking $q=2$,
we find that we have constructed two unipotent characters:
$\triv_{GL_2}=\chi^{(2)}_{q=2}$ of dimension $1$, and
$\St_2=\chi^{(1,1)}_{q=2}$ of dimension $q=2$.
This lets us identify the unique cuspidal character
$\rho_2$ of $GL_2(\FF_2)$, using knowledge of
the character table of $GL_2(\FF_2) \cong \Symm_3$:

\vskip.1in
\begin{tabular}{|c||c|c|c|c|}\hline
 & & & & \\
&
& $\left[ \begin{matrix} 1 & 0 \\ 0 & 1 \end{matrix} \right]$
&
$\left[ \begin{matrix} 1 & 1 \\ 0 & 1 \end{matrix} \right]$,
$\left[ \begin{matrix} 1 & 0 \\ 1 & 1 \end{matrix} \right]$,
$\left[ \begin{matrix} 0 & 1 \\ 1 & 0 \end{matrix} \right]$
&
$\left[ \begin{matrix} 1 & 1 \\ 1 & 0 \end{matrix} \right]$,
$\left[ \begin{matrix} 0 & 1 \\ 1 & 1 \end{matrix} \right]$\\
 & & & & \\\hline\hline
$\triv_{GL_2}=\chi^{(2)}_{q=2}$ &\text{unipotent}& $1$ & $1$ & $1$ \\\hline
$\St_2=\chi^{(1,1)}_{q=2}$ &\text{unipotent}& $2$ & $0$ & $-1$ \\\hline
$\rho_2$ &\text{cuspidal}& $1$ & $-1$ & $1$ \\\hline
\end{tabular}

\vskip.1in
\noindent
In other words, the cuspidal character $\rho_2$ of $GL_2(\FF_2)$
corresponds under the isomorphism $GL_2(\FF_2) \cong \Symm_3$
to the sign character $\sgn_{\Symm_3}$.
\end{example}

\begin{example}
Continuing the previous example to $q=2$ and $n=3$ lets us analyze
the irreducible characters of $GL_3(\FF_2)$.
Recalling our labelling
$\rho_1,\rho_2,\rho_3,\rho_3'$ from Example~\ref{low-degree-cuspidal-labelling-example} of the cuspidal characters of $GL_n(\FF_2)$ for $n=1,2,3$,
Zelevinsky's Theorem~\ref{Zelevinsky-decomposition-theorem}
tells us that the $GL_3(\FF_2)$-irreducible
characters should be labelled by functions
$
\{\rho_1,\rho_2,\rho_3,\rho_3'\}
\overset{\underline{\lambda}}{\longrightarrow} \Par
$
for which
\[
1 \cdot |\underline{\lambda}(\rho_1)| +
2 \cdot |\underline{\lambda}(\rho_2)| +
3 \cdot |\underline{\lambda}(\rho_3)| +
3 \cdot |\underline{\lambda}(\rho'_3)| = 3 .
\]
We will label such an irreducible character
$\chi^{\underline{\lambda}}=\chi^{(\underline{\lambda}(\rho_1),\underline{\lambda}(\rho_2),\underline{\lambda}(\rho_3),\underline{\lambda}(\rho'_3))}$.

Three of these irreducibles will be the unipotent characters,
mapping under the isomorphism
from Proposition~\ref{Sym-to-Steinberg-isomorphism} as follows:
\begin{enumerate}
\item[$\bullet$]
$
s_{(3)} =h_3 \longmapsto
\chi^{((3),\varnothing,\varnothing,\varnothing)} = \triv_{GL_3}
$
of dimension $1$.
\item[$\bullet$]
\[
s_{(2,1)}=\det \left[
\begin{matrix}
h_2 & h_3 \\
1  & h_1
\end{matrix} \right]
=h_2h_1-h_3
\longmapsto \chi^{((2,1),\varnothing,\varnothing,\varnothing)}
            =\Ind_{P_{2,1}}^{GL_3} \triv_{P_{2,1}} -   \triv_{GL_3},
\]
of dimension $\qbin{3}{2,1}{q}-\qbin{3}{3}{q}=[3]_q-1=q^2+q
\overset{q=2}{\rightsquigarrow} 6$.
\item[$\bullet$]
Lastly,
\begin{align*}
s_{(1,1,1)}&=\det \left[
\begin{matrix}
h_1 & h_2 & h_3 \\
1  & h_1 & h_2 \\
0  & 1  & h_1
\end{matrix} \right]
=h_1^3-h_2 h_1 - h_1 h_2 + h_3 \\
&\longmapsto \St_3=\chi^{((1,1,1),\varnothing,\varnothing,\varnothing)}
=\Ind_{B}^{GL_3} \triv_{B} -
\Ind_{P_{2,1}}^{GL_3} \triv_{P_{2,1}} -
\Ind_{P_{1,2}}^{GL_3} \triv_{P_{1,2}} +   \triv_{GL_3}
\end{align*}
of dimension
\begin{align*}
&\qbin{3}{1,1,1}{q}-\qbin{3}{2,1}{q}-
\qbin{3}{1,2}{q}+\qbin{3}{3}{q} \\
&=[3]!_q - [3]_q-[3]_q + 1=q^3
\overset{q=2}{\rightsquigarrow} 8.
\end{align*}
\end{enumerate}

There should also be one non-unipotent, non-cuspidal character, namely
\[
\chi^{((1),(1),\varnothing,\varnothing)}
= \rho_1 \rho_2
= \Ind_{P_{1,2}}^{GL_3} \Inf_{GL_1 \times GL_2}^{P_{1,2}}
      \left( \triv_{GL_1} \otimes \rho_2 \right)
\]
having dimension $\qbin{3}{1,2}{q} \cdot 1 \cdot 1=[3]_q
\overset{q=2}{\rightsquigarrow} 7$.

Finally, we expect cuspidal characters
$\rho_3=\chi^{(\varnothing,\varnothing,(1),\varnothing)},
\rho_3'=\chi^{(\varnothing,\varnothing,\varnothing,(1))}$,
whose dimensions $d_3, d_3'$ can be deduced from the equation
\[
1^2+6^2+8^2+7^2+d_3^2+(d_3')^2=|GL_3(\FF_2)|
=\left[(q^3-q^0)(q^3-q^1)(q^3-q^2)\right]_{q=2}=168.
\]
This forces $d_3^2+(d_3')^2=18$,
whose only solution in positive integers
is $d_3=d_3'=3$.

We can check our predictions of the dimensions
for the various  $GL_3(\FF_2)$-irreducible
characters since $GL_3(\FF_2)$ is the finite simple group of
order $168$ (also isomorphic to $PSL_2(\FF_7)$), with known
character table (see James and Liebeck \cite[p. 318]{JamesLiebeck}):\\
\begin{tabular}{|c||c||c|c|c|c|c|c|}\hline
&\text{centralizer order} & $168$ & $8$ & $4$ & $3$  & $7$ & $7$ \\\hline
 &unipotent?/cuspidal? &   &  &  &   &  &\\\hline\hline
$\triv_{GL_3}=\chi^{((3),\varnothing,\varnothing,\varnothing)}$
  &\text{unipotent}& $1$ & $1$ & $1$ & $1$ & $1$ & $1$ \\\hline
$\phantom{\St_3=}\chi^{((2,1),\varnothing,\varnothing,\varnothing)}$
  &\text{unipotent}& $6$ & $2$ & $0$ & $0$ & $-1$ & $-1$ \\\hline
$\St_3=\chi^{((1,1,1),\varnothing,\varnothing,\varnothing)}$
  &\text{unipotent}& $8$ & $0$ & $0$ & $-1$ & $1$ & $1$ \\\hline
$\phantom{\St_3=}\chi^{((1),(1),\varnothing,\varnothing)}$
  &\text{ }& $7$ & $-1$ & $-1$ & $1$ & $0$ & $0$ \\\hline
$\rho_3=\chi^{(\varnothing,\varnothing,(1),\varnothing)}$
  &\text{cuspidal}& $3$ & $-1$ & $1$ & $0$ & $\alpha$ & $\overline{\alpha}$ \\\hline
$\rho_3'=\chi^{(\varnothing,\varnothing,\varnothing,(1))}$
  &\text{cuspidal}& $3$ & $-1$ & $1$ & $0$ & $\overline{\alpha}$ & $\alpha$ \\\hline
\end{tabular}\\
Here $\alpha:=-1/2+i\sqrt{7}/2$.

\begin{remark}
It is known (see e.g. Bump \cite[Cor. 7.4]{Bump})
that, for $n \geq 2$, the dimension of any cuspidal irreducible character
$\rho$ of $GL_n(\FF_q)$ is
\[
(q^{n-1}-1)(q^{n-2}-1)\cdots(q^2-1)(q-1).
\]
Note that when $q=2$,
\begin{enumerate}
\item[$\bullet$]
for $n=2$ this gives $2^1-1=1$ for the dimension of $\rho_2$, and
\item[$\bullet$]
for $n=3$ it gives $(2^2-1)(2-1)=3$ for the dimensions of $\rho_3,\rho_3'$,
\end{enumerate}
agreeing with our calculations above.
Much more is known about the character table of $GL_n(\FF_q)$; see
Remark~\ref{Green-polynomial-remark} below,
Zelevinsky \cite[Chap. 11]{Zelevinsky}, and
Macdonald \cite[Chap. IV]{Macdonald}.
\end{remark}

\end{example}

\subsection{The Hall algebra}
\label{Hall-algebra-section}

There is another interesting Hopf subalgebra (and quotient Hopf algebra) of $A(GL)$,
related to unipotent conjugacy classes in $GL_n(\FF_q)$.

\begin{definition}
Say that an element $g$ in $GL_n(\FF_q)$ is
\emph{unipotent}\index{unipotent element of $GL_n$} if its eigenvalues are
all equal to $1$. Equivalently, $g \in GL_n(\FF_q)$ is unipotent if and only
if $g - \id_{\FF_q^n}$ is nilpotent. A conjugacy class in $GL_n(\FF_q)$ is
\emph{unipotent}\index{unipotent conjugacy class in $GL_n$}
if its elements are unipotent.

Denote by $\HHH_n$ the $\CC$-subspace of $R_\CC(GL_n)$ consisting
of those class functions which are supported only on unipotent conjugacy classes,
and let $\HHH=\bigoplus_{n \geq 0} \HHH_n$\index{$\HHH$} as a  $\CC$-subspace of
$A_\CC(GL)=\bigoplus_{n \geq 0} R_\CC(GL_n)$.
\end{definition}

% [DG][v29] Added second and third sentence of this definition.

\begin{proposition}
\label{Hall-algebra-is-Hopf-sub-quotient}
The subspace $\HHH$ is
a Hopf subalgebra of $A_\CC(GL)$, which
is graded, connected, and of finite type, and self-dual with
respect to the inner product on class functions inherited
from $A_\CC(GL)$.
It is also a quotient Hopf algebra of $A_\CC(GL)$,
as the $\CC$-linear surjection
$
A_\CC(GL) \twoheadrightarrow \HHH
$
restricting class functions to unipotent classes is a
Hopf algebra homomorphism. This surjection
has kernel $\HHH^\perp$, which is both an ideal
and a two-sided coideal.
\end{proposition}

% [DG][v29] Made it more explicit that the surjection is a Hopf
% morphism (also in the proof).

\begin{proof}
It is immediately clear that $\HHH^\perp$ is a graded
$\CC$-vector subspace of $A_\CC\left(GL\right)$, whose
$n$-th homogeneous component consists of those class
functions on $GL_n$ whose values on all unipotent
classes are $0$.
(This holds no matter whether the perpendicular space
is taken with respect to the Hermitian form
$\left(\cdot,\cdot\right)_G$ or with respect to
the bilinear form $\left<\cdot,\cdot\right>_G$.)
In other words, $\HHH^\perp$ is the kernel of the
surjection $A_\CC(GL) \twoheadrightarrow \HHH$ defined
in the proposition.

Given two class functions $\chi_i, \chi_j$ on $GL_i, GL_j$
and $g$ in $GL_{i+j}$, one has
\begin{equation}
\label{Hall-algebra-product}
\left( \chi_i \cdot \chi_j \right)(g) =
\frac{1}{|P_{i,j}|}
  \sum\limits_{\substack{h \in GL_{i+j}:\\
    h^{-1} g h=\left[ \begin{matrix} g_i & * \\0 & g_j\end{matrix} \right] \in P_{i,j}}}
  \chi_i(g_i) \chi_j(g_j).
\end{equation}
Since $g$ is unipotent if and only if $h^{-1}gh$ is unipotent if and only if both
$g_i,g_j$ are unipotent, the formula \eqref{Hall-algebra-product} shows
both that $\HHH$ is a subalgebra\footnote{Indeed,
if $\chi_i$ and $\chi_j$ are both supported only on unipotent classes,
then the same holds for $\chi_i \cdot \chi_j$.
% [DG][v29] The previous sentence was an "if and only if",
% but the converse is false if one of $\chi_i$ and
% $\chi_j$ is zero, and I am not sure if it is true otherwise.
% (It is also not needed.)
}
and that $\HHH^\perp$ is a two-sided ideal\footnote{In
fact, if one of $\chi_i$ and $\chi_j$ annihilates all
unipotent classes, then so does $\chi_i \cdot \chi_j$.}.
It also shows that
the surjection $A_\CC(GL) \twoheadrightarrow \HHH$
restricting every class function to unipotent classes is an
algebra homomorphism\footnote{because if $g$ is unipotent, then the
only values of $\chi_i$ and $\chi_j$ appearing on the
right hand side of \eqref{Hall-algebra-product} are those
on unipotent elements}.

Similarly, for class functions $\chi$ on $GL_n$ and $(g_i, g_j)$ in $GL_{i,j}=GL_i \times GL_j$, one has
\[
\Delta(\chi)(g_i, g_j) = \frac{1}{q^{ij}} \sum_{k \in \FF_q^{i \times j}}
\chi\left[ \begin{matrix} g_i & k \\0 & g_j\end{matrix} \right]
\]
using \eqref{deflation-character-formula-v2}.
This shows both that $\HHH$ is a sub-coalgebra
of $A=A_\CC(GL)$ (that is, it satisfies
$\Delta \HHH \subset \HHH \otimes \HHH$)
and that $\HHH^\perp$ is a two-sided coideal
(that is, we have
$\Delta(\HHH^\perp) \subset \HHH^\perp \otimes A + A \otimes \HHH^\perp$),
since it shows that if
$\chi$ is supported only on unipotent classes,
then $\Delta(\chi)$ vanishes on
$(g_1,g_2)$ that have either $g_1$ or $g_2$
non-unipotent.
% [DG][v29] Again, removed one direction of an "if and only if".
% The other direction was true (because $\Delta(\chi)$ has
% a $GL_0 \times GL_n$-component which basically is a copy of
% $\chi$), but it was again useless and only confusing.
It also
shows that the surjection $A_\CC(GL) \twoheadrightarrow \HHH$
restricting every class function to unipotent classes is a
coalgebra homomorphism.
The rest follows.
\end{proof}

The subspace $\HHH$ is called the \dfn{Hall algebra}. It has
an obvious orthogonal $\CC$-basis, with interesting structure constants.

\begin{definition}
\label{def.hall-alg.hall-coeffs}
Given a partition $\lambda$ of $n$, let
$J_{\lambda}$ denote the $GL_n$-conjugacy class of
unipotent matrices whose \dfn{Jordan type} (that is, the list of the
sizes of the Jordan blocks, in decreasing order) is given by
$\lambda$.
Furthermore, let $z_\lambda(q)$ denote the size of the centralizer
of any element of this conjugacy class $J_{\lambda}$.

The indicator class functions\footnote{Here we use the following notation:
Whenever $P$ is a subset of a group $G$, we denote by $\triv_P$ the map
$G \to \CC$ which sends every element of $P$ to $1$ and all remaining
elements of $G$ to $0$. This is not in conflict with the notation
$\triv_G$ for the trivial character of $G$, since $\triv_P = \triv_G$
for $P = G$. Note that $\triv_P$ is a class function when $P$ is a union
of conjugacy classes of $G$.}
 $\{ \triv_{J_\lambda} \}_{\lambda \in \Par}$
form a $\CC$-basis for $\HHH$ whose multiplicative structure constants
are called the \dfn{Hall coefficients} $g^{\lambda}_{\mu,\nu}(q)$:
\[
\triv_{J_\mu} \triv_{J_\nu} = \sum_\lambda g^{\lambda}_{\mu,\nu}(q) \,\, \triv_{J_\lambda}.
\]
\end{definition}

% [DG][v19] Added the above footnote.

\noindent
Because the dual basis to $\{\triv_{J_\lambda}\}$ is
$\{ z_\lambda(q) \triv_{J_\lambda} \}$,
self-duality of $\HHH$ shows that the Hall coefficients are
(essentially) also structure constants for the comultiplication:
\[
\Delta \triv_{J_\lambda}
= \sum_{\mu,\nu} g^{\lambda}_{\mu,\nu}(q)
\frac{z_{\mu}(q) z_{\nu}(q)}{z_{\lambda}(q)} \cdot
\triv_{J_\mu} \otimes \triv_{J_\nu}.
\]

% [DG][v29] Fixed the definition of $z_\lambda(q)$ and the
% formulas referencing it.

The Hall coefficient $g^{\lambda}_{\mu,\nu}(q)$ has the following interpretation.

\begin{proposition}
\label{prop.hall-alg.product}
Fix any $g$ in $GL_n(\FF_q)$
acting unipotently on $\FF_q^n$ with Jordan type $\lambda$.
Then $g^{\lambda}_{\mu,\nu}(q)$ counts the $g$-stable $\FF_q$-subspaces $V \subset \FF_q^n$
for which the restriction
$g|V$ acts with Jordan type $\mu$, and the induced map $\bar{g}$
on the quotient space $\FF_q^n/V$ has Jordan type $\nu$.
\end{proposition}
\begin{proof}
Given $\mu,\nu$ partitions of $i,j$ with $i+j=n$, taking $\chi_i, \chi_j$ equal
to $\triv_{J_\mu}, \triv_{J_\nu}$ in
\eqref{Hall-algebra-product}
shows that for any $g$ in $GL_n$, the value of
$\left( \triv_{J_\mu} \cdot \triv_{J_{\nu}} \right)(g)$ is given by
\begin{equation}
\label{Hall-coefficient-crux}
\frac{1}{|P_{i,j}|}
  \left|\left\{ h \in GL_n: h^{-1} g h=\left[ \begin{matrix} g_i & * \\0 & g_j\end{matrix} \right] \text{ with }g_i \in J_{\mu},g_j \in J_{\nu} \right\}\right| .
\end{equation}
Let $S$ denote the set appearing in \eqref{Hall-coefficient-crux},
and let $\FF_q^i$ denote the $i$-dimensional
subspace of $\FF_q^n$ spanned by the first $i$ standard basis vectors.
Note that the condition on an element $h$ in $S$ saying that $h^{-1}gh$ is in block upper-triangular form can be re-expressed by saying that the subspace $V:=h(\FF_q^i)$ is $g$-stable.
One then sees that the map $h \overset{\varphi}{\longmapsto} V=h(\FF_q^i)$ surjects
$S$ onto the set of $i$-dimensional $g$-stable subspaces $V$ of $\FF_q^n$ for which
$g|V$ and $\bar{g}$ are unipotent of types $\mu,\nu$, respectively.  Furthermore, for
any particular such $V$,
its fiber $\varphi^{-1}(V)$ in $S$ is a coset of the stabilizer within
% [DG][v29] Added "a coset of" on the previous line.
$GL_n$ of $V$, which is conjugate to $P_{i,j}$, and hence has cardinality
$|\varphi^{-1}(V)|=|P_{i,j}|$.  This proves the assertion of the proposition.
\end{proof}

The Hall algebra $\HHH$ will turn out to be isomorphic to the ring $\Lambda_\CC$ of symmetric
functions with $\CC$ coefficients, via a composite $\varphi$ of three maps
\[
\Lambda_\CC \longrightarrow
A(GL)(\iota)_\CC \longrightarrow
A(GL)_\CC \longrightarrow
\HHH
\]
in which the first map is the isomorphism from
Proposition~\ref{Sym-to-Steinberg-isomorphism}, the second is inclusion,
and the third is the quotient map from
Proposition~\ref{Hall-algebra-is-Hopf-sub-quotient}.

\begin{theorem}
\label{Sym-to-Hall-isomorphism-thm}
The above composite $\varphi$ is a Hopf algebra isomorphism, sending
\[
\begin{array}{rcl}
h_n &\longmapsto & \sum_{\lambda \in \Par_n} \triv_{J_\lambda}, \\
e_n &\longmapsto & q^{\binom{n}{2}} \triv_{J_{(1^n)}}, \\
p_n &\longmapsto & \sum_{\lambda \in \Par_n} (q;q)_{\ell(\lambda) - 1}\triv_{J_\lambda} \qquad \left(\text{for } n > 0\right),
\end{array}
\]
where we are using the notation
\[
(x;q)_{m}
:= (1-x)(1-qx)(1-q^2x) \cdots (1-q^{m-1}x)
  \qquad \text{ for all } m \in \NN \text{ and } x \text{ in any ring}.
\]
\end{theorem}

\begin{proof}
That $\varphi$ is a graded Hopf morphism follows because it is a composite
of three such morphisms.  We claim that once one shows the formula for the
(nonzero) image $\varphi(p_n)$ given above is correct,
then this will already show $\varphi$ is an isomorphism,
by the following argument.  Note first that $\Lambda_\CC$ and $\HHH$
both have dimension $|\Par_n|$ for their $n$-th homogeneous
components, so it suffices to show that the graded map $\varphi$ is injective.
On the other hand, both $\Lambda_\CC$ and $\HHH$
are (graded, connected, finite type) \emph{self-dual} Hopf algebras
(although with respect to a sesquilinear form),
so Theorem~\ref{general-structure-theorem} says that each is the symmetric algebra
on its space of primitive elements.  Thus it suffices to check that
$\varphi$ is injective when restricted to their subspaces of
primitives.\footnote{An alternative way to see that it suffices to
check this is by recalling Exercise~\ref{exe.I=0}(c).}
For $\Lambda_\CC$, by Corollary~\ref{primitives-in-Lambda} the primitives
are spanned by $\{p_1,p_2,\ldots\}$, with only
one basis element in each degree $n\geq 1$.  Hence
$\varphi$ is injective on the subspace of primitives if and only if
it does not annihilate any $p_n$.

% [DG][v14] There is a subtlety here that even I don't want to deal
% with: The Hall algebra $\HHH$ is self-dual with respect to a
% *Hermitian*, not bilinear, form; this isn't something you have
% defined. I guess the Zelevinsky results still work in this case,
% and I also guess it is easy to identify a $\QQ$-form of $A(GL)_\CC$
% in which the elements $h_n$, $e_n$, $p_n$ of $\HHH$ are contained,
% and thus get around this distinction.
% Anyway I don't care much, because Exercise~\ref{exe.I=0} makes it
% possible to completely avoid using PSH theory here. Would you
% object to me adding a reference to this exercise as an alternative
% proof?

% [DG][v23] Added footnote. I don't really like either argument, to be
% honest. What about the following alternative, which I suggest adding
% as an exercise: check that
% $\varphi(e_\lambda)$ equals $(nonzero scalar) \triv_{J_\lambda^t}$ plus
% dominance-smaller terms? (Whence the map is unitriangular.)

% [DG][v29] Added the exercise which follows just after this proof.

Thus it only remains to show the above formulas for the images of
$h_n,e_n,p_n$ under $\varphi$.
This is clear for $h_n$, since
Proposition~\ref{Sym-to-Steinberg-isomorphism} shows
that it maps under the first two composites to the indicator function
$\triv_{GL_n}$ which then restricts to the sum of indicators
$\sum_{\lambda \in \Par_n} \triv_{J_\lambda}$ in $\HHH$.
For $e_n, p_n$, we resort to generating functions.
Let $\tilde{h}_n,\tilde{e}_n,\tilde{p}_n$ denote the
three putative images in $\HHH$ of $h_n,e_n,p_n$,
appearing on the right side in the theorem, and define
generating functions
\[
\tilde{H}(t):=\sum_{n \geq 0} \tilde{h}_n t^n, \quad
\tilde{E}(t):=\sum_{n \geq 0} \tilde{e}_n t^n, \quad
\tilde{P}(t):=\sum_{n \geq 0} \tilde{p}_{n+1} t^n
\qquad \qquad \text{ in } \HHH[[t]] .
\]
We wish to show that the map
$\varphi[[t]] : \Lambda_\CC [[t]] \to \HHH[[t]]$
(induced by $\varphi$) maps $H(t),E(t),P(t)$ in $\Lambda[[t]]$
to these three generating functions\footnote{See
\eqref{h-generating-function},
\eqref{e-generating-function},
\eqref{E(t)-H(t)-P(t)-relation} for the definitions
of $H(t),E(t),P(t)$.}.  Since we have already shown this is
correct for $H(t)$, by \eqref{E(t)-H(t)-relation},
\eqref{E(t)-H(t)-P(t)-relation}, it suffices to check that in $\HHH[[t]]$ one has
\[
\begin{array}{rccl}
\tilde{H}(t) \tilde{E}(-t) = 1,& \text{ or equivalently, }
  &\sum_{k=0}^n (-1)^k \tilde{e}_k \tilde{h}_{n-k} = \delta_{0,n} ;\\
\tilde{H}'(t) \tilde{E}(-t) = \tilde{P}(t),& \text{  or equivalently, }
  &\sum_{k=0}^n (-1)^k (n-k) \tilde{e}_k \tilde{h}_{n-k} = \tilde{p}_n.
\end{array}
\]
Thus it would be helpful to evaluate the class function
$\tilde{e}_k \tilde{h}_{n-k}$.  Note that a unipotent $g$ in $GL_n$
having $\ell$ Jordan blocks has an $\ell$-dimensional $1$-eigenspace, so that
the number of $k$-dimensional $g$-stable $\FF_q$-subspaces of $\FF_q^n$ on
which $g$ has Jordan type $\left(1^k\right)$ (that is, on which $g$
acts as the identity) is the \dfn{$q$-binomial coefficient}
\[
\qbin{\ell}{k}{q} = \frac{(q;q)_\ell}{(q;q)_k (q;q)_{\ell-k}} ,
\]
counting $k$-dimensional $\FF_q$-subspaces $V$ of
an $\ell$-dimensional $\FF_q$-vector space; see, e.g., \cite[\S 1.7]{Stanley}.
Hence, for a unipotent $g$ in $GL_n$
having $\ell$ Jordan blocks, we have
\[
(\tilde{e}_k \tilde{h}_{n-k})(g)
= q^{\binom{k}{2}} \cdot \left(\triv_{J_{(1^k)}} \cdot \tilde{h}_{n-k}\right) (g)
= q^{\binom{k}{2}} \cdot \sum_{\nu \in \Par_{n-k}} \left(\triv_{J_{(1^k)}} \cdot \triv_{J_\nu}\right) (g)
=  q^{\binom{k}{2}} \qbin{\ell}{k}{q}
\]
(by Proposition \ref{prop.hall-alg.product}).
Thus one needs for $\ell \geq 1$ that
\begin{align}
\sum_{k=0}^\ell (-1)^k
q^{\binom{k}{2}} \qbin{\ell}{k}{q}
&= 0,\label{first-q-binomial-identity}\\
\sum_{k=0}^\ell (-1)^k (n-k)
q^{\binom{k}{2}} \qbin{\ell}{k}{q}
&= (q;q)_{\ell - 1}.\label{second-q-binomial-identity}
\end{align}
Identity \eqref{first-q-binomial-identity}
comes from setting $x=1$ in the \dfn{$q$-binomial theorem} \cite[Exer. 3.119]{Stanley}:
\begin{equation}
\sum_{k=0}^\ell (-1)^k q^{\binom{k}{2}} \qbin{\ell}{k}{q} x^{\ell-k}
=(x-1)(x-q)(x-q^2)\cdots (x-q^{\ell-1}).
\label{q-binomial-theorem}
\end{equation}
Identity \eqref{second-q-binomial-identity}
comes from applying $\frac{d}{dx}$ to \eqref{q-binomial-theorem},
then setting $x=1$, and finally
adding $(n-\ell)$ times \eqref{first-q-binomial-identity}.
\end{proof}

% [DG][v29] Added some more detail.

\begin{exercise}
\label{exe.hall-alg.dominance}
Fix a prime power $q$.
For any $k \in \NN$, and any $k$ partitions
$\lambda^{(1)}, \lambda^{(2)}, \ldots, \lambda^{(k)}$, we define
a family
$\left(g_{\lambda^{(1)}, \lambda^{(2)}, \ldots, \lambda^{(k)}}^\lambda \left(q\right) \right) _ {\lambda \in \Par}$
of elements of $\CC$ by the equation
\[
\triv_{J_{\lambda^{(1)}}} \triv_{J_{\lambda^{(2)}}} \cdots \triv_{J_{\lambda^{(k)}}}
= \sum_{\lambda \in \Par} g_{\lambda^{(1)}, \lambda^{(2)}, \ldots, \lambda^{(k)}}^\lambda \left(q\right) \triv_{J_\lambda}
\]
in $\HHH$.
This notation generalizes the notation $g^\lambda_{\mu, \nu} \left(q\right)$ we introduced in Definition~\ref{def.hall-alg.hall-coeffs}. Note that $g^\lambda_\mu \left(q\right) = \delta_{\lambda, \mu}$ for any two partitions $\lambda$ and $\mu$, and that $g^\lambda \left(q\right) = \delta_{\lambda, \varnothing}$ for any partition $\lambda$ (where $g^\lambda \left(q\right)$ is to be understood as $g_{\lambda^{(1)}, \lambda^{(2)}, \ldots, \lambda^{(k)}}^\lambda \left(q\right)$ for $k = 0$).

\begin{itemize}

\item[(a)] Let $\lambda \in \Par$, and let $n = \left| \lambda \right|$. Let $V$ be an $n$-dimensional $\FF_q$-vector space, and let $g$ be a unipotent endomorphism of $V$ having Jordan type $\lambda$. Let $k \in \NN$, and let $\lambda^{(1)}, \lambda^{(2)}, \ldots, \lambda^{(k)}$ be $k$ partitions. A \emph{$\left(\lambda^{(1)}, \lambda^{(2)}, \ldots, \lambda^{(k)}\right)$-compatible $g$-flag}\index{compatible flag} will mean a sequence $0 = V_0 \subset V_1 \subset V_2 \subset \cdots \subset V_k = V$ of $g$-invariant $\FF_q$-vector subspaces $V_i$ of $V$ such that for every $i \in \left\{1, 2, \ldots, k\right\}$, the endomorphism of $V_i / V_{i-1}$ induced by $g$ \ \ \ \ \footnote{This is well-defined. In fact, both $V_i$ and $V_{i-1}$ are $g$-invariant, so that $g$ restricts to an endomorphism of $V_i$, which further restricts to an endomorphism of $V_{i-1}$, and thus gives rise to an endomorphism of $V_i / V_{i-1}$.} has Jordan type $\lambda^{(i)}$.

Show that
$g_{\lambda^{(1)}, \lambda^{(2)}, \ldots, \lambda^{(k)}}^\lambda \left(q\right)$
is the number of $\left(\lambda^{(1)}, \lambda^{(2)}, \ldots, \lambda^{(k)}\right)$-compatible $g$-flags.\footnote{This can be seen as a generalization of Proposition~\ref{prop.hall-alg.product}. In fact, if $\mu$ and $\nu$ are two partitions, then a $\left(\mu, \nu\right)$-compatible $g$-flag is a sequence $0 = V_0 \subset V_1 \subset V_2 = V$ of $g$-invariant $\FF_q$-vector subspaces $V_i$ of $V$ such that the endomorphism of $V_1 / V_0 \cong V_1$ induced by $g$ has Jordan type $\mu$, and the endomorphism of $V_2 / V_1 \cong V / V_1$ induced by $g$ has Jordan type $\nu$. Choosing such a sequence amounts to choosing $V_1$ (since there is only one choice for each of $V_0$ and $V_2$), and the conditions on this $V_1$ are precisely the conditions on $V$ in Proposition~\ref{prop.hall-alg.product}.}

\item[(b)] Let $\lambda \in \Par$. Let $k \in \NN$, and let $\lambda^{(1)}, \lambda^{(2)}, \ldots, \lambda^{(k)}$ be $k$ partitions. Show that
$g_{\lambda^{(1)}, \lambda^{(2)}, \ldots, \lambda^{(k)}}^\lambda \left(q\right) = 0$
unless
$\left| \lambda^{(1)} \right| + \left| \lambda^{(2)} \right| + \cdots + \left| \lambda^{(k)} \right| = \left| \lambda \right|$
and
$\lambda^{(1)} + \lambda^{(2)} + \cdots + \lambda^{(k)} \triangleright \lambda$.
(Here and in the following, we are using the notations of Exercise~\ref{exe.liri.simple}).

\item[(c)] Let $\lambda \in \Par$, and let us write the transpose partition $\lambda^t$ as $\lambda^t = \left(\left(\lambda^t\right)_1, \left(\lambda^t\right)_2, \ldots, \left(\lambda^t\right)_\ell\right)$. Show that
$g_{\left(1^{\left(\lambda^t\right)_1}\right), \left(1^{\left(\lambda^t\right)_2}\right), \ldots, \left(1^{\left(\lambda^t\right)_\ell}\right)}^\lambda \left(q\right) \neq 0$.

\item[(d)] Let $n \in \NN$ and $\lambda \in \Par_n$. Show that
\[
\varphi\left(e_\lambda\right)
= \sum_{\mu \in \Par_n ; \  \lambda^t \triangleright \mu}
  \alpha_{\lambda, \mu} \triv_{J_\mu}
\]
for some coefficients $\alpha_{\lambda, \mu} \in \CC$ satisfying $\alpha_{\lambda, \lambda^t} \neq 0$.

\item[(e)] Give another proof of the fact that the map $\varphi$ is injective.

\end{itemize}

[\textbf{Hint:} For (b), use Exercise~\ref{exe.liri.jordan}(b).]

\end{exercise}

% [DG][v31] Added the hint (which I had intended to add in v29 but forgot). This
% should clear up the issue from your last email?

We next indicate, without proof, how $\HHH$ relates to the
classical Hall algebra.

\begin{definition}
Let $p$ be a prime.
The usual \dfn{Hall algebra}, or what Schiffmann \cite[\S 2.3]{Schiffmann}
calls \dfn{Steinitz's classical Hall algebra}
(see also Macdonald \cite[Chap. II]{Macdonald}),
has $\ZZ$-basis elements
$\{ u_\lambda \}_{\lambda \in \Par}$, with the multiplicative structure constants
$g^{\lambda}_{\mu,\nu}(p)$ in
\[
u_\mu u_\nu = \sum_\lambda
g^{\lambda}_{\mu,\nu}(p)\,\, u_\lambda
\]
defined as follows: fix a finite abelian $p$-group $L$ of
\emph{type $\lambda$}\index{type of an abelian $p$-group},
meaning that
\[
L \cong \bigoplus_{i=1}^{\ell(\lambda)} \ZZ/p^{\lambda_i}\ZZ,
\]
and let $g^{\lambda}_{\mu,\nu}(p)$ be the number of subgroups $M$ of $L$ of
type $\mu$, for which the quotient $N:=L/M$ is of type $\nu$.  In other words,
$g^{\lambda}_{\mu,\nu}(p)$ counts, for a fixed abelian
$p$-group $L$ of type $\lambda$,  the number of short exact sequences
$
0 \rightarrow M
\rightarrow L
\rightarrow N
\rightarrow 0
$
in which $M,N$ have types $\mu,\nu$, respectively (modulo isomorphism
of short exact sequences restricting to the identity on $L$).
\end{definition}

% [DG][v23] Added "Let $p$ be a prime" to the above definition -- or was
% this in the context?

We claim that when one takes the finite field $\FF_q$ of order
$q=p$ a \emph{prime}, the $\ZZ$-linear map
\begin{equation}
\label{classical-Hall-to-Zelevinsky-Hall}
u_\lambda \longmapsto \triv_{J_\lambda}
\end{equation}
gives an isomorphism from this classical Hall algebra
to the $\ZZ$-algebra $\HHH_\ZZ \subset \HHH$.
The key point is \dfn{Hall's Theorem}, a non-obvious
statement for which Macdonald includes two proofs in
\cite[Chap. II]{Macdonald}, one of them due to
Zelevinsky\footnote{See also \cite[Thm. 2.6, Prop. 2.7]{Schiffmann}
for quick proofs of part of it, similar to Zelevinsky's.
Another proof, based on a recent category-theoretical
paradigm, can be found in \cite[Theorem 3.53]{Dyckerhoff}.}.
To state it, we first recall some notions about discrete valuation rings.

% [DG][v52] Added {Dyckerhoff} reference. This looks fascinating and
% readable -- I'll probably start reading it soon.

\begin{definition}
A \dfn{discrete valuation ring} (short \dfn{DVR}) $\ooo$ is
a principal ideal domain having only one maximal ideal $\mmm \neq 0$,
with quotient $k=\ooo/\mmm$ called its \dfn{residue field}.

% [DG][v29] Added nonzeroness condition on $\mmm$.

The structure theorem for finitely generated modules
over a PID implies that an $\ooo$-module $L$
with finite composition series of composition length $n$
must have $L \cong \bigoplus_{i=1}^{\ell(\lambda)} \ooo/\mmm^{\lambda_i}$
for some partition $\lambda$ of $n$; say $L$ has
\emph{type}\index{type of a module over a DVR} $\lambda$
in this situation.
\end{definition}

Here are the two crucial examples for us.

\begin{example}
\label{DVR-example-F[[t]]}
For any field $\FF$, the power series ring $\ooo=\FF[[t]]$ is a DVR
with maximal ideal $\mmm=(t)$ and
residue field $k=\ooo/\mmm=\FF[[t]]/(t)\cong\FF$.
An $\ooo$-module $L$ of type $\lambda$ is an $\FF$-vector space together
with an $\FF$-linear transformation $T \in \End L$ that acts
on $L$ nilpotently
(so that $g:=T+1$ acts unipotently, where $1 = \id_L$) with
Jordan blocks of sizes given by $\lambda$:  each summand
$\ooo/\mmm^{\lambda_i}=\FF[[t]]/(t^{\lambda_i})$ of $L$ has
an $\FF$-basis $\{1,t,t^2,\ldots,t^{\lambda_i-1}\}$ on
which the map $T$ that multiplies by $t$ acts as a nilpotent Jordan block of
size $\lambda_i$.
Note also that, in this setting, $\ooo$-submodules are the same as
$T$-stable (or $g$-stable) $\FF$-subspaces.
\end{example}

\begin{example}
\label{DVR-example-Z_p}
The ring of $p$-adic integers $\ooo=\ZZ_p$ is a DVR with
maximal ideal $\mmm=(p)$ and
residue field $k=\ooo/\mmm=\ZZ_p/p\ZZ_p \cong \ZZ/p\ZZ$.
An $\ooo$-module $L$ of type $\lambda$ is an abelian $p$-group of
type $\lambda$:  for each summand,
$\ooo/\mmm^{\lambda_i}=\ZZ_p/p^{\lambda_i}\ZZ_p \cong \ZZ/p^{\lambda_i}\ZZ$.
Note also that, in this setting, $\ooo$-submodules are the same as
subgroups.
\end{example}

\noindent
One last notation: $n(\lambda):=\sum_{i \geq 1} (i-1)\lambda_i$,
for $\lambda$ in $\Par$.  Hall's Theorem is as follows.

\begin{theorem}
\label{Hall's-theorem}
Assume $\ooo$ is a DVR with maximal ideal $\mmm$, and that its residue field
$k=\ooo/\mmm$ is finite of cardinality $q$.
Fix an $\ooo$-module $L$ of type $\lambda$.
Then the number of $\ooo$-submodules $M$ of
type $\mu$ for which the quotient $N=L/M$ is of type $\nu$ can be
written as the specialization
\[
[g^{\lambda}_{\mu,\nu}(t)]_{t=q}
\]
of a polynomial $g^{\lambda}_{\mu,\nu}(t)$ in $\ZZ[t]$, called the
\dfn{Hall polynomial}.

Furthermore, the Hall polynomial $g^{\lambda}_{\mu,\nu}(t)$
has degree at most $n(\lambda)-(n(\mu)+n(\nu))$, and its
coefficient of $t^{n(\lambda)-(n(\mu)+n(\nu))}$ is the Littlewood-Richardson
coefficient $c^{\lambda}_{\mu,\nu}$.
\end{theorem}

\noindent
Comparing what Hall's Theorem says in
Examples~\ref{DVR-example-F[[t]]}
and \ref{DVR-example-Z_p}, shows that
the map \eqref{classical-Hall-to-Zelevinsky-Hall}
gives the desired isomorphism from the classical Hall algebra
to $\HHH_\ZZ$.

We close this section with some remarks on the vast literature
on Hall algebras that we will \emph{not} discuss here.

\begin{remark}
Macdonald's version of Hall's Theorem \cite[(4.3)]{Macdonald} is stronger
than Theorem~\ref{Hall's-theorem}, and useful for certain applications:
he shows that $g^{\lambda}_{\mu,\nu}(t)$ is the zero polynomial whenever
the Littlewood-Richardson coefficient $c^{\lambda}_{\mu,\nu}$ is zero.
\end{remark}

\begin{remark}
In general, not all coefficients of the Hall polynomials
$g^\lambda_{\mu, \nu}(t)$ are nonnegative (see Butler/Hales
\cite{ButlerHales} for a study of when they are); it often
happens that $g^\lambda_{\mu, \nu}(1) = 0$ despite
$g^\lambda_{\mu, \nu}(t)$ not being the zero
polynomial\footnote{Actually, Butler/Hales show in
\cite[proof of Prop. 2.4]{ButlerHales} that the values
$g^\lambda_{\mu, \nu}(1)$ are the structure constants of
the ring $\Lambda$ with respect to its basis
$\left(m_\lambda\right)_{\lambda \in \Par}$: we have
\[
m_\mu m_\nu = \sum_{\lambda \in \Par} g^\lambda_{\mu, \nu} (1) m_\lambda
\]
for all partitions $\mu$ and $\nu$.}. However,
in \cite[Thm. 4.2]{Klein-hall}, Klein showed that the polynomial
values $g^\lambda_{\mu, \nu}\left(p\right)$ for $p$ prime
are always positive when $c^\lambda_{\mu, \nu} \neq 0$. (This
easily yields the same result for $p$ a prime power.)
\end{remark}

% [DG][v29] Added above remark.

% [DG][v32] Added footnote. I suspect that the unexpected
% proneness of $g^\lambda_{\mu, \nu}(1)$ to vanish is the reason
% why the converse of Exercise~\ref{exe.liri.jordan}(b) cannot
% be proven by a construction of the kind "consider the free
% vector space with basis the set of cells of $\lambda$, and map
% every cell to some other cell or $0$" (viz., the converse of
% Exercise~\ref{exe.liri.jordan}(b) does not hold over the
% "field" $\FF_1$). I have not checked this in detail.

\begin{remark}
\label{Green-polynomial-remark}
Zelevinsky in \cite[Chaps 10, 11]{Zelevinsky} uses
the isomorphism $\Lambda_\CC \rightarrow \HHH$ to derive J. Green's formula
for the value of any irreducible character $\chi$ of $GL_n$ on any unipotent
class $J_\lambda$.  The answer involves values of irreducible characters
of $\Symm_n$ along with \emph{Green's polynomials}
$Q^\lambda_\mu(q)$
(see Macdonald \cite[\S III.7]{Macdonald};  they are
denoted $Q(\lambda,\mu)$ by Zelevinsky),
which express the images under the isomorphism
of Theorem~\ref{Sym-to-Hall-isomorphism-thm}
of the symmetric function basis $\{p_\mu\}$ in terms of the
basis $\{\triv_{J_\lambda}\}$.
\end{remark}

\begin{remark}
The Hall polynomials $g^{\lambda}_{\mu,\nu}(t)$ also essentially give the
multiplicative structure constants for $\Lambda(\xx)[t]$ with
respect to its basis of \emph{Hall-Littlewood symmetric functions}
$P_\lambda=P_{\lambda}(\xx;t)$:
\[
P_{\mu} P_{\nu}
= \sum_{\lambda} t^{n(\lambda)-(n(\mu)+n(\nu))} g^{\lambda}_{\mu,\nu}(t^{-1}) P_{\lambda}.
\]
See Macdonald \cite[\S III.3]{Macdonald}.
\end{remark}

\begin{remark}
Schiffmann \cite{Schiffmann} discusses self-dual Hopf algebras which
vastly generalize the classical Hall algebra called \emph{Ringel-Hall algebras},
associated to abelian categories which are hereditary.  Examples
come from categories of nilpotent representations of quivers;
the quiver having exactly one node and one arc recovers the
classical Hall algebra $\HHH_\ZZ$ discussed above.
\end{remark}

\begin{remark}
The general linear groups $GL_n(\FF_q)$ are one of four families
of so-called \emph{classical groups}. Progress has been made on
extending Zelevinsky's PSH theory to the other families:

(a) Work of Thiem and Vinroot \cite{ThiemVinroot} shows that
the tower $\{G_*\}$ of
\emph{finite unitary groups} $U_n(\FF_{q^2})$
give rise to another positive self-dual Hopf algebra
$A=\bigoplus_{n \geq 0}R( U_n(\FF_{q^2}) )$,
in which the role of Harish-Chandra induction
is played by \emph{Deligne-Lusztig induction}.  In this theory,
character and degree formulas for $U_n(\FF_{q^2})$ are related
to those of $GL_n(\FF_q)$ by substituting $q \mapsto -q$, along
with appropriate scalings by $\pm 1$, a phenomenon sometimes called
\emph{Ennola duality}.  See also \cite[\S 4]{RStantonThiem}.

(b) van Leeuwen \cite{Leeuwen-hopf}
has studied $\bigoplus_{n\geq0}R\left(  Sp_{2n}\left(
\FF_{q}\right)  \right)  $, $\bigoplus_{n\geq0}R\left(
O_{2n}\left(  \FF_{q}\right)  \right)  $ and $\bigoplus_{n\geq
0}R\left(  U_{n}\left(  \FF_{q^{2}}\right)  \right)  $ not as Hopf
algebras, but rather as so-called \emph{twisted PSH-modules} over the
PSH $A(GL)$ (a ``deformed'' version of the older notion of
Hopf modules). He classified these PSH-modules axiomatically
similarly to Zelevinsky's above classification of PSH's.

(c) In a recent honors thesis \cite{ShelleyAbrahamson-honors},
Shelley-Abrahamson defined yet another variation of the concept of
Hopf modules, named \emph{$2$-compatible Hopf modules}, and
identified $\bigoplus_{n\geq0}R\left(  Sp_{2n}\left(
\FF_{q}\right)  \right)  $ and $\bigoplus_{n\geq0}R\left(
O_{2n+1}\left(  \FF_{q}\right)  \right)  $ as such modules over
$A(GL)$.
\end{remark}

% [DG][v3] Here's hoping I haven't misrepresented either of (b) and (c)
% -- I have not really read any of them... I am wondering if
% (b) and (c) are related, particularly for the $Sp_{2n}$ family
% they have in common. The definitions of van Leeuwen's and
% Shelley-Abrahamson's Hopf modules don't completely seem to
% match, although they are similar.

\newpage

%%%%%%%%%%%%%%%%%%%%%%%%%%%%%%%%%%%%%
\section{Quasisymmetric functions and $P$-partitions}
\label{Qsym-section}
%%%%%%%%%%%%%%%%%%%%%%%%%%%%%%%%%%%%%

We discuss here our next important example of a Hopf algebra arising in
combinatorics:  the \emph{quasisymmetric functions} of Gessel \cite{Gessel},
with roots in work of Stanley \cite{Stanley-thesis} on $P$-partitions.
Other treatments of quasisymmetric functions
can be found in \cite[Section 7.19]{Stanley} and \cite[Chapter 8]{Sagan2020}
(with focus on their enumerative applications rather than on their Hopf
structure) and in \cite[Chapter 6]{Meliot} (with a focus on their
representation-theoretical meaning).
Quasisymmetric functions have found applications in combinatorial
enumeration (\cite[Chapter 8]{Sagan2020}, \cite[Section 7.19]{Stanley}),
topology (\cite{BakerRichter}) and algebraic geometry
(\cite{Morava}, \cite{Oesinghaus}).

\subsection{Definitions, and Hopf structure}

The definitions of quasisymmetric functions require a totally ordered variable set.
Usually we will use a variable set denoted $\xx=(x_1,x_2,\ldots)$ with the
usual ordering $x_1 < x_2 < \cdots$.  However, it is good to have some flexibility
in changing the ordering, which is why we make the following definition.

\begin{definition}
Given any totally ordered set $I$, create a totally ordered
variable set $\{x_i\}_{i\in I}$, and then let $R(\{x_i\}_{i\in I})$
denote the power series
of bounded degree in $\{x_i\}_{i \in I}$
having coefficients in $\kk$.

The \dfn{ring of quasisymmetric functions}\index{quasisymmetric function}
\dfn{$\Qsym(\{x_i\}_{i \in I})$}\index{$\Qsym$}
\emph{over the alphabet $\{x_i\}_{i \in I}$}
will be the $\kk$-submodule consisting of the elements $f$ in
$R(\{x_i\}_{i\in I})$ that have the same coefficient on
the monomials $x_{i_1}^{\alpha_1} \cdots x_{i_\ell}^{\alpha_\ell}$
and $x_{j_1}^{\alpha_1} \cdots x_{j_\ell}^{\alpha_\ell}$
whenever both $i_1 < \cdots < i_\ell$ and $j_1 < \cdots < j_\ell$
in the total order on $I$. We write
$\Qsym_\kk (\{x_i\}_{i \in I})$ instead of
$\Qsym (\{x_i\}_{i \in I})$
to stress the choice of base ring $\kk$.
\end{definition}

% [DG][v17] I restructured this section, making the distinction between
% $\Qsym(\{x_i\}_{i \in I})$ and the "universal" $\Qsym$ more explicit.
% (They are identified only for $I = \{1,2,3,\ldots\}$ now.)

It immediately follows from this definition that
$\Qsym(\{x_i\}_{i \in I})$ is a free $\kk$-submodule of
$R(\{x_i\}_{i\in I})$, having as $\kk$-basis elements the
\emph{monomial quasisymmetric functions}\index{monomial quasisymmetric function}\index{$M_\alpha$}
\[
M_\alpha(\{x_i\}_{i\in I}):=
\sum_{i_1 < \cdots < i_\ell\text{ in }I}
x_{i_1}^{\alpha_1} \cdots x_{i_\ell}^{\alpha_\ell}
\]
for all compositions\footnote{Recall that compositions
were defined in Definition~\ref{def.composition-almost-composition},
along with related concepts such as length and size.}
$\alpha$ satisfying $\ell(\alpha) \leq \left|I\right|$.
When $I$ is infinite, this means that the $M_\alpha$ for all
compositions $\alpha$ form a basis of $\Qsym(\{x_i\}_{i \in I})$.

Note that
$\Qsym (\{x_i\}_{i \in I})
 = \bigoplus_{n \geq 0} \Qsym_n (\{x_i\}_{i \in I})$ is a
graded $\kk$-module of finite type, where
$\Qsym_n (\{x_i\}_{i \in I})$ is the $\kk$-submodule
of quasisymmetric functions which are homogeneous of degree $n$.
Letting \dfn{$\Comp$} denote the set of all compositions $\alpha$, and
\dfn{$\Comp_n$} the compositions $\alpha$ of $n$ (that is, compositions whose
parts sum to $n$), the subset
$\{ M_\alpha \}_{\alpha \in \Comp_n;\ \ell(\alpha) \leq \left|I\right|}$
gives a $\kk$-basis for $\Qsym_n (\{x_i\}_{i \in I})$.

\begin{example}
Taking the variable set $\xx=(x_1 < x_2 < \cdots)$ to define $\Qsym(\xx)$,
for $n=0,1,2,3$, one has these basis elements in $\Qsym_n(\xx)$:
\[
\begin{array}{rll}
M_{()} =M_\varnothing&= 1 ,&\\
& & \\
M_{(1)} &= x_1+x_2+x_3+ \cdots &= m_{(1)}=s_{(1)}=e_1=h_1=p_1 ,\\
& & \\
M_{(2)} &= x_1^2+x_2^2+x_3^2+ \cdots &=m_{(2)}=p_2 ,\\
M_{(1,1)} &= x_1x_2+x_1x_3+x_2x_3+ \cdots&=m_{(1,1)}=e_2 ,\\
& & \\
M_{(3)} &= x_1^3+x_2^3+x_3^3+ \cdots &=m_{(3)}=p_3 ,\\
M_{(2,1)} &= x_1^2 x_2 + x_1^2 x_3 + x_2^2 x_3 + \cdots ,&\\
M_{(1,2)} &= x_1 x_2^2 + x_1 x_3^2 + x_2 x_3^2 + \cdots ,&\\
M_{(1,1,1)} &= x_1 x_2 x_3 +x_1 x_2 x_4 + x_1 x_3 x_4 +  \cdots &=m_{(1,1,1)}=e_3 .
\end{array}
\]
\end{example}

It is not obvious that $\Qsym(\xx)$ is a subalgebra of $R(\xx)$,
but we will show this momentarily.  For example,
\begin{align*}
M_{(a)}M_{(b, c)}
&=(x_1^a + x_2^a + x_3^a+\cdots)
(x_1^b x_2^c + x_1^b x_3^c + x_2^b x_3^c + \cdots)\\
&=x_1^{a+b} x_2^c + \cdots +x_1^b x_3^{a+c} + \cdots
  + x_1^a x_2^b x_3^c + \cdots
  +x_1^b x_2^a x_3^c + \cdots
  + x_1^b x_2^c x_3^a + \cdots\\
&=M_{(a+b, c)} + M_{(b, a+c)} + M_{(a, b, c)} + M_{(b, a, c)} + M_{(b, c, a)} .
\end{align*}

\begin{proposition}
\label{QSym-multiplication}
For any infinite totally ordered set $I$, one has that
$\Qsym(\{x_i\}_{i \in I})$  is a $\kk$-subalgebra of
$R(\{x_i\}_{i\in I})$, with multiplication in the
$\{M_\alpha\}$-basis as follows:
Fix three disjoint chain posets
$(i_1 < \cdots <i_\ell)$, $(j_1 < \cdots < j_m)$ and
$(k_1 < k_2 < \cdots)$. Now, if
$\alpha=(\alpha_1,\ldots,\alpha_\ell)$
and $\beta=(\beta_1,\ldots,\beta_m)$
are two compositions,
then
\begin{equation}
\label{Qsym-product-on-monomials}
M_\alpha M_\beta = \sum_f M_{\wt (f)}
\end{equation}
in which the sum is over all $p \in \NN$ and all maps $f$ from the
disjoint union of two chains to a chain
\begin{equation}
\label{map-from-two-chains-to-one}
(i_1 < \cdots <i_\ell) \sqcup (j_1 < \cdots < j_m)
\overset{f}{\longrightarrow}
(k_1 < \cdots <k_p)
\end{equation}
which are both surjective and strictly order-preserving (that is,
if $x$ and $y$ are two elements in the domain satisfying $x < y$,
then $f(x) < f(y)$),
and where the composition $\wt(f):=(\wt_1(f),\ldots,\wt_p(f))$ is defined by
$
\wt_s(f):= \sum_{i_u \in f^{-1}(k_s)} \alpha_u + \sum_{j_v \in f^{-1}(k_s)} \beta_v.
$
\end{proposition}

\begin{example}
\label{exa.QSym-multiplication}
For this example, set $\alpha = \left(2,1\right)$ and
$\beta = \left(3,4,2\right)$. Let us compute $M_\alpha M_\beta$ using
\eqref{Qsym-product-on-monomials}. Indeed, the length of $\alpha$ is
$\ell = 2$, and the length of $\beta$ is $m = 3$, so the sum on the
right hand side of \eqref{Qsym-product-on-monomials} is a sum over all
$p \in \NN$ and all surjective strictly order-preserving maps $f$
from the disjoint union $(i_1 < i_2) \sqcup (j_1 < j_2 < j_3)$ of two
chains to the chain $(k_1 < k_2 < \cdots < k_p)$. Such maps can exist
only when $p \leq 5$ (due to having to be surjective) and only for
$p \geq 3$ (since, being strictly order-preserving, they have to
be injective when restricted to $(j_1 < j_2 < j_3)$). Hence,
enumerating them is a finite problem. The reader can check
that the value obtained fo $M_\alpha M_\beta$ is
\begin{align*}
& M_{\left(2,1,3,4,2\right)} + M_{\left(2,3,1,4,2\right)} + M_{\left(2,3,4,1,2\right)} + M_{\left(2,3,4,2,1\right)}
+ M_{\left(3,2,1,4,2\right)} \\
& + M_{\left(3,2,4,1,2\right)} + M_{\left(3,2,4,2,1\right)}
+ M_{\left(3,4,2,1,2\right)} + M_{\left(3,4,2,2,1\right)}
+ M_{\left(3,4,2,2,1\right)} \\
& + M_{\left(2, 3, 4, 3\right)} + M_{\left(2, 3, 5, 2\right)} + M_{\left(2, 4, 4, 2\right)} + M_{\left(3, 2, 4, 3\right)} + M_{\left(3, 2, 5, 2\right)} + M_{\left(3, 4, 2, 3\right)} \\
& + M_{\left(3, 4, 4, 1\right)} + M_{\left(3, 6, 1, 2\right)} + M_{\left(3, 6, 2, 1\right)} + M_{\left(5, 1, 4, 2\right)} + M_{\left(5, 4, 1, 2\right)} + M_{\left(5, 4, 2, 1\right)} \\
& + M_{\left(5, 4, 3\right)} + M_{\left(5, 5, 2\right)} + M_{\left(3, 6, 3\right)}.
\end{align*}
Here, we have listed the addends corresponding to $p=5$ on the first
two rows, the addends corresponding to $p=4$ on the next two rows,
and those corresponding to $p=3$ on the fifth row. The reader might
notice that the first two rows (i.e., the addends with $p=5$) are
basically a list of shuffles of $\alpha$ and $\beta$: In general,
the maps \eqref{map-from-two-chains-to-one} for $p = \ell + m$ are in
bijection with the elements of $\operatorname{Sh}_{\ell, m}$
\ \ \ \ \footnote{The bijection takes a map $f$ to the inverse of
the permutation $\sigma \in \Symm_p$ which sends every
$x \in \left\{1, 2, \ldots, \ell\right\}$ to the index $y$ satisfying
$f\left(i_x\right) = k_y$, and sends every
$x \in \left\{\ell + 1, \ell + 2, \ldots, \ell + m\right\}$ to the
index $y$ satisfying $f\left(j_{x-\ell}\right) = k_y$.}, and the
corresponding compositions $\wt (f)$ are the shuffles of $\alpha$
and $\beta$. Therefore the name ``overlapping shuffle product''.
\end{example}

% [DG][v25] Added the above example, as the overlapping shuffle rule
% isn't the easiest thing to understand. I fear my example doesn't
% help it, though...

% [DG][v26] Shortened the above example.

\begin{proof}[Proof of Proposition~\ref{QSym-multiplication}]
It clearly suffices to prove the formula
\eqref{Qsym-product-on-monomials}.
Let $\alpha = \left(\alpha_1, \ldots, \alpha_\ell\right)$
and $\beta = \left(\beta_1, \ldots, \beta_m\right)$ be two
compositions.
Fix three disjoint chain posets
$(i_1 < \cdots <i_\ell)$, $(j_1 < \cdots < j_m)$ and
$(k_1 < k_2 < \cdots)$.

Thus, multiplying
$M_\alpha = \sum_{u_1 < \cdots < u_\ell} x_{u_1}^{\alpha_1} \cdots x_{u_\ell}^{\alpha_\ell}$
with
$M_\beta = \sum_{v_1 < \cdots < v_m} x_{v_1}^{\beta_1} \cdots x_{v_m}^{\beta_m}$,
we obtain
\begin{align}
M_\alpha M_\beta
&= \sum_{u_1 < \cdots < u_\ell} \sum_{v_1 < \cdots < v_m}
\left( x_{u_1}^{\alpha_1} \cdots x_{u_\ell}^{\alpha_\ell} \right)
\left( x_{v_1}^{\beta_1} \cdots x_{v_m}^{\beta_m} \right)
\nonumber \\
&= \sum_{\gamma = \left(\gamma_1, \ldots, \gamma_p\right) \in \Comp}
  \sum_{w_1 < \cdots < w_p \text{ in } I}
  N^\gamma_{w_1, \ldots, w_p}
  x_{w_1}^{\gamma_1} \cdots x_{w_p}^{\gamma_p} ,
\label{pf.QSym-multiplication.1}
\end{align}
where $N^\gamma_{w_1, \ldots, w_p}$ is the number of all pairs
\begin{align}
\left(\left( u_1 < \cdots < u_\ell \right),
       \left( v_1 < \cdots < v_m \right) \right)
 \in I^\ell \times I^m
\label{pf.QSym-multiplication.2a}
\end{align}
of two strictly increasing tuples satisfying
\begin{align}
\left( x_{u_1}^{\alpha_1} \cdots x_{u_\ell}^{\alpha_\ell} \right)
\left( x_{v_1}^{\beta_1} \cdots x_{v_m}^{\beta_m} \right)
= x_{w_1}^{\gamma_1} \cdots x_{w_p}^{\gamma_p} .
\label{pf.QSym-multiplication.2b}
\end{align}
\footnote{In the second equality in \eqref{pf.QSym-multiplication.1},
  we have used the fact that
  each monomial can be uniquely written in the form
  $x_{w_1}^{\gamma_1} \cdots x_{w_p}^{\gamma_p}$ for some
  composition
  $\gamma = \left(\gamma_1, \ldots, \gamma_p\right) \in \Comp$
  and some strictly increasing tuple
  $\left( w_1 < \cdots < w_p \right) \in I^p$.}
Thus, we need to show that $N^\gamma_{w_1, \ldots, w_p}$
(for a given
$\gamma = \left(\gamma_1, \ldots, \gamma_p\right) \in \Comp$
and a given
$\left( w_1 < \cdots < w_p \right) \in I^p$)
is also
the number of all surjective strictly order-preserving maps
\begin{align}
(i_1 < \cdots <i_\ell) \sqcup (j_1 < \cdots < j_m)
\overset{f}{\longrightarrow}
(k_1 < \cdots <k_p)
\text{ satisfying } \wt(f) = \gamma
\label{pf.QSym-multiplication.3}
\end{align}
(because then, \eqref{pf.QSym-multiplication.1} will simplify
to \eqref{Qsym-product-on-monomials}).

In order to show this, it suffices to construct a bijection from
the set of all pairs \eqref{pf.QSym-multiplication.2a}
satisfying \eqref{pf.QSym-multiplication.2b}
to the set of all surjective strictly order-preserving maps
\eqref{pf.QSym-multiplication.3}.
This bijection is easy to construct:
Given a pair \eqref{pf.QSym-multiplication.2a}
satisfying \eqref{pf.QSym-multiplication.2b},
the bijection sends it to the map \eqref{pf.QSym-multiplication.3}
determined by:
\begin{align*}
i_g &\overset{f}{\mapsto} k_h, \text{ where } h \text{ is chosen such that }
                        u_g = w_h; \\
j_g &\overset{f}{\mapsto} k_h, \text{ where } h \text{ is chosen such that }
                        v_g = w_h.
\end{align*}
Proving that this bijection is well-defined and bijective
is straightforward\footnote{
The inverse of this bijection sends each
map \eqref{pf.QSym-multiplication.3} to the
pair \eqref{pf.QSym-multiplication.2a} determined by
\begin{align*}
u_g = w_h, \text{ where } h \text{ is chosen such that } f \left(i_g\right) = k_h; \\
v_g = w_h, \text{ where } h \text{ is chosen such that } f \left(j_g\right) = k_h.
\end{align*}
}.
\end{proof}

% [DG][v62] Made the proof above a lot more detailed.

The multiplication rule \eqref{Qsym-product-on-monomials} shows that the
$\kk$-algebra $\Qsym(\{x_i\}_{i \in I})$ does not depend much on
$I$, as long as $I$ is infinite. More precisely, all such $\kk$-algebras
are mutually isomorphic. We can use this to define a $\kk$-algebra
of quasisymmetric functions without any reference to $I$:

\begin{definition}
\label{def.Qsym}
Let $\Qsym$ be the $\kk$-algebra defined as having $\kk$-basis
$\{M_\alpha\}_{\alpha \in \Comp}$ and with multiplication defined
$\kk$-linearly by \eqref{Qsym-product-on-monomials}. This is called
the \dfn{$\kk$-algebra of quasisymmetric functions}. We write
$\Qsym_\kk$ instead of $\Qsym$ to stress the choice of base ring $\kk$.

The $\kk$-algebra $\Qsym$ is graded, and its $n$-th graded component
$\Qsym_n$ has $\kk$-basis $\{ M_\alpha \}_{\alpha \in \Comp_n}$.

For every infinite totally ordered set $I$, the $\kk$-algebra
$\Qsym$ is isomorphic to the $\kk$-algebra
$\Qsym(\{x_i\}_{i \in I})$.  The isomorphism sends
$M_\alpha \longmapsto M_\alpha(\{x_i\}_{i\in I})$.

In particular, we obtain the isomorphism
$\Qsym \cong \Qsym\left(\xx\right)$ for
$\xx$ being the infinite chain $\left(x_1 < x_2 < x_3 < \cdots\right)$.
We will identify $\Qsym$ with $\Qsym\left(\xx\right)$ along this
isomorphism. This allows us to regard quasisymmetric functions either
as power series in a specific set of variables (``alphabet''), or as
formal linear combinations of $M_\alpha$'s, whatever is more
convenient.

For any infinite alphabet
$\left\{x_i\right\}_{i \in I}$ and any $f \in \Qsym$, we denote
by $f\left(\left\{x_i\right\}_{i \in I}\right)$ the image of $f$
under the algebra isomorphism
$\Qsym \to \Qsym\left(\left\{x_i\right\}_{i \in I}\right)$ defined in
Definition~\ref{def.Qsym}.
\end{definition}

The comultiplication of $\Qsym$ will extend the one that we defined
for $\Lambda$, but we need to take care about the order of the
variables this time.
We consider the linear order from
\eqref{two-variable-set-ordering} on two sets of variables
$
(\xx,\yy)=(x_1 < x_2 < \cdots < y_1 < y_2 < \cdots),
$
and we embed the $\kk$-algebra $\Qsym(\xx) \otimes \Qsym(\yy)$ into
the $\kk$-algebra $R(\xx, \yy)$ by identifying every $f \otimes g \in
\Qsym(\xx) \otimes \Qsym(\yy)$ with $fg \in R(\xx, \yy)$
(this embedding is indeed injective\footnote{This is because it
sends the basis elements $M_\beta(\xx) \otimes M_\gamma(\yy)$ of
the former $\kk$-algebra to the linearly independent power series
$M_\beta(\xx) M_\gamma(\yy)$.}). It can then be seen that
\[
\Qsym(\xx,\yy) \subset \Qsym(\xx) \otimes \Qsym(\yy)
\]
(where the right hand side is viewed as $\kk$-subalgebra of
$R(\xx, \yy)$ via said embedding)\footnote{This is not completely
obvious, but can be easily checked by verifying that
$M_\alpha(\xx, \yy)
= \sum\limits_{\substack{(\beta,\gamma):\\ \beta \cdot \gamma=\alpha}} M_\beta (\xx) \otimes M_\gamma (\yy)$
for every composition $\alpha$ (see the proof of
Proposition~\ref{Qsym-coproduct-on-monomials} for why this holds).},
so that one can define $\Qsym \overset{\Delta}{\longrightarrow} \Qsym\otimes \Qsym$
%\begin{equation}
%\label{Qsym-coproduct-defintion}
%\begin{array}{rcl}
%\Qsym=\Qsym(\xx) & \overset{\Delta}{\longrightarrow}
%   & \Qsym(\xx,\yy) \hookrightarrow \Qsym\otimes \Qsym\\
%f(\xx)=f(x_1,x_2,\ldots) & \longmapsto & f(\xx,\yy)=f(x_1,x_2,\ldots,y_1,y_2,\ldots).
%\end{array}
%\end{equation}
as the  composite of the maps in the bottom row here:
\begin{equation}
\label{Qsym-coproduct-definition}
\begin{array}{rccccl}
      &                 & R(\xx,\yy)     & = & R(\xx, \yy) \\
      &                 & \cup           &       & \cup\\

\Qsym & \cong & \Qsym(\xx,\yy) & \hookrightarrow &
\Qsym(\xx) \otimes \Qsym(\yy) & \cong \Qsym \otimes \Qsym , \\
f&\longmapsto & f(\xx,\yy)=f(x_1,x_2,\ldots,y_1,y_2,\ldots) .& &
\end{array}
\end{equation}
(Recall that $f(\xx,\yy)$ is formally defined as the image of $f$ under the
algebra isomorphism $\Qsym \to \Qsym(\xx,\yy)$ defined in
Definition~\ref{def.Qsym}.)

% [DG][v17] There used to be a $R(\xx,\yy)$ in \eqref{Qsym-coproduct-defintion},
% but it was claiming that $R(\xx,\yy) \cong R(\xx) \otimes R(\yy)$, which
% is false ($\sum_i x_i \otimes y_i$ is in the former but not in the latter).
% When $\kk$ is not a field, there is not even a guarantee that the latter
% embeds into the former
% ( http://projecteuclid.org/DPubS?verb=Display&version=1.0&service=UI&handle=euclid.pjm/1102959646 ).
% 
% I still believe it would be better to define the coproduct on the basis.
% How do you check that
% $\Qsym(\xx,\yy) \subset \Qsym(\xx) \otimes \Qsym(\yy)$
% if not by computing it on $M_\alpha$ ?

% [DG][v19] Edited the definition of the coproduct again, this time
% explaining every step. I might have overswung the pendulum to the
% pedantic side this time...

\begin{example}
For example,
\begin{align*}
\Delta M_{(a,b,c)} &= M_{(a,b,c)}(x_1,x_2,\ldots,y_1,y_2,\ldots) \\
&= x_1^a x_2^b x_3^c+ x_1^a x_2^b x_4^c + \cdots \\
 &\quad + x_1^a x_2^b \cdot y_1^c+ x_1^a x_2^b \cdot y_2^c + \cdots \\
 &\quad + x_1^a \cdot y_1^b y_2^c + x_1^a \cdot y_1^b y_3^c + \cdots \\
 &\quad + y_1^a y_2^b y_3^c + y_1^a y_2^b y_4^c+ \cdots\\
&= M_{(a,b,c)}(\xx) + M_{(a,b)}(\xx) M_{(c)}(\yy) +  M_{(a)}(\xx) M_{(b,c)}(\yy)+ M_{(a,b,c)}(\yy) \\
&= M_{(a,b,c)} \otimes \one + M_{(a,b)} \otimes M_{(c)}
 + M_{(a)} \otimes M_{(b,c)} + \one \otimes M_{(a,b,c)} .
\end{align*}
\end{example}

\noindent
Defining the \emph{concatenation}\index{concatenation of compositions}
$\beta \cdot \gamma$ of
two compositions $\beta=(\beta_1,\ldots,\beta_r),
\gamma=(\gamma_1,\ldots,\gamma_s)$ to be the composition
$(\beta_1,\ldots,\beta_r,\gamma_1,\ldots,\gamma_s)$, one has the following
description of the coproduct in the $\{M_\alpha\}$ basis.

\begin{proposition}
\label{Qsym-coproduct-on-monomials}
For a composition $\alpha=(\alpha_1,\ldots,\alpha_\ell)$,
one has
\[
\Delta M_\alpha = \sum_{k=0}^\ell M_{(\alpha_1,\ldots,\alpha_k)} \otimes M_{(\alpha_{k+1},\ldots,\alpha_\ell)} = \sum\limits_{\substack{(\beta,\gamma):\\ \beta \cdot \gamma=\alpha}} M_\beta \otimes M_\gamma .
\]
\end{proposition}
\begin{proof}
We work with the infinite totally ordered set $I = \left\{1<2<3<\cdots\right\}$.
The definition of $\Delta$ yields
\begin{equation}
\Delta M_\alpha = M_\alpha(\xx,\yy) = \sum_{p_1 < p_2 < \cdots < p_\ell \text{ in } (\xx, \yy)} p_1^{\alpha_1} p_2^{\alpha_2} \cdots p_\ell^{\alpha_\ell} ,
\label{pf.Qsym-coproduct-on-monomials.1}
\end{equation}
where the sum runs over strictly increasing $\ell$-tuples $\left(p_1 < p_2 < \cdots < p_\ell\right)$ of variables in the variable set $(\xx, \yy)$.
But every such $\ell$-tuple $\left(p_1 < p_2 < \cdots < p_\ell\right)$ can be expressed uniquely in the form $\left(x_{i_1}, \ldots, x_{i_k}, y_{j_1}, \ldots, y_{j_{\ell-k}}\right)$ for some $k \in \left\{ 0,1,\ldots,\ell \right\}$ and some subscripts $i_1 < \cdots < i_k$ and $j_1 < \cdots < j_{\ell-k}$ in $I$. The corresponding monomial $p_1^{\alpha_1} p_2^{\alpha_2} \cdots p_\ell^{\alpha_\ell}$ then rewrites as $x_{i_1}^{\alpha_1} \cdots x_{i_k}^{\alpha_k} \cdot y_{j_1}^{\alpha_{k+1}} \cdots y_{j_{\ell-k}}^{\alpha_\ell}$. Thus, the sum on the right hand side of \eqref{pf.Qsym-coproduct-on-monomials.1} rewrites as
\begin{align*}
& \sum_{k=0}^{\ell} \ \ \sum_{i_1 < \cdots < i_k} \ \ \sum_{j_1 < \cdots < j_{\ell-k}} x_{i_1}^{\alpha_1} \cdots x_{i_k}^{\alpha_k} \cdot
y_{j_1}^{\alpha_{k+1}} \cdots y_{j_{\ell-k}}^{\alpha_\ell} \\
&= \sum_{k=0}^{\ell} \underbrace{\left(\sum_{i_1 < \cdots < i_k} x_{i_1}^{\alpha_1} \cdots x_{i_k}^{\alpha_k} \right)}_{= M_{(\alpha_1,\ldots,\alpha_k)}(\xx)} \cdot
\underbrace{\left(\sum_{j_1 < \cdots < j_{\ell-k}} y_{j_1}^{\alpha_{k+1}} \cdots y_{j_{\ell-k}}^{\alpha_\ell}\right)}_{= M_{(\alpha_{k+1},\ldots,\alpha_\ell)}(\yy)}
\\
&= \sum_{k=0}^{\ell} M_{(\alpha_1,\ldots,\alpha_k)}(\xx) M_{(\alpha_{k+1},\ldots,\alpha_\ell)}(\yy) .
\end{align*}
Thus, \eqref{pf.Qsym-coproduct-on-monomials.1} becomes
\begin{align*}
\Delta M_\alpha = \sum_{p_1 < p_2 < \cdots < p_\ell \text{ in } (\xx, \yy)} p_1^{\alpha_1} p_2^{\alpha_2} \cdots p_\ell^{\alpha_\ell}
&= \sum_{k=0}^{\ell} M_{(\alpha_1,\ldots,\alpha_k)}(\xx) M_{(\alpha_{k+1},\ldots,\alpha_\ell)}(\yy) \\
&= \sum_{k=0}^{\ell} M_{(\alpha_1,\ldots,\alpha_k)} \otimes M_{(\alpha_{k+1},\ldots,\alpha_\ell)}
= \sum\limits_{\substack{(\beta,\gamma):\\ \beta \cdot \gamma=\alpha}} M_\beta \otimes M_\gamma .
\end{align*}
% Old proof:
% This comes from expressing a monomial in $\Delta M_\alpha = M_\alpha(\xx,\yy)$
% uniquely in the form
% $x_{i_1}^{\alpha_1} \cdots x_{i_k}^{\alpha_k} \cdot
% y_{j_1}^{\alpha_{k+1}} \cdots y_{j_{\ell-k}}^{\alpha_\ell}$
% for some $k \in \left\{ 0,1,\ldots,\ell \right\}$ and some subscripts
% $i_1 < \cdots < i_k$ and $j_1 < \cdots < j_{\ell-k}$.
\end{proof}

% [DG][v72] More details in the above proof (and 1 less typo).

\begin{proposition}
The quasisymmetric functions $\Qsym$ form a connected graded Hopf algebra of finite
type, which is commutative, and contains
the symmetric functions $\Lambda$ as a Hopf subalgebra.
\end{proposition}
\begin{proof}
To prove coassociativity of $\Delta$, we need to be slightly careful.
It seems reasonable to argue by
$
(\Delta \otimes \id) \circ \Delta f =
f(\xx,\yy,\zz)=
(\id \otimes \Delta) \circ \Delta f
$
as in the case of $\Lambda$, but this would now require further
justification, as terms like $f(\xx,\yy)$ and $f(\xx,\yy,\zz)$ are
no longer directly defined as evaluations of $f$ on some sequences
(but rather are defined as images of $f$ under certain
homomorphisms). However, it is very easy to see that $\Delta$ is
coassociative by checking $(\Delta \otimes \id) \circ \Delta
= (\id \otimes \Delta) \circ \Delta $ on the $\left\{M_\alpha\right\}$
basis: Proposition~\ref{Qsym-coproduct-on-monomials} yields
\begin{align*}
\left((\Delta \otimes \id) \circ \Delta\right) M_{\alpha}
&= \sum_{k=0}^\ell \Delta(M_{(\alpha_1,\ldots,\alpha_k)}) \otimes M_{(\alpha_{k+1},\ldots,\alpha_\ell)} \\
&= \sum_{k=0}^\ell \left(\sum_{i=0}^k M_{(\alpha_1,\ldots,\alpha_i)} \otimes M_{(\alpha_{i+1},\ldots,\alpha_k)}\right) \otimes M_{(\alpha_{k+1},\ldots,\alpha_\ell)} \\
&= \sum_{k=0}^\ell \sum_{i=0}^k M_{(\alpha_1,\ldots,\alpha_i)} \otimes M_{(\alpha_{i+1},\ldots,\alpha_k)} \otimes M_{(\alpha_{k+1},\ldots,\alpha_\ell)}
\end{align*}
and the same expression for
$\left((\id \otimes \Delta) \circ \Delta\right) M_{\alpha}$.

The coproduct $\Delta$ of $\Qsym$ is an algebra morphism because
it is defined as a composite of algebra morphisms
in the bottom row of \eqref{Qsym-coproduct-definition}.
To prove that the restriction of
$\Delta$ to the subring $\Lambda$ of $\Qsym$ is the comultiplication
of $\Lambda$, it thus is enough to check that it sends the elementary
symmetric function $e_n$ to $\sum_{i=0}^n e_i \otimes e_{n-i}$ for
every $n\in\NN$. This again follows from
Proposition~\ref{Qsym-coproduct-on-monomials}, since
$e_n = M_{(1,1,\ldots,1)}$
(with $n$ times $1$).

% [DG][v3] The above might be shortenable.

The counit is as usual for a connected graded coalgebra,
and just as in the case of $\Lambda$,
sends a quasisymmetric function
$f(\xx)$ to its constant term $f(0,0,\ldots)$.
This is an evaluation, and hence an algebra morphism.
Hence $\Qsym$ forms a bialgebra, and as it is
graded and connected, also
a Hopf algebra by Proposition~\ref{graded-connected-bialgebras-have-antipodes}.
It is clearly of finite type and contains
$\Lambda$ as a Hopf subalgebra.
\end{proof}

We will identify the antipode in $\Qsym$ shortly, but we first
deal with another slightly subtle issue.
In addition to the counit evaluation $\epsilon(f) = f(0,0,\ldots)$,
starting in Section~\ref{ABS-part-I-section}, we will
want to specialize elements in $\Qsym(\xx)$ by making other
variable substitutions, in which all but a finite list of variables
are set to zero.  We justify this here.

\begin{proposition}
\label{prop.Qsym.eval}
Fix a totally ordered set $I$, a commutative $\kk$-algebra $A$,
a finite list of variables $x_{i_1}, \ldots,x_{i_m}$, say
with $i_1 < \cdots < i_m$ in $I$,
and an ordered list of elements $(a_1,\ldots,a_m) \in A^m$.

Then there is a well-defined evaluation homomorphism
\[
\begin{array}{rcl}
\Qsym(\{x_i\}_{i \in I}) &\longrightarrow &A ,\\
f &\longmapsto&
\left[ f
\right]_{\substack{x_{i_1}=a_1,\ldots,x_{i_m}=a_m\\
x_j=0 \text{ for }j \not\in \{i_1,\ldots,i_m\}}}.
\end{array}
\]
Furthermore, this homomorphism
depends only upon the list $(a_1,\ldots,a_m)$, as it coincides
with the following:
\[
\begin{array}{rccl}
\Qsym(\{x_i\}_{i \in I}) &\cong
\Qsym(x_1,x_2,\ldots) &\longrightarrow &A ,\\
&f(x_1,x_2,\ldots)& \longmapsto& f(a_1,\ldots,a_m,0,0\ldots).
\end{array}
\]
(This latter statement is stated for the case when $I$ is infinite; otherwise,
read ``$x_1,x_2,\ldots,x_{\left|I\right|}$'' for
``$x_1,x_2,\ldots$'', and interpret
$(a_1,\ldots,a_m,0,0\ldots)$ as an $\left|I\right|$-tuple.)
\end{proposition}
\begin{proof}
One already can make sense of evaluating
$x_{i_1}=a_1,\ldots,x_{i_m}=a_m$ and $x_j=0$ for $j \not \in \{i_1,\ldots,i_m\}$
in the ambient ring $R(\{x_i\}_{i \in I})$ containing $\Qsym(\{x_i\}_{i \in I})$,
since a power series $f$ of bounded degree will have finitely
many monomials that only involve the variables
$x_{i_1},\ldots,x_{i_m}$.  The last assertion follows from quasisymmetry
of $f$, and is perhaps checked most easily when $f=M_\alpha(\{x_i\}_{i \in I})$
for some $\alpha$.
\end{proof}

The antipode in $\Qsym$ has a reasonably simple expression
in the $\{M_\alpha\}$ basis, but requiring a definition.

\begin{definition}
\label{def.QSym.D}
For $\alpha,\beta$ in $\Comp_n$, say that $\alpha$
\emph{refines}\index{refinement of compositions}\index{refining a composition}
$\beta$
or $\beta$
\emph{coarsens}\index{coarsening of compositions}\index{coarsening a composition}
$\alpha$ if, informally, one can obtain
$\beta$ from $\alpha$ by combining some of its adjacent parts. 
Alternatively, this can be defined as follows:
One has a bijection $\Comp_n \rightarrow 2^{[n-1]}$ where
$[n-1]:=\{1,2,\ldots,n-1\}$ which sends
$\alpha=(\alpha_1,\ldots,\alpha_\ell)$ having
length $\ell(\alpha)=\ell$ to its subset of partial sums\index{$D(\alpha)$}
\[
D(\alpha):= \left\{ \alpha_1,\alpha_1+\alpha_2,\ldots,\alpha_1+\cdots+\alpha_{\ell-1} \right\},
\]
and this sends the refinement ordering to the inclusion ordering
on the Boolean algebra $2^{[n-1]}$ (to be more precise: a composition
$\alpha \in \Comp_n$ refines a composition $\beta \in \Comp_n$
if and only if $D(\alpha) \supset D(\beta)$).

There is also a bijection sending every composition
$\alpha$ to its \dfn{ribbon}\index{ribbon diagram}\index{$\Rib(\alpha)$}
diagram $\Rib\left(\alpha\right)$:
the skew diagram $\lambda/\mu$
having rows of sizes $\alpha_1,\ldots,\alpha_\ell$ read from bottom to top
with exactly one column of overlap between adjacent rows.
These bijections and the refinement partial order are illustrated here for $n=4$:
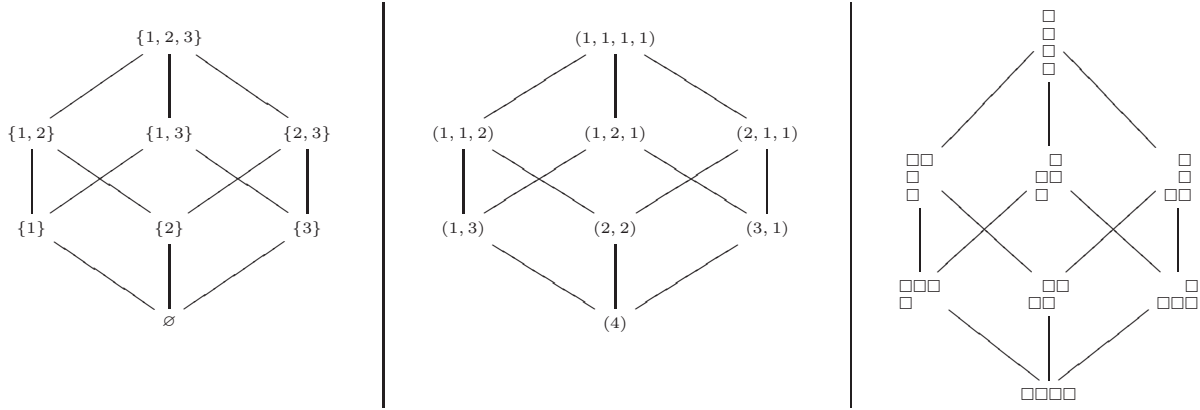
\[
\begin{array}{c|c|c}
\tiny
\xymatrix{
 & \{1,2,3\} & \\
\{1,2\} \ar@{-}[ur]& \{1,3\} \ar@{-}[u] & \{2,3\} \ar@{-}[ul] \\
\{1\} \ar@{-}[ur]\ar@{-}[u]& \{2\}\ar@{-}[ur]\ar@{-}[ul] & \{3\} \ar@{-}[ul]\ar@{-}[u]\\
&\varnothing \ar@{-}[ur]\ar@{-}[u]\ar@{-}[ul]&
}
\quad \phantom{a} & \quad
\tiny \xymatrix{
 & (1,1,1,1)& \\
(1,1,2) \ar@{-}[ur]& (1,2,1) \ar@{-}[u] & (2,1,1) \ar@{-}[ul] \\
(1,3) \ar@{-}[ur]\ar@{-}[u]& (2,2) \ar@{-}[ur]\ar@{-}[ul] & (3,1) \ar@{-}[ul]\ar@{-}[u]\\
&(4) \ar@{-}[ur]\ar@{-}[u]\ar@{-}[ul]&
}
\quad \phantom{a} & \quad
\tiny \xymatrix{
 &\txt{$\sq$ \\ $\sq$ \\ $\sq$ \\  $\sq$} & \\
  \txt{$\sq \sq$\\$\sq\phantom{\sq}$ \\$\sq\phantom{\sq}$} \ar@{-}[ur]
 & \txt{$\phantom{\sq}\sq$ \\  $\sq \sq$ \\$\sq\phantom{\sq}$ } \ar@{-}[u]
 & \txt{$\phantom{\sq}\sq$ \\ $\phantom{\sq}\sq$ \\  $\sq \sq$} \ar@{-}[ul] \\
 \txt{$\sq \sq \sq$\\$\sq\phantom{\sq}\phantom{\sq}$} \ar@{-}[ur]\ar@{-}[u]
 & \txt{$\phantom{\sq}\sq\sq$\\$\sq \sq \phantom{\sq}$} \ar@{-}[ur]\ar@{-}[ul]
 & \txt{$\phantom{\sq}\phantom{\sq}\sq$\\$\sq \sq \sq$} \ar@{-}[ul]\ar@{-}[u]\\
&\txt{$\sq\sq\sq\sq$} \ar@{-}[ur]\ar@{-}[u]\ar@{-}[ul]&
}
\end{array}
\]
(where we have drawn each ribbon diagram with its boxes spaced
out).

Given $\alpha=(\alpha_1,\ldots,\alpha_\ell)$,
its \emph{reverse}\index{reverse composition}\index{$\rev(\alpha)$}
composition is
$\rev(\alpha)=(\alpha_{\ell},\alpha_{\ell-1},\ldots,\alpha_2,\alpha_1)$.
Note that $\alpha \mapsto \rev(\alpha)$ is a poset automorphism of
$\Comp_n$ for the refinement ordering.
\end{definition}

% [DG][v3] I have replaced "from bottom-to-top" by "from bottom to top".
% Done the same to "from left-to-right".

% [DG][v78] Introduced the notation "\Rib" for the ribbon
% corresponding to a composition (so we no longer identify
% it with the composition itself, thus avoiding a possible
% confusion when the composition is a partition).

\begin{theorem}
\label{Qsym-antipode-on-monomials}
For any composition $\alpha$ in $\Comp$,
\[
S(M_\alpha)= (-1)^{\ell(\alpha)} \sum\limits_{\substack{\gamma \in \Comp:\\ \gamma \text{ coarsens }\rev(\alpha)}} M_{\gamma} .
\]
\end{theorem}

For example,
\[
S(M_{(a,b,c)})=-\left( M_{(c,b,a)} + M_{(b+c,a)} + M_{(c,a+b)} + M_{(a+b+c)} \right) .
\]

\begin{proof}
We give Ehrenborg's proof\footnote{A different proof was
given by Malvenuto and Reutenauer \cite[Cor. 2.3]{MalvenutoReutenauer}, and is
sketched in Remark~\ref{second-antipode-proof} below.} \cite[Prop. 3.4]{Ehrenborg}
via induction on $\ell=\ell(\alpha)$.
One has easy base cases when
$\ell(\alpha)=0$, where
$
S(M_\varnothing) = S(1)=1=(-1)^0 M_{\rev(\varnothing)},
$
and when
$\ell(\alpha)=1$, where $M_{(n)}$ is primitive by
Proposition~\ref{Qsym-coproduct-on-monomials}, so Proposition~\ref{antipode-on-primitives} shows
$
S(M_{(n)})=-M_{(n)}=(-1)^1 M_{\rev((n))}.
$

For the inductive step, apply the inductive definition of $S$
from the proof of
Proposition~\ref{graded-connected-bialgebras-have-antipodes}:
\begin{align*}
S(M_{(\alpha_1,\ldots,\alpha_\ell)})
&=-\sum_{i=0}^{\ell-1} S(M_{(\alpha_1,\ldots,\alpha_i)})
M_{(\alpha_{i+1},\ldots,\alpha_\ell)}  \\
&=\sum_{i=0}^{\ell-1}
\sum\limits_{\substack{ \beta\text{ coarsening } \\ (\alpha_i,\alpha_{i-1},\ldots,\alpha_1)}}
(-1)^{i+1} M_{\beta} M_{(\alpha_{i+1},\ldots,\alpha_\ell)} .
\end{align*}
The idea will be to cancel terms of opposite sign
that appear in the expansions of the products
$M_\beta  M_{(\alpha_{i+1},\ldots,\alpha_\ell)}$.  Note that
each composition $\beta$ appearing above has
first part $\beta_1$ of the form $\alpha_i+\alpha_{i-1}+\cdots+\alpha_h$ for some
$h \leq i$ (unless $\beta = \varnothing$), and hence each term $M_\gamma$ in the
expansion of the product
$
M_\beta  M_{(\alpha_{i+1},\ldots,\alpha_\ell)}
$
has $\gamma_1$ (that is, the first entry of $\gamma$) a sum that
can take one of these three forms:
\begin{enumerate}
\item[$\bullet$]
$\alpha_i+\alpha_{i-1}+\cdots+\alpha_h$,
\item[$\bullet$]
$\alpha_{i+1}+(\alpha_i+\alpha_{i-1}+\cdots+\alpha_h)$,
\item[$\bullet$]
$\alpha_{i+1}$.
\end{enumerate}
Say that the \emph{type} of $\gamma$
is $i$ in the first case, and $i+1$ in the second two cases\footnote{We
imagine that we label the terms obtained by expanding
$M_\beta M_{(\alpha_{i+1},\ldots,\alpha_\ell)}$ by distinct labels,
so that each term knows how exactly it was created (i.e., which $i$, which $\beta$
and which map $f$ as in \eqref{map-from-two-chains-to-one} gave
rise to it). Strictly speaking, it is these triples $\left(i,\beta,f\right)$
that we should be assigning types to, not terms.};
in other words, the type is the largest subscript
$k$ on a part $\alpha_k$ which was combined in the sum $\gamma_1$.
It is not hard to see that a given $\gamma$ for which the type $k$
is strictly smaller than $\ell$ arises from exactly two
pairs $(\beta,\gamma),(\beta',\gamma)$, having opposite signs $(-1)^k$
and $(-1)^{k+1}$ in the above sum\footnote{Strictly speaking, this
means that we have an involution on the set of our $\left(i,\beta,f\right)$
triples having type smaller than $\ell$, and this involution switches the
sign of $(-1)^i M_{\wt (f)}$.}.
For example, if $\alpha=(\alpha_1,\ldots,\alpha_8)$,
then the composition
$\gamma=(\alpha_6+\alpha_5+\alpha_4,\alpha_3,\alpha_7,
          \alpha_8+\alpha_2+\alpha_1)$
of type $6$
can arise from either of
\begin{align*}
\beta&=(\alpha_6+\alpha_5+\alpha_4,\alpha_3,\alpha_2+\alpha_1)
\text{ with }i=6\text{ and sign }(-1)^7 ,\\
\beta'&=(\alpha_5+\alpha_4,\alpha_3,\alpha_2+\alpha_1)
\text{ with }i=5 \text{ and sign }(-1)^6.
\end{align*}
Similarly,
$\gamma=(\alpha_6,\alpha_5+\alpha_4,\alpha_3,\alpha_7,
          \alpha_8+\alpha_2+\alpha_1)$ can arise from either of
\begin{align*}
\beta&=(\alpha_6,\alpha_5+\alpha_4,\alpha_3,\alpha_2+\alpha_1)
\text{ with }i=6\text{ and sign }(-1)^7,\\
\beta'&=(\alpha_5+\alpha_4,\alpha_3,\alpha_2+\alpha_1)
\text{ with }i=5\text{ and sign }(-1)^6.
\end{align*}
Thus one can cancel almost all the terms, excepting those
with $\gamma$ of type $\ell$ among the terms $M_\gamma$
in the expansion of the last ($i=\ell-1$) summand $M_\beta M_{(\alpha_\ell)}$.
A bit of thought shows that these are the
$\gamma$ coarsening $\rev(\alpha)$, and all have sign $(-1)^\ell$.
\end{proof}

% [DG][v3] I added the last two footnotes in the proof above.

\subsection{The fundamental basis and $P$-partitions}

There is a second important basis for $\Qsym$ which arose originally
in Stanley's $P$-partition theory \cite{Stanley-thesis}.\footnote{See
\cite{Gessel-Ppart} for a history of $P$-partitions; our notations,
however, strongly differ from those in \cite{Gessel-Ppart}.}

\begin{definition}
A \dfn{labelled poset} will here mean
a partially ordered set $P$ whose underlying set
is some finite subset of the integers.  A \dfn{$P$-partition} is a function
$P \overset{f}{\rightarrow} \{1,2,\ldots\}$ with the following
two properties:
\begin{enumerate}
\item[$\bullet$]
If $i \in P$ and $j \in P$ satisfy
$i <_P j$ and $i <_\ZZ j$,
then $f(i) \leq f(j)$.
\item[$\bullet$]
If $i \in P$ and $j \in P$ satisfy
$i <_P j$ and $i >_\ZZ j$, then
$f(i) < f(j)$.
\end{enumerate}
Denote by \dfn{$\AAA(P)$} the set of all $P$-partitions $f$, and let
$
F_P(\xx):=\sum_{f \in \AAA(P)} \xx_f
$
where $\xx_f:=\prod_{i \in P} x_{f(i)}$. This $F_P(\xx)$
is an element of $\kk\left[\left[\xx\right]\right]
:= \kk\left[\left[x_1,x_2,\ldots\right]\right]$.
\end{definition}

% [DG][v72] Made the above definition a bit less jargonish by
% replacing "implies" by proper "if/then"-clauses. Similar
% changes made everywhere the word "implies"/"imply" appeared.

\begin{example}
\label{labelled-poset-example}
Depicted is a labelled poset $P$, along with the relations among the
four values $f=(f(1),f(2),f(3),f(4))$ that define its $P$-partitions $f$:
\[
\xymatrix{
 & &2\\
4& &1\ar@{-}[u]\\
 &3\ar@{-}[ur]\ar@{-}[ul]&
} \qquad
\xymatrix{
    & &f(2)\\
f(4)& &f(1)\ar@{-}[u]^{\leq}\\
    &f(3)\ar@{-}[ul]^{\leq} \ar@{-}[ur]^{<}&
}
\]
\end{example}

\begin{remark}
Stanley's treatment of $P$-partitions in \cite[\S 3.15 and \S 7.19]{Stanley} uses
a language different from ours. First, Stanley works not with labelled posets
$P$, but with pairs $\left(P, \omega\right)$ of a poset $P$ and a bijective
labelling $\omega : P \to [n]$. Thus, the relation $<_\ZZ$ is not given on $P$
a priori, but has to be pulled back from $[n]$ using $\omega$ (and it depends
on $\omega$, whence Stanley speaks of ``$\left(P, \omega\right)$-partitions'').
Furthermore, what we call ``$P$-partition'' is called a ``reverse $P$-partition''
in \cite{Stanley}.
Finally, Stanley uses the notations $F_P$ and $F_{P, \omega}$ for something
different from what we denote by $F_P$, whereas what we call $F_P$ is dubbed
$K_{P, \omega}$ in \cite[\S 7.19]{Stanley}.
\end{remark}

% [DG][v25] Added preceding remark on the babylon in this subject.

The so-called \emph{fundamental quasisymmetric functions}
are an important special case of the $F_P(\xx)$.
We shall first define them directly and then see how they are obtained as
$P$-partition enumerators $F_P(\xx)$ for some special labelled posets $P$.

\begin{definition}
\label{fundamental-qsym-def}
Let $n\in\NN$ and $\alpha\in\Comp_n$. We define the
\dfn{fundamental quasisymmetric function}\index{$L_\alpha$}
$L_\alpha = L_\alpha(\xx) \in \Qsym$
by
\begin{equation}
\label{fundamental-qsym-fn-defn}
L_\alpha := \sum\limits_{\substack{\beta \in \Comp_n:\\
                               \beta \text{ refines }\alpha}} M_\beta.
\end{equation}
\end{definition}

\begin{example}
\label{the-two-extreme-fundamentals-example}
The extreme cases for $\alpha$ in $\Comp_n$ give quasisymmetric functions
$L_\alpha$ which are symmetric:
\begin{align*}
L_{(1^n)}&=M_{(1^n)}=e_n,\\
L_{(n)}&=\sum_{\alpha \in \Comp_n} M_{\alpha}= h_n .\\
\end{align*}
\end{example}

Before studying the $L_\alpha$ in earnest, we recall a basic fact
about finite sets, which is sometimes known as the ``principle of
inclusion and exclusion'' (although it is more general than the
formula for the size of a union of sets that commonly
goes by this name):

\begin{lemma}
\label{lem.moebius-bool.subsets}
Let $G$ be a finite set.
Let $V$ be a $\kk$-module.
For each subset $A$ of $G$, we let $f_A$ and $g_A$ be two
elements of $V$.

\begin{enumerate}
\item[(a)] If
\[
\text{every } A \subset G \text{ satisfies } g_A = \sum_{B \subset A} f_B ,
\]
then
\[
\text{every } A \subset G \text{ satisfies } f_A = \sum_{B \subset A} \left(-1\right)^{\left|A \setminus B\right|} g_B .
\]

\item[(b)] If
\[
\text{every } A \subset G \text{ satisfies } g_A = \sum_{B \subset G; \  B \supset A} f_B ,
\]
then
\[
\text{every } A \subset G \text{ satisfies } f_A = \sum_{B \subset G; \  B \supset A} \left(-1\right)^{\left|B \setminus A\right|} g_B .
\]
\end{enumerate}
\end{lemma}

\begin{proof}
This can be proven by elementary arguments (easy exercise).
Alternatively, Lemma~\ref{lem.moebius-bool.subsets} can be
viewed as a particular case of the M\"obius inversion principle
(see, e.g., \cite[Propositions 3.7.1 and 3.7.2]{Stanley})
applied to the Boolean lattice $2^G$
(whose M\"obius function is very simple: see
\cite[Example 3.8.3]{Stanley}).
(This is spelled out in \cite[Example 4.52]{Loehr-bij}, for example.)
\end{proof}

Lemma~\ref{lem.moebius-bool.subsets} can be translated into the language
of compositions:

\begin{lemma}
\label{lem.moebius-bool.comps}
Let $n \in \NN$.
Let $V$ be a $\kk$-module.
For each $\alpha \in \Comp_n$, we let $f_\alpha$ and $g_\alpha$ be two
elements of $V$.

\begin{enumerate}
\item[(a)] If
\[
\text{every } \alpha \in \Comp_n \text{ satisfies } g_\alpha = \sum_{\beta \text{ coarsens } \alpha} f_\beta ,
\]
then
\[
\text{every } \alpha \in \Comp_n \text{ satisfies } f_\alpha = \sum_{\beta \text{ coarsens } \alpha} \left(-1\right)^{\ell\left(\alpha\right) - \ell\left(\beta\right)} g_\beta .
\]

\item[(b)] If
\[
\text{every } \alpha \in \Comp_n \text{ satisfies } g_\alpha = \sum_{\beta \text{ refines } \alpha} f_\beta ,
\]
then
\[
\text{every } \alpha \in \Comp_n \text{ satisfies } f_\alpha = \sum_{\beta \text{ refines } \alpha} \left(-1\right)^{\ell\left(\beta\right) - \ell\left(\alpha\right)} g_\beta .
\]
\end{enumerate}
\end{lemma}

\begin{proof}%[Proof of Lemma~\ref{lem.moebius-bool.comps}.]
Set $\left[n-1\right] = \left\{1,2,\ldots,n-1\right\}$.
Recall (from Definition~\ref{def.QSym.D}) that there is a bijection
$D : \Comp_n \rightarrow 2^{\left[n-1\right]}$ that sends each
$\alpha \in \Comp_n$ to $D \left(\alpha\right) \subset \left[n-1\right]$.
This bijection $D$ has the properties that:
\begin{itemize}
\item a composition $\beta$ refines a composition $\alpha$ if and only if
      $D\left(\beta\right) \supset D\left(\alpha\right)$;
\item a composition $\beta$ coarsens a composition $\alpha$ if and only if
      $D\left(\beta\right) \subset D\left(\alpha\right)$;
\item any composition $\alpha \in \Comp_n$ satisfies
      $\left|D\left(\alpha\right)\right| = \ell\left(\alpha\right)-1$
      (unless $n=0$), and thus
\item any compositions $\alpha$ and $\beta$ in $\Comp_n$ satisfy
      $\left|D\left(\alpha\right)\right| - \left|D\left(\beta\right)\right|
      = \ell\left(\alpha\right) - \ell\left(\beta\right)$.
\end{itemize}
This creates a dictionary between compositions in $\Comp_n$ and
subsets of $\left[n-1\right]$.
Now, apply Lemma~\ref{lem.moebius-bool.subsets} to
$G = \left[n-1\right]$, $f_A = f_{D^{-1}\left(A\right)}$ and
$g_A = g_{D^{-1}\left(A\right)}$, and translate using the
dictionary.
\end{proof}

% [DG][v70] Added the above two lemmas, and replaced the references
% to "inclusion-exclusion" by references to these lemmas.
% The reason is that "inclusion-exclusion" means a rather different
% thing for beginners in combinatorics.

Now, we can see the following about the fundamental quasisymmetric
functions:

\begin{proposition}
\label{prop.QSym.L-basis}
The family $\{L_\alpha\}_{\alpha \in \Comp}$
is a $\kk$-basis for $\Qsym$, and
each $n \in \NN$ and $\alpha \in \Comp_n$ satisfy
\begin{equation}
M_\alpha = \sum\limits_{\substack{\beta \in \Comp_n:\\
           \beta \text{ refines }\alpha}} (-1)^{\ell(\beta)-\ell(\alpha)} L_\beta.
\label{eq.QSym.M-through-L}
\end{equation}
\end{proposition}

\begin{proof}%[Proof of Proposition~\ref{prop.QSym.L-basis}.]
Fix $n \in \NN$.
Recall the equality \eqref{fundamental-qsym-fn-defn}.
Thus, Lemma~\ref{lem.moebius-bool.comps}(b) (applied to
$V = \Qsym$, $f_\alpha = M_\alpha$ and $g_\alpha = L_\alpha$)
yields \eqref{eq.QSym.M-through-L}.

Recall that the family $\left(M_\alpha\right)_{\alpha \in \Comp_n}$
is a basis of the $\kk$-module $\Qsym_n$.
The equality \eqref{fundamental-qsym-fn-defn} shows that the
family $\left(L_\alpha\right)_{\alpha \in \Comp_n}$ expands
invertibly
triangularly\footnote{See Section~\ref{sect.STmat} for a definition
of this concept.} with respect to the family
$\left(M_\alpha\right)_{\alpha \in \Comp_n}$
(where $\Comp_n$ is equipped with the refinement
order).\footnote{In fact, it expands unitriangularly with respect
to the latter family.}
Thus, Corollary~\ref{cor.STmat.expansion-tria-inv}(e) (applied to
$\Qsym_n$, $\Comp_n$, $\left(M_\alpha\right)_{\alpha \in \Comp_n}$
and $\left(L_\alpha\right)_{\alpha \in \Comp_n}$ instead of
$M$, $S$, $\left(e_s\right)_{s\in S}$ and
$\left(f_s\right)_{s\in S}$) shows that the
family $\left(L_\alpha\right)_{\alpha \in \Comp_n}$ is a basis
of the $\kk$-module $\Qsym_n$.
Combining this fact for all $n \in \NN$, we conclude that
the family $\left(L_\alpha\right)_{\alpha \in \Comp}$ is a basis
of the $\kk$-module $\Qsym$.
This completes the proof of Proposition~\ref{prop.QSym.L-basis}.
\end{proof}

% [DG][v70] Promoted the above to a proposition, with an
% explicit proof.

\begin{proposition}
\label{Qsym-technical-lemma-Lxy}
Let $n\in\NN$. Let $\alpha$ be a composition of $n$. Let $I$ be an
infinite totally ordered set. Then,
\[
L_\alpha \left(  \left\{  x_{i}\right\}_{i\in I}\right)
=\sum\limits_{\substack{i_1 \leq i_2 \leq \cdots \leq i_n \text{ in } I; \\ 
i_{j} < i_{j+1} \text{ if } j\in D\left(  \alpha\right) }}
x_{i_1} x_{i_2} \cdots x_{i_n},
\]
where $L_\alpha\left(  \left\{  x_i \right\}_{i\in I}\right)$ is
defined as the image of $L_\alpha$ under the isomorphism
$\Qsym \to \Qsym\left( \left\{  x_i \right\}_{i\in I}\right)$ obtained
in Definition~\ref{def.Qsym}.
In particular, for the standard (totally ordered) variable set
$\xx = \left(x_1 < x_2 < \cdots\right)$, we obtain
\begin{equation}
\label{fundamental-qsym-fn-defn-alternative}
L_\alpha = L_\alpha\left(\xx\right)
=\sum\limits_{\substack{(1 \leq ) i_1 \leq i_2 \leq \cdots \leq i_n; \\
i_{j} < i_{j+1} \text{ if } j\in D\left(  \alpha\right) }}
x_{i_1} x_{i_2} \cdots x_{i_n}.
\end{equation}
\end{proposition}

\begin{proof}
Every composition $\beta =
\left(\beta_1,\ldots,\beta_\ell\right)$ of $n$ satisfies
%% begin pedantry
%\footnote{For the second equality of \eqref{Qsym-technical-lemma-Lxy-proof}, we use the bijection between:
%\begin{enumerate}
%\item[$\bullet$] the sequences $k_1 < \cdots < k_\ell$ in $I$,
%\item[$\bullet$] the sequences $i_{1}\leq i_{2}\leq\cdots \leq i_{n}$ in $I$
%which satisfy $i_j < i_{j+1}$ precisely for those $j$
%satisfying $j\in D\left(  \beta\right)$.
%\end{enumerate}
%This bijection sends $k_1 < \cdots < k_\ell$ to
%$\underbrace{k_1 \leq \cdots \leq k_1}_{\beta_1 \text{ times}} \leq
%\underbrace{k_2 \leq \cdots \leq k_2}_{\beta_2 \text{ times}} \leq \cdots \leq
%\underbrace{k_\ell \leq \cdots \leq k_\ell}_{\beta_\ell \text{ times}}$,
%and thus replaces the monomial $x_{k_1}^{\beta_1} \cdots x_{k_\ell}^{\beta_\ell}$
%by $x_{i_{1}}x_{i_{2}}\cdots x_{i_{n}}$.}
%% end pedantry
\begin{equation}
\label{Qsym-technical-lemma-Lxy-proof}
M_\beta\left(\{x_i\}_{i\in I}\right)
= \sum_{k_1 < \cdots < k_\ell\text{ in }I}
x_{k_1}^{\beta_1} \cdots x_{k_\ell}^{\beta_\ell}
=\sum\limits_{\substack{i_{1}\leq
i_{2}\leq\cdots \leq i_{n}\text{ in } I;\\
i_{j}<i_{j+1}\text{ if and only if }j\in D\left(  \beta\right)
}}x_{i_{1}}x_{i_{2}}\cdots x_{i_{n}}.
\end{equation}
Applying the ring homomorphism
$\Qsym \to \Qsym\left( \left\{  x_{i}\right\}_{i\in I}\right)$
to \eqref{fundamental-qsym-fn-defn}, we obtain
\begin{align*}
L_{\alpha}\left(  \left\{  x_{i}\right\}_{i\in I}\right)
&= \sum\limits_{\substack{\beta \in \Comp_n:\\
                               \beta \text{ refines }\alpha}} M_\beta \left(  \left\{  x_{i}\right\}_{i\in I}\right) 
\overset{\eqref{Qsym-technical-lemma-Lxy-proof}}{=} \sum\limits_{\substack{\beta \in \Comp_n:\\
                               \beta \text{ refines }\alpha}} 
\quad \sum\limits_{\substack{i_{1}\leq
i_{2}\leq\cdots \leq i_{n}\text{ in } I;\\
i_{j}<i_{j+1}\text{ if and only if }j\in D\left(  \beta\right)
}}x_{i_{1}}x_{i_{2}}\cdots x_{i_{n}} \\
&= \sum\limits_{\substack{\beta \in \Comp_n:\\    D(\alpha) \subset D(\beta)}} \quad
\sum\limits_{\substack{i_{1}\leq
i_{2}\leq\cdots \leq i_{n}\text{ in } I;\\
i_{j}<i_{j+1}\text{ if and only if }j\in D\left(  \beta\right)
}}x_{i_{1}}x_{i_{2}}\cdots x_{i_{n}} \\
&= \sum\limits_{\substack{Z \subset [n-1]:\\
                               D(\alpha) \subset Z}} 
\quad \sum\limits_{\substack{i_{1}\leq
i_{2}\leq\cdots \leq i_{n}\text{ in } I;\\
i_{j}<i_{j+1}\text{ if and only if }j\in Z
}}x_{i_{1}}x_{i_{2}}\cdots x_{i_{n}}
= \sum\limits_{\substack{i_1 \leq i_2 \leq \cdots \leq i_n \text{ in } I;\\
i_{j}<i_{j+1} \text{ if } j\in D\left(  \alpha\right)  }}
x_{i_{1}}x_{i_{2}}\cdots x_{i_{n}}.
\end{align*}
\end{proof}

\begin{proposition}
\label{fundamental-P-partition-prop}
Assume that the labelled poset $P$ is a total or linear order
$w=(w_1 < \cdots < w_n)$
(that is, $P = \left\{w_1,w_2,\ldots,w_n\right\}$ as sets,
and the order $<_P$ is given by $w_1 <_P w_2 <_P \cdots <_P w_n$).
Let $\Des(w)$ be the
\emph{descent set}\index{descent set of a permutation} of $w$,
defined by
\[
\Des(w) := \{ i : w_i >_{\ZZ} w_{i+1} \} \subset \{1,2,\ldots,n-1\} .
\]
Let $\alpha \in \Comp_n$ be the unique composition in $\Comp_n$ having
partial sums $D(\alpha)=\Des(w)$.
Then, the generating function
$F_w(\xx)$ equals the fundamental quasisymmetric function
$L_\alpha$.
In particular, $F_w(\xx)$ depends only upon the descent set $\Des(w)$.
\end{proposition}
% [DG][v72] Restated the above proposition so it's no longer one
% long convoluted sentence; also, added the "that is" parenthetical
% to clarify what the $<$ signs in "$(w_1 < \cdots < w_n)$" mean.
\noindent
E.g., total order $w=35142$ has $\Des(w)=\{2,4\}$ and
composition $\alpha=(2,2,1)$, so
\begin{align*}
F_{35142}(\xx)
&=\sum_{f(3) \leq f(5) < f(1) \leq f(4) < f(2)}
         x_{f(3)}x_{f(5)}x_{f(1)}x_{f(4)}x_{f(2)}\\
&=\sum_{i_1 \leq i_2 < i_3 \leq i_4 < i_5}
         x_{i_1}x_{i_2}x_{i_3}x_{i_4}x_{i_5}\\
&= L_{(2,2,1)} = M_{(2,2,1)}+M_{(2,1,1,1)}
                               +M_{(1,1,2,1)}+M_{(1,1,1,1,1)}.
\end{align*}

\begin{proof}[Proof of Proposition~\ref{fundamental-P-partition-prop}.]
Write $F_w(\xx)$ as a sum of monomials $x_{f(w_1)} \cdots x_{f(w_n)}$
over all $w$-partitions $f$. These $w$-partitions are exactly the
maps $f : w \to \left\{1,2,3,\ldots\right\}$ satisfying
$f(w_1) \leq \cdots \leq f(w_n)$ and having
strict inequalities $f(w_i) < f(w_{i+1})$ whenever $i$ is in $\Des(w)$
(because if two elements $w_a$ and $w_b$ of $w$ satisfy $w_a <_w w_b$
and $w_a >_{\ZZ} w_b$, then they must satisfy $a < b$ and
$i \in \Des(w)$ for some $i \in \left\{a, a+1, \ldots, b-1\right\}$; thus, the
conditions ``$f(w_1) \leq \cdots \leq f(w_n)$'' and
``$f(w_i) < f(w_{i+1})$ whenever $i$ is in $\Des(w)$''
ensure that $f\left(w_a\right) < f\left(w_b\right)$ in this case).
Therefore,
they are in bijection with the weakly increasing sequences
$\left(i_1 \leq i_2 \leq \cdots \leq i_n\right)$ of positive integers
having strict inequalities $i_j < i_{j+1}$ whenever $i \in \Des(w)$
(namely, the bijection sends any $w$-partition $f$ to the sequence
$\left(f\left(w_1\right) \leq f\left(w_2\right) \leq \cdots \leq f\left(w_n\right)\right)$).
Hence,
\[
F_w(\xx)
= \sum_{f \in \mathcal{A}(w)} \xx_f
= \sum\limits_{\substack{(1 \leq ) i_1 \leq i_2 \leq \cdots \leq i_n; \\
i_{j} < i_{j+1} \text{ if } j\in \Des(w) }}
x_{i_1} x_{i_2} \cdots x_{i_n}
= \sum\limits_{\substack{(1 \leq ) i_{1} \leq i_{2} \leq \cdots \leq i_n; \\
i_{j} < i_{j+1} \text{ if } j\in D\left(  \alpha\right) }}
x_{i_1} x_{i_2} \cdots x_{i_n}
\]
(since $\Des(w) = D(\alpha)$).
Comparing this with
\eqref{fundamental-qsym-fn-defn-alternative}, we conclude that
$F_w(\xx) = L_\alpha$.
%% Old argument:
%The underlying set $\{f(w_i)\}_{i=1}^n$ will equal $\{j_1 < \cdots < j_\ell\}$
%with indices in increasing order having
%a multiplicity sequence $\beta=(\beta_1,\ldots,\beta_\ell)$
%that gives a composition $\beta$ refining $\alpha$.
\end{proof}

% [DG][v72] Added more details to the above proof.

The next proposition (\cite[Cor. 7.19.5]{Stanley},
\cite[Cor. 3.3.24]{LuotoMykytiukvanWilligenburg}) is an algebraic shadow of
Stanley's main lemma \cite[Thm. 7.19.4]{Stanley} in $P$-partition theory.
  It expands any $F_P(\xx)$ in the $\{L_\alpha\}$ basis,
as a sum over the set \dfn{$\LLL(P)$} of all \emph{linear extensions} $w$
of $P$\ \ \ \ \footnote{Let us explain what we mean by linear extensions
and how we represent them.
\par
If $\mathbf{P}$ is a finite poset, then a \dfn{linear extension} of
$\mathbf{P}$ denotes a total order $w$ on the set $\mathbf{P}$ having the
property that every two elements $i$ and $j$ of $\mathbf{P}$ satisfying
$i <_{\mathbf{P}} j$ satisfy $i <_w j$.
(In other words, it is a linear order on
the ground set $\mathbf{P}$ which extends $\mathbf{P}$ as a poset;
therefore the name.) We identify such a total order $w$ with the % (unique)
list $\left(\mathbf{p}_1, \mathbf{p}_2, \ldots, \mathbf{p}_n\right)$
containing all elements of $\mathbf{P}$ in $w$-increasing order (that
is, $\mathbf{p}_1 <_w \mathbf{p}_2 <_w \cdots <_w \mathbf{p}_n$).
\par
(Stanley, in \cite[\S 3.5]{Stanley}, defines linear extensions in a
slightly different way: For him, a linear extension of a finite
poset $\mathbf{P}$ is an order-preserving bijection from $\mathbf{P}$
to the subposet $\left\{ 1, 2, \ldots, \left|\mathbf{P}\right| \right\}$
of $\ZZ$. But this is equivalent to our definition, since a bijection
like this can be used to transport the order relation of
$\left\{ 1, 2, \ldots, \left|\mathbf{P}\right| \right\}$
back to $\mathbf{P}$, thus resulting in a total order on $\mathbf{P}$
which is a linear extension of $\mathbf{P}$ in our sense.)
}. %, that is, the set of all extensions of $P$ to a linear order.
E.g., the poset $P$ from Example~\ref{labelled-poset-example} has
$\LLL(P) = \{3124,3142,3412\}$.

% [DG][v41] Added footnote about the definition of linear extension.
% Is it any good? (I don't want to leave this lingo undefined, but I am
% not good at defining things :/ .)

% [DG][v42] Shortened as discussed in email.

\begin{theorem}
\label{Stanley's-P-partition-theorem}
For any labelled poset $P$,
\[
F_P(\xx) = \sum_{w \in \LLL(P)} F_w(\xx).
\]
\end{theorem}
\begin{proof}
We give Gessel's proof \cite[Thm. 1]{Gessel}, via induction on the
number of pairs $i,j$ which are incomparable in $P$.  When this quantity
is $0$, then $P$ is itself a linear order $w$, so that $\LLL(P)=\{w\}$
and there is nothing to prove.

In the inductive step,
let $i,j$ be incomparable elements.  Consider the two
posets $P_{i<j}$ and $P_{j<i}$ which are obtained
from $P$ by adding in an order relation between $i$ and $j$, and then taking
the transitive closure;  it is not hard to see that these transitive closures
cannot contain a cycle, so that these really do define two posets.
The result then follows by induction applied to  $P_{i<j}, P_{j<i}$, once
one notices that
$
\LLL(P) = \LLL(P_{i<j}) \sqcup  \LLL(P_{j<i})
$
since every linear extension $w$ of $P$ either has $i$ before $j$ or vice-versa, and
$
\AAA(P) = \AAA(P_{i<j}) \sqcup \AAA(P_{j<i})
$
since, assuming that $i <_\ZZ j$ without loss of generality, every
$f$ in $\AAA(P)$ either satisfies $f(i) \leq f(j)$ or $f(i) > f(j)$.
\end{proof}

\begin{example}
To illustrate the induction in the above proof, consider the poset
$P$ from Example~\ref{labelled-poset-example},
having $\LLL(P)=\{3124,3142,3412\}$.  Then choosing as incomparable pair
$(i,j)=(1,4)$, one has
\[
\xymatrix{
        &4& &2\\
P_{i<j}=& &1\ar@{-}[ur]\ar@{-}[ul]& \\
        & &3\ar@{-}[u]&}
\qquad
\xymatrix{
f(4)& &f(2)\\
 &f(1)\ar@{-}[ur]^{\leq} \ar@{-}[ul]^{\leq}& \\
 &f(3)\ar@{-}[u]^{<}&}
\qquad
\xymatrix{ \\  \text{, thus } \LLL(P_{i<j})=\{3124,3142\}\\ }
\]
and % \vskip.1in
\[
\xymatrix{
 &2\\
P_{j<i}=&1\ar@{-}[u]\\
 &4\ar@{-}[u]\\
 &3\ar@{-}[u]}
\qquad \qquad \qquad \qquad \qquad %\qquad
\xymatrix{
  f(2)\\
  f(1)\ar@{-}[u]^{\leq} \\
  f(4)\ar@{-}[u]^{<} \\
  f(3)\ar@{-}[u]^{\leq}}
\qquad
\xymatrix{ \\  \text{, thus } \LLL(P_{j<i})=\{3412\} .\\ }
\]
\end{example}

\begin{exercise}
\label{exe.P-partition.stanley.altproof}
Give an alternative proof for
Theorem~\ref{Stanley's-P-partition-theorem}.

[\textbf{Hint:} For every
$f : P \rightarrow \left\{ 1, 2, 3, \ldots \right\}$, we
can define a binary relation $\prec_{f}$ on the set $P$ by
letting $i \prec_{f} j$ hold if and only if
\[
\left( f \left( i \right) < f \left( j \right)
\text{ or }
\left( f \left( i \right) = f\left( j \right)
  \text{ and } i <_{\ZZ} j \right) \right)  .
\]
Show that this binary relation $\prec_{f}$ is (the smaller
relation of) a total order. When $f$ is a $P$-partition, then
endowing the set $P$ with this total order yields a linear
extension of $P$. Use this to show that the set
$\mathcal{A}\left( P \right)$ is the union of its
disjoint subsets $\mathcal{A}\left( w \right)$ with
$w \in \mathcal{L}\left( P \right)$.]
\end{exercise}

% [DG][v50] Added above exercise. It is an alternative to the
% "induction over the number of pairs of incomparable elements"
% argument that you and Gessel use.

Various other properties of the quasisymmetric functions
$F_P\left(\xx\right)$ are studied, e.g., in
\cite{McNamaraWard}.

We next wish to describe the
structure maps for the Hopf algebra $\Qsym$ in the basis $\{L_\alpha\}$
of fundamental quasisymmetric functions.  For this purpose, two more
definitions are useful.

\begin{definition}
\label{def.compositions.near-conc}
Given two nonempty compositions
$\alpha=(\alpha_1,\ldots,\alpha_\ell)$ and
$\beta=(\beta_1,\ldots,\beta_m)$, their
\emph{near-concatenation}\index{near-concatenation of compositions}
is
\[
\alpha \odot \beta:=
(\alpha_1,\ldots,\alpha_{\ell-1},\alpha_\ell+\beta_1,\beta_2,\ldots,\beta_m) .
\]
For example, the figure below depicts for $\alpha=(1,3,3)$
(black squares) and $\beta=(4,2)$ (white squares) the concatenation and
near-concatenation as ribbons:\footnote{The ribbons are drawn with their
boxes spaced out in order to facilitate counting.}
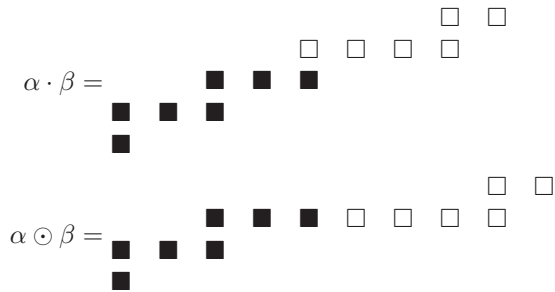
\begin{align*}
\Rib\left(\alpha \cdot \beta\right)&=\ %
\begingroup
\setlength\arraycolsep{0.2pc} % Make the spacing equal in rows and in columns.
\begin{matrix}
   &   &   &   &   &   &   &\sq&\sq\\
   &   &   &   &\sq&\sq&\sq&\sq&   \\
   &   &\bsq&\bsq&\bsq&   &   &   &   \\
\bsq&\bsq&\bsq&   &   &   &   &   &   \\
\bsq&   &   &   &   &   &   &   &
\end{matrix}
\endgroup
\\
\Rib\left(\alpha \odot \beta\right)&=\ %
\begingroup
\setlength\arraycolsep{0.2pc} % Make the spacing equal in rows and in columns.
\begin{matrix}
   &   &   &   &   &   &   &   &\sq&\sq\\
   &   &\bsq&\bsq&\bsq&\sq&\sq&\sq&\sq&\\
\bsq&\bsq&\bsq&   &   &   &   &   &   &\\
\bsq&   &   &   &   &   &   &   &   &
\end{matrix}
\endgroup
\end{align*}

Lastly, given $\alpha$ in $\Comp_n$, let $\omega(\alpha)$ be the unique composition in $\Comp_n$ whose partial sums $D(\omega(\alpha))$
form the complementary set within $[n-1]$ to the partial sums $D(\rev(\alpha))$;  alternatively, one can check this means that the ribbon for $\omega(\alpha)$ is obtained from that of $\alpha$
by conjugation or transposing, that is, if $\Rib\left(\alpha\right) =\lambda/\mu$ then
$\Rib\left(\omega(\alpha)\right) =\lambda^t/\mu^t$.
E.g. if $\alpha=(4,2,2)$ so that $n=8$,
then $\rev(\alpha)=(2,2,4)$ has $D(\rev(\alpha))=\{2,4\} \subset [7]$,
complementary to the set $\{1,3,5,6,7\}$ which are the partial sums
for $\omega(\alpha)=(1,2,2,1,1,1)$, and the ribbon diagrams of $\alpha$
and $\omega(\alpha)$ are
\[
\begin{matrix}
\Rib\left(\alpha\right) =\ &
\begingroup
\setlength\arraycolsep{0.2pc} % Make the spacing equal in rows and in columns.
\begin{matrix}
   &   &   &   &\sq&\sq\\
   &   &   &\sq&\sq& \\
\sq&\sq&\sq&\sq&   & \\
\end{matrix}
\endgroup
\qquad \text{ and } \qquad
\Rib\left(\omega(\alpha)\right) =\ &
\begingroup
\setlength\arraycolsep{0.2pc} % Make the spacing equal in rows and in columns.
\begin{matrix}
   &   &\sq\\
   &   &\sq\\
   &   &\sq\\
   &\sq&\sq\\
\sq&\sq&   \\
\sq&   &   \\
\end{matrix}
\endgroup
\end{matrix}
\]

\end{definition}

\begin{proposition}
\label{fundamental-basis-structure-maps}
The structure maps for the Hopf algebra $\Qsym$ in the basis $\{L_\alpha\}$
of fundamental quasisymmetric functions are as follows:
\begin{align}
\Delta L_\alpha &= \sum\limits_{\substack{(\beta,\gamma):\\
\beta \cdot \gamma=\alpha \text{ or }
\beta \odot \gamma=\alpha}} L_\beta \otimes L_\gamma ,
\label{Qsym-coproduct-on-fundamentals}\\
L_\alpha L_\beta &= \sum_{w \in \shuf{w_\alpha}{w_\beta}} L_{\gamma(w)} ,
\label{Qsym-product-on-fundamentals}\\
S(L_\alpha)&= (-1)^{|\alpha|} L_{\omega(\alpha)} .
\label{Qsym-antipode-on-fundamentals}
\end{align}
Here we are making use of the following notations in
\eqref{Qsym-product-on-fundamentals} (recall also
Definition~\ref{shuffles}):
\begin{enumerate}
\item[$\bullet$] A \dfn{labelled linear order}
will mean a labelled poset $P$ whose order
$<_{P}$ is a total order. We will identify any labelled linear order $P$ with
the word (over the alphabet $\ZZ$) obtained by
writing down the elements of $P$ in increasing order (with respect to the
total order $<_{P}$). This way, every word (over the alphabet $\ZZ$)
which has no two equal letters becomes identified with a labelled linear order.
\item[$\bullet$] $w_{\alpha}$ is any labelled linear order with underlying set
$\left\{  1,2,\ldots,\left\vert \alpha\right\vert \right\}  $ such that
$\Des\left(  w_{\alpha}\right)  =D\left(  \alpha\right)  $.
\item[$\bullet$] $w_{\beta}$ is any labelled linear order with underlying set $\left\{
\left\vert \alpha\right\vert +1,\left\vert \alpha\right\vert +2,\ldots,\left\vert
\alpha\right\vert +\left\vert \beta\right\vert \right\}  $ such that
$\Des\left(  w_{\beta}\right)  =D\left(  \beta\right)  $.
\item[$\bullet$] \dfn{$\gamma(w)$} is the unique composition of
$\left|\alpha\right| + \left|\beta\right|$ with $D(\gamma(w))=\Des(w)$.
\end{enumerate}
(The right hand side of \eqref{Qsym-product-on-fundamentals} is
to be read as a sum over all $w$, for a fixed choice of $w_\alpha$
and $w_\beta$.)
\end{proposition}

% [DG][v63] Added the preceding parenthetical sentence.

\noindent
At first glance the formula \eqref{Qsym-coproduct-on-fundamentals} for $\Delta L_\alpha$
might seem more complicated than the formula of
Proposition~\ref{Qsym-coproduct-on-monomials}
for $\Delta M_\alpha$.  However, it is equally simple
when viewed in terms of ribbon
diagrams:  it cuts the ribbon diagram $\Rib\left(\alpha\right)$
into two smaller ribbons $\Rib\left(\beta\right)$ and
$\Rib\left(\gamma\right)$, in all $|\alpha|+1$ possible ways,
via \emph{horizontal} cuts ($\beta \cdot \gamma = \alpha$)
or \emph{vertical} cuts ($\beta \odot \gamma = \alpha$).
For example,
\[
\begingroup
\setlength\arraycolsep{0.2pc} % Make the spacing equal in rows and in columns.
\begin{array}{rcccccl}
& \Delta L_{(3,2)} \\
&=1 \otimes L_{(3,2)}
&+L_{(1)} \otimes L_{(2,2)}
&+L_{(2)} \otimes L_{(1,2)}
&+L_{(3)} \otimes L_{(2)}
&+L_{(3,1)} \otimes L_{(1)}
&+L_{(3,2)} \otimes 1 .\\
    &\begin{matrix}  & & \sq & \sq \\ \underline{\sq} & \sq & \sq & \end{matrix}
    &\begin{matrix}  & & \sq & \sq \\ \sq & |\sq & \sq & \end{matrix}
    &\begin{matrix}  & & \sq & \sq \\ \sq & \sq & |\sq & \end{matrix}
    &\begin{matrix}  & & \underline{\sq} & \sq \\ \sq & \sq & \sq & \end{matrix}
    &\begin{matrix}  & & \sq & |\sq \\ \sq & \sq & \sq & \end{matrix}
    &\begin{matrix}  & & \sq & \overline{\sq} \\ \sq & \sq & \sq & \end{matrix}
\end{array}
\endgroup
\]

% [DG] I have edited the first and the last of these pictures so as
% to show that $1 \otimes L_{(3,2)}$ and $L_{(3,2)} \otimes 1$ correspond
% to horizontal (rather than vertical) cuts.

\begin{example}
To multiply $L_{(1,1)} L_{(2)}$, one could pick $w_\alpha = 21$ and $w_\beta=34$,
and then
\[
\begingroup
\setlength\arraycolsep{0.2pc} % Smaller spacing.
\begin{array}{rcccccccccccl}
L_{(1,1)} L_{(2)}
=\sum\limits_{w \in \shuf{21}{34}} L_{\gamma(w)}
&=&L_{\gamma(2134)}
&+&L_{\gamma(2314)}
&+&L_{\gamma(3214)}
&+&L_{\gamma(2341)}
&+&L_{\gamma(3241)}
&+&L_{\gamma(3421)}\\
&=&L_{(1,3)}
&+&L_{(2,2)}
&+&L_{(1,1,2)}
&+&L_{(3,1)}
&+&L_{(1,2,1)}
&+&L_{(2,1,1)}.
\end{array}
\endgroup
\]
\end{example}

Before we prove Proposition~\ref{fundamental-basis-structure-maps},
we state a simple lemma:

\begin{lemma}
\label{lem.QSym.ppar.djun}
Let $Q$ and $R$ be two labelled posets whose underlying
sets are disjoint.
Let $Q \sqcup R$ be the disjoint union of these posets
$Q$ and $R$; this is again a labelled poset.
Then,
\[
F_Q \left( \xx \right) F_R \left( \xx \right)
= F_{Q \sqcup R} \left( \xx \right) .
\]
\end{lemma}

\begin{proof} %[Proof of Lemma~\ref{lem.QSym.ppar.djun}.]
We identify the underlying set of $Q \sqcup R$ with $Q \cup R$
(since the sets $Q$ and $R$ are already disjoint).
If $f : Q \sqcup R \to \left\{1, 2, 3, \ldots\right\}$ is a
$Q \sqcup R$-partition, then its restrictions $f\mid_Q$ and
$f\mid_R$ are a $Q$-partition and an $R$-partition,
respectively.
Conversely, any pair of a $Q$-partition and an $R$-partition
can be combined to form a $Q \sqcup R$-partition.
Thus, there is a bijective correspondence between the addends
in the expanded sum
$F_Q \left( \xx \right) F_R \left( \xx \right)$ and the
addends in $F_{Q \sqcup R} \left( \xx \right)$.
\end{proof}

% [DG][v70] Added the above lemma (mostly in order to be
% able to cite it separately from the proof in which it
% appears).

\begin{proof}[Proof of Proposition~\ref{fundamental-basis-structure-maps}.]
To prove formula \eqref{Qsym-coproduct-on-fundamentals}
for $\alpha$ in $\Comp_n$, note that
\begin{equation}
\Delta L_\alpha = L_\alpha(\xx,\yy)
= \sum_{k=0}^n
\sum\limits_{\substack{
1 \leq i_1 \leq \cdots \leq i_k,\\
1 \leq i_{k+1} \leq \cdots \leq i_n:\\
i_r < i_{r+1} \text{ for }r \in D(\alpha) \setminus \left\{k\right\}
}}
 x_{i_1} \cdots x_{i_k} \cdot y_{i_{k+1}} \cdots y_{i_n}
\label{pf.fundamental-basis-structure-maps.DeltaLalpha}
\end{equation}
by Proposition~\ref{Qsym-technical-lemma-Lxy}
(where we identify $\Qsym \otimes \Qsym$ with a
$\kk$-subalgebra of $R\left(\xx, \yy\right)$ by means of
the embedding
$\Qsym \otimes \Qsym \overset{\cong}{\rightarrow}
\Qsym\left(\xx\right) \otimes \Qsym\left(\yy\right)
\hookrightarrow R\left(\xx, \yy\right)$ as in the definition
of the comultiplication on $\Qsym$).
One then realizes that the inner sums corresponding
to values of $k$ that lie (resp. do not lie) in $D(\alpha) \cup \left\{0,n\right\}$
correspond to the terms $L_\beta(\xx) L_\gamma(\yy)$
for pairs $(\beta,\gamma)$ in which
$\beta \cdot \gamma =\alpha$ (resp. $\beta \odot \gamma =\alpha$).

For formula \eqref{Qsym-product-on-fundamentals}, let
$P$ be the labelled poset which is the disjoint union of linear orders $w_\alpha, w_\beta$.
Then
\[
L_\alpha L_\beta
=F_{w_\alpha}(\xx) F_{w_\beta}(\xx)
=F_P(\xx) = \sum_{w \in \LLL(P)} F_w(\xx)
=\sum_{w \in \shuf{w_\alpha}{w_\beta}} L_{\gamma(w)}
\]
where the first equality used Proposition~\ref{fundamental-P-partition-prop},
the second equality comes from Lemma~\ref{lem.QSym.ppar.djun},
%the definition of a $P$-partition,
the third equality from Theorem~\ref{Stanley's-P-partition-theorem},
and the fourth from the equality $\LLL(P)=\shuf{w_\alpha}{w_\beta}$.

To prove formula \eqref{Qsym-antipode-on-fundamentals}, compute
using Theorem~\ref{Qsym-antipode-on-monomials} that
\[
%\begin{aligned}
S(L_\alpha)
=\sum_{\beta \text{ refining } \alpha} S(M_\beta)
=\sum\limits_{\substack{(\beta,\gamma):\\
         \beta \text{ refines } \alpha,\\
         \gamma \text{ coarsens }\rev(\beta)}} (-1)^{\ell(\beta)} M_\gamma
=\sum_{\gamma} M_\gamma
     \sum_{\beta}  (-1)^{\ell(\beta)}
%&=(-1)^{|\alpha|}
%    \sum_{\gamma: D(\alpha) \cap D(\rev(\gamma))= \varnothing} M_\gamma \\
%&=(-1)^{|\alpha|} \sum_{\gamma: D(\gamma) \supset D(\omega(\alpha))} M_\gamma% \\
%&=(-1)^{|\alpha|} L_{\omega(\alpha)} \\
%\end{aligned}
\]
in which the last inner sum is over $\beta$ for which
\[
D(\beta) \supset D(\alpha) \cup D(\rev(\gamma)).
\]
The alternating signs make such inner sums vanish
unless they have only the single
term where $D(\beta)=[n-1]$ (that is, $\beta=(1^n)$).  This happens
exactly when
$
D(\rev(\gamma)) \cup D(\alpha)=[n-1]
$
or equivalently, when $D(\rev(\gamma))$ contains the complement
of $D(\alpha)$, that is, when $D(\gamma)$ contains the complement of
$D(\rev(\alpha))$, that is, when
$\gamma$ refines $\omega(\alpha)$.  Thus
\[
S(L_\alpha)=
\sum\limits_{\substack{\gamma \in \Comp_n:\\\gamma \text{ refines }\omega(\alpha)}}
M_\gamma  \cdot (-1)^n
=(-1)^{|\alpha|} L_{\omega(\alpha)}.
\]
\end{proof}

The antipode formula \eqref{Qsym-antipode-on-fundamentals}
for $L_\alpha$ leads to a general interpretation for the antipode
of $\Qsym$ acting on $P$-partition enumerators $F_P(\xx)$.

\begin{definition}
Given a labelled poset $P$ on $\{1,2,\ldots,n\}$, let the
\emph{opposite}\index{opposite labelled poset}
or \emph{dual}\index{dual labelled poset}
labelled poset \dfn{$P^{\opp}$}
be the labelled poset on $\{1,2,\ldots,n\}$ that has
$i <_{P^{\opp}} j$ if and only if $j <_P i$.
\end{definition}
For example,
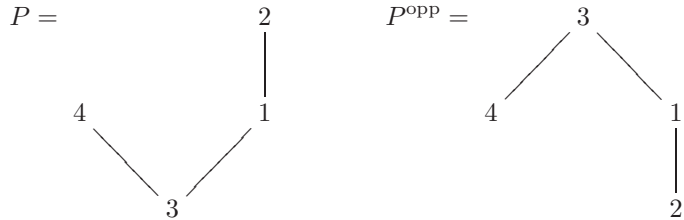
\[
P=
\xymatrix{
 & &2\\
4& &1\ar@{-}[u]\\
 &3\ar@{-}[ur]\ar@{-}[ul]&
} \qquad \qquad
P^{\opp}=
\xymatrix{
    &3\ar@{-}[dl]\ar@{-}[dr]&\\
  4& & 1\ar@{-}[d]\\
    & &2
}
\]

The following observation is straightforward.

\begin{proposition}
\label{rotating-ribbons-is-reversing-permutations}
When $P$ is a linear order
corresponding to some permutation $w=(w_1,\ldots,w_n)$ in $\Symm_n$,
then $w^{\opp}=w w_0$ where $w_0 \in \Symm_n$ is the permutation that
swaps $i \leftrightarrow n+1-i$ (this is the so-called \dfn{longest
permutation}, thus named due to it having the highest ``Coxeter
length'' among all permutations in $\Symm_n$).  Furthermore, in this situation
one has $F_w(\xx)=L_\alpha$, that is, $\Des(w) = D(\alpha)$ if and only if
$\Des(w^{\opp})= D(\omega(\alpha))$, that is $F_{w^{\opp}}(\xx) = L_{\omega(\alpha)}$.
Thus,
\[
S(F_w(\xx)) =(-1)^n F_{w^{\opp}}(\xx).
\]
\end{proposition}

For example, given the compositions considered earlier:
\[
\alpha =(4,2,2)=
\begingroup
\setlength\arraycolsep{0.2pc} % Make the spacing equal in rows and in columns.
\begin{matrix}
   &   &   &   &\sq&\sq\\
   &   &   &\sq&\sq& \\
\sq&\sq&\sq&\sq&   & \\
\end{matrix}
\endgroup
\qquad \text{ and } \qquad
\omega(\alpha) =(1,2,2,1,1,1)=
\begingroup
\setlength\arraycolsep{0.2pc} % Make the spacing equal in rows and in columns.
\begin{matrix}
   &   &\sq\\
   &   &\sq\\
   &   &\sq\\
   &\sq&\sq\\
\sq&\sq&   \\
\sq&   &   \\
\end{matrix}
\endgroup
\]
if one picks $w=1235\cdot47\cdot68$
(with descent positions marked by dots) having
$\Des(w)=\{4,6\}=D(\alpha)$, then
$w^{\opp}=ww_0=8\cdot67\cdot45\cdot3\cdot2\cdot1$
has $\Des(w^{\opp})=\{1,3,5,6,7\}=D(\omega(\alpha))$.

\begin{corollary}
\label{antipode-in-Qsym-as-reciprocity-corollary}
For any labelled poset $P$ on $\{1,2,\ldots,n\}$, one has
\[
S \left( F_P(\xx) \right)= (-1)^n F_{P^{\opp}}(\xx).
\]
\end{corollary}
\begin{proof}
Since $S$ is linear, one can apply Theorem~\ref{Stanley's-P-partition-theorem}
and Proposition~\ref{rotating-ribbons-is-reversing-permutations},
obtaining
\[
S \left( F_P(\xx) \right)=
\sum_{w \in \LLL(P)} S( F_w(\xx) )
=\sum_{w \in \LLL(P)} (-1)^{n} F_{w^{\opp}}(\xx)=
(-1)^n F_{P^{\opp}}(\xx),
\]
as $\LLL(P^{\opp}) = \{w^{\opp}: w \in \LLL(P)\}$.
\end{proof}

\begin{remark}
Malvenuto and Reutenauer, in \cite[Theorem 3.1]{MalvenutoReutenauer-plethysm},
prove an even more general antipode formula, which encompasses our
Corollary~\ref{antipode-in-Qsym-as-reciprocity-corollary},
Proposition~\ref{rotating-ribbons-is-reversing-permutations},
Theorem~\ref{Qsym-antipode-on-monomials}
and \eqref{Qsym-antipode-on-fundamentals}.
See \cite[Theorem 4.2]{Grinberg-dp} for a restatement and
a self-contained proof of this theorem (and
\cite[Theorem 4.7]{Grinberg-dp} for an even further
generalization).
\end{remark}

We remark on a special case of Corollary~\ref{antipode-in-Qsym-as-reciprocity-corollary}
to which we alluded earlier, related to skew Schur functions.

\begin{corollary}
\label{antipode-on-skew-Schur-corollary}
In $\Lambda$, the action of $\omega$ and
the antipode $S$ on skew Schur functions $s_{\lambda/\mu}$
are as follows:
\begin{align}
\omega(s_{\lambda/\mu})&=s_{\lambda^t/\mu^t} ,\label{omega-on-skew-Schur}\\
S(s_{\lambda/\mu})&=(-1)^{|\lambda/\mu|} s_{\lambda^t/\mu^t} \label{antipode-on-skew-Schur}.
\end{align}
\end{corollary}
\begin{proof}
Given a skew shape $\lambda/\mu$, one can always create a labelled
poset $P$ which is its \dfn{skew Ferrers poset}, together with one of many
\emph{column-strict labellings}\index{column-strict labelling},
in such a way that
$F_P(\xx)=s_{\lambda/\mu}(\xx)$.  An example is shown here
for $\lambda/\mu=(4,4,2)/(1,1,0)$:
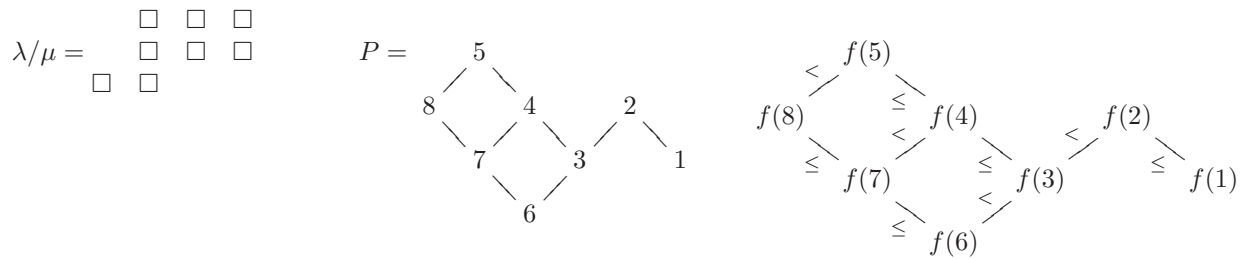
\[
\lambda/\mu=
\begingroup
\setlength\arraycolsep{0.2pc} % Make the spacing equal in rows and in columns.
\begin{matrix}
   &\sq&\sq&\sq\\
   &\sq&\sq&\sq\\
\sq&\sq&  &
\end{matrix}
\endgroup
\qquad \qquad
P=\xymatrix@C=8pt@R=8pt{
                      &5&                         &       &       & \\
8\ar@{-}[ur]&       &4\ar@{-}[ul]&       &2& \\
                      &7\ar@{-}[ul]\ar@{-}[ur]&       &3\ar@{-}[ul]\ar@{-}[ur]&       &1\ar@{-}[ul]\\
       &       &6\ar@{-}[ul]\ar@{-}[ur]&       &       &   \\
}
\qquad
\xymatrix@C=8pt@R=8pt{
                      &f(5)&                         &       &       & \\
f(8)\ar@{-}[ur]^{<}&       &f(4)\ar@{-}[ul]^{\leq}&       &f(2)& \\
                      &f(7)\ar@{-}[ul]^{\leq}\ar@{-}[ur]^{<}&       &f(3)\ar@{-}[ul]^{\leq}\ar@{-}[ur]^{<}&       &f(1)\ar@{-}[ul]^{\leq}\\
       &       &f(6)\ar@{-}[ul]^{\leq}\ar@{-}[ur]^{<}&       &       &   \\
}
\]
The general definition is as follows: Let $P$ be the set of all boxes of the skew diagram
$\lambda / \mu$. Label these boxes by the numbers $1,2,\ldots,n$ (where $n=\left\vert
\lambda / \mu\right\vert $) row by row from bottom to top (reading every
row from left to right), and then define an order relation $<_{P}$ on $P$ by
requiring that every box be smaller (in $P$) than its right neighbor and smaller (in
$P$) than its lower neighbor.
It is not hard to see that in this situation, $F_{P^{\opp}}(\xx)=\sum_T \xx^{\cont(T)}$ as $T$ ranges over all
\emph{reverse semistandard tableaux}\index{reverse semistandard tableau}
or \emph{column-strict plane partitions}\index{column-strict plane partition}
of $\lambda^t/\mu^t$:
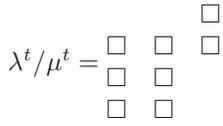
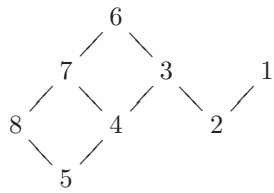
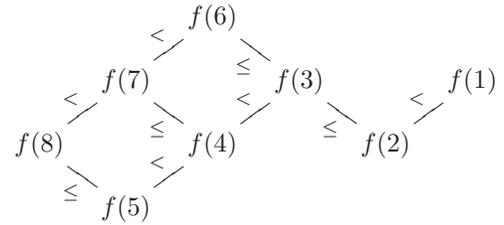
\[
\lambda^t/\mu^t=
\begingroup
\setlength\arraycolsep{0.2pc} % Make the spacing equal in rows and in columns.
\begin{matrix}
   &   &\sq\\
\sq&\sq&\sq&\\
\sq&\sq&  &\\
\sq&\sq&  &
\end{matrix}
\endgroup
\qquad \qquad
P^{\opp}=\xymatrix@C=8pt@R=8pt{
 & &6& & &\\
 &7\ar@{-}[ur]& &3\ar@{-}[ul]& &1\\
8\ar@{-}[ur]& &4\ar@{-}[ur]\ar@{-}[ul]& &2\ar@{-}[ur]\ar@{-}[ul]& \\
 &5\ar@{-}[ur]\ar@{-}[ul]& & & &
}
\qquad
\xymatrix@C=8pt@R=8pt{
 & &f(6)& & &\\
 &f(7)\ar@{-}[ur]^{<}& &f(3)\ar@{-}[ul]^{\leq}& &f(1)\\
f(8)\ar@{-}[ur]^{<}& &f(4)\ar@{-}[ur]^{<}\ar@{-}[ul]^{\leq}& &f(2)\ar@{-}[ur]^{<}\ar@{-}[ul]^{\leq}& \\
 &f(5)\ar@{-}[ur]^{<}\ar@{-}[ul]^{\leq}& & & &
}
\]
But this means that $F_{P^{\opp}}(\xx)=s_{\lambda^t/\mu^t}(\xx)$, since
the fact that skew Schur functions lie in $\Lambda$ implies that they
can be defined either as generating functions for
column-strict tableaux or reverse semistandard tableaux;
see Remark~\ref{other-linear-orders-remark} above,
or \cite[Prop. 7.10.4]{Stanley}.

% [DG] I have added a semi-formal definition of the poset $P$. It is
% somewhat clumsy.

% [DG][v3] Replaced "reverse column-strict tableaux" by "reverse
% semistandard tableaux" to match the notations in a footnote in Chapter 2.
% Maybe I should have done it the other way round...

Thus we have
\begin{align*}
F_P(\xx)&=s_{\lambda/\mu}(\xx) ,\\
F_{P^{\opp}}(\xx)&=s_{\lambda^t/\mu^t}(\xx).
\end{align*}
Corollary~\ref{cor.bialg-mor-is-Hopf}
tells us that the antipode for $\Qsym$ must specialize
to the antipode for $\Lambda$ (see also Remark~\ref{antipode-comparison-remark} below), so \eqref{antipode-on-skew-Schur}
is a special case of
Corollary~\ref{antipode-in-Qsym-as-reciprocity-corollary}.
Then \eqref{omega-on-skew-Schur}
follows from the relation \eqref{antipode-versus-omega-in-Sym} that
$S(f)=(-1)^n \omega(f)$ for $f$ in $\Lambda_n$.
\end{proof}

\begin{remark}
Before leaving $P$-partitions temporarily, we mention two
open questions about them.

The first is a conjecture of Stanley from his thesis \cite{Stanley-thesis}.
As mentioned in the proof of Corollary~\ref{antipode-on-skew-Schur-corollary},
each skew Schur function $s_{\lambda/\mu}(\xx)$ is a
special instance of $P$-partition enumerator $F_P(\xx)$.

\begin{conjecture}
A labelled poset $P$ has $F_P(\xx)$ symmetric, and not just quasisymmetric,
\emph{if and only if} $P$ is a column-strict labelling of some
skew Ferrers poset $\lambda/\mu$.
\end{conjecture}
\noindent
A somewhat weaker result in this direction
was proven by Malvenuto in her thesis \cite[Thm. 6.4]{Malvenuto}, showing
that if a labelled poset $P$ has the stronger property that its set of linear extensions
$\LLL(P)$ is a union of \emph{plactic} or \emph{Knuth equivalence classes}, then
$P$ must be a column-strict labelling of a skew Ferrers poset.

The next question is due to P. McNamara, and is suggested by the
obvious factorizations of $P$-partition enumerators
$F_{P_1 \sqcup P_2}(\xx) = F_{P_1}(\xx) F_{P_2}(\xx)$ (Lemma~\ref{lem.QSym.ppar.djun}).

\begin{question}
\label{McNamara's-question}
If $\kk$ is a field,
does a \emph{connected} labelled poset $P$ always have
$F_P(\xx)$ \emph{irreducible} within the ring $\Qsym$?
\end{question}

The phrasing of this question requires further comment.
It is assumed here that $\xx=(x_1,x_2,\ldots)$ is infinite;
for example when $P$ is a $2$-element chain labelled ``against
the grain'' (i.e., the bigger element of the chain has the
smaller label),
then $F_P(\xx)=e_2(\xx)$ is irreducible, but its specialization
to two variables $\xx=(x_1,x_2)$ is $e_2(x_1,x_2)=x_1 x_2$,
which is reducible.  If one wishes to work in finitely many variables
$\xx=(x_1,\ldots,x_m)$ one can perhaps assume that $m$ is at least
$|P|+1$.

% [DG] I specified the labelling in the example.

When working in $\Qsym=\Qsym(\xx)$ in
infinitely many variables, it is perhaps not so clear where factorizations occur.
For example, if $f$ lies in $\Qsym$ and factors
$f=g \cdot h$ with $g, h$ in $R(\xx)$, does this imply that $g, h$ also lie in $\Qsym$?
The answer is ``Yes'' (for $\kk = \ZZ$), but this is not obvious, and
was proven by P. Pylyavskyy in \cite[Chap. 11]{Pylyavskyy-thesis}.

One also might wonder whether $\Qsym_\ZZ$ is a unique factorization domain, but
this follows from the result of M. Hazewinkel (\cite{Hazewinkel}
and \cite[Thm. 6.7.5]{HazewinkelGubareniKirichenko}, and
Theorem~\ref{thm.QSym.lyndon} further below) who
proved a conjecture of Ditters that $\Qsym_\ZZ$
is a polynomial algebra;  earlier
Malvenuto and Reutenauer \cite[Cor. 2.2]{MalvenutoReutenauer} had shown
that $\Qsym_\QQ$ is a polynomial algebra.  In fact,
one can find polynomial generators $\{P_{\alpha}\}$ for $\Qsym_\QQ$ as a
subset of the dual basis to the $\QQ$-basis $\{\xi_\alpha\}$ for $\Nsym_\QQ$
which comes from taking
products $\xi_{\alpha}:= \xi_{\alpha_1} \cdots \xi_{\alpha_\ell}$
of the elements $\{\xi_n\}$ defined in Remark~\ref{second-antipode-proof} below.  Specifically,
one takes those $P_\alpha$ for which the composition $\alpha$ is a \emph{Lyndon composition};
see the First proof of Proposition~\ref{prop.QSym.lyndon} for a mild variation
on this construction.

Hazewinkel's proof \cite[Thm. 6.7.5]{HazewinkelGubareniKirichenko}
of the polynomiality of $\Qsym_\ZZ$ also shows that $\Qsym$ is a
polynomial ring over $\Lambda$ (see
Corollary~\ref{cor.QSym.lyndon.free}); in particular, this yields that
$\Qsym$ is a free $\Lambda$-module.\footnote{The latter statement
has an analogue in finitely many indeterminates, proven by Lauve
and Mason in \cite[Corollary 13]{LauveMason-stable}: The
quasisymmetric functions $\Qsym\left(\left\{x_i\right\}_{i\in I}\right)$
are free as a $\Lambda\left(\left\{x_i\right\}_{i\in I}\right)$-module
for any totally ordered set $I$, infinite or not. In the case of
finite $I$, this cannot be derived by Hazewinkel's arguments, as
the ring $\Qsym\left(\left\{x_i\right\}_{i\in I}\right)$ is not
in general a polynomial ring (e.g., when $\kk = \QQ$ and
$I = \left\{1, 2\right\}$, this ring is not even a UFD, as
witnessed by $\left(x_1^2 x_2\right) \cdot \left(x_1 x_2^2\right)
= \left(x_1 x_2\right)^3$).}

% [DG][v34] Added above paragraph on Lauve-Mason work.

% [DG][v36] Updated with references.

An affirmative answer to Question~\ref{McNamara's-question}
is known at least in the special case where $P$ is a connected column-strict labelling
of a skew Ferrers diagram, that is, when $F_P(\xx)=s_{\lambda/\mu}(\xx)$
for some connected skew diagram $\lambda/\mu$;  see \cite{BarekatRWilligenburg}.
\end{remark}

\subsection{\label{subsect.Lalpha-std}Standardization of $n$-tuples and the
fundamental basis}

Another equivalent description of the fundamental quasisymmetric functions
$L_{\alpha}$ (Lemma~\ref{lem.Lalpha-std.std-L} below) relies on the concept of
words and of their standardizations. We shall study words in detail in Chapter
\ref{sect.QSym.lyndon}; at this point, we merely introduce the few notions
that we will need:

\begin{definition}
\label{def.words.1}
We fix a totally ordered set $\mathfrak{A}$, which we call
the \dfn{alphabet}.

We recall that a \emph{word over $\mathfrak{A}$}\index{word}
is just a (finite) tuple of
elements of $\mathfrak{A}$. A word $\left(  w_{1},w_{2},\ldots,w_{n}\right)  $
can be written as $w_{1}w_{2}\cdots w_{n}$ when this incurs no ambiguity.

If $w\in\mathfrak{A}^{n}$ is a word and $i\in\left\{  1,2,\ldots,n\right\}  $,
then the \dfn{$i$-th letter} of $w$ means the $i$-th entry of the $n$-tuple
$w$. This $i$-th letter will be denoted by \dfn{$w_i$}.
\end{definition}

Our next definition relies on a simple fact about permutations and
words:\footnote{See Exercise \ref{exe.Lalpha-std.proofs} below for a proof of
Proposition \ref{prop.words.std-exists}.}

\begin{proposition}
\label{prop.words.std-exists}Let $w=\left(  w_{1},w_{2},\ldots,w_{n}\right)
\in\mathfrak{A}^{n}$ be any word. Then, there exists a unique permutation
$\sigma\in\Symm_n$ such that for every two elements $a$ and $b$ of
$\left\{  1,2,\ldots,n\right\}  $ satisfying $a<b$, we have $\left(
\sigma\left(  a\right)  <\sigma\left(  b\right)  \text{ if and only if }%
w_{a}\leq w_{b}\right)  $.
\end{proposition}

\begin{definition}
\label{def.words.std} Let $w\in\mathfrak{A}^{n}$ be any word. The unique
permutation $\sigma\in\Symm_n$ defined in Proposition
\ref{prop.words.std-exists} is called the \dfn{standardization} of $w$, and
is denoted by \dfn{$\operatorname*{std}w$}.
\end{definition}

\begin{example}
If $\mathfrak{A}$ is the alphabet $\left\{  1<2<3<\cdots\right\}  $, then
$\operatorname*{std}\left(  41211424\right)  $ is the permutation which is
written (in one-line notation) as $61423758$.
\end{example}

A simple method to compute the standardization of a word $w\in\mathfrak{A}%
^{n}$ is the following: Replace all occurrences of the smallest letter
appearing in $w$ by the numbers $1,2,\ldots,m_{1}$ (where $m_{1}$ is the
number of these occurrences); then replace all occurrences of the
second-smallest letter appearing in $w$ by the numbers $m_{1}+1,m_{1}%
+2,\ldots,m_{1}+m_{2}$ (where $m_{2}$ is the number of these occurrences), and
so on, until all letters are replaced by numbers.\footnote{Here, a number is
not considered to be a letter; thus, a number that replaces a letter will
always be left in peace afterwards.} The result is the standardization of $w$,
in one-line notation.

Another method to compute the standardization $\operatorname*{std}w$ of a word
$w=\left(  w_{1},w_{2},\ldots,w_{n}\right)  \in\mathfrak{A}^{n}$ is based on
sorting. Namely, consider the total order on the set $\mathfrak{A}%
\times \ZZ$ given by%
\[
\left(  a,i\right)  \leq\left(  b,j\right)  \text{ if and only if }\left(
\text{either }a<b\text{ or }\left(  a=b\text{ and }i\leq j\right)  \right)  .
\]
(In other words, two pairs in $\mathfrak{A}\times \ZZ$ are compared by
first comparing their first entries, and then, in the case of a tie, using the
second entries as tiebreakers.) Now, in order to compute $\operatorname*{std}%
w$, we sort the $n$-tuple $\left(  \left(  w_{1},1\right)  ,\left(
w_{2},2\right)  ,\ldots,\left(  w_{n},n\right)  \right)  \in\left(
\mathfrak{A}\times \ZZ\right)  ^{n}$ into increasing order (with respect
to the total order just described), thus obtaining a new $n$-tuple of the form
$\left(  \left(  w_{\tau\left(  1\right)  },\tau\left(  1\right)  \right)
,\left(  w_{\tau\left(  2\right)  },\tau\left(  2\right)  \right)
,\ldots,\left(  w_{\tau\left(  n\right)  },\tau\left(  n\right)  \right)
\right)  $ for some $\tau\in\Symm_n$; the standardization
$\operatorname*{std}w$ is then $\tau^{-1}$.

\begin{definition}
\label{def.QSym.gamma(sigma)} Let $n\in \NN$. Let $\sigma\in
\Symm_n$. Define a subset $\Des \sigma$ of $\left\{
1,2,\ldots,n-1\right\}  $ by
\[
\Des \sigma = \left\{  i\in\left\{  1,2,\ldots,n-1\right\}
\ \mid\ \sigma\left(  i\right)  >\sigma\left(  i+1\right)  \right\}  .
\]
(This is a particular case of the definition of $\Des w$ in
Exercise \ref{exe.shareshian-wachs.path}, if we identify $\sigma$ with the
$n$-tuple $\left(  \sigma\left(  1\right)  ,\sigma\left(  2\right)
,\ldots,\sigma\left(  n\right)  \right)  $. It is also a particular case of
the definition of $\Des w$ in Proposition
\ref{fundamental-P-partition-prop}, if we identify $\sigma$ with the total
order $\left(  \sigma\left(  1\right)  <\sigma\left(  2\right)  <\cdots
<\sigma\left(  n\right)  \right)  $ on the set $\left\{  1,2,\ldots,n\right\}
$.)

There is a unique composition $\alpha$ of $n$ satisfying $D\left(
\alpha\right)  = \Des \sigma$ (where $D\left(  \alpha\right)  $
is defined as in Definition \ref{def.QSym.D}). This composition will be
denoted by $\gamma\left(  \sigma\right)  $.
\end{definition}

The following lemma (equivalent to \cite[Lemma 9.39]{Reutenauer}) yields
another description of the fundamental quasisymmetric functions:

\begin{lemma}
\label{lem.Lalpha-std.std-L}Let $\mathfrak{A}$ denote the totally ordered set
$\left\{  1<2<3<\cdots\right\}  $ of positive integers. For each word
$w=\left(  w_{1},w_{2},\ldots,w_{n}\right)  \in\mathfrak{A}^{n}$, we define a
monomial $\xx_{w}$ in $\kk\left[  \left[  \xx\right]
\right]  $ by $\xx_{w}=x_{w_{1}}x_{w_{2}}\cdots x_{w_{n}}$.

Let $n\in \NN$ and $\sigma\in\Symm_n$. Then,%
\[
L_{\gamma\left(  \sigma\right)  }=\sum\limits_{\substack{w\in\mathfrak{A}%
^{n};\\\operatorname*{std}w=\sigma^{-1}}}\xx_{w}.
\]

\end{lemma}

\begin{exercise}
\label{exe.Lalpha-std.proofs}Prove Proposition \ref{prop.words.std-exists} and
Lemma \ref{lem.Lalpha-std.std-L}.
\end{exercise}

\subsection{The Hopf algebra $\Nsym$ dual to $\Qsym$}

We introduce here the (graded) dual Hopf algebra to $\Qsym$. This is
well-defined, as $\Qsym$ is connected graded of finite type.

\begin{definition}
Let $\Nsym:=\Qsym^o$\index{$\Nsym$}, with dual pairing
$\Nsym \otimes \Qsym \overset{(\cdot,\cdot)}{\longrightarrow} \kk$.
Let $\{H_\alpha\}$\index{$H_\alpha$}
be the $\kk$-basis of $\Nsym$ dual to the $\kk$-basis
$\{M_\alpha\}$ of $\Qsym$, so that
\[
(H_\alpha,M_\beta)=\delta_{\alpha,\beta}.
\]
When the base ring $\kk$ is not clear from the context, we write
$\Nsym_\kk$ in lieu of $\Nsym$.

The Hopf algebra $\Nsym$ is known as the \dfn{Hopf algebra of
noncommutative symmetric functions}\index{noncommutative symmetric
function}. Its study goes back to
\cite{GelfandKrobLascouxLeclercRetakhThibon}.
\end{definition}

% [DG][v7] Added the previous two sentences.

% [DG][v3] Here and in the following, all $\ZZ$s replaced by $\kk$s
% (except when necessary in the proof).

\begin{theorem}
\label{thm.NSym.free}
Letting $H_n:=H_{(n)}$ for $n=0,1,2,\ldots$, with $H_0=1$,
one has that
\begin{equation}
\Nsym \cong \kk\langle H_1,H_2,\ldots \rangle,
\label{eq.thm.NSym.free}
\end{equation}
the free associative (but not commutative) algebra on
generators $\{H_1,H_2,\ldots\}$ with coproduct determined
by\footnote{The abbreviated summation indexing
$\sum_{i+j=n} t_{i,j}$ used here is intended to mean
\[
\sum_{ \substack{ (i,j) \in \NN^2; \\ i+j=n } } t_{i,j}.
\] }
\begin{equation}
\label{Nsym-coproduct-on-H}
\Delta H_n = \sum_{i+j=n} H_i \otimes H_j.
\end{equation}
\end{theorem}
\begin{proof}
Since Proposition~\ref{Qsym-coproduct-on-monomials} asserts that
$\Delta M_\alpha
      = \sum_{(\beta,\gamma): \beta \cdot \gamma=\alpha}
       M_\beta \otimes M_\gamma$, and since $\{H_\alpha\}$ are
dual to $\{M_\alpha\}$, one concludes that for any compositions $\beta, \gamma$, one has
\[
H_\beta H_\gamma= H_{\beta \cdot \gamma} .
\]
Iterating this gives
\begin{equation}
\label{eq.Nsym.Halpha}
H_{\alpha}=H_{(\alpha_1,\ldots,\alpha_\ell)}
    =H_{\alpha_1} \cdots H_{\alpha_\ell}.
\end{equation}
Since the $H_\alpha$ are a $\kk$-basis for $\Nsym$, this shows
$\Nsym \cong \kk\langle H_1,H_2,\ldots \rangle$.

Note that $H_n=H_{(n)}$ is dual to $M_{(n)}$, so to
understand $\Delta H_n$, one should understand how  $M_{(n)}$ can appear
as a term in the product $M_{\alpha} M_{\beta}$.
By \eqref{Qsym-product-on-monomials} this occurs only
if $\alpha=(i), \beta=(j)$ where $i+j=n$, where
\[
M_{(i)} M_{(j)}= M_{(i+j)}+M_{(i,j)}+M_{(j,i)}
\]
(where the $M_{(i,j)}$ and $M_{(j,i)}$ addends have to be
disregarded if one of $i$ and $j$ is $0$).
By duality, this implies the formula \eqref{Nsym-coproduct-on-H}.
\end{proof}

\begin{corollary}
\label{cor.NSym.pi}
The algebra homomorphism defined by
\[
\begin{array}{rcl}
\Nsym &\overset{\pi}{\longrightarrow} &\Lambda ,\\
H_n &\longmapsto &h_n
\end{array}
\]
is a Hopf algebra surjection, and adjoint to the
inclusion $\Lambda \overset{i}{\hookrightarrow} \Qsym$
(with respect to the dual pairing
$\Nsym \otimes \Qsym \overset{(\cdot,\cdot)}{\longrightarrow} \kk$).
\end{corollary}

% [DG][v43] Added
% "(with respect to the dual pairing
% $\Nsym \otimes \Qsym \overset{(\cdot,\cdot)}{\longrightarrow} \kk$)."

\begin{proof}
As an algebra morphism, $\pi$ may be identified with the surjection
$T(V) \rightarrow \Sym(V)$ from the tensor algebra on a graded
free $\kk$-module $V$ with basis $\{H_1, H_2,\ldots\}$ to the
symmetric algebra on $V$, since
\begin{align*}
\Nsym &\cong \kk\langle H_1,H_2,\ldots\rangle ,\\
\Lambda &\cong \kk[h_1,h_2,\ldots] .
\end{align*}
As \eqref{Nsym-coproduct-on-H} and Proposition~\ref{symm-comultiplication-formulas-prop}(iii)
assert that
\begin{align*}
\Delta H_n &= \sum_{i+j=n} H_i \otimes H_j ,\\
\Delta h_n &= \sum_{i+j=n} h_i \otimes h_j ,\\
\end{align*}
this map $\pi$ is also a bialgebra morphism, and hence a Hopf morphism by
Corollary~\ref{cor.bialg-mor-is-Hopf}.

To check $\pi$ is adjoint to $i$, let $\lambda(\alpha)$
denote the partition which is the weakly decreasing rearrangement of the composition $\alpha$,
and note that the bases $\{H_\alpha\}$
of $\Nsym$ and $\{m_\lambda\}$ of $\Lambda$ satisfy
\[
(\pi(H_\alpha),m_\lambda)
=(h_{\lambda(\alpha)},m_\lambda)
=\left\{ \begin{matrix}
1 & \text{ if }\lambda(\alpha)=\lambda \\
0 & \text{ otherwise }
\end{matrix} \right\}
=\left( H_{\alpha},\sum_{\beta: \lambda(\beta)=\lambda} M_\beta \right)
=(H_\alpha, i(m_\lambda)).
\]
\end{proof}

\begin{remark}
\label{second-antipode-proof}
For those who prefer generating functions to sign-reversing involutions,
we sketch here  Malvenuto and Reutenauer's elegant proof
\cite[Cor. 2.3]{MalvenutoReutenauer} of the antipode
formula (Theorem~\ref{Qsym-antipode-on-monomials}).  One needs to
know that when $\QQ$ is a subring of $\kk$, and $A$ is a $\kk$-algebra
(possibly noncommutative),
in the ring of power series $A[[t]]$ where $t$ commutes with
all of $A$, one still has familiar facts, such as
\[
a(t) = \log b(t) \quad \text{ if and only if }\quad b(t) = \exp a(t)
\]
and whenever $a(t), b(t)$ commute in $A[[t]]$, one has
\begin{align}
\exp \left( a(t)+b(t) \right) &= \exp a(t) \exp b(t) , \label{exp-of-sum} \\
\log \left( a(t) b(t) \right) &= \log a(t) + \log b(t) . \label{log-of-product}
\end{align}
Start by assuming WLOG that $\kk=\ZZ$ (as $\Nsym_\kk
= \Nsym_\ZZ \otimes_\ZZ \kk$ in the general case).
Now, define in $\Nsym_\QQ=\Nsym \otimes_\ZZ \QQ$ the elements
$\{\xi_1,\xi_2,\ldots\}$\index{$\xi_n$} via generating functions
in $\Nsym_\QQ[[t]]$:
\begin{equation}
\begin{aligned}
\label{relation-definining-pi's}
\widetilde{H}(t)&:=\sum_{n \geq 0} H_n t^n,\\
\xi(t)&:=\sum_{n \geq 1} \xi_n t^n=\log \widetilde{H}(t) .
\end{aligned}
\end{equation}
One first checks that this makes each $\xi_n$ primitive, via
a computation in the ring $(\Nsym_\QQ \otimes \Nsym_\QQ)[[t]]$ (into which
we ``embed'' the ring $(\Nsym_{\QQ}[[t]]) \otimes_{\QQ[[t]]} (\Nsym_{\QQ}[[t]])$
via the canonical ring homomorphism from the latter into the former
\footnote{This ring homomorphism might fail to be injective, whence the ``embed''
stands in quotation marks. This does not need to worry us, since we will
not draw any conclusions in
$(\Nsym_{\QQ}[[t]]) \otimes_{\QQ[[t]]} (\Nsym_{\QQ}[[t]])$
from our computation.
\par
We are also somewhat cavalier with the notation $\Delta$: we use it both
for the comultiplication $\Delta : \Nsym_\QQ \to \Nsym_\QQ \otimes \Nsym_\QQ$
of the Hopf algebra $\Nsym_\QQ$ and for the continuous $\kk$-algebra
homomorphism $\Nsym_\QQ\left[\left[t\right]\right] \to
\left(\Nsym_\QQ \otimes \Nsym_\QQ\right) \left[\left[t\right]\right]$ it
induces.}):
\begin{align*}
\Delta \xi(t) &= \Delta \left( \log \sum_{n \geq 0} H_n t^n\right)
= \log \sum_{n \geq 0} \Delta(H_n) t^n
= \log \sum_{n \geq 0} \left( \sum_{i+j=n} H_i \otimes H_j \right) t^n\\
&= \log \left( \left(\sum_{i \geq 0} H_i t^i\right)  \otimes
                \left(\sum_{j \geq 0} H_j t^j\right) \right)
= \log \left( \left(\sum_{i \geq 0} H_i t^i \otimes \one  \right)
                \left( \one \otimes \sum_{j \geq 0} H_j t^j\right)  \right)\\
&\overset{\eqref{log-of-product}}{=}\log \widetilde{H}(t) \otimes \one + \one \otimes \log \widetilde{H}(t)
= \xi(t) \otimes \one + \one \otimes \xi(t).
\end{align*}
Comparing coefficients in this equality yields
$\Delta(\xi_n) = \xi_n \otimes \one + \one \otimes \xi_n$.
Thus $S(\xi_n) = -\xi_n$, by Proposition~\ref{antipode-on-primitives}.
This allows one to determine $S(H_n)$ and $S(H_\alpha)$, after
one first inverts the relation \eqref{relation-definining-pi's}
to get that $\widetilde{H}(t) = \exp \xi(t)$, and hence
\begin{align*}
S(\widetilde{H}(t)) &= S(\exp \xi(t))
=\exp S(\xi(t))
=\exp\left( -\xi(t) \right)
\overset{\eqref{exp-of-sum}}{=}(\exp \xi(t)) ^{-1} \\
&=
\widetilde{H}(t)^{-1}=\left( 1 + H_1 t + H_2 t^2 + \cdots \right)^{-1}.
\end{align*}
Upon expanding the right side, and comparing coefficients of $t^n$, this gives
\[
S(H_n)= \sum_{\beta \in \Comp_n} (-1)^{\ell(\beta)} H_\beta
\]
and hence
\[
S(H_\alpha) = S(H_{\alpha_\ell}) \cdots S(H_{\alpha_2}) S(H_{\alpha_1})
=\sum\limits_{\substack{\gamma: \\ \gamma \text{ refines }\rev(\alpha)}}
(-1)^{\ell(\gamma)} H_\gamma
=\sum\limits_{\substack{\gamma: \\ \rev(\gamma) \text{ refines }\alpha}}
(-1)^{\ell(\gamma)} H_\gamma
\]
(because if $\mu$ and $\nu$ are two compositions, then $\mu$ refines $\nu$
if and only if $\rev(\mu)$ refines $\rev(\nu)$).
As $S_{\Nsym}, S_{\Qsym}$ are adjoint, and $\{H_\alpha\},
\{M_\alpha\}$ are dual bases, this is equivalent to
saying that
\[
S(M_\alpha)= (-1)^{\ell(\alpha)} \sum\limits_{\substack{\gamma:\\ \rev(\alpha) \text{ refines } \gamma}} M_{\gamma}
\qquad \text{for all $\alpha \in \Comp$.}
\]
But this is precisely the claim of
Theorem~\ref{Qsym-antipode-on-monomials}.
Thus,
Theorem~\ref{Qsym-antipode-on-monomials} is proven once again.

Let us say a bit more about the elements $\xi_n$ defined in
\eqref{relation-definining-pi's} above. The elements $n\xi_n$
are noncommutative
analogues of the power sum symmetric functions $p_n$ (and, indeed, are
lifts of the latter to $\Nsym$, as Exercise~\ref{exe.xi-lifts-p} below
shows). They are called the \dfn{noncommutative power sums of the second
kind} in \cite{GelfandKrobLascouxLeclercRetakhThibon}\footnote{See
Exercise~\ref{exe.Psi} for the ones of the first kind.},
and their products form a basis of $\Nsym$. They are furthermore
useful in studying the so-called \emph{Eulerian idempotent} of a
cocommutative Hopf algebra, as shown in Exercise~\ref{exe.eulerian-idp}
below.
\end{remark}

% [DG][v43] Renamed $H(t)$ as $\widetilde{H}(t)$ to avoid conflict with
% the $H(t)$ over the symmetric functions. Also expanded the footnote.

\begin{exercise}
\label{exe.xi-lifts-p}
Assume that $\QQ$ is a subring of $\kk$.
Define a sequence
of elements $\xi_1,\xi_2,\xi_3,\ldots$ of $\Nsym = \Nsym_\kk$ by
\eqref{relation-definining-pi's}.
\begin{itemize}
\item[(a)] For every $n \geq 1$, show that $\xi_n$ is a primitive
homogeneous element of $\Nsym$ of degree $n$.
\item[(b)] For every $n \geq 1$, show that $\pi\left(n\xi_n\right)$
is the $n$-th power sum symmetric function $p_n \in \Lambda$.
\item[(c)] For every $n \geq 1$, show that
\begin{equation}
\label{eq.exe.xi-lifts-p.c.1}
\xi_n = \sum_{\alpha \in \Comp_n}
 \left(-1\right)^{\ell\left(\alpha\right)-1}
 \frac{1}{\ell\left(\alpha\right)} H_\alpha .
\end{equation}
\item[(d)] For every
composition $\alpha$, define an element $\xi_\alpha$ of $\Nsym$
by $\xi_\alpha = \xi_{\alpha_1} \xi_{\alpha_2} \cdots
\xi_{\alpha_\ell}$, where $\alpha$ is written in the form
$\alpha = \left(\alpha_1, \alpha_2, \ldots, \alpha_\ell\right)$
with $\ell = \ell\left(\alpha\right)$.
Show that
\begin{equation}
\label{eq.exe.xi-lifts-p.d.1}
H_n = \sum_{\alpha \in \Comp_n} \frac{1}{\ell\left(\alpha\right)!} \xi_\alpha
\end{equation}
for every $n \in \NN$.

Use this to prove that
$\left(\xi_\alpha\right)_{\alpha \in \Comp_n}$ is a $\kk$-basis
of $\Nsym_n$ for every $n \in \NN$.
\end{itemize}
\end{exercise}

% [DG][v19] Added parts (c) and (d) to the above exercise.

\begin{exercise}
\label{exe.eulerian-idp}
Assume that $\QQ$ is a subring of $\kk$. 
Let $A$ be a cocommutative connected graded $\kk$-bialgebra.
Let $A = \bigoplus_{n\geq 0} A_n$ be the decomposition of $A$
into homogeneous components.
If $f$ is any $\kk$-linear map $A \to A$ annihilating $A_0$, then $f$
is locally $\star$-nilpotent\footnote{See the proof of
Proposition~\ref{Takeuchi-prop} for what this means.}, and so the
sum $\log^\star \left(f + u\epsilon\right)
:= \sum_{n \geq 1} \left(-1\right)^{n-1} \frac{1}{n} f^{\star n}$
is a well-defined endomorphism of $A$\ \ \ \ \footnote{%
This definition of $\log^{\star}\left(  f+u\epsilon\right)  $ is actually a
particular case of Definition~\ref{def.exp-log.log}. This can be seen as follows:

We have $f\left(  A_{0}\right)  =0$. Thus,
Proposition~\ref{prop.convolution-series}(h) (applied to $C=A$) yields
$f\in\mathfrak{n}\left(  A,A\right)  $ (where $\mathfrak{n}\left(  A,A\right)
$ is defined as in Section~\ref{subsect.leray}), so that $\left(
f+u\epsilon\right)  -u\epsilon=f\in\mathfrak{n}\left(  A,A\right)  $.
Therefore, Definition~\ref{def.exp-log.log} defines a map $\log^{\star}\left(
f+u\epsilon\right)  \in\mathfrak{n}\left(  A,A\right)  $. This map is
identical to the map $\log^{\star}\left(  f+u\epsilon\right)  :=
\sum_{n\geq1}\left(  -1\right)  ^{n-1}\frac{1}{n}f^{\star n}$ we have just
defined, because Proposition~\ref{prop.exp-log}(f) (applied to $C=A$) shows
that the map $\log^{\star}\left(  f+u\epsilon\right)  $ defined using
Definition~\ref{def.exp-log.log} satisfies%
\[
\log^{\star}\left(  f+u\epsilon\right)  =\sum_{n\geq1}\dfrac{\left(
-1\right)  ^{n-1}}{n}f^{\star n}=\sum_{n\geq1}\left(  -1\right)  ^{n-1}%
\frac{1}{n}f^{\star n}.
\]
}. Let $\mathfrak{e}$ denote the
endomorphism $\log^\star \left(\id_A\right)$ of $A$ (obtained by setting
$f = \id_A - u\epsilon : A \to A$). Show that $\mathfrak{e}$
is a projection from $A$ to the $\kk$-submodule
$\Liep$ of all primitive elements of $A$ (and thus, in particular, is
idempotent).

\textbf{Hint:}
For every $n \geq 0$, let $\pi_n : A \to A$ be the projection
onto the $n$-th homogeneous component $A_n$.
Since $\Nsym$ is the free $\kk$-algebra with generators
$H_1, H_2, H_3, \ldots$, we can define a $\kk$-algebra homomorphism
$\mathfrak{W} : \Nsym \to \left(\End A, \star\right)$ by sending
$H_n$ to $\pi_n$. Show that:
\begin{itemize}
\item[(a)] The map $\mathfrak{e} : A \to A$ is graded.
For every $n \geq 0$, we will denote the map
$\pi_n \circ \mathfrak{e} = \mathfrak{e} \circ \pi_n : A \to A$ by
$\mathfrak{e}_n$.
\item[(b)] We have $\mathfrak{W}\left(\xi_n\right) = \mathfrak{e}_n$
for all $n \geq 1$, where $\xi_n$ is defined as in
Exercise~\ref{exe.xi-lifts-p}.
\item[(c)] If $w$ is an element of $\Nsym$, and if we write
$\Delta\left(w\right) = \sum_{(w)} w_1 \otimes w_2$ using the Sweedler
notation, then $\Delta \circ \left(\mathfrak{W}\left(w\right)\right)
= \left(\sum_{(w)} \mathfrak{W}\left(w_1\right) \otimes
\mathfrak{W}\left(w_2\right) \right) \circ \Delta$.
\item[(d)] We have $\mathfrak{e}_n\left(A\right) \subset \Liep$ for every
$n \geq 0$.
\item[(e)] We have $\mathfrak{e}\left(A\right) \subset \Liep$.
\item[(f)] The map $\mathfrak{e}$ fixes any element of $\Liep$.
\end{itemize}
\end{exercise}

\begin{remark}
The endomorphism $\mathfrak{e}$ of Exercise~\ref{exe.eulerian-idp}
is known as the \dfn{Eulerian idempotent} of $A$, and can be contrasted
with the Dynkin idempotent of Remark~\ref{rmk.dynkin}. It has been
studied in \cite{Patras}, \cite{PatrasReutenauer-hli},
\cite{Burgunder} and \cite{DuchampMinhTolluChienNghia}, and relates to
the Hochschild cohomology of commutative algebras \cite[\S 4.5.2]{Loday}.
\end{remark}

% [DG][v14] Added the two exercises above and this remark. I also
% renamed the $\pi_n$ as $\xi_n$ to avoid confusion with the $\pi$
% that projects $\Nsym$ on $\Lambda$ and with the projections $\pi_m$
% in the Aguiar-Bergeron-Sottile universal property (which also is
% related to this all...).
% 
% Exercise \ref{exe.eulerian-idp} originates in my feeling that the
% argument proving the primitivity of the $\xi_n$ in Remark
% \ref{second-antipode-proof} looked similar to the proof of the
% properties of the Eulerian idempotent in
% \cite{DuchampMinhTolluChienNghia}. I am likely not the first to
% have this feeling (can it be that one of the proofs inspired the
% other?), but it is not clear who should be cited for it.
% 
% The map $\mathfrak{W}$ constructed in the hint to the exercise
% becomes more important if one starts studying the internal product
% on $\Nsym$; it transports knowledge both ways (from $\Nsym$ to
% the endomorphism algebra of a cocommutative connected graded
% bialgebra and back).
% 
% Also, I should really stop using fraktur...

\begin{exercise}
\label{exe.eulerian-idp.2}
Assume that $\QQ$ is a subring of $\kk$. Let $A$, $A_n$ and
$\mathfrak{e}$ be as in Exercise~\ref{exe.eulerian-idp}.

\begin{itemize}
\item[(a)] Show that
$\mathfrak{e}^{\star n} \circ \mathfrak{e}^{\star m}
= n! \delta_{n,m} \mathfrak{e}^{\star n}$
for all $n \in \NN$ and $m \in \NN$.

\item[(b)] Show that
$\mathfrak{e}^{\star n} \circ \id_A^{\star m}
= \id_{A}^{\star m} \circ \mathfrak{e}^{\star n}
= m^n \mathfrak{e}^{\star n}$
for all $n \in \NN$ and $m \in \NN$.
\end{itemize}
\end{exercise}

% [DG][v28] Added exercise above.

We next explore the basis for $\Nsym$ dual to the $\{L_\alpha\}$ in $\Qsym$.
\begin{definition}
Define the
\emph{noncommutative ribbon functions}\index{noncommutative ribbon function}\index{$R_\alpha$}
$\{R_\alpha\}_{\alpha \in \Comp}$ to be the
$\kk$-basis of $\Nsym$ dual to
the fundamental basis $\{L_\alpha\}_{\alpha \in \Comp}$
of $\Qsym$,
so that
\[
(R_\alpha,L_\beta)=\delta_{\alpha,\beta}
\qquad \text{ for all $\alpha, \beta \in \Comp$.}
\]
\end{definition}

\begin{theorem} \phantomsection
\label{Nsym-structure-on-ribbons-theorem}
\begin{itemize}
\item[(a)] One has that
\begin{align}
H_\alpha&=\sum_{\beta \text{ coarsens }\alpha} R_\beta ;
\label{H-as-sum-of-ribbons}\\
R_\alpha&=\sum_{\beta \text{ coarsens }\alpha} (-1)^{\ell(\beta)-\ell(\alpha)} H_\beta .
\label{ribbon-as-sum-of-Hs}
\end{align}
\item[(b)] The surjection $\Nsym \overset{\pi}{\longrightarrow} \Lambda$
sends $R_\alpha \longmapsto s_{\Rib\left(\alpha\right)}$,
the skew Schur function associated
to the ribbon $\Rib\left(\alpha\right)$.
\item[(c)] Furthermore,
\begin{align}
R_\alpha R_\beta &= R_{\alpha \cdot \beta} + R_{\alpha \odot \beta} \qquad \text{ if }\alpha\text{ and }\beta\text{ are nonempty} ;
\label{Nsym-product-in-ribbons}\\
S(R_\alpha) &= (-1)^{|\alpha|} R_{\omega(\alpha)} .
\label{Nsym-antipode-in-ribbons}
\end{align}
Finally, $R_\varnothing$ is the multiplicative identity of $\Nsym$.
\end{itemize}
\end{theorem}

\begin{proof}
(a) For \eqref{H-as-sum-of-ribbons}, note that
\[
H_\alpha
= \sum_{\beta} (H_\alpha,L_\beta) R_\beta
= \sum_{\beta}
 \left( H_\alpha, \sum\limits_{\substack{\gamma:\\\gamma \text{ refines }\beta}} M_\gamma \right)
    R_\beta
=  \sum\limits_{\substack{\beta:\\ \beta \text{ coarsens }\alpha}}     R_\beta.
\]
The equality \eqref{ribbon-as-sum-of-Hs} follows from
\eqref{H-as-sum-of-ribbons}
by Lemma~\ref{lem.moebius-bool.comps}(a).

(b) Write $\alpha$ as $\left(\alpha_1,\ldots,\alpha_\ell\right)$.
To show that $\pi(R_\alpha)=s_{\Rib\left(\alpha\right)}$,
we instead examine $\pi(H_\alpha)$:
\[
\pi(H_\alpha)
=\pi(H_{\alpha_1} \cdots H_{\alpha_\ell})
=h_{\alpha_1} \cdots h_{\alpha_\ell}
= s_{(\alpha_1)} \cdots s_{(\alpha_\ell)}
= s_{(\alpha_1) \oplus \cdots \oplus (\alpha_\ell)}
\]
where $(\alpha_1) \oplus \cdots \oplus (\alpha_\ell)$ is some
skew shape which is a horizontal strip having rows of lengths $\alpha_1,\ldots,\alpha_\ell$
from bottom to top.  We claim
\[
s_{(\alpha_1) \oplus \cdots \oplus (\alpha_\ell)}
=\sum\limits_{\substack{\beta:\\\beta\text{ coarsens }\alpha}}
s_{\Rib\left(\beta\right)} ,
\]
because column-strict tableaux $T$ of shape
$(\alpha_1) \oplus \cdots \oplus (\alpha_\ell)$ biject to column-strict
tableaux $T'$ of some ribbon $\Rib\left(\beta\right)$
with $\beta$ coarsening $\alpha$, as follows:
Let $a_i,b_i$ denote the leftmost, rightmost entries of the $i$-th row from the
bottom in $T$, of length $\alpha_i$, and
\begin{enumerate}
\item[$\bullet$] if $b_i \leq a_{i+1}$,
merge parts $\alpha_i,\alpha_{i+1}$ in $\beta$,
and concatenate the rows of length $\alpha_i,\alpha_{i+1}$ in $T'$, or
\item[$\bullet$] if $b_i > a_{i+1}$,
do not merge parts $\alpha_i,\alpha_{i+1}$ in $\beta$,
and let these two rows overlap in one column in $T'$.
\end{enumerate}
E.g., if $\alpha=(3,3,2,3,2)$, then
\begin{align*}
\text{the tableau } &
T =
\begin{matrix}
 & & & & & & & & & & &3&4 \\
 & & & & & & & &4&4&5& &\\
 & & & & & &4&4& & & & &\\
 & & &2&2&3& & & & & & &\\
1&1&3& & & & & & & & & &
\end{matrix}
\text{ of shape }
(\alpha_1) \oplus \cdots \oplus (\alpha_\ell) \\
\\
\text{maps to the tableau } &
T' =
\begin{matrix}
 & & & & & & & & &3&4 \\
 & &2&2&3&4&4&4&4&5& &\\
1&1&3& & & & & & & & & &
\end{matrix}
\text{ of shape }
\Rib\left(\beta\right)
\text{ for } \beta=(3,8,2) .
\end{align*}
% [DG][v80] Reorganized the above example.
The reverse bijection breaks the rows of $T'$ into the rows of $T$ of lengths dictated by
the parts of $\alpha$.  Having shown
$\pi(H_\alpha) =\sum_{\beta:\\\beta\text{ coarsens }\alpha} s_{\Rib\left(\beta\right)}$,
we can now apply Lemma~\ref{lem.moebius-bool.comps}(a) to obtain
\[
s_{\Rib\left(\alpha\right)}
= \sum_{\beta:\\\beta\text{ coarsens }\alpha} \left(-1\right)^{\ell\left(\alpha\right)-\ell\left(\beta\right)} \pi\left(H_\beta\right)
= \pi\left(R_\alpha\right)
\qquad \left(\text{by \eqref{ribbon-as-sum-of-Hs}}\right) ;
\]
thus, $\pi(R_\alpha) = s_{\Rib\left(\alpha\right)}$ is proven.

(c) Finally, \eqref{Nsym-product-in-ribbons} and \eqref{Nsym-antipode-in-ribbons}
follow from \eqref{Qsym-coproduct-on-fundamentals} and
\eqref{Qsym-antipode-on-fundamentals} by duality.
\end{proof}

\begin{remark}
\label{antipode-comparison-remark}
Since the maps
\[
\xymatrix{
\Nsym \arsurj[dr]_\pi & & \Qsym \\
      &\Lambda \arinjrev[ur]_i&
}
\]
are Hopf morphisms, they must respect the antipodes $S_{\Lambda}, S_{\Qsym}, S_{\Nsym}$,
but it is interesting to compare them explicitly using the fundamental basis for
$\Qsym$ and the ribbon basis for $\Nsym$.

On one hand \eqref{Qsym-antipode-on-fundamentals}
shows that $S_{\Qsym}(L_\alpha)= (-1)^{|\alpha|} L_{\omega(\alpha)}$
extends the map $S_{\Lambda}$ since $L_{(1^n)} = e_n$ and $L_{(n)}=h_n$, as observed in
Example~\ref{the-two-extreme-fundamentals-example}, and $\omega((n))=(1^n)$.

On the other hand, \eqref{Nsym-antipode-in-ribbons} shows that
$S_{\Nsym}(R_\alpha)= (-1)^{|\alpha|} R_{\omega(\alpha)}$  lifts the map $S_{\Lambda}$ to $S_{\Nsym}$:
Theorem~\ref{Nsym-structure-on-ribbons-theorem}(b) showed that
$R_\alpha$ lifts the skew Schur function $s_{\Rib\left(\alpha\right)}$,
while \eqref{skew-Schur-antipodes} asserted that
$S(s_{\lambda/\mu}) = (-1)^{|\lambda/\mu|} s_{\lambda^t/\mu^t}$,
and a ribbon $\Rib\left(\alpha\right)=\lambda/\mu$ has
$\Rib\left(\omega(\alpha)\right)=\lambda^t/\mu^t$.
\end{remark}

\begin{exercise} \phantomsection
\label{exe.Psi}
\begin{itemize}
\item[(a)] Show that any integers $n$ and $i$ with $0\leq
i<n$ satisfy
\[
R_{\left(  1^{i},n-i\right)  }=\sum_{j=0}^{i}\left(  -1\right)  ^{i-j}
R_{\left(  1^{j}\right)  }H_{n-j}.
\]
(Here, as usual, $1^i$ stands for the number $1$ repeated $i$ times.)

\item[(b)] Show that any integers $n$ and $i$ with $0\leq i<n$ satisfy
\[
\left(  -1\right)  ^{i}R_{\left(  1^{i},n-i\right)  }=\sum_{j=0}^{i}S\left(
H_{j}\right)  H_{n-j}.
\]

\item[(c)] For every positive integer $n$, define an element $\Psi_{n}$ of
$\Nsym$ by
\[
\Psi_{n}=\sum_{i=0}^{n-1}\left(  -1\right)  ^{i}R_{\left(  1^{i},n-i\right)
}.
\]
Show that $\Psi_{n}=\left(  S\star E\right)  \left(  H_{n}\right)  $, where
the map $E : \Nsym \to \Nsym$ is defined as in
Exercise~\ref{exe.dynkin} (for $A=\Nsym$).
Conclude that $\Psi_{n}$ is primitive.

\item[(d)] Prove that
\[
\sum_{k=0}^{n-1}H_{k}\Psi_{n-k} = nH_{n}
\]
for every $n\in\NN$.

\item[(e)] Define two power series $\psi\left(  t\right)  $ and
$\widetilde{H}\left(  t\right)  $ in $\Nsym\left[  \left[
t\right]  \right]  $ by
\begin{align*}
\psi\left(  t\right)   &  =\sum_{n\geq 1}\Psi_{n}t^{n-1};\\
\widetilde{H}\left(  t\right)   &  =\sum_{n\geq0}H_{n}t^{n}.
\end{align*}
Show that\footnote{The derivative $\frac{d}{dt}Q\left(t\right)$ of
a power series $Q\left(t\right) \in R\left[\left[t\right]\right]$
over a noncommutative ring $R$ is defined just as in the case of
$R$ commutative: by setting
$\frac{d}{dt}Q\left(t\right) = \sum_{i \geq 1} iq_i t^{i-1}$,
where $Q\left(t\right)$ is written in the form
$Q\left(t\right) = \sum_{i \geq 0} q_i t^i$.}
$\dfrac{d}{dt}\widetilde{H}\left(  t\right)  =\widetilde{H}\left(
t\right)  \cdot\psi\left(  t\right)  $.

(The functions $\Psi_n$ are called \dfn{noncommutative power sums of the
first kind}; they are studied in
\cite{GelfandKrobLascouxLeclercRetakhThibon}. The power sums of the second
kind are the $n\xi_n$ in Remark~\ref{second-antipode-proof}.)

\item[(f)] Show that $\pi\left(  \Psi_n \right)  $ equals the power sum
symmetric function $p_n$ for every positive integer $n$.

\item[(g)] Show that every positive integer $n$ satisfies
\[
p_n = \sum\limits_{i=0}^{n-1} \left(-1\right)^i s_{\left(n-i, 1^i\right)}
  \qquad \text{ in }\Lambda .
\]

\item[(h)] For every nonempty
composition $\alpha$, define a positive integer
$\operatorname{lp}\left(\alpha\right)$ by
$\operatorname{lp}\left(\alpha\right) = \alpha_\ell$,
where $\alpha$ is written in the form
$\alpha = \left(\alpha_1, \alpha_2, \ldots, \alpha_\ell\right)$
with $\ell = \ell\left(\alpha\right)$.
(Thus, $\operatorname{lp}\left(\alpha\right)$ is the last
part of $\alpha$.) Show that every positive integer $n$
satisfies
\begin{equation}
\label{eq.exe.Psi.h.0}
\Psi_n = \sum_{\alpha \in \Comp_n} \left(-1\right)^{\ell\left(\alpha\right)-1} \operatorname{lp}\left(\alpha\right) H_\alpha .
\end{equation}

\item[(i)] Assume that $\QQ$ is a subring of $\kk$. For every
composition $\alpha$, define an element $\Psi_\alpha$ of $\Nsym$
by $\Psi_\alpha = \Psi_{\alpha_1} \Psi_{\alpha_2} \cdots
\Psi_{\alpha_\ell}$, where $\alpha$ is written in the form
$\alpha = \left(\alpha_1, \alpha_2, \ldots, \alpha_\ell\right)$
with $\ell = \ell\left(\alpha\right)$.
For every composition
$\alpha$, define $\pi_u\left(\alpha\right)$ to be the positive
integer
$\alpha_1 \left(\alpha_1 + \alpha_2\right) \cdots \left(\alpha_1 + \alpha_2 + \cdots + \alpha_\ell\right)$,
where $\alpha$ is written in the form
$\alpha = \left(\alpha_1, \alpha_2, \ldots, \alpha_\ell\right)$
with $\ell = \ell\left(\alpha\right)$.
Show that
\begin{equation}
\label{eq.exe.Psi.h.1}
H_n = \sum_{\alpha \in \Comp_n} \frac{1}{\pi_u\left(\alpha\right)} \Psi_\alpha 
\end{equation}
for every $n \in \NN$.

Use this to prove that
$\left(\Psi_\alpha\right)_{\alpha \in \Comp_n}$ is a $\kk$-basis
of $\Nsym_n$ for every $n \in \NN$.

\item[(j)] Assume that $\QQ$ is a subring of $\kk$.
Let $V$ be the
free $\kk$-module with basis
$\left(\mathfrak{b}_n\right)_{n \in \left\{1, 2, 3, \ldots\right\}}$.
Define a $\kk$-module homomorphism
$f : V \to \Nsym$ by requiring that
$f\left(\mathfrak{b}_n\right) = \Psi_n$ for every
$n \in \left\{1, 2, 3, \ldots\right\}$.
Let $F$ be the $\kk$-algebra homomorphism
$T\left(V\right) \to \Nsym$ induced by this $f$ (using the
universal property of the tensor algebra $T\left(V\right)$).
Show that $F$ is a Hopf algebra isomorphism (where the
Hopf algebra structure on $T\left(V\right)$ is as in
Example~\ref{exa.tensor-alg.antipode}).

\item[(k)] Assume that $\QQ$ is a subring of $\kk$. Let $V$ be as
in Exercise~\ref{exe.Psi}(j). Show that $\Qsym$ is isomorphic to
the shuffle algebra $\operatorname{Sh}\left(V\right)$ (defined as
in Proposition~\ref{prop.shuffle-alg}) as Hopf algebras.

\item[(l)] Solve parts (a) and (b) of Exercise~\ref{exe.hookschur}
again using the ribbon basis functions $R_\alpha$.
\end{itemize}
\end{exercise}

% [DG][v7] Added exercise above.

% [DG][v14] Added part (g), which is an application of $\Nsym$
% to $\Lambda$. Of course, it is not very exciting if one knows
% Murnaghan-Nakayama...

% [DG][v19] Added parts (h) and (i). Copypasted parts of the solution
% from the previous exercise because I was tired.

% [DG][v30] Added part (j).

% [DG][v31] Added parts (j) and (k), moving the old (j) to place (l).

One might wonder whether the Frobenius endomorphisms of $\Lambda$ (defined in
Exercise~\ref{exe.witt.ghost-app}) and the Verschiebung endomorphisms of
$\Lambda$ (defined in Exercise~\ref{exe.Lambda.verschiebung}) generalize to
analogous operators on either $\Qsym$ or $\Nsym$. The next two
exercises (whose claims mostly come from \cite[\S 13]{Hazewinkel2}) answer
this question: The Frobenius endomorphisms extend to $\Qsym$, and
the Verschiebung ones lift to $\Nsym$.

\begin{exercise}
\label{exe.QSym.frobenius}
For every $n \in \left\{ 1, 2, 3, \ldots \right\}$,
define a map $\mathbf{F}_n : \Qsym \to \Qsym$ by setting
\[
\mathbf{F}_n \left( a \right)
= a \left( x_1^n, x_2^n, x_3^n, \ldots \right)
\qquad\qquad\text{for every }a \in \Qsym .
\]
(So what $\mathbf{F}_n$ does to a quasi-symmetric function is replacing all
variables $x_1, x_2, x_3, \ldots$ by their $n$-th powers.)

\begin{enumerate}
\item[(a)] Show that $\mathbf{F}_n : \Qsym \to \Qsym$ is a $\kk$-algebra
homomorphism for every $n \in \left\{ 1, 2, 3, \ldots \right\}$.

\item[(b)] Show that $\mathbf{F}_n \circ \mathbf{F}_m = \mathbf{F}_{nm}$
for any two positive integers $n$ and $m$.

\item[(c)] Show that $\mathbf{F}_1 = \id$.

\item[(d)] Prove that
$\mathbf{F}_n
\left( M_{\left( \beta_1, \beta_2, \ldots, \beta_s \right)} \right)
= M_{\left( n\beta_1, n\beta_2, \ldots, n\beta_s \right)}$
for every $n \in \left\{ 1, 2, 3, \ldots \right\}$
and $\left( \beta_1, \beta_2, \ldots, \beta_s \right) \in \Comp$.

\item[(e)] Prove that $\mathbf{F}_n : \Qsym \to \Qsym$ is a Hopf
algebra homomorphism for every
$n \in \left\{ 1, 2, 3, \ldots \right\}$.

\item[(f)] Consider the maps $\mathbf{f}_n : \Lambda \to \Lambda$
defined in Exercise~\ref{exe.witt.ghost-app}. Show that
$\mathbf{F}_n \mid_{\Lambda} = \mathbf{f}_n$ for every
$n \in \left\{ 1, 2, 3, \ldots \right\}$.

\item[(g)] Assume that $\kk = \ZZ$. Prove that $\mathbf{f}_p
\left( a \right) \equiv a^p \operatorname{mod} p\Qsym$ for
every $a \in \Qsym$ and every prime number $p$.

\item[(h)] Give a new solution to
Exercise~\ref{exe.witt.ghost-app}(d).
\end{enumerate}
\end{exercise}

\begin{exercise}
\label{exe.NSym.verschiebung}
For every $n \in \left\{ 1, 2, 3, \ldots \right\}$,
define a $\kk$-algebra homomorphism
$\mathbf{V}_n : \Nsym \rightarrow \Nsym$ by
\[
\mathbf{V}_n \left( H_m \right)
= \begin{cases}
H_{m/n},& \text{if } n\mid m;\\
0,& \text{if } n\nmid m
\end{cases}
\qquad\qquad\text{for every positive integer }m
\]
\footnote{This is well-defined, since $\Nsym$ is (isomorphic
to) the free associative algebra with generators $H_1, H_2, H_3, \ldots$
(according to \eqref{eq.thm.NSym.free}).}.

\begin{enumerate}
\item[(a)] Show that any positive integers $n$ and $m$ satisfy
\[
\mathbf{V}_n \left( \Psi_m \right)
= \begin{cases}
n\Psi_{m/n},& \text{if } n\mid m;\\
0,& \text{if } n\nmid m
\end{cases} \quad ,
\]
where the elements $\Psi_m$ and $\Psi_{m/n}$ of $\Nsym$ are
as defined in Exercise~\ref{exe.Psi}(c).

\item[(b)] Show that if $\QQ$ is a subring of $\kk$, then any
positive integers $n$ and $m$ satisfy
\[
\mathbf{V}_n \left( \xi_m \right)
= \begin{cases}
\xi_{m/n},& \text{if } n\mid m;\\
0,& \text{if } n\nmid m
\end{cases} \quad ,
\]
where the elements $\xi_m$ and $\xi_{m/n}$ of $\Nsym$ are as
defined in Exercise~\ref{exe.xi-lifts-p}.

\item[(c)] Prove that
$\mathbf{V}_n \circ \mathbf{V}_m = \mathbf{V}_{nm}$ for
any two positive integers $n$ and $m$.

\item[(d)] Prove that $\mathbf{V}_1 = \id$.

\item[(e)] Prove that $\mathbf{V}_n : \Nsym \rightarrow
 \Nsym$ is a Hopf algebra homomorphism for every
$n \in \left\{ 1, 2, 3, \ldots \right\}$.
\end{enumerate}

Now, consider also the maps $\mathbf{F}_n : \Qsym
\rightarrow \Qsym$ defined in Exercise~\ref{exe.witt.ghost-app}.
Fix a positive integer $n$.

\begin{enumerate}
\item[(f)] Prove that the maps $\mathbf{F}_n : \Qsym
\rightarrow \Qsym$ and $\mathbf{V}_n : \Nsym
\rightarrow \Nsym$ are adjoint with respect to the dual pairing
$\Nsym \otimes \Qsym
\overset{\left(  \cdot , \cdot\right)  }{\longrightarrow} \kk$.

\item[(g)] Consider the maps
$\mathbf{v}_n : \Lambda \rightarrow \Lambda$
defined in Exercise~\ref{exe.Lambda.verschiebung}. Show that
the surjection $\pi : \Nsym \rightarrow \Lambda$ satisfies
$\mathbf{v}_n \circ \pi = \pi \circ \mathbf{V}_n$ for every
$n \in \left\{ 1, 2, 3, \ldots \right\}$.

\item[(h)] Give a new solution to
Exercise~\ref{exe.Lambda.verschiebung}(f).
\end{enumerate}
\end{exercise}

% [DG][v43] Added the above two exercises.

\newpage

%%%%%%%%%%%%%%%%%%%%%%%%%%%%%%%%%%%%%%%%%%%%%%%%%%%%%%%%%%%%%%%%%%%%%%%%%%%%
\section{Polynomial generators for $\Qsym$ and Lyndon words}
\label{sect.QSym.lyndon}
%%%%%%%%%%%%%%%%%%%%%%%%%%%%%%%%%%%%%%%%%%%%%%%%%%%%%%%%%%%%%%%%%%%%%%%%%%%%

In this chapter, we shall construct an algebraically independent
generating set for $\Qsym$ as a $\kk$-algebra, thus showing that
$\Qsym$ is a polynomial ring over $\kk$. This has been
done by Malvenuto \cite[Cor. 4.19]{Malvenuto} when $\kk$
is a field of characteristic $0$, and by Hazewinkel
\cite{Hazewinkel} in the general case. We will begin by
introducing the notion of \emph{Lyndon words}
(Section~\ref{subsect.lyndon}), on which both of these
constructions rely; we will then
(Section~\ref{subsect.lyndon.shuffles}) elucidate the connection
of Lyndon words with shuffles, and afterwards
(Section~\ref{subsect.shuffle.radford}) apply it to prove
\emph{Radford's theorem} stating that the shuffle algebra of
a free $\kk$-module over a commutative $\QQ$-algebra is a
polynomial algebra (Theorem~\ref{thm.shuffle.radford}). The
shuffle algebra is not yet $\Qsym$, but Radford's theorem on the
shuffle algebra serves
as a natural stepping stone for the study of the more complicated
algebra $\Qsym$. We will prove -- in two ways -- that $\Qsym$
is a polynomial algebra when $\QQ$ is a subring of $\kk$ in
Section~\ref{subsect.shuffle.QSym.1}, and then we will finally
prove the general case in Section~\ref{subsect.shuffle.QSym.2}.
In Section~\ref{subsect.lyndon.gr}, we will explore a different
aspect of the combinatorics of words: the notion of necklaces
(which are in bijection with Lyndon words, as
Exercise~\ref{exe.words.necklaces} will show) and the
\emph{Gessel-Reutenauer bijection}, which help us define and
understand the \emph{Gessel-Reutenauer symmetric functions}.
This will rely on Section~\ref{subsect.lyndon}, but not on any
of the other sections of Chapter~\ref{sect.QSym.lyndon}.

Strictly speaking, this whole Chapter~\ref{sect.QSym.lyndon}
is a digression, as it involves almost no coalgebraic or
Hopf-algebraic structures, and its results will not be used in
further chapters (which means it can be skipped if so desired).
However, it sheds additional light on both quasisymmetric and
symmetric functions,
and serves as an excuse to study Lyndon words, which are a
combinatorial object of independent interest (and are involved
in the study of free algebras and Hopf algebras, apart from
$\Qsym$ -- see \cite{Radford-shuffle} and
\cite{Reutenauer}\footnote{They also are involved in indexing
basis elements of combinatorial Hopf algebras other than
$\Qsym$. See Bergeron/Zabrocki \cite{BergeronZabrocki}.}).

We will take a scenic route to the proof of Hazewinkel's
theorem. A reader only interested in the proof proper can
restrict themselves to reading only the following:

\begin{itemize}

\item from Section~\ref{subsect.lyndon}, everything up to
Corollary~\ref{cor.words.uv=vu.3}, then from
Definition~\ref{def.words.lyndon} up to 
Proposition~\ref{prop.words.lyndon.preorder}, then from
Definition~\ref{def.words.CFL} up to
Lemma~\ref{lem.words.CFLlemma}, and finally
Theorem~\ref{thm.words.lyndon.std}.
(Proposition~\ref{prop.words.lyndon.cutoff} and
Theorem~\ref{thm.words.lyndon.equiv} are also relevant
if one wants to use a different definition of Lyndon words,
as they prove the equivalence of most such definitions.)

\item from Section~\ref{subsect.lyndon.shuffles},
everything except for
Exercise~\ref{exe.shuffle.lyndon.cflcrit}.

\item from Section~\ref{subsect.shuffle.radford},
Definition~\ref{def.shuffle.polyalg},
Lemma~\ref{lem.lyndonbases}, and
Lemma~\ref{lem.shuffle.preserves-lex}.

\item from Section~\ref{subsect.shuffle.QSym.1},
Definition~\ref{def.lyndon-composition},
Theorem~\ref{thm.QSym.lyndon}, then from
Proposition~\ref{prop.QSym.lyndon.MaMb} up to
Definition~\ref{def.wll}, and
Lemma~\ref{lem.QSym.lyndon.triang.wll-lemma}.

\item all of Section~\ref{subsect.shuffle.QSym.2}.

\end{itemize}

\noindent Likewise, Section~\ref{subsect.lyndon.gr} can be
read immediately after Section~\ref{subsect.lyndon}.

\subsection{\label{subsect.lyndon}Lyndon words}

Lyndon words have been independently defined by Shirshov \cite{Shirshov},
Lyndon \cite{Lyndon}, Radford \cite[\S 2]{Radford-shuffle} and de
Bruijn/Klarner \cite{deBruijnKlarner} (though using different and sometimes
incompatible notations). They have since been surfacing in various places in
noncommutative algebra (particularly the study of free Lie algebras);
expositions of their theory can be found in \cite[\S 5]{Lothaire},
\cite[\S 5.1]{Reutenauer} and \cite[\S 1]{Laue} (in German).
We will follow our own approach to the properties
of Lyndon words that we need.

\begin{definition}
\label{def.words}We fix a totally ordered set $\mathfrak{A}$, which we call
the \dfn{alphabet}. Throughout Section~\ref{subsect.lyndon} and
Section~\ref{subsect.lyndon.shuffles}, we will understand
``word'' to mean a word over $\mathfrak{A}$.

We recall that a \dfn{word} is just a (finite) tuple of elements of
$\mathfrak{A}$. In other words, a word is an element of the set $\bigsqcup
_{n\geq0}\mathfrak{A}^{n}$. We denote this set by
\dfn{$\mathfrak{A}^{\ast}$}.

The \dfn{empty word} is the unique tuple with $0$ elements. It is denoted
by $\varnothing$. If $w\in\mathfrak{A}^{n}$ is a word and $i\in\left\{
1,2,\ldots,n\right\}  $, then the \dfn{$i$-th letter} of $w$ means the
$i$-th entry of the $n$-tuple $w$. This $i$-th letter will be denoted by
$w_{i}$.

The \emph{length}\index{length of a word}
$\ell\left(  w\right)  $ of a word $w\in\bigsqcup_{n\geq
0}\mathfrak{A}^{n}$ is defined to be the $n\in \NN$ satisfying
$w\in\mathfrak{A}^{n}$. Thus, $w=\left(  w_{1},w_{2},\ldots,w_{\ell\left(
w\right)  }\right)  $ for every word $w$.

Given two words $u$ and $v$, we say that $u$ is \dfn{longer} than $v$
(or, equivalently, $v$ is \dfn{shorter} than $u$) if and only if
$\ell\left(u\right) > \ell\left(v\right)$.

The \dfn{concatenation} of two words $u$ and $v$ is defined to be the word
$\left(  u_{1},u_{2},\ldots,u_{\ell\left(  u\right)  },v_{1},v_{2}%
,\ldots,v_{\ell\left(  v\right)  }\right)  $. This concatenation is denoted by
$uv$ or $u\cdot v$. The set $\mathfrak{A}^{\ast}$ of all words is a monoid
with respect to concatenation, with neutral element $\varnothing$. It is
precisely the free monoid on generators $\mathfrak{A}$. If $u$ is a word and
$i\in \NN$, we will understand $u^{i}$ to mean the $i$-th power of $u$
in this monoid (that is, the word $\underbrace{uu \cdots u}_{i\text{ times}}$).

The elements of $\mathfrak{A}$ are called \emph{letters}\index{letter},
and will be
identified with elements of $\mathfrak{A}^{1}\subset\bigsqcup_{n\geq
0}\mathfrak{A}^{n}=\mathfrak{A}^{\ast}$. This identification equates every
letter $u\in\mathfrak{A}$ with the one-letter word $\left(  u\right)
\in\mathfrak{A}^{1}$. Thus, every word $\left(  u_{1},u_{2},\ldots
,u_{n}\right)  \in\mathfrak{A}^{\ast}$ equals the concatenation $u_{1}u_{2}
\cdots u_{n}$ of letters, hence allowing us to use $u_{1}u_{2} \cdots u_{n}$
as a brief notation for the word $\left(  u_{1},u_{2},\ldots,u_{n}\right)  $.

If $w$ is a word, then:

\begin{itemize}
\item a \dfn{prefix} of $w$ means a word of the form $\left(  w_{1}%
,w_{2},\ldots,w_{i}\right)  $ for some $i\in\left\{  0,1,\ldots,\ell\left(
w\right)  \right\}  $;

\item a \dfn{suffix} of $w$ means a word of the form $\left(
w_{i+1},w_{i+2},\ldots,w_{\ell\left(  w\right)  }\right)  $ for some
$i\in\left\{  0,1,\ldots,\ell\left(  w\right)  \right\}  $;

\item a \dfn{proper suffix} of $w$ means a word of the form $\left(
w_{i+1},w_{i+2},\ldots,w_{\ell\left(  w\right)  }\right)  $ for some
$i\in\left\{  1,2,\ldots,\ell\left(  w\right)  \right\}  $.
\end{itemize}

\noindent In other words,

\begin{itemize}
\item a \dfn{prefix} of $w\in\mathfrak{A}^{\ast}$ is a word $u\in
\mathfrak{A}^{\ast}$ such that there exists a $v\in\mathfrak{A}^{\ast}$
satisfying $w=uv$;

\item a \dfn{suffix} of $w\in\mathfrak{A}^{\ast}$ is a word $v\in
\mathfrak{A}^{\ast}$ such that there exists a $u\in\mathfrak{A}^{\ast}$
satisfying $w=uv$;

\item a \dfn{proper suffix} of $w\in\mathfrak{A}^{\ast}$ is a word
$v\in\mathfrak{A}^{\ast}$ such that there exists a nonempty $u\in
\mathfrak{A}^{\ast}$ satisfying $w=uv$.
\end{itemize}

Clearly, any proper suffix of $w\in\mathfrak{A}^{\ast}$ is a suffix
of $w$. Moreover, if $w \in \mathfrak{A}^{\ast}$ is any word, then
a proper suffix of $w$ is the same thing as a suffix of $w$ distinct
from $w$.

We define a relation $\leq$ on the set $\mathfrak{A}^{\ast}$ as follows: For
two words $u\in\mathfrak{A}^{\ast}$ and $v\in\mathfrak{A}^{\ast}$, we set
$u\leq v$ to hold if and only if%
\begin{align*}
&  \text{\textbf{either }}\text{there exists an }i\in\left\{  1,2,\ldots
,\min\left\{  \ell\left(  u\right)  ,\ell\left(  v\right)  \right\}  \right\}
\\
&  \ \ \ \ \ \ \ \ \ \ \text{such that }\left(  u_{i}<v_{i}\text{, and every
}j\in\left\{  1,2,\ldots,i-1\right\}  \text{ satisfies }u_{j}=v_{j}\right)
,\\
&  \text{\textbf{or} the word }u\text{ is a prefix of }v.
\end{align*}
This order relation (taken as the smaller-or-equal relation) makes
$\mathfrak{A}^{\ast}$ into a poset (by Proposition \ref{prop.words.lex}(a)
below), and we will always be regarding $\mathfrak{A}^{\ast}$ as endowed with
this poset structure (thus, notations such as $<$, $\leq$, $>$ and $\geq$ will
be referring to this poset structure). This poset is actually totally ordered
(see Proposition \ref{prop.words.lex}(a)).
\end{definition}

Here are some examples of words compared by the relation $\leq$:%
\begin{align*}
113  &  \leq114,\ \ \ \ \ \ \ \ \ \ 113\leq132,\ \ \ \ \ \ \ \ \ \ 19\leq
195,\ \ \ \ \ \ \ \ \ \ 41\leq412,\\
41  &  \leq421,\ \ \ \ \ \ \ \ \ \ 539\leq54,\ \ \ \ \ \ \ \ \ \ \varnothing
\leq21,\ \ \ \ \ \ \ \ \ \ \varnothing\leq\varnothing
\end{align*}
(where $\mathfrak{A}$ is the alphabet $\left\{  1<2<3<\cdots\right\}  $).

Notice that if $u$ and $v$ are two words of the same length (i.e., we have
$u,v\in\mathfrak{A}^{n}$ for one and the same $n$), then $u\leq v$ holds if
and only if $u$ is lexicographically smaller-or-equal to $v$. In other words,
the relation $\leq$ is an extension of the lexicographic order on every
$\mathfrak{A}^{n}$ to $\mathfrak{A}^{\ast}$. This is the reason why this
relation $\leq$ is usually called the
\emph{lexicographic order}\index{lexicographic order on words}
on $\mathfrak{A}^{\ast}$.
In particular, we will be using this name.\footnote{The
relation $\leq$ is also known as the \dfn{dictionary order}, due to the
fact that it is the order in which words appear in a dictionary.} However,
unlike the lexicographic order on $\mathfrak{A}^{n}$, it does not always
respect concatenation from the right: It can happen that $u,v,w\in
\mathfrak{A}^{\ast}$ satisfy $u\leq v$ but not $uw\leq vw$. (For example,
$u=1$, $v=13$ and $w=4$, again with $\mathfrak{A}=\left\{  1<2<3<\cdots
\right\}  $.) We will see in Proposition \ref{prop.words.lex} that this is
rather an exception than the rule and the relation $\leq$ still behaves mostly
predictably with respect to concatenation.

Some basic properties of the order relation $\leq$ just defined are collected
in the following proposition:

\begin{proposition}
\phantomsection
\label{prop.words.lex}

\begin{itemize}
\item[(a)] The order relation $\leq$ is (the smaller-or-equal relation of) a
total order on the set $\mathfrak{A}^{\ast}$.

\item[(b)] If $a,c,d\in\mathfrak{A}^{\ast}$ satisfy $c\leq d$, then $ac\leq
ad$.

\item[(c)] If $a,c,d\in\mathfrak{A}^{\ast}$ satisfy $ac\leq ad$, then $c\leq
d$.

\item[(d)] If $a,b,c,d\in\mathfrak{A}^{\ast}$ satisfy $a\leq c$, then either
we have $ab\leq cd$ or the word $a$ is a prefix of $c$.

\item[(e)] If $a,b,c,d\in\mathfrak{A}^{\ast}$ satisfy $ab\leq cd$, then either
we have $a\leq c$ or the word $c$ is a prefix of $a$.

\item[(f)] If $a,b,c,d\in\mathfrak{A}^{\ast}$ satisfy $ab\leq cd$ and
$\ell\left(  a\right)  \leq\ell\left(  c\right)  $, then $a\leq c$.

\item[(g)] If $a,b,c\in\mathfrak{A}^{\ast}$ satisfy $a\leq b\leq ac$, then $a$
is a prefix of $b$.

\item[(h)] If $a\in\mathfrak{A}^{\ast}$ is a prefix of $b\in\mathfrak{A}%
^{\ast}$, then $a\leq b$.

\item[(i)] If $a$ and $b$ are two prefixes of $c\in\mathfrak{A}^{\ast}$, then
either $a$ is a prefix of $b$, or $b$ is a prefix of $a$.

\item[(j)] If $a,b,c\in\mathfrak{A}^{\ast}$ are such that $a\leq b$ and
$\ell\left(  a\right)  \geq\ell\left(  b\right)  $, then $ac\leq bc$.

\item[(k)] If $a\in\mathfrak{A}^{\ast}$ and $b\in\mathfrak{A}^{\ast}$ are such
that $b$ is nonempty, then $a<ab$.
\end{itemize}
\end{proposition}

\begin{exercise}
\label{exe.words.lex}Prove Proposition \ref{prop.words.lex}.

[\textbf{Hint:} No part of Proposition \ref{prop.words.lex} requires more than
straightforward case analysis. However, the proof of (a) can be simplified by
identifying the order relation $\leq$ on $\mathfrak{A}^{\ast}$ as a
restriction of the lexicographic order on the set $\mathfrak{B}^{\infty}$,
where $\mathfrak{B}$ is a suitable extension of the alphabet $\mathfrak{A}$.
What is this extension, and how to embed $\mathfrak{A}^{\ast}$ into
$\mathfrak{B}^{\infty}$ ?]
\end{exercise}

Proposition \ref{prop.words.lex} provides a set of tools for working with the
lexicographic order without having to refer to its definition; we shall use it
extensively. Proposition \ref{prop.words.lex}(h) (and its equivalent form
stating that $a\leq ac$ for every $a\in\mathfrak{A}^{\ast}$ and $c\in
\mathfrak{A}^{\ast}$) and Proposition \ref{prop.words.lex}(k) will often be
used without explicit mention.

Before we define Lyndon words, let us show two more facts about words which
will be used later. First, when do words commute?

\begin{proposition}
\label{prop.words.uv=vu.2}Let $u,v\in\mathfrak{A}^{\ast}$ satisfy $uv=vu$.
Then, there exist a $t\in\mathfrak{A}^{\ast}$ and two nonnegative integers $n
$ and $m$ such that $u=t^{n}$ and $v=t^{m}$.
\end{proposition}

\begin{proof}
We prove this by strong induction on $\ell\left(  u\right)  +\ell\left(
v\right)  $. We assume WLOG that $\ell\left(  u\right)  $ and $\ell\left(
v\right)  $ are positive (because otherwise, one of $u$ and $v$ is the empty
word, and everything is trivial). It is easy to see that either $u$ is a
prefix of $v$, or $v$ is a prefix of $u$\ \ \ \ \footnote{\textit{Proof.} The
word $u$ is a prefix of $uv$. But the word $v$ is also a prefix of $uv$ (since
$uv=vu$). Hence, Proposition \ref{prop.words.lex}(i) (applied to $a=u$, $b=v$
and $c=uv$) yields that either $u$ is a prefix of $v$, or $v$ is a prefix of
$u$, qed.}. We assume WLOG that $u$ is a prefix of $v$ (since our situation is
symmetric). Thus, we can write $v$ in the form $v=uw$ for some $w\in
\mathfrak{A}^{\ast}$. Consider this $w$. Clearly, $\ell\left(  u\right)
+\ell\left(  w\right)  =\ell\left(  \underbrace{uw}_{=v}\right)  =\ell\left(
v\right)  <\ell\left(  u\right)  +\ell\left(  v\right)  $ (since $\ell\left(
v\right)  $ is positive). Since $v=uw$, the equality $uv=vu$ becomes
$uuw=uwu$. Cancelling $u$ from this equality, we obtain $uw=wu$. Now, we can
apply Proposition \ref{prop.words.uv=vu.2} to $w$ instead of $v$ (by the
induction assumption, since $\ell\left(  u\right)  +\ell\left(  w\right)
<\ell\left(  u\right)  +\ell\left(  v\right)  $), and obtain that there exist
a $t\in\mathfrak{A}^{\ast}$ and two nonnegative integers $n$ and $m$ such that
$u=t^{n}$ and $w=t^{m}$. Consider this $t$ and these $n$ and $m$. Of course,
$u=t^{n}$ and $v=\underbrace{u}_{=t^{n}}\underbrace{w}_{=t^{m}}=t^{n}%
t^{m}=t^{n+m}$. So the induction step is complete, and Proposition
\ref{prop.words.uv=vu.2} is proven.
\end{proof}

\begin{proposition}
\label{prop.words.uv=vu.3}Let $u,v,w\in\mathfrak{A}^{\ast}$ be nonempty words
satisfying $uv\geq vu$, $vw\geq wv$ and $wu\geq uw$. Then, there exist a
$t\in\mathfrak{A}^{\ast}$ and three nonnegative integers $n$, $m$ and $p$ such
that $u=t^{n}$, $v=t^{m}$ and $w=t^{p}$.
\end{proposition}

\begin{proof}
We prove this by strong induction on $\ell\left(  u\right)  +\ell\left(
v\right)  +\ell\left(  w\right)  $. Clearly, $\ell\left(  u\right)  $,
$\ell\left(  v\right)  $ and $\ell\left(  w\right)  $ are positive (since $u$,
$v$ and $w$ are nonempty). We assume WLOG that $\ell\left(  u\right)
=\min\left\{  \ell\left(  u\right)  ,\ell\left(  v\right)  ,\ell\left(
w\right)  \right\}  $ (because there is a cyclic symmetry in our situation).
Thus, $\ell\left(  u\right)  \leq\ell\left(  v\right)  $ and $\ell\left(
u\right)  \leq\ell\left(  w\right)  $. But $vu\leq uv$. Hence, Proposition
\ref{prop.words.lex}(e) (applied to $a=v$, $b=u$, $c=u$ and $d=v$) yields that
either we have $v\leq u$ or the word $u$ is a prefix of $v$. But Proposition
\ref{prop.words.lex}(f) (applied to $a=u$, $b=w$, $c=w$ and $d=u$) yields
$u\leq w$ (since $uw\leq wu$ and $\ell\left(  u\right)  \leq\ell\left(
w\right)  $). Furthermore, $wv\leq vw$. Hence, Proposition
\ref{prop.words.lex}(e) (applied to $a=w$, $b=v$, $c=v$ and $d=w$) yields that
either we have $w\leq v$ or the word $v$ is a prefix of $w$.

From what we have found so far, it is easy to see that $u$ is a prefix of
$v$\ \ \ \ \footnote{\textit{Proof.} Assume the contrary. Then, $u$ is not a
prefix of $v$. Hence, we must have $v\leq u$ (since either we have $v\leq u$
or the word $u$ is a prefix of $v$), and in fact $v<u$ (because $v=u$ would
contradict to $u$ not being a prefix of $v$). Thus, $v<u\leq w$. But recall
that either we have $w\leq v$ or the word $v$ is a prefix of $w$. Thus, $v$
must be a prefix of $w$ (because $v<w$ rules out $w\leq v$). In other words,
there exists a $q\in\mathfrak{A}^{\ast}$ such that $w=vq$. Consider this $q$.
We have $v<u\leq w=vq$. Thus, Proposition \ref{prop.words.lex}(g) (applied to
$a=v$, $b=u$ and $c=q$) yields that $v$ is a prefix of $u$. In light of
$\ell\left(  u\right)  \leq\ell\left(  v\right)  $, this is only possible if
$v=u$, but this contradicts $v<u$. This contradiction completes this proof.}.
In other words, there exists a $v^{\prime}\in\mathfrak{A}^{\ast}$ such that
$v=uv^{\prime}$. Consider this $v^{\prime}$.

If the word $v^{\prime}$ is empty, then the statement of Proposition
\ref{prop.words.uv=vu.3} can be easily deduced from Proposition
\ref{prop.words.uv=vu.2}\footnote{\textit{Proof.} Assume that the word
$v^{\prime}$ is empty. Then, $v=uv^{\prime}$ becomes $v=u$. Therefore, $vw\geq
wv$ becomes $uw\geq wu$. Combined with $wu\geq uw$, this yields $uw=wu $.
Hence, Proposition \ref{prop.words.uv=vu.2} (applied to $w$ instead of $v$)
yields that there exist a $t\in\mathfrak{A}^{\ast}$ and two nonnegative
integers $n$ and $m$ such that $u=t^{n}$ and $w=t^{m}$. Clearly, $v=u=t^{n}$
as well, and so the statement of Proposition \ref{prop.words.uv=vu.3} is
true.}. Thus, we assume WLOG that this is not the case. Hence, $v^{\prime}$ is nonempty.

Using $v=uv^{\prime}$, we can rewrite $uv\geq vu$ as $uuv^{\prime}\geq
uv^{\prime}u$. That is, $uv^{\prime}u\leq uuv^{\prime}$, so that $v^{\prime
}u\leq uv^{\prime}$ (by Proposition \ref{prop.words.lex}(c), applied to $a=u$,
$c=v^{\prime}u$ and $d=uv^{\prime}$). That is, $uv^{\prime}\geq v^{\prime}u$.
But $\ell\left(  uw\right)  =\ell\left(  u\right)  +\ell\left(  w\right)
=\ell\left(  w\right)  +\ell\left(  u\right)  =\ell\left(  wu\right)  \geq
\ell\left(  wu\right)  $. Hence, Proposition \ref{prop.words.lex}(i) (applied
to $a=uw$, $b=wu$ and $c=v^{\prime}$) yields $uwv^{\prime}\leq wuv^{\prime}$
(since $uw\leq wu$). Now, $\underbrace{uv^{\prime}}_{=v}w=vw\geq
w\underbrace{v}_{=uv^{\prime}}=wuv^{\prime}\geq uwv^{\prime}$ (since
$uwv^{\prime}\leq wuv^{\prime}$), so that $uwv^{\prime}\leq uv^{\prime}w$.
Hence, $wv^{\prime}\leq v^{\prime}w$ (by Proposition \ref{prop.words.lex}(c),
applied to $a=u$, $c=wv^{\prime}$ and $d=v^{\prime}w$), so that $v^{\prime
}w\geq wv^{\prime}$. Now, we can apply Proposition \ref{prop.words.uv=vu.3} to
$v^{\prime}$ instead of $v$ (by the induction hypothesis, because
$\underbrace{\ell\left(  u\right)  +\ell\left(  v^{\prime}\right)
}_{\substack{=\ell\left(  uv^{\prime}\right)  =\ell\left(  v\right)
\\\text{(since }uv^{\prime}=v\text{)}}}+\ell\left(  w\right)  =\ell\left(
v\right)  +\ell\left(  w\right)  <\ell\left(  u\right)  +\ell\left(  v\right)
+\ell\left(  w\right)  $). As a result, we see that there exist a
$t\in\mathfrak{A}^{\ast}$ and three nonnegative integers $n$, $m$ and $p$ such
that $u=t^{n}$, $v^{\prime}=t^{m}$ and $w=t^{p}$. Clearly, this $t$ and these
$n,m,p$ satisfy $v=\underbrace{u}_{=t^{n}}\underbrace{v^{\prime}}_{=t^{m}%
}=t^{n}t^{m}=t^{n+m}$, and so the statement of Proposition
\ref{prop.words.uv=vu.3} is satisfied. The induction step is thus complete.
\end{proof}

\begin{corollary}
\label{cor.words.uv=vu.3}Let $u,v,w\in\mathfrak{A}^{\ast}$ be words satisfying
$uv\geq vu$ and $vw\geq wv$. Assume that $v$ is nonempty. Then, $uw\geq wu$.
\end{corollary}

\begin{proof}
Assume the contrary. Thus, $uw<wu$, so that $wu\geq uw$.

If $u$ or $w$ is empty, then everything is obvious. We thus WLOG assume that
$u$ and $w$ are nonempty. Thus, Proposition \ref{prop.words.uv=vu.3} shows
that there exist a $t\in\mathfrak{A}^{\ast}$ and three nonnegative integers
$n$, $m$ and $p$ such that $u=t^{n}$, $v=t^{m}$ and $w=t^{p}$. But this yields
$wu=t^{p}t^{n}=t^{p+n}=t^{n+p}=\underbrace{t^{n}}_{=u}\underbrace{t^{p}}%
_{=w}=uw$, contradicting $uw<wu$. This contradiction finishes the proof.
\end{proof}

\begin{exercise}
\label{exe.words.uw-wu}Find an alternative proof of Corollary
\ref{cor.words.uv=vu.3} which does not use Proposition
\ref{prop.words.uv=vu.3}.
\end{exercise}

The above results have a curious consequence, which we are not going to use:

\begin{corollary}
\label{cor.words.uv=vu.partord}We can define a preorder on the set
$\mathfrak{A}^{\ast}\setminus\left\{  \varnothing\right\}  $ of all nonempty
words by defining a nonempty word $u$ to be greater-or-equal to a nonempty
word $v$ (with respect to this preorder) if and only if $uv\geq vu$. Two
nonempty words $u,v$ are equivalent with respect to the equivalence relation
induced by this preorder if and only if there exist a $t\in\mathfrak{A}^{\ast
}$ and two nonnegative integers $n$ and $m$ such that $u=t^{n}$ and $v=t^{m}$.
\end{corollary}

\begin{proof}
The alleged preorder is transitive (by Corollary \ref{cor.words.uv=vu.3}) and
reflexive (obviously), and hence is really a preorder. The claim in the second
sentence follows from Proposition \ref{prop.words.uv=vu.2}.
\end{proof}

As another consequence of Proposition \ref{prop.words.uv=vu.3}, we obtain a
classical property of words \cite[Proposition 1.3.1]{Lothaire}:

\begin{exercise}
\label{exe.words.u=pn=qm}Let $u$ and $v$ be words and $n$ and $m$ be positive
integers such that $u^{n}=v^{m}$. Prove that there exists a word $t$ and
positive integers $i$ and $j$ such that $u=t^{i}$ and $v=t^{j}$.
\end{exercise}

Here is another application of Corollary \ref{cor.words.uv=vu.3}:

\begin{exercise}
\label{exe.words.unvm.vmun}Let $n$ and $m$ be positive integers. Let
$u\in\mathfrak{A}^{\ast}$ and $v\in\mathfrak{A}^{\ast}$ be two words. Prove
that $uv\geq vu$ holds if and only if $u^{n}v^{m}\geq v^{m}u^{n}$ holds.
\end{exercise}

\begin{exercise}
\label{exe.words.un.vm}Let $n$ and $m$ be positive integers. Let
$u\in\mathfrak{A}^{\ast}$ and $v\in\mathfrak{A}^{\ast}$ be two words
satisfying $n\ell\left(  u\right)  =m\ell\left(  v\right)  $. Prove that
$uv\geq vu$ holds if and only if $u^{n}\geq v^{m}$ holds.
\end{exercise}

We can also generalize Propositions \ref{prop.words.uv=vu.2} and
\ref{prop.words.uv=vu.3}:

\begin{exercise}
\label{exe.words.uv=vu.k}Let $u_{1},u_{2},\ldots,u_{k}$ be nonempty words such
that every $i\in\left\{  1,2,\ldots,k\right\}  $ satisfies $u_{i}u_{i+1}\geq
u_{i+1}u_{i}$, where $u_{k+1}$ means $u_{1}$. Show that there exist a word $t
$ and nonnegative integers $n_{1},n_{2},\ldots,n_{k}$ such that $u_{1}%
=t^{n_{1}} $, $u_{2}=t^{n_{2}}$, $\ldots$, $u_{k}=t^{n_{k}}$.
\end{exercise}

Now, we define the notion of a Lyndon word. There are several definitions in
literature, some of which will be proven equivalent in Theorem
\ref{thm.words.lyndon.equiv}.

\begin{definition}
\label{def.words.lyndon}A word $w\in\mathfrak{A}^{\ast}$ is said to be
\emph{Lyndon}\index{Lyndon word} if it is nonempty and satisfies the
following property: Every
nonempty proper suffix $v$ of $w$ satisfies $v>w$.
\end{definition}

For example, the word $113$ is Lyndon (because its nonempty proper suffixes
are $13$ and $3$, and these are both $>113$), and the word $242427$ is Lyndon
(its nonempty proper suffixes are $42427$, $2427$, $427$, $27$ and $7$, and
again these are each $>242427$). The words $2424$ and $35346$ are not Lyndon
(the word $2424$ has a nonempty proper suffix $24\leq2424$, and the word
$35346$ has a nonempty proper suffix $346\leq35346$). Every word of length $1$
is Lyndon (since it has no nonempty proper suffixes). A word $w=\left(
w_{1},w_{2}\right)  $ with two letters is Lyndon if and only if $w_{1}<w_{2}$.
A word $w=\left(  w_{1},w_{2},w_{3}\right)  $ of length $3$ is Lyndon if and
only if $w_{1}<w_{3}$ and $w_{1}\leq w_{2}$. A four-letter word $w=\left(
w_{1},w_{2},w_{3},w_{4}\right)  $ is Lyndon if and only if $w_{1}<w_{4}$,
$w_{1}\leq w_{3}$, $w_{1}\leq w_{2}$ and $\left(  \text{if }w_{1}=w_{3}\text{
then }w_{2}<w_{4}\right)  $. (These rules only get more complicated as the
words grow longer.)

We will show several properties of Lyndon words now. We begin with
trivialities which will make some arguments a bit shorter:

\begin{proposition}
\label{prop.words.lyndon.trivia}Let $w$ be a Lyndon word. Let $u$ and $v$ be
words such that $w=uv$.

\begin{itemize}
\item[(a)] If $v$ is nonempty, then $v\geq w$.

\item[(b)] If $v$ is nonempty, then $v>u$.

\item[(c)] If $u$ and $v$ are nonempty, then $vu>uv$.

\item[(d)] We have $vu\geq uv$.
\end{itemize}
\end{proposition}

\begin{proof}
(a) Assume that $v$ is nonempty. Clearly, $v$ is a suffix of $w$ (since
$w=uv$). If $v$ is a proper suffix of $w$, then the definition of a Lyndon
word yields that $v>w$ (since $w$ is a Lyndon word); otherwise, $v$ must be
$w$ itself. In either case, we have $v\geq w$. Hence, Proposition
\ref{prop.words.lyndon.trivia}(a) is proven.

(b) Assume that $v$ is nonempty. From Proposition
\ref{prop.words.lyndon.trivia}(a), we obtain $v\geq w=uv>u$ (since $v$ is
nonempty). This proves Proposition \ref{prop.words.lyndon.trivia}(b).

(c) Assume that $u$ and $v$ are nonempty. Since $u$ is nonempty, we have
$vu>v\geq w$ (by Proposition \ref{prop.words.lyndon.trivia}(a)). Since $w=uv$,
this becomes $vu>uv$. This proves Proposition \ref{prop.words.lyndon.trivia}(c).

(d) We need to prove that $vu\geq uv$. If either $u$ or $v$ is empty, $vu$ and
$uv$ are obviously equal, and thus $vu\geq uv$ is true in this case. Hence, we
can WLOG assume that $u$ and $v$ are nonempty. Assume this. Then, $vu\geq uv$
follows from Proposition \ref{prop.words.lyndon.trivia}(c). This proves
Proposition \ref{prop.words.lyndon.trivia}(d).
\end{proof}

\begin{corollary}
\label{cor.words.lyndon.trivia}Let $w$ be a Lyndon word. Let $v$ be a nonempty
suffix of $w$. Then, $v\geq w$.
\end{corollary}

\begin{proof}
Since $v$ is a nonempty suffix of $w$, there exists $u\in\mathfrak{A}^{\ast}$
such that $w=uv$. Thus, $v\geq w$ follows from Proposition
\ref{prop.words.lyndon.trivia}(a).
\end{proof}

Our next proposition is \cite[Lemma 6.5.4]{HazewinkelGubareniKirichenko}; its
part (a) is also \cite[(5.1.2)]{Reutenauer}:

\begin{proposition}
\label{prop.words.lyndon.concat}Let $u$ and $v$ be two Lyndon words such that
$u<v$. Then:

\begin{itemize}
\item[(a)] The word $uv$ is Lyndon.

\item[(b)] We have $uv<v$.
\end{itemize}
\end{proposition}

\begin{proof}
(b) The word $u$ is Lyndon and thus nonempty. Hence, $uv\neq v$%
\ \ \ \ \footnote{\textit{Proof.} Assume the contrary. Then, $uv=v$. Thus,
$uv=v=\varnothing v$. Cancelling $v$ from this equation, we obtain
$u=\varnothing$. That is, $u$ is empty. This contradicts the fact that $u$ is
nonempty. This contradiction proves that our assumption was wrong, qed.}. If
$uv\leq v\varnothing$, then Proposition \ref{prop.words.lyndon.concat}(b)
easily follows\footnote{\textit{Proof.} Assume that $uv\leq v\varnothing$.
Thus, $uv\leq v\varnothing=v$. Since $uv\neq v$, this becomes $uv<v$, so that
Proposition \ref{prop.words.lyndon.concat}(b) is proven.}. Hence, for the rest
of this proof, we can WLOG assume that we don't have $uv\leq v\varnothing$.
Assume this.

We have $u<v$. Hence, Proposition \ref{prop.words.lex}(d) (applied to $a=u$,
$b=v$, $c=v$ and $d=\varnothing$) yields that either we have $uv\leq
v\varnothing$ or the word $u$ is a prefix of $v$. Since we don't have $uv\leq
v\varnothing$, we thus see that the word $u$ is a prefix of $v$. In other
words, there exists a $t\in\mathfrak{A}^{\ast}$ satisfying $v=ut$. Consider
this $t$. Then, $t$ is nonempty (else we would have $v=u\underbrace{t}%
_{=\varnothing}=u$ in contradiction to $u<v$).

Now, $v=ut$. Hence, $t$ is a proper suffix of $v$ (proper because $u$ is
nonempty). Thus, $t$ is a nonempty proper suffix of $v$. Since every nonempty
proper suffix of $v$ is $>v$ (because $v$ is Lyndon), this shows that $t>v$.
Hence, $v\leq t$. Thus, Proposition \ref{prop.words.lex}(b) (applied to $a=u$,
$c=v$ and $d=t$) yields $uv\leq ut=v$. Combined with $uv\neq v$, this yields
$uv<v$. Hence, Proposition \ref{prop.words.lyndon.concat}(b) is proven.

(a) The word $v$ is nonempty (since it is Lyndon). Hence, $uv$ is nonempty. It
thus remains to check that every nonempty proper suffix $p$ of $uv$ satisfies
$p>uv$.

So let $p$ be a nonempty proper suffix of $uv$. We must show that $p>uv$.
Since $p$ is a nonempty proper suffix of $uv$, we must be in one of the
following two cases (depending on whether this suffix begins before the suffix
$v$ of $uv$ begins or afterwards):

\textit{Case 1:} The word $p$ is a nonempty suffix of $v$. (Note that $p=v$ is allowed.)

\textit{Case 2:} The word $p$ has the form $qv$ where $q$ is a nonempty proper
suffix of $u$.

Let us first handle Case 1. In this case, $p$ is a nonempty suffix of $v$.
Since $v$ is Lyndon, this yields that $p\geq v$ (by Corollary
\ref{cor.words.lyndon.trivia}, applied to $v$ and $p$ instead of $w$ and $v$).
But Proposition \ref{prop.words.lyndon.concat}(b) yields $uv<v$, thus $v>uv$.
Hence, $p\geq v>uv$. We thus have proven $p>uv$ in Case 1.

Let us now consider Case 2. In this case, $p$ has the form $qv$ where $q$ is a
nonempty proper suffix of $u$. Consider this $q$. Clearly, $q>u$ (since $u$ is
Lyndon and since $q$ is a nonempty proper suffix of $u$), so that $u\leq q$.
Thus, Proposition \ref{prop.words.lex}(d) (applied to $a=u$, $b=v$, $c=q$ and
$d=v$) yields that either we have $uv\leq qv$ or the word $u$ is a prefix of
$q$. Since $u$ being a prefix of $q$ is impossible (in fact, $q$ is a proper
suffix of $u$, thus shorter than $u$), we thus must have $uv\leq qv$. Since
$uv\neq qv$ (because otherwise we would have $uv=qv$, thus $u=q$ (because we
can cancel $v$ from the equality $uv=qv$), contradicting $q>u$), this can be
strengthened to $uv<qv=p$. Thus, $p>uv$ is proven in Case 2 as well.

Now that $p>uv$ is shown to hold in both cases, we conclude that $p>uv$ always holds.

Now, let us forget that we fixed $p$. We have thus shown that every nonempty
proper suffix $p$ of $uv$ satisfies $p>uv$. Since $uv$ is nonempty, this
yields that $uv$ is Lyndon (by the definition of a Lyndon word). Thus, the
proof of Proposition \ref{prop.words.lyndon.concat}(a) is complete.
\end{proof}

\begin{commentedout}

\begin{proof}
[Alternative Proof of Proposition~\ref{prop.words.lyndon.concat}%
.][\textbf{Note:} The following proof of
Proposition~\ref{prop.words.lyndon.concat} is the first proof that came into
my mind. It is similar to the proof above, but organized in a less clean way
(and not as easily generalizable to poset alphabets).]

(a) It is clear that $uv$ is nonempty. It thus remains to check that every
nonempty proper suffix $p$ of $uv$ satisfies $p>uv$.

So let $p$ be a nonempty proper suffix of $uv$. We must show that $p>uv$.
Since $p$ is a nonempty proper suffix of $uv$, we must be in one of the
following two cases (depending on whether this suffix begins before the suffix
$v$ of $uv$ begins or afterwards):

\textit{Case 1:} The word $p$ is a nonempty suffix of $v$. (Note that $p=v$ is allowed.)

\textit{Case 2:} The word $p$ has the form $qv$ where $q$ is a nonempty proper
suffix of $u$.

Let us first handle Case 1. In this case, $p$ is a nonempty suffix of $v$.
Since $v$ is Lyndon, this yields that $p\geq v$ (by Corollary
\ref{cor.words.lyndon.trivia}, applied to $v$ and $p$ instead of $w$ and $v$).
We now assume (for the sake of contradiction) that $p\leq uv$. Then, $u<v\leq
p\leq uv$. Proposition \ref{prop.words.lex}(g) (applied to $a=u$, $b=p$ and
$c=v$) thus yields $u$ is a prefix of $p$. Thus, there exists a $t\in
\mathfrak{A}^{\ast}$ satisfying $p=ut$. Consider this $t$. Then, $t$ is
nonempty (else we would have $p=u\underbrace{t}_{=\varnothing}=u$ in
contradiction to $u<p$).

Now, $v=r\underbrace{p}_{=ut}=rut=\left(  ru\right)  t$. Hence, $t$ is a
nonempty proper suffix of $v$ (proper because $u$ is nonempty). Since every
nonempty proper suffix of $v$ is $>v$ (because $v$ is Lyndon), this shows that
$t>v$. Hence, $v\leq t$. Thus, Proposition \ref{prop.words.lex}(b) (applied to
$a=u$, $c=v$ and $d=t$) yields $uv\leq ut=p$. Combined with $p\leq uv$, this
yields $uv=p=ut$. Cancelling $u$ in this equality, we obtain $v=t$,
contradicting $t>v$. This contradiction dispels our assumption that $p\leq
uv$. Hence, $p>uv$. We thus have proven $p>uv$ in Case 1.

Let us now consider Case 2. In this case, $p$ has the form $qv$ where $q$ is a
nonempty proper suffix of $u$. Consider this $q$. Clearly, $q>u$ (since $u$ is
Lyndon), so that $u\leq q$. Thus, Proposition \ref{prop.words.lex}(d) (applied
to $a=u$, $b=v$, $c=q$ and $d=v$) yields that either we have $uv\leq qv$ or
the word $u$ is a prefix of $q$. Since $u$ being a prefix of $q $ is
impossible (in fact, $q$ is a proper suffix of $u$, thus shorter than $u $),
we thus must have $uv\leq qv$. Since $uv\neq qv$ (because otherwise we would
have $uv=qv$, thus $u=q$ (because we can cancel $v$ from the equality
$uv=qv$), contradicting $q>u$), this can be strengthened to $uv<qv=p $. Thus,
$p>uv$ is proven in Case 2 as well.

Now that $p>uv$ is shown to hold in both cases, the proof of Proposition
\ref{prop.words.lyndon.concat}(a) is complete.

(b) We know from Proposition \ref{prop.words.lyndon.concat}(a) that $uv$ is
Lyndon. Hence, Proposition \ref{prop.words.lyndon.trivia}(a) (applied to
$w=uv$) yields $v\geq uv$. Since $v\neq uv$, this yields $v>uv$. This proves
Proposition \ref{prop.words.lyndon.concat}(b).
\end{proof}
\end{commentedout}

Proposition \ref{prop.words.lyndon.concat}(b), combined with Corollary
\ref{cor.words.uv=vu.3}, leads to a technical result which we will find good
use for later:

\begin{corollary}
\label{cor.words.lyndon.uv=vu}Let $u$ and $v$ be two Lyndon words such that
$u<v$. Let $z$ be a word such that $zv\geq vz$ and $uz\geq zu$. Then, $z$ is
the empty word.
\end{corollary}

\begin{proof}
Assume the contrary. Then, $z$ is nonempty. Thus, Corollary
\ref{cor.words.uv=vu.3} (applied to $z$ and $v$ instead of $v$ and $w$) yields
$uv\geq vu$. But Proposition \ref{prop.words.lyndon.concat}(b) yields
$uv<v\leq vu$, contradicting $uv\geq vu$. This contradiction completes our proof.
\end{proof}

We notice that the preorder of Corollary \ref{cor.words.uv=vu.partord} becomes
particularly simple on Lyndon words:

\begin{proposition}
\label{prop.words.lyndon.preorder}Let $u$ and $v$ be two Lyndon words. Then,
$u\geq v$ if and only if $uv\geq vu$.
\end{proposition}

\begin{proof}
We distinguish between three cases:

\textit{Case 1:} We have $u<v$.

\textit{Case 2:} We have $u=v$.

\textit{Case 3:} We have $u>v$.

Let us consider Case 1. In this case, we have $u<v$. Thus,%
\begin{align*}
uv  &  <v\ \ \ \ \ \ \ \ \ \ \left(  \text{by Proposition
\ref{prop.words.lyndon.concat}(b)}\right) \\
&  \leq vu.
\end{align*}
Hence, we have neither $u\geq v$ nor $uv\geq vu$ (because we have $u<v$ and
$uv<vu$). Thus, Proposition \ref{prop.words.lyndon.preorder} is proven in Case 1.

In Case 2, we have $u=v$. Therefore, in Case 2, both inequalities $u\geq v$
and $uv\geq vu$ hold (and actually are equalities). Thus, Proposition
\ref{prop.words.lyndon.preorder} is proven in Case 2 as well.

Let us finally consider Case 3. In this case, we have $u>v$. In other words,
$v<u$. Thus,%
\begin{align*}
vu  &  <u\ \ \ \ \ \ \ \ \ \ \left(  \text{by Proposition
\ref{prop.words.lyndon.concat}(b), applied to }v\text{ and }u\text{ instead of
}u\text{ and }v\right) \\
&  \leq uv.
\end{align*}
Hence, we have both $u\geq v$ and $uv\geq vu$ (because we have $v<u$ and
$vu<uv$). Thus, Proposition \ref{prop.words.lyndon.preorder} is proven in Case 3.

Proposition \ref{prop.words.lyndon.preorder} is now proven in all three
possible cases.
\end{proof}

\begin{proposition}
\label{prop.words.lyndon.cutoff}Let $w$ be a nonempty word. Let $v$ be the
(lexicographically) smallest nonempty suffix of $w$. Then:

\begin{itemize}
\item[(a)] The word $v$ is a Lyndon word.

\item[(b)] Assume that $w$ is not a Lyndon word. Then there exists a nonempty
$u\in\mathfrak{A}^{\ast}$ such that $w=uv$, $u\geq v$ and $uv\geq vu$.
\end{itemize}
\end{proposition}

\begin{proof}
(a) Every nonempty proper suffix of $v$ is $\geq v$ (since every nonempty
proper suffix of $v$ is a nonempty suffix of $w$, but $v$ is the smallest such
suffix) and therefore $>v$ (since a proper suffix of $v$ cannot be $=v$).
Combined with the fact that $v$ is nonempty, this yields that $v$ is Lyndon.
Proposition \ref{prop.words.lyndon.cutoff}(a) is proven.

(b) Assume that $w$ is not a Lyndon word. Then, $w\neq v$ (since $v$ is Lyndon
(by Proposition \ref{prop.words.lyndon.cutoff}(a)) while $w$ is not). Now, $v$
is a suffix of $w$. Thus, there exists an $u\in\mathfrak{A}^{\ast}$ such that
$w=uv$. Consider this $u$. Clearly, $u$ is nonempty (since $uv=w\neq v$).
Assume (for the sake of contradiction) that $u<v$. Let $v^{\prime}$ be the
(lexicographically) smallest nonempty suffix of $u$. Then, $v^{\prime}$ is a
Lyndon word (by Proposition \ref{prop.words.lyndon.cutoff}(a), applied to $u$
and $v^{\prime}$ instead of $w$ and $v$) and satisfies $v^{\prime}\leq u$
(since $u$ is a nonempty suffix of $u$, whereas $v^{\prime}$ is the smallest
such suffix). Thus, $v^{\prime}$ and $v$ are Lyndon words such that
$v^{\prime}\leq u<v$. Proposition \ref{prop.words.lyndon.concat}(a) (applied
to $v^{\prime}$ instead of $u$) now yields that the word $v^{\prime}v$ is
Lyndon. Hence, every nonempty proper suffix of $v^{\prime}v$ is $>v^{\prime}%
v$. Since $v$ is a nonempty proper suffix of $v^{\prime}v$, this yields that
$v>v^{\prime}v$.

But $v^{\prime}$ is a nonempty suffix of $u$, so that $v^{\prime}v$ is a
nonempty suffix of $uv=w$. Since $v$ is the smallest such suffix, this yields
that $v^{\prime}v\geq v$. This contradicts $v>v^{\prime}v$. Our assumption
(that $u<v$) therefore falls. We conclude that $u\geq v$.

It remains to prove that $uv\geq vu$. Assume the contrary. Then, $uv<vu$.
Thus, there exists at least one suffix $t$ of $u$ such that $tv<vt$ (namely,
$t=u$). Let $p$ be the \textbf{minimum-length} such suffix. Then, $pv<vp$.
Thus, $p$ is nonempty.

Since $p$ is a suffix of $u$, it is clear that $pv$ is a suffix of $uv=w$. So
we know that $pv$ is a nonempty suffix of $w$. Since $v$ is the smallest such
suffix, this yields that $v\leq pv<vp$. Thus, Proposition \ref{prop.words.lex}%
(g) (applied to $a=v$, $b=pv$ and $c=p$) yields that $v$ is a prefix of $pv$.
In other words, there exists a $q\in\mathfrak{A}^{\ast}$ such that $pv=vq$.
Consider this $q$. This $q$ is nonempty (because otherwise we would have
$pv=v\underbrace{q}_{=\varnothing}=v$, contradicting the fact that $p$ is
nonempty). From $vq=pv<vp$, we obtain $q\leq p$ (by Proposition
\ref{prop.words.lex}(c), applied to $a=v$, $c=q$ and $d=p$).

We know that $q$ is a suffix of $pv$ (since $vq=pv$), whereas $pv$ is a suffix
of $w$. Thus, $q$ is a suffix of $w$. So $q$ is a nonempty suffix of $w$.
Since $v$ is the smallest such suffix, this yields that $v\leq q$. We now have
$v\leq q\leq p\leq pv<vp$. Hence, $v$ is a prefix of $p$ (by Proposition
\ref{prop.words.lex}(g), applied to $a=v$, $b=p$ and $c=p$). In other words,
there exists an $r\in\mathfrak{A}^{\ast}$ such that $p=vr$. Consider this $r$.
Clearly, $r$ is a suffix of $p$, while $p$ is a suffix of $u$; therefore, $r$
is a suffix of $u$. Also, $pv<vp$ rewrites as $vrv<vvr$ (because $p=vr$).
Thus, Proposition \ref{prop.words.lex}(c) (applied to $a=v$, $c=rv$ and
$d=vr$) yields $rv\leq vr$. Since $rv\neq vr$ (because otherwise, we would
have $rv=vr$, thus $v\underbrace{rv}_{=vr}=vvr$, contradicting $vrv<vvr$),
this becomes $rv<vr$.

Now, $r$ is a suffix of $u$ such that $rv<vr$. Since $p$ is the minimum-length
such suffix, this yields $\ell\left(  r\right)  \geq\ell\left(  p\right)  $.
But this contradicts the fact that $\ell\left(  \underbrace{p}_{=vr}\right)
=\ell\left(  vr\right)  =\underbrace{\ell\left(  v\right)  }_{>0}+\ell\left(
r\right)  >\ell\left(  r\right)  $. This contradiction proves our assumption
wrong; thus, we have shown that $uv\geq vu$. Proposition
\ref{prop.words.lyndon.cutoff}(b) is proven.
\end{proof}

\begin{theorem}
\label{thm.words.lyndon.equiv}Let $w$ be a nonempty word. The following four
assertions are equivalent:

\begin{itemize}
\item \textit{Assertion $\mathcal{A}$:} The word $w$ is Lyndon.

\item \textit{Assertion $\mathcal{B}$:} Any nonempty words $u$ and
$v$ satisfying $w=uv$ satisfy $v>w$.

\item \textit{Assertion $\mathcal{C}$:} Any nonempty words $u$ and
$v$ satisfying $w=uv$ satisfy $v>u$.

\item \textit{Assertion $\mathcal{D}$:} Any nonempty words $u$ and
$v$ satisfying $w=uv$ satisfy $vu>uv$.
\end{itemize}
\end{theorem}

\begin{proof}
\textit{Proof of the implication $\mathcal{A}\Longrightarrow\mathcal{B}%
$:} If Assertion $\mathcal{A}$ holds, then Assertion $\mathcal{B}$
clearly holds (in fact, whenever $u$ and $v$ are nonempty words satisfying
$w=uv$, then $v$ is a nonempty proper suffix of $w$, and therefore $>w$ by the
definition of a Lyndon word).

\textit{Proof of the implication $\mathcal{A}\Longrightarrow\mathcal{C}%
$:} This implication follows from Proposition
\ref{prop.words.lyndon.trivia}(b).

\textit{Proof of the implication $\mathcal{A}\Longrightarrow\mathcal{D}%
$:} This implication follows from Proposition
\ref{prop.words.lyndon.trivia}(c).

\textit{Proof of the implication $\mathcal{B}\Longrightarrow\mathcal{A}%
$:} Assume that Assertion $\mathcal{B}$ holds. If $v$ is a nonempty
proper suffix of $w$, then there exists an $u\in\mathfrak{A}^{\ast}$
satisfying $w=uv$. This $u$ is nonempty because $v$ is a proper suffix, and
thus Assertion $\mathcal{B}$ yields $v>w$. Hence, every nonempty proper suffix
$v$ of $w$ satisfies $v>w$. By the definition of a Lyndon word, this yields
that $w$ is Lyndon, so that Assertion $\mathcal{A}$ holds.

\textit{Proof of the implication $\mathcal{C}\Longrightarrow\mathcal{A}%
$:} Assume that Assertion $\mathcal{C}$ holds. If $w$ was not Lyndon,
then Proposition \ref{prop.words.lyndon.cutoff}(b) would yield nonempty words
$u$ and $v$ such that $w=uv$ and $u\geq v$; this would contradict Assertion
$\mathcal{C}$. Thus, $w$ is Lyndon, and Assertion $\mathcal{A}$ holds.

\textit{Proof of the implication $\mathcal{D}\Longrightarrow\mathcal{A}%
$:} Assume that Assertion $\mathcal{D}$ holds. If $w$ was not Lyndon,
then Proposition \ref{prop.words.lyndon.cutoff}(b) would yield nonempty words
$u$ and $v$ such that $w=uv$ and $uv\geq vu$; this would contradict Assertion
$\mathcal{D}$. Thus, $w$ is Lyndon, and Assertion $\mathcal{A}$ holds.

Now we have proven enough implications to conclude the equivalence of all four assertions.
\end{proof}

Theorem \ref{thm.words.lyndon.equiv} connects our definition of Lyndon words
with some of the definitions appearing in literature. For example, Lothaire
\cite[\S 5.1]{Lothaire}, Shirshov \cite{Shirshov} and de Bruijn/Klarner
\cite[\S 4]{deBruijnKlarner} define Lyndon words using Assertion $\mathcal{D}$
(note, however, that Shirshov takes $<$ instead of $>$ and calls Lyndon words
``regular words''; also, de Bruijn/Klarner
call Lyndon words ``normal words'').
Chen-Fox-Lyndon \cite[\S 1]{ChenFoxLyndon}, Reutenauer \cite{Reutenauer} and
Radford \cite{Radford-shuffle} use our definition (but Chen-Fox-Lyndon call
the Lyndon words ``standard sequences'', and
Radford calls them ``primes'' and uses $<$
instead of $>$).

Theorem \ref{thm.words.lyndon.equiv} appears (with different notations) in
Zhou-Lu \cite[Proposition 1.4]{ZhouLu}. The equivalence $\mathcal{D}%
\Longleftrightarrow\mathcal{A}$ of our Theorem \ref{thm.words.lyndon.equiv} is
equivalent to \cite[Proposition 5.12]{Lothaire} and to \cite[$\mathfrak{A}%
^{\prime\prime}=\mathfrak{A}^{\prime\prime\prime}$]{ChenFoxLyndon}.

The following exercise provides a different (laborious) approach to Theorem
\ref{thm.words.lyndon.equiv}:

\begin{exercise}
\phantomsection\label{exe.words.lyndon.equiv2}

\begin{enumerate}
\item[(a)] Prove that if $u\in\mathfrak{A}^{\ast}$ and $v\in\mathfrak{A}%
^{\ast}$ are two words satisfying $uv<vu$, then there exists a nonempty suffix
$s$ of $u$ satisfying $sv<v$.

\item[(b)] Give a new proof of Theorem \ref{thm.words.lyndon.equiv} (avoiding
the use of Proposition \ref{prop.words.lyndon.cutoff}).
\end{enumerate}

[\textbf{Hint:} For (a), perform strong induction on $\ell\left(  u\right)
+\ell\left(  v\right)  $, assume the contrary, and distinguish between the
case when $u\leq v$ and the case when $v$ is a prefix of $u$. For (b), use
part (a) in proving the implication $\mathcal{D}\Longrightarrow\mathcal{B}$,
and factor $v$ as $v=u^{m}v^{\prime}$ with $m$ maximal in the proof of the
implication $\mathcal{C}\Longrightarrow\mathcal{B}$.]
\end{exercise}

The following two exercises are taken from \cite{Hazewinkel-lh}%
\footnote{Exercise \ref{exe.words.lyndon.power-criterion} is more or less
\cite[Lemma 4.3]{Hazewinkel-lh} with a converse added; Exercise
\ref{exe.words.lyndon.concatn} is \cite[Lemma 4.2]{Hazewinkel-lh}.}.

\begin{exercise}
\label{exe.words.lyndon.power-criterion}Let $w$ be a nonempty word. Prove that
$w$ is Lyndon if and only if every nonempty word $t$ and every positive
integer $n$ satisfy $\left(  \text{if }w\leq t^{n}\text{, then }w\leq
t\right)  $.
\end{exercise}

\begin{exercise}
\label{exe.words.lyndon.concatn}Let $w_{1}$, $w_{2}$, $\ldots$, $w_{n}$ be $n$
Lyndon words, where $n$ is a positive integer. Assume that $w_{1}\leq
w_{2}\leq\cdots\leq w_{n}$ and $w_{1}<w_{n}$. Show that $w_{1}w_{2}\cdots
w_{n}$ is a Lyndon word.
\end{exercise}

The following exercise is a generalization (albeit not in an obvious way) of
Exercise \ref{exe.words.lyndon.concatn}:

\begin{exercise}
\label{exe.words.lyndon.concatn-lyndon}Let $w_{1}$, $w_{2}$, $\ldots$, $w_{n}$
be $n$ Lyndon words, where $n$ is a positive integer. Assume that
$w_{i}w_{i+1}\cdots w_{n}\geq w_{1}w_{2}\cdots w_{n}$ for every $i\in\left\{
1,2,\ldots,n\right\}  $. Show that $w_{1}w_{2}\cdots w_{n}$ is a Lyndon word.
\end{exercise}

We are now ready to meet one of the most important features of Lyndon
words: a bijection between all words and multisets of Lyndon
words\footnote{And it is not even the only such bijection: we will see
another in Subsection~\ref{subsub.lyndon.gr.bij}.}; it is clear that
such a bijection is vital for constructing polynomial generating sets
of commutative algebras with bases indexed by words, such as $\Qsym$ or
shuffle algebras. This bijection is given by the
\emph{Chen-Fox-Lyndon factorization}:

\begin{definition}
\label{def.words.CFL}Let $w$ be a word. A \dfn{Chen-Fox-Lyndon
factorization} (in short, \dfn{CFL factorization}) of $w$ means a tuple
$\left(  a_{1},a_{2},\ldots,a_{k}\right)  $ of Lyndon words satisfying
$w=a_{1}a_{2} \cdots a_{k}$ and $a_{1}\geq a_{2}\geq\cdots\geq a_{k}$.
\end{definition}

\begin{example}
The tuple $\left(  23,2,14,13323,13,12,12,1\right)  $ is a CFL factorization
of the word $23214133231312121$ over the alphabet $\left\{  1,2,3,\ldots
\right\}  $ (ordered by $1<2<3<\cdots$), since $23$, $2$, $14$, $13323$, $13$,
$12$, $12$ and $1$ are Lyndon words satisfying $23214133231312121=23\cdot
2\cdot14\cdot13323\cdot13\cdot12\cdot12\cdot1$ and $23\geq2\geq14\geq
13323\geq13\geq12\geq12\geq1$.
\end{example}

The bijection is given by the following \dfn{Chen-Fox-Lyndon theorem}
(\cite[Theorem 6.5.5]{HazewinkelGubareniKirichenko}, \cite[Thm. 5.1.5]%
{Lothaire}, \cite[part of Thm. 2.1.4]{Radford-shuffle}):

\begin{theorem}
\label{thm.words.CFL}Let $w$ be a word. Then, there exists a unique CFL
factorization of $w$.
\end{theorem}

Before we prove this, we need to state and prove a lemma
(which is \cite[Proposition 5.1.6]{Lothaire}):

\begin{lemma}
\label{lem.words.CFLlemma}Let $\left( a_1, a_2, \ldots, a_k \right)  $ be a
CFL factorization of a nonempty word $w$. Let $p$ be a nonempty suffix of $w$.
Then, $p\geq a_k$.
\end{lemma}

\begin{proof}%[Proof of Lemma \ref{lem.words.CFLlemma}]
We will prove Lemma \ref{lem.words.CFLlemma} by induction over the
(obviously) positive integer $k$.

\textit{Induction base:} Assume that $k=1$. Thus, $\left(  a_{1},a_{2}%
,\ldots,a_{k}\right)  =\left(  a_{1}\right)  $ is a tuple of Lyndon words
satisfying $w=a_{1}a_{2}\cdots a_{k}$. We have $w=a_{1}a_{2}\cdots a_{k}%
=a_{1}$ (since $k=1$), so that $w$ is a Lyndon word (since $a_{1}$ is a Lyndon
word). Thus, Corollary \ref{cor.words.lyndon.trivia} (applied to $v=p$) yields
$p\geq w=a_{1}=a_{k}$ (since $1=k$). Thus, Lemma \ref{lem.words.CFLlemma} is
proven in the case $k=1$. The induction base is complete.

\textit{Induction step:} Let $K$ be a positive integer. Assume (as the
induction hypothesis) that Lemma \ref{lem.words.CFLlemma} is proven for $k=K$.
We now need to show that Lemma \ref{lem.words.CFLlemma} holds for $k=K+1$.

So let $\left(  a_{1},a_{2},\ldots,a_{K+1}\right)  $ be a CFL factorization of
a nonempty word $w$. Let $p$ be a nonempty suffix of $w$. We need to prove
that $p\geq a_{K+1}$.

By the definition of a CFL factorization, $\left(  a_{1},a_{2},\ldots
,a_{K+1}\right)  $ is a tuple of Lyndon words satisfying $w=a_{1}a_{2}\cdots
a_{K+1}$ and $a_{1}\geq a_{2}\geq\cdots\geq a_{K+1}$. Let $w^{\prime}%
=a_{2}a_{3}\cdots a_{K+1}$; then, $w=a_{1}a_{2}\cdots a_{K+1}=a_{1}%
\underbrace{\left(  a_{2}a_{3}\cdots a_{K+1}\right)  }_{=w^{\prime}}%
=a_{1}w^{\prime}$. Hence, every nonempty suffix of $w$ is either a nonempty
suffix of $w^{\prime}$, or has the form $qw^{\prime}$ for a nonempty suffix
$q$ of $a_{1}$. Since $p$ is a nonempty suffix of $w$, we thus must be in one
of the following two cases:

\textit{Case 1:} The word $p$ is a nonempty suffix of $w^{\prime}$.

\textit{Case 2:} The word $p$ has the form $qw^{\prime}$ for a nonempty suffix
$q$ of $a_{1}$.

Let us first consider Case 1. In this case, $p$ is a nonempty suffix of
$w^{\prime}$. The $K$-tuple $\left(  a_{2},a_{3},\ldots,a_{K+1}\right)  $ of
Lyndon words satisfies $w^{\prime}=a_{2}a_{3} \cdots a_{K+1}$ and $a_{2}\geq
a_{3}\geq\cdots\geq a_{K+1}$; therefore, $\left(  a_{2},a_{3},\ldots
,a_{K+1}\right)  $ is a CFL factorization of $w^{\prime}$. We can thus apply
Lemma \ref{lem.words.CFLlemma} to $K$, $w^{\prime}$ and $\left(  a_{2}%
,a_{3},\ldots,a_{K+1}\right)  $ instead of $k$, $w$ and $\left(  a_{1}%
,a_{2},\ldots,a_{k}\right)  $ (because we assumed that Lemma
\ref{lem.words.CFLlemma} is proven for $k=K$). As a result, we obtain that
$p\geq a_{K+1}$. Thus, $p\geq a_{K+1}$ is proven in Case 1.

Let us now consider Case 2. In this case, $p$ has the form $qw^{\prime}$ for a
nonempty suffix $q$ of $a_{1}$. Consider this $q$. Since $a_{1}$ is a Lyndon
word, we have $q\geq a_{1}$ (by Corollary \ref{cor.words.lyndon.trivia},
applied to $a_{1}$ and $q$ instead of $w$ and $v$). Thus, $q\geq a_{1}\geq
a_{2}\geq\cdots\geq a_{K+1}$, so that $p=qw^{\prime}\geq q\geq a_{K+1}$. Thus,
$p\geq a_{K+1}$ is proven in Case 2.

We have now proven $p\geq a_{K+1}$ in all cases. This proves that Lemma
\ref{lem.words.CFLlemma} holds for $k=K+1$. The induction step is thus
finished, and with it the proof of Lemma \ref{lem.words.CFLlemma}.
\end{proof}

\begin{proof}
[Proof of Theorem \ref{thm.words.CFL}]Let us first prove that there exists a
CFL factorization of $w$.

Indeed, there clearly exists a tuple $\left(  a_{1},a_{2},\ldots,a_{k}\right)
$ of Lyndon words satisfying $w=a_{1}a_{2} \cdots a_{k}$\ \ \ \ \footnote{For
instance, the tuple $\left(  w_{1},w_{2},\ldots,w_{\ell\left(  w\right)
}\right)  $ of one-letter words is a valid example (recall that one-letter
words are always Lyndon).}. Fix such a tuple with \textbf{minimum} $k$. We
claim that $a_{1}\geq a_{2}\geq\cdots\geq a_{k}$.

Indeed, if some $i\in\left\{  1,2,\ldots,k-1\right\}  $ would satisfy
$a_{i}<a_{i+1}$, then the word $a_{i}a_{i+1}$ would be Lyndon (by Proposition
\ref{prop.words.lyndon.concat}(a), applied to $u=a_{i}$ and $v=a_{i+1}$),
whence $\left(  a_{1},a_{2},\ldots,a_{i-1},a_{i}a_{i+1},a_{i+2},a_{i+3}%
,\ldots,a_{k}\right)  $ would also be a tuple of Lyndon words satisfying
$w=a_{1}a_{2} \cdots a_{i-1}\left(  a_{i}a_{i+1}\right)  a_{i+2}a_{i+3} \cdots
a_{k}$ but having length $k-1<k$, contradicting the fact that $k$ is the
minimum length of such a tuple. Hence, no $i\in\left\{  1,2,\ldots
,k-1\right\}  $ can satisfy $a_{i}<a_{i+1}$. In other words, every
$i\in\left\{  1,2,\ldots,k-1\right\}  $ satisfies $a_{i}\geq a_{i+1}$. In
other words, $a_{1}\geq a_{2}\geq\cdots\geq a_{k}$. Thus, $\left(  a_{1}%
,a_{2},\ldots,a_{k}\right)  $ is a CFL factorization of $w$, so we have shown
that such a CFL factorization exists.

It remains to show that there exists at most one CFL factorization of $w$. We
shall prove this by induction over $\ell\left(  w\right)  $. Thus, we fix a
word $w$ and assume that%
\begin{equation}
\text{for every word }v\text{ with }\ell\left(  v\right)  <\ell\left(
w\right)  \text{, there exists at most one CFL factorization of }v.
\label{pf.thm.words.CFL.uni.1}%
\end{equation}
We now have to prove that there exists at most one CFL factorization of $w$.

Indeed, let $\left(  a_{1},a_{2},\ldots,a_{k}\right)  $ and $\left(
b_{1},b_{2},\ldots,b_{m}\right)  $ be two CFL factorizations of $w$. We need
to prove that $\left(  a_{1},a_{2},\ldots,a_{k}\right)  =\left(  b_{1}%
,b_{2},\ldots,b_{m}\right)  $. If $w$ is empty, then this is obvious, so we
WLOG assume that it is not; thus, $k>0$ and $m>0$.

Since $\left(  b_{1},b_{2},\ldots,b_{m}\right)  $ is a CFL factorization of
$w$, we have $w=b_{1}b_{2} \cdots b_{m}$, and thus $b_{m}$ is a nonempty
suffix of $w$. Thus, Lemma \ref{lem.words.CFLlemma} (applied to $p=b_{m}$)
yields $b_{m}\geq a_{k}$. The same argument (but with the roles of $\left(
a_{1},a_{2},\ldots,a_{k}\right)  $ and $\left(  b_{1},b_{2},\ldots
,b_{m}\right)  $ switched) shows that $a_{k}\geq b_{m}$. Combined with
$b_{m}\geq a_{k}$, this yields $a_{k}=b_{m}$. Now let $v=a_{1}a_{2} \cdots
a_{k-1}$. Then, $\left(  a_{1},a_{2},\ldots,a_{k-1}\right)  $ is a CFL
factorization of $v$ (since $a_{1}\geq a_{2}\geq\cdots\geq a_{k-1}$).

Since $\left(  a_{1},a_{2},\ldots,a_{k}\right)  $ is a CFL factorization of
$w$, we have $w=a_{1}a_{2} \cdots a_{k}=\underbrace{a_{1}a_{2} \cdots a_{k-1}%
}_{=v}\underbrace{a_{k}}_{=b_{m}}=vb_{m}$, so that%
\[
vb_{m}=w=b_{1}b_{2} \cdots b_{m}=b_{1}b_{2} \cdots b_{m-1}b_{m}.
\]
Cancelling $b_{m}$ yields $v=b_{1}b_{2} \cdots b_{m-1}$. Thus, $\left(
b_{1},b_{2},\ldots,b_{m-1}\right)  $ is a CFL factorization of $v$ (since
$b_{1}\geq b_{2}\geq\cdots\geq b_{m-1}$). Since $\ell\left(  v\right)
<\ell\left(  w\right)  $ (because $v=a_{1}a_{2} \cdots a_{k-1}$ is shorter
than $w=a_{1}a_{2} \cdots a_{k}$), we can apply (\ref{pf.thm.words.CFL.uni.1})
to obtain that there exists at most one CFL factorization of $v$. But we
already know two such CFL factorizations: $\left(  a_{1},a_{2},\ldots
,a_{k-1}\right)  $ and $\left(  b_{1},b_{2},\ldots,b_{m-1}\right)  $. Thus,
$\left(  a_{1},a_{2},\ldots,a_{k-1}\right)  =\left(  b_{1},b_{2}%
,\ldots,b_{m-1}\right)  $, which, combined with $a_{k}=b_{m}$, leads to
$\left(  a_{1},a_{2},\ldots,a_{k}\right)  =\left(  b_{1},b_{2},\ldots
,b_{m}\right)  $. This is exactly what we needed to prove. So we have shown
(by induction) that there exists at most one CFL factorization of $w$. This
completes the proof of Theorem \ref{thm.words.CFL}.
\end{proof}

The CFL factorization allows us to count all Lyndon words of a given length if
$\mathfrak{A}$ is finite:

\begin{exercise}
\label{exe.words.lyndon.count}Assume that the alphabet $\mathfrak{A}$ is
finite. Let $q=\left\vert \mathfrak{A}\right\vert $.
Let $\mu$ be the number-theoretic M\"{o}bius
function (defined as in Exercise~\ref{exe.witt.ghost-equiv}).
Show that the number of
Lyndon words of length $n$ equals $\dfrac{1}{n}\sum\limits_{d\mid n}\mu\left(
d\right)  q^{n/d}$ for every positive integer $n$ (where ``%
$\sum\limits_{d\mid n}$'' means a sum over all positive
divisors of $n$).\ \ \ \ %
\footnote{In particular, $\dfrac{1}{n}\sum\limits_{d\mid n}%
\mu\left(  d\right)  q^{n/d}$ is an integer.}
\end{exercise}

Exercise \ref{exe.words.lyndon.count} is a well-known result and appears,
e.g., in \cite[Theorem 1.5]{ChenFoxLyndon} or in
\cite[Section 5.1]{Lothaire}.

We will now study another kind of factorization: not of an arbitrary word into
Lyndon words, but of a Lyndon word into two smaller Lyndon words. This
factorization is called \dfn{standard factorization} (\cite[\S 5.1]%
{Lothaire}) or \dfn{canonical factorization} (\cite[Lemma 6.5.33]%
{HazewinkelGubareniKirichenko}); we only introduce it from the viewpoint we
are interested in, namely its providing a way to do induction over Lyndon
words\footnote{e.g., allowing to solve Exercise
\ref{exe.words.lyndon.concatn-lyndon} in a simpler way}. Here is what we need
to know:

\begin{theorem}
\label{thm.words.lyndon.std}Let $w$ be a Lyndon word of length $>1$. Let $v$
be the (lexicographically) smallest nonempty \textbf{proper} suffix of $w$.
Since $v$ is a proper suffix of $w$, there exists a nonempty $u\in
\mathfrak{A}^{\ast}$ such that $w=uv$. Consider this $u$. Then:

\begin{itemize}
\item[(a)] The words $u$ and $v$ are Lyndon.

\item[(b)] We have $u<w<v$.
\end{itemize}
\end{theorem}

\begin{proof}
Every nonempty proper suffix of $v$ is $\geq v$ (since every nonempty proper
suffix of $v$ is a nonempty proper suffix of $w$, but $v$ is the smallest such
suffix) and therefore $>v$ (since a proper suffix of $v$ cannot be $=v$).
Combined with the fact that $v$ is nonempty, this yields that $v$ is Lyndon.

Since $w$ is Lyndon, we know that every nonempty proper suffix of $w$ is $>w$.
Applied to the nonempty proper suffix $v$ of $w$, this yields that $v>w$.
Hence, $w<v$. Since $v$ is nonempty, we have $u<uv=w<v$. This proves Theorem
\ref{thm.words.lyndon.std}(b).

Let $p$ be a nonempty proper suffix of $u$. Then, $pv$ is a nonempty proper
suffix of $uv=w$. Thus, $pv>w$ (since every nonempty proper suffix of $w$ is
$>w$). Thus, $pv>w=uv$, so that $uv<pv$. Thus, Proposition
\ref{prop.words.lex}(e) (applied to $a=u$, $b=v$, $c=p$ and $d=v$) yields that
either we have $u\leq p$ or the word $p$ is a prefix of $u$.

Let us assume (for the sake of contradiction) that $p\leq u$. Then, $p<u$
(because $p$ is a proper suffix of $u$, and therefore $p\neq u$). Hence, we
cannot have $u\leq p$. Thus, the word $p$ is a prefix of $u$ (since either we
have $u\leq p$ or the word $p$ is a prefix of $u$). In other words, there
exists a $q\in\mathfrak{A}^{\ast}$ such that $u=pq$. Consider this $q$. We
have $w=\underbrace{u}_{=pq}v=pqv=p\left(  qv\right)  $, and thus $qv$ is a
proper suffix of $w$ (proper because $p$ is nonempty). Moreover, $qv$ is
nonempty (since $v$ is nonempty). Hence, $qv$ is a nonempty proper suffix of
$w$. Since $v$ is the smallest such suffix, this entails that $v\leq qv$.
Proposition \ref{prop.words.lex}(b) (applied to $a=p$, $c=v$ and $d=qv$) thus
yields $pv\leq pqv$. Hence, $pv\leq pqv=w$, which contradicts $pv>w$. This
contradiction shows that our assumption (that $p\leq u$) was false. We thus
have $p>u$.

We now have shown that $p>u$ whenever $p$ is a nonempty proper suffix of $u$.
Combined with the fact that $u$ is nonempty, this shows that $u$ is a Lyndon
word. This completes the proof of Theorem \ref{thm.words.lyndon.std}(a).
\end{proof}

Another approach to the standard factorization is given in the following exercise:

\begin{exercise}
\label{exe.words.lyndon.std2}Let $w$ be a Lyndon word of length $>1$. Let $v$
be the longest proper suffix of $w$ such that $v$ is Lyndon\footnote{This is
well-defined, because there exists at least one proper suffix $v$ of $w$ such
that $v$ is Lyndon. (Indeed, the last letter of $w$ forms such a suffix,
because it is a proper suffix of $w$ (since $w$ has length $>1$) and is Lyndon
(since it is a one-letter word, and since every one-letter word is Lyndon).)}.
Since $v$ is a proper suffix of $w$, there exists a nonempty $u\in
\mathfrak{A}^{\ast}$ such that $w=uv$. Consider this $u$. Prove that:

\begin{enumerate}
\item[(a)] The words $u$ and $v$ are Lyndon.

\item[(b)] We have $u<w<v$.

\item[(c)] The words $u$ and $v$ are precisely the words $u$ and $v$
constructed in Theorem \ref{thm.words.lyndon.std}.
\end{enumerate}
\end{exercise}

Notice that a well-known recursive characterization of Lyndon words
\cite[$\mathfrak{A}^{\prime}=\mathfrak{A}^{\prime\prime}$]{ChenFoxLyndon} can
be easily derived from Theorem \ref{thm.words.lyndon.std} and Proposition
\ref{prop.words.lyndon.concat}(a). We will not dwell on it.

The following exercise surveys some variations on the characterizations of
Lyndon words\footnote{Compare this with
\cite[\S 7.2.11, Theorem Q]{Knuth-TAoCP4a}.}:

\begin{exercise}
\label{exe.words.powerlyndon}Let $w$ be a nonempty word. Consider the
following nine assertions:

\begin{itemize}
\item \textit{Assertion $\mathcal{A}^{\prime}$:} The word $w$ is a
power of a Lyndon word.

\item \textit{Assertion $\mathcal{B}^{\prime}$:} If $u$ and $v$ are
nonempty words satisfying $w=uv$, then either we have $v\geq w$ or the word
$v$ is a prefix of $w$.

\item \textit{Assertion $\mathcal{C}^{\prime}$:} If $u$ and $v$ are
nonempty words satisfying $w=uv$, then either we have $v\geq u$ or the word
$v$ is a prefix of $u$.

\item \textit{Assertion $\mathcal{D}^{\prime}$:} If $u$ and $v$ are
nonempty words satisfying $w=uv$, then we have $vu\geq uv$.

\item \textit{Assertion $\mathcal{E}^{\prime}$:} If $u$ and $v$ are
nonempty words satisfying $w=uv$, then either we have $v\geq u$ or the word
$v$ is a prefix of $w$.

\item \textit{Assertion $\mathcal{F}^{\prime}$:} The word $w$ is a
prefix of a Lyndon word in $\mathfrak{A}^{\ast}$.

\item \textit{Assertion $\mathcal{F}^{\prime\prime}$:} Let $m$ be an
object not in the alphabet $\mathfrak{A}$. Let us equip the set $\mathfrak{A}%
\cup\left\{  m\right\}  $ with a total order which extends the total order on
the alphabet $\mathfrak{A}$ and which satisfies $\left(  a<m\text{ for every
}a\in\mathfrak{A}\right)  $. Then, the word $wm\in\left(  \mathfrak{A}%
\cup\left\{  m\right\}  \right)  ^{\ast}$ (the concatenation of the word $w$
with the one-letter word $m$) is a Lyndon word.

\item \textit{Assertion $\mathcal{G}^{\prime}$:} There exists a
Lyndon word $t\in\mathfrak{A}^{\ast}$, a positive integer $\ell$ and a prefix
$p$ of $t$ (possibly empty) such that $w=t^{\ell}p$.

\item \textit{Assertion $\mathcal{H}^{\prime}$:} There exists a
Lyndon word $t\in\mathfrak{A}^{\ast}$, a nonnegative integer $\ell$ and a
prefix $p$ of $t$ (possibly empty) such that $w=t^{\ell}p$.

\item[(a)] Prove the equivalence $\mathcal{A}^{\prime}\Longleftrightarrow
\mathcal{D}^{\prime}$.

\item[(b)] Prove the equivalence $\mathcal{B}^{\prime}\Longleftrightarrow
\mathcal{C}^{\prime}\Longleftrightarrow\mathcal{E}^{\prime}\Longleftrightarrow
\mathcal{F}^{\prime\prime}\Longleftrightarrow\mathcal{G}^{\prime
}\Longleftrightarrow\mathcal{H}^{\prime}$.

\item[(c)] Prove the implication $\mathcal{F}^{\prime}\Longrightarrow
\mathcal{B}^{\prime}$.

\item[(d)] Prove the implication $\mathcal{D}^{\prime}\Longrightarrow
\mathcal{B}^{\prime}$. (The implication $\mathcal{B}^{\prime}\Longrightarrow
\mathcal{D}^{\prime}$ is false, as witnessed by the word $11211$.)

\item[(e)] Prove that if there exists a letter $\mu\in\mathfrak{A}$ such that
$\left(  \mu>a\text{ for every letter }a\text{ of }w\right)  $, then the
equivalence $\mathcal{F}^{\prime}\Longleftrightarrow\mathcal{F}^{\prime\prime
}$ holds.

\item[(f)] Prove that if there exists a letter $\mu\in\mathfrak{A}$ such that
$\left(  \mu>a\text{ for some letter }a\text{ of }w\right)  $, then the
equivalence $\mathcal{F}^{\prime}\Longleftrightarrow\mathcal{F}^{\prime\prime
}$ holds.
\end{itemize}
\end{exercise}

The next exercise (based on work of Hazewinkel \cite{Hazewinkel-poset})
extends some of the above properties of Lyndon words (and words in general) to
a more general setting, in which the alphabet $\mathfrak{A}$ is no longer
required to be totally ordered, but only needs to be a poset:

\begin{exercise}
\label{exe.words.poset}In this exercise, we shall loosen the requirement that
the alphabet $\mathfrak{A}$ be a totally ordered set: Instead, we will only
require $\mathfrak{A}$ to be a poset. The resulting more general setting will
be called the
\emph{partial-order setting}\index{partial-order setting for words},
to distinguish it from the
\emph{total-order setting}\index{total-order setting for words}
in which $\mathfrak{A}$ is required to be a
totally ordered set. All results in Chapter \ref{sect.QSym.lyndon} so far
address the total-order setting. In this exercise, we will generalize some of
them to the partial-order setting.

All notions that we have defined in the total-order setting (the notion of a
word, the relation $\leq$, the notion of a Lyndon word, etc.) are defined in
precisely the same way in the partial-order setting. However, the poset
$\mathfrak{A}^{\ast}$ is no longer totally ordered in the partial-order setting.

\begin{enumerate}
\item[(a)] Prove that Proposition \ref{prop.words.lex} holds in the
partial-order setting, as long as one replaces ``a total
order'' by ``a partial order'' in part (a) of this Proposition.

\item[(b)] Prove (in the partial-order setting) that if $a,b,c,d\in
\mathfrak{A}^{\ast}$ are four words such that the words $ab$ and $cd$ are
comparable (with respect to the partial order $\leq$), then the words $a$ and
$c$ are comparable.

\item[(c)] Prove that Proposition \ref{prop.words.uv=vu.2}, Proposition
\ref{prop.words.uv=vu.3}, Corollary \ref{cor.words.uv=vu.3}, Corollary
\ref{cor.words.uv=vu.partord}, Exercise \ref{exe.words.u=pn=qm}, Exercise
\ref{exe.words.unvm.vmun}, Exercise \ref{exe.words.un.vm}, Exercise
\ref{exe.words.uv=vu.k}, Proposition \ref{prop.words.lyndon.trivia}, Corollary
\ref{cor.words.lyndon.trivia}, Proposition \ref{prop.words.lyndon.concat},
Corollary \ref{cor.words.lyndon.uv=vu}, Proposition
\ref{prop.words.lyndon.preorder}, Theorem \ref{thm.words.lyndon.equiv},
Exercise \ref{exe.words.lyndon.equiv2}(a), Exercise
\ref{exe.words.lyndon.concatn}, Exercise \ref{exe.words.lyndon.concatn-lyndon}%
, Exercise \ref{exe.words.lyndon.std2}(a) and Exercise
\ref{exe.words.lyndon.std2}(b) still hold in the partial-order setting.

\item[(d)] Find a counterexample to Exercise
\ref{exe.words.lyndon.power-criterion} in the partial-order setting.

\item[(e)] Salvage Exercise \ref{exe.words.lyndon.power-criterion} in the
partial-order setting (i.e., find a statement which is easily equivalent to
this exercise in the total-order setting, yet true in the partial-order setting).

\item[(f)] In the partial-order setting, a \dfn{Hazewinkel-CFL
factorization} of a word $w$ will mean a tuple $\left(  a_{1},a_{2}%
,\ldots,a_{k}\right)  $ of Lyndon words such that $w=a_{1}a_{2}\cdots a_{k}$
and such that no $i\in\left\{  1,2,\ldots,k-1\right\}  $ satisfies
$a_{i}<a_{i+1}$. Prove that every word $w$ has a unique Hazewinkel-CFL
factorization (in the partial-order setting).\footnote{This result, as well as
the validity of Proposition \ref{prop.words.lyndon.concat} in the
partial-order setting, are due to Hazewinkel \cite{Hazewinkel-poset}.}

\item[(g)] Prove that Exercise~\ref{exe.words.powerlyndon} still holds in the
partial-order setting.
\end{enumerate}
\end{exercise}

% [DG][v45] Added part (g) to this exercise.

The reader is invited to try extending other results to the partial-order
setting (it seems that no research has been done on this except for
Hazewinkel's \cite{Hazewinkel-poset}). We shall now, however, return to
the total-order setting (which has the most known applications).

Another extension of the notion of Lyndon words has been introduced
in 2018 by Dolce, Restivo and Reutenauer
\cite{DolceRestivoReutenauer-oglw}; it is based on a generalized version
of the lexicographic order, in which different letters are compared
differently depending on their positions in the word
(i.e., there is one total order for comparing
first letters, another for comparing second letters, etc.).

Lyndon words are related to various other objects in mathematics, such as
free Lie algebras (Subsection \ref{subsect.lyndon.freelie} below),
shuffles and shuffle algebras (Sections
\ref{subsect.lyndon.shuffles} and \ref{subsect.shuffle.radford} below),
$\Qsym$ (Sections \ref{subsect.shuffle.QSym.1} and
\ref{subsect.shuffle.QSym.2}), Markov chains on combinatorial Hopf algebras
(\cite{DiaconisPangRam}), de Bruijn sequences (\cite{FredricksenMaiorana},
\cite{Moreno}, \cite{Moreno-corr}, \cite[\S 7.2.11, Algorithm F]{Knuth-TAoCP4a}),
symmetric functions (specifically, the transition matrices
between the bases $\left(h_\lambda\right)_{\lambda \in \Par}$,
$\left(e_\lambda\right)_{\lambda \in \Par}$ and
$\left(m_\lambda\right)_{\lambda \in \Par}$; see
\cite{KulikauskasRemmel} for this),
and the Burrows-Wheeler algorithm for data compression
(see Remark~\ref{rmk.lyndon.gr.bw} below for a quick
idea, and \cite{CrochemoreDesarmenienPerrin},
\cite{GesselRestivoReutenauer}, \cite{Kufleitner} for more).
They are also connected to \emph{necklaces} (in the
combinatorial sense) -- a combinatorial object that also happens to be related
to a lot of algebra (\cite[Chapter 5]{Rota}, \cite{Dayton-ghosts}). Let us
survey the basics of this latter classical connection in an exercise:

% [DG][v52] Added reference to \cite{KulikauskasRemmel}, which seems
% very appropriate seeing how much space we have dedicated to
% symmetric functions. (No, I have not read that paper.)

\begin{exercise}
\label{exe.words.necklaces}
Let $\mathfrak{A}$ be any set (not necessarily
totally ordered).
Let \dfn{$C$} denote the infinite cyclic group, written multiplicatively.
Fix a generator \dfn{$c$} of $C$.\ \ \ \ \footnote{So $C$ is a
group isomorphic to $\left( \ZZ ,+\right)  $, and the isomorphism
$\left( \ZZ ,+\right)  \rightarrow C$ sends every $n\in \ZZ$ to
$c^n$. (Recall that we write the binary operation of $C$ as $\cdot$ instead
of $+$.)} Fix a positive integer $n$. The group $C$ acts on $\mathfrak{A}^{n}$
from the left according to the rule
\[
c\cdot\left(  a_{1},a_{2},\ldots,a_{n}\right)  =\left(  a_{2},a_{3}%
,\ldots,a_{n},a_{1}\right)  \ \ \ \ \ \ \ \ \ \ \text{for all }\left(
a_{1},a_{2},\ldots,a_{n}\right)  \in \mathfrak{A}^{n}.
\]
\footnote{In other words, $c$ rotates any $n$-tuple of elements of
$\mathfrak{A}$ cyclically to the left. Thus, $c^{n}\in C$ acts trivially on
$\mathfrak{A}^{n}$, and so this action of $C$ on $\mathfrak{A}^{n}$ factors
through $C/\left\langle c^{n}\right\rangle $ (a cyclic group of order $n$).}
The orbits of this $C$-action will be called
\emph{$n$-necklaces}\index{$n$-necklace}\footnote{Classically, one
visualizes them as necklaces of $n$ beads of
$\left\vert \mathfrak{A} \right\vert$ colors. (The colors are the elements of
$\mathfrak{A}$.) For example, the necklace containing an $n$-tuple
$\left( w_1, w_2, \ldots, w_n \right)$ is visualized as follows:
\[
\xymatrix@R=2.0pc@C=.8pc{
& w_1 \ar@/^1ex/[rr] & & w_2 \ar@/^1ex/[dr] \\
w_n \ar@/^1ex/[ur] & & \;\; & & w_3 \ar@/^1ex/[dl] \\
& w_{n-1} \ar@/^1ex/[ul] & & \iddots \ar@/^1ex/[ll]
}
\]
with $w_1, w_2, \ldots, w_n$ being the colors of the respective beads.
The intuition behind this is that a necklace is an object
that doesn't really change when we rotate it in its plane. However, to make this
intuition match the definition, we need to think of a necklace as being stuck
in its (fixed) plane, so that we cannot lift it up and turn it around,
dropping it back to its plane in a reflected state.}; they form a set
partition of the set $\mathfrak{A}^{n}$.

The $n$-necklace containing a given $n$-tuple $w\in\mathfrak{A}^{n}$ will be
denoted by $\left[  w\right]  $.

\begin{enumerate}
\item[(a)] Prove that
every $n$-necklace $N$ is a finite nonempty set and satisfies
$\left\vert N\right\vert \mid n$. (Recall that $N$ is an orbit, thus a set; as
usual, $\left\vert N\right\vert $ denotes the cardinality of this set.)
\end{enumerate}

The \emph{period}\index{period of an $n$-necklace}
of an $n$-necklace $N$ is defined as the positive integer
$\left\vert N\right\vert $. (This $\left\vert N\right\vert $ is indeed a
positive integer, since $N$ is a finite nonempty set.)\footnote{For example,
the $6$-necklace $\left[ 232232 \right]$ -- or, visually,
\[
\xymatrix@R=2.0pc@C=.8pc{
& 2 \ar@/^1ex/[rr] & & 3 \ar@/^1ex/[dr] \\
2 \ar@/^1ex/[ur] & & \;\; & & 2 \ar@/^1ex/[dl] \\
& 3 \ar@/^1ex/[ul] & & 2 \ar@/^1ex/[ll]
}
\]
-- has period $3$, as it is a set of size $3$
(with elements $232232$, $322322$ and $223223$).
The word ``period'' hints at the geometric meaning: If an $n$-necklace $N$
is represented by coloring the vertices of a regular $n$-gon, then
its period is the smallest positive integer $d$ such that the colors
are preserved when the $n$-gon is rotated by $2\pi d/n$.}

An $n$-necklace is said to be
\emph{aperiodic}\index{aperiodic $n$-necklace} if its period is $n$.

\begin{enumerate}
\item[(b)] Given any $n$-tuple $w = \left(w_1, w_2, \ldots, w_n\right)
\in \mathfrak{A}^n$, prove that the $n$-necklace $\left[w\right]$ is
aperiodic if and only if every $k \in \left\{ 1, 2, \ldots, n-1 \right\}$
satisfies $\left(w_{k+1}, w_{k+2}, \ldots, w_n, w_1, w_2, \ldots, w_k\right)
\neq w$.
\end{enumerate}

From now on, we assume that the set $\mathfrak{A}$ is totally ordered. We use
$\mathfrak{A}$ as our alphabet to define the notions of words, the
lexicographic order, and Lyndon words. All notations that we introduced for
words will thus be used for elements of $\mathfrak{A}^{n}$.

\begin{enumerate}
\item[(c)] Prove that every aperiodic $n$-necklace contains exactly one Lyndon word.

\item[(d)] If $N$ is an $n$-necklace which is not aperiodic, then prove that
$N$ contains no Lyndon word.

\item[(e)] Show that the aperiodic $n$-necklaces are in bijection with Lyndon
words of length $n$.
\end{enumerate}

From now on, we assume that the set $\mathfrak{A}$ is finite. Define the
number-theoretic M\"{o}bius function $\mu$ and the Euler totient function
$\phi$ as in Exercise \ref{exe.witt.ghost-equiv}.

\begin{enumerate}
\item[(f)] Prove that the number of all aperiodic $n$-necklaces is%
\[
\dfrac{1}{n}\sum_{d\mid n}\mu\left(  d\right)  \left\vert \mathfrak{A}%
\right\vert ^{n/d}.
\]

\item[(g)] Prove that the number of all $n$-necklaces is%
\[
\dfrac{1}{n}\sum_{d\mid n}\phi\left(  d\right)  \left\vert \mathfrak{A}%
\right\vert ^{n/d}.
\]

\item[(h)] Solve Exercise~\ref{exe.words.lyndon.count} again.

\item[(i)] Forget that we fixed $\mathfrak{A}$. Show that every $q\in
 \ZZ$ satisfies $n\mid\sum_{d\mid n}\mu\left(  d\right)  q^{n/d}$ and
$n\mid\sum_{d\mid n}\phi\left(  d\right)  q^{n/d}$.
\end{enumerate}

[\textbf{Hint:} For (c), use Theorem~\ref{thm.words.lyndon.equiv}. For (i),
either use parts (f) and (g) and a trick to extend to $q$ negative; or recall
Exercise~\ref{exe.witt.ghost-exa}.]
\end{exercise}

% [DG][v37] Added this exercise (and made some minor corrections to my
% section).

We will pick up the topic of necklaces again in
Section~\ref{subsect.lyndon.gr}, where we will connect it back
to symmetric functions.

\subsubsection{\label{subsect.lyndon.freelie}Free Lie algebras}

In this brief subsection, we shall review the connection between Lyndon words
and free Lie algebras (following \cite[Kap. 4]{Laue}, but avoiding the
generality of Hall sets in favor of just using Lyndon words). None of this
material shall be used in the rest of these notes. We will only prove some
basic results; for more thorough and comprehensive treatments of free Lie
algebras, see \cite{Reutenauer}, \cite[Chapter 2]{Bourbaki-Lie2-3-fr} and
\cite[Kap. 4]{Laue}.

We begin with some properties of Lyndon words.

\begin{exercise}
\label{exe.words.lyndon.stf.props1}Let $w\in\mathfrak{A}^{\ast}$ be a nonempty
word. Let $v$ be the longest Lyndon suffix of $w$\ \ \ \ \footnote{Of course,
a Lyndon suffix of $w$ just means a suffix $p$ of $w$ such that $p$ is
Lyndon.}. Let $t$ be a Lyndon word. Then, $t$ is the longest Lyndon suffix of
$wt$ if and only if we do not have $v<t$.
\end{exercise}

(We have written ``we do not have $v<t$'' instead of ``$v\geq t$'' in Exercise
\ref{exe.words.lyndon.stf.props1} for reasons of generalizability: This way,
Exercise \ref{exe.words.lyndon.stf.props1} generalizes to the partial-order
setting introduced in Exercise \ref{exe.words.poset}, whereas the version with
``$v\geq t$'' does not.)

\begin{exercise}
\label{exe.words.lyndon.stf.props2}Let $w\in\mathfrak{A}^{\ast}$ be a word of
length $>1$. Let $v$ be the longest Lyndon proper suffix of $w$%
\ \ \ \ \footnote{Of course, a Lyndon proper suffix of $w$ just means a proper
suffix $p$ of $w$ such that $p$ is Lyndon.}. Let $t$ be a Lyndon word. Then,
$t$ is the longest Lyndon proper suffix of $wt$ if and only if we do not have
$v<t$.
\end{exercise}

(Exercise \ref{exe.words.lyndon.stf.props2}, while being a trivial consequence
of Exercise \ref{exe.words.lyndon.stf.props1}, is rather useful in the study
of free Lie algebras. It generalizes both \cite[Lemma (1.6)]{ChenFoxLyndon}
(which is obtained by taking $w=c$, $v=b$ and $t=d$) and \cite[Proposition
5.1.4]{Lothaire} (which is obtained by taking $v=m$ and $t=n$).)

\begin{definition}
For the rest of Subsection \ref{subsect.lyndon.freelie}, we let
\dfn{$\mathfrak{L}$} be the
set of all Lyndon words (over the alphabet $\mathfrak{A}$).
\end{definition}

\begin{definition}
\label{def.words.lyndon.stf}Let $w$ be a Lyndon word of length $>1$. Let $v$
be the longest proper suffix of $w$ such that $v$ is Lyndon. (This is
well-defined, as we know from Exercise \ref{exe.words.lyndon.std2}.) Since $v$
is a proper suffix of $w$, there exists a nonempty $u\in\mathfrak{A}^{\ast}$
such that $w=uv$. Consider this $u$. (Clearly, this $u$ is unique.) Theorem
\ref{thm.words.lyndon.std}(a) shows that the words $u$ and $v$ are Lyndon. In
other words, $u\in\mathfrak{L}$ and $v\in\mathfrak{L}$. Hence, $\left(
u,v\right)  \in\mathfrak{L}\times\mathfrak{L}$. The pair $\left(  u,v\right)
\in\mathfrak{L}\times\mathfrak{L}$ is called the \dfn{standard
factorization} of $w$, and is denoted by \dfn{$\operatorname*{stf}w$}.
\end{definition}

For the sake of easier reference, we gather a few basic properties of the
standard factorization:

\begin{exercise}
\label{exe.words.lyndon.stf.basic}Let $w$ be a Lyndon word of length $>1$. Let
$\left(  g,h\right)  =\operatorname*{stf}w$. Prove the following:

\begin{enumerate}
\item[(a)] The word $h$ is the longest Lyndon proper suffix of $w$.

\item[(b)] We have $w=gh$.

\item[(c)] We have $g<gh<h$.

\item[(d)] The word $g$ is Lyndon.

\item[(e)] We have $g\in\mathfrak{L}$, $h\in\mathfrak{L}$, $\ell\left(
g\right)  <\ell\left(  w\right)  $ and $\ell\left(  h\right)  <\ell\left(
w\right)  $.

\item[(f)] Let $t$ be a Lyndon word. Then, $t$ is the longest Lyndon proper
suffix of $wt$ if and only if we do not have $h<t$.
\end{enumerate}
\end{exercise}

\begin{exercise}
\label{exe.freelie.1}Let $\mathfrak{g}$ be a Lie algebra. For every Lyndon
word $w$, let $b_{w}$ be an element of $\mathfrak{g}$. Assume that for every
Lyndon word $w$ of length $>1$, we have%
\begin{equation}
b_{w}=\left[  b_{u},b_{v}\right]  ,\ \ \ \ \ \ \ \ \ \ \text{where }\left(
u,v\right)  =\operatorname*{stf}w. \label{eq.exe.freelie.1.bw}%
\end{equation}

Let $B$ be the $\kk$-submodule of $\mathfrak{g}$ spanned by the family
$\left(  b_{w}\right)  _{w\in\mathfrak{L}}$.

\begin{enumerate}
\item[(a)] Prove that $B$ is a Lie subalgebra of $\mathfrak{g}$.

\item[(b)] Let $\mathfrak{h}$ be a $\kk$-Lie algebra. Let
$f:B\rightarrow\mathfrak{h}$ be a $\kk$-module homomorphism. Assume
that whenever $w$ is a Lyndon word of length $>1$, we have%
\begin{equation}
f\left(  \left[  b_{u},b_{v}\right]  \right)  =\left[  f\left(  b_{u}\right)
,f\left(  b_{v}\right)  \right]  ,\ \ \ \ \ \ \ \ \ \ \text{where }\left(
u,v\right)  =\operatorname*{stf}w. \label{eq.exe.freelie.1.bcond}%
\end{equation}
Prove that $f$ is a Lie algebra homomorphism.
\end{enumerate}

[\textbf{Hint:} Given two words $w$ and $w^{\prime}$, write $w\sim w^{\prime}$
if and only if $w^{\prime}$ is a permutation of $w$. Part (a) follows from the
fact that for any $\left(  p,q\right)  \in\mathfrak{L}\times\mathfrak{L}$
satisfying $p<q$, we have $\left[  b_{p},b_{q}\right]  \in B_{pq,q}$, where
$B_{h,s}$ denotes the $\kk$-linear span of $\left\{  b_{w}\ \mid
\ w\in\mathfrak{L}\text{, }w\sim h\text{ and }w<s\right\}  $ for any two words
$h$ and $s$. Prove this fact by a double induction, first inducting over
$\ell\left(  pq\right)  $, and then (for fixed $\ell\left(  pq\right)  $)
inducting over the rank of $q$ in lexicographic order (i.e., assume that the
fact is already proven for every $q^{\prime}<q$ instead of $q$). In the
induction step, assume that $\left(  p,q\right)  \neq\operatorname*{stf}%
\left(  pq\right)  $ (since otherwise the claim is rather obvious) and
conclude that $p$ has length $>1$; thus, set $\left(  u,v\right)
=\operatorname*{stf}p$, so that $\left[  \underbrace{b_{p}}_{=\left[
b_{u},b_{v}\right]  },b_{q}\right]  =\left[  \left[  b_{u},b_{v}\right]
,b_{q}\right]  =\left[  \left[  b_{u},b_{q}\right]  ,b_{v}\right]  -\left[
\left[  b_{v},b_{q}\right]  ,b_{u}\right]  $, and use Exercise
\ref{exe.words.lyndon.stf.props2} to obtain $v<q$.

The proof of (b) proceeds by a similar induction, piggybacking on the $\left[
b_{p},b_{q}\right]  \in B_{pq,q}$ claim.]
\end{exercise}

\begin{exercise}
\label{exe.freelie.2}Let $V$ be the free $\kk$-module with basis
$\left(  x_{a}\right)  _{a\in\mathfrak{A}}$. For every word $w\in
\mathfrak{A}^{\ast}$, let $x_{w}$ be the tensor $x_{w_{1}}\otimes x_{w_{2}%
}\otimes\cdots\otimes x_{w_{\ell\left(  w\right)  }}$. As we know from Example
\ref{exa.tensor-alg}, the tensor algebra $T\left(  V\right)  $ is a free
$\kk$-module with basis $\left(  x_{w}\right)  _{w\in\mathfrak{A}%
^{\ast}}$. We regard $V$ as a $\kk$-submodule of $T\left(  V\right)  $.

The tensor algebra $T\left(  V\right)  $ becomes a Lie algebra via the
commutator (i.e., its Lie bracket is defined by $\left[  \alpha,\beta\right]
=\alpha\beta-\beta\alpha$ for all $\alpha\in T\left(  V\right)  $ and
$\beta\in T\left(  V\right)  $).

We define a sequence $\left(  \mathfrak{g}_{1},\mathfrak{g}_{2},\mathfrak{g}%
_{3},\ldots\right)  $ of $\kk$-submodules of $T\left(  V\right)  $ as
follows: Recursively, we set $\mathfrak{g}_{1}=V$, and for every $i\in\left\{
2,3,4,\ldots\right\}  $, we set $\mathfrak{g}_{i}=\left[  V,\mathfrak{g}%
_{i-1}\right]  $. Let $\mathfrak{g}$ be the $\kk$-submodule
$\mathfrak{g}_{1}+\mathfrak{g}_{2}+\mathfrak{g}_{3}+\cdots$ of $T\left(
V\right)  $.

Prove the following:

\begin{enumerate}
\item[(a)] The $\kk$-submodule $\mathfrak{g}$ is a Lie subalgebra of
$T\left(  V\right)  $.

\item[(b)] If $\mathfrak{k}$ is any Lie subalgebra of $T\left(  V\right)  $
satisfying $V\subset\mathfrak{k}$, then $\mathfrak{g}\subset\mathfrak{k}$.
\end{enumerate}

Now, for every $w\in\mathfrak{L}$, we define an element $b_{w}$ of $T\left(
V\right)  $ as follows: We define $b_{w}$ by recursion on the length of $w$.
If the length of $w$ is $1$\ \ \ \ \footnote{The length of any $w\in
\mathfrak{L}$ must be at least $1$. (Indeed, if $w\in\mathfrak{L}$, then the
word $w$ is Lyndon and thus nonempty, and hence its length must be at least
$1$.)}, then we have $w=\left(  a\right)  $ for some letter $a\in\mathfrak{A}%
$, and we set $b_{w}=x_{a}$ for this letter $a$. If the length of $w$ is $>1$,
then we set $b_{w}=\left[  b_{u},b_{v}\right]  $, where $\left(  u,v\right)
=\operatorname*{stf}w$\ \ \ \ \footnote{This is well-defined, because $b_{u}$
and $b_{v}$ have already been defined. [\textit{Proof.} Let $\left(
u,v\right)  =\operatorname*{stf}w$. Then, Exercise
\ref{exe.words.lyndon.stf.basic}(e) (applied to $\left(  g,h\right)  =\left(
u,v\right)  $) shows that $u\in\mathfrak{L}$, $v\in\mathfrak{L}$, $\ell\left(
u\right)  <\ell\left(  w\right)  $ and $\ell\left(  v\right)  <\ell\left(
w\right)  $. Recall that we are defining $b_{w}$ by recursion on the length of
$w$. Hence, $b_{p}$ is already defined for every $p\in\mathfrak{L}$ satisfying
$\ell\left(  p\right)  <\ell\left(  w\right)  $. Applying this to $p=u$, we
see that $b_{u}$ is already defined (since $u\in\mathfrak{L}$ and $\ell\left(
u\right)  <\ell\left(  w\right)  $). The same argument (but applied to $v$
instead of $u$) shows that $b_{v}$ is already defined. Hence, $b_{u}$ and
$b_{v}$ have already been defined. Thus, $b_{w}$ is well-defined by
$b_{w}=\left[  b_{u},b_{v}\right]  $, qed.]}.

Prove the following:

\begin{enumerate}
\item[(c)] For every $w\in\mathfrak{L}$, we have%
\[
b_{w}\in x_{w}+\sum\limits_{\substack{v\in\mathfrak{A}^{\ell\left(  w\right)
};\\v>w}}\kk x_{v}.
\]

\item[(d)] The family $\left(  b_{w}\right)  _{w\in\mathfrak{L}}$ is a basis
of the $\kk$-module $\mathfrak{g}$.

\item[(e)] Let $\mathfrak{h}$ be any $\kk$-Lie algebra. Let
$\xi:\mathfrak{A}\rightarrow\mathfrak{h}$ be any map. Then, there exists a
unique Lie algebra homomorphism $\Xi:\mathfrak{g}\rightarrow\mathfrak{h}$ such
that every $a\in\mathfrak{A}$ satisfies $\Xi\left(  x_{a}\right)  =\xi\left(
a\right)  $.
\end{enumerate}
\end{exercise}

\begin{remark}
\label{rmk.freelie.2.rmk}Let $V$ and $\mathfrak{g}$ be as in Exercise
\ref{exe.freelie.2}. In the language of universal algebra, the statement of
Exercise \ref{exe.freelie.2}(e) says that $\mathfrak{g}$ (or, to be more
precise, the pair $\left(  \mathfrak{g},f\right)  $, where $f:\mathfrak{A}%
\rightarrow\mathfrak{g}$ is the map sending each $a\in\mathfrak{A}$ to
$x_{a}\in\mathfrak{g}$) satisfies the universal property of the free Lie
algebra on the set $\mathfrak{A}$. Thus, this exercise allows us to call
$\mathfrak{g}$ the \dfn{free Lie algebra} on $\mathfrak{A}$. Most authors
define the free Lie algebra differently, but all reasonable definitions of a
free Lie algebra\footnote{Here, we call a definition ``reasonable'' if the
``free Lie algebra'' it defines satisfies the universal property.} lead
to isomorphic Lie algebras (because the universal property determines the free
Lie algebra uniquely up to canonical isomorphism).

Notice that the Lie algebra $\mathfrak{g}$ does not depend on the total order
on the alphabet $\mathfrak{A}$, but the basis $\left(  b_{w}\right)
_{w\in\mathfrak{L}}$ constructed in Exercise \ref{exe.freelie.2}(d) does.
There is no known basis of $\mathfrak{g}$ defined without ordering
$\mathfrak{A}$.

It is worth noticing that our construction of $\mathfrak{g}$ proves not only
that the free Lie algebra on $\mathfrak{A}$ exists, but also that this free
Lie algebra can be realized as a Lie subalgebra of the (associative) algebra
$T\left(  V\right)  $. Therefore, if we want to prove that a certain identity
holds in every Lie algebra, we only need to check that this identity holds in
every associative algebra (if all Lie brackets are replaced by commutators);
the universal property of the free Lie algebra (i.e., Exercise
\ref{exe.freelie.2}(e)) will then ensure that this identity also holds in
every Lie algebra $\mathfrak{h}$.

There is much more to say about free Lie algebras than what we have said
here; in particular, there are connections to symmetric functions,
necklaces, representations of symmetric groups and $\Nsym$. See
\cite[\S 5.3]{Lothaire}, \cite{Reutenauer},
\cite[Chapter 2]{Bourbaki-Lie2-3-fr}, \cite[\S 4]{Laue} and
\cite{BlessenohlLaue} for further developments\footnote{The claim made
in \cite[page 2]{BlessenohlLaue} that ``$\left\{x_1, \ldots, x_n\right\}$
generates freely a Lie subalgebra of $A_R$'' is essentially our
Exercise~\ref{exe.freelie.2}(e).}.
\end{remark}

% [DG][v54] Added this subsection on free Lie algebras.

\subsection{\label{subsect.lyndon.shuffles}Shuffles and Lyndon words}

We will now connect the theory of Lyndon words with the notion of shuffle
products. We have already introduced the latter notion in Definition
\ref{shuffles}, but we will now study it more closely and introduce some more
convenient notations (e.g., we will need a notation for single shuffles, not
just the whole multiset).\footnote{Parts (a) and (c) of the below Definition
\ref{def.shuffle.sigma} define notions which have already been introduced in
Definition \ref{shuffles}. Of course, the definitions of these notions are
equivalent; however, the variables are differently labelled in the two
definitions (for example, the variables $u$, $v$, $w$ and $\sigma$ of
Definition \ref{def.shuffle.sigma}(c) correspond to the variables $a$, $b$,
$c$ and $w$ of Definition \ref{shuffles}). The labels in Definition
\ref{def.shuffle.sigma} have been chosen to match with the rest of Section
\ref{subsect.lyndon.shuffles}.}

\begin{definition}
\phantomsection
\label{def.shuffle.sigma}

\begin{itemize}
\item[(a)] Let $n\in \NN$ and $m\in \NN$. Then,
\dfn{$\operatorname{Sh}_{n,m}$} denotes the subset
\[
\left\{  \sigma\in\Symm_{n+m}\ :\ \sigma^{-1}\left(  1\right)
<\sigma^{-1}\left(  2\right)  <\cdots<\sigma^{-1}\left(  n\right)
;\ \sigma^{-1}\left(  n+1\right)  <\sigma^{-1}\left(  n+2\right)
<\cdots<\sigma^{-1}\left(  n+m\right)  \right\}
\]
of the symmetric group $\Symm_{n+m}$.

\item[(b)] Let $u=\left(  u_{1},u_{2},\ldots,u_{n}\right)  $ and $v=\left(
v_{1},v_{2},\ldots,v_{m}\right)  $ be two words. If $\sigma\in
\operatorname{Sh}_{n,m}$, then, $u\underset{\sigma}{\shuffle}v$
will denote
the word $\left(  w_{\sigma\left(  1\right)  },w_{\sigma\left(  2\right)
},\ldots,w_{\sigma\left(  n+m\right)  }\right)  $, where $\left(  w_{1}%
,w_{2},\ldots,w_{n+m}\right)  $ is the concatenation $u\cdot v=\left(
u_{1},u_{2},\ldots,u_{n},v_{1},v_{2},\ldots,v_{m}\right)  $. We notice that
the multiset of all letters of $u\underset{\sigma}{\shuffle}v$ is the disjoint
union of the multiset of all letters of $u$ with the multiset of all letters
of $v$. As a consequence, $\ell\left(  u\underset{\sigma}{\shuffle}v\right)
=\ell\left(  u\right)  +\ell\left(  v\right)  $.

\item[(c)] Let $u=\left(  u_{1},u_{2},\ldots,u_{n}\right)  $ and $v=\left(
v_{1},v_{2},\ldots,v_{m}\right)  $ be two words. The \dfn{multiset of
shuffles of $u$ and $v$}\index{shuffle of words}
is defined as the multiset $\left\{  \left(
w_{\sigma\left(  1\right)  },w_{\sigma\left(  2\right)  },\ldots
,w_{\sigma\left(  n+m\right)  }\right)  \ :\ \sigma\in
\operatorname{Sh}_{n,m}\right\}  _{\text{multiset}}$, where $\left(  w_{1}%
,w_{2},\ldots,w_{n+m}\right)  $ is the concatenation $u\cdot v=\left(
u_{1},u_{2},\ldots,u_{n},v_{1},v_{2},\ldots,v_{m}\right)  $. In other words,
the multiset of shuffles of $u$ and $v$ is the multiset
\[
\left\{  u\underset{\sigma}{\shuffle}v\ :\ \sigma\in\operatorname{Sh}%
_{n,m}\right\}  _{\text{multiset}}.
\]
It is denoted by \dfn{$u\shuffle v$}.
\end{itemize}
\end{definition}

The next fact provides the main connection between Lyndon words and shuffles:

\begin{theorem}
\label{thm.shuffle.lyndon.used}Let $u$ and $v$ be two words.

Let $\left(  a_{1},a_{2},\ldots,a_{p}\right)  $ be the CFL factorization of
$u$. Let $\left(  b_{1},b_{2},\ldots,b_{q}\right)  $ be the CFL factorization
of $v$.

\begin{itemize}
\item[(a)] Let $\left(  c_{1},c_{2},\ldots,c_{p+q}\right)  $ be the result of
sorting the list $\left(  a_{1},a_{2},\ldots,a_{p},b_{1},b_{2},\ldots
,b_{q}\right)  $ in decreasing order\footnote{with respect to the total order
on $\mathfrak{A}^{\ast}$ whose greater-or-equal relation is $\geq$}. Then, the
lexicographically highest element of the multiset $u\shuffle v$ is $c_{1}c_{2}
\cdots c_{p+q}$ (and $\left(  c_{1},c_{2},\ldots,c_{p+q}\right)  $ is the CFL
factorization of this element).

\item[(b)] Let $\mathfrak{L}$ denote the set of all Lyndon words. If $w$ is a
Lyndon word and $z$ is any word, let $\operatorname{mult}_{w}z$
denote the number of terms in the CFL factorization of $z$ which are equal to
$w$. The multiplicity with which the lexicographically highest element of the
multiset $u\shuffle v$ appears in the multiset $u\shuffle v$ is
$\prod_{w\in\mathfrak{L}}
  \dbinom{\operatorname{mult}_{w} u + \operatorname{mult}_{w}v}
         {\operatorname{mult}_{w}u}$.
(This product is well-defined because almost all of its factors are $1$.)

\item[(c)] If $a_{i}\geq b_{j}$ for every $i\in\left\{  1,2,\ldots,p\right\}
$ and $j\in\left\{  1,2,\ldots,q\right\}  $, then the lexicographically
highest element of the multiset $u\shuffle v$ is $uv$.

\item[(d)] If $a_{i} > b_{j}$ for every $i\in\left\{  1,2,\ldots,p\right\}  $
and $j\in\left\{  1,2,\ldots,q\right\}  $, then the multiplicity with which
the word $uv$ appears in the multiset $u\shuffle v$ is $1$.

\item[(e)] Assume that $u$ is a Lyndon word. Also, assume that $u\geq b_{j}$
for every $j\in\left\{  1,2,\ldots,q\right\}  $. Then, the lexicographically
highest element of the multiset $u\shuffle v$ is $uv$, and the multiplicity
with which this word $uv$ appears in the multiset $u\shuffle v$ is
$\operatorname{mult}_{u}v+1$.
\end{itemize}
\end{theorem}

\begin{comment}
Here is an old (equivalent, but messy) version of the statement of
Theorem~\ref{thm.shuffle.lyndon.used}(a):
(a) Let $C$ be the disjoint union of the multisets $\left\{  a_{1}%
,a_{2},\ldots,a_{p}\right\}  _{\text{multiset}}$ and $\left\{  b_{1}%
,b_{2},\ldots,b_{q}\right\}  _{\text{multiset}}$ (where ``disjoint union''
means that the multiplicities of each
element are added). Let $c_{1},c_{2},\ldots,c_{p+q}$ be the elements of $C$
listed in decreasing order (as often as they appear in $C$). Then, the
lexicographically highest element of the multiset $u\shuffle v$ is $c_{1}%
c_{2} \cdots c_{p+q}$ (and $\left(  c_{1},c_{2},\ldots,c_{p+q}\right)  $ is the CFL
factorization of this element).
\end{comment}

\begin{example}
For this example, let $u$ and $v$ be the words $u=23232$ and $v=323221$ over
the alphabet $\mathfrak{A}=\left\{  1,2,3, \ldots\right\}  $ with total order
given by $1<2<3< \cdots$. The CFL factorizations of $u$ and $v$ are $\left(
23,23,2\right)  $ and $\left(  3,23,2,2,1\right)  $, respectively. Thus, using
the notations of Theorem \ref{thm.shuffle.lyndon.used}, we have $p=3$,
$\left(  a_{1},a_{2},\ldots,a_{p}\right)  =\left(  23,23,2\right)  $, $q=5$
and $\left(  b_{1},b_{2},\ldots,b_{q}\right)  =\left(  3,23,2,2,1\right)  $.
Thus, Theorem \ref{thm.shuffle.lyndon.used}(a) predicts that the
lexicographically highest element of the multiset $u\shuffle v$ is $c_{1}%
c_{2}c_{3}c_{4}c_{5}c_{6}c_{7}c_{8}$, where $c_{1},c_{2},c_{3},c_{4}%
,c_{5},c_{6},c_{7},c_{8}$ are the words $23,23,2,3,23,2,2,1$ listed in
decreasing order (in other words, $\left(  c_{1},c_{2},c_{3},c_{4},c_{5}%
,c_{6},c_{7},c_{8}\right)  =\left(  3,23,23,23,2,2,2,1\right)  $). In other
words, Theorem \ref{thm.shuffle.lyndon.used}(a) predicts that the
lexicographically highest element of the multiset $u\shuffle v$ is
$32323232221$. We could verify this by brute force, but this would be
laborious since the multiset $u\shuffle v$ has $\dbinom{5+6}{5}=462$ elements
(with multiplicities). Theorem \ref{thm.shuffle.lyndon.used}(b) predicts that
this lexicographically highest element $32323232221$ appears in the multiset
$u\shuffle v$ with a multiplicity of $\prod_{w\in\mathfrak{L}}\dbinom
{\operatorname{mult}_{w}u+\operatorname{mult}_{w}%
v}{\operatorname{mult}_{w}u}$. This product $\prod_{w\in
\mathfrak{L}}\dbinom{\operatorname{mult}_{w}u+
\operatorname{mult}_{w}v}{\operatorname{mult}_{w}u}$ is infinite, but all but
finitely many of its factors are $1$ and therefore can be omitted; the only
factors which are not $1$ are those corresponding to Lyndon words $w$ which
appear both in the CFL factorization of $u$ and in the CFL factorization of
$v$ (since for any other factor, at least one of the numbers
$\operatorname{mult}_{w}u$ or $\operatorname{mult}_{w}v$
equals $0$, and therefore the binomial coefficient $\dbinom
{\operatorname{mult}_{w}u+\operatorname{mult}_{w}%
v}{\operatorname{mult}_{w}u}$ equals $1$). Thus, in order to compute
the product $\prod_{w\in\mathfrak{L}}\dbinom{
\operatorname{mult}_{w}u+\operatorname{mult}_{w}v}{
\operatorname{mult}_{w}u}$, we only need to multiply these factors. In our example,
these are the factors for $w=23$ and for $w=2$ (these are the only Lyndon
words which appear both in the CFL factorization $\left(  23,23,2\right)  $ of
$u$ and in the CFL factorization $\left(  3,23,2,2,1\right)  $ of $v$). So we
have%
\[
\prod_{w\in\mathfrak{L}}\dbinom{\operatorname{mult}_{w}%
u+\operatorname{mult}_{w}v}{\operatorname{mult}_{w}%
u}=\underbrace{\dbinom{\operatorname{mult}_{23}%
u+\operatorname{mult}_{23}v}{\operatorname{mult}_{23}u}%
}_{=\dbinom{2+1}{2}=3}\underbrace{\dbinom{\operatorname{mult}_{2}%
u+\operatorname{mult}_{2}v}{\operatorname{mult}_{2}u}%
}_{=\dbinom{1+2}{1}=3}=3\cdot3=9.
\]
The word $32323232221$ must thus appear in the multiset $u\shuffle v$ with a
multiplicity of $9$. This, too, could be checked by brute force.
\end{example}

Theorem \ref{thm.shuffle.lyndon.used} (and Theorem
\ref{thm.shuffle.lyndon.full} further below, which describes more precisely
how the lexicographically highest element of $u\shuffle v$ emerges by
shuffling $u$ and $v$) is fairly close to \cite[Theorem 2.2.2]%
{Radford-shuffle} (and will be used for the same purposes), the main
difference being that we are talking about the shuffle product of two (not
necessarily Lyndon) words, while Radford (and most other authors) study the
shuffle product of many Lyndon words.

In order to prove Theorem \ref{thm.shuffle.lyndon.used}, we will need to make
some stronger statements, for which we first have to introduce some more notation:

\begin{definition}
\begin{itemize}

\item[(a)] If $p$ and $q$ are two integers, then $\left[  p:q\right]  ^{+}$
denotes the interval $\left\{  p+1,p+2,\ldots,q\right\}  $ of $\ZZ$.
Note that $\left\vert \left[  p:q\right]  ^{+}\right\vert =q-p$ if $q\geq p$.

\item[(b)] If $I$ and $J$ are two nonempty intervals of $\ZZ$, then we
say that $I<J$ if and only if every $i\in I$ and $j\in J$ satisfy $i<j$. This
defines a partial order on the set of nonempty intervals of $\ZZ$.
(Roughly speaking, $I<J$ if the interval $I$ ends before $J$ begins.)

\item[(c)] If $w$ is a word with $n$ letters (for some $n\in \NN$), and
$I$ is an interval of $\ZZ$ such that $I\subset\left[  0:n\right]
^{+}$, then $w\left[  I\right]  $ will denote the word $\left(  w_{p+1}%
,w_{p+2},\ldots,w_{q}\right)  $, where $I$ is written in the form $I=\left[
p:q\right]  ^{+}$ with $q \geq p$. Obviously, $\ell\left(  w\left[  I\right]
\right)  =\left\vert I\right\vert =q-p$. A word of the form $w\left[
I\right]  $ for an interval $I\subset\left[  0:n\right]  ^{+}$ (equivalently,
a word which is a prefix of a suffix of $w$) is called a
\emph{factor}\index{factor of a word} of $w$.

\item[(d)] Let $\alpha$ be a composition. Then, we define a tuple
\dfn{$\operatorname*{intsys}\alpha$}
of intervals of $\ZZ$ as follows: Write
$\alpha$ in the form $\left(  \alpha_{1},\alpha_{2},\ldots,\alpha_{\ell
}\right)  $ (so that $\ell= \ell\left(  \alpha\right)  $). Then, set
$\operatorname*{intsys}\alpha=\left(  I_{1},I_{2},\ldots,I_{\ell}\right)  $,
where
\[
I_{i}=\left[  \sum_{k=1}^{i-1}\alpha_{k}:\sum_{k=1}^{i}\alpha_{k}\right]
^{+}\ \ \ \ \ \ \ \ \ \ \text{for every }i\in\left\{  1,2,\ldots,\ell\right\}
.
\]
This $\ell$-tuple $\operatorname*{intsys}\alpha$ is a tuple of nonempty
intervals of $\ZZ$. This tuple $\operatorname{intsys}\alpha$ is called
the \dfn{interval system corresponding to $\alpha$}. (This is precisely the
$\ell$-tuple $\left(  I_{1},I_{2},\ldots,I_{\ell}\right)  $ constructed in
Definition \ref{def.compositions.intervaldec}.) The length of the tuple
$\operatorname*{intsys}\alpha$ is $\ell\left(  \alpha\right)  $.
\end{itemize}
\end{definition}

\begin{example}
\begin{itemize}

\item[(a)] We have $\left[  2:4\right]  ^{+}=\left\{  3,4\right\}  $ and
$\left[  3:3\right]  ^{+}=\varnothing$.

\item[(b)] We have $\left[  2:4\right]  ^{+}<\left[  4:5\right]  ^{+}<\left[
6:8\right]  ^{+}$, but we have neither $\left[  2:4\right]  ^{+}<\left[
3:5\right]  ^{+}$ nor $\left[  3:5\right]  ^{+}<\left[  2:4\right]  ^{+}$.

\item[(c)] If $w$ is the word $915352$, then $w\left[  \left[  0:3\right]
^{+}\right]  =\left(  w_{1},w_{2},w_{3}\right)  =915$ and $w\left[  \left[
2:4\right]  ^{+}\right]  =\left(  w_{3},w_{4}\right)  =53$.

\item[(d)] If $\alpha$ is the composition $\left(  4,1,4,2,3\right)  $, then
the interval system corresponding to $\alpha$ is
\begin{align*}
\operatorname*{intsys}\alpha &  =\left(  \left[  0:4\right]  ^{+},\left[
4:5\right]  ^{+},\left[  5:9\right]  ^{+},\left[  9:11\right]  ^{+},\left[
11:14\right]  ^{+}\right) \\
&  =\left(  \left\{  1,2,3,4\right\}  ,\left\{  5\right\}  ,\left\{
6,7,8,9\right\}  ,\left\{  10,11\right\}  ,\left\{  12,13,14\right\}  \right)
.
\end{align*}

\end{itemize}
\end{example}

The following properties of the notions introduced in the preceding definition
are easy to check:

\begin{remark}
\phantomsection
\label{rmk.intsys.trivia}

\begin{itemize}
\item[(a)] If $I$ and $J$ are two nonempty intervals of $\ZZ$
satisfying $I < J$, then $I$ and $J$ are disjoint.

\item[(b)] If $I$ and $J$ are two disjoint nonempty intervals of
$\ZZ$, then either $I < J$ or $J < I$.

\item[(c)] Let $\alpha$ be a composition. Write $\alpha$ in the form $\left(
\alpha_{1},\alpha_{2},\ldots,\alpha_{\ell}\right)  $ (so that $\ell=
\ell\left(  \alpha\right)  $). The interval system $\operatorname{intsys}%
\alpha$ can be described as the unique $\ell$-tuple $\left(  I_{1}%
,I_{2},\ldots,I_{\ell}\right)  $ of nonempty intervals of $\ZZ$
satisfying the following three properties:

\begin{itemize}
\item The intervals $I_{1}$, $I_{2}$, $\ldots$, $I_{\ell}$ form a set
partition of the set $\left[  0:n\right]  ^{+}$, where $n=\left\vert
\alpha\right\vert $.

\item We have $I_{1}<I_{2}< \cdots<I_{\ell}$.

\item We have $\left\vert I_{i}\right\vert =\alpha_{i}$ for every
$i\in\left\{  1,2,\ldots,\ell\right\}  $.
\end{itemize}
\end{itemize}
\end{remark}

\begin{exercise}
\label{exe.rmk.intsys.trivia}Prove Remark \ref{rmk.intsys.trivia}.
\end{exercise}

The following two lemmas are collections of more or less trivial consequences
of what it means to be an element of $\operatorname{Sh}_{n,m}$ and
what it means to be a shuffle:

\begin{lemma}
\label{lem.shuffle.lyndon.easystep}Let $n\in \NN$ and $m\in \NN$.
Let $\sigma\in\operatorname{Sh}_{n,m}$.

\begin{itemize}
\item[(a)] If $I$ is an interval of $\ZZ$ such that $I\subset\left[
0:n+m\right]  ^{+}$, then $\sigma\left(  I\right)  \cap\left[  0:n\right]
^{+}$ and $\sigma\left(  I\right)  \cap\left[  n:n+m\right]  ^{+}$ are intervals.

\item[(b)] Let $K$ and $L$ be nonempty intervals of $\ZZ$ such that
$K\subset\left[  0:n\right]  ^{+}$ and $L\subset\left[  0:n\right]  ^{+}$ and
$K<L$ and such that $K\cup L$ is an interval. Assume that $\sigma^{-1}\left(
K\right)  $ and $\sigma^{-1}\left(  L\right)  $ are intervals, but
$\sigma^{-1}\left(  K\right)  \cup\sigma^{-1}\left(  L\right)  $ is not an
interval. Then, there exists a nonempty interval $P\subset\left[
n:n+m\right]  ^{+}$ such that $\sigma^{-1}\left(  P\right)  $, $\sigma
^{-1}\left(  K\right)  \cup\sigma^{-1}\left(  P\right)  $ and $\sigma
^{-1}\left(  P\right)  \cup\sigma^{-1}\left(  L\right)  $ are intervals and
such that $\sigma^{-1}\left(  K\right)  <\sigma^{-1}\left(  P\right)
<\sigma^{-1}\left(  L\right)  $.

\item[(c)] Lemma \ref{lem.shuffle.lyndon.easystep}(b) remains valid if
``$K\subset\left[  0:n\right]  ^{+}$ and
$L\subset\left[  0:n\right]  ^{+}$''
and ``$P\subset\left[ n:n+m\right]  ^{+}$''
are replaced by ``$K\subset\left[  n:n+m\right]  ^{+}$ and
$L\subset\left[  n:n+m\right]  ^{+}$'' and
``$P\subset\left[  0:n\right]  ^{+}$'', respectively.
\end{itemize}
\end{lemma}

\begin{exercise}
\label{exe.lem.shuffle.lyndon.easystep}Prove Lemma
\ref{lem.shuffle.lyndon.easystep}.
\end{exercise}

\begin{lemma}
\label{lem.shuffle.lyndon.step}Let $u$ and $v$ be two words. Let
$n=\ell\left(  u\right)  $ and $m=\ell\left(  v\right)  $. Let $\sigma
\in\operatorname{Sh}_{n,m}$.

\begin{itemize}
\item[(a)] If $I$ is an interval of $\ZZ$ satisfying either
$I\subset\left[  0:n\right]  ^{+}$ or $I\subset\left[  n:n+m\right]  ^{+}$,
and if $\sigma^{-1}\left(  I\right)  $ is an interval, then%
\begin{equation}
\left(  u\underset{\sigma}{\shuffle}v\right)  \left[  \sigma^{-1}\left(
I\right)  \right]  =\left(  uv\right)  \left[  I\right]  .
\label{pf.lem.shuffle.lyndon.step.gamma}%
\end{equation}

\item[(b)] Assume that $u\underset{\sigma}{\shuffle}v$ is the
lexicographically highest element of the multiset $u\shuffle v$. Let
$I\subset\left[  0:n\right]  ^{+}$ and $J\subset\left[  n:n+m\right]  ^{+}$ be
two nonempty intervals. Assume that $\sigma^{-1}\left(  I\right)  $ and
$\sigma^{-1}\left(  J\right)  $ are also intervals, that $\sigma^{-1}\left(
I\right)  <\sigma^{-1}\left(  J\right)  $, and that $\sigma^{-1}\left(
I\right)  \cup\sigma^{-1}\left(  J\right)  $ is an interval as well. Then,
$\left(  uv\right)  \left[  I\right]  \cdot\left(  uv\right)  \left[
J\right]  \geq\left(  uv\right)  \left[  J\right]  \cdot\left(  uv\right)
\left[  I\right]  $.

\item[(c)] Lemma \ref{lem.shuffle.lyndon.step}(b) remains valid if
``$I\subset\left[  0:n\right]  ^{+}$ and $J\subset\left[
n:n+m\right]  ^{+}$'' is replaced by
``$I\subset\left[  n:n+m\right]  ^{+}$ and
$J\subset\left[  0:n\right]  ^{+}$''.
\end{itemize}
\end{lemma}

\begin{exercise}
\label{exe.lem.shuffle.lyndon.step}Prove Lemma \ref{lem.shuffle.lyndon.step}.

[\textbf{Hint:} For (b), show that there exists a $\tau\in
\operatorname{Sh}_{n,m}$ such that $u\underset{\tau}{\shuffle}v$ differs from
$u\underset{\sigma}{\shuffle}v$ only in the order of the subwords $\left(
uv\right)  \left[  I\right]  $ and $\left(  uv\right)  \left[  J\right]  $.]
\end{exercise}

We are still a few steps away from stating our results in a way that allows
comfortably proving Theorem \ref{thm.shuffle.lyndon.used}. For the latter aim,
we introduce the notion of \emph{$\alpha$-clumping permutations}, and
characterize them in two ways:

\begin{definition}
\label{def.intervalperm} Let $n\in{ \NN }$. Let $\alpha$ be a composition
of $n$. Let $\ell=\ell\left(  \alpha\right)  $.

\begin{itemize}
\item[(a)] For every set $S$ of positive integers, let $\overrightarrow{S}$
denote the list of all elements of $S$ in increasing order (with each element
appearing exactly once). Notice that this list $\overrightarrow{S}$ is a word
over the set of positive integers.

\item[(b)] For every $\tau\in\Symm_{\ell}$, we define a permutation
$\operatorname{iper}\left(  \alpha,\tau\right)  \in\Symm_n$ as follows:

The interval system corresponding to $\alpha$ is an $\ell$-tuple of intervals
(since $\ell\left(  \alpha\right)  =\ell$); denote this $\ell$-tuple by
$\left(  I_{1},I_{2},\ldots,I_{\ell}\right)  $. Now, define
$\operatorname{iper}\left(  \alpha,\tau\right)  $ to be the permutation in
$\Symm_n$ which (in one-line notation) is the word
$\overrightarrow{I_{\tau\left(  1\right)  }}\overrightarrow{I_{\tau\left(
2\right)  }}\cdots\overrightarrow{I_{\tau\left(  \ell\right)  }}$ (a
concatenation of $\ell$ words). This is well-defined\footnote{In fact, from
the properties of interval systems, we know that the intervals $I_{1}$,
$I_{2}$, $\ldots$, $I_{\ell}$ form a set partition of the set $\left[
0:n\right]  ^{+}$. Hence, the intervals $I_{\tau\left(  1\right)  }$,
$I_{\tau\left(  2\right)  }$, $\ldots$, $I_{\tau\left(  \ell\right)  }$ form a
set partition of the set $\left[  0:n\right]  ^{+}$. As a consequence, the
word $\overrightarrow{I_{\tau\left(  1\right)  }}\overrightarrow{I_{\tau
\left(  2\right)  }}\cdots\overrightarrow{I_{\tau\left(  \ell\right)  }}$ is a
permutation of the word $12\ldots n$, and so there exists a permutation in
$\Symm_n$ which (in one-line notation) is this word, qed.}; hence,
$\operatorname{iper}\left(  \alpha,\tau\right)  \in\Symm_n$ is defined.

\item[(c)] The interval system corresponding to $\alpha$ is an $\ell$-tuple of
intervals (since $\ell\left(  \alpha\right)  =\ell$); denote this $\ell$-tuple
by $\left(  I_{1},I_{2},\ldots,I_{\ell}\right)  $.

A permutation $\sigma\in\Symm_n$ is said to be
\emph{$\alpha$-clumping}\index{$\alpha$-clumping permutation}\index{clumping}
if every $i\in\left\{  1,2,\ldots,\ell\right\}  $ has the two
properties that:

\begin{itemize}
\item the set $\sigma^{-1}\left(  I_{i}\right)  $ is an interval;

\item the restriction of the map $\sigma^{-1}$ to the interval $I_{i}$ is increasing.
\end{itemize}
\end{itemize}
\end{definition}

\begin{example}
For this example, let $n=7$ and $\alpha=\left(  2,1,3,1\right)  $. Then,
$\ell=\ell\left(  \alpha\right)  =4$ and $\left(  I_{1},I_{2},I_{3}%
,I_{4}\right)  =\left(  \left\{  1,2\right\}  ,\left\{  3\right\}  ,\left\{
4,5,6\right\}  ,\left\{  7\right\}  \right)  $ (where we are using the
notations of Definition \ref{def.intervalperm}). Hence, $\overrightarrow{I_{1}%
}=12$, $\overrightarrow{I_{2}}=3$, $\overrightarrow{I_{3}}=456$ and
$\overrightarrow{I_{4}}=7$.

\begin{itemize}
\item[(a)] If $\tau\in\Symm_{\ell}=\Symm_4$ is the permutation
$\left(  2,3,1,4\right)  $, then $\operatorname*{iper}\left(  \alpha
,\tau\right)  $ is the permutation in $\Symm_7$ which (in one-line
notation) is the word $\overrightarrow{I_{\tau\left(  1\right)  }%
}\overrightarrow{I_{\tau\left(  2\right)  }}\overrightarrow{I_{\tau\left(
3\right)  }}\overrightarrow{I_{\tau\left(  4\right)  }}=\overrightarrow{I_{2}%
}\overrightarrow{I_{3}}\overrightarrow{I_{1}}\overrightarrow{I_{4}}=3456127$.

If $\tau\in\Symm_{\ell}=\Symm_4$ is the permutation $\left(
3,1,4,2\right)  $, then $\operatorname*{iper}\left(  \alpha,\tau\right)  $ is
the permutation in $\Symm_7$ which (in one-line notation) is the word
$\overrightarrow{I_{\tau\left(  1\right)  }}\overrightarrow{I_{\tau\left(
2\right)  }}\overrightarrow{I_{\tau\left(  3\right)  }}\overrightarrow{I_{\tau
\left(  4\right)  }}=\overrightarrow{I_{3}}\overrightarrow{I_{1}%
}\overrightarrow{I_{4}}\overrightarrow{I_{2}}=4561273$.

\item[(b)] The permutation $\sigma=\left(  3,7,4,5,6,1,2\right)
\in\Symm_7$ (given here in one-line notation) is $\alpha$-clumping, because:

\begin{itemize}
\item every $i\in\left\{  1,2,\ldots,\ell\right\}  =\left\{  1,2,3,4\right\}
$ has the property that $\sigma^{-1}\left(  I_{i}\right)  $ is an interval
(namely, $\sigma^{-1}\left(  I_{1}\right)  =\sigma^{-1}\left(  \left\{
1,2\right\}  \right)  =\left\{  6,7\right\}  $, $\sigma^{-1}\left(
I_{2}\right)  =\sigma^{-1}\left(  \left\{  3\right\}  \right)  =\left\{
1\right\}  $, $\sigma^{-1}\left(  I_{3}\right)  =\sigma^{-1}\left(  \left\{
4,5,6\right\}  \right)  =\left\{  3,4,5\right\}  $ and $\sigma^{-1}\left(
I_{4}\right)  =\sigma^{-1}\left(  \left\{  7\right\}  \right)  =\left\{
2\right\}  $), and

\item the restrictions of the map $\sigma^{-1}$ to the intervals $I_{i}$ are
increasing (this means that $\sigma^{-1}\left(  1\right)  <\sigma^{-1}\left(
2\right)  $ and $\sigma^{-1}\left(  4\right)  <\sigma^{-1}\left(  5\right)
<\sigma^{-1}\left(  6\right)  $, since the one-element intervals $I_{2}$ and
$I_{4}$ do not contribute anything to this condition).
\end{itemize}
\end{itemize}
\end{example}

Here is a more or less trivial observation:

\begin{proposition}
\label{prop.iper.intinverse}Let $n\in{ \NN }$. Let $\alpha$ be a
composition of $n$. Let $\ell=\ell\left(  \alpha\right)  $. Write $\alpha$ in
the form $\left(  \alpha_{1},\alpha_{2},\ldots,\alpha_{\ell}\right)  $. The
interval system corresponding to $\alpha$ is an $\ell$-tuple of intervals
(since $\ell\left(  \alpha\right)  =\ell$); denote this $\ell$-tuple by
$\left(  I_{1},I_{2},\ldots,I_{\ell}\right)  $. Let $\tau\in\Symm_{\ell}$.
Set $\sigma=\operatorname*{iper}\left(  \alpha,\tau\right)  $.

\begin{itemize}
\item[(a)] We have $\sigma^{-1}\left(  I_{\tau\left(  j\right)  }\right)
=\left[  \sum_{k=1}^{j-1}\alpha_{\tau\left(  k\right)  }:\sum_{k=1}^{j}%
\alpha_{\tau\left(  k\right)  }\right]  ^{+}$ for every $j\in\left\{
1,2,\ldots,\ell\right\}  $.

\item[(b)] For every $j\in\left\{  1,2,\ldots,\ell\right\}  $, the restriction
of the map $\sigma^{-1}$ to the interval $I_{\tau\left(  j\right)  }$ is increasing.

\item[(c)] The permutation $\operatorname*{iper}\left(  \alpha,\tau\right)  $
is $\alpha$-clumping.

\item[(d)] Let $i\in\left\{  1,2,\ldots,\ell-1\right\}  $. Then, the sets
$\sigma^{-1}\left(  I_{\tau\left(  i\right)  }\right)  $, $\sigma^{-1}\left(
I_{\tau\left(  i+1\right)  }\right)  $ and $\sigma^{-1}\left(  I_{\tau\left(
i\right)  }\right)  \cup\sigma^{-1}\left(  I_{\tau\left(  i+1\right)
}\right)  $ are nonempty intervals. Also, $\sigma^{-1}\left(  I_{\tau\left(
i\right)  }\right)  <\sigma^{-1}\left(  I_{\tau\left(  i+1\right)  }\right)  $.
\end{itemize}
\end{proposition}

\begin{exercise}
\label{exe.prop.iper.intinverse}Prove Proposition \ref{prop.iper.intinverse}.
\end{exercise}

\begin{proposition}
\label{prop.iper.biject} Let $n\in{ \NN }$. Let $\alpha$ be a composition
of $n$. Let $\ell=\ell\left(  \alpha\right)  $.

\begin{itemize}
\item[(a)] Define a map
\begin{align*}
\operatorname{iper}_{\alpha}:\Symm_{\ell}  &  \longrightarrow\left\{
\omega\in\Symm_n\mid\omega\text{ is }\alpha\text{-clumping}\right\}
,\\
\tau &  \longmapsto\operatorname{iper}\left(  \alpha,\tau\right)
\end{align*}
\footnote{This map is well-defined because for every $\tau\in\Symm_{\ell}$,
the permutation $\operatorname{iper}\left(  \alpha,\tau\right)  $ is
$\alpha$-clumping (according to Proposition \ref{prop.iper.intinverse}(c)).}.
This map $\operatorname{iper}_{\alpha}$ is bijective.

\item[(b)] Let $\sigma\in\Symm_n$ be an $\alpha$-clumping
permutation. Then, there exists a unique $\tau\in\Symm_{\ell}$
satisfying $\sigma=\operatorname{iper}\left(  \alpha,\tau\right)  $.
\end{itemize}
\end{proposition}

\begin{exercise}
\label{exe.prop.iper.biject}Prove Proposition \ref{prop.iper.biject}.
\end{exercise}

Next, we recall that the concatenation $\alpha\cdot\beta$ of two compositions
$\alpha$ and $\beta$ is defined in the same way as the concatenation of two
words; if we regard compositions as words over the alphabet $\left\{
1,2,3,\ldots\right\}  $, then the concatenation $\alpha\cdot\beta$ of two
compositions $\alpha$ and $\beta$ \textbf{is} the concatenation $\alpha\beta$
of the words $\alpha$ and $\beta$. Thus, we are going to write $\alpha\beta$
for the concatenation $\alpha\cdot\beta$ of two compositions $\alpha$ and
$\beta$ from now on.

\begin{proposition}
\label{prop.iper.shuffle}Let $n\in \NN$ and $m\in \NN$. Let
$\alpha$ be a composition of $n$, and $\beta$ be a composition of $m$. Let
$p=\ell\left(  \alpha\right)  $ and $q=\ell\left(  \beta\right)  $. Let
$\tau\in\Symm_{p+q}$. Notice that $\operatorname*{iper}\left(
\alpha\beta,\tau\right)  \in\Symm_{n+m}$ (since $\alpha\beta$ is a
composition of $n+m$ having length $\ell\left(  \alpha\beta\right)
=\ell\left(  \alpha\right)  +\ell\left(  \beta\right)  =p+q$). Then, $\tau
\in\operatorname{Sh}_{p,q}$ if and only if $\operatorname*{iper}%
\left(  \alpha\beta,\tau\right)  \in\operatorname{Sh}_{n,m}$.
\end{proposition}

\begin{exercise}
\label{exe.prop.iper.shuffle}Prove Proposition \ref{prop.iper.shuffle}.
\end{exercise}

Here is one more simple fact:

\begin{lemma}
\label{lem.shuffle.lyndon.step2}Let $u$ and $v$ be two words. Let
$n=\ell\left(  u\right)  $ and $m=\ell\left(  v\right)  $. Let $\alpha$ be a
composition of $n$, and let $\beta$ be a composition of $m$. Let
$p=\ell\left(  \alpha\right)  $ and $q=\ell\left(  \beta\right)  $. The
concatenation $\alpha\beta$ is a composition of $n+m$ having length
$\ell\left(  \alpha\beta\right)  =\ell\left(  \alpha\right)  +\ell\left(
\beta\right)  =p+q$. Thus, the interval system corresponding to $\alpha\beta$
is a $\left(  p+q\right)  $-tuple of intervals which covers $\left[
0:n+m\right]  ^{+}$. Denote this $\left(  p+q\right)  $-tuple by $\left(
I_{1},I_{2},\ldots,I_{p+q}\right)  $.

Let $\tau\in\operatorname{Sh}_{p,q}$. Set $\sigma
=\operatorname*{iper}\left(  \alpha\beta,\tau\right)  $. Then,%
\[
u\underset{\sigma}{\shuffle}v=\left(  uv\right)  \left[  I_{\tau\left(
1\right)  }\right]  \cdot\left(  uv\right)  \left[  I_{\tau\left(  2\right)
}\right]  \cdot\cdots\cdot\left(  uv\right)  \left[  I_{\tau\left(
p+q\right)  }\right]  .
\]

\end{lemma}

\begin{exercise}
\label{exe.lem.shuffle.lyndon.step2}Prove Lemma \ref{lem.shuffle.lyndon.step2}.
\end{exercise}

Having these notations and trivialities in place, we can say a bit more about
the lexicographically highest element of a shuffle product than what was said
in Theorem \ref{thm.shuffle.lyndon.used}:

\begin{theorem}
\label{thm.shuffle.lyndon.full}Let $u$ and $v$ be two words. Let
$n=\ell\left(  u\right)  $ and $m=\ell\left(  v\right)  $.

Let $\left(  a_{1},a_{2},\ldots,a_{p}\right)  $ be the CFL factorization of
$u$. Let $\left(  b_{1},b_{2},\ldots,b_{q}\right)  $ be the CFL factorization
of $v$.

Let $\alpha$ be the $p$-tuple $\left(  \ell\left(  a_{1}\right)  ,\ell\left(
a_{2}\right)  ,\ldots,\ell\left(  a_{p}\right)  \right)  $. Then, $\alpha$ is
a composition\footnote{since Lyndon words are nonempty, and thus $\ell\left(
a_{i}\right)  >0$ for every $i$} of length $p$ and size $\sum_{k=1}^{p}%
\ell\left(  a_{k}\right)  =\ell\left(  \underbrace{a_{1}a_{2} \cdots a_{p}%
}_{=u}\right)  =\ell\left(  u\right)  =n$.

Let $\beta$ be the $q$-tuple $\left(  \ell\left(  b_{1}\right)  ,\ell\left(
b_{2}\right)  ,\ldots,\ell\left(  b_{q}\right)  \right)  $. Then, $\beta$ is a
composition of length $q$ and size $\sum_{k=1}^{q}\ell\left(  b_{k}\right)
=m$.\ \ \ \ \footnote{The proof of this is the same as the proof of the fact
that $\alpha$ is a composition of length $p$ and size $\sum_{k=1}^{p}%
\ell\left(  \alpha_{k}\right)  =n$.}

Now, $\alpha$ is a composition of length $p$ and size $n$, and $\beta$ is a
composition of length $q$ and size $m$. Thus, the concatenation $\alpha\beta$
of these two tuples is a composition of length $p+q$ and size $n+m$. The
interval system corresponding to this composition $\alpha\beta$ is a $\left(
p+q\right)  $-tuple (since said composition has length $p+q$); denote this
$\left(  p+q\right)  $-tuple by $\left(  I_{1},I_{2},\ldots,I_{p+q}\right)  $.

\begin{itemize}
\item[(a)] If $\tau\in\operatorname{Sh}_{p,q}$ satisfies $\left(
uv\right)  \left[  I_{\tau\left(  1\right)  }\right]  \geq\left(  uv\right)
\left[  I_{\tau\left(  2\right)  }\right]  \geq\cdots\geq\left(  uv\right)
\left[  I_{\tau\left(  p+q\right)  }\right]  $, and if we set $\sigma
=\operatorname*{iper}\left(  \alpha\beta,\tau\right)  $, then $\sigma
\in\operatorname{Sh}_{n,m}$, and the word $u\underset{\sigma
}{\shuffle}v$ is the lexicographically highest element of the multiset
$u\shuffle v$.

\item[(b)] Let $\sigma\in\operatorname{Sh}_{n,m}$ be a permutation
such that $u\underset{\sigma}{\shuffle}v$ is the lexicographically highest
element of the multiset $u\shuffle v$. Then, there exists a unique permutation
$\tau\in\operatorname{Sh}_{p,q}$ satisfying $\left(  uv\right)
\left[  I_{\tau\left(  1\right)  }\right]  \geq\left(  uv\right)  \left[
I_{\tau\left(  2\right)  }\right]  \geq\cdots\geq\left(  uv\right)  \left[
I_{\tau\left(  p+q\right)  }\right]  $ and $\sigma=\operatorname*{iper}\left(
\alpha\beta,\tau\right)  $.
\end{itemize}
\end{theorem}

\begin{proof}
Before we step to the actual proof, we need to make some preparation. First of
all, $\left(  I_{1},I_{2},\ldots,I_{p+q}\right)  $ is the interval system
corresponding to the composition $\alpha\beta$. In other words,%
\begin{equation}
\left(  I_{1},I_{2},\ldots,I_{p+q}\right)  =\operatorname*{intsys}\left(
\alpha\beta\right)  . \label{pf.thm.shuffle.lyndon.full.0}%
\end{equation}
But since $\alpha=\left(  \ell\left(  a_{1}\right)  ,\ell\left(  a_{2}\right)
,\ldots,\ell\left(  a_{p}\right)  \right)  $ and $\beta=\left(  \ell\left(
b_{1}\right)  ,\ell\left(  b_{2}\right)  ,\ldots,\ell\left(  b_{q}\right)
\right)  $, we have%
\[
\alpha\beta=\left(  \ell\left(  a_{1}\right)  ,\ell\left(  a_{2}\right)
,\ldots,\ell\left(  a_{p}\right)  ,\ell\left(  b_{1}\right)  ,\ell\left(
b_{2}\right)  ,\ldots,\ell\left(  b_{q}\right)  \right)  .
\]
Thus, (\ref{pf.thm.shuffle.lyndon.full.0}) rewrites as%
\[
\left(  I_{1},I_{2},\ldots,I_{p+q}\right)  =\operatorname*{intsys}\left(
\ell\left(  a_{1}\right)  ,\ell\left(  a_{2}\right)  ,\ldots,\ell\left(
a_{p}\right)  ,\ell\left(  b_{1}\right)  ,\ell\left(  b_{2}\right)
,\ldots,\ell\left(  b_{q}\right)  \right)  .
\]
By the definition of $\operatorname*{intsys}\left(  \ell\left(  a_{1}\right)
,\ell\left(  a_{2}\right)  ,\ldots,\ell\left(  a_{p}\right)  ,\ell\left(
b_{1}\right)  ,\ell\left(  b_{2}\right)  ,\ldots,\ell\left(  b_{q}\right)
\right)  $, we thus have%
\[
I_{i}=\left[  \sum_{k=1}^{i-1}\ell\left(  a_{k}\right)  :\sum_{k=1}^{i}%
\ell\left(  a_{k}\right)  \right]  ^{+}\ \ \ \ \ \ \ \ \ \ \text{for every
}i\in\left\{  1,2,\ldots,p\right\}  ,
\]
and besides
\[
I_{p+j}=\left[  n+\sum_{k=1}^{j-1}\ell\left(  b_{k}\right)  :n+\sum_{k=1}%
^{j}\ell\left(  b_{k}\right)  \right]  ^{+}\ \ \ \ \ \ \ \ \ \ \text{for every
}j\in\left\{  1,2,\ldots,q\right\}
\]
(since $\sum_{k=1}^{p}\ell\left(  a_{k}\right)  =n$). Moreover,
Remark~\ref{rmk.intsys.trivia}(c) shows that $\left(  I_{1},I_{2}%
,\ldots,I_{p+q}\right)  $ is a $\left(  p+q\right)  $-tuple of nonempty
intervals of $\ZZ$ and satisfies the following three properties:

\begin{itemize}
\item The intervals $I_{1}$, $I_{2}$, $\ldots$, $I_{p+q}$ form a set partition
of the set $\left[  0:n+m\right]  ^{+}$.

\item We have $I_{1}<I_{2}< \cdots<I_{p+q}$.

\item We have $\left\vert I_{i}\right\vert =\ell\left(  a_{i}\right)  $ for
every $i\in\left\{  1,2,\ldots,p\right\}  $ and $\left\vert I_{p+j}\right\vert
=\ell\left(  b_{j}\right)  $ for every $j\in\left\{  1,2,\ldots,q\right\}  $.
\end{itemize}

Of course, every $i\in\left\{  1,2,\ldots,p\right\}  $ satisfies%
\begin{equation}
I_{i}\subset\left[  0:n\right]  ^{+}\ \ \ \ \ \ \ \ \ \ \text{and}%
\ \ \ \ \ \ \ \ \ \ \left(  uv\right)  \left[  I_{i}\right]  =u\left[
I_{i}\right]  =a_{i}. \label{pf.thm.shuffle.lyndon.full.prep1}%
\end{equation}
Meanwhile, every $i\in\left\{  p+1,p+2,\ldots,p+q\right\}  $ satisfies%
\begin{equation}
I_{i}\subset\left[  n:n+m\right]  ^{+}\ \ \ \ \ \ \ \ \ \ \text{and}%
\ \ \ \ \ \ \ \ \ \ \left(  uv\right)  \left[  I_{i}\right]  =v\left[
I_{i}-n\right]  =b_{i-p} \label{pf.thm.shuffle.lyndon.full.prep1'}%
\end{equation}
(where $I_{i}-n$ denotes the interval $\left\{  k-n\ \mid\ k\in I_{i}\right\}
$). We thus see that%
\begin{equation}
\left(  uv\right)  \left[  I_{i}\right]  \text{ is a Lyndon word}%
\ \ \ \ \ \ \ \ \ \ \text{for every }i\in\left\{  1,2,\ldots,p+q\right\}
\label{pf.thm.shuffle.lyndon.full.prep1lynd}%
\end{equation}
\footnote{Indeed, when $i\leq p$, this follows from
(\ref{pf.thm.shuffle.lyndon.full.prep1}) and the fact that $a_{i}$ is Lyndon;
whereas in the other case, this follows from
(\ref{pf.thm.shuffle.lyndon.full.prep1'}) and the fact that $b_{i-p}$ is
Lyndon.}.

By the definition of a CFL factorization, we have $a_{1}\geq a_{2}\geq
\cdots\geq a_{p}$ and $b_{1}\geq b_{2}\geq\cdots\geq b_{q}$.

We have $\sigma\in\operatorname{Sh}_{n,m}$, so that $\sigma
^{-1}\left(  1\right)  <\sigma^{-1}\left(  2\right)  <\cdots<\sigma
^{-1}\left(  n\right)  $ and $\sigma^{-1}\left(  n+1\right)  <\sigma
^{-1}\left(  n+2\right)  <\cdots<\sigma^{-1}\left(  n+m\right)  $. In other
words, the restriction of the map $\sigma^{-1}$ to the interval $\left[
0:n\right]  ^{+}$ is strictly increasing, and so is the restriction of the map
$\sigma^{-1}$ to the interval $\left[  n:n+m\right]  ^{+}$.

(b) We will first show that%
\begin{equation}
\text{if }J\subset\left[  0:n\right]  ^{+}\text{ is an interval such that the
word }\left(  uv\right)  \left[  J\right]  \text{ is Lyndon, then }\sigma
^{-1}\left(  J\right)  \text{ is an interval.}
\label{pf.thm.shuffle.lyndon.full.1}%
\end{equation}

\textit{Proof of (\ref{pf.thm.shuffle.lyndon.full.1}):} We will prove
(\ref{pf.thm.shuffle.lyndon.full.1}) by strong induction over $\left\vert
J\right\vert $.

So, fix some $N\in \NN$. Assume (as the induction hypothesis) that
(\ref{pf.thm.shuffle.lyndon.full.1}) has been proven whenever $\left\vert
J\right\vert <N$. We now need to prove (\ref{pf.thm.shuffle.lyndon.full.1})
when $\left\vert J\right\vert =N$.

Let $J\subset\left[  0:n\right]  ^{+}$ be an interval such that the word
$\left(  uv\right)  \left[  J\right]  $ is Lyndon and such that $\left\vert
J\right\vert =N$. We have to prove that $\sigma^{-1}\left(  J\right)  $ is an
interval. This is obvious if $\left\vert J\right\vert =1$ (because in this
case, $\sigma^{-1}\left(  J\right)  $ is a one-element set, thus trivially an
interval). Hence, we WLOG assume that we don't have $\left\vert J\right\vert
=1$. We also don't have $\left\vert J\right\vert =0$, because $\left(
uv\right)  \left[  J\right]  $ has to be Lyndon (and the empty word is not).
So we have $\left\vert J\right\vert >1$. Now, $\ell\left(  \left(  uv\right)
\left[  J\right]  \right)  =\left\vert J\right\vert >1$, and thus $\left(
uv\right)  \left[  J\right]  $ is a Lyndon word of length $>1$. Let
$v^{\prime}$ be the (lexicographically) smallest nonempty \textbf{proper}
suffix of $\left(  uv\right)  \left[  J\right]  $. Since $v^{\prime}$ is a
proper suffix of $w$, there exists a nonempty $u^{\prime}\in\mathfrak{A}%
^{\ast}$ such that $\left(  uv\right)  \left[  J\right]  =u^{\prime}v^{\prime
}$. Consider this $u^{\prime}$.

Now, Theorem \ref{thm.words.lyndon.std}(a) (applied to $\left(  uv\right)
\left[  J\right]  $, $u^{\prime}$ and $v^{\prime}$ instead of $w$, $u$ and
$v$) yields that the words $u^{\prime}$ and $v^{\prime}$ are Lyndon. Also,
Theorem \ref{thm.words.lyndon.std}(b) (applied to $\left(  uv\right)  \left[
J\right]  $, $u^{\prime}$ and $v^{\prime}$ instead of $w$, $u$ and $v$) yields
that $u^{\prime} < \left(uv\right) \left[ J \right] < v^{\prime}$.

But from the fact that $\left(  uv\right)  \left[  J\right]  =u^{\prime
}v^{\prime}$ with $u^{\prime}$ and $v^{\prime}$ both being nonempty, it
becomes immediately clear that we can write $J$ as a union of two disjoint
nonempty intervals $K$ and $L$ such that $K<L$, $u^{\prime}=\left(  uv\right)
\left[  K\right]  $ and $v^{\prime}=\left(  uv\right)  \left[  L\right]  $.
Consider these $K$ and $L$. The intervals $K$ and $L$ are nonempty and have
their sizes add up to $\left\vert J\right\vert$
(since they are disjoint and their union is $J$),
and hence both must have size smaller than $\left\vert J\right\vert =N$. So
$K\subset\left[  0:n\right]  ^{+}$ is an interval of size $\left\vert
K\right\vert <N$ having the property that $\left(  uv\right)  \left[
K\right]  $ is Lyndon (since $\left(  uv\right)  \left[  K\right]  =u^{\prime
}$ is Lyndon). Thus, we can apply (\ref{pf.thm.shuffle.lyndon.full.1}) to $K$
instead of $J$ (because of the induction hypothesis). As a result, we conclude
that $\sigma^{-1}\left(  K\right)  $ is an interval. Similarly, we can apply
(\ref{pf.thm.shuffle.lyndon.full.1}) to $L$ instead of $J$ (we know that
$\left(  uv\right)  \left[  L\right]  $ is Lyndon since $\left(  uv\right)
\left[  L\right]  =v^{\prime}$), and learn that $\sigma^{-1}\left(  L\right)
$ is an interval. The intervals $\sigma^{-1}\left(  K\right)  $ and
$\sigma^{-1}\left(  L\right)  $ are both nonempty (since $K$ and $L$ are
nonempty), and their union is $\sigma^{-1}\left(  J\right)  $ (because the
union of $K$ and $L$ is $J$). The nonempty intervals $K$ and $L$ both are
subsets of $\left[  0:n\right]  ^{+}$ (since their union is $J\subset\left[
0:n\right]  ^{+}$), and their union $K\cup L$ is an interval (since their
union $K\cup L$ is $J$, and we know that $J$ is an interval).

Now, assume (for the sake of contradiction) that $\sigma^{-1}\left(  J\right)
$ is not an interval. Since $J$ is the union of $K$ and $L$, we have $J=K\cup
L$ and thus $\sigma^{-1}\left(  J\right)  =\sigma^{-1}\left(  K\cup L\right)
=\sigma^{-1}\left(  K\right)  \cup\sigma^{-1}\left(  L\right)  $ (since
$\sigma$ is a bijection). Therefore, $\sigma^{-1}\left(  K\right)  \cup
\sigma^{-1}\left(  L\right)  $ is not an interval (since $\sigma^{-1}\left(
J\right)  $ is not an interval). Thus, Lemma \ref{lem.shuffle.lyndon.easystep}%
(b) yields that there exists a nonempty interval $P\subset\left[
n:n+m\right]  ^{+}$ such that $\sigma^{-1}\left(  P\right)  $, $\sigma
^{-1}\left(  K\right)  \cup\sigma^{-1}\left(  P\right)  $ and $\sigma
^{-1}\left(  P\right)  \cup\sigma^{-1}\left(  L\right)  $ are intervals and
such that $\sigma^{-1}\left(  K\right)  <\sigma^{-1}\left(  P\right)
<\sigma^{-1}\left(  L\right)  $. Consider this $P$. Since $P$ is nonempty, we
have $\left\vert P\right\vert \neq0$.

Lemma \ref{lem.shuffle.lyndon.step}(b) (applied to $K$ and $P$ instead of $I$
and $J$) yields
\begin{equation}
\left(  uv\right)  \left[  K\right]  \cdot\left(  uv\right)  \left[  P\right]
\geq\left(  uv\right)  \left[  P\right]  \cdot\left(  uv\right)  \left[
K\right]  . \label{pf.thm.shuffle.lyndon.full.5a}%
\end{equation}
Since $\left(  uv\right)  \left[  K\right]  =u^{\prime}$, this rewrites as%
\begin{equation}
u^{\prime}\cdot\left(  uv\right)  \left[  P\right]  \geq\left(  uv\right)
\left[  P\right]  \cdot u^{\prime}. \label{pf.thm.shuffle.lyndon.full.5a1}%
\end{equation}

But Lemma \ref{lem.shuffle.lyndon.step}(c) (applied to $P$ and $L$ instead of
$I$ and $J$) yields
\begin{equation}
\left(  uv\right)  \left[  P\right]  \cdot\left(  uv\right)  \left[  L\right]
\geq\left(  uv\right)  \left[  L\right]  \cdot\left(  uv\right)  \left[
P\right]  . \label{pf.thm.shuffle.lyndon.full.5b}%
\end{equation}
Since $\left(  uv\right)  \left[  L\right]  =v^{\prime}$, this rewrites as%
\begin{equation}
\left(  uv\right)  \left[  P\right]  \cdot v^{\prime}\geq v^{\prime}%
\cdot\left(  uv\right)  \left[  P\right]  .
\label{pf.thm.shuffle.lyndon.full.5b1}%
\end{equation}

Recall also that $u^{\prime}<v^{\prime}$, and that both words $u^{\prime}$ and
$v^{\prime}$ are Lyndon. Now, Corollary \ref{cor.words.lyndon.uv=vu} (applied
to $u^{\prime}$, $v^{\prime}$ and $\left(  uv\right)  \left[  P\right]  $
instead of $u$, $v$ and $z$) yields that $\left(  uv\right)  \left[  P\right]
$ is the empty word (because of (\ref{pf.thm.shuffle.lyndon.full.5a1}) and
(\ref{pf.thm.shuffle.lyndon.full.5b1})), so that $\ell\left(  \left(
uv\right)  \left[  P\right]  \right)  =0$. This contradicts $\ell\left(
\left(  uv\right)  \left[  P\right]  \right)  =\left\vert P\right\vert \neq0$.
This contradiction shows that our assumption (that $\sigma^{-1}\left(
J\right)  $ is not an interval) was wrong. Hence, $\sigma^{-1}\left(
J\right)  $ is an interval. This completes the induction step, and thus
(\ref{pf.thm.shuffle.lyndon.full.1}) is proven.

Similarly to (\ref{pf.thm.shuffle.lyndon.full.1}), we can show that%
\begin{equation}
\text{if }J\subset\left[  n:n+m\right]  ^{+}\text{ is an interval such that
the word }\left(  uv\right)  \left[  J\right]  \text{ is Lyndon, then }%
\sigma^{-1}\left(  J\right)  \text{ is an interval.}
\label{pf.thm.shuffle.lyndon.full.1'}%
\end{equation}

Now, let $i\in\left\{  1,2,\ldots,p+q\right\}  $ be arbitrary. We are going to
prove that%
\begin{equation}
\sigma^{-1}\left(  I_{i}\right)  \text{ is an interval.}
\label{pf.thm.shuffle.lyndon.full.21}%
\end{equation}

\textit{Proof of (\ref{pf.thm.shuffle.lyndon.full.21}):} We must be in one of
the following two cases:

\textit{Case 1:} We have $i\in\left\{  1,2,\ldots,p\right\}  $.

\textit{Case 2:} We have $i\in\left\{  p+1,p+2,\ldots,p+q\right\}  $.

Let us first consider Case 1. In this case, we have $i\in\left\{
1,2,\ldots,p\right\}  $. Thus, $I_{i}\subset\left[  0:n\right]  ^{+}$ (by
(\ref{pf.thm.shuffle.lyndon.full.prep1})). Also,
(\ref{pf.thm.shuffle.lyndon.full.prep1}) yields that $\left(  uv\right)
\left[  I_{i}\right]  =a_{i}$ is a Lyndon word. Hence,
(\ref{pf.thm.shuffle.lyndon.full.1}) (applied to $J=I_{i}$) yields that
$\sigma^{-1}\left(  I_{i}\right)  $ is an interval. Thus,
(\ref{pf.thm.shuffle.lyndon.full.21}) is proven in Case 1.

Similarly, we can prove (\ref{pf.thm.shuffle.lyndon.full.21}) in Case 2, using
(\ref{pf.thm.shuffle.lyndon.full.prep1'}) and
(\ref{pf.thm.shuffle.lyndon.full.1'}) instead of
(\ref{pf.thm.shuffle.lyndon.full.prep1}) and
(\ref{pf.thm.shuffle.lyndon.full.1}), respectively. Hence,
(\ref{pf.thm.shuffle.lyndon.full.21}) is proven.

So we know that $\sigma^{-1}\left(  I_{i}\right)  $ is an interval. But we
also know that either $I_{i}\subset\left[  0:n\right]  ^{+}$ or $I_{i}%
\subset\left[  n:n+m\right]  ^{+}$ (depending on whether $i\leq p$ or $i>p$).
As a consequence, the restriction of the map $\sigma^{-1}$ to the interval
$I_{i}$ is increasing (because the restriction of the map $\sigma^{-1}$ to the
interval $\left[  0:n\right]  ^{+}$ is strictly increasing, and so is the
restriction of the map $\sigma^{-1}$ to the interval $\left[  n:n+m\right]
^{+}$).

Now, let us forget that we fixed $i$. We thus have shown that every
$i\in\left\{  1,2,\ldots,p+q\right\}  $ has the two properties that:

\begin{itemize}
\item the set $\sigma^{-1}\left(  I_{i}\right)  $ is an interval;

\item the restriction of the map $\sigma^{-1}$ to the interval $I_{i}$ is increasing.
\end{itemize}

\noindent In other words, the permutation $\sigma$ is
$\left(  \alpha\beta\right)$-clumping
(since $\left(  I_{1},I_{2},\ldots,I_{p+q}\right)  $ is the
interval system corresponding to the composition $\alpha\beta$). Hence,
Proposition \ref{prop.iper.biject}(b) (applied to $n+m$, $\alpha\beta$ and
$p+q$ instead of $n$, $\alpha$ and $\ell$) shows that there exists a unique
$\tau\in\Symm_{p+q}$ satisfying $\sigma=\operatorname*{iper}\left(
\alpha\beta,\tau\right)  $. Thus, the uniqueness part of Theorem
\ref{thm.shuffle.lyndon.full}(b) (i.e., the claim that the $\tau$ in Theorem
\ref{thm.shuffle.lyndon.full}(b) is unique if it exists) is proven.

It now remains to prove the existence part of Theorem
\ref{thm.shuffle.lyndon.full}(b), i.e., to prove that there exists at least
one permutation $\tau\in\operatorname{Sh}_{p,q}$ satisfying $\left(
uv\right)  \left[  I_{\tau\left(  1\right)  }\right]  \geq\left(  uv\right)
\left[  I_{\tau\left(  2\right)  }\right]  \geq\cdots\geq\left(  uv\right)
\left[  I_{\tau\left(  p+q\right)  }\right]  $ and $\sigma
=\operatorname*{iper}\left(  \alpha\beta,\tau\right)  $. We already know that
there exists a unique $\tau\in\Symm_{p+q}$ satisfying $\sigma
=\operatorname*{iper}\left(  \alpha\beta,\tau\right)  $. Consider this $\tau$.
We will now prove that $\left(  uv\right)  \left[  I_{\tau\left(  1\right)
}\right]  \geq\left(  uv\right)  \left[  I_{\tau\left(  2\right)  }\right]
\geq\cdots\geq\left(  uv\right)  \left[  I_{\tau\left(  p+q\right)  }\right]
$ and $\tau\in\operatorname{Sh}_{p,q}$. Once this is done, the
existence part of Theorem \ref{thm.shuffle.lyndon.full}(b) will be proven, and
thus the proof of Theorem \ref{thm.shuffle.lyndon.full}(b) will be complete.

Proposition \ref{prop.iper.shuffle} yields that $\tau\in
\operatorname{Sh}_{p,q}$ if and only if $\operatorname*{iper}\left(  \alpha\beta
,\tau\right)  \in\operatorname{Sh}_{n,m}$. Since we know that
$\operatorname*{iper}\left(  \alpha\beta,\tau\right)  =\sigma\in
\operatorname{Sh}_{n,m}$, we thus conclude that $\tau\in
\operatorname{Sh}_{p,q}$. The only thing that remains to be proven
now is that%
\begin{equation}
\left(  uv\right)  \left[  I_{\tau\left(  1\right)  }\right]  \geq\left(
uv\right)  \left[  I_{\tau\left(  2\right)  }\right]  \geq\cdots\geq\left(
uv\right)  \left[  I_{\tau\left(  p+q\right)  }\right]  .
\label{pf.thm.shuffle.lyndon.full.25}%
\end{equation}

\textit{Proof of (\ref{pf.thm.shuffle.lyndon.full.25}):} We have $\tau
\in\operatorname{Sh}_{p,q}$. In other words, $\tau^{-1}\left(
1\right)  <\tau^{-1}\left(  2\right)  <\cdots<\tau^{-1}\left(  p\right)  $ and
$\tau^{-1}\left(  p+1\right)  <\tau^{-1}\left(  p+2\right)  <\cdots<\tau
^{-1}\left(  p+q\right)  $. In other words, the restriction of the map
$\tau^{-1}$ to the interval $\left[  0:p\right]  ^{+}$ is strictly increasing,
and so is the restriction of the map $\tau^{-1}$ to the interval $\left[
p:p+q\right]  ^{+}$.

Let $i\in\left\{  1,2,\ldots,p+q-1\right\}  $. We will show that
\begin{equation}
\left(  uv\right)  \left[  I_{\tau\left(  i\right)  }\right]  \geq\left(
uv\right)  \left[  I_{\tau\left(  i+1\right)  }\right]  .
\label{pf.thm.shuffle.lyndon.full.26}%
\end{equation}

Clearly, both $\tau\left(  i\right)  $ and $\tau\left(  i+1\right)  $ belong
to $\left\{  1,2,\ldots,p+q\right\}  =\left\{  1,2,\ldots,p\right\}
\cup\left\{  p+1,p+2,\ldots,p+q\right\}  $. Thus, we must be in one of the
following four cases:

\textit{Case 1:} We have $\tau\left(  i\right)  \in\left\{  1,2,\ldots
,p\right\}  $ and $\tau\left(  i+1\right)  \in\left\{  1,2,\ldots,p\right\}  $.

\textit{Case 2:} We have $\tau\left(  i\right)  \in\left\{  1,2,\ldots
,p\right\}  $ and $\tau\left(  i+1\right)  \in\left\{  p+1,p+2,\ldots
,p+q\right\}  $.

\textit{Case 3:} We have $\tau\left(  i\right)  \in\left\{  p+1,p+2,\ldots
,p+q\right\}  $ and $\tau\left(  i+1\right)  \in\left\{  1,2,\ldots,p\right\}
$.

\textit{Case 4:} We have $\tau\left(  i\right)  \in\left\{  p+1,p+2,\ldots
,p+q\right\}  $ and $\tau\left(  i+1\right)  \in\left\{  p+1,p+2,\ldots
,p+q\right\}  $.

Let us consider Case 1 first. In this case, we have $\tau\left(  i\right)
\in\left\{  1,2,\ldots,p\right\}  $ and $\tau\left(  i+1\right)  \in\left\{
1,2,\ldots,p\right\}  $. From the fact that the restriction of the map
$\tau^{-1}$ to the interval $\left[  0:p\right]  ^{+}$ is strictly increasing,
we can easily deduce $\tau\left(  i\right)  <\tau\left(  i+1\right)
$\ \ \ \ \footnote{\textit{Proof.} Assume the contrary. Then, $\tau\left(
i\right)  \geq\tau\left(  i+1\right)  $. Since both $\tau\left(  i\right)  $
and $\tau\left(  i+1\right)  $ belong to $\left\{  1,2,\ldots,p\right\}
=\left[  0:p\right]  ^{+}$, this yields $\tau^{-1}\left(  \tau\left(
i\right)  \right)  \geq\tau^{-1}\left(  \tau\left(  i+1\right)  \right)  $
(since the restriction of the map $\tau^{-1}$ to the interval $\left[
0:p\right]  ^{+}$ is strictly increasing), which contradicts $\tau^{-1}\left(
\tau\left(  i\right)  \right)  =i<i+1=\tau^{-1}\left(  \tau\left(  i+1\right)
\right)  $. This contradiction proves the assumption wrong, qed.}. Therefore,
$a_{\tau\left(  i\right)  }\geq a_{\tau\left(  i+1\right)  }$ (since
$a_{1}\geq a_{2}\geq\cdots\geq a_{p}$).

But $\left(  uv\right)  \left[  I_{\tau\left(  i\right)  }\right]
=a_{\tau\left(  i\right)  }$ (by (\ref{pf.thm.shuffle.lyndon.full.prep1}),
applied to $\tau\left(  i\right)  $ instead of $i$) and $\left(  uv\right)
\left[  I_{\tau\left(  i+1\right)  }\right]  =a_{\tau\left(  i+1\right)  }$
(similarly). In view of these equalities, the inequality $a_{\tau\left(
i\right)  }\geq a_{\tau\left(  i+1\right)  }$ rewrites as $\left(  uv\right)
\left[  I_{\tau\left(  i\right)  }\right]  \geq\left(  uv\right)  \left[
I_{\tau\left(  i+1\right)  }\right]  $. Thus,
(\ref{pf.thm.shuffle.lyndon.full.26}) is proven in Case 1.

Similarly, we can show (\ref{pf.thm.shuffle.lyndon.full.26}) in Case 4
(observing that $\left(  uv\right)  \left[  I_{\tau\left(  i\right)  }\right]
=b_{\tau\left(  i\right)  -p}$ and $\left(  uv\right)  \left[  I_{\tau\left(
i+1\right)  }\right]  =b_{\tau\left(  i+1\right)  -p}$ in this case).

Let us now consider Case 2. In this case, we have $\tau\left(  i\right)
\in\left\{  1,2,\ldots,p\right\}  $ and $\tau\left(  i+1\right)  \in\left\{
p+1,p+2,\ldots,p+q\right\}  $. From $\tau\left(  i\right)  \in\left\{
1,2,\ldots,p\right\}  $, we conclude that $I_{\tau\left(  i\right)  }%
\subset\left[  0:n\right]  ^{+}$. From $\tau\left(  i+1\right)  \in\left\{
p+1,p+2,\ldots,p+q\right\}  $, we conclude that $I_{\tau\left(  i+1\right)
}\subset\left[  n:n+m\right]  ^{+}$. The intervals $I_{\tau\left(  i\right)
}$ and $I_{\tau\left(  i+1\right)  }$ are clearly nonempty.

Proposition \ref{prop.iper.intinverse}(d) (applied to $n+m$, $\alpha\beta$,
$p+q$ and $\left(  I_{1},I_{2},\ldots,I_{p+q}\right)  $ instead of $n$,
$\alpha$, $\ell$ and $\left(  I_{1},I_{2},\ldots,I_{\ell}\right)  $) yields
that the sets $\sigma^{-1}\left(  I_{\tau\left(  i\right)  }\right)  $,
$\sigma^{-1}\left(  I_{\tau\left(  i+1\right)  }\right)  $ and $\sigma
^{-1}\left(  I_{\tau\left(  i\right)  }\right)  \cup\sigma^{-1}\left(
I_{\tau\left(  i+1\right)  }\right)  $ are nonempty intervals, and that we
have $\sigma^{-1}\left(  I_{\tau\left(  i\right)  }\right)  <\sigma
^{-1}\left(  I_{\tau\left(  i+1\right)  }\right)  $. Hence, Lemma
\ref{lem.shuffle.lyndon.step}(b) (applied to $I=I_{\tau\left(  i\right)  }$
and $J=I_{\tau\left(  i+1\right)  }$) yields%
\[
\left(  uv\right)  \left[  I_{\tau\left(  i\right)  }\right]  \cdot\left(
uv\right)  \left[  I_{\tau\left(  i+1\right)  }\right]  \geq\left(  uv\right)
\left[  I_{\tau\left(  i+1\right)  }\right]  \cdot\left(  uv\right)  \left[
I_{\tau\left(  i\right)  }\right]  .
\]
But $\left(  uv\right)  \left[  I_{\tau\left(  i\right)  }\right]  $ and
$\left(  uv\right)  \left[  I_{\tau\left(  i+1\right)  }\right]  $ are Lyndon
words (as a consequence of (\ref{pf.thm.shuffle.lyndon.full.prep1lynd})).
Thus, Proposition \ref{prop.words.lyndon.preorder} (applied to $\left(
uv\right)  \left[  I_{\tau\left(  i\right)  }\right]  $ and $\left(
uv\right)  \left[  I_{\tau\left(  i+1\right)  }\right]  $ instead of $u$ and
$v$) shows that $\left(  uv\right)  \left[  I_{\tau\left(  i\right)  }\right]
\geq\left(  uv\right)  \left[  I_{\tau\left(  i+1\right)  }\right]  $ if and
only if $\left(  uv\right)  \left[  I_{\tau\left(  i\right)  }\right]
\cdot\left(  uv\right)  \left[  I_{\tau\left(  i+1\right)  }\right]
\geq\left(  uv\right)  \left[  I_{\tau\left(  i+1\right)  }\right]
\cdot\left(  uv\right)  \left[  I_{\tau\left(  i\right)  }\right]  $. Since we
know that $\left(  uv\right)  \left[  I_{\tau\left(  i\right)  }\right]
\cdot\left(  uv\right)  \left[  I_{\tau\left(  i+1\right)  }\right]
\geq\left(  uv\right)  \left[  I_{\tau\left(  i+1\right)  }\right]
\cdot\left(  uv\right)  \left[  I_{\tau\left(  i\right)  }\right]  $ holds, we
thus conclude that $\left(  uv\right)  \left[  I_{\tau\left(  i\right)
}\right]  \geq\left(  uv\right)  \left[  I_{\tau\left(  i+1\right)  }\right]
$. Thus, (\ref{pf.thm.shuffle.lyndon.full.26}) is proven in Case 2.

The proof of (\ref{pf.thm.shuffle.lyndon.full.26}) in Case 3 is analogous to
that in Case 2 (the main difference being that Lemma
\ref{lem.shuffle.lyndon.step}(c) is used in lieu of Lemma
\ref{lem.shuffle.lyndon.step}(b)).

Thus, (\ref{pf.thm.shuffle.lyndon.full.26}) is proven in all possible cases.
So we always have (\ref{pf.thm.shuffle.lyndon.full.26}). In other words,
$\left(  uv\right)  \left[  I_{\tau\left(  i\right)  }\right]  \geq\left(
uv\right)  \left[  I_{\tau\left(  i+1\right)  }\right]  $.

Now, forget that we fixed $i$. We hence have shown that $\left(  uv\right)
\left[  I_{\tau\left(  i\right)  }\right]  \geq\left(  uv\right)  \left[
I_{\tau\left(  i+1\right)  }\right]  $ for all $i\in\left\{  1,2,\ldots
,p+q-1\right\}  $. This proves (\ref{pf.thm.shuffle.lyndon.full.25}), and thus
completes our proof of Theorem \ref{thm.shuffle.lyndon.full}(b).

(a) Let $\tau\in\operatorname{Sh}_{p,q}$ be such that
\begin{equation}
\left(  uv\right)  \left[  I_{\tau\left(  1\right)  }\right]  \geq\left(
uv\right)  \left[  I_{\tau\left(  2\right)  }\right]  \geq\cdots\geq\left(
uv\right)  \left[  I_{\tau\left(  p+q\right)  }\right]  .
\label{pf.thm.shuffle.lyndon.full.a.1}%
\end{equation}
Set $\sigma=\operatorname*{iper}\left(  \alpha\beta,\tau\right)  $. Then,
Proposition \ref{prop.iper.shuffle} yields that $\tau\in
\operatorname{Sh}_{p,q}$ if and only if $\operatorname*{iper}\left(  \alpha\beta
,\tau\right)  \in\operatorname{Sh}_{n,m}$. Since we know that
$\tau\in\operatorname{Sh}_{p,q}$, we can deduce from this that
$\operatorname*{iper}\left(  \alpha\beta,\tau\right)  \in
\operatorname{Sh}_{n,m}$, so that $\sigma=\operatorname*{iper}\left(  \alpha\beta
,\tau\right)  \in\operatorname{Sh}_{n,m}$.

It remains to prove that the word $u\underset{\sigma}{\shuffle}v$ is the
lexicographically highest element of the multiset $u\shuffle v$.

It is clear that the multiset $u\shuffle v$ has \textbf{some}
lexicographically highest element. This element has the form
$u\underset{\widetilde{\sigma}}{\shuffle}v$ for some $\widetilde{\sigma}%
\in\operatorname{Sh}_{n,m}$ (because any element of this multiset
has such a form). Consider this $\widetilde{\sigma}$. Theorem
\ref{thm.shuffle.lyndon.full}(b) (applied to $\widetilde{\sigma}$ instead of
$\sigma$) yields that there exists a unique permutation $\widetilde{\tau}%
\in\operatorname{Sh}_{p,q}$ satisfying $\left(  uv\right)  \left[
I_{\widetilde{\tau}\left(  1\right)  }\right]  \geq\left(  uv\right)  \left[
I_{\widetilde{\tau}\left(  2\right)  }\right]  \geq\cdots\geq\left(
uv\right)  \left[  I_{\widetilde{\tau}\left(  p+q\right)  }\right]  $ and
$\widetilde{\sigma}=\operatorname*{iper}\left(  \alpha\beta,\widetilde{\tau
}\right)  $. (What we call $\widetilde{\tau}$ here is what has been called
$\tau$ in Theorem \ref{thm.shuffle.lyndon.full}(b).)

Now, the chain of inequalities $\left(  uv\right)  \left[  I_{\widetilde{\tau
}\left(  1\right)  }\right]  \geq\left(  uv\right)  \left[  I_{\widetilde{\tau
}\left(  2\right)  }\right]  \geq\cdots\geq\left(  uv\right)  \left[
I_{\widetilde{\tau}\left(  p+q\right)  }\right]  $ shows that the list
$\left(  \left(  uv\right)  \left[  I_{\widetilde{\tau}\left(  1\right)
}\right]  ,\left(  uv\right)  \left[  I_{\widetilde{\tau}\left(  2\right)
}\right]  ,\ldots,\left(  uv\right)  \left[  I_{\widetilde{\tau}\left(
p+q\right)  }\right]  \right)  $ is the result of sorting the list $\left(
\left(  uv\right)  \left[  I_{1}\right]  ,\left(  uv\right)  \left[
I_{2}\right]  ,\ldots,\left(  uv\right)  \left[  I_{p+q}\right]  \right)  $ in
decreasing order. But the chain of inequalities
(\ref{pf.thm.shuffle.lyndon.full.a.1}) shows that the list \newline$\left(
\left(  uv\right)  \left[  I_{\tau\left(  1\right)  }\right]  ,\left(
uv\right)  \left[  I_{\tau\left(  2\right)  }\right]  ,\ldots,\left(
uv\right)  \left[  I_{\tau\left(  p+q\right)  }\right]  \right)  $ is the
result of sorting the same list $\left(  \left(  uv\right)  \left[
I_{1}\right]  ,\left(  uv\right)  \left[  I_{2}\right]  ,\ldots,\left(
uv\right)  \left[  I_{p+q}\right]  \right)  $ in decreasing order. So each of
the two lists $\left(  \left(  uv\right)  \left[  I_{\widetilde{\tau}\left(
1\right)  }\right]  ,\left(  uv\right)  \left[  I_{\widetilde{\tau}\left(
2\right)  }\right]  ,\ldots,\left(  uv\right)  \left[  I_{\widetilde{\tau
}\left(  p+q\right)  }\right]  \right)  $ and \newline$\left(  \left(
uv\right)  \left[  I_{\tau\left(  1\right)  }\right]  ,\left(  uv\right)
\left[  I_{\tau\left(  2\right)  }\right]  ,\ldots,\left(  uv\right)  \left[
I_{\tau\left(  p+q\right)  }\right]  \right)  $ is the result of sorting one
and the same list \newline$\left(  \left(  uv\right)  \left[  I_{1}\right]
,\left(  uv\right)  \left[  I_{2}\right]  ,\ldots,\left(  uv\right)  \left[
I_{p+q}\right]  \right)  $ in decreasing order. Since the result of sorting a
given list in decreasing order is unique, this yields%
\[
\left(  \left(  uv\right)  \left[  I_{\widetilde{\tau}\left(  1\right)
}\right]  ,\left(  uv\right)  \left[  I_{\widetilde{\tau}\left(  2\right)
}\right]  ,\ldots,\left(  uv\right)  \left[  I_{\widetilde{\tau}\left(
p+q\right)  }\right]  \right)  =\left(  \left(  uv\right)  \left[
I_{\tau\left(  1\right)  }\right]  ,\left(  uv\right)  \left[  I_{\tau\left(
2\right)  }\right]  ,\ldots,\left(  uv\right)  \left[  I_{\tau\left(
p+q\right)  }\right]  \right)  .
\]
Hence,%
\begin{equation}
\left(  uv\right)  \left[  I_{\widetilde{\tau}\left(  1\right)  }\right]
\cdot\left(  uv\right)  \left[  I_{\widetilde{\tau}\left(  2\right)  }\right]
\cdot\cdots\cdot\left(  uv\right)  \left[  I_{\widetilde{\tau}\left(
p+q\right)  }\right]  =\left(  uv\right)  \left[  I_{\tau\left(  1\right)
}\right]  \cdot\left(  uv\right)  \left[  I_{\tau\left(  2\right)  }\right]
\cdot\cdots\cdot\left(  uv\right)  \left[  I_{\tau\left(  p+q\right)
}\right]  . \label{pf.thm.shuffle.lyndon.full.a.5}%
\end{equation}
But Lemma \ref{lem.shuffle.lyndon.step2} yields
\begin{equation}
u\underset{\sigma}{\shuffle}v=\left(  uv\right)  \left[  I_{\tau\left(
1\right)  }\right]  \cdot\left(  uv\right)  \left[  I_{\tau\left(  2\right)
}\right]  \cdot\cdots\cdot\left(  uv\right)  \left[  I_{\tau\left(
p+q\right)  }\right]  . \label{pf.thm.shuffle.lyndon.full.a.6}%
\end{equation}
Meanwhile, Lemma \ref{lem.shuffle.lyndon.step2} (applied to $\widetilde{\tau}$
and $\widetilde{\sigma}$ instead of $\tau$ and $\sigma$) yields
\begin{align*}
u\underset{\widetilde{\sigma}}{\shuffle}v  &  =\left(  uv\right)  \left[
I_{\widetilde{\tau}\left(  1\right)  }\right]  \cdot\left(  uv\right)  \left[
I_{\widetilde{\tau}\left(  2\right)  }\right]  \cdot\cdots\cdot\left(
uv\right)  \left[  I_{\widetilde{\tau}\left(  p+q\right)  }\right] \\
&  =\left(  uv\right)  \left[  I_{\tau\left(  1\right)  }\right]  \cdot\left(
uv\right)  \left[  I_{\tau\left(  2\right)  }\right]  \cdot\cdots\cdot\left(
uv\right)  \left[  I_{\tau\left(  p+q\right)  }\right]
\ \ \ \ \ \ \ \ \ \ \left(  \text{by (\ref{pf.thm.shuffle.lyndon.full.a.5}%
)}\right) \\
&  =u\underset{\sigma}{\shuffle}v\ \ \ \ \ \ \ \ \ \ \left(  \text{by
(\ref{pf.thm.shuffle.lyndon.full.a.6})}\right)  .
\end{align*}
Thus, $u\underset{\sigma}{\shuffle}v$ is the lexicographically highest element
of the multiset $u\shuffle v$ (since we know that
$u\underset{\widetilde{\sigma}}{\shuffle}v$ is the lexicographically highest
element of the multiset $u\shuffle v$). This proves Theorem
\ref{thm.shuffle.lyndon.full}(a).
\end{proof}

Now, in order to prove Theorem \ref{thm.shuffle.lyndon.used}, we record a very
simple fact about counting shuffles:

\begin{proposition}
\label{prop.shuffle.decrmax}Let $p\in \NN$ and $q\in \NN$. Let
$\mathfrak{W}$ be a totally ordered set, and let $h:\left\{  1,2,\ldots
,p+q\right\}  \rightarrow\mathfrak{W}$ be a map. Assume that $h\left(
1\right)  \geq h\left(  2\right)  \geq\cdots\geq h\left(  p\right)  $ and
$h\left(  p+1\right)  \geq h\left(  p+2\right)  \geq\cdots\geq h\left(
p+q\right)  $.

For every $w\in\mathfrak{W}$, let $\mathfrak{a}\left(  w\right)  $ denote the
number of all $i\in\left\{  1,2,\ldots,p\right\}  $ satisfying $h\left(
i\right)  =w$, and let $\mathfrak{b}\left(  w\right)  $ denote the number of
all $i\in\left\{  p+1,p+2,\ldots,p+q\right\}  $ satisfying $h\left(  i\right)
=w$.

Then, the number of $\tau\in\operatorname{Sh}_{p,q}$ satisfying
$h\left(  \tau\left(  1\right)  \right)  \geq h\left(  \tau\left(  2\right)
\right)  \geq\cdots\geq h\left(  \tau\left(  p+q\right)  \right)  $ is
$\prod_{w\in\mathfrak{W}}\dbinom{\mathfrak{a}\left(  w\right)  +\mathfrak{b}%
\left(  w\right)  }{\mathfrak{a}\left(  w\right)  }$. (Of course, all but
finitely many factors of this product are $1$.)
\end{proposition}

\begin{exercise}
\label{exe.prop.shuffle.decrmax}Prove Proposition \ref{prop.shuffle.decrmax}.
\end{exercise}

\begin{proof}
[Proof of Theorem \ref{thm.shuffle.lyndon.used}.]Let $n=\ell\left(  u\right)
$ and $m=\ell\left(  v\right)  $. Define $\alpha$, $\beta$ and $\left(
I_{1},I_{2},\ldots,I_{p+q}\right)  $ as in Theorem
\ref{thm.shuffle.lyndon.full}.

Since $\left(  a_{1},a_{2},\ldots,a_{p}\right)  $ is the CFL factorization of
$u$, we have $a_{1}\geq a_{2}\geq\cdots\geq a_{p}$ and $a_{1}a_{2} \cdots
a_{p}=u$. Similarly, $b_{1}\geq b_{2}\geq\cdots\geq b_{q}$ and $b_{1}b_{2}
\cdots b_{q}=v$.

From (\ref{pf.thm.shuffle.lyndon.full.prep1}), we see that $\left(  uv\right)
\left[  I_{i}\right]  =a_{i}$ for every $i\in\left\{  1,2,\ldots,p\right\}  $.
From (\ref{pf.thm.shuffle.lyndon.full.prep1'}), we see that $\left(
uv\right)  \left[  I_{i}\right]  =b_{i-p}$ for every $i\in\left\{
p+1,p+2,\ldots,p+q\right\}  $. Combining these two equalities, we obtain%
\begin{equation}
\left(  uv\right)  \left[  I_{i}\right]  = \begin{cases}
a_{i},& \text{if }i\leq p;\\
b_{i-p},& \text{if }i>p
\end{cases}
\ \ \ \ \ \ \ \ \ \ \text{for every }i\in\left\{  1,2,\ldots
,p+q\right\}  .
\label{pf.thm.shuffle.lyndon.used.uvIi}
\end{equation}
In other words,
\begin{equation}
\left(  \left(  uv\right)  \left[  I_{1}\right]  ,\left(  uv\right)  \left[
I_{2}\right]  ,\ldots,\left(  uv\right)  \left[  I_{p+q}\right]  \right)
=\left(  a_{1},a_{2},\ldots,a_{p},b_{1},b_{2},\ldots,b_{q}\right)  .
\label{pf.thm.shuffle.lyndon.used.a.1}%
\end{equation}

(a) Let $z$ be the lexicographically highest element of the multiset
$u\shuffle v$. We must prove that $z=c_{1}c_{2} \cdots c_{p+q}$.

Since $z\in u\shuffle v$, we can write $z$ in the form $u\underset{\sigma
}{\shuffle}v$ for some $\sigma\in\operatorname{Sh}_{n,m}$ (since we
can write any element of $u\shuffle v$ in this form). Consider this $\sigma$.
Then, $u\underset{\sigma}{\shuffle}v=z$ is the lexicographically highest
element of the multiset $u\shuffle v$. Hence, Theorem
\ref{thm.shuffle.lyndon.full}(b) yields that there exists a unique permutation
$\tau\in\operatorname{Sh}_{p,q}$ satisfying $\left(  uv\right)
\left[  I_{\tau\left(  1\right)  }\right]  \geq\left(  uv\right)  \left[
I_{\tau\left(  2\right)  }\right]  \geq\cdots\geq\left(  uv\right)  \left[
I_{\tau\left(  p+q\right)  }\right]  $ and $\sigma=\operatorname*{iper}\left(
\alpha\beta,\tau\right)  $. Consider this $\tau$.

Now, $\tau\in\operatorname{Sh}_{p,q}\subset\Symm_{p+q}$ is a
permutation, and thus the list $\left(  \left(  uv\right)  \left[
I_{\tau\left(  1\right)  }\right]  ,\left(  uv\right)  \left[  I_{\tau\left(
2\right)  }\right]  ,\ldots,\left(  uv\right)  \left[  I_{\tau\left(
p+q\right)  }\right]  \right)  $ is a rearrangement of the list $\left(
\left(  uv\right)  \left[  I_{1}\right]  ,\left(  uv\right)  \left[
I_{2}\right]  ,\ldots,\left(  uv\right)  \left[  I_{p+q}\right]  \right)  $.
Due to (\ref{pf.thm.shuffle.lyndon.used.a.1}), this rewrites as follows: The
list $\left(  \left(  uv\right)  \left[  I_{\tau\left(  1\right)  }\right]
,\left(  uv\right)  \left[  I_{\tau\left(  2\right)  }\right]  ,\ldots,\left(
uv\right)  \left[  I_{\tau\left(  p+q\right)  }\right]  \right)  $ is a
rearrangement of the list $\left(  a_{1},a_{2},\ldots,a_{p},b_{1},b_{2}%
,\ldots,b_{q}\right)  $. Hence, $\left(  \left(  uv\right)  \left[
I_{\tau\left(  1\right)  }\right]  ,\left(  uv\right)  \left[  I_{\tau\left(
2\right)  }\right]  ,\ldots,\left(  uv\right)  \left[  I_{\tau\left(
p+q\right)  }\right]  \right)  $ is the result of sorting the list $\left(
a_{1},a_{2},\ldots,a_{p},b_{1},b_{2},\ldots,b_{q}\right)  $ in decreasing
order (since $\left(  uv\right)  \left[  I_{\tau\left(  1\right)  }\right]
\geq\left(  uv\right)  \left[  I_{\tau\left(  2\right)  }\right]  \geq
\cdots\geq\left(  uv\right)  \left[  I_{\tau\left(  p+q\right)  }\right]  $).
But since the result of sorting the list $\left(  a_{1},a_{2},\ldots
,a_{p},b_{1},b_{2},\ldots,b_{q}\right)  $ in decreasing order is $\left(
c_{1},c_{2},\ldots,c_{p+q}\right)  $, this becomes%
\[
\left(  \left(  uv\right)  \left[  I_{\tau\left(  1\right)  }\right]  ,\left(
uv\right)  \left[  I_{\tau\left(  2\right)  }\right]  ,\ldots,\left(
uv\right)  \left[  I_{\tau\left(  p+q\right)  }\right]  \right)  =\left(
c_{1},c_{2},\ldots,c_{p+q}\right)  .
\]
Hence,%
\[
\left(  uv\right)  \left[  I_{\tau\left(  1\right)  }\right]  \cdot\left(
uv\right)  \left[  I_{\tau\left(  2\right)  }\right]  \cdot\cdots\cdot\left(
uv\right)  \left[  I_{\tau\left(  p+q\right)  }\right]  =c_{1}\cdot c_{2}%
\cdot\cdots\cdot c_{p+q}.
\]
But Lemma \ref{lem.shuffle.lyndon.step2} yields
\[
u\underset{\sigma}{\shuffle}v=\left(  uv\right)  \left[  I_{\tau\left(
1\right)  }\right]  \cdot\left(  uv\right)  \left[  I_{\tau\left(  2\right)
}\right]  \cdot\cdots\cdot\left(  uv\right)  \left[  I_{\tau\left(
p+q\right)  }\right]  .
\]
Altogether, we have
\begin{align*}
z = u\underset{\sigma}{\shuffle}v=\left(  uv\right)  \left[  I_{\tau\left(
1\right)  }\right]  \cdot\left(  uv\right)  \left[  I_{\tau\left(  2\right)
}\right]  \cdot\cdots\cdot\left(  uv\right)  \left[  I_{\tau\left(
p+q\right)  }\right]
= c_{1}\cdot c_{2}\cdot\cdots\cdot c_{p+q}=c_{1}c_{2} \cdots c_{p+q}.
\end{align*}
This proves Theorem \ref{thm.shuffle.lyndon.used}(a).

(b) Recall that $u\shuffle v=\left\{  u\underset{\sigma}{\shuffle}
v\ :\ \sigma\in\operatorname{Sh}_{n,m}\right\}  _{\text{multiset}}$. Hence,
\begin{align*}
&  \left(  \text{the multiplicity with which the lexicographically highest
element of the multiset}\right. \\
&  \ \ \ \ \ \ \ \ \ \ \left.  u\shuffle v\text{ appears in the multiset
}u\shuffle v\right) \\
&  =\left(  \text{the number of all }\sigma\in
\operatorname{Sh}_{n,m}\text{ such that }u\underset{\sigma}{\shuffle}v\text{ is
the}\right. \\
&  \ \ \ \ \ \ \ \ \ \ \left.  \text{lexicographically highest element of the
multiset }u\shuffle v\right)  .
\end{align*}

However, for a given $\sigma\in\operatorname{Sh}_{n,m}$, we know that
$u\underset{\sigma}{\shuffle}v$ is the lexicographically highest element of
the multiset $u\shuffle v$ if and only if $\sigma$ can be written in the form
$\sigma=\operatorname*{iper}\left(  \alpha\beta,\tau\right)  $ for some
$\tau\in\operatorname{Sh}_{p,q}$ satisfying $\left(  uv\right)
\left[  I_{\tau\left(  1\right)  }\right]  \geq\left(  uv\right)  \left[
I_{\tau\left(  2\right)  }\right]  \geq\cdots\geq\left(  uv\right)  \left[
I_{\tau\left(  p+q\right)  }\right]  $.\ \ \ \ \footnote{In fact, the
``if'' part of this assertion follows from
Theorem \ref{thm.shuffle.lyndon.full}(a), whereas its ``only
if'' part follows from Theorem \ref{thm.shuffle.lyndon.full}(b).}
Hence,
\begin{align*}
&  \left(  \text{the number of all }\sigma\in\operatorname{Sh}_{n,m}%
\text{ such that }u\underset{\sigma}{\shuffle}v\text{ is the}\right. \\
&  \ \ \ \ \ \ \ \ \ \ \left.  \text{lexicographically highest element of the
multiset }u\shuffle v\right) \\
&  =\left(  \text{the number of all }\sigma\in
\operatorname{Sh}_{n,m}\text{ which can be written in the form }\sigma
=\operatorname*{iper}\left(  \alpha\beta,\tau\right)  \right. \\
&  \ \ \ \ \ \ \ \ \ \ \left.  \text{for some }\tau\in
\operatorname{Sh}_{p,q}\text{ satisfying }\left(  uv\right)  \left[  I_{\tau\left(
1\right)  }\right]  \geq\left(  uv\right)  \left[  I_{\tau\left(  2\right)
}\right]  \geq\cdots\geq\left(  uv\right)  \left[  I_{\tau\left(  p+q\right)
}\right]  \right) \\
&  =\left(  \text{the number of all }\tau\in\operatorname{Sh}_{p,q}%
\text{ satisfying }\left(  uv\right)  \left[  I_{\tau\left(  1\right)
}\right]  \geq\left(  uv\right)  \left[  I_{\tau\left(  2\right)  }\right]
\geq\cdots\geq\left(  uv\right)  \left[  I_{\tau\left(  p+q\right)  }\right]
\right)
\end{align*}
(because if a $\sigma\in\operatorname{Sh}_{n,m}$ can be written in
the form $\sigma=\operatorname*{iper}\left(  \alpha\beta,\tau\right)  $ for
some $\tau\in\operatorname{Sh}_{p,q}$ satisfying $\left(  uv\right)
\left[  I_{\tau\left(  1\right)  }\right]  \geq\left(  uv\right)  \left[
I_{\tau\left(  2\right)  }\right]  \geq\cdots\geq\left(  uv\right)  \left[
I_{\tau\left(  p+q\right)  }\right]  $, then $\sigma$ can be written
\textbf{uniquely} in this form\footnote{\textit{Proof.} Let $\sigma
\in\operatorname{Sh}_{n,m}$ be such that $\sigma$ can be written in
the form $\sigma=\operatorname*{iper}\left(  \alpha\beta,\tau\right)  $ for
some $\tau\in\operatorname{Sh}_{p,q}$ satisfying $\left(  uv\right)
\left[  I_{\tau\left(  1\right)  }\right]  \geq\left(  uv\right)  \left[
I_{\tau\left(  2\right)  }\right]  \geq\cdots\geq\left(  uv\right)  \left[
I_{\tau\left(  p+q\right)  }\right]  $. Then, the word $u\underset{\sigma
}{\shuffle}v$ is the lexicographically highest element of the multiset
$u\shuffle v$ (according to Theorem \ref{thm.shuffle.lyndon.full}(a)). Hence,
there exists a unique permutation $\tau\in\operatorname{Sh}_{p,q}$
satisfying $\left(  uv\right)  \left[  I_{\tau\left(  1\right)  }\right]
\geq\left(  uv\right)  \left[  I_{\tau\left(  2\right)  }\right]  \geq
\cdots\geq\left(  uv\right)  \left[  I_{\tau\left(  p+q\right)  }\right]  $
and $\sigma=\operatorname*{iper}\left(  \alpha\beta,\tau\right)  $ (according
to Theorem \ref{thm.shuffle.lyndon.full}(b)). In other words, $\sigma$ can be
written \textbf{uniquely} in the form $\sigma=\operatorname*{iper}\left(
\alpha\beta,\tau\right)  $ for some $\tau\in\operatorname{Sh}_{p,q}$
satisfying $\left(  uv\right)  \left[  I_{\tau\left(  1\right)  }\right]
\geq\left(  uv\right)  \left[  I_{\tau\left(  2\right)  }\right]  \geq
\cdots\geq\left(  uv\right)  \left[  I_{\tau\left(  p+q\right)  }\right]  $,
qed.}). Thus,%
\begin{align}
&  \left(  \text{the multiplicity with which the lexicographically highest
element of the multiset}\right. \nonumber\\
&  \ \ \ \ \ \ \ \ \ \ \left.  u\shuffle v\text{ appears in the multiset
}u\shuffle v\right) \nonumber\\
&  =\left(  \text{the number of all }\sigma\in
\operatorname{Sh}_{n,m}\text{ such that }u\underset{\sigma}{\shuffle}v\text{ is
the}\right. \nonumber\\
&  \ \ \ \ \ \ \ \ \ \ \left.  \text{lexicographically highest element of the
multiset }u\shuffle v\right) \nonumber\\
&  =\left(  \text{the number of all }\tau\in\operatorname{Sh}_{p,q}%
\text{ satisfying }\left(  uv\right)  \left[  I_{\tau\left(  1\right)
}\right]  \geq\left(  uv\right)  \left[  I_{\tau\left(  2\right)  }\right]
\geq\cdots\geq\left(  uv\right)  \left[  I_{\tau\left(  p+q\right)  }\right]
\right)  . \label{pf.thm.shuffle.lyndon.used.b.3}%
\end{align}

Now, define a map $h:\left\{  1,2,\ldots,p+q\right\}  \rightarrow\mathfrak{L}$
by%
\[
h\left(  i\right)  = \begin{cases}
a_{i},& \text{if }i\leq p;\\
b_{i-p},& \text{if }i>p
\end{cases}
\ \ \ \ \ \ \ \ \ \ \text{for every }i\in\left\{  1,2,\ldots
,p+q\right\}  .
\]
Then, $h\left(  1\right)  \geq h\left(  2\right)  \geq\cdots\geq h\left(
p\right)  $ (because this is just a rewriting of $a_{1}\geq a_{2}\geq
\cdots\geq a_{p}$) and $h\left(  p+1\right)  \geq h\left(  p+2\right)
\geq\cdots\geq h\left(  p+q\right)  $ (since this is just a rewriting of
$b_{1}\geq b_{2}\geq\cdots\geq b_{q}$). For every $w\in\mathfrak{L}$, the
number of all $i\in\left\{  1,2,\ldots,p\right\}  $ satisfying $h\left(
i\right)  =w$ is
\begin{align*}
&  \left\vert \left\{  i\in\left\{  1,2,\ldots,p\right\}  \ \mid
\ \underbrace{h\left(  i\right)  }_{=a_{i}}=w\right\}  \right\vert \\
&  =\left\vert \left\{  i\in\left\{  1,2,\ldots,p\right\}  \ \mid
\ a_{i}=w\right\}  \right\vert \\
&  =\left(  \text{the number of terms in the list }\left(  a_{1},a_{2}%
,\ldots,a_{p}\right)  \text{ which are equal to }w\right) \\
&  =\left(  \text{the number of terms in the CFL factorization of }u\text{
which are equal to }w\right) \\
&  \ \ \ \ \ \ \ \ \ \ \left(  \text{since the list }\left(  a_{1}%
,a_{2},\ldots,a_{p}\right)  \text{ is the CFL factorization of }u\right) \\
&  =\operatorname{mult}_{w}u
\end{align*}
(because $\operatorname{mult}_{w}u$ is defined as the number of
terms in the CFL factorization of $u$ which are equal to $w$). Similarly, for
every $w\in\mathfrak{L}$, the number of all $i\in\left\{  p+1,p+2,\ldots
,p+q\right\}  $ satisfying $h\left(  i\right)  =w$ equals
$\operatorname{mult}_{w}v$. Thus, we can apply Proposition
\ref{prop.shuffle.decrmax} to $\mathfrak{W}=\mathfrak{L}$, $\mathfrak{a}%
\left(  w\right)  =\operatorname{mult}_{w}u$ and $\mathfrak{b}%
\left(  w\right)  =\operatorname{mult}_{w}v$. As a result, we see
that the number of $\tau\in\operatorname{Sh}_{p,q}$ satisfying
$h\left(  \tau\left(  1\right)  \right)  \geq h\left(  \tau\left(  2\right)
\right)  \geq\cdots\geq h\left(  \tau\left(  p+q\right)  \right)  $ is
$\prod_{w\in\mathfrak{L}}\dbinom{\operatorname{mult}_{w}%
u+\operatorname{mult}_{w}v}{\operatorname{mult}_{w}u}$. In
other words,%
\begin{align}
&  \left(  \text{the number of all }\tau\in\operatorname{Sh}_{p,q}%
\text{ satisfying }h\left(  \tau\left(  1\right)  \right)  \geq h\left(
\tau\left(  2\right)  \right)  \geq\cdots\geq h\left(  \tau\left(  p+q\right)
\right)  \right) \nonumber\\
&  =\prod_{w\in\mathfrak{L}}\dbinom{\operatorname{mult}_{w}%
u+\operatorname{mult}_{w}v}{\operatorname{mult}_{w}u}.
\label{pf.thm.shuffle.lyndon.used.b.8}%
\end{align}

However, for every $i\in\left\{  1,2,\ldots,p+q\right\}  $, we have%
\[
h\left(  i\right)  =  \begin{cases}
a_{i},& \text{if }i\leq p;\\
b_{i-p},& \text{if }i>p
\end{cases}
\quad =\left(  uv\right)  \left[  I_{i}\right]  \ \ \ \ \ \ \ \ \ \ \left(
\text{by (\ref{pf.thm.shuffle.lyndon.used.uvIi})}\right)  .
\]
Hence, for any $\tau\in\operatorname{Sh}_{p,q}$, the condition
$h\left(  \tau\left(  1\right)  \right)  \geq h\left(  \tau\left(  2\right)
\right)  \geq\cdots\geq h\left(  \tau\left(  p+q\right)  \right)  $ is
equivalent to $\left(  uv\right)  \left[  I_{\tau\left(  1\right)  }\right]
\geq\left(  uv\right)  \left[  I_{\tau\left(  2\right)  }\right]  \geq
\cdots\geq\left(  uv\right)  \left[  I_{\tau\left(  p+q\right)  }\right]  $.
Thus,
\begin{align*}
&  \left(  \text{the number of all }\tau\in\operatorname{Sh}_{p,q}%
\text{ satisfying }\underbrace{h\left(  \tau\left(  1\right)  \right)  \geq
h\left(  \tau\left(  2\right)  \right)  \geq\cdots\geq h\left(  \tau\left(
p+q\right)  \right)  }_{\substack{\text{this is equivalent to}\\\left(
uv\right)  \left[  I_{\tau\left(  1\right)  }\right]  \geq\left(  uv\right)
\left[  I_{\tau\left(  2\right)  }\right]  \geq\cdots\geq\left(  uv\right)
\left[  I_{\tau\left(  p+q\right)  }\right]  }}\right) \\
&  =\left(  \text{the number of all }\tau\in\operatorname{Sh}_{p,q}%
\text{ satisfying }\left(  uv\right)  \left[  I_{\tau\left(  1\right)
}\right]  \geq\left(  uv\right)  \left[  I_{\tau\left(  2\right)  }\right]
\geq\cdots\geq\left(  uv\right)  \left[  I_{\tau\left(  p+q\right)  }\right]
\right) \\
&  =\left(  \text{the multiplicity with which the lexicographically highest
element of the multiset}\right. \\
&  \ \ \ \ \ \ \ \ \ \ \left.  u\shuffle v\text{ appears in the multiset
}u\shuffle v\right)
\end{align*}
(by (\ref{pf.thm.shuffle.lyndon.used.b.3})). Compared with
(\ref{pf.thm.shuffle.lyndon.used.b.8}), this yields%
\begin{align*}
&  \left(  \text{the multiplicity with which the lexicographically highest
element of the multiset}\right. \\
&  \ \ \ \ \ \ \ \ \ \ \left.  u\shuffle v\text{ appears in the multiset
}u\shuffle v\right) \\
&  =\prod_{w\in\mathfrak{L}}\dbinom{\operatorname{mult}_{w}%
u+\operatorname{mult}_{w}v}{\operatorname{mult}_{w}u}.
\end{align*}
This proves Theorem \ref{thm.shuffle.lyndon.used}(b).

(c) We shall use the notations of Theorem \ref{thm.shuffle.lyndon.used}(a) and
Theorem \ref{thm.shuffle.lyndon.used}(b).

Assume that $a_{i}\geq b_{j}$ for every $i\in\left\{  1,2,\ldots,p\right\}  $
and $j\in\left\{  1,2,\ldots,q\right\}  $. This, combined with $a_{1}\geq
a_{2}\geq\cdots\geq a_{p}$ and $b_{1}\geq b_{2}\geq\cdots\geq b_{q}$, yields
that $a_{1}\geq a_{2}\geq\cdots\geq a_{p}\geq b_{1}\geq b_{2}\geq\cdots\geq
b_{q}$. Thus, the list $\left(  a_{1},a_{2},\ldots,a_{p},b_{1},b_{2}%
,\ldots,b_{q}\right)  $ is weakly decreasing. Thus, the result of sorting the
list $\left(  a_{1},a_{2},\ldots,a_{p},b_{1},b_{2},\ldots,b_{q}\right)  $ in
decreasing order is the list $\left(  a_{1},a_{2},\ldots,a_{p},b_{1}%
,b_{2},\ldots,b_{q}\right)  $ itself. But since this result is $\left(
c_{1},c_{2},\ldots,c_{p+q}\right)  $, this shows that $\left(  c_{1}%
,c_{2},\ldots,c_{p+q}\right)  =\left(  a_{1},a_{2},\ldots,a_{p},b_{1}%
,b_{2},\ldots,b_{q}\right)  $. Hence, $c_{1}c_{2} \cdots c_{p+q}%
=\underbrace{a_{1}a_{2} \cdots a_{p}}_{=u}\underbrace{b_{1}b_{2} \cdots b_{q}%
}_{=v}=uv$. Now, Theorem \ref{thm.shuffle.lyndon.used}(a) yields that the
lexicographically highest element of the multiset $u\shuffle v$ is $c_{1}c_{2}
\cdots c_{p+q}=uv$. This proves Theorem \ref{thm.shuffle.lyndon.used}(c).

(d) We shall use the notations of Theorem \ref{thm.shuffle.lyndon.used}(a) and
Theorem \ref{thm.shuffle.lyndon.used}(b).

Assume that $a_{i}>b_{j}$ for every $i\in\left\{  1,2,\ldots,p\right\}  $ and
$j\in\left\{  1,2,\ldots,q\right\}  $. Thus, $a_{i}\geq b_{j}$ for every
$i\in\left\{  1,2,\ldots,p\right\}  $ and $j\in\left\{  1,2,\ldots,q\right\}
$. Hence, Theorem \ref{thm.shuffle.lyndon.used}(c) yields that the
lexicographically highest element of the multiset $u\shuffle v$ is $uv$.
Therefore, Theorem \ref{thm.shuffle.lyndon.used}(b) shows that the
multiplicity with which this word $uv$ appears in the multiset $u\shuffle v$
is $\prod_{w\in\mathfrak{L}}\dbinom{\operatorname{mult}_{w}%
u+\operatorname{mult}_{w}v}{\operatorname{mult}_{w}u}$.

Now, every $w\in\mathfrak{L}$ satisfies $\dbinom{
\operatorname{mult}_{w}u+\operatorname{mult}_{w}v}{
\operatorname{mult}_{w}u}=1$\ \ \ $\ $\footnote{\textit{Proof.} Assume the contrary.
Then, there exists at least one $w\in\mathfrak{L}$ such that $\dbinom
{\operatorname{mult}_{w}u+\operatorname{mult}_{w}%
v}{\operatorname{mult}_{w}u}\neq1$. Consider this $w$. Both
$\operatorname{mult}_{w}u$ and $\operatorname{mult}_{w}v$
must be positive (since $\dbinom{\operatorname{mult}_{w}%
u+\operatorname{mult}_{w}v}{\operatorname{mult}_{w}u}%
\neq1$). Since $\operatorname{mult}_{w}u$ is positive, there must be
at least one term in the CFL factorization of $u$ which is equal to $w$. In
other words, there is at least one $i\in\left\{  1,2,\ldots,p\right\}  $
satisfying $a_{i}=w$ (since $\left(  a_{1},a_{2},\ldots,a_{p}\right)  $ is the
CFL factorization of $u$). Similarly, there is at least one $j\in\left\{
1,2,\ldots,q\right\}  $ satisfying $b_{j}=w$. These $i$ and $j$ satisfy
$a_{i}=w=b_{j}$, which contradicts $a_{i}>b_{j}$. This contradiction shows
that our assumption was false, qed.}. Thus, as we know, the multiplicity with
which this word $uv$ appears in the multiset $u\shuffle v$ is $\prod
_{w\in\mathfrak{L}}\underbrace{\dbinom{\operatorname{mult}_{w}%
u+\operatorname{mult}_{w}v}{\operatorname{mult}_{w}u}%
}_{=1}=\prod_{w\in\mathfrak{L}}1=1$. This proves Theorem
\ref{thm.shuffle.lyndon.used}(d).

(e) We shall use the notations of Theorem \ref{thm.shuffle.lyndon.used}(a) and
Theorem \ref{thm.shuffle.lyndon.used}(b).

Since $u$ is a Lyndon word, the $1$-tuple $\left(  u\right)  $ is the CFL
factorization of $u$. Hence, we can apply Theorem
\ref{thm.shuffle.lyndon.used}(c) to $1$ and $\left(  u\right)  $ instead of
$p$ and $\left(  a_{1},a_{2},\ldots,a_{p}\right)  $. As a result, we conclude
that the lexicographically highest element of the multiset $u\shuffle v$ is
$uv$. It remains to prove that the multiplicity with which this word $uv$
appears in the multiset $u\shuffle v$ is $\operatorname{mult}_{u}%
v+1$.

For every $w\in\mathfrak{L}$ satisfying $w\neq u$, we have
\begin{equation}
\operatorname{mult}_{w}u=0 \label{pf.thm.shuffle.lyndon.used.e.0}%
\end{equation}
\footnote{\textit{Proof of (\ref{pf.thm.shuffle.lyndon.used.e.0}):} Let
$w\in\mathfrak{L}$ be such that $w\neq u$. Then, the number of terms in the
list $\left(  u\right)  $ which are equal to $w$ is $0$. Since $\left(
u\right)  $ is the CFL factorization of $u$, this rewrites as follows: The
number of terms in the CFL factorization of $u$ which are equal to $w$ is $0$.
In other words, $\operatorname{mult}_{w}u=0$. This proves
(\ref{pf.thm.shuffle.lyndon.used.e.0}).}. Also, $
\operatorname{mult}_{u}u=1$ (for a similar reason). But $uv$ is the lexicographically
highest element of the multiset $u\shuffle v$. Hence, the multiplicity with
which the word $uv$ appears in the multiset $u\shuffle v$ is the multiplicity
with which the lexicographically highest element of the multiset $u\shuffle v$
appears in the multiset $u\shuffle v$. According to Theorem
\ref{thm.shuffle.lyndon.used}(b), the latter multiplicity is%
\begin{align*}
&  \prod_{w\in\mathfrak{L}}\dbinom{\operatorname{mult}_{w}%
u+\operatorname{mult}_{w}v}{\operatorname{mult}_{w}u}\\
&  =\underbrace{\dbinom{\operatorname{mult}_{u}%
u+\operatorname{mult}_{u}v}{\operatorname{mult}_{u}u}%
}_{\substack{=\dbinom{1+\operatorname{mult}_{u}v}{1}\\\text{(since
}\operatorname{mult}_{u}u=1\text{)}}}\cdot\prod\limits_{\substack{w\in
\mathfrak{L;}\\w\neq u}}\underbrace{\dbinom{
\operatorname{mult}_{w}u+\operatorname{mult}_{w}v}{
\operatorname{mult}_{w}u}}_{\substack{=\dbinom{0+\operatorname{mult}_{w}v}%
{0}\\\text{(since }\operatorname{mult}_{w}u=0\text{ (by
(\ref{pf.thm.shuffle.lyndon.used.e.0})))}}}\ \ \ \ \ \ \ \ \ \ \left(
\text{since }u\in\mathfrak{L}\right) \\
&  =\underbrace{\dbinom{1+\operatorname{mult}_{u}v}{1}%
}_{=1+\operatorname{mult}_{u}v=\operatorname{mult}_{u}%
v+1}\cdot\prod\limits_{\substack{w\in\mathfrak{L;}\\w\neq u}}\underbrace{\dbinom
{0+\operatorname{mult}_{w}v}{0}}_{=1}=\left(  
\operatorname{mult}_{u}v+1\right)  \cdot\underbrace{\prod\limits_{\substack{w\in\mathfrak{L;}%
\\w\neq u}}1}_{=1}=\operatorname{mult}_{u}v+1.
\end{align*}
This proves Theorem \ref{thm.shuffle.lyndon.used}(e).
\end{proof}

As an application of our preceding results, we can prove a further necessary
and sufficient criterion for a word to be Lyndon; this criterion is due to
Chen/Fox/Lyndon \cite[$\mathfrak{A}^{\prime\prime}=\mathfrak{A}^{\prime
\prime\prime\prime}$]{ChenFoxLyndon}:

\begin{exercise}
\label{exe.shuffle.lyndon.cflcrit}Let $w\in\mathfrak{A}^{\ast}$ be a nonempty
word. Prove that $w$ is Lyndon if and only if for any two nonempty words
$u\in\mathfrak{A}^{\ast}$ and $v\in\mathfrak{A}^{\ast}$ satisfying $w=uv$,
there exists at least one $s\in u\shuffle v$ satisfying $s>w$.
\end{exercise}

\subsection{\label{subsect.shuffle.radford}Radford's theorem on the shuffle
algebra}

We recall that our goal in Chapter \ref{sect.QSym.lyndon} is to exhibit an
algebraically independent generating set of the $\kk$-algebra
$\Qsym$. Having the notion of Lyndon words -- which will, to
some extent, but not literally, parametrize this generating set -- in place,
we could start the construction of this generating set immediately. However,
it might come off as rather unmotivated this way, and so we begin with some
warmups. First, we shall prove Radford's theorem on the shuffle algebra.

\begin{definition}
\label{def.shuffle.polyalg}
A \dfn{polynomial algebra} will mean a $\kk$-algebra which is
isomorphic to the polynomial ring $\kk \left[  x_{i}\ \mid\ i\in
I\right]  $ as a $\kk$-algebra (for some indexing set $I$). Note that
$I$ need not be finite.

Equivalently, a polynomial algebra can be defined as a $\kk$-algebra
which has an algebraically independent (over $\kk$) generating set. Yet
equivalently, a polynomial algebra can be defined as a $\kk$-algebra
which is isomorphic to the symmetric algebra of a free $\kk$-module.
\end{definition}

Keep in mind that when we say that a certain bialgebra $A$ is a polynomial
algebra, we are making no statement about the coalgebra structure on $A$.
The isomorphism from $A$ to the symmetric algebra of a free $\kk$-module
need not be a coalgebra isomorphism, and the algebraically independent
generating set of $A$ need not consist of primitives. Thus, showing that
a bialgebra $A$ is a polynomial algebra does not trivialize the study of
its bialgebraic structure.

\begin{remark}
\label{rmk.shuffle.shufflealg}
Let $V$ be a $\kk$-module, and let
$\mathfrak{A}$ be a totally ordered set. Let $b_{a}$ be an element of $V$ for
every $a\in\mathfrak{A}$. Consider the shuffle algebra $\operatorname*{Sh}%
\left(  V\right)  $ (defined in Definition \ref{prop.shuffle-alg}).

For every word $w\in\mathfrak{A}^{\ast}$ over the alphabet $\mathfrak{A}$, let
us define an element $b_{w}$ of $\operatorname*{Sh}\left(  V\right)  $ by
$b_{w}=b_{w_{1}}b_{w_{2}} \cdots b_{w_{\ell}}$, where $\ell$ is the length of
$w$. (The multiplication used here is that of $T\left(  V\right)  $, not that
of $\operatorname*{Sh}\left(  V\right)  $; the latter is denoted by
$\shufmult$.)

Let $u\in\mathfrak{A}^{\ast}$ and $v\in\mathfrak{A}^{\ast}$ be two words over
the alphabet $\mathfrak{A}$. Let $n=\ell\left(  u\right)  $ and $m=\ell\left(
v\right)  $. Then,
\[
b_{u} \shufmult b_{v}=\sum_{\sigma\in
\operatorname{Sh}_{n,m}}b_{u\underset{\sigma}{\shuffle}v}.
\]

\end{remark}

\begin{exercise}
\label{exe.rmk.shuffle.shufflealg}Prove Remark \ref{rmk.shuffle.shufflealg}.

[\textbf{Hint:} This follows from the definition of $\shufmult$.]
\end{exercise}

We can now state Radford's theorem \cite[Theorem 3.1.1(e)]{Radford-shuffle}:

\begin{theorem}
\label{thm.shuffle.radford}Assume that $\QQ$ is a subring of
$\kk$. Let $V$ be a free $\kk$-module with a basis $\left(
b_{a}\right)  _{a\in\mathfrak{A}}$, where $\mathfrak{A}$ is a totally ordered
set. Then, the shuffle algebra $\operatorname*{Sh}\left(  V\right)  $ (defined
in Definition \ref{prop.shuffle-alg}) is a polynomial $\kk$-algebra. An
algebraically independent generating set of $\operatorname*{Sh}\left(
V\right)  $ can be constructed as follows:

For every word $w\in\mathfrak{A}^{\ast}$ over the alphabet $\mathfrak{A}$, let
us define an element $b_{w}$ of $\operatorname*{Sh}\left(  V\right)  $ by
$b_{w}=b_{w_{1}}b_{w_{2}} \cdots b_{w_{\ell}}$, where $\ell$ is the length of
$w$. (The multiplication used here is that of $T\left(  V\right)  $, not that
of $\operatorname*{Sh}\left(  V\right)  $; the latter is denoted by
$\shufmult$.) Let $\mathfrak{L}$ denote the set of all Lyndon words
over the alphabet $\mathfrak{A}$. Then, $\left(  b_{w}\right)  _{w\in
\mathfrak{L}}$ is an algebraically independent generating set of the
$\kk$-algebra $\operatorname*{Sh}\left(  V\right)  $.
\end{theorem}

\begin{example}
For this example, let $\mathfrak{A}$ be the alphabet $\left\{  1,2,3,\ldots
\right\}  $ with total order given by $1<2<3< \cdots$, and assume that
$\QQ$ is a subring of $\kk$. Let $V$ be the free $\kk$-module with basis
$\left(  b_{a}\right)  _{a\in\mathfrak{A}}$. We use the
notations of Theorem \ref{thm.shuffle.radford}. Then, Theorem
\ref{thm.shuffle.radford} yields that $\left(  b_{w}\right)  _{w\in
\mathfrak{L}}$ is an algebraically independent generating set of the
$\kk$-algebra $\operatorname*{Sh}\left(  V\right)  $. Here are some
examples of elements of $\operatorname*{Sh}\left(  V\right)  $ written as
polynomials in this generating set:
\begin{align*}
b_{12}  &  =b_{12}\ \ \ \ \ \ \ \ \ \ \left(  \text{the word }12\text{ itself
is Lyndon}\right)  ;\\
b_{21}  &  =b_{1} \shufmult b_{2}-b_{12};\\
b_{11}  &  =\dfrac{1}{2}b_{1} \shufmult b_{1};\\
b_{123}  &  =b_{123}\ \ \ \ \ \ \ \ \ \ \left(  \text{the word }123\text{
itself is Lyndon}\right)  ;\\
b_{132}  &  =b_{132}\ \ \ \ \ \ \ \ \ \ \left(  \text{the word }132\text{
itself is Lyndon}\right)  ;\\
b_{213}  &  =b_{2} \shufmult b_{13}-b_{123}-b_{132};\\
b_{231}  &  =b_{23} \shufmult b_{1}-b_{2} \shufmult
b_{13}+b_{132};\\
b_{312}  &  =b_{3} \shufmult b_{12}-b_{123}-b_{132};\\
b_{321}  &  =b_{1} \shufmult b_{2} \shufmult b_{3}%
-b_{23} \shufmult b_{1}-b_{3} \shufmult b_{12}+b_{123};\\
b_{112}  &  =b_{112}\ \ \ \ \ \ \ \ \ \ \left(  \text{the word }112\text{
itself is Lyndon}\right)  ;\\
b_{121}  &  =b_{12} \shufmult b_{1}-2b_{112};\\
b_{1212}  &  =\dfrac{1}{2}b_{12} \shufmult b_{12}-2b_{1122};\\
b_{4321}  &  =b_{1} \shufmult b_{2} \shufmult
b_{3} \shufmult b_{4}-b_{1} \shufmult b_{2}
\shufmult b_{34}-b_{1} \shufmult b_{23} \shufmult
b_{4}-b_{12} \shufmult b_{3} \shufmult b_{4}\\
&  \ \ \ \ \ \ \ \ \ \ +b_{1} \shufmult b_{234}+b_{12}
\shufmult b_{34}+b_{123} \shufmult b_{4}-b_{1234}.
\end{align*}
\footnote{A pattern emerges in the formulas for $b_{21}$, $b_{321}$ and
$b_{4321}$: for every $n\in \NN$, we have
\[
b_{\left(  n,n-1,\ldots,1\right)  }=\sum_{\alpha\in
\Comp_{n}}\left(  -1\right)  ^{n-\ell\left(  \alpha\right)  }%
b_{\mathbf{d}_{1}\left(  \alpha\right)  } \shufmult b_{\mathbf{d}%
_{2}\left(  \alpha\right)  } \shufmult \cdots \shufmult%
b_{\mathbf{d}_{\ell\left(  \alpha\right)  }\left(  \alpha\right)  },
\]
where $\left(  \mathbf{d}_{1}\left(  \alpha\right)  \right)  \cdot\left(
\mathbf{d}_{2}\left(  \alpha\right)  \right)  \cdot\cdots\cdot\left(
\mathbf{d}_{\ell\left(  \alpha\right)  }\left(  \alpha\right)  \right)  $ is
the factorization of the word $\left(  1,2,\ldots,n\right)  $ into factors of
length $\alpha_{1}$, $\alpha_{2}$, $\ldots$, $\alpha_{\ell}$ (where
$\alpha=\left(  \alpha_{1},\alpha_{2},\ldots,\alpha_{\ell}\right)  $). This
can be proved by an application of Lemma~\ref{lem.moebius-bool.comps}(a)
(as it is easy to see that for any composition $\alpha$ of $n$, we have
\begin{align*}
b_{\mathbf{d}_{1}\left(  \alpha\right)  } \shufmult b_{\mathbf{d}%
_{2}\left(  \alpha\right)  } \shufmult \cdots \shufmult%
b_{\mathbf{d}_{\ell\left(  \alpha\right)  }\left(  \alpha\right)  }
&= \left(\text{the sum of } b_\pi \text{ for all words }
     \pi \in \Symm_n \text{ satisfying }
     \operatorname{Des}\left(  \pi^{-1}\right)  \subset D\left(  \alpha\right)
   \right) \\
&= \sum\limits_{\substack{\beta \in \Comp_n; \\  \beta \text{ coarsens } \alpha}}
   \ \ %
   \sum\limits_{\substack{\pi \in \Symm_n; \\  \gamma\left(  \pi^{-1}\right) = \beta}}
   b_\pi ,
\end{align*}
where $\gamma\left(\pi^{-1}\right)$ denotes the composition $\tau$
of $n$ satisfying $D \left( \tau \right) = \operatorname{Des}
\left( \pi^{-1} \right)$).}
\end{example}

Note that Theorem \ref{thm.shuffle.radford} cannot survive without the
condition that $\QQ$ be a subring of $\kk$. For instance, for
any $v\in V$, we have $v \shufmult v=2vv$ in $\operatorname*{Sh}%
\left(  V\right)  $, which vanishes if $2=0$ in $\kk$; this stands in
contrast to the fact that polynomial $\kk$-algebras are integral
domains when $\kk$ itself is one. We will see that
$\Qsym$ is less sensitive towards the base ring in this regard
(although proving that $\Qsym$ is a polynomial algebra is much
easier when $\QQ$ is a subring of $\kk$).

\begin{remark}
Theorem~\ref{thm.shuffle.radford} can be contrasted with the
following fact:
If $\QQ$ is a subring of $\kk$,
then the shuffle algebra $\operatorname{Sh}\left(V\right)$
of \textbf{any} $\kk$-module $V$ (not necessarily free!)
is isomorphic (as a $\kk$-algebra) to the symmetric algebra
$\Sym \left( \left(\ker \epsilon\right)
             / \left(\ker \epsilon\right)^2 \right)$
(by Theorem~\ref{thm.leray.leray-e}(e), applied to
$A = \operatorname{Sh}\left(V\right)$).
This fact is closely related to
Theorem~\ref{thm.shuffle.radford}, but neither follows from
it (since Theorem~\ref{thm.shuffle.radford} only considers the
case of free $\kk$-modules $V$) nor yields it (since this
fact does not provide explicit generators for the
$\kk$-module
$\left(\ker \epsilon\right) / \left(\ker \epsilon\right)^2$
and thus for the $\kk$-algebra $\operatorname{Sh}\left(V\right)$).
\end{remark}

In our proof of Theorem \ref{thm.shuffle.radford} (but not only there), we
will use part (a) of the following lemma\footnote{And in a later proof, we
will also use its part (c) (which is tailored for application to
$\Qsym$).}, which makes proving that certain families indexed
by Lyndon words generate certain $\kk$-algebras more comfortable:

\begin{lemma}
\label{lem.lyndonbases}Let $A$ be a commutative $\kk$-algebra. Let
$\mathfrak{A}$ be a totally ordered set. Let $\mathfrak{L}$ be the set of all
Lyndon words over the alphabet $\mathfrak{A}$. Let $b_{w}$ be an element of
$A$ for every $w\in\mathfrak{L}$. For every word $u\in\mathfrak{A}^{\ast}$,
define an element $\mathbf{b}_{u}$ of $A$ by $\mathbf{b}_{u}=b_{a_{1}}%
b_{a_{2}} \cdots b_{a_{p}}$, where $\left(  a_{1},a_{2},\ldots,a_{p}\right)  $
is the CFL factorization of $u$.

\begin{itemize}
\item[(a)] The family $\left(  b_{w}\right)  _{w\in\mathfrak{L}}$ is an
algebraically independent generating set of the $\kk$-algebra $A$ if
and only if the family $\left(  \mathbf{b}_{u}\right)  _{u\in\mathfrak{A}%
^{\ast}}$ is a basis of the $\kk$-module $A$.

\item[(b)] The family $\left(  b_{w}\right)  _{w\in\mathfrak{L}}$ generates
the $\kk$-algebra $A$ if and only if the family $\left(  \mathbf{b}%
_{u}\right)  _{u\in\mathfrak{A}^{\ast}}$ spans the $\kk$-module $A$.

\item[(c)] Assume that the $\kk$-algebra $A$ is graded. Let
$\operatorname*{wt}:\mathfrak{A}\rightarrow\left\{  1,2,3,\ldots\right\}  $ be
any map such that for every $N\in\left\{  1,2,3,\ldots\right\}  $, the set
$\operatorname{wt}^{-1}\left(  N\right)  $ is finite.

For every word $w\in\mathfrak{A}^{\ast}$, define an element
$\operatorname*{Wt}\left(  w\right)  \in \NN$ by $\operatorname*{Wt}%
\left(  w\right)  =\operatorname*{wt}\left(  w_{1}\right)  +\operatorname*{wt}%
\left(  w_{2}\right)  + \cdots+\operatorname*{wt}\left(  w_{k}\right)  $,
where $k$ is the length of $w$.

Assume that for every $w\in\mathfrak{L}$, the element $b_{w}$ of $A$ is
homogeneous of degree $\operatorname*{Wt}\left(  w\right)  $.

Assume further that the $\kk$-module $A$ has a basis $\left(
g_{u}\right)  _{u\in\mathfrak{A}^{\ast}}$ having the property that for every
$u\in\mathfrak{A}^{\ast}$, the element $g_{u}$ of $A$ is homogeneous of degree
$\operatorname*{Wt}\left(  u\right)  $.

Assume also that the family $\left(  b_{w}\right)  _{w\in\mathfrak{L}}$
generates the $\kk$-algebra $A$.

Then, this family $\left(  b_{w}\right)  _{w\in\mathfrak{L}}$ is an
algebraically independent generating set of the $\kk$-algebra $A$.
\end{itemize}
\end{lemma}

\begin{exercise}
\label{exe.lem.lyndonbases}Prove Lemma \ref{lem.lyndonbases}.

[\textbf{Hint:} For (a) and (b), notice that the $\mathbf{b}_{u}$ are the
``monomials'' in the $b_{w}$. For (c), use
Exercise \ref{exe.linalg.sur-bij}(b) in every homogeneous component of $A$.]
\end{exercise}

The main workhorse of our proof of Theorem \ref{thm.shuffle.radford} will be
the following consequence of Theorem \ref{thm.shuffle.lyndon.used}(c):

\begin{proposition}
\label{prop.shuffle.radford.triang}Let $V$ be a free $\kk$-module with
a basis $\left(  b_{a}\right)  _{a\in\mathfrak{A}}$, where $\mathfrak{A}$ is a
totally ordered set.

For every word $w\in\mathfrak{A}^{\ast}$ over the alphabet $\mathfrak{A}$, let
us define an element $b_{w}$ of $\operatorname*{Sh}\left(  V\right)  $ by
$b_{w}=b_{w_{1}}b_{w_{2}} \cdots b_{w_{\ell}}$, where $\ell$ is the length of
$w$. (The multiplication used here is that of $T\left(  V\right)  $, not that
of $\operatorname*{Sh}\left(  V\right)  $; the latter is denoted by
$\shufmult$.)

For every word $u\in\mathfrak{A}^{\ast}$, define an element $\mathbf{b}_{u}$
by $\mathbf{b}_{u}=b_{a_{1}} \shufmult b_{a_{2}} \shufmult
\cdots \shufmult b_{a_{p}}$, where $\left(  a_{1},a_{2},\ldots
,a_{p}\right)  $ is the CFL factorization of $u$.

If $\ell\in \NN$ and if $x\in\mathfrak{A}^{\ell}$ is a word, then there
is a family $\left(  \eta_{x,y}\right)  _{y\in\mathfrak{A}^{\ell}}%
\in \NN ^{\mathfrak{A}^{\ell}}$ of elements of $\NN$ satisfying
\[
\mathbf{b}_{x}=\sum\limits_{\substack{y\in\mathfrak{A}^{\ell};\\y\leq x}}\eta
_{x,y}b_{y}
\]
and $\eta_{x,x}\neq0$ (in $\NN$).
\end{proposition}

Before we prove this, let us show a very simple lemma:

\begin{lemma}
\label{lem.shuffle.preserves-lex}Let $\mathfrak{A}$ be a totally ordered set.
Let $n\in \NN$ and $m\in \NN$. Let $\sigma\in
\operatorname{Sh}_{n,m}$.

(a) If $u$, $v$ and $v^{\prime}$ are three words satisfying $\ell\left(
u\right)  =n$, $\ell\left(  v\right)  =m$, $\ell\left(  v^{\prime}\right)  =m$
and $v^{\prime}<v$, then $u\underset{\sigma}{\shuffle}v^{\prime}%
<u\underset{\sigma}{\shuffle}v$.

(b) If $u$, $u^{\prime}$ and $v$ are three words satisfying $\ell\left(
u\right)  =n$, $\ell\left(  u^{\prime}\right)  =n$, $\ell\left(  v\right)  =m$
and $u^{\prime}<u$, then $u^{\prime}\underset{\sigma}{\shuffle}%
v<u\underset{\sigma}{\shuffle}v$.

(c) If $u$, $v$ and $v^{\prime}$ are three words satisfying $\ell\left(
u\right)  =n$, $\ell\left(  v\right)  =m$, $\ell\left(  v^{\prime}\right)  =m$
and $v^{\prime}\leq v$, then $u\underset{\sigma}{\shuffle}v^{\prime}\leq
u\underset{\sigma}{\shuffle}v$.
\end{lemma}

\begin{exercise}
\label{exe.shuffle.preserves-lex.pf}Prove Lemma
\ref{lem.shuffle.preserves-lex}.
\end{exercise}

\begin{exercise}
\label{exe.prop.shuffle.radford.triang} Prove Proposition
\ref{prop.shuffle.radford.triang}.

[\textbf{Hint:} Proceed by induction over $\ell$. In the induction step, apply
Theorem \ref{thm.shuffle.lyndon.used}(c)\footnote{Or Theorem
\ref{thm.shuffle.lyndon.used}(e), if you prefer.} to $u = a_{1}$ and $v =
a_{2} a_{3} \cdots a_{p}$, where $\left(  a_{1}, a_{2}, \ldots, a_{p}\right)
$ is the CFL factorization of $x$. Use Lemma \ref{lem.shuffle.preserves-lex}
to get rid of smaller terms.]
\end{exercise}

\begin{exercise}
\label{exe.thm.shuffle.radford} Prove Theorem \ref{thm.shuffle.radford}.

[\textbf{Hint:} According to Lemma \ref{lem.lyndonbases}(a), it suffices to
show that the family $\left(  \mathbf{b}_{u}\right)  _{u\in\mathfrak{A}^{\ast
}}$ defined in Proposition \ref{prop.shuffle.radford.triang} is a basis of the
$\kk$-module $\operatorname*{Sh}\left(  V\right)  $. When
$\mathfrak{A}$ is finite, the latter can be proven by triangularity using
Proposition \ref{prop.shuffle.radford.triang}. Reduce the general case to that
of finite $\mathfrak{A}$.]
\end{exercise}

\subsection{\label{subsect.shuffle.QSym.1}Polynomial freeness of
$\Qsym$: statement and easy parts}

\begin{definition}
\label{def.lyndon-composition}
For the rest of Section \ref{subsect.shuffle.QSym.1} and
for Section \ref{subsect.shuffle.QSym.2}, we introduce the following
notations: We let $\mathfrak{A}$ be the totally ordered set $\left\{
1,2,3,\ldots\right\}  $ with its natural order (that is, $1<2<3<\cdots$.)
Thus, the words over $\mathfrak{A}$ are precisely the compositions. That is,
$\mathfrak{A}^{\ast}= \Comp $. We let $\mathfrak{L}$ denote the
set of all Lyndon words over $\mathfrak{A}$. These Lyndon words are also
called \emph{Lyndon compositions}\index{Lyndon composition}.
\end{definition}

A natural question is how many Lyndon compositions of a given size exist.
While we will not use the answer, we nevertheless record it:

\begin{exercise}
\label{exe.words.lyndon.compositions.count}Show that the number of Lyndon
compositions of size $n$ equals
\[
\dfrac{1}{n}\sum\limits_{d\mid n}\mu\left(  d\right)  \left(  2^{n/d}%
-1\right)  =\dfrac{1}{n}\sum\limits_{d\mid n}\mu\left(  d\right)
2^{n/d}-\delta_{n,1}
\]
for every positive integer $n$ (where ``$\sum\limits_{d\mid n}$''
means a sum over all positive divisors of $n$, and where
$\mu$ is the number-theoretic M\"{o}bius function).

[\textbf{Hint:} One solution is similar to the solution of Exercise
\ref{exe.words.lyndon.count} using CFL factorization. Another proceeds by
defining a bijection between Lyndon compositions and Lyndon words over a
two-letter alphabet $\left\{  \mathbf{0},\mathbf{1}\right\}  $ (with
$\mathbf{0}<\mathbf{1}$) which are $\neq\mathbf{1}$.\ \ \ \ \footnote{This
bijection is obtained by restricting the bijection
\begin{align*}
\Comp  &  \rightarrow
\left\{ w \in \left\{  \mathbf{0},\mathbf{1}\right\} ^{\ast}
\mid w \text{ does not start with } \mathbf{1} \right\}
,\\
\left(  \alpha_{1},\alpha_{2},\ldots,\alpha_{\ell}\right)   &  \mapsto
\mathbf{01}^{\alpha_{1}-1}\mathbf{01}^{\alpha_{2}-1}\cdots\mathbf{01}%
^{\alpha_{\ell}-1}%
\end{align*}
(where $\mathbf{01}^{k}$ is to be read as $\mathbf{0}\left(  \mathbf{1}%
^{k}\right)  $, not as $\left(  \mathbf{01}\right)  ^{k}$) to the set of
Lyndon compositions. The idea behind this bijection is well-known in the
Grothendieck-Teichm\"{u}ller community: see, e.g., \cite[\S 3.1]{Henderson}
(and see \cite[Note 5.16]{GelfandKrobLascouxLeclercRetakhThibon} for a
different appearance of this idea).}]
\end{exercise}

Let us now state Hazewinkel's result (\cite[Theorem 8.1]{Hazewinkel},
\cite[\S 6.7]{HazewinkelGubareniKirichenko}) which is the main goal of Chapter
\ref{sect.QSym.lyndon}:

\begin{theorem}
\label{thm.QSym.lyndon}The $\kk$-algebra $\Qsym$ is a
polynomial algebra. It is isomorphic, as a graded $\kk$-algebra, to the
$\kk$-algebra $\kk \left[  x_{w}\ \mid\ w\in\mathfrak{L}\right]
$. Here, the grading on $\kk \left[  x_{w}\ \mid\ w\in\mathfrak{L}%
\right]  $ is defined by setting $\deg\left(  x_{w}\right)  =\sum_{i=1}%
^{\ell\left(  w\right)  }w_{i}$ for every $w\in\mathfrak{L}$.
\end{theorem}

We shall prove Theorem \ref{thm.QSym.lyndon} in the next section
(Section \ref{subsect.shuffle.QSym.2}). But the particular case of Theorem
\ref{thm.QSym.lyndon} when $\QQ$ is a subring of $\kk$ can be
proven more easily; we state it as a proposition:

\begin{proposition}
\label{prop.QSym.lyndon}Assume that $\QQ$ is a subring of $\kk$.
Then, Theorem \ref{thm.QSym.lyndon} holds.
\end{proposition}

We will give two proofs of Proposition \ref{prop.QSym.lyndon} in this
Section \ref{subsect.shuffle.QSym.1}; a third proof of Proposition
\ref{prop.QSym.lyndon} will immediately result from the proof of Theorem
\ref{thm.QSym.lyndon} in Section \ref{subsect.shuffle.QSym.2}. (There
\textbf{is} virtue in giving three different proofs, as they all construct
different isomorphisms $\kk \left[  x_{w}\ \mid\ w\in\mathfrak{L}%
\right]  \rightarrow\Qsym$.)

Our first proof -- originating in Malvenuto's \cite[Corollaire 4.20]%
{Malvenuto} -- can be given right away; it relies on Exercise \ref{exe.Psi}:

\begin{proof}
[First proof of Proposition \ref{prop.QSym.lyndon}.]Let $V$ be the free
$\kk$-module with basis $\left(  \mathfrak{b}_{n}\right)
_{n\in\left\{  1,2,3,\ldots\right\}  }$. Endow the $\kk$-module $V$
with a grading by assigning to each basis vector $\mathfrak{b}_{n}$ the degree
$n$. Exercise \ref{exe.Psi}(k) shows that $\Qsym$ is isomorphic
to the shuffle algebra $\operatorname{Sh}\left(  V\right)  $ (defined as in
Proposition~\ref{prop.shuffle-alg}) as Hopf algebras. By being a bit more
careful, we can obtain the slightly stronger result that $\Qsym$
is isomorphic to the shuffle algebra $\operatorname{Sh}\left(  V\right)  $ as
\textbf{graded} Hopf algebras\footnote{\textit{Proof.} In the solution of
Exercise \ref{exe.Psi}(k), we have shown that $\Qsym\cong
T\left(  V\right)  ^{o}$ as graded Hopf algebras. But Remark
\ref{rmk.shuffle-alg.dual}(b) shows that the Hopf algebra $T\left(  V\right)
^{o}$ is naturally isomorphic to the shuffle algebra $\operatorname*{Sh}%
\left(  V^{o}\right)  $ as Hopf algebras; it is easy to see that the natural
isomorphism $T\left(  V\right)  ^{o}\rightarrow\operatorname*{Sh}\left(
V^{o}\right)  $ is graded (because it is the direct sum of the isomorphisms
$\left(  V^{\otimes n}\right)  ^{o}\rightarrow\left(  V^{o}\right)  ^{\otimes
n}$ over all $n\in \NN$, and each of these isomorphisms is graded).
Hence, $T\left(  V\right)  ^{o}\cong\operatorname*{Sh}\left(  V^{o}\right)  $
as graded Hopf algebras. But $V^{o}\cong V$ as graded $\kk$-modules
(since $V$ is of finite type), and thus $\operatorname*{Sh}\left(
V^{o}\right)  \cong\operatorname*{Sh}\left(  V\right)  $ as graded Hopf
algebras. Altogether, we obtain $\Qsym\cong T\left(  V\right)
^{o}\cong\operatorname*{Sh}\left(  V^{o}\right)  \cong\operatorname*{Sh}%
\left(  V\right)  $ as graded Hopf algebras, qed.}. In particular,
$\Qsym\cong\operatorname*{Sh}\left(  V\right)  $ as graded
$\kk$-algebras.

Theorem \ref{thm.shuffle.radford} (applied to $b_{a}=\mathfrak{b}_{a}$) yields
that the shuffle algebra $\operatorname*{Sh}\left(  V\right)  $ is a
polynomial $\kk$-algebra, and that an algebraically independent
generating set of $\operatorname*{Sh}\left(  V\right)  $ can be constructed as follows:

For every word $w\in\mathfrak{A}^{\ast}$ over the alphabet $\mathfrak{A}$, let
us define an element $\mathfrak{b}_{w}$ of $\operatorname*{Sh}\left(
V\right)  $ by $\mathfrak{b}_{w}=\mathfrak{b}_{w_{1}}\mathfrak{b}_{w_{2}}
\cdots\mathfrak{b}_{w_{\ell}}$, where $\ell$ is the length of $w$. (The
multiplication used here is that of $T\left(  V\right)  $, not that of
$\operatorname*{Sh}\left(  V\right)  $; the latter is denoted by
$\shufmult$.) Then, $\left(  \mathfrak{b}_{w}\right)
_{w\in\mathfrak{L}}$ is an algebraically independent generating set of the
$\kk$-algebra $\operatorname*{Sh}\left(  V\right)  $.

For every $w\in\mathfrak{A}^{\ast}$, we have $\mathfrak{b}_{w}=\mathfrak{b}%
_{w_{1}}\mathfrak{b}_{w_{2}} \cdots\mathfrak{b}_{w_{\ell\left(  w\right)  }}$
(by the definition of $\mathfrak{b}_{w}$). For every $w\in\mathfrak{A}^{\ast}%
$, the element $\mathfrak{b}_{w}=\mathfrak{b}_{w_{1}}\mathfrak{b}_{w_{2}}
\cdots\mathfrak{b}_{w_{\ell\left(  w\right)  }}$ of $\operatorname*{Sh}\left(
V\right)  $ is homogeneous of degree $\sum_{i=1}^{\ell\left(  w\right)
}\underbrace{\deg\left(  \mathfrak{b}_{w_{i}}\right)  }_{=w_{i}}=
\sum_{i=1}^{\ell\left(  w\right)  }w_{i}$.

Now, define a grading on the $\kk$-algebra $\kk \left[
x_{w}\ \mid\ w\in\mathfrak{L}\right]  $ by setting $\deg\left(  x_{w}\right)
=\sum_{i=1}^{\ell\left(  w\right)  }w_{i}$ for every $w\in\mathfrak{L}$. By
the universal property of the polynomial algebra $\kk \left[
x_{w}\ \mid\ w\in\mathfrak{L}\right]  $, we can define a $\kk$-algebra
homomorphism $\Phi: \kk \left[  x_{w}\ \mid\ w\in\mathfrak{L}\right]
\rightarrow\operatorname*{Sh}\left(  V\right)  $ by setting%
\[
\Phi\left(  x_{w}\right)  =\mathfrak{b}_{w}\ \ \ \ \ \ \ \ \ \ \text{for every
}w\in\mathfrak{L}.
\]
This homomorphism $\Phi$ is a $\kk$-algebra isomorphism (since $\left(
\mathfrak{b}_{w}\right)  _{w\in\mathfrak{L}}$ is an algebraically independent
generating set of the $\kk$-algebra $\operatorname*{Sh}\left(
V\right)  $) and is graded (because for every $w\in\mathfrak{L}$, the element
$\mathfrak{b}_{w}$ of $\operatorname*{Sh}\left(  V\right)  $ is homogeneous of
degree $\sum_{i=1}^{\ell\left(  w\right)  }w_{i}=\deg\left(  x_{w}\right)  $).
Thus, $\Phi$ is an isomorphism of graded $\kk$-algebras. Hence,
$\operatorname*{Sh}\left(  V\right)  \cong \kk \left[  x_{w}\ \mid
\ w\in\mathfrak{L}\right]  $ as graded $\kk$-algebras. Altogether,
$\Qsym\cong\operatorname*{Sh}\left(  V\right)  \cong%
 \kk \left[  x_{w}\ \mid\ w\in\mathfrak{L}\right]  $ as graded
$\kk$-algebras. Thus, $\Qsym$ is a polynomial algebra.
This proves Theorem \ref{thm.QSym.lyndon} under the assumption that
$\QQ$ be a subring of $\kk$. In other words, this proves
Proposition \ref{prop.QSym.lyndon}.
\end{proof}

Our second proof of Proposition \ref{prop.QSym.lyndon} comes from
Hazewinkel/Gubareni/Kirichenko \cite{HazewinkelGubareniKirichenko} (where
Proposition \ref{prop.QSym.lyndon} appears as \cite[Theorem 6.5.13]%
{HazewinkelGubareniKirichenko}). This proof will construct an explicit
algebraically independent family generating the $\kk$-algebra
$\Qsym$.\ \ \ \ \footnote{We could, of course, obtain such a
family from our above proof as well (this is done by Malvenuto in
\cite[Corollaire 4.20]{Malvenuto}), but it won't be a very simple one.} The
generating set will be very unsophisticated: it will be $\left(  M_{\alpha
}\right)  _{\alpha\in\mathfrak{L}}$, where $\mathfrak{A}$ and $\mathfrak{L}$
are as in Theorem \ref{thm.QSym.lyndon}. Here, we are using the fact that
words over the alphabet $\left\{  1,2,3, \ldots\right\}$ are the same thing
as compositions, so, in particular, a monomial quasisymmetric function
$M_{\alpha}$ is defined for every such word $\alpha$.

It takes a bit of work to show that this family indeed fits the bill. We begin
with a corollary of Proposition \ref{QSym-multiplication} that is essentially
obtained by throwing away all non-bijective maps $f$:

\begin{proposition}
\label{prop.QSym.lyndon.MaMb}Let $\alpha\in\mathfrak{A}^{\ast}$ and $\beta
\in\mathfrak{A}^{\ast}$. Then,%
\begin{align*}
&  M_{\alpha}M_{\beta}\\
&  =\sum_{\gamma\in\alpha\shuffle\beta}M_{\gamma}+\left(  \text{a sum of terms
of the form }M_{\delta}\text{ with }\delta\in\mathfrak{A}^{\ast}\text{
satisfying }\ell\left(  \delta\right)  <\ell\left(  \alpha\right)
+\ell\left(  \beta\right)  \right)  .
\end{align*}
\footnote{The sum $\sum_{\gamma\in\alpha\shuffle\beta}M_{\gamma}$ ranges over
the \textbf{multiset} $\alpha\shuffle\beta$; if an element appears several
times in $\alpha\shuffle\beta$, then it has accordingly many addends
corresponding to it.}
\end{proposition}

\begin{exercise}
\label{exe.prop.QSym.lyndon.MaMb}Prove Proposition \ref{prop.QSym.lyndon.MaMb}.

[\textbf{Hint:} Recall what was said about the $p=\ell+m$ case in Example
\ref{exa.QSym-multiplication}.]
\end{exercise}

\begin{corollary}
\label{cor.QSym.lyndon.MaMb.leq}Let $\alpha\in\mathfrak{A}^{\ast}$ and
$\beta\in\mathfrak{A}^{\ast}$. Then, $M_{\alpha}M_{\beta}$ is a sum of terms
of the form $M_{\delta}$ with $\delta\in\mathfrak{A}^{\ast}$ satisfying
$\ell\left(  \delta\right)  \leq\ell\left(  \alpha\right)  +\ell\left(
\beta\right)  $.
\end{corollary}

\begin{exercise}
\label{exe.cor.QSym.lyndon.MaMb.leq}Prove Corollary
\ref{cor.QSym.lyndon.MaMb.leq}.
\end{exercise}

We now define a partial order on the compositions of a given nonnegative integer:

\begin{definition}
\label{def.wll}Let $n\in \NN$. We define a binary relation
$\underset{\operatorname*{wll}}{\leq}$ on the set $
\Comp_{n}$ as follows: For two compositions $\alpha$ and $\beta$ in
$\Comp_{n}$, we set $\alpha
\underset{\operatorname*{wll}}{\leq}\beta$ if and only if
\[
\text{either }\ell\left(  \alpha\right)  <\ell\left(  \beta\right)  \text{ or
}\left(  \ell\left(  \alpha\right)  =\ell\left(  \beta\right)  \text{ and
}\alpha\leq\beta\text{ in lexicographic order}\right)  \text{.}%
\]
This binary relation $\underset{\operatorname*{wll}}{\leq}$ is the
smaller-or-equal relation of a total order on $
\Comp_{n}$; we refer to said total order as the \dfn{wll-order} on
$\Comp_{n}$, and we denote by
$\underset{\operatorname*{wll}}{<}$ the smaller relation of this total order.
\end{definition}

Notice that if $\alpha$ and $\beta$ are two compositions satisfying
$\ell\left(  \alpha\right)  =\ell\left(  \beta\right)  $, then $\alpha
\leq\beta$ in lexicographic order if and only if $\alpha\leq\beta$ with
respect to the relation $\leq$ defined in Definition \ref{def.words}.

A remark about the name ``wll-order'' is in
order. We have taken this notation from \cite[Definition 6.7.14]{Hazewinkel},
where it is used for an extension of this order to the whole set
$\Comp$. We will never use this extension, as we will only ever
compare two compositions of the same integer.\footnote{In \cite[Definition
6.7.14]{Hazewinkel}, the name ``wll-order'' is introduced as an
abbreviation for ``\textbf{w}eight first,
then \textbf{l}ength, then \textbf{l}exicographic'' (in the
sense that two compositions are first compared by their weights, then, if the
weights are equal, by their lengths, and finally, if the lengths are also
equal, by the lexicographic order). For us, the alternative explanation
``\textbf{w}ord \textbf{l}ength, then \textbf{l}exicographic'' serves
just as well.}

We now state a fact which is similar (and plays a similar role) to Proposition
\ref{prop.shuffle.radford.triang}:

\begin{proposition}
\label{prop.QSym.lyndon.triang}For every composition $u\in \Comp %
=\mathfrak{A}^{\ast}$, define an element $\mathbf{M}_{u}\in
\Qsym$ by $\mathbf{M}_{u}=M_{a_{1}}M_{a_{2}} \cdots M_{a_{p}}$,
where $\left(  a_{1},a_{2},\ldots,a_{p}\right)  $ is the CFL factorization of
the word $u$.

If $n\in \NN$ and if $x\in\Comp_{n}$, then there
is a family $\left(  \eta_{x,y}\right)  _{y\in
\Comp_{n}}\in \NN ^{\Comp_{n}}$ of elements
of $\NN$ satisfying
\[
\mathbf{M}_{x}=\sum\limits_{\substack{y\in\Comp_{n};\\
y\underset{\operatorname*{wll}}{\leq}x}}\eta_{x,y}M_{y}
\]
and $\eta_{x,x}\neq0$ (in $\NN$).
\end{proposition}

Before we prove it, let us show the following lemma:

\begin{lemma}
\label{lem.QSym.lyndon.triang.wll-lemma}Let $n\in \NN$ and
$m\in \NN$. Let $u\in\Comp_{n}$ and
$v\in\Comp_{m}$. Let $z$ be the lexicographically
highest element of the multiset $u\shuffle v$.

(a) We have $z\in\Comp_{n+m}$.

(b) There exists a positive integer $h$ such that
\[
M_{u}M_{v}=hM_{z}+\left(  \text{a sum of terms of the form }M_{w}\text{ with
}w\in\Comp_{n+m}\text{ satisfying }%
w\underset{\operatorname*{wll}}{<}z\right)  .
\]

(c) Let $v^{\prime}\in\Comp_{m}$ be such that
$v^{\prime}\underset{\operatorname*{wll}}{<}v$. Then,%
\[
M_{u}M_{v^{\prime}}=\left(  \text{a sum of terms of the form }M_{w}\text{ with
}w\in\Comp_{n+m}\text{ satisfying }%
w\underset{\operatorname*{wll}}{<}z\right)  .
\]

\end{lemma}

\begin{exercise}
\label{exe.lem.QSym.lyndon.triang.wll-lemma}Prove Lemma
\ref{lem.QSym.lyndon.triang.wll-lemma}.

[\textbf{Hint:} For (b), set $h$ to be the multiplicity
with which the word $z$ appears in the multiset $u\shuffle v$,
then use Proposition \ref{prop.QSym.lyndon.MaMb} and
notice that $M_{u}M_{v}$ is homogeneous of degree $n+m$.
For (c), use (b) for
$v^{\prime}$ instead of $v$ and notice that Lemma
\ref{lem.shuffle.preserves-lex}(a) shows that the lexicographically highest
element of the multiset $u\shuffle v^{\prime}$ is
$\underset{\operatorname*{wll}}{<}z$.]
\end{exercise}

\begin{exercise}
\label{exe.prop.QSym.lyndon.triang}Prove Proposition
\ref{prop.QSym.lyndon.triang}.

[\textbf{Hint:} Proceed by strong induction over $n$. In the induction step,
let $\left(  a_{1},a_{2},\ldots,a_{p}\right)  $ be the CFL factorization of
$x$, and set $u=a_{1}$ and $v=a_{2}a_{3}\cdots a_{p}$; then apply Proposition
\ref{prop.QSym.lyndon.triang} to $v$ instead of $x$, and multiply the
resulting equality $\mathbf{M}_{v}=\sum\limits_{\substack{y\in
\Comp_{\left\vert v\right\vert };\\y\underset{\operatorname*{wll}}{\leq}%
v}}\eta_{v,y}M_{y}$ with $M_{u}$ to obtain an expression for $M_{u}%
\mathbf{M}_{v}=\mathbf{M}_{x}$. Use Lemma
\ref{lem.QSym.lyndon.triang.wll-lemma} to show that this expression has the
form $\sum\limits_{\substack{y\in\Comp_{n}%
;\\y\underset{\operatorname*{wll}}{\leq}x}}\eta_{x,y}M_{y}$ with $\eta
_{x,x}\neq0$; here it helps to remember that the lexicographically highest
element of the multiset $u\shuffle v$ is $uv=x$ (by Theorem
\ref{thm.shuffle.lyndon.used}(c)).]
\end{exercise}

We are almost ready to give our second proof of Proposition
\ref{prop.QSym.lyndon}; our last step is the following proposition:

\begin{proposition}
\label{prop.QSym.lyndon.M-genset}Assume that $\QQ$ is a subring of
$\kk$. Then, $\left(  M_{w}\right)  _{w\in\mathfrak{L}}$ is an
algebraically independent generating set of the $\kk$-algebra
$\Qsym$.
\end{proposition}

\begin{exercise}
\label{exe.prop.QSym.lyndon.M-genset}Prove Proposition
\ref{prop.QSym.lyndon.M-genset}.

[\textbf{Hint:} Define $\mathbf{M}_{u}$ for every $u\in \Comp $
as in Proposition \ref{prop.QSym.lyndon.triang}. Conclude from Proposition
\ref{prop.QSym.lyndon.triang} that, for every $n\in \NN$, the family
$\left(  \mathbf{M}_{u}\right)  _{u\in\Comp_{n}}$
expands invertibly triangularly\footnote{See
Definition~\ref{def.STmat.expansion-tria}(b) for the
meaning of this.} (with respect to the total order
$\underset{\operatorname*{wll}}{\leq}$ on $\Comp_{n}$)
with respect to the basis $\left(  M_{u}\right)  _{u\in
\Comp_{n}}$ of $\Qsym_{n}$. Conclude that this family $\left(
\mathbf{M}_{u}\right)  _{u\in\Comp_{n}}$ is a basis of
$\Qsym_{n}$ itself, and so the whole family $\left(
\mathbf{M}_{u}\right)  _{u\in \Comp }$ is a basis of
$\Qsym$. Conclude using Lemma \ref{lem.lyndonbases}(a).]
\end{exercise}

\begin{proof}
[Second proof of Proposition \ref{prop.QSym.lyndon}.]Proposition
\ref{prop.QSym.lyndon.M-genset} yields that $\left(  M_{w}\right)
_{w\in\mathfrak{L}}$ is an algebraically independent generating set of the
$\kk$-algebra $\Qsym$.

Define a grading on the $\kk$-algebra $\kk \left[  x_{w}%
\ \mid\ w\in\mathfrak{L}\right]  $ by setting $\deg\left(  x_{w}\right)
=\sum_{i=1}^{\ell\left(  w\right)  }w_{i}$ for every $w\in\mathfrak{L}$. By
the universal property of the polynomial algebra $\kk \left[
x_{w}\ \mid\ w\in\mathfrak{L}\right]  $, we can define a $\kk$-algebra
homomorphism $\Phi: \kk \left[  x_{w}\ \mid\ w\in\mathfrak{L}\right]
\rightarrow\Qsym$ by setting%
\[
\Phi\left(  x_{w}\right)  =M_{w}\ \ \ \ \ \ \ \ \ \ \text{for every }%
w\in\mathfrak{L}.
\]
This homomorphism $\Phi$ is a $\kk$-algebra isomorphism (since $\left(
M_{w}\right)  _{w\in\mathfrak{L}}$ is an algebraically independent generating
set of the $\kk$-algebra $\Qsym$) and is graded (because
for every $w\in\mathfrak{L}$, the element $M_{w}$ of $\Qsym$ is
homogeneous of degree $\left\vert w\right\vert =\sum_{i=1}^{\ell\left(
w\right)  }w_{i}=\deg\left(  x_{w}\right)  $). Thus, $\Phi$ is an isomorphism
of graded $\kk$-algebras. Hence, $\Qsym\cong%
 \kk \left[  x_{w}\ \mid\ w\in\mathfrak{L}\right]  $ as graded
$\kk$-algebras. In particular, this shows that $\Qsym$
is a polynomial algebra. This proves Theorem \ref{thm.QSym.lyndon} under the
assumption that $\QQ$ be a subring of $\kk$. Proposition
\ref{prop.QSym.lyndon} is thus proven again.
\end{proof}

\subsection{\label{subsect.shuffle.QSym.2}Polynomial freeness of
$\Qsym$: the general case}

We now will prepare for proving Theorem \ref{thm.QSym.lyndon} without any
assumptions on $\kk$. In our proof, we follow \cite{Hazewinkel} and
\cite[\S 6.7]{HazewinkelGubareniKirichenko}, but without using the language of
plethysm and Frobenius maps. We start with the following definition:

\begin{definition}
\label{def.QSym.Malphas}Let $\alpha$ be a composition. Write $\alpha$ in the
form $\alpha=\left(  \alpha_{1},\alpha_{2},\ldots,\alpha_{\ell}\right)  $ with
$\ell=\ell\left(  \alpha\right)  $.

(a) Let $\operatorname*{SIS}\left(  \ell\right)  $ denote the set of all
strictly increasing $\ell$-tuples $\left(  i_{1},i_{2},\ldots,i_{\ell}\right)
$ of positive integers.\footnote{``Strictly
increasing'' means that $i_{1}<i_{2}<\cdots<i_{\ell}$ here.
Of course, the elements of $\operatorname*{SIS}\left(  \ell\right)  $ are in
1-to-1 correspondence with $\ell$-element subsets of $\left\{  1,2,3,\ldots
\right\}  $.} For every $\ell$-tuple $\mathbf{i}=\left(  i_{1},i_{2}%
,\ldots,i_{\ell}\right)  \in\operatorname*{SIS}\left(  \ell\right)  $, we
denote the monomial $x_{i_{1}}^{\alpha_{1}}x_{i_{2}}^{\alpha_{2}}\cdots
x_{i_{\ell}}^{\alpha_{\ell}}$ by $\xx_{\mathbf{i}}^{\alpha}$. This
$\xx_{\mathbf{i}}^{\alpha}$ is a monomial of degree $\alpha_{1}%
+\alpha_{2}+\cdots+\alpha_{\ell}=\left\vert \alpha\right\vert $. Then,
\begin{equation}
M_{\alpha}=\sum_{\mathbf{i}\in\operatorname*{SIS}\left(  \ell\right)
}\xx_{\mathbf{i}}^{\alpha}. \label{eq.def.QSym.Malphas.a.1}%
\end{equation}
\footnote{\textit{Proof of (\ref{eq.def.QSym.Malphas.a.1}):} By the definition
of $M_{\alpha}$, we have
\begin{align*}
M_{\alpha}  &  =\underbrace{\sum_{i_{1}<i_{2}<\cdots<i_{\ell}\text{ in
}\left\{  1,2,3,\ldots\right\}  }}_{=\sum_{\left(  i_{1},i_{2},\ldots,i_{\ell
}\right)  \in\operatorname*{SIS}\left(  \ell\right)  }}x_{i_{1}}^{\alpha_{1}%
}x_{i_{2}}^{\alpha_{2}}\cdots x_{i_{\ell}}^{\alpha_{\ell}}=\sum_{\left(
i_{1},i_{2},\ldots,i_{\ell}\right)  \in\operatorname*{SIS}\left(  \ell\right)
}x_{i_{1}}^{\alpha_{1}}x_{i_{2}}^{\alpha_{2}}\cdots x_{i_{\ell}}^{\alpha
_{\ell}}=\sum_{\mathbf{i}=\left(  i_{1},i_{2},\ldots,i_{\ell}\right)
\in\operatorname*{SIS}\left(  \ell\right)  }\underbrace{x_{i_{1}}^{\alpha_{1}%
}x_{i_{2}}^{\alpha_{2}}\cdots x_{i_{\ell}}^{\alpha_{\ell}}}%
_{\substack{=\xx_{\mathbf{i}}^{\alpha}\\\text{(by the definition of
}\xx_{\mathbf{i}}^{\alpha}\text{)}}}\\
&  =\sum_{\mathbf{i}=\left(  i_{1},i_{2},\ldots,i_{\ell}\right)
\in\operatorname*{SIS}\left(  \ell\right)  }\xx_{\mathbf{i}}^{\alpha
}=\sum_{\mathbf{i}\in\operatorname*{SIS}\left(  \ell\right)  }
\xx_{\mathbf{i}}^{\alpha},
\end{align*}
qed.}

(b) Consider the ring $\kk \left[  \left[  \xx\right]  \right]  $
endowed with the coefficientwise topology\footnote{This topology is defined as
follows:
\par
We endow the ring $\kk$ with the discrete topology. Then, we can regard
the $\kk$-module $\kk \left[  \left[  \xx\right]  \right]  $ as a direct
product of infinitely many copies of $\kk$ (by identifying every power
series in $\kk \left[  \left[  \xx\right]  \right]  $ with the
family of its coefficients). Hence, the product topology is a well-defined
topology on $\kk \left[  \left[  \xx\right]  \right]  $; this
topology is denoted as the \dfn{coefficientwise topology}. A sequence
$\left(  a_{n}\right)  _{n\in \NN }$ of power series converges to a power
series $a$ with respect to this topology if and only if for every monomial
$\mathfrak{m}$, all sufficiently high $n\in \NN$ satisfy
\[
\left(  \text{the coefficient of }\mathfrak{m}\text{ in }a_{n}\right)
=\left(  \text{the coefficient of }\mathfrak{m}\text{ in }a\right)  .
\]
Note that this is \textbf{not} the topology obtained by taking the completion
of $\kk \left[  x_{1},x_{2},x_{3},\ldots\right]  $ with respect to the
standard grading (in which all $x_{i}$ have degree $1$). (The latter
completion is actually a smaller ring than $\kk\left[\left[\xx\right]\right]$.)}.
The family $\left(
\xx_{\mathbf{i}}^{\alpha}\right)  _{\mathbf{i}\in\operatorname*{SIS}%
\left(  \ell\right)  }$ of elements of $\kk \left[  \left[
\xx\right]  \right]  $ is power-summable\footnote{Let us define what
``power-summable'' means for us:
\par
A family $\left(  n_{\mathbf{i}}\right)  _{\mathbf{i}\in\mathbf{I}}%
\in \NN ^{\mathbf{I}}$ (where $\mathbf{I}$ is some set) is said to be
\dfn{finitely supported} if all but finitely many $\mathbf{i}\in\mathbf{I}$
satisfy $n_{\mathbf{i}}=0$.
\par
If $\left(  n_{\mathbf{i}}\right)  _{\mathbf{i}\in\mathbf{I}}\in
 \NN ^{\mathbf{I}}$ is a finitely supported family, then
$\sum_{\mathbf{i}\in\mathbf{I}}n_{\mathbf{i}}$ is a well-defined element of
$\NN$. If $N\in \NN$, then a family $\left(  n_{\mathbf{i}%
}\right)  _{\mathbf{i}\in\mathbf{I}}\in \NN ^{\mathbf{I}}$ will be called
\emph{$\left(  \leq N\right)  $-supported} if it is finitely supported and
satisfies $\sum_{\mathbf{i}\in\mathbf{I}}n_{\mathbf{i}}\leq N$.
\par
We say that a family $\left(  s_{\mathbf{i}}\right)  _{\mathbf{i}\in
\mathbf{I}}\in R^{\mathbf{I}}$ of elements of a topological commutative
$\kk$-algebra $R$ is \emph{power-summable}\index{power-summable family}
if it satisfies the following property: For every $N\in \NN$, the sum
\[
\sum\limits_{\substack{\left(  n_{\mathbf{i}}\right)  _{\mathbf{i}\in\mathbf{I}}%
\in \NN ^{\mathbf{I}};\\\left(  n_{\mathbf{i}}\right)  _{\mathbf{i}%
\in\mathbf{I}}\text{ is }\left(  \leq N\right)  \text{-supported }}%
}\alpha_{\left(  n_{\mathbf{i}}\right)  _{\mathbf{i}\in\mathbf{I}}}%
\prod_{\mathbf{i}\in\mathbf{I}}s_{\mathbf{i}}^{n_{\mathbf{i}}}%
\]
converges in the topology on $R$ for every choice of scalars $\alpha_{\left(
n_{\mathbf{i}}\right)  _{\mathbf{i}\in\mathbf{I}}}\in \kk$ corresponding
to all $\left(  \leq N\right)  $-supported $\left(  n_{\mathbf{i}}\right)
_{\mathbf{i}\in\mathbf{I}}\in \NN ^{\mathbf{I}}$. In our specific case,
we consider $\kk \left[  \left[  \xx\right]  \right]  $ as a
topological commutative $\kk$-algebra,
where the topology is the coefficientwise
topology. The fact that the family $\left(  \xx_{\mathbf{i}}^{\alpha
}\right)  _{\mathbf{i}\in\operatorname*{SIS}\left(  \ell\right)  }$ is
power-summable then can be proven as follows:
\par
\begin{itemize}
\item If $\alpha\neq\varnothing$, then this fact follows from the
(easily-verified) observation that every given monomial in the variables
$x_{1},x_{2},x_{3},\ldots$ can be written as a product of monomials of the
form $\xx_{\mathbf{i}}^{\alpha}$ (with $\mathbf{i}\in
\operatorname*{SIS}\left(  \ell\right)  $) in only finitely many ways.
\par
\item If $\alpha=\varnothing$, then this fact follows by noticing that
$\left(  \xx_{\mathbf{i}}^{\alpha}\right)  _{\mathbf{i}\in
\operatorname*{SIS}\left(  \ell\right)  }$ is a finite family (indeed,
$\operatorname*{SIS}\left(  \ell\right)  =\operatorname*{SIS}\left(  0\right)
=\left\{  \left(  \right)  \right\}  $), and every finite family is
power-summable.
\end{itemize}
}. Hence, for every $f\in\Lambda$, there is a well-defined power series
$f\left(  \left(  \xx_{\mathbf{i}}^{\alpha}\right)  _{\mathbf{i}%
\in\operatorname*{SIS}\left(  \ell\right)  }\right)  \in \kk \left[
\left[  \xx\right]  \right]  $ obtained by ``evaluating'' $f$
at $\left(  \xx_{\mathbf{i}}^{\alpha
}\right)  _{\mathbf{i}\in\operatorname*{SIS}\left(  \ell\right)  }%
$\ \ \ \ \footnote{Here is how this power series $f\left(  \left(
\xx_{\mathbf{i}}^{\alpha}\right)  _{\mathbf{i}\in\operatorname*{SIS}%
\left(  \ell\right)  }\right)  $ is formally defined:
\par
Let $R$ be any topological commutative $\kk$-algebra, and let $\left(
s_{\mathbf{i}}\right)  _{\mathbf{i}\in\mathbf{I}}\in R^{\mathbf{I}}$ be any
power-summable family of elements of $R$. Assume that the indexing set
$\mathbf{I}$ is countably infinite, and fix a bijection $\mathfrak{j}:\left\{
1,2,3,\ldots\right\}  \rightarrow\mathbf{I}$. Let $g\in R\left(
\xx\right)  $ be arbitrary. Then, we can substitute $s_{\mathfrak{j}%
\left(  1\right)  }$, $s_{\mathfrak{j}\left(  2\right)  }$, $s_{\mathfrak{j}%
\left(  3\right)  }$, $\ldots$ for the variables $x_{1}$, $x_{2}$, $x_{3}$,
$\ldots$ in $g$, thus obtaining an infinite sum which converges in $R$ (in
fact, its convergence follows from the fact that the family $\left(
s_{\mathbf{i}}\right)  _{\mathbf{i}\in\mathbf{I}}\in R^{\mathbf{I}}$ is
power-summable). The value of this sum will be denoted by $g\left(  \left(
s_{\mathbf{i}}\right)  _{\mathbf{i}\in\mathbf{I}}\right)  $. In general, this
value depends on the choice of the bijection $\mathfrak{j}$, so the notation
$g\left(  \left(  s_{\mathbf{i}}\right)  _{\mathbf{i}\in\mathbf{I}}\right)  $
is unambiguous only if this bijection $\mathfrak{j}$ is chosen once and for
all. However, when $g\in\Lambda$, one can easily see that the choice of
$\mathfrak{j}$ has no effect on $g\left(  \left(  s_{\mathbf{i}}\right)
_{\mathbf{i}\in\mathbf{I}}\right)  $.
\par
We can still define $g\left(  \left(  s_{\mathbf{i}}\right)  _{\mathbf{i}%
\in\mathbf{I}}\right)  $ when the set $\mathbf{I}$ is finite instead of being
countably infinite. In this case, we only need to modify our above definition
as follows: Instead of fixing a bijection $\mathfrak{j}:\left\{
1,2,3,\ldots\right\}  \rightarrow\mathbf{I}$, we now fix a bijection
$\mathfrak{j}:\left\{  1,2,\ldots,\left\vert \mathbf{I}\right\vert \right\}
\rightarrow\mathbf{I}$, and instead of substituting $s_{\mathfrak{j}\left(
1\right)  }$, $s_{\mathfrak{j}\left(  2\right)  }$, $s_{\mathfrak{j}\left(
3\right)  }$, $\ldots$ for the variables $x_{1}$, $x_{2}$, $x_{3}$, $\ldots$
in $g$, we now substitute $s_{\mathfrak{j}\left(  1\right)  }$,
$s_{\mathfrak{j}\left(  2\right)  }$, $\ldots$, $s_{\mathfrak{j}\left(
\left\vert \mathbf{I}\right\vert \right)  }$, $0$, $0$, $0$, $\ldots$ for the
variables $x_{1}$, $x_{2}$, $x_{3}$, $\ldots$ in $g$. Again, the same
observations hold as before: $g\left(  \left(  s_{\mathbf{i}}\right)
_{\mathbf{i}\in\mathbf{I}}\right)  $ is independent on $\mathfrak{j}$ if
$g\in\Lambda$.
\par
Hence, $g\left(  \left(  s_{\mathbf{i}}\right)  _{\mathbf{i}\in\mathbf{I}%
}\right)  $ is well-defined for every $g\in R\left(  \xx\right)  $,
every countable (i.e., finite or countably infinite) set $\mathbf{I}$, every
topological commutative $\kk$-algebra $R$ and every power-summable
family $\left(  s_{\mathbf{i}}\right)  _{\mathbf{i}\in\mathbf{I}}\in
R^{\mathbf{I}}$ of elements of $R$, as long as a bijection $\mathfrak{j}$ is
chosen. In particular, we can apply this to $g=f$, $\mathbf{I}%
=\operatorname*{SIS}\left(  \ell\right)  $, $R= \kk \left[  \left[
\xx\right]  \right]  $ and $\left(  s_{\mathbf{i}}\right)
_{\mathbf{i}\in\mathbf{I}}=\left(  \xx_{\mathbf{i}}^{\alpha}\right)
_{\mathbf{i}\in\operatorname*{SIS}\left(  \ell\right)  }$, choosing
$\mathfrak{j}$ to be the bijection which sends every positive integer $k$ to
the $k$-th smallest element of $\operatorname*{SIS}\left(  \ell\right)  $ in
the lexicographic order. (Of course, since $f\in\Lambda$, the choice of
$\mathfrak{j}$ is irrelevant.)}. In particular, for every $s\in \ZZ$, we
can evaluate the symmetric function $e_s \in \Lambda$\ \ \ \ \footnote{Recall
that $e_0 = 1$, and that $e_s = 0$ for $s<0$.} at $\left(
\xx_{\mathbf{i}}^{\alpha}\right)  _{\mathbf{i}\in\operatorname*{SIS}
\left(  \ell\right)  }$. The resulting power series $e_{s}\left(  \left(
\xx_{\mathbf{i}}^{\alpha}\right)  _{\mathbf{i}\in\operatorname*{SIS}%
\left(  \ell\right)  }\right)  \in \kk \left[  \left[  \xx\right]
\right]  $ will be denoted $M_{\alpha}^{\left\langle s\right\rangle }$. Thus,%
\[
M_{\alpha}^{\left\langle s\right\rangle }
= e_{s}\left(  \left(  \xx_{\mathbf{i}}^{\alpha}
\right)  _{\mathbf{i}\in\operatorname*{SIS}\left(
\ell\right)  }\right)  .
\]

\end{definition}

The power series $M_{\alpha}^{\left\langle s\right\rangle }$ are the power
series $e_{s}\left(  \alpha\right)  $ in \cite{HazewinkelGubareniKirichenko}.
We will shortly (in Corollary \ref{cor.QSym.Malphas.QSym}(a)) see that
$M_{\alpha}^{\left\langle s\right\rangle }\in\Qsym$ (although
this is also easy to prove by inspection). Here are some examples of
$M_{\alpha}^{\left\langle s\right\rangle }$:

\begin{example}
If $\alpha$ is a composition and $\ell$ denotes its length
$\ell\left(\alpha\right)$, then
\[
M_{\alpha}^{\left\langle 0\right\rangle }=\underbrace{e_{0}}_{=1}\left(
\left(  \xx_{\mathbf{i}}^{\alpha}\right)  _{\mathbf{i}\in
\operatorname*{SIS}\left(  \ell\right)  }\right)  =1\left(  \left(
\xx_{\mathbf{i}}^{\alpha}\right)  _{\mathbf{i}\in\operatorname*{SIS}%
\left(  \ell\right)  }\right)  =1
\]
and%
\[
M_{\alpha}^{\left\langle 1\right\rangle }=e_{1}\left(  \left(
\xx_{\mathbf{i}}^{\alpha}\right)  _{\mathbf{i}\in\operatorname*{SIS}\left(
\ell\right)  }\right)  =\sum_{\mathbf{i}\in\operatorname*{SIS}\left(
\ell\right)  }\xx_{\mathbf{i}}^{\alpha}=M_{\alpha}%
\ \ \ \ \ \ \ \ \ \ \left(  \text{by (\ref{eq.def.QSym.Malphas.a.1})}\right)
\]
and\footnote{This is not completely
obvious, but easy to check (see Exercise~\ref{exe.QSym.Malphas.ps}(b)).}
\[
M_{\alpha}^{\left\langle 2\right\rangle }=e_{2}\left(  \left(
\xx_{\mathbf{i}}^{\alpha}\right)  _{\mathbf{i}\in\operatorname*{SIS}
\left(  \ell\right)  }\right)
=\sum\limits_{\substack{\mathbf{i}\in\operatorname*{SIS}%
\left(  \ell\right)  ,\ \mathbf{j}\in\operatorname*{SIS}\left(  \ell\right)
;\\\mathbf{i}<\mathbf{j}}}\xx_{\mathbf{i}}^{\alpha}
\xx_{\mathbf{j}}^{\alpha}
\]
(where the notation ``$\mathbf{i}<\mathbf{j}$'' should be interpreted
with respect to an arbitrary but fixed total order on
the set $\operatorname*{SIS}\left(  \ell\right)  $ -- for example, the
lexicographic order). Applying the last of these three equalities to
$\alpha=\left(  2,1\right)  $, we obtain
\begin{align*}
M_{\left(  2,1\right)  }^{\left\langle 2\right\rangle }  &  =
\sum_{\substack{\mathbf{i}\in\operatorname*{SIS}\left(  2\right)  ,
\ \mathbf{j} \in\operatorname*{SIS}\left(  2\right)  ,\\
\mathbf{i}<\mathbf{j}}}
\xx_{\mathbf{i}}^{\left(  2,1\right)  }\xx_{\mathbf{j}%
}^{\left(  2,1\right)  }=\sum\limits_{\substack{\left(  i_{1},i_{2}\right)
\in\operatorname*{SIS}\left(  2\right)  ,\ \left(  j_{1},j_{2}\right)
\in\operatorname*{SIS}\left(  2\right)  ;\\\left(  i_{1},i_{2}\right)
<\left(  j_{1},j_{2}\right)  }}\underbrace{\xx_{\left(  i_{1}%
,i_{2}\right)  }^{\left(  2,1\right)  }}_{=x_{i_{1}}^{2}x_{i_{2}}^{1}%
}\underbrace{\xx_{\left(  j_{1},j_{2}\right)  }^{\left(  2,1\right)  }%
}_{=x_{j_{1}}^{2}x_{j_{2}}^{1}}\\
&  =\sum\limits_{\substack{\left(  i_{1},i_{2}\right)  \in\operatorname*{SIS}\left(
2\right)  ,\ \left(  j_{1},j_{2}\right)  \in\operatorname*{SIS}\left(
2\right)  ;\\\left(  i_{1},i_{2}\right)  <\left(  j_{1},j_{2}\right)
}}x_{i_{1}}^{2}x_{i_{2}}^{1}x_{j_{1}}^{2}x_{j_{2}}^{1}\\
&  =\underbrace{\sum\limits_{\substack{i_{1}<i_{2};\ j_{1}<j_{2};\\i_{1}<j_{1}%
}}x_{i_{1}}^{2}x_{i_{2}}^{1}x_{j_{1}}^{2}x_{j_{2}}^{1}}_{=M_{\left(
2,1,2,1\right)  }+M_{\left(  2,3,1\right)  }+2M_{\left(  2,2,1,1\right)
}+M_{\left(  2,2,2\right)  }}+\underbrace{\sum\limits_{\substack{i_{1}<i_{2}%
;\ j_{1}<j_{2};\\i_{1}=j_{1};\ i_{2}<j_{2}}}x_{i_{1}}^{2}x_{i_{2}}^{1}%
x_{j_{1}}^{2}x_{j_{2}}^{1}}_{=M_{\left(  4,1,1\right)  }}\\
&  \ \ \ \ \ \ \ \ \ \ \left(  \text{here, we have WLOG assumed that the order
on }\operatorname*{SIS}\left(  2\right)  \text{ is lexicographic}\right) \\
&  =M_{\left(  2,1,2,1\right)  }+M_{\left(  2,3,1\right)  }+2M_{\left(
2,2,1,1\right)  }+M_{\left(  2,2,2\right)  }+M_{\left(  4,1,1\right)  }.
\end{align*}

Of course, every negative integer $s$ satisfies $M_{\alpha}^{\left\langle
s\right\rangle }=\underbrace{e_{s}}_{=0}\left(  \left(
\xx_{\mathbf{i}}^{\alpha}\right)  _{\mathbf{i}\in\operatorname*{SIS}
\left(  \ell\right)  }\right)  =0$.
\end{example}

There is a determinantal formula for the $s!M_{\alpha}^{\left\langle
s\right\rangle }$ (and thus also for $M_{\alpha}^{\left\langle s\right\rangle
}$ when $s!$ is invertible in $\kk$), but in order to state it, we need
to introduce one more notation:

\begin{definition}
Let $\alpha=\left(  \alpha_{1},\alpha_{2},\ldots,\alpha_{\ell}\right)  $ be a
composition, and let $k$ be a positive integer. Then, $\alpha\left\{
k\right\}  $ will denote the composition $\left(  k\alpha_{1},k\alpha
_{2},\ldots,k\alpha_{\ell}\right)  $. Clearly, $\ell\left(  \alpha\left\{
k\right\}  \right)  =\ell\left(  \alpha\right)  $ and $\left\vert
\alpha\left\{  k\right\}  \right\vert =k\left\vert \alpha\right\vert $.
\end{definition}

\begin{exercise}
\label{exe.QSym.Malphas.ps}Let $\alpha$ be a composition. Write the
composition $\alpha$ in the form $\alpha=\left(  \alpha_{1},\alpha_{2}%
,\ldots,\alpha_{\ell}\right)  $ with $\ell=\ell\left(  \alpha\right)  $.

\begin{enumerate}
\item[(a)] Show that the $s$-th power-sum symmetric function $p_{s}\in\Lambda$
satisfies%
\[
p_{s}\left(  \left(  \xx_{\mathbf{i}}^{\alpha}\right)  _{\mathbf{i}%
\in\operatorname*{SIS}\left(  \ell\right)  }\right)  =M_{\alpha\left\{
s\right\}  }%
\]
for every positive integer $s$.

\item[(b)] Let us fix a total order on the set $\operatorname*{SIS}\left(
\ell\right)  $ (for example, the lexicographic order). Show that the $s$-th
elementary symmetric function $e_{s}\in\Lambda$ satisfies%
\[
M_{\alpha}^{\left\langle s\right\rangle }=e_{s}\left(  \left(  \xx%
_{\mathbf{i}}^{\alpha}\right)  _{\mathbf{i}\in\operatorname*{SIS}\left(
\ell\right)  }\right)  =\sum\limits_{\substack{\left(  \mathbf{i}_{1},\mathbf{i}%
_{2},\ldots,\mathbf{i}_{s}\right)  \in\left(  \operatorname*{SIS}\left(
\ell\right)  \right)  ^{s};\\\mathbf{i}_{1}<\mathbf{i}_{2}<\cdots
<\mathbf{i}_{s}}}\xx_{\mathbf{i}_{1}}^{\alpha}\xx_{\mathbf{i}%
_{2}}^{\alpha}\cdots\xx_{\mathbf{i}_{s}}^{\alpha}%
\]
for every $s\in \NN$.

\item[(c)] Let $s\in \NN$, and let $n$ be a positive integer. Let
$e_{s}^{\left\langle n\right\rangle }$ be the symmetric function
$\sum_{i_{1}<i_{2}<\cdots<i_{s}}
x_{i_{1}}^{n}x_{i_{2}}^{n}\cdots x_{i_{s}}^{n} \in\Lambda$.
Then, show that
\[
M_{\alpha\left\{  n\right\}  }^{\left\langle s\right\rangle }=e_{s}%
^{\left\langle n\right\rangle }\left(  \left(  \xx_{\mathbf{i}}%
^{\alpha}\right)  _{\mathbf{i}\in\operatorname*{SIS}\left(  \ell\right)
}\right)  .
\]

\item[(d)] Let $s\in \NN$, and let $n$ be a positive integer. Prove that
there exists a polynomial $P\in \kk \left[  z_{1},z_{2},z_{3}%
,\ldots\right]  $ such that $M_{\alpha\left\{  n\right\}  }^{\left\langle
s\right\rangle }=P\left(  M_{\alpha}^{\left\langle 1\right\rangle },M_{\alpha
}^{\left\langle 2\right\rangle },M_{\alpha}^{\left\langle 3\right\rangle
},\ldots\right)  $.
\end{enumerate}

[\textbf{Hint:} For (a), (b) and (c), apply the definition of $f\left(
\left(  \xx_{\mathbf{i}}^{\alpha}\right)  _{\mathbf{i}\in
\operatorname*{SIS}\left(  \ell\right)  }\right)  $ with $f$ a symmetric
function\footnote{There are two subtleties that need to be addressed:
\par
\begin{itemize}
\item the fact that the definition of $f\left(  \left(  \xx%
_{\mathbf{i}}^{\alpha}\right)  _{\mathbf{i}\in\operatorname*{SIS}\left(
\ell\right)  }\right)  $ distinguishes between two cases depending on whether
or not $\operatorname*{SIS}\left(  \ell\right)  $ is finite;
\par
\item the fact that the total order on the set $\left\{  1,2,3,\ldots\right\}
$ (which appears in the summation subscript in the equality $e_{s}%
=\sum\limits_{\substack{\left(  i_{1},i_{2},\ldots,i_{s}\right)  \in\left\{
1,2,3,\ldots\right\}  ^{s};\\i_{1}<i_{2}<\cdots<i_{s}}}x_{i_{1}}x_{i_{2}%
}\cdots x_{i_{s}}$) has nothing to do with the total order on the set
$\operatorname*{SIS}\left(  \ell\right)  $ (which appears in the summation
subscript in $\sum\limits_{\substack{\left(  \mathbf{i}_{1},\mathbf{i}_{2}%
,\ldots,\mathbf{i}_{s}\right)  \in\left(  \operatorname*{SIS}\left(
\ell\right)  \right)  ^{s};\\\mathbf{i}_{1}<\mathbf{i}_{2}<\cdots
<\mathbf{i}_{s}}}\xx_{\mathbf{i}_{1}}^{\alpha}\xx_{\mathbf{i}%
_{2}}^{\alpha}\cdots\xx_{\mathbf{i}_{s}}^{\alpha}$). For instance, the
former total order is well-founded, whereas the latter may and may not be. So
there is (generally) no bijection between $\left\{  1,2,3,\ldots\right\}  $
and $\operatorname*{SIS}\left(  \ell\right)  $ preserving these orders (even
if $\operatorname*{SIS}\left(  \ell\right)  $ is infinite). Fortunately, this
does not matter much, because the total order is only being used to ensure
that every product of $s$ distinct elements appears exactly once in
the sum.
\end{itemize}
}. For (d), recall that $\Lambda$ is generated by $e_{1},e_{2},e_{3},\ldots$.]
\end{exercise}

\begin{exercise}
\label{exe.QSym.Malphas.M1s}Let $s\in \NN$. Show that the composition
$\left(  1\right)  $ satisfies $M_{\left(  1\right)  }^{\left\langle
s\right\rangle }=e_{s}$.
\end{exercise}

\begin{proposition}
\label{prop.QSym.Malphas.dets}Let $\alpha=\left(  \alpha_{1},\alpha_{2}%
,\ldots,\alpha_{\ell}\right)  $ be a composition.

(a) Let $n\in \NN$. Define a matrix $A_{n}^{\left\langle \alpha
\right\rangle }=\left(  a_{i,j}^{\left\langle \alpha\right\rangle }\right)
_{i,j=1,2,\ldots,n}$ by
\[
a_{i,j}^{\left\langle \alpha\right\rangle }=
\begin{cases}
M_{\alpha\left\{  i-j+1\right\}  }, & \text{if }i\geq j;\\
i, & \text{if }i=j-1;\\
0, & \text{if }i<j-1
\end{cases}
\qquad\qquad\text{ for all }\left(  i,j\right)  \in\left\{  1,2,\ldots
,n\right\}  ^{2}.
\]
This matrix $A_{n}^{\left\langle \alpha\right\rangle }$ looks as follows:
\[
A_{n}^{\left\langle \alpha\right\rangle }=\left(
\begin{array}[c]{cccccc}
M_{\alpha\left\{  1\right\}  } & 1 & 0 & \cdots & 0 & 0\\
M_{\alpha\left\{  2\right\}  } & M_{\alpha\left\{  1\right\}  } & 2 & \cdots &
0 & 0\\
M_{\alpha\left\{  3\right\}  } & M_{\alpha\left\{  2\right\}  } &
M_{\alpha\left\{  1\right\}  } & \cdots & 0 & 0\\
\vdots & \vdots & \vdots & \ddots & \vdots & \vdots\\
M_{\alpha\left\{  n-1\right\}  } & M_{\alpha\left\{  n-2\right\}  } &
M_{\alpha\left\{  n-3\right\}  } & \cdots & M_{\alpha\left\{  1\right\}  } &
n-1\\
M_{\alpha\left\{  n\right\}  } & M_{\alpha\left\{  n-1\right\}  } &
M_{\alpha\left\{  n-2\right\}  } & \cdots & M_{\alpha\left\{  2\right\}  } &
M_{\alpha\left\{  1\right\}  }%
\end{array}
\right)  .
\]
Then, $\det\left(  A_{n}^{\left\langle \alpha\right\rangle }\right)
=n!M_{\alpha}^{\left\langle n\right\rangle }$.

(b) Let $n$ be a positive integer. Define a matrix $B_{n}^{\left\langle
\alpha\right\rangle }=\left(  b_{i,j}^{\left\langle \alpha\right\rangle
}\right)  _{i,j=1,2,\ldots,n}$ by
\[
b_{i,j}^{\left\langle \alpha\right\rangle }=
\begin{cases}
iM_{\alpha}^{\left\langle i\right\rangle }, & \text{if }j=1;\\
M_{\alpha}^{\left\langle i-j+1\right\rangle }, & \text{if }j>1
\end{cases}
\qquad\qquad\text{ for all }\left(  i,j\right)  \in\left\{  1,2,\ldots
,n\right\}  ^{2}.
\]
The matrix $B_{n}^{\left\langle \alpha\right\rangle }$ looks as follows:
\begin{align*}
B_{n}^{\left\langle \alpha\right\rangle }  &  =\left(
\begin{array}[c]{cccccc}
M_{\alpha}^{\left\langle 1\right\rangle } & M_{\alpha}^{\left\langle
0\right\rangle } & M_{\alpha}^{\left\langle -1\right\rangle } & \cdots &
M_{\alpha}^{\left\langle -n+3\right\rangle } & M_{\alpha}^{\left\langle
-n+2\right\rangle }\\
2M_{\alpha}^{\left\langle 2\right\rangle } & M_{\alpha}^{\left\langle
1\right\rangle } & M_{\alpha}^{\left\langle 0\right\rangle } & \cdots &
M_{\alpha}^{\left\langle -n+4\right\rangle } & M_{\alpha}^{\left\langle
-n+3\right\rangle }\\
3M_{\alpha}^{\left\langle 3\right\rangle } & M_{\alpha}^{\left\langle
2\right\rangle } & M_{\alpha}^{\left\langle 1\right\rangle } & \cdots &
M_{\alpha}^{\left\langle -n+5\right\rangle } & M_{\alpha}^{\left\langle
-n+4\right\rangle }\\
\vdots & \vdots & \vdots & \ddots & \vdots & \vdots\\
\left(  n-1\right)  M_{\alpha}^{\left\langle n-1\right\rangle } & M_{\alpha
}^{\left\langle n-2\right\rangle } & M_{\alpha}^{\left\langle n-3\right\rangle
} & \cdots & M_{\alpha}^{\left\langle 1\right\rangle } & M_{\alpha
}^{\left\langle 0\right\rangle }\\
nM_{\alpha}^{\left\langle n\right\rangle } & M_{\alpha}^{\left\langle
n-1\right\rangle } & M_{\alpha}^{\left\langle n-2\right\rangle } & \cdots &
M_{\alpha}^{\left\langle 2\right\rangle } & M_{\alpha}^{\left\langle
1\right\rangle }%
\end{array}
\right) \\
&  =\left(
\begin{array}[c]{cccccc}
M_{\alpha}^{\left\langle 1\right\rangle } & 1 & 0 & \cdots & 0 & 0\\
2M_{\alpha}^{\left\langle 2\right\rangle } & M_{\alpha}^{\left\langle
1\right\rangle } & 1 & \cdots & 0 & 0\\
3M_{\alpha}^{\left\langle 3\right\rangle } & M_{\alpha}^{\left\langle
2\right\rangle } & M_{\alpha}^{\left\langle 1\right\rangle } & \cdots & 0 &
0\\
\vdots & \vdots & \vdots & \ddots & \vdots & \vdots\\
\left(  n-1\right)  M_{\alpha}^{\left\langle n-1\right\rangle } & M_{\alpha
}^{\left\langle n-2\right\rangle } & M_{\alpha}^{\left\langle n-3\right\rangle
} & \cdots & M_{\alpha}^{\left\langle 1\right\rangle } & 1\\
nM_{\alpha}^{\left\langle n\right\rangle } & M_{\alpha}^{\left\langle
n-1\right\rangle } & M_{\alpha}^{\left\langle n-2\right\rangle } & \cdots &
M_{\alpha}^{\left\langle 2\right\rangle } & M_{\alpha}^{\left\langle
1\right\rangle }
\end{array}
\right)  .
\end{align*}
Then, $\det\left(  B_{n}^{\left\langle \alpha\right\rangle }\right)
=M_{\alpha\left\{  n\right\}  }$.
\end{proposition}

\begin{exercise}
\label{exe.prop.QSym.Malphas.dets}Prove Proposition
\ref{prop.QSym.Malphas.dets}.

[\textbf{Hint:} Substitute $\left(  \xx_{\mathbf{i}}^{\alpha}\right)
_{\mathbf{i}\in\operatorname*{SIS}\left(  \ell\right)  }$ for the variable set
in Exercise \ref{exe.pe-ep-determinants}, and recall Exercise
\ref{exe.QSym.Malphas.ps}(a).]
\end{exercise}

\begin{corollary}
\label{cor.QSym.Malphas.QSym}Let $\alpha$ be a composition. Let $s\in
 \ZZ$.

\begin{enumerate}
\item[(a)] We have $M_{\alpha}^{\left\langle s\right\rangle }\in
\Qsym$.

\item[(b)] We have $M_{\alpha}^{\left\langle s\right\rangle }\in
\Qsym_{s\left\vert \alpha\right\vert }$.
\end{enumerate}
\end{corollary}

\begin{exercise}
\label{exe.cor.QSym.Malphas.QSym}Prove Corollary \ref{cor.QSym.Malphas.QSym}.
\end{exercise}

We make one further definition:

\begin{definition}
\label{def.QSym.lyndon.red}Let $\alpha$ be a nonempty composition. Then, we
denote by $\gcd\alpha$ the greatest common divisor of the parts of $\alpha$.
(For instance, $\gcd\left(  8,6,4\right)  =2$.) We also define
$\operatorname*{red}\alpha$ to be the composition $\left(  \dfrac{\alpha_{1}%
}{\gcd\alpha},\dfrac{\alpha_{2}}{\gcd\alpha},\ldots,\dfrac{\alpha_{\ell}}%
{\gcd\alpha}\right)  $, where $\alpha$ is written in the form $\left(
\alpha_{1},\alpha_{2},\ldots,\alpha_{\ell}\right)  $.

We say that a nonempty composition $\alpha$ is
\emph{reduced}\index{reduced composition} if $\gcd \alpha = 1$.

We define $\mathfrak{RL}$ to be the set of all reduced Lyndon compositions. In
other words, $\mathfrak{RL}=\left\{  w\in\mathfrak{L}\ \mid\ w\text{ is
reduced}\right\}  $ (since $\mathfrak{L}$ is the set of all Lyndon compositions).
\end{definition}

Hazewinkel, in
\cite[proof of Thm. 6.7.5]{HazewinkelGubareniKirichenko},
denotes $\mathfrak{RL}$ by $eLYN$, calling reduced Lyndon compositions
``elementary Lyndon words''.

\begin{remark}
\label{rmk.QSym.lyndon.red}Let $\alpha$ be a nonempty composition.

(a) We have $\alpha=\left(  \operatorname*{red}\alpha\right)  \left\{
\gcd\alpha\right\}  $.

(b) The composition $\alpha$ is Lyndon if and only if the composition
$\operatorname*{red}\alpha$ is Lyndon.

(c) The composition $\operatorname*{red}\alpha$ is reduced.

(d) If $\alpha$ is reduced, then $\operatorname*{red}\alpha=\alpha$.

(e) If $s\in\left\{  1,2,3,\ldots\right\}  $, then the composition
$\alpha\left\{  s\right\}  $ is nonempty and satisfies $\operatorname*{red}%
\left(  \alpha\left\{  s\right\}  \right)  =\operatorname*{red}\alpha$ and
$\gcd\left(  \alpha\left\{  s\right\}  \right)  =s\gcd\alpha$.

(f) We have $\left(  \gcd\alpha\right)  \left\vert \operatorname*{red}%
\alpha\right\vert =\left\vert \alpha\right\vert $.
\end{remark}

\begin{exercise}
\label{exe.rmk.QSym.lyndon.red}Prove Remark \ref{rmk.QSym.lyndon.red}.
\end{exercise}

Our goal in this section is now to prove the following result of Hazewinkel:

\begin{theorem}
\label{thm.QSym.lyndon.Z-genset}The family $\left(  M_{w}^{\left\langle
s\right\rangle }\right)  _{\left(  w,s\right)  \in\mathfrak{RL}\times\left\{
1,2,3,\ldots\right\}  }$ is an algebraically independent generating set of the
$\kk$-algebra $\Qsym$.
\end{theorem}

This will (almost) immediately yield Theorem \ref{thm.QSym.lyndon}.

Our first step towards proving Theorem \ref{thm.QSym.lyndon.Z-genset} is the
following observation:

\begin{lemma}
\label{lem.QSym.lyndon.Z-genset.reindex}The family $\left(  M_{w}%
^{\left\langle s\right\rangle }\right)  _{\left(  w,s\right)  \in
\mathfrak{RL}\times\left\{  1,2,3,\ldots\right\}  }$ is a reindexing of the
family $\left(  M_{\operatorname*{red}\alpha}^{\left\langle \gcd
\alpha\right\rangle }\right)  _{\alpha\in\mathfrak{L}}$.
\end{lemma}

\begin{exercise}
\label{exe.lem.QSym.lyndon.Z-genset.reindex}Prove Lemma
\ref{lem.QSym.lyndon.Z-genset.reindex}.
\end{exercise}

Next, we show a lemma:

\begin{lemma}
\label{lem.QSym.lyndon.Malphas1}Let $\alpha$ be a nonempty composition. Let
$s\in \NN$. Then,%
\begin{equation}
s!M_{\alpha}^{\left\langle s\right\rangle }-M_{\alpha}^{s}\in
\sum\limits_{\substack{\beta\in\Comp_{s\left\vert \alpha
\right\vert };\\\ell\left(  \beta\right)  \leq\left(  s-1\right)  \ell\left(
\alpha\right)  }} \kk M_{\beta}.
\label{eq.lem.QSym.lyndon.Malphas1}
\end{equation}
(That is, $s!M_{\alpha}^{\left\langle s\right\rangle }-M_{\alpha}^{s}$ is a
$\kk$-linear combination of terms of the form $M_{\beta}$ with $\beta$
ranging over the compositions of $s\left\vert \alpha\right\vert $ satisfying
$\ell\left(  \beta\right)  \leq\left(  s-1\right)  \ell\left(  \alpha\right)
$.)
\end{lemma}

\begin{exercise}
\label{exe.lem.QSym.lyndon.Malphas1}Prove Lemma~\ref{lem.QSym.lyndon.Malphas1}.

[\textbf{Hint:} There are two approaches: One is to apply Proposition
\ref{prop.QSym.Malphas.dets}(a) and expand the determinant; the other is to
argue which monomials can appear in $s!M_{\alpha}^{\left\langle s\right\rangle
}-M_{\alpha}^{s}$.]
\end{exercise}

We now return to studying products of monomial quasisymmetric functions:

\begin{lemma}
\label{lem.QSym.lyndon.triangl.wll-lemma.precise}Let $n\in \NN$ and
$m\in \NN$. Let $u\in\Comp_{n}$ and
$v\in\Comp_{m}$. Let $z$ be the lexicographically
highest element of the multiset $u\shuffle v$. Let $h$ be the multiplicity
with which the word $z$ appears in the multiset $u\shuffle v$.
Then,\footnote{The following equality makes sense because we have
$z\in\Comp_{n+m}$ (by
Lemma~\ref{lem.QSym.lyndon.triang.wll-lemma}(a)).}
\[
M_{u}M_{v}=hM_{z}+\left(  \text{a sum of terms of the form }M_{w}\text{ with
}w\in\Comp_{n+m}\text{ satisfying }%
w\underset{\operatorname*{wll}}{<}z\right)  .
\]

\end{lemma}

\begin{proof}
[Proof of Lemma \ref{lem.QSym.lyndon.triangl.wll-lemma.precise}.]Lemma
\ref{lem.QSym.lyndon.triangl.wll-lemma.precise} was shown during the proof of
Lemma \ref{lem.QSym.lyndon.triang.wll-lemma}(b).
\end{proof}

\begin{corollary}
\label{cor.QSym.lyndon.triangl.wll-lemma.lyn-cfl}Let $n\in \NN$ and
$m\in \NN$. Let $u\in\Comp_{n}$ and
$v\in\Comp_{m}$. Regard $u$ and $v$ as words in
$\mathfrak{A}^{\ast}$. Assume that $u$ is a Lyndon word. Let $\left(
b_{1},b_{2},\ldots,b_{q}\right)  $ be the CFL factorization of the word $v$.

Assume that $u\geq b_{j}$ for every $j\in\left\{  1,2,\ldots,q\right\}  $. Let%
\[
h=1+\left\vert \left\{  j\in\left\{  1,2,\ldots,q\right\}  \ \mid
\ b_{j}=u\right\}  \right\vert .
\]
Then,%
\[
M_{u}M_{v}=hM_{uv}+\left(  \text{a sum of terms of the form }M_{w}\text{ with
}w\in\Comp_{n+m}\text{ satisfying }%
w\underset{\operatorname*{wll}}{<}uv\right)  .
\]

\end{corollary}

\begin{exercise}
\label{exe.cor.QSym.lyndon.triangl.wll-lemma.lyn-cfl}Prove Corollary
\ref{cor.QSym.lyndon.triangl.wll-lemma.lyn-cfl}.

[\textbf{Hint:} Apply Lemma \ref{lem.QSym.lyndon.triangl.wll-lemma.precise},
and notice that $uv$ is the lexicographically highest element of the multiset
$u\shuffle v$ (by Theorem \ref{thm.shuffle.lyndon.used}(e)), and that $h$ is
the multiplicity with which this word $uv$ appears in the multiset
$u\shuffle
v$ (this is a rewriting of Theorem \ref{thm.shuffle.lyndon.used}(e)).]
\end{exercise}

\begin{corollary}
\label{cor.QSym.lyndon.triangl.wll-lemma.lyn-xs}Let $k\in \NN$ and
$s\in \NN$. Let $x\in\Comp_{k}$ be such that $x$
is a Lyndon word. Then:

\begin{enumerate}
\item[(a)] The lexicographically highest element of the multiset $x\shuffle
x^{s}$ is $x^{s+1}$.

\item[(b)] We have%
\[
M_{x}M_{x^{s}}=\left(  s+1\right)  M_{x^{s+1}}+\left(  \text{a sum of terms of
the form }M_{w}\text{ with }w\in\Comp_{\left(
s+1\right)  k}\text{ satisfying }w\underset{\operatorname*{wll}}{<}%
x^{s+1}\right)  .
\]

\item[(c)] Let $t\in\Comp_{sk}$ be such that
$t\underset{\operatorname*{wll}}{<}x^{s}$. Then,%
\[
M_{x}M_{t}=\left(  \text{a sum of terms of the form }M_{w}\text{ with }%
w\in\Comp_{\left(  s+1\right)  k}\text{ satisfying
}w\underset{\operatorname*{wll}}{<}x^{s+1}\right)  .
\]

\end{enumerate}
\end{corollary}

\begin{exercise}
\label{exe.cor.QSym.lyndon.triangl.wll-lemma.lyn-xs}Prove Corollary
\ref{cor.QSym.lyndon.triangl.wll-lemma.lyn-xs}.

[\textbf{Hint:} Notice that $\left(  \underbrace{x,x,\ldots,x}_{s\text{
times}}\right)  $ is the CFL factorization of the word $x^{s}$. Now, part (a)
of Corollary \ref{cor.QSym.lyndon.triangl.wll-lemma.lyn-xs} follows from
Theorem \ref{thm.shuffle.lyndon.used}(c), part (b) follows from Corollary
\ref{cor.QSym.lyndon.triangl.wll-lemma.lyn-cfl}, and part (c) from Lemma
\ref{lem.QSym.lyndon.triang.wll-lemma}(c) (using part (a)).]
\end{exercise}

\begin{corollary}
\label{cor.QSym.lyndon.triangl.wll-lemma.cfl-cfl.strict}Let $n\in \NN$
and $m\in \NN$. Let $u\in\Comp_{n}$ and
$v\in\Comp_{m}$. Regard $u$ and $v$ as words in
$\mathfrak{A}^{\ast}$. Let $\left(  a_{1},a_{2},\ldots,a_{p}\right)  $ be the
CFL factorization of $u$. Let $\left(  b_{1},b_{2},\ldots,b_{q}\right)  $ be
the CFL factorization of the word $v$. Assume that $a_{i}>b_{j}$ for every
$i\in\left\{  1,2,\ldots,p\right\}  $ and $j\in\left\{  1,2,\ldots,q\right\}
$. Then,%
\[
M_{u}M_{v}=M_{uv}+\left(  \text{a sum of terms of the form }M_{w}\text{ with
}w\in\Comp_{n+m}\text{ satisfying }%
w\underset{\operatorname*{wll}}{<}uv\right)  .
\]

\end{corollary}

\begin{exercise}
\label{exe.cor.QSym.lyndon.triangl.wll-lemma.cfl-cfl.strict}Prove Corollary
\ref{cor.QSym.lyndon.triangl.wll-lemma.cfl-cfl.strict}.

[\textbf{Hint:} Combine Lemma \ref{lem.QSym.lyndon.triangl.wll-lemma.precise}
with the parts (c) and (d) of Theorem \ref{thm.shuffle.lyndon.used}.]
\end{exercise}

\begin{corollary}
\label{cor.QSym.lyndon.triangl.wll-lemma.BC}Let $n\in \NN$. Let
$u\in\Comp_{n}$ be a nonempty composition. Regard $u$
as a word in $\mathfrak{A}^{\ast}$. Let $\left(  a_{1},a_{2},\ldots
,a_{p}\right)  $ be the CFL factorization of $u$. Let $k\in\left\{
1,2,\ldots,p-1\right\}  $ be such that $a_{k}>a_{k+1}$. Let $x$ be the word
$a_{1}a_{2}\cdots a_{k}$, and let $y$ be the word $a_{k+1}a_{k+2}\cdots a_{p}%
$. Then,%
\[
M_{u}=M_{x}M_{y}-\left(  \text{a sum of terms of the form }M_{w}\text{ with
}w\in\Comp_{n}\text{ satisfying }%
w\underset{\operatorname*{wll}}{<}u\right)  .
\]

\end{corollary}

\begin{exercise}
\label{exe.cor.QSym.lyndon.triangl.wll-lemma.BC}Prove Corollary
\ref{cor.QSym.lyndon.triangl.wll-lemma.BC}.

[\textbf{Hint:} Apply Corollary
\ref{cor.QSym.lyndon.triangl.wll-lemma.cfl-cfl.strict} to $x$, $y$,
$\left\vert x\right\vert $, $\left\vert y\right\vert $, $k$, $p-k$, $\left(
a_{1},a_{2},\ldots,a_{k}\right)  $ and $\left(  a_{k+1},a_{k+2},\ldots
,a_{p}\right)  $ instead of $u$, $v$, $n$, $m$, $p$, $q$, $\left(  a_{1}%
,a_{2},\ldots,a_{p}\right)  $ and $\left(  b_{1},b_{2},\ldots,b_{q}\right)  $;
then, notice that $xy=u$ and $\left\vert x\right\vert +\left\vert y\right\vert
=n$.]
\end{exercise}

\begin{corollary}
\label{cor.QSym.lyndon.triangl.wll-lemma.As}Let $k\in \NN$. Let
$x\in\Comp_{k}$ be a composition. Assume that $x$ is a
Lyndon word. Let $s\in \NN$. Then,%
\[
M_{x}^{s}-s!M_{x^{s}}\in\sum\limits_{\substack{w\in
\Comp_{sk};\\w\underset{\operatorname*{wll}}{<}x^{s}}} \kk M_{w}.
\]
(Recall that $x^{s}$ is defined to be the word $\underbrace{xx\cdots
x}_{s\text{ times}}$.)
\end{corollary}

\begin{exercise}
\label{exe.cor.QSym.lyndon.triangl.wll-lemma.As}Prove Corollary
\ref{cor.QSym.lyndon.triangl.wll-lemma.As}.

[\textbf{Hint:} Rewrite the claim of Corollary
\ref{cor.QSym.lyndon.triangl.wll-lemma.As} in the form $M_{x}^{s}\in
s!M_{x^{s}}+\sum\limits_{\substack{w\in\Comp_{sk}%
;\\w\underset{\operatorname*{wll}}{<}x^{s}}} \kk M_{w}$. This can be
proven by induction over $s$, where in the induction step we need the
following two observations:

\begin{enumerate}
\item We have $M_{x}M_{x^{s}}\in\left(  s+1\right)  M_{x^{s+1}}+
\sum\limits_{\substack{w\in\Comp_{\left(  s+1\right)
k};\\w\underset{\operatorname*{wll}}{<}x^{s+1}}} \kk M_{w}$.

\item For every $t\in\Comp_{sk}$ satisfying
$t\underset{\operatorname*{wll}}{<}x^{s}$, we have $M_{x}M_{t}\in
\sum\limits_{\substack{w\in\Comp_{\left(  s+1\right)
k};\\w\underset{\operatorname*{wll}}{<}x^{s+1}}} \kk M_{w}$.
\end{enumerate}

These two observations follow from parts (b) and (c) of Corollary
\ref{cor.QSym.lyndon.triangl.wll-lemma.lyn-xs}.]
\end{exercise}

\begin{corollary}
\label{cor.QSym.lyndon.triangl.wll-lemma.As2}Let $k\in \NN$. Let
$x\in\Comp_{k}$ be a composition. Assume that $x$ is a
Lyndon word. Let $s\in \NN$. Then,%
\[
M_{x}^{\left\langle s\right\rangle }-M_{x^{s}}\in\sum\limits_{\substack{w\in
\Comp_{sk};\\w\underset{\operatorname*{wll}}{<}x^{s}%
}} \kk M_{w}.
\]
(Recall that $x^{s}$ is defined to be the word $\underbrace{xx\cdots
x}_{s\text{ times}}$.)
\end{corollary}

\begin{exercise}
\label{exe.cor.QSym.lyndon.triangl.wll-lemma.As2}Prove Corollary
\ref{cor.QSym.lyndon.triangl.wll-lemma.As2}.

[\textbf{Hint:} Lemma \ref{lem.QSym.lyndon.Malphas1} (applied to $\alpha=x$)
yields%
\[
s!M_{x}^{\left\langle s\right\rangle }-M_{x}^{s}\in\sum\limits_{\substack{\beta
\in\Comp_{sk};\\\ell\left(  \beta\right)  \leq\left(
s-1\right)  \ell\left(  x\right)  }}M_{\beta}=\sum\limits_{\substack{w\in
\Comp_{sk};\\\ell\left(  w\right)  \leq\left(
s-1\right)  \ell\left(  x\right)  }} \kk M_{w}\subset\sum\limits_{\substack{w\in
\Comp_{sk};\\w\underset{\operatorname*{wll}}{<}x^{s}%
}} \kk M_{w}%
\]
\footnote{since every $w\in\Comp_{sk}$ with the
property that $\ell\left(  w\right)  \leq\left(  s-1\right)  \ell\left(
x\right)  $ must satisfy $w\underset{\operatorname*{wll}}{<}x^{s}$}. Adding
this to the claim of Corollary \ref{cor.QSym.lyndon.triangl.wll-lemma.As},
obtain $s!M_{x}^{\left\langle s\right\rangle }-s!M_{x^{s}}\in
\sum\limits_{\substack{w\in\Comp_{sk}%
;\\w\underset{\operatorname*{wll}}{<}x^{s}}} \kk M_{w}$, that is,
$s!\left(  M_{x}^{\left\langle s\right\rangle }-M_{x^{s}}\right)  \in
\sum\limits_{\substack{w\in\Comp_{sk}%
;\\w\underset{\operatorname*{wll}}{<}x^{s}}} \kk M_{w}$. It remains to
get rid of the $s!$ on the left hand side. Assume WLOG that $\kk %
= \ZZ$, and argue that every $f\in\Qsym$ satisfying
$s!\cdot f\in\sum\limits_{\substack{w\in\Comp_{sk}%
;\\w\underset{\operatorname*{wll}}{<}x^{s}}} \kk M_{w}$ must itself lie
in $\sum\limits_{\substack{w\in\Comp_{sk}%
;\\w\underset{\operatorname*{wll}}{<}x^{s}}} \kk M_{w}$.]
\end{exercise}

We are now ready to prove Theorem \ref{thm.QSym.lyndon.Z-genset}:

\begin{exercise}
\label{exe.thm.QSym.lyndon.Z-genset}Prove Theorem
\ref{thm.QSym.lyndon.Z-genset}.

[\textbf{Hint:} Lemma \ref{lem.QSym.lyndon.Z-genset.reindex} yields that the
family $\left(  M_{w}^{\left\langle s\right\rangle }\right)  _{\left(
w,s\right)  \in\mathfrak{RL}\times\left\{  1,2,3,\ldots\right\}  }$ is a
reindexing of the family $\left(  M_{\operatorname*{red}w}^{\left\langle \gcd
w\right\rangle }\right)  _{w\in\mathfrak{L}}$. Hence, it is enough to prove
that the family $\left(  M_{\operatorname*{red}w}^{\left\langle \gcd
w\right\rangle }\right)  _{w\in\mathfrak{L}}$ is an algebraically independent
generating set of the $\kk$-algebra $\Qsym$. The latter
claim, in turn, will follow from Lemma \ref{lem.lyndonbases}%
(c)\footnote{applied to $A=\Qsym$, $b_{w}%
=M_{\operatorname*{red}w}^{\left\langle \gcd w\right\rangle }$,
$\operatorname*{wt}\left(  N\right)  =N$ and $g_{u}=M_{u}$} once it is proven
that the family $\left(  M_{\operatorname*{red}w}^{\left\langle \gcd
w\right\rangle }\right)  _{w\in\mathfrak{L}}$ generates the
$\kk$-algebra $\Qsym$. So it remains to show that the family
$\left(  M_{\operatorname*{red}w}^{\left\langle \gcd w\right\rangle }\right)
_{w\in\mathfrak{L}}$ generates the $\kk$-algebra $\Qsym$.

Let $U$ denote the $\kk$-subalgebra of $\Qsym$ generated
by $\left(  M_{\operatorname*{red}w}^{\left\langle \gcd w\right\rangle
}\right)  _{w\in\mathfrak{L}}$. It then suffices to prove that
$U=\Qsym$. To this purpose, it is enough to prove that
\begin{equation}
M_{\beta}\in U\ \ \ \ \ \ \ \ \ \ \text{for every composition }\beta.
\label{exe.thm.QSym.lyndon.Z-genset.hint1}%
\end{equation}

For every reduced Lyndon composition $\alpha$ and every $j\in\left\{
1,2,3,\ldots\right\}  $, the quasisymmetric function $M_{\alpha}^{\left\langle
j\right\rangle }$ is an element of the family $\left(  M_{\operatorname*{red}%
w}^{\left\langle \gcd w\right\rangle }\right)  _{w\in\mathfrak{L}}$ and thus
belongs to $U$. Combine this with Exercise \ref{exe.QSym.Malphas.ps}(d) to see
that%
\begin{equation}
M_{\beta}^{\left\langle s\right\rangle }\in U\ \ \ \ \ \ \ \ \ \ \text{for
every Lyndon composition }\beta\text{ and every }s\in\left\{  1,2,3,\ldots
\right\}  \label{exe.thm.QSym.lyndon.Z-genset.hint2}%
\end{equation}
(because every Lyndon composition $\beta$ can be written as $\alpha\left\{
n\right\}  $ for a reduced Lyndon composition $\alpha$ and an $n\in\left\{
1,2,3,\ldots\right\}  $). Now, prove (\ref{exe.thm.QSym.lyndon.Z-genset.hint1}%
) by strong induction: first, induct on $\left\vert \beta\right\vert $, and
then, for fixed $\left\vert \beta\right\vert $, induct on $\beta$ in the
wll-order. The induction step looks as follows: Fix some composition $\alpha$,
and assume (as induction hypothesis) that:

\begin{itemize}
\item (\ref{exe.thm.QSym.lyndon.Z-genset.hint1}) holds for every composition
$\beta$ satisfying $\left\vert \beta\right\vert <\left\vert \alpha\right\vert
$;

\item (\ref{exe.thm.QSym.lyndon.Z-genset.hint1}) holds for every composition
$\beta$ satisfying $\left\vert \beta\right\vert =\left\vert \alpha\right\vert
$ and $\beta\underset{\operatorname*{wll}}{<}\alpha$.
\end{itemize}

\noindent It remains to prove that (\ref{exe.thm.QSym.lyndon.Z-genset.hint1})
holds for $\beta=\alpha$.
In other words, it remains to prove that $M_{\alpha}\in U$.

Let $\left(  a_{1},a_{2},\ldots,a_{p}\right)  $ be the CFL factorization of
the word $\alpha$. Assume WLOG that $p\neq0$ (else, all is trivial). We are in
one of the following two cases:

\textit{Case 1:} All of the words $a_{1}$, $a_{2}$, $\ldots$, $a_{p}$ are equal.

\textit{Case 2:} Not all of the words $a_{1}$, $a_{2}$, $\ldots$, $a_{p}$ are equal.

In Case 2, there exists a $k\in\left\{  1,2,\ldots,p-1\right\}  $ satisfying
$a_{k}>a_{k+1}$ (since $a_{1}\geq a_{2}\geq\cdots\geq a_{p}$), and thus
Corollary \ref{cor.QSym.lyndon.triangl.wll-lemma.BC} (applied to $u=\alpha$,
$n=\left\vert \alpha\right\vert $, $x=a_{1}a_{2}\cdots a_{k}$ and
$y=a_{k+1}a_{k+2}\cdots a_{p}$) shows that
\begin{align*}
M_{\alpha}  &  =\underbrace{M_{a_{1}a_{2}\cdots a_{k}}}_{\substack{\in
U\\\text{(by the induction}\\\text{hypothesis)}}}
\underbrace{M_{a_{k+1}a_{k+2}\cdots a_{p}}}_{\substack{\in U\\
\text{(by the induction}\\\text{hypothesis)}}}\\
&  \ \ \ \ \ \ \ \ \ \ -\left(  \text{a sum of terms of the form
}\underbrace{M_{w}}_{\substack{\in U\\\text{(by the induction}%
\\\text{hypothesis)}}}\text{ with }w\in
\Comp_{\left\vert \alpha\right\vert }\text{ satisfying }%
w\underset{\operatorname*{wll}}{<}\alpha\right) \\
&  \in UU-\left(  \text{a sum of terms in }U\right)  \subset U.
\end{align*}

Hence, it only remains to deal with Case 1. In this case, set $x=a_{1}%
=a_{2}=\cdots=a_{p}$. Thus, $\alpha=a_{1}a_{2}\cdots a_{p}=x^{p}$, whence
$\left\vert \alpha\right\vert =p\left\vert x\right\vert $. But Corollary
\ref{cor.QSym.lyndon.triangl.wll-lemma.As2} (applied to $s=p$ and
$k=\left\vert x\right\vert $) yields%
\begin{align*}
M_{x}^{\left\langle p\right\rangle }-M_{x^{p}}  &  \in\sum\limits_{\substack{w\in
\Comp_{p\left\vert x\right\vert }%
;\\w\underset{\operatorname*{wll}}{<}x^{p}}} \kk M_{w}=
\sum\limits_{\substack{w\in\Comp_{\left\vert \alpha\right\vert
};\\w\underset{\operatorname*{wll}}{<}\alpha}} \kk \underbrace{M_{w}%
}_{\substack{\in U\\\text{(by the induction}\\\text{hypothesis)}%
}}\ \ \ \ \ \ \ \ \ \ \left(  \text{since }p\left\vert x\right\vert
=\left\vert \alpha\right\vert \text{ and }x^{p}=\alpha\right) \\
&  \subset\sum\limits_{\substack{w\in\Comp_{N}%
;\\w\underset{\operatorname*{wll}}{<}\alpha}} \kk U\subset U,
\end{align*}
so that $M_{x^{p}}\in\underbrace{M_{x}^{\left\langle p\right\rangle }%
}_{\substack{\in U\\\text{(by (\ref{exe.thm.QSym.lyndon.Z-genset.hint2}))}%
}}-U\subset U-U\subset U$. This rewrites as $M_{\alpha}\in U$ (since
$\alpha=x^{p}$). So $M_{\alpha}\in U$ is proven in both Cases 1 and 2, and
thus the induction proof of (\ref{exe.thm.QSym.lyndon.Z-genset.hint1}) is finished.]
\end{exercise}

\begin{exercise}
\label{exe.thm.QSym.lyndon}Prove Theorem \ref{thm.QSym.lyndon}.
\end{exercise}

Of course, this proof of Theorem \ref{thm.QSym.lyndon} yields a new (third)
proof for Proposition \ref{prop.QSym.lyndon}.

We notice the following corollary of our approach to Theorem
\ref{thm.QSym.lyndon}:

\begin{corollary}
\label{cor.QSym.lyndon.free}The $\Lambda$-algebra $\Qsym$ is a
polynomial algebra (over $\Lambda$).
\end{corollary}

\begin{exercise}
\label{exe.cor.QSym.lyndon.free}Prove Corollary \ref{cor.QSym.lyndon.free}.

[\textbf{Hint:} The algebraically independent generating set $\left(
M_{w}^{\left\langle s\right\rangle }\right)  _{\left(  w,s\right)
\in\mathfrak{RL}\times\left\{  1,2,3,\ldots\right\}  }$ of
$\Qsym$ contains the elements $M_{\left(  1\right)
}^{\left\langle s\right\rangle }=e_{s}\in\Lambda$ for all $s\in\left\{
1,2,3,\ldots\right\}  $.]
\end{exercise}

\subsection{\label{subsect.lyndon.gr}The Gessel-Reutenauer bijection and
symmetric functions}

In this section, we shall discuss the Gessel-Reutenauer bijection between
words and multisets of aperiodic necklaces, and use it to study another family
of symmetric functions.

The Gessel-Reutenauer bijection was studied in \cite{GesselReutenauer}, where
it was applied to various enumeration problems (e.g., counting permutations in
$\Symm_{n}$ with given descent set and given cycle type); it is also
closely related to the Burrows-Wheeler bijection used in data compression
(\cite{CrochemoreDesarmenienPerrin}), and to the structure of free Lie algebras
(\cite{GesselRestivoReutenauer}, \cite{Reutenauer}).
We shall first introduce the Gessel-Reutenauer
bijection and study it combinatorially in Subsection
\ref{subsub.lyndon.gr.bij}; then, in the following Subsection
\ref{subsub.lyndon.grfun}, we shall apply it to symmetric functions.

\subsubsection{\label{subsub.lyndon.gr.bij}Necklaces and the Gessel-Reutenauer
bijection}

We begin with definitions, some of which have already been made in Exercise
\ref{exe.words.necklaces}:

\begin{definition}
\label{def.lyndon.gr.cyclic}Throughout Section \ref{subsect.lyndon.gr}, we
shall freely use Definition \ref{def.words} and Definition
\ref{def.words.lyndon}. We fix a totally ordered alphabet $\mathfrak{A}$.
(This alphabet can be arbitrary, although most examples will use
$\mathfrak{A}=\left\{  1<2<3<\cdots\right\}  $.)

Let \dfn{$C$} denote the infinite cyclic group, written multiplicatively.
Fix a generator \dfn{$c$} of $C$.\ \ \ \ \footnote{So $C$ is a group
isomorphic to $\left(  \ZZ ,+\right)  $, and the isomorphism
$\left(  \ZZ ,+\right) \rightarrow C$ sends every $n\in \ZZ$ to $c^{n}$.
(Recall that we write the binary operation of $C$ as $\cdot$ instead of
$+$.)}

For any positive integer $n$, the group $C$ acts on $\mathfrak{A}^{n}$ from
the left according to the rule
\[
c\cdot\left(  a_{1},a_{2},\ldots,a_{n}\right)  =\left(  a_{2},a_{3}%
,\ldots,a_{n},a_{1}\right)  \ \ \ \ \ \ \ \ \ \ \text{for all }\left(
a_{1},a_{2},\ldots,a_{n}\right)  \in\mathfrak{A}^{n}.
\]
\footnote{In other words, $c$ rotates any $n$-tuple of elements of
$\mathfrak{A}$ cyclically to the left. Thus, $c^{n}\in C$ acts trivially on
$\mathfrak{A}^{n}$, and so this action of $C$ on $\mathfrak{A}^{n}$ factors
through $C/\left\langle c^{n}\right\rangle $ (a cyclic group of order $n$).}
The orbits of this $C$-action will be called
\emph{$n$-necklaces}\index{$n$-necklace}%
\footnote{See Exercise \ref{exe.words.necklaces} for the motivation behind
this word.
\par
Notice that there are no $0$-necklaces, because we required $n$ to be positive
in the definition of a necklace. This is intentional.}; they form a set
partition of the set $\mathfrak{A}^{n}$.

The $n$-necklace containing a given $n$-tuple $w\in\mathfrak{A}^{n}$ will be
denoted by $\left[  w\right]  $.

A \dfn{necklace} shall mean an $n$-necklace for some positive integer $n$.
Thus, for each nonempty word $w$, there is a well-defined necklace $\left[
w\right]  $ (namely, $\left[  w\right]  $ is an $n$-necklace, where
$n=\ell\left(  w\right)  $).

The \emph{period}\index{period of a necklace}
of a necklace $N$ is defined as the positive integer
$\left\vert N\right\vert $. (This $\left\vert N\right\vert $ is indeed a
positive integer, since $N$ is a finite nonempty set\footnote{by
Exercise \ref{exe.words.necklaces}(a), because $N$ is an $n$-necklace
for some positive integer $n$}.)

An $n$-necklace is said to be \emph{aperiodic}\index{aperiodic $n$-necklace}
if its period is $n$.
\end{definition}

\begin{example}
Let $\mathfrak{A}$ be the alphabet $\left\{  1<2<3<\cdots\right\}  $. The
orbit of the word $223$ under the $C$-action is the $3$-necklace $\left\{
223,232,322\right\}  $; it is an aperiodic $3$-necklace. The orbit of the word
$223223$ under the $C$-action is the $6$-necklace $\left\{
223223,232232,322322\right\}  $; it is not aperiodic (since it has period
$3$). The orbit of any nonempty word $w=\left(  w_{1},w_{2},\ldots
,w_{n}\right)  \in\mathfrak{A}^{n}$ is the $n$-necklace
\[
\left\{  \left(  w_{i},w_{i+1},\ldots,w_{n},w_{1},w_{2},\ldots,w_{i-1}\right)
\ \mid\ i\in\left\{  1,2,\ldots,n\right\}  \right\}  .
\]
We can draw this $n$-necklace on the plane as follows:
\[
\xymatrix@R=2.0pc@C=.8pc{
& w_1 \ar@/^1ex/[rr] & & w_2 \ar@/^1ex/[dr] \\
w_n \ar@/^1ex/[ur] & & \;\; & & w_3 \ar@/^1ex/[dl] \\
& w_{n-1} \ar@/^1ex/[ul] & & \iddots \ar@/^1ex/[ll]
}
\]
\end{example}

It is easy to see that the notion of an ``aperiodic
necklace'' we just defined is equivalent to the notion of a
``primitive necklace'' used in Exercise
\ref{exe.polynomials-and-necklaces}(b).

Exercise \ref{exe.words.necklaces}(a) shows that any $n$-necklace for any
positive integer $n$ is a finite nonempty set.
In other words, any necklace is a finite nonempty set.

Let us next introduce some notations regarding words and permutations. We
recall that a cycle of a permutation $\tau\in\Symm_{n}$ is an orbit
under the action of $\tau$ on $\left\{  1,2,\ldots,n\right\}  $. (This orbit
can be a $1$-element set, when $\tau$ has fixed points.) We begin with a basic definition:

\begin{definition}
\label{def.words.perm-ord}Let $\tau\in\Symm_{n}$ be a permutation. Let
$h\in\left\{  1,2,\ldots,n\right\}  $.

\begin{enumerate}
\item[(a)] We let $\operatorname*{ord}\nolimits_{\tau}\left(  h\right)  $
denote the smallest positive integer $i$ such that $\tau^{i}\left(  h\right)
=h$. (Basic properties of permutations show that this $i$ exists.)

\item[(b)] Let $w=\left(  w_{1},w_{2},\ldots,w_{n}\right)  \in\mathfrak{A}%
^{n}$ be a word. Then, $w_{\tau,h}$ shall denote the word $w_{\tau^{1}\left(
h\right)  }w_{\tau^{2}\left(  h\right)  }\cdots w_{\tau^{k}\left(  h\right)
}$, where $k=\operatorname*{ord}\nolimits_{\tau}\left(  h\right)  $.
\end{enumerate}
\end{definition}

\begin{example}
Let $\tau$ be the permutation $3142765\in\Symm_{7}$ (in one-line
notation). Then, $\operatorname*{ord}\nolimits_{\tau}\left(  1\right)  =4$
(since $\tau^{4}\left(  1\right)  =1$, but $\tau^{i}\left(  1\right)  \neq1$
for every positive integer $i<4$). Likewise, $\operatorname*{ord}%
\nolimits_{\tau}\left(  2\right)  =4$ and $\operatorname*{ord}\nolimits_{\tau
}\left(  3\right)  =4$ and $\operatorname*{ord}\nolimits_{\tau}\left(
4\right)  =4$ and $\operatorname*{ord}\nolimits_{\tau}\left(  5\right)  =2$
and $\operatorname*{ord}\nolimits_{\tau}\left(  6\right)  =1$ and
$\operatorname*{ord}\nolimits_{\tau}\left(  7\right)  =2$.

Now, let $w$ be the word $4112524\in\mathfrak{A}^{7}$. Then,%
\begin{align*}
w_{\tau,3}  &  =w_{\tau^{1}\left(  3\right)  }w_{\tau^{2}\left(  3\right)
}w_{\tau^{3}\left(  3\right)  }w_{\tau^{4}\left(  3\right)  }%
\ \ \ \ \ \ \ \ \ \ \left(  \text{since }\operatorname*{ord}\nolimits_{\tau
}\left(  3\right)  =4\right) \\
&  =w_{4}w_{2}w_{1}w_{3}\\
&  \ \ \ \ \ \ \ \ \ \ \left(  \text{since }\tau^{1}\left(  3\right)  =4\text{
and }\tau^{2}\left(  3\right)  =\tau\left(  4\right)  =2\text{ and }\tau
^{3}\left(  3\right)  =\tau\left(  2\right)  =1\text{ and }\tau^{4}\left(
3\right)  =\tau\left(  1\right)  =3\right) \\
&  =2141.
\end{align*}
Likewise, we can check that $w_{\tau,1}=w_{3}w_{4}w_{2}w_{1}=1214$ and
$w_{\tau,5}=w_{7}w_{5}=45$ and $w_{\tau,6}=w_{6}=2$.
\end{example}

We begin the study of the words $w_{\tau,h}$ by stating some of their simplest
properties:\footnote{See Exercise \ref{exe.lyndon.gr.basics} below for the
proof of Proposition \ref{prop.words.read1}, as well as for the proofs of all
other propositions stated before Exercise \ref{exe.lyndon.gr.basics}.}

\begin{proposition}
\label{prop.words.read1}Let $w=\left(  w_{1},w_{2},\ldots,w_{n}\right)
\in\mathfrak{A}^{n}$ be a word. Let $\tau\in\Symm_{n}$. Let
$h\in\left\{  1,2,\ldots,n\right\}  $. Then:

\begin{enumerate}
\item[(a)] The word $w_{\tau,h}$ is nonempty and has length
$\operatorname*{ord}\nolimits_{\tau}\left(  h\right)  $.

\item[(b)] The first letter of the word $w_{\tau,h}$ is $w_{\tau\left(
h\right)  }$.

\item[(c)] The last letter of the word $w_{\tau,h}$ is $w_{h}$.

\item[(d)] We have $w_{\tau,\tau\left(  h\right)  }=c\cdot w_{\tau,h}$.

\item[(e)] We have $w_{\tau,\tau^{i}\left(  h\right)  }=c^{i}\cdot w_{\tau,h}$
for each $i\in \ZZ $.
\end{enumerate}
\end{proposition}

Recall that if $n\in \NN$ and if $w\in\mathfrak{A}^{n}$ is a word, then
a permutation $\operatorname*{std}w\in\Symm_{n}$ was defined in
Definition \ref{def.words.std}. The words $w_{\tau,h}$ have particularly nice
properties when $\tau=\left(  \operatorname*{std}w\right)  ^{-1}$:

\begin{lemma}
\label{lem.sol.lyndon.gr.basics.1}Let $w=\left(  w_{1},w_{2},\ldots
,w_{n}\right)  \in\mathfrak{A}^{n}$ be a word. Let $\tau$ be the permutation
$\left(  \operatorname*{std}w\right)  ^{-1}\in\Symm_{n}$. Let $\alpha$
and $\beta$ be two elements of $\left\{  1,2,\ldots,n\right\}  $ such that
$\alpha<\beta$. Then:

\begin{enumerate}
\item[(a)] If $\tau^{-1}\left(  \alpha\right)  <\tau^{-1}\left(  \beta\right)
$, then $w_{\alpha}\leq w_{\beta}$.

\item[(b)] If $\tau^{-1}\left(  \alpha\right)  \geq\tau^{-1}\left(
\beta\right)  $, then $w_{\alpha}>w_{\beta}$.

\item[(c)] We have $w_{\tau\left(  \alpha\right)  }\leq w_{\tau\left(
\beta\right)  }$.

\item[(d)] If $\tau\left(  \alpha\right)  \geq\tau\left(  \beta\right)  $,
then $w_{\tau\left(  \alpha\right)  }<w_{\tau\left(  \beta\right)  }$.

\item[(e)] If $w_{\tau\left(  \alpha\right)  }=w_{\tau\left(  \beta\right)  }%
$, then $\tau\left(  \alpha\right)  <\tau\left(  \beta\right)  $.

\item[(f)] If $w_{\tau,\alpha}=w_{\tau,\beta}$, then $\tau\left(
\alpha\right)  <\tau\left(  \beta\right)  $ and $w_{\tau,\tau\left(
\alpha\right)  }=w_{\tau,\tau\left(  \beta\right)  }$.

\item[(g)] If $w_{\tau,\alpha}=w_{\tau,\beta}$, then $\tau^{i}\left(
\alpha\right)  <\tau^{i}\left(  \beta\right)  $ for each $i\in \NN $.

\item[(h)] Let $j\in \NN $ be such that every $i\in\left\{
0,1,\ldots,j-1\right\}  $ satisfies $w_{\tau^{i+1}\left(  \alpha\right)
}=w_{\tau^{i+1}\left(  \beta\right)  }$. Then, $w_{\tau^{j+1}\left(
\alpha\right)  }\leq w_{\tau^{j+1}\left(  \beta\right)  }$.
\end{enumerate}
\end{lemma}

\begin{proposition}
\label{prop.words.read-std2}Let $w\in\mathfrak{A}^{n}$ be a word. Let $\tau$
be the permutation $\left(  \operatorname*{std}w\right)  ^{-1}\in
\Symm_{n}$. Let $z$ be a cycle of $\tau$. Then:

\begin{enumerate}
\item[(a)] For each $h\in z$, we have $\left[  w_{\tau,h}\right]  =\left\{
w_{\tau,i}\ \mid\ i\in z\right\}  $.

\item[(b)] If $\alpha$ and $\beta$ are two distinct elements of $z$, then
$w_{\tau,\alpha}\neq w_{\tau,\beta}$.

\item[(c)] We have $\left|  \left\{  w_{\tau,i}\ \mid\ i\in z\right\}
\right|  = \left|  z\right|  $.

\item[(d)] The set $\left\{  w_{\tau,i}\ \mid\ i\in z\right\}  $ is an
aperiodic necklace.
\end{enumerate}
\end{proposition}

\begin{exercise}
\label{exe.lyndon.gr.basics}Prove Proposition \ref{prop.words.read1}, Lemma
\ref{lem.sol.lyndon.gr.basics.1} and Proposition \ref{prop.words.read-std2}.
\end{exercise}

\begin{definition}
Let $w \in\mathfrak{A}^{n}$ be a
word. Let $\tau$ be the permutation $\left(  \operatorname*{std}w\right)
^{-1}\in\Symm_{n}$. Let $z$ be a cycle of $\tau$. Then, we define an
aperiodic necklace $\left[  w\right]  _{z}$ by $\left[  w\right]
_{z}=\left\{  w_{\tau,i}\ \mid\ i\in z\right\}  $. (This is indeed an
aperiodic necklace, according to Proposition \ref{prop.words.read-std2}(d).)
\end{definition}

\begin{example}
Let $\mathfrak{A}$ be the alphabet $\left\{ 1<2<3<\cdots \right\}$, and let
$w$ be the word $2511321\in\mathfrak{A}^{7}$. Let $\tau$ be the
permutation $\left(  \operatorname*{std}w\right)  ^{-1}\in\Symm_{7}$;
this is the permutation $3471652$ (in one-line notation). One cycle of $\tau$
is $z=\left\{  1,3,7,2,4\right\}  $. The corresponding aperiodic necklace
$\left[  w\right]  _{z}$ is%
\begin{align*}
\left[  w\right]  _{z}  &  =\left\{  w_{\tau,i}\ \mid\ i\in z\right\}
=\left\{  w_{\tau,1},w_{\tau,3},w_{\tau,7},w_{\tau,2},w_{\tau,4}\right\}
\ \ \ \ \ \ \ \ \ \ \left(  \text{since }z=\left\{  1,3,7,2,4\right\}  \right)
\\
&  =\left\{  11512,15121,51211,12115,21151\right\}  =\left[  11512\right]  .
\end{align*}

\end{example}

\begin{definition}
We let \dfn{$\mathfrak{N}$} be the set of all necklaces.
We let \dfn{$\mathfrak{N}^{\mathfrak{a}}$} be the set of all aperiodic
necklaces.
We let \dfn{$\mathfrak{MN}^{\mathfrak{a}}$} be the set of all finite
multisets of aperiodic necklaces.
\end{definition}

\begin{definition}
\label{def.lyndon.gr.GR}We define a map $\operatorname*{GR}:\mathfrak{A}%
^{\ast}\rightarrow\mathfrak{MN}^{\mathfrak{a}}$\index{$\operatorname{GR}$}
as follows:

Let $w\in\mathfrak{A}^{\ast}$. Let $n=\ell\left(  w\right)  $ (so that
$w\in\mathfrak{A}^{n}$). Let $\tau$ be the permutation $\left(
\operatorname*{std}w\right)  ^{-1}\in\Symm_{n}$. Then, we define the
multiset $\operatorname*{GR}w\in\mathfrak{MN}^{\mathfrak{a}}$ by setting%
\[
\operatorname*{GR}w=\left\{  \left[  w\right]  _{z}\ \mid\ z\text{ is a cycle
of }\tau\right\}  _{\operatorname*{multiset}}.
\]
(This multiset $\operatorname*{GR}w$ is indeed a finite multiset of aperiodic
necklaces\footnote{Indeed, this multiset $\operatorname*{GR}w$ is finite
(since $\tau$ has only finitely many cycles), and its elements $\left[
w\right]  _{z}$ are aperiodic necklaces (as we have seen in the definition of
$\left[  w\right]  _{z}$).}, and thus belongs to $\mathfrak{MN}^{\mathfrak{a}%
}$.)
\end{definition}

\begin{example}
\label{exa.lyndon.gr.GR}Let $\mathfrak{A}$ be the alphabet $\left\{
1<2<3<\cdots\right\}  $, and let $w=33232112\in\mathfrak{A}^{8}$.

To compute $\operatorname*{GR}w$, we first notice that $\operatorname*{std}%
w=67384125$ (in one-line notation). Hence, the
permutation $\tau$ from Definition \ref{def.lyndon.gr.GR} satisfies
$\tau=\left(  \operatorname*{std}w\right)  ^{-1}=
67358124$. The cycles of $\tau$ are $\left\{  1,6\right\}  $,
$\left\{  2,7\right\}  $, $\left\{  3\right\}  $ and $\left\{  4,5,8\right\}
$. Thus,%
\begin{align*}
\operatorname*{GR}w  &  =\left\{  \left[  w\right]  _{z}\ \mid\ z\text{ is a
cycle of }\tau\right\}  _{\operatorname*{multiset}}=\left\{  \left[  w\right]
_{\left\{  1,6\right\}  },\left[  w\right]  _{\left\{  2,7\right\}  },\left[
w\right]  _{\left\{  3\right\}  },\left[  w\right]  _{\left\{  4,5,8\right\}
}\right\}  _{\operatorname*{multiset}}\\
&  =\left\{  \left[  31\right]  ,\left[  31\right]  ,\left[  2\right]
,\left[  322\right]  \right\}  _{\operatorname*{multiset}}=\left\{  \left[
13\right]  ,\left[  13\right]  ,\left[  2\right]  ,\left[  223\right]
\right\}  _{\operatorname*{multiset}}%
\end{align*}
(since $\left[  31\right]  =\left[  13\right]  $ and $\left[  322\right]
=\left[  223\right]  $ as necklaces).
Drawn on the plane, the necklaces in $\operatorname*{GR}w$ look as follows:
\[
\begin{array}{cccc}
\xymatrix@R=1.4pc@C=.8pc{
\\
1 \ar@/^3ex/[rr] & & 3 \ar@/^3ex/[ll]
}
&
\qquad
\xymatrix@R=1.4pc@C=.8pc{
\\
1 \ar@/^3ex/[rr] & & 3 \ar@/^3ex/[ll]
}
&\quad\qquad
\xymatrix@R=1.4pc@C=.8pc{
\\ 2
}
&\quad\qquad
\xymatrix@R=1.4pc@C=1.0pc{
2 \ar@/^2ex/[drr] \\
& & 2 \ar@/^2ex/[dll] \\
3 \ar@/^4ex/[uu]
}
\end{array}
\]
\end{example}

The map $\operatorname*{GR}$ is called the \dfn{Gessel-Reutenauer
bijection}. In order to show that it indeed is a bijection, we shall construct
its inverse. First, we introduce some further objects.

\begin{definition}
A nonempty word $w$ is said to be \emph{aperiodic}\index{aperiodic word}
if there exist no
$m\geq2$ and $u\in\mathfrak{A}^{\ast}$ satisfying $w=u^{m}$.

Let \dfn{$\mathfrak{A}^{\mathfrak{a}}$} be the set of all aperiodic words in
$\mathfrak{A}^{\ast}$.
\end{definition}

For example, the word $132231$ is aperiodic, but the word $132132$ is not
(since $132132 = u^m$ for $u = 132$ and $m = 2$).

Aperiodic words are directly connected to aperiodic necklaces, as the
following facts show:\footnote{See Exercise \ref{exe.lyndon.gr.aper-rel} for
the proofs of all unproved statements made until Exercise
\ref{exe.lyndon.gr.aper-rel}.}

\begin{proposition}
\label{prop.words.aper-neck}Let $w\in\mathfrak{A}^{\ast}$ be a nonempty word.
Then, the word $w$ is aperiodic if and only if the necklace $\left[  w\right]
$ is aperiodic.
\end{proposition}

\begin{corollary}
\label{cor.words.aper-neck-cw}Let $w\in\mathfrak{A}^{\ast}$ be an aperiodic
word. Then, the word $c\cdot w$ is aperiodic.\footnote{See Definition
\ref{def.lyndon.gr.cyclic} for the definition of $c$ and its
action on words.}
\end{corollary}

\begin{corollary}
\label{cor.words.aper-neck-ms}Each aperiodic necklace is a set of aperiodic words.
\end{corollary}

Let us now introduce a new total order on the set $\mathfrak{A}^{\mathfrak{a}%
}$ of all aperiodic words:

\begin{definition}
Let $u$ and $v$ be two aperiodic words. Then, we write $u\leq_{\omega}v$ if
and only if $uv\leq vu$. Thus, we have defined a binary relation
\dfn{$\leq_{\omega}$}
on the set $\mathfrak{A}^{\mathfrak{a}}$ of all aperiodic words.
\end{definition}

\begin{proposition}
\label{prop.words.aper-rel.tord}The relation $\leq_{\omega}$ on the set
$\mathfrak{A}^{\mathfrak{a}}$ is the smaller-or-equal relation of a total order.
\end{proposition}

For the next proposition, we should recall Definition
\ref{def.lyndon.gr.cyclic} (and, in particular, the meaning of $c$ and its
action on words).

\begin{proposition}
\label{prop.words.aper-rel.1}Let $u$ and $v$ be two aperiodic words.

\begin{enumerate}
\item[(a)] We have $u\leq_{\omega}v$ if and only if either $u_{1}<v_{1}$ or
$\left(  u_{1}=v_{1}\text{ and }c\cdot u\leq_{\omega}c\cdot v\right)
$.\ \ \ \footnote{The relation ``$c\cdot u\leq_{\omega}c\cdot
v$'' here makes sense because the words $c\cdot u$ and
$c\cdot v$ are aperiodic (by Corollary \ref{cor.words.aper-neck-cw}).}

\item[(b)] If $u\neq v$, then there exists some $i\in \NN $ satisfying
$\left(  c^{i}\cdot u\right)  _{1}\neq\left(  c^{i}\cdot v\right)  _{1}$.

\item[(c)] We have $u\leq_{\omega}v$ if and only if the smallest
$i\in \NN $ satisfying $\left(  c^{i}\cdot u\right)  _{1}\neq\left(
c^{i}\cdot v\right)  _{1}$ \textbf{either} does not exist \textbf{or}
satisfies $\left(  c^{i}\cdot u\right)  _{1}<\left(  c^{i}\cdot v\right)
_{1}$.

\item[(d)] Let $n$ and $m$ be positive integers such that $n\ell\left(
u\right)  =m\ell\left(  v\right)  $. We have $u\leq_{\omega}v$ if and only if
$u^{n}\leq v^{m}$.
\end{enumerate}
\end{proposition}

\begin{remark}
We are avoiding the use of infinite words here; if we didn't, we could restate
the relation $\leq_{\omega}$ in a simpler way (which is easily seen to be
equivalent to Proposition \ref{prop.words.aper-rel.1}(c)): Two aperiodic words
$u$ and $v$ satisfy $u\leq_{\omega}v$ if and only if $u^{\infty}\leq
v^{\infty}$. Here, for any nonempty word $w$, we are letting $w^{\infty}$
denote the infinite word
\[
\left(  w_{1},w_{2},\ldots,w_{\ell\left(  w\right)  },w_{1},w_{2}%
,\ldots,w_{\ell\left(  w\right)  },w_{1},w_{2},\ldots,w_{\ell\left(  w\right)
},\ldots\right)
\]
(that is, the word $w$ repeated endlessly), and the symbol ``$\leq$'' in
``$u^{\infty}\leq v^{\infty}$'' refers to the lexicographic order on
$\mathfrak{A}^{\infty}$.

Other equivalent descriptions of the relation $\leq_{\omega}$ (or, more
precisely, of the ``strictly less'' relation corresponding to it) can be
found in \cite[Corollary 11]{DolceRestivoReutenauer-svl}.
\end{remark}

\begin{proposition}
\label{prop.words.read-std3}Let $w
\in\mathfrak{A}^{n}$ be a word. Let $\tau$ be the permutation $\left(
\operatorname*{std}w\right)  ^{-1}\in\Symm_{n}$. Then:

\begin{enumerate}
\item[(a)] The words $w_{\tau,1},w_{\tau,2},\ldots,w_{\tau,n}$ are aperiodic.

\item[(b)] We have $w_{\tau,1}\leq_{\omega}w_{\tau,2}\leq_{\omega}\cdots
\leq_{\omega}w_{\tau,n}$.
\end{enumerate}
\end{proposition}

\begin{exercise}
\label{exe.lyndon.gr.aper-rel}Prove Proposition \ref{prop.words.aper-neck},
Corollary \ref{cor.words.aper-neck-cw}, Corollary \ref{cor.words.aper-neck-ms}%
, Proposition \ref{prop.words.aper-rel.tord}, Proposition
\ref{prop.words.aper-rel.1} and Proposition \ref{prop.words.read-std3}.
\end{exercise}

We need two more notations about multisets:

\begin{definition}
Let $T$ be a totally ordered set, and let $\leq_T$ be the smaller-or-equal
relation of $T$. Let $M$ be a finite multiset of elements of $T$.
Then, there is a unique list $\left( m_1, m_2, \ldots, m_n \right)$
such that
\[
\left\{ m_1, m_2, \ldots, m_n \right\}_{\operatorname*{multiset}}
= M
\ \ \ \ \ \ \ \ \ \ \text{and}\ \ \ \ \ \ \ \ \ \ %
m_1 \leq_T m_2 \leq_T \cdots \leq_T m_n .
\]
This list $\left( m_1, m_2, \ldots, m_n \right)$ is obtained by listing
all elements of $M$ (with their multiplicities) in increasing order
(increasing with respect to $\leq_T$). We shall refer to this list
$\left( m_1, m_2, \ldots, m_n \right)$ as the
\dfn{$\leq_T$-increasing list}\index{increasing list of a multiset}
of $M$.

(For example, the $\leq_{\ZZ}$-increasing list of $\left\{
1,2,3,2,1\right\}  _{\operatorname*{multiset}}$ is
$\left( 1,1,2,2,3 \right) $.)
\end{definition}

\begin{definition}
Let $S$ be a finite multiset.

\begin{enumerate}
\item[(a)] The \emph{support}\index{support of a multiset}
\dfn{$\operatorname*{Supp}S$} is defined to be the
set of all elements of $S$.
Thus, if
$S = \left\{ m_1, m_2, \ldots, m_n \right\}_{\operatorname*{multiset}}$,
then
$\operatorname*{Supp} S = \left\{ m_1, m_2, \ldots, m_n \right\}  $.

\item[(b)] For each $s\in S$, let $M_{s}$ be a finite multiset. Then, we
define the \dfn{multiset union} \dfn{$\biguplus_{s\in S}M_{s}$} to be
the finite multiset $M$ with the following property:
For any object $x$, we have
\begin{align*}
\left(  \text{multiplicity of }x\text{ in }M\right)
&  =\sum_{s\in\operatorname*{Supp}S}\left(  \text{multiplicity of }s\text{ in
}S\right)  \cdot\left(  \text{multiplicity of }x\text{ in }M_{s}\right)  .
\end{align*}

For example:

\begin{itemize}
\item If $S=\left\{  1,2,3\right\}  _{\operatorname*{multiset}}$ and
$M_{s}=\left\{  s,s+1\right\}  _{\operatorname*{multiset}}$ for each
$s\in\operatorname*{Supp}S$, then $\biguplus_{s\in S}M_{s}=\left\{
1,2,2,3,3,4\right\}  _{\operatorname*{multiset}}$.

\item If $S=\left\{  1,1,2\right\}  _{\operatorname*{multiset}}$ and
$M_{s}=\left\{  s,s+1\right\}  _{\operatorname*{multiset}}$ for each
$s\in\operatorname*{Supp}S$, then $\biguplus_{s\in S}M_{s}=\left\{
1,1,2,2,2,3\right\}  _{\operatorname*{multiset}}$.
\end{itemize}

We regard each set as a multiset; thus, the multiset union $\biguplus_{s\in
S}M_{s}$ is also defined when the $M_{s}$ are sets.
\end{enumerate}
\end{definition}

Now, we can construct the inverse of the Gessel-Reutenauer bijection:

\begin{definition}
\label{def.lyndon.gr.RG}We define a map $\operatorname*{RG}:\mathfrak{MN}%
^{\mathfrak{a}}\rightarrow\mathfrak{A}^{\ast}$\index{$\operatorname{RG}$}
as follows:

Let $M\in\mathfrak{MN}^{\mathfrak{a}}$ be a finite multiset of aperiodic
necklaces. Let $M^{\prime}=\biguplus_{N\in M}N$. (We are here using the fact
that each necklace $N\in M$ is a finite set, thus a finite multiset.) Notice
that $M^{\prime}$ is a finite multiset of aperiodic words\footnote{Indeed:
\par
\begin{itemize}
\item Each $N\in M$ is an aperiodic necklace (since $M$ is a multiset
of aperiodic necklaces), and thus (by Corollary \ref{cor.words.aper-neck-ms})
a set of aperiodic words.
Therefore, $\biguplus_{N\in M}N$ is a multiset of aperiodic
words.
\par
\item Each $N\in M$ is a necklace, and thus is a finite set (since any
necklace is a finite set). Since the multiset $M$ is also finite, this shows
that $\biguplus_{N\in M}N$ is finite.
\end{itemize}
\par
Thus, $\biguplus_{N\in M}N$ is a finite multiset of aperiodic words. In other
words, $M^{\prime}$ is a finite multiset of aperiodic words (since $M^{\prime
}=\biguplus_{N\in M}N$).}. Let $\left(  m_{1},m_{2},\ldots,m_{n}\right)  $ be
the $\leq_{\omega}$-increasing list of $M^{\prime}$. For each $i\in\left\{
1,2,\ldots,n\right\}  $, let $\ell_{i}$ be the last letter of the nonempty
word $m_{i}$. Then, $\operatorname*{RG}\left(  M\right)  $ is defined to be
the word $\left(  \ell_{1},\ell_{2},\ldots,\ell_{n}\right)  \in\mathfrak{A}%
^{\ast}$.
\end{definition}

\begin{example}
\label{exa.lyndon.gr.RG}Let $\mathfrak{A}$ be the alphabet $\left\{
1<2<3<\cdots\right\}  $, and let $M=\left\{  \left[  13\right]  ,\left[
13\right]  ,\left[  2\right]  ,\left[  223\right]  \right\}
_{\operatorname*{multiset}}$. Clearly, $M\in\mathfrak{M}\mathfrak{N}%
^{\mathfrak{a}}$ (since $M$ is a finite multiset of aperiodic necklaces).
(Actually, $M$ is the multiset of aperiodic necklaces drawn in
Example~\ref{exa.lyndon.gr.GR}.)
In order to
compute the word $\operatorname*{RG}\left(  M\right)  $, let us first compute
the multiset $M^{\prime}$ from Definition \ref{def.lyndon.gr.RG}. Indeed, the
definition of $M^{\prime}$ yields%
\begin{align*}
M^{\prime}  &  =\biguplus_{N\in M}N=\underbrace{\left[  13\right]
}_{=\left\{  13,31\right\}  }\uplus\underbrace{\left[  13\right]  }_{=\left\{
13,31\right\}  }\uplus\underbrace{\left[  2\right]  }_{=\left\{  2\right\}
}\uplus\underbrace{\left[  223\right]  }_{=\left\{  223,232,322\right\}  }\\
&  \ \ \ \ \ \ \ \ \ \ \left(
\text{where we are using the notation }M_{1}\uplus M_{2}\uplus\cdots\uplus
M_{k}
\text{ for a multiset union }\biguplus\limits_{s\in\left\{  1,2,\ldots
,k\right\}  }M_{s}
\right) \\
&  =\left\{  13,31\right\}  \uplus\left\{  13,31\right\}  \uplus\left\{
2\right\}  \uplus\left\{  223,232,322\right\} \\
&  =\left\{  13,31,13,31,2,223,232,322\right\}  _{\operatorname*{multiset}}.
\end{align*}
Hence, the $\leq_{\omega}$-increasing list of $M^{\prime}$ is $\left(
13,13,2,223,232,31,31,322\right)  $ (since $13\leq_{\omega}13\leq_{\omega
}2\leq_{\omega}223\leq_{\omega}232\leq_{\omega}31\leq_{\omega}31\leq_{\omega
}322$). The last letters of the words in this list are $3,3,2,3,2,1,1,2$ (in
this order). Hence, Definition \ref{def.lyndon.gr.RG} shows that
\[
\operatorname*{RG}\left(  M\right)  =\left(  3,3,2,3,2,1,1,2\right)
=33232112.
\]

\end{example}

\begin{remark}
The $\leq_{\omega}$-increasing list of a multiset $M^{\prime}$ of aperiodic
words is not always the same as its $\leq$-increasing list. For example, the
$\leq_{\omega}$-increasing list of $\left\{  2,21\right\}  $ is $\left(
21,2\right)  $ (since $21\leq_{\omega}2$), whereas its $\leq$-increasing list
is $\left(  2,21\right)  $ (since $2\leq21$).
\end{remark}

A comparison of Examples \ref{exa.lyndon.gr.GR} and \ref{exa.lyndon.gr.RG}
suggests that the maps $\operatorname*{GR}$ and $\operatorname*{RG}$ undo
one another. This is indeed true, as the following theorem (due to Gessel
and Reutenauer \cite[Lemma 3.4 and Example 3.5]{GesselReutenauer}; also
proved in \cite[Theorem 7.20]{Reutenauer}, \cite[Theorem 3.1 and
Proposition 3.1]{DiaconisMcGrathPitman} and
\cite[\S 2]{GesselRestivoReutenauer}) shows:

\begin{theorem}
\label{thm.lyndon.gr.bij}The maps $\operatorname*{GR}:\mathfrak{A}^{\ast
}\rightarrow\mathfrak{MN}^{\mathfrak{a}}$ and $\operatorname*{RG}%
:\mathfrak{MN}^{\mathfrak{a}}\rightarrow\mathfrak{A}^{\ast}$ are mutually
inverse bijections.
\end{theorem}

\begin{exercise}
\label{exe.lyndon.gr.bij}Prove Theorem \ref{thm.lyndon.gr.bij}.

[\textbf{Hint:} First, use Proposition \ref{prop.words.read-std3} to show that
$\operatorname*{RG}\circ\operatorname*{GR} = \id$. Then recall
the fact that any injective map between two finite sets of the same sizes is a
bijection. This does not directly apply here, since the sets $\mathfrak{A}%
^{\ast}$ and $\mathfrak{MN}^{\mathfrak{a}}$ are usually not finite. However,
$\operatorname*{GR}$ can be restricted to a map between two appropriate finite
subsets, obtained by focussing on a finite sub-alphabet of $\mathfrak{A}$
and fixing the length of the words; these subsets can be shown to have equal
size using the Chen-Fox-Lyndon factorization (see the following paragraph for
the connection).\footnote{This argument roughly follows
\cite{GesselRestivoReutenauer}.}]
\end{exercise}

Theorem \ref{thm.lyndon.gr.bij} shows that the sets $\mathfrak{A}^{\ast}$ and
$\mathfrak{MN}^{\mathfrak{a}}$ are in bijection. This bijection is in some
sense similar to the Chen-Fox-Lyndon factorization\footnote{The
Chen-Fox-Lyndon factorization (Theorem \ref{thm.words.CFL}) provides a
bijection between words in $\mathfrak{A}^{\ast}$ and multisets of Lyndon words
(because the factors in the CFL factorization of a word $w\in\mathfrak{A}%
^{\ast}$ can be stored in a multiset), whereas the Gessel-Reutenauer bijection
$\operatorname*{GR}:\mathfrak{A}^{\ast}\rightarrow\mathfrak{MN}^{\mathfrak{a}%
}$ is a bijection between words in $\mathfrak{A}^{\ast}$ and multisets of
aperiodic necklaces. Since the Lyndon words are in bijection with the
aperiodic necklaces (by Exercise \ref{exe.words.necklaces}(e)), we can thus
view the two bijections as having the same targets (and the same domains).
That said, they are not the same bijection.}, and preserves various quantities
(for example, the number of times a given letter $a$ appears in a word
$w\in\mathfrak{A}^{\ast}$ equals the number of times this letter $a$ appears
in the words in the corresponding multiset $\operatorname*{GR}w\in
\mathfrak{MN}^{\mathfrak{a}}$, provided that we pick one representative
of each necklace in $\operatorname*{GR}w$),
and predictably affects other quantities (for
example, the cycles of the standardization $\operatorname*{std}w$ of a word
$w\in\mathfrak{A}^{\ast}$ have the same lengths as the aperiodic necklaces in
the corresponding multiset $\operatorname*{GR}w\in\mathfrak{MN}^{\mathfrak{a}%
}$); these properties have ample applications to enumerative questions
(discussed in \cite{GesselReutenauer}).

\begin{remark}
\label{rmk.lyndon.gr.bw}
The Gessel-Reutenauer bijection relates to the \dfn{Burrows-Wheeler
transformation} (e.g., \cite[\S 2]{CrochemoreDesarmenienPerrin}). Indeed, the
latter sends an aperiodic word $w\in\mathfrak{A}^{\mathfrak{a}}$ to the word
$\operatorname*{RG}\left(  \left\{  \left[  w\right]  \right\}
_{\operatorname*{multiset}}\right)  $ obtained by applying $\operatorname*{RG}%
$ to the multiset consisting of the single aperiodic necklace $\left[
w\right]  $. This transformation is occasionally applied in (lossless) data
compression, as the word $\operatorname*{RG}\left(  \left\{  \left[  w\right]
\right\}  _{\operatorname*{multiset}}\right)  $ tends to have many strings of
consecutive equal letters when $w$ has substrings occurring multiple times
(for example, if $\mathfrak{A}=\left\{  a<b<c<d<\cdots\right\}  $ and
$w=bananaban$, then $\operatorname*{RG}\left(  \left\{  \left[  w\right]
\right\}  _{\operatorname*{multiset}}\right)  =nnbbnaaaa$), and strings of
consecutive equal letters can easily be compressed. (In order to guarantee
that $w$ can be recovered from the result, one can add a new letter $\zeta$ --
called a ``sentinel symbol'' -- to the
alphabet $\mathfrak{A}$, and apply the Burrows-Wheeler transformation to the
word $w\zeta$ instead of $w$. This also ensures that $w\zeta$ is an aperiodic
word, so the Burrows-Wheeler transformation can be applied to $w\zeta$ even
if it cannot be applied to $w$.)

Kufleitner, in \cite[\S 4]{Kufleitner}, suggests a bijective variant of the
Burrows-Wheeler transformation. In our notations, it sends a word
$w\in\mathfrak{A}^{\ast}$ to the word $\operatorname*{RG}\left(  \left\{
\left[  a_{1}\right]  ,\left[  a_{2}\right]  ,\ldots,\left[  a_{k}\right]
\right\}  _{\operatorname*{multiset}}\right)  $, where $\left(  a_{1}%
,a_{2},\ldots,a_{k}\right)  $ is the CFL factorization of $w$.
\end{remark}

For variants and generalizations of the Gessel-Reutenauer bijection,
see \cite{Kufleitner}, \cite{Steinhardt}, \cite{ShareshianWachs-euler},
\cite{DressSiebeneicher-diophant} and \cite{Reiner-signedperm}.

\subsubsection{\label{subsub.lyndon.grfun}The Gessel-Reutenauer symmetric
functions}

In this subsection, we shall study a certain family of symmetric functions.
First, we recall that every word $w\in\mathfrak{A}^{\ast}$ has a unique CFL
factorization (see Theorem \ref{thm.words.CFL}). Based on this fact, we can
make the following definition:

\begin{definition}
For the rest of Subsection \ref{subsub.lyndon.grfun}, we let $\mathfrak{A}$ be
the alphabet $\left\{  1<2<3<\cdots\right\}  $.

Let $w\in\mathfrak{A}^{\ast}$ be a word. The
\emph{CFL type}\index{CFL type of a word} of $w$ is
defined to be the partition whose parts are the positive integers $\ell\left(
a_{1}\right)  ,\ell\left(  a_{2}\right)  ,\ldots,\ell\left(  a_{k}\right)  $
(listed in decreasing order), where $\left(  a_{1},a_{2},\ldots,a_{k}\right)
$ is the CFL factorization of $w$. This CFL type is denoted by
$\operatorname*{CFLtype}w$.
\end{definition}

\begin{example}
Let $w$ be the word $213212412112$. Then, the tuple $\left(
2,132,124,12,112\right)  $ is the CFL factorization of $w$. Hence, the CFL
type of $w$ is the partition whose parts are the positive integers
$\ell\left(  2\right)  ,\ell\left(  132\right)  ,\ell\left(  124\right)
,\ell\left(  12\right)  ,\ell\left(  112\right)  $ (listed in decreasing
order). In other words, the CFL type of $w$ is the partition $\left(
3,3,3,2,1\right)  $ (since the positive integers $\ell\left(  2\right)
,\ell\left(  132\right)  ,\ell\left(  124\right)  ,\ell\left(  12\right)
,\ell\left(  112\right)  $ are $1,3,3,2,3$).
\end{example}

\begin{definition}
\label{def.lyndon.grfun.GRfun}For each word $w=\left(  w_{1},w_{2}%
,\ldots,w_{n}\right)  \in\mathfrak{A}^{\ast}$, we define a monomial
$\xx_{w}$ in $\kk \left[  \left[  \xx\right]  \right]  $
by setting $\xx_{w}=x_{w_{1}}x_{w_{2}}\cdots x_{w_{n}}$. (For example,
$\xx_{\left(  1,3,2,1\right)  }=x_{1}x_{3}x_{2}x_{1}=x_{1}^{2}%
x_{2}x_{3}$.)

For any partition $\lambda$, we define a power series
$\mathbf{GR}_{\lambda}\in \kk\left[  \left[  \xx\right]  \right]  $%
\index{$\mathbf{GR}_{\lambda}$} by
\[
\mathbf{GR}_{\lambda}=\sum\limits_{\substack{w\in\mathfrak{A}^{\ast
};\\\operatorname*{CFLtype}w=\lambda}}\xx_{w}.
\]

\end{definition}

\begin{example}
Let us compute $\mathbf{GR}_{\left(  2,1\right)  }$. Indeed, the words
$w\in\mathfrak{A}^{\ast}$ satisfying $\operatorname*{CFLtype}w=\left(
2,1\right)  $ are the words whose CFL factorization consists of two words, one
of which has length $1$ and the other has length $2$. In other words, these
words $w\in\mathfrak{A}^{\ast}$ must have the form $w=a_{1}a_{2}$ for two
Lyndon words $a_{1}$ and $a_{2}$ satisfying $a_{1}\geq a_{2}$ and $\left(
\ell\left(  a_{1}\right)  ,\ell\left(  a_{2}\right)  \right)  \in\left\{
\left(  1,2\right)  ,\left(  2,1\right)  \right\}  $. A straightforward
analysis of possibilities reveals that these are precisely the $3$-letter
words $w=\left(  w_{1},w_{2},w_{3}\right)  $ satisfying either $\left(
w_{1}<w_{2}\text{ and }w_{1}\geq w_{3}\right)  $ or $\left(  w_{1}>w_{2}\text{
and }w_{2}<w_{3}\right)  $. Hence,%
\begin{align*}
\mathbf{GR}_{\left(  2,1\right)  }  &  =\sum\limits_{\substack{w\in
\mathfrak{A}^{\ast};\\\operatorname*{CFLtype}w=\left(  2,1\right)
}}\xx_{w}=\sum\limits_{\substack{w\in\mathfrak{A}^{\ast};\\w_{1}%
<w_{2}\text{ and }w_{1}\geq w_{3}}}\xx_{w}+\sum\limits_{\substack{w\in
\mathfrak{A}^{\ast};\\w_{1}>w_{2}\text{ and }w_{2}<w_{3}}}\xx_{w}\\
&  =\sum\limits_{\substack{w\in\mathfrak{A}^{\ast};\\w_{1}<w_{2}\text{ and
}w_{1}\geq w_{3}}}\xx_{w}+\sum\limits_{\substack{w\in\mathfrak{A}%
^{\ast};\\w_{1}>w_{2}\text{ and }w_{2}<w_{3}\text{ and }w_{1}\leq w_{3}%
}}\xx_{w}+\sum\limits_{\substack{w\in\mathfrak{A}^{\ast};\\w_{1}%
>w_{2}\text{ and }w_{2}<w_{3}\text{ and }w_{1}>w_{3}}}\xx_{w}\\
&  \ \ \ \ \ \ \ \ \ \ \left(
\text{here, we have split the second sum according to the relation between }w_{1}\text{ and }w_{3}
\right) \\
&  =\sum\limits_{\substack{w\in\mathfrak{A}^{\ast};\\w_{3}\leq w_{1}<w_{2}%
}}\xx_{w}+\sum\limits_{\substack{w\in\mathfrak{A}^{\ast};\\w_{2}%
<w_{1}\leq w_{3}}}\xx_{w}+\sum\limits_{\substack{w\in\mathfrak{A}%
^{\ast};\\w_{2}<w_{3}<w_{1}}}\xx_{w}%
\end{align*}
(here, we rewrote the conditions under the summation signs). The three sums on
the right hand side are clearly quasisymmetric functions. Using
(\ref{fundamental-qsym-fn-defn-alternative}), we can rewrite them as
$L_{\left(  2,1\right)  }$, $L_{\left(  1,2\right)  }$ and $L_{\left(
1,1,1\right)  }$, respectively. Thus, we obtain%
\begin{align*}
\mathbf{GR}_{\left(  2,1\right)  }  &  =L_{\left(  2,1\right)  }+L_{\left(
1,2\right)  }+L_{\left(  1,1,1\right)  }=3M_{\left(  1,1,1\right)
}+M_{\left(  1,2\right)  }+M_{\left(  2,1\right)  }\\
&  =3m_{\left(  1,1,1\right)  }+m_{\left(  2,1\right)  }.
\end{align*}
Thus, $\mathbf{GR}_{\left(  2,1\right)  }$ is actually a symmetric function!
We shall soon (in Proposition \ref{prop.lyndon.grfun.symm}) see that this is
not a coincidence.
\end{example}

We shall now state various properties of the power series $\mathbf{GR}%
_{\lambda}$; their proofs are all part of Exercise \ref{exe.lyndon.grfun.all}.

\begin{proposition}
\label{prop.lyndon.grfun.n}Let $n$ be a positive integer. Then:

\begin{enumerate}
\item[(a)] The partition $\left(  n\right)  $ satisfies%
\[
\mathbf{GR}_{\left(  n\right)  }=\sum\limits_{\substack{w\in\mathfrak{A}%
^{n};\\w\text{ is Lyndon}}}\xx_{w}.
\]

\item[(b)] Assume that $\kk$ is a $\QQ$-algebra. Then,%
\[
\mathbf{GR}_{\left(  n\right)  }=\dfrac{1}{n}\sum_{d\mid n}\mu\left(
d\right)  p_{d}^{n/d}.
\]
Here, $\mu$ denotes the number-theoretical M\"{o}bius function (defined as in
Exercise \ref{exe.witt.ghost-equiv}), and the summation sign
``$\sum_{d\mid n}$'' is understood to range over all
\textbf{positive} divisors $d$ of $n$.
\end{enumerate}
\end{proposition}

\begin{proposition}
\label{prop.lyndon.grfun.symm}Let $\lambda$ be a partition. Then, the power
series $\mathbf{GR}_{\lambda}$ belongs to $\Lambda$.
\end{proposition}

Thus, $\left(  \mathbf{GR}_{\lambda}\right)  _{\lambda\in \Par}$
is a family of symmetric functions.\footnote{Several sources, including
\cite{GesselReutenauer}, \cite[Exercise 7.89]{Stanley} and
\cite{ElizaldeTroyka},
write $L_\lambda$ for what we call $\mathbf{GR}_{\lambda}$.
(So would we if $L_\alpha$ didn't already have another meaning here.)}
Unlike many other such families we have
studied, it is not a basis of $\Lambda$; it is not linearly independent (e.g.,
it satisfies $\mathbf{GR}_{\left(  2,1,1\right)  }=\mathbf{GR}_{\left(
4\right)  }$). Nevertheless, it satisfies a Cauchy-kernel-like
identity\footnote{Recall that $\mathfrak{L}$ denotes the set of Lyndon words
in $\mathfrak{A}^{\ast}$.}:

\begin{proposition}
\label{prop.lyndon.grfun.cauchy}Consider two countable sets of indeterminates
$\xx=\left(  x_{1},x_{2},x_{3},\ldots\right)  $ and $\mathbf{y}=\left(
y_{1},y_{2},y_{3},\ldots\right)  $.

\begin{enumerate}
\item[(a)] In the power series ring $\kk \left[  \left[  \xx%
,\mathbf{y}\right]  \right]  = \kk \left[  \left[  x_{1},x_{2}%
,x_{3},\ldots,y_{1},y_{2},y_{3},\ldots\right]  \right]  $, we have%
\[
\sum_{\lambda\in \Par}\mathbf{GR}_{\lambda}\left(
\xx\right)  p_{\lambda}\left(  \mathbf{y}\right)  =\sum_{\lambda
\in \Par}p_{\lambda}\left(  \xx\right)  \mathbf{GR}%
_{\lambda}\left(  \mathbf{y}\right)  .
\]

\item[(b)] For each word $w=\left(  w_{1},w_{2},\ldots,w_{n}\right)
\in\mathfrak{A}^{\ast}$, we define a monomial $\mathbf{y}_{w}$ in $\kk
\left[  \left[  \mathbf{y}\right]  \right]  $ by setting
$\mathbf{y}_{w}=y_{w_{1}}y_{w_{2}}\cdots y_{w_{n}}$. Then,%
\begin{align*}
\sum_{\lambda\in \Par}\mathbf{GR}_{\lambda}\left(
\xx\right)  p_{\lambda}\left(  \mathbf{y}\right)   &  =\sum
_{w\in\mathfrak{A}^{\ast}}\xx_{w}p_{\operatorname*{CFLtype}w}\left(
\mathbf{y}\right)  =\prod_{w\in\mathfrak{L}}\prod_{u\in\mathfrak{L}}\dfrac
{1}{1-\xx_{w}^{\ell\left(  u\right)  }\mathbf{y}_{u}^{\ell\left(
w\right)  }}\\
&  =\sum_{\lambda\in \Par}p_{\lambda}\left(  \xx\right)
\mathbf{GR}_{\lambda}\left(  \mathbf{y}\right)  .
\end{align*}

\end{enumerate}
\end{proposition}

The proof of this proposition rests upon the following simple
equality\footnote{Recall that $\mathfrak{L}$ denotes the set of Lyndon words
in $\mathfrak{A}^{\ast}$. Also, recall that $\mathfrak{A}=\left\{
1<2<3<\cdots\right\}  $. Thus, $p_{1}=\sum_{i\geq1}x_{i}=\sum_{a\in
\mathfrak{A}}x_{a}$.}:

\begin{proposition}
\label{prop.lyndon.grfun.CFLps}In the power series ring $\left(
\kk \left[  \left[  \xx\right]  \right]  \right)  \left[  \left[
t\right]  \right]  $, we have%
\[
\dfrac{1}{1-p_{1}t}=\prod_{w\in\mathfrak{L}}\dfrac{1}{1-\xx_{w}%
t^{\ell\left(  w\right)  }}.
\]

\end{proposition}

We can furthermore represent the symmetric functions $\mathbf{GR}_{\lambda}$
in terms of the fundamental basis $\left(  L_{\alpha}\right)  _{\alpha
\in\operatorname*{Comp}}$ of $\operatorname*{QSym}$; here, the
Gessel-Reutenauer bijection from Theorem \ref{thm.lyndon.gr.bij} reveals its
usefulness. We will use Definition~\ref{def.QSym.gamma(sigma)}.

\begin{proposition}
\label{prop.lyndon.grfun.F}Let $\lambda$ be a partition. Let $n=\left\vert
\lambda\right\vert $. Then,%
\[
\mathbf{GR}_{\lambda}=\sum\limits_{\substack{\sigma\in\Symm_n;
\\\sigma\text{ has cycle type }\lambda}}L_{\gamma\left(  \sigma\right)
}.
\]

\end{proposition}

The proof of this relies on Lemma~\ref{lem.Lalpha-std.std-L} (see
Exercise~\ref{exe.lyndon.grfun.all} below for the details).

\begin{definition}
\label{def.lyndon.grfun.type}
Let $\Symm = \bigsqcup_{n\in \NN }\Symm_{n}$ (an external
disjoint union). For each $\sigma \in \Symm$, we let
$\operatorname*{type}\sigma$ denote the cycle type of $\sigma$.
\end{definition}

\begin{proposition}
\label{prop.lyndon.grfun.cauchy2}Consider two countable sets of indeterminates
$\xx=\left(  x_{1},x_{2},x_{3},\ldots\right)  $ and $\mathbf{y}=\left(
y_{1},y_{2},y_{3},\ldots\right)  $.

In the power
series ring $\kk \left[  \left[  \xx,\mathbf{y}\right]  \right]
$, we have%
\[
\sum_{\lambda\in \Par}\mathbf{GR}_{\lambda}\left(
\xx\right)  p_{\lambda}\left(  \mathbf{y}\right)  =\sum_{\lambda
\in \Par}p_{\lambda}\left(  \xx\right)  \mathbf{GR}%
_{\lambda}\left(  \mathbf{y}\right)  =\sum_{\sigma\in\Symm}
L_{\gamma\left(  \sigma\right)  }\left(  \xx\right)
p_{\operatorname*{type}\sigma}\left(  \mathbf{y}\right)  .
\]

\end{proposition}

Let us finally give two alternative descriptions of the $\mathbf{GR}_{\lambda
}$ that do not rely on the notion of CFL factorization. First, we state a fact
that is essentially trivial:

\begin{proposition}
\label{prop.lyndon.grfun.xN}Let $N$ be a necklace. Let $w$ and $w^{\prime}$ be
two elements of $N$. Then:

\begin{enumerate}
\item[(a)] There exist words $u$ and $v$ such that $w=uv$ and $w^{\prime}=vu$.

\item[(b)] We have $\xx_{w}=\xx_{w^{\prime}}$.
\end{enumerate}
\end{proposition}

\begin{definition}
Let $N\in\mathfrak{N}$ be a necklace. Then, we define a monomial
$\xx_{N}$ in $\kk \left[  \left[  \xx\right]  \right]  $
by setting $\xx_{N}=\xx_{w}$, where $w$ is any element of $N$.
(This is well-defined, because Proposition \ref{prop.lyndon.grfun.xN}(b) shows
that $\xx_{w}$ does not depend on the choice of $w$.)
\end{definition}

\begin{definition}
Let $M$ be a finite multiset of necklaces. Then, we define a monomial
$\xx_{M}$ in $\kk \left[  \left[  \xx\right]  \right]  $
by setting $\xx_{M}=\xx_{N_{1}}\xx_{N_{2}}\cdots
\xx_{N_{k}}$, where $M$ is written in the form $M=\left\{  N_{1}%
,N_{2},\ldots,N_{k}\right\}  _{\operatorname*{multiset}}$.
\end{definition}

\begin{definition}
Let $M$ be a finite multiset of necklaces. Then, we can obtain a partition by
listing the sizes of the necklaces in $M$ in decreasing order. This partition
will be called the \emph{type}\index{type of a finite multiset of necklaces}
of $M$, and will be denoted by \dfn{$\operatorname*{type}M$}.
\end{definition}

\begin{example}
If $M=\left\{  \left[  13\right]  ,\left[  13\right]  ,\left[  2\right]
,\left[  223\right]  \right\}  _{\operatorname*{multiset}}$, then the type of
$M$ is $\left(  3,2,2,1\right)  $ (because the sizes of the necklaces in $M$
are $2,2,1,3$).
\end{example}

\begin{proposition}
\label{prop.lyndon.grfun.alt1}Let $\lambda$ be a partition. Then,%
\[
\mathbf{GR}_{\lambda}=\sum\limits_{\substack{M\in\mathfrak{MN}^{\mathfrak{a}%
};\\\operatorname*{type}M=\lambda}}\xx_{M}.
\]

\end{proposition}

This was our first alternative description of $\mathbf{GR}_{\lambda}$. Note
that it is used as a definition of $\mathbf{GR}_{\lambda}$ in \cite[(2.1)]%
{GesselReutenauer} (where $\mathbf{GR}_{\lambda}$ is denoted by $L_{\lambda}%
$). Using the Gessel-Reutenauer bijection, we can restate it as follows:

\begin{proposition}
\label{prop.lyndon.grfun.alt2}Let $\lambda$ be a partition. Then,%
\[
\mathbf{GR}_{\lambda}=\sum\limits_{\substack{w\in\mathfrak{A}^{\ast
};\\\operatorname*{type}\left(  \operatorname*{GR}w\right)  =\lambda
}}\xx_{w}.
\]

\end{proposition}

Let us finally give a second alternative description of $\mathbf{GR}_{\lambda
}$:

\begin{proposition}
\label{prop.lyndon.grfun.alt3}Let $\lambda$ be a partition. Then,%
\[
\mathbf{GR}_{\lambda}=\sum\limits_{\substack{w\in\mathfrak{A}^{\ast
};\\\operatorname*{type}\left(  \operatorname*{std}w\right)  =\lambda
}}\xx_{w}.
\]

\end{proposition}

\begin{exercise}
\label{exe.lyndon.grfun.all}Prove all statements made in Subsection
\ref{subsub.lyndon.grfun}.

[\textbf{Hint:} Here is one way to proceed:

\begin{itemize}
\item First prove Proposition \ref{prop.lyndon.grfun.CFLps}, by using the CFL
factorization to argue that both sides equal $\sum_{w\in\mathfrak{A}^{\ast}%
}\xx_{w}t^{\ell\left(  w\right)  }$.

\item Use a similar argument to derive Proposition
\ref{prop.lyndon.grfun.cauchy} (starting with part (b)).

\item Proposition \ref{prop.lyndon.grfun.xN} is almost trivial.

\item Derive Proposition \ref{prop.lyndon.grfun.alt1} from the definition of
$\mathbf{GR}_{\lambda}$ using the uniqueness of the CFL factorization.

\item Derive Proposition \ref{prop.lyndon.grfun.alt2} from Proposition
\ref{prop.lyndon.grfun.alt1} using the bijectivity of $\operatorname*{GR}$.

\item Derive Proposition \ref{prop.lyndon.grfun.alt3} from Proposition
\ref{prop.lyndon.grfun.alt2}.

\item Obtain Proposition \ref{prop.lyndon.grfun.F} by combining Proposition
\ref{prop.lyndon.grfun.alt3} with Lemma \ref{lem.Lalpha-std.std-L}.

\item Derive Proposition \ref{prop.lyndon.grfun.cauchy2} from Propositions
\ref{prop.lyndon.grfun.F} and \ref{prop.lyndon.grfun.cauchy}.

\item Derive Proposition \ref{prop.lyndon.grfun.symm} either from Proposition
\ref{prop.lyndon.grfun.alt1} or from Proposition
\ref{prop.lyndon.grfun.cauchy}. (In the latter case, make sure to work with
$\kk = \QQ$ first, and then extend to all other $\kk$, as
the proof will rely on the $\kk$-linear independence of $\left(
p_{\lambda}\right)  _{\lambda\in \Par}$, which doesn't hold for
all $\kk$.)

\item Prove Proposition \ref{prop.lyndon.grfun.n}(a) directly using the
definition of $\mathbf{GR}_{\left(  n\right)  }$.

\item Show that each positive integer $n$ satisfies $p_{1}^{n}=\sum_{d\mid
n}d\cdot\mathbf{GR}_{\left(  d\right)  }\left(  x_{1}^{n/d},x_{2}^{n/d}%
,x_{3}^{n/d},\ldots\right)  $ by taking logarithms in Proposition
\ref{prop.lyndon.grfun.CFLps}. Use this and
(\ref{eq.exe.witt.ghost-equiv.mu*1}) to prove Proposition
\ref{prop.lyndon.grfun.n}(b) recursively.
\end{itemize}

\noindent Other approaches are, of course, possible.]
\end{exercise}

\begin{remark}
\label{rmk.lyndon.grfun.char}
Let $n$ be a positive integer. The symmetric function $\mathbf{GR}_{\left(
n\right)  }$ has a few more properties:

\begin{enumerate}
\item[(a)] It is an $\NN$-linear combination of Schur functions. To
state the precise rule, we need a few more notations: A \dfn{standard
tableau} can be defined as a column-strict tableau $T$ with
$\cont \left(  T\right)  =\left(  1^{m}\right)  $, where $m$ is
the number of boxes of $T$. (That is, each of the numbers $1,2,\ldots,m$
appears exactly once in $T$, and no other numbers appear.) If $T$ is a
standard tableau with $m$ boxes, then a
\emph{descent}\index{descent of a standard tableau} of $T$ means an
$i\in\left\{  1,2,\ldots,m-1\right\}  $ such that the entry $i+1$ appears in
$T$ in a row further down than $i$ does. The
\emph{major index $\operatorname*{maj}T$}\index{major index of a standard tableau}
of a standard tableau $T$ is defined to be the sum
of its descents.\footnote{For example, the tableau
\[
\begin{array}[c]{cccc}
1 & 3 & 4 & 8\\
2 & 5 & 6 & 9\\
7 &  &  &
\end{array}
\]
is standard and has descents $1,4,6,8$ and major index $1+4+6+8=19$.} Now,
\[
\mathbf{GR}_{\left(  n\right)  }
=\sum_{\lambda\in \Par_{n}} a_{\lambda,1}s_{\lambda},
\]
where $a_{\lambda,1}$ is the number of standard tableaux $T$ of shape
$\lambda$ satisfying $\operatorname*{maj}T\equiv1\operatorname{mod}n$. (See
\cite[Exercise 7.89 (c)]{Stanley}.)

\item[(b)] Assume that $\kk = \CC$.
Recall the map $\ch : A\left(  \Symm \right) \rightarrow \Lambda$
from Theorem \ref{symmetric-group-Frobenius-map-theorem}.
Embed the cyclic group $C_n = \ZZ / n \ZZ$ as a subgroup in the
symmetric group $\Symm_n$ by identifying some generator $g$ of
$C_n$ with some $n$-cycle in $\Symm_n$. Let $\omega$ be a primitive
$n$-th root of unity in $\CC$ (for instance, $\exp\left(2\pi i/n\right)$).
Let $\gamma : C_n \rightarrow \CC$ be the character of $C_n$ that sends
each $g^i \in C_n$ to $\omega^i$. Then,
\[
\mathbf{GR}_{\left( n \right)}
= \ch \left( \Ind_{C_n}^{\Symm_n} \gamma \right)  .
\]
(See \cite[Exercise 7.89 (b)]{Stanley}.)

\item[(c)] The character $\Ind_{C_n}^{\Symm_n} \gamma$ of $\Symm_n$
is actually the character of a representation.
To construct it, set $\kk = \CC$, and
recall the notations from Exercise~\ref{exe.freelie.2}
(while keeping $\mathfrak{A} = \left\{1,2,3,\ldots\right\}$).
Let $\mathfrak{m}_n$ be the $\CC$-vector subspace of
$T\left(V\right)$ spanned by the products
$x_{\sigma\left(1\right)} x_{\sigma\left(2\right)}
 \cdots x_{\sigma\left(n\right)}$
with $\sigma \in \Symm_n$.
The symmetric group $\Symm_n$ acts on $T\left(V\right)$ by
algebra homomorphisms, with $\sigma \in \Symm_n$
sending each $x_i$ to $x_{\sigma\left(i\right)}$ when $i \leq n$
and to $x_i$ otherwise.
Both $\mathfrak{g}_n$ and $\mathfrak{m}_n$ are
$\CC \Symm_n$-submodules of $T\left(V\right)$.
Thus, so is the intersection $\mathfrak{g}_n \cap \mathfrak{m}_n$.
It is not hard to see that this intersection is spanned by all
``nested commutators''
$\left[x_{\sigma\left(1\right)},
       \left[x_{\sigma\left(2\right)},
             \left[x_{\sigma\left(3\right)},
                   \ldots
             \right]
       \right]
 \right]$ (in $T\left(V\right)$) with $\sigma \in \Symm_n$.
The character of this $\CC \Symm_n$-module
$\mathfrak{g}_n \cap \mathfrak{m}_n$ is precisely the
$\Ind_{C_n}^{\Symm_n} \gamma$ from Remark~\ref{rmk.lyndon.grfun.char}(b),
so applying the Frobenius characteristic map $\ch$ to it yields
the symmetric function $\mathbf{GR}_{\left(n\right)}$.
(See \cite[Theorem 9.41(i)]{Reutenauer}. There are similar
ways to obtain $\mathbf{GR}_\lambda$ for all $\lambda \in \Par$.)
\end{enumerate}
\end{remark}

\begin{exercise}
\label{exe.lyndon.grfun.ch-cyclic}
Prove the claim of Remark \ref{rmk.lyndon.grfun.char}(b).

[\textbf{Hint:} It helps to recall (or prove) that for any positive integer
$m$, the sum of all primitive $m$-th roots of unity in $\mathbb{C}$ is
$\mu\left(  m\right)  $.]
\end{exercise}

% [DG][v79] Added the above exercise.

The symmetric functions $\mathbf{GR}_{\lambda}$ for more general partitions
$\lambda$ can be expressed in terms of the symmetric functions
$\mathbf{GR}_{\left( n \right)}$ (which, as we recall from
Proposition \ref{prop.lyndon.grfun.n}(b), have a simple expression
in terms of the $p_m$) using the concept of
\emph{plethysm}; see \cite[Theorem 3.6]{GesselReutenauer}.

In \cite{GesselReutenauer}, Gessel and Reutenauer apply the symmetric
functions $\mathbf{GR}_{\lambda}$ to questions of permutation enumeration
via the following result%
\footnote{Proposition~\ref{prop.lyndon.grfun.count-perms}(a) is
\cite[Corollary 2.2]{GesselReutenauer};
Proposition~\ref{prop.lyndon.grfun.count-perms}(b) is
\cite[Theorem 2.1]{GesselReutenauer}.}:

\begin{proposition}
\label{prop.lyndon.grfun.count-perms}
Let $n \in \NN$. Let $\lambda \in \Par_n$ and
$\beta = \left( \beta_1, \beta_2, \ldots, \beta_k \right) \in \Comp_n$.
We shall use the notations introduced in Definition~\ref{def.QSym.D}.
Definition~\ref{def.QSym.gamma(sigma)} and
Definition~\ref{def.lyndon.grfun.type}.

\begin{enumerate}
\item[(a)] Let $\mu \in \Par_n$ be the partition
obtained by sorting the entries of $\beta$ into decreasing order. Then,%
\begin{align*}
&  \left(  \text{the number of permutations }\sigma\in\Symm_n \text{
satisfying }\operatorname*{type}\sigma=\lambda\text{ such that }\beta\text{
refines }\gamma\left(  \sigma\right)  \right) \\
&  =\left(  \text{the number of permutations }\sigma\in\Symm_n \text{
satisfying }\operatorname*{type}\sigma=\lambda\text{ and }
\Des \sigma \subset D\left(  \beta\right)  \right) \\
&  =\left(  \text{the coefficient of }x_{1}^{\beta_{1}}x_{2}^{\beta_{2}}\cdots
x_{k}^{\beta_{k}}\text{ in }\mathbf{GR}_{\lambda}\right)  =\left(  \text{the
coefficient of }\xx^{\mu}\text{ in }\mathbf{GR}_{\lambda}\right) \\
&  =\left(  \mathbf{GR}_{\lambda},h_{\mu}\right)  \ \ \ \ \ \ \ \ \ \ \left(
\text{this is the Hall inner product of }\mathbf{GR}_{\lambda}\in\Lambda\text{
and }h_{\mu}\in\Lambda\right)  .
\end{align*}

\item[(b)] Recall the ribbon diagram $\Rib\left(\beta\right)$
corresponding to the composition $\beta$
(defined as in Definition~\ref{def.QSym.D}). Then,
\begin{align*}
&  \left(  \text{the number of permutations }\sigma\in\Symm_n \text{
satisfying }\operatorname*{type}\sigma=\lambda\text{ and }\beta=\gamma\left(
\sigma\right)  \right) \\
&  =\left(  \text{the number of permutations }\sigma\in\Symm_n \text{
satisfying }\operatorname*{type}\sigma=\lambda\text{ and }
\Des \sigma =D\left(  \beta\right)  \right) \\
&  =\left(  \mathbf{GR}_{\lambda},s_{\Rib\left(\beta\right)}\right)
\ \ \ \ \ \ \ \ \ \ \left(
\text{this is the Hall inner product of }\mathbf{GR}_{\lambda}\in\Lambda\text{
and }s_{\Rib\left(\beta\right)}\in\Lambda\right)  .
\end{align*}
\end{enumerate}
\end{proposition}

\begin{exercise}
\label{exe.lyndon.grfun.count-perms}
Prove Proposition~\ref{prop.lyndon.grfun.count-perms}.

[\textbf{Hint:} Use Proposition \ref{prop.lyndon.grfun.F}, Theorem
\ref{Nsym-structure-on-ribbons-theorem}, the equality \eqref{eq.Nsym.Halpha}
and the adjointness between $\pi$ and $i$ in Corollary \ref{cor.NSym.pi}.]
\end{exercise}

By strategic application of Proposition~\ref{prop.lyndon.grfun.count-perms},
Gessel and Reutenauer arrive at several enumerative consequences, such as
the following:

\begin{itemize}
\item (\cite[Theorem 8.3]{GesselReutenauer}) If $A$ is a proper subset of
$\left\{  1,2,\ldots,n-1\right\}  $, then
\begin{align*}
&  \left(  \text{the number of permutations }\sigma\in\Symm_{n}\text{
satisfying }\left\vert \operatorname*{Fix}\sigma\right\vert =0
\text{ and } \Des \sigma =A\right) \\
&  =\left(  \text{the number of permutations }\sigma\in\Symm_{n}\text{
satisfying }\left\vert \operatorname*{Fix}\sigma\right\vert =1
\text{ and } \Des \sigma =A\right)  ,
\end{align*}
where $\operatorname*{Fix}\sigma$ denotes the set of all fixed points of a
permutation $\sigma$.
This can also be proved bijectively; such a bijective proof can be
obtained by combining \cite[Theorems 5.1 and 6.1]{DesarmenienWachs}.

\item (\cite[Theorem 9.4]{GesselReutenauer}) If $i\in\left\{  1,2,\ldots
,n-1\right\}  $, then
\[
\left(  \text{the number of }n\text{-cycles }\sigma\in\Symm_n \text{
satisfying } \Des \sigma = \left\{  i\right\}  \right)
=\sum_{d\mid\gcd\left(  n,i\right)  }\mu\left(  d\right)  \dbinom{n/d}{i/d}.
\]
Note that this also equals the number of necklaces $\left[  \left(
w_{1},w_{2},\ldots,w_{n}\right)  \right]  $ (or, equivalently, Lyndon words
$\left(  w_{1},w_{2},\ldots,w_{n}\right)  $) with $w_{1},w_{2},\ldots,w_{n}%
\in\left\{  0,1\right\}  $ and $w_{1}+w_{2}+\cdots+w_{n}=i$. This suggests
that there should be a bijection between
$\left\{ n\text{-cycles }\sigma\in\Symm_n \text{
satisfying } \Des \sigma = \left\{  i\right\} \right\}$ and the set
of such necklaces; and indeed, such a bijection can be found in
\cite[Theorem 1]{CrochemoreDesarmenienPerrin}.
\end{itemize}

\noindent See \cite{GesselReutenauer} and \cite{ElizaldeTroyka} for more
such applications.

% [DG][v78] Added the above section.

% \subsection{Application:  Multiple zeta values and
% Hoffman's stuffle conjecture}
%
% Perhaps to be filled in later....

% [DG][v36] Commented out.

\newpage

%%%%%%%%%%%%%%%%%%%%%%%%%%%%%%%%%%%%%%%%%%%%%%%%%%%%%%%%%%%%%%%%%%%%%%%%%%%%
\section{Aguiar-Bergeron-Sottile character theory Part I: $\Qsym$ as a terminal object}
\label{sect.ABS}
%%%%%%%%%%%%%%%%%%%%%%%%%%%%%%%%%%%%%%%%%%%%%%%%%%%%%%%%%%%%%%%%%%%%%%%%%%%%

It turns out that the
universal mapping property of $\Nsym$ as a free associative algebra leads
via duality to a universal property for its dual $\Qsym$, elegantly explaining several
combinatorial invariants that take the form of quasisymmetric or symmetric functions:

\begin{enumerate}
\item[$\bullet$]  Ehrenborg's quasisymmetric function of a \emph{ranked poset} \cite{Ehrenborg},
\item[$\bullet$]  Stanley's \emph{chromatic} symmetric function of a \emph{graph} \cite{Stanley-chromatic-qsym-fn},
\item[$\bullet$]  the quasisymmetric function of a \emph{matroid} considered in \cite{BilleraJiaR}.
\end{enumerate}

\subsection{Characters and the universal property}
\label{ABS-part-I-section}

\begin{definition}
Given a Hopf algebra $A$ over $\kk$, a
\emph{character}\index{character of a Hopf algebra} is an algebra morphism
$A \overset{\zeta}{\longrightarrow} \kk$, that is,
\begin{enumerate}
\item[$\bullet$]
$\zeta(\one_A)=\one_\kk$,
\item[$\bullet$]
$\zeta$ is $\kk$-linear, and
\item[$\bullet$]
$\zeta(a b)=\zeta(a)\zeta(b)$ for $a,b$ in $A$.
\end{enumerate}
\end{definition}

\begin{example}
\label{Qsym-zeta-character-example}
A particularly important character for $A=\Qsym$ is defined as follows:\footnote{We
are using the notation of Proposition~\ref{prop.Qsym.eval} here, and we are still
identifying $\Qsym$ with $\Qsym\left(\xx\right)$, where $\xx$
denotes the infinite chain $\left(x_1 < x_2 < \cdots\right)$.}
\[
\begin{array}{rcl}
\Qsym & \overset{\zeta_Q}{\longrightarrow} & \kk , \\
f(\xx)& \longmapsto & f(1,0,0,\ldots) = \left[ f(\xx) \right]_{x_1=1,x_2=x_3=\cdots=0}.
\end{array}
\]
Hence,
\[
\zeta_Q(M_\alpha)=
\zeta_Q(L_\alpha)=
\begin{cases}
1, & \text{ if } \alpha=(n) \text{ for some }n;\\
0, & \text{ otherwise.}
\end{cases}
\]
In other words, the restriction $\zeta_Q|_{\Qsym_n}$ coincides with the
functional $H_n$ in $\Nsym_n=\Hom_\kk(\Qsym_n,\kk)$:  one has for $f$ in $\Qsym_n$ that
\begin{equation}
\label{Qsym-zeta-character-example.zetaQf}
\zeta_Q(f)=(H_n,f).
\end{equation}

% [DG][v17] Added footnote mostly as a reminder.

It is worth remarking that there is nothing special about setting $x_1=1$
and $x_2=x_3=\cdots=0$:  for quasisymmetric $f$, we could have defined the same character
$\zeta_Q$ by picking any variable, say $x_n$, and sending
\[
f(\xx) \longmapsto \left[ f(\xx) \right]_{\substack{x_n=1,\text{ and}\\x_m=0 \text{ for }m\neq n}}.
\]
\end{example}

This character  $\Qsym \overset{\zeta_Q}{\longrightarrow} \kk$
has a certain universal property, known as the
\dfn{Aguiar-Bergeron-Sottile universality theorem}
(part of \cite[Theorem 4.1]{AguiarBergeronSottile}):

\begin{theorem}
\label{Qsym-as-terminal-object-theorem}
Let $A$ be a connected graded Hopf algebra, and let
$A \overset{\zeta}{\longrightarrow} \kk$ be a character.
Then, there is a unique graded Hopf morphism
$A \overset{\Psi}{\longrightarrow} \Qsym$
making the following diagram commute:
\begin{equation}
\label{ABS-terminal-morphism-diagram}
\xymatrix{
A\ar[rr]^{\Psi} \ar[dr]_{\zeta} & & \Qsym \ar[dl]^{\zeta_Q}\\
  &\kk&
}
\end{equation}
Furthermore, $\Psi$ is given by the following formula on
homogeneous elements:
\begin{equation}
\label{ABS-terminal-morphism-formula}
\Psi(a)=\sum_{\alpha \in \Comp_n} \zeta_\alpha(a) M_\alpha
\qquad \text{ for all $n \in \NN$ and $a \in A_n$,}
\end{equation}
where for $\alpha=(\alpha_1,\ldots,\alpha_\ell)$, the map $\zeta_\alpha$ is the composite
\[
A_n \overset{\Delta^{(\ell-1)}}{\longrightarrow} A^{\otimes \ell}
    \overset{\pi_\alpha}{\longrightarrow} A_{\alpha_1} \otimes \cdots \otimes A_{\alpha_\ell}
     \overset{\zeta^{\otimes \ell}}{\longrightarrow} \kk
\]
in which
$
A^{\otimes \ell}
\overset{\pi_\alpha}{\longrightarrow}
A_{\alpha_1} \otimes \cdots \otimes A_{\alpha_\ell}
$
is the canonical projection.
\end{theorem}

% [DG] Replaced "is a projection" by "is the canonical projection" so that it
% doesn't sound like a hidden existence quantor.

% [DG][v14] Removed "finite type" condition and added a footnote about this
% in the proof. It turns out that not just the proposition but also its
% proof has always worked in this generality, provided one handles
% dualization correctly (Exercise \ref{exe.dual-algebra}(e)).
%
% I also removed "finite type" from some corollaries below.

\begin{proof}
One argues that $\Psi$ is unique, and has formula \eqref{ABS-terminal-morphism-formula},
using only that $\zeta$ is $\kk$-linear and sends $\one$ to $\one$ and that
$\Psi$ is a graded \emph{$\kk$-coalgebra} map making \eqref{ABS-terminal-morphism-diagram} commute.
Equivalently, consider the adjoint \emph{$\kk$-algebra} map\footnote{Here we
are using the fact that there is a 1-to-1 correspondence between
graded $\kk$-linear maps $A \to \Qsym$ and graded $\kk$-linear maps
$\Qsym^o \to A^o$ given by $f \mapsto f^*$, and this correspondence
has the property that a given graded map $f : A \to \Qsym$ is a $\kk$-coalgebra morphism if and only if
$f^*$ is a $\kk$-algebra morphism. This is a particular case of
Exercise~\ref{exe.dual-algebra}(f).}
\[
\Nsym = \Qsym^o \overset{\Psi^*}{\longrightarrow} A^o.
\]
Commutativity of \eqref{ABS-terminal-morphism-diagram} implies that for $a$ in $A_n$,
\[
(\Psi^*(H_n),a) = (H_n,\Psi(a)) \overset{\eqref{Qsym-zeta-character-example.zetaQf}}{=} \zeta_Q( \Psi(a) ) = \zeta(a),
\]
whereas gradedness of $\Psi^*$ yields that $(\Psi^*(H_m),a) = 0$ whenever $a \in A_n$ and
$m \neq n$.
In other words, $\Psi^*(H_n)$ is the element of $A^o$ defined as the following functional on $A$:
\begin{equation}
\label{image-of-Hn-as-functional}
\Psi^*(H_n)(a) =
\begin{cases}
\zeta(a), & \text{ if } a \in A_n;\\
0, & \text{ if }a \in A_m\text{ for some }m \neq n.
\end{cases}
\end{equation}
By the universal property for $\Nsym \cong \kk\langle H_1,H_2,\ldots \rangle$
as free associative $\kk$-algebra, we see that
any choice of a $\kk$-linear map $A \overset{\zeta}{\rightarrow} \kk$
uniquely produces a
$\kk$-algebra morphism $\Psi^* : \Qsym^o \to A^o$
which satisfies \eqref{image-of-Hn-as-functional} for all $n \geq 1$.
It is easy to see that this $\Psi^*$ then automatically satisfies
\eqref{image-of-Hn-as-functional} for $n = 0$ as well if $\zeta$ sends
$\one$ to $\one$ (it is here that
we use $\zeta(\one)=\one$ and the connectedness of $A$). Hence, any
given $\kk$-linear map $A \overset{\zeta}{\rightarrow} \kk$
sending $\one$ to $\one$ uniquely produces a
$\kk$-algebra morphism $\Psi^* : \Qsym^o \to A^o$
which satisfies \eqref{image-of-Hn-as-functional} for all $n \geq 0$.
Formula
\eqref{ABS-terminal-morphism-formula} follows as
\[
\Psi(a) = \sum_{\alpha \in \Comp} (H_\alpha, \Psi(a)) \,\, M_\alpha
\]
and for a composition $\alpha=(\alpha_1,\ldots,\alpha_\ell)$, one has
\begin{align*}
(H_\alpha, \Psi(a))
=(\Psi^*(H_\alpha),a)
&=\left( \Psi^*(H_{\alpha_1}) \cdots \Psi^*(H_{\alpha_\ell}), a \right)\\
%&=( \Psi^*(H_{\alpha_1}) \otimes \cdots \otimes \Psi^*(H_{\alpha_\ell}),( \pi_\alpha \circ \Delta^{(\ell-1)})(a))
&= \left( \Psi^*(H_{\alpha_1}) \otimes \cdots \otimes \Psi^*(H_{\alpha_\ell}), \Delta^{(\ell-1)}(a) \right) \\
&\overset{\eqref{image-of-Hn-as-functional}}{=} \left(\zeta^{\otimes \ell} \circ \pi_\alpha \right) \left(\Delta^{(\ell-1)}(a)\right)
=\zeta_{\alpha}(a),
\end{align*}
where the definition of $\zeta_\alpha$ was used in the last equality.

% [DG][v14] I sprayed some more detail into this proof.

We wish to show that if, in addition, $A$ is a Hopf algebra and $A \overset{\zeta}{\longrightarrow} \kk$ is a character (i.e., an algebra morphism), then  $A \overset{\Psi}{\longrightarrow} \Qsym$ will be an algebra morphism, that is, the two maps $A \otimes A \longrightarrow \Qsym$
given by $\Psi \circ m$ and $m \circ (\Psi \otimes \Psi)$ coincide.  To see this, consider
these two diagrams having the two maps in question as the composites of their top rows:
\begin{equation}
\xymatrix{
A \otimes A \ar[r]^-{m}\ar[dr]_{\zeta \otimes \zeta} & A\ar[r]^\Psi \ar[d]_{\zeta} &\Qsym \ar[dl]^{\zeta_Q}\\
  &\kk&
}
\qquad
\xymatrix@C=55pt{
A \otimes A \ar[r]^{\Psi \otimes \Psi} \ar[dr]_{\zeta \otimes \zeta}
  & \Qsym^{\otimes 2} \ar[r]^-{m} \ar[d]_{\zeta_Q \otimes \zeta_Q}
    &\Qsym \ar[dl]^{\zeta_Q}\\
  &\kk&
}
\end{equation}
The fact that $\zeta, \zeta_Q$ are algebra morphisms makes the above diagrams commute, so that
applying the uniqueness in the first part of the proof to the
character $A \otimes A \overset{\zeta \otimes \zeta}{\longrightarrow} \kk$ proves
the desired equality $\Psi \circ m=m \circ (\Psi \otimes \Psi)$.
\end{proof}

\begin{remark}
\label{Sym-as-terminal-object-remark}
When one assumes in addition that $A$ is cocommutative,
it follows that the image of $\Psi$ will lie
in the subalgebra $\Lambda \subset \Qsym$, e.g. from the explicit
formula \eqref{ABS-terminal-morphism-formula} and the fact
that one will have $\zeta_\alpha = \zeta_\beta$ whenever $\beta$
is a rearrangement of $\alpha$.  In other words, the character
$\Lambda \overset{\zeta_\Lambda}{\longrightarrow} \kk$ defined by
restricting $\zeta_Q$ to $\Lambda$, or by
\[
\zeta_\Lambda(m_\lambda)=
\begin{cases}
1, &\text{ if }\lambda=(n) \text{ for some }n;\\
0, &\text{ otherwise,}
\end{cases}
\]
has a universal property as terminal object
with respect to characters on cocommutative %co- or 
Hopf algebras.
\end{remark}

% [DG][v78] Removed "co- or" from the previous remark.
% Restated the theorem to uncrowd the first sentence,
% while adding the fiollowing sentence in order to
% fix the "induced" terminology:

The graded Hopf morphism $\Psi$ in
Theorem~\ref{Qsym-as-terminal-object-theorem}
will be called the \dfn{map $A \to \Qsym$ induced
by the character $\zeta$}.

We close this section by discussing a well-known polynomiality and reciprocity
phenomenon; see, e.g., Humpert and Martin \cite[Prop. 2.2]{HumpertMartin},
Stanley \cite[\S 4]{Stanley-chromatic-qsym-fn}.

\begin{definition}
The \dfn{binomial Hopf algebra} (over the commutative ring $\kk$)
is the polynomial algebra $\kk\left[  m\right]  $ in a single
variable $m$, with a Hopf algebra structure transported from the symmetric
algebra $\Sym\left(  \kk^1\right)  $ (which is a Hopf
algebra by virtue of Example~\ref{symmetric-algebra-as-bialgebra-example},
applied to $V=\kk^1$) along the
isomorphism $\Sym\left(  \kk^1\right)  \rightarrow
\kk\left[  m\right]  $ which sends the standard basis element of
$\kk^1$ to $m$. Thus the element $m$ is primitive; that is,
$\Delta m = \one \otimes m + m \otimes \one$ and $S(m) = -m$.
As $S$ is an algebra anti-endomorphism by
Proposition~\ref{antipodes-are-antiendomorphisms} and $\kk[m]$ is commutative,
one has $S(g)(m) = g(-m)$ for all polynomials $g(m)$ in $\kk[m]$.
\end{definition}

% [DG][v3] Not the best way to state it, I fear.

\begin{definition}
\label{principal-specialization-at-1-defn}
For an element $f(\xx)$ in $\Qsym$ and a nonnegative integer $m$,
let $\ps^1(f)(m)$ denote the element of $\kk$
obtained by \dfn{principal specialization at $q=1$}
\begin{align*}
\ps^1(f)(m)&=\left[ f(\xx)
\right]_{\substack{x_1=x_2=\cdots =x_m=1,\\x_{m+1}=x_{m+2}= \cdots =0}} \\
&=f(\underbrace{1,1,\ldots,1}_{m\text{ ones}},0,0,\ldots).
\end{align*}
\end{definition}

\begin{proposition}
\label{zeta-polynomiality-proposition}
Assume that $\QQ$ is a subring of $\kk$.
The map $\ps^1$ has the following properties.
\begin{enumerate}
\item[(i)]
Let $f \in \Qsym$.
There is a unique
polynomial in $\kk[m]$ which agrees for each nonnegative integer $m$
with $\ps^1(f)(m)$, and which, by abuse of notation, we will also denote
$\ps^1(f)(m)$.  If $f$ lies in $\Qsym_n$, then $\ps^1(f)(m)$ is a polynomial
of degree at most $n$, taking these values on $M_\alpha, L_\alpha$
for $\alpha=(\alpha_1,\ldots,\alpha_\ell)$ in $\Comp_n$:
\begin{align*}
\ps^1(M_\alpha)(m) &=\binom{m}{\ell}, \\
\ps^1(L_\alpha)(m) &=\binom{m-\ell+n}{n}.
\end{align*}
\item[(ii)] The map $\Qsym \overset{\ps^1}{\longrightarrow} \kk[m]$
is a Hopf morphism into the binomial Hopf algebra.
\item[(iii)] For all $m$ in $\ZZ$ and $f$ in $\Qsym$ one has
\[
\zeta_Q^{\star m}(f) = \ps^1(f)(m).
\]
In particular, one also has
\[
\zeta_Q^{\star (-m)}(f) = \ps^1(S(f))(m) =\ps^1(f)(-m).
\]
\item[(iv)]
For a graded Hopf algebra $A$ with a character
$A \overset{\zeta}{\longrightarrow} \kk$, and any element $a$ in $A_n$,
the polynomial $\ps^1(\Psi(a))(m)$ in $\kk[m]$ has degree at most $n$,
and when specialized to $m$ in $\ZZ$ satisfies
\[
\zeta^{\star m}(a) = \ps^1(\Psi(a))(m).
\]
\end{enumerate}
\end{proposition}

% [DG] I required characteristic $0$ to make sense of the polynomial
% $\binom{m}{\ell}$.

\begin{proof}
To prove assertion (i), note that one has
\begin{align*}
\ps^1(M_\alpha)(m)
=M_\alpha(1,1,\dots,1,0,0,\ldots)
&=\sum_{1 \leq i_1 < \cdots < i_\ell \leq m}
 \left[ x_{i_1}^{\alpha_1} \cdots x_{i_\ell}^{\alpha_\ell} \right]_{x_j=1}
= \binom{m}{\ell},\\
\ps^1(L_\alpha)(m)
=L_\alpha(1,1,\dots,1,0,0,\ldots)
&=\sum\limits_{\substack{1 \leq i_1 \leq \cdots \leq i_n \leq m : \\ i_k < i_{k+1} \text{ if }k \in D(\alpha)}}  \left[ x_{i_1} \cdots x_{i_n} \right]_{x_j=1}\\
&=|\{1 \leq j_1 \leq j_2 \leq \cdots \leq j_n \leq m-\ell+1\}|
=\binom{m-\ell+n}{n}.
\end{align*}
As $\{M_\alpha\}_{\alpha \in \Comp_n}$ form a basis for $\Qsym_n$, and
$\binom{m}{\ell}$ is a polynomial function in $m$ of degree $\ell (\leq n)$,
one concludes that for $f$ in $\Qsym_n$ one has that $\ps^1(f)(m)$ is a polynomial function in $m$ of
degree at most $n$. The polynomial giving rise to this function is unique,
since infinitely many of its values are fixed.

To prove assertion (ii), note that $\ps^1$ is an algebra morphism because it
is an evaluation homomorphism.
To check that it is a coalgebra morphism, it suffices to check
$\Delta \circ \ps^1 = (\ps^1 \otimes \ps^1) \circ \Delta$ on each $M_\alpha$ for
$\alpha=(\alpha_1,\ldots,\alpha_\ell)$ in $\Comp_n$. Using the  Vandermonde
summation $\binom{A+B}{\ell} = \sum_{k} \binom{A}{k} \binom{B}{\ell-k}$, one has
\[
(\Delta \circ \ps^1)(M_\alpha)
= \Delta \binom{m}{\ell}
= \binom{ m \otimes \one + \one \otimes m }{\ell}
=\sum_{k=0}^\ell \binom{ m \otimes \one }{k}  \binom{ \one \otimes m }{\ell-k}
=\sum_{k=0}^\ell \binom{ m }{k} \otimes \binom{  m }{\ell-k}
\]
while at the same time
\[
\left((\ps^1 \otimes \ps^1) \circ \Delta\right)(M_\alpha)
= \sum_{k=0}^{\ell} \ps^1(M_{(\alpha_1,\ldots,\alpha_k)})
                     \otimes \ps^1(M_{(\alpha_{k+1},\ldots,\alpha_\ell)})
=\sum_{k=0}^{\ell} \binom{m}{k} \otimes \binom{m}{\ell-k}.
\]
Thus $\ps^1$ is a bialgebra morphism, and hence also a Hopf morphism, by
Corollary~\ref{cor.bialg-mor-is-Hopf}.

For assertion (iii), first assume $m$ lies in $\{0,1,2,\ldots\}$.
Since $\zeta_Q(f)=f(1,0,0,\ldots)$, one has
\begin{align*}
\zeta_Q^{\star m}(f)
&= \zeta_Q ^{\otimes m}\circ \Delta^{(m-1)} f(\xx)
=  \zeta_Q ^{\otimes m} \left( f(\xx^{(1)},\xx^{(2)},\ldots,\xx^{(m)}) \right) \\
&=  \left[ f(\xx^{(1)},\xx^{(2)},\ldots,\xx^{(m)})
         \right]_{\substack{x^{(1)}_1=x^{(2)}_1=\cdots=x^{(m)}_1=1,
                    \\x^{(j)}_2=x^{(j)}_3=\cdots=0 \text{ for all }j}}\\
&=f(1,0,0,\ldots,1,0,0,\ldots,\cdots,1,0,0,\ldots)
=f(\underbrace{1,1,\ldots,1}_{m\text{ ones}},0,0,\ldots) = \ps^1(f)(m).
\end{align*}
\footnote{See Exercise~\ref{exe.zeta-polynomiality.altpf} for an alternative
way to prove this, requiring less thought to verify its soundness.}
But then Proposition~\ref{antipodes-give-convolution-inverses}(a) also implies
\begin{align*}
\zeta_Q^{\star (-m)}(f)
&=\left( \zeta_Q^{\star (-1)} \right)^{\star m}(f)
= \left( \zeta_Q \circ S \right)^{\star m}(f)
= \zeta_Q^{\star m}(S(f))\\
&= \ps^1(S(f))(m)
= S(\ps^1(f))(m)
=\ps^1(f)(-m).
\end{align*}
For assertion (iv), note that
\[
\zeta^{\star m}(a)
= (\zeta_Q \circ \Psi)^{\star m}(a)
= (\zeta_Q^{\star m})(\Psi(a))
= \ps^1(\Psi(a))(m),
\]
where the three equalities come from
\eqref{ABS-terminal-morphism-diagram},
Proposition~\ref{antipodes-give-convolution-inverses}(a),
and assertion (iii) above, respectively.
\end{proof}

% [DG][v17] I added the alternative proof in the footnote to avoid
% formalizing the argument you gave (which is rather tricky imho).
% Maybe better as an exercise?

% [DG][v19] Now it's an exercise.

\begin{remark}
Aguiar, Bergeron and Sottile give
a very cute (third) proof
of the $\Qsym$ antipode formula
Theorem~\ref{Qsym-antipode-on-monomials},
via Theorem~\ref{Qsym-as-terminal-object-theorem},
in \cite[Example 4.8]{AguiarBergeronSottile}.
They apply Theorem~\ref{Qsym-as-terminal-object-theorem}
to the \emph{coopposite coalgebra} $\Qsym^{cop}$ and its
character $\zeta_Q^{\star (-1)}$.  One can show that the map
$\Qsym^{cop} \overset{\Psi}{\rightarrow} \Qsym$ induced by $\zeta_Q^{\star (-1)}$
is $\Psi=S$, the antipode of $\Qsym$, because
$S : \Qsym \to \Qsym$ is
a coalgebra anti-endomorphism (by Exercise~\ref{exe.coalg-anti})
satisfying
$\zeta_Q^{\star (-1)} = \zeta_Q \circ S$.  They then use
the formula \eqref{ABS-terminal-morphism-formula}
for $\Psi=S$ (together with the polynomiality
Proposition~\ref{zeta-polynomiality-proposition})
to derive Theorem~\ref{Qsym-antipode-on-monomials}.
\end{remark}

% [DG][v14] Replaced
% $\Qsym \overset{\Psi}{\rightarrow} \Qsym^{cop}$
% by
% $\Qsym^{cop} \overset{\Psi}{\rightarrow} \Qsym$.
% Replaced ^{\star ^{-1}} by ^{\star (-1)}.

\begin{exercise}
\label{exe.zeta-polynomiality.altpf}
Show that
$\zeta_Q^{\star m}(f) = \ps^1(f)(m)$ for all $f \in \Qsym$ and
$m \in \left\{0,1,2,\ldots\right\}$. (This was already
proven in Proposition~\ref{zeta-polynomiality-proposition}(iii);
give an alternative proof using
Proposition~\ref{Qsym-coproduct-on-monomials}.)
\end{exercise}

\subsection{Example: Ehrenborg's quasisymmetric function of a ranked poset}

Here we consider incidence algebras, coalgebras and Hopf algebras generally,
and then particularize to the case of graded posets, to recover Ehrenborg's interesting quasisymmetric function invariant
via Theorem~\ref{Qsym-as-terminal-object-theorem}.

\subsubsection{Incidence algebras, coalgebras, Hopf algebras}

\begin{definition}
Given a family $\PPP$ of finite partially ordered sets $P$,
let $\kk[\PPP]$ denote the free $\kk$-module
whose basis consists of symbols $[P]$ corresponding
to isomorphism classes of posets $P$ in $\PPP$.

We will assume throughout that each $P$ in $\PPP$ is
\emph{bounded}\index{bounded poset},
that is, it has a unique minimal element $\hat{0}:=\hat{0}_P$ and a unique
maximal element $\hat{1}:=\hat{1}_P$.  In particular, $P \neq \varnothing$,
although it is allowed that $|P|=1$, so that $\hat{0}=\hat{1}$;
denote this isomorphism class of posets with one element by $[o]$.

If $\PPP$ is closed under taking intervals
\[
[x,y]:=[x,y]_P:=\{ z \in P: x \leq_P z \leq_P y \} ,
\]
then one can easily see that the following coproduct and counit
endow $\kk[\PPP]$ with the structure of a coalgebra,
called the
\emph{(reduced) incidence coalgebra}\index{reduced incidence coalgebra}:
\begin{align*}
\Delta [P]&:= \sum_{x \in P} [\hat{0},x] \otimes [x,\hat{1}],\\
\epsilon [P]&:=
\begin{cases}
1, & \text{ if }|P|=1 ; \\
0, & \text{ otherwise.}
\end{cases}
\end{align*}
The dual algebra $\kk[\PPP]^\ast$ is generally called the
\emph{reduced incidence algebra (modulo isomorphism)}
for the family $\PPP$ (see, e.g., \cite{Schmitt-antipodes}).
It contains the important element
$\kk[\PPP] \overset{\zeta}{\longrightarrow} \kk$,
called the
\emph{$\zeta$-function}\index{$\zeta$-function in an incidence algebra}
that takes the value $\zeta[P]=1$ for all $P$.

% [DG][v17] Wait: haven't you defined the graded dual only for graded
% Hopf algebras? This one is only filtered. Should I add a definition
% of filtered Hopf algebras? Or do you want the full dual? Or
% Sweedler's zero-dual?

% [DG][v19] Correction made.

If $\PPP$ (is not empty and) satisfies the further property of being
\emph{hereditary}\index{hereditary class of posets}
in the sense that for every $P_1,P_2$ in $\PPP$, the
\dfn{Cartesian product poset} $P_1 \times P_2$ with componentwise
partial order is also in $\PPP$, then one can check that
the following product and unit endow $\kk[\PPP]$ with
the structure of a (commutative) algebra:
\begin{align*}
[P_1] \cdot [P_2]&:=m([P_1] \otimes [P_2]) := [P_1 \times P_2],\\
\one_{\kk[\PPP]}&:=[o].
\end{align*}
\end{definition}

\begin{proposition}
\label{incidence-coalgebras-are-Hopf}
For any hereditary family $\PPP$ of finite posets,
$\kk[\PPP]$ is a bialgebra, and even a Hopf algebra
with antipode $S$ given as in \eqref{Takeuchi-formula} (Takeuchi's formula):
\[
S[P]=\sum_{k \geq 0} (-1)^k
     \sum_{\hat{0}=x_0 < \cdots <x_k=\hat{1}}
        [x_0,x_1] \cdots [x_{k-1},x_k].
\]
\end{proposition}
\begin{proof}
Checking the commutativity of the pentagonal diagram in
\eqref{bialgebra-diagrams}
amounts to the fact that, for any
$(x_1,x_2)  <_{P_1 \times P_2} (y_1,y_2)$,
one has a poset isomorphism
\[
\left[ (x_1,x_2) \,\, , \,\, (y_1,y_2) \right]_{P_1 \times P_2}
\cong [x_1,y_1]_{P_1} \times [x_2,y_2]_{P_2}.
\]
Commutativity of the remaining diagrams in
\eqref{bialgebra-diagrams} is straightforward, and so $\kk[\PPP]$
is a bialgebra.  But then Remark~\ref{non-graded-Takeuchi-remark}
implies that it is a Hopf algebra, with antipode $S$ as in \eqref{Takeuchi-formula},
because the map $f:=\id_{\kk[\PPP]} - u\epsilon$ (sending the class $[o]$ to $0$,
and fixing all other $[P]$) is locally $\star$-nilpotent:
\[
f^{\star k}[P]
=\sum_{\hat{0}=x_0 < \cdots <x_{k}=\hat{1}}
[x_0,x_1] \cdots [x_{k-1},x_{k}]
\]
will vanish due to an empty sum whenever $k$ exceeds the maximum length of
a chain in the finite poset $P$.
\end{proof}

It is perhaps worth remarking how this generalizes the M\"obius function
formula of P. Hall.  Note that the zeta function
$\kk[\PPP] \overset{\zeta}{\longrightarrow} \kk$ is a \emph{character},
that is, an algebra morphism.
Proposition~\ref{antipodes-give-convolution-inverses}(a)
then tells us that $\zeta$ should have a convolutional
inverse $\kk[\PPP] \overset{\mu=\zeta^{\star -1}}{\longrightarrow} \kk$,
traditionally called the
\emph{M\"obius function}\index{M\"obius function of a poset}, with
the formula
$\mu = \zeta^{\star -1} = \zeta \circ S$.  Rewriting this via
the antipode formula for $S$ given
in Proposition~\ref{incidence-coalgebras-are-Hopf} yields
P. Hall's formula.

\begin{corollary}
\label{P.Hall-formula}
For a finite bounded poset $P$, one has
\[
\mu[P] =
\sum_{k \geq 0} (-1)^k
|\{ \text{chains }\hat{0}=x_0 < \cdots <x_k=\hat{1}\text{ in }P\}|.
\]
\end{corollary}

We can also notice that $S$ is an algebra anti-endomorphism
(by Proposition~\ref{antipodes-are-antiendomorphisms}), thus an
algebra endomorphism (since $\kk[\PPP]$ is commutative, so
Exercise~\ref{exe.comm-cocomm.anti}(a) shows that the
algebra anti-endomorphisms of $\kk[\PPP]$ are
the same as the algebra endomorphisms of $\kk[\PPP]$). Hence,
$\mu = \zeta \circ S$ is a composition of two algebra
homomorphisms, thus an algebra homomorphism itself. We
therefore obtain the following classical fact:

\begin{corollary}
For two finite bounded posets $P$ and $Q$, we have
$\mu[P \times Q] = \mu[P] \cdot \mu[Q]$.
\end{corollary}

\subsubsection{The incidence Hopf algebras for ranked posets and Ehrenborg's function}

\begin{definition}
Take $\PPP$ to be the class of bounded
\emph{ranked}\index{ranked poset} finite
posets $P$, that is, those for which all maximal chains from
$\hat{0}$ to $\hat{1}$ have the same length $r(P)$.  This is a hereditary
class, as it implies that any interval is $[x,y]_P$ is also ranked,
and the product of two bounded ranked posets is also bounded and ranked.
It also uniquely defines a
\emph{rank function}\index{rank function of a poset}
$P \overset{r}{\longrightarrow} \NN$
in which $r(\hat{0})=0$ and $r(x)$ is the length of any maximal chain
from $\hat{0}$ to $x$.
\end{definition}

\begin{example}
\label{polytope-example}
Consider a pyramid with apex vertex $a$ over a square base with vertices $b,c,d,e$:
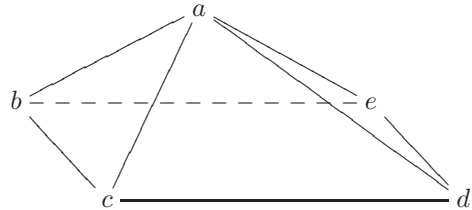
\[
\xymatrix{
 & &a\ar@{-}[dll]\ar@{-}[ddl]\ar@{-}[drr]\ar@{-}[ddrrr]&&& & \\
b\ar@{--}[rrrr]\ar@{-}[dr]& & &&e\ar@{-}[dr]& \\
 &c\ar@{-}[rrrr]& & &&d
}
\]
Ordering its faces by inclusion gives a bounded ranked poset $P$,
where the rank of an element is one more than the dimension
of the face it represents:
\[
\xymatrix{
 & & & & & & & & \text{rank:}\\
   &   &   &abcd\ar@{-}[drrrr]\ar@{-}[d]\ar@{-}[dl]\ar@{-}[dll]\ar@{-}[dlll]&  &  &  & & 4\\
abc\ar@{-}[d]\ar@{-}[dr]&acd\ar@{-}[d]\ar@{-}[dr]&ade\ar@{-}[d]\ar@{-}[dr]&abe\ar@{-}[d]\ar@{-}[dlll] &  &  &  &bcde\ar@{-}[d]\ar@{-}[dl]\ar@{-}[dll]\ar@{-}[dlll]& 3\\
ab &ac &ad &ae  &be\ar@{-}[ul]&bc\ar@{-}[ulllll]&cd\ar@{-}[ulllll]&de\ar@{-}[ulllll]& 2\\
a\ar@{-}[u]\ar@{-}[ur]\ar@{-}[urr]\ar@{-}[urrr]  &   &   &    &b\ar@{-}[u]\ar@{-}[ur]\ar@{-}[ullll] &c\ar@{-}[u]\ar@{-}[ur]\ar@{-}[ullll] &d\ar@{-}[u]\ar@{-}[ur]\ar@{-}[ullll] &e\ar@{-}[u]\ar@{-}[ulll]\ar@{-}[ullll] & 1\\
   &   &   &    &\varnothing\ar@{-}[ullll]\ar@{-}[u]\ar@{-}[ur]\ar@{-}[urr]\ar@{-}[urrr]& & & & 0\\
}
\]
\end{example}

\begin{definition}
\dfn{Ehrenborg's quasisymmetric function} $\Psi[P]$ for a bounded ranked poset $P$
is the image of $[P]$ under the map
$\kk[\PPP] \overset{\Psi}{\longrightarrow} \Qsym$
induced by the zeta function
$\kk[\PPP] \overset{\zeta}{\longrightarrow} \kk$ as
a character, via Theorem~\ref{Qsym-as-terminal-object-theorem}.
\end{definition}

The quasisymmetric function $\Psi[P]$ captures several interesting
combinatorial invariants of $P$; see Stanley \cite[Chap. 3]{Stanley}
for more background on these notions.

\begin{definition}
Let $P$ be a bounded ranked poset $P$ of rank $r(P):=r(\hat{1})$.
Define its \emph{rank-generating function}\index{rank-generating function of a poset}\index{$RGF(P,q)$}
\[
RGF(P,q):=\sum_{p \in P} q^{r(p)} \in \ZZ\left[q\right],
\]
its \emph{characteristic polynomial}\index{characteristic polynomial of a poset}\index{$\chi(P,q)$}
\[
\chi(P,q):=\sum_{p \in P} \mu(\hat{0},p)q^{r(p)} \in \ZZ\left[q\right]
\]
(where $\mu(u,v)$ is shorthand for $\mu([u,v])$),
and its \emph{zeta polynomial}\index{zeta polynomial of a poset}\index{$Z(P,m)$}
\begin{align}
\label{zeta-polynomial-def}
Z(P,m)&=|\{\text{multichains }
\hat{0} \leq_P  p_1 \leq_P \cdots \leq_P p_{m-1}  \leq_P \hat{1} \} |\\
 &= \sum_{s=0}^{r(P)-1} \binom{m}{s+1}
          |\{\text{chains }\hat{0} < p_1 < \cdots < p_s < \hat{1}\}| \in \QQ\left[m\right]
\label{zeta-polynomial-as-chain-sum}
\end{align}
\footnote{Actually, \eqref{zeta-polynomial-as-chain-sum} is false
if $\left|P\right| = 1$ (but only then). We use
\eqref{zeta-polynomial-def} to define $Z(P,m)$ in this case.}.
Also, for each subset $S \subset \{1,2,\ldots,r(P)-1\}$,
define
the \emph{flag number}\index{flag number of a poset}\index{$f_S$}
$f_S$ of $P$ by
\[
f_S=|\{\text{chains } \hat{0} <_P p_1 <_P \cdots <_P p_s  <_P \hat{1}
\text{ with }\{r(p_1),\ldots,r(p_s)\}=S\}|.
\]
These flag numbers are the components of the
\emph{flag $f$-vector}\index{flag $f$-vector of a poset}
$(f_S)_{S \subset [r-1]}$ of $P$.
Further define the \emph{flag $h$-vector}\index{flag $h$-vector of a poset}
$(h_T)_{T \subset [r-1]}$ of $P$,
whose entries $h_T$ are given by
$f_S=\sum_{ T \subset S} h_T$, or, equivalently%
\footnote{The equivalence follows from inclusion-exclusion
(more specifically, from
the converse of Lemma~\ref{lem.moebius-bool.subsets}(a)).},
by
$h_S=\sum_{ T \subset S} (-1)^{|S \setminus T|} f_T$.
\end{definition}

% [DG][v13] Replaced "r" by "r(P)" in
% "and for a subset  $S \subset \{1,2,\ldots,r-1\}$,".

% [DG][v19] Added footnote about border case when $P$ has
% one element only.

% [DG][v80] Replaced run-on sentence by several shorter ones.

\begin{example}
For the poset $P$ in Example~\ref{polytope-example}, one
has $RGF(P,q)=1+5q+8q^2+5q^3+q^4$.  Since $P$ is the poset of
faces of a polytope, the M\"obius function values for its
intervals are easily predicted:  $\mu(x,y)=(-1)^{r[x,y]}$,
that is, $P$ is an \dfn{Eulerian ranked poset}; see Stanley
\cite[\S 3.16]{Stanley}.  Hence
its characteristic polynomial is trivially related to the
rank generating function, sending $q \mapsto -q$, that is,
\[
\chi(P,q)= RGF(P,-q)=1-5q+8q^2-5q^3+q^4.
\]
Its flag $f$-vector and $h$-vector entries are given in the following
table.
\vskip.1in
\begin{tabular}{|c|c|cl|}\hline
$S$           & $f_S$ &       &$h_S$  \\ \hline\hline
$\varnothing$ & $1$   &       &$1$ \\ \hline
$\{1\}$       & $5$   & $5-1=$&$4$ \\ \hline
$\{2\}$       & $8$   & $8-1=$&$7$ \\ \hline
$\{3\}$       & $5$   & $5-1=$&$4$ \\ \hline
$\{1,2\}$     & $16$  & $16-(5+8)+1=$&$4$ \\ \hline
$\{1,3\}$     & $16$  & $16-(5+5)+1=$&$7$ \\ \hline
$\{2,3\}$     & $16$  & $16-(5+8)+1=$&$4$ \\ \hline
$\{1,2,3\}$   & $32$  & $32-(16+16+16)+(5+8+5)-1=$&$1$ \\ \hline
\end{tabular}
\vskip.1in
\noindent
and using \eqref{zeta-polynomial-as-chain-sum}, its zeta polynomial is
\[
Z(P,m)=1 \binom{m}{1} + (5+8+5) \binom{m}{2} + (16+16+16) \binom{m}{3} +
32 \binom{m}{4}=\frac{m^2(2m-1)(2m+1)}{3}.
\]
\end{example}

\begin{theorem}
\label{Ehrenborg's-theorem}
Assume that $\QQ$ is a subring of $\kk$.
Ehrenborg's quasisymmetric function $\Psi[P]$ for a bounded ranked poset $P$
encodes
\begin{enumerate}
\item[(i)]
the flag $f$-vector entries $f_S$ and flag $h$-vector entries $h_S$
as its $M_\alpha$ and $L_\alpha$ expansion coefficients\footnote{In fact,
Ehrenborg \emph{defined} $\Psi[P]$  in \cite[Defn. 4.1]{Ehrenborg}
via this $M_\alpha$ expansion, and then showed that it gave a Hopf morphism.} :
\[
\Psi[P] = \sum_\alpha f_{D(\alpha)}(P) \,\, M_\alpha
= \sum_\alpha h_{D(\alpha)}(P) \,\, L_\alpha ,
\]
\item[(ii)]
the zeta polynomial as the specialization from
Definition~\ref{principal-specialization-at-1-defn}
\[
Z(P,m)=\ps^1(\Psi[P])(m)
=\left[ \Psi[P] \right]_{\substack{x_1=x_2=\cdots =x_m=1,\\ x_{m+1}=x_{m+2}= \cdots =0}} ,
\]
\item[(iii)]
the rank-generating function as the specialization
\[
RGF(P,q)
=\left[ \Psi[P] \right]_{\substack{x_1=q,x_2=1,\\x_3=x_4= \cdots =0}} ,
\]
\item[(iv)]
the characteristic polynomial as the convolution
\[
\chi(P,q)= ((\psi_q \circ S) \star \zeta_Q) \circ \Psi[P] ,
\]
where $\Qsym \overset{\psi_q}{\longrightarrow} \kk[q]$
maps $f(\xx) \longmapsto f(q,0,0,\ldots)$.
\end{enumerate}
\end{theorem}
\begin{proof}
In assertion (i), the
expansion $\Psi[P]=\sum_\alpha f_{D(\alpha)}(P) \, M_\alpha$ is
\eqref{ABS-terminal-morphism-formula}, since
$\zeta_\alpha[P]=f_{D(\alpha)}(P)$.  The $L_\alpha$ expansion
follows from this, as
$L_\alpha=\sum_{\beta: D(\beta) \supset D(\alpha)} M_\beta$
and $f_S(P)=\sum_{T \subset S} h_T$.

Assertion (ii) is immediate from
Proposition~\ref{zeta-polynomiality-proposition}(iv),
since $Z(P,m)=\zeta^{\star m}[P]$.

Assertion (iii) can be deduced from assertion (i),
but it is perhaps more fun and in the spirit of things to proceed
as follows. Note that
$\psi_q(M_\alpha)=q^{n}$ for $\alpha=(n)$, and $\psi_q(M_\alpha)$ vanishes for
all other $\alpha \neq (n)$ in $\Comp_n$.  Hence
for a bounded ranked poset $P$ one has
\begin{equation}
\label{q-rank-functional-equation}
(\psi_q \circ \Psi)[P]= q^{r(P)}.
\end{equation}
But if we treat $\zeta_Q : \Qsym \to \kk$ as a map
$\Qsym \to \kk\left[q\right]$, then
\eqref{pre-and-post-composition-in-convolution}
(applied to $\kk\left[\PPP\right]$, $\Qsym$, $\kk\left[q\right]$,
$\kk\left[q\right]$, $\Psi$, $\id_{\kk\left[q\right]}$,
$\psi_q$ and $\zeta_Q$
instead of $C$, $C'$, $A$, $A'$, $\gamma$, $\alpha$,
$f$ and $g$) shows that
\begin{equation}
\left( \psi_q \star \zeta_Q \right) \circ \Psi
= \left( \psi_q \circ \Psi \right) \star \left( \zeta_Q \circ \Psi \right) ,
\label{pf.Ehrenborg's-theorem.compose-star}
\end{equation}
since $\Psi : \kk\left[\PPP\right] \to \Qsym$ is a $\kk$-coalgebra
homomorphism.
Consequently, one can compute
\begin{align*}
RGF(P,q) & = \sum_{p \in P} q^{r(p)} \cdot 1 = \sum_{p \in P} q^{r([\hat{0},p])} \cdot \zeta[p,\hat{1}]
\overset{\substack{\eqref{q-rank-functional-equation}, \\ \eqref{ABS-terminal-morphism-diagram}}}{=}
\sum_{p \in P} (\psi_q \circ \Psi)[\hat{0},p] \cdot
               (\zeta_Q \circ \Psi)[p,\hat{1}] \\
&= \left( \left( \psi_q \circ \Psi \right) \star \left( \zeta_Q \circ \Psi \right) \right) [P]
\overset{\eqref{pf.Ehrenborg's-theorem.compose-star}}{=}(\psi_q \star \zeta_Q)(\Psi[P])
=(\psi_q \otimes \zeta_Q)\left( \Delta \Psi[P]\right) \\
&=\left[ \Psi[P](\xx,\yy) \right]_{\substack{x_1=q,x_2=x_3=\cdots=0\\y_1=1,y_2=y_3=\cdots=0}}
=\left[ \Psi[P](\xx) \right]_{\substack{x_1=q,x_2=1,\\x_3=x_4=\cdots=0}} .
\end{align*}

%Assertion (ii) follows from observing the evaluation for $\alpha$ in $\Comp_n$
%\[
%\left[ M_\alpha\right]_{\substack{x_1=q,x_2=1,\\x_3=x_4= \cdots =0}}
%=\begin{cases}
%1+q^n & \text{ if }\alpha=(n),\\
%q^k & \text{ if } \alpha=(k,n-k), 0 < k < n\\
%0 & \text{ if }\ell(\alpha) \geq 3.
%\end{cases}
%\]
%This implies that
%\[
%\left[ \Psi[P] \right]_{\substack{x_1=q,x_2=1,\\x_3=x_4= \cdots =0}}
%= 1+q^{r(P)} + \sum_{r=1}^{r(P)-1} f_{(k,n-k)}(P) q^k
%=RGF(P,q).
%\]

Similarly, for assertion (iv) first note  that
\begin{equation}
\left( \left( \psi_q \circ S \right) \star \zeta_Q \right) \circ \Psi
= \left( \psi_q \circ S \circ \Psi \right) \star \left( \zeta_Q \circ \Psi \right) ,
\label{pf.Ehrenborg's-theorem.compose-star2}
\end{equation}
(this is proven similarly to
\eqref{pf.Ehrenborg's-theorem.compose-star}, but now
using the map $\psi_q \circ S$ instead of $\psi_q$).
Now,
Proposition~\ref{incidence-coalgebras-are-Hopf} and
Corollary~\ref{P.Hall-formula} let one calculate that
\begin{align*}
(\psi_q \circ \Psi \circ S)[P]
&= \sum_k (-1)^k \sum_{\hat{0}=x_0 < \cdots < x_k=\hat{1}}
            (\psi_q \circ \Psi)([x_0,x_1]) \cdots (\psi_q \circ \Psi)([x_{k-1},x_k]) \\
&
\overset{\eqref{q-rank-functional-equation}}{=} \sum_k (-1)^k \sum_{\hat{0}=x_0 < \cdots < x_k=\hat{1}} q^{r(P)}
= \mu(\hat{0},\hat{1}) q^{r(P)}.
\end{align*}
This is used in the penultimate equality here:
\begin{align*}
((\psi_q \circ S) \star \zeta_Q) \circ \Psi[P]
%&= m \circ ((\psi_q \circ S) \otimes \zeta_Q) \circ \Delta \circ \Psi[P] \\
%&= m \circ ((\psi_q \circ S) \otimes \zeta_Q) \circ (\Psi \otimes \Psi) \circ\Delta [P] \\
&\overset{\eqref{pf.Ehrenborg's-theorem.compose-star2}}{=}
  ((\psi_q \circ S \circ \Psi) \star (\zeta_Q \circ \Psi)) [P]
= ((\psi_q \circ \Psi \circ S) \star \zeta) [P] \\
&= \sum_{p \in P}  (\psi_q \circ \Psi \circ S)[\hat{0},p] \cdot \zeta[p,\hat{1}]
=\sum_{p \in P} \mu[\hat{0},p] q^{r(p)} = \chi(P,q).
\end{align*}
\end{proof}

\subsection{Example: Stanley's chromatic symmetric function of a graph}

We introduce the \emph{chromatic Hopf algebra of graphs} and an associated
character $\zeta$ so that the map $\Psi$ from
Theorem~\ref{Qsym-as-terminal-object-theorem}
sends a graph $G$ to Stanley's \emph{chromatic symmetric function}
of $G$.  Then principal specialization $\ps^1$ sends this
to the \emph{chromatic polynomial} of the graph.

\subsubsection{The chromatic Hopf algebra of graphs}

\begin{definition}
The \dfn{chromatic Hopf algebra} (see Schmitt \cite[\S 3.2]{Schmitt-structures})
\dfn{$\GGG$} is a free $\kk$-module whose
$\kk$-basis elements $[G]$ are indexed by isomorphism classes of
(finite) simple graphs $G=(V,E)$.
Define for $G_1=(V_1,E_1), G_2=(V_2,E_2)$ the multiplication
\[
[G_1] \cdot [G_2]:=[G_1 \sqcup G_2]
\]
where $[G_1 \sqcup G_2]$ denote the isomorphism class of the
disjoint union, on vertex set $V=V_1 \sqcup V_2$ which
is a disjoint union of copies of their vertex sets $V_1,V_2$, with
edge set $E=E_1 \sqcup E_2$.  For example,
\[
\left[
\xymatrix@C=4pt@R=4pt{
\bullet& &\bullet\\
       &\bullet\ar@{-}[ul]\ar@{-}[ur]&
}
\right] \cdot
\left[
\xymatrix@C=4pt@R=4pt{
\bullet\\
\bullet\ar@{-}[u]
}
\right]
=
\left[
\xymatrix@C=4pt@R=4pt{
\bullet& &\bullet &\bullet\\
       &\bullet\ar@{-}[ul]\ar@{-}[ur]& &\bullet\ar@{-}[u]
}
\right]
\]
  Thus the class $[\varnothing]$
of the empty graph $\varnothing$ having $V=\varnothing, E=\varnothing$
is a unit element.

Given a graph $G = \left(V, E\right)$ and a subset $V' \subset V$,
the \dfn{subgraph induced on vertex set $V'$} is defined as
the graph $G|_{V'}:=(V',E')$ with
edge set $E'=\{e \in E: e= \left\{v_1,v_2\right\} \subset V'\}$.
This lets one define a comultiplication
$\Delta : \GGG \to \GGG \otimes \GGG$ by setting
\[
\Delta[G]:=\sum_{(V_1,V_2): V_1 \sqcup V_2=V} [G|_{V_1}] \otimes [G|_{V_2}].
\]
Define a counit $\epsilon : \GGG \to \kk$ by
\[
\epsilon [G]:=
\begin{cases}
1, & \text{ if }G=\varnothing \text{;} \\
0, & \text{ otherwise.}
\end{cases}
\]
\end{definition}

\begin{proposition}
\label{prop.GGG.hopf}
The above maps endow $\GGG$ with the structure of a connected graded
finite type Hopf algebra over $\kk$, which is both commutative and cocommutative.
\end{proposition}

\begin{example}
Here are some examples of these structure maps:
\begin{align*}
\left[
\xymatrix@C=4pt@R=4pt{
\bullet& &\bullet\\
       &\bullet\ar@{-}[ul]\ar@{-}[ur]&
}
\right] \cdot
\left[
\xymatrix@C=4pt@R=4pt{
\bullet\\
\bullet\ar@{-}[u]
}
\right]
&=
\left[
\xymatrix@C=4pt@R=4pt{
\bullet& &\bullet &\bullet\\
       &\bullet\ar@{-}[ul]\ar@{-}[ur]& &\bullet\ar@{-}[u]
}
\right] ;\\
%%%%%%%%%%%%%%
\Delta \left[
\xymatrix@C=4pt@R=4pt{
\bullet& &\bullet\\
       &\bullet\ar@{-}[ul]\ar@{-}[ur]&
}
\right]
&= \one
\otimes
\left[
\xymatrix@C=4pt@R=4pt{
\bullet& &\bullet\\
       &\bullet\ar@{-}[ul]\ar@{-}[ur]&
}
\right]
+
2
\left[
\xymatrix@C=4pt@R=4pt{\bullet}
\right]
\otimes
\left[
\xymatrix@C=4pt@R=4pt{
\bullet\\
\bullet\ar@{-}[u]
}
\right]
+
2
\left[
\xymatrix@C=4pt@R=4pt{
\bullet\\
\bullet\ar@{-}[u]
}
\right]
\otimes
\left[
\xymatrix@C=4pt@R=4pt{\bullet}
\right]
+
\left[
\xymatrix@C=4pt@R=4pt{\bullet & \bullet}
\right]
\otimes
\left[
\xymatrix@C=4pt@R=4pt{\bullet}
\right] \\
& \qquad
+
\left[
\xymatrix@C=4pt@R=4pt{\bullet}
\right]
\otimes
\left[
\xymatrix@C=4pt@R=4pt{\bullet & \bullet}
\right]
+
\left[
\xymatrix@C=4pt@R=4pt{
\bullet& &\bullet\\
       &\bullet\ar@{-}[ul]\ar@{-}[ur]&
}
\right]
\otimes \one
\end{align*}
\end{example}

\begin{proof}[Proof of Proposition~\ref{prop.GGG.hopf}.]
The associativity of the multiplication and comultiplication
should be clear as
\begin{align*}
m^{(2)}([G_1] \otimes [G_2] \otimes [G_3])
&=[G_1 \sqcup G_2 \sqcup G_3] ,\\
\Delta^{(2)}[G]
&=\sum\limits_{\substack{(V_1,V_2,V_3):\\V=V_1 \sqcup V_2 \sqcup V_3}}
       [G|_{V_1}] \otimes [G|_{V_2}] \otimes [G|_{V_3}].
\end{align*}
Checking the unit and counit conditions are straightforward.
Commutativity of the pentagonal bialgebra diagram
in \eqref{bialgebra-diagrams} comes down to check that,
given graphs $G_1,G_2$ on disjoint vertex
sets $V_1,V_2$ , when one applies to $[G_1] \otimes [G_2]$
either the composite $\Delta \circ m$ or the composite
$(m \otimes m) \circ (\id \otimes T \otimes \id) \circ (\Delta \otimes \Delta)$,
the result is the same:
\[
\sum\limits_{\substack{(V_{11},V_{12},V_{21},V_{22}):\\
        V_1=V_{11} \sqcup V_{12}\\
        V_2=V_{21} \sqcup V_{22}}}
[G_1|_{V_{11}} \sqcup G_2|_{V_{21}}] \otimes
[G_1|_{V_{12}} \sqcup G_2|_{V_{22}}].
\]
Letting $\GGG_n$ be the $\kk$-span of $[G]$ having $n$ vertices makes $\GGG$
a bialgebra which is graded and connected, and hence also a
Hopf algebra by Proposition~\ref{graded-connected-bialgebras-have-antipodes}.
Cocommutativity should be clear, and commutativity follows from the
graph isomorphism $G_1 \sqcup G_2 \cong G_2 \sqcup G_1$.
Finally, $\GGG$ is of finite type since there are only finitely
many isomorphism classes of simple graphs on $n$ vertices for
every given $n$.
\end{proof}

\begin{remark}
\label{Humpert-Martin-antipode-remark}
Humpert and Martin \cite[Theorem 3.1]{HumpertMartin}
gave the following expansion for the antipode in the
chromatic Hopf algebra, containing fewer terms than Takeuchi's general
formula \eqref{Takeuchi-formula}:  given a graph $G=(V,E)$, one has
\begin{equation}
\label{eq.Humpert-Martin-antipode-remark.S(G)}
S[G]=\sum_{F} (-1)^{|V|-\rank(F)} \acyc(G/F) [G_{V,F}].
\end{equation}
Here $F$ runs over all subsets of edges that form
\emph{flats}\index{flat in a graphic matroid}
in the graphic matroid for $G$, meaning that if $e=\{v,v'\}$ is an
edge in $E$ for which one has a path of edges in $F$ connecting $v$ to $v'$,
then $e$ also lies in $F$.  Here $G/F$ denotes the quotient graph in which all of the
edges of $F$ have been \emph{contracted}, while $\acyc(G/F)$ denotes its number
of \emph{acyclic orientations}, and $G_{V,F}:=(V,F)$ as a simple
graph.\footnote{The notation $\rank(F)$ denotes the
\emph{rank}\index{rank in a graphic matroid}
of $F$ in the graphic matroid of $G$. We can define it without
reference to matroid theory as the maximum cardinality of a
subset $F'$ of $F$ such that the graph $G_{V,F'}$ is acyclic.
Equivalently, $\rank(F)$ is $|V|-c(F)$, where $c(F)$ denotes the
number of connected components of the graph $G_{V,F}$. Thus,
the equality \eqref{eq.Humpert-Martin-antipode-remark.S(G)}
can be rewritten as
$S[G]=\sum_{F} (-1)^{c(F)} \acyc(G/F) [G_{V,F}]$.
In this form, this equality is also proven in
\cite[Thm. 7.1]{BenedettiSagan}.}
\end{remark}

% [DG][v44] Added footnote about the meaning of $\rank(F)$, and also
% referenced \cite{BenedettiSagan} for another proof.

\begin{remark}
In \cite{BenedettiHallamMachacek}, Benedetti, Hallam and Machacek
define a Hopf algebra of simplicial complexes, which contains
$\GGG$ as a Hopf subalgebra (and also has $\GGG$ as a quotient Hopf
algebra). They compute a formula for its antipode similar to
(and generalizing) \eqref{eq.Humpert-Martin-antipode-remark.S(G)}.
\end{remark}

% [DG][v50] Added above remark.

\begin{remark}
The chromatic Hopf algebra $\GGG$ is used in \cite{Lando} and
\cite[\S 14.4]{ChmutovDuzhinMostovoy} to study
\emph{Vassiliev invariants of knots}. In fact, a certain quotient
of $\GGG$ (named $\mathcal{F}$ in \cite{Lando} and $\mathcal{L}$ in
\cite[\S 14.4]{ChmutovDuzhinMostovoy}) is shown to naturally
host invariants of \emph{chord diagrams} and therefore Vassiliev
invariants of knots.
\end{remark}

% [DG][v41] Above remark added.

\begin{remark}
\label{rmk.GGG.free}
The $\kk$-algebra $\GGG$ is isomorphic to a polynomial algebra (in
infinitely many indeterminates) over $\kk$. Indeed, every finite
graph can be uniquely written as a disjoint union of finitely
many connected finite graphs (up to order). Therefore, the basis
elements $\left[G\right]$ of $\GGG$ corresponding to connected
finite graphs $G$ are algebraically independent in $\GGG$ and
generate the whole $\kk$-algebra $\GGG$ (indeed, the disjoint
unions of connected finite graphs are precisely the monomials in
these elements). Thus, $\GGG$ is isomorphic to a polynomial
$\kk$-algebra with countably many generators (one for each
isomorphism class of connected finite graphs). As a consequence,
for example, we see that $\GGG$ is an integral domain if $\kk$ is
an integral domain.
\end{remark}

% [DG][v51] Added above remark.

\subsubsection{A ``ribbon basis'' for $\GGG$ and self-duality}

In this subsection, we shall explore a second basis of $\GGG$
and a bilinear form on $\GGG$. This material will not be used
in the rest of these notes (except in
Exercise~\ref{exe.GGG.chromatic-sharp}), but it is of some
interest and provides an example of how a commutative
cocommutative Hopf algebra can be studied.

First, let us define a second basis of $\GGG$, which is obtained
by M\"obius inversion (in an appropriate sense) from the standard
basis
$\left( \left[ G \right] \right) _{\left[ G \right]
\text{ is an isomorphism class of finite graphs}}$:

\begin{definition}
\label{def.GGG.moebius-basis}
For every finite graph $G = \left(V, E\right)$, set
\[
\left[G\right]^\sharp = \sum\limits_{\substack{H = \left(V, E^\prime\right);\\
                                        E^\prime \supset E^c}}
                        \left(-1\right)^{\left|E^\prime \setminus E^c\right|}
                        \left[H\right]
                      \in \GGG,
\]
where $E^c$ denotes the complement of the subset $E$ in the set of
all two-element subsets of $V$. Clearly, $\left[G\right]^\sharp$
depends only on the isomorphism class $\left[G\right]$ of $G$, not on
$G$ itself.
\end{definition}

\begin{proposition} \phantomsection
\label{prop.GGG.moebius-basis}
\begin{itemize}
\item[(a)] Every finite graph $G = \left(V, E\right)$
satisfies
\[
\left[G\right] = \sum\limits_{\substack{H = \left(V, E^\prime\right);\\
                                        E^\prime \cap E = \varnothing}}
                        \left[H\right]^\sharp .
\]

\item[(b)] The elements $\left[G\right]^\sharp$, where
$\left[G\right]$ ranges over all isomorphism classes of finite graphs,
form a basis of the $\kk$-module $\GGG$.

\item[(c)] For any graph $H = \left(V, E\right)$, we have
\begin{equation}
\label{eq.prop.GGG.moebius-basis.c}
\Delta\left[H\right]^\sharp
=\sum\limits_{\substack{\left(V_1,V_2\right);\\ V=V_1 \sqcup V_2 ;\\ H=H|_{V_1} \sqcup H|_{V_2}}}
\left[H|_{V_1}\right]^\sharp \otimes \left[H|_{V_2}\right]^\sharp .
\end{equation}

\item[(d)] For any two graphs $H_1 = \left(V_1, E_1\right)$ and
$H_2 = \left(V_2, E_2\right)$, we have
\begin{equation}
\label{eq.prop.GGG.moebius-basis.d}
\left[H_1\right]^\sharp \left[H_2\right]^\sharp
 = \sum\limits_{\substack{H=\left(V_1 \sqcup V_2,E\right) ;\\H|_{V_1}=H_1;\\ H|_{V_2}=H_2}} \left[H\right]^\sharp.
\end{equation}
\end{itemize}
\end{proposition}

For example,
\begin{align*}
\left[
\xymatrix@C=4pt@R=4pt{
\bullet& &\bullet\\
       &\bullet\ar@{-}[ul]\ar@{-}[ur]&
}
\right]^\sharp
&=
\left[
\xymatrix@C=4pt@R=4pt{
\bullet\ar@{-}[rr]& &\bullet\\
       &\bullet\ar@{-}[ul]\ar@{-}[ur]&
}
\right]
-
\left[
\xymatrix@C=4pt@R=4pt{
\bullet\ar@{-}[rr]& &\bullet\\
       &\bullet\ar@{-}[ur]&
}
\right]
- \left[
\xymatrix@C=4pt@R=4pt{
\bullet\ar@{-}[rr]& &\bullet\\
       &\bullet\ar@{-}[ul]&
}
\right]
+ \left[
\xymatrix@C=4pt@R=4pt{
\bullet\ar@{-}[rr]& &\bullet\\
       &\bullet&
}
\right]
\\
&=
\left[
\xymatrix@C=4pt@R=4pt{
\bullet\ar@{-}[rr]& &\bullet\\
       &\bullet\ar@{-}[ul]\ar@{-}[ur]&
}
\right]
-
2 \left[
\xymatrix@C=4pt@R=4pt{
\bullet\ar@{-}[rr]& &\bullet\\
       &\bullet\ar@{-}[ur]&
}
\right]
+ \left[
\xymatrix@C=4pt@R=4pt{
\bullet\ar@{-}[rr]& &\bullet\\
       &\bullet&
}
\right] .
\end{align*}

Proving Proposition~\ref{prop.GGG.moebius-basis} is part of
Exercise~\ref{exe.GGG.moebius-basis} further below.

The equalities that express the elements $\left[G\right]^\sharp$
in terms of the elements $\left[H\right]$ (as in
Definition~\ref{def.GGG.moebius-basis}), and vice versa
(Proposition~\ref{prop.GGG.moebius-basis}(a)), are reminiscent
of the relations \eqref{ribbon-as-sum-of-Hs} and
\eqref{H-as-sum-of-ribbons} between the bases
$\left(R_\alpha\right)$ and $\left(H_\alpha\right)$ of $\Nsym$.
In this sense, we can call the basis of $\GGG$ formed by the
$\left[G\right]^\sharp$ a ``ribbon basis'' of $\GGG$.

We now define a $\kk$-bilinear form on $\GGG$:

\begin{definition}
\label{def.GGG.bilform}
For any two graphs $G$ and $H$, let
$\operatorname*{Iso} \left( G, H \right)$
denote the set of all isomorphisms from $G$ to
$H$\ \ \ \ \footnote{We recall that if $G = \left( V, E \right)$ and
$H = \left( W, F \right)$ are two graphs, then an
\emph{isomorphism}\index{isomorphism of graphs}
from $G$ to $H$ means a bijection $\varphi : V \rightarrow W$ such
that $\varphi_{\ast} \left( E \right) = F$. Here, $\varphi_{\ast}$
denotes the map from the powerset of $V$ to the powerset of $W$ which
sends every $T \subset V$ to $\varphi \left( T \right) \subset W$.}.
Let us now define a $\kk$-bilinear form $\left( \cdot, \cdot \right)
: \GGG \times \GGG \rightarrow \kk$ on $\GGG$ by
setting
\[
\left( \left[ G \right] ^{\sharp}, \left[ H \right] \right)
= \left\vert \operatorname*{Iso} \left( G, H \right) \right\vert .
\]
\footnote{This is well-defined, because:
\par
\begin{itemize}
\item the number $\left\vert \operatorname*{Iso} \left( G, H \right)
\right\vert$ depends only on the isomorphism classes $\left[ G \right]$
and $\left[ H \right]$ of $G$ and $H$, but not on $G$ and $H$
themselves;
\par
\item the elements $\left[ G \right] ^{\sharp}$, where
$\left[ G \right]$ ranges over all isomorphism classes of finite
graphs, form a basis of the $\kk$-module $\GGG$ (because of
Proposition~\ref{prop.GGG.moebius-basis}(b));
\par
\item the elements $\left[ G \right]$, where
$\left[ G \right]$ ranges over all isomorphism classes of finite
graphs, form a basis of the $\kk$-module $\GGG$.
\end{itemize}
}
\end{definition}

\begin{proposition}
\label{prop.GGG.bilform.sym}
The form $\left( \cdot, \cdot \right)
: \GGG \times \GGG \rightarrow \kk$ is symmetric.
\end{proposition}

Again, we refer to Exercise~\ref{exe.GGG.moebius-basis} for a proof
of Proposition~\ref{prop.GGG.bilform.sym}.

The basis of $\GGG$ constructed in
Proposition~\ref{prop.GGG.moebius-basis}(b) and the bilinear form
$\left(\cdot, \cdot\right)$ defined in
Definition~\ref{def.GGG.bilform} can be used to construct a Hopf
algebra homomorphism from $\GGG$ to its graded dual $\GGG^o$:

\begin{definition}
\label{def.GGG.selfdual}
For any finite graph $G$, let $\operatorname*{aut} \left( G \right)$
denote the number
$\left\vert \operatorname*{Iso} \left( G, G \right) \right\vert$.
Notice that this is a positive integer, since the set
$\operatorname*{Iso} \left( G, G \right)$ is nonempty (it contains
$\id_G$).

Now, recall that the Hopf algebra $\GGG$ is a connected graded Hopf
algebra of finite type. The $n$-th homogeneous component is spanned by
the $\left[ G \right]$ where $G$ ranges over the graphs with $n$
vertices. Since $\GGG$ is of finite type, its graded dual $\GGG^{o}$
is defined.
Let $\left( \left[ G \right] ^{\ast} \right) _{\left[ G \right]
\text{ is an isomorphism class of finite graphs}}$
be the basis of $\GGG^{o}$ dual to the basis
$\left( \left[ G \right] \right) _{\left[ G \right]
\text{ is an isomorphism class of finite graphs}}$ of $\GGG$. Define
a $\kk$-linear map $\psi : \GGG \rightarrow \GGG^{o}$ by
\[
\psi \left( \left[ G \right] ^{\sharp} \right)
= \operatorname*{aut} \left( G \right) \cdot \left[ G \right] ^{\ast}
\ \ \ \ \ \ \ \ \ \ \text{for every finite graph } G \text{.}
\]
\footnote{This is well-defined, since
$\left( \left[ G \right] ^{\sharp} \right) _{\left[ G \right]
\text{ is an isomorphism class of finite graphs}}$
is a basis of the $\kk$-module $\GGG$ (because of
Proposition~\ref{prop.GGG.moebius-basis}(b)).}
\end{definition}

\begin{proposition}
\label{prop.GGG.selfdual}
Consider the map $\psi : \GGG \to \GGG^{o}$ defined in
Definition~\ref{def.GGG.selfdual}.

\begin{itemize}
\item[(a)] This map $\psi$ satisfies
$\left( \psi \left( a \right) \right) \left( b \right)
= \left( a, b \right)$ for all $a \in \GGG$ and $b \in \GGG$.

\item[(b)] The map $\psi : \GGG \rightarrow \GGG^{o}$ is a Hopf
algebra homomorphism.

\item[(c)] If $\QQ$ is a subring of $\kk$, then the map $\psi$
is a Hopf algebra isomorphism $\GGG \rightarrow \GGG^{o}$.
\end{itemize}
\end{proposition}

\begin{exercise}
\label{exe.GGG.moebius-basis}
Prove Proposition~\ref{prop.GGG.moebius-basis},
Proposition~\ref{prop.GGG.bilform.sym} and
Proposition~\ref{prop.GGG.selfdual}.
\end{exercise}

\begin{remark}
Proposition~\ref{prop.GGG.selfdual}(c) shows that the Hopf algebra
$\GGG$ is self-dual when $\QQ$ is a subring of $\kk$. On
the other hand, if $\kk$ is a field of positive characteristic, then
$\GGG$ is never self-dual. Here is a quick way to see this: The
elements $\left[ G \right] ^{\ast}$ of $\GGG^{o}$ defined in
Definition~\ref{def.GGG.selfdual} have the property that
\[
\left( \left[ \circ \right] ^{\ast} \right) ^{n}
= n! \cdot \sum\limits_{\substack{\left[ G \right] \text{ is an isomorphism}
\\\text{class of finite graphs on}\\ n \text{ vertices}}}
\left[ G \right] ^{\ast}
\]
for every $n \in \NN$, where $\circ$ denotes the graph with one
vertex.\footnote{To see this, observe that the tensor
$\left[ \circ \right] ^{\otimes n}$ appears in the iterated coproduct
$\Delta^{\left( n-1 \right)} \left( \left[ G \right] \right)$
exactly $n!$ times whenever $G$ is a graph on $n$ vertices.} Thus,
if $p$ is a prime and $\kk$ is a field of characteristic $p$, then
$\left( \left[ \circ \right] ^{\ast} \right) ^{p} = 0$. Hence, the
$\kk$-algebra $\GGG^{o}$ has nilpotents in this situation. However,
the $\kk$-algebra $\GGG$ does not (indeed, Remark~\ref{rmk.GGG.free}
shows that it is an integral domain whenever $\kk$ is an integral
domain). Thus, when $\kk$ is a field of characteristic $p$, then
$\GGG$ and $\GGG^{o}$ are not isomorphic as $\kk$-algebras (let alone
as Hopf algebras).
\end{remark}

\begin{commentedout}
\textbf{Update:} The following claims (specifically, the equations
\eqref{dual-chromatic-Hopf-structure-maps-1} and
\eqref{dual-chromatic-Hopf-structure-maps-2} and
Proposition \ref{prop.GGG.selfdual.OLD-AND-FALSE}) are FALSE.

It turns out that the chromatic Hopf algebra $\GGG$ is self-dual.
In the dual Hopf algebra $\GGG^o$,
let $\{[G]^*\}$ denote the dual basis elements, so that $([H]^*,[G])=\delta_{[H],[G]}$.
To describe the structure maps in $\GGG^o$ explicitly, for any
graph $H = \left(V, E\right)$ one has
\begin{equation}
\label{dual-chromatic-Hopf-structure-maps-1}
\Delta[H]^*
=\sum\limits_{\substack{(V_1,V_2):\\ V=V_1 \sqcup V_2\\H=H|_{V_1} \sqcup H|_{V_2}}}
[H|_{V_1}]^* \otimes [H|_{V_2}]^*,
\end{equation}
and for any two graphs $H_1 = \left(V_1, E_1\right)$ and
$H_2 = \left(V_2, E_2\right)$ one has
\begin{equation}
\label{dual-chromatic-Hopf-structure-maps-2}
[H_1]^* [H_2]^*
 = \sum\limits_{\substack{H=(V_1 \sqcup V_2,E)\\H|_{V_1}=H_1\\H|_{V_2}=H_2}} [H]^*.
\end{equation}

% [DG][v14] I have split the two aligned equations, because $V_1$
% means a different thing in the first one than it means in the
% second.

\begin{proposition}
\label{prop.GGG.selfdual.OLD-AND-FALSE}
One has a Hopf isomorphism
$\GGG \overset{\varphi}{\longrightarrow} \GGG^o$
defined by
\[
[G] \longmapsto \sum\limits_{\substack{H=(V,E'):\\ E' \cap E = \varnothing}} [H]^*
\qquad \text{ for all graphs } G = \left(V,E\right).
\]
\end{proposition}
For example, this isomorphism maps
\[
\varphi \left[
\xymatrix@C=4pt@R=4pt{
\bullet& &\bullet\\
       &\bullet\ar@{-}[ul]\ar@{-}[ur]&
}
\right]
=
\left[
\xymatrix@C=4pt@R=4pt{
\bullet\ar@{-}[rr]& &\bullet\\
       &\bullet&
}
\right]^*
+
\left[
\xymatrix@C=4pt@R=4pt{
\bullet& &\bullet\\
       &\bullet&
}
\right]^*.
\]

\begin{proof}[Proof of Proposition~\ref{prop.GGG.selfdual.OLD-AND-FALSE}.]
First note that $\varphi$ is a $\kk$-module
isomorphism via triangularity:  one has $H=(V,E')$
with $E' \cap E =\varnothing$ if and only
if $H$ is an edge subgraph of the
\emph{complementary graph} $\overline{G}$
to $G$ on the same vertex set $V$.

One can then check that for graphs $G_1,G_2$ on vertex sets $V_1, V_2$,
the fact that $\varphi([G_1][G_2])=
\varphi[G_1]\varphi[G_2]$
amounts, using
\eqref{dual-chromatic-Hopf-structure-maps-2},
to both being a sum of $[H]^*$
over graphs $H$ on $V_1 \sqcup V_2$ that share no
edges with $G_1 \sqcup G_2$.

Similarly, one can check that the fact
that $\Delta \varphi[G]=(\varphi \otimes \varphi)(\Delta[G])$
amounts, using
\eqref{dual-chromatic-Hopf-structure-maps-1},
to both being a sum of $[H_1]^* \otimes [H_2]^*$
over triples $(H,H_1,H_2)$ where
$H$ is a graph on the same vertex
set $V$ as $G$ but sharing no edges with $G$, and
with $H=H_1 \sqcup H_2$.
\end{proof}
\end{commentedout}

% [DG][v50] Replaced false claim of selfduality by a proposition giving
% a new basis.

% [DG][v51] Extended the proposition and added a remark. Together,
% they now show that $\GGG$ is selfdual when $\kk$ is a $\QQ$-algebra,
% and not selfdual when $\kk$ is a field of characteristic $p$.

% [DG][v53] Restructured the statements and proofs around the
% ribbon basis and the bilinear form on $\GGG$.

\subsubsection{Stanley's chromatic symmetric function of a graph}

\begin{definition}
\label{def.GGG.chromatic}
\dfn{Stanley's chromatic symmetric function}\index{chromatic symmetric function}
$\Psi[G]$ for a simple graph $G=(V,E)$
is the image of $[G]$ under the map
$\GGG \overset{\Psi}{\longrightarrow} \Qsym$
induced via Theorem~\ref{Qsym-as-terminal-object-theorem}
from the \dfn{edge-free character}
$\GGG \overset{\zeta}{\longrightarrow} \kk$ defined by
\begin{equation}
\label{edge-free-character-defn}
\zeta[G]=\begin{cases}
1, & \text{ if }G\text{ has no edges, that is, }G\text{ is an independent/stable set of vertices};\\
0, & \text{ otherwise.}
\end{cases}
\end{equation}
Note that, because $\GGG$ is cocommutative, $\Psi[G]$ is symmetric and not just quasisymmetric;
see Remark~\ref{Sym-as-terminal-object-remark}.
\end{definition}

Recall that for a graph $G=(V,E)$,
a (vertex-)coloring $f: V \rightarrow \{1,2,\ldots\}$
is called \emph{proper}\index{proper coloring}
if no edge $e=\{v,v'\}$ in $E$ has $f(v)=f(v')$.

\begin{proposition}
\label{prop.GGG.chromatic}
For a graph $G=(V,E)$, the symmetric function $\Psi[G]$ has the expansion
\footnote{In fact, Stanley \emph{defined} $\Psi[G]$  in
\cite[Defn. 2.1]{Stanley-chromatic-qsym-fn}
via this expansion.}
\[
\Psi[G]=\sum\limits_{\substack{\text{proper colorings}\\f:V \rightarrow \{1,2,\ldots\}}}
         \xx_f
\]
where $\xx_f:=\prod_{v \in V} x_{f(v)}$.
In particular, its specialization from Proposition~\ref{principal-specialization-at-1-defn}
gives the chromatic polynomial of $G$:
\[
\ps^1 \Psi[G](m) = \chi_G(m)
= \left|\{\text{proper colorings }f:V \rightarrow \{1,2,\ldots,m\}\}\right| .
\]
\end{proposition}
\begin{proof}
The iterated coproduct
$\GGG \overset{\Delta^{(\ell-1)}}{\longrightarrow} \GGG^{\otimes \ell}$
sends
\[
[G] \longmapsto \sum\limits_{\substack{(V_1,\ldots,V_{\ell}):\\ V=V_1 \sqcup \cdots \sqcup V_\ell}}
[G|_{V_1}] \otimes \cdots \otimes [G|_{V_\ell}]
\]
and the map $\zeta^{\otimes \ell}$ sends each addend on the right to
$1$ or $0$, depending upon whether each $V_i \subset V$ is a stable set or not,
that is, whether the assignment of color $i$ to the vertices in $V_i$ gives
a proper coloring of $G$.  Thus formula \eqref{ABS-terminal-morphism-formula}
shows that the coefficient $\zeta_\alpha$ of $x_1^{\alpha_1} \cdots x_\ell^{\alpha_\ell}$ in $\Psi[G]$ counts the
proper colorings $f$ in which $|f^{-1}(i)|=\alpha_i$ for each $i$.
\end{proof}

% [DG][v25] Replaced "the element on the right" by "each
% addend on the right".

\begin{example}
For the complete graph $K_n$ on $n$ vertices,
one has
\begin{align*}
\Psi[K_n]&=n! e_n , \qquad \text{thus} \\
\ps^1 (\Psi[K_n])(m)
&=n! e_n(\underbrace{1,1,\ldots,1}_{m\text{ ones}})
=n! \binom{m}{n}\\
& =m(m-1) \cdots (m-(n-1))=\chi_{K_n}(m).
\end{align*}
In particular, the single vertex graph $K_1$ has
$\Psi[K_1]=e_1$, and since the Hopf morphism $\Psi$ is in particular
an algebra morphism, a graph $K_1^{\sqcup n}$ having $n$ isolated
vertices and no edges will have $\Psi[K_1^{\sqcup n}]=e_1^n$.

As a slightly more interesting example, the graph $P_3$ which is
a path having three vertices and two edges will have
\[
\Psi[P_3]=m_{(2,1)}+6m_{(1,1,1)}=e_2 e_1 + 3 e_3 .
\]
\end{example}

One might wonder, based on the previous examples, when $\Psi[G]$ is
\dfn{$e$-positive}, that is, when does its unique expansion in the $\{e_\lambda\}$
basis for $\Lambda$ have nonnegative coefficients?  This is an even stronger
assertion than \emph{$s$-positivity}\index{$s$-positive},
that is, having nonnegative coefficients for the expansion in terms
of Schur functions $\{s_\lambda\}$, since each $e_\lambda$ is $s$-positive.
This weaker property fails, starting with the
\dfn{claw graph $K_{3,1}$}, which has
\[
\Psi[K_{3,1}]=s_{(3,1)}-s_{(2,2)}+5s_{(2,1,1)}+8s_{(1,1,1,1)}.
\]
On the other hand, a result of Gasharov \cite[Theorem 2]{Gasharov} shows that one at
least has $s$-positivity for $\Psi[\inc(P)]$ where
$\inc(P)$ is the \emph{incomparability graph} of a poset
which is $(\mathbf{3}+\mathbf{1})$-free;  we refer the
reader to Stanley \cite[\S 5]{Stanley-chromatic-qsym-fn} for a discussion
of the following conjecture, which remains open\footnote{A recent
refinement for incomparability graphs of posets which are both $(\mathbf{3}+\mathbf{1})$- and
$(\mathbf{2}+\mathbf{2})$-free, also known as \emph{unit interval orders}
is discussed by Shareshian and Wachs \cite{ShareshianWachs}.}:

\begin{conjecture}
For any $(\mathbf{3}+\mathbf{1})$-free poset $P$, the incomparability graph
$\inc(P)$ has $\Psi[\inc(P)]$ an $e$-positive symmetric function.
\end{conjecture}

Here is another question about $\Psi[G]$:  how well does it distinguish nonisomorphic
graphs? Stanley gave this example of two graphs $G_1, G_2$ having $\Psi[G_1]=\Psi[G_2]$:
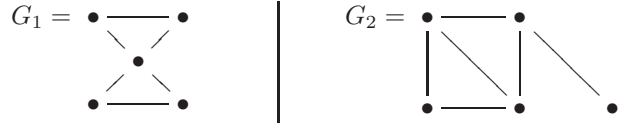
\[
\begin{array}{c|c}
G_1=
\xymatrix@C=6pt@R=6pt{
\bullet\ar@{-}[rr] &  &\bullet\\
        &\bullet\ar@{-}[ul]\ar@{-}[ur]\ar@{-}[dr]\ar@{-}[dl]& \\
\bullet\ar@{-}[rr] &  &\bullet
}
\qquad \phantom{a} & \qquad
G_2=
\xymatrix{
\bullet\ar@{-}[r] \ar@{-}[dr]\ar@{-}[d]&\bullet\ar@{-}[d]\ar@{-}[dr]  &\\
\bullet\ar@{-}[r] &\bullet &\bullet
}
\end{array}
\]
\noindent
At least $\Psi[G]$ appears to do better at distinguishing \emph{trees}, much
better than its specialization, the chromatic polynomial $\chi_G(m)$, which takes
the same value $m(m-1)^{n-1}$ on all trees with $n$ vertices.
\begin{question}
Does the chromatic symmetric function (for $\kk = \ZZ$) distinguish trees?
\end{question}
\noindent
It has been checked that the answer is affirmative for trees on $23$ vertices or less.
There are also interesting partial results on this question
by Martin, Morin and Wagner \cite{MartinMorinWagner}.

\vskip.1in
We close this section with a few other properties of $\Psi[G]$ proven by Stanley
which follow easily from the theory we have developed.
For example, his work makes
no explicit mention of the chromatic Hopf algebra $\GGG$, and the fact that $\Psi$ is
a Hopf morphism (although he certainly notes the trivial algebra morphism property
$\Psi[G_1 \sqcup G_2]=\Psi[G_1] \Psi[G_2]$).  One property he proves is implicitly related to
$\Psi$ as a coalgebra morphism:  he considers (in the case when $\QQ$ is a subring
of $\kk$) the effect on $\Psi$ of the operator
$
\frac{\partial}{\partial p_1}:
\Lambda_\QQ \longrightarrow \Lambda_\QQ
$
which acts by first expressing a symmetric function
$f \in \Lambda_\QQ$ as a polynomial in the power sums $\{ p_n \}$, and then applies
the partial derivative operator $\frac{\partial}{\partial p_1}$ of the
polynomial ring $\QQ \left[p_1, p_2, p_3, \ldots \right]$.
It is not hard to see that  $\frac{\partial}{\partial p_1}$
is the same as the skewing operator $s_{(1)}^\perp=p_1^\perp$:  both act as  derivations
on $\Lambda_\QQ=\QQ[p_1,p_2,\ldots]$
(since $p_1 \in \Lambda_\QQ$ is primitive),
and agree in their effect on each $p_n$, in that
both send $p_1 \mapsto  1$, and both annihilate $p_2,p_3,\ldots$.

\begin{proposition}(Stanley \cite[Cor. 2.12(a)]{Stanley-chromatic-qsym-fn})
For any graph $G=(V,E)$, one has
\[
\frac{\partial}{\partial p_1} \Psi[G] = \sum_{v \in V} \Psi[G|_{V\setminus v}].
\]
\end{proposition}
\begin{proof}
Since $\Psi$ is a coalgebra homomorphism, we have
\[
\Delta \Psi[G] =
 (\Psi \otimes \Psi) \Delta[G]
=\sum\limits_{\substack{(V_1,V_2):\\V=V_1 \sqcup V_2}} \Psi[G|_{V_1}] \otimes \Psi[G|_{V_2}].
\]
Using this expansion (and the equality
$\frac{\partial}{\partial p_1} = s_{\left(1\right)}^\perp$), we now compute
\[
\frac{\partial}{\partial p_1} \Psi[G] = s_{(1)}^\perp \Psi[G] =
\sum\limits_{\substack{(V_1,V_2):\\V=V_1 \sqcup V_2}} (s_{(1)},\Psi[G|_{V_1}])  \cdot \Psi[G|_{V_2}]
=\sum_{v \in V} \Psi[G|_{V\setminus v}]
\]
(since degree considerations force $(s_{(1)}, \Psi[G|_{V_1}]) = 0$ unless $|V_1|=1$,
in which case $\Psi[G|_{V_1}]=s_{(1)}$).
\end{proof}

\begin{definition}
Given a graph $G=(V,E)$, an acyclic orientation $\Omega$ of the edges $E$
(that is, an orientation of each edge such that the resulting directed graph
has no cycles),
and a vertex-coloring $f: V \rightarrow \{1,2,\ldots\}$,
say that the pair $(\Omega,f)$ are \dfn{weakly compatible}
if whenever $\Omega$ orients an edge $\{v,v'\}$ in $E$ as $v \rightarrow v'$,
one has $f(v) \leq f(v')$.   Note that a
\emph{proper} vertex-coloring $f$ of a graph $G=(V,E)$
is weakly compatible with a unique acyclic orientation $\Omega$.
\end{definition}

% [DG] I added a parenthesized definition of acyclic orientations, and
% replaced "compatible" by "weakly compatible" here and in the proof (since
% you don't define just "compatible").

\begin{proposition}(Stanley \cite[Prop. 4.1, Thm. 4.2]{Stanley-chromatic-qsym-fn})
\label{Stanley-chromatic-reciprocity}
The involution $\omega$ of $\Lambda$ sends $\Psi[G]$ to
$
\omega \left( \Psi[G] \right) = \sum_{(\Omega,f)} \xx_f
$
in which the sum runs over weakly compatible
pairs $(\Omega,f)$ of an acyclic orientation $\Omega$ and
vertex-coloring $f$.

Furthermore, the chromatic polynomial $\chi_G(m)$ has the property that
$(-1)^{|V|} \chi_G(-m)$ counts all such weakly compatible pairs $(\Omega,f)$
in which $f: V \rightarrow \{1,2,\ldots,m\}$ is a vertex-$m$-coloring.
\end{proposition}
\begin{proof}
As observed above, a proper coloring $f$ is weakly compatible with a unique
acyclic orientation $\Omega$ of $G$.  Denote by $P_\Omega$ the
poset on $V$ which is the transitive closure of $\Omega$, endowed
with a \dfn{strict labelling} by integers, that is,
every $i \in P_\Omega$ and $j \in P_\Omega$ satisfying $i<_{P_\Omega} j$ must satisfy $i >_\ZZ j$.
Then proper colorings $f$ that induce $\Omega$ are the same as $P_{\Omega}$-partitions,
so that
\begin{equation}
\label{chromatic-symm-fn-as-sum-of-posets}
\Psi[G] = \sum_{\Omega} F_{P_{\Omega}}(\xx).
\end{equation}
Applying the antipode $S$ and
using Corollary~\ref{antipode-in-Qsym-as-reciprocity-corollary} gives
\[
\omega \left( \Psi[G] \right)
=(-1)^{|V|} S \left( \Psi[G] \right)
= \sum_{\Omega}  F_{P_{\Omega}^{\opp}}(\xx)
= \sum_{(\Omega,f)} \xx_{f}
\]
where in the last line one sums over weakly compatible pairs as in the proposition.  The last equality comes from the fact that since each $P_\Omega$ has been given a strict labelling, $P_\Omega^{\opp}$ acquires a
\emph{weak (or natural) labelling}\index{weak labelling}\index{natural labelling},
that is,
every $i \in P_\Omega$ and $j \in P_\omega$ satisfying $i <_{P_\Omega^{\opp}} j$ must satisfy $i <_\ZZ j$.

The last assertion follows from Proposition~\ref{zeta-polynomiality-proposition}(iii).
\end{proof}

\begin{remark}
The interpretation
of $\chi_G(-m)$ in Proposition~\ref{Stanley-chromatic-reciprocity}
is a much older result of Stanley \cite{Stanley-negative-colors}.
The special case interpreting $\chi_G(-1)$ as $(-1)^{|V|}$ times
the number of acyclic orientations of $G$
has sometimes been called Stanley's \dfn{(-1)-color theorem}.
It also follows (via Proposition~\ref{zeta-polynomiality-proposition})
from Humpert and Martin's antipode formula for $\GGG$ discussed in
Remark~\ref{Humpert-Martin-antipode-remark}: taking $\zeta$ to be
the character of $\GGG$ given in \eqref{edge-free-character-defn},
\[
\chi_G(-1)=\zeta^{\star (-1)}[G] = \zeta(S[G])
=\sum_{F} (-1)^{|V|-\rank(F)} \acyc(G/F) \zeta[G_{V,F}]
=(-1)^{|V|}\acyc(G)
\]
where the last equality uses the vanishing of $\zeta$
on graphs that have edges, so only the $F=\varnothing$ term survives.
\end{remark}

\begin{exercise}
\label{exe.GGG.chromatic-sharp}
If $V$ and $X$ are two sets, and if $f : V \rightarrow X$ is any map,
then $\operatorname*{eqs} f$ will denote the set
\[
\left\{ \left\{ u, u^{\prime} \right\} \ \mid\ u \in V,
\ u^{\prime} \in V, \ u \neq u^{\prime} \text{ and }
f \left( u \right) = f \left( u^{\prime} \right) \right\} .
\]
This is a subset of the set of all two-element subsets of $V$.

If $G = \left( V, E \right)$ is a finite graph, then show that the
map $\Psi$ introduced in Definition~\ref{def.GGG.chromatic} satisfies
\[
\Psi \left( \left[ G \right] ^{\sharp} \right)
= \sum\limits_{\substack{f : V \rightarrow \left\{ 1, 2, 3, \ldots \right\}
; \\ \operatorname*{eqs} f = E}} \xx_{f},
\]
where $\xx_{f} := \prod_{v \in V} x_{f \left( v \right)}$.
Here, $\left[ G\right] ^{\sharp}$ is defined as in
Definition~\ref{def.GGG.moebius-basis}.
\end{exercise}

% [DG][v51] Added the above exercise.

\subsection{Example: The quasisymmetric function of a matroid}

We introduce the \emph{matroid-minor Hopf algebra} of Schmitt \cite{Schmitt},
and studied extensively by
Crapo and Schmitt \cite{CrapoSchmitt1, CrapoSchmitt2, CrapoSchmitt3}.   A very simple character
$\zeta$ on this Hopf algebra will then give rise, via the map $\Psi$ from
Theorem~\ref{Qsym-as-terminal-object-theorem}, to the quasisymmetric function
invariant of matroids from the work of Billera, Jia and the second
author \cite{BilleraJiaR}.

\subsubsection{The matroid-minor Hopf algebra}

We begin by reviewing some notions from matroid theory; see Oxley \cite{Oxley} for background,
undefined terms and unproven facts.

\begin{definition}
A \dfn{matroid} $M$ of rank $r$ on a (finite) ground set $E$
is specified by a nonempty collection
\dfn{$\BBB(M)$} of $r$-element subsets of $E$ with the following
\dfn{exchange property}:
\begin{statement}
For any $B, B'$ in $\BBB(M)$ and $b$ in $B$,
there exists $b'$ in $B'$ with $(B \setminus \{b\}) \cup \{b'\}$ in
$\BBB(M)$.
\end{statement}
The elements of $\BBB\left(M\right)$ are called the
\emph{bases}\index{basis of a matroid} of the matroid $M$.
\end{definition}

\begin{example}
\label{vector-matroid-example}
\begin{comment}
\textbf{Old version:}
A matroid $M$ is \emph{represented} by a collection of vectors $E=\{e_1,\ldots,e_n\}$
in a vector space if $\BBB(M)$ is the collection of subsets $B=\{e_{i_1},\ldots,e_{i_r}\}$
having the property that $B$ forms a basis for the span of all of the vectors in $E$.

For example, if $E=\{a,b,c,d\}$ are the four vectors $a=(1,0),b=(1,1),c=(0,1)=d$ in $\RR^2$
depicted here
\[
\xymatrix{
c,d & b \\
  \ar[r] \arsurj[u] \ar[ur]     & a
}
\]
then $\BBB(M)=\{\{a,b\},\{a,c\},\{a,d\},\{b,c\},\{b,d\}\}$.
\end{comment}
A matroid $M$ with ground set $E$ is
\emph{represented}\index{represented matroid} by
a family of vectors $S = \left(v_e\right)_{e \in E}$
in a vector space if $\BBB\left(M\right)$ is the collection
of subsets $B \subset E$ having the property that
the subfamily $\left(v_e\right)_{e \in B}$ is a basis for
the span of all of the vectors in $S$.

For example, if $M$ is the matroid with
$\BBB(M)=\{\{a,b\},\{a,c\},\{a,d\},\{b,c\},\{b,d\}\}$ on
the ground set $E=\{a,b,c,d\}$, then $M$ is represented by
the family $S = \left(v_a, v_b, v_c, v_d\right)$ of
the four vectors $v_a=(1,0), v_b=(1,1), v_c=(0,1)=v_d$ in
$\RR^2$ depicted here
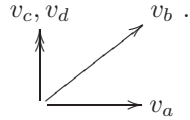
\[
\xymatrix{
v_c,v_d & v_b \\
  \ar[r] \arsurj[u] \ar[ur]     & v_a
} .
\]

% [DG][v29] Rewrote the above example to distinguish between the elements
% of the matroids and the vectors representing them (the correspondence
% is not 1-to-1, so I would rather not identify them).

Conversely, whenever $E$ is a finite set and
$S = \left(v_e\right)_{e \in E}$ is a family of vectors
in a vector space, then the set
\[
\left\{ B \subset E :
\text{the subfamily } \left(v_e\right)_{e \in B}
\text{ is a basis for the span of all of the vectors in } S \right\}
\]
is a matroid on the ground set $E$.

A matroid is said to be \emph{linear}\index{linear matroid}
if there exists a family of vectors in a vector space
representing it. Not all matroids are linear, but many
important ones are.

% [DG][v29] Added preceding two paragraphs.
\end{example}

\begin{example}
\label{graphic-matroid-example}
A special case of matroids $M$ represented by vectors are
\emph{graphic matroids}\index{graphic matroid}, coming
from a graph $G=(V,E)$, with parallel edges and self-loops allowed.  One represents
these by vectors in $\RR^V$ with standard basis $\{ \epsilon_v\}_{v \in V}$ by
associating the vector $\epsilon_v -\epsilon_{v'}$ to any edge connecting a vertex
$v$ with a vertex $v'$.
% [DG][v29] Replaced "to the edge $e=\{v,v'\}$" by "to any edge connecting a vertex
% $v$ with a vertex $v'$" to fit the model of a graph allowing multiple edges.
One can check (or see \cite[\S 1.2]{Oxley}) that the bases $B$ in $\BBB(M)$ correspond to the edge
sets of \emph{spanning forests}\index{spanning forest}
for $G$, that is, edge sets which are acyclic and contain
one spanning tree for each connected component of $G$.
%For example, the graph $G=(V,E)$ shown below has the same matroid $\BBB(M)$
%as the one represented by the vectors in Example~\ref{vector-matroid-example}:
For example, the matroid $\BBB(M)$ corresponding to the graph $G=(V,E)$ shown below:
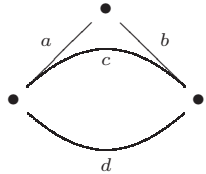
\[
\xymatrix{
&\bullet \ar@{-}[dl]_{a}\ar@{-}[dr]^{b}& \\
\bullet \ar@(ur,ul)@{-}[rr]_{c}  \ar@(dr,dl)@{-}[rr]_{d}& &\bullet
}
\]
%whose spanning trees are the edge sets $\BBB(M)=\{ \{a,b\},\{a,c\},\{a,d\},\{b,c\},\{b,d\}\}$.
is exactly the matroid represented by the vectors in
Example~\ref{vector-matroid-example}; indeed, the spanning forests of
this graph $G$ are the edge sets
$\{a,b\},\{a,c\},\{a,d\},\{b,c\},\{b,d\}$.
% [DG][v29] Rewrote preceding sentence.
(In this example, spanning forests are the same as spanning trees,
since $G$ is connected.)
\end{example}

To define the matroid-minor Hopf algebra one needs the basic matroid
operations of \emph{deletion} and \emph{contraction}.
These  model the operations of deleting or contracting an edge in a graph.
For configurations of vectors they model the deletion of a vector, or the passage to
images in the quotient space modulo the span of a vector.

\begin{definition}
Given a matroid $M$ of rank $r$ and an element $e$ of its ground set $E$,  say that
$e$ is \emph{loop}\index{loop of a matroid}
(resp. \emph{coloop}\index{coloop of a matroid})
of $M$ if $e$ lies in no basis (resp. every basis) $B$ in $\BBB(M)$.
If $e$ is not a coloop, the
\emph{deletion}\index{deletion in a matroid}\index{$M \setminus e$}
$M \setminus e$ is a matroid of rank $r$
on ground set $E \setminus \{e\}$ having bases
\begin{equation}
\label{deletion-bases}
\BBB(M \setminus e):=\{ B \in \BBB(M): e \not\in B\}.
\end{equation}
If $e$ is not a loop, the
\emph{contraction}\index{contraction in a matroid}\index{$M / e$}
$M / e$ is a matroid of rank $r-1$
on ground set $E \setminus \{e\}$ having bases
\begin{equation}
\label{contraction-bases}
\BBB(M/ e):=\{ B \setminus \{e\}: e \in B \in \BBB(M)\}.
\end{equation}
When $e$ is a loop of $M$, then $M/e$ has rank $r$ instead of $r-1$ and one defines its
bases as in \eqref{deletion-bases} rather than \eqref{contraction-bases}; similarly, if $e$ is a coloop
of $M$ then $M \setminus e$ has rank $r-1$ instead of $r$ and one defines its
bases as in \eqref{contraction-bases} rather than \eqref{deletion-bases}.
\end{definition}

\begin{example}
Starting with the graph $G$ and its graphic matroid $M$ from Example ~\ref{graphic-matroid-example}, the deletion $M\setminus a$ and contraction $M/c$
correspond to the graphs $G \setminus a$ and $G/c$ shown here:
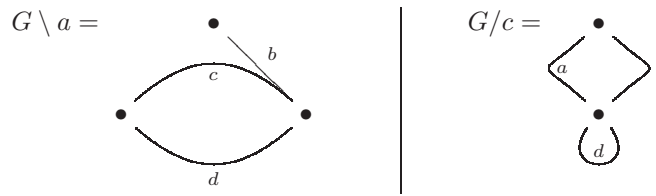
\[
\begin{array}{c|c}
G \setminus a=
\xymatrix{
&\bullet \ar@{-}[dr]^{b}& \\
\bullet \ar@(ur,ul)@{-}[rr]_{c}  \ar@(dr,dl)@{-}[rr]_{d}& &\bullet
}
\qquad \phantom{a}
&
\qquad
G/c=
\xymatrix{
\bullet \ar@(dl,ul)@{-}[d]^{a}\ar@(dr,ur)@{-}[d]^{b}\\
\bullet\ar@(dr,dl) @{-}_{d}
}
\end{array}
\]
One has
\begin{enumerate}
\item[$\bullet$]
$\BBB(M \setminus a)=\{\{b,c\},\{b,d\}\}$,
so that $b$ has become a coloop in $M \setminus a$, and
\item[$\bullet$]
$\BBB(M/c)=\{\{a\},\{b\}\}$,
so that $d$ has become a loop in $M/c$.
\end{enumerate}
\end{example}

\begin{definition}
Deletions and contractions commute with each other. Thus,
given a matroid $M$ with ground set $E$, and a subset
$A \subset E$, two well-defined matroids can be constructed:
\begin{enumerate}
\item[$\bullet$]
the \emph{restriction}\index{restriction of a matroid}\index{$M\mid_A$}
$M|_A$, which is a matroid on ground set $A$,
obtained from $M$ by deleting all
$e \in E \setminus A$ in any order, and
\item[$\bullet$]
the \emph{quotient/contraction}\index{quotient of a matroid}\index{contraction of a matroid}\index{$M/A$}
$M/A$, which is a matroid
on ground set $E \setminus A$,
obtained from $M$ by contracting all $e \in A$ in any order.
\end{enumerate}

% [DG][v59] Reworded the above definition to make it clearer
% (or so I like to think).

We will also need the
\emph{direct sum}\index{direct sum of matroids} \dfn{$M_1 \oplus M_2$}
of two matroids $M_1$ and $M_2$. This is the matroid
whose ground set $E=E_1 \sqcup E_2$ is the disjoint union of a copy of the
ground sets $E_1,E_2$ for $M_1,M_2$, and whose bases
are
\[
\BBB(M_1\oplus M_2):=\{B_1 \sqcup B_2: B_i \in \BBB(M_i) \text{ for }i=1,2\}.
\]

Lastly, say that two matroids $M_1, M_2$ are
\dfn{isomorphic}\index{isomorphic matroids} if there is a bijection
of their ground sets $E_1 \overset{\varphi}{\longrightarrow} E_2$ having the property
that $\varphi \BBB(M_1) = \BBB(M_2)$.
\end{definition}

Now one can define the matroid-minor Hopf algebra, originally introduced by Schmitt \cite[\S 15]{Schmitt}, and studied further by Crapo and Schmitt \cite{CrapoSchmitt1, CrapoSchmitt2, CrapoSchmitt3}.
\begin{definition}
Let \dfn{$\MMM$} have $\kk$-basis  elements $[M]$ indexed by isomorphism classes of
matroids. Define the multiplication via
\[
[M_1] \cdot [M_2]:=[M_1 \oplus M_2] ,
\]
so that the class $[\varnothing]$ of the \dfn{empty matroid} $\varnothing$ having
empty ground set gives a unit.  Define the comultiplication for $M$ a matroid on ground set $E$ via
\[
\Delta[M]:=\sum_{A \subset E} [M|_A] \otimes [M/A],
\]
and a counit
\[
\epsilon [M]:=
\begin{cases}
1, & \text{ if }M=\varnothing \text{;} \\
0, & \text{ otherwise.}
\end{cases}
\]
\end{definition}

\begin{proposition}
The above maps endow $\MMM$ with the structure of a connected graded
finite type Hopf algebra over $\kk$, which is commutative.
\end{proposition}
\begin{proof}
Checking the unit and counit conditions are straightforward.
Associativity and commutativity of the multiplication follow because the direct sum operation
$\oplus$ for matroids is associative and commutative up to isomorphism.  Coassociativity
follows because for a matroid $M$ on ground set $E$,
one has the following equality between the two candidates
for $\Delta^{(2)}[M]$:
\begin{align*}
&\sum_{\varnothing \subset A_1 \subset A_2 \subset E}
[M|_{A_1}] \otimes [(M|_{A_2})/A_1] \otimes [M/A_2] \\
&=
\sum_{\varnothing \subset A_1 \subset A_2 \subset E}[M|_{A_1}] \otimes [(M/A_1)|_{A_2 \setminus A_1}] \otimes [M/A_2]
\end{align*}
due to the matroid isomorphism $(M|_{A_2})/A_1 \cong (M/A_1)|_{A_2 \setminus A_1}$.
Commutativity of the bialgebra diagram in \eqref{bialgebra-diagrams} amounts
to the fact that for a pair of matroids $M_1,M_2$ and subsets $A_1,A_2$ of their
(disjoint) ground sets $E_1,E_2$, one has isomorphisms
\begin{align*}
M_1|_{A_1} \oplus M_2|_{A_2}
  & \cong \left( M_1 \oplus M_2\right)|_{A_1 \sqcup A_2}, \\
M_1/A_1 \oplus M_2/A_2
  & \cong \left( M_1 \oplus M_2\right)/(A_1 \sqcup A_2). \\
\end{align*}
Letting $\MMM_n$ be the $\kk$-span of $[M]$ for matroids whose ground set $E$
has cardinality $|E|=n$, one can then easily check that $\MMM$
becomes a bialgebra which is
graded, connected, and of finite type, hence also a
Hopf algebra by Proposition~\ref{graded-connected-bialgebras-have-antipodes}.
\end{proof}

See \cite{DuchampHoangNghiaKrajewskiTanasa} for an application
of $\MMM$ (and the operator $\exp^{\star}$ from
Section~\ref{subsect.leray}) to proving the \emph{Tutte recipe
theorem}, a ``universal'' property of the Tutte polynomial
of a matroid.

\subsubsection{A quasisymmetric function for matroids}

\begin{definition}
Define  a character $\MMM \overset{\zeta}{\longrightarrow} \kk$ by
\[
\zeta[M]=\begin{cases}
1, & \text{ if }M\text{ has only one basis};\\
0, & \text{ otherwise.}
\end{cases}
\]
It is easily checked that this is a character, that is, an algebra morphism
$\MMM \overset{\zeta}{\longrightarrow} \kk$.
Note that if $M$ has only one basis, say $\BBB(M)=\{B\}$, then
$B:=\coloops(M)$ is the set of coloops of $M$, and
$E\setminus B=\loops(M)$ is the set of loops of $M$.  Equivalently,
$M = \bigoplus_{e \in E} M|_{\{e\}}$ is the direct sum of matroids each
having one element, each a coloop or loop.

Define $\Psi[M]$ for a matroid $M$ to be
the image of $[M]$ under the map
$\MMM\overset{\Psi}{\longrightarrow} \Qsym$
induced via Theorem~\ref{Qsym-as-terminal-object-theorem}
from the above character $\zeta$.
\end{definition}

It turns out that $\Psi[M]$ is intimately related with greedy algorithms and
finding minimum cost bases.
A fundamental property of matroids (and one that characterizes them, in fact; see \cite[\S 1.8]{Oxley})
is that no matter how one assigns costs $f: E \rightarrow \RR$ to the elements of $E$,
the following \emph{greedy algorithm} (generalizing \emph{Kruskal's algorithm} for finding
minimum cost spanning trees) always succeeds
in finding one basis $B$ in $\BBB(M)$ achieving the minimum \emph{total cost} $f(B):=\sum_{b \in B} f(b)$:
\begin{algorithm}
\label{greedy-algorithm}
Start with the empty subset $I_0=\varnothing$ of $E$.
For $j=1,2,\ldots,r$, having already defined the set $I_{j-1}$, let
$e$ be the element of $E \setminus I_{j-1}$ having the lowest cost $f(e)$ among all those
for which $I_{j-1} \cup \{e\}$ is \emph{independent}, that is, still a subset of at least one
basis $B$ in $\BBB(M)$.  Then define $I_j:=I_{j-1} \cup \{e\}$.
Repeat this until $j=r$, and $B=I_r$ will be among the bases that achieve the minimum cost.
\end{algorithm}

\begin{definition}
Say that a cost function $f: E \rightarrow \{1,2,\ldots\}$ is
\emph{$M$-generic}\index{$M$-generic cost function}
if there is a \emph{unique} basis $B$ in $\BBB(M)$
achieving the minimum cost $f(B)$.
\end{definition}

\begin{example}
For the graphic matroid $M$ of Example~\ref{graphic-matroid-example},
this cost function $f_1: E \rightarrow \{1,2,\ldots\}$
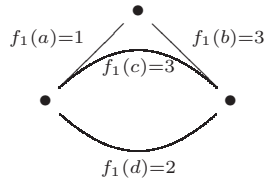
\[
\xymatrix{
&\bullet \ar@{-}[dl]_{f_1(a)=1}\ar@{-}[dr]^{f_1(b)=3}& \\
\bullet \ar@(ur,ul)@{-}[rr]_{f_1(c)={3}}  \ar@(dr,dl)@{-}[rr]_{f_1(d)=2}& &\bullet
}
\]
is $M$-generic, as it minimizes uniquely on the basis $\{a,d\}$,
whereas this cost function $f_2: E \rightarrow \{1,2,\ldots\}$
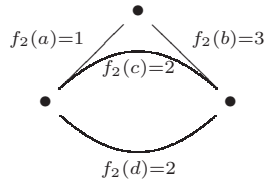
\[
\xymatrix{
&\bullet \ar@{-}[dl]_{f_2(a)=1}\ar@{-}[dr]^{f_2(b)=3}& \\
\bullet \ar@(ur,ul)@{-}[rr]_{f_2(c)=2}  \ar@(dr,dl)@{-}[rr]_{f_2(d)=2}& &\bullet
}
\]
is \emph{not} $M$-generic, as it achieves its minimum value on the two
bases $\{a,c\}, \{a,d\}$.
\end{example}

\begin{proposition}
\label{matroid-qsym-expansion}
For a matroid $M$ on ground set $E$,
one has this expansion\footnote{In fact,
this expansion was the original definition of $\Psi[M]$  in
\cite[Defn. 1.1]{BilleraJiaR}.}
\[
\Psi[M]=\sum\limits_{\substack{M \text{-generic}\\f:E \rightarrow \{1,2,\ldots\}}}
         \xx_f
\]
where $\xx_f:=\prod_{e \in E} x_{f(e)}$.
In particular, for $m \geq 0$, its specialization $ps^1$
from Definition~\ref{principal-specialization-at-1-defn}
has this interpretation:
\[
\ps^1 \Psi[M](m)
=|\{M\text{-generic }f:E \rightarrow \{1,2,\ldots,m\}\}|.
\]
\end{proposition}
\begin{proof}
The iterated coproduct
$\MMM \overset{\Delta^{(\ell-1)}}{\longrightarrow} \MMM^{\otimes \ell}$
sends
\[
[M] \longmapsto
\sum
[M|_{A_1}] \otimes [(M|_{A_2})/A_1] \otimes \cdots \otimes [(M|_{A_{\ell}})/A_{\ell-1}]
\]
where the sum is over flags of nested subsets
\begin{equation}
\label{flag-of-nested-subsets}
\varnothing = A_0 \subset A_1 \subset \cdots
\subset A_{\ell-1} \subset A_{\ell}=E.
\end{equation}
The map $\zeta^{\otimes \ell}$ sends each summand to
$1$ or $0$, depending upon whether each $(M|_{A_j})/A_{j-1}$ has
a unique basis or not.  Thus formula \eqref{ABS-terminal-morphism-formula}
shows that the coefficient $\zeta_\alpha$ of
$x_{i_1}^{\alpha_1} \cdots x_{i_\ell}^{\alpha_\ell}$ in $\Psi[M]$
counts the flags of subsets in \eqref{flag-of-nested-subsets} for
which $|A_j \setminus A_{j-1}|=\alpha_j$ and $(M|_{A_j})/A_{j-1}$
has a unique basis, for each $j$.

Given a flag as in  \eqref{flag-of-nested-subsets},
associate the cost function $f: E \rightarrow \{1,2,\ldots\}$
whose value on each element of $A_j \setminus A_{j-1}$ is $i_j$;  conversely, given
any cost function $f$, say whose distinct values are $i_1 < \cdots < i_\ell$,
one associates the flag having $A_j \setminus A_{j-1}=f^{-1}(i_j)$ for each $j$.

Now, apply the greedy algorithm (Algorithm~\ref{greedy-algorithm})
to find a minimum-cost basis of $M$ for such a cost function $f$.
At each step of the greedy
algorithm, one new element is added to the independent set; these
elements weakly increase in cost as the algorithm
progresses\footnote{\textit{Proof.} Let $e$ be the element added
at step $i$, and let $e'$ be the element added at step $i+1$.
We want to show that $f\left(e\right) \leq f\left(e'\right)$.
But the element $e'$ could already have been added at step $i$.
Since it wasn't, we thus conclude that the element $e$ that was
added instead must have been cheaper or equally expensive.
In other words,
$f\left(e\right) \leq f\left(e'\right)$, qed.}. Thus,
the algorithm first adds some elements of cost $i_1$, then
adds some elements of cost $i_2$, then adds some elements of
cost $i_3$, and so on. We can therefore subdivide the execution
of the algorithm into phases $1, 2, \ldots, \ell$, where each
phase consists of some finite number of steps, such that all
elements added in phase $k$ have cost $i_k$. (A phase may be
empty.) For each $k \in \left\{1, 2, \ldots, \ell\right\}$, we
let $\beta_k$ be the number of steps in phase $k$; in other words,
$\beta_k$ is the number of elements of elements of cost $i_k$
added during the algorithm.

We will prove below, using induction on $s=0,1,2,\ldots,\ell$
the following \textbf{claim}:
After having completed phases $1, 2, \ldots, s$ in the
greedy algorithm (Algorithm~\ref{greedy-algorithm}), there is
\emph{a unique choice}
for the independent set produced thus far,
namely
\begin{equation}
\label{t-step-unique-greedy-basis}
I_{\beta_1+\beta_2+\cdots+\beta_s}=\bigsqcup_{j=1}^{s} \coloops( (M|_{A_j})/A_{j-1} ),
\end{equation}
if and only if each of the matroids $(M|_{A_j})/A_{j-1}$ for $j=1,2,\ldots,s$
\emph{has a unique basis}.

The case $s=\ell$ in this claim would show what we want, namely that
$f$ is $M$-generic, minimizing uniquely on the basis
shown in \eqref{t-step-unique-greedy-basis} with $s=\ell$,
if and only if each $(M|_{A_j})/A_{j-1}$ has a unique basis.

The assertion of the claim is trivially true for $s=0$.
In the inductive step, one may assume that
\begin{enumerate}
\item[$\bullet$]
the independent set $I_{\beta_1+\beta_2+\cdots+\beta_{s-1}}$ takes the form in
\eqref{t-step-unique-greedy-basis}, replacing $s$ by $s-1$,
\item[$\bullet$]
it is the unique $f$-minimizing basis
for $M|_{A_{s-1}}$, and
\item[$\bullet$]
$(M|_{A_j})/A_{j-1}$ has a unique basis for $j=1,2,\ldots,s-1$.
\end{enumerate}
Since $A_{s-1}$ exactly consists of all of the elements $e$ of $E$ whose
costs $f(e)$ lie in the range
$\{i_1,i_2,\ldots,i_{s-1}\}$, in phase $s$ the algorithm
will work in the quotient matroid $M/A_{s-1}$ and
attempt to augment $I_{\beta_1+\beta_2+\cdots+\beta_{s-1}}$
using the next-cheapest elements, namely the
elements of $A_{s} \setminus A_{s-1}$, which all have cost $f$ equal to $i_s$.
Thus the algorithm will have no choices about how to do this augmentation
if and only if $(M|_{A_s})/A_{s-1}$ has a unique basis, namely its set of
coloops, in which case the algorithm will choose to add all of these coloops,
giving $I_{\beta_1+\beta_2+\cdots+\beta_s}$ as described in \eqref{t-step-unique-greedy-basis}.
This completes the induction.

% [DG][v19] Replaced
% "will work in the quotient matroid $M/A_{j-1}$ and" by
% "will work in the quotient matroid $M/A_{s-1}$ and".

% [DG][v59] Corrected the proof, introducing the concept of phases
% and the $\beta_k$. The old version used $\alpha_k$ instead of
% $\beta_k$; but this would mean that the independent set would
% eventually grow to size $\left|E\right|$, which is clearly false.

The last assertion follows from Proposition~\ref{zeta-polynomiality-proposition}.
\end{proof}

\begin{example}
If $M$ has one basis then every function $f:E \rightarrow \{1,2,\ldots\}$ is $M$-generic, and
\[
\Psi[M]=\sum_{f: E \rightarrow \{1,2,\ldots\}} \xx_f = (x_1+x_2+ \cdots)^{|E|} =
M_{(1)}^{|E|}.
\]
\end{example}
\begin{example}
Let $U_{r,n}$ denote the \dfn{uniform matroid} of rank $r$ on $n$ elements $E$, having
$\BBB(U_{r,n})$ equal to all of the $r$-element subsets of $E$.

As $U_{1,2}$ has $E=\{1,2\}$ and $\BBB=\{\{1\},\{2\}\}$,
genericity means $f(1) \neq f(2)$, so
\[
\Psi[U_{1,2}]=\sum\limits_{\substack{(f(1),f(2)):\\f(1) \neq f(2)}} x_{f(1)} x_{f(2)}
=x_1 x_2 + x_2 x_1 + x_1 x_3+x_3 x_1 + \cdots = 2M_{(1,1)}.
\]

Similarly $U_{1,3}$ has $E=\{1,2,3\}$ with $\BBB=\{\{1\},\{2\},\{3\}\}$,
and genericity means either that $f(1),f(2),f(3)$
are all distinct, or that
two of them are the same and the third is smaller.
This shows
\begin{align*}
\Psi[U_{1,3}]&=3\sum_{i<j} x_i x_j^2+6\sum_{i<j<k} x_i x_j x_k\\
&=3M_{(1,2)}+6M_{(1,1,1)} ;\\
\ps^1 \Psi[U_{1,3}](m) &= 3 \binom{m}{2} + 6 \binom{m}{3}=\frac{m(m-1)(2m-1)}{2} .
\end{align*}
One can similarly analyze $U_{2,3}$ and check that
\begin{align*}
\Psi[U_{2,3}]&=3M_{(2,1)}+6M_{(1,1,1)} ;\\
\ps^1 \Psi[U_{2,3}](m) &= 3 \binom{m}{2} + 6 \binom{m}{3}=\frac{m(m-1)(2m-1)}{2} .
\end{align*}
\end{example}

These last examples illustrate the behavior of $\Psi$ under
the duality operation on matroids.

\begin{definition}
Given a matroid $M$ of rank $r$ on ground set $E$, its
\emph{dual}\index{dual matroid} or \dfn{orthogonal matroid}
$M^\perp$ is a matroid of
rank $|E|-r$ on the same ground set $E$, having
\[
\BBB(M^\perp):=\{ E \setminus B \}_{B \in \BBB(M)}.
\]
\end{definition}

See \cite[Theorem 2.1.1]{Oxley} or \cite[Section 4]{Cameron-matroids}
for a proof of the fact that this is well-defined (i.e., that
the collection $\{ E \setminus B \}_{B \in \BBB(M)}$
really satisfies the exchange property).
Here are a few examples of dual matroids.

\begin{example}
The dual of a uniform matroid is another uniform matroid:
\[
U_{r,n}^\perp = U_{n-r,n}.
\]
\end{example}

\begin{example}
If $M$ is matroid of rank $r$ represented by family of vectors $\{e_1,\ldots,e_n\}$ in a vector space over some field $\kk$,
one can find a family of vectors $\{e_1^\perp,\ldots,e_n^\perp\}$ that represent $M^\perp$ in the following way.  Pick a basis for the span of the vectors $\{e_i\}_{i=1}^n$, and create a matrix $A$ in $\kk^{r \times n}$ whose columns express the $e_i$ in terms of this basis.  Then pick any matrix $A^\perp$ whose row space is the null space of $A$, and one finds that the columns $\{e_i^\perp\}_{i=1}^n$ of $A^\perp$ represent $M^\perp$.  See Oxley \cite[\S 2.2]{Oxley}.
\end{example}

\begin{example}
\label{graphical-duality}
Let $G=(V,E)$ be a graph embedded in the plane with edge set $E$, giving rise to a graphic matroid $M$ on ground set $E$.  Let $G^\perp$ be a planar dual of $G$, so that, in particular, for each edge $e$ in $E$, the graph
$G^\perp$ has one edge $e^\perp$, crossing $e$ transversely.
Then the graphic matroid of $G^\perp$ is $M^\perp$.  See Oxley \cite[\S 2.3]{Oxley}.
\end{example}

\begin{proposition}
If $\Psi[M]=\sum_\alpha c_\alpha M_\alpha$ then
$\Psi[M^\perp] = \sum_\alpha c_\alpha M_{\rev(\alpha)}$.

Consequently, $\ps^1 \Psi[M](m) = \ps^1 \Psi[M^\perp](m)$.
\end{proposition}
\begin{proof}
First, let us prove that if
$\Psi[M]=\sum_\alpha c_\alpha M_\alpha$ then
$\Psi[M^\perp] = \sum_\alpha c_\alpha M_{\rev(\alpha)}$.
In other words, let us show that for any given
composition $\alpha$, the coefficient of $M_\alpha$
in $\Psi[M]$ (when $\Psi[M]$ is expanded in the basis
$\left(M_\beta\right)_{\beta \in \Comp}$ of $\Qsym$)
equals the coefficient of $M_{\rev(\alpha)}$ in $\Psi[M^\perp]$.
This amounts to showing that
for any composition $\alpha=(\alpha_1,\ldots,\alpha_\ell)$,
the cardinality of the set of $M$-generic $f$ having $\xx_f=\xx^{\alpha}$
is the same as the cardinality of the set of
$M^\perp$-generic $f^\perp$ having $\xx_{f^\perp}=\xx^{\rev(\alpha)}$.
We claim that the map $f \longmapsto f^\perp$ in which
$f^\perp(e)=\ell+1-f(e)$ gives a bijection between these sets.
To see this, note that any
basis $B$ of $M$ satisfies
\begin{align}
f(B) + f(E \setminus B) &= \sum_{e\in E}f(e) , \label{first-dual-base-relation}\\
f(E \setminus B) + f^\perp(E \setminus B) & = (\ell+1)(|E|-r),
\label{second-dual-base-relation}
\end{align}
where $r$ denotes the rank of $M$.
Thus $B$ is $f$-minimizing
if and only if $E\setminus B$ is $f$-maximizing
(by \eqref{first-dual-base-relation})
if and only if $E \setminus B$ is $f^\perp$-minimizing
(by \eqref{second-dual-base-relation}).
Consequently $f$ is $M$-generic if and only if $f^\perp$ is $M^\perp$-generic.

The last assertion follows, for example, from the calculation in
Proposition~\ref{zeta-polynomiality-proposition}(i) that
$\ps^1 (M_\alpha)(m) = \binom{m}{\ell(\alpha)}$ together with the fact that
$\ell(\rev(\alpha))= \ell(\alpha)$.
\end{proof}

Just as \eqref{chromatic-symm-fn-as-sum-of-posets} showed that
Stanley's chromatic symmetric function of a graph has an expansion as a sum
of $P$-partition enumerators for certain strictly labelled posets\footnote{%
A labelled poset $P$ is said to be
\emph{strictly labelled}\index{strictly labelled poset}
if every two
elements $i$ and $j$ of $P$ satisfying $i <_P j$ satisfy $i >_{\ZZ} j$.}
$P$, the same
holds for $\Psi[M]$.

\begin{definition}
Given a matroid $M$ on ground set $E$, and a basis $B$ in $\BBB(M)$,
define the \dfn{base-cobase poset} \dfn{$P_B$} to have $b < b'$ whenever
$b$ lies in $B$ and $b'$ lies in $E \setminus B$ and
$\left(B  \setminus \{b\}\right) \cup \{b'\}$ is in $\BBB(M)$.
\end{definition}

\begin{proposition}
\label{matroid-qsym-fn-as-sum-of-posets-prop}
For any matroid $M$, one has
$\Psi[M]=\sum_{B \in \BBB(M)} F_{(P_B,\strict)}(\xx)$
where $F_{(P_,\strict)}(\xx)$ for a poset $P$ means the $P$-partition
enumerator for any strict labelling of $P$, i.e.
a labelling such that
the $P$-partitions satisfy $f(i) < f(j)$ whenever $i<_P j$.

In particular, $\Psi[M]$ expands nonnegatively in the
$\{L_\alpha\}$ basis.
\end{proposition}
\begin{proof}
A basic result about matroids, due to Edmonds \cite{Edmonds}, describes
the \emph{edges} in the \emph{matroid base polytope} which is
the convex hull of all vectors $\{ \sum_{b \in B} \epsilon_b\}_{B \in \BBB(M)}$
inside $\RR^E$ with standard basis $\{\epsilon_e\}_{e \in E}$.  He shows
that all such edges connect two bases $B, B'$ that differ by a single
\emph{basis exchange}, that is, $B'=
\left(B  \setminus \{b\}\right) \cup \{b'\}$ for some $b$ in $B$
and $b'$ in $E \setminus B$.

Polyhedral theory then says that
a cost function $f$ on $E$ will minimize uniquely at $B$ if and only if
one has a strict increase
$f(B)<f(B')$ along each such edge $B \rightarrow B'$ emanating from $B$,
that is, if and only if $f(b) < f(b')$ whenever $b<_{P_B} b'$ in the
base-cobase poset $P_B$, that is, $f$ lies in $\AAA(P_B,\strict)$.
\end{proof}

\begin{example}
\label{graphic-matroid-example-continued}
The graphic matroid from Example~\ref{graphic-matroid-example} has
this matroid base polytope, with the bases $B$ in $\BBB(M)$ labelling the
vertices:
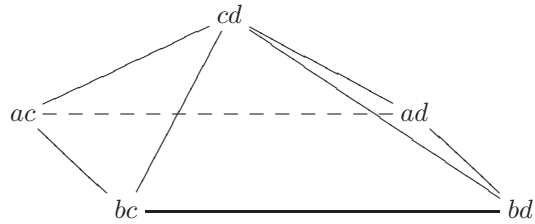
\[
\xymatrix{
 & &cd\ar@{-}[dll]\ar@{-}[ddl]\ar@{-}[drr]\ar@{-}[ddrrr]&&& & \\
ac\ar@{--}[rrrr]\ar@{-}[dr]& & &&ad\ar@{-}[dr]& \\
 &bc\ar@{-}[rrrr]& & &&bd
}
\]
The base-cobase posets $P_B$ for its five vertices $B$ are as follows:
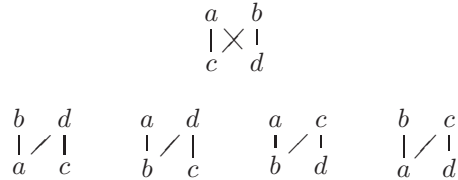
\[
\xymatrix@C=6pt@R=6pt{
a\ar@{-}[d]\ar@{-}[dr] & b\ar@{-}[d] \ar@{-}[dl] \\
c & d
}
\]
\[
\xymatrix@C=6pt@R=6pt{
b\ar@{-}[d]& d\ar@{-}[d] \ar@{-}[dl] \\
a & c
}\qquad
\xymatrix@C=6pt@R=6pt{
a\ar@{-}[d]& d\ar@{-}[d] \ar@{-}[dl] \\
b & c
}\qquad
\xymatrix@C=6pt@R=6pt{
a\ar@{-}[d]& c\ar@{-}[d] \ar@{-}[dl] \\
b & d
}\qquad
\xymatrix@C=6pt@R=6pt{
b\ar@{-}[d]& c\ar@{-}[d] \ar@{-}[dl] \\
a & d
}
\]
One can label the first of these five strictly as
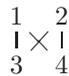
\[
\xymatrix@C=6pt@R=6pt{
1\ar@{-}[d]\ar@{-}[dr] & 2\ar@{-}[d] \ar@{-}[dl] \\
3 & 4
}
\]
and compute its strict $P$-partition enumerator from the linear
extensions $\{3412,3421,4312,4321\}$ as
\[
L_{(2,2)}+L_{(2,1,1)}+L_{(1,1,2)}+L_{(1,1,1,1)} ,
\]
while any of the last four can be labelled strictly as
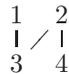
\[
\xymatrix@C=6pt@R=6pt{
1\ar@{-}[d]& 2\ar@{-}[d] \ar@{-}[dl] \\
3 & 4
}
\]
and they each have an extra linear extension $3142$ giving their
strict $P$-partition enumerators as
\[
L_{(2,2)}+L_{(2,1,1)}+L_{(1,1,2)}+L_{(1,1,1,1)}+L_{(1,2,1)}.
\]
Hence one has
\[
\Psi[M]=5L_{(2,2)}+5L_{(1,1,2)}+4L_{(1,2,1)}+5L_{(2,1,1)}+5L_{(1,1,1,1)}.
\]
As $M$ is a graphic matroid for a self-dual planar graph, one has a matroid isomorphism $M \cong M^\perp$ (see Example~\ref{graphical-duality}), reflected in the fact that $\Psi[M]$ is invariant under the symmetry
swapping $M_\alpha \leftrightarrow M_{\rev(\alpha)}$ (and simultaneously
swapping  $L_\alpha \leftrightarrow L_{\rev(\alpha)}$).
\end{example}

This $P$-partition expansion for $\Psi[M]$ also allows us to identify
its image under the antipode of $\Qsym$.

\begin{proposition}
\label{matroid-qsym-reciprocity}
For a matroid $M$ on ground set $E$, one has
\[
S(\Psi[M]) = (-1)^{|E|} \sum_{f: E \rightarrow \{1,2,\ldots\}}
                  |\{f\text{-maximizing bases }B\}| \cdot \xx_f
\]
and
\[
\ps^1\Psi[M](-m) = (-1)^{|E|} \sum_{f: E \rightarrow \{1,2,\ldots,m\}}
                  |\{f\text{-maximizing bases }B\}|.
\]
In particular, the expected number of $f$-maximizing bases among
all cost functions $f:E\rightarrow\left\{  1,2, \ldots ,m\right\}  $
is $(-m)^{-|E|} \ps^1\Psi[M](-m)$.
\end{proposition}

% [DG][v3] I replaced "all cost functions that take at most $m$ values"
% by "all cost functions $f:E\rightarrow\left\{  1,2, \ldots ,m\right\}  $".
% The new wording is clearer, but less expressive (the old one
% presumably implied that the expected number does not depend on the
% choice of $m$-element set to which the cost functions are required
% to go).

\begin{proof}
Corollary~\ref{antipode-in-Qsym-as-reciprocity-corollary} implies
\[
S(\Psi[M])
=\sum_{B \in \BBB(M)} S(F_{(P_B,\strict)}(\xx)) \\
=(-1)^{|E|} \sum_{B \in \BBB(M)} F_{(P_B^{\opp},\weak)}(\xx) ,
\]
where $F_{(P,\weak)}(\xx)$ is the enumerator
for $P$-partitions in which $P$ has been \emph{naturally} labelled,
so that they satisfy $f(i) \leq f(j)$ whenever $i <_P j$.
When $P=P_B^{\opp}$, this is exactly the condition for $f$ to achieve
its maximum value at $f(B)$ (possibly not uniquely), that is,
for $f$ to lie in the \emph{closed} normal cone to the vertex indexed
by $B$ in the matroid base polytope;  compare this with the
discussion in the proof of
Proposition~\ref{matroid-qsym-fn-as-sum-of-posets-prop}.  Thus one has
\[
S(\Psi[M])=(-1)^{|E|} \sum\limits_{\substack{(B,f):\\B \in \BBB(M)\\f\text{ maximizing at }B}}
\xx_f ,
\]
which agrees with the statement of the proposition, after reversing
the order of the summation.

The rest follows from Proposition~\ref{zeta-polynomiality-proposition}.
\end{proof}

\begin{example}
We saw in Example~\ref{graphic-matroid-example-continued}
that the matroid $M$ from Example~\ref{graphic-matroid-example} has
\[
\Psi[M]=5L_{(2,2)}+5L_{(1,1,2)}+4L_{(1,2,1)}+5L_{(2,1,1)}+5L_{(1,1,1,1)},
\]
and therefore will have
\[
%\begin{aligned}
\ps^1 \Psi[M](m)= 5 \binom{m-2+4}{4} + (5+4+5) \binom{m-3+4}{4} + 5 \binom{m-4+4}{4}
            = \frac{m(m-1)(2m^2-2m+1)}{2}
%\end{aligned}
\]
using $\ps^1(L_\alpha)(m) =\binom{m-\ell+|\alpha|}{|\alpha|}$ from
Proposition~\ref{zeta-polynomiality-proposition} (i).
Let us first do a reality-check on a few of its values with
$m \geq 0$ using Proposition~\ref{matroid-qsym-expansion},
and for negative $m$ using Proposition~\ref{matroid-qsym-reciprocity}:
\vskip.1in
\begin{tabular}{|c||c|c|c|c|}\hline
$m$ & $-1$ & $0$ & $1$ & $2$ \\\hline
$\ps^1 \Psi[M](m)$ & $5$ & $0$ & $0$ & $5$\\ \hline
\end{tabular}
\vskip.1in
When $m=0$, interpreting the set of cost functions $f:E \rightarrow \{1,2,\ldots,m\}$
as being empty explains why the value shown is $0$.
When $m=1$, there is only one function $f:E \rightarrow \{1\}$, and it is
not $M$-generic; any of the $5$ bases in $\BBB(M)$ will minimize $f(B)$, explaining
both why the value for $m=1$ is $0$, but also explaining the
value of $5$ for $m=-1$.
The value of $5$ for $m=2$ counts these
$M$-generic cost functions $f:E \rightarrow \{1,2\}$:
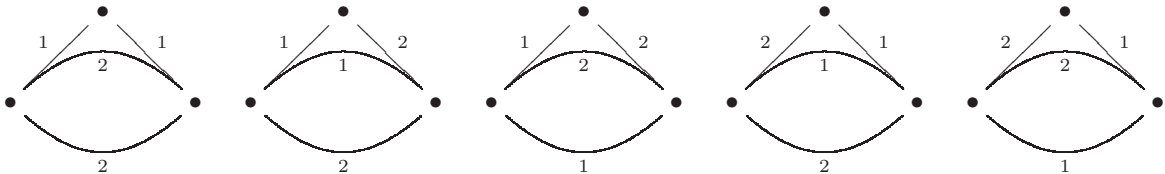
\[
\xymatrix{
&\bullet \ar@{-}[dl]_{1}\ar@{-}[dr]^{1}& \\
\bullet \ar@(ur,ul)@{-}[rr]_{2}  \ar@(dr,dl)@{-}[rr]_{2}& &\bullet
}\quad
\xymatrix{
&\bullet \ar@{-}[dl]_{1}\ar@{-}[dr]^{2}& \\
\bullet \ar@(ur,ul)@{-}[rr]_{1}  \ar@(dr,dl)@{-}[rr]_{2}& &\bullet
}\quad
\xymatrix{
&\bullet \ar@{-}[dl]_{1}\ar@{-}[dr]^{2}& \\
\bullet \ar@(ur,ul)@{-}[rr]_{2}  \ar@(dr,dl)@{-}[rr]_{1}& &\bullet
}\quad
\xymatrix{
&\bullet \ar@{-}[dl]_{2}\ar@{-}[dr]^{1}& \\
\bullet \ar@(ur,ul)@{-}[rr]_{1}  \ar@(dr,dl)@{-}[rr]_{2}& &\bullet
}\quad
\xymatrix{
&\bullet \ar@{-}[dl]_{2}\ar@{-}[dr]^{1}& \\
\bullet \ar@(ur,ul)@{-}[rr]_{2}  \ar@(dr,dl)@{-}[rr]_{1}& &\bullet
}
\]
\noindent
Lastly, Proposition~\ref{matroid-qsym-reciprocity} predicts
the expected number of $f$-minimizing bases for
$f: E \rightarrow \{1,2,\ldots,m\}$ as
\[
(-m)^{-|E|}\ps^1 \Psi[M](-m)= (-m)^{-4}\frac{m(m+1)(2m^2+2m+1)}{2}
=\frac{(m+1)(2m^2+2m+1)}{2m^3},
\]
whose limit as $m \rightarrow \infty$ is $1$, consistent with the notion
that ``most'' cost functions should be generic with respect to the bases of
$M$, and maximize/minimize on a unique basis.
\end{example}

\begin{remark}
It is not coincidental that there is a
similarity of results for Stanley's chromatic symmetric
function of a graph $\Psi[G]$ and for the matroid quasisymmetric
function $\Psi[M]$, such as the $P$-partition expansions
\eqref{chromatic-symm-fn-as-sum-of-posets} versus Proposition~\ref{matroid-qsym-fn-as-sum-of-posets-prop}, and the reciprocity results
Proposition~\ref{Stanley-chromatic-reciprocity} versus
Proposition~\ref{matroid-qsym-reciprocity}.
It was noted in \cite[\S 9]{BilleraJiaR} that one can associate a similar
quasisymmetric function invariant to any \emph{generalized permutohedra}
in the sense of Postnikov \cite{Postnikov}.  Furthermore, recent work of
Ardila and Aguiar \cite{AguiarArdila} has shown that there is a Hopf algebra
of such generalized permutohedra, arising from a \emph{Hopf monoid} in the sense of
Aguiar and Mahajan \cite{AguiarMahajan}.  This Hopf algebra
generalizes the chromatic Hopf algebra of
graphs\footnote{Aguiar and Ardila actually work with a larger
Hopf algebra of graphs.
Namely, their concept of graphs allows parallel edges,
and it also allows ``half-edges'', which have only one endpoint.
If $G = \left(V, E\right)$ is such a graph (where $E$ is the
set of its edges and its half-edges), and if $V'$ is a subset
of $V$, then they define $G/_{V'}$ to be
the graph on vertex set $V'$ obtained from $G$ by
\begin{itemize}
\item removing all vertices that are not in $V'$,
\item removing all edges that have no endpoint in $V'$, and
all half-edges that have no endpoint in $V'$, and
\item replacing all edges that have
only one endpoint in $V'$ by half-edges.
\end{itemize}
(This is to be contrasted with the induced subgraph $G\mid_{V'}$,
which is constructed in the same way but with the edges that
have only one endpoint in $V'$ getting removed as well.)
The comultiplication they define on the Hopf algebra of such
graphs sends the isomorphism class $\left[G\right]$ of a
graph $G = \left(V, E\right)$ to
$\sum\limits_{\left(V_1, V_2\right) : V_1 \sqcup V_2 = V}
\left[G\mid_{V_1}\right] \otimes \left[G/_{V_2}\right]$.
This is no longer a cocommutative Hopf algebra;
our Hopf algebra $\GGG$ is a quotient of it.
In \cite[Corollary 13.10]{AguiarArdila}, Ardila and Aguiar
compute the antipode of the Hopf monoid of such graphs;
this immediately leads to a formula for the antipode of the
corresponding Hopf algebra, because what they call the
Fock functor $\overline{\mathcal{K}}$ preserves antipodes
\cite[Theorem 2.18]{AguiarArdila}.}
and the matroid-minor Hopf algebra, and its quasisymmetric function
invariant derives as usual from Theorem~\ref{Qsym-as-terminal-object-theorem}.
Their work \cite{AguiarArdila} also provides a generalization of the chromatic Hopf algebra
antipode formula of Humpert and Martin \cite{HumpertMartin} discussed
in Remark~\ref{Humpert-Martin-antipode-remark} above.
\end{remark}

\newpage

%%%%%%%%%%%%%%%%%%%%%%%%%%%%%%%%%%%%%%%%%%%%%%%%%%%%%%%%%%%%%%%%%%%%%%%%%%%%%%%
\section{The Malvenuto-Reutenauer Hopf algebra of permutations}
\label{sect.MR}
%%%%%%%%%%%%%%%%%%%%%%%%%%%%%%%%%%%%%%%%%%%%%%%%%%%%%%%%%%%%%%%%%%%%%%%%%%%%%%%

Like so many Hopf algebras we have seen, the \dfn{Malvenuto-Reutenauer
Hopf algebra} $\FQsym$ can be
thought of fruitfully in more than one way.  One is that it gives
a natural noncommutative lift of the quasisymmetric
$P$-partition enumerators and the
fundamental basis $\{L_\alpha\}$ of $\Qsym$, rendering their
product and coproduct formulas even more natural.

\subsection{Definition and Hopf structure}

\begin{definition}
We shall regard permutations as words (over the alphabet
$\left\{1,2,3,\ldots\right\}$) by identifying
every permutation $\pi \in \Symm_n$
with the word $\left(\pi(1), \pi(2), \ldots, \pi(n)\right)$.

Define $\FQsym=\bigoplus_{n \geq 0} \FQsym_n$\index{$\FQsym$}
to be a graded $\kk$-module
in which $\FQsym_n$ has $\kk$-basis $\{F_w\}_{w \in \Symm_n}$ indexed by
the permutations $w=(w_1,\ldots,w_n)$ in $\Symm_n$.
\end{definition}

We first attempt to lift the product and coproduct formulas
\eqref{Qsym-product-on-fundamentals},
\eqref{Qsym-coproduct-on-fundamentals}
in the $\{L_\alpha\}$ basis of $\Qsym$.
We attempt to define a product for $u \in \Symm_k$ and
$v \in \Symm_{\ell}$ as follows\footnote{Recall that we regard
permutations as words.}:
\begin{equation}
\label{FQsym-product-on-F's}
F_u F_v:= \sum_{w \in \shuf{u}{v[k]}} F_w ,
\end{equation}
where for any word $v=(v_1,\ldots,v_\ell)$
we set $v[k]:=(k+v_1,\ldots,k+v_\ell)$. Note that the multiset
$\shuf{u}{v[k]}$ is an actual set in this situation (i.e., has
each element appear only once) and is a subset of $\Symm_{k+\ell}$.

% [DG][v14] Added the previous sentence. Also added "as follows:"
% before the equation.

% [DG][v37] Added "where we regard permutations as words..." part.

The coproduct will be defined using the notation of standardization
of $\std(w)$ a word $w$ in some linearly ordered alphabet (see
Definition~\ref{def.words.std}).

\begin{example}
\label{standardization-example}
Considering words in the Roman alphabet $a<b<c< \cdots$, we have
\[
\begin{array}{rccccccccccl}
\std(b&a&c&c &b&a&a&b&a&c &b)\\
=   (5&1&9&10&6&2&3&7&4&11&8).
\end{array}
\]
\end{example}

Using this, define for $w=(w_1,\ldots,w_n)$ in $\Symm_n$
the element $\Delta F_w \in \FQsym \otimes \FQsym$ by
\begin{equation}
\label{FQsym-coproduct-on-F's}
\Delta F_w := \sum_{k=0}^n
F_{\std(w_1,w_2,\ldots,w_k)} \otimes F_{\std(w_{k+1},w_{k+2},\ldots,w_n)}.
\end{equation}

It is possible to check directly that the maps defined in
\eqref{FQsym-product-on-F's} and
\eqref{FQsym-coproduct-on-F's}
endow $\FQsym$ with the structure of a connected graded finite type Hopf
algebra; see Hazewinkel, Gubareni, Kirichenko \cite[Thm. 7.1.8]{HazewinkelGubareniKirichenko}.
However in justifying this here, we
will follow the approach of Duchamp, Hivert and Thibon
\cite[\S 3]{NCSF6}, which exhibits $\FQsym$ as a subalgebra of a larger
ring of (noncommutative) power series of bounded degree
in a totally ordered alphabet.

\begin{definition}
Given a totally ordered set $I$, create a totally ordered variable
set $\{X_i\}_{i \in I}$, and the ring
\emph{$R\langle \{X_i\}_{i\in I} \rangle$}\index{$R<\{X_i\}_{i \in I}>$}
of \dfn{noncommutative power series of bounded degree}
in this alphabet\footnote{Let us recall the definition of
$R\langle \{X_i\}_{i\in I} \rangle$. \par
Let $N$ denote the \idx{free monoid}
on the alphabet $\{X_i\}_{i\in I}$; it
consists of words $X_{i_1} X_{i_2} \cdots X_{i_k}$. We define a
topological $\kk$-module $\ncps{\kk}{\{X_i\}_{i\in I}}$ to be the
Cartesian product $\kk^N$ (equipped with the product topology),
but we identify its element $\left(\delta_{w, u}\right)_{u \in N}$
with the word $w$ for every $w \in N$. Thus, every element
$\left(\lambda_w\right)_{w \in N} \in \kk^N =
\ncps{\kk}{\{X_i\}_{i\in I}}$ can be rewritten as the convergent sum
$\sum_{w \in N} \lambda_w w$. We call
$\lambda_w$ the \emph{coefficient of $w$} in this element
(or the \emph{coefficient of this element before $w$}). The
elements of $\ncps{\kk}{\{X_i\}_{i\in I}}$ will be referred to as
\dfn{noncommutative power series}.
We define
a multiplication on $\ncps{\kk}{\{X_i\}_{i\in I}}$ by the formula
\[
 \left(\sum_{w \in N} \lambda_w w\right)
 \left(\sum_{w \in N} \mu_w w\right)
 = \sum_{w \in N}
   \left(\sum_{\left(u, v\right) \in N^2; \  w = uv} \lambda_u \mu_v\right)
   w .
\]
(This is well-defined thanks to
the fact that, for each $w \in N$, there are only finitely many
$\left(u, v\right) \in N^2$ satisfying $w = uv$.)
Thus, $\ncps{\kk}{\{X_i\}_{i\in I}}$ becomes a $\kk$-algebra with unity
$1$ (the empty word). (It is similar to the monoid algebra $\kk N$
of $N$ over $\kk$, with the only difference that infinite sums are
allowed.)
\par
Now, we define $R\langle \{X_i\}_{i\in I} \rangle$ to be the
$\kk$-subalgebra of $\ncps{\kk}{\{X_i\}_{i\in I}}$ consisting of all
noncommutative power series
$\sum_{w \in N} \lambda_w w \in \ncps{\kk}{\{X_i\}_{i\in I}}$
\dfn{of bounded degree} (i.e., such that all words $w \in N$ of
sufficiently high length satisfy $\lambda_w = 0$).}.  Many times,
we will use a variable
set $\XX:=(X_1 < X_2 < \cdots)$, and call the ring $R \langle \XX \rangle$.
\end{definition}

% [DG][v50] Added footnote defining noncommutative power series.
% Noone likes to do this...

We first identify the algebra structure for $\FQsym$ as the subalgebra
of finite type within $R\langle \{X_i\}_{i\in I} \rangle$ spanned
by the elements
\begin{equation}
F_w = F_w( \{X_i\}_{i\in I}):= \sum\limits_{\substack{\ii=(i_1,\ldots,i_n):\\ \std(\ii) = w^{-1}}} \XX_\ii ,
\label{eq.FQSym.Fw-as-series}
\end{equation}
where $\XX_{\ii}:=X_{i_1} \cdots X_{i_n}$,
as $w$ ranges over $\bigcup_{n \geq 0} \Symm_n$ .

\begin{example}
\label{F's-example}
For the alphabet $\XX=(X_1 < X_2 < \cdots)$, in $R\langle \XX \rangle$ one has
\begin{align*}
F_{1}&=\sum_{1 \leq i} X_i=X_1 + X_2 + \cdots , \\
F_{12}&=\sum_{1 \leq i \leq j} X_i X_j=X_1^2 + X_2^2 + \cdots + X_1 X_2 + X_1 X_3 + X_2 X_3 + X_1 X_4+\cdots , \\
F_{21}&=\sum_{1 \leq i < j} X_j X_i =X_2 X_1 + X_3 X_1 + X_3 X_2 + X_4 X_1 + \cdots , \\
F_{312}&=\sum_{\ii: \std(\ii)=231} \XX_\ii
        =\sum_{1 \leq i < j \leq k} X_j X_k X_i \\
        &=X_2^2 X_1 + X_3^2 X_1 + X_3^2 X_2 + \cdots
         + X_2 X_3 X_1 + X_2 X_4 X_1 + \cdots .
\end{align*}
\end{example}

\begin{proposition}
\label{prop.FQsym.iso}
For any totally ordered infinite set $I$,
the elements $\{F_w\}$ as $w$ ranges over $\bigcup_{n \geq 0} \Symm_n$
form a $\kk$-basis for a subalgebra $\FQsym( \{X_i\}_{i\in I})$
of $R \langle \XX \rangle$, which is connected graded and of finite type,
having multiplication defined $\kk$-linearly by
\eqref{FQsym-product-on-F's}.

Consequently all such algebras are isomorphic
to a single algebra $\FQsym$, having basis $\{F_w\}$ and
multiplication given by the rule \eqref{FQsym-product-on-F's},
with the isomorphism mapping $F_w \longmapsto F_w (\{X_i\}_{i\in I})$.
\end{proposition}

For example,
\begin{align*}
F_1 F_{21}
&= (X_1 + X_2 + X_3 + \cdots)
    (X_2 X_1 + X_3 X_1 + X_3 X_2 + X_4 X_1 + \cdots) \\
&= X_1 \cdot X_3 X_2 + X_1 \cdot X_4 X_2 + \cdots
 +X_1 \cdot X_2 X_1 + X_2 \cdot X_3 X_2 + X_2 \cdot X_4 X_2 + \cdots \\
&\qquad   +X_2 \cdot X_3 X_1 + X_2 \cdot X_4 X_1 + \cdots
   +X_2 \cdot X_2 X_1 + X_3 \cdot X_3 X_1 + X_3 \cdot X_3 X_2 + \cdots \\
&\qquad   +X_3 \cdot X_2 X_1 + X_4 \cdot X_2 X_1 + \cdots \\
& = \sum_{\ii:\std(\ii)=132} \XX_\ii
+\sum_{\ii:\std(\ii)=231} \XX_\ii
+\sum_{\ii:\std(\ii)=321} \XX_\ii
= F_{132}+F_{312}+F_{321} = \sum_{w \in \shuf{1}{32}} F_w .
\end{align*}

\begin{proof}[Proof of Proposition~\ref{prop.FQsym.iso}]
The elements $\{F_w(\{X_i\}_{i\in I})\}$ are linearly independent as
they are supported on disjoint monomials, and so
form a $\kk$-basis for their span.  The fact that they multiply
via rule \eqref{FQsym-product-on-F's} is the equivalence of
conditions (i) and (iii) in the following Lemma~\ref{standardization-lemma},
from which all the remaining assertions follow.
\end{proof}

\begin{lemma}
\label{standardization-lemma}
For a triple of permutations
\begin{align*}
u&=(u_1,\ldots,u_k) \text{ in }\Symm_k,\\
v&=(v_1,\ldots,v_{n-k}) \text{ in }\Symm_{n-k},\\
w&=(w_1,\ldots,w_n) \text{ in }\Symm_n,
\end{align*}
the following conditions are equivalent:
\begin{enumerate}
\item[(i)]
$w^{-1}$ lies in the set $\shuf{u^{-1}}{v^{-1}[k]}$.
\item[(ii)]
$u=\std(w_1,\ldots,w_k)$ and
$v=\std(w_{k+1},\ldots,w_n)$,
\item[(iii)]
for some word $\ii=(i_1,\ldots,i_n)$
with $\std(\ii)=w$ one has
$u=\std(i_1,\ldots,i_k)$ and
$v=\std(i_{k+1},\ldots,i_n)$.
\end{enumerate}
\end{lemma}
\begin{proof} %[Proof of Lemma~\ref{standardization-lemma}.]
The implication (ii) $\Rightarrow$ (iii) is clear since $\std(w)=w$.
The reverse implication (iii) $\Rightarrow$ (ii) is best illustrated
by example, e.g. considering Example~\ref{standardization-example}
as concatenated, with $n=11$ and $k=6$ and $n-k=5$:
\[
\begin{tabu}{rccccccrrcccccl}
%\begin{array}{rccccccccccccccl}
w=\std&(b&a&c&c&b&a& | &            &a&b&a&c &b)\\
 =   &(5&1&9&10&6&2& | &            &3&7&4&11&8)\\
       & & &  & & &          &    & & & &  & \\
\hline
u=\std&(5&1&9&10 &6&2)& \multicolumn{1}{c||}{\qquad} & v = \std &(3&7&4&11&8)\\
  =&  (3&1&5&6 &4&2)& \multicolumn{1}{c||}{\qquad} & =    &(1&3&2&5 &4)\\
 =\std&(b&a&c&c &b&a)& \multicolumn{1}{c||}{\qquad} & = \std& (a&b&a&c &b)
%\end{array}
\end{tabu}
\]

% [DG][v14] I changed this to a tabu environment (using the "tabu"
% package) and modified the layout. Is it any better? (My main issue
% was that the equalities for $u$ and $v$ seemed to form one common
% equality chain.)

The equivalence of (i) and (ii) is a fairly standard consequence
of unique parabolic factorization $W=W^J W_J$
where $W=\Symm_n$ and $W_J=\Symm_k \times \Symm_{n-k}$,
so that $W^J$ are the  minimum-length coset representatives for cosets
$xW_J$ (that is, the permutations $x \in \Symm_n$ satisfying
$x_1 < \cdots < x_k$ and $x_{k+1} < \cdots < x_n$).
One can uniquely express any $w$ in $W$ as
$w=x y$ with $x$ in $W^J$ and $y$ in $W_J$, which here means
that $y=u \cdot v[k]=v[k] \cdot u$ for some $u$ in $\Symm_k$
and $v$ in $\Symm_{n-k}$.  Therefore $w=x u v[k]$,
if and only if $w^{-1}=u^{-1} v^{-1}[k] x^{-1}$, which means that
$w^{-1}$ is the shuffle of the sequences $u^{-1}$ in positions
$\{x_1,\ldots,x_k\}$ and $v^{-1}[k]$ in positions $\{x_{k+1},\ldots,x_n\}$.
\end{proof}

\begin{example}
To illustrate the equivalence of (i) and (ii) and the parabolic
factorization in the preceding proof, let $n=9$ and $k=5$ with
\begin{align*}
w&=\left(
\begin{matrix}
1 & 2& 3 & 4  & 5 & | & 6 & 7 & 8 & 9 \\
4 & 9 & 6 & 1 & 5 & | & 8 & 2 & 3 & 7
\end{matrix}
\right)\\
&=
\left(
\begin{matrix}
1 & 2& 3 & 4  & 5 & | & 6 & 7 & 8 & 9 \\
1 & 4 & 5 & 6 & 9 & | & 2 & 3 & 7 & 8
\end{matrix}
\right)
\left(\begin{matrix}
1 & 2& 3 & 4  & 5 \\
2 & 5 & 4 & 1 & 3
\end{matrix}
\right)
\left(\begin{matrix}
6 & 7& 8 & 9 \\
9 & 6 & 7 & 8
\end{matrix}
\right)\\
&=x \cdot u \cdot v[k] ;\\
\text{then} &  \\
w^{-1}&=\left(
\begin{matrix}
1 & 2& 3 & 4  & 5  & 6 & 7 & 8 & 9 \\
4 & 9 & 6 & 1 & 5  & 8 & 2 & 3 & 7
\end{matrix}
\right)\\
&=
\left(\begin{matrix}
1 & 2& 3 & 4  & 5 \\
4 & 1 & 5 & 3 & 2
\end{matrix}
\right)
\left(\begin{matrix}
6 & 7& 8 & 9 \\
7 & 8 & 9 & 6
\end{matrix}
\right)
\left(
\begin{matrix}
1 & 2& 3 & 4  & 5 & 6 & 7 & 8 & 9 \\
\underline{1} & \underline{\underline{6}} & \underline{\underline{7}} & \underline{2} & \underline{3}  & \underline{4} & \underline{\underline{8}} & \underline{\underline{9}} & \underline{5}
\end{matrix}
\right) \\
&= u^{-1} \cdot v^{-1}[k] \cdot x^{-1} .
\end{align*}
\end{example}

Proposition~\ref{prop.FQsym.iso} yields that
$\FQsym$ is isomorphic to the $\kk$-subalgebra
$\FQsym\left(\XX\right)$ of the $\kk$-algebra $R\left<\XX\right>$
when $\XX$ is the variable set
$\left(X_1 < X_2 < \cdots\right)$. We identify
$\FQsym$ with $\FQsym\left(\XX\right)$ along this isomorphism.
% [DG][v17] Captain Obvious was here and inserted the above paragraph.
For any infinite alphabet
$\left\{X_i\right\}_{i \in I}$ and any $f \in \FQsym$, we denote
by $f\left(\left\{X_i\right\}_{i \in I}\right)$ the image of $f$
under the algebra isomorphism
$\FQsym \to \FQsym\left(\left\{X_i\right\}_{i \in I}\right)$ defined in
Proposition~\ref{prop.FQsym.iso}.

% [DG][v19] This is copypasted from something I just added to the
% QSym section. I need it below to define the coproduct.

%[VR][v20] Starting here, Vic began rewriting the alphabet-doubling argument...

One can now use this to define a coalgebra structure on $\FQsym$.
Roughly speaking, one wants to first 
evaluate an element
$f$ in $\FQsym \cong \FQsym\left(\XX\right) \cong \FQsym\left(\XX,\YY\right)$
as $f(\XX,\YY)$, using the linearly ordered 
variable set $(\XX,\YY) :=(X_1 < X_2 < \cdots < Y_1 < Y_2 < \cdots )$.
Then one should take the image of $f(\XX,\YY)$ after
imposing the partial commutativity relations 
\begin{equation}
\label{partial-commutativity-relations}
X_i Y_j = Y_j X_i \text{ for every pair }(X_i,Y_j) \in \XX \times \YY,
\end{equation}
and hope that this image lies
in a subalgebra isomorphic to 
\[
\FQsym\left(\XX\right) 
\otimes 
\FQsym\left(\YY\right) 
\cong \FQsym \otimes \FQsym.
\]
We argue this somewhat carefully.   Start by 
considering the canonical monoid epimorphism 
\begin{equation}
\label{semigroup-epimorphism}
F\langle \XX,\YY \rangle \overset{\rho}{\twoheadrightarrow} M,
\end{equation}
where $F\langle \XX,\YY \rangle$ denotes the \emph{free monoid} on 
the alphabet $(\XX,\YY)$ and $M$ denotes its 
quotient monoid imposing the partial commutativity relations 
\eqref{partial-commutativity-relations}.
Let $\kk^M$ denote the $\kk$-module of all functions
$f: M \rightarrow \kk$, with pointwise addition and scalar multiplication;
similarly define $\kk^{F\langle \XX,\YY \rangle}$.
As both monoids $F\langle \XX,\YY \rangle$ and $M$
enjoy the property that an element 
$m$ has only finitely many factorizations as $m=m_1 m_2$,
one can define a convolution algebra structure on 
both $\kk^{F\langle \XX,\YY \rangle}$ and $\kk^M$ via
\[
(f_1 \star f_2)(m) 
= \sum\limits_{\substack{(m_1,m_2) \in N \times N:\\ m=m_1m_2}}
f_1(m_1)f_2(m_2),
\]
where $N$ is respectively $F\langle \XX,\YY \rangle$ or $M$.
As fibers of the map $\rho$ 
in \eqref{semigroup-epimorphism} are finite, it induces a map
of convolution algebras,
which we also call $\rho$:
\begin{equation}
\label{semigroup-algebra-epimorphism}
\kk^{F\langle \XX,\YY\rangle} \overset{\rho}{\twoheadrightarrow} \kk^M.
\end{equation}
Now recall that $R\langle \XX\rangle$ denotes the algebra of
noncommutative formal power
series in the variable set $\XX$, of bounded degree, with coefficients in $\kk$.  One similarly has
the ring $R\langle \XX,\YY \rangle$, which can be identified with
the subalgebra of $\kk^{F\langle \XX,\YY\rangle}$ consisting of the
functions $f: F\langle \XX,\YY \rangle \rightarrow \kk$ having a 
bound on the length of the words
in their support (the value of $f$ on a word in $(\XX,\YY)$ gives
its power series coefficient corresponding to said word). 
We let $R\langle M\rangle$ denote
the analogous subalgebra of $\kk^M$; this can be thought of as the
algebra of bounded degree ``partially commutative power series''
in the variable sets $\XX$ and $\YY$.
Note that $\rho$ restricts to a map 
\begin{equation}
\label{bounded-power-series-morphism}
R\langle \XX,\YY \rangle \overset{\rho}{\rightarrow} R\langle M\rangle.
\end{equation}
Finally, we claim (and see Proposition~\ref{prop.FQsym.coproduct.wd} below 
for a proof) that this further restricts to a map
\begin{equation}
\label{FQsym-coproduct-defining-morphism}
\FQsym\left(\XX,\YY\right) \overset{\rho}{\rightarrow} 
\FQsym\left(\XX\right) \otimes \FQsym\left(\YY\right)
\end{equation}
in which the target is identified with its image under the
(injective\footnote{as images of the basis 
$F_u(\XX) \otimes F_v(\YY)$ of
$\FQsym(\XX) \otimes \FQsym(\YY)$ are supported on disjoint monomials in
$R\langle M\rangle$, so linearly independent.}) multiplication map 
\[
\begin{array}{rcl}
\FQsym\left(\XX\right) \otimes \FQsym\left(\YY\right) 
& \hookrightarrow & R\langle M\rangle , \\
f(\XX) \otimes g(\YY) &\mapsto& f(\XX)g(\YY).
\end{array}
\]
Using the identification of $\FQsym$ with all three of
$\FQsym\left(\XX \right), \FQsym\left(\YY\right), \FQsym\left(\XX,\YY\right)$,
the map $\rho$ in \eqref{FQsym-coproduct-defining-morphism} will then define
a coproduct structure on $\FQsym$.  Abusing notation, for $f$ in $\FQsym$, we will simply write $\Delta(f) = f(\XX,\YY)$ instead of $\rho(f(\XX,\YY))$.

\begin{example}
Recall from Example~\ref{F's-example} that one has
\[
F_{312}=\sum_{\ii: \std(\ii)=231} \XX_\ii
        =\sum_{1 \leq i < j \leq k} X_j X_k X_i ,
\]
and therefore its coproduct is
\begin{align*}
\Delta F_{312}
%&=\overline{\rho}\left(F_{312}(X_1,X_2,\ldots,Y_1,Y_2,\ldots)\right) \\
&=F_{312}(X_1,X_2,\ldots,Y_1,Y_2,\ldots) \qquad \qquad \left(\text{by our abuse of notation}\right)\\
&=\sum_{i < j \leq k} X_j X_k X_i
+\sum\limits_{\substack{i < j,\\ k}} X_j Y_k X_i
+\sum\limits_{\substack{i,\\ j \leq k}} Y_j Y_k X_i
+\sum_{i < j \leq k} Y_j Y_k Y_i\\
&= \sum_{i < j \leq k} X_j X_k X_i \cdot 1
+\sum\limits_{\substack{i < j,\\ k}} X_j X_i \cdot Y_k
+\sum\limits_{\substack{i,\\ j \leq k}} X_i \cdot Y_j Y_k
+\sum_{i < j \leq k} 1 \cdot Y_j Y_k Y_i \\
&=F_{312}(\XX)\cdot 1
+ F_{21}(\XX) \cdot F_1(\YY)
+ F_1(\XX) \cdot F_{12}(\YY)
+ 1 \cdot F_{312}(\YY) \\
&=F_{312} \otimes \one
+ F_{21} \otimes F_1
+ F_1 \otimes F_{12}
+ \one \otimes F_{312}.
\end{align*}
\end{example}

\begin{proposition}
\label{prop.FQsym.coproduct.wd}

\begin{comment}
The image of the $\kk$-algebra homomorphism
$\FQsym \to \powser{\kk}{M_{\XX \sim \YY}}$ sending every $f \in \FQsym$
to $\overline{\rho}\left(f\left(\XX, \YY\right)\right)
\in \powser{\kk}{M_{\XX \sim \YY}}$
lies in $\FQsym \otimes \FQsym$, giving rise to a coproduct
\begin{align*}
\nonumber
\FQsym & \overset{\Delta}{\longrightarrow} \FQsym \otimes \FQsym , \\
f & \longmapsto \overline{\rho}\left(f\left(\XX, \YY\right)\right)
\in \FQsym \otimes \FQsym \subset \powser{\kk}{M_{\XX \sim \YY}}.
\end{align*}
\end{comment}

The map $\rho$ in \eqref{bounded-power-series-morphism}
does restrict as claimed to a map 
as in \eqref{FQsym-coproduct-defining-morphism},
and hence defines a coproduct on $\FQsym$,
acting on the $\{F_w\}$ basis by
the rule \eqref{FQsym-coproduct-on-F's}.  
This endows $\FQsym$ with
the structure of a connected graded finite type Hopf algebra.
\end{proposition}
\begin{proof}

\begin{comment}
The $\kk$-algebra homomorphism
$\FQsym \to \powser{\kk}{M_{\XX \sim \YY}}$ sending every $f \in \FQsym$
to $\overline{\rho}\left(f\left(\XX, \YY\right)\right)
\in \powser{\kk}{M_{\XX \sim \YY}}$ sends every $F_w$ to
\end{comment}

Let $I$ be the totally ordered set
$\left\{1 < 2 < 3 < \cdots\right\}$.
Let $J$ be the totally ordered set \newline
$\left\{1 < 2 < 3 < \cdots < \widetilde 1 < \widetilde 2 < \widetilde 3 < \cdots \right\}$.
We set $X_{\widetilde i} = Y_i$ for every positive integer $i$. Then,
the alphabet $(\XX, \YY)$ can be written as $\left\{X_i\right\}_{i \in J}$.

If $\ii$ is a word over the alphabet $I = \left\{1 < 2 < 3 < \cdots\right\}$,
then we denote by $\widetilde{\ii}$ the word over $J$
obtained from $\ii$ by replacing every letter $i$ by $\widetilde i$.

For the first assertion of Proposition~\ref{prop.FQsym.coproduct.wd},
it suffices to check that $F_w$ indeed has the image under $\Delta$
claimed in \eqref{FQsym-coproduct-on-F's}. Let $n \in \NN$ and
$w \in \Symm_n$. Then,
\begin{align}
\Delta F_w
% \overline{\rho}\left(F_w(\XX,\YY)\right)
 &= F_w\left(\XX, \YY\right)  \qquad \qquad \left(\text{by our abuse of notation}\right) \nonumber \\
 &= \sum_{\ii \in J^n : \std(\ii)=w^{-1} } (\XX,\YY)_{\ii} = \sum_{\mathbf{t} \in J^n : \std(\mathbf{t})=w^{-1} } (\XX,\YY)_{\mathbf{t}} \nonumber \\
 &= \sum_{k=0}^n \sum_{(\ii,\jj) \in I^k \times I^{n-k}} \sum\limits_{\substack{\mathbf{t} \in J^n : \\ \std(\mathbf{t})=w^{-1} ; \\ \mathbf{t} \in \shuf{\ii}{\widetilde{\jj}}}} (\XX,\YY)_{\mathbf{t}}
\label{pf.prop.FQsym.coproduct.wd.1}
\end{align}
(since for every $\mathbf{t} \in J^n$, there exists exactly one choice of $k \in \left\{0,1,\ldots,n\right\}$ and $\left(\ii,\jj\right) \in I^k \times I^{n-k}$ satisfying $\mathbf{t} \in \shuf{\ii}{\widetilde{\jj}}$; namely, $\ii$ is the restriction of $\mathbf{t}$ to the subalphabet $I$ of $J$, whereas $\jj$ is the restriction of $\mathbf{t}$ to $J \setminus I$, and $k$ is the length of $\ii$).

We now fix $k$ and $\left(\ii, \jj\right)$, and try to simplify the inner sum
$\sum\limits_{\substack{\mathbf{t} \in J^n : \\ \std(\mathbf{t})=w^{-1} ; \\ \mathbf{t} \in \shuf{\ii}{\widetilde{\jj}}}} (\XX,\YY)_{\mathbf{t}}$ on the right hand side of \eqref{pf.prop.FQsym.coproduct.wd.1}. First we notice that this sum is nonempty if and only if there exists some $\mathbf{t} \in \shuf{\ii}{\widetilde{\jj}}$ satisfying $\std(\mathbf{t})=w^{-1}$. This existence is easily seen to be equivalent to $w^{-1} \in \shuf{\std(\ii)}{\std(\jj)[k]}$ (since the standardization of any shuffle in $\shuf{\ii}{\widetilde{\jj}}$ is the corresponding shuffle in $\shuf{\std(\ii)}{\std(\jj)[k]}$). This, in turn, is equivalent to $\std(\ii)=(\std(w_1,\ldots,w_k))^{-1}$ and $\std(\jj)=(\std(w_{k+1},\ldots,w_n))^{-1}$ (according to the equivalence (i) $\Longleftrightarrow $ (ii) in Lemma~\ref{standardization-lemma}). Hence, the inner sum on the right hand side of \eqref{pf.prop.FQsym.coproduct.wd.1} is nonempty if and only if $\std(\ii)=(\std(w_1,\ldots,w_k))^{-1}$ and $\std(\jj)=(\std(w_{k+1},\ldots,w_n))^{-1}$. When it is nonempty, it has only one addend\footnote{In fact, the elements $\std\left(\mathbf{t}\right)$ for $\mathbf{t} \in \shuf{\ii}{\widetilde{\jj}}$ are distinct, and thus only one of them can equal $w^{-1}$.}, and this addend is $\left(\XX,\YY\right)_{\mathbf{t}} = \XX_\ii \YY_\jj$ (since $\mathbf{t} \in \shuf{\ii}{\widetilde{\jj}}$). Summarizing, we see that the inner sum on the right hand side of \eqref{pf.prop.FQsym.coproduct.wd.1} equals $\XX_\ii \YY_\jj$ when $\std(\ii)=(\std(w_1,\ldots,w_k))^{-1}$ and $\std(\jj)=(\std(w_{k+1},\ldots,w_n))^{-1}$, and is empty otherwise. Thus, \eqref{pf.prop.FQsym.coproduct.wd.1} simplifies to 
\begin{align*}
\Delta F_w
 &= \sum_{k=0}^n \sum\limits_{\substack{(\ii,\jj) \in I^k \times I^{n-k}:\\
          \std(\ii)=(\std(w_1,\ldots,w_k))^{-1}\\
            \std(\jj)=(\std(w_{k+1},\ldots,w_n))^{-1}}} \XX_\ii \YY_\jj \\
 &= \sum_{k=0}^n F_{\std(w_1,\ldots,w_k)}(\XX) F_{\std(w_{k+1},\ldots,w_n)}(\YY) \\
 &= \sum_{k=0}^n F_{\std(w_1,\ldots,w_k)} \otimes F_{\std(w_{k+1},\ldots,w_n)} \in \FQsym \otimes \FQsym .
\end{align*}
This proves \eqref{FQsym-coproduct-on-F's}, and thus the first assertion
of Proposition~\ref{prop.FQsym.coproduct.wd}.

% [DG][v28] Expanded the above argument (I need the details in front
% of my eyes since I am checking an analogous proof for the
% $m$-free quasisymmetric functions).

\begin{comment}
This proves the first statement of the proposition, in particularly
showing that $\Delta$ is well-defined. 
Of course, we now have
\[
\Delta F_w = \overline{\rho}\left(F_w(\XX,\YY)\right) = \sum_{k=0}^n F_{\std(w_1,\ldots,w_k)} \otimes F_{\std(w_{k+1},\ldots,w_n)},
\]
so that the action of $\Delta$ is indeed given by \eqref{FQsym-coproduct-on-F's}.
\end{comment}

From this, it is easy to derive that $\Delta$ satisfies coassociativity
(i.e., the diagram \eqref{coassociativity-diagram} holds for
$C = \FQsym$). (Alternatively, one can obtain this from the
associativity of multiplication using Corollary~\ref{FQsym-self-dual}.)
We have already verified the rule \eqref{FQsym-coproduct-on-F's}.
The connected graded structure on $\FQsym$ gives a counit and an antipode
for free.
\end{proof}

%[VR][v20] This ends the area where Vic made changes to the alphabet-doubling arguments. 

\begin{exercise}
\label{exe.FQSym.free}
We say that a permutation $w\in \Symm_{n}$ is
\emph{connected}\index{connected permutation}
if $n$ is a positive integer and if there exists no
$i\in\left\{  1,2,\ldots,n-1\right\}  $ satisfying $f\left(  \left\{
1,2,\ldots,i\right\}  \right)  =\left\{  1,2,\ldots,i\right\}  $. Let $\mathfrak{CS}%
$ denote the set of all connected permutations of all $n\in \NN$. Show
that $ \FQsym $ is a free (noncommutative) $\kk$-algebra
with generators $\left(  F_{w}\right)  _{w\in\mathfrak{CS}}$. (This statement
means that $\left(  F_{w_{1}}F_{w_{2}}\cdots F_{w_{k}}\right)  _{k\in
 \NN ;\ \left(  w_{1},w_{2},\ldots,w_{k}\right)  \in\mathfrak{CS}^{k}}$
is a basis of the $\kk$-module $\FQsym$.)

[\textbf{Hint:} This is a result of Poirier and Reutenauer \cite[Theorem
2.1]{PoirierReutenauer}; it is much easier than the similar Theorem
\ref{thm.QSym.lyndon}.]
\end{exercise}

% [DG][v37] Added above exercise.

\begin{corollary}
\label{FQsym-self-dual}
The Hopf algebra $\FQsym$ is self-dual:
Let $\{G_w\}$ be the dual $\kk$-basis
to the $\kk$-basis $\{F_w\}$ for $\FQsym$. Then,
the $\kk$-linear map sending $G_w \longmapsto F_{w^{-1}}$
is a Hopf algebra isomorphism $\FQsym^{o} \longrightarrow \FQsym$.
\end{corollary}
\begin{proof}
For any $0 \leq k \leq n$, any $u \in \Symm_k$ and any
$v \in \Symm_{n-k}$, one has
\[
F_{u^{-1}} F_{v^{-1}}
= \sum_{w^{-1} \in \shuf{u^{-1}}{v^{-1}[k]}} F_{w^{-1}}
= \sum\limits_{\substack{w \in \Symm_n:\\
           \std(w_1,\ldots,w_k)=u \\
           \std(w_{k+1},\ldots,w_n)=v }} F_{w^{-1}}
\]
via the equivalence of (i) and (ii) in Lemma~\ref{standardization-lemma}.
On the other hand, in $\FQsym^o$, the dual $\kk$-basis $\{G_w\}$
to the $\kk$-basis $\{F_w\}$ for $\FQsym$ should have product
formula
\[
G_u G_v = \sum\limits_{\substack{w \in \Symm_n:\\
           \std(w_1,\ldots,w_k)=u \\
           \std(w_{k+1},\ldots,w_n)=v }} G_w
\]
coming from the coproduct formula \eqref{FQsym-coproduct-on-F's}
for $\FQsym$ in the $\{F_w\}$-basis. Comparing these equalities,
we see that the $\kk$-linear map $\tau$ sending $G_w \longmapsto F_{w^{-1}}$
is an isomorphism $\FQsym^{o} \longrightarrow \FQsym$ of
$\kk$-algebras. Hence, the adjoint $\tau^\ast : \FQsym^{o}
\to \left(\FQsym^{o}\right)^{o}$ of this map is an isomorphism
of $\kk$-coalgebras. But identifying $\left(\FQsym^{o}\right)^{o}$
with $\FQsym$ in the natural way (since $\FQsym$ is of finite type),
we easily see that $\tau^\ast = \tau$, whence $\tau$ itself is
an isomorphism of both $\kk$-algebras and $\kk$-coalgebras, hence
of $\kk$-bialgebras, hence of Hopf algebras.
\end{proof}

% [DG][v25] Revamped the above corollary. It used to speak of a map
% $\FQsym \longrightarrow \FQsym^{o}$ sending $F_w$ to $F_{w^{-1}}$,
% which was more than confusing as there was no identification of
% $\FQsym$ with its dual given a priori. Added some more details to
% the proof.

We can now be a bit more precise about the relations between
the various algebras
\[
\Lambda, \Qsym, \Nsym, \FQsym, R \langle \XX \rangle, R(\xx).
\]
Not only does $\FQsym$ allow one to \emph{lift} the Hopf structure
of $\Qsym$, it dually allows one to \emph{extend} the Hopf structure
of $\Nsym$.  To set up this duality,
note that Corollary~\ref{FQsym-self-dual}
motivates the choice of an inner product on $\FQsym$ in which
\[
(F_u, F_v):=\delta_{u^{-1},v}.
\]
We wish to identify the images of the ribbon basis
$\{R_\alpha\}$ of $\Nsym$ when included in $\FQsym$.

\begin{definition}
For any composition $\alpha$, define an element
\dfn{$\RRR_\alpha$} of $\FQsym$ by
% the \dfn{free quasi-ribbon function}
\[
\RRR_\alpha
:=\sum_{\substack{w \in \Symm_{\left|\alpha\right|}: \\ \Des(w)=D(\alpha) }} F_{w^{-1}}
=\sum\limits_{\substack{(w,\ii): \\ w \in \Symm_{\left|\alpha\right|}; \\ \Des(w)=D(\alpha);\\ \std(\ii)=w}} \XX_\ii
=\sum_{\ii: \Des(\ii)=D(\alpha)} \XX_\ii ,
\]
% where the $w$ in the sums are supposed to belong to $\Symm_{\left|\alpha\right|}$, and
where the \emph{descent set}\index{descent set of a sequence}
of a sequence $\ii=(i_1,\ldots,i_n)$
is defined by
\[
\Des(\ii):=\{j \in \{1,2,\ldots,n-1\}: i_j > i_{j+1} \}=
\Des(\std(\ii)).
\]
Alternatively,
\begin{equation}
\RRR_\alpha= \sum_T \XX_T
\label{eq.FQSym.ribbon.tableau-sum}
\end{equation}
in which the sum is over column-strict tableaux of
the ribbon skew shape $\Rib\left(\alpha\right)$,
and $\XX_T=\XX_{\ii}$
in which $\ii$ is the sequence of entries of $T$ read in
order from the southwest toward the northeast.

\end{definition}

% [DG][v14] I commented out "the \dfn{free quasi-ribbon function}" because
% it seems to me that this notation is used for the $F_w$ in NCSF6 and in
% other sources, but never for the $\RRR_\alpha$.

\begin{example}
Taking $\alpha=(1,3,2)$, with ribbon shape and column-strict fillings $T$ as
shown:
\[
\Rib\left(\alpha\right) =
\begingroup % keep the change local
\setlength\arraycolsep{0.2pc} % Make the spacing equal in rows and in columns.
\begin{matrix}
   &   &\sq&\sq\\
\sq&\sq&\sq& \\
\sq&   &   &
\end{matrix}
\endgroup
\qquad\qquad \text{and} \qquad\qquad
T=\quad
\begin{matrix}
   &   &    &   &i_5&\leq&i_6\\
   &   &    &   &\wedge&    &   \\
i_2&\leq&i_3&\leq&i_4& \\
\wedge&   &   &     &  & \\
i_1&   &   &     &  &
\end{matrix}
\]
one has that
\[
\RRR_{(1,3,2)}=
\sum\limits_{\substack{ \ii=(i_1,i_2,i_3,i_4,i_5,i_6):\\
                   \Des(\ii)=D(\alpha)=\{1,4\}}} \XX_\ii
=\sum_{i_1 > i_2 \leq i_3 \leq i_4 > i_5 \leq i_6}
   X_{i_1}  X_{i_2} X_{i_3}  X_{i_4}  X_{i_5} X_{i_6}
=\sum_T \XX_T .
\]
\end{example}

\begin{corollary}
\label{cor.FQsym.pi-iota}
For every $n \in \NN$ and $w \in \Symm_n$, we let \dfn{$\gamma(w)$}
denote the unique composition $\alpha$ of $n$ satisfying
$D\left(\alpha\right) = \Des\left(w\right)$.
% [DG][v50] Added the preceding sentence; the $\gamma(w)$ notation
% is used a few lines below.
% Also, slightly reorganized the next sentence.
\begin{enumerate}
\item[(a)]
The $\kk$-linear map\index{$\pi : \FQsym \to \Qsym$}
\[
\begin{array}{rcl}
\FQsym &\overset{\pi}{\twoheadrightarrow}& \Qsym , \\
F_w & \longmapsto &L_{\gamma(w)} \\
\end{array}
\]
is a surjective Hopf algebra homomorphism.
\item[(b)]
The $\kk$-linear map\index{$\iota : \Nsym \to \FQsym$}
\[
\begin{array}{rcl}
\Nsym &\overset{\iota}{\hookrightarrow}& \FQsym , \\
R_\alpha & \longmapsto & \RRR_\alpha
\end{array}
\]
is an injective Hopf algebra homomorphism.
\item[(c)]
The linear maps $\pi$ and $\iota$ are adjoint maps with respect to the
above choice of inner product on $\FQsym$
and the usual dual pairing between $\Nsym$ and $\Qsym$.
\end{enumerate}
Now, consider the abelianization map
$\operatorname{ab} : R\langle \XX \rangle \twoheadrightarrow R(\xx)$
defined as the continuous $\kk$-algebra homomorphism
sending the noncommutative variable $X_i$ to the commutative $x_i$.
\begin{enumerate}
\item[(d)]
The map $\pi$ is a restriction of $\operatorname{ab}$.
\item[(e)]
The map $\iota$ lets one factor
the surjection $\Nsym \twoheadrightarrow \Lambda$ as follows:
\[
\begin{array}{rcccl}
\Nsym & \rightarrow & \FQsym \hookrightarrow R\langle \XX \rangle & \overset{\operatorname{ab}}{\rightarrow} & R(\xx) , \\
 R_\alpha &\longmapsto & \RRR_\alpha & \longmapsto & s_{\Rib\left(\alpha\right)}(\xx) .
\end{array}
\]
\end{enumerate}
\end{corollary}

\begin{proof}
Given $n \in \NN$, each composition $\alpha$ of $n$ can be written in
the form $\gamma\left(w\right)$ for some $w \in \Symm_n$.
\ \ \ \ \footnote{Indeed, write our composition $\alpha$
as $\left(\alpha_1, \alpha_2, \ldots, \alpha_k\right)$.
Then, we can pick $w$ to be the permutation whose first
$\alpha_1$ entries are the largest $\alpha_1$ elements of
$\left\{1,2,\ldots,n\right\}$ in increasing order;
whose next $\alpha_2$ entries are the next-largest $\alpha_2$
elements of $\left\{1,2,\ldots,n\right\}$ in increasing order;
and so on.
This permutation $w$ will satisfy
$\Des\left(w\right) = \left\{\alpha_1, \alpha_1+\alpha_2, \ldots, \alpha_1+\alpha_2+\cdots+\alpha_{k-1}\right\}
= D\left(\alpha\right)$ and thus
$\gamma\left(w\right) = \alpha$.}
Hence, each fundamental quasisymmetric function $L_\alpha$ lies
in the image of $\pi$. Thus, $\pi$ is surjective.

Also, for each $n \in \NN$ and $\alpha \in \Comp_n$,
the element $\RRR_\alpha$ is a nonempty sum of noncommutative
monomials (nonempty because $\alpha$ can be written in
the form $\gamma\left(w\right)$ for some $w \in \Symm_n$).
Moreover, the elements $\RRR_\alpha$ for varying $n$ and
$\alpha$ are supported on disjoint monomials.
Thus, these elements are linearly independent.
Hence, the map $\iota$ is injective.

(d) Let $\mathfrak{A}$ denote the totally ordered set
$\left\{  1<2<3<\cdots\right\}  $ of positive integers.
For each word
$w = \left( w_1, w_2, \ldots, w_n \right) \in
\mathfrak{A}^{n}$, we define a monomial
$\xx_w$ in $\kk \left[\left[\xx\right]\right]$
by $\xx_w = x_{w_1} x_{w_2} \cdots x_{w_n}$.

Let $n \in \NN$ and $\sigma \in \Symm_n$. Then,
\[
L_{\gamma\left( \sigma \right)  }
= \sum\limits_{\substack{w\in\mathfrak{A}^{n};
\\ \std w=\sigma^{-1}}} \xx_w
\]
(by Lemma~\ref{lem.Lalpha-std.std-L}).
But \eqref{eq.FQSym.Fw-as-series} (applied to $w = \sigma$) yields
\begin{align*}
F_{\sigma}
&= \sum\limits_{\substack{\ii=(i_1,\ldots,i_n):\\ \std(\ii) = \sigma^{-1}}} \XX_\ii
= \sum\limits_{\substack{w\in\mathfrak{A}^{n};\\ \std w=\sigma^{-1}}}\XX_{w}
\end{align*}
and thus
\begin{align*}
\operatorname{ab}\left(F_{\sigma}\right)
&= \operatorname{ab}\left(\sum\limits_{\substack{w\in\mathfrak{A}^{n};
\\ \std w=\sigma^{-1}}}\XX_{w}\right)
= \sum\limits_{\substack{w\in\mathfrak{A}^{n};\\ \std w=\sigma^{-1}}} \underbrace{\operatorname{ab}\left(\XX_w\right)}_{=\xx_w}
= \sum\limits_{\substack{w\in\mathfrak{A}^{n}; \\ \std w=\sigma^{-1}}} \xx_w
= L_{\gamma\left( \sigma \right)  }
= \pi\left(F_\sigma\right) .
\end{align*}
We have shown this for all $n \in \NN$ and $\sigma \in \Symm_n$.
Thus, $\pi$ is a restriction of $\operatorname{ab}$.
This proves Corollary~\ref{cor.FQsym.pi-iota}(d).

(a) Let $n \in \NN$ and $w = \left(w_1, w_2, \ldots, w_n\right) \in \Symm_n$.
Let $\alpha$ be the composition $\gamma\left(w\right)$ of $n$.
Thus, the definition of $\pi$ yields
$\pi\left(F_w\right) = L_\alpha$.
But applying the map $\pi \otimes \pi$ to the equality
\eqref{FQsym-coproduct-on-F's}, we obtain
\begin{align}
\left(\pi \otimes \pi\right) \left(\Delta F_w\right)
&= \left(\pi \otimes \pi\right)
    \left( \sum_{k=0}^n F_{\std(w_1,w_2,\ldots,w_k)}
                \otimes F_{\std(w_{k+1},w_{k+2},\ldots,w_n)}
   \right) \nonumber\\
&= \sum_{k=0}^n \pi\left(F_{\std(w_1,w_2,\ldots,w_k)}\right) \otimes \pi\left(F_{\std(w_{k+1},w_{k+2},\ldots,w_n)}\right) \nonumber\\
&= \sum_{k=0}^n L_{\gamma\left(\std(w_1,w_2,\ldots,w_k)\right)} \otimes L_{\gamma\left(\std(w_{k+1},w_{k+2},\ldots,w_n)\right)}
\label{pf.cor.FQsym.pi-iota.a.1}
\end{align}
(by the definition of $\pi$).
Now, for each $k \in \left\{0,1,\ldots, n\right\}$, the two
compositions $\gamma\left(\std(w_1,w_2,\ldots,w_k)\right)$
$\gamma\left(\std(w_{k+1},w_{k+2},\ldots,w_n)\right)$ form
a pair $\left(\beta, \gamma\right)$ of compositions
satisfying\footnote{See Definition~\ref{def.compositions.near-conc}
for the notation we are using.}
either $\beta \cdot \gamma=\alpha$ or $\beta \odot \gamma=\alpha$,
and in fact they form the only such pair satisfying
$\left|\beta\right| = k$ and $\left|\gamma\right| = n-k$.
Thus, the right hand side of \eqref{pf.cor.FQsym.pi-iota.a.1}
can be rewritten as
\[
\sum\limits_{\substack{(\beta,\gamma):\\
\beta \cdot \gamma=\alpha \text{ or }
\beta \odot \gamma=\alpha}} L_\beta \otimes L_\gamma .\]
But this sum is $\Delta L_\alpha$, as we know from
\eqref{Qsym-coproduct-on-fundamentals}.
Hence, \eqref{pf.cor.FQsym.pi-iota.a.1} becomes
\[
\left(\pi \otimes \pi\right) \left(\Delta F_w\right)
= \Delta L_\alpha = \Delta \left(\pi\left(F_w\right)\right)
\qquad\qquad \left(\text{since $L_\alpha = \pi\left(F_w\right)$}\right).
\]

We have proven this for each $n \in \NN$ and $w \in \Symm_n$.
Thus, we have proven that
$\left(\pi \otimes \pi\right) \circ \Delta_{\FQsym}
= \Delta_{\Qsym} \circ \pi$.
Combined with $\epsilon_{\FQsym} = \epsilon_{\Qsym} \circ \pi$
(which is easy to check), this shows that $\pi$ is a
coalgebra homomorphism.

We can similarly see that $\pi$ is an algebra homomorphism
by checking that it respects the product
(compare \eqref{Qsym-product-on-fundamentals}
and \eqref{FQsym-product-on-F's}).
However, this also follows trivially from
Corollary~\ref{cor.FQsym.pi-iota}(d).

Thus, $\pi$ is a bialgebra morphism, and therefore a Hopf
algebra morphism (by
Corollary~\ref{cor.bialg-mor-is-Hopf}).
This proves Corollary~\ref{cor.FQsym.pi-iota}(a).

(c) For any composition $\alpha$ and any $w \in \Symm$,
we have
\begin{align*}
(\iota(R_\alpha), F_w)
&=(\RRR_\alpha,F_w)
=\sum_{u: \Des(u)=D(\alpha)} (F_{u^{-1}}, F_w)
=
\begin{cases}
1, &\text{ if }\Des(w)=D(\alpha) ; \\
0, &\text{ otherwise }
\end{cases}
\quad
=
\begin{cases}
1, &\text{ if }\gamma(w)=\alpha ; \\
0, &\text{ otherwise }
\end{cases}
\\
&=(R_\alpha, L_{\gamma(w)})
=(R_\alpha,\pi(F_w)).
\end{align*}
Thus, the maps $\pi$ and $\iota$ are adjoint.
This proves Corollary~\ref{cor.FQsym.pi-iota}(c).

(b) Again, there are several ways to prove this.
Here is one:

First, note that $\iota\left(1\right) = 1$
(because $R_{\varnothing} = 1$ and $\RRR_{\varnothing} = 1$).
Next, let
$\alpha$ and $\beta$ be two nonempty compositions.
Let $m = \left|\alpha\right|$ and $n = \left|\beta\right|$.
Then,
$R_\alpha R_\beta = R_{\alpha \cdot \beta} + R_{\alpha \odot \beta}$
(by \eqref{Nsym-product-in-ribbons}) and thus
\begin{align}
\iota\left(R_\alpha R_\beta\right)
 &= \iota\left(R_{\alpha \cdot \beta} + R_{\alpha \odot \beta}\right)
 = \underbrace{\iota\left(R_{\alpha \cdot \beta}\right)}_{= \RRR_{\alpha \cdot \beta} = \sum_{\ii: \Des(\ii)=D(\alpha \cdot \beta)} \XX_\ii}
  + \underbrace{\iota\left(R_{\alpha \odot \beta}\right)}_{= \RRR_{\alpha \odot \beta} = \sum_{\ii: \Des(\ii)=D(\alpha \odot \beta)} \XX_\ii} \nonumber\\
 &= \sum_{\ii: \Des(\ii)=D(\alpha \cdot \beta)} \XX_\ii + \sum_{\ii: \Des(\ii)=D(\alpha \odot \beta)} \XX_\ii
 = \sum_{\ii: \Des(\ii)=D(\alpha \cdot \beta) \text{ or } \Des(\ii)=D(\alpha \odot \beta)} \XX_\ii \nonumber\\
 &= \sum\limits_{\substack{\ii = \left(i_1, i_2, \ldots, i_{m+n}\right): \\ \Des\left(i_1, i_2, \ldots, i_m\right) = D\left(\alpha\right) \text{ and} \\ \Des\left(i_{m+1}, i_{m+2}, \ldots, i_{m+n}\right) = D\left(\beta\right)}} \XX_\ii
\label{pf.cor.FQsym.pi-iota.b.1}
\end{align}
(since the words $\ii$ of length $m+n$ satisfying
$\Des(\ii)=D(\alpha \cdot \beta)$ or $\Des(\ii)=D(\alpha \odot \beta)$
are precisely the words $\ii = \left(i_1, i_2, \ldots, i_{m+n}\right)$
satisfying
$\Des\left(i_1, i_2, \ldots, i_m\right) = D\left(\alpha\right)$ and
$\Des\left(i_{m+1}, i_{m+2}, \ldots, i_{m+n}\right) = D\left(\beta\right)$).
But choosing a word $\ii = \left(i_1, i_2, \ldots, i_{m+n}\right)$
satisfying $\Des\left(i_1, i_2, \ldots, i_m\right) = D\left(\alpha\right)$
and
$\Des\left(i_{m+1}, i_{m+2}, \ldots, i_{m+n}\right) = D\left(\beta\right)$
is tantamount to choosing a pair $\left(\mathbf{u}, \mathbf{v}\right)$
of a word $\mathbf{u} = \left(i_1, i_2, \ldots, i_m\right)$ satisfying
$\Des\mathbf{u} = D\left(\alpha\right)$ and a word
$\mathbf{v} = \left(i_{m+1}, i_{m+2}, \ldots, i_{m+n}\right)$ satisfying
$\Des\mathbf{v} = D\left(\beta\right)$.
Thus, \eqref{pf.cor.FQsym.pi-iota.b.1} becomes
\begin{align*}
\iota\left(R_\alpha R_\beta\right)
 &= \sum\limits_{\substack{\ii = \left(i_1, i_2, \ldots, i_{m+n}\right): \\ \Des\left(i_1, i_2, \ldots, i_m\right) = D\left(\alpha\right) \text{ and} \\ \Des\left(i_{m+1}, i_{m+2}, \ldots, i_{m+n}\right) = D\left(\beta\right)}} \XX_\ii
 =      \sum_{\mathbf{u} : \Des\mathbf{u} = D\left(\alpha\right)}
    \ \ \sum_{\mathbf{v} : \Des\mathbf{v} = D\left(\beta\right)}
        \XX_{\mathbf{u}} \XX_{\mathbf{v}} \\
 &= \underbrace{\left(\sum_{\mathbf{u} : \Des\mathbf{u} = D\left(\alpha\right)} \XX_{\mathbf{u}}\right)}_{= \RRR_\alpha = \iota\left(R_\alpha\right)}
    \underbrace{\left(\sum_{\mathbf{v} : \Des\mathbf{v} = D\left(\beta\right)} \XX_{\mathbf{v}}\right)}_{= \RRR_\beta = \iota\left(R_\beta\right)}
 = \iota\left(R_\alpha\right) \iota\left(R_\beta\right) .
\end{align*}
Thus, we have proven the equality
$\iota\left(R_\alpha R_\beta\right)
= \iota\left(R_\alpha\right) \iota\left(R_\beta\right)$
whenever $\alpha$ and $\beta$ are two nonempty compositions.
It also holds if we drop the ``nonempty'' requirement
(since $R_{\varnothing} = 1$ and $\iota\left(1\right) = 1$).
Thus, the $\kk$-linear map $\iota$ respects the multiplication.
Since $\iota\left(1\right) = 1$, this shows that $\iota$
is a $\kk$-algebra homomorphism.

For each $n \in \NN$, we let $\id_n$ be the identity
permutation in $\Symm_n$.
Next, we observe that each $n \in \NN$ satisfies
$H_n = R_{\left(n\right)}$ (this follows, e.g., from
\eqref{H-as-sum-of-ribbons}, because the composition
$\left(n\right)$ is coarsened only by itself).
Hence, each $n \in \NN$ satisfies
\begin{align}
\iota\left(H_n\right)
&= \iota\left(R_{\left(n\right)}\right)
= \RRR_{\left(n\right)}
= \sum_{\substack{w \in \Symm_n:\\ \Des(w)=D\left(\left(n\right)\right) }} F_{w^{-1}} \nonumber\\
&= F_{\id_n^{-1}}
\qquad \left(\text{since the only $w \in \Symm_n$ satisfying $\Des(w)=D\left(\left(n\right)\right)$ is $\id_n$}\right) \nonumber\\
&= F_{\id_n} .
\label{pf.cor.FQsym.pi-iota.b.iotaH}
\end{align}

In order to show that $\iota$ is a $\kk$-coalgebra
homomorphism, it suffices to check the equalities
$\left(\iota \otimes \iota\right) \circ \Delta_{\Nsym}
= \Delta_{\FQsym} \circ \iota$
and
$\epsilon_{\Nsym} = \epsilon_{\FQsym} \circ \iota$.
We shall only prove the first one, since the second is easy.
Since $\iota$, $\Delta_{\Nsym}$ and $\Delta_{\FQsym}$
are $\kk$-algebra homomorphisms, it suffices to check it
on the generators $H_1, H_2, H_3, \ldots$ of $\Nsym$.
But on these generators, it follows from comparing
\begin{align*}
\left(\left(\iota \otimes \iota\right) \circ \Delta_{\Nsym}\right) \left(H_n\right)
&= \left(\iota \otimes \iota\right) \left( \Delta_{\Nsym} H_n \right)
= \left(\iota \otimes \iota\right) \left( \sum_{i+j=n} H_i \otimes H_j \right)
\qquad \left(\text{by \eqref{Nsym-coproduct-on-H}}\right) \\
&= \sum_{i+j=n} \underbrace{\iota\left(H_i\right)}_{\substack{= F_{\id_i} \\ \text{(by \eqref{pf.cor.FQsym.pi-iota.b.iotaH})}}}
  \otimes \underbrace{\iota\left(H_j\right)}_{\substack{= F_{\id_j} \\ \text{(by \eqref{pf.cor.FQsym.pi-iota.b.iotaH})}}}
= \sum_{i+j=n} F_{\id_i} \otimes F_{\id_j}
= \sum_{k=0}^n F_{\id_k} \otimes F_{\id_{n-k}}
\end{align*}
with
\begin{align*}
\left(\Delta_{\FQsym} \circ \iota\right) \left(H_n\right)
&= \Delta_{\FQsym} \left(\iota\left(H_n\right)\right)
= \Delta_{\FQsym} \left( F_{\id_n} \right)
\qquad \left(\text{by \eqref{pf.cor.FQsym.pi-iota.b.iotaH}}\right) \\
&= \sum_{k=0}^n F_{\id_k} \otimes F_{\id_{n-k}}
\qquad \left(\text{by \eqref{FQsym-coproduct-on-F's}}\right) .
\end{align*}

Thus, we know that $\iota$ is a $\kk$-algebra homomorphism
and a $\kk$-coalgebra homomorphism.
Hence, $\iota$ is a bialgebra morphism, and therefore a Hopf
algebra morphism
(by Corollary~\ref{cor.bialg-mor-is-Hopf}).
This proves Corollary~\ref{cor.FQsym.pi-iota}(b).

An alternative proof of Corollary~\ref{cor.FQsym.pi-iota}(b)
can be obtained by adjointness from
Corollary~\ref{cor.FQsym.pi-iota}(a).
Both the inner product on $\FQsym$ and the dual pairing
$\left(\cdot, \cdot\right) : \Nsym \otimes \Qsym \to \kk$
respect the Hopf structures (i.e., the maps
$\Delta_{\Nsym}$ and $m_{\Qsym}$ are mutually adjoint with
respect to these forms, and so are the maps
$m_{\Nsym}$ and $\Delta_{\Qsym}$, and the maps
$\Delta_{\FQsym}$ and $m_{\FQsym}$, and so on).
Corollary~\ref{cor.FQsym.pi-iota}(c) shows that the map
$\iota$ is adjoint to the map $\pi$ with respect to
these two bilinear forms.
Hence, we have a commutative diagram
\[
\xymatrix{
\Nsym \arinjrev[r]^{\iota} \ar[d]^{\cong} & \FQsym \ar[d]^{\cong} \\
\Qsym^o \ar[r]_{\pi^{\ast}} & \FQsym^o
}
\]
of Hopf algebras (where the two vertical arrows are the isomorphisms
induced by the two bilinear forms).
Thus, Corollary~\ref{cor.FQsym.pi-iota}(b) follows from
Corollary~\ref{cor.FQsym.pi-iota}(a) by duality.

(e) For each composition $\alpha$,
the abelianization map $\operatorname{ab}$
sends the noncommutative tableau monomial $\XX_T$ to the commutative tableau
monomial $\xx_T$
whenever $T$ is a tableau of ribbon shape $\Rib\left(\alpha\right)$.
Thus, $\operatorname{ab}$ sends $\RRR_\alpha$ to
$s_{\Rib\left(\alpha\right)}(\xx)$ (because of the formula
\eqref{eq.FQSym.ribbon.tableau-sum}).
Hence, the composition
$\Nsym \rightarrow \FQsym \hookrightarrow R\langle \XX \rangle \overset{\operatorname{ab}}{\rightarrow} R(\xx)$
does indeed send $R_\alpha$ to $s_{\Rib\left(\alpha\right)}(\xx)$.
But so does the projection $\pi : \Nsym \to \Lambda$,
according to Theorem~\ref{Nsym-structure-on-ribbons-theorem}(b).
Hence, the composition factors the projection.
This proves Corollary~\ref{cor.FQsym.pi-iota}(e).
\end{proof}

We summarize some of this picture as follows:
\[
\xymatrix@M+2pt{
\FQsym \ar@{.}[rr]^{\text{dual}} &         & \FQsym \arsurj[d]^{\pi}\\
\Nsym \ar@{.}[rr]^{\text{dual}} \arsurj[d]^{\pi} \arinjrev[u]_{\iota} &         & \Qsym \\
\Lambda \ar@{.}[rr]^{\text{dual}} & & \Lambda \arinjrev[u]
}
\]
Furthermore, if we denote by $\iota$ the canonical inclusion
$\Lambda \to \Qsym$ as well, then the diagram
\[
\xymatrix{
& \FQsym \arsurj[dr]^{\pi}\\
\Nsym \arsurj[dr]^{\pi} \arinjrev[ur]_{\iota} &         & \Qsym \\
       &\Lambda\arinjrev[ur]_{\iota}& 
}
\]
is commutative (according to
Corollary~\ref{cor.FQsym.pi-iota}(e)).

\begin{remark}
Different notations for $\FQsym$ appear in the literature. In the book
\cite{BlessenohlLaue} (which presents an unusual approach to the
character theory of the symmetric group using $\FQsym$), the Hopf
algebra $\FQsym$ is called $\mathcal{P}$, and its basis that we call
$\left\{G_w\right\}_{w\in\Symm_n}$ is denoted
$\left\{w\right\}_{w\in\Symm_n}$. In
\cite[Chapter 7]{HazewinkelGubareniKirichenko}, the Hopf algebra
$\FQsym$ and its basis $\left\{F_w\right\}_{w\in\Symm_n}$ are denoted
$MPR$ and $\left\{w\right\}_{w\in\Symm_n}$, respectively.
\end{remark}

% [DG][v58] Added the above remark.

% [DG][v25] What do you think about moving all the following (skeletal)
% sections inside "Further topics" before arXiving (lowering them by 1
% level, so they become subsections)?

%%%%%%%%%%%%%%%%%%%%%%%%%%%%%%%%%%%%%
\section{Further topics}
\label{further-topics-section}
%%%%%%%%%%%%%%%%%%%%%%%%%%%%%%%%%%%%%

The following is a list of topics that were, at one point, planned to
be touched in class, but did not make the cut. They might get
elaborated upon in a future version of these notes.

% [DG][v58] Added the above paragraph, replacing the drafty
% "Some of these we may touch on in class, others are appropriate for student talks.".
% Also, added references for the topics below.

%%%%%%%%%%%%%%%%%%%%%%%%%%%%%%%%%%%%%
\subsubsection{$0$-Hecke algebras}
\label{0-Hecke-section}
%%%%%%%%%%%%%%%%%%%%%%%%%%%%%%%%%%%%%

\begin{itemize}

\item \textbf{Review of representation theory of finite-dimensional algebras.}

Review the notions of indecomposables, simples, projectives, along
with the theorems of Krull-Remak-Schmidt, of Jordan-H\"older,
and the two kinds of Grothendieck groups dual to each other.

\item \textbf{$0$-Hecke algebra representation theory.}

Describe  the simples and projectives, following
Denton, Hivert, Schilling, Thiery \cite{DentonHivertSchillingThiery}
on $\JJJ$-trivial monoids.

\item \textbf{Nsym and Qsym as Grothendieck groups.}

Give Krob and Thibon's interpretation (see \cite[\S 5]{Thibon} for a
brief summary) of
\begin{enumerate}
\item[$\bullet$]
$\Qsym$ and the Grothendieck group of
composition series, and
\item[$\bullet$]
$\Nsym$ and the Grothendieck group of projectives.
\end{enumerate}

\begin{remark}
Mention P. McNamara's interpretation, in the case of \emph{supersolvable lattices}, of the Ehrenborg quasisymmetric function as the composition series enumerator for an $H_n(0)$-action on the maximal chains
\end{remark}

\end{itemize}

%%%%%%%%%%%%%%%%%%%%%%%%%%%%%%%%%%%%%
\subsubsection{Aguiar-Bergeron-Sottile character theory Part II: Odd and even
characters, subalgebras}
\label{ABS-section}
%%%%%%%%%%%%%%%%%%%%%%%%%%%%%%%%%%%%%

%%%%%%%%%%%%%%%%%%%%%%%%%%%%%%%%%%%%%
\subsubsection{Face enumeration, Eulerian posets, and cd-indices}
\label{cd-index-section}
%%%%%%%%%%%%%%%%%%%%%%%%%%%%%%%%%%%%%

Borrowing from Billera's ICM notes \cite{Billera}.

\begin{itemize}

\item f-vectors, h-vectors

\item flag f-vectors, flag h-vectors

\item ab-indices and cd-indices

\end{itemize}

% [DG][v37] I've put the unwritten parts here, into the "Further
% topics" list.

%%%%%%%%%%%%%%%%%%%%%%%%%%%%%%%%%%%%%
\subsubsection{Other topics}
%%%%%%%%%%%%%%%%%%%%%%%%%%%%%%%%%%%%%

\begin{itemize}
\item Loday-Ronco Hopf algebra of planar binary trees \cite{LodayRonco-trees}
\item Poirier-Reutenauer Hopf algebra of tableaux
\item Reading Hopf algebra of Baxter permutations
\item Hopf monoids, e.g. of Hopf algebra of generalized permutohedra, of matroids, of graphs, Stanley chromatic symmetric functions and Tutte polynomials
\item Lam-Pylyavskyy Hopf algebra of set-valued tableaux
\item Connes-Kreimer Hopf algebra and renormalization
\item Noncommutative symmetric functions and $\Omega \Sigma \CC P^\infty$
\item Maschke's theorem and ``integrals'' for Hopf algebras
\item Nichols-Zoeller structure theorem and group-like elements
\item Cartier-Milnor-Moore structure theorem and primitive elements
\item Quasi-triangular Hopf algebras and quantum groups
\item The Steenrod algebra, its dual, and tree Hopf algebras
\item Ringel-Hall algebras of quivers
\item Ellis-Khovanov odd symmetric function Hopf algebras \cite{EllisKhovanov}
(see also Lauda-Russell \cite{LaudaRussell})
\end{itemize}

Student talks given in class were:
\begin{enumerate}
\item Al Garver, on Maschke's theorem for finite-dimensional Hopf algebras
\item Jonathan Hahn, on the paper by Humpert and Martin.
\item Emily Gunawan, on the paper by Lam, Lauve and Sottile.
\item Jonas Karlsson, on the paper by Connes and Kreimer
\item Thomas McConville, on Butcher's group and generalized Runge-Kutta methods.
\item Cihan Bahran, on universal enveloping algebras and the Poincar\'e-Birkhoff-Witt theorem.
\item Theodosios Douvropolos, on the Cartier-Milnor-Moore theorem.
\item Alex Csar, on the Loday-Ronco Hopf algebra of binary trees
\item Kevin Dilks, on Reading's Hopf algebra of (twisted) Baxter permutations
\item Becky Patrias, on the paper by Lam and Pylyavskyy
\item Meng Wu, on multiple zeta values and Hoffman's homomorphism from $\Qsym$
\end{enumerate}

\section{Some open problems and conjectures}
\begin{itemize}
\item Is there a proof of the Assaf-McNamara skew Pieri rule that
gives a resolution of Specht or Schur/Weyl modules whose character
corresponds to $s_{\lambda/\mu} h_n$, whose terms model their alternating sum?

\item Explicit antipodes in the Lam-Pylyavskyy Hopf algebras? (Answered by Patrias in \cite{Patrias-antipode}.)

\item P. McNamara's question \cite[Question 7.1]{McNamaraWard}:  are $P$-partition enumerators irreducible for connected posets $P$?

\item Stanley's question:  are the only $P$-partition enumerators
which are symmetric (not just quasisymmetric) those for which $P$ is
a skew shape with a column-strict labelling?

\item Does Stanley's chromatic symmetric function distinguish trees?

\item Hoffman's stuffle conjecture

\item Billera-Brenti's nonnegativity conjecture for
the total $cd$-index of  Bruhat intervals
(\cite[Conjecture 6.1]{BilleraBrenti})

\end{itemize}

\newpage

\section{\label{chp.appendix}Appendix: Some basics}

In this appendix, we briefly discuss some basic notions from linear algebra
and elementary combinatorics that are used in these notes.

\subsection{\label{sect.STmat}Linear expansions and triangularity}

In this Section, we shall recall some fundamental results from linear algebra
(most importantly, the notions of a change-of-basis matrix and of a
unitriangular matrix), but in greater generality than how it is usually done
in textbooks. We shall use these results later when studying bases of
combinatorial Hopf algebras; but per se, this section has nothing to do with
Hopf algebras.

\subsubsection{Matrices}

Let us first define the notion of a matrix whose rows and columns are indexed
by arbitrary objects (as opposed to numbers):\footnote{As before, $\kk$
denotes a commutative ring.}

\begin{definition}
\label{def.STmat.STmat}
Let $S$ and $T$ be two sets. An \dfn{$S\times T$-matrix over $\kk$}
shall mean a family
$\left(  a_{s,t} \right)  _{\left(  s,t\right)  \in S\times T}
\in \kk^{S\times T}$ of
elements of $\kk$ indexed by elements of $S\times T$. Thus, the set of
all $S\times T$-matrices over $\kk$ is $\kk^{S\times T}$.

We shall abbreviate ``$S\times T$-matrix over $\kk$'' by
``$S\times T$-matrix'' when the value of $\kk$ is clear from the context.
\end{definition}

This definition of $S\times T$-matrices generalizes the usual notion of
matrices (i.e., the notion of $n\times m$-matrices): Namely, if $n\in
 \NN$ and $m\in \NN$, then the $\left\{  1,2,\ldots,n\right\}
\times\left\{  1,2,\ldots,m\right\}  $-matrices are precisely the $n\times
m$-matrices (in the usual meaning of this word). We shall often use the word
``\dfn{matrix}'' for both the usual notion
of matrices and for the more general notion of $S\times T$-matrices.

Various concepts defined for $n\times m$-matrices (such as addition and
multiplication of matrices, or the notion of a row) can be generalized to
$S\times T$-matrices in a straightforward way. The following four definitions
are examples of such generalizations:

\begin{definition}
\label{def.STmat.+}Let $S$ and $T$ be two sets.

\begin{enumerate}
\item[(a)] The sum of two $S\times T$-matrices is defined by $\left(
a_{s,t}\right)  _{\left(  s,t\right)  \in S\times T}+\left(  b_{s,t}\right)
_{\left(  s,t\right)  \in S\times T}=\left(  a_{s,t}+b_{s,t}\right)  _{\left(
s,t\right)  \in S\times T}$.

\item[(b)] If $u\in\kk$ and if $\left(  a_{s,t}\right)  _{\left(
s,t\right)  \in S\times T}\in\kk^{S\times T}$, then we define $u\left(
a_{s,t}\right)  _{\left(  s,t\right)  \in S\times T}$ to be the $S\times
T$-matrix $\left(  ua_{s,t}\right)  _{\left(  s,t\right)  \in S\times T}$.

\item[(c)] Let $A=\left(  a_{s,t}\right)  _{\left(  s,t\right)  \in S\times
T}$ be an $S\times T$-matrix. For every $s\in S$, we define the
\emph{$s$-th row of $A$}\index{$s$-th row of an $S \times T$-matrix}
to be the $\left\{  1\right\}  \times T$-matrix $\left(
a_{s,t}\right)  _{\left(  i,t\right)  \in\left\{  1\right\}  \times T}$.
(Notice that $\left\{  1\right\}  \times T$-matrices are a generalization of
row vectors.) Similarly, for every $t\in T$, we define the
\emph{$t$-th column of $A$}\index{$t$-th column of an $S \times T$-matrix}
to be the $S\times\left\{  1\right\}  $-matrix $\left(
a_{s,t}\right)  _{\left(  s,i\right)  \in S\times\left\{  1\right\}  }$.
\end{enumerate}
\end{definition}

\begin{definition}
\label{def.STmat.1}Let $S$ be a set.

\begin{enumerate}
\item[(a)] The $S\times S$ \dfn{identity matrix} is defined to be the
$S\times S$-matrix $\left(  \delta_{s,t}\right)  _{\left(  s,t\right)  \in
S\times S}$. This $S\times S$-matrix is denoted by $I_{S}$. (Notice that the
$n\times n$ identity matrix $I_{n}$ is $I_{\left\{  1,2,\ldots,n\right\}  }$
for each $n\in \NN$.)

\item[(b)] An $S\times S$-matrix $\left(  a_{s,t}\right)  _{\left(
s,t\right)  \in S\times S}$ is said to be
\emph{diagonal}\index{diagonal $S\times S$-matrix}
if every $\left(  s,t\right)  \in S\times T$ satisfying
$s\neq t$ satisfies $a_{s,t}=0$.

\item[(c)] Let $A=\left(  a_{s,t}\right)  _{\left(  s,t\right)  \in S\times
S}$ be an $S\times S$-matrix. The
\emph{diagonal}\index{diagonal of an $S\times S$-matrix}
of $A$ means the family
$\left(  a_{s,s}\right)  _{s\in S}$. The
\emph{diagonal entries}\index{diagonal entries of an $S\times S$-matrix}
of $A$ are
the entries of this diagonal $\left(  a_{s,s}\right)  _{s\in S}$.
\end{enumerate}
\end{definition}

\begin{definition}
\label{def.STmat.*}Let $S$, $T$ and $U$ be three sets. Let $A=\left(
a_{s,t}\right)  _{\left(  s,t\right)  \in S\times T}$ be an $S\times
T$-matrix, and let $B=\left(  b_{t,u}\right)  _{\left(  t,u\right)  \in
T\times U}$ be a $T\times U$-matrix. Assume that the sum $\sum_{t\in T}%
a_{s,t}b_{t,u}$ is well-defined for every $\left(  s,u\right)  \in S\times U$.
(For example, this is guaranteed to hold if the set $T$ is finite. For
infinite $T$, it may and may not hold.) Then, the $S\times U$-matrix $AB$ is
defined by%
\[
AB=\left(  \sum_{t\in T}a_{s,t}b_{t,u}\right)  _{\left(  s,u\right)  \in
S\times U}.
\]

\end{definition}

\begin{definition}
\label{def.STmat.inv}Let $S$ and $T$ be two finite sets. We say that an
$S\times T$-matrix $A$ is
\emph{invertible}\index{invertible $S \times T$-matrix}
if and only if there exists a
$T\times S$-matrix $B$ satisfying $AB=I_{S}$ and $BA=I_{T}$. In this case,
this matrix $B$ is unique; it is denoted by $A^{-1}$ and is called the
\emph{inverse}\index{inverse of an $S \times T$-matrix} of $A$.
\end{definition}

The definitions that we have just given are straightforward generalizations of
the analogous definitions for $n\times m$-matrices; thus, unsurprisingly, many
properties of $n\times m$-matrices still hold for $S\times T$-matrices. For example:

\begin{proposition}
\phantomsection\label{prop.STmat.ass}

\begin{enumerate}
\item[(a)] Let $S$ and $T$ be two sets. Let $A$ be an $S\times T$-matrix.
Then, $I_{S}A=A$ and $AI_{T}=A$.

\item[(b)] Let $S$, $T$ and $U$ be three sets such that $T$ is finite. Let $A$
and $B$ be two $S\times T$-matrices. Let $C$ be a $T\times U$-matrix. Then,
$\left(  A+B\right)  C=AC+BC$.

\item[(c)] Let $S$, $T$, $U$ and $V$ be four sets such that $T$ and $U$ are
finite. Let $A$ be an $S\times T$-matrix. Let $B$ be a $T\times U$-matrix. Let
$C$ be a $U\times V$-matrix. Then, $\left(  AB\right)  C=A\left(  BC\right)  $.
\end{enumerate}
\end{proposition}

The proof of Proposition \ref{prop.STmat.ass} (and of similar properties that
will be left unstated) is analogous to the proofs of the corresponding
properties of $n\times m$-matrices.\footnote{A little \textbf{warning}: In
Proposition \ref{prop.STmat.ass}(c), the condition that $T$ and $U$ be finite
can be loosened (we leave this to the interested reader), but cannot be
completely disposed of. It can happen that both $\left(  AB\right)  C$ and
$A\left(  BC\right)  $ are defined, but $\left(  AB\right)  C=A\left(
BC\right)  $ does not hold (if we remove this condition). For example, this
happens if $S=\ZZ$, $T=\ZZ$, $U=\ZZ$, $V=\ZZ$,
$A=\left(
\begin{cases}
1, & \text{if }i\geq j;\\
0, & \text{if }i<j
\end{cases}
\right)  _{\left(  i,j\right)  \in \ZZ \times \ZZ}$, $B=\left(
\delta_{i,j}-\delta_{i,j+1}\right)  _{\left(  i,j\right)  \in \ZZ%
\times \ZZ}$ and $C=\left(
\begin{cases}
0, & \text{if }i\geq j;\\
1, & \text{if }i<j
\end{cases}
\right)  _{\left(  i,j\right)  \in \ZZ \times \ZZ}$. (Indeed, in
this example, it is easy to check that $AB=I_{\ZZ}$ and
$BC=-I_{\ZZ}$ and thus $\underbrace{\left(  AB\right)  }%
_{=I_{\ZZ}}C=I_{\ZZ}C=C\neq-A=A\underbrace{\left(
-I_{\ZZ}\right)  }_{=BC}=A\left(  BC\right)  $.)
\par
This seeming paradox is due to the subtleties of rearranging infinite sums
(similarly to how a conditionally convergent series of real numbers can change
its value when its entries are rearranged).} As a consequence of these
properties, it is easy to see that if $S$ is any finite set, then
$\kk^{S\times S}$ is a $\kk$-algebra.

In general, $S\times T$-matrices (unlike $n\times m$-matrices) do not have a
predefined order on their rows and their columns. Thus, the classical notion
of a triangular $n\times n$-matrix cannot be generalized to a notion of a
``triangular $S\times S$-matrix'' when $S$ is
just a set with no additional structure. However, when $S$ is a poset, such a
generalization can be made:

\begin{definition}
\label{def.STmat.tria}Let $S$ be a poset. Let $A=\left(  a_{s,t}\right)
_{\left(  s,t\right)  \in S\times S}$ be an $S\times S$-matrix.

\begin{enumerate}
\item[(a)] The matrix $A$ is said to be
\emph{triangular}\index{triangular $S\times S$-matrix}
if and only if
every $\left(  s,t\right)  \in S\times S$ which does not satisfy $t\leq s$
must satisfy $a_{s,t}=0$. (Here, $\leq$ denotes the smaller-or-equal relation
of the poset $S$.)

\item[(b)] The matrix $A$ is said to be
\emph{unitriangular}\index{unitriangular $S\times S$-matrix}
if and only if
$A$ is triangular and has the further property that, for every $s\in S$, we
have $a_{s,s}=1$.

\item[(c)] The matrix $A$ is said to be
\emph{invertibly triangular}\index{invertibly triangular $S\times S$-matrix}
if and
only if $A$ is triangular and has the further property that, for every $s\in
S$, the element $a_{s,s}$ of $\kk$ is invertible.
\end{enumerate}

Of course, all three notions of ``triangular'', ``unitriangular'' and
``invertibly triangular'' depend on the partial order on $S$.

Clearly, every invertibly triangular $S\times S$-matrix is triangular. Also,
every unitriangular $S\times S$-matrix is invertibly triangular (because the
element $1$ of $\kk$ is invertible).

We can restate the definition of ``invertibly
triangular'' as follows: The matrix $A$ is said to be
\emph{invertibly triangular} if and only if it is triangular and its
diagonal entries are invertible. Similarly, we can restate the definition of
``unitriangular'' as follows: The matrix $A$
is said to be \emph{unitriangular} if and only if it is triangular and all
its diagonal entries equal $1$.
\end{definition}

Definition \ref{def.STmat.tria}(a) generalizes both the notion of
upper-triangular matrices and the notion of lower-triangular matrices. To wit:

\begin{example}
Let $n\in \NN$. Let $N_{1}$ be the poset whose ground set is $\left\{
1,2,\ldots,n\right\}  $ and whose smaller-or-equal relation $\leq_{1}$ is
given by%
\[
s\leq_{1}t\ \Longleftrightarrow\ s\leq t\text{ (as integers).}%
\]
(This is the usual order relation on this set.) Let $N_{2}$ be the poset whose
ground set is $\left\{  1,2,\ldots,n\right\}  $ and whose order relation
$\leq_{2}$ is given by%
\[
s\leq_{2}t\ \Longleftrightarrow\ s\geq t\text{ (as integers).}%
\]
Let $A\in\kk^{n\times n}$.

\begin{enumerate}
\item[(a)] The matrix $A$ is upper-triangular if and only if $A$ is triangular
when regarded as an $N_{1}\times N_{1}$-matrix.

\item[(b)] The matrix $A$ is lower-triangular if and only if $A$ is triangular
when regarded as an $N_{2}\times N_{2}$-matrix.
\end{enumerate}
\end{example}

More interesting examples of triangular matrices are obtained when the order
on $S$ is not a total order:

\begin{example}
Let $S$ be the poset whose ground set is $\left\{  1,2,3\right\}  $ and whose
smaller relation $<_{S}$ is given by $1<_{S}2$ and $3<_{S}2$. Then, the
triangular $S\times S$-matrices are precisely the $3\times3$-matrices of the
form $\left(
\begin{array}
[c]{ccc}%
a_{1,1} & 0 & 0\\
a_{2,1} & a_{2,2} & a_{2,3}\\
0 & 0 & a_{3,3}%
\end{array}
\right)  $ with $a_{1,1},a_{2,1},a_{2,2},a_{2,3},a_{3,3}\in\kk$.
\end{example}

We shall now state some basic properties of triangular matrices:

\begin{proposition}
\label{prop.STmat.triangular}Let $S$ be a finite poset.

\begin{enumerate}
\item[(a)] The triangular $S\times S$-matrices form a subalgebra of the
$\kk$-algebra $\kk^{S\times S}$.

\item[(b)] The invertibly triangular $S\times S$-matrices form a group with
respect to multiplication.

\item[(c)] The unitriangular $S\times S$-matrices form a group with respect to multiplication.

\item[(d)] Any invertibly triangular $S\times S$-matrix is invertible, and its
inverse is again invertibly triangular.

\item[(e)] Any unitriangular $S\times S$-matrix is invertible, and its inverse
is again unitriangular.
\end{enumerate}
\end{proposition}

\begin{exercise}
\label{exe.STmat.mat}Prove Proposition~\ref{prop.STmat.triangular}.
\end{exercise}

\subsubsection{Expansion of a family in another}

We will often study situations where two families $\left(  e_{s}\right)
_{s\in S}$ and $\left(  f_{t}\right)  _{t\in T}$ of vectors in a $\kk%
$-module $M$ are given, and the vectors $e_{s}$ can be written as linear
combinations of the vectors $f_{t}$. In such situations, we can form an
$S\times T$-matrix out of the coefficients of these linear combinations; this
is one of the ways how matrices arise in the theory of modules. Let us define
the notations we are going to use in such situations:

\begin{definition}
\label{def.STmat.expansion}Let $M$ be a $\kk$-module. Let $\left(
e_{s}\right)  _{s\in S}$ and $\left(  f_{t}\right)  _{t\in T}$ be two families
of elements of $M$. (The sets $S$ and $T$ may and may not be finite.)

Let $A=\left(  a_{s,t}\right)  _{\left(  s,t\right)  \in S\times T}$ be an
$S\times T$-matrix. Assume that, for every $s\in S$, all but finitely many
$t\in T$ satisfy $a_{s,t}=0$. (This assumption is automatically satisfied if
$T$ is finite.)

We say that the family $\left(  e_{s}\right)  _{s\in S}$ \dfn{expands in
the family $\left(  f_{t}\right)  _{t\in T}$ through the matrix $A$}
if
\begin{equation}
\text{every }s\in S\text{ satisfies }e_{s}=\sum_{t\in T}a_{s,t}f_{t}.
\label{eq.def.STmat.expansion.es=}%
\end{equation}
In this case, we furthermore say that the matrix $A$ is a
\dfn{change-of-basis matrix} (or \dfn{transition matrix}) from the
family $\left(  e_{s}\right)  _{s\in S}$ to the family $\left(  f_{t}\right)
_{t\in T}$.
\end{definition}

\begin{remark}
The notation in Definition \ref{def.STmat.expansion} is not really standard;
even we ourselves will occasionally deviate in its use. In the formulation
``the family $\left(  e_{s}\right)  _{s\in S}$ expands in the
family $\left(  f_{t}\right)  _{t\in T}$ through the matrix $A$'', the word
``in'' can be replaced by ``with respect to'', and the word
``through'' can be replaced by ``using''.

The notion of a ``change-of-basis matrix'' is
slightly misleading, because neither of the families $\left(  e_{s}\right)
_{s\in S}$ and $\left(  f_{t}\right)  _{t\in T}$ has to be a basis. Our use of
the words ``transition matrix'' should not be
confused with the different meaning that these words have in the theory of
Markov chains. The indefinite article in ``a change-of-basis
matrix'' is due to the fact that, for given families $\left(
e_{s}\right)  _{s\in S}$ and $\left(  f_{t}\right)  _{t\in T}$, there might be
more than one change-of-basis matrix from $\left(  e_{s}\right)  _{s\in S}$ to
$\left(  f_{t}\right)  _{t\in T}$. (There also might be no such matrix.) When
$\left(  e_{s}\right)  _{s\in S}$ and $\left(  f_{t}\right)  _{t\in T}$ are
bases of the $\kk$-module $M$, there exists precisely one
change-of-basis matrix from $\left(  e_{s}\right)  _{s\in S}$ to $\left(
f_{t}\right)  _{t\in T}$.
\end{remark}

So a change-of-basis matrix $A=\left(  a_{s,t}\right)  _{\left(  s,t\right)
\in S\times T}$ from one family $\left(  e_{s}\right)  _{s\in S}$ to another
family $\left(  f_{t}\right)  _{t\in T}$ allows us to write the elements of
the former family as linear combinations of the elements of the latter (using
(\ref{eq.def.STmat.expansion.es=})). When such a matrix $A$ is invertible (and
the sets $S$ and $T$ are finite\footnote{We are requiring the finiteness of
$S$ and $T$ mainly for the sake of simplicity. We could allow $S$ and $T$ to
be infinite, but then we would have to make some finiteness requirements on
$A$ and $A^{-1}$.}), it also (indirectly) allows us to do the opposite: i.e.,
to write the elements of the latter family as linear combinations of the
elements of the former. This is because if $A$ is an invertible
change-of-basis matrix from $\left(  e_{s}\right)  _{s\in S}$ to $\left(
f_{t}\right)  _{t\in T}$, then $A^{-1}$ is a change-of-basis matrix from
$\left(  f_{t}\right)  _{t\in T}$ to $\left(  e_{s}\right)  _{s\in S}$. This
is part (a) of the following theorem:

\begin{theorem}
\label{thm.STmat.expansion-inv}Let $M$ be a $\kk$-module. Let $S$ and
$T$ be two finite sets. Let $\left(  e_{s}\right)  _{s\in S}$ and $\left(
f_{t}\right)  _{t\in T}$ be two families of elements of $M$.

Let $A$ be an invertible $S\times T$-matrix. Thus, $A^{-1}$ is a $T\times S$-matrix.

Assume that the family $\left(  e_{s}\right)  _{s\in S}$ expands in the family
$\left(  f_{t}\right)  _{t\in T}$ through the matrix $A$. Then:

\begin{enumerate}
\item[(a)] The family $\left(  f_{t}\right)  _{t\in T}$ expands in the family
$\left(  e_{s}\right)  _{s\in S}$ through the matrix $A^{-1}$.

\item[(b)] The $\kk$-submodule of $M$ spanned by the family $\left(
e_{s}\right)  _{s\in S}$ is the $\kk$-submodule of $M$ spanned by the
family $\left(  f_{t}\right)  _{t\in T}$.

\item[(c)] The family $\left(  e_{s}\right)  _{s\in S}$ spans the $\kk%
$-module $M$ if and only if the family $\left(  f_{t}\right)  _{t\in T}$ spans
the $\kk$-module $M$.

\item[(d)] The family $\left(  e_{s}\right)  _{s\in S}$ is $\kk%
$-linearly independent if and only if the family $\left(  f_{t}\right)  _{t\in
T}$ is $\kk$-linearly independent.

\item[(e)] The family $\left(  e_{s}\right)  _{s\in S}$ is a basis of the
$\kk$-module $M$ if and only if the family $\left(  f_{t}\right)
_{t\in T}$ is a basis of the $\kk$-module $M$.
\end{enumerate}
\end{theorem}

\begin{exercise}
\label{exe.STmat.expansion-inv}Prove Theorem~\ref{thm.STmat.expansion-inv}.
\end{exercise}

\begin{definition}
\label{def.STmat.expansion-tria}Let $M$ be a $\kk$-module. Let $S$ be a
finite poset. Let $\left(  e_{s}\right)  _{s\in S}$ and $\left(  f_{s}\right)
_{s\in S}$ be two families of elements of $M$.

\begin{enumerate}
\item[(a)] We say that the family $\left(  e_{s}\right)  _{s\in S}$
\dfn{expands triangularly in the family $\left(  f_{s}\right)  _{s\in S}$}
if and only if there exists a triangular $S\times S$-matrix $A$ such that the
family $\left(  e_{s}\right)  _{s\in S}$ expands in the family $\left(
f_{s}\right)  _{s\in S}$ through the matrix $A$.

\item[(b)] We say that the family $\left(  e_{s}\right)  _{s\in S}$
\dfn{expands invertibly triangularly in the family $\left(  f_{s}\right)
_{s\in S}$} if and only if there exists an invertibly triangular $S\times
S$-matrix $A$ such that the family $\left(  e_{s}\right)  _{s\in S}$ expands
in the family $\left(  f_{s}\right)  _{s\in S}$ through the matrix $A$.

\item[(c)] We say that the family $\left(  e_{s}\right)  _{s\in S}$
\dfn{expands unitriangularly in the family $\left(  f_{s}\right)  _{s\in
S}$} if and only if there exists a unitriangular $S\times S$-matrix $A$
such that the family $\left(  e_{s}\right)  _{s\in S}$ expands in the
family $\left(  f_{s}\right)  _{s\in S}$ through the matrix $A$.
\end{enumerate}

Clearly, if the family $\left(  e_{s}\right)  _{s\in S}$ expands
unitriangularly in the family $\left(  f_{s}\right)  _{s\in S}$, then it also
expands invertibly triangularly in the family $\left(  f_{s}\right)  _{s\in
S}$ (because any unitriangular matrix is an invertibly triangular matrix).
\end{definition}

We notice that in Definition \ref{def.STmat.expansion-tria}, the two families
$\left(  e_{s}\right)  _{s\in S}$ and $\left(  f_{s}\right)  _{s\in S}$ must
be indexed by one and the same set $S$.

The concepts of ``expanding triangularly'', ``expanding invertibly
triangularly'' and ``expanding unitriangularly'' can also be characterized
without referring to matrices, as follows:

\begin{remark}
\label{rmk.STmat.expansion-tria-shortcuts}
Let $M$ be a $\kk$-module.
Let $S$ be a finite poset. Let $\left( e_s \right)_{s \in S}$ and
$\left( f_s \right)_{s \in S}$ be two families of elements of $M$.
Let $<$ denote the smaller relation of the poset $S$, and let
$\leq$ denote the smaller-or-equal relation of the poset $S$.
Then:

\begin{enumerate}

\item[(a)] The family $\left( e_s \right)_{s \in S}$ expands triangularly
in the family $\left( f_s \right)_{s \in S}$ if and only if every $s\in S$
satisfies
\[
e_{s}=\left(  \text{a }\kk\text{-linear combination of the elements
}f_{t}\text{ for }t\in S\text{ satisfying }t\leq s\right)  .
\]

\item[(b)] The family $\left( e_s \right)_{s \in S}$ expands invertibly
triangularly in the family $\left( f_s \right)_{s \in S}$ if and only if
every $s\in S$ satisfies
\[
e_{s}=\alpha_{s}f_{s}+\left(  \text{a }\kk\text{-linear combination of
the elements }f_{t}\text{ for }t\in S\text{ satisfying }t<s\right)
\]
for some invertible $\alpha_{s}\in\kk$.

\item[(c)] The family $\left( e_s \right)_{s \in S}$ expands
unitriangularly in the family $\left( f_s \right)_{s \in S}$ if and only if
every $s\in S$ satisfies
\[
e_{s}=f_{s}+\left(  \text{a }\kk\text{-linear combination of the
elements }f_{t}\text{ for }t\in S\text{ satisfying }t<s\right)  .
\]

\end{enumerate}
\end{remark}

All three parts of Remark~\ref{rmk.STmat.expansion-tria-shortcuts} follow
easily from the definitions.

\begin{example}
Let $n\in \NN$. For this example, let $S$ be the poset $\left\{
1,2,\ldots,n\right\}  $ (with its usual order). Let $M$ be a $\kk%
$-module, and let $\left(  e_{s}\right)  _{s\in S}$ and $\left(  f_{s}\right)
_{s\in S}$ be two families of elements of $M$. We shall identify these
families $\left(  e_{s}\right)  _{s\in S}$ and $\left(  f_{s}\right)  _{s\in
S}$ with the $n$-tuples $\left(  e_{1},e_{2},\ldots,e_{n}\right)  $ and
$\left(  f_{1},f_{2},\ldots,f_{n}\right)  $. Then, the family $\left(
e_{s}\right)  _{s\in S}=\left(  e_{1},e_{2},\ldots,e_{n}\right)  $ expands
triangularly in the family $\left(  f_{s}\right)  _{s\in S}=\left(
f_{1},f_{2},\ldots,f_{n}\right)  $ if and only if, for every $s\in\left\{
1,2,\ldots,n\right\}  $, the vector $e_{s}$ is a $\kk$-linear
combination of $f_{1},f_{2},\ldots,f_{s}$. Moreover, the family $\left(
e_{s}\right)  _{s\in S}=\left(  e_{1},e_{2},\ldots,e_{n}\right)  $ expands
unitriangularly in the family $\left(  f_{s}\right)  _{s\in S}=\left(
f_{1},f_{2},\ldots,f_{n}\right)  $ if and only if, for every $s\in\left\{
1,2,\ldots,n\right\}  $, the vector $e_{s}$ is a sum of $f_{s}$ with a
$\kk$-linear combination of $f_{1},f_{2},\ldots,f_{s-1}$.
\end{example}

\begin{corollary}
\label{cor.STmat.expansion-tria-inv}Let $M$ be a $\kk$-module. Let $S$
be a finite poset. Let $\left(  e_{s}\right)  _{s\in S}$ and $\left(
f_{s}\right)  _{s\in S}$ be two families of elements of $M$. Assume that the
family $\left(  e_{s}\right)  _{s\in S}$ expands invertibly triangularly in
the family $\left(  f_{s}\right)  _{s\in S}$. Then:

\begin{enumerate}
\item[(a)] The family $\left(  f_{s}\right)  _{s\in S}$ expands invertibly
triangularly in the family $\left(  e_{s}\right)  _{s\in S}$.

\item[(b)] The $\kk$-submodule of $M$ spanned by the family $\left(
e_{s}\right)  _{s\in S}$ is the $\kk$-submodule of $M$ spanned by the
family $\left(  f_{s}\right)  _{s\in S}$.

\item[(c)] The family $\left(  e_{s}\right)  _{s\in S}$ spans the $\kk%
$-module $M$ if and only if the family $\left(  f_{s}\right)  _{s\in S}$ spans
the $\kk$-module $M$.

\item[(d)] The family $\left(  e_{s}\right)  _{s\in S}$ is $\kk%
$-linearly independent if and only if the family $\left(  f_{s}\right)  _{s\in
S}$ is $\kk$-linearly independent.

\item[(e)] The family $\left(  e_{s}\right)  _{s\in S}$ is a basis of the
$\kk$-module $M$ if and only if the family $\left(  f_{s}\right)
_{s\in S}$ is a basis of the $\kk$-module $M$.
\end{enumerate}
\end{corollary}

\begin{exercise}
\label{exe.STmat.expansion-tria-inv}Prove
Remark~\ref{rmk.STmat.expansion-tria-shortcuts} and
Corollary~\ref{cor.STmat.expansion-tria-inv}.
\end{exercise}

An analogue of Corollary \ref{cor.STmat.expansion-tria-inv} can be stated for
unitriangular expansions, but we leave this to the reader.

% [DG][v55] Added the above appendix, and referenced it in the text
% accordingly.

\newpage

\section{Further hints to the exercises (work in progress)}

The following pages contain hints to (some of\footnote{Currently only
the ones from Chapter~\ref{Hopf-intro-section}.}) the exercises in the text
(beyond the hints occasionally included in the exercises themselves). Some of
the hints rise to the level of outlined solutions.

\newcommand\hintsboilerplate{Detailed solutions (sometimes different
from the solutions hinted at) can be found in Chapter \ref{chp.solutions}.}

\begin{nosolutions}
\renewcommand\hintsboilerplate{Note that there is also a version of
this text that contains detailed
solutions for all the exercises; this version can be downloaded from
\url{http://www.cip.ifi.lmu.de/~grinberg/algebra/HopfComb-sols.pdf} (or
compiled from the sourcecode of the text).}
\end{nosolutions}

\hintsboilerplate

\textbf{Warning:} The hints below are new and have never been proofread. Typos
(or worse) are likely. In case of doubt, consult the detailed solutions.

\subsection{Hints for Chapter \ref{Hopf-intro-section}}

\textit{Hint to Exercise \ref{exe.counit.unique}.} The claim of the exercise
is dual to the classical fact that if $A$ is a $\kk$-module and
$m:A\otimes A\rightarrow A$ is a $\kk$-linear map, then there exists
\textit{at most one} $\kk$-linear map $u:\mathbf{k}\rightarrow A$ such
that the diagram \eqref{unit-diagram} commutes\footnote{This fact is just the
linearization of the known fact that any binary operation has at most one
neutral element.}. Take any proof of this latter fact, rewrite it in an
``element-free'' fashion\footnote{This means
rewriting it completely in terms of linear maps rather than elements. For
example, instead of talking about $m\left(  m\left(  a\otimes b\right)
\otimes c\right)  $ for three elements $a,b,c\in A$, you should talk about the
map $m\circ\left(  m\otimes \id_{A}\right)  :A\otimes
A\otimes A\rightarrow A$ (which is, of course, the map that sends each
$a\otimes b\otimes c$ to $m\left(  m\left(  a\otimes b\right)  \otimes
c\right)  $). Instead of computing with elements, you should compute with maps
(and commutative diagrams).}, and ``reverse all
arrows''. This will yield a solution to Exercise
\ref{exe.counit.unique}.

For an alternative solution, use Sweedler notation (as in
\eqref{eq.sweedler-not.def1}) as follows: The commutativity of the diagram
\eqref{counit-diagram} says that%
\[
c=\sum_{\left(  c\right)  }\epsilon\left(  c_{1}\right)  c_{2}=\sum_{\left(
c\right)  }\epsilon\left(  c_{2}\right)  c_{1}\qquad\text{for each }c\in C.
\]
Thus, if $\epsilon_{1}$ and $\epsilon_{2}$ are two $\kk$-linear maps
$\epsilon:C\rightarrow\mathbf{k}$ such that the diagram \eqref{counit-diagram}
commutes, then each $c\in C$ satisfies%
\[
c=\sum_{\left(  c\right)  }\epsilon_{1}\left(  c_{1}\right)  c_{2}%
=\sum_{\left(  c\right)  }\epsilon_{1}\left(  c_{2}\right)  c_{1}%
\]
and%
\[
c=\sum_{\left(  c\right)  }\epsilon_{2}\left(  c_{1}\right)  c_{2}%
=\sum_{\left(  c\right)  }\epsilon_{2}\left(  c_{2}\right)  c_{1}.
\]
Apply $\epsilon_{2}$ to both sides of the equality $c=\sum_{\left(  c\right)
}\epsilon_{1}\left(  c_{2}\right)  c_{1}$, and apply $\epsilon_{1}$ to both
sides of the equality $c=\sum_{\left(  c\right)  }\epsilon_{2}\left(
c_{1}\right)  c_{2}$. Compare the results, and conclude that $\epsilon
_{1}=\epsilon_{2}$.

\bigskip

\textit{Hint to Exercise \ref{exe.tensor-prod-algs}.} Part (a) is well-known,
and part (b) is dual to part (a). So the trick is (again) to rewrite the
classical proof of part (a) in an ``element-free'' way, and then
``reversing all arrows''. Alternatively, part (b) can be solved using
Sweedler notation.

\bigskip

\textit{Hint to Exercise \ref{exe.tensorprod.morphism}.} Same method as for
Exercise \ref{exe.tensor-prod-algs} above.

\bigskip

\textit{Hint to Exercise \ref{exe.kernel-is-coideal}.} (a) Use the following
fact from linear algebra: If $U$, $V$, $U^{\prime}$ and $V^{\prime}$ are four
$\kk$-modules, and $\phi:U\rightarrow U^{\prime}$ and $\psi
:V\rightarrow V^{\prime}$ are two surjective $\kk$-linear maps, then
the kernel of $\phi\otimes\psi:U\otimes V\rightarrow U^{\prime}\otimes
V^{\prime}$ is
\[
\ker\left(  \phi\otimes\psi\right)  =\left(  \ker\phi\right)  \otimes
V+U\otimes\left(  \ker\psi\right)  .
\]

(b) The fact just mentioned also holds if we no longer require $\phi$ and
$\psi$ to be surjective, but instead require $\kk$ to be a field.

\bigskip

\textit{Hint to Exercise \ref{exe.graded.iso}.} Let $f : V \to W$
be an invertible graded $\kk$-linear map. Let $n \in \NN$
and $w \in W_n$. Show that the $n$-th homogeneous component of
$f^{-1}\left(w\right)$ is also a preimage of $w$ under $f$, and
thus must equal $f^{-1}\left(w\right)$. Therefore,
$f^{-1}\left(w\right) \in W_n$.

\bigskip

\textit{Hint to Exercise \ref{exe.primitives.graded-coideal}.} (a) Define the
$\kk$-linear map $\widetilde{\Delta}:A\rightarrow A\otimes A$ by
$\widetilde{\Delta}\left(  x\right)  =\Delta\left(  x\right)  -\left(
x\otimes1+1\otimes x\right)  $. Argue that $\widetilde{\Delta}$ is graded, so
its kernel $\ker\widetilde{\Delta}$ is a graded $\kk$-submodule of $A$.
But this kernel is precisely $\mathfrak{p}$.

(b) The hard part is to show that $\epsilon\left(  \mathfrak{p}\right)  =0$.
To do so, consider any $x\in\mathfrak{p}$, and apply the map $\epsilon
\otimes \id$ to both sides of the equality $\Delta\left(
x\right)  =x\otimes1+1\otimes x$. The result simplifies to $x=\epsilon\left(
x\right)  \cdot1_{A}+x$. Thus, $\epsilon\left(  x\right)  \cdot1_{A}=0$. Now
apply $\epsilon$ to this, thus obtaining $\epsilon\left(  x\right)  =0$.

\bigskip

\textit{Hint to Exercise \ref{graded-connected-exercise}.} (a) This follows
from $1_{A}\in A_{0}$, which is part of what it means for $A$ to be a graded
$\kk$-algebra.

(b) Let $\epsilon^{\prime}:A_{0}\rightarrow\mathbf{k}$ be the restriction of
the map $\epsilon$ to $A_{0}$. We know that $\epsilon^{\prime}$ is surjective
(since $\epsilon^{\prime}\left(  1_{A}\right)  =1_{\mathbf{k}}$), and that
both $A_{0}$ and $\kk$ are free $\kk$-modules of rank $1$ (since
connectedness of $A$ means $A_{0}\cong\mathbf{k}$ as $\kk$-modules). It
is an an easy exercise in linear algebra to conclude from these facts that
$\epsilon^{\prime}$ is an isomorphism. Since $\epsilon^{\prime}\circ
u= \id_{\kk}$, we thus conclude that
$u:\mathbf{k}\rightarrow A_{0}$ is an isomorphism as well (from $\kk$
to $A_{0}$).

(c) This follows from part (b).

(e) This follows from how we solved part (b).

(d) Since the bialgebra $A$ is graded, the map $\epsilon$ must be graded.
Thus, for each positive integer $n$, we have $\epsilon\left(  A_{n}\right)
\subset\mathbf{k}_{n}=0$. This quickly yields $\epsilon\left(  I\right)  =0$
(where $I=\bigoplus_{n>0}A_{n}$), hence $I\subset\ker\epsilon$. On the other
hand, $\ker\epsilon\subset I$ can be shown as follows: Let $a\in\ker\epsilon$;
write $a$ in the form $a=a^{\prime}+a^{\prime\prime}$ for some $a^{\prime}\in
A_{0}$ and some $a^{\prime\prime}\in I$, and then argue that $0=\epsilon
\left(  a\right)  =\epsilon\left(  a^{\prime}+a^{\prime\prime}\right)
=\epsilon\left(  a^{\prime}\right)  +\underbrace{\epsilon\left(
a^{\prime\prime}\right)  }_{\substack{=0\\\text{(since }a^{\prime\prime}\in
I\subset\ker\epsilon\text{)}}}=\epsilon\left(  a^{\prime}\right)  $, so that
$a^{\prime}=0$ by part (e) and therefore $a\in I$.

(f) This is most intuitive with Sweedler notation: Let $x\in A$. Then,
$\Delta\left(  x\right)  =\sum_{\left(  x\right)  }x_{1}\otimes x_{2}$.
Applying $\id \otimes\epsilon$ and recalling the commutativity
of \eqref{counit-diagram}, we thus get $x=\sum_{\left(  x\right)  }%
\epsilon\left(  x_{2}\right)  x_{1}$. Thus,%
\begin{align*}
\underbrace{\Delta\left(  x\right)  }_{=\sum_{\left(  x\right)  }x_{1}\otimes
x_{2}}-\underbrace{x}_{=\sum_{\left(  x\right)  }\epsilon\left(  x_{2}\right)
x_{1}}\otimes1  &  =\sum_{\left(  x\right)  }x_{1}\otimes x_{2}-\sum_{\left(
x\right)  }\epsilon\left(  x_{2}\right)  x_{1}\otimes1\\
&  =\sum_{\left(  x\right)  }\underbrace{x_{1}}_{\in A}\otimes
\underbrace{\left(  x_{2}-\epsilon\left(  x_{2}\right)  \cdot1\right)
}_{\substack{\in\ker\epsilon=I\\\text{(by part (d))}}}\in A\otimes I.
\end{align*}

(g) Let $x\in I$. Proceeding similarly to part (f), show that%
\[
\Delta\left(  x\right)  -1\otimes x-x\otimes1+\epsilon\left(  x\right)
1\otimes1=\sum_{\left(  x\right)  }\underbrace{\left(  x_{1}-\epsilon\left(
x_{1}\right)  \cdot1\right)  }_{\substack{\in\ker\epsilon=I\\\text{(by part
(d))}}}\otimes\underbrace{\left(  x_{2}-\epsilon\left(  x_{2}\right)
\cdot1\right)  }_{\substack{\in\ker\epsilon=I\\\text{(by part (d))}}}\in
I\otimes I.
\]
Since $x\in I=\ker\epsilon$, the $\epsilon\left(  x\right)  1\otimes1$ term on
the left hand side vanishes.

(h) This follows from part (g), since a simple homogeneity argument shows that
$\left(  I\otimes I\right)  _{n}=\sum_{k=1}^{n-1}A_{k}\otimes A_{n-k}$.

\bigskip

\textit{Hint to Exercise \ref{exe.graded.Dq}.} We need to check the four
equalities $D_{q} \circ m=m\circ\left(  D_{q} \otimes D_{q} \right)  $ and
$D_{q} \circ u=u$ and $\left(  D_{q} \otimes D_{q} \right)  \circ\Delta
=\Delta\circ D_{q} $ and $\epsilon\circ D_{q} =\epsilon$. This can easily be
done by hand (just check everything on homogeneous elements); a more erudite
proof proceeds as follows: Generalize the map $D_{q}$ to a map $D_{q,V}%
:V\rightarrow V$ defined (in the same way as $D_{q}$) for every graded
$\kk$-module $V$, and show that these maps $D_{q,V}$ are functorial
(i.e., if $f:V\rightarrow W$ is a graded $\kk$-linear map between two
graded $\kk$-modules $V$ and $W$, then $D_{q,W}\circ f=f\circ D_{q,V}$)
and ``respect tensor products'' (i.e., we
have $D_{q,V\otimes W}=D_{q,V}\otimes D_{q,W}$ for any two graded $\mathbf{k}%
$-modules $V$ and $W$). The four equalities are then easily obtained from
these two facts, without having to introduce elements.

\bigskip

\textit{Hint to Exercise \ref{exe.bialg.tensor}.} (a) Our definition of the
$\kk$-coalgebra $A\otimes B$ yields%
\[
\Delta_{A\otimes B}=\left(  \id_{A}\otimes T_{A,B}%
\otimes \id_{B}\right)  \circ\left(  \Delta_{A}\otimes\Delta
_{B}\right)  \ \ \ \ \ \ \ \ \ \ \text{and}\ \ \ \ \ \ \ \ \ \ \epsilon
_{A\otimes B}=\theta\circ\left(  \epsilon_{A}\otimes\epsilon_{B}\right)  ,
\]
where $\theta$ is the canonical $\kk$-module isomorphism $\mathbf{k}%
\otimes\mathbf{k}\rightarrow\mathbf{k}$. All maps on the right hand sides are
$\kk$-algebra homomorphisms (see Exercise \ref{exe.tensorprod.morphism}%
(a)); thus, so are $\Delta_{A\otimes B}$ and $\epsilon_{A\otimes B}$.

(b) Straightforward.

\bigskip

\textit{Hint to Exercise \ref{exe.convolution.assoc}.} Simple computation
(either element-free or with Sweedler notation).

\bigskip

\textit{Hint to Exercise \ref{exe.convolution.tensor}.} Simple computation
(either element-free or with Sweedler notation).

\bigskip

\textit{Hint to Exercise \ref{exe.convolution.tensor-curry}.} Straightforward
computation, best done using Sweedler notation.

\bigskip

\textit{Hint to Exercise \ref{exe.hopf.tensor}.} Use Exercise
\ref{exe.convolution.assoc}.

\bigskip

\textit{Hint to Exercise \ref{exe.iterated.m}.} The following is more context
than hint (see the last paragraph for an actual hint).

It is easiest to prove this by calculating with elements. To wit, in order to
prove that two $\kk$-linear maps from $A^{\otimes\left(  k+1\right)  }$
are identical, it suffices to show that they agree on all pure tensors
$a_{1}\otimes a_{2}\otimes\cdots\otimes a_{k+1}\in A^{\otimes\left(
k+1\right)  }$. But the recursive definition of $m^{\left(  k\right)  }$ shows
that%
\begin{equation}
m^{\left(  k\right)  }\left(  a_{1}\otimes a_{2}\otimes\cdots\otimes
a_{k+1}\right)  =a_{1}\left(  a_{2}\left(  a_{3}\left(  \cdots\left(
a_{k}a_{k+1}\right)  \cdots\right)  \right)  \right)
\label{hint.iterated.m.1}%
\end{equation}
for all $a_{1},a_{2},\ldots,a_{k+1}\in A$. Now, the ``general
associativity'' law (a fundamental result in abstract
algebra, commonly used without mention) says that, because the multiplication
of $A$ is associative, the parentheses in the product $a_{1}\left(
a_{2}\left(  a_{3}\left(  \cdots\left(  a_{k}a_{k+1}\right)  \cdots\right)
\right)  \right)  $ can be omitted without making it ambiguous -- i.e., any
two ways of parenthesizing the product $a_{1}a_{2}\cdots a_{k+1}$ evaluate to
the same result. (For example, for $k=4$, this says that
\[
a_{1}\left(  a_{2}\left(  a_{3}a_{4}\right)  \right)  =a_{1}\left(  \left(
a_{2}a_{3}\right)  a_{4}\right)  =\left(  a_{1}a_{2}\right)  \left(
a_{3}a_{4}\right)  =\left(  a_{1}\left(  a_{2}a_{3}\right)  \right)
a_{4}=\left(  \left(  a_{1}a_{2}\right)  a_{3}\right)  a_{4}%
\]
for all $a_{1},a_{2},a_{3},a_{4}\in A$.) Thus, we can rewrite
\eqref{hint.iterated.m.1} as
\[
m^{\left(  k\right)  }\left(  a_{1}\otimes a_{2}\otimes\cdots\otimes
a_{k+1}\right)  =a_{1}a_{2}\cdots a_{k+1}.
\]
Using this formula, all four parts of the exercise become trivial: For
example, part (a) simply says that
\[
a_{1}a_{2}\cdots a_{k+1}=\left(  a_{1}a_{2}\cdots a_{i+1}\right)  \left(
a_{i+2}a_{i+3}\cdots a_{k+1}\right)
\]
for all $a_{1},a_{2},\ldots,a_{k+1}\in A$, because we have
\[
\left(  m\circ\left(  m^{\left(  i\right)  }\otimes m^{\left(  k-1-i\right)
}\right)  \right)  \left(  a_{1}\otimes a_{2}\otimes\cdots\otimes
a_{k+1}\right)  =\left(  a_{1}a_{2}\cdots a_{i+1}\right)  \left(
a_{i+2}a_{i+3}\cdots a_{k+1}\right)  .
\]
Likewise, part (c) simply says that%
\[
a_{1}a_{2}\cdots a_{k+1}=a_{1}a_{2}\cdots a_{i}\left(  a_{i+1}a_{i+2}\right)
a_{i+3}a_{i+4}\cdots a_{k+1}%
\]
for all $a_{1},a_{2},\ldots,a_{k+1}\in A$. Parts (b) and (d) are particular
cases of parts (a) and (c), respectively.

Of course, in order for this to be a complete solution, you have to prove the
``general associativity'' law used above. It
turns out that doing so is not much easier than solving the exercise from
scratch (in fact, part (a) of the exercise is an equivalent form of the
``general associativity'' law). So we can
just as well start from scratch and solve part (a) directly by induction on
$k$, then derive part (b) as its particular case, then solve part (c) by
induction on $k$ using the result of part (b), then derive part (d) as a
particular case of (c).

\bigskip

\textit{Hint to Exercise~\ref{exe.iterated.Delta}.} If you have solved
Exercise \ref{exe.iterated.m} in an ``element-free'' way, then you can
reverse all arrows in said
solution and thus obtain a solution to Exercise~\ref{exe.iterated.Delta}.

\bigskip

\textit{Hint to Exercise~\ref{exe.iterated.patras}.} (a) Induction on $k$,
using Exercise \ref{exe.tensorprod.morphism}(b).

(b) This is dual to (a).

(d) For every $\kk$-coalgebra $C$, consider the map $\Delta
_{C}^{\left(  k\right)  }:C\rightarrow C^{\otimes\left(  k+1\right)  }$ (this
is the map $\Delta^{\left(  k\right)  }$ defined in
Exercise~\ref{exe.iterated.Delta}). This map $\Delta_{C}^{\left(  k\right)  }$
is clearly functorial in $C$. By this we mean that if $C$ and $D$ are any two
$\kk$-coalgebras, and $f:C\rightarrow D$ is any $\kk$-coalgebra
homomorphism, then the diagram
\[
\xymatrixcolsep{4pc}
\xymatrix{
C \ar[r]^{f} \ar[d]^{\Delta^{\left(k\right)}_C} & D \ar[d]^{\Delta
^{\left(k\right)}_D} \\
C^{\otimes\left(k+1\right)} \ar[r]_{f^{\otimes\left(k+1\right)}}
& D^{\otimes\left(k+1\right)}
}
\]
commutes. Now, apply this to $C=H^{\otimes\left(  \ell+1\right)  }$, $D=H$ and
$f=m_{H}^{\left(  \ell\right)  }$ (using part (a)).

(c) This is dual to (d).

\bigskip

\textit{Hint to Exercise~\ref{exe.convolution.k}.} Induction on $k$.

\bigskip

\textit{Hint to Exercise~\ref{exe.coalg-anti}.} This is dual to
Proposition~\ref{antipodes-are-antiendomorphisms}, so the usual strategy
(viz., rewriting element-free and reversing all arrows) applies.

\bigskip

\textit{Hint to Exercise~\ref{exe.coalg-anti.gen}.} (a) A straightforward
generalization of the proof of
Proposition~\ref{antipodes-are-antiendomorphisms} (which corresponds to the
particular case when $C=A$ and $r = \id$) does the trick.

(b) This is dual to (a).

(c) Easy.

(d) Apply Exercise~\ref{exe.coalg-anti.gen}(a) to $C=A$ and
$r = \id_{A}$; then, apply
Proposition~\ref{antipodes-give-convolution-inverses}(a) to $H=A$ and
$\alpha=S$.

(e) Let $s:C\rightarrow A$ be the $\kk$-linear map that sends every
homogeneous element $c\in C_{n}$ (for every $n\in \NN$) to the $n$-th
homogeneous component of $r^{\star\left(  -1\right)  }\left(  c\right)  $.
Then, $s$ is graded, and (this takes some work) is also a $\star$-inverse to
$r$. But $r$ has only one $\star$-inverse.

\bigskip

\textit{Hint to Exercise~\ref{exe.antipode.inverse}.} (a) Rewrite the
assumption as $m\circ\left(  P\otimes \id \right)  \circ
T\circ\Delta=u\circ\epsilon$, where $T$ is the twist map $T_{A,A}$.
Proposition~\ref{antipodes-are-antiendomorphisms} leads to $m\circ\left(
S\otimes S\right)  =S\circ m\circ T$ and $u=S\circ u$.
Exercise~\ref{exe.coalg-anti} leads to $\left(  S\otimes S\right)  \circ
\Delta=T\circ\Delta\circ S$ and $\epsilon\circ S=\epsilon$. Use these to show
that $\left(  P\circ S\right)  \star S=u\circ\epsilon$, so that $P\circ
S= \id$. Also, show that $S\star\left(  S\circ P\right)
=u\circ\epsilon$, so that $S\circ P= \id$.

(b) Similar to (a).

(c) Let $A$ be a connected graded Hopf algebra. Just as a left $\star$-inverse
$S$ to $\id_{A}$ has been constructed in the proof of
Proposition~\ref{graded-connected-bialgebras-have-antipodes}, we could
construct a $\kk$-linear map $P:A\rightarrow A$ such that every $a\in
A$ satisfies $\sum_{\left(  a\right)  }P\left(  a_{2}\right)  \cdot
a_{1}=u\left(  \epsilon\left(  a\right)  \right)  $. Now apply part (a).

\bigskip

\textit{Hint to Exercise~\ref{exe.subcoalgebra.summand}.} Since $D$ is a
direct summand of $C$, we can identify the tensor products $D\otimes C$,
$C\otimes D$ and $D\otimes D$ with their canonical images inside $C\otimes C$.
Now, we can show that $\Delta\left(  D\right)  \subset D\otimes D$ as follows:
Let $p:C\rightarrow D$ be the canonical projection from $C$ onto its direct
summand $D$; then, $\Delta\left(  D\right)  \subset D\otimes C$ shows that
$\left(  p\otimes \id \right)  \circ\Delta=\Delta$, and
$\Delta\left(  D\right)  \subset C\otimes D$ shows that $\left(
\id \otimes p\right)  \circ\Delta=\Delta$. Hence,
\[
\underbrace{\left(  p\otimes p\right)  }_{=\left(  p\otimes\id%
\right)  \circ\left(  \id \otimes p\right)  }\circ\Delta=\left(
p\otimes \id \right)  \circ\underbrace{\left(  \id \otimes p
\right)  \circ\Delta}_{=\Delta}=\left(  p \otimes \id
\right)  \circ\Delta=\Delta.
\]
This yields $\Delta\left(  D\right)  \subset D\otimes D$. Hence, we get a map
$\Delta_{D}:D\rightarrow D\otimes D$ by restricting $\Delta$. Obviously, the
map $\epsilon:C\rightarrow\mathbf{k}$ restricts to a map $\epsilon
_{D}:D\rightarrow\mathbf{k}$ as well. It remains to check the commutativity of
the diagrams \eqref{coassociativity-diagram} and
\eqref{counit-diagram} for $D$ instead of $C$; but this is inherited
from $C$.

\bigskip

\textit{Hint to Exercise~\ref{exe.subcoalgebra.ABS}.} (a) Let $\widetilde{f}%
=\left(  \id_{C}\otimes f\otimes \id_{C}\right)
\circ\Delta^{\left(  2\right)  }:C\rightarrow C\otimes U\otimes C$; then,
$K=\ker\widetilde{f}$. Show (by manipulation of maps, using
Exercise~\ref{exe.iterated.Delta}(b)) that $\left(  \id_{C}
\otimes \id_{U}\otimes\Delta\right)  \circ\widetilde{f}%
=\left(  \widetilde{f}\otimes \id_{C}\right)  \circ\Delta$. Now,%
\[
K=\ker\widetilde{f}\subset\ker\left(  \underbrace{\left(  \id_{C}
\otimes \id_{U}\otimes\Delta\right)  \circ\widetilde{f}%
}_{=\left(  \widetilde{f}\otimes \id_{C}\right)  \circ\Delta
}\right)  =\ker\left(  \left(  \widetilde{f}\otimes \id_{C}
\right)  \circ\Delta\right)  =\Delta^{-1}\left(  \ker\left(
\widetilde{f}\otimes \id_{C}\right)  \right)
\]
and therefore%
\begin{align*}
\Delta\left(  K\right)   &  \subset\ker\left(  \widetilde{f}\otimes
\id_{C}\right)  =\underbrace{\left(  \ker\widetilde{f}\right)
}_{=K}\otimes C\ \ \ \ \ \ \ \ \ \ \left(  \text{since tensoring over a field
is left-exact}\right) \\
&  =K\otimes C.
\end{align*}
Similarly, $\Delta\left(  K\right)  \subset C\otimes K$. Now, apply
Exercise~\ref{exe.subcoalgebra.summand} to $D=K$.

(b) Let $E$ be a $\kk$-subcoalgebra of $C$ which is a subset of $\ker
f$. Then, $\Delta^{\left(  2\right)  }\left(  E\right)  \subset E\otimes
E\otimes E$ (since $E$ is a subcoalgebra) and $f\left(  E\right)  =0$ (since
$E\subset\ker f$). Now,
\begin{align*}
\left(  \left(  \id_{C}\otimes f\otimes \id_{C}\right)  \circ
\Delta^{\left(  2\right)  }\right)  \left(  E\right)   &
=\left(  \id_{C}\otimes f\otimes \id_{C}\right)
\left(  \underbrace{\Delta^{\left(  2\right)  }\left(  E\right)  }_{\subset
E\otimes E\otimes E}\right) \\
&  \subset\left(  \id_{C}\otimes f\otimes \id_{C}\right)
\left(  E\otimes E\otimes E\right) \\
&  = \id_{C}\left(  E\right)  \otimes\underbrace{f\left(
E\right)  }_{=0}\otimes \id_{C}\left(  E\right)  =0.
\end{align*}
Hence, $E\subset\ker\left(  \left(  \id_{C}\otimes
f\otimes \id_{C}\right)  \circ\Delta^{\left(  2\right)  }\right)
=K$.

[\textit{Remark:} Exercise~\ref{exe.subcoalgebra.ABS}(a) would not hold if we
allowed $\kk$ to be an arbitrary commutative ring rather than a field.]

\bigskip

\textit{Hint to Exercise~\ref{exe.antipode.connected.general}.} (a) Here is
Takeuchi's argument: We know that the map $h\mid_{C_{0}}\in\Hom %
\left(  C_{0},A\right)  $ is $\star$-invertible; let $\widetilde{g}$ be its
$\star$-inverse. Extend $\widetilde{g}$ to a $\kk$-linear map
$g:C\rightarrow A$ by defining it as $0$ on every $C_{n}$ for $n>0$. It is
then easy to see that $\left(  h\star g\right)  \mid_{C_{0}}=\left(  g\star
h\right)  \mid_{C_{0}}=\left(  u\epsilon\right)  \mid_{C_{0}}$. This allows us
to assume WLOG that $h\mid_{C_{0}}=\left(  u\epsilon\right)  \mid_{C_{0}}$
(because once we know that $h\star g$ and $g\star h$ are $\star$-invertible,
it follows that so is $h$). Assuming this, we conclude that $h-u\epsilon$
annihilates $C_{0}$. Define $f$ as $h-u\epsilon$. Now, we can proceed as in
the proof of Proposition~\ref{Takeuchi-prop} to show that $\sum_{k\geq
0}\left(  -1\right)  ^{k}f^{\star k}$ is a well-defined linear map
$C\rightarrow A$ and a two-sided $\star$-inverse for $h$. Thus, $h$ is $\star
$-invertible, and part (a) of the exercise is proven. (An alternative proof
proceeds by mimicking the proof of
Proposition~\ref{graded-connected-bialgebras-have-antipodes}, again by first
assuming WLOG that $h\mid_{C_{0}}=\left(  u\epsilon\right)  \mid_{C_{0}}$.)

(b) Apply part (a) to $C=A$ and the map $\id_{A}:A\rightarrow A$.

(c) Applying part (b), we see that $A$ is a Hopf algebra (since $A_{0}%
=\mathbf{k}$ is a Hopf algebra) in the setting of
Proposition~\ref{graded-connected-bialgebras-have-antipodes}. This yields the
existence of the antipode. Its uniqueness is trivial, and its gradedness
follows from Exercise \ref{exe.coalg-anti.gen}(e).

\bigskip

\textit{Hint to Exercise~\ref{exe.I=0}.} (a) Let $I$ be a two-sided coideal of
$A$ such that $I\cap\mathfrak{p}=0$ and such that $I=\bigoplus_{n\geq0}\left(
I\cap A_{n}\right)  $. Let $I_{n}=I\cap A_{n}$ for every $n\in \NN$.
Then, $I=\bigoplus_{n\geq0}I_{n}$. Since $I$ is a two-sided coideal, we have
$\epsilon\left(  I\right)  =0$.

We want to prove that $I=0$. It clearly suffices to show that every
$n\in \NN$ satisfies $I_{n}=0$ (since $I=\bigoplus_{n\geq0}I_{n}$). We
shall show this by strong induction: We fix an $N\in \NN$, and we assume
(as induction hypothesis) that $I_{n}=0$ for all $n<N$. We must prove that
$I_{N}=0$.

Fix $i\in I_{N}$; we aim to show that $i=0$. We have $i\in I_{N}\subset A_{N}$
and thus $\Delta\left(  i\right)  \in\left(  A\otimes A\right)  _{N}$ (since
$\Delta$ is a graded map). On the other hand, from $i\in I_{N}\subset I$, we
obtain
\begin{align*}
\Delta\left(  i\right)   &  \in\Delta\left(  I\right)  \subset\underbrace{I}%
_{=\bigoplus_{n\geq0}I_{n}}\otimes\underbrace{A}_{=\bigoplus_{m\geq0}A_{m}%
}+\underbrace{A}_{=\bigoplus_{m\geq0}A_{m}}\otimes\underbrace{I}%
_{=\bigoplus_{n\geq0}I_{n}}\ \ \ \ \ \ \ \ \ \ \left(  \text{since }I\text{ is
a two-sided coideal}\right) \\
&  =\sum\limits_{\left(  m,n\right)  \in \NN^{2}}I_{n}\otimes
A_{m}+\sum\limits_{\left(  m,n\right)  \in \NN^{2}}A_{m}\otimes I_{n}.
\end{align*}
Combining this with $\Delta\left(  i\right)  \in\left(  A\otimes A\right)
_{N}$, we obtain%
\begin{align*}
\Delta\left(  i\right)   &  \in\sum\limits_{\substack{\left(  m,n\right)
\in \NN^{2};\\m+n=N}}I_{n}\otimes A_{m}+\sum\limits_{\substack{\left(
m,n\right)  \in \NN^{2};\\m+n=N}}A_{m}\otimes I_{n}%
\ \ \ \ \ \ \ \ \ \ \left(  \text{since }I_{n}\otimes A_{m}\text{ and }%
A_{m}\otimes I_{n}\text{ are subsets of }\left(  A\otimes A\right)
_{n+m}\right) \\
&  =\sum_{n=0}^{N}I_{n}\otimes A_{N-n}+\sum_{n=0}^{N}A_{N-n}\otimes I_{n}\\
&  =I_{N}\otimes\underbrace{A_{0}}_{=\mathbf{k}\cdot1_{A}}+\sum_{n=0}%
^{N-1}\underbrace{I_{n}}_{\substack{=0\\\text{(by the induction}%
\\\text{hypothesis)}}}\otimes A_{N-n}+\underbrace{A_{0}}_{=\mathbf{k}%
\cdot1_{A}}\otimes I_{N}+\sum_{n=0}^{N-1}A_{N-n}\otimes\underbrace{I_{n}%
}_{\substack{=0\\\text{(by the induction}\\\text{hypothesis)}}}\\
&  =I_{N}\otimes\left(  \mathbf{k}\cdot1_{A}\right)  +\left(  \mathbf{k}%
\cdot1_{A}\right)  \otimes I_{N}.
\end{align*}
In other words,
\begin{equation}
\Delta\left(  i\right)  =j\otimes1_{A}+1_{A}\otimes k \label{hint.I=0.deltai}%
\end{equation}
for some $j,k\in I_{N}$. By applying $\epsilon\otimes \id$ to
both sides of this equality, and recalling the commutativity of
\eqref{counit-diagram}, we obtain $i=\epsilon\left(  j\right)  1_{A}+k$. But
$\epsilon\left(  j\right)  =0$ (since $j\in I_{N}\subset I$, so $\epsilon
\left(  j\right)  \in\epsilon\left(  I\right)  =0$), so this simplifies to
$i=k$. Similarly, $i=j$. Hence, \eqref{hint.I=0.deltai} rewrites as
$\Delta\left(  i\right)  =i\otimes1_{A}+1_{A}\otimes i$, which shows that
$i\in\mathfrak{p}$, hence $i\in I\cap\mathfrak{p}=0$ and thus $i=0$. This was
for proved for each $i\in I_{N}$, so we obtain $I_{N}=0$. This completes the
induction step, and so part (a) is solved.

(b) Exercise \ref{exe.kernel-is-coideal}(a) shows that $\ker f$ is a two-sided
coideal of $C$. If $f\mid_{\mathfrak{p}}$ is injective, then $\left(  \ker
f\right)  \cap\mathfrak{p}=0$. Now, apply part (a) of the current exercise to
$I=\ker f$.

(c) Proceed as in part (b), but use Exercise~\ref{exe.kernel-is-coideal}(b)
instead of Exercise~\ref{exe.kernel-is-coideal}(a).

\bigskip

\textit{Hint to Exercise~\ref{exe.liealg}.} (a) Straightforward (if slightly
laborious) computations.

(b) Direct verification (the hard part of which has been done in
\eqref{commutator-of-primitives} already).

(c) For every subset $S$ of a $\kk$-module $U$, we let $\left\langle
S\right\rangle $ denote the $\kk$-submodule of $U$ spanned by $S$. Our
definition of $J$ thus becomes%
\begin{equation}
J=T\left(  \mathfrak{p}\right)  \cdot C\cdot T\left(  \mathfrak{p}\right)  ,
\label{hint.liealg.1}%
\end{equation}
where $C=\left\langle xy-yx-\left[  x,y\right]  \ \mid\ x,y\in\mathfrak{p}%
\right\rangle $. A simple computation shows that each element of $C$ is
primitive. Hence,
\[
\Delta\left(  C\right)  \subset C\otimes T\left(  \mathfrak{p}\right)
+T\left(  \mathfrak{p}\right)  \otimes C.
\]
Applying $\Delta$ to both sides of \eqref{hint.liealg.1}, and recalling that
$\Delta$ is a $\kk$-algebra homomorphism, we find%
\begin{align*}
\Delta\left(  J\right)   &  =\underbrace{\Delta\left(  T\left(  \mathfrak{p}%
\right)  \right)  }_{\subset T\left(  \mathfrak{p}\right)  \otimes T\left(
\mathfrak{p}\right)  }\cdot\underbrace{\Delta\left(  C\right)  }_{\subset
C\otimes T\left(  \mathfrak{p}\right)  +T\left(  \mathfrak{p}\right)  \otimes
C}\cdot\underbrace{\Delta\left(  T\left(  \mathfrak{p}\right)  \right)
}_{\subset T\left(  \mathfrak{p}\right)  \otimes T\left(  \mathfrak{p}\right)
}\\
&  \subset\left(  T\left(  \mathfrak{p}\right)  \otimes T\left(
\mathfrak{p}\right)  \right)  \cdot\left(  C\otimes T\left(  \mathfrak{p}%
\right)  +T\left(  \mathfrak{p}\right)  \otimes C\right)  \cdot\left(
T\left(  \mathfrak{p}\right)  \otimes T\left(  \mathfrak{p}\right)  \right) \\
&  =J\otimes T\left(  \mathfrak{p}\right)  +T\left(  \mathfrak{p}\right)
\otimes J.
\end{align*}
A similar (but simpler) argument shows $\epsilon\left(  J\right)  =0$. Thus,
$J$ is a two-sided coideal of $T\left(  \mathfrak{p}\right)  $. This yields
that $T\left(  \mathfrak{p}\right)  /J$ is a $\kk$-bialgebra.

(d) We need to show that $S\left(  J\right)  \subset J$. This can be done in a
similar way as we proved $\Delta\left(  J\right)  \subset J\otimes T\left(
\mathfrak{p}\right)  +T\left(  \mathfrak{p}\right)  \otimes J$ in part (c),
once you know (from Proposition~\ref{antipodes-are-antiendomorphisms}) that
the antipode $S$ of $T\left(  \mathfrak{p}\right)  $ is a $\kk$-algebra anti-homomorphism.

\bigskip

\textit{Hint to Exercise~\ref{exe.commutative.convolution}.} Straightforward
and easy verification.

\bigskip

\textit{Hint to Exercise~\ref{exe.cocommutative.coalghom}.} Straightforward
and easy verification. Parts (a) and (b) are dual, of course.

\bigskip

\textit{Hint to Exercise~\ref{exe.comm-cocomm.anti}.} (a) Straightforward and
easy verification.

(b) The dual says the following: Let $A$ and $B$ be two $\mathbf{k}%
$-coalgebras, at least one of which is cocommutative. Prove that the
$\kk$-coalgebra anti-homomorphisms from $A$ to $B$ are the same as the
$\kk$-coalgebra homomorphisms from $A$ to $B$.

\bigskip

\textit{Hint to Exercise~\ref{exe.iterated.comm}.} For every $1\leq i<j\leq
k$, let $t_{i,j}$ be the transposition in $\mathfrak{S}_{k}$ which transposes
$i$ with $j$. It is well-known that the symmetric group $\mathfrak{S}_{k}$ is
generated by the transpositions $t_{i,i+1}$ with $i$ ranging over $\left\{
1,2,\ldots,k-1\right\}  $. However, we have $\left(  \rho\left(  \pi\right)
\right)  \circ\left(  \rho\left(  \psi\right)  \right)  =\rho\left(  \pi
\psi\right)  $ for any two elements $\pi$ and $\psi$ of $\mathfrak{S}_{k}$.
Thus, it suffices to check that
\[
m^{\left(  k-1\right)  }\circ\left(  \rho\left(  t_{i,i+1}\right)  \right)
=m^{\left(  k-1\right)  }\qquad\qquad\text{for all }i\in\left\{
1,2,\ldots,k-1\right\}  .
\]
But this is not hard to check using $m^{\left(  k-1\right)  }=m^{\left(
k-2\right)  }\circ\left(  \id_{A^{\otimes\left(  i-1\right)  }%
}\otimes m\otimes \id_{A^{\otimes\left(  k-1-i\right)  }}\right)
$ (a consequence of Exercise~\ref{exe.iterated.m}(c)) and $m\circ T=m$.

\bigskip

\textit{Hint to Exercise~\ref{exe.iterated.cocomm}.} Here is the dual
statement: Let $C$ be a cocommutative $\kk$-coalgebra, and let $k
\in \NN$. The symmetric group $\mathfrak{S}_{k}$ acts on the $k$-fold
tensor power $C^{\otimes k}$ by permuting the tensor factors: $\sigma\left(
v_{1} \otimes v_{2} \otimes\cdots\otimes v_{k}\right)  = v_{\sigma^{-1}\left(
1\right)  } \otimes v_{\sigma^{-1}\left(  2\right)  } \otimes\cdots\otimes
v_{\sigma^{-1}\left(  k\right)  }$ for all $v_{1}, v_{2}, \ldots, v_{k} \in C$
and $\sigma\in\mathfrak{S}_{k}$. For every $\pi\in\mathfrak{S}_{k}$, denote by
$\rho\left(  \pi\right)  $ the action of $\pi$ on $C^{\otimes k}$ (this is an
endomorphism of $C^{\otimes k}$). Show that every $\pi\in\mathfrak{S}_{k}$
satisfies $\left(  \rho\left(  \pi\right)  \right)  \circ\Delta^{\left(
k-1\right)  } = \Delta^{\left(  k-1\right)  }$. (Recall that $\Delta^{\left(
k-1\right)  } : C \to C^{\otimes k}$ is defined as in
Exercise~\ref{exe.iterated.Delta} for $k \geq1$, and by $\Delta^{\left(
-1\right)  } = \epsilon: C \to\mathbf{k}$ for $k = 0$.)

\bigskip

\textit{Hint to Exercise~\ref{exe.alghoms.comm}.} (a) Use
Exercise~\ref{exe.cocommutative.coalghom}(b) and
Exercise~\ref{exe.tensorprod.morphism}(a) to represent $f\star g$ as a
composition of three $\kk$-algebra homomorphisms.

(b) Induction on $k$, using part (a).

(c) Use Proposition~\ref{antipodes-are-antiendomorphisms}, Proposition
\ref{antipodes-give-convolution-inverses}(a) and the easy fact that a
composition of a $\kk$-algebra homomorphism with a $\kk$-algebra
anti-homomorphism (in either order) always is a $\kk$-algebra anti-homomorphism.

(d) Use Exercise~\ref{exe.cocommutative.coalghom}(b). Then, proceed by
induction on $k$ as in the solution of Exercise~\ref{exe.iterated.patras}(a).

(e) Use Proposition~\ref{prop.convolution.functor}.

(f) Let $H$ be a commutative $\kk$-bialgebra. Let $k$ and $\ell$ be two
nonnegative integers. Then, Exercise~\ref{exe.alghoms.comm}(b) (applied to
$A=H$ and $f_{i}= \id_{H}$) yields that $\id_{H}^{\star k}$ is a
$\kk$-algebra homomorphism $H\rightarrow H$. Now,
apply Exercise~\ref{exe.alghoms.comm}(e) to $H$, $H$, $H$, $H$, $\ell$,
$\id_{H}$, $\id_{H}^{\star k}$ and
$\id_{H}$ instead of $C$, $C^{\prime}$, $A$, $A^{\prime}$, $k$,
$f_{i}$, $\alpha$ and $\gamma$.

(g) This is an exercise in bootstrapping. First, let $k\in \NN$. Then,
part (b) of this exercise shows that $\id_{H}^{\star k}$ is a
$\kk$-algebra homomorphism. Use this together with part (c) to conclude
that $\id_{H}^{\star k}\circ S$ is again a $\kk$-algebra
homomorphism and a $\star$-inverse to $\id_{H}^{\star k}$; thus,
$\id_{H}^{\star k}\circ S
=\left(  \id_{H}^{\star k}\right)  ^{\star\left(  -1\right)  }
=\id_{H}^{\star\left(  -k\right)  }$, and this map
$\id_{H}^{\star\left(  -k\right)  }$ is a $\kk$-algebra homomorphism.

Now forget that we fixed $k$. We thus have shown that
$\id_{H}^{\star k}$ and $\id_{H}^{\star\left(  -k\right)  }$ are
$\kk$-algebra homomorphisms for each $k\in \NN$. In other words,
\begin{equation}
\id_{H}^{\star k}\text{ is a }\mathbf{k}\text{-algebra
homomorphism}\qquad\qquad\text{ for every }k\in \ZZ.
\label{hint.alghoms.comm.g.alghom}%
\end{equation}

Furthermore, we have proved the equality $\id_{H}^{\star k}\circ
S= \id_{H}^{\star\left(  -k\right)  }$ for each $k\in \NN$.
Repeating the proof of this, but now taking $k\in \ZZ$ instead of
$k\in \NN$, we conclude that it also holds for each $k\in \ZZ$
(since we already have proved \eqref{hint.alghoms.comm.g.alghom}). In other
words,%
\begin{equation}
\id_{H}^{\star\left(  -k\right)  }= \id_{H}^{\star
k}\circ S\qquad\qquad\text{ for every }k\in \ZZ.
\label{hint.alghoms.comm.g.starinv}%
\end{equation}

Now, fix two integers $k$ and $\ell$. From \eqref{hint.alghoms.comm.g.alghom},
we know that $\id_{H}^{\star k}$ is a $\kk$-algebra
homomorphism. Hence, if $\ell$ is nonnegative, then we can prove
$\id_{H}^{\star k}\circ \id_{H}^{\star\ell
}= \id_{H}^{\star\left(  k\ell\right)  }$ just as we did in the
solution to Exercise~\ref{exe.alghoms.comm}(f). But the case when $\ell$ is
negative can be reduced to the previous case by applying
\eqref{hint.alghoms.comm.g.starinv} (once to $-\ell$ instead of $k$, and once
again to $-k\ell$ instead of $k$). Thus, in each case, we obtain
$\id_{H}^{\star k}\circ \id_{H}^{\star\ell}
= \id_{H}^{\star\left(  k\ell\right)  }$.

(h) The dual of Exercise~\ref{exe.alghoms.comm}(a) is the following exercise:

\begin{statement}
If $H$ is a $\kk$-bialgebra and $C$ is a cocommutative $\mathbf{k}%
$-coalgebra, and if $f$ and $g$ are two $\kk$-coalgebra homomorphisms
$C \to H$, then prove that $f \star g$ also is a $\kk$-coalgebra
homomorphism $C \to H$.
\end{statement}

The dual of Exercise~\ref{exe.alghoms.comm}(b) is the following exercise:

\begin{statement}
If $H$ is a $\kk$-bialgebra and $C$ is a cocommutative $\mathbf{k}%
$-coalgebra, and if $f_{1}, f_{2}, \ldots, f_{k}$ are several $\mathbf{k}%
$-coalgebra homomorphisms $C \to H$, then prove that $f_{1} \star f_{2}
\star\cdots\star f_{k}$ also is a $\kk$-coalgebra homomorphism $C \to
H$.
\end{statement}

The dual of Exercise~\ref{exe.alghoms.comm}(c) is the following exercise:

\begin{statement}
If $H$ is a Hopf algebra and $C$ is a cocommutative $\kk$-coalgebra,
and if $f : C \to H$ is a $\kk$-coalgebra homomorphism, then prove that
$S \circ f : C \to H$ (where $S$ is the antipode of $H$) is again a
$\kk$-coalgebra homomorphism, and is a $\star$-inverse to $f$.
\end{statement}

The dual of Exercise~\ref{exe.alghoms.comm}(d) is the following exercise:

\begin{statement}
If $C$ is a cocommutative $\kk$-coalgebra, then show that
$\Delta^{\left(  k\right)  }$ is a $\kk$-coalgebra homomorphism for
every $k \in \NN$. (The map $\Delta^{\left(  k\right)  } : C \to
C^{\otimes\left(  k+1\right)  }$ is defined as in
Exercise~\ref{exe.iterated.Delta}.)
\end{statement}

The dual of Exercise~\ref{exe.alghoms.comm}(e) is
Exercise~\ref{exe.alghoms.comm}(e) itself (up to renaming objects and maps).

The dual of Exercise~\ref{exe.alghoms.comm}(f) is the following exercise:

\begin{statement}
If $H$ is a cocommutative $\kk$-bialgebra, and $k$ and $\ell$ are two
nonnegative integers, then prove that $\id_{H}^{\star\ell}
\circ \id_{H}^{\star k} = \id_{H}^{\star\left(  \ell k\right)  }$.
\end{statement}

The dual of Exercise~\ref{exe.alghoms.comm}(g) is the following exercise:

\begin{statement}
If $H$ is a cocommutative $\kk$-Hopf algebra, and $k$ and $\ell$ are
two integers, then prove that $\id_{H}^{\star\ell}
\circ \id_{H}^{\star k} = \id_{H}^{\star\left(  \ell k\right)  }$.
\end{statement}

\bigskip

\textit{Hint to Exercise~\ref{exe.cocomm-antipode-S2}.} This is dual to
Corollary \ref{commutative-implies-involutive-antipode-cor} (but can also
easily be shown using Exercise~\ref{exe.coalg-anti.gen}(b),
Exercise~\ref{exe.comm-cocomm.anti}(b) and
Proposition~\ref{antipodes-give-convolution-inverses}(b)).

\bigskip

\textit{Hint to Exercise~\ref{exe.dynkin}.} (a) This can be proved
computationally (using Sweedler notation), but there is a nicer argument as well:

A \emph{coderivation} of a $\kk$-coalgebra $\left(  C,\Delta
,\epsilon\right)  $ is defined as a $\kk$-linear map $F:C\rightarrow C$
such that $\Delta\circ F=\left(  F\otimes \id+ \id\otimes F\right)
\circ\Delta$. (The reader can check that this axiom is the
result of writing the axiom for a derivation in element-free terms and
reversing all arrows. Nothing less should be expected.) It is easy to see that
$E$ is a coderivation. Hence, it will be enough to check that $\left(  S\star
f\right)  \left(  a\right)  $ and $\left(  f\star S\right)  \left(  a\right)
$ are primitive whenever $f:A\rightarrow A$ is a coderivation and $a\in A$. So
fix a coderivation $f:A\rightarrow A$. Notice that the antipode $S$ of $A$ is
a coalgebra anti-endomorphism (by Exercise~\ref{exe.coalg-anti}), thus a
coalgebra endomorphism (by Exercise~\ref{exe.comm-cocomm.anti}(b)). Thus,
$\Delta\circ S=\left(  S\otimes S\right)  \circ\Delta$. Moreover,
$\Delta:A\rightarrow A\otimes A$ is a coalgebra homomorphism (by
Exercise~\ref{exe.cocommutative.coalghom}(a)) and an algebra homomorphism
(since $A$ is a bialgebra). Applying
\eqref{pre-and-post-composition-in-convolution} to $A\otimes A$, $A$, $A$,
$\Delta$, $\id_{A}$, $S$ and $f$ instead of $A^{\prime}$, $C$,
$C^{\prime}$, $\alpha$, $\gamma$, $f$ and $g$, we obtain
\begin{align*}
\Delta\circ\left(  S\star f\right)   &  =\underbrace{\left(  \Delta\circ
S\right)  }_{=\left(  S\otimes S\right)  \circ\Delta}\star\underbrace{\left(
\Delta\circ f\right)  }_{\substack{=\left(  f\otimes \id
+ \id \otimes f\right)  \circ\Delta\\\text{(since }f\text{ is a
coderivation)}}}\\
&  =\left(  \left(  S\otimes S\right)  \circ\Delta\right)  \star\left(
\left(  f\otimes \id + \id \otimes f\right)
\circ\Delta\right)  =\left(  \left(  S\otimes S\right)  \star\left(
f\otimes \id + \id \otimes f\right)  \right)
\circ\Delta\\
&  =\underbrace{\left(  \left(  S\otimes S\right)  \star\left(  f\otimes
\id \right)  \right)  }_{\substack{=\left(  S\star f\right)
\otimes\left(  S\star \id \right)  \\\text{(by
Exercise~\ref{exe.convolution.tensor}(a))}}}\circ\Delta+\underbrace{\left(
\left(  S\otimes S\right)  \star\left(  \id \otimes f\right)
\right)  }_{\substack{=\left(  S\star \id \right)  \otimes\left(
S\star f\right)  \\\text{(by Exercise~\ref{exe.convolution.tensor}(a))}}%
}\circ\Delta\\
&  =\left(  \left(  S\star f\right)  \otimes\underbrace{\left(  S\star
\id \right)  }_{=u\epsilon}\right)  \circ\Delta+\left(
\underbrace{\left(  S\star \id \right)  }_{=u\epsilon}%
\otimes\left(  S\star f\right)  \right)  \circ\Delta\\
&  =\left(  \left(  S\star f\right)  \otimes u\epsilon\right)  \circ
\Delta+\left(  u\epsilon\otimes\left(  S\star f\right)  \right)  \circ\Delta.
\end{align*}
Hence, every $a\in A$ satisfies
\begin{align*}
\left(  \Delta\circ\left(  S\star f\right)  \right)  \left(  a\right)   &
=\left(  \left(  \left(  S\star f\right)  \otimes u\epsilon\right)
\circ\Delta+\left(  u\epsilon\otimes\left(  S\star f\right)  \right)
\circ\Delta\right)  \left(  a\right) \\
&  =\left(  S\star f\right)  \left(  a\right)  \otimes1+1\otimes\left(  S\star
f\right)  \left(  a\right)
\end{align*}
(after some brief computations using \eqref{counit-diagram}). In other words,
for every $a\in A$, the element $\left(  S\star f\right)  \left(  a\right)  $
is primitive. Similarly the same can be shown for $\left(  f\star S\right)
\left(  a\right)  $, and so we are done.

(b) is a very simple computation. (Alternatively, the $\left(  S\star
E\right)  \left(  p\right)  =E\left(  p\right)  $ part follows from applying
part (c) to $a=1$, and similarly one can show $\left(  E\star S\right)
\left(  p\right)  =E\left(  p\right)  $.)

(c) This is another computation, using
Proposition~\ref{antipode-on-primitives} and the (easy) observation that $E$
is a derivation of the algebra $A$.

(d) Assume that the graded algebra $A=\bigoplus_{n\geq0}A_{n}$ is connected
and that ${\mathbb{Q}}$ is a subring of $\kk$. Let $B$ be the
$\kk$-subalgebra of $A$ generated by $\mathfrak{p}$. In order to prove
part (d), we need to show that $A\subset B$. Clearly, it suffices to show that
$A_{n}\subset B$ for every $n\in \NN$. We prove this by strong
induction on $n$; thus, we fix some $n\in \NN$, and assume as induction
hypothesis that $A_{m}\subset B$ for every $m<n$. Our goal is then to show
that $A_{n}\subset B$. This being trivial for $n=0$ (since $A$ is connected),
we WLOG assume that $n>0$. Let $a\in A_{n}$. Part (a) of this exercise yields
$\left(  S\star E\right)  \left(  a\right)  \in\mathfrak{p}\subset B$. On the
other hand, Exercise \ref{graded-connected-exercise}(h) (applied to $x=a$)
yields%
\[
\Delta\left(  a\right)  \in1\otimes a+a\otimes1+\sum_{k=1}^{n-1}A_{k}\otimes
A_{n-k}.
\]
Hence, from the definition of convolution, we obtain%
\begin{align*}
\left(  S\star E\right)  \left(  a\right)   &  \in\underbrace{S\left(
1\right)  }_{=1}E\left(  a\right)  +S\left(  a\right)  \underbrace{E\left(
1\right)  }_{=0}+\underbrace{\left(  m\circ\left(  S\otimes E\right)  \right)
\left(  \sum_{k=1}^{n-1}A_{k}\otimes A_{n-k}\right)  }_{=\sum_{k=1}%
^{n-1}S\left(  A_{k}\right)  E\left(  A_{n-k}\right)  }\\
&  =E\left(  a\right)  +\sum_{k=1}^{n-1}\underbrace{S\left(  A_{k}\right)
}_{\substack{\subset A_{k}\\\text{(since }S\text{ is graded)}}%
}\underbrace{E\left(  A_{n-k}\right)  }_{\substack{\subset A_{n-k}\subset
B\\\text{(by the induction}\\\text{hypothesis)}}}\subset E\left(  a\right)
+\sum_{k=1}^{n-1}\underbrace{A_{k}}_{\substack{\subset B\\\text{(by the
induction}\\\text{hypothesis)}}}B\subset E\left(  a\right)  +B
\end{align*}
(since $B$ is a subalgebra). Hence, $E\left(  a\right)  \in\left(  S\star
E\right)  \left(  a\right)  +B=B$ (since $\left(  S\star E\right)  \left(
a\right)  \in B$). Since $E\left(  a\right)  =na$, this becomes $na\in B$,
thus $a\in B$ (since $\mathbb{Q}$ is a subring of $\kk$). Since we have
shown this for each $a\in A_{n}$, we thus obtain $A_{n}\subset B$, and our
induction is complete.

This solution of part (d) is not the most generalizable one -- for instance,
(d) also holds if $A$ is connected filtered instead of connected graded, and
then a different argument is necessary. This is a part of the
Cartier-Milnor-Moore theorem, and appears e.g. in \cite[\S 3.2]%
{DuchampMinhTolluChienNghia}.

(e) If $a\in T\left(  V\right)  $ is homogeneous of positive degree and $p\in
V$, then part (c) quickly yields $\left(  S\star E\right)  \left(  ap\right)
=\left[  \left(  S\star E\right)  \left(  a\right)  ,p\right]  $. This allows
proving (e) by induction over $n$, with the induction base $n=1$ being a
consequence of part (b).

\bigskip

\textit{Hint to Exercise~\ref{exe.dual-algebra}.} (a) This can be done by
diagram chasing. For example, if $\mathfrak{m}$ denotes the map $\Delta
_{C}^{\ast}\circ\rho_{C,C}:C^{\ast}\otimes C^{\ast}\rightarrow C^{\ast}$, then
the diagram
\[
\xymatrix{
& & C^* \otimes C^* \otimes C^* \ar@/_5pc/[ddll]^{\mathfrak{m}\otimes\id}
\ar[dl]^{\rho_{C,C}\otimes\id} \ar[dr]^{\id\otimes\rho_{C,C}} \ar@
/^5pc/[ddrr]^{\id\otimes\mathfrak{m}} & & \\
& \left(C\otimes C\right)^* \otimes C^* \ar[dl]^{\Delta_C^*\otimes\id}
\ar[dr]^{\rho_{C\otimes C,C}} & & C^* \otimes\left(C\otimes C\right
)^* \ar[dr]_{\id\otimes\Delta_C^*} \ar[dl]^{\rho_{C,C\otimes C}} & \\
C^* \otimes C^* \ar[dr]^{\rho_{C,C}} \ar@/_5pc/[ddrr]^{\mathfrak{m}}
& & \left(C\otimes C\otimes C\right)^* \ar[dl]^{\left(\Delta_C\otimes\id
\right)^*} \ar[dr]^{\left(\id\otimes\Delta_C\right)^*} & & C^* \otimes
C^* \ar[dl]_{\rho_{C,C}} \ar@/^5pc/[ddll]^{\mathfrak{m}} \\
& \left(C\otimes C\right)^* \ar[dr]^{\Delta_C^*} & & \left(C\otimes
C\right)^* \ar[dl]^{\Delta_C^*} & \\
& & C^* & &
}
\]
is commutative (since each of its little triangles and squares is); thus,
$\mathfrak{m}\circ\left(  \mathfrak{m}\otimes \id \right)
=\mathfrak{m}\circ\left(  \id \otimes\mathfrak{m}\right)  $ for
$\mathfrak{m}$. This proves that the diagram \eqref{associativity-diagram}
commutes for our algebra $C^{\ast}$. The commutativity of \eqref{unit-diagram}
is obtained similarly.

Alternatively, we could also solve part (a) trivially by first solving part
(b) and then recalling Exercise \ref{exe.convolution.assoc}.

(b) Straightforward verification on pure tensors.

(c) Let $C=\bigoplus_{n\geq0}C_{n}$ be a graded $\kk$-coalgebra. For
every $n\in \NN$, we identify $\left(  C_{n}\right)  ^{\ast}$ with a
$\kk$-submodule of $C^{\ast}$, namely with the $\kk$-submodule
$\left\{  f\in C^{\ast}\ \mid\ f\left(  C_{p}\right)  =0\text{ for all }%
p\in \NN \text{ satisfying }p\neq n\right\}  $. By the definition of
$C^{o}$, we have $C^{o}=\bigoplus_{n\geq0}\left(  C_{n}\right)  ^{\ast}$.
Hence, it remains to show that $\left(  C_{a}\right)  ^{\ast}\left(
C_{b}\right)  ^{\ast}\subset\left(  C_{a+b}\right)  ^{\ast}$ for all
$a,b\in \NN$, and that $1_{C^{\ast}}\in\left(  C_{0}\right)  ^{\ast}$.
But this is straightforward using the gradedness of $\Delta$ and $\epsilon$.

(d) Diagram chasing or simple element-wise verification.

(e) Simple linear algebra (no Hopf algebras involved here).

(f) The ``only if'' direction is proved in
the same way as part (d) (or as a corollary of part (d), since $D^{\circ}$ and
$C^{\circ}$ are subalgebras of $D^{\ast}$ and $C^{\ast}$). It remains to prove
the ``if'' direction.

Assume that $f^{\ast}:D^{o}\rightarrow C^{o}$ is a $\kk$-algebra
morphism. We want to show that $f:C\rightarrow D$ is a $\kk$-coalgebra
morphism. In other words, we want to show that the two diagrams
\begin{equation}
\xymatrix{ C \ar[d]_{\Delta_C} \ar[r]^{f}
& D \ar[d]^{\Delta_D} \\ C \otimes C \ar[r]^{f \otimes f} & D \otimes
D }\ \ \ \ \ \ \ \ \ \ \text{and}
\ \ \ \ \ \ \ \ \ \ \xymatrix{ C \ar[dr]_{\epsilon_C} \ar[rr]^{f}%
& & D \ar[dl]^{\epsilon_D} \\ &\kk& }
\label{hint.dual-algebra.e.1}
\end{equation}
commute. Let us start with the left one of these diagrams. The graded
$\kk$-module $D$ is of finite type, and therefore the map $\rho
_{D,D}:D^{o}\otimes D^{o}\rightarrow\left(  D\otimes D\right)  ^{o}$ (a
restriction of the map $\rho_{D,D}:D^{\ast}\otimes D^{\ast}\rightarrow\left(
D\otimes D\right)  ^{\ast}$) is an isomorphism. Its inverse $\rho_{D,D}%
^{-1}:\left(  D\otimes D\right)  ^{o}\rightarrow D^{o}\otimes D^{o}$ is
therefore well-defined\footnote{Beware: we don't have an inverse of the
non-restricted map $\rho_{D,D}:D^{\ast}\otimes D^{\ast}\rightarrow\left(
D\otimes D\right)  ^{\ast}$.}. We can thus form the (asymmetric!) diagram
\begin{equation}
{\xymatrixcolsep{5pc}
\xymatrix{ D^o \ar[rrr]^{f^*}
& & & C^o \\ & D^o \otimes D^o \ar[ul]^{m_{D^*}} \ar[r]^{f^* \otimes f^*}
& C^o \otimes C^o \ar[ur]^{m_{C^*}} \ar[dr]^{\rho_{C,C}} & \\ \left(D\otimes
D\right)^o \ar[uu]^{\Delta_D^*} \ar[ur]^{\rho_{D,D}^{-1}} \ar[rrr]_{\left
(f\otimes f\right)^*} & & & \left(C\otimes C\right)^o \ar[uu]^{\Delta_C^*} }%
}.
\label{hint.dual-algebra.e.2}
\end{equation}
(The arrows labelled $m_{C^{\ast}}$ and $m_{D^{\ast}}$ could just as well have
been labelled $m_{C^{o}}$ and $m_{D^{o}}$, since the multiplication maps
$m_{C^{o}}$ and $m_{D^{o}}$ are restrictions of $m_{C^{\ast}}$ and
$m_{D^{\ast}}$.) Argue that the diagram \eqref{hint.dual-algebra.e.2}
commutes. Thus, $f^{\ast}\circ\Delta_{D}^{\ast}=\Delta_{C}^{\ast}\circ\left(
f\otimes f\right)  ^{\ast}$ as maps from $\left(  D\otimes D\right)  ^{o}$ to
$C^{o}$. In other words, $\left(  \Delta_{D}\circ f\right)  ^{\ast}=\left(
\left(  f\otimes f\right)  \circ\Delta_{C}\right)  ^{\ast}$ as maps from
$\left(  D\otimes D\right)  ^{o}$ to $C^{o}$. But a general linear-algebraic
fact states that if $U$ and $V$ are two graded $\kk$-modules such that
$V$ is of finite type, and if $\alpha$ and $\beta$ are two graded $\mathbf{k}%
$-linear maps $U\rightarrow V$ such that $\alpha^{\ast}=\beta^{\ast}$ as maps
from $V^{o}$ to $U^{o}$, then $\alpha=\beta$\ \ \ \ \footnote{This follows
immediately from Exercise~\ref{exe.dual-algebra} (e).}. Hence, $\left(
\Delta_{D}\circ f\right)  ^{\ast}=\left(  \left(  f\otimes f\right)
\circ\Delta_{C}\right)  ^{\ast}$ leads to $\Delta_{D}\circ f=\left(  f\otimes
f\right)  \circ\Delta_{C}$. In other words, the first diagram in
\eqref{hint.dual-algebra.e.1} commutes. The second is similar but easier.
Thus, $f$ is a $\kk$-coalgebra morphism, and the ``if'' direction is proved.

\bigskip

\textit{Hint to Exercise~\ref{exe.Sym.1dim}.} Straightforward computations.
For part (d), first show (independently of whether $\kk$ is a field and
its characteristic) that $\left(  f^{\left(  1\right)  }\right)
^{m}=m!f^{\left(  m\right)  }$ for every $m\in \NN$.

\bigskip

\textit{Hint to Exercise~\ref{exe.Sym.wedge}.} It is best to solve parts (c)
and (d) before approaching (b).

(a) Both maps $\Delta_{\operatorname*{Sym}V}$ and%
\[%
\begin{array}
[c]{ccc}%
\mathbf{k}[{\mathbf{x}}] & \overset{\Delta}{\longrightarrow} & \mathbf{k}%
[{\mathbf{x}},{\mathbf{y}}],\\
f(x_{1},\ldots,x_{n}) & \longmapsto & f(x_{1}+y_{1},\ldots,x_{n}+y_{n})
\end{array}
\]
are $\kk$-algebra homomorphisms. Thus, in order to check that they are
equal, it suffices to verify that they agree on $V$ (since $V$ generates
$\operatorname*{Sym}V$).

(c) This is a straightforward computation unless you get confused with the
topologist's sign convention. The latter convention affects the twist map
$T=T_{T\left(  V\right)  ,T\left(  V\right)  }:T\left(  V\right)  \otimes
T\left(  V\right)  \rightarrow T\left(  V\right)  \otimes T\left(  V\right)  $
(in particular, we now have $T\left(  x\otimes x\right)  =-x\otimes x$ instead
of $T\left(  x\otimes x\right)  =x\otimes x$), and thus also affects the
multiplication in the $\kk$-algebra $T\left(  V\right)  \otimes
T\left(  V\right)  $, because this multiplication is given by%
\[
m_{T\left(  V\right)  \otimes T\left(  V\right)  }=\left(  m_{T\left(
V\right)  }\otimes m_{T\left(  V\right)  }\right)  \circ\left(
\id \otimes T \otimes \id \right)  .
\]
Make sure you understand why this leads to $\left(  1\otimes x\right)
\cdot\left(  x\otimes1\right)  =-x\otimes x$ (whereas $\left(  x\otimes
1\right)  \cdot\left(  1\otimes x\right)  =x\otimes x$).

(d) The trickiest part is showing that $J$ is a graded $\kk$-submodule
of $T\left(  V\right)  $. It suffices to check that $J$ is generated (as a
two-sided ideal) by homogeneous elements\footnote{Make sure you understand
why.}; however, this is not completely trivial, as the designated generators
$x^{2}$ for $x\in V$ need not be homogeneous. However, it helps to observe
that $J$ is also the two-sided ideal generated by the set%
\[
\left\{  x\otimes x\right\}  _{x\in V\text{ is homogeneous}}\cup\left\{
x\otimes y+y\otimes x\right\}  _{x,y\in V\text{ are homogeneous}}%
\]
(why?), which set does consist of homogeneous elements. Thus, $J$ is a graded
$\kk$-submodule of $T\left(  V\right)  $. From part (c), it is easy to
observe that $J$ is a two-sided coideal of $T\left(  V\right)  $ as well.
Hence, $T\left(  V\right)  /J$ inherits a graded $\kk$-bialgebra
structure from $T\left(  V\right)  $. The rest is easy.

(b) is now a consequence of what has been done in (d).

\bigskip

\textit{Hint to Exercise~\ref{exe.convolution.dual}.} Easy and straightforward.

\bigskip

\textit{Hint to Exercise~\ref{exe.shuffle-alg}.} The hint after the exercise
shows the way; here are a few more pointers. The solution proceeds in two steps:

\begin{itemize}
\item \textit{Step 1:} Show that Proposition~\ref{prop.shuffle-alg} holds when
$V$ is a finite free $\kk$-module.

\item \textit{Step 2:} Use this to conclude that
Proposition~\ref{prop.shuffle-alg} always holds.
\end{itemize}

The trick to Step 1 is to reduce the proof to Example \ref{exa.shuffle-alg}.
In a bit more detail: If $V$ is a finite free $\kk$-module with basis
$\left(  v_{1},v_{2},\ldots,v_{n}\right)  $, then we know from Example
\ref{exa.shuffle-alg} that the graded dual $A^{o}$ of its tensor algebra
$A:=T\left(  V\right)  $ is a Hopf algebra whose basis $\left\{  y_{\left(
i_{1},i_{2},\ldots,i_{\ell}\right)  }\right\}  $ is indexed by words in the
alphabet $I:=\left\{  1,2,\ldots,n\right\}  $. This allows us to define a
$\kk$-linear map $\phi:A^{o}\rightarrow T\left(  V\right)  $ by
setting
\[
\phi\left(  y_{\left(  i_{1},i_{2},\ldots,i_{\ell}\right)  }\right)
=v_{i_{1}}v_{i_{2}}\cdots v_{i_{\ell}}\ \ \ \ \ \ \ \ \ \ \text{for every
}\ell\in \NN \text{ and }\left(  i_{1},i_{2},\ldots,i_{\ell}\right)
\in I^{\ell}.
\]
This $\kk$-linear map $\phi$ then is an isomorphism from the Hopf
algebra $A^{o}$ to the putative Hopf algebra $\left(  \operatorname*{Sh}%
\left(  V\right)  ,\shufmult,1_{T\left(  V\right)  },\Delta_{\shuffle}%
,\epsilon,S\right)  $, in the sense that it is invertible (since it sends a
basis to a basis) and satisfies the five equalities%
\begin{align*}
\phi\circ m_{A^{o}}  &  =m_{\shuffle}\circ\left(  \phi\otimes\phi\right)  ,\\
\phi\circ u_{A^{o}}  &  =u,\\
\left(  \phi\otimes\phi\right)  \circ\Delta_{A^{o}}  &  =\Delta_{\shuffle}%
\circ\phi,\\
\epsilon_{A^{o}}  &  =\epsilon\circ\phi,\\
\phi\circ S_{A^{o}}  &  =S\circ\phi
\end{align*}
(check all these -- for instance, the first of these equalities follows by
comparing \eqref{eq.exa.shuffle-alg.m3} with the definition of $\shufmult$).
Thus, the latter putative Hopf algebra is an actual Hopf algebra (since the
former is). This proves Proposition~\ref{prop.shuffle-alg} for our finite free
$V$, and thus completes Step 1.

Step 2 demonstrates the power of functoriality. We want to prove
Proposition~\ref{prop.shuffle-alg} in the general case, knowing that it holds
when $V$ is finite free. So let $V$ be an arbitrary $\kk$-module. For
the sake of brevity, we shall write $\mathbf{V}$ for $T\left(  V\right)  $.
Let $m_{\shuffle}$ denote the $\kk$-linear map $\mathbf{V}%
\otimes\mathbf{V}\rightarrow\mathbf{V}$ which sends every $a\otimes b$ to
$a\shufmult b$. One of the things that need to be shown is the commutativity
of the diagram
\begin{equation}
\xymatrix{
& \mathbf{V} \otimes\mathbf{V} \ar[dl]_{\Delta_\shuffle\otimes\Delta_\shuffle}
\ar[ddr]_{m_\shuffle} & \\
\mathbf{V} \otimes\mathbf{V} \otimes\mathbf{V} \otimes\mathbf{V}
\ar[dd]_{\id\otimes T \otimes\id} & & \\
& & \mathbf{V} \ar[ddl]_{\Delta_\shuffle}\\
\mathbf{V} \otimes\mathbf{V} \otimes\mathbf{V} \otimes\mathbf{V}
\ar[dr]_{m_\shuffle\otimes m_\shuffle} & & \\
& \mathbf{V} \otimes\mathbf{V} &
}
,
\label{hint.shuffle-alg.diag-to-prove}
\end{equation}
where $T$ is the twist map $T_{\mathbf{V},\mathbf{V}}$. By linearity, it is
clearly enough to verify this only on the pure tensors; that is, it is enough
to check that every $a\in\mathbf{V}$ and $b\in\mathbf{V}$ satisfy
\begin{equation}
\left(  \left(  m_{\shuffle}\otimes m_{\shuffle}\right)  \circ\left(
\id \otimes T \otimes \id \right)  \circ\left(
\Delta_{\shuffle}\otimes\Delta_{\shuffle}\right)  \right)  \left(  a\otimes
b\right)  =\left(  \Delta_{\shuffle}\circ m_{\shuffle}\right)  \left(
a\otimes b\right)  . \label{hint.shuffle-alg.diag-toprove}%
\end{equation}
So let $a,b\in\mathbf{V}$ be arbitrary. WLOG assume that $a=v_{1}v_{2}\cdots
v_{p}$ and $b=v_{p+1}v_{p+2}\cdots v_{p+q}$ for some $p,q\in \NN$ and
$v_{1},v_{2},\ldots,v_{p+q}\in V$. Define $W$ to be the free $\mathbf{k}%
$-module with basis $\left(  x_{1},x_{2},\ldots,x_{p+q}\right)  $, and let
$\mathbf{W}$ be its tensor algebra $T\left(  W\right)  $. Then, $W$ is a
finite free $\kk$-module, and so we know from Step 1 that
Proposition~\ref{prop.shuffle-alg} holds for $W$ instead of $V$. But we can
define a $\kk$-linear map $f:W\rightarrow V$ that sends $x_{1}%
,x_{2},\ldots,x_{p+q}$ to $v_{1},v_{2},\ldots,v_{p+q}$, respectively. This map
$f:W\rightarrow V$ clearly induces a $\kk$-algebra homomorphism
$\mathbf{f}:=T\left(  f\right)  :\mathbf{W}\rightarrow\mathbf{V}$ that
respects all relevant shuffle-algebraic structure (i.e., it satisfies
$\mathbf{f}\circ m_{\shuffle}=m_{\shuffle}\circ\left(  \mathbf{f}%
\otimes\mathbf{f}\right)  $ and $\left(  \mathbf{f}\otimes\mathbf{f}\right)
\circ\Delta_{\shuffle}=\Delta_{\shuffle}\circ\mathbf{f}$ and so on), simply
because this structure has been defined canonically in terms of each of $V$
and $W$. Thus, in the diagram%
\[
\xymatrixcolsep{3.5pc}
\xymatrix{
& &  \mathbf{W} \otimes\mathbf{W} \ar[ddll]_{\Delta_\shuffle\otimes
\Delta_\shuffle}
\ar[dddrr]^{m_\shuffle} \ar[d]_{\mathbf{f}\otimes\mathbf{f}} & \\
& & \mathbf{V} \otimes\mathbf{V} \ar[dl]_{\Delta_\shuffle\otimes
\Delta_\shuffle} \ar[ddr]_{m_\shuffle} & \\
\mathbf{W} \otimes\mathbf{W} \otimes\mathbf{W} \otimes\mathbf{W} \ar
[dd]_{\id\otimes T \otimes\id}
\ar[r]_-{\mathbf{f}\otimes\mathbf{f}\otimes\mathbf{f}\otimes\mathbf{f}} &
\mathbf{V} \otimes\mathbf{V} \otimes\mathbf{V} \otimes\mathbf{V}
\ar[dd]_{\id\otimes T \otimes\id}
& & \\
& & &  \mathbf{V} \ar[ddl]_{\Delta_\shuffle} &  \mathbf{W} \ar[dddll]^{\Delta
_\shuffle}
\ar[l]^{\mathbf{f}} \\
\mathbf{W} \otimes\mathbf{W} \otimes\mathbf{W} \otimes\mathbf{W}
\ar[ddrr]_{m_\shuffle\otimes m_\shuffle}
\ar[r]^-{\mathbf{f}\otimes\mathbf{f}\otimes\mathbf{f}\otimes\mathbf{f}} &
\mathbf{V} \otimes\mathbf{V} \otimes\mathbf{V} \otimes\mathbf{V}
\ar[dr]_{m_\shuffle\otimes m_\shuffle}
& & \\
& &  \mathbf{V} \otimes\mathbf{V} & \\
& &  \mathbf{W} \otimes\mathbf{W} \ar[u]^{\mathbf{f}\otimes\mathbf{f}} &
}
,
\]
all the little quadrilaterals commute. The outer pentagon also commutes, since
Proposition~\ref{prop.shuffle-alg} holds for $W$ instead of $V$. If
$\mathbf{f}$ was surjective, then we would be able to conclude that the inner
pentagon also commutes, so we would immediately get the commutativity of
\eqref{hint.shuffle-alg.diag-to-prove}. But even if $\mathbf{f}$ is not
surjective, we are almost there: The inner pentagon commutes on the image of
the map $\mathbf{f}\otimes\mathbf{f}:\mathbf{W}\otimes\mathbf{W}%
\rightarrow\mathbf{V}\otimes\mathbf{V}$ (because when we start at
$\mathbf{W}\otimes\mathbf{W}$, we can walk around the outer pentagon instead,
which is known to commute), but this image contains $a\otimes b$ (since
$a=v_{1}v_{2}\cdots v_{p}=\mathbf{f}\left(  x_{1}x_{2}\cdots x_{p}\right)  $
and similarly $b=\mathbf{f}\left(  x_{p+1}x_{p+2}\cdots x_{p+q}\right)  $), so
we conclude that \eqref{hint.shuffle-alg.diag-toprove} holds, as we wanted to show.

This is only one of the diagrams we need to prove in order to prove
Proposition~\ref{prop.shuffle-alg}, but the other diagrams are done in the
exact same way.

\bigskip

\textit{Hint to Exercise~\ref{exe.pw-fin-supp.P1-5}.} Straightforward
reasoning using facts like ``a union of finitely many finite
sets is finite'' and ``a tensor is a sum of
finitely many pure tensors''.

\bigskip

\textit{Hint to Exercise~\ref{exe.convolution-series}.} Parts (a), (b), (d)
and (e) of Proposition~\ref{prop.convolution-series} are easy. (In proving
\eqref{eq.exe.convolution-series.b.2} and later, it helps to first establish
an extension of \eqref{eq.exe.convolution-series.b.1} to infinite
sums\footnote{Namely: Let $\left(  r_{q}\right)  _{q\in Q}\in\left(
\mathbf{k}\left[  \left[  T\right]  \right]  \right)  ^{Q}$ be a family of
power series such that the (possibly infinite) sum $\sum_{q\in Q}r_{q}$
converges in $\mathbf{k}\left[  \left[  T\right]  \right]  $. Let
$f\in\mathfrak{n}\left(  C,A\right)  $. Then, the family $\left(  \left(
r_{q}\right)  ^{\star}\left(  f\right)  \right)  _{q\in Q}\in\left(
\Hom \left(  C,A\right)  \right)  ^{Q}$ is pointwise finitely
supported and satisfies $\left(  \sum_{q\in Q}r_{q}\right)  ^{\star}\left(
f\right)  =\sum_{q\in Q}\left(  r_{q}\right)  ^{\star}\left(  f\right)  $.}.)
For part (c), recall that the binomial formula $\left(  a+b\right)  ^{n}%
=\sum_{k=0}^{n}\dbinom{n}{k}a^{k}b^{n-k}$ holds for any two commuting elements
$a$ and $b$ of any ring (such as $f$ and $g$ in the convolution algebra
$\Hom \left(  C,A\right)  $). Part (f) follows from (e) using
\eqref{eq.exe.convolution-series.b.2}. Part (g) is best proved in two steps:
First, use induction to prove part (g) in the case when $u=T^{k}$ for some
$k\in \NN$ (this relies on \eqref{eq.exe.convolution-series.b.2}); then,
notice that both sides of \eqref{eq.exe.convolution-series.g.1} depend
$\kk$-linearly on $u$, whence the general case follows (up to some
mudfighting with infinite sums). Part (h) is an instance of the
``local $\star$-nilpotence'' already observed
in the proof of Proposition \ref{Takeuchi-formula}. Part (j) follows from (h).
Part (i) follows from Proposition \ref{prop.convolution.functor} (applied to
$C^{\prime}=C$, $A^{\prime}=B$, $\gamma= \id_{C}$ and
$\alpha=s$) in a similar way as part (g) followed from \eqref{eq.exe.convolution-series.b.2}.

\bigskip

\textit{Hint to Exercise~\ref{exe.exp-log}.}
Proposition~\ref{prop.exp-log-inv} is a classical result, often proved by a
lazy reference to the mythical complex analysis class the reader has surely
seen it in. Here is a do-it-yourself purely algebraic proof:

\begin{itemize}
\item \textit{Step 1:} If $u,v\in\mathbf{k}\left[  \left[  T\right]  \right]
$ are two power series having the same constant term and satisfying $\dfrac
{d}{dT}u=\dfrac{d}{dT}v$, then $u=v$. This simple lemma (whose analogue for
differentiable functions is a fundamental fact of real analysis) is easily
proved by comparing coefficients in $\dfrac{d}{dT}u=\dfrac{d}{dT}v$ and
recalling that $\kk$ is a $\mathbb{Q}$-algebra (so $1,2,3,\ldots$ are
invertible in $\kk$).

\item \textit{Step 2:} If $u,v\in\mathbf{k}\left[  \left[  T\right]  \right]
$ are two power series having constant term $1$ and satisfying $\left(
\dfrac{d}{dT}u\right)  \cdot v=\left(  \dfrac{d}{dT}v\right)  \cdot u$, then
$u=v$. This can be proved by applying Step 1 to $uv^{-1}$ and $1$ instead of
$u$ and $v$.

\item \textit{Step 3:} The power series $\overline{\log}\left[  \overline
{\exp}\right]  $ and $\overline{\exp}\left[  \overline{\log}\right]  $ are
well-defined and have constant term $0$. (Easy.)

\item \textit{Step 4:} If $w\in\mathbf{k}\left[  \left[  T\right]  \right]  $
is a power series having constant term $0$, then%
\begin{align*}
\dfrac{d}{dT}\left(  \overline{\exp}\left[  w\right]  \right)   &  =\left(
\dfrac{d}{dT}w\right)  \cdot\exp\left[  w\right]
\ \ \ \ \ \ \ \ \ \ \text{and}\\
\dfrac{d}{dT}\left(  \overline{\log}\left[  w\right]  \right)   &  =\left(
\dfrac{d}{dT}w\right)  \cdot\dfrac{1}{1+w}.
\end{align*}
These formulas can be derived from the chain rule, or more directly from
$\overline{\exp}\left[  w\right]  =\sum_{n\geq1}\dfrac{1}{n!}w^{n}$ and
$\overline{\log}\left[  w\right]  =\sum_{n\geq1}\dfrac{\left(  -1\right)
^{n-1}}{n}w^{n}$.

\item \textit{Step 5:} Show $\overline{\exp}\left[  \overline{\log}\right]
=T$ by applying Step 2 to $u=\exp\left[  \overline{\log}\right]  $ and $v=1+T$.

\item \textit{Step 6:} Show $\overline{\log}\left[  \overline{\exp}\right]
=T$ by applying Step 1 to $u=\overline{\log}\left[  \overline{\exp}\right]  $
and $v=T$.
\end{itemize}

Lemma~\ref{lem.exp-log-log} easily follows from
Proposition~\ref{prop.convolution-series}(f).

Remains to prove Proposition~\ref{prop.exp-log}. It is easy to see that
$\log^{\star}\left(  \exp^{\star}f\right)  =\overline{\log}^{\star}\left(
\overline{\exp}^{\star}f\right)  $ for each $f\in\mathfrak{n}\left(
C,A\right)  $; thus, Proposition~\ref{prop.exp-log}(a) follows from
(\ref{eq.exe.convolution-series.g.1}) using Proposition~\ref{prop.exp-log-inv}
and Proposition~\ref{prop.convolution-series}(f) (since $T^{\star}\left(
f\right)  =f$). A similar argument yields Proposition~\ref{prop.exp-log}(b)
(this time, we need to observe that $\exp^{\star}\left(  \log^{\star}g\right)
=\overline{\exp}^{\star}\left(  \overline{\log}^{\star}\left(  g-u_{A}%
\epsilon_{C}\right)  \right)  +u_{A}\epsilon_{C}$ first). To prove
Proposition~\ref{prop.exp-log}(c), first use
Proposition~\ref{prop.convolution-series}(c) to show that $\exp^{\star}\left(
f+g\right)  $ is well-defined; then, apply the well-known fact that
$\exp\left(  x+y\right)  =\exp x\cdot\exp y$ for any two commuting elements
$x$ and $y$ of a ring (provided the exponentials are well-defined; some
yak-shaving is required here to convince oneself that the infinite sums behave
well)\footnote{If you have not seen this well-known fact, prove it by a quick
computation using the binomial formula.}. Part (d) is trivial. Part (e) is an
induction on $n$. Part (f) is a rehash of the definition of $\log^{\star
}\left(  f+u_{A}\epsilon_{C}\right)  =\overline{\log}^{\star}f$.

\bigskip

\textit{Hint to Exercise~\ref{exe.leray.exp-alg}.}
Proposition~\ref{prop.convolution.nat}(a) is easily proved by unpacking the
definition of convolution (just like
Proposition~\ref{prop.convolution.functor}). Part (b) follows from (a) by induction.

The trick to Proposition~\ref{prop.leray.f*nxy} is to realize that if
$f\in\Hom \left(  C,A\right)  $ is as in
Proposition~\ref{prop.leray.f*nxy}, then every $x,y\in C$ satisfy%
\begin{equation}
f\left(  xy\right)  =\epsilon\left(  y\right)  f\left(  x\right)
+\epsilon\left(  x\right)  f\left(  y\right)  ,
\label{hint.leray.exp-alg.f*nxy.1}%
\end{equation}
because $xy-\epsilon\left(  x\right)  y-\epsilon\left(  y\right)
x=\epsilon\left(  x\right)  \epsilon\left(  y\right)  \cdot
1+\underbrace{\left(  x-\epsilon\left(  x\right)  \right)  }_{\in\ker\epsilon
}\underbrace{\left(  y-\epsilon\left(  y\right)  \right)  }_{\in\ker\epsilon}$
is annihilated by $f$. Once this equality is known, it is not hard to prove
Proposition~\ref{prop.leray.f*nxy} ``by hand'' by induction on $n$
(using Sweedler notation). Alternatively, for a cleaner
proof, the equality \eqref{hint.leray.exp-alg.f*nxy.1} can be restated in an
element-free way as
\[
f\circ m_{C}=m_{A}\circ\left(  f\otimes\mathfrak{i}+\mathfrak{i}\otimes
f\right)  ,
\]
where $\mathfrak{i}=u_{A}\circ\epsilon_{C}$ is the unity of the $\mathbf{k}%
$-algebra $\left(  \Hom \left(  C,A\right)  ,\star\right)  $;
then, an application of Proposition~\ref{prop.convolution.nat}(b) shows that
every $n\in \NN$ satisfies%
\begin{align*}
f^{\star n}\circ m_{C}  &  =m_{A}\circ\underbrace{\left(  f\otimes
\mathfrak{i}+\mathfrak{i}\otimes f\right)  ^{\star n}}_{\substack{=\sum
_{i=0}^{n}\dbinom{n}{i}\left(  f\otimes\mathfrak{i}\right)  ^{\star i}%
\star\left(  \mathfrak{i}\otimes f\right)  ^{\star\left(  n-i\right)
}\\\text{(by the binomial formula,}\\\text{since }f\otimes\mathfrak{i}\text{
and }\mathfrak{i}\otimes f\text{ commute in}\\\text{the convolution algebra
}\Hom \left(  C\otimes C,A\otimes A\right)  \text{)}}%
}=m_{A}\circ\left(  \sum_{i=0}^{n}\dbinom{n}{i}\underbrace{\left(
f\otimes\mathfrak{i}\right)  ^{\star i}\star\left(  \mathfrak{i}\otimes
f\right)  ^{\star\left(  n-i\right)  }}_{\substack{=f^{\star i}\otimes
f^{\star\left(  n-i\right)  }\\\text{(by repeated application of
Exercise~\ref{exe.convolution.tensor}(a))}}}\right) \\
&  =m_{A}\circ\left(  \sum_{i=0}^{n}\dbinom{n}{i}f^{\star i}\otimes
f^{\star\left(  n-i\right)  }\right)  ,
\end{align*}
which is precisely Proposition~\ref{prop.leray.f*nxy} (restated in an
element-free way).

Proposition~\ref{prop.leray.exp-algh1} is an easy consequence of
Proposition~\ref{prop.leray.f*nxy}, since $\left(  \exp^{\star}f\right)
\left(  xy\right)  =\sum_{n\in \NN}\dfrac{1}{n!}f^{\star n}\left(
xy\right)  $. (Again, fighting infinite sums is probably the most laborious
part of the proof.)

Lemma \ref{lem.sol.eulerian-idp.2.poly1} can be reduced to the fact that the
matrix $\left(  i^{N+1-j}\right)  _{i,j=1,2,\ldots,N+1}\in{\mathbb{Q}%
}^{\left(  N+1\right)  \times\left(  N+1\right)  }$ is invertible (since its
determinant is the Vandermonde determinant $\prod_{1\leq i<j\leq
N+1}\underbrace{\left(  i-j\right)  }_{\neq0}\neq0$) and thus has trivial
kernel (not just over $\mathbb{Q}$, but on any torsionfree abelian group).

Lemma~\ref{lem.leray.poly1-inf} follows from Lemma
\ref{lem.sol.eulerian-idp.2.poly1}, because a finitely supported family
indexed by nonnegative integers must become all zeroes from some point on.

The proof of Proposition~\ref{prop.leray.exp-algh2} is rather surprising: It
suffices to show that $f\left(  xy\right)  =0$ for all $x,y\in\ker\epsilon$.
So let us fix $x,y\in\ker\epsilon$. Proposition~\ref{prop.convolution-series}%
(h) yields $f\in\mathfrak{n}\left(  C,A\right)  $. Let $t\in \NN$ be
arbitrary. Then, Proposition~\ref{prop.exp-log}(e) (applied to $n=t$) shows
that $tf\in\mathfrak{n}\left(  C,A\right)  $ and $\exp^{\star}\left(
tf\right)  =\left(  \exp^{\star}f\right)  ^{\star t}$. But
Exercise~\ref{exe.alghoms.comm}(b) shows that $\left(  \exp^{\star}f\right)
^{\star t}$ is a $\kk$-algebra homomorphism $C\rightarrow A$. Hence,
$\left(  \exp^{\star}f\right)  ^{\star t}\left(  xy\right)  =\left(
\exp^{\star}f\right)  ^{\star t}\left(  x\right)  \cdot\left(  \exp^{\star
}f\right)  ^{\star t}\left(  y\right)  $. Rewriting $\left(  \exp^{\star
}f\right)  ^{\star t}$ as $\exp^{\star}\left(  tf\right)  =\sum_{n\in
 \NN}\dfrac{1}{n!}f^{\star n}t^{n}$ on both sides, and multiplying out
the right hand side, we can rewrite this as%
\[
\sum_{k\in \NN}\dfrac{1}{k!}f^{\star k}\left(  xy\right)  t^{k}%
=\sum_{k\in \NN}\left(  \sum_{i=0}^{k}\dfrac{f^{\star i}\left(
x\right)  }{i!}\cdot\dfrac{f^{\star\left(  k-i\right)  }\left(  y\right)
}{\left(  k-i\right)  !}\right)  t^{k}.
\]
In other words,%
\[
\sum_{k\in \NN} w_k t^{k}=0,
\ \ \ \ \ \ \ \ \ \ \text{where we set }
w_k = \dfrac{1}{k!}f^{\star k}\left(  xy\right)
-\sum_{i=0}^{k}\dfrac{f^{\star i}\left(  x\right)  }{i!}\cdot\dfrac
{f^{\star\left(  k-i\right)  }\left(  y\right)  }{\left(  k-i\right)
!} .
\]
But we have proved this for all $t\in \NN$. Thus,
Lemma~\ref{lem.leray.poly1-inf} shows that%
\[
w_k = 0\ \ \ \ \ \ \ \ \ \ \text{for every }%
k\in \NN.
\]
Applying this to $k=1$ and simplifying, we obtain $f\left(  xy\right)
-\epsilon\left(  x\right)  f\left(  y\right)  -f\left(  x\right)
\epsilon\left(  y\right)  =0$. Since $x,y\in\ker\epsilon$, this simplifies
even further to $f\left(  xy\right)  =0$, which proves
Proposition~\ref{prop.leray.exp-algh2}.

Finally, we need to prove Proposition~\ref{prop.leray.exp-algh3}. Set
$F=\exp^{\star}f$ and $\widetilde{F}=F-u_{A}\epsilon_{C}$, so that
$\widetilde{F}\in\mathfrak{n}\left(  C,A\right)  $. Then,
Proposition~\ref{prop.leray.exp-algh1} shows that $F:C\rightarrow A$ is a
$\kk$-algebra homomorphism, so it remains to show that $F$ is
surjective. But it is easy to see using Proposition~\ref{prop.exp-log}(a) that
$f=\overline{\log}^{\star}\widetilde{F}$.

Define $\widetilde{\id} \in \mathfrak{n}\left(  C,C\right)  $ by
$\widetilde{\id} = \id_{C}-u_{C}\epsilon_{C}$.
Then, it is not hard to see that $F\circ
\widetilde{\id}=\widetilde{F}$. Hence, $f=\overline{\log}%
^{\star}\underbrace{\widetilde{F}}_{=F\circ\widetilde{\id}}
=\overline{\log}^{\star}\left(  F\circ\widetilde{\id}\right)
=F\circ\left(  \overline{\log}^{\star}\left(  \widetilde{\id}
\right)  \right)  $ (by Proposition~\ref{prop.convolution-series}(i), since
$F$ is a $\kk$-algebra homomorphism). Therefore, $f\left(  C\right)
\subset F\left(  C\right)  $. Since $F$ is a $\kk$-algebra
homomorphism, this entails that $F\left(  C\right)  $ is a $\mathbf{k}%
$-subalgebra of $A$ that contains $f\left(  C\right)  $ as a subset. But this
causes $F\left(  C\right)  $ to be the whole $A$ (since $f\left(  C\right)  $
generates $A$). Thus, $F$ is surjective, so
Proposition~\ref{prop.leray.exp-algh3} is proven.

\bigskip

\textit{Hint to Exercise~\ref{exe.leray.leray-e}.} We must prove
Theorem~\ref{thm.leray.leray-e}. Part (a) is easy. For the remainder of the
proof, we set $\widetilde{\id}= \id_{A}%
-u_{A}\epsilon_{A}\in\operatorname*{End}A$, and equip ourselves with some
simple lemmas:

\begin{itemize}
\item The kernel $\ker\epsilon$ is an ideal of $A$.

\item We have $\widetilde{\id}\in\mathfrak{n}\left(
A,A\right)  $ and $\ker\widetilde{\id} = \kk \cdot1_{A}$ and
$\widetilde{\id}\left(  A\right)  =\ker\epsilon$.

\item We have $A/\left(  \mathbf{k}\cdot1_{A}+\left(  \ker\epsilon\right)
^{2}\right)  \cong\left(  \ker\epsilon\right)  /\left(  \ker\epsilon\right)
^{2}$ as $\kk$-modules.
\end{itemize}

Now, to the proof of Theorem~\ref{thm.leray.leray-e}(b). Using $\mathfrak{e}%
=\log^{\star}\left(  \id_{A}\right)  =\overline{\log
}^{\star}\widetilde{\id}$ and $\widetilde{\id}
\left(  1_{A}\right)  =0$, it is easy to see that $\mathfrak{e}\left(
1_{A}\right)  =0$. Hence, $\mathfrak{e}\left(  A_{0}\right)  =0$ since $A$ is
connected. Thus, Proposition~\ref{prop.leray.exp-algh2} shows that
$\mathfrak{e}\left(  \left(  \ker\epsilon\right)  ^{2}\right)  =0$ (since
$\exp^{\star}\mathfrak{e}= \id_{A}$ is a $\mathbf{k}%
$-algebra homomorphism). Combined with $\mathfrak{e}\left(  1_{A}\right)  =0$,
this yields $\mathbf{k}\cdot1_{A}+\left(  \ker\epsilon\right)  ^{2}\subset
\ker\mathfrak{e}$. But this inclusion is actually an equality, as we can show
by the following computation: We have $\mathfrak{e}=\overline{\log}^{\star
}\widetilde{\id}=\sum_{n\geq1}\dfrac{\left(  -1\right)  ^{n-1}%
}{n}\widetilde{\id}^{\star n}$, and therefore each $x\in A$
satisfies%
\begin{align*}
\mathfrak{e}\left(  x\right)   &  =\sum_{n\geq1}\dfrac{\left(  -1\right)
^{n-1}}{n}\widetilde{\id}^{\star n}\left(  x\right)
=\underbrace{\widetilde{\id}\left(  x\right)  }%
_{\substack{=x-\epsilon\left(  x\right)  1_{A}\\\text{(by the definition of
}\widetilde{\id}\text{)}}}+\sum_{n\geq2}\dfrac{\left(
-1\right)  ^{n-1}}{n}\underbrace{\widetilde{\id}^{\star
n}\left(  x\right)  }_{\substack{\in\left(  \widetilde{\id}
\left(  A\right)  \right)  ^{n}\\\text{(by induction on }n\text{,}%
\\\text{using the definition}\\\text{of convolution)}}}\\
&  \in x-\epsilon\left(  x\right)  1_{A}+\sum_{n\geq2}\dfrac{\left(
-1\right)  ^{n-1}}{n}\left(  \underbrace{\widetilde{\id}\left(
A\right)  }_{=\ker\epsilon}\right)  ^{n}=x-\underbrace{\epsilon\left(
x\right)  }_{\in\mathbf{k}}1_{A}+\underbrace{\sum_{n\geq2}\dfrac{\left(
-1\right)  ^{n-1}}{n}\left(  \ker\epsilon\right)  ^{n}}_{\subset\left(
\ker\epsilon\right)  ^{2}}\subset x-\mathbf{k}\cdot1_{A}+\left(  \ker
\epsilon\right)  ^{2},
\end{align*}
so that
\begin{equation}
x-\mathfrak{e}\left(  x\right)  \in\mathbf{k}\cdot1_{A}+\left(  \ker
\epsilon\right)  ^{2}. \label{hint.leray.leray-e.b.2}%
\end{equation}
If $x\in\ker\mathfrak{e}$, then this simplifies to $x\in\mathbf{k}\cdot
1_{A}+\left(  \ker\epsilon\right)  ^{2}$. Thus, $\ker\mathfrak{e}%
\subset\mathbf{k}\cdot1_{A}+\left(  \ker\epsilon\right)  ^{2}$. Combining this
with $\mathbf{k}\cdot1_{A}+\left(  \ker\epsilon\right)  ^{2}\subset
\ker\mathfrak{e}$, we obtain $\ker\mathfrak{e}=\mathbf{k}\cdot1_{A}+\left(
\ker\epsilon\right)  ^{2}$. But the homomorphism theorem yields%
\[
\mathfrak{e}\left(  A\right)  \cong A/\underbrace{\ker\mathfrak{e}%
}_{=\mathbf{k}\cdot1_{A}+\left(  \ker\epsilon\right)  ^{2}}=A/\left(
\mathbf{k}\cdot1_{A}+\left(  \ker\epsilon\right)  ^{2}\right)  \cong\left(
\ker\epsilon\right)  /\left(  \ker\epsilon\right)  ^{2}%
\ \ \ \ \ \ \ \ \ \ \left(  \text{as seen above}\right)
\]
as $\kk$-modules. This completes the proof of
Theorem~\ref{thm.leray.leray-e}(b).

Theorem~\ref{thm.leray.leray-e}(c) just requires showing that $\mathfrak{q}%
\left(  A_{0}\right)  =0$, which is a consequence of $\mathfrak{e}\left(
A_{0}\right)  =0$.

Next, we shall prove Theorem~\ref{thm.leray.leray-e}(d). We have
$\mathfrak{q}\in\mathfrak{n}\left(  A,\operatorname{Sym}\left(  \mathfrak{e}%
\left(  A\right)  \right)  \right)  $. Furthermore, $\mathfrak{q}\left(
A\right)  $ generates the $\kk$-algebra $\operatorname{Sym}\left(
\mathfrak{e}\left(  A\right)  \right)  $ (since $\mathfrak{q}\left(  A\right)
=\operatorname{Sym}^{1}\left(  \mathfrak{e}\left(  A\right)  \right)  $). From
Theorem~\ref{thm.leray.leray-e}(b), we get $\ker\mathfrak{e}=\mathbf{k}%
\cdot1_{A}+\left(  \ker\epsilon\right)  ^{2}$, from which we easily obtain
$\mathfrak{q}\left(  1_{A}\right)  =0$ and $\mathfrak{q}\left(  \left(
\ker\epsilon\right)  ^{2}\right)  =0$. Thus,
Proposition~\ref{prop.leray.exp-algh3} (applied to $A$, $\operatorname{Sym}%
\left(  \mathfrak{e}\left(  A\right)  \right)  $ and $\mathfrak{q}$ instead of
$C$, $A$ and $f$) shows that $\exp^{\star}\mathfrak{q}:A\rightarrow
\operatorname{Sym}\left(  \mathfrak{e}\left(  A\right)  \right)  $ is a
surjective $\kk$-algebra homomorphism. But $\mathfrak{s}$ is a
$\kk$-algebra homomorphism $\operatorname{Sym}\left(  \mathfrak{e}%
\left(  A\right)  \right)  \rightarrow A$ and satisfies $\mathbf{i}%
=\mathfrak{s}\circ\iota_{\mathfrak{e}\left(  A\right)  }$ (by its definition).
Thus, Proposition~\ref{prop.convolution-series}(i) (applied to $A$,
$\operatorname{Sym}\left(  \mathfrak{e}\left(  A\right)  \right)  $, $A$,
$\mathfrak{s}$, $\exp$ and $\mathfrak{q}$ instead of $C$, $A$, $B$, $s$, $u$
and $f$) shows that $\mathfrak{s}\circ\mathfrak{q}\in\mathfrak{n}\left(
A,A\right)  $ and $\exp^{\star}\left(  \mathfrak{s}\circ\mathfrak{q}\right)
=\mathfrak{s}\circ\left(  \exp^{\star}\mathfrak{q}\right)  $. However, it is
easy to see that $\mathfrak{s}\circ\mathfrak{q}=\mathfrak{e}$ (since
$\mathbf{i}=\mathfrak{s}\circ\iota_{\mathfrak{e}\left(  A\right)  }$); this
lets us rewrite the equality $\exp^{\star}\left(  \mathfrak{s}\circ
\mathfrak{q}\right)  =\mathfrak{s}\circ\left(  \exp^{\star}\mathfrak{q}%
\right)  $ as $\exp^{\star}\mathfrak{e}=\mathfrak{s}\circ\left(  \exp^{\star
}\mathfrak{q}\right)  $. Comparing this with $\exp^{\star}\mathfrak{e}%
= \id_{A}$, we obtain $\mathfrak{s}\circ\left(
\exp^{\star}\mathfrak{q}\right)  = \id_{A}$. Since $\exp^{\star
}\mathfrak{q}$ is surjective, this entails that the maps $\exp^{\star
}\mathfrak{q}$ and $\mathfrak{s}$ are mutually inverse. This proves
Theorem~\ref{thm.leray.leray-e}(d).

Theorem~\ref{thm.leray.leray-e}(d) shows that $A\cong\operatorname{Sym}\left(
\mathfrak{e}\left(  A\right)  \right)  $ as $\kk$-algebras, but
Theorem~\ref{thm.leray.leray-e}(b) shows that $\mathfrak{e}\left(  A\right)
\cong\left(  \ker\epsilon\right)  /\left(  \ker\epsilon\right)  ^{2}$ as
$\kk$-modules. Combining these, we obtain
Theorem~\ref{thm.leray.leray-e}(e).

Finally, to prove Theorem~\ref{thm.leray.leray-e}(f), we notice that each
$x\in A$ satisfies%
\begin{align*}
x-\mathfrak{e}\left(  x\right)   &  \in\mathbf{k}\cdot1_{A}+\left(
\ker\epsilon\right)  ^{2}\ \ \ \ \ \ \ \ \ \ \left(  \text{by
\eqref{hint.leray.leray-e.b.2}}\right) \\
&  =\ker\mathfrak{e}\ \ \ \ \ \ \ \ \ \ \left(  \text{by
Theorem~\ref{thm.leray.leray-e}(b)}\right)
\end{align*}
and thus $0=\mathfrak{e}\left(  x-\mathfrak{e}\left(  x\right)  \right)
=\mathfrak{e}\left(  x\right)  -\left(  \mathfrak{e}\circ\mathfrak{e}\right)
\left(  x\right)  $.

% [DG][v80] Added the above chapter (mostly pulling from the
% solutions, but omitting all but the most important milestones).

%%%%%%%%%%%%%%%%%%%%%%%%%%%%%%%%%%%%%%%%%%%%%%%%%%%%%%%%%
\section*{Acknowledgements}
The authors thank the following for helpful comments and/or teaching
them about Hopf algebras:
Marcelo Aguiar, Federico Ardila, Lou Billera, Richard Ehrenborg,
Mark Haiman, Florent Hivert,
Christophe Hohlweg, Jia Huang, Jang Soo Kim, Aaron Lauve,
Dominique Manchon, John Palmieri, Alexander Postnikov,
Margie Readdy, Nathan Reading, Christophe Reutenauer,
Hans-J\"urgen Schneider,
Richard Stanley, Josh Swanson, Muge Taskin, Jean-Yves Thibon.

Parts of this text have been written during stays at the
Mathematisches Forschungsinstitut Oberwolfach (2019 and 2020)%
\footnote{This research was supported through the programme
``Oberwolfach Leibniz Fellows'' by the
Mathematisches Forschungsinstitut Oberwolfach in 2019 and 2020.}
and at the Institut Mittag--Leffler Djursholm (Spring 2020,
supported by the Swedish Research Council under grant
no. 2016-06596);
DG thanks both for their hospitality.
%%%%%%%%%%%%%%%%%%%%%%%%%%%%%%%%%%%%%%%%%%%%%%%%%%%%%%%%%

% [DG][v29] Added Stanley to the list.

\horiline

\IfFileExists{./HopfComb-v80.ind}{\input{HopfComb-v80.ind}}{\newpage An index goes here!}

\begin{nosolutions}
\end{document}